\documentclass[12pt]{book}
\usepackage{amsmath, amscd, amssymb, amsthm}
\usepackage[subnum]{cases}
\usepackage{changepage}
\usepackage{marginnote}
\usepackage{hyperref}
\usepackage{cleveref}
\usepackage[all]{xy}
\usepackage{tabularx}
\usepackage{ltablex}
\usepackage[T1]{fontenc}
\usepackage{setspace}
\usepackage{fancybox}
\usepackage{ragged2e}
\usepackage[margin=0.8in]{geometry}
\hypersetup{colorlinks,linkcolor={black},citecolor={black}}
\usepackage{times}
\usepackage{newtxmath}

\theoremstyle{definition}
\newtheorem{theorem}{Theorem}[section]

\newtheorem{proposition}[theorem]{Proposition}

\newtheorem{definition}[theorem]{Definition}
\newtheorem{remark}[theorem]{Remark}
\numberwithin{equation}{section}

\usepackage{fancyhdr}
\begin{document}

\normalfont

\title{$\infty$-Categorical Approaches to Hodge-Iwasawa Theory II: $\infty$-Categorical and Derived Hodge-Iwasawa Modules}
\author{Xin Tong}
\date{}

\maketitle

\newpage

\subsection*{Abstract}
This paper is a discussion on $\infty$-categorical approaches to Hodge-Iwasawa Theory, which was initiated in our project on the $\infty$-categorical approaches to Hodge-Iwasawa Theory. The theory aims at the serious  unification of $p$-adic Hodge Theory and $p$-adic Iwasawa Theory, by taking deformation of Hodge-theoretic constructions along some consideration in Iwasawa Theory beyond the Iwasawa deformation of certain motives in the general sense. The Hodge modules in our current consideration will be essentially within $\infty$-categorical derived categories of inductive Banach modules and $\infty$-categorical derived categories of condensed solidification of certain topological modules.

\newpage

\tableofcontents

\newpage

\chapter{Introduction}

\section{Results and Notations}

\subsection{Introduction to the Main Ideas}

\indent We discuss in this article the corresponding $\infty$-categorical and homotopicalization of some our serious consideration in \cite{10T2}, with some goal in mind to generalize \cite{10KL1} and \cite{10KL2} along some deformation point of view. We hope to remind the readers of the fact that Hodge-Iwasawa theory is a combination of consideration in both Hodge theory and Iwasawa theory with certain moduli stack consideration closely after  
\cite{10BF1}, \cite{10BF2}, \cite{10FK}, \cite{10FS}, \cite{10He}, \cite{10KL1}, \cite{10KL2}, \cite{10KP}, \cite{10KPX}, \cite{10Na1}, \cite{10Na2}, \cite{10PR}, \cite{10RZ}, \cite{10Sch2}, \cite{10SW}, \cite{10Wit}. The moduli stack consideration is following those in \cite{102He1}, \cite{102HH}, \cite{102HV}, \cite{102PR}, \cite{102RZ}, \cite{102SW}, \cite{102FS}, \cite{102Sch}, \cite{102Laff}, \cite{102GL}, \cite{102Ked}, \cite{102Har1}, \cite{102Har2}, \cite{102EG}, \cite{10GEH}, \cite{10HHS}, \cite{10D}, \cite{10L}, \cite{10W}\footnote{The consideration will be essentially after \cite{102EG}, \cite{10GEH}, \cite{10HHS}. See \cite[Conjecture 5.1.18, Section 5.2, Theorem 5.2.4]{10GEH} for the detail of Emerton-Gee-Hellmann conjecture on the moduli stack of $(\varphi,\Gamma)$-modules.}. We are also inspired by \cite{103AB1}, \cite{103AB2}, \cite{103AB3}, \cite{103AI1}, \cite{103AI2}, \cite{103AI3}, \cite{103BMS1}, \cite{103BMS2}, \cite{103BS1}, \cite{103BS2}, \cite{103Fa1}, \cite{103Fa2}, \cite{103Fa3}, \cite{103Fon1}, \cite{103Fon2}, \cite{103Fon3}, \cite{103Fon4}, \cite{103Fon5}, \cite{103Iwa}, \cite{103Ka1}, \cite{103Ka2}, \cite{103KL1}, \cite{103KL2}, \cite{103KL3}, \cite{103KL4}, \cite{103KP}, \cite{103Lu1}, \cite{103Lu2}, \cite{103Lu3}, \cite{103Mann}, \cite{103Sch1}, \cite{103Sch3}, \cite{103Wi1}, \cite{103Wi2}, \cite{103Wi3}, \cite{10HP}. We consider the multivariatization as well from \cite{10CKZ}, \cite{10PZ} and \cite{10BCM}. The deformation is so important that we need to establish the corresponding results even in the situation where the spaces are relatively complicated to study. For instance the importance to look at $\infty$-categorical and homotopical aspects took the deep roots from the sheafiness of the topological period rings we are considering. However the sheafiness issue could be weakened after we look at some development from \cite{10BK}, \cite{10BBBK}, \cite{10BBK},  \cite{10KKM}, \cite{10BBM}, \cite{10CS1}, \cite{10CS2}. Generalizing the sheafiness to derived sheafiness, one could even achieve more such as the corresponding $\infty$-descent in quite general sense. The following two examples of results give the main ideas of the picture along the discussion above after \cite{10BK}, \cite{10BBBK}, \cite{10BBK},  \cite{10KKM}, \cite{10BBM}, \cite{10CS1}, \cite{10CS2}.\\

\begin{proposition}
There is a well-defined functor from the $\infty$-category 
\begin{align}
\mathrm{Quasicoherentpresheaves,Condensed}_{*}	
\end{align}
where $*$ is one of the following spaces:
\begin{align}
&\mathrm{Spec}^\mathrm{CS}\widetilde{\Phi}_{\psi,\Gamma,A}/\mathrm{Fro}^\mathbb{Z},	\\
\end{align}
\begin{align}
&\mathrm{Spec}^\mathrm{CS}\breve{\Phi}_{\psi,\Gamma,A}/\mathrm{Fro}^\mathbb{Z},	\\
\end{align}
\begin{align}
&\mathrm{Spec}^\mathrm{CS}{\Phi}_{\psi,\Gamma,A}/\mathrm{Fro}^\mathbb{Z},	
\end{align}
to the $\infty$-category of $\mathrm{Fro}$-equivariant quasicoherent presheaves over similar spaces above correspondingly without the $\mathrm{Fro}$-quotients, and to the $\infty$-category of $\mathrm{Fro}$-equivariant quasicoherent modules over global sections of the structure $\infty$-sheaves of the similar spaces above correspondingly without the $\mathrm{Fro}$-quotients. Here for those space without notation related to the radius and the corresponding interval we consider the total unions $\bigcap_r,\bigcup_I$ in order to achieve the whole spaces to achieve the analogues of the corresponding FF curves from \cite{10KL1}, \cite{10KL2}, \cite{10FF} for
\[
\xymatrix@R+0pc@C+0pc{
\underset{r}{\mathrm{homotopylimit}}~\mathrm{Spec}^\mathrm{CS}\widetilde{\Phi}^r_{\psi,\Gamma,A},\underset{I}{\mathrm{homotopycolimit}}~\mathrm{Spec}^\mathrm{CS}\widetilde{\Phi}^I_{\psi,\Gamma,A},	\\
}
\]
\[
\xymatrix@R+0pc@C+0pc{
\underset{r}{\mathrm{homotopylimit}}~\mathrm{Spec}^\mathrm{CS}\breve{\Phi}^r_{\psi,\Gamma,A},\underset{I}{\mathrm{homotopycolimit}}~\mathrm{Spec}^\mathrm{CS}\breve{\Phi}^I_{\psi,\Gamma,A},	\\
}
\]
\[
\xymatrix@R+0pc@C+0pc{
\underset{r}{\mathrm{homotopylimit}}~\mathrm{Spec}^\mathrm{CS}{\Phi}^r_{\psi,\Gamma,A},\underset{I}{\mathrm{homotopycolimit}}~\mathrm{Spec}^\mathrm{CS}{\Phi}^I_{\psi,\Gamma,A}.	
}
\]
\[ 
\xymatrix@R+0pc@C+0pc{
\underset{r}{\mathrm{homotopylimit}}~\mathrm{Spec}^\mathrm{CS}\widetilde{\Phi}^r_{\psi,\Gamma,A}/\mathrm{Fro}^\mathbb{Z},\underset{I}{\mathrm{homotopycolimit}}~\mathrm{Spec}^\mathrm{CS}\widetilde{\Phi}^I_{\psi,\Gamma,A}/\mathrm{Fro}^\mathbb{Z},	\\
}
\]
\[ 
\xymatrix@R+0pc@C+0pc{
\underset{r}{\mathrm{homotopylimit}}~\mathrm{Spec}^\mathrm{CS}\breve{\Phi}^r_{\psi,\Gamma,A}/\mathrm{Fro}^\mathbb{Z},\underset{I}{\mathrm{homotopycolimit}}~\breve{\Phi}^I_{\psi,\Gamma,A}/\mathrm{Fro}^\mathbb{Z},	\\
}
\]
\[ 
\xymatrix@R+0pc@C+0pc{
\underset{r}{\mathrm{homotopylimit}}~\mathrm{Spec}^\mathrm{CS}{\Phi}^r_{\psi,\Gamma,A}/\mathrm{Fro}^\mathbb{Z},\underset{I}{\mathrm{homotopycolimit}}~\mathrm{Spec}^\mathrm{CS}{\Phi}^I_{\psi,\Gamma,A}/\mathrm{Fro}^\mathbb{Z}.	
}
\]	
In this situation we will have the target category being family parametrized by $r$ or $I$ in compatible glueing sense as in \cite[Definition 5.4.10]{10KL2}. In this situation for modules parametrized by the intervals we have the equivalence of $\infty$-categories by using \cite[Proposition 13.8]{10CS2}. Here the corresponding quasicoherent Frobenius modules are defined to be the corresponding homotopy colimits and limits of Frobenius modules:
\begin{align}
\underset{r}{\mathrm{homotopycolimit}}~M_r,\\
\underset{I}{\mathrm{homotopylimit}}~M_I,	
\end{align}
where each $M_r$ is a Frobenius-equivariant module over the period ring with respect to some radius $r$ while each $M_I$ is a Frobenius-equivariant module over the period ring with respect to some interval $I$\footnote{This means that we can descend to some ring or space with respect to some radius $r$, since we are talking about the corresponding projective limits of stacks (carrying inductive limits of $\infty$-sheaves of rings). This is in some sense (especially when we consider finally the corresponding possibly very non-quasicompact space $X$ instead of a ring $A$) deviating from the situations above over full Robba rings where when deforming over $X$ for instance noncompact we will achieve families of the parameters ${r}$ going to $\infty$. This is not actually hard to understand since we are taking two combined limits in different orders:
\begin{align}
\underset{r}{\mathrm{homotopycolimit}}~\underset{Y\subset X}{\mathrm{homotopycolimit}}\neq \underset{Y\subset X}{\mathrm{homotopycolimit}}~\underset{r}{\mathrm{homotopycolimit}}.
\end{align}
}.\\
\end{proposition}

\begin{proposition}
Similar proposition holds for 
\begin{align}
\mathrm{Quasicoherentsheaves,IndBanach}_{*}.	
\end{align}	
\end{proposition}

\

\begin{proposition}
There is a well-defined functor from the $\infty$-category 
\begin{align}
\mathrm{Quasicoherentpresheaves,Perfectcomplex,Condensed}_{*}	
\end{align}
where $*$ is one of the following spaces:
\begin{align}
&\mathrm{Spec}^\mathrm{CS}\widetilde{\Phi}_{\psi,\Gamma,A}/\mathrm{Fro}^\mathbb{Z},	\\
\end{align}
\begin{align}
&\mathrm{Spec}^\mathrm{CS}\breve{\Phi}_{\psi,\Gamma,A}/\mathrm{Fro}^\mathbb{Z},	\\
\end{align}
\begin{align}
&\mathrm{Spec}^\mathrm{CS}{\Phi}_{\psi,\Gamma,A}/\mathrm{Fro}^\mathbb{Z},	
\end{align}
to the $\infty$-category of $\mathrm{Fro}$-equivariant quasicoherent presheaves over similar spaces above correspondingly without the $\mathrm{Fro}$-quotients, and to the $\infty$-category of $\mathrm{Fro}$-equivariant quasicoherent modules over global sections of the structure $\infty$-sheaves of the similar spaces above correspondingly without the $\mathrm{Fro}$-quotients. Here for those space without notation related to the radius and the corresponding interval we consider the total unions $\bigcap_r,\bigcup_I$ in order to achieve the whole spaces to achieve the analogues of the corresponding FF curves from \cite{10KL1}, \cite{10KL2}, \cite{10FF} for
\[
\xymatrix@R+0pc@C+0pc{
\underset{r}{\mathrm{homotopylimit}}~\mathrm{Spec}^\mathrm{CS}\widetilde{\Phi}^r_{\psi,\Gamma,A},\underset{I}{\mathrm{homotopycolimit}}~\mathrm{Spec}^\mathrm{CS}\widetilde{\Phi}^I_{\psi,\Gamma,A},	\\
}
\]
\[
\xymatrix@R+0pc@C+0pc{
\underset{r}{\mathrm{homotopylimit}}~\mathrm{Spec}^\mathrm{CS}\breve{\Phi}^r_{\psi,\Gamma,A},\underset{I}{\mathrm{homotopycolimit}}~\mathrm{Spec}^\mathrm{CS}\breve{\Phi}^I_{\psi,\Gamma,A},	\\
}
\]
\[
\xymatrix@R+0pc@C+0pc{
\underset{r}{\mathrm{homotopylimit}}~\mathrm{Spec}^\mathrm{CS}{\Phi}^r_{\psi,\Gamma,A},\underset{I}{\mathrm{homotopycolimit}}~\mathrm{Spec}^\mathrm{CS}{\Phi}^I_{\psi,\Gamma,A}.	
}
\]
\[ 
\xymatrix@R+0pc@C+0pc{
\underset{r}{\mathrm{homotopylimit}}~\mathrm{Spec}^\mathrm{CS}\widetilde{\Phi}^r_{\psi,\Gamma,A}/\mathrm{Fro}^\mathbb{Z},\underset{I}{\mathrm{homotopycolimit}}~\mathrm{Spec}^\mathrm{CS}\widetilde{\Phi}^I_{\psi,\Gamma,A}/\mathrm{Fro}^\mathbb{Z},	\\
}
\]
\[ 
\xymatrix@R+0pc@C+0pc{
\underset{r}{\mathrm{homotopylimit}}~\mathrm{Spec}^\mathrm{CS}\breve{\Phi}^r_{\psi,\Gamma,A}/\mathrm{Fro}^\mathbb{Z},\underset{I}{\mathrm{homotopycolimit}}~\breve{\Phi}^I_{\psi,\Gamma,A}/\mathrm{Fro}^\mathbb{Z},	\\
}
\]
\[ 
\xymatrix@R+0pc@C+0pc{
\underset{r}{\mathrm{homotopylimit}}~\mathrm{Spec}^\mathrm{CS}{\Phi}^r_{\psi,\Gamma,A}/\mathrm{Fro}^\mathbb{Z},\underset{I}{\mathrm{homotopycolimit}}~\mathrm{Spec}^\mathrm{CS}{\Phi}^I_{\psi,\Gamma,A}/\mathrm{Fro}^\mathbb{Z}.	
}
\]	
In this situation we will have the target category being family parametrized by $r$ or $I$ in compatible glueing sense as in \cite[Definition 5.4.10]{10KL2}. In this situation for modules parametrized by the intervals we have the equivalence of $\infty$-categories by using \cite[Proposition 12.18]{10CS2}. Here the corresponding quasicoherent Frobenius modules are defined to be the corresponding homotopy colimits and limits of Frobenius modules:
\begin{align}
\underset{r}{\mathrm{homotopycolimit}}~M_r,\\
\underset{I}{\mathrm{homotopylimit}}~M_I,	
\end{align}
where each $M_r$ is a Frobenius-equivariant module over the period ring with respect to some radius $r$ while each $M_I$ is a Frobenius-equivariant module over the period ring with respect to some interval $I$.\\
\end{proposition}

\begin{proposition}
Similar proposition holds for 
\begin{align}
\mathrm{Quasicoherentsheaves,Perfectcomplex,IndBanach}_{*}.	
\end{align}	
\end{proposition}

\begin{proposition}
There is a well-defined functor from the $\infty$-category 
\begin{align}
\mathrm{Quasicoherentpresheaves,Condensed}_{*}	
\end{align}
where $*$ is one of the following spaces:
\begin{align}
&\underset{\mathrm{Spec}}{\mathcal{O}}^\mathrm{CS}\widetilde{\Phi}_{\psi,\Gamma,A}/\mathrm{Fro}^\mathbb{Z},	\\
\end{align}
\begin{align}
&\underset{\mathrm{Spec}}{\mathcal{O}}^\mathrm{CS}\breve{\Phi}_{\psi,\Gamma,A}/\mathrm{Fro}^\mathbb{Z},	\\
\end{align}
\begin{align}
&\underset{\mathrm{Spec}}{\mathcal{O}}^\mathrm{CS}{\Phi}_{\psi,\Gamma,A}/\mathrm{Fro}^\mathbb{Z},	
\end{align}
to the $\infty$-category of $\mathrm{Fro}$-equivariant quasicoherent presheaves over similar spaces above correspondingly without the $\mathrm{Fro}$-quotients, and to the $\infty$-category of $\mathrm{Fro}$-equivariant quasicoherent modules over global sections of the structure $\infty$-sheaves of the similar spaces above correspondingly without the $\mathrm{Fro}$-quotients. Here for those space without notation related to the radius and the corresponding interval we consider the total unions $\bigcap_r,\bigcup_I$ in order to achieve the whole spaces to achieve the analogues of the corresponding FF curves from \cite{10KL1}, \cite{10KL2}, \cite{10FF} for
\[
\xymatrix@R+0pc@C+0pc{
\underset{r}{\mathrm{homotopycolimit}}~\underset{\mathrm{Spec}}{\mathcal{O}}^\mathrm{CS}\widetilde{\Phi}^r_{\psi,\Gamma,A},\underset{I}{\mathrm{homotopylimit}}~\underset{\mathrm{Spec}}{\mathcal{O}}^\mathrm{CS}\widetilde{\Phi}^I_{\psi,\Gamma,A},	\\
}
\]
\[
\xymatrix@R+0pc@C+0pc{
\underset{r}{\mathrm{homotopycolimit}}~\underset{\mathrm{Spec}}{\mathcal{O}}^\mathrm{CS}\breve{\Phi}^r_{\psi,\Gamma,A},\underset{I}{\mathrm{homotopylimit}}~\underset{\mathrm{Spec}}{\mathcal{O}}^\mathrm{CS}\breve{\Phi}^I_{\psi,\Gamma,A},	\\
}
\]
\[
\xymatrix@R+0pc@C+0pc{
\underset{r}{\mathrm{homotopycolimit}}~\underset{\mathrm{Spec}}{\mathcal{O}}^\mathrm{CS}{\Phi}^r_{\psi,\Gamma,A},\underset{I}{\mathrm{homotopylimit}}~\underset{\mathrm{Spec}}{\mathcal{O}}^\mathrm{CS}{\Phi}^I_{\psi,\Gamma,A}.	
}
\]
\[ 
\xymatrix@R+0pc@C+0pc{
\underset{r}{\mathrm{homotopycolimit}}~\underset{\mathrm{Spec}}{\mathcal{O}}^\mathrm{CS}\widetilde{\Phi}^r_{\psi,\Gamma,A}/\mathrm{Fro}^\mathbb{Z},\underset{I}{\mathrm{homotopylimit}}~\underset{\mathrm{Spec}}{\mathcal{O}}^\mathrm{CS}\widetilde{\Phi}^I_{\psi,\Gamma,A}/\mathrm{Fro}^\mathbb{Z},	\\
}
\]
\[ 
\xymatrix@R+0pc@C+0pc{
\underset{r}{\mathrm{homotopycolimit}}~\underset{\mathrm{Spec}}{\mathcal{O}}^\mathrm{CS}\breve{\Phi}^r_{\psi,\Gamma,A}/\mathrm{Fro}^\mathbb{Z},\underset{I}{\mathrm{homotopylimit}}~\breve{\Phi}^I_{\psi,\Gamma,A}/\mathrm{Fro}^\mathbb{Z},	\\
}
\]
\[ 
\xymatrix@R+0pc@C+0pc{
\underset{r}{\mathrm{homotopycolimit}}~\underset{\mathrm{Spec}}{\mathcal{O}}^\mathrm{CS}{\Phi}^r_{\psi,\Gamma,A}/\mathrm{Fro}^\mathbb{Z},\underset{I}{\mathrm{homotopylimit}}~\underset{\mathrm{Spec}}{\mathcal{O}}^\mathrm{CS}{\Phi}^I_{\psi,\Gamma,A}/\mathrm{Fro}^\mathbb{Z}.	
}
\]	
In this situation we will have the target category being family parametrized by $r$ or $I$ in compatible glueing sense as in \cite[Definition 5.4.10]{10KL2}. In this situation for modules parametrized by the intervals we have the equivalence of $\infty$-categories by using \cite[Proposition 13.8]{10CS2}. Here the corresponding quasicoherent Frobenius modules are defined to be the corresponding homotopy colimits and limits of Frobenius modules:
\begin{align}
\underset{r}{\mathrm{homotopycolimit}}~M_r,\\
\underset{I}{\mathrm{homotopylimit}}~M_I,	
\end{align}
where each $M_r$ is a Frobenius-equivariant module over the period ring with respect to some radius $r$ while each $M_I$ is a Frobenius-equivariant module over the period ring with respect to some interval $I$.\\
\end{proposition}

\begin{proposition}
Similar proposition holds for 
\begin{align}
\mathrm{Quasicoherentsheaves,IndBanach}_{*}.	
\end{align}	
\end{proposition}

\

\begin{proposition}
There is a well-defined functor from the $\infty$-category 
\begin{align}
\mathrm{Quasicoherentpresheaves,Perfectcomplex,Condensed}_{*}	
\end{align}
where $*$ is one of the following spaces:
\begin{align}
&\underset{\mathrm{Spec}}{\mathcal{O}}^\mathrm{CS}\widetilde{\Phi}_{\psi,\Gamma,A}/\mathrm{Fro}^\mathbb{Z},	\\
\end{align}
\begin{align}
&\underset{\mathrm{Spec}}{\mathcal{O}}^\mathrm{CS}\breve{\Phi}_{\psi,\Gamma,A}/\mathrm{Fro}^\mathbb{Z},	\\
\end{align}
\begin{align}
&\underset{\mathrm{Spec}}{\mathcal{O}}^\mathrm{CS}{\Phi}_{\psi,\Gamma,A}/\mathrm{Fro}^\mathbb{Z},	
\end{align}
to the $\infty$-category of $\mathrm{Fro}$-equivariant quasicoherent presheaves over similar spaces above correspondingly without the $\mathrm{Fro}$-quotients, and to the $\infty$-category of $\mathrm{Fro}$-equivariant quasicoherent modules over global sections of the structure $\infty$-sheaves of the similar spaces above correspondingly without the $\mathrm{Fro}$-quotients. Here for those space without notation related to the radius and the corresponding interval we consider the total unions $\bigcap_r,\bigcup_I$ in order to achieve the whole spaces to achieve the analogues of the corresponding FF curves from \cite{10KL1}, \cite{10KL2}, \cite{10FF} for
\[
\xymatrix@R+0pc@C+0pc{
\underset{r}{\mathrm{homotopycolimit}}~\underset{\mathrm{Spec}}{\mathcal{O}}^\mathrm{CS}\widetilde{\Phi}^r_{\psi,\Gamma,A},\underset{I}{\mathrm{homotopylimit}}~\underset{\mathrm{Spec}}{\mathcal{O}}^\mathrm{CS}\widetilde{\Phi}^I_{\psi,\Gamma,A},	\\
}
\]
\[
\xymatrix@R+0pc@C+0pc{
\underset{r}{\mathrm{homotopycolimit}}~\underset{\mathrm{Spec}}{\mathcal{O}}^\mathrm{CS}\breve{\Phi}^r_{\psi,\Gamma,A},\underset{I}{\mathrm{homotopylimit}}~\underset{\mathrm{Spec}}{\mathcal{O}}^\mathrm{CS}\breve{\Phi}^I_{\psi,\Gamma,A},	\\
}
\]
\[
\xymatrix@R+0pc@C+0pc{
\underset{r}{\mathrm{homotopycolimit}}~\underset{\mathrm{Spec}}{\mathcal{O}}^\mathrm{CS}{\Phi}^r_{\psi,\Gamma,A},\underset{I}{\mathrm{homotopylimit}}~\underset{\mathrm{Spec}}{\mathcal{O}}^\mathrm{CS}{\Phi}^I_{\psi,\Gamma,A}.	
}
\]
\[ 
\xymatrix@R+0pc@C+0pc{
\underset{r}{\mathrm{homotopycolimit}}~\underset{\mathrm{Spec}}{\mathcal{O}}^\mathrm{CS}\widetilde{\Phi}^r_{\psi,\Gamma,A}/\mathrm{Fro}^\mathbb{Z},\underset{I}{\mathrm{homotopylimit}}~\underset{\mathrm{Spec}}{\mathcal{O}}^\mathrm{CS}\widetilde{\Phi}^I_{\psi,\Gamma,A}/\mathrm{Fro}^\mathbb{Z},	\\
}
\]
\[ 
\xymatrix@R+0pc@C+0pc{
\underset{r}{\mathrm{homotopycolimit}}~\underset{\mathrm{Spec}}{\mathcal{O}}^\mathrm{CS}\breve{\Phi}^r_{\psi,\Gamma,A}/\mathrm{Fro}^\mathbb{Z},\underset{I}{\mathrm{homotopylimit}}~\breve{\Phi}^I_{\psi,\Gamma,A}/\mathrm{Fro}^\mathbb{Z},	\\
}
\]
\[ 
\xymatrix@R+0pc@C+0pc{
\underset{r}{\mathrm{homotopycolimit}}~\underset{\mathrm{Spec}}{\mathcal{O}}^\mathrm{CS}{\Phi}^r_{\psi,\Gamma,A}/\mathrm{Fro}^\mathbb{Z},\underset{I}{\mathrm{homotopylimit}}~\underset{\mathrm{Spec}}{\mathcal{O}}^\mathrm{CS}{\Phi}^I_{\psi,\Gamma,A}/\mathrm{Fro}^\mathbb{Z}.	
}
\]	
In this situation we will have the target category being family parametrized by $r$ or $I$ in compatible glueing sense as in \cite[Definition 5.4.10]{10KL2}. In this situation for modules parametrized by the intervals we have the equivalence of $\infty$-categories by using \cite[Proposition 12.18]{10CS2}. Here the corresponding quasicoherent Frobenius modules are defined to be the corresponding homotopy colimits and limits of Frobenius modules:
\begin{align}
\underset{r}{\mathrm{homotopycolimit}}~M_r,\\
\underset{I}{\mathrm{homotopylimit}}~M_I,	
\end{align}
where each $M_r$ is a Frobenius-equivariant module over the period ring with respect to some radius $r$ while each $M_I$ is a Frobenius-equivariant module over the period ring with respect to some interval $I$.\\
\end{proposition}

\begin{proposition}
Similar proposition holds for 
\begin{align}
\mathrm{Quasicoherentsheaves,Perfectcomplex,IndBanach}_{*}.	
\end{align}	
\end{proposition}

\newpage

\subsection{Some Notations}

\

\noindent [$?$] $\mathrm{Spec}^\mathrm{BK}$: Bambozzi-Kremnizer $\infty$-topos.\\
\noindent [$?$] $\mathrm{Spec}^\mathrm{CS}$: Clausen-Scholze $\infty$-topos.\\
\noindent [$?$] $\underset{\mathrm{Spec}}{\mathcal{O}}^\mathrm{BK}$: $\infty$-ringed structure sheaf of Bambozzi-Kremnizer $\infty$-topos.\\
\noindent [$?$] $\underset{\mathrm{Spec}}{\mathcal{O}}^\mathrm{CS}$: $\infty$-ringed structure sheaf of Clausen-Scholze $\infty$-topos.\\
\noindent [$?$] $\mathrm{Quasicoherentsheaves,IndBanach}_{*}$: $(\infty,1)$-category of $(\infty,1)$-quasicoherent sheaves of $(\infty,1)$-Banach modules.\\
\noindent [$?$] $\mathrm{Quasicoherentsheaves,Perfectcomplex,IndBanach}_{*}$: $(\infty,1)$-category of $(\infty,1)$-quasicoh\\
erent sheaves of $(\infty,1)$-Banach modules with certain requirement on the perfectness. \\
\noindent [$?$] $\mathrm{Quasicoherentpresheaves,Condensed}_{*}$:$(\infty,1)$-category of $(\infty,1)$-quasicoherent sheaves of $(\infty,1)$-condensed solid modules.\\
\noindent [$?$] $\mathrm{Quasicoherentpresheaves,Perfectcomplex,Condensed}_{*}$:$(\infty,1)$-category of $(\infty,1)$-quasi\\
coherent sheaves of $(\infty,1)$-condensed solid modules with certain requirement on the perfectness. \\
\noindent [$?$] $A$: Banach Rings. \\
\noindent [$?$] $-,\circ$: Deformations of Banach Rings as Functors.\\
\noindent [$?$] $X$: Preadic Spaces. \\
\noindent [$?$] $\square,X_\square$: Colimits of Rings and Stacks.\\

\section{Multivariate Hodge Iwasawa Modules}

This chapter follows closely \cite{10T1}, \cite{10T2}, \cite{10T3}, \cite{10KPX}, \cite{10KP}, \cite{10KL1}, \cite{10KL2}, \cite{10BK}, \cite{10BBBK}, \cite{10BBM}, \cite{10KKM}, \cite{10CS1}, \cite{10CS2}, \cite{10CKZ}, \cite{10PZ}, \cite{10BCM}, \cite{10LBV}\footnote{Note that all our constructions in this article are motivated by certain deformation needed in our project on the Hodge-Iwasawa theory closely after \cite{10BF1}, \cite{10BF2}, \cite{10FK}, \cite{10FS}, \cite{10He}, \cite{10KL1}, \cite{10KL2}, \cite{10KP}, \cite{10KPX},  \cite{10Na1}, \cite{10Na2}, \cite{10PR}, \cite{10RZ}, \cite{10Sch2}, \cite{10SW}, \cite{10Wit}. Also in our context one can study the deformation version of the construction in \cite{10LBV} by taking the corresponding de Rham complex of rigid motives over the spaces we consider in this paper.}.

\subsection{Frobenius Quasicoherent Modules I}

\begin{definition}
Let $\psi$ be a toric tower over $\mathbb{Q}_p$ as in \cite[Chapter 7]{10KL2} with base $\mathbb{Q}_p\left<X_1^{\pm 1},...,X_k^{\pm 1}\right>$. Then from \cite{10KL1} and \cite[Definition 5.2.1]{10KL2} we have the following class of Kedlaya-Liu rings (with the following replacement: $\Delta$ stands for $A$, $\nabla$ stands for $B$, while $\Phi$ stands for $C$) by taking product in the sense of self $\Gamma$-th power:

\[
\xymatrix@R+0pc@C+0pc{
\widetilde{\Delta}_{\psi,\Gamma},\widetilde{\nabla}_{\psi,\Gamma},\widetilde{\Phi}_{\psi,\Gamma},\widetilde{\Delta}^+_{\psi,\Gamma},\widetilde{\nabla}^+_{\psi,\Gamma},\widetilde{\Delta}^\dagger_{\psi,\Gamma},\widetilde{\nabla}^\dagger_{\psi,\Gamma},\widetilde{\Phi}^r_{\psi,\Gamma},\widetilde{\Phi}^I_{\psi,\Gamma}, 
}
\]

\[
\xymatrix@R+0pc@C+0pc{
\breve{\Delta}_{\psi,\Gamma},\breve{\nabla}_{\psi,\Gamma},\breve{\Phi}_{\psi,\Gamma},\breve{\Delta}^+_{\psi,\Gamma},\breve{\nabla}^+_{\psi,\Gamma},\breve{\Delta}^\dagger_{\psi,\Gamma},\breve{\nabla}^\dagger_{\psi,\Gamma},\breve{\Phi}^r_{\psi,\Gamma},\breve{\Phi}^I_{\psi,\Gamma},	
}
\]

\[
\xymatrix@R+0pc@C+0pc{
{\Delta}_{\psi,\Gamma},{\nabla}_{\psi,\Gamma},{\Phi}_{\psi,\Gamma},{\Delta}^+_{\psi,\Gamma},{\nabla}^+_{\psi,\Gamma},{\Delta}^\dagger_{\psi,\Gamma},{\nabla}^\dagger_{\psi,\Gamma},{\Phi}^r_{\psi,\Gamma},{\Phi}^I_{\psi,\Gamma}.	
}
\]
We now consider the following rings with $A$ being a Banach ring over $\mathbb{Q}_p$. Taking the product we have:
\[
\xymatrix@R+0pc@C+0pc{
\widetilde{\Phi}_{\psi,\Gamma,A},\widetilde{\Phi}^r_{\psi,\Gamma,A},\widetilde{\Phi}^I_{\psi,\Gamma,A},	
}
\]
\[
\xymatrix@R+0pc@C+0pc{
\breve{\Phi}_{\psi,\Gamma,A},\breve{\Phi}^r_{\psi,\Gamma,A},\breve{\Phi}^I_{\psi,\Gamma,A},	
}
\]
\[
\xymatrix@R+0pc@C+0pc{
{\Phi}_{\psi,\Gamma,A},{\Phi}^r_{\psi,\Gamma,A},{\Phi}^I_{\psi,\Gamma,A}.	
}
\]
They carry multi Frobenius action $\varphi_\Gamma$ and multi $\mathrm{Lie}_\Gamma:=\mathbb{Z}_p^{\times\Gamma}$ action. In our current situation after \cite{10CKZ} and \cite{10PZ} we consider the following $(\infty,1)$-categories of $(\infty,1)$-modules.\\
\end{definition}

\begin{definition}
First we consider the Bambozzi-Kremnizer spectrum $\mathrm{Spec}^\mathrm{BK}(*)$ attached to any of those in the above from \cite{10BK} by taking derived rational localization:
\begin{align}
&\mathrm{Spec}^\mathrm{BK}\widetilde{\Phi}_{\psi,\Gamma,A},\mathrm{Spec}^\mathrm{BK}\widetilde{\Phi}^r_{\psi,\Gamma,A},\mathrm{Spec}^\mathrm{BK}\widetilde{\Phi}^I_{\psi,\Gamma,A},	
\end{align}
\begin{align}
&\mathrm{Spec}^\mathrm{BK}\breve{\Phi}_{\psi,\Gamma,A},\mathrm{Spec}^\mathrm{BK}\breve{\Phi}^r_{\psi,\Gamma,A},\mathrm{Spec}^\mathrm{BK}\breve{\Phi}^I_{\psi,\Gamma,A},	
\end{align}
\begin{align}
&\mathrm{Spec}^\mathrm{BK}{\Phi}_{\psi,\Gamma,A},
\mathrm{Spec}^\mathrm{BK}{\Phi}^r_{\psi,\Gamma,A},\mathrm{Spec}^\mathrm{BK}{\Phi}^I_{\psi,\Gamma,A}.	
\end{align}
Then we take the corresponding quotients by using the corresponding Frobenius operators:
\begin{align}
&\mathrm{Spec}^\mathrm{BK}\widetilde{\Phi}_{\psi,\Gamma,A}/\mathrm{Fro}^\mathbb{Z},	\\
\end{align}
\begin{align}
&\mathrm{Spec}^\mathrm{BK}\breve{\Phi}_{\psi,\Gamma,A}/\mathrm{Fro}^\mathbb{Z},	\\
\end{align}
\begin{align}
&\mathrm{Spec}^\mathrm{BK}{\Phi}_{\psi,\Gamma,A}/\mathrm{Fro}^\mathbb{Z}.	
\end{align}
Here for those space without notation related to the radius and the corresponding interval we consider the total unions $\bigcap_r,\bigcup_I$ in order to achieve the whole spaces to achieve the analogues of the corresponding FF curves from \cite{10KL1}, \cite{10KL2}, \cite{10FF} for
\[
\xymatrix@R+0pc@C+0pc{
\underset{r}{\mathrm{homotopylimit}}~\mathrm{Spec}^\mathrm{BK}\widetilde{\Phi}^r_{\psi,\Gamma,A},\underset{I}{\mathrm{homotopycolimit}}~\mathrm{Spec}^\mathrm{BK}\widetilde{\Phi}^I_{\psi,\Gamma,A},	\\
}
\]
\[
\xymatrix@R+0pc@C+0pc{
\underset{r}{\mathrm{homotopylimit}}~\mathrm{Spec}^\mathrm{BK}\breve{\Phi}^r_{\psi,\Gamma,A},\underset{I}{\mathrm{homotopycolimit}}~\mathrm{Spec}^\mathrm{BK}\breve{\Phi}^I_{\psi,\Gamma,A},	\\
}
\]
\[
\xymatrix@R+0pc@C+0pc{
\underset{r}{\mathrm{homotopylimit}}~\mathrm{Spec}^\mathrm{BK}{\Phi}^r_{\psi,\Gamma,A},\underset{I}{\mathrm{homotopycolimit}}~\mathrm{Spec}^\mathrm{BK}{\Phi}^I_{\psi,\Gamma,A}.	
}
\]
\[  
\xymatrix@R+0pc@C+0pc{
\underset{r}{\mathrm{homotopylimit}}~\mathrm{Spec}^\mathrm{BK}\widetilde{\Phi}^r_{\psi,\Gamma,A}/\mathrm{Fro}^\mathbb{Z},\underset{I}{\mathrm{homotopycolimit}}~\mathrm{Spec}^\mathrm{BK}\widetilde{\Phi}^I_{\psi,\Gamma,A}/\mathrm{Fro}^\mathbb{Z},	\\
}
\]
\[ 
\xymatrix@R+0pc@C+0pc{
\underset{r}{\mathrm{homotopylimit}}~\mathrm{Spec}^\mathrm{BK}\breve{\Phi}^r_{\psi,\Gamma,A}/\mathrm{Fro}^\mathbb{Z},\underset{I}{\mathrm{homotopycolimit}}~\mathrm{Spec}^\mathrm{BK}\breve{\Phi}^I_{\psi,\Gamma,A}/\mathrm{Fro}^\mathbb{Z},	\\
}
\]
\[ 
\xymatrix@R+0pc@C+0pc{
\underset{r}{\mathrm{homotopylimit}}~\mathrm{Spec}^\mathrm{BK}{\Phi}^r_{\psi,\Gamma,A}/\mathrm{Fro}^\mathbb{Z},\underset{I}{\mathrm{homotopycolimit}}~\mathrm{Spec}^\mathrm{BK}{\Phi}^I_{\psi,\Gamma,A}/\mathrm{Fro}^\mathbb{Z}.	
}
\]

\end{definition}

\indent Meanwhile we have the corresponding Clausen-Scholze analytic stacks from \cite{10CS2}, therefore applying their construction we have:

\begin{definition}
Here we define the following products by using the solidified tensor product from \cite{10CS1} and \cite{10CS2}. Namely $A$ will still as above as a Banach ring over $\mathbb{Q}_p$. Then we take solidified tensor product $\overset{\blacksquare}{\otimes}$ of any of the following
\[
\xymatrix@R+0pc@C+0pc{
\widetilde{\Delta}_{\psi,\Gamma},\widetilde{\nabla}_{\psi,\Gamma},\widetilde{\Phi}_{\psi,\Gamma},\widetilde{\Delta}^+_{\psi,\Gamma},\widetilde{\nabla}^+_{\psi,\Gamma},\widetilde{\Delta}^\dagger_{\psi,\Gamma},\widetilde{\nabla}^\dagger_{\psi,\Gamma},\widetilde{\Phi}^r_{\psi,\Gamma},\widetilde{\Phi}^I_{\psi,\Gamma}, 
}
\]

\[
\xymatrix@R+0pc@C+0pc{
\breve{\Delta}_{\psi,\Gamma},\breve{\nabla}_{\psi,\Gamma},\breve{\Phi}_{\psi,\Gamma},\breve{\Delta}^+_{\psi,\Gamma},\breve{\nabla}^+_{\psi,\Gamma},\breve{\Delta}^\dagger_{\psi,\Gamma},\breve{\nabla}^\dagger_{\psi,\Gamma},\breve{\Phi}^r_{\psi,\Gamma},\breve{\Phi}^I_{\psi,\Gamma},	
}
\]

\[
\xymatrix@R+0pc@C+0pc{
{\Delta}_{\psi,\Gamma},{\nabla}_{\psi,\Gamma},{\Phi}_{\psi,\Gamma},{\Delta}^+_{\psi,\Gamma},{\nabla}^+_{\psi,\Gamma},{\Delta}^\dagger_{\psi,\Gamma},{\nabla}^\dagger_{\psi,\Gamma},{\Phi}^r_{\psi,\Gamma},{\Phi}^I_{\psi,\Gamma},	
}
\]  	
with $A$. Then we have the notations:
\[
\xymatrix@R+0pc@C+0pc{
\widetilde{\Delta}_{\psi,\Gamma,A},\widetilde{\nabla}_{\psi,\Gamma,A},\widetilde{\Phi}_{\psi,\Gamma,A},\widetilde{\Delta}^+_{\psi,\Gamma,A},\widetilde{\nabla}^+_{\psi,\Gamma,A},\widetilde{\Delta}^\dagger_{\psi,\Gamma,A},\widetilde{\nabla}^\dagger_{\psi,\Gamma,A},\widetilde{\Phi}^r_{\psi,\Gamma,A},\widetilde{\Phi}^I_{\psi,\Gamma,A}, 
}
\]

\[
\xymatrix@R+0pc@C+0pc{
\breve{\Delta}_{\psi,\Gamma,A},\breve{\nabla}_{\psi,\Gamma,A},\breve{\Phi}_{\psi,\Gamma,A},\breve{\Delta}^+_{\psi,\Gamma,A},\breve{\nabla}^+_{\psi,\Gamma,A},\breve{\Delta}^\dagger_{\psi,\Gamma,A},\breve{\nabla}^\dagger_{\psi,\Gamma,A},\breve{\Phi}^r_{\psi,\Gamma,A},\breve{\Phi}^I_{\psi,\Gamma,A},	
}
\]

\[
\xymatrix@R+0pc@C+0pc{
{\Delta}_{\psi,\Gamma,A},{\nabla}_{\psi,\Gamma,A},{\Phi}_{\psi,\Gamma,A},{\Delta}^+_{\psi,\Gamma,A},{\nabla}^+_{\psi,\Gamma,A},{\Delta}^\dagger_{\psi,\Gamma,A},{\nabla}^\dagger_{\psi,\Gamma,A},{\Phi}^r_{\psi,\Gamma,A},{\Phi}^I_{\psi,\Gamma,A}.	
}
\]
\end{definition}

\begin{definition}
First we consider the Clausen-Scholze spectrum $\mathrm{Spec}^\mathrm{CS}(*)$ attached to any of those in the above from \cite{10CS2} by taking derived rational localization:
\begin{align}
\mathrm{Spec}^\mathrm{CS}\widetilde{\Delta}_{\psi,\Gamma,A},\mathrm{Spec}^\mathrm{CS}\widetilde{\nabla}_{\psi,\Gamma,A},\mathrm{Spec}^\mathrm{CS}\widetilde{\Phi}_{\psi,\Gamma,A},\mathrm{Spec}^\mathrm{CS}\widetilde{\Delta}^+_{\psi,\Gamma,A},\mathrm{Spec}^\mathrm{CS}\widetilde{\nabla}^+_{\psi,\Gamma,A},\\
\mathrm{Spec}^\mathrm{CS}\widetilde{\Delta}^\dagger_{\psi,\Gamma,A},\mathrm{Spec}^\mathrm{CS}\widetilde{\nabla}^\dagger_{\psi,\Gamma,A},\mathrm{Spec}^\mathrm{CS}\widetilde{\Phi}^r_{\psi,\Gamma,A},\mathrm{Spec}^\mathrm{CS}\widetilde{\Phi}^I_{\psi,\Gamma,A},	\\
\end{align}
\begin{align}
\mathrm{Spec}^\mathrm{CS}\breve{\Delta}_{\psi,\Gamma,A},\breve{\nabla}_{\psi,\Gamma,A},\mathrm{Spec}^\mathrm{CS}\breve{\Phi}_{\psi,\Gamma,A},\mathrm{Spec}^\mathrm{CS}\breve{\Delta}^+_{\psi,\Gamma,A},\mathrm{Spec}^\mathrm{CS}\breve{\nabla}^+_{\psi,\Gamma,A},\\
\mathrm{Spec}^\mathrm{CS}\breve{\Delta}^\dagger_{\psi,\Gamma,A},\mathrm{Spec}^\mathrm{CS}\breve{\nabla}^\dagger_{\psi,\Gamma,A},\mathrm{Spec}^\mathrm{CS}\breve{\Phi}^r_{\psi,\Gamma,A},\breve{\Phi}^I_{\psi,\Gamma,A},	\\
\end{align}
\begin{align}
\mathrm{Spec}^\mathrm{CS}{\Delta}_{\psi,\Gamma,A},\mathrm{Spec}^\mathrm{CS}{\nabla}_{\psi,\Gamma,A},\mathrm{Spec}^\mathrm{CS}{\Phi}_{\psi,\Gamma,A},\mathrm{Spec}^\mathrm{CS}{\Delta}^+_{\psi,\Gamma,A},\mathrm{Spec}^\mathrm{CS}{\nabla}^+_{\psi,\Gamma,A},\\
\mathrm{Spec}^\mathrm{CS}{\Delta}^\dagger_{\psi,\Gamma,A},\mathrm{Spec}^\mathrm{CS}{\nabla}^\dagger_{\psi,\Gamma,A},\mathrm{Spec}^\mathrm{CS}{\Phi}^r_{\psi,\Gamma,A},\mathrm{Spec}^\mathrm{CS}{\Phi}^I_{\psi,\Gamma,A}.	
\end{align}

Then we take the corresponding quotients by using the corresponding Frobenius operators:
\begin{align}
&\mathrm{Spec}^\mathrm{CS}\widetilde{\Delta}_{\psi,\Gamma,A}/\mathrm{Fro}^\mathbb{Z},\mathrm{Spec}^\mathrm{CS}\widetilde{\nabla}_{\psi,\Gamma,A}/\mathrm{Fro}^\mathbb{Z},\mathrm{Spec}^\mathrm{CS}\widetilde{\Phi}_{\psi,\Gamma,A}/\mathrm{Fro}^\mathbb{Z},\mathrm{Spec}^\mathrm{CS}\widetilde{\Delta}^+_{\psi,\Gamma,A}/\mathrm{Fro}^\mathbb{Z},\\
&\mathrm{Spec}^\mathrm{CS}\widetilde{\nabla}^+_{\psi,\Gamma,A}/\mathrm{Fro}^\mathbb{Z}, \mathrm{Spec}^\mathrm{CS}\widetilde{\Delta}^\dagger_{\psi,\Gamma,A}/\mathrm{Fro}^\mathbb{Z},\mathrm{Spec}^\mathrm{CS}\widetilde{\nabla}^\dagger_{\psi,\Gamma,A}/\mathrm{Fro}^\mathbb{Z},	\\
\end{align}
\begin{align}
&\mathrm{Spec}^\mathrm{CS}\breve{\Delta}_{\psi,\Gamma,A}/\mathrm{Fro}^\mathbb{Z},\breve{\nabla}_{\psi,\Gamma,A}/\mathrm{Fro}^\mathbb{Z},\mathrm{Spec}^\mathrm{CS}\breve{\Phi}_{\psi,\Gamma,A}/\mathrm{Fro}^\mathbb{Z},\mathrm{Spec}^\mathrm{CS}\breve{\Delta}^+_{\psi,\Gamma,A}/\mathrm{Fro}^\mathbb{Z},\\
&\mathrm{Spec}^\mathrm{CS}\breve{\nabla}^+_{\psi,\Gamma,A}/\mathrm{Fro}^\mathbb{Z}, \mathrm{Spec}^\mathrm{CS}\breve{\Delta}^\dagger_{\psi,\Gamma,A}/\mathrm{Fro}^\mathbb{Z},\mathrm{Spec}^\mathrm{CS}\breve{\nabla}^\dagger_{\psi,\Gamma,A}/\mathrm{Fro}^\mathbb{Z},	\\
\end{align}
\begin{align}
&\mathrm{Spec}^\mathrm{CS}{\Delta}_{\psi,\Gamma,A}/\mathrm{Fro}^\mathbb{Z},\mathrm{Spec}^\mathrm{CS}{\nabla}_{\psi,\Gamma,A}/\mathrm{Fro}^\mathbb{Z},\mathrm{Spec}^\mathrm{CS}{\Phi}_{\psi,\Gamma,A}/\mathrm{Fro}^\mathbb{Z},\mathrm{Spec}^\mathrm{CS}{\Delta}^+_{\psi,\Gamma,A}/\mathrm{Fro}^\mathbb{Z},\\
&\mathrm{Spec}^\mathrm{CS}{\nabla}^+_{\psi,\Gamma,A}/\mathrm{Fro}^\mathbb{Z}, \mathrm{Spec}^\mathrm{CS}{\Delta}^\dagger_{\psi,\Gamma,A}/\mathrm{Fro}^\mathbb{Z},\mathrm{Spec}^\mathrm{CS}{\nabla}^\dagger_{\psi,\Gamma,A}/\mathrm{Fro}^\mathbb{Z}.	
\end{align}
Here for those space with notations related to the radius and the corresponding interval we consider the total unions $\bigcap_r,\bigcup_I$ in order to achieve the whole spaces to achieve the analogues of the corresponding FF curves from \cite{10KL1}, \cite{10KL2}, \cite{10FF} for
\[
\xymatrix@R+0pc@C+0pc{
\underset{r}{\mathrm{homotopylimit}}~\mathrm{Spec}^\mathrm{CS}\widetilde{\Phi}^r_{\psi,\Gamma,A},\underset{I}{\mathrm{homotopycolimit}}~\mathrm{Spec}^\mathrm{CS}\widetilde{\Phi}^I_{\psi,\Gamma,A},	\\
}
\]
\[
\xymatrix@R+0pc@C+0pc{
\underset{r}{\mathrm{homotopylimit}}~\mathrm{Spec}^\mathrm{CS}\breve{\Phi}^r_{\psi,\Gamma,A},\underset{I}{\mathrm{homotopycolimit}}~\mathrm{Spec}^\mathrm{CS}\breve{\Phi}^I_{\psi,\Gamma,A},	\\
}
\]
\[
\xymatrix@R+0pc@C+0pc{
\underset{r}{\mathrm{homotopylimit}}~\mathrm{Spec}^\mathrm{CS}{\Phi}^r_{\psi,\Gamma,A},\underset{I}{\mathrm{homotopycolimit}}~\mathrm{Spec}^\mathrm{CS}{\Phi}^I_{\psi,\Gamma,A}.	
}
\]
\[ 
\xymatrix@R+0pc@C+0pc{
\underset{r}{\mathrm{homotopylimit}}~\mathrm{Spec}^\mathrm{CS}\widetilde{\Phi}^r_{\psi,\Gamma,A}/\mathrm{Fro}^\mathbb{Z},\underset{I}{\mathrm{homotopycolimit}}~\mathrm{Spec}^\mathrm{CS}\widetilde{\Phi}^I_{\psi,\Gamma,A}/\mathrm{Fro}^\mathbb{Z},	\\
}
\]
\[ 
\xymatrix@R+0pc@C+0pc{
\underset{r}{\mathrm{homotopylimit}}~\mathrm{Spec}^\mathrm{CS}\breve{\Phi}^r_{\psi,\Gamma,A}/\mathrm{Fro}^\mathbb{Z},\underset{I}{\mathrm{homotopycolimit}}~\breve{\Phi}^I_{\psi,\Gamma,A}/\mathrm{Fro}^\mathbb{Z},	\\
}
\]
\[ 
\xymatrix@R+0pc@C+0pc{
\underset{r}{\mathrm{homotopylimit}}~\mathrm{Spec}^\mathrm{CS}{\Phi}^r_{\psi,\Gamma,A}/\mathrm{Fro}^\mathbb{Z},\underset{I}{\mathrm{homotopycolimit}}~\mathrm{Spec}^\mathrm{CS}{\Phi}^I_{\psi,\Gamma,A}/\mathrm{Fro}^\mathbb{Z}.	
}
\]

\end{definition}

\

\begin{definition}
We then consider the corresponding quasipresheaves of the corresponding ind-Banach or monomorphic ind-Banach modules from \cite{10BBK}, \cite{10KKM}:
\begin{align}
\mathrm{Quasicoherentpresheaves,IndBanach}_{*}	
\end{align}
where $*$ is one of the following spaces:
\begin{align}
&\mathrm{Spec}^\mathrm{BK}\widetilde{\Phi}_{\psi,\Gamma,A}/\mathrm{Fro}^\mathbb{Z},	\\
\end{align}
\begin{align}
&\mathrm{Spec}^\mathrm{BK}\breve{\Phi}_{\psi,\Gamma,A}/\mathrm{Fro}^\mathbb{Z},	\\
\end{align}
\begin{align}
&\mathrm{Spec}^\mathrm{BK}{\Phi}_{\psi,\Gamma,A}/\mathrm{Fro}^\mathbb{Z}.	
\end{align}
Here for those space without notation related to the radius and the corresponding interval we consider the total unions $\bigcap_r,\bigcup_I$ in order to achieve the whole spaces to achieve the analogues of the corresponding FF curves from \cite{10KL1}, \cite{10KL2}, \cite{10FF} for
\[
\xymatrix@R+0pc@C+0pc{
\underset{r}{\mathrm{homotopylimit}}~\mathrm{Spec}^\mathrm{BK}\widetilde{\Phi}^r_{\psi,\Gamma,A},\underset{I}{\mathrm{homotopycolimit}}~\mathrm{Spec}^\mathrm{BK}\widetilde{\Phi}^I_{\psi,\Gamma,A},	\\
}
\]
\[
\xymatrix@R+0pc@C+0pc{
\underset{r}{\mathrm{homotopylimit}}~\mathrm{Spec}^\mathrm{BK}\breve{\Phi}^r_{\psi,\Gamma,A},\underset{I}{\mathrm{homotopycolimit}}~\mathrm{Spec}^\mathrm{BK}\breve{\Phi}^I_{\psi,\Gamma,A},	\\
}
\]
\[
\xymatrix@R+0pc@C+0pc{
\underset{r}{\mathrm{homotopylimit}}~\mathrm{Spec}^\mathrm{BK}{\Phi}^r_{\psi,\Gamma,A},\underset{I}{\mathrm{homotopycolimit}}~\mathrm{Spec}^\mathrm{BK}{\Phi}^I_{\psi,\Gamma,A}.	
}
\]
\[  
\xymatrix@R+0pc@C+0pc{
\underset{r}{\mathrm{homotopylimit}}~\mathrm{Spec}^\mathrm{BK}\widetilde{\Phi}^r_{\psi,\Gamma,A}/\mathrm{Fro}^\mathbb{Z},\underset{I}{\mathrm{homotopycolimit}}~\mathrm{Spec}^\mathrm{BK}\widetilde{\Phi}^I_{\psi,\Gamma,A}/\mathrm{Fro}^\mathbb{Z},	\\
}
\]
\[ 
\xymatrix@R+0pc@C+0pc{
\underset{r}{\mathrm{homotopylimit}}~\mathrm{Spec}^\mathrm{BK}\breve{\Phi}^r_{\psi,\Gamma,A}/\mathrm{Fro}^\mathbb{Z},\underset{I}{\mathrm{homotopycolimit}}~\mathrm{Spec}^\mathrm{BK}\breve{\Phi}^I_{\psi,\Gamma,A}/\mathrm{Fro}^\mathbb{Z},	\\
}
\]
\[ 
\xymatrix@R+0pc@C+0pc{
\underset{r}{\mathrm{homotopylimit}}~\mathrm{Spec}^\mathrm{BK}{\Phi}^r_{\psi,\Gamma,A}/\mathrm{Fro}^\mathbb{Z},\underset{I}{\mathrm{homotopycolimit}}~\mathrm{Spec}^\mathrm{BK}{\Phi}^I_{\psi,\Gamma,A}/\mathrm{Fro}^\mathbb{Z}.	
}
\]

\end{definition}

\begin{definition}
We then consider the corresponding quasisheaves of the corresponding condensed solid topological modules from \cite{10CS2}:
\begin{align}
\mathrm{Quasicoherentsheaves, Condensed}_{*}	
\end{align}
where $*$ is one of the following spaces:
\begin{align}
&\mathrm{Spec}^\mathrm{CS}\widetilde{\Delta}_{\psi,\Gamma,A}/\mathrm{Fro}^\mathbb{Z},\mathrm{Spec}^\mathrm{CS}\widetilde{\nabla}_{\psi,\Gamma,A}/\mathrm{Fro}^\mathbb{Z},\mathrm{Spec}^\mathrm{CS}\widetilde{\Phi}_{\psi,\Gamma,A}/\mathrm{Fro}^\mathbb{Z},\mathrm{Spec}^\mathrm{CS}\widetilde{\Delta}^+_{\psi,\Gamma,A}/\mathrm{Fro}^\mathbb{Z},\\
&\mathrm{Spec}^\mathrm{CS}\widetilde{\nabla}^+_{\psi,\Gamma,A}/\mathrm{Fro}^\mathbb{Z},\mathrm{Spec}^\mathrm{CS}\widetilde{\Delta}^\dagger_{\psi,\Gamma,A}/\mathrm{Fro}^\mathbb{Z},\mathrm{Spec}^\mathrm{CS}\widetilde{\nabla}^\dagger_{\psi,\Gamma,A}/\mathrm{Fro}^\mathbb{Z},	\\
\end{align}
\begin{align}
&\mathrm{Spec}^\mathrm{CS}\breve{\Delta}_{\psi,\Gamma,A}/\mathrm{Fro}^\mathbb{Z},\breve{\nabla}_{\psi,\Gamma,A}/\mathrm{Fro}^\mathbb{Z},\mathrm{Spec}^\mathrm{CS}\breve{\Phi}_{\psi,\Gamma,A}/\mathrm{Fro}^\mathbb{Z},\mathrm{Spec}^\mathrm{CS}\breve{\Delta}^+_{\psi,\Gamma,A}/\mathrm{Fro}^\mathbb{Z},\\
&\mathrm{Spec}^\mathrm{CS}\breve{\nabla}^+_{\psi,\Gamma,A}/\mathrm{Fro}^\mathbb{Z},\mathrm{Spec}^\mathrm{CS}\breve{\Delta}^\dagger_{\psi,\Gamma,A}/\mathrm{Fro}^\mathbb{Z},\mathrm{Spec}^\mathrm{CS}\breve{\nabla}^\dagger_{\psi,\Gamma,A}/\mathrm{Fro}^\mathbb{Z},	\\
\end{align}
\begin{align}
&\mathrm{Spec}^\mathrm{CS}{\Delta}_{\psi,\Gamma,A}/\mathrm{Fro}^\mathbb{Z},\mathrm{Spec}^\mathrm{CS}{\nabla}_{\psi,\Gamma,A}/\mathrm{Fro}^\mathbb{Z},\mathrm{Spec}^\mathrm{CS}{\Phi}_{\psi,\Gamma,A}/\mathrm{Fro}^\mathbb{Z},\mathrm{Spec}^\mathrm{CS}{\Delta}^+_{\psi,\Gamma,A}/\mathrm{Fro}^\mathbb{Z},\\
&\mathrm{Spec}^\mathrm{CS}{\nabla}^+_{\psi,\Gamma,A}/\mathrm{Fro}^\mathbb{Z}, \mathrm{Spec}^\mathrm{CS}{\Delta}^\dagger_{\psi,\Gamma,A}/\mathrm{Fro}^\mathbb{Z},\mathrm{Spec}^\mathrm{CS}{\nabla}^\dagger_{\psi,\Gamma,A}/\mathrm{Fro}^\mathbb{Z}.	
\end{align}
Here for those space with notations related to the radius and the corresponding interval we consider the total unions $\bigcap_r,\bigcup_I$ in order to achieve the whole spaces to achieve the analogues of the corresponding FF curves from \cite{10KL1}, \cite{10KL2}, \cite{10FF} for
\[
\xymatrix@R+0pc@C+0pc{
\underset{r}{\mathrm{homotopylimit}}~\mathrm{Spec}^\mathrm{CS}\widetilde{\Phi}^r_{\psi,\Gamma,A},\underset{I}{\mathrm{homotopycolimit}}~\mathrm{Spec}^\mathrm{CS}\widetilde{\Phi}^I_{\psi,\Gamma,A},	\\
}
\]
\[
\xymatrix@R+0pc@C+0pc{
\underset{r}{\mathrm{homotopylimit}}~\mathrm{Spec}^\mathrm{CS}\breve{\Phi}^r_{\psi,\Gamma,A},\underset{I}{\mathrm{homotopycolimit}}~\mathrm{Spec}^\mathrm{CS}\breve{\Phi}^I_{\psi,\Gamma,A},	\\
}
\]
\[
\xymatrix@R+0pc@C+0pc{
\underset{r}{\mathrm{homotopylimit}}~\mathrm{Spec}^\mathrm{CS}{\Phi}^r_{\psi,\Gamma,A},\underset{I}{\mathrm{homotopycolimit}}~\mathrm{Spec}^\mathrm{CS}{\Phi}^I_{\psi,\Gamma,A}.	
}
\]
\[ 
\xymatrix@R+0pc@C+0pc{
\underset{r}{\mathrm{homotopylimit}}~\mathrm{Spec}^\mathrm{CS}\widetilde{\Phi}^r_{\psi,\Gamma,A}/\mathrm{Fro}^\mathbb{Z},\underset{I}{\mathrm{homotopycolimit}}~\mathrm{Spec}^\mathrm{CS}\widetilde{\Phi}^I_{\psi,\Gamma,A}/\mathrm{Fro}^\mathbb{Z},	\\
}
\]
\[ 
\xymatrix@R+0pc@C+0pc{
\underset{r}{\mathrm{homotopylimit}}~\mathrm{Spec}^\mathrm{CS}\breve{\Phi}^r_{\psi,\Gamma,A}/\mathrm{Fro}^\mathbb{Z},\underset{I}{\mathrm{homotopycolimit}}~\breve{\Phi}^I_{\psi,\Gamma,A}/\mathrm{Fro}^\mathbb{Z},	\\
}
\]
\[ 
\xymatrix@R+0pc@C+0pc{
\underset{r}{\mathrm{homotopylimit}}~\mathrm{Spec}^\mathrm{CS}{\Phi}^r_{\psi,\Gamma,A}/\mathrm{Fro}^\mathbb{Z},\underset{I}{\mathrm{homotopycolimit}}~\mathrm{Spec}^\mathrm{CS}{\Phi}^I_{\psi,\Gamma,A}/\mathrm{Fro}^\mathbb{Z}.	
}
\]

\end{definition}

\

\begin{proposition}
There is a well-defined functor from the $\infty$-category 
\begin{align}
\mathrm{Quasicoherentpresheaves,Condensed}_{*}	
\end{align}
where $*$ is one of the following spaces:
\begin{align}
&\mathrm{Spec}^\mathrm{CS}\widetilde{\Phi}_{\psi,\Gamma,A}/\mathrm{Fro}^\mathbb{Z},	\\
\end{align}
\begin{align}
&\mathrm{Spec}^\mathrm{CS}\breve{\Phi}_{\psi,\Gamma,A}/\mathrm{Fro}^\mathbb{Z},	\\
\end{align}
\begin{align}
&\mathrm{Spec}^\mathrm{CS}{\Phi}_{\psi,\Gamma,A}/\mathrm{Fro}^\mathbb{Z},	
\end{align}
to the $\infty$-category of $\mathrm{Fro}$-equivariant quasicoherent presheaves over similar spaces above correspondingly without the $\mathrm{Fro}$-quotients, and to the $\infty$-category of $\mathrm{Fro}$-equivariant quasicoherent modules over global sections of the structure $\infty$-sheaves of the similar spaces above correspondingly without the $\mathrm{Fro}$-quotients. Here for those space without notation related to the radius and the corresponding interval we consider the total unions $\bigcap_r,\bigcup_I$ in order to achieve the whole spaces to achieve the analogues of the corresponding FF curves from \cite{10KL1}, \cite{10KL2}, \cite{10FF} for
\[
\xymatrix@R+0pc@C+0pc{
\underset{r}{\mathrm{homotopylimit}}~\mathrm{Spec}^\mathrm{CS}\widetilde{\Phi}^r_{\psi,\Gamma,A},\underset{I}{\mathrm{homotopycolimit}}~\mathrm{Spec}^\mathrm{CS}\widetilde{\Phi}^I_{\psi,\Gamma,A},	\\
}
\]
\[
\xymatrix@R+0pc@C+0pc{
\underset{r}{\mathrm{homotopylimit}}~\mathrm{Spec}^\mathrm{CS}\breve{\Phi}^r_{\psi,\Gamma,A},\underset{I}{\mathrm{homotopycolimit}}~\mathrm{Spec}^\mathrm{CS}\breve{\Phi}^I_{\psi,\Gamma,A},	\\
}
\]
\[
\xymatrix@R+0pc@C+0pc{
\underset{r}{\mathrm{homotopylimit}}~\mathrm{Spec}^\mathrm{CS}{\Phi}^r_{\psi,\Gamma,A},\underset{I}{\mathrm{homotopycolimit}}~\mathrm{Spec}^\mathrm{CS}{\Phi}^I_{\psi,\Gamma,A}.	
}
\]
\[ 
\xymatrix@R+0pc@C+0pc{
\underset{r}{\mathrm{homotopylimit}}~\mathrm{Spec}^\mathrm{CS}\widetilde{\Phi}^r_{\psi,\Gamma,A}/\mathrm{Fro}^\mathbb{Z},\underset{I}{\mathrm{homotopycolimit}}~\mathrm{Spec}^\mathrm{CS}\widetilde{\Phi}^I_{\psi,\Gamma,A}/\mathrm{Fro}^\mathbb{Z},	\\
}
\]
\[ 
\xymatrix@R+0pc@C+0pc{
\underset{r}{\mathrm{homotopylimit}}~\mathrm{Spec}^\mathrm{CS}\breve{\Phi}^r_{\psi,\Gamma,A}/\mathrm{Fro}^\mathbb{Z},\underset{I}{\mathrm{homotopycolimit}}~\breve{\Phi}^I_{\psi,\Gamma,A}/\mathrm{Fro}^\mathbb{Z},	\\
}
\]
\[ 
\xymatrix@R+0pc@C+0pc{
\underset{r}{\mathrm{homotopylimit}}~\mathrm{Spec}^\mathrm{CS}{\Phi}^r_{\psi,\Gamma,A}/\mathrm{Fro}^\mathbb{Z},\underset{I}{\mathrm{homotopycolimit}}~\mathrm{Spec}^\mathrm{CS}{\Phi}^I_{\psi,\Gamma,A}/\mathrm{Fro}^\mathbb{Z}.	
}
\]	
In this situation we will have the target category being family parametrized by $r$ or $I$ in compatible glueing sense as in \cite[Definition 5.4.10]{10KL2}. In this situation for modules parametrized by the intervals we have the equivalence of $\infty$-categories by using \cite[Proposition 13.8]{10CS2}. Here the corresponding quasicoherent Frobenius modules are defined to be the corresponding homotopy colimits and limits of Frobenius modules:
\begin{align}
\underset{r}{\mathrm{homotopycolimit}}~M_r,\\
\underset{I}{\mathrm{homotopylimit}}~M_I,	
\end{align}
where each $M_r$ is a Frobenius-equivariant module over the period ring with respect to some radius $r$ while each $M_I$ is a Frobenius-equivariant module over the period ring with respect to some interval $I$.\\
\end{proposition}

\begin{proposition}
Similar proposition holds for 
\begin{align}
\mathrm{Quasicoherentsheaves,IndBanach}_{*}.	
\end{align}	
\end{proposition}

\

\begin{definition}
We then consider the corresponding quasipresheaves of perfect complexes the corresponding ind-Banach or monomorphic ind-Banach modules from \cite{10BBK}, \cite{10KKM}:
\begin{align}
\mathrm{Quasicoherentpresheaves,Perfectcomplex,IndBanach}_{*}	
\end{align}
where $*$ is one of the following spaces:
\begin{align}
&\mathrm{Spec}^\mathrm{BK}\widetilde{\Phi}_{\psi,\Gamma,A}/\mathrm{Fro}^\mathbb{Z},	\\
\end{align}
\begin{align}
&\mathrm{Spec}^\mathrm{BK}\breve{\Phi}_{\psi,\Gamma,A}/\mathrm{Fro}^\mathbb{Z},	\\
\end{align}
\begin{align}
&\mathrm{Spec}^\mathrm{BK}{\Phi}_{\psi,\Gamma,A}/\mathrm{Fro}^\mathbb{Z}.	
\end{align}
Here for those space without notation related to the radius and the corresponding interval we consider the total unions $\bigcap_r,\bigcup_I$ in order to achieve the whole spaces to achieve the analogues of the corresponding FF curves from \cite{10KL1}, \cite{10KL2}, \cite{10FF} for
\[
\xymatrix@R+0pc@C+0pc{
\underset{r}{\mathrm{homotopylimit}}~\mathrm{Spec}^\mathrm{BK}\widetilde{\Phi}^r_{\psi,\Gamma,A},\underset{I}{\mathrm{homotopycolimit}}~\mathrm{Spec}^\mathrm{BK}\widetilde{\Phi}^I_{\psi,\Gamma,A},	\\
}
\]
\[
\xymatrix@R+0pc@C+0pc{
\underset{r}{\mathrm{homotopylimit}}~\mathrm{Spec}^\mathrm{BK}\breve{\Phi}^r_{\psi,\Gamma,A},\underset{I}{\mathrm{homotopycolimit}}~\mathrm{Spec}^\mathrm{BK}\breve{\Phi}^I_{\psi,\Gamma,A},	\\
}
\]
\[
\xymatrix@R+0pc@C+0pc{
\underset{r}{\mathrm{homotopylimit}}~\mathrm{Spec}^\mathrm{BK}{\Phi}^r_{\psi,\Gamma,A},\underset{I}{\mathrm{homotopycolimit}}~\mathrm{Spec}^\mathrm{BK}{\Phi}^I_{\psi,\Gamma,A}.	
}
\]
\[  
\xymatrix@R+0pc@C+0pc{
\underset{r}{\mathrm{homotopylimit}}~\mathrm{Spec}^\mathrm{BK}\widetilde{\Phi}^r_{\psi,\Gamma,A}/\mathrm{Fro}^\mathbb{Z},\underset{I}{\mathrm{homotopycolimit}}~\mathrm{Spec}^\mathrm{BK}\widetilde{\Phi}^I_{\psi,\Gamma,A}/\mathrm{Fro}^\mathbb{Z},	\\
}
\]
\[ 
\xymatrix@R+0pc@C+0pc{
\underset{r}{\mathrm{homotopylimit}}~\mathrm{Spec}^\mathrm{BK}\breve{\Phi}^r_{\psi,\Gamma,A}/\mathrm{Fro}^\mathbb{Z},\underset{I}{\mathrm{homotopycolimit}}~\mathrm{Spec}^\mathrm{BK}\breve{\Phi}^I_{\psi,\Gamma,A}/\mathrm{Fro}^\mathbb{Z},	\\
}
\]
\[ 
\xymatrix@R+0pc@C+0pc{
\underset{r}{\mathrm{homotopylimit}}~\mathrm{Spec}^\mathrm{BK}{\Phi}^r_{\psi,\Gamma,A}/\mathrm{Fro}^\mathbb{Z},\underset{I}{\mathrm{homotopycolimit}}~\mathrm{Spec}^\mathrm{BK}{\Phi}^I_{\psi,\Gamma,A}/\mathrm{Fro}^\mathbb{Z}.	
}
\]

\end{definition}

\begin{definition}
We then consider the corresponding quasisheaves of perfect complexes of the corresponding condensed solid topological modules from \cite{10CS2}:
\begin{align}
\mathrm{Quasicoherentsheaves, Perfectcomplex, Condensed}_{*}	
\end{align}
where $*$ is one of the following spaces:
\begin{align}
&\mathrm{Spec}^\mathrm{CS}\widetilde{\Delta}_{\psi,\Gamma,A}/\mathrm{Fro}^\mathbb{Z},\mathrm{Spec}^\mathrm{CS}\widetilde{\nabla}_{\psi,\Gamma,A}/\mathrm{Fro}^\mathbb{Z},\mathrm{Spec}^\mathrm{CS}\widetilde{\Phi}_{\psi,\Gamma,A}/\mathrm{Fro}^\mathbb{Z},\mathrm{Spec}^\mathrm{CS}\widetilde{\Delta}^+_{\psi,\Gamma,A}/\mathrm{Fro}^\mathbb{Z},\\
&\mathrm{Spec}^\mathrm{CS}\widetilde{\nabla}^+_{\psi,\Gamma,A}/\mathrm{Fro}^\mathbb{Z},\mathrm{Spec}^\mathrm{CS}\widetilde{\Delta}^\dagger_{\psi,\Gamma,A}/\mathrm{Fro}^\mathbb{Z},\mathrm{Spec}^\mathrm{CS}\widetilde{\nabla}^\dagger_{\psi,\Gamma,A}/\mathrm{Fro}^\mathbb{Z},	\\
\end{align}
\begin{align}
&\mathrm{Spec}^\mathrm{CS}\breve{\Delta}_{\psi,\Gamma,A}/\mathrm{Fro}^\mathbb{Z},\breve{\nabla}_{\psi,\Gamma,A}/\mathrm{Fro}^\mathbb{Z},\mathrm{Spec}^\mathrm{CS}\breve{\Phi}_{\psi,\Gamma,A}/\mathrm{Fro}^\mathbb{Z},\mathrm{Spec}^\mathrm{CS}\breve{\Delta}^+_{\psi,\Gamma,A}/\mathrm{Fro}^\mathbb{Z},\\
&\mathrm{Spec}^\mathrm{CS}\breve{\nabla}^+_{\psi,\Gamma,A}/\mathrm{Fro}^\mathbb{Z},\mathrm{Spec}^\mathrm{CS}\breve{\Delta}^\dagger_{\psi,\Gamma,A}/\mathrm{Fro}^\mathbb{Z},\mathrm{Spec}^\mathrm{CS}\breve{\nabla}^\dagger_{\psi,\Gamma,A}/\mathrm{Fro}^\mathbb{Z},	\\
\end{align}
\begin{align}
&\mathrm{Spec}^\mathrm{CS}{\Delta}_{\psi,\Gamma,A}/\mathrm{Fro}^\mathbb{Z},\mathrm{Spec}^\mathrm{CS}{\nabla}_{\psi,\Gamma,A}/\mathrm{Fro}^\mathbb{Z},\mathrm{Spec}^\mathrm{CS}{\Phi}_{\psi,\Gamma,A}/\mathrm{Fro}^\mathbb{Z},\mathrm{Spec}^\mathrm{CS}{\Delta}^+_{\psi,\Gamma,A}/\mathrm{Fro}^\mathbb{Z},\\
&\mathrm{Spec}^\mathrm{CS}{\nabla}^+_{\psi,\Gamma,A}/\mathrm{Fro}^\mathbb{Z}, \mathrm{Spec}^\mathrm{CS}{\Delta}^\dagger_{\psi,\Gamma,A}/\mathrm{Fro}^\mathbb{Z},\mathrm{Spec}^\mathrm{CS}{\nabla}^\dagger_{\psi,\Gamma,A}/\mathrm{Fro}^\mathbb{Z}.	
\end{align}
Here for those space with notations related to the radius and the corresponding interval we consider the total unions $\bigcap_r,\bigcup_I$ in order to achieve the whole spaces to achieve the analogues of the corresponding FF curves from \cite{10KL1}, \cite{10KL2}, \cite{10FF} for
\[
\xymatrix@R+0pc@C+0pc{
\underset{r}{\mathrm{homotopylimit}}~\mathrm{Spec}^\mathrm{CS}\widetilde{\Phi}^r_{\psi,\Gamma,A},\underset{I}{\mathrm{homotopycolimit}}~\mathrm{Spec}^\mathrm{CS}\widetilde{\Phi}^I_{\psi,\Gamma,A},	\\
}
\]
\[
\xymatrix@R+0pc@C+0pc{
\underset{r}{\mathrm{homotopylimit}}~\mathrm{Spec}^\mathrm{CS}\breve{\Phi}^r_{\psi,\Gamma,A},\underset{I}{\mathrm{homotopycolimit}}~\mathrm{Spec}^\mathrm{CS}\breve{\Phi}^I_{\psi,\Gamma,A},	\\
}
\]
\[
\xymatrix@R+0pc@C+0pc{
\underset{r}{\mathrm{homotopylimit}}~\mathrm{Spec}^\mathrm{CS}{\Phi}^r_{\psi,\Gamma,A},\underset{I}{\mathrm{homotopycolimit}}~\mathrm{Spec}^\mathrm{CS}{\Phi}^I_{\psi,\Gamma,A}.	
}
\]
\[ 
\xymatrix@R+0pc@C+0pc{
\underset{r}{\mathrm{homotopylimit}}~\mathrm{Spec}^\mathrm{CS}\widetilde{\Phi}^r_{\psi,\Gamma,A}/\mathrm{Fro}^\mathbb{Z},\underset{I}{\mathrm{homotopycolimit}}~\mathrm{Spec}^\mathrm{CS}\widetilde{\Phi}^I_{\psi,\Gamma,A}/\mathrm{Fro}^\mathbb{Z},	\\
}
\]
\[ 
\xymatrix@R+0pc@C+0pc{
\underset{r}{\mathrm{homotopylimit}}~\mathrm{Spec}^\mathrm{CS}\breve{\Phi}^r_{\psi,\Gamma,A}/\mathrm{Fro}^\mathbb{Z},\underset{I}{\mathrm{homotopycolimit}}~\breve{\Phi}^I_{\psi,\Gamma,A}/\mathrm{Fro}^\mathbb{Z},	\\
}
\]
\[ 
\xymatrix@R+0pc@C+0pc{
\underset{r}{\mathrm{homotopylimit}}~\mathrm{Spec}^\mathrm{CS}{\Phi}^r_{\psi,\Gamma,A}/\mathrm{Fro}^\mathbb{Z},\underset{I}{\mathrm{homotopycolimit}}~\mathrm{Spec}^\mathrm{CS}{\Phi}^I_{\psi,\Gamma,A}/\mathrm{Fro}^\mathbb{Z}.	
}
\]

\end{definition}

\begin{proposition}
There is a well-defined functor from the $\infty$-category 
\begin{align}
\mathrm{Quasicoherentpresheaves,Perfectcomplex,Condensed}_{*}	
\end{align}
where $*$ is one of the following spaces:
\begin{align}
&\mathrm{Spec}^\mathrm{CS}\widetilde{\Phi}_{\psi,\Gamma,A}/\mathrm{Fro}^\mathbb{Z},	\\
\end{align}
\begin{align}
&\mathrm{Spec}^\mathrm{CS}\breve{\Phi}_{\psi,\Gamma,A}/\mathrm{Fro}^\mathbb{Z},	\\
\end{align}
\begin{align}
&\mathrm{Spec}^\mathrm{CS}{\Phi}_{\psi,\Gamma,A}/\mathrm{Fro}^\mathbb{Z},	
\end{align}
to the $\infty$-category of $\mathrm{Fro}$-equivariant quasicoherent presheaves over similar spaces above correspondingly without the $\mathrm{Fro}$-quotients, and to the $\infty$-category of $\mathrm{Fro}$-equivariant quasicoherent modules over global sections of the structure $\infty$-sheaves of the similar spaces above correspondingly without the $\mathrm{Fro}$-quotients. Here for those space without notation related to the radius and the corresponding interval we consider the total unions $\bigcap_r,\bigcup_I$ in order to achieve the whole spaces to achieve the analogues of the corresponding FF curves from \cite{10KL1}, \cite{10KL2}, \cite{10FF} for
\[
\xymatrix@R+0pc@C+0pc{
\underset{r}{\mathrm{homotopylimit}}~\mathrm{Spec}^\mathrm{CS}\widetilde{\Phi}^r_{\psi,\Gamma,A},\underset{I}{\mathrm{homotopycolimit}}~\mathrm{Spec}^\mathrm{CS}\widetilde{\Phi}^I_{\psi,\Gamma,A},	\\
}
\]
\[
\xymatrix@R+0pc@C+0pc{
\underset{r}{\mathrm{homotopylimit}}~\mathrm{Spec}^\mathrm{CS}\breve{\Phi}^r_{\psi,\Gamma,A},\underset{I}{\mathrm{homotopycolimit}}~\mathrm{Spec}^\mathrm{CS}\breve{\Phi}^I_{\psi,\Gamma,A},	\\
}
\]
\[
\xymatrix@R+0pc@C+0pc{
\underset{r}{\mathrm{homotopylimit}}~\mathrm{Spec}^\mathrm{CS}{\Phi}^r_{\psi,\Gamma,A},\underset{I}{\mathrm{homotopycolimit}}~\mathrm{Spec}^\mathrm{CS}{\Phi}^I_{\psi,\Gamma,A}.	
}
\]
\[ 
\xymatrix@R+0pc@C+0pc{
\underset{r}{\mathrm{homotopylimit}}~\mathrm{Spec}^\mathrm{CS}\widetilde{\Phi}^r_{\psi,\Gamma,A}/\mathrm{Fro}^\mathbb{Z},\underset{I}{\mathrm{homotopycolimit}}~\mathrm{Spec}^\mathrm{CS}\widetilde{\Phi}^I_{\psi,\Gamma,A}/\mathrm{Fro}^\mathbb{Z},	\\
}
\]
\[ 
\xymatrix@R+0pc@C+0pc{
\underset{r}{\mathrm{homotopylimit}}~\mathrm{Spec}^\mathrm{CS}\breve{\Phi}^r_{\psi,\Gamma,A}/\mathrm{Fro}^\mathbb{Z},\underset{I}{\mathrm{homotopycolimit}}~\breve{\Phi}^I_{\psi,\Gamma,A}/\mathrm{Fro}^\mathbb{Z},	\\
}
\]
\[ 
\xymatrix@R+0pc@C+0pc{
\underset{r}{\mathrm{homotopylimit}}~\mathrm{Spec}^\mathrm{CS}{\Phi}^r_{\psi,\Gamma,A}/\mathrm{Fro}^\mathbb{Z},\underset{I}{\mathrm{homotopycolimit}}~\mathrm{Spec}^\mathrm{CS}{\Phi}^I_{\psi,\Gamma,A}/\mathrm{Fro}^\mathbb{Z}.	
}
\]	
In this situation we will have the target category being family parametrized by $r$ or $I$ in compatible glueing sense as in \cite[Definition 5.4.10]{10KL2}. In this situation for modules parametrized by the intervals we have the equivalence of $\infty$-categories by using \cite[Proposition 12.18]{10CS2}. Here the corresponding quasicoherent Frobenius modules are defined to be the corresponding homotopy colimits and limits of Frobenius modules:
\begin{align}
\underset{r}{\mathrm{homotopycolimit}}~M_r,\\
\underset{I}{\mathrm{homotopylimit}}~M_I,	
\end{align}
where each $M_r$ is a Frobenius-equivariant module over the period ring with respect to some radius $r$ while each $M_I$ is a Frobenius-equivariant module over the period ring with respect to some interval $I$.\\
\end{proposition}

\begin{proposition}
Similar proposition holds for 
\begin{align}
\mathrm{Quasicoherentsheaves,Perfectcomplex,IndBanach}_{*}.	
\end{align}	
\end{proposition}

\newpage
\subsection{Frobenius Quasicoherent Modules II: Deformation in Banach Rings}

\begin{definition}
Let $\psi$ be a toric tower over $\mathbb{Q}_p$ as in \cite[Chapter 7]{10KL2} with base $\mathbb{Q}_p\left<X_1^{\pm 1},...,X_k^{\pm 1}\right>$. Then from \cite{10KL1} and \cite[Definition 5.2.1]{10KL2} we have the following class of Kedlaya-Liu rings (with the following replacement: $\Delta$ stands for $A$, $\nabla$ stands for $B$, while $\Phi$ stands for $C$) by taking product in the sense of self $\Gamma$-th power:

\[
\xymatrix@R+0pc@C+0pc{
\widetilde{\Delta}_{\psi,\Gamma},\widetilde{\nabla}_{\psi,\Gamma},\widetilde{\Phi}_{\psi,\Gamma},\widetilde{\Delta}^+_{\psi,\Gamma},\widetilde{\nabla}^+_{\psi,\Gamma},\widetilde{\Delta}^\dagger_{\psi,\Gamma},\widetilde{\nabla}^\dagger_{\psi,\Gamma},\widetilde{\Phi}^r_{\psi,\Gamma},\widetilde{\Phi}^I_{\psi,\Gamma}, 
}
\]

\[
\xymatrix@R+0pc@C+0pc{
\breve{\Delta}_{\psi,\Gamma},\breve{\nabla}_{\psi,\Gamma},\breve{\Phi}_{\psi,\Gamma},\breve{\Delta}^+_{\psi,\Gamma},\breve{\nabla}^+_{\psi,\Gamma},\breve{\Delta}^\dagger_{\psi,\Gamma},\breve{\nabla}^\dagger_{\psi,\Gamma},\breve{\Phi}^r_{\psi,\Gamma},\breve{\Phi}^I_{\psi,\Gamma},	
}
\]

\[
\xymatrix@R+0pc@C+0pc{
{\Delta}_{\psi,\Gamma},{\nabla}_{\psi,\Gamma},{\Phi}_{\psi,\Gamma},{\Delta}^+_{\psi,\Gamma},{\nabla}^+_{\psi,\Gamma},{\Delta}^\dagger_{\psi,\Gamma},{\nabla}^\dagger_{\psi,\Gamma},{\Phi}^r_{\psi,\Gamma},{\Phi}^I_{\psi,\Gamma}.	
}
\]
We now consider the following rings with $-$ being any deforming Banach ring over $\mathbb{Q}_p$. Taking the product we have:
\[
\xymatrix@R+0pc@C+0pc{
\widetilde{\Phi}_{\psi,\Gamma,-},\widetilde{\Phi}^r_{\psi,\Gamma,-},\widetilde{\Phi}^I_{\psi,\Gamma,-},	
}
\]
\[
\xymatrix@R+0pc@C+0pc{
\breve{\Phi}_{\psi,\Gamma,-},\breve{\Phi}^r_{\psi,\Gamma,-},\breve{\Phi}^I_{\psi,\Gamma,-},	
}
\]
\[
\xymatrix@R+0pc@C+0pc{
{\Phi}_{\psi,\Gamma,-},{\Phi}^r_{\psi,\Gamma,-},{\Phi}^I_{\psi,\Gamma,-}.	
}
\]
They carry multi Frobenius action $\varphi_\Gamma$ and multi $\mathrm{Lie}_\Gamma:=\mathbb{Z}_p^{\times\Gamma}$ action. In our current situation after \cite{10CKZ} and \cite{10PZ} we consider the following $(\infty,1)$-categories of $(\infty,1)$-modules.\\
\end{definition}

\begin{definition}
First we consider the Bambozzi-Kremnizer spectrum $\mathrm{Spec}^\mathrm{BK}(*)$ attached to any of those in the above from \cite{10BK} by taking derived rational localization:
\begin{align}
&\mathrm{Spec}^\mathrm{BK}\widetilde{\Phi}_{\psi,\Gamma,-},\mathrm{Spec}^\mathrm{BK}\widetilde{\Phi}^r_{\psi,\Gamma,-},\mathrm{Spec}^\mathrm{BK}\widetilde{\Phi}^I_{\psi,\Gamma,-},	
\end{align}
\begin{align}
&\mathrm{Spec}^\mathrm{BK}\breve{\Phi}_{\psi,\Gamma,-},\mathrm{Spec}^\mathrm{BK}\breve{\Phi}^r_{\psi,\Gamma,-},\mathrm{Spec}^\mathrm{BK}\breve{\Phi}^I_{\psi,\Gamma,-},	
\end{align}
\begin{align}
&\mathrm{Spec}^\mathrm{BK}{\Phi}_{\psi,\Gamma,-},
\mathrm{Spec}^\mathrm{BK}{\Phi}^r_{\psi,\Gamma,-},\mathrm{Spec}^\mathrm{BK}{\Phi}^I_{\psi,\Gamma,-}.	
\end{align}

Then we take the corresponding quotients by using the corresponding Frobenius operators:
\begin{align}
&\mathrm{Spec}^\mathrm{BK}\widetilde{\Phi}_{\psi,\Gamma,-}/\mathrm{Fro}^\mathbb{Z},	\\
\end{align}
\begin{align}
&\mathrm{Spec}^\mathrm{BK}\breve{\Phi}_{\psi,\Gamma,-}/\mathrm{Fro}^\mathbb{Z},	\\
\end{align}
\begin{align}
&\mathrm{Spec}^\mathrm{BK}{\Phi}_{\psi,\Gamma,-}/\mathrm{Fro}^\mathbb{Z}.	
\end{align}
Here for those space without notation related to the radius and the corresponding interval we consider the total unions $\bigcap_r,\bigcup_I$ in order to achieve the whole spaces to achieve the analogues of the corresponding FF curves from \cite{10KL1}, \cite{10KL2}, \cite{10FF} for
\[
\xymatrix@R+0pc@C+0pc{
\underset{r}{\mathrm{homotopylimit}}~\mathrm{Spec}^\mathrm{BK}\widetilde{\Phi}^r_{\psi,\Gamma,-},\underset{I}{\mathrm{homotopycolimit}}~\mathrm{Spec}^\mathrm{BK}\widetilde{\Phi}^I_{\psi,\Gamma,-},	\\
}
\]
\[
\xymatrix@R+0pc@C+0pc{
\underset{r}{\mathrm{homotopylimit}}~\mathrm{Spec}^\mathrm{BK}\breve{\Phi}^r_{\psi,\Gamma,-},\underset{I}{\mathrm{homotopycolimit}}~\mathrm{Spec}^\mathrm{BK}\breve{\Phi}^I_{\psi,\Gamma,-},	\\
}
\]
\[
\xymatrix@R+0pc@C+0pc{
\underset{r}{\mathrm{homotopylimit}}~\mathrm{Spec}^\mathrm{BK}{\Phi}^r_{\psi,\Gamma,-},\underset{I}{\mathrm{homotopycolimit}}~\mathrm{Spec}^\mathrm{BK}{\Phi}^I_{\psi,\Gamma,-}.	
}
\]
\[  
\xymatrix@R+0pc@C+0pc{
\underset{r}{\mathrm{homotopylimit}}~\mathrm{Spec}^\mathrm{BK}\widetilde{\Phi}^r_{\psi,\Gamma,-}/\mathrm{Fro}^\mathbb{Z},\underset{I}{\mathrm{homotopycolimit}}~\mathrm{Spec}^\mathrm{BK}\widetilde{\Phi}^I_{\psi,\Gamma,-}/\mathrm{Fro}^\mathbb{Z},	\\
}
\]
\[ 
\xymatrix@R+0pc@C+0pc{
\underset{r}{\mathrm{homotopylimit}}~\mathrm{Spec}^\mathrm{BK}\breve{\Phi}^r_{\psi,\Gamma,-}/\mathrm{Fro}^\mathbb{Z},\underset{I}{\mathrm{homotopycolimit}}~\mathrm{Spec}^\mathrm{BK}\breve{\Phi}^I_{\psi,\Gamma,-}/\mathrm{Fro}^\mathbb{Z},	\\
}
\]
\[ 
\xymatrix@R+0pc@C+0pc{
\underset{r}{\mathrm{homotopylimit}}~\mathrm{Spec}^\mathrm{BK}{\Phi}^r_{\psi,\Gamma,-}/\mathrm{Fro}^\mathbb{Z},\underset{I}{\mathrm{homotopycolimit}}~\mathrm{Spec}^\mathrm{BK}{\Phi}^I_{\psi,\Gamma,-}/\mathrm{Fro}^\mathbb{Z}.	
}
\]

\end{definition}

\indent Meanwhile we have the corresponding Clausen-Scholze analytic stacks from \cite{10CS2}, therefore applying their construction we have:

\begin{definition}
Here we define the following products by using the solidified tensor product from \cite{10CS1} and \cite{10CS2}. Namely $A$ will still as above as a Banach ring over $\mathbb{Q}_p$. Then we take solidified tensor product $\overset{\blacksquare}{\otimes}$ of any of the following
\[
\xymatrix@R+0pc@C+0pc{
\widetilde{\Delta}_{\psi,\Gamma},\widetilde{\nabla}_{\psi,\Gamma},\widetilde{\Phi}_{\psi,\Gamma},\widetilde{\Delta}^+_{\psi,\Gamma},\widetilde{\nabla}^+_{\psi,\Gamma},\widetilde{\Delta}^\dagger_{\psi,\Gamma},\widetilde{\nabla}^\dagger_{\psi,\Gamma},\widetilde{\Phi}^r_{\psi,\Gamma},\widetilde{\Phi}^I_{\psi,\Gamma}, 
}
\]

\[
\xymatrix@R+0pc@C+0pc{
\breve{\Delta}_{\psi,\Gamma},\breve{\nabla}_{\psi,\Gamma},\breve{\Phi}_{\psi,\Gamma},\breve{\Delta}^+_{\psi,\Gamma},\breve{\nabla}^+_{\psi,\Gamma},\breve{\Delta}^\dagger_{\psi,\Gamma},\breve{\nabla}^\dagger_{\psi,\Gamma},\breve{\Phi}^r_{\psi,\Gamma},\breve{\Phi}^I_{\psi,\Gamma},	
}
\]

\[
\xymatrix@R+0pc@C+0pc{
{\Delta}_{\psi,\Gamma},{\nabla}_{\psi,\Gamma},{\Phi}_{\psi,\Gamma},{\Delta}^+_{\psi,\Gamma},{\nabla}^+_{\psi,\Gamma},{\Delta}^\dagger_{\psi,\Gamma},{\nabla}^\dagger_{\psi,\Gamma},{\Phi}^r_{\psi,\Gamma},{\Phi}^I_{\psi,\Gamma},	
}
\]  	
with $A$. Then we have the notations:
\[
\xymatrix@R+0pc@C+0pc{
\widetilde{\Delta}_{\psi,\Gamma,-},\widetilde{\nabla}_{\psi,\Gamma,-},\widetilde{\Phi}_{\psi,\Gamma,-},\widetilde{\Delta}^+_{\psi,\Gamma,-},\widetilde{\nabla}^+_{\psi,\Gamma,-},\widetilde{\Delta}^\dagger_{\psi,\Gamma,-},\widetilde{\nabla}^\dagger_{\psi,\Gamma,-},\widetilde{\Phi}^r_{\psi,\Gamma,-},\widetilde{\Phi}^I_{\psi,\Gamma,-}, 
}
\]

\[
\xymatrix@R+0pc@C+0pc{
\breve{\Delta}_{\psi,\Gamma,-},\breve{\nabla}_{\psi,\Gamma,-},\breve{\Phi}_{\psi,\Gamma,-},\breve{\Delta}^+_{\psi,\Gamma,-},\breve{\nabla}^+_{\psi,\Gamma,-},\breve{\Delta}^\dagger_{\psi,\Gamma,-},\breve{\nabla}^\dagger_{\psi,\Gamma,-},\breve{\Phi}^r_{\psi,\Gamma,-},\breve{\Phi}^I_{\psi,\Gamma,-},	
}
\]

\[
\xymatrix@R+0pc@C+0pc{
{\Delta}_{\psi,\Gamma,-},{\nabla}_{\psi,\Gamma,-},{\Phi}_{\psi,\Gamma,-},{\Delta}^+_{\psi,\Gamma,-},{\nabla}^+_{\psi,\Gamma,-},{\Delta}^\dagger_{\psi,\Gamma,-},{\nabla}^\dagger_{\psi,\Gamma,-},{\Phi}^r_{\psi,\Gamma,-},{\Phi}^I_{\psi,\Gamma,-}.	
}
\]
\end{definition}

\begin{definition}
First we consider the Clausen-Scholze spectrum $\mathrm{Spec}^\mathrm{CS}(*)$ attached to any of those in the above from \cite{10CS2} by taking derived rational localization:
\begin{align}
\mathrm{Spec}^\mathrm{CS}\widetilde{\Delta}_{\psi,\Gamma,-},\mathrm{Spec}^\mathrm{CS}\widetilde{\nabla}_{\psi,\Gamma,-},\mathrm{Spec}^\mathrm{CS}\widetilde{\Phi}_{\psi,\Gamma,-},\mathrm{Spec}^\mathrm{CS}\widetilde{\Delta}^+_{\psi,\Gamma,-},\mathrm{Spec}^\mathrm{CS}\widetilde{\nabla}^+_{\psi,\Gamma,-},\\
\mathrm{Spec}^\mathrm{CS}\widetilde{\Delta}^\dagger_{\psi,\Gamma,-},\mathrm{Spec}^\mathrm{CS}\widetilde{\nabla}^\dagger_{\psi,\Gamma,-},\mathrm{Spec}^\mathrm{CS}\widetilde{\Phi}^r_{\psi,\Gamma,-},\mathrm{Spec}^\mathrm{CS}\widetilde{\Phi}^I_{\psi,\Gamma,-},	\\
\end{align}
\begin{align}
\mathrm{Spec}^\mathrm{CS}\breve{\Delta}_{\psi,\Gamma,-},\breve{\nabla}_{\psi,\Gamma,-},\mathrm{Spec}^\mathrm{CS}\breve{\Phi}_{\psi,\Gamma,-},\mathrm{Spec}^\mathrm{CS}\breve{\Delta}^+_{\psi,\Gamma,-},\mathrm{Spec}^\mathrm{CS}\breve{\nabla}^+_{\psi,\Gamma,-},\\
\mathrm{Spec}^\mathrm{CS}\breve{\Delta}^\dagger_{\psi,\Gamma,-},\mathrm{Spec}^\mathrm{CS}\breve{\nabla}^\dagger_{\psi,\Gamma,-},\mathrm{Spec}^\mathrm{CS}\breve{\Phi}^r_{\psi,\Gamma,-},\breve{\Phi}^I_{\psi,\Gamma,-},	\\
\end{align}
\begin{align}
\mathrm{Spec}^\mathrm{CS}{\Delta}_{\psi,\Gamma,-},\mathrm{Spec}^\mathrm{CS}{\nabla}_{\psi,\Gamma,-},\mathrm{Spec}^\mathrm{CS}{\Phi}_{\psi,\Gamma,-},\mathrm{Spec}^\mathrm{CS}{\Delta}^+_{\psi,\Gamma,-},\mathrm{Spec}^\mathrm{CS}{\nabla}^+_{\psi,\Gamma,-},\\
\mathrm{Spec}^\mathrm{CS}{\Delta}^\dagger_{\psi,\Gamma,-},\mathrm{Spec}^\mathrm{CS}{\nabla}^\dagger_{\psi,\Gamma,-},\mathrm{Spec}^\mathrm{CS}{\Phi}^r_{\psi,\Gamma,-},\mathrm{Spec}^\mathrm{CS}{\Phi}^I_{\psi,\Gamma,-}.	
\end{align}

Then we take the corresponding quotients by using the corresponding Frobenius operators:
\begin{align}
&\mathrm{Spec}^\mathrm{CS}\widetilde{\Delta}_{\psi,\Gamma,-}/\mathrm{Fro}^\mathbb{Z},\mathrm{Spec}^\mathrm{CS}\widetilde{\nabla}_{\psi,\Gamma,-}/\mathrm{Fro}^\mathbb{Z},\mathrm{Spec}^\mathrm{CS}\widetilde{\Phi}_{\psi,\Gamma,-}/\mathrm{Fro}^\mathbb{Z},\mathrm{Spec}^\mathrm{CS}\widetilde{\Delta}^+_{\psi,\Gamma,-}/\mathrm{Fro}^\mathbb{Z},\\
&\mathrm{Spec}^\mathrm{CS}\widetilde{\nabla}^+_{\psi,\Gamma,-}/\mathrm{Fro}^\mathbb{Z}, \mathrm{Spec}^\mathrm{CS}\widetilde{\Delta}^\dagger_{\psi,\Gamma,-}/\mathrm{Fro}^\mathbb{Z},\mathrm{Spec}^\mathrm{CS}\widetilde{\nabla}^\dagger_{\psi,\Gamma,-}/\mathrm{Fro}^\mathbb{Z},	\\
\end{align}
\begin{align}
&\mathrm{Spec}^\mathrm{CS}\breve{\Delta}_{\psi,\Gamma,-}/\mathrm{Fro}^\mathbb{Z},\breve{\nabla}_{\psi,\Gamma,-}/\mathrm{Fro}^\mathbb{Z},\mathrm{Spec}^\mathrm{CS}\breve{\Phi}_{\psi,\Gamma,-}/\mathrm{Fro}^\mathbb{Z},\mathrm{Spec}^\mathrm{CS}\breve{\Delta}^+_{\psi,\Gamma,-}/\mathrm{Fro}^\mathbb{Z},\\
&\mathrm{Spec}^\mathrm{CS}\breve{\nabla}^+_{\psi,\Gamma,-}/\mathrm{Fro}^\mathbb{Z}, \mathrm{Spec}^\mathrm{CS}\breve{\Delta}^\dagger_{\psi,\Gamma,-}/\mathrm{Fro}^\mathbb{Z},\mathrm{Spec}^\mathrm{CS}\breve{\nabla}^\dagger_{\psi,\Gamma,-}/\mathrm{Fro}^\mathbb{Z},	\\
\end{align}
\begin{align}
&\mathrm{Spec}^\mathrm{CS}{\Delta}_{\psi,\Gamma,-}/\mathrm{Fro}^\mathbb{Z},\mathrm{Spec}^\mathrm{CS}{\nabla}_{\psi,\Gamma,-}/\mathrm{Fro}^\mathbb{Z},\mathrm{Spec}^\mathrm{CS}{\Phi}_{\psi,\Gamma,-}/\mathrm{Fro}^\mathbb{Z},\mathrm{Spec}^\mathrm{CS}{\Delta}^+_{\psi,\Gamma,-}/\mathrm{Fro}^\mathbb{Z},\\
&\mathrm{Spec}^\mathrm{CS}{\nabla}^+_{\psi,\Gamma,-}/\mathrm{Fro}^\mathbb{Z}, \mathrm{Spec}^\mathrm{CS}{\Delta}^\dagger_{\psi,\Gamma,-}/\mathrm{Fro}^\mathbb{Z},\mathrm{Spec}^\mathrm{CS}{\nabla}^\dagger_{\psi,\Gamma,-}/\mathrm{Fro}^\mathbb{Z}.	
\end{align}
Here for those space with notations related to the radius and the corresponding interval we consider the total unions $\bigcap_r,\bigcup_I$ in order to achieve the whole spaces to achieve the analogues of the corresponding FF curves from \cite{10KL1}, \cite{10KL2}, \cite{10FF} for
\[
\xymatrix@R+0pc@C+0pc{
\underset{r}{\mathrm{homotopylimit}}~\mathrm{Spec}^\mathrm{CS}\widetilde{\Phi}^r_{\psi,\Gamma,-},\underset{I}{\mathrm{homotopycolimit}}~\mathrm{Spec}^\mathrm{CS}\widetilde{\Phi}^I_{\psi,\Gamma,-},	\\
}
\]
\[
\xymatrix@R+0pc@C+0pc{
\underset{r}{\mathrm{homotopylimit}}~\mathrm{Spec}^\mathrm{CS}\breve{\Phi}^r_{\psi,\Gamma,-},\underset{I}{\mathrm{homotopycolimit}}~\mathrm{Spec}^\mathrm{CS}\breve{\Phi}^I_{\psi,\Gamma,-},	\\
}
\]
\[
\xymatrix@R+0pc@C+0pc{
\underset{r}{\mathrm{homotopylimit}}~\mathrm{Spec}^\mathrm{CS}{\Phi}^r_{\psi,\Gamma,-},\underset{I}{\mathrm{homotopycolimit}}~\mathrm{Spec}^\mathrm{CS}{\Phi}^I_{\psi,\Gamma,-}.	
}
\]
\[ 
\xymatrix@R+0pc@C+0pc{
\underset{r}{\mathrm{homotopylimit}}~\mathrm{Spec}^\mathrm{CS}\widetilde{\Phi}^r_{\psi,\Gamma,-}/\mathrm{Fro}^\mathbb{Z},\underset{I}{\mathrm{homotopycolimit}}~\mathrm{Spec}^\mathrm{CS}\widetilde{\Phi}^I_{\psi,\Gamma,-}/\mathrm{Fro}^\mathbb{Z},	\\
}
\]
\[ 
\xymatrix@R+0pc@C+0pc{
\underset{r}{\mathrm{homotopylimit}}~\mathrm{Spec}^\mathrm{CS}\breve{\Phi}^r_{\psi,\Gamma,-}/\mathrm{Fro}^\mathbb{Z},\underset{I}{\mathrm{homotopycolimit}}~\breve{\Phi}^I_{\psi,\Gamma,-}/\mathrm{Fro}^\mathbb{Z},	\\
}
\]
\[ 
\xymatrix@R+0pc@C+0pc{
\underset{r}{\mathrm{homotopylimit}}~\mathrm{Spec}^\mathrm{CS}{\Phi}^r_{\psi,\Gamma,-}/\mathrm{Fro}^\mathbb{Z},\underset{I}{\mathrm{homotopycolimit}}~\mathrm{Spec}^\mathrm{CS}{\Phi}^I_{\psi,\Gamma,-}/\mathrm{Fro}^\mathbb{Z}.	
}
\]

\end{definition}

\

\begin{definition}
We then consider the corresponding quasipresheaves of the corresponding ind-Banach or monomorphic ind-Banach modules from \cite{10BBK}, \cite{10KKM}:
\begin{align}
\mathrm{Quasicoherentpresheaves,IndBanach}_{*}	
\end{align}
where $*$ is one of the following spaces:
\begin{align}
&\mathrm{Spec}^\mathrm{BK}\widetilde{\Phi}_{\psi,\Gamma,-}/\mathrm{Fro}^\mathbb{Z},	\\
\end{align}
\begin{align}
&\mathrm{Spec}^\mathrm{BK}\breve{\Phi}_{\psi,\Gamma,-}/\mathrm{Fro}^\mathbb{Z},	\\
\end{align}
\begin{align}
&\mathrm{Spec}^\mathrm{BK}{\Phi}_{\psi,\Gamma,-}/\mathrm{Fro}^\mathbb{Z}.	
\end{align}
Here for those space without notation related to the radius and the corresponding interval we consider the total unions $\bigcap_r,\bigcup_I$ in order to achieve the whole spaces to achieve the analogues of the corresponding FF curves from \cite{10KL1}, \cite{10KL2}, \cite{10FF} for
\[
\xymatrix@R+0pc@C+0pc{
\underset{r}{\mathrm{homotopylimit}}~\mathrm{Spec}^\mathrm{BK}\widetilde{\Phi}^r_{\psi,\Gamma,-},\underset{I}{\mathrm{homotopycolimit}}~\mathrm{Spec}^\mathrm{BK}\widetilde{\Phi}^I_{\psi,\Gamma,-},	\\
}
\]
\[
\xymatrix@R+0pc@C+0pc{
\underset{r}{\mathrm{homotopylimit}}~\mathrm{Spec}^\mathrm{BK}\breve{\Phi}^r_{\psi,\Gamma,-},\underset{I}{\mathrm{homotopycolimit}}~\mathrm{Spec}^\mathrm{BK}\breve{\Phi}^I_{\psi,\Gamma,-},	\\
}
\]
\[
\xymatrix@R+0pc@C+0pc{
\underset{r}{\mathrm{homotopylimit}}~\mathrm{Spec}^\mathrm{BK}{\Phi}^r_{\psi,\Gamma,-},\underset{I}{\mathrm{homotopycolimit}}~\mathrm{Spec}^\mathrm{BK}{\Phi}^I_{\psi,\Gamma,-}.	
}
\]
\[  
\xymatrix@R+0pc@C+0pc{
\underset{r}{\mathrm{homotopylimit}}~\mathrm{Spec}^\mathrm{BK}\widetilde{\Phi}^r_{\psi,\Gamma,-}/\mathrm{Fro}^\mathbb{Z},\underset{I}{\mathrm{homotopycolimit}}~\mathrm{Spec}^\mathrm{BK}\widetilde{\Phi}^I_{\psi,\Gamma,-}/\mathrm{Fro}^\mathbb{Z},	\\
}
\]
\[ 
\xymatrix@R+0pc@C+0pc{
\underset{r}{\mathrm{homotopylimit}}~\mathrm{Spec}^\mathrm{BK}\breve{\Phi}^r_{\psi,\Gamma,-}/\mathrm{Fro}^\mathbb{Z},\underset{I}{\mathrm{homotopycolimit}}~\mathrm{Spec}^\mathrm{BK}\breve{\Phi}^I_{\psi,\Gamma,-}/\mathrm{Fro}^\mathbb{Z},	\\
}
\]
\[ 
\xymatrix@R+0pc@C+0pc{
\underset{r}{\mathrm{homotopylimit}}~\mathrm{Spec}^\mathrm{BK}{\Phi}^r_{\psi,\Gamma,-}/\mathrm{Fro}^\mathbb{Z},\underset{I}{\mathrm{homotopycolimit}}~\mathrm{Spec}^\mathrm{BK}{\Phi}^I_{\psi,\Gamma,-}/\mathrm{Fro}^\mathbb{Z}.	
}
\]

\end{definition}

\begin{definition}
We then consider the corresponding quasisheaves of the corresponding condensed solid topological modules from \cite{10CS2}:
\begin{align}
\mathrm{Quasicoherentsheaves, Condensed}_{*}	
\end{align}
where $*$ is one of the following spaces:
\begin{align}
&\mathrm{Spec}^\mathrm{CS}\widetilde{\Delta}_{\psi,\Gamma,-}/\mathrm{Fro}^\mathbb{Z},\mathrm{Spec}^\mathrm{CS}\widetilde{\nabla}_{\psi,\Gamma,-}/\mathrm{Fro}^\mathbb{Z},\mathrm{Spec}^\mathrm{CS}\widetilde{\Phi}_{\psi,\Gamma,-}/\mathrm{Fro}^\mathbb{Z},\mathrm{Spec}^\mathrm{CS}\widetilde{\Delta}^+_{\psi,\Gamma,-}/\mathrm{Fro}^\mathbb{Z},\\
&\mathrm{Spec}^\mathrm{CS}\widetilde{\nabla}^+_{\psi,\Gamma,-}/\mathrm{Fro}^\mathbb{Z},\mathrm{Spec}^\mathrm{CS}\widetilde{\Delta}^\dagger_{\psi,\Gamma,-}/\mathrm{Fro}^\mathbb{Z},\mathrm{Spec}^\mathrm{CS}\widetilde{\nabla}^\dagger_{\psi,\Gamma,-}/\mathrm{Fro}^\mathbb{Z},	\\
\end{align}
\begin{align}
&\mathrm{Spec}^\mathrm{CS}\breve{\Delta}_{\psi,\Gamma,-}/\mathrm{Fro}^\mathbb{Z},\breve{\nabla}_{\psi,\Gamma,-}/\mathrm{Fro}^\mathbb{Z},\mathrm{Spec}^\mathrm{CS}\breve{\Phi}_{\psi,\Gamma,-}/\mathrm{Fro}^\mathbb{Z},\mathrm{Spec}^\mathrm{CS}\breve{\Delta}^+_{\psi,\Gamma,-}/\mathrm{Fro}^\mathbb{Z},\\
&\mathrm{Spec}^\mathrm{CS}\breve{\nabla}^+_{\psi,\Gamma,-}/\mathrm{Fro}^\mathbb{Z},\mathrm{Spec}^\mathrm{CS}\breve{\Delta}^\dagger_{\psi,\Gamma,-}/\mathrm{Fro}^\mathbb{Z},\mathrm{Spec}^\mathrm{CS}\breve{\nabla}^\dagger_{\psi,\Gamma,-}/\mathrm{Fro}^\mathbb{Z},	\\
\end{align}
\begin{align}
&\mathrm{Spec}^\mathrm{CS}{\Delta}_{\psi,\Gamma,-}/\mathrm{Fro}^\mathbb{Z},\mathrm{Spec}^\mathrm{CS}{\nabla}_{\psi,\Gamma,-}/\mathrm{Fro}^\mathbb{Z},\mathrm{Spec}^\mathrm{CS}{\Phi}_{\psi,\Gamma,-}/\mathrm{Fro}^\mathbb{Z},\mathrm{Spec}^\mathrm{CS}{\Delta}^+_{\psi,\Gamma,-}/\mathrm{Fro}^\mathbb{Z},\\
&\mathrm{Spec}^\mathrm{CS}{\nabla}^+_{\psi,\Gamma,-}/\mathrm{Fro}^\mathbb{Z}, \mathrm{Spec}^\mathrm{CS}{\Delta}^\dagger_{\psi,\Gamma,-}/\mathrm{Fro}^\mathbb{Z},\mathrm{Spec}^\mathrm{CS}{\nabla}^\dagger_{\psi,\Gamma,-}/\mathrm{Fro}^\mathbb{Z}.	
\end{align}
Here for those space with notations related to the radius and the corresponding interval we consider the total unions $\bigcap_r,\bigcup_I$ in order to achieve the whole spaces to achieve the analogues of the corresponding FF curves from \cite{10KL1}, \cite{10KL2}, \cite{10FF} for
\[
\xymatrix@R+0pc@C+0pc{
\underset{r}{\mathrm{homotopylimit}}~\mathrm{Spec}^\mathrm{CS}\widetilde{\Phi}^r_{\psi,\Gamma,-},\underset{I}{\mathrm{homotopycolimit}}~\mathrm{Spec}^\mathrm{CS}\widetilde{\Phi}^I_{\psi,\Gamma,-},	\\
}
\]
\[
\xymatrix@R+0pc@C+0pc{
\underset{r}{\mathrm{homotopylimit}}~\mathrm{Spec}^\mathrm{CS}\breve{\Phi}^r_{\psi,\Gamma,-},\underset{I}{\mathrm{homotopycolimit}}~\mathrm{Spec}^\mathrm{CS}\breve{\Phi}^I_{\psi,\Gamma,-},	\\
}
\]
\[
\xymatrix@R+0pc@C+0pc{
\underset{r}{\mathrm{homotopylimit}}~\mathrm{Spec}^\mathrm{CS}{\Phi}^r_{\psi,\Gamma,-},\underset{I}{\mathrm{homotopycolimit}}~\mathrm{Spec}^\mathrm{CS}{\Phi}^I_{\psi,\Gamma,-}.	
}
\]
\[ 
\xymatrix@R+0pc@C+0pc{
\underset{r}{\mathrm{homotopylimit}}~\mathrm{Spec}^\mathrm{CS}\widetilde{\Phi}^r_{\psi,\Gamma,-}/\mathrm{Fro}^\mathbb{Z},\underset{I}{\mathrm{homotopycolimit}}~\mathrm{Spec}^\mathrm{CS}\widetilde{\Phi}^I_{\psi,\Gamma,-}/\mathrm{Fro}^\mathbb{Z},	\\
}
\]
\[ 
\xymatrix@R+0pc@C+0pc{
\underset{r}{\mathrm{homotopylimit}}~\mathrm{Spec}^\mathrm{CS}\breve{\Phi}^r_{\psi,\Gamma,-}/\mathrm{Fro}^\mathbb{Z},\underset{I}{\mathrm{homotopycolimit}}~\breve{\Phi}^I_{\psi,\Gamma,-}/\mathrm{Fro}^\mathbb{Z},	\\
}
\]
\[ 
\xymatrix@R+0pc@C+0pc{
\underset{r}{\mathrm{homotopylimit}}~\mathrm{Spec}^\mathrm{CS}{\Phi}^r_{\psi,\Gamma,-}/\mathrm{Fro}^\mathbb{Z},\underset{I}{\mathrm{homotopycolimit}}~\mathrm{Spec}^\mathrm{CS}{\Phi}^I_{\psi,\Gamma,-}/\mathrm{Fro}^\mathbb{Z}.	
}
\]

\end{definition}

\

\begin{proposition}
There is a well-defined functor from the $\infty$-category 
\begin{align}
\mathrm{Quasicoherentpresheaves,Condensed}_{*}	
\end{align}
where $*$ is one of the following spaces:
\begin{align}
&\mathrm{Spec}^\mathrm{CS}\widetilde{\Phi}_{\psi,\Gamma,-}/\mathrm{Fro}^\mathbb{Z},	\\
\end{align}
\begin{align}
&\mathrm{Spec}^\mathrm{CS}\breve{\Phi}_{\psi,\Gamma,-}/\mathrm{Fro}^\mathbb{Z},	\\
\end{align}
\begin{align}
&\mathrm{Spec}^\mathrm{CS}{\Phi}_{\psi,\Gamma,-}/\mathrm{Fro}^\mathbb{Z},	
\end{align}
to the $\infty$-category of $\mathrm{Fro}$-equivariant quasicoherent presheaves over similar spaces above correspondingly without the $\mathrm{Fro}$-quotients, and to the $\infty$-category of $\mathrm{Fro}$-equivariant quasicoherent modules over global sections of the structure $\infty$-sheaves of the similar spaces above correspondingly without the $\mathrm{Fro}$-quotients. Here for those space without notation related to the radius and the corresponding interval we consider the total unions $\bigcap_r,\bigcup_I$ in order to achieve the whole spaces to achieve the analogues of the corresponding FF curves from \cite{10KL1}, \cite{10KL2}, \cite{10FF} for
\[
\xymatrix@R+0pc@C+0pc{
\underset{r}{\mathrm{homotopylimit}}~\mathrm{Spec}^\mathrm{CS}\widetilde{\Phi}^r_{\psi,\Gamma,-},\underset{I}{\mathrm{homotopycolimit}}~\mathrm{Spec}^\mathrm{CS}\widetilde{\Phi}^I_{\psi,\Gamma,-},	\\
}
\]
\[
\xymatrix@R+0pc@C+0pc{
\underset{r}{\mathrm{homotopylimit}}~\mathrm{Spec}^\mathrm{CS}\breve{\Phi}^r_{\psi,\Gamma,-},\underset{I}{\mathrm{homotopycolimit}}~\mathrm{Spec}^\mathrm{CS}\breve{\Phi}^I_{\psi,\Gamma,-},	\\
}
\]
\[
\xymatrix@R+0pc@C+0pc{
\underset{r}{\mathrm{homotopylimit}}~\mathrm{Spec}^\mathrm{CS}{\Phi}^r_{\psi,\Gamma,-},\underset{I}{\mathrm{homotopycolimit}}~\mathrm{Spec}^\mathrm{CS}{\Phi}^I_{\psi,\Gamma,-}.	
}
\]
\[ 
\xymatrix@R+0pc@C+0pc{
\underset{r}{\mathrm{homotopylimit}}~\mathrm{Spec}^\mathrm{CS}\widetilde{\Phi}^r_{\psi,\Gamma,-}/\mathrm{Fro}^\mathbb{Z},\underset{I}{\mathrm{homotopycolimit}}~\mathrm{Spec}^\mathrm{CS}\widetilde{\Phi}^I_{\psi,\Gamma,-}/\mathrm{Fro}^\mathbb{Z},	\\
}
\]
\[ 
\xymatrix@R+0pc@C+0pc{
\underset{r}{\mathrm{homotopylimit}}~\mathrm{Spec}^\mathrm{CS}\breve{\Phi}^r_{\psi,\Gamma,-}/\mathrm{Fro}^\mathbb{Z},\underset{I}{\mathrm{homotopycolimit}}~\breve{\Phi}^I_{\psi,\Gamma,-}/\mathrm{Fro}^\mathbb{Z},	\\
}
\]
\[ 
\xymatrix@R+0pc@C+0pc{
\underset{r}{\mathrm{homotopylimit}}~\mathrm{Spec}^\mathrm{CS}{\Phi}^r_{\psi,\Gamma,-}/\mathrm{Fro}^\mathbb{Z},\underset{I}{\mathrm{homotopycolimit}}~\mathrm{Spec}^\mathrm{CS}{\Phi}^I_{\psi,\Gamma,-}/\mathrm{Fro}^\mathbb{Z}.	
}
\]	
In this situation we will have the target category being family parametrized by $r$ or $I$ in compatible glueing sense as in \cite[Definition 5.4.10]{10KL2}. In this situation for modules parametrized by the intervals we have the equivalence of $\infty$-categories by using \cite[Proposition 13.8]{10CS2}. Here the corresponding quasicoherent Frobenius modules are defined to be the corresponding homotopy colimits and limits of Frobenius modules:
\begin{align}
\underset{r}{\mathrm{homotopycolimit}}~M_r,\\
\underset{I}{\mathrm{homotopylimit}}~M_I,	
\end{align}
where each $M_r$ is a Frobenius-equivariant module over the period ring with respect to some radius $r$ while each $M_I$ is a Frobenius-equivariant module over the period ring with respect to some interval $I$.\\
\end{proposition}

\begin{proposition}
Similar proposition holds for 
\begin{align}
\mathrm{Quasicoherentsheaves,IndBanach}_{*}.	
\end{align}	
\end{proposition}

\

\begin{definition}
We then consider the corresponding quasipresheaves of perfect complexes the corresponding ind-Banach or monomorphic ind-Banach modules from \cite{10BBK}, \cite{10KKM}:
\begin{align}
\mathrm{Quasicoherentpresheaves,Perfectcomplex,IndBanach}_{*}	
\end{align}
where $*$ is one of the following spaces:
\begin{align}
&\mathrm{Spec}^\mathrm{BK}\widetilde{\Phi}_{\psi,\Gamma,-}/\mathrm{Fro}^\mathbb{Z},	\\
\end{align}
\begin{align}
&\mathrm{Spec}^\mathrm{BK}\breve{\Phi}_{\psi,\Gamma,-}/\mathrm{Fro}^\mathbb{Z},	\\
\end{align}
\begin{align}
&\mathrm{Spec}^\mathrm{BK}{\Phi}_{\psi,\Gamma,-}/\mathrm{Fro}^\mathbb{Z}.	
\end{align}
Here for those space without notation related to the radius and the corresponding interval we consider the total unions $\bigcap_r,\bigcup_I$ in order to achieve the whole spaces to achieve the analogues of the corresponding FF curves from \cite{10KL1}, \cite{10KL2}, \cite{10FF} for
\[
\xymatrix@R+0pc@C+0pc{
\underset{r}{\mathrm{homotopylimit}}~\mathrm{Spec}^\mathrm{BK}\widetilde{\Phi}^r_{\psi,\Gamma,-},\underset{I}{\mathrm{homotopycolimit}}~\mathrm{Spec}^\mathrm{BK}\widetilde{\Phi}^I_{\psi,\Gamma,-},	\\
}
\]
\[
\xymatrix@R+0pc@C+0pc{
\underset{r}{\mathrm{homotopylimit}}~\mathrm{Spec}^\mathrm{BK}\breve{\Phi}^r_{\psi,\Gamma,-},\underset{I}{\mathrm{homotopycolimit}}~\mathrm{Spec}^\mathrm{BK}\breve{\Phi}^I_{\psi,\Gamma,-},	\\
}
\]
\[
\xymatrix@R+0pc@C+0pc{
\underset{r}{\mathrm{homotopylimit}}~\mathrm{Spec}^\mathrm{BK}{\Phi}^r_{\psi,\Gamma,-},\underset{I}{\mathrm{homotopycolimit}}~\mathrm{Spec}^\mathrm{BK}{\Phi}^I_{\psi,\Gamma,-}.	
}
\]
\[  
\xymatrix@R+0pc@C+0pc{
\underset{r}{\mathrm{homotopylimit}}~\mathrm{Spec}^\mathrm{BK}\widetilde{\Phi}^r_{\psi,\Gamma,-}/\mathrm{Fro}^\mathbb{Z},\underset{I}{\mathrm{homotopycolimit}}~\mathrm{Spec}^\mathrm{BK}\widetilde{\Phi}^I_{\psi,\Gamma,-}/\mathrm{Fro}^\mathbb{Z},	\\
}
\]
\[ 
\xymatrix@R+0pc@C+0pc{
\underset{r}{\mathrm{homotopylimit}}~\mathrm{Spec}^\mathrm{BK}\breve{\Phi}^r_{\psi,\Gamma,-}/\mathrm{Fro}^\mathbb{Z},\underset{I}{\mathrm{homotopycolimit}}~\mathrm{Spec}^\mathrm{BK}\breve{\Phi}^I_{\psi,\Gamma,-}/\mathrm{Fro}^\mathbb{Z},	\\
}
\]
\[ 
\xymatrix@R+0pc@C+0pc{
\underset{r}{\mathrm{homotopylimit}}~\mathrm{Spec}^\mathrm{BK}{\Phi}^r_{\psi,\Gamma,-}/\mathrm{Fro}^\mathbb{Z},\underset{I}{\mathrm{homotopycolimit}}~\mathrm{Spec}^\mathrm{BK}{\Phi}^I_{\psi,\Gamma,-}/\mathrm{Fro}^\mathbb{Z}.	
}
\]

\end{definition}

\begin{definition}
We then consider the corresponding quasisheaves of perfect complexes of the corresponding condensed solid topological modules from \cite{10CS2}:
\begin{align}
\mathrm{Quasicoherentsheaves, Perfectcomplex, Condensed}_{*}	
\end{align}
where $*$ is one of the following spaces:
\begin{align}
&\mathrm{Spec}^\mathrm{CS}\widetilde{\Delta}_{\psi,\Gamma,-}/\mathrm{Fro}^\mathbb{Z},\mathrm{Spec}^\mathrm{CS}\widetilde{\nabla}_{\psi,\Gamma,-}/\mathrm{Fro}^\mathbb{Z},\mathrm{Spec}^\mathrm{CS}\widetilde{\Phi}_{\psi,\Gamma,-}/\mathrm{Fro}^\mathbb{Z},\mathrm{Spec}^\mathrm{CS}\widetilde{\Delta}^+_{\psi,\Gamma,-}/\mathrm{Fro}^\mathbb{Z},\\
&\mathrm{Spec}^\mathrm{CS}\widetilde{\nabla}^+_{\psi,\Gamma,-}/\mathrm{Fro}^\mathbb{Z},\mathrm{Spec}^\mathrm{CS}\widetilde{\Delta}^\dagger_{\psi,\Gamma,-}/\mathrm{Fro}^\mathbb{Z},\mathrm{Spec}^\mathrm{CS}\widetilde{\nabla}^\dagger_{\psi,\Gamma,-}/\mathrm{Fro}^\mathbb{Z},	\\
\end{align}
\begin{align}
&\mathrm{Spec}^\mathrm{CS}\breve{\Delta}_{\psi,\Gamma,-}/\mathrm{Fro}^\mathbb{Z},\breve{\nabla}_{\psi,\Gamma,-}/\mathrm{Fro}^\mathbb{Z},\mathrm{Spec}^\mathrm{CS}\breve{\Phi}_{\psi,\Gamma,-}/\mathrm{Fro}^\mathbb{Z},\mathrm{Spec}^\mathrm{CS}\breve{\Delta}^+_{\psi,\Gamma,-}/\mathrm{Fro}^\mathbb{Z},\\
&\mathrm{Spec}^\mathrm{CS}\breve{\nabla}^+_{\psi,\Gamma,-}/\mathrm{Fro}^\mathbb{Z},\mathrm{Spec}^\mathrm{CS}\breve{\Delta}^\dagger_{\psi,\Gamma,-}/\mathrm{Fro}^\mathbb{Z},\mathrm{Spec}^\mathrm{CS}\breve{\nabla}^\dagger_{\psi,\Gamma,-}/\mathrm{Fro}^\mathbb{Z},	\\
\end{align}
\begin{align}
&\mathrm{Spec}^\mathrm{CS}{\Delta}_{\psi,\Gamma,-}/\mathrm{Fro}^\mathbb{Z},\mathrm{Spec}^\mathrm{CS}{\nabla}_{\psi,\Gamma,-}/\mathrm{Fro}^\mathbb{Z},\mathrm{Spec}^\mathrm{CS}{\Phi}_{\psi,\Gamma,-}/\mathrm{Fro}^\mathbb{Z},\mathrm{Spec}^\mathrm{CS}{\Delta}^+_{\psi,\Gamma,-}/\mathrm{Fro}^\mathbb{Z},\\
&\mathrm{Spec}^\mathrm{CS}{\nabla}^+_{\psi,\Gamma,-}/\mathrm{Fro}^\mathbb{Z}, \mathrm{Spec}^\mathrm{CS}{\Delta}^\dagger_{\psi,\Gamma,-}/\mathrm{Fro}^\mathbb{Z},\mathrm{Spec}^\mathrm{CS}{\nabla}^\dagger_{\psi,\Gamma,-}/\mathrm{Fro}^\mathbb{Z}.	
\end{align}
Here for those space with notations related to the radius and the corresponding interval we consider the total unions $\bigcap_r,\bigcup_I$ in order to achieve the whole spaces to achieve the analogues of the corresponding FF curves from \cite{10KL1}, \cite{10KL2}, \cite{10FF} for
\[
\xymatrix@R+0pc@C+0pc{
\underset{r}{\mathrm{homotopylimit}}~\mathrm{Spec}^\mathrm{CS}\widetilde{\Phi}^r_{\psi,\Gamma,-},\underset{I}{\mathrm{homotopycolimit}}~\mathrm{Spec}^\mathrm{CS}\widetilde{\Phi}^I_{\psi,\Gamma,-},	\\
}
\]
\[
\xymatrix@R+0pc@C+0pc{
\underset{r}{\mathrm{homotopylimit}}~\mathrm{Spec}^\mathrm{CS}\breve{\Phi}^r_{\psi,\Gamma,-},\underset{I}{\mathrm{homotopycolimit}}~\mathrm{Spec}^\mathrm{CS}\breve{\Phi}^I_{\psi,\Gamma,-},	\\
}
\]
\[
\xymatrix@R+0pc@C+0pc{
\underset{r}{\mathrm{homotopylimit}}~\mathrm{Spec}^\mathrm{CS}{\Phi}^r_{\psi,\Gamma,-},\underset{I}{\mathrm{homotopycolimit}}~\mathrm{Spec}^\mathrm{CS}{\Phi}^I_{\psi,\Gamma,-}.	
}
\]
\[ 
\xymatrix@R+0pc@C+0pc{
\underset{r}{\mathrm{homotopylimit}}~\mathrm{Spec}^\mathrm{CS}\widetilde{\Phi}^r_{\psi,\Gamma,-}/\mathrm{Fro}^\mathbb{Z},\underset{I}{\mathrm{homotopycolimit}}~\mathrm{Spec}^\mathrm{CS}\widetilde{\Phi}^I_{\psi,\Gamma,-}/\mathrm{Fro}^\mathbb{Z},	\\
}
\]
\[ 
\xymatrix@R+0pc@C+0pc{
\underset{r}{\mathrm{homotopylimit}}~\mathrm{Spec}^\mathrm{CS}\breve{\Phi}^r_{\psi,\Gamma,-}/\mathrm{Fro}^\mathbb{Z},\underset{I}{\mathrm{homotopycolimit}}~\breve{\Phi}^I_{\psi,\Gamma,-}/\mathrm{Fro}^\mathbb{Z},	\\
}
\]
\[ 
\xymatrix@R+0pc@C+0pc{
\underset{r}{\mathrm{homotopylimit}}~\mathrm{Spec}^\mathrm{CS}{\Phi}^r_{\psi,\Gamma,-}/\mathrm{Fro}^\mathbb{Z},\underset{I}{\mathrm{homotopycolimit}}~\mathrm{Spec}^\mathrm{CS}{\Phi}^I_{\psi,\Gamma,-}/\mathrm{Fro}^\mathbb{Z}.	
}
\]

\end{definition}

\begin{proposition}
There is a well-defined functor from the $\infty$-category 
\begin{align}
\mathrm{Quasicoherentpresheaves,Perfectcomplex,Condensed}_{*}	
\end{align}
where $*$ is one of the following spaces:
\begin{align}
&\mathrm{Spec}^\mathrm{CS}\widetilde{\Phi}_{\psi,\Gamma,-}/\mathrm{Fro}^\mathbb{Z},	\\
\end{align}
\begin{align}
&\mathrm{Spec}^\mathrm{CS}\breve{\Phi}_{\psi,\Gamma,-}/\mathrm{Fro}^\mathbb{Z},	\\
\end{align}
\begin{align}
&\mathrm{Spec}^\mathrm{CS}{\Phi}_{\psi,\Gamma,-}/\mathrm{Fro}^\mathbb{Z},	
\end{align}
to the $\infty$-category of $\mathrm{Fro}$-equivariant quasicoherent presheaves over similar spaces above correspondingly without the $\mathrm{Fro}$-quotients, and to the $\infty$-category of $\mathrm{Fro}$-equivariant quasicoherent modules over global sections of the structure $\infty$-sheaves of the similar spaces above correspondingly without the $\mathrm{Fro}$-quotients. Here for those space without notation related to the radius and the corresponding interval we consider the total unions $\bigcap_r,\bigcup_I$ in order to achieve the whole spaces to achieve the analogues of the corresponding FF curves from \cite{10KL1}, \cite{10KL2}, \cite{10FF} for
\[
\xymatrix@R+0pc@C+0pc{
\underset{r}{\mathrm{homotopylimit}}~\mathrm{Spec}^\mathrm{CS}\widetilde{\Phi}^r_{\psi,\Gamma,-},\underset{I}{\mathrm{homotopycolimit}}~\mathrm{Spec}^\mathrm{CS}\widetilde{\Phi}^I_{\psi,\Gamma,-},	\\
}
\]
\[
\xymatrix@R+0pc@C+0pc{
\underset{r}{\mathrm{homotopylimit}}~\mathrm{Spec}^\mathrm{CS}\breve{\Phi}^r_{\psi,\Gamma,-},\underset{I}{\mathrm{homotopycolimit}}~\mathrm{Spec}^\mathrm{CS}\breve{\Phi}^I_{\psi,\Gamma,-},	\\
}
\]
\[
\xymatrix@R+0pc@C+0pc{
\underset{r}{\mathrm{homotopylimit}}~\mathrm{Spec}^\mathrm{CS}{\Phi}^r_{\psi,\Gamma,-},\underset{I}{\mathrm{homotopycolimit}}~\mathrm{Spec}^\mathrm{CS}{\Phi}^I_{\psi,\Gamma,-}.	
}
\]
\[ 
\xymatrix@R+0pc@C+0pc{
\underset{r}{\mathrm{homotopylimit}}~\mathrm{Spec}^\mathrm{CS}\widetilde{\Phi}^r_{\psi,\Gamma,-}/\mathrm{Fro}^\mathbb{Z},\underset{I}{\mathrm{homotopycolimit}}~\mathrm{Spec}^\mathrm{CS}\widetilde{\Phi}^I_{\psi,\Gamma,-}/\mathrm{Fro}^\mathbb{Z},	\\
}
\]
\[ 
\xymatrix@R+0pc@C+0pc{
\underset{r}{\mathrm{homotopylimit}}~\mathrm{Spec}^\mathrm{CS}\breve{\Phi}^r_{\psi,\Gamma,-}/\mathrm{Fro}^\mathbb{Z},\underset{I}{\mathrm{homotopycolimit}}~\breve{\Phi}^I_{\psi,\Gamma,-}/\mathrm{Fro}^\mathbb{Z},	\\
}
\]
\[ 
\xymatrix@R+0pc@C+0pc{
\underset{r}{\mathrm{homotopylimit}}~\mathrm{Spec}^\mathrm{CS}{\Phi}^r_{\psi,\Gamma,-}/\mathrm{Fro}^\mathbb{Z},\underset{I}{\mathrm{homotopycolimit}}~\mathrm{Spec}^\mathrm{CS}{\Phi}^I_{\psi,\Gamma,-}/\mathrm{Fro}^\mathbb{Z}.	
}
\]	
In this situation we will have the target category being family parametrized by $r$ or $I$ in compatible glueing sense as in \cite[Definition 5.4.10]{10KL2}. In this situation for modules parametrized by the intervals we have the equivalence of $\infty$-categories by using \cite[Proposition 12.18]{10CS2}. Here the corresponding quasicoherent Frobenius modules are defined to be the corresponding homotopy colimits and limits of Frobenius modules:
\begin{align}
\underset{r}{\mathrm{homotopycolimit}}~M_r,\\
\underset{I}{\mathrm{homotopylimit}}~M_I,	
\end{align}
where each $M_r$ is a Frobenius-equivariant module over the period ring with respect to some radius $r$ while each $M_I$ is a Frobenius-equivariant module over the period ring with respect to some interval $I$.\\
\end{proposition}

\begin{proposition}
Similar proposition holds for 
\begin{align}
\mathrm{Quasicoherentsheaves,Perfectcomplex,IndBanach}_{*}.	
\end{align}	
\end{proposition}

\indent Then we have the following functoriality results:

\begin{proposition}
We have the following commutative diagram:
\[
\xymatrix@R+5pc@C+0pc{
\mathrm{Quasicoherentsheaves,Condensed}_{*}\ar[d]\ar[d]\ar[d]\ar[r]\ar[r]\ar[r] & \mathrm{Quasicoherentsheaves,Condensed}_{.}\ar[d]\ar[d]\ar[d]	\\
\mathrm{Quasicoherentsheaves,Condensed}_{*_{\psi_0}}\ar[r]\ar[r]\ar[r] & \mathrm{Quasicoherentsheaves,Condensed}_{._{\psi_0}}.	\\
}
\]	
\end{proposition}

\begin{proposition}
We have the following commutative diagram:
\[
\xymatrix@R+5pc@C+0pc{
\mathrm{Quasicoherentsheaves,IndBanach}_{*}\ar[d]\ar[d]\ar[d]\ar[r]\ar[r]\ar[r] & \mathrm{Quasicoherentsheaves,IndBanach}_{.}\ar[d]\ar[d]\ar[d]	\\
\mathrm{Quasicoherentsheaves,IndBanach}_{*_{\psi_0}}\ar[r]\ar[r]\ar[r] & \mathrm{Quasicoherentsheaves,IndBanach}_{._{\psi_0}}.	\\
}
\]	
\end{proposition}

\begin{proposition}
We have the following commutative diagram:
\[\tiny
\xymatrix@R+5pc@C+0pc{
\mathrm{Quasicoherentsheaves,Perfectcomplex,Condensed}_{*}\ar[d]\ar[d]\ar[d]\ar[r]\ar[r]\ar[r] & \mathrm{Quasicoherentsheaves,Perfectcomplex,Condensed}_{.}\ar[d]\ar[d]\ar[d]	\\
\mathrm{Quasicoherentsheaves,Perfectcomplex,Condensed}_{*_{\psi_0}}\ar[r]\ar[r]\ar[r] & \mathrm{Quasicoherentsheaves,Perfectcomplex,Condensed}_{._{\psi_0}}.	\\
}
\]	
\end{proposition}

\begin{proposition}
We have the following commutative diagram:
\[\tiny
\xymatrix@R+5pc@C+0pc{
\mathrm{Quasicoherentsheaves,Perfectcomplex,IndBanach}_{*}\ar[d]\ar[d]\ar[d]\ar[r]\ar[r]\ar[r] & \mathrm{Quasicoherentsheaves,Perfectcomplex,IndBanach}_{.}\ar[d]\ar[d]\ar[d]	\\
\mathrm{Quasicoherentsheaves,Perfectcomplex,IndBanach}_{*_{\psi_0}}\ar[r]\ar[r]\ar[r] & \mathrm{Quasicoherentsheaves,Perfectcomplex,IndBanach}_{._{\psi_0}}.	\\
}
\]	
\end{proposition}

\newpage
\subsection{Frobenius Quasicoherent Modules III: Deformation in $(\infty,1)$-Ind-Banach Rings}

\begin{definition}
Let $\psi$ be a toric tower over $\mathbb{Q}_p$ as in \cite[Chapter 7]{10KL2} with base $\mathbb{Q}_p\left<X_1^{\pm 1},...,X_k^{\pm 1}\right>$. Then from \cite{10KL1} and \cite[Definition 5.2.1]{10KL2} we have the following class of Kedlaya-Liu rings (with the following replacement: $\Delta$ stands for $A$, $\nabla$ stands for $B$, while $\Phi$ stands for $C$) by taking product in the sense of self $\Gamma$-th power:

\[
\xymatrix@R+0pc@C+0pc{
\widetilde{\Delta}_{\psi,\Gamma},\widetilde{\nabla}_{\psi,\Gamma},\widetilde{\Phi}_{\psi,\Gamma},\widetilde{\Delta}^+_{\psi,\Gamma},\widetilde{\nabla}^+_{\psi,\Gamma},\widetilde{\Delta}^\dagger_{\psi,\Gamma},\widetilde{\nabla}^\dagger_{\psi,\Gamma},\widetilde{\Phi}^r_{\psi,\Gamma},\widetilde{\Phi}^I_{\psi,\Gamma}, 
}
\]

\[
\xymatrix@R+0pc@C+0pc{
\breve{\Delta}_{\psi,\Gamma},\breve{\nabla}_{\psi,\Gamma},\breve{\Phi}_{\psi,\Gamma},\breve{\Delta}^+_{\psi,\Gamma},\breve{\nabla}^+_{\psi,\Gamma},\breve{\Delta}^\dagger_{\psi,\Gamma},\breve{\nabla}^\dagger_{\psi,\Gamma},\breve{\Phi}^r_{\psi,\Gamma},\breve{\Phi}^I_{\psi,\Gamma},	
}
\]

\[
\xymatrix@R+0pc@C+0pc{
{\Delta}_{\psi,\Gamma},{\nabla}_{\psi,\Gamma},{\Phi}_{\psi,\Gamma},{\Delta}^+_{\psi,\Gamma},{\nabla}^+_{\psi,\Gamma},{\Delta}^\dagger_{\psi,\Gamma},{\nabla}^\dagger_{\psi,\Gamma},{\Phi}^r_{\psi,\Gamma},{\Phi}^I_{\psi,\Gamma}.	
}
\]
We now consider the following rings with $\square$ being a homotopy colimit
\begin{align}
 \underset{i}{\mathrm{homotopycolimit}}\square_i
 \end{align}
 of $\mathbb{Q}_p\left<Y_1,...,Y_i\right>,i=1,2,...$ in $\infty$-categories of simplicial ind-Banach rings in \cite{10BBBK}
 \begin{align}
  \mathrm{SimplicialInd-BanachRings}_{\mathbb{Q}_p}
\end{align}  
or animated analytic condensed commutative algebras in \cite{10CS2} 
\begin{align}   
\mathrm{SimplicialAnalyticCondensed}_{\mathbb{Q}_p}.
\end{align}   
Taking the product we have:
\[
\xymatrix@R+0pc@C+0pc{
\widetilde{\Phi}_{\psi,\Gamma,\square},\widetilde{\Phi}^r_{\psi,\Gamma,\square},\widetilde{\Phi}^I_{\psi,\Gamma,\square},	
}
\]
\[
\xymatrix@R+0pc@C+0pc{
\breve{\Phi}_{\psi,\Gamma,\square},\breve{\Phi}^r_{\psi,\Gamma,\square},\breve{\Phi}^I_{\psi,\Gamma,\square},	
}
\]
\[
\xymatrix@R+0pc@C+0pc{
{\Phi}_{\psi,\Gamma,\square},{\Phi}^r_{\psi,\Gamma,\square},{\Phi}^I_{\psi,\Gamma,\square}.	
}
\]
They carry multi Frobenius action $\varphi_\Gamma$ and multi $\mathrm{Lie}_\Gamma:=\mathbb{Z}_p^{\times\Gamma}$ action. In our current situation after \cite{10CKZ} and \cite{10PZ} we consider the following $(\infty,1)$-categories of $(\infty,1)$-modules.\\
\end{definition}

\begin{definition}
First we consider the Bambozzi-Kremnizer spectrum $\mathrm{Spec}^\mathrm{BK}(*)$ attached to any of those in the above from \cite{10BK} by taking derived rational localization:
\begin{align}
&\mathrm{Spec}^\mathrm{BK}\widetilde{\Phi}_{\psi,\Gamma,\square},\mathrm{Spec}^\mathrm{BK}\widetilde{\Phi}^r_{\psi,\Gamma,\square},\mathrm{Spec}^\mathrm{BK}\widetilde{\Phi}^I_{\psi,\Gamma,\square},	
\end{align}
\begin{align}
&\mathrm{Spec}^\mathrm{BK}\breve{\Phi}_{\psi,\Gamma,\square},\mathrm{Spec}^\mathrm{BK}\breve{\Phi}^r_{\psi,\Gamma,\square},\mathrm{Spec}^\mathrm{BK}\breve{\Phi}^I_{\psi,\Gamma,\square},	
\end{align}
\begin{align}
&\mathrm{Spec}^\mathrm{BK}{\Phi}_{\psi,\Gamma,\square},
\mathrm{Spec}^\mathrm{BK}{\Phi}^r_{\psi,\Gamma,\square},\mathrm{Spec}^\mathrm{BK}{\Phi}^I_{\psi,\Gamma,\square}.	
\end{align}

Then we take the corresponding quotients by using the corresponding Frobenius operators:
\begin{align}
&\mathrm{Spec}^\mathrm{BK}\widetilde{\Phi}_{\psi,\Gamma,\square}/\mathrm{Fro}^\mathbb{Z},	\\
\end{align}
\begin{align}
&\mathrm{Spec}^\mathrm{BK}\breve{\Phi}_{\psi,\Gamma,\square}/\mathrm{Fro}^\mathbb{Z},	\\
\end{align}
\begin{align}
&\mathrm{Spec}^\mathrm{BK}{\Phi}_{\psi,\Gamma,\square}/\mathrm{Fro}^\mathbb{Z}.	
\end{align}
Here for those space without notation related to the radius and the corresponding interval we consider the total unions $\bigcap_r,\bigcup_I$ in order to achieve the whole spaces to achieve the analogues of the corresponding FF curves from \cite{10KL1}, \cite{10KL2}, \cite{10FF} for
\[
\xymatrix@R+0pc@C+0pc{
\underset{r}{\mathrm{homotopylimit}}~\mathrm{Spec}^\mathrm{BK}\widetilde{\Phi}^r_{\psi,\Gamma,\square},\underset{I}{\mathrm{homotopycolimit}}~\mathrm{Spec}^\mathrm{BK}\widetilde{\Phi}^I_{\psi,\Gamma,\square},	\\
}
\]
\[
\xymatrix@R+0pc@C+0pc{
\underset{r}{\mathrm{homotopylimit}}~\mathrm{Spec}^\mathrm{BK}\breve{\Phi}^r_{\psi,\Gamma,\square},\underset{I}{\mathrm{homotopycolimit}}~\mathrm{Spec}^\mathrm{BK}\breve{\Phi}^I_{\psi,\Gamma,\square},	\\
}
\]
\[
\xymatrix@R+0pc@C+0pc{
\underset{r}{\mathrm{homotopylimit}}~\mathrm{Spec}^\mathrm{BK}{\Phi}^r_{\psi,\Gamma,\square},\underset{I}{\mathrm{homotopycolimit}}~\mathrm{Spec}^\mathrm{BK}{\Phi}^I_{\psi,\Gamma,\square}.	
}
\]
\[  
\xymatrix@R+0pc@C+0pc{
\underset{r}{\mathrm{homotopylimit}}~\mathrm{Spec}^\mathrm{BK}\widetilde{\Phi}^r_{\psi,\Gamma,\square}/\mathrm{Fro}^\mathbb{Z},\underset{I}{\mathrm{homotopycolimit}}~\mathrm{Spec}^\mathrm{BK}\widetilde{\Phi}^I_{\psi,\Gamma,\square}/\mathrm{Fro}^\mathbb{Z},	\\
}
\]
\[ 
\xymatrix@R+0pc@C+0pc{
\underset{r}{\mathrm{homotopylimit}}~\mathrm{Spec}^\mathrm{BK}\breve{\Phi}^r_{\psi,\Gamma,\square}/\mathrm{Fro}^\mathbb{Z},\underset{I}{\mathrm{homotopycolimit}}~\mathrm{Spec}^\mathrm{BK}\breve{\Phi}^I_{\psi,\Gamma,\square}/\mathrm{Fro}^\mathbb{Z},	\\
}
\]
\[ 
\xymatrix@R+0pc@C+0pc{
\underset{r}{\mathrm{homotopylimit}}~\mathrm{Spec}^\mathrm{BK}{\Phi}^r_{\psi,\Gamma,\square}/\mathrm{Fro}^\mathbb{Z},\underset{I}{\mathrm{homotopycolimit}}~\mathrm{Spec}^\mathrm{BK}{\Phi}^I_{\psi,\Gamma,\square}/\mathrm{Fro}^\mathbb{Z}.	
}
\]

\end{definition}

\indent Meanwhile we have the corresponding Clausen-Scholze analytic stacks from \cite{10CS2}, therefore applying their construction we have:

\begin{definition}
Here we define the following products by using the solidified tensor product from \cite{10CS1} and \cite{10CS2}. Namely $A$ will still as above as a Banach ring over $\mathbb{Q}_p$. Then we take solidified tensor product $\overset{\blacksquare}{\otimes}$ of any of the following
\[
\xymatrix@R+0pc@C+0pc{
\widetilde{\Delta}_{\psi,\Gamma},\widetilde{\nabla}_{\psi,\Gamma},\widetilde{\Phi}_{\psi,\Gamma},\widetilde{\Delta}^+_{\psi,\Gamma},\widetilde{\nabla}^+_{\psi,\Gamma},\widetilde{\Delta}^\dagger_{\psi,\Gamma},\widetilde{\nabla}^\dagger_{\psi,\Gamma},\widetilde{\Phi}^r_{\psi,\Gamma},\widetilde{\Phi}^I_{\psi,\Gamma}, 
}
\]

\[
\xymatrix@R+0pc@C+0pc{
\breve{\Delta}_{\psi,\Gamma},\breve{\nabla}_{\psi,\Gamma},\breve{\Phi}_{\psi,\Gamma},\breve{\Delta}^+_{\psi,\Gamma},\breve{\nabla}^+_{\psi,\Gamma},\breve{\Delta}^\dagger_{\psi,\Gamma},\breve{\nabla}^\dagger_{\psi,\Gamma},\breve{\Phi}^r_{\psi,\Gamma},\breve{\Phi}^I_{\psi,\Gamma},	
}
\]

\[
\xymatrix@R+0pc@C+0pc{
{\Delta}_{\psi,\Gamma},{\nabla}_{\psi,\Gamma},{\Phi}_{\psi,\Gamma},{\Delta}^+_{\psi,\Gamma},{\nabla}^+_{\psi,\Gamma},{\Delta}^\dagger_{\psi,\Gamma},{\nabla}^\dagger_{\psi,\Gamma},{\Phi}^r_{\psi,\Gamma},{\Phi}^I_{\psi,\Gamma},	
}
\]  	
with $A$. Then we have the notations:
\[
\xymatrix@R+0pc@C+0pc{
\widetilde{\Delta}_{\psi,\Gamma,\square},\widetilde{\nabla}_{\psi,\Gamma,\square},\widetilde{\Phi}_{\psi,\Gamma,\square},\widetilde{\Delta}^+_{\psi,\Gamma,\square},\widetilde{\nabla}^+_{\psi,\Gamma,\square},\widetilde{\Delta}^\dagger_{\psi,\Gamma,\square},\widetilde{\nabla}^\dagger_{\psi,\Gamma,\square},\widetilde{\Phi}^r_{\psi,\Gamma,\square},\widetilde{\Phi}^I_{\psi,\Gamma,\square}, 
}
\]

\[
\xymatrix@R+0pc@C+0pc{
\breve{\Delta}_{\psi,\Gamma,\square},\breve{\nabla}_{\psi,\Gamma,\square},\breve{\Phi}_{\psi,\Gamma,\square},\breve{\Delta}^+_{\psi,\Gamma,\square},\breve{\nabla}^+_{\psi,\Gamma,\square},\breve{\Delta}^\dagger_{\psi,\Gamma,\square},\breve{\nabla}^\dagger_{\psi,\Gamma,\square},\breve{\Phi}^r_{\psi,\Gamma,\square},\breve{\Phi}^I_{\psi,\Gamma,\square},	
}
\]

\[
\xymatrix@R+0pc@C+0pc{
{\Delta}_{\psi,\Gamma,\square},{\nabla}_{\psi,\Gamma,\square},{\Phi}_{\psi,\Gamma,\square},{\Delta}^+_{\psi,\Gamma,\square},{\nabla}^+_{\psi,\Gamma,\square},{\Delta}^\dagger_{\psi,\Gamma,\square},{\nabla}^\dagger_{\psi,\Gamma,\square},{\Phi}^r_{\psi,\Gamma,\square},{\Phi}^I_{\psi,\Gamma,\square}.	
}
\]
\end{definition}

\begin{definition}
First we consider the Clausen-Scholze spectrum $\mathrm{Spec}^\mathrm{CS}(*)$ attached to any of those in the above from \cite{10CS2} by taking derived rational localization:
\begin{align}
\mathrm{Spec}^\mathrm{CS}\widetilde{\Delta}_{\psi,\Gamma,\square},\mathrm{Spec}^\mathrm{CS}\widetilde{\nabla}_{\psi,\Gamma,\square},\mathrm{Spec}^\mathrm{CS}\widetilde{\Phi}_{\psi,\Gamma,\square},\mathrm{Spec}^\mathrm{CS}\widetilde{\Delta}^+_{\psi,\Gamma,\square},\mathrm{Spec}^\mathrm{CS}\widetilde{\nabla}^+_{\psi,\Gamma,\square},\\
\mathrm{Spec}^\mathrm{CS}\widetilde{\Delta}^\dagger_{\psi,\Gamma,\square},\mathrm{Spec}^\mathrm{CS}\widetilde{\nabla}^\dagger_{\psi,\Gamma,\square},\mathrm{Spec}^\mathrm{CS}\widetilde{\Phi}^r_{\psi,\Gamma,\square},\mathrm{Spec}^\mathrm{CS}\widetilde{\Phi}^I_{\psi,\Gamma,\square},	\\
\end{align}
\begin{align}
\mathrm{Spec}^\mathrm{CS}\breve{\Delta}_{\psi,\Gamma,\square},\breve{\nabla}_{\psi,\Gamma,\square},\mathrm{Spec}^\mathrm{CS}\breve{\Phi}_{\psi,\Gamma,\square},\mathrm{Spec}^\mathrm{CS}\breve{\Delta}^+_{\psi,\Gamma,\square},\mathrm{Spec}^\mathrm{CS}\breve{\nabla}^+_{\psi,\Gamma,\square},\\
\mathrm{Spec}^\mathrm{CS}\breve{\Delta}^\dagger_{\psi,\Gamma,\square},\mathrm{Spec}^\mathrm{CS}\breve{\nabla}^\dagger_{\psi,\Gamma,\square},\mathrm{Spec}^\mathrm{CS}\breve{\Phi}^r_{\psi,\Gamma,\square},\breve{\Phi}^I_{\psi,\Gamma,\square},	\\
\end{align}
\begin{align}
\mathrm{Spec}^\mathrm{CS}{\Delta}_{\psi,\Gamma,\square},\mathrm{Spec}^\mathrm{CS}{\nabla}_{\psi,\Gamma,\square},\mathrm{Spec}^\mathrm{CS}{\Phi}_{\psi,\Gamma,\square},\mathrm{Spec}^\mathrm{CS}{\Delta}^+_{\psi,\Gamma,\square},\mathrm{Spec}^\mathrm{CS}{\nabla}^+_{\psi,\Gamma,\square},\\
\mathrm{Spec}^\mathrm{CS}{\Delta}^\dagger_{\psi,\Gamma,\square},\mathrm{Spec}^\mathrm{CS}{\nabla}^\dagger_{\psi,\Gamma,\square},\mathrm{Spec}^\mathrm{CS}{\Phi}^r_{\psi,\Gamma,\square},\mathrm{Spec}^\mathrm{CS}{\Phi}^I_{\psi,\Gamma,\square}.	
\end{align}

Then we take the corresponding quotients by using the corresponding Frobenius operators:
\begin{align}
&\mathrm{Spec}^\mathrm{CS}\widetilde{\Delta}_{\psi,\Gamma,\square}/\mathrm{Fro}^\mathbb{Z},\mathrm{Spec}^\mathrm{CS}\widetilde{\nabla}_{\psi,\Gamma,\square}/\mathrm{Fro}^\mathbb{Z},\mathrm{Spec}^\mathrm{CS}\widetilde{\Phi}_{\psi,\Gamma,\square}/\mathrm{Fro}^\mathbb{Z},\mathrm{Spec}^\mathrm{CS}\widetilde{\Delta}^+_{\psi,\Gamma,\square}/\mathrm{Fro}^\mathbb{Z},\\
&\mathrm{Spec}^\mathrm{CS}\widetilde{\nabla}^+_{\psi,\Gamma,\square}/\mathrm{Fro}^\mathbb{Z}, \mathrm{Spec}^\mathrm{CS}\widetilde{\Delta}^\dagger_{\psi,\Gamma,\square}/\mathrm{Fro}^\mathbb{Z},\mathrm{Spec}^\mathrm{CS}\widetilde{\nabla}^\dagger_{\psi,\Gamma,\square}/\mathrm{Fro}^\mathbb{Z},	\\
\end{align}
\begin{align}
&\mathrm{Spec}^\mathrm{CS}\breve{\Delta}_{\psi,\Gamma,\square}/\mathrm{Fro}^\mathbb{Z},\breve{\nabla}_{\psi,\Gamma,\square}/\mathrm{Fro}^\mathbb{Z},\mathrm{Spec}^\mathrm{CS}\breve{\Phi}_{\psi,\Gamma,\square}/\mathrm{Fro}^\mathbb{Z},\mathrm{Spec}^\mathrm{CS}\breve{\Delta}^+_{\psi,\Gamma,\square}/\mathrm{Fro}^\mathbb{Z},\\
&\mathrm{Spec}^\mathrm{CS}\breve{\nabla}^+_{\psi,\Gamma,\square}/\mathrm{Fro}^\mathbb{Z}, \mathrm{Spec}^\mathrm{CS}\breve{\Delta}^\dagger_{\psi,\Gamma,\square}/\mathrm{Fro}^\mathbb{Z},\mathrm{Spec}^\mathrm{CS}\breve{\nabla}^\dagger_{\psi,\Gamma,\square}/\mathrm{Fro}^\mathbb{Z},	\\
\end{align}
\begin{align}
&\mathrm{Spec}^\mathrm{CS}{\Delta}_{\psi,\Gamma,\square}/\mathrm{Fro}^\mathbb{Z},\mathrm{Spec}^\mathrm{CS}{\nabla}_{\psi,\Gamma,\square}/\mathrm{Fro}^\mathbb{Z},\mathrm{Spec}^\mathrm{CS}{\Phi}_{\psi,\Gamma,\square}/\mathrm{Fro}^\mathbb{Z},\mathrm{Spec}^\mathrm{CS}{\Delta}^+_{\psi,\Gamma,\square}/\mathrm{Fro}^\mathbb{Z},\\
&\mathrm{Spec}^\mathrm{CS}{\nabla}^+_{\psi,\Gamma,\square}/\mathrm{Fro}^\mathbb{Z}, \mathrm{Spec}^\mathrm{CS}{\Delta}^\dagger_{\psi,\Gamma,\square}/\mathrm{Fro}^\mathbb{Z},\mathrm{Spec}^\mathrm{CS}{\nabla}^\dagger_{\psi,\Gamma,\square}/\mathrm{Fro}^\mathbb{Z}.	
\end{align}
Here for those space with notations related to the radius and the corresponding interval we consider the total unions $\bigcap_r,\bigcup_I$ in order to achieve the whole spaces to achieve the analogues of the corresponding FF curves from \cite{10KL1}, \cite{10KL2}, \cite{10FF} for
\[
\xymatrix@R+0pc@C+0pc{
\underset{r}{\mathrm{homotopylimit}}~\mathrm{Spec}^\mathrm{CS}\widetilde{\Phi}^r_{\psi,\Gamma,\square},\underset{I}{\mathrm{homotopycolimit}}~\mathrm{Spec}^\mathrm{CS}\widetilde{\Phi}^I_{\psi,\Gamma,\square},	\\
}
\]
\[
\xymatrix@R+0pc@C+0pc{
\underset{r}{\mathrm{homotopylimit}}~\mathrm{Spec}^\mathrm{CS}\breve{\Phi}^r_{\psi,\Gamma,\square},\underset{I}{\mathrm{homotopycolimit}}~\mathrm{Spec}^\mathrm{CS}\breve{\Phi}^I_{\psi,\Gamma,\square},	\\
}
\]
\[
\xymatrix@R+0pc@C+0pc{
\underset{r}{\mathrm{homotopylimit}}~\mathrm{Spec}^\mathrm{CS}{\Phi}^r_{\psi,\Gamma,\square},\underset{I}{\mathrm{homotopycolimit}}~\mathrm{Spec}^\mathrm{CS}{\Phi}^I_{\psi,\Gamma,\square}.	
}
\]
\[ 
\xymatrix@R+0pc@C+0pc{
\underset{r}{\mathrm{homotopylimit}}~\mathrm{Spec}^\mathrm{CS}\widetilde{\Phi}^r_{\psi,\Gamma,\square}/\mathrm{Fro}^\mathbb{Z},\underset{I}{\mathrm{homotopycolimit}}~\mathrm{Spec}^\mathrm{CS}\widetilde{\Phi}^I_{\psi,\Gamma,\square}/\mathrm{Fro}^\mathbb{Z},	\\
}
\]
\[ 
\xymatrix@R+0pc@C+0pc{
\underset{r}{\mathrm{homotopylimit}}~\mathrm{Spec}^\mathrm{CS}\breve{\Phi}^r_{\psi,\Gamma,\square}/\mathrm{Fro}^\mathbb{Z},\underset{I}{\mathrm{homotopycolimit}}~\breve{\Phi}^I_{\psi,\Gamma,\square}/\mathrm{Fro}^\mathbb{Z},	\\
}
\]
\[ 
\xymatrix@R+0pc@C+0pc{
\underset{r}{\mathrm{homotopylimit}}~\mathrm{Spec}^\mathrm{CS}{\Phi}^r_{\psi,\Gamma,\square}/\mathrm{Fro}^\mathbb{Z},\underset{I}{\mathrm{homotopycolimit}}~\mathrm{Spec}^\mathrm{CS}{\Phi}^I_{\psi,\Gamma,\square}/\mathrm{Fro}^\mathbb{Z}.	
}
\]

\end{definition}

\

\begin{definition}
We then consider the corresponding quasipresheaves of the corresponding ind-Banach or monomorphic ind-Banach modules from \cite{10BBK}, \cite{10KKM}\footnote{Here the categories are defined to be the corresponding homotopy colimits of the corresponding categories with respect to each $\square_i$.}:
\begin{align}
\mathrm{Quasicoherentpresheaves,IndBanach}_{*}	
\end{align}
where $*$ is one of the following spaces:
\begin{align}
&\mathrm{Spec}^\mathrm{BK}\widetilde{\Phi}_{\psi,\Gamma,\square}/\mathrm{Fro}^\mathbb{Z},	\\
\end{align}
\begin{align}
&\mathrm{Spec}^\mathrm{BK}\breve{\Phi}_{\psi,\Gamma,\square}/\mathrm{Fro}^\mathbb{Z},	\\
\end{align}
\begin{align}
&\mathrm{Spec}^\mathrm{BK}{\Phi}_{\psi,\Gamma,\square}/\mathrm{Fro}^\mathbb{Z}.	
\end{align}
Here for those space without notation related to the radius and the corresponding interval we consider the total unions $\bigcap_r,\bigcup_I$ in order to achieve the whole spaces to achieve the analogues of the corresponding FF curves from \cite{10KL1}, \cite{10KL2}, \cite{10FF} for
\[
\xymatrix@R+0pc@C+0pc{
\underset{r}{\mathrm{homotopylimit}}~\mathrm{Spec}^\mathrm{BK}\widetilde{\Phi}^r_{\psi,\Gamma,\square},\underset{I}{\mathrm{homotopycolimit}}~\mathrm{Spec}^\mathrm{BK}\widetilde{\Phi}^I_{\psi,\Gamma,\square},	\\
}
\]
\[
\xymatrix@R+0pc@C+0pc{
\underset{r}{\mathrm{homotopylimit}}~\mathrm{Spec}^\mathrm{BK}\breve{\Phi}^r_{\psi,\Gamma,\square},\underset{I}{\mathrm{homotopycolimit}}~\mathrm{Spec}^\mathrm{BK}\breve{\Phi}^I_{\psi,\Gamma,\square},	\\
}
\]
\[
\xymatrix@R+0pc@C+0pc{
\underset{r}{\mathrm{homotopylimit}}~\mathrm{Spec}^\mathrm{BK}{\Phi}^r_{\psi,\Gamma,\square},\underset{I}{\mathrm{homotopycolimit}}~\mathrm{Spec}^\mathrm{BK}{\Phi}^I_{\psi,\Gamma,\square}.	
}
\]
\[  
\xymatrix@R+0pc@C+0pc{
\underset{r}{\mathrm{homotopylimit}}~\mathrm{Spec}^\mathrm{BK}\widetilde{\Phi}^r_{\psi,\Gamma,\square}/\mathrm{Fro}^\mathbb{Z},\underset{I}{\mathrm{homotopycolimit}}~\mathrm{Spec}^\mathrm{BK}\widetilde{\Phi}^I_{\psi,\Gamma,\square}/\mathrm{Fro}^\mathbb{Z},	\\
}
\]
\[ 
\xymatrix@R+0pc@C+0pc{
\underset{r}{\mathrm{homotopylimit}}~\mathrm{Spec}^\mathrm{BK}\breve{\Phi}^r_{\psi,\Gamma,\square}/\mathrm{Fro}^\mathbb{Z},\underset{I}{\mathrm{homotopycolimit}}~\mathrm{Spec}^\mathrm{BK}\breve{\Phi}^I_{\psi,\Gamma,\square}/\mathrm{Fro}^\mathbb{Z},	\\
}
\]
\[ 
\xymatrix@R+0pc@C+0pc{
\underset{r}{\mathrm{homotopylimit}}~\mathrm{Spec}^\mathrm{BK}{\Phi}^r_{\psi,\Gamma,\square}/\mathrm{Fro}^\mathbb{Z},\underset{I}{\mathrm{homotopycolimit}}~\mathrm{Spec}^\mathrm{BK}{\Phi}^I_{\psi,\Gamma,\square}/\mathrm{Fro}^\mathbb{Z}.	
}
\]

\end{definition}

\begin{definition}
We then consider the corresponding quasisheaves of the corresponding condensed solid topological modules from \cite{10CS2}:
\begin{align}
\mathrm{Quasicoherentsheaves, Condensed}_{*}	
\end{align}
where $*$ is one of the following spaces:
\begin{align}
&\mathrm{Spec}^\mathrm{CS}\widetilde{\Delta}_{\psi,\Gamma,\square}/\mathrm{Fro}^\mathbb{Z},\mathrm{Spec}^\mathrm{CS}\widetilde{\nabla}_{\psi,\Gamma,\square}/\mathrm{Fro}^\mathbb{Z},\mathrm{Spec}^\mathrm{CS}\widetilde{\Phi}_{\psi,\Gamma,\square}/\mathrm{Fro}^\mathbb{Z},\mathrm{Spec}^\mathrm{CS}\widetilde{\Delta}^+_{\psi,\Gamma,\square}/\mathrm{Fro}^\mathbb{Z},\\
&\mathrm{Spec}^\mathrm{CS}\widetilde{\nabla}^+_{\psi,\Gamma,\square}/\mathrm{Fro}^\mathbb{Z},\mathrm{Spec}^\mathrm{CS}\widetilde{\Delta}^\dagger_{\psi,\Gamma,\square}/\mathrm{Fro}^\mathbb{Z},\mathrm{Spec}^\mathrm{CS}\widetilde{\nabla}^\dagger_{\psi,\Gamma,\square}/\mathrm{Fro}^\mathbb{Z},	\\
\end{align}
\begin{align}
&\mathrm{Spec}^\mathrm{CS}\breve{\Delta}_{\psi,\Gamma,\square}/\mathrm{Fro}^\mathbb{Z},\breve{\nabla}_{\psi,\Gamma,\square}/\mathrm{Fro}^\mathbb{Z},\mathrm{Spec}^\mathrm{CS}\breve{\Phi}_{\psi,\Gamma,\square}/\mathrm{Fro}^\mathbb{Z},\mathrm{Spec}^\mathrm{CS}\breve{\Delta}^+_{\psi,\Gamma,\square}/\mathrm{Fro}^\mathbb{Z},\\
&\mathrm{Spec}^\mathrm{CS}\breve{\nabla}^+_{\psi,\Gamma,\square}/\mathrm{Fro}^\mathbb{Z},\mathrm{Spec}^\mathrm{CS}\breve{\Delta}^\dagger_{\psi,\Gamma,\square}/\mathrm{Fro}^\mathbb{Z},\mathrm{Spec}^\mathrm{CS}\breve{\nabla}^\dagger_{\psi,\Gamma,\square}/\mathrm{Fro}^\mathbb{Z},	\\
\end{align}
\begin{align}
&\mathrm{Spec}^\mathrm{CS}{\Delta}_{\psi,\Gamma,\square}/\mathrm{Fro}^\mathbb{Z},\mathrm{Spec}^\mathrm{CS}{\nabla}_{\psi,\Gamma,\square}/\mathrm{Fro}^\mathbb{Z},\mathrm{Spec}^\mathrm{CS}{\Phi}_{\psi,\Gamma,\square}/\mathrm{Fro}^\mathbb{Z},\mathrm{Spec}^\mathrm{CS}{\Delta}^+_{\psi,\Gamma,\square}/\mathrm{Fro}^\mathbb{Z},\\
&\mathrm{Spec}^\mathrm{CS}{\nabla}^+_{\psi,\Gamma,\square}/\mathrm{Fro}^\mathbb{Z}, \mathrm{Spec}^\mathrm{CS}{\Delta}^\dagger_{\psi,\Gamma,\square}/\mathrm{Fro}^\mathbb{Z},\mathrm{Spec}^\mathrm{CS}{\nabla}^\dagger_{\psi,\Gamma,\square}/\mathrm{Fro}^\mathbb{Z}.	
\end{align}
Here for those space with notations related to the radius and the corresponding interval we consider the total unions $\bigcap_r,\bigcup_I$ in order to achieve the whole spaces to achieve the analogues of the corresponding FF curves from \cite{10KL1}, \cite{10KL2}, \cite{10FF} for
\[
\xymatrix@R+0pc@C+0pc{
\underset{r}{\mathrm{homotopylimit}}~\mathrm{Spec}^\mathrm{CS}\widetilde{\Phi}^r_{\psi,\Gamma,\square},\underset{I}{\mathrm{homotopycolimit}}~\mathrm{Spec}^\mathrm{CS}\widetilde{\Phi}^I_{\psi,\Gamma,\square},	\\
}
\]
\[
\xymatrix@R+0pc@C+0pc{
\underset{r}{\mathrm{homotopylimit}}~\mathrm{Spec}^\mathrm{CS}\breve{\Phi}^r_{\psi,\Gamma,\square},\underset{I}{\mathrm{homotopycolimit}}~\mathrm{Spec}^\mathrm{CS}\breve{\Phi}^I_{\psi,\Gamma,\square},	\\
}
\]
\[
\xymatrix@R+0pc@C+0pc{
\underset{r}{\mathrm{homotopylimit}}~\mathrm{Spec}^\mathrm{CS}{\Phi}^r_{\psi,\Gamma,\square},\underset{I}{\mathrm{homotopycolimit}}~\mathrm{Spec}^\mathrm{CS}{\Phi}^I_{\psi,\Gamma,\square}.	
}
\]
\[ 
\xymatrix@R+0pc@C+0pc{
\underset{r}{\mathrm{homotopylimit}}~\mathrm{Spec}^\mathrm{CS}\widetilde{\Phi}^r_{\psi,\Gamma,\square}/\mathrm{Fro}^\mathbb{Z},\underset{I}{\mathrm{homotopycolimit}}~\mathrm{Spec}^\mathrm{CS}\widetilde{\Phi}^I_{\psi,\Gamma,\square}/\mathrm{Fro}^\mathbb{Z},	\\
}
\]
\[ 
\xymatrix@R+0pc@C+0pc{
\underset{r}{\mathrm{homotopylimit}}~\mathrm{Spec}^\mathrm{CS}\breve{\Phi}^r_{\psi,\Gamma,\square}/\mathrm{Fro}^\mathbb{Z},\underset{I}{\mathrm{homotopycolimit}}~\breve{\Phi}^I_{\psi,\Gamma,\square}/\mathrm{Fro}^\mathbb{Z},	\\
}
\]
\[ 
\xymatrix@R+0pc@C+0pc{
\underset{r}{\mathrm{homotopylimit}}~\mathrm{Spec}^\mathrm{CS}{\Phi}^r_{\psi,\Gamma,\square}/\mathrm{Fro}^\mathbb{Z},\underset{I}{\mathrm{homotopycolimit}}~\mathrm{Spec}^\mathrm{CS}{\Phi}^I_{\psi,\Gamma,\square}/\mathrm{Fro}^\mathbb{Z}.	
}
\]

\end{definition}

\

\begin{proposition}
There is a well-defined functor from the $\infty$-category 
\begin{align}
\mathrm{Quasicoherentpresheaves,Condensed}_{*}	
\end{align}
where $*$ is one of the following spaces:
\begin{align}
&\mathrm{Spec}^\mathrm{CS}\widetilde{\Phi}_{\psi,\Gamma,\square}/\mathrm{Fro}^\mathbb{Z},	\\
\end{align}
\begin{align}
&\mathrm{Spec}^\mathrm{CS}\breve{\Phi}_{\psi,\Gamma,\square}/\mathrm{Fro}^\mathbb{Z},	\\
\end{align}
\begin{align}
&\mathrm{Spec}^\mathrm{CS}{\Phi}_{\psi,\Gamma,\square}/\mathrm{Fro}^\mathbb{Z},	
\end{align}
to the $\infty$-category of $\mathrm{Fro}$-equivariant quasicoherent presheaves over similar spaces above correspondingly without the $\mathrm{Fro}$-quotients, and to the $\infty$-category of $\mathrm{Fro}$-equivariant quasicoherent modules over global sections of the structure $\infty$-sheaves of the similar spaces above correspondingly without the $\mathrm{Fro}$-quotients. Here for those space without notation related to the radius and the corresponding interval we consider the total unions $\bigcap_r,\bigcup_I$ in order to achieve the whole spaces to achieve the analogues of the corresponding FF curves from \cite{10KL1}, \cite{10KL2}, \cite{10FF} for
\[
\xymatrix@R+0pc@C+0pc{
\underset{r}{\mathrm{homotopylimit}}~\mathrm{Spec}^\mathrm{CS}\widetilde{\Phi}^r_{\psi,\Gamma,\square},\underset{I}{\mathrm{homotopycolimit}}~\mathrm{Spec}^\mathrm{CS}\widetilde{\Phi}^I_{\psi,\Gamma,\square},	\\
}
\]
\[
\xymatrix@R+0pc@C+0pc{
\underset{r}{\mathrm{homotopylimit}}~\mathrm{Spec}^\mathrm{CS}\breve{\Phi}^r_{\psi,\Gamma,\square},\underset{I}{\mathrm{homotopycolimit}}~\mathrm{Spec}^\mathrm{CS}\breve{\Phi}^I_{\psi,\Gamma,\square},	\\
}
\]
\[
\xymatrix@R+0pc@C+0pc{
\underset{r}{\mathrm{homotopylimit}}~\mathrm{Spec}^\mathrm{CS}{\Phi}^r_{\psi,\Gamma,\square},\underset{I}{\mathrm{homotopycolimit}}~\mathrm{Spec}^\mathrm{CS}{\Phi}^I_{\psi,\Gamma,\square}.	
}
\]
\[ 
\xymatrix@R+0pc@C+0pc{
\underset{r}{\mathrm{homotopylimit}}~\mathrm{Spec}^\mathrm{CS}\widetilde{\Phi}^r_{\psi,\Gamma,\square}/\mathrm{Fro}^\mathbb{Z},\underset{I}{\mathrm{homotopycolimit}}~\mathrm{Spec}^\mathrm{CS}\widetilde{\Phi}^I_{\psi,\Gamma,\square}/\mathrm{Fro}^\mathbb{Z},	\\
}
\]
\[ 
\xymatrix@R+0pc@C+0pc{
\underset{r}{\mathrm{homotopylimit}}~\mathrm{Spec}^\mathrm{CS}\breve{\Phi}^r_{\psi,\Gamma,\square}/\mathrm{Fro}^\mathbb{Z},\underset{I}{\mathrm{homotopycolimit}}~\breve{\Phi}^I_{\psi,\Gamma,\square}/\mathrm{Fro}^\mathbb{Z},	\\
}
\]
\[ 
\xymatrix@R+0pc@C+0pc{
\underset{r}{\mathrm{homotopylimit}}~\mathrm{Spec}^\mathrm{CS}{\Phi}^r_{\psi,\Gamma,\square}/\mathrm{Fro}^\mathbb{Z},\underset{I}{\mathrm{homotopycolimit}}~\mathrm{Spec}^\mathrm{CS}{\Phi}^I_{\psi,\Gamma,\square}/\mathrm{Fro}^\mathbb{Z}.	
}
\]	
In this situation we will have the target category being family parametrized by $r$ or $I$ in compatible glueing sense as in \cite[Definition 5.4.10]{10KL2}. In this situation for modules parametrized by the intervals we have the equivalence of $\infty$-categories by using \cite[Proposition 13.8]{10CS2}. Here the corresponding quasicoherent Frobenius modules are defined to be the corresponding homotopy colimits and limits of Frobenius modules:
\begin{align}
\underset{r}{\mathrm{homotopycolimit}}~M_r,\\
\underset{I}{\mathrm{homotopylimit}}~M_I,	
\end{align}
where each $M_r$ is a Frobenius-equivariant module over the period ring with respect to some radius $r$ while each $M_I$ is a Frobenius-equivariant module over the period ring with respect to some interval $I$.\\
\end{proposition}

\begin{proposition}
Similar proposition holds for 
\begin{align}
\mathrm{Quasicoherentsheaves,IndBanach}_{*}.	
\end{align}	
\end{proposition}

\

\begin{definition}
We then consider the corresponding quasipresheaves of perfect complexes the corresponding ind-Banach or monomorphic ind-Banach modules from \cite{10BBK}, \cite{10KKM}:
\begin{align}
\mathrm{Quasicoherentpresheaves,Perfectcomplex,IndBanach}_{*}	
\end{align}
where $*$ is one of the following spaces:
\begin{align}
&\mathrm{Spec}^\mathrm{BK}\widetilde{\Phi}_{\psi,\Gamma,\square}/\mathrm{Fro}^\mathbb{Z},	\\
\end{align}
\begin{align}
&\mathrm{Spec}^\mathrm{BK}\breve{\Phi}_{\psi,\Gamma,\square}/\mathrm{Fro}^\mathbb{Z},	\\
\end{align}
\begin{align}
&\mathrm{Spec}^\mathrm{BK}{\Phi}_{\psi,\Gamma,\square}/\mathrm{Fro}^\mathbb{Z}.	
\end{align}
Here for those space without notation related to the radius and the corresponding interval we consider the total unions $\bigcap_r,\bigcup_I$ in order to achieve the whole spaces to achieve the analogues of the corresponding FF curves from \cite{10KL1}, \cite{10KL2}, \cite{10FF} for
\[
\xymatrix@R+0pc@C+0pc{
\underset{r}{\mathrm{homotopylimit}}~\mathrm{Spec}^\mathrm{BK}\widetilde{\Phi}^r_{\psi,\Gamma,\square},\underset{I}{\mathrm{homotopycolimit}}~\mathrm{Spec}^\mathrm{BK}\widetilde{\Phi}^I_{\psi,\Gamma,\square},	\\
}
\]
\[
\xymatrix@R+0pc@C+0pc{
\underset{r}{\mathrm{homotopylimit}}~\mathrm{Spec}^\mathrm{BK}\breve{\Phi}^r_{\psi,\Gamma,\square},\underset{I}{\mathrm{homotopycolimit}}~\mathrm{Spec}^\mathrm{BK}\breve{\Phi}^I_{\psi,\Gamma,\square},	\\
}
\]
\[
\xymatrix@R+0pc@C+0pc{
\underset{r}{\mathrm{homotopylimit}}~\mathrm{Spec}^\mathrm{BK}{\Phi}^r_{\psi,\Gamma,\square},\underset{I}{\mathrm{homotopycolimit}}~\mathrm{Spec}^\mathrm{BK}{\Phi}^I_{\psi,\Gamma,\square}.	
}
\]
\[  
\xymatrix@R+0pc@C+0pc{
\underset{r}{\mathrm{homotopylimit}}~\mathrm{Spec}^\mathrm{BK}\widetilde{\Phi}^r_{\psi,\Gamma,\square}/\mathrm{Fro}^\mathbb{Z},\underset{I}{\mathrm{homotopycolimit}}~\mathrm{Spec}^\mathrm{BK}\widetilde{\Phi}^I_{\psi,\Gamma,\square}/\mathrm{Fro}^\mathbb{Z},	\\
}
\]
\[ 
\xymatrix@R+0pc@C+0pc{
\underset{r}{\mathrm{homotopylimit}}~\mathrm{Spec}^\mathrm{BK}\breve{\Phi}^r_{\psi,\Gamma,\square}/\mathrm{Fro}^\mathbb{Z},\underset{I}{\mathrm{homotopycolimit}}~\mathrm{Spec}^\mathrm{BK}\breve{\Phi}^I_{\psi,\Gamma,\square}/\mathrm{Fro}^\mathbb{Z},	\\
}
\]
\[ 
\xymatrix@R+0pc@C+0pc{
\underset{r}{\mathrm{homotopylimit}}~\mathrm{Spec}^\mathrm{BK}{\Phi}^r_{\psi,\Gamma,\square}/\mathrm{Fro}^\mathbb{Z},\underset{I}{\mathrm{homotopycolimit}}~\mathrm{Spec}^\mathrm{BK}{\Phi}^I_{\psi,\Gamma,\square}/\mathrm{Fro}^\mathbb{Z}.	
}
\]

\end{definition}

\begin{definition}
We then consider the corresponding quasisheaves of perfect complexes of the corresponding condensed solid topological modules from \cite{10CS2}:
\begin{align}
\mathrm{Quasicoherentsheaves, Perfectcomplex, Condensed}_{*}	
\end{align}
where $*$ is one of the following spaces:
\begin{align}
&\mathrm{Spec}^\mathrm{CS}\widetilde{\Delta}_{\psi,\Gamma,\square}/\mathrm{Fro}^\mathbb{Z},\mathrm{Spec}^\mathrm{CS}\widetilde{\nabla}_{\psi,\Gamma,\square}/\mathrm{Fro}^\mathbb{Z},\mathrm{Spec}^\mathrm{CS}\widetilde{\Phi}_{\psi,\Gamma,\square}/\mathrm{Fro}^\mathbb{Z},\mathrm{Spec}^\mathrm{CS}\widetilde{\Delta}^+_{\psi,\Gamma,\square}/\mathrm{Fro}^\mathbb{Z},\\
&\mathrm{Spec}^\mathrm{CS}\widetilde{\nabla}^+_{\psi,\Gamma,\square}/\mathrm{Fro}^\mathbb{Z},\mathrm{Spec}^\mathrm{CS}\widetilde{\Delta}^\dagger_{\psi,\Gamma,\square}/\mathrm{Fro}^\mathbb{Z},\mathrm{Spec}^\mathrm{CS}\widetilde{\nabla}^\dagger_{\psi,\Gamma,\square}/\mathrm{Fro}^\mathbb{Z},	\\
\end{align}
\begin{align}
&\mathrm{Spec}^\mathrm{CS}\breve{\Delta}_{\psi,\Gamma,\square}/\mathrm{Fro}^\mathbb{Z},\breve{\nabla}_{\psi,\Gamma,\square}/\mathrm{Fro}^\mathbb{Z},\mathrm{Spec}^\mathrm{CS}\breve{\Phi}_{\psi,\Gamma,\square}/\mathrm{Fro}^\mathbb{Z},\mathrm{Spec}^\mathrm{CS}\breve{\Delta}^+_{\psi,\Gamma,\square}/\mathrm{Fro}^\mathbb{Z},\\
&\mathrm{Spec}^\mathrm{CS}\breve{\nabla}^+_{\psi,\Gamma,\square}/\mathrm{Fro}^\mathbb{Z},\mathrm{Spec}^\mathrm{CS}\breve{\Delta}^\dagger_{\psi,\Gamma,\square}/\mathrm{Fro}^\mathbb{Z},\mathrm{Spec}^\mathrm{CS}\breve{\nabla}^\dagger_{\psi,\Gamma,\square}/\mathrm{Fro}^\mathbb{Z},	\\
\end{align}
\begin{align}
&\mathrm{Spec}^\mathrm{CS}{\Delta}_{\psi,\Gamma,\square}/\mathrm{Fro}^\mathbb{Z},\mathrm{Spec}^\mathrm{CS}{\nabla}_{\psi,\Gamma,\square}/\mathrm{Fro}^\mathbb{Z},\mathrm{Spec}^\mathrm{CS}{\Phi}_{\psi,\Gamma,\square}/\mathrm{Fro}^\mathbb{Z},\mathrm{Spec}^\mathrm{CS}{\Delta}^+_{\psi,\Gamma,\square}/\mathrm{Fro}^\mathbb{Z},\\
&\mathrm{Spec}^\mathrm{CS}{\nabla}^+_{\psi,\Gamma,\square}/\mathrm{Fro}^\mathbb{Z}, \mathrm{Spec}^\mathrm{CS}{\Delta}^\dagger_{\psi,\Gamma,\square}/\mathrm{Fro}^\mathbb{Z},\mathrm{Spec}^\mathrm{CS}{\nabla}^\dagger_{\psi,\Gamma,\square}/\mathrm{Fro}^\mathbb{Z}.	
\end{align}
Here for those space with notations related to the radius and the corresponding interval we consider the total unions $\bigcap_r,\bigcup_I$ in order to achieve the whole spaces to achieve the analogues of the corresponding FF curves from \cite{10KL1}, \cite{10KL2}, \cite{10FF} for
\[
\xymatrix@R+0pc@C+0pc{
\underset{r}{\mathrm{homotopylimit}}~\mathrm{Spec}^\mathrm{CS}\widetilde{\Phi}^r_{\psi,\Gamma,\square},\underset{I}{\mathrm{homotopycolimit}}~\mathrm{Spec}^\mathrm{CS}\widetilde{\Phi}^I_{\psi,\Gamma,\square},	\\
}
\]
\[
\xymatrix@R+0pc@C+0pc{
\underset{r}{\mathrm{homotopylimit}}~\mathrm{Spec}^\mathrm{CS}\breve{\Phi}^r_{\psi,\Gamma,\square},\underset{I}{\mathrm{homotopycolimit}}~\mathrm{Spec}^\mathrm{CS}\breve{\Phi}^I_{\psi,\Gamma,\square},	\\
}
\]
\[
\xymatrix@R+0pc@C+0pc{
\underset{r}{\mathrm{homotopylimit}}~\mathrm{Spec}^\mathrm{CS}{\Phi}^r_{\psi,\Gamma,\square},\underset{I}{\mathrm{homotopycolimit}}~\mathrm{Spec}^\mathrm{CS}{\Phi}^I_{\psi,\Gamma,\square}.	
}
\]
\[ 
\xymatrix@R+0pc@C+0pc{
\underset{r}{\mathrm{homotopylimit}}~\mathrm{Spec}^\mathrm{CS}\widetilde{\Phi}^r_{\psi,\Gamma,\square}/\mathrm{Fro}^\mathbb{Z},\underset{I}{\mathrm{homotopycolimit}}~\mathrm{Spec}^\mathrm{CS}\widetilde{\Phi}^I_{\psi,\Gamma,\square}/\mathrm{Fro}^\mathbb{Z},	\\
}
\]
\[ 
\xymatrix@R+0pc@C+0pc{
\underset{r}{\mathrm{homotopylimit}}~\mathrm{Spec}^\mathrm{CS}\breve{\Phi}^r_{\psi,\Gamma,\square}/\mathrm{Fro}^\mathbb{Z},\underset{I}{\mathrm{homotopycolimit}}~\breve{\Phi}^I_{\psi,\Gamma,\square}/\mathrm{Fro}^\mathbb{Z},	\\
}
\]
\[ 
\xymatrix@R+0pc@C+0pc{
\underset{r}{\mathrm{homotopylimit}}~\mathrm{Spec}^\mathrm{CS}{\Phi}^r_{\psi,\Gamma,\square}/\mathrm{Fro}^\mathbb{Z},\underset{I}{\mathrm{homotopycolimit}}~\mathrm{Spec}^\mathrm{CS}{\Phi}^I_{\psi,\Gamma,\square}/\mathrm{Fro}^\mathbb{Z}.	
}
\]

\end{definition}

\begin{proposition}
There is a well-defined functor from the $\infty$-category 
\begin{align}
\mathrm{Quasicoherentpresheaves,Perfectcomplex,Condensed}_{*}	
\end{align}
where $*$ is one of the following spaces:
\begin{align}
&\mathrm{Spec}^\mathrm{CS}\widetilde{\Phi}_{\psi,\Gamma,\square}/\mathrm{Fro}^\mathbb{Z},	\\
\end{align}
\begin{align}
&\mathrm{Spec}^\mathrm{CS}\breve{\Phi}_{\psi,\Gamma,\square}/\mathrm{Fro}^\mathbb{Z},	\\
\end{align}
\begin{align}
&\mathrm{Spec}^\mathrm{CS}{\Phi}_{\psi,\Gamma,\square}/\mathrm{Fro}^\mathbb{Z},	
\end{align}
to the $\infty$-category of $\mathrm{Fro}$-equivariant quasicoherent presheaves over similar spaces above correspondingly without the $\mathrm{Fro}$-quotients, and to the $\infty$-category of $\mathrm{Fro}$-equivariant quasicoherent modules over global sections of the structure $\infty$-sheaves of the similar spaces above correspondingly without the $\mathrm{Fro}$-quotients. Here for those space without notation related to the radius and the corresponding interval we consider the total unions $\bigcap_r,\bigcup_I$ in order to achieve the whole spaces to achieve the analogues of the corresponding FF curves from \cite{10KL1}, \cite{10KL2}, \cite{10FF} for
\[
\xymatrix@R+0pc@C+0pc{
\underset{r}{\mathrm{homotopylimit}}~\mathrm{Spec}^\mathrm{CS}\widetilde{\Phi}^r_{\psi,\Gamma,\square},\underset{I}{\mathrm{homotopycolimit}}~\mathrm{Spec}^\mathrm{CS}\widetilde{\Phi}^I_{\psi,\Gamma,\square},	\\
}
\]
\[
\xymatrix@R+0pc@C+0pc{
\underset{r}{\mathrm{homotopylimit}}~\mathrm{Spec}^\mathrm{CS}\breve{\Phi}^r_{\psi,\Gamma,\square},\underset{I}{\mathrm{homotopycolimit}}~\mathrm{Spec}^\mathrm{CS}\breve{\Phi}^I_{\psi,\Gamma,\square},	\\
}
\]
\[
\xymatrix@R+0pc@C+0pc{
\underset{r}{\mathrm{homotopylimit}}~\mathrm{Spec}^\mathrm{CS}{\Phi}^r_{\psi,\Gamma,\square},\underset{I}{\mathrm{homotopycolimit}}~\mathrm{Spec}^\mathrm{CS}{\Phi}^I_{\psi,\Gamma,\square}.	
}
\]
\[ 
\xymatrix@R+0pc@C+0pc{
\underset{r}{\mathrm{homotopylimit}}~\mathrm{Spec}^\mathrm{CS}\widetilde{\Phi}^r_{\psi,\Gamma,\square}/\mathrm{Fro}^\mathbb{Z},\underset{I}{\mathrm{homotopycolimit}}~\mathrm{Spec}^\mathrm{CS}\widetilde{\Phi}^I_{\psi,\Gamma,\square}/\mathrm{Fro}^\mathbb{Z},	\\
}
\]
\[ 
\xymatrix@R+0pc@C+0pc{
\underset{r}{\mathrm{homotopylimit}}~\mathrm{Spec}^\mathrm{CS}\breve{\Phi}^r_{\psi,\Gamma,\square}/\mathrm{Fro}^\mathbb{Z},\underset{I}{\mathrm{homotopycolimit}}~\breve{\Phi}^I_{\psi,\Gamma,\square}/\mathrm{Fro}^\mathbb{Z},	\\
}
\]
\[ 
\xymatrix@R+0pc@C+0pc{
\underset{r}{\mathrm{homotopylimit}}~\mathrm{Spec}^\mathrm{CS}{\Phi}^r_{\psi,\Gamma,\square}/\mathrm{Fro}^\mathbb{Z},\underset{I}{\mathrm{homotopycolimit}}~\mathrm{Spec}^\mathrm{CS}{\Phi}^I_{\psi,\Gamma,\square}/\mathrm{Fro}^\mathbb{Z}.	
}
\]	
In this situation we will have the target category being family parametrized by $r$ or $I$ in compatible glueing sense as in \cite[Definition 5.4.10]{10KL2}. In this situation for modules parametrized by the intervals we have the equivalence of $\infty$-categories by using \cite[Proposition 12.18]{10CS2}. Here the corresponding quasicoherent Frobenius modules are defined to be the corresponding homotopy colimits and limits of Frobenius modules:
\begin{align}
\underset{r}{\mathrm{homotopycolimit}}~M_r,\\
\underset{I}{\mathrm{homotopylimit}}~M_I,	
\end{align}
where each $M_r$ is a Frobenius-equivariant module over the period ring with respect to some radius $r$ while each $M_I$ is a Frobenius-equivariant module over the period ring with respect to some interval $I$.\\
\end{proposition}

\begin{proposition}
Similar proposition holds for 
\begin{align}
\mathrm{Quasicoherentsheaves,Perfectcomplex,IndBanach}_{*}.	
\end{align}	
\end{proposition}

\section{Univariate Hodge Iwasawa Modules}

This chapter follows closely \cite{10T1}, \cite{10T2}, \cite{10T3}, \cite{10KPX}, \cite{10KP}, \cite{10KL1}, \cite{10KL2}, \cite{10BK}, \cite{10BBBK}, \cite{10BBM}, \cite{10KKM}, \cite{10CS1}, \cite{10CS2}, \cite{10LBV}. 

\subsection{Frobenius Quasicoherent Modules I}

\begin{definition}
Let $\psi$ be a toric tower over $\mathbb{Q}_p$ as in \cite[Chapter 7]{10KL2} with base $\mathbb{Q}_p\left<X_1^{\pm 1},...,X_k^{\pm 1}\right>$. Then from \cite{10KL1} and \cite[Definition 5.2.1]{10KL2} we have the following class of Kedlaya-Liu rings (with the following replacement: $\Delta$ stands for $A$, $\nabla$ stands for $B$, while $\Phi$ stands for $C$) by taking product in the sense of self $\Gamma$-th power\footnote{Here $|\Gamma|=1$.}:

\[
\xymatrix@R+0pc@C+0pc{
\widetilde{\Delta}_{\psi},\widetilde{\nabla}_{\psi},\widetilde{\Phi}_{\psi},\widetilde{\Delta}^+_{\psi},\widetilde{\nabla}^+_{\psi},\widetilde{\Delta}^\dagger_{\psi},\widetilde{\nabla}^\dagger_{\psi},\widetilde{\Phi}^r_{\psi},\widetilde{\Phi}^I_{\psi}, 
}
\]

\[
\xymatrix@R+0pc@C+0pc{
\breve{\Delta}_{\psi},\breve{\nabla}_{\psi},\breve{\Phi}_{\psi},\breve{\Delta}^+_{\psi},\breve{\nabla}^+_{\psi},\breve{\Delta}^\dagger_{\psi},\breve{\nabla}^\dagger_{\psi},\breve{\Phi}^r_{\psi},\breve{\Phi}^I_{\psi},	
}
\]

\[
\xymatrix@R+0pc@C+0pc{
{\Delta}_{\psi},{\nabla}_{\psi},{\Phi}_{\psi},{\Delta}^+_{\psi},{\nabla}^+_{\psi},{\Delta}^\dagger_{\psi},{\nabla}^\dagger_{\psi},{\Phi}^r_{\psi},{\Phi}^I_{\psi}.	
}
\]
We now consider the following rings with $A$ being a Banach ring over $\mathbb{Q}_p$. Taking the product we have:
\[
\xymatrix@R+0pc@C+0pc{
\widetilde{\Phi}_{\psi,A},\widetilde{\Phi}^r_{\psi,A},\widetilde{\Phi}^I_{\psi,A},	
}
\]
\[
\xymatrix@R+0pc@C+0pc{
\breve{\Phi}_{\psi,A},\breve{\Phi}^r_{\psi,A},\breve{\Phi}^I_{\psi,A},	
}
\]
\[
\xymatrix@R+0pc@C+0pc{
{\Phi}_{\psi,A},{\Phi}^r_{\psi,A},{\Phi}^I_{\psi,A}.	
}
\]
They carry multi Frobenius action $\varphi_\Gamma$ and multi $\mathrm{Lie}_\Gamma:=\mathbb{Z}_p^{\times\Gamma}$ action. In our current situation after \cite{10CKZ} and \cite{10PZ} we consider the following $(\infty,1)$-categories of $(\infty,1)$-modules.\\
\end{definition}

\begin{definition}
First we consider the Bambozzi-Kremnizer spectrum $\mathrm{Spec}^\mathrm{BK}(*)$ attached to any of those in the above from \cite{10BK} by taking derived rational localization:
\begin{align}
&\mathrm{Spec}^\mathrm{BK}\widetilde{\Phi}_{\psi,A},\mathrm{Spec}^\mathrm{BK}\widetilde{\Phi}^r_{\psi,A},\mathrm{Spec}^\mathrm{BK}\widetilde{\Phi}^I_{\psi,A},	
\end{align}
\begin{align}
&\mathrm{Spec}^\mathrm{BK}\breve{\Phi}_{\psi,A},\mathrm{Spec}^\mathrm{BK}\breve{\Phi}^r_{\psi,A},\mathrm{Spec}^\mathrm{BK}\breve{\Phi}^I_{\psi,A},	
\end{align}
\begin{align}
&\mathrm{Spec}^\mathrm{BK}{\Phi}_{\psi,A},
\mathrm{Spec}^\mathrm{BK}{\Phi}^r_{\psi,A},\mathrm{Spec}^\mathrm{BK}{\Phi}^I_{\psi,A}.	
\end{align}

Then we take the corresponding quotients by using the corresponding Frobenius operators:
\begin{align}
&\mathrm{Spec}^\mathrm{BK}\widetilde{\Phi}_{\psi,A}/\mathrm{Fro}^\mathbb{Z},	\\
\end{align}
\begin{align}
&\mathrm{Spec}^\mathrm{BK}\breve{\Phi}_{\psi,A}/\mathrm{Fro}^\mathbb{Z},	\\
\end{align}
\begin{align}
&\mathrm{Spec}^\mathrm{BK}{\Phi}_{\psi,A}/\mathrm{Fro}^\mathbb{Z}.	
\end{align}
Here for those space without notation related to the radius and the corresponding interval we consider the total unions $\bigcap_r,\bigcup_I$ in order to achieve the whole spaces to achieve the analogues of the corresponding FF curves from \cite{10KL1}, \cite{10KL2}, \cite{10FF} for
\[
\xymatrix@R+0pc@C+0pc{
\underset{r}{\mathrm{homotopylimit}}~\mathrm{Spec}^\mathrm{BK}\widetilde{\Phi}^r_{\psi,A},\underset{I}{\mathrm{homotopycolimit}}~\mathrm{Spec}^\mathrm{BK}\widetilde{\Phi}^I_{\psi,A},	\\
}
\]
\[
\xymatrix@R+0pc@C+0pc{
\underset{r}{\mathrm{homotopylimit}}~\mathrm{Spec}^\mathrm{BK}\breve{\Phi}^r_{\psi,A},\underset{I}{\mathrm{homotopycolimit}}~\mathrm{Spec}^\mathrm{BK}\breve{\Phi}^I_{\psi,A},	\\
}
\]
\[
\xymatrix@R+0pc@C+0pc{
\underset{r}{\mathrm{homotopylimit}}~\mathrm{Spec}^\mathrm{BK}{\Phi}^r_{\psi,A},\underset{I}{\mathrm{homotopycolimit}}~\mathrm{Spec}^\mathrm{BK}{\Phi}^I_{\psi,A}.	
}
\]
\[  
\xymatrix@R+0pc@C+0pc{
\underset{r}{\mathrm{homotopylimit}}~\mathrm{Spec}^\mathrm{BK}\widetilde{\Phi}^r_{\psi,A}/\mathrm{Fro}^\mathbb{Z},\underset{I}{\mathrm{homotopycolimit}}~\mathrm{Spec}^\mathrm{BK}\widetilde{\Phi}^I_{\psi,A}/\mathrm{Fro}^\mathbb{Z},	\\
}
\]
\[ 
\xymatrix@R+0pc@C+0pc{
\underset{r}{\mathrm{homotopylimit}}~\mathrm{Spec}^\mathrm{BK}\breve{\Phi}^r_{\psi,A}/\mathrm{Fro}^\mathbb{Z},\underset{I}{\mathrm{homotopycolimit}}~\mathrm{Spec}^\mathrm{BK}\breve{\Phi}^I_{\psi,A}/\mathrm{Fro}^\mathbb{Z},	\\
}
\]
\[ 
\xymatrix@R+0pc@C+0pc{
\underset{r}{\mathrm{homotopylimit}}~\mathrm{Spec}^\mathrm{BK}{\Phi}^r_{\psi,A}/\mathrm{Fro}^\mathbb{Z},\underset{I}{\mathrm{homotopycolimit}}~\mathrm{Spec}^\mathrm{BK}{\Phi}^I_{\psi,A}/\mathrm{Fro}^\mathbb{Z}.	
}
\]

\end{definition}

\indent Meanwhile we have the corresponding Clausen-Scholze analytic stacks from \cite{10CS2}, therefore applying their construction we have:

\begin{definition}
Here we define the following products by using the solidified tensor product from \cite{10CS1} and \cite{10CS2}. Namely $A$ will still as above as a Banach ring over $\mathbb{Q}_p$. Then we take solidified tensor product $\overset{\blacksquare}{\otimes}$ of any of the following
\[
\xymatrix@R+0pc@C+0pc{
\widetilde{\Delta}_{\psi},\widetilde{\nabla}_{\psi},\widetilde{\Phi}_{\psi},\widetilde{\Delta}^+_{\psi},\widetilde{\nabla}^+_{\psi},\widetilde{\Delta}^\dagger_{\psi},\widetilde{\nabla}^\dagger_{\psi},\widetilde{\Phi}^r_{\psi},\widetilde{\Phi}^I_{\psi}, 
}
\]

\[
\xymatrix@R+0pc@C+0pc{
\breve{\Delta}_{\psi},\breve{\nabla}_{\psi},\breve{\Phi}_{\psi},\breve{\Delta}^+_{\psi},\breve{\nabla}^+_{\psi},\breve{\Delta}^\dagger_{\psi},\breve{\nabla}^\dagger_{\psi},\breve{\Phi}^r_{\psi},\breve{\Phi}^I_{\psi},	
}
\]

\[
\xymatrix@R+0pc@C+0pc{
{\Delta}_{\psi},{\nabla}_{\psi},{\Phi}_{\psi},{\Delta}^+_{\psi},{\nabla}^+_{\psi},{\Delta}^\dagger_{\psi},{\nabla}^\dagger_{\psi},{\Phi}^r_{\psi},{\Phi}^I_{\psi},	
}
\]  	
with $A$. Then we have the notations:
\[
\xymatrix@R+0pc@C+0pc{
\widetilde{\Delta}_{\psi,A},\widetilde{\nabla}_{\psi,A},\widetilde{\Phi}_{\psi,A},\widetilde{\Delta}^+_{\psi,A},\widetilde{\nabla}^+_{\psi,A},\widetilde{\Delta}^\dagger_{\psi,A},\widetilde{\nabla}^\dagger_{\psi,A},\widetilde{\Phi}^r_{\psi,A},\widetilde{\Phi}^I_{\psi,A}, 
}
\]

\[
\xymatrix@R+0pc@C+0pc{
\breve{\Delta}_{\psi,A},\breve{\nabla}_{\psi,A},\breve{\Phi}_{\psi,A},\breve{\Delta}^+_{\psi,A},\breve{\nabla}^+_{\psi,A},\breve{\Delta}^\dagger_{\psi,A},\breve{\nabla}^\dagger_{\psi,A},\breve{\Phi}^r_{\psi,A},\breve{\Phi}^I_{\psi,A},	
}
\]

\[
\xymatrix@R+0pc@C+0pc{
{\Delta}_{\psi,A},{\nabla}_{\psi,A},{\Phi}_{\psi,A},{\Delta}^+_{\psi,A},{\nabla}^+_{\psi,A},{\Delta}^\dagger_{\psi,A},{\nabla}^\dagger_{\psi,A},{\Phi}^r_{\psi,A},{\Phi}^I_{\psi,A}.	
}
\]
\end{definition}

\begin{definition}
First we consider the Clausen-Scholze spectrum $\mathrm{Spec}^\mathrm{CS}(*)$ attached to any of those in the above from \cite{10CS2} by taking derived rational localization:
\begin{align}
\mathrm{Spec}^\mathrm{CS}\widetilde{\Delta}_{\psi,A},\mathrm{Spec}^\mathrm{CS}\widetilde{\nabla}_{\psi,A},\mathrm{Spec}^\mathrm{CS}\widetilde{\Phi}_{\psi,A},\mathrm{Spec}^\mathrm{CS}\widetilde{\Delta}^+_{\psi,A},\mathrm{Spec}^\mathrm{CS}\widetilde{\nabla}^+_{\psi,A},\\
\mathrm{Spec}^\mathrm{CS}\widetilde{\Delta}^\dagger_{\psi,A},\mathrm{Spec}^\mathrm{CS}\widetilde{\nabla}^\dagger_{\psi,A},\mathrm{Spec}^\mathrm{CS}\widetilde{\Phi}^r_{\psi,A},\mathrm{Spec}^\mathrm{CS}\widetilde{\Phi}^I_{\psi,A},	\\
\end{align}
\begin{align}
\mathrm{Spec}^\mathrm{CS}\breve{\Delta}_{\psi,A},\breve{\nabla}_{\psi,A},\mathrm{Spec}^\mathrm{CS}\breve{\Phi}_{\psi,A},\mathrm{Spec}^\mathrm{CS}\breve{\Delta}^+_{\psi,A},\mathrm{Spec}^\mathrm{CS}\breve{\nabla}^+_{\psi,A},\\
\mathrm{Spec}^\mathrm{CS}\breve{\Delta}^\dagger_{\psi,A},\mathrm{Spec}^\mathrm{CS}\breve{\nabla}^\dagger_{\psi,A},\mathrm{Spec}^\mathrm{CS}\breve{\Phi}^r_{\psi,A},\breve{\Phi}^I_{\psi,A},	\\
\end{align}
\begin{align}
\mathrm{Spec}^\mathrm{CS}{\Delta}_{\psi,A},\mathrm{Spec}^\mathrm{CS}{\nabla}_{\psi,A},\mathrm{Spec}^\mathrm{CS}{\Phi}_{\psi,A},\mathrm{Spec}^\mathrm{CS}{\Delta}^+_{\psi,A},\mathrm{Spec}^\mathrm{CS}{\nabla}^+_{\psi,A},\\
\mathrm{Spec}^\mathrm{CS}{\Delta}^\dagger_{\psi,A},\mathrm{Spec}^\mathrm{CS}{\nabla}^\dagger_{\psi,A},\mathrm{Spec}^\mathrm{CS}{\Phi}^r_{\psi,A},\mathrm{Spec}^\mathrm{CS}{\Phi}^I_{\psi,A}.	
\end{align}

Then we take the corresponding quotients by using the corresponding Frobenius operators:
\begin{align}
&\mathrm{Spec}^\mathrm{CS}\widetilde{\Delta}_{\psi,A}/\mathrm{Fro}^\mathbb{Z},\mathrm{Spec}^\mathrm{CS}\widetilde{\nabla}_{\psi,A}/\mathrm{Fro}^\mathbb{Z},\mathrm{Spec}^\mathrm{CS}\widetilde{\Phi}_{\psi,A}/\mathrm{Fro}^\mathbb{Z},\mathrm{Spec}^\mathrm{CS}\widetilde{\Delta}^+_{\psi,A}/\mathrm{Fro}^\mathbb{Z},\\
&\mathrm{Spec}^\mathrm{CS}\widetilde{\nabla}^+_{\psi,A}/\mathrm{Fro}^\mathbb{Z}, \mathrm{Spec}^\mathrm{CS}\widetilde{\Delta}^\dagger_{\psi,A}/\mathrm{Fro}^\mathbb{Z},\mathrm{Spec}^\mathrm{CS}\widetilde{\nabla}^\dagger_{\psi,A}/\mathrm{Fro}^\mathbb{Z},	\\
\end{align}
\begin{align}
&\mathrm{Spec}^\mathrm{CS}\breve{\Delta}_{\psi,A}/\mathrm{Fro}^\mathbb{Z},\breve{\nabla}_{\psi,A}/\mathrm{Fro}^\mathbb{Z},\mathrm{Spec}^\mathrm{CS}\breve{\Phi}_{\psi,A}/\mathrm{Fro}^\mathbb{Z},\mathrm{Spec}^\mathrm{CS}\breve{\Delta}^+_{\psi,A}/\mathrm{Fro}^\mathbb{Z},\\
&\mathrm{Spec}^\mathrm{CS}\breve{\nabla}^+_{\psi,A}/\mathrm{Fro}^\mathbb{Z}, \mathrm{Spec}^\mathrm{CS}\breve{\Delta}^\dagger_{\psi,A}/\mathrm{Fro}^\mathbb{Z},\mathrm{Spec}^\mathrm{CS}\breve{\nabla}^\dagger_{\psi,A}/\mathrm{Fro}^\mathbb{Z},	\\
\end{align}
\begin{align}
&\mathrm{Spec}^\mathrm{CS}{\Delta}_{\psi,A}/\mathrm{Fro}^\mathbb{Z},\mathrm{Spec}^\mathrm{CS}{\nabla}_{\psi,A}/\mathrm{Fro}^\mathbb{Z},\mathrm{Spec}^\mathrm{CS}{\Phi}_{\psi,A}/\mathrm{Fro}^\mathbb{Z},\mathrm{Spec}^\mathrm{CS}{\Delta}^+_{\psi,A}/\mathrm{Fro}^\mathbb{Z},\\
&\mathrm{Spec}^\mathrm{CS}{\nabla}^+_{\psi,A}/\mathrm{Fro}^\mathbb{Z}, \mathrm{Spec}^\mathrm{CS}{\Delta}^\dagger_{\psi,A}/\mathrm{Fro}^\mathbb{Z},\mathrm{Spec}^\mathrm{CS}{\nabla}^\dagger_{\psi,A}/\mathrm{Fro}^\mathbb{Z}.	
\end{align}
Here for those space with notations related to the radius and the corresponding interval we consider the total unions $\bigcap_r,\bigcup_I$ in order to achieve the whole spaces to achieve the analogues of the corresponding FF curves from \cite{10KL1}, \cite{10KL2}, \cite{10FF} for
\[
\xymatrix@R+0pc@C+0pc{
\underset{r}{\mathrm{homotopylimit}}~\mathrm{Spec}^\mathrm{CS}\widetilde{\Phi}^r_{\psi,A},\underset{I}{\mathrm{homotopycolimit}}~\mathrm{Spec}^\mathrm{CS}\widetilde{\Phi}^I_{\psi,A},	\\
}
\]
\[
\xymatrix@R+0pc@C+0pc{
\underset{r}{\mathrm{homotopylimit}}~\mathrm{Spec}^\mathrm{CS}\breve{\Phi}^r_{\psi,A},\underset{I}{\mathrm{homotopycolimit}}~\mathrm{Spec}^\mathrm{CS}\breve{\Phi}^I_{\psi,A},	\\
}
\]
\[
\xymatrix@R+0pc@C+0pc{
\underset{r}{\mathrm{homotopylimit}}~\mathrm{Spec}^\mathrm{CS}{\Phi}^r_{\psi,A},\underset{I}{\mathrm{homotopycolimit}}~\mathrm{Spec}^\mathrm{CS}{\Phi}^I_{\psi,A}.	
}
\]
\[ 
\xymatrix@R+0pc@C+0pc{
\underset{r}{\mathrm{homotopylimit}}~\mathrm{Spec}^\mathrm{CS}\widetilde{\Phi}^r_{\psi,A}/\mathrm{Fro}^\mathbb{Z},\underset{I}{\mathrm{homotopycolimit}}~\mathrm{Spec}^\mathrm{CS}\widetilde{\Phi}^I_{\psi,A}/\mathrm{Fro}^\mathbb{Z},	\\
}
\]
\[ 
\xymatrix@R+0pc@C+0pc{
\underset{r}{\mathrm{homotopylimit}}~\mathrm{Spec}^\mathrm{CS}\breve{\Phi}^r_{\psi,A}/\mathrm{Fro}^\mathbb{Z},\underset{I}{\mathrm{homotopycolimit}}~\breve{\Phi}^I_{\psi,A}/\mathrm{Fro}^\mathbb{Z},	\\
}
\]
\[ 
\xymatrix@R+0pc@C+0pc{
\underset{r}{\mathrm{homotopylimit}}~\mathrm{Spec}^\mathrm{CS}{\Phi}^r_{\psi,A}/\mathrm{Fro}^\mathbb{Z},\underset{I}{\mathrm{homotopycolimit}}~\mathrm{Spec}^\mathrm{CS}{\Phi}^I_{\psi,A}/\mathrm{Fro}^\mathbb{Z}.	
}
\]

\end{definition}

\

\begin{definition}
We then consider the corresponding quasipresheaves of the corresponding ind-Banach or monomorphic ind-Banach modules from \cite{10BBK}, \cite{10KKM}:
\begin{align}
\mathrm{Quasicoherentpresheaves,IndBanach}_{*}	
\end{align}
where $*$ is one of the following spaces:
\begin{align}
&\mathrm{Spec}^\mathrm{BK}\widetilde{\Phi}_{\psi,A}/\mathrm{Fro}^\mathbb{Z},	\\
\end{align}
\begin{align}
&\mathrm{Spec}^\mathrm{BK}\breve{\Phi}_{\psi,A}/\mathrm{Fro}^\mathbb{Z},	\\
\end{align}
\begin{align}
&\mathrm{Spec}^\mathrm{BK}{\Phi}_{\psi,A}/\mathrm{Fro}^\mathbb{Z}.	
\end{align}
Here for those space without notation related to the radius and the corresponding interval we consider the total unions $\bigcap_r,\bigcup_I$ in order to achieve the whole spaces to achieve the analogues of the corresponding FF curves from \cite{10KL1}, \cite{10KL2}, \cite{10FF} for
\[
\xymatrix@R+0pc@C+0pc{
\underset{r}{\mathrm{homotopylimit}}~\mathrm{Spec}^\mathrm{BK}\widetilde{\Phi}^r_{\psi,A},\underset{I}{\mathrm{homotopycolimit}}~\mathrm{Spec}^\mathrm{BK}\widetilde{\Phi}^I_{\psi,A},	\\
}
\]
\[
\xymatrix@R+0pc@C+0pc{
\underset{r}{\mathrm{homotopylimit}}~\mathrm{Spec}^\mathrm{BK}\breve{\Phi}^r_{\psi,A},\underset{I}{\mathrm{homotopycolimit}}~\mathrm{Spec}^\mathrm{BK}\breve{\Phi}^I_{\psi,A},	\\
}
\]
\[
\xymatrix@R+0pc@C+0pc{
\underset{r}{\mathrm{homotopylimit}}~\mathrm{Spec}^\mathrm{BK}{\Phi}^r_{\psi,A},\underset{I}{\mathrm{homotopycolimit}}~\mathrm{Spec}^\mathrm{BK}{\Phi}^I_{\psi,A}.	
}
\]
\[  
\xymatrix@R+0pc@C+0pc{
\underset{r}{\mathrm{homotopylimit}}~\mathrm{Spec}^\mathrm{BK}\widetilde{\Phi}^r_{\psi,A}/\mathrm{Fro}^\mathbb{Z},\underset{I}{\mathrm{homotopycolimit}}~\mathrm{Spec}^\mathrm{BK}\widetilde{\Phi}^I_{\psi,A}/\mathrm{Fro}^\mathbb{Z},	\\
}
\]
\[ 
\xymatrix@R+0pc@C+0pc{
\underset{r}{\mathrm{homotopylimit}}~\mathrm{Spec}^\mathrm{BK}\breve{\Phi}^r_{\psi,A}/\mathrm{Fro}^\mathbb{Z},\underset{I}{\mathrm{homotopycolimit}}~\mathrm{Spec}^\mathrm{BK}\breve{\Phi}^I_{\psi,A}/\mathrm{Fro}^\mathbb{Z},	\\
}
\]
\[ 
\xymatrix@R+0pc@C+0pc{
\underset{r}{\mathrm{homotopylimit}}~\mathrm{Spec}^\mathrm{BK}{\Phi}^r_{\psi,A}/\mathrm{Fro}^\mathbb{Z},\underset{I}{\mathrm{homotopycolimit}}~\mathrm{Spec}^\mathrm{BK}{\Phi}^I_{\psi,A}/\mathrm{Fro}^\mathbb{Z}.	
}
\]

\end{definition}

\begin{definition}
We then consider the corresponding quasisheaves of the corresponding condensed solid topological modules from \cite{10CS2}:
\begin{align}
\mathrm{Quasicoherentsheaves, Condensed}_{*}	
\end{align}
where $*$ is one of the following spaces:
\begin{align}
&\mathrm{Spec}^\mathrm{CS}\widetilde{\Delta}_{\psi,A}/\mathrm{Fro}^\mathbb{Z},\mathrm{Spec}^\mathrm{CS}\widetilde{\nabla}_{\psi,A}/\mathrm{Fro}^\mathbb{Z},\mathrm{Spec}^\mathrm{CS}\widetilde{\Phi}_{\psi,A}/\mathrm{Fro}^\mathbb{Z},\mathrm{Spec}^\mathrm{CS}\widetilde{\Delta}^+_{\psi,A}/\mathrm{Fro}^\mathbb{Z},\\
&\mathrm{Spec}^\mathrm{CS}\widetilde{\nabla}^+_{\psi,A}/\mathrm{Fro}^\mathbb{Z},\mathrm{Spec}^\mathrm{CS}\widetilde{\Delta}^\dagger_{\psi,A}/\mathrm{Fro}^\mathbb{Z},\mathrm{Spec}^\mathrm{CS}\widetilde{\nabla}^\dagger_{\psi,A}/\mathrm{Fro}^\mathbb{Z},	\\
\end{align}
\begin{align}
&\mathrm{Spec}^\mathrm{CS}\breve{\Delta}_{\psi,A}/\mathrm{Fro}^\mathbb{Z},\breve{\nabla}_{\psi,A}/\mathrm{Fro}^\mathbb{Z},\mathrm{Spec}^\mathrm{CS}\breve{\Phi}_{\psi,A}/\mathrm{Fro}^\mathbb{Z},\mathrm{Spec}^\mathrm{CS}\breve{\Delta}^+_{\psi,A}/\mathrm{Fro}^\mathbb{Z},\\
&\mathrm{Spec}^\mathrm{CS}\breve{\nabla}^+_{\psi,A}/\mathrm{Fro}^\mathbb{Z},\mathrm{Spec}^\mathrm{CS}\breve{\Delta}^\dagger_{\psi,A}/\mathrm{Fro}^\mathbb{Z},\mathrm{Spec}^\mathrm{CS}\breve{\nabla}^\dagger_{\psi,A}/\mathrm{Fro}^\mathbb{Z},	\\
\end{align}
\begin{align}
&\mathrm{Spec}^\mathrm{CS}{\Delta}_{\psi,A}/\mathrm{Fro}^\mathbb{Z},\mathrm{Spec}^\mathrm{CS}{\nabla}_{\psi,A}/\mathrm{Fro}^\mathbb{Z},\mathrm{Spec}^\mathrm{CS}{\Phi}_{\psi,A}/\mathrm{Fro}^\mathbb{Z},\mathrm{Spec}^\mathrm{CS}{\Delta}^+_{\psi,A}/\mathrm{Fro}^\mathbb{Z},\\
&\mathrm{Spec}^\mathrm{CS}{\nabla}^+_{\psi,A}/\mathrm{Fro}^\mathbb{Z}, \mathrm{Spec}^\mathrm{CS}{\Delta}^\dagger_{\psi,A}/\mathrm{Fro}^\mathbb{Z},\mathrm{Spec}^\mathrm{CS}{\nabla}^\dagger_{\psi,A}/\mathrm{Fro}^\mathbb{Z}.	
\end{align}
Here for those space with notations related to the radius and the corresponding interval we consider the total unions $\bigcap_r,\bigcup_I$ in order to achieve the whole spaces to achieve the analogues of the corresponding FF curves from \cite{10KL1}, \cite{10KL2}, \cite{10FF} for
\[
\xymatrix@R+0pc@C+0pc{
\underset{r}{\mathrm{homotopylimit}}~\mathrm{Spec}^\mathrm{CS}\widetilde{\Phi}^r_{\psi,A},\underset{I}{\mathrm{homotopycolimit}}~\mathrm{Spec}^\mathrm{CS}\widetilde{\Phi}^I_{\psi,A},	\\
}
\]
\[
\xymatrix@R+0pc@C+0pc{
\underset{r}{\mathrm{homotopylimit}}~\mathrm{Spec}^\mathrm{CS}\breve{\Phi}^r_{\psi,A},\underset{I}{\mathrm{homotopycolimit}}~\mathrm{Spec}^\mathrm{CS}\breve{\Phi}^I_{\psi,A},	\\
}
\]
\[
\xymatrix@R+0pc@C+0pc{
\underset{r}{\mathrm{homotopylimit}}~\mathrm{Spec}^\mathrm{CS}{\Phi}^r_{\psi,A},\underset{I}{\mathrm{homotopycolimit}}~\mathrm{Spec}^\mathrm{CS}{\Phi}^I_{\psi,A}.	
}
\]
\[ 
\xymatrix@R+0pc@C+0pc{
\underset{r}{\mathrm{homotopylimit}}~\mathrm{Spec}^\mathrm{CS}\widetilde{\Phi}^r_{\psi,A}/\mathrm{Fro}^\mathbb{Z},\underset{I}{\mathrm{homotopycolimit}}~\mathrm{Spec}^\mathrm{CS}\widetilde{\Phi}^I_{\psi,A}/\mathrm{Fro}^\mathbb{Z},	\\
}
\]
\[ 
\xymatrix@R+0pc@C+0pc{
\underset{r}{\mathrm{homotopylimit}}~\mathrm{Spec}^\mathrm{CS}\breve{\Phi}^r_{\psi,A}/\mathrm{Fro}^\mathbb{Z},\underset{I}{\mathrm{homotopycolimit}}~\breve{\Phi}^I_{\psi,A}/\mathrm{Fro}^\mathbb{Z},	\\
}
\]
\[ 
\xymatrix@R+0pc@C+0pc{
\underset{r}{\mathrm{homotopylimit}}~\mathrm{Spec}^\mathrm{CS}{\Phi}^r_{\psi,A}/\mathrm{Fro}^\mathbb{Z},\underset{I}{\mathrm{homotopycolimit}}~\mathrm{Spec}^\mathrm{CS}{\Phi}^I_{\psi,A}/\mathrm{Fro}^\mathbb{Z}.	
}
\]

\end{definition}

\

\begin{proposition}
There is a well-defined functor from the $\infty$-category 
\begin{align}
\mathrm{Quasicoherentpresheaves,Condensed}_{*}	
\end{align}
where $*$ is one of the following spaces:
\begin{align}
&\mathrm{Spec}^\mathrm{CS}\widetilde{\Phi}_{\psi,A}/\mathrm{Fro}^\mathbb{Z},	\\
\end{align}
\begin{align}
&\mathrm{Spec}^\mathrm{CS}\breve{\Phi}_{\psi,A}/\mathrm{Fro}^\mathbb{Z},	\\
\end{align}
\begin{align}
&\mathrm{Spec}^\mathrm{CS}{\Phi}_{\psi,A}/\mathrm{Fro}^\mathbb{Z},	
\end{align}
to the $\infty$-category of $\mathrm{Fro}$-equivariant quasicoherent presheaves over similar spaces above correspondingly without the $\mathrm{Fro}$-quotients, and to the $\infty$-category of $\mathrm{Fro}$-equivariant quasicoherent modules over global sections of the structure $\infty$-sheaves of the similar spaces above correspondingly without the $\mathrm{Fro}$-quotients. Here for those space without notation related to the radius and the corresponding interval we consider the total unions $\bigcap_r,\bigcup_I$ in order to achieve the whole spaces to achieve the analogues of the corresponding FF curves from \cite{10KL1}, \cite{10KL2}, \cite{10FF} for
\[
\xymatrix@R+0pc@C+0pc{
\underset{r}{\mathrm{homotopylimit}}~\mathrm{Spec}^\mathrm{CS}\widetilde{\Phi}^r_{\psi,A},\underset{I}{\mathrm{homotopycolimit}}~\mathrm{Spec}^\mathrm{CS}\widetilde{\Phi}^I_{\psi,A},	\\
}
\]
\[
\xymatrix@R+0pc@C+0pc{
\underset{r}{\mathrm{homotopylimit}}~\mathrm{Spec}^\mathrm{CS}\breve{\Phi}^r_{\psi,A},\underset{I}{\mathrm{homotopycolimit}}~\mathrm{Spec}^\mathrm{CS}\breve{\Phi}^I_{\psi,A},	\\
}
\]
\[
\xymatrix@R+0pc@C+0pc{
\underset{r}{\mathrm{homotopylimit}}~\mathrm{Spec}^\mathrm{CS}{\Phi}^r_{\psi,A},\underset{I}{\mathrm{homotopycolimit}}~\mathrm{Spec}^\mathrm{CS}{\Phi}^I_{\psi,A}.	
}
\]
\[ 
\xymatrix@R+0pc@C+0pc{
\underset{r}{\mathrm{homotopylimit}}~\mathrm{Spec}^\mathrm{CS}\widetilde{\Phi}^r_{\psi,A}/\mathrm{Fro}^\mathbb{Z},\underset{I}{\mathrm{homotopycolimit}}~\mathrm{Spec}^\mathrm{CS}\widetilde{\Phi}^I_{\psi,A}/\mathrm{Fro}^\mathbb{Z},	\\
}
\]
\[ 
\xymatrix@R+0pc@C+0pc{
\underset{r}{\mathrm{homotopylimit}}~\mathrm{Spec}^\mathrm{CS}\breve{\Phi}^r_{\psi,A}/\mathrm{Fro}^\mathbb{Z},\underset{I}{\mathrm{homotopycolimit}}~\breve{\Phi}^I_{\psi,A}/\mathrm{Fro}^\mathbb{Z},	\\
}
\]
\[ 
\xymatrix@R+0pc@C+0pc{
\underset{r}{\mathrm{homotopylimit}}~\mathrm{Spec}^\mathrm{CS}{\Phi}^r_{\psi,A}/\mathrm{Fro}^\mathbb{Z},\underset{I}{\mathrm{homotopycolimit}}~\mathrm{Spec}^\mathrm{CS}{\Phi}^I_{\psi,A}/\mathrm{Fro}^\mathbb{Z}.	
}
\]	
In this situation we will have the target category being family parametrized by $r$ or $I$ in compatible glueing sense as in \cite[Definition 5.4.10]{10KL2}. In this situation for modules parametrized by the intervals we have the equivalence of $\infty$-categories by using \cite[Proposition 13.8]{10CS2}. Here the corresponding quasicoherent Frobenius modules are defined to be the corresponding homotopy colimits and limits of Frobenius modules:
\begin{align}
\underset{r}{\mathrm{homotopycolimit}}~M_r,\\
\underset{I}{\mathrm{homotopylimit}}~M_I,	
\end{align}
where each $M_r$ is a Frobenius-equivariant module over the period ring with respect to some radius $r$ while each $M_I$ is a Frobenius-equivariant module over the period ring with respect to some interval $I$.\\
\end{proposition}

\begin{proposition}
Similar proposition holds for 
\begin{align}
\mathrm{Quasicoherentsheaves,IndBanach}_{*}.	
\end{align}	
\end{proposition}

\

\begin{definition}
We then consider the corresponding quasipresheaves of perfect complexes the corresponding ind-Banach or monomorphic ind-Banach modules from \cite{10BBK}, \cite{10KKM}:
\begin{align}
\mathrm{Quasicoherentpresheaves,Perfectcomplex,IndBanach}_{*}	
\end{align}
where $*$ is one of the following spaces:
\begin{align}
&\mathrm{Spec}^\mathrm{BK}\widetilde{\Phi}_{\psi,A}/\mathrm{Fro}^\mathbb{Z},	\\
\end{align}
\begin{align}
&\mathrm{Spec}^\mathrm{BK}\breve{\Phi}_{\psi,A}/\mathrm{Fro}^\mathbb{Z},	\\
\end{align}
\begin{align}
&\mathrm{Spec}^\mathrm{BK}{\Phi}_{\psi,A}/\mathrm{Fro}^\mathbb{Z}.	
\end{align}
Here for those space without notation related to the radius and the corresponding interval we consider the total unions $\bigcap_r,\bigcup_I$ in order to achieve the whole spaces to achieve the analogues of the corresponding FF curves from \cite{10KL1}, \cite{10KL2}, \cite{10FF} for
\[
\xymatrix@R+0pc@C+0pc{
\underset{r}{\mathrm{homotopylimit}}~\mathrm{Spec}^\mathrm{BK}\widetilde{\Phi}^r_{\psi,A},\underset{I}{\mathrm{homotopycolimit}}~\mathrm{Spec}^\mathrm{BK}\widetilde{\Phi}^I_{\psi,A},	\\
}
\]
\[
\xymatrix@R+0pc@C+0pc{
\underset{r}{\mathrm{homotopylimit}}~\mathrm{Spec}^\mathrm{BK}\breve{\Phi}^r_{\psi,A},\underset{I}{\mathrm{homotopycolimit}}~\mathrm{Spec}^\mathrm{BK}\breve{\Phi}^I_{\psi,A},	\\
}
\]
\[
\xymatrix@R+0pc@C+0pc{
\underset{r}{\mathrm{homotopylimit}}~\mathrm{Spec}^\mathrm{BK}{\Phi}^r_{\psi,A},\underset{I}{\mathrm{homotopycolimit}}~\mathrm{Spec}^\mathrm{BK}{\Phi}^I_{\psi,A}.	
}
\]
\[  
\xymatrix@R+0pc@C+0pc{
\underset{r}{\mathrm{homotopylimit}}~\mathrm{Spec}^\mathrm{BK}\widetilde{\Phi}^r_{\psi,A}/\mathrm{Fro}^\mathbb{Z},\underset{I}{\mathrm{homotopycolimit}}~\mathrm{Spec}^\mathrm{BK}\widetilde{\Phi}^I_{\psi,A}/\mathrm{Fro}^\mathbb{Z},	\\
}
\]
\[ 
\xymatrix@R+0pc@C+0pc{
\underset{r}{\mathrm{homotopylimit}}~\mathrm{Spec}^\mathrm{BK}\breve{\Phi}^r_{\psi,A}/\mathrm{Fro}^\mathbb{Z},\underset{I}{\mathrm{homotopycolimit}}~\mathrm{Spec}^\mathrm{BK}\breve{\Phi}^I_{\psi,A}/\mathrm{Fro}^\mathbb{Z},	\\
}
\]
\[ 
\xymatrix@R+0pc@C+0pc{
\underset{r}{\mathrm{homotopylimit}}~\mathrm{Spec}^\mathrm{BK}{\Phi}^r_{\psi,A}/\mathrm{Fro}^\mathbb{Z},\underset{I}{\mathrm{homotopycolimit}}~\mathrm{Spec}^\mathrm{BK}{\Phi}^I_{\psi,A}/\mathrm{Fro}^\mathbb{Z}.	
}
\]

\end{definition}

\begin{definition}
We then consider the corresponding quasisheaves of perfect complexes of the corresponding condensed solid topological modules from \cite{10CS2}:
\begin{align}
\mathrm{Quasicoherentsheaves, Perfectcomplex, Condensed}_{*}	
\end{align}
where $*$ is one of the following spaces:
\begin{align}
&\mathrm{Spec}^\mathrm{CS}\widetilde{\Delta}_{\psi,A}/\mathrm{Fro}^\mathbb{Z},\mathrm{Spec}^\mathrm{CS}\widetilde{\nabla}_{\psi,A}/\mathrm{Fro}^\mathbb{Z},\mathrm{Spec}^\mathrm{CS}\widetilde{\Phi}_{\psi,A}/\mathrm{Fro}^\mathbb{Z},\mathrm{Spec}^\mathrm{CS}\widetilde{\Delta}^+_{\psi,A}/\mathrm{Fro}^\mathbb{Z},\\
&\mathrm{Spec}^\mathrm{CS}\widetilde{\nabla}^+_{\psi,A}/\mathrm{Fro}^\mathbb{Z},\mathrm{Spec}^\mathrm{CS}\widetilde{\Delta}^\dagger_{\psi,A}/\mathrm{Fro}^\mathbb{Z},\mathrm{Spec}^\mathrm{CS}\widetilde{\nabla}^\dagger_{\psi,A}/\mathrm{Fro}^\mathbb{Z},	\\
\end{align}
\begin{align}
&\mathrm{Spec}^\mathrm{CS}\breve{\Delta}_{\psi,A}/\mathrm{Fro}^\mathbb{Z},\breve{\nabla}_{\psi,A}/\mathrm{Fro}^\mathbb{Z},\mathrm{Spec}^\mathrm{CS}\breve{\Phi}_{\psi,A}/\mathrm{Fro}^\mathbb{Z},\mathrm{Spec}^\mathrm{CS}\breve{\Delta}^+_{\psi,A}/\mathrm{Fro}^\mathbb{Z},\\
&\mathrm{Spec}^\mathrm{CS}\breve{\nabla}^+_{\psi,A}/\mathrm{Fro}^\mathbb{Z},\mathrm{Spec}^\mathrm{CS}\breve{\Delta}^\dagger_{\psi,A}/\mathrm{Fro}^\mathbb{Z},\mathrm{Spec}^\mathrm{CS}\breve{\nabla}^\dagger_{\psi,A}/\mathrm{Fro}^\mathbb{Z},	\\
\end{align}
\begin{align}
&\mathrm{Spec}^\mathrm{CS}{\Delta}_{\psi,A}/\mathrm{Fro}^\mathbb{Z},\mathrm{Spec}^\mathrm{CS}{\nabla}_{\psi,A}/\mathrm{Fro}^\mathbb{Z},\mathrm{Spec}^\mathrm{CS}{\Phi}_{\psi,A}/\mathrm{Fro}^\mathbb{Z},\mathrm{Spec}^\mathrm{CS}{\Delta}^+_{\psi,A}/\mathrm{Fro}^\mathbb{Z},\\
&\mathrm{Spec}^\mathrm{CS}{\nabla}^+_{\psi,A}/\mathrm{Fro}^\mathbb{Z}, \mathrm{Spec}^\mathrm{CS}{\Delta}^\dagger_{\psi,A}/\mathrm{Fro}^\mathbb{Z},\mathrm{Spec}^\mathrm{CS}{\nabla}^\dagger_{\psi,A}/\mathrm{Fro}^\mathbb{Z}.	
\end{align}
Here for those space with notations related to the radius and the corresponding interval we consider the total unions $\bigcap_r,\bigcup_I$ in order to achieve the whole spaces to achieve the analogues of the corresponding FF curves from \cite{10KL1}, \cite{10KL2}, \cite{10FF} for
\[
\xymatrix@R+0pc@C+0pc{
\underset{r}{\mathrm{homotopylimit}}~\mathrm{Spec}^\mathrm{CS}\widetilde{\Phi}^r_{\psi,A},\underset{I}{\mathrm{homotopycolimit}}~\mathrm{Spec}^\mathrm{CS}\widetilde{\Phi}^I_{\psi,A},	\\
}
\]
\[
\xymatrix@R+0pc@C+0pc{
\underset{r}{\mathrm{homotopylimit}}~\mathrm{Spec}^\mathrm{CS}\breve{\Phi}^r_{\psi,A},\underset{I}{\mathrm{homotopycolimit}}~\mathrm{Spec}^\mathrm{CS}\breve{\Phi}^I_{\psi,A},	\\
}
\]
\[
\xymatrix@R+0pc@C+0pc{
\underset{r}{\mathrm{homotopylimit}}~\mathrm{Spec}^\mathrm{CS}{\Phi}^r_{\psi,A},\underset{I}{\mathrm{homotopycolimit}}~\mathrm{Spec}^\mathrm{CS}{\Phi}^I_{\psi,A}.	
}
\]
\[ 
\xymatrix@R+0pc@C+0pc{
\underset{r}{\mathrm{homotopylimit}}~\mathrm{Spec}^\mathrm{CS}\widetilde{\Phi}^r_{\psi,A}/\mathrm{Fro}^\mathbb{Z},\underset{I}{\mathrm{homotopycolimit}}~\mathrm{Spec}^\mathrm{CS}\widetilde{\Phi}^I_{\psi,A}/\mathrm{Fro}^\mathbb{Z},	\\
}
\]
\[ 
\xymatrix@R+0pc@C+0pc{
\underset{r}{\mathrm{homotopylimit}}~\mathrm{Spec}^\mathrm{CS}\breve{\Phi}^r_{\psi,A}/\mathrm{Fro}^\mathbb{Z},\underset{I}{\mathrm{homotopycolimit}}~\breve{\Phi}^I_{\psi,A}/\mathrm{Fro}^\mathbb{Z},	\\
}
\]
\[ 
\xymatrix@R+0pc@C+0pc{
\underset{r}{\mathrm{homotopylimit}}~\mathrm{Spec}^\mathrm{CS}{\Phi}^r_{\psi,A}/\mathrm{Fro}^\mathbb{Z},\underset{I}{\mathrm{homotopycolimit}}~\mathrm{Spec}^\mathrm{CS}{\Phi}^I_{\psi,A}/\mathrm{Fro}^\mathbb{Z}.	
}
\]

\end{definition}

\begin{proposition}
There is a well-defined functor from the $\infty$-category 
\begin{align}
\mathrm{Quasicoherentpresheaves,Perfectcomplex,Condensed}_{*}	
\end{align}
where $*$ is one of the following spaces:
\begin{align}
&\mathrm{Spec}^\mathrm{CS}\widetilde{\Phi}_{\psi,A}/\mathrm{Fro}^\mathbb{Z},	\\
\end{align}
\begin{align}
&\mathrm{Spec}^\mathrm{CS}\breve{\Phi}_{\psi,A}/\mathrm{Fro}^\mathbb{Z},	\\
\end{align}
\begin{align}
&\mathrm{Spec}^\mathrm{CS}{\Phi}_{\psi,A}/\mathrm{Fro}^\mathbb{Z},	
\end{align}
to the $\infty$-category of $\mathrm{Fro}$-equivariant quasicoherent presheaves over similar spaces above correspondingly without the $\mathrm{Fro}$-quotients, and to the $\infty$-category of $\mathrm{Fro}$-equivariant quasicoherent modules over global sections of the structure $\infty$-sheaves of the similar spaces above correspondingly without the $\mathrm{Fro}$-quotients. Here for those space without notation related to the radius and the corresponding interval we consider the total unions $\bigcap_r,\bigcup_I$ in order to achieve the whole spaces to achieve the analogues of the corresponding FF curves from \cite{10KL1}, \cite{10KL2}, \cite{10FF} for
\[
\xymatrix@R+0pc@C+0pc{
\underset{r}{\mathrm{homotopylimit}}~\mathrm{Spec}^\mathrm{CS}\widetilde{\Phi}^r_{\psi,A},\underset{I}{\mathrm{homotopycolimit}}~\mathrm{Spec}^\mathrm{CS}\widetilde{\Phi}^I_{\psi,A},	\\
}
\]
\[
\xymatrix@R+0pc@C+0pc{
\underset{r}{\mathrm{homotopylimit}}~\mathrm{Spec}^\mathrm{CS}\breve{\Phi}^r_{\psi,A},\underset{I}{\mathrm{homotopycolimit}}~\mathrm{Spec}^\mathrm{CS}\breve{\Phi}^I_{\psi,A},	\\
}
\]
\[
\xymatrix@R+0pc@C+0pc{
\underset{r}{\mathrm{homotopylimit}}~\mathrm{Spec}^\mathrm{CS}{\Phi}^r_{\psi,A},\underset{I}{\mathrm{homotopycolimit}}~\mathrm{Spec}^\mathrm{CS}{\Phi}^I_{\psi,A}.	
}
\]
\[ 
\xymatrix@R+0pc@C+0pc{
\underset{r}{\mathrm{homotopylimit}}~\mathrm{Spec}^\mathrm{CS}\widetilde{\Phi}^r_{\psi,A}/\mathrm{Fro}^\mathbb{Z},\underset{I}{\mathrm{homotopycolimit}}~\mathrm{Spec}^\mathrm{CS}\widetilde{\Phi}^I_{\psi,A}/\mathrm{Fro}^\mathbb{Z},	\\
}
\]
\[ 
\xymatrix@R+0pc@C+0pc{
\underset{r}{\mathrm{homotopylimit}}~\mathrm{Spec}^\mathrm{CS}\breve{\Phi}^r_{\psi,A}/\mathrm{Fro}^\mathbb{Z},\underset{I}{\mathrm{homotopycolimit}}~\breve{\Phi}^I_{\psi,A}/\mathrm{Fro}^\mathbb{Z},	\\
}
\]
\[ 
\xymatrix@R+0pc@C+0pc{
\underset{r}{\mathrm{homotopylimit}}~\mathrm{Spec}^\mathrm{CS}{\Phi}^r_{\psi,A}/\mathrm{Fro}^\mathbb{Z},\underset{I}{\mathrm{homotopycolimit}}~\mathrm{Spec}^\mathrm{CS}{\Phi}^I_{\psi,A}/\mathrm{Fro}^\mathbb{Z}.	
}
\]	
In this situation we will have the target category being family parametrized by $r$ or $I$ in compatible glueing sense as in \cite[Definition 5.4.10]{10KL2}. In this situation for modules parametrized by the intervals we have the equivalence of $\infty$-categories by using \cite[Proposition 12.18]{10CS2}. Here the corresponding quasicoherent Frobenius modules are defined to be the corresponding homotopy colimits and limits of Frobenius modules:
\begin{align}
\underset{r}{\mathrm{homotopycolimit}}~M_r,\\
\underset{I}{\mathrm{homotopylimit}}~M_I,	
\end{align}
where each $M_r$ is a Frobenius-equivariant module over the period ring with respect to some radius $r$ while each $M_I$ is a Frobenius-equivariant module over the period ring with respect to some interval $I$.\\
\end{proposition}

\begin{proposition}
Similar proposition holds for 
\begin{align}
\mathrm{Quasicoherentsheaves,Perfectcomplex,IndBanach}_{*}.	
\end{align}	
\end{proposition}

\newpage
\subsection{Frobenius Quasicoherent Modules II: Deformation in Banach Rings}

\begin{definition}
Let $\psi$ be a toric tower over $\mathbb{Q}_p$ as in \cite[Chapter 7]{10KL2} with base $\mathbb{Q}_p\left<X_1^{\pm 1},...,X_k^{\pm 1}\right>$. Then from \cite{10KL1} and \cite[Definition 5.2.1]{10KL2} we have the following class of Kedlaya-Liu rings (with the following replacement: $\Delta$ stands for $A$, $\nabla$ stands for $B$, while $\Phi$ stands for $C$) by taking product in the sense of self $\Gamma$-th power\footnote{Here $|\Gamma|=1$.}:

\[
\xymatrix@R+0pc@C+0pc{
\widetilde{\Delta}_{\psi},\widetilde{\nabla}_{\psi},\widetilde{\Phi}_{\psi},\widetilde{\Delta}^+_{\psi},\widetilde{\nabla}^+_{\psi},\widetilde{\Delta}^\dagger_{\psi},\widetilde{\nabla}^\dagger_{\psi},\widetilde{\Phi}^r_{\psi},\widetilde{\Phi}^I_{\psi}, 
}
\]

\[
\xymatrix@R+0pc@C+0pc{
\breve{\Delta}_{\psi},\breve{\nabla}_{\psi},\breve{\Phi}_{\psi},\breve{\Delta}^+_{\psi},\breve{\nabla}^+_{\psi},\breve{\Delta}^\dagger_{\psi},\breve{\nabla}^\dagger_{\psi},\breve{\Phi}^r_{\psi},\breve{\Phi}^I_{\psi},	
}
\]

\[
\xymatrix@R+0pc@C+0pc{
{\Delta}_{\psi},{\nabla}_{\psi},{\Phi}_{\psi},{\Delta}^+_{\psi},{\nabla}^+_{\psi},{\Delta}^\dagger_{\psi},{\nabla}^\dagger_{\psi},{\Phi}^r_{\psi},{\Phi}^I_{\psi}.	
}
\]
We now consider the following rings with $-$ being any deforming Banach ring over $\mathbb{Q}_p$. Taking the product we have:
\[
\xymatrix@R+0pc@C+0pc{
\widetilde{\Phi}_{\psi,-},\widetilde{\Phi}^r_{\psi,-},\widetilde{\Phi}^I_{\psi,-},	
}
\]
\[
\xymatrix@R+0pc@C+0pc{
\breve{\Phi}_{\psi,-},\breve{\Phi}^r_{\psi,-},\breve{\Phi}^I_{\psi,-},	
}
\]
\[
\xymatrix@R+0pc@C+0pc{
{\Phi}_{\psi,-},{\Phi}^r_{\psi,-},{\Phi}^I_{\psi,-}.	
}
\]
They carry multi Frobenius action $\varphi_\Gamma$ and multi $\mathrm{Lie}_\Gamma:=\mathbb{Z}_p^{\times\Gamma}$ action. In our current situation after \cite{10CKZ} and \cite{10PZ} we consider the following $(\infty,1)$-categories of $(\infty,1)$-modules.\\
\end{definition}

\begin{definition}
First we consider the Bambozzi-Kremnizer spectrum $\mathrm{Spec}^\mathrm{BK}(*)$ attached to any of those in the above from \cite{10BK} by taking derived rational localization:
\begin{align}
&\mathrm{Spec}^\mathrm{BK}\widetilde{\Phi}_{\psi,-},\mathrm{Spec}^\mathrm{BK}\widetilde{\Phi}^r_{\psi,-},\mathrm{Spec}^\mathrm{BK}\widetilde{\Phi}^I_{\psi,-},	
\end{align}
\begin{align}
&\mathrm{Spec}^\mathrm{BK}\breve{\Phi}_{\psi,-},\mathrm{Spec}^\mathrm{BK}\breve{\Phi}^r_{\psi,-},\mathrm{Spec}^\mathrm{BK}\breve{\Phi}^I_{\psi,-},	
\end{align}
\begin{align}
&\mathrm{Spec}^\mathrm{BK}{\Phi}_{\psi,-},
\mathrm{Spec}^\mathrm{BK}{\Phi}^r_{\psi,-},\mathrm{Spec}^\mathrm{BK}{\Phi}^I_{\psi,-}.	
\end{align}

Then we take the corresponding quotients by using the corresponding Frobenius operators:
\begin{align}
&\mathrm{Spec}^\mathrm{BK}\widetilde{\Phi}_{\psi,-}/\mathrm{Fro}^\mathbb{Z},	\\
\end{align}
\begin{align}
&\mathrm{Spec}^\mathrm{BK}\breve{\Phi}_{\psi,-}/\mathrm{Fro}^\mathbb{Z},	\\
\end{align}
\begin{align}
&\mathrm{Spec}^\mathrm{BK}{\Phi}_{\psi,-}/\mathrm{Fro}^\mathbb{Z}.	
\end{align}
Here for those space without notation related to the radius and the corresponding interval we consider the total unions $\bigcap_r,\bigcup_I$ in order to achieve the whole spaces to achieve the analogues of the corresponding FF curves from \cite{10KL1}, \cite{10KL2}, \cite{10FF} for
\[
\xymatrix@R+0pc@C+0pc{
\underset{r}{\mathrm{homotopylimit}}~\mathrm{Spec}^\mathrm{BK}\widetilde{\Phi}^r_{\psi,-},\underset{I}{\mathrm{homotopycolimit}}~\mathrm{Spec}^\mathrm{BK}\widetilde{\Phi}^I_{\psi,-},	\\
}
\]
\[
\xymatrix@R+0pc@C+0pc{
\underset{r}{\mathrm{homotopylimit}}~\mathrm{Spec}^\mathrm{BK}\breve{\Phi}^r_{\psi,-},\underset{I}{\mathrm{homotopycolimit}}~\mathrm{Spec}^\mathrm{BK}\breve{\Phi}^I_{\psi,-},	\\
}
\]
\[
\xymatrix@R+0pc@C+0pc{
\underset{r}{\mathrm{homotopylimit}}~\mathrm{Spec}^\mathrm{BK}{\Phi}^r_{\psi,-},\underset{I}{\mathrm{homotopycolimit}}~\mathrm{Spec}^\mathrm{BK}{\Phi}^I_{\psi,-}.	
}
\]
\[  
\xymatrix@R+0pc@C+0pc{
\underset{r}{\mathrm{homotopylimit}}~\mathrm{Spec}^\mathrm{BK}\widetilde{\Phi}^r_{\psi,-}/\mathrm{Fro}^\mathbb{Z},\underset{I}{\mathrm{homotopycolimit}}~\mathrm{Spec}^\mathrm{BK}\widetilde{\Phi}^I_{\psi,-}/\mathrm{Fro}^\mathbb{Z},	\\
}
\]
\[ 
\xymatrix@R+0pc@C+0pc{
\underset{r}{\mathrm{homotopylimit}}~\mathrm{Spec}^\mathrm{BK}\breve{\Phi}^r_{\psi,-}/\mathrm{Fro}^\mathbb{Z},\underset{I}{\mathrm{homotopycolimit}}~\mathrm{Spec}^\mathrm{BK}\breve{\Phi}^I_{\psi,-}/\mathrm{Fro}^\mathbb{Z},	\\
}
\]
\[ 
\xymatrix@R+0pc@C+0pc{
\underset{r}{\mathrm{homotopylimit}}~\mathrm{Spec}^\mathrm{BK}{\Phi}^r_{\psi,-}/\mathrm{Fro}^\mathbb{Z},\underset{I}{\mathrm{homotopycolimit}}~\mathrm{Spec}^\mathrm{BK}{\Phi}^I_{\psi,-}/\mathrm{Fro}^\mathbb{Z}.	
}
\]

\end{definition}

\indent Meanwhile we have the corresponding Clausen-Scholze analytic stacks from \cite{10CS2}, therefore applying their construction we have:

\begin{definition}
Here we define the following products by using the solidified tensor product from \cite{10CS1} and \cite{10CS2}. Namely $A$ will still as above as a Banach ring over $\mathbb{Q}_p$. Then we take solidified tensor product $\overset{\blacksquare}{\otimes}$ of any of the following
\[
\xymatrix@R+0pc@C+0pc{
\widetilde{\Delta}_{\psi},\widetilde{\nabla}_{\psi},\widetilde{\Phi}_{\psi},\widetilde{\Delta}^+_{\psi},\widetilde{\nabla}^+_{\psi},\widetilde{\Delta}^\dagger_{\psi},\widetilde{\nabla}^\dagger_{\psi},\widetilde{\Phi}^r_{\psi},\widetilde{\Phi}^I_{\psi}, 
}
\]

\[
\xymatrix@R+0pc@C+0pc{
\breve{\Delta}_{\psi},\breve{\nabla}_{\psi},\breve{\Phi}_{\psi},\breve{\Delta}^+_{\psi},\breve{\nabla}^+_{\psi},\breve{\Delta}^\dagger_{\psi},\breve{\nabla}^\dagger_{\psi},\breve{\Phi}^r_{\psi},\breve{\Phi}^I_{\psi},	
}
\]

\[
\xymatrix@R+0pc@C+0pc{
{\Delta}_{\psi},{\nabla}_{\psi},{\Phi}_{\psi},{\Delta}^+_{\psi},{\nabla}^+_{\psi},{\Delta}^\dagger_{\psi},{\nabla}^\dagger_{\psi},{\Phi}^r_{\psi},{\Phi}^I_{\psi},	
}
\]  	
with $A$. Then we have the notations:
\[
\xymatrix@R+0pc@C+0pc{
\widetilde{\Delta}_{\psi,-},\widetilde{\nabla}_{\psi,-},\widetilde{\Phi}_{\psi,-},\widetilde{\Delta}^+_{\psi,-},\widetilde{\nabla}^+_{\psi,-},\widetilde{\Delta}^\dagger_{\psi,-},\widetilde{\nabla}^\dagger_{\psi,-},\widetilde{\Phi}^r_{\psi,-},\widetilde{\Phi}^I_{\psi,-}, 
}
\]

\[
\xymatrix@R+0pc@C+0pc{
\breve{\Delta}_{\psi,-},\breve{\nabla}_{\psi,-},\breve{\Phi}_{\psi,-},\breve{\Delta}^+_{\psi,-},\breve{\nabla}^+_{\psi,-},\breve{\Delta}^\dagger_{\psi,-},\breve{\nabla}^\dagger_{\psi,-},\breve{\Phi}^r_{\psi,-},\breve{\Phi}^I_{\psi,-},	
}
\]

\[
\xymatrix@R+0pc@C+0pc{
{\Delta}_{\psi,-},{\nabla}_{\psi,-},{\Phi}_{\psi,-},{\Delta}^+_{\psi,-},{\nabla}^+_{\psi,-},{\Delta}^\dagger_{\psi,-},{\nabla}^\dagger_{\psi,-},{\Phi}^r_{\psi,-},{\Phi}^I_{\psi,-}.	
}
\]
\end{definition}

\begin{definition}
First we consider the Clausen-Scholze spectrum $\mathrm{Spec}^\mathrm{CS}(*)$ attached to any of those in the above from \cite{10CS2} by taking derived rational localization:
\begin{align}
\mathrm{Spec}^\mathrm{CS}\widetilde{\Delta}_{\psi,-},\mathrm{Spec}^\mathrm{CS}\widetilde{\nabla}_{\psi,-},\mathrm{Spec}^\mathrm{CS}\widetilde{\Phi}_{\psi,-},\mathrm{Spec}^\mathrm{CS}\widetilde{\Delta}^+_{\psi,-},\mathrm{Spec}^\mathrm{CS}\widetilde{\nabla}^+_{\psi,-},\\
\mathrm{Spec}^\mathrm{CS}\widetilde{\Delta}^\dagger_{\psi,-},\mathrm{Spec}^\mathrm{CS}\widetilde{\nabla}^\dagger_{\psi,-},\mathrm{Spec}^\mathrm{CS}\widetilde{\Phi}^r_{\psi,-},\mathrm{Spec}^\mathrm{CS}\widetilde{\Phi}^I_{\psi,-},	\\
\end{align}
\begin{align}
\mathrm{Spec}^\mathrm{CS}\breve{\Delta}_{\psi,-},\breve{\nabla}_{\psi,-},\mathrm{Spec}^\mathrm{CS}\breve{\Phi}_{\psi,-},\mathrm{Spec}^\mathrm{CS}\breve{\Delta}^+_{\psi,-},\mathrm{Spec}^\mathrm{CS}\breve{\nabla}^+_{\psi,-},\\
\mathrm{Spec}^\mathrm{CS}\breve{\Delta}^\dagger_{\psi,-},\mathrm{Spec}^\mathrm{CS}\breve{\nabla}^\dagger_{\psi,-},\mathrm{Spec}^\mathrm{CS}\breve{\Phi}^r_{\psi,-},\breve{\Phi}^I_{\psi,-},	\\
\end{align}
\begin{align}
\mathrm{Spec}^\mathrm{CS}{\Delta}_{\psi,-},\mathrm{Spec}^\mathrm{CS}{\nabla}_{\psi,-},\mathrm{Spec}^\mathrm{CS}{\Phi}_{\psi,-},\mathrm{Spec}^\mathrm{CS}{\Delta}^+_{\psi,-},\mathrm{Spec}^\mathrm{CS}{\nabla}^+_{\psi,-},\\
\mathrm{Spec}^\mathrm{CS}{\Delta}^\dagger_{\psi,-},\mathrm{Spec}^\mathrm{CS}{\nabla}^\dagger_{\psi,-},\mathrm{Spec}^\mathrm{CS}{\Phi}^r_{\psi,-},\mathrm{Spec}^\mathrm{CS}{\Phi}^I_{\psi,-}.	
\end{align}

Then we take the corresponding quotients by using the corresponding Frobenius operators:
\begin{align}
&\mathrm{Spec}^\mathrm{CS}\widetilde{\Delta}_{\psi,-}/\mathrm{Fro}^\mathbb{Z},\mathrm{Spec}^\mathrm{CS}\widetilde{\nabla}_{\psi,-}/\mathrm{Fro}^\mathbb{Z},\mathrm{Spec}^\mathrm{CS}\widetilde{\Phi}_{\psi,-}/\mathrm{Fro}^\mathbb{Z},\mathrm{Spec}^\mathrm{CS}\widetilde{\Delta}^+_{\psi,-}/\mathrm{Fro}^\mathbb{Z},\\
&\mathrm{Spec}^\mathrm{CS}\widetilde{\nabla}^+_{\psi,-}/\mathrm{Fro}^\mathbb{Z}, \mathrm{Spec}^\mathrm{CS}\widetilde{\Delta}^\dagger_{\psi,-}/\mathrm{Fro}^\mathbb{Z},\mathrm{Spec}^\mathrm{CS}\widetilde{\nabla}^\dagger_{\psi,-}/\mathrm{Fro}^\mathbb{Z},	\\
\end{align}
\begin{align}
&\mathrm{Spec}^\mathrm{CS}\breve{\Delta}_{\psi,-}/\mathrm{Fro}^\mathbb{Z},\breve{\nabla}_{\psi,-}/\mathrm{Fro}^\mathbb{Z},\mathrm{Spec}^\mathrm{CS}\breve{\Phi}_{\psi,-}/\mathrm{Fro}^\mathbb{Z},\mathrm{Spec}^\mathrm{CS}\breve{\Delta}^+_{\psi,-}/\mathrm{Fro}^\mathbb{Z},\\
&\mathrm{Spec}^\mathrm{CS}\breve{\nabla}^+_{\psi,-}/\mathrm{Fro}^\mathbb{Z}, \mathrm{Spec}^\mathrm{CS}\breve{\Delta}^\dagger_{\psi,-}/\mathrm{Fro}^\mathbb{Z},\mathrm{Spec}^\mathrm{CS}\breve{\nabla}^\dagger_{\psi,-}/\mathrm{Fro}^\mathbb{Z},	\\
\end{align}
\begin{align}
&\mathrm{Spec}^\mathrm{CS}{\Delta}_{\psi,-}/\mathrm{Fro}^\mathbb{Z},\mathrm{Spec}^\mathrm{CS}{\nabla}_{\psi,-}/\mathrm{Fro}^\mathbb{Z},\mathrm{Spec}^\mathrm{CS}{\Phi}_{\psi,-}/\mathrm{Fro}^\mathbb{Z},\mathrm{Spec}^\mathrm{CS}{\Delta}^+_{\psi,-}/\mathrm{Fro}^\mathbb{Z},\\
&\mathrm{Spec}^\mathrm{CS}{\nabla}^+_{\psi,-}/\mathrm{Fro}^\mathbb{Z}, \mathrm{Spec}^\mathrm{CS}{\Delta}^\dagger_{\psi,-}/\mathrm{Fro}^\mathbb{Z},\mathrm{Spec}^\mathrm{CS}{\nabla}^\dagger_{\psi,-}/\mathrm{Fro}^\mathbb{Z}.	
\end{align}
Here for those space with notations related to the radius and the corresponding interval we consider the total unions $\bigcap_r,\bigcup_I$ in order to achieve the whole spaces to achieve the analogues of the corresponding FF curves from \cite{10KL1}, \cite{10KL2}, \cite{10FF} for
\[
\xymatrix@R+0pc@C+0pc{
\underset{r}{\mathrm{homotopylimit}}~\mathrm{Spec}^\mathrm{CS}\widetilde{\Phi}^r_{\psi,-},\underset{I}{\mathrm{homotopycolimit}}~\mathrm{Spec}^\mathrm{CS}\widetilde{\Phi}^I_{\psi,-},	\\
}
\]
\[
\xymatrix@R+0pc@C+0pc{
\underset{r}{\mathrm{homotopylimit}}~\mathrm{Spec}^\mathrm{CS}\breve{\Phi}^r_{\psi,-},\underset{I}{\mathrm{homotopycolimit}}~\mathrm{Spec}^\mathrm{CS}\breve{\Phi}^I_{\psi,-},	\\
}
\]
\[
\xymatrix@R+0pc@C+0pc{
\underset{r}{\mathrm{homotopylimit}}~\mathrm{Spec}^\mathrm{CS}{\Phi}^r_{\psi,-},\underset{I}{\mathrm{homotopycolimit}}~\mathrm{Spec}^\mathrm{CS}{\Phi}^I_{\psi,-}.	
}
\]
\[ 
\xymatrix@R+0pc@C+0pc{
\underset{r}{\mathrm{homotopylimit}}~\mathrm{Spec}^\mathrm{CS}\widetilde{\Phi}^r_{\psi,-}/\mathrm{Fro}^\mathbb{Z},\underset{I}{\mathrm{homotopycolimit}}~\mathrm{Spec}^\mathrm{CS}\widetilde{\Phi}^I_{\psi,-}/\mathrm{Fro}^\mathbb{Z},	\\
}
\]
\[ 
\xymatrix@R+0pc@C+0pc{
\underset{r}{\mathrm{homotopylimit}}~\mathrm{Spec}^\mathrm{CS}\breve{\Phi}^r_{\psi,-}/\mathrm{Fro}^\mathbb{Z},\underset{I}{\mathrm{homotopycolimit}}~\breve{\Phi}^I_{\psi,-}/\mathrm{Fro}^\mathbb{Z},	\\
}
\]
\[ 
\xymatrix@R+0pc@C+0pc{
\underset{r}{\mathrm{homotopylimit}}~\mathrm{Spec}^\mathrm{CS}{\Phi}^r_{\psi,-}/\mathrm{Fro}^\mathbb{Z},\underset{I}{\mathrm{homotopycolimit}}~\mathrm{Spec}^\mathrm{CS}{\Phi}^I_{\psi,-}/\mathrm{Fro}^\mathbb{Z}.	
}
\]

\end{definition}

\

\begin{definition}
We then consider the corresponding quasipresheaves of the corresponding ind-Banach or monomorphic ind-Banach modules from \cite{10BBK}, \cite{10KKM}:
\begin{align}
\mathrm{Quasicoherentpresheaves,IndBanach}_{*}	
\end{align}
where $*$ is one of the following spaces:
\begin{align}
&\mathrm{Spec}^\mathrm{BK}\widetilde{\Phi}_{\psi,-}/\mathrm{Fro}^\mathbb{Z},	\\
\end{align}
\begin{align}
&\mathrm{Spec}^\mathrm{BK}\breve{\Phi}_{\psi,-}/\mathrm{Fro}^\mathbb{Z},	\\
\end{align}
\begin{align}
&\mathrm{Spec}^\mathrm{BK}{\Phi}_{\psi,-}/\mathrm{Fro}^\mathbb{Z}.	
\end{align}
Here for those space without notation related to the radius and the corresponding interval we consider the total unions $\bigcap_r,\bigcup_I$ in order to achieve the whole spaces to achieve the analogues of the corresponding FF curves from \cite{10KL1}, \cite{10KL2}, \cite{10FF} for
\[
\xymatrix@R+0pc@C+0pc{
\underset{r}{\mathrm{homotopylimit}}~\mathrm{Spec}^\mathrm{BK}\widetilde{\Phi}^r_{\psi,-},\underset{I}{\mathrm{homotopycolimit}}~\mathrm{Spec}^\mathrm{BK}\widetilde{\Phi}^I_{\psi,-},	\\
}
\]
\[
\xymatrix@R+0pc@C+0pc{
\underset{r}{\mathrm{homotopylimit}}~\mathrm{Spec}^\mathrm{BK}\breve{\Phi}^r_{\psi,-},\underset{I}{\mathrm{homotopycolimit}}~\mathrm{Spec}^\mathrm{BK}\breve{\Phi}^I_{\psi,-},	\\
}
\]
\[
\xymatrix@R+0pc@C+0pc{
\underset{r}{\mathrm{homotopylimit}}~\mathrm{Spec}^\mathrm{BK}{\Phi}^r_{\psi,-},\underset{I}{\mathrm{homotopycolimit}}~\mathrm{Spec}^\mathrm{BK}{\Phi}^I_{\psi,-}.	
}
\]
\[  
\xymatrix@R+0pc@C+0pc{
\underset{r}{\mathrm{homotopylimit}}~\mathrm{Spec}^\mathrm{BK}\widetilde{\Phi}^r_{\psi,-}/\mathrm{Fro}^\mathbb{Z},\underset{I}{\mathrm{homotopycolimit}}~\mathrm{Spec}^\mathrm{BK}\widetilde{\Phi}^I_{\psi,-}/\mathrm{Fro}^\mathbb{Z},	\\
}
\]
\[ 
\xymatrix@R+0pc@C+0pc{
\underset{r}{\mathrm{homotopylimit}}~\mathrm{Spec}^\mathrm{BK}\breve{\Phi}^r_{\psi,-}/\mathrm{Fro}^\mathbb{Z},\underset{I}{\mathrm{homotopycolimit}}~\mathrm{Spec}^\mathrm{BK}\breve{\Phi}^I_{\psi,-}/\mathrm{Fro}^\mathbb{Z},	\\
}
\]
\[ 
\xymatrix@R+0pc@C+0pc{
\underset{r}{\mathrm{homotopylimit}}~\mathrm{Spec}^\mathrm{BK}{\Phi}^r_{\psi,-}/\mathrm{Fro}^\mathbb{Z},\underset{I}{\mathrm{homotopycolimit}}~\mathrm{Spec}^\mathrm{BK}{\Phi}^I_{\psi,-}/\mathrm{Fro}^\mathbb{Z}.	
}
\]

\end{definition}

\begin{definition}
We then consider the corresponding quasisheaves of the corresponding condensed solid topological modules from \cite{10CS2}:
\begin{align}
\mathrm{Quasicoherentsheaves, Condensed}_{*}	
\end{align}
where $*$ is one of the following spaces:
\begin{align}
&\mathrm{Spec}^\mathrm{CS}\widetilde{\Delta}_{\psi,-}/\mathrm{Fro}^\mathbb{Z},\mathrm{Spec}^\mathrm{CS}\widetilde{\nabla}_{\psi,-}/\mathrm{Fro}^\mathbb{Z},\mathrm{Spec}^\mathrm{CS}\widetilde{\Phi}_{\psi,-}/\mathrm{Fro}^\mathbb{Z},\mathrm{Spec}^\mathrm{CS}\widetilde{\Delta}^+_{\psi,-}/\mathrm{Fro}^\mathbb{Z},\\
&\mathrm{Spec}^\mathrm{CS}\widetilde{\nabla}^+_{\psi,-}/\mathrm{Fro}^\mathbb{Z},\mathrm{Spec}^\mathrm{CS}\widetilde{\Delta}^\dagger_{\psi,-}/\mathrm{Fro}^\mathbb{Z},\mathrm{Spec}^\mathrm{CS}\widetilde{\nabla}^\dagger_{\psi,-}/\mathrm{Fro}^\mathbb{Z},	\\
\end{align}
\begin{align}
&\mathrm{Spec}^\mathrm{CS}\breve{\Delta}_{\psi,-}/\mathrm{Fro}^\mathbb{Z},\breve{\nabla}_{\psi,-}/\mathrm{Fro}^\mathbb{Z},\mathrm{Spec}^\mathrm{CS}\breve{\Phi}_{\psi,-}/\mathrm{Fro}^\mathbb{Z},\mathrm{Spec}^\mathrm{CS}\breve{\Delta}^+_{\psi,-}/\mathrm{Fro}^\mathbb{Z},\\
&\mathrm{Spec}^\mathrm{CS}\breve{\nabla}^+_{\psi,-}/\mathrm{Fro}^\mathbb{Z},\mathrm{Spec}^\mathrm{CS}\breve{\Delta}^\dagger_{\psi,-}/\mathrm{Fro}^\mathbb{Z},\mathrm{Spec}^\mathrm{CS}\breve{\nabla}^\dagger_{\psi,-}/\mathrm{Fro}^\mathbb{Z},	\\
\end{align}
\begin{align}
&\mathrm{Spec}^\mathrm{CS}{\Delta}_{\psi,-}/\mathrm{Fro}^\mathbb{Z},\mathrm{Spec}^\mathrm{CS}{\nabla}_{\psi,-}/\mathrm{Fro}^\mathbb{Z},\mathrm{Spec}^\mathrm{CS}{\Phi}_{\psi,-}/\mathrm{Fro}^\mathbb{Z},\mathrm{Spec}^\mathrm{CS}{\Delta}^+_{\psi,-}/\mathrm{Fro}^\mathbb{Z},\\
&\mathrm{Spec}^\mathrm{CS}{\nabla}^+_{\psi,-}/\mathrm{Fro}^\mathbb{Z}, \mathrm{Spec}^\mathrm{CS}{\Delta}^\dagger_{\psi,-}/\mathrm{Fro}^\mathbb{Z},\mathrm{Spec}^\mathrm{CS}{\nabla}^\dagger_{\psi,-}/\mathrm{Fro}^\mathbb{Z}.	
\end{align}
Here for those space with notations related to the radius and the corresponding interval we consider the total unions $\bigcap_r,\bigcup_I$ in order to achieve the whole spaces to achieve the analogues of the corresponding FF curves from \cite{10KL1}, \cite{10KL2}, \cite{10FF} for
\[
\xymatrix@R+0pc@C+0pc{
\underset{r}{\mathrm{homotopylimit}}~\mathrm{Spec}^\mathrm{CS}\widetilde{\Phi}^r_{\psi,-},\underset{I}{\mathrm{homotopycolimit}}~\mathrm{Spec}^\mathrm{CS}\widetilde{\Phi}^I_{\psi,-},	\\
}
\]
\[
\xymatrix@R+0pc@C+0pc{
\underset{r}{\mathrm{homotopylimit}}~\mathrm{Spec}^\mathrm{CS}\breve{\Phi}^r_{\psi,-},\underset{I}{\mathrm{homotopycolimit}}~\mathrm{Spec}^\mathrm{CS}\breve{\Phi}^I_{\psi,-},	\\
}
\]
\[
\xymatrix@R+0pc@C+0pc{
\underset{r}{\mathrm{homotopylimit}}~\mathrm{Spec}^\mathrm{CS}{\Phi}^r_{\psi,-},\underset{I}{\mathrm{homotopycolimit}}~\mathrm{Spec}^\mathrm{CS}{\Phi}^I_{\psi,-}.	
}
\]
\[ 
\xymatrix@R+0pc@C+0pc{
\underset{r}{\mathrm{homotopylimit}}~\mathrm{Spec}^\mathrm{CS}\widetilde{\Phi}^r_{\psi,-}/\mathrm{Fro}^\mathbb{Z},\underset{I}{\mathrm{homotopycolimit}}~\mathrm{Spec}^\mathrm{CS}\widetilde{\Phi}^I_{\psi,-}/\mathrm{Fro}^\mathbb{Z},	\\
}
\]
\[ 
\xymatrix@R+0pc@C+0pc{
\underset{r}{\mathrm{homotopylimit}}~\mathrm{Spec}^\mathrm{CS}\breve{\Phi}^r_{\psi,-}/\mathrm{Fro}^\mathbb{Z},\underset{I}{\mathrm{homotopycolimit}}~\breve{\Phi}^I_{\psi,-}/\mathrm{Fro}^\mathbb{Z},	\\
}
\]
\[ 
\xymatrix@R+0pc@C+0pc{
\underset{r}{\mathrm{homotopylimit}}~\mathrm{Spec}^\mathrm{CS}{\Phi}^r_{\psi,-}/\mathrm{Fro}^\mathbb{Z},\underset{I}{\mathrm{homotopycolimit}}~\mathrm{Spec}^\mathrm{CS}{\Phi}^I_{\psi,-}/\mathrm{Fro}^\mathbb{Z}.	
}
\]

\end{definition}

\

\begin{proposition}
There is a well-defined functor from the $\infty$-category 
\begin{align}
\mathrm{Quasicoherentpresheaves,Condensed}_{*}	
\end{align}
where $*$ is one of the following spaces:
\begin{align}
&\mathrm{Spec}^\mathrm{CS}\widetilde{\Phi}_{\psi,-}/\mathrm{Fro}^\mathbb{Z},	\\
\end{align}
\begin{align}
&\mathrm{Spec}^\mathrm{CS}\breve{\Phi}_{\psi,-}/\mathrm{Fro}^\mathbb{Z},	\\
\end{align}
\begin{align}
&\mathrm{Spec}^\mathrm{CS}{\Phi}_{\psi,-}/\mathrm{Fro}^\mathbb{Z},	
\end{align}
to the $\infty$-category of $\mathrm{Fro}$-equivariant quasicoherent presheaves over similar spaces above correspondingly without the $\mathrm{Fro}$-quotients, and to the $\infty$-category of $\mathrm{Fro}$-equivariant quasicoherent modules over global sections of the structure $\infty$-sheaves of the similar spaces above correspondingly without the $\mathrm{Fro}$-quotients. Here for those space without notation related to the radius and the corresponding interval we consider the total unions $\bigcap_r,\bigcup_I$ in order to achieve the whole spaces to achieve the analogues of the corresponding FF curves from \cite{10KL1}, \cite{10KL2}, \cite{10FF} for
\[
\xymatrix@R+0pc@C+0pc{
\underset{r}{\mathrm{homotopylimit}}~\mathrm{Spec}^\mathrm{CS}\widetilde{\Phi}^r_{\psi,-},\underset{I}{\mathrm{homotopycolimit}}~\mathrm{Spec}^\mathrm{CS}\widetilde{\Phi}^I_{\psi,-},	\\
}
\]
\[
\xymatrix@R+0pc@C+0pc{
\underset{r}{\mathrm{homotopylimit}}~\mathrm{Spec}^\mathrm{CS}\breve{\Phi}^r_{\psi,-},\underset{I}{\mathrm{homotopycolimit}}~\mathrm{Spec}^\mathrm{CS}\breve{\Phi}^I_{\psi,-},	\\
}
\]
\[
\xymatrix@R+0pc@C+0pc{
\underset{r}{\mathrm{homotopylimit}}~\mathrm{Spec}^\mathrm{CS}{\Phi}^r_{\psi,-},\underset{I}{\mathrm{homotopycolimit}}~\mathrm{Spec}^\mathrm{CS}{\Phi}^I_{\psi,-}.	
}
\]
\[ 
\xymatrix@R+0pc@C+0pc{
\underset{r}{\mathrm{homotopylimit}}~\mathrm{Spec}^\mathrm{CS}\widetilde{\Phi}^r_{\psi,-}/\mathrm{Fro}^\mathbb{Z},\underset{I}{\mathrm{homotopycolimit}}~\mathrm{Spec}^\mathrm{CS}\widetilde{\Phi}^I_{\psi,-}/\mathrm{Fro}^\mathbb{Z},	\\
}
\]
\[ 
\xymatrix@R+0pc@C+0pc{
\underset{r}{\mathrm{homotopylimit}}~\mathrm{Spec}^\mathrm{CS}\breve{\Phi}^r_{\psi,-}/\mathrm{Fro}^\mathbb{Z},\underset{I}{\mathrm{homotopycolimit}}~\breve{\Phi}^I_{\psi,-}/\mathrm{Fro}^\mathbb{Z},	\\
}
\]
\[ 
\xymatrix@R+0pc@C+0pc{
\underset{r}{\mathrm{homotopylimit}}~\mathrm{Spec}^\mathrm{CS}{\Phi}^r_{\psi,-}/\mathrm{Fro}^\mathbb{Z},\underset{I}{\mathrm{homotopycolimit}}~\mathrm{Spec}^\mathrm{CS}{\Phi}^I_{\psi,-}/\mathrm{Fro}^\mathbb{Z}.	
}
\]	
In this situation we will have the target category being family parametrized by $r$ or $I$ in compatible glueing sense as in \cite[Definition 5.4.10]{10KL2}. In this situation for modules parametrized by the intervals we have the equivalence of $\infty$-categories by using \cite[Proposition 13.8]{10CS2}. Here the corresponding quasicoherent Frobenius modules are defined to be the corresponding homotopy colimits and limits of Frobenius modules:
\begin{align}
\underset{r}{\mathrm{homotopycolimit}}~M_r,\\
\underset{I}{\mathrm{homotopylimit}}~M_I,	
\end{align}
where each $M_r$ is a Frobenius-equivariant module over the period ring with respect to some radius $r$ while each $M_I$ is a Frobenius-equivariant module over the period ring with respect to some interval $I$.\\
\end{proposition}

\begin{proposition}
Similar proposition holds for 
\begin{align}
\mathrm{Quasicoherentsheaves,IndBanach}_{*}.	
\end{align}	
\end{proposition}

\

\begin{definition}
We then consider the corresponding quasipresheaves of perfect complexes the corresponding ind-Banach or monomorphic ind-Banach modules from \cite{10BBK}, \cite{10KKM}:
\begin{align}
\mathrm{Quasicoherentpresheaves,Perfectcomplex,IndBanach}_{*}	
\end{align}
where $*$ is one of the following spaces:
\begin{align}
&\mathrm{Spec}^\mathrm{BK}\widetilde{\Phi}_{\psi,-}/\mathrm{Fro}^\mathbb{Z},	\\
\end{align}
\begin{align}
&\mathrm{Spec}^\mathrm{BK}\breve{\Phi}_{\psi,-}/\mathrm{Fro}^\mathbb{Z},	\\
\end{align}
\begin{align}
&\mathrm{Spec}^\mathrm{BK}{\Phi}_{\psi,-}/\mathrm{Fro}^\mathbb{Z}.	
\end{align}
Here for those space without notation related to the radius and the corresponding interval we consider the total unions $\bigcap_r,\bigcup_I$ in order to achieve the whole spaces to achieve the analogues of the corresponding FF curves from \cite{10KL1}, \cite{10KL2}, \cite{10FF} for
\[
\xymatrix@R+0pc@C+0pc{
\underset{r}{\mathrm{homotopylimit}}~\mathrm{Spec}^\mathrm{BK}\widetilde{\Phi}^r_{\psi,-},\underset{I}{\mathrm{homotopycolimit}}~\mathrm{Spec}^\mathrm{BK}\widetilde{\Phi}^I_{\psi,-},	\\
}
\]
\[
\xymatrix@R+0pc@C+0pc{
\underset{r}{\mathrm{homotopylimit}}~\mathrm{Spec}^\mathrm{BK}\breve{\Phi}^r_{\psi,-},\underset{I}{\mathrm{homotopycolimit}}~\mathrm{Spec}^\mathrm{BK}\breve{\Phi}^I_{\psi,-},	\\
}
\]
\[
\xymatrix@R+0pc@C+0pc{
\underset{r}{\mathrm{homotopylimit}}~\mathrm{Spec}^\mathrm{BK}{\Phi}^r_{\psi,-},\underset{I}{\mathrm{homotopycolimit}}~\mathrm{Spec}^\mathrm{BK}{\Phi}^I_{\psi,-}.	
}
\]
\[  
\xymatrix@R+0pc@C+0pc{
\underset{r}{\mathrm{homotopylimit}}~\mathrm{Spec}^\mathrm{BK}\widetilde{\Phi}^r_{\psi,-}/\mathrm{Fro}^\mathbb{Z},\underset{I}{\mathrm{homotopycolimit}}~\mathrm{Spec}^\mathrm{BK}\widetilde{\Phi}^I_{\psi,-}/\mathrm{Fro}^\mathbb{Z},	\\
}
\]
\[ 
\xymatrix@R+0pc@C+0pc{
\underset{r}{\mathrm{homotopylimit}}~\mathrm{Spec}^\mathrm{BK}\breve{\Phi}^r_{\psi,-}/\mathrm{Fro}^\mathbb{Z},\underset{I}{\mathrm{homotopycolimit}}~\mathrm{Spec}^\mathrm{BK}\breve{\Phi}^I_{\psi,-}/\mathrm{Fro}^\mathbb{Z},	\\
}
\]
\[ 
\xymatrix@R+0pc@C+0pc{
\underset{r}{\mathrm{homotopylimit}}~\mathrm{Spec}^\mathrm{BK}{\Phi}^r_{\psi,-}/\mathrm{Fro}^\mathbb{Z},\underset{I}{\mathrm{homotopycolimit}}~\mathrm{Spec}^\mathrm{BK}{\Phi}^I_{\psi,-}/\mathrm{Fro}^\mathbb{Z}.	
}
\]

\end{definition}

\begin{definition}
We then consider the corresponding quasisheaves of perfect complexes of the corresponding condensed solid topological modules from \cite{10CS2}:
\begin{align}
\mathrm{Quasicoherentsheaves, Perfectcomplex, Condensed}_{*}	
\end{align}
where $*$ is one of the following spaces:
\begin{align}
&\mathrm{Spec}^\mathrm{CS}\widetilde{\Delta}_{\psi,-}/\mathrm{Fro}^\mathbb{Z},\mathrm{Spec}^\mathrm{CS}\widetilde{\nabla}_{\psi,-}/\mathrm{Fro}^\mathbb{Z},\mathrm{Spec}^\mathrm{CS}\widetilde{\Phi}_{\psi,-}/\mathrm{Fro}^\mathbb{Z},\mathrm{Spec}^\mathrm{CS}\widetilde{\Delta}^+_{\psi,-}/\mathrm{Fro}^\mathbb{Z},\\
&\mathrm{Spec}^\mathrm{CS}\widetilde{\nabla}^+_{\psi,-}/\mathrm{Fro}^\mathbb{Z},\mathrm{Spec}^\mathrm{CS}\widetilde{\Delta}^\dagger_{\psi,-}/\mathrm{Fro}^\mathbb{Z},\mathrm{Spec}^\mathrm{CS}\widetilde{\nabla}^\dagger_{\psi,-}/\mathrm{Fro}^\mathbb{Z},	\\
\end{align}
\begin{align}
&\mathrm{Spec}^\mathrm{CS}\breve{\Delta}_{\psi,-}/\mathrm{Fro}^\mathbb{Z},\breve{\nabla}_{\psi,-}/\mathrm{Fro}^\mathbb{Z},\mathrm{Spec}^\mathrm{CS}\breve{\Phi}_{\psi,-}/\mathrm{Fro}^\mathbb{Z},\mathrm{Spec}^\mathrm{CS}\breve{\Delta}^+_{\psi,-}/\mathrm{Fro}^\mathbb{Z},\\
&\mathrm{Spec}^\mathrm{CS}\breve{\nabla}^+_{\psi,-}/\mathrm{Fro}^\mathbb{Z},\mathrm{Spec}^\mathrm{CS}\breve{\Delta}^\dagger_{\psi,-}/\mathrm{Fro}^\mathbb{Z},\mathrm{Spec}^\mathrm{CS}\breve{\nabla}^\dagger_{\psi,-}/\mathrm{Fro}^\mathbb{Z},	\\
\end{align}
\begin{align}
&\mathrm{Spec}^\mathrm{CS}{\Delta}_{\psi,-}/\mathrm{Fro}^\mathbb{Z},\mathrm{Spec}^\mathrm{CS}{\nabla}_{\psi,-}/\mathrm{Fro}^\mathbb{Z},\mathrm{Spec}^\mathrm{CS}{\Phi}_{\psi,-}/\mathrm{Fro}^\mathbb{Z},\mathrm{Spec}^\mathrm{CS}{\Delta}^+_{\psi,-}/\mathrm{Fro}^\mathbb{Z},\\
&\mathrm{Spec}^\mathrm{CS}{\nabla}^+_{\psi,-}/\mathrm{Fro}^\mathbb{Z}, \mathrm{Spec}^\mathrm{CS}{\Delta}^\dagger_{\psi,-}/\mathrm{Fro}^\mathbb{Z},\mathrm{Spec}^\mathrm{CS}{\nabla}^\dagger_{\psi,-}/\mathrm{Fro}^\mathbb{Z}.	
\end{align}
Here for those space with notations related to the radius and the corresponding interval we consider the total unions $\bigcap_r,\bigcup_I$ in order to achieve the whole spaces to achieve the analogues of the corresponding FF curves from \cite{10KL1}, \cite{10KL2}, \cite{10FF} for
\[
\xymatrix@R+0pc@C+0pc{
\underset{r}{\mathrm{homotopylimit}}~\mathrm{Spec}^\mathrm{CS}\widetilde{\Phi}^r_{\psi,-},\underset{I}{\mathrm{homotopycolimit}}~\mathrm{Spec}^\mathrm{CS}\widetilde{\Phi}^I_{\psi,-},	\\
}
\]
\[
\xymatrix@R+0pc@C+0pc{
\underset{r}{\mathrm{homotopylimit}}~\mathrm{Spec}^\mathrm{CS}\breve{\Phi}^r_{\psi,-},\underset{I}{\mathrm{homotopycolimit}}~\mathrm{Spec}^\mathrm{CS}\breve{\Phi}^I_{\psi,-},	\\
}
\]
\[
\xymatrix@R+0pc@C+0pc{
\underset{r}{\mathrm{homotopylimit}}~\mathrm{Spec}^\mathrm{CS}{\Phi}^r_{\psi,-},\underset{I}{\mathrm{homotopycolimit}}~\mathrm{Spec}^\mathrm{CS}{\Phi}^I_{\psi,-}.	
}
\]
\[ 
\xymatrix@R+0pc@C+0pc{
\underset{r}{\mathrm{homotopylimit}}~\mathrm{Spec}^\mathrm{CS}\widetilde{\Phi}^r_{\psi,-}/\mathrm{Fro}^\mathbb{Z},\underset{I}{\mathrm{homotopycolimit}}~\mathrm{Spec}^\mathrm{CS}\widetilde{\Phi}^I_{\psi,-}/\mathrm{Fro}^\mathbb{Z},	\\
}
\]
\[ 
\xymatrix@R+0pc@C+0pc{
\underset{r}{\mathrm{homotopylimit}}~\mathrm{Spec}^\mathrm{CS}\breve{\Phi}^r_{\psi,-}/\mathrm{Fro}^\mathbb{Z},\underset{I}{\mathrm{homotopycolimit}}~\breve{\Phi}^I_{\psi,-}/\mathrm{Fro}^\mathbb{Z},	\\
}
\]
\[ 
\xymatrix@R+0pc@C+0pc{
\underset{r}{\mathrm{homotopylimit}}~\mathrm{Spec}^\mathrm{CS}{\Phi}^r_{\psi,-}/\mathrm{Fro}^\mathbb{Z},\underset{I}{\mathrm{homotopycolimit}}~\mathrm{Spec}^\mathrm{CS}{\Phi}^I_{\psi,-}/\mathrm{Fro}^\mathbb{Z}.	
}
\]

\end{definition}

\begin{proposition}
There is a well-defined functor from the $\infty$-category 
\begin{align}
\mathrm{Quasicoherentpresheaves,Perfectcomplex,Condensed}_{*}	
\end{align}
where $*$ is one of the following spaces:
\begin{align}
&\mathrm{Spec}^\mathrm{CS}\widetilde{\Phi}_{\psi,-}/\mathrm{Fro}^\mathbb{Z},	\\
\end{align}
\begin{align}
&\mathrm{Spec}^\mathrm{CS}\breve{\Phi}_{\psi,-}/\mathrm{Fro}^\mathbb{Z},	\\
\end{align}
\begin{align}
&\mathrm{Spec}^\mathrm{CS}{\Phi}_{\psi,-}/\mathrm{Fro}^\mathbb{Z},	
\end{align}
to the $\infty$-category of $\mathrm{Fro}$-equivariant quasicoherent presheaves over similar spaces above correspondingly without the $\mathrm{Fro}$-quotients, and to the $\infty$-category of $\mathrm{Fro}$-equivariant quasicoherent modules over global sections of the structure $\infty$-sheaves of the similar spaces above correspondingly without the $\mathrm{Fro}$-quotients. Here for those space without notation related to the radius and the corresponding interval we consider the total unions $\bigcap_r,\bigcup_I$ in order to achieve the whole spaces to achieve the analogues of the corresponding FF curves from \cite{10KL1}, \cite{10KL2}, \cite{10FF} for
\[
\xymatrix@R+0pc@C+0pc{
\underset{r}{\mathrm{homotopylimit}}~\mathrm{Spec}^\mathrm{CS}\widetilde{\Phi}^r_{\psi,-},\underset{I}{\mathrm{homotopycolimit}}~\mathrm{Spec}^\mathrm{CS}\widetilde{\Phi}^I_{\psi,-},	\\
}
\]
\[
\xymatrix@R+0pc@C+0pc{
\underset{r}{\mathrm{homotopylimit}}~\mathrm{Spec}^\mathrm{CS}\breve{\Phi}^r_{\psi,-},\underset{I}{\mathrm{homotopycolimit}}~\mathrm{Spec}^\mathrm{CS}\breve{\Phi}^I_{\psi,-},	\\
}
\]
\[
\xymatrix@R+0pc@C+0pc{
\underset{r}{\mathrm{homotopylimit}}~\mathrm{Spec}^\mathrm{CS}{\Phi}^r_{\psi,-},\underset{I}{\mathrm{homotopycolimit}}~\mathrm{Spec}^\mathrm{CS}{\Phi}^I_{\psi,-}.	
}
\]
\[ 
\xymatrix@R+0pc@C+0pc{
\underset{r}{\mathrm{homotopylimit}}~\mathrm{Spec}^\mathrm{CS}\widetilde{\Phi}^r_{\psi,-}/\mathrm{Fro}^\mathbb{Z},\underset{I}{\mathrm{homotopycolimit}}~\mathrm{Spec}^\mathrm{CS}\widetilde{\Phi}^I_{\psi,-}/\mathrm{Fro}^\mathbb{Z},	\\
}
\]
\[ 
\xymatrix@R+0pc@C+0pc{
\underset{r}{\mathrm{homotopylimit}}~\mathrm{Spec}^\mathrm{CS}\breve{\Phi}^r_{\psi,-}/\mathrm{Fro}^\mathbb{Z},\underset{I}{\mathrm{homotopycolimit}}~\breve{\Phi}^I_{\psi,-}/\mathrm{Fro}^\mathbb{Z},	\\
}
\]
\[ 
\xymatrix@R+0pc@C+0pc{
\underset{r}{\mathrm{homotopylimit}}~\mathrm{Spec}^\mathrm{CS}{\Phi}^r_{\psi,-}/\mathrm{Fro}^\mathbb{Z},\underset{I}{\mathrm{homotopycolimit}}~\mathrm{Spec}^\mathrm{CS}{\Phi}^I_{\psi,-}/\mathrm{Fro}^\mathbb{Z}.	
}
\]	
In this situation we will have the target category being family parametrized by $r$ or $I$ in compatible glueing sense as in \cite[Definition 5.4.10]{10KL2}. In this situation for modules parametrized by the intervals we have the equivalence of $\infty$-categories by using \cite[Proposition 12.18]{10CS2}. Here the corresponding quasicoherent Frobenius modules are defined to be the corresponding homotopy colimits and limits of Frobenius modules:
\begin{align}
\underset{r}{\mathrm{homotopycolimit}}~M_r,\\
\underset{I}{\mathrm{homotopylimit}}~M_I,	
\end{align}
where each $M_r$ is a Frobenius-equivariant module over the period ring with respect to some radius $r$ while each $M_I$ is a Frobenius-equivariant module over the period ring with respect to some interval $I$.\\
\end{proposition}

\begin{proposition}
Similar proposition holds for 
\begin{align}
\mathrm{Quasicoherentsheaves,Perfectcomplex,IndBanach}_{*}.	
\end{align}	
\end{proposition}

\section{Multivariate Hodge Iwasawa Prestacks}

This chapter follows closely \cite{10T1}, \cite{10T2}, \cite{10T3}, \cite{10KPX}, \cite{10KP}, \cite{10KL1}, \cite{10KL2}, \cite{10BK}, \cite{10BBBK}, \cite{10BBM}, \cite{10KKM}, \cite{10CS1}, \cite{10CS2}, \cite{10CKZ}, \cite{10PZ}, \cite{10BCM}, \cite{10LBV}.

\subsection{Frobenius Quasicoherent Prestacks I}

\begin{definition}
We now consider the pro-\'etale site of $\mathrm{Spa}\mathbb{Q}_p\left<X_1^{\pm 1},...,X_k^{\pm 1}\right>$ from \cite{10Sch}, denote that by $*$. To be more accurate we replace one component for $\Gamma$ with the pro-\'etale site of $\mathrm{Spa}\mathbb{Q}_p\left<X_1^{\pm 1},...,X_k^{\pm 1}\right>$. And we treat then all the functor to be prestacks for this site\footnote{Here for those imperfect rings, the notation will mean that the specific component forming the pro-\'etale site will be the perfect version of the corresponding ring. Certainly if we have $|\Gamma|=1$ then we have that all the rings are perfect in \cite{10KL1} and \cite{10KL2}.}. Then from \cite{10KL1} and \cite[Definition 5.2.1]{10KL2} we have the following class of Kedlaya-Liu rings (with the following replacement: $\Delta$ stands for $A$, $\nabla$ stands for $B$, while $\Phi$ stands for $C$) by taking product in the sense of self $\Gamma$-th power:

\[
\xymatrix@R+0pc@C+0pc{
\widetilde{\Delta}_{*,\Gamma},\widetilde{\nabla}_{*,\Gamma},\widetilde{\Phi}_{*,\Gamma},\widetilde{\Delta}^+_{*,\Gamma},\widetilde{\nabla}^+_{*,\Gamma},\widetilde{\Delta}^\dagger_{*,\Gamma},\widetilde{\nabla}^\dagger_{*,\Gamma},\widetilde{\Phi}^r_{*,\Gamma},\widetilde{\Phi}^I_{*,\Gamma}, 
}
\]

\[
\xymatrix@R+0pc@C+0pc{
\breve{\Delta}_{*,\Gamma},\breve{\nabla}_{*,\Gamma},\breve{\Phi}_{*,\Gamma},\breve{\Delta}^+_{*,\Gamma},\breve{\nabla}^+_{*,\Gamma},\breve{\Delta}^\dagger_{*,\Gamma},\breve{\nabla}^\dagger_{*,\Gamma},\breve{\Phi}^r_{*,\Gamma},\breve{\Phi}^I_{*,\Gamma},	
}
\]

\[
\xymatrix@R+0pc@C+0pc{
{\Delta}_{*,\Gamma},{\nabla}_{*,\Gamma},{\Phi}_{*,\Gamma},{\Delta}^+_{*,\Gamma},{\nabla}^+_{*,\Gamma},{\Delta}^\dagger_{*,\Gamma},{\nabla}^\dagger_{*,\Gamma},{\Phi}^r_{*,\Gamma},{\Phi}^I_{*,\Gamma}.	
}
\]
We now consider the following rings with $A$ being a Banach ring over $\mathbb{Q}_p$. Taking the product we have:
\[
\xymatrix@R+0pc@C+0pc{
\widetilde{\Phi}_{*,\Gamma,A},\widetilde{\Phi}^r_{*,\Gamma,A},\widetilde{\Phi}^I_{*,\Gamma,A},	
}
\]
\[
\xymatrix@R+0pc@C+0pc{
\breve{\Phi}_{*,\Gamma,A},\breve{\Phi}^r_{*,\Gamma,A},\breve{\Phi}^I_{*,\Gamma,A},	
}
\]
\[
\xymatrix@R+0pc@C+0pc{
{\Phi}_{*,\Gamma,A},{\Phi}^r_{*,\Gamma,A},{\Phi}^I_{*,\Gamma,A}.	
}
\]
They carry multi Frobenius action $\varphi_\Gamma$ and multi $\mathrm{Lie}_\Gamma:=\mathbb{Z}_p^{\times\Gamma}$ action. In our current situation after \cite{10CKZ} and \cite{10PZ} we consider the following $(\infty,1)$-categories of $(\infty,1)$-modules.\\
\end{definition}

\begin{definition}
First we consider the Bambozzi-Kremnizer spectrum $\mathrm{Spec}^\mathrm{BK}(*)$ attached to any of those in the above from \cite{10BK} by taking derived rational localization:
\begin{align}
&\mathrm{Spec}^\mathrm{BK}\widetilde{\Phi}_{*,\Gamma,A},\mathrm{Spec}^\mathrm{BK}\widetilde{\Phi}^r_{*,\Gamma,A},\mathrm{Spec}^\mathrm{BK}\widetilde{\Phi}^I_{*,\Gamma,A},	
\end{align}
\begin{align}
&\mathrm{Spec}^\mathrm{BK}\breve{\Phi}_{*,\Gamma,A},\mathrm{Spec}^\mathrm{BK}\breve{\Phi}^r_{*,\Gamma,A},\mathrm{Spec}^\mathrm{BK}\breve{\Phi}^I_{*,\Gamma,A},	
\end{align}
\begin{align}
&\mathrm{Spec}^\mathrm{BK}{\Phi}_{*,\Gamma,A},
\mathrm{Spec}^\mathrm{BK}{\Phi}^r_{*,\Gamma,A},\mathrm{Spec}^\mathrm{BK}{\Phi}^I_{*,\Gamma,A}.	
\end{align}

Then we take the corresponding quotients by using the corresponding Frobenius operators:
\begin{align}
&\mathrm{Spec}^\mathrm{BK}\widetilde{\Phi}_{*,\Gamma,A}/\mathrm{Fro}^\mathbb{Z},	\\
\end{align}
\begin{align}
&\mathrm{Spec}^\mathrm{BK}\breve{\Phi}_{*,\Gamma,A}/\mathrm{Fro}^\mathbb{Z},	\\
\end{align}
\begin{align}
&\mathrm{Spec}^\mathrm{BK}{\Phi}_{*,\Gamma,A}/\mathrm{Fro}^\mathbb{Z}.	
\end{align}
Here for those space without notation related to the radius and the corresponding interval we consider the total unions $\bigcap_r,\bigcup_I$ in order to achieve the whole spaces to achieve the analogues of the corresponding FF curves from \cite{10KL1}, \cite{10KL2}, \cite{10FF} for
\[
\xymatrix@R+0pc@C+0pc{
\underset{r}{\mathrm{homotopylimit}}~\mathrm{Spec}^\mathrm{BK}\widetilde{\Phi}^r_{*,\Gamma,A},\underset{I}{\mathrm{homotopycolimit}}~\mathrm{Spec}^\mathrm{BK}\widetilde{\Phi}^I_{*,\Gamma,A},	\\
}
\]
\[
\xymatrix@R+0pc@C+0pc{
\underset{r}{\mathrm{homotopylimit}}~\mathrm{Spec}^\mathrm{BK}\breve{\Phi}^r_{*,\Gamma,A},\underset{I}{\mathrm{homotopycolimit}}~\mathrm{Spec}^\mathrm{BK}\breve{\Phi}^I_{*,\Gamma,A},	\\
}
\]
\[
\xymatrix@R+0pc@C+0pc{
\underset{r}{\mathrm{homotopylimit}}~\mathrm{Spec}^\mathrm{BK}{\Phi}^r_{*,\Gamma,A},\underset{I}{\mathrm{homotopycolimit}}~\mathrm{Spec}^\mathrm{BK}{\Phi}^I_{*,\Gamma,A}.	
}
\]
\[  
\xymatrix@R+0pc@C+0pc{
\underset{r}{\mathrm{homotopylimit}}~\mathrm{Spec}^\mathrm{BK}\widetilde{\Phi}^r_{*,\Gamma,A}/\mathrm{Fro}^\mathbb{Z},\underset{I}{\mathrm{homotopycolimit}}~\mathrm{Spec}^\mathrm{BK}\widetilde{\Phi}^I_{*,\Gamma,A}/\mathrm{Fro}^\mathbb{Z},	\\
}
\]
\[ 
\xymatrix@R+0pc@C+0pc{
\underset{r}{\mathrm{homotopylimit}}~\mathrm{Spec}^\mathrm{BK}\breve{\Phi}^r_{*,\Gamma,A}/\mathrm{Fro}^\mathbb{Z},\underset{I}{\mathrm{homotopycolimit}}~\mathrm{Spec}^\mathrm{BK}\breve{\Phi}^I_{*,\Gamma,A}/\mathrm{Fro}^\mathbb{Z},	\\
}
\]
\[ 
\xymatrix@R+0pc@C+0pc{
\underset{r}{\mathrm{homotopylimit}}~\mathrm{Spec}^\mathrm{BK}{\Phi}^r_{*,\Gamma,A}/\mathrm{Fro}^\mathbb{Z},\underset{I}{\mathrm{homotopycolimit}}~\mathrm{Spec}^\mathrm{BK}{\Phi}^I_{*,\Gamma,A}/\mathrm{Fro}^\mathbb{Z}.	
}
\]

\end{definition}

\indent Meanwhile we have the corresponding Clausen-Scholze analytic stacks from \cite{10CS2}, therefore applying their construction we have:

\begin{definition}
Here we define the following products by using the solidified tensor product from \cite{10CS1} and \cite{10CS2}. Namely $A$ will still as above as a Banach ring over $\mathbb{Q}_p$. Then we take solidified tensor product $\overset{\blacksquare}{\otimes}$ of any of the following
\[
\xymatrix@R+0pc@C+0pc{
\widetilde{\Delta}_{*,\Gamma},\widetilde{\nabla}_{*,\Gamma},\widetilde{\Phi}_{*,\Gamma},\widetilde{\Delta}^+_{*,\Gamma},\widetilde{\nabla}^+_{*,\Gamma},\widetilde{\Delta}^\dagger_{*,\Gamma},\widetilde{\nabla}^\dagger_{*,\Gamma},\widetilde{\Phi}^r_{*,\Gamma},\widetilde{\Phi}^I_{*,\Gamma}, 
}
\]

\[
\xymatrix@R+0pc@C+0pc{
\breve{\Delta}_{*,\Gamma},\breve{\nabla}_{*,\Gamma},\breve{\Phi}_{*,\Gamma},\breve{\Delta}^+_{*,\Gamma},\breve{\nabla}^+_{*,\Gamma},\breve{\Delta}^\dagger_{*,\Gamma},\breve{\nabla}^\dagger_{*,\Gamma},\breve{\Phi}^r_{*,\Gamma},\breve{\Phi}^I_{*,\Gamma},	
}
\]

\[
\xymatrix@R+0pc@C+0pc{
{\Delta}_{*,\Gamma},{\nabla}_{*,\Gamma},{\Phi}_{*,\Gamma},{\Delta}^+_{*,\Gamma},{\nabla}^+_{*,\Gamma},{\Delta}^\dagger_{*,\Gamma},{\nabla}^\dagger_{*,\Gamma},{\Phi}^r_{*,\Gamma},{\Phi}^I_{*,\Gamma},	
}
\]  	
with $A$. Then we have the notations:
\[
\xymatrix@R+0pc@C+0pc{
\widetilde{\Delta}_{*,\Gamma,A},\widetilde{\nabla}_{*,\Gamma,A},\widetilde{\Phi}_{*,\Gamma,A},\widetilde{\Delta}^+_{*,\Gamma,A},\widetilde{\nabla}^+_{*,\Gamma,A},\widetilde{\Delta}^\dagger_{*,\Gamma,A},\widetilde{\nabla}^\dagger_{*,\Gamma,A},\widetilde{\Phi}^r_{*,\Gamma,A},\widetilde{\Phi}^I_{*,\Gamma,A}, 
}
\]

\[
\xymatrix@R+0pc@C+0pc{
\breve{\Delta}_{*,\Gamma,A},\breve{\nabla}_{*,\Gamma,A},\breve{\Phi}_{*,\Gamma,A},\breve{\Delta}^+_{*,\Gamma,A},\breve{\nabla}^+_{*,\Gamma,A},\breve{\Delta}^\dagger_{*,\Gamma,A},\breve{\nabla}^\dagger_{*,\Gamma,A},\breve{\Phi}^r_{*,\Gamma,A},\breve{\Phi}^I_{*,\Gamma,A},	
}
\]

\[
\xymatrix@R+0pc@C+0pc{
{\Delta}_{*,\Gamma,A},{\nabla}_{*,\Gamma,A},{\Phi}_{*,\Gamma,A},{\Delta}^+_{*,\Gamma,A},{\nabla}^+_{*,\Gamma,A},{\Delta}^\dagger_{*,\Gamma,A},{\nabla}^\dagger_{*,\Gamma,A},{\Phi}^r_{*,\Gamma,A},{\Phi}^I_{*,\Gamma,A}.	
}
\]
\end{definition}

\begin{definition}
First we consider the Clausen-Scholze spectrum $\mathrm{Spec}^\mathrm{CS}(*)$ attached to any of those in the above from \cite{10CS2} by taking derived rational localization:
\begin{align}
\mathrm{Spec}^\mathrm{CS}\widetilde{\Delta}_{*,\Gamma,A},\mathrm{Spec}^\mathrm{CS}\widetilde{\nabla}_{*,\Gamma,A},\mathrm{Spec}^\mathrm{CS}\widetilde{\Phi}_{*,\Gamma,A},\mathrm{Spec}^\mathrm{CS}\widetilde{\Delta}^+_{*,\Gamma,A},\mathrm{Spec}^\mathrm{CS}\widetilde{\nabla}^+_{*,\Gamma,A},\\
\mathrm{Spec}^\mathrm{CS}\widetilde{\Delta}^\dagger_{*,\Gamma,A},\mathrm{Spec}^\mathrm{CS}\widetilde{\nabla}^\dagger_{*,\Gamma,A},\mathrm{Spec}^\mathrm{CS}\widetilde{\Phi}^r_{*,\Gamma,A},\mathrm{Spec}^\mathrm{CS}\widetilde{\Phi}^I_{*,\Gamma,A},	\\
\end{align}
\begin{align}
\mathrm{Spec}^\mathrm{CS}\breve{\Delta}_{*,\Gamma,A},\breve{\nabla}_{*,\Gamma,A},\mathrm{Spec}^\mathrm{CS}\breve{\Phi}_{*,\Gamma,A},\mathrm{Spec}^\mathrm{CS}\breve{\Delta}^+_{*,\Gamma,A},\mathrm{Spec}^\mathrm{CS}\breve{\nabla}^+_{*,\Gamma,A},\\
\mathrm{Spec}^\mathrm{CS}\breve{\Delta}^\dagger_{*,\Gamma,A},\mathrm{Spec}^\mathrm{CS}\breve{\nabla}^\dagger_{*,\Gamma,A},\mathrm{Spec}^\mathrm{CS}\breve{\Phi}^r_{*,\Gamma,A},\breve{\Phi}^I_{*,\Gamma,A},	\\
\end{align}
\begin{align}
\mathrm{Spec}^\mathrm{CS}{\Delta}_{*,\Gamma,A},\mathrm{Spec}^\mathrm{CS}{\nabla}_{*,\Gamma,A},\mathrm{Spec}^\mathrm{CS}{\Phi}_{*,\Gamma,A},\mathrm{Spec}^\mathrm{CS}{\Delta}^+_{*,\Gamma,A},\mathrm{Spec}^\mathrm{CS}{\nabla}^+_{*,\Gamma,A},\\
\mathrm{Spec}^\mathrm{CS}{\Delta}^\dagger_{*,\Gamma,A},\mathrm{Spec}^\mathrm{CS}{\nabla}^\dagger_{*,\Gamma,A},\mathrm{Spec}^\mathrm{CS}{\Phi}^r_{*,\Gamma,A},\mathrm{Spec}^\mathrm{CS}{\Phi}^I_{*,\Gamma,A}.	
\end{align}

Then we take the corresponding quotients by using the corresponding Frobenius operators:
\begin{align}
&\mathrm{Spec}^\mathrm{CS}\widetilde{\Delta}_{*,\Gamma,A}/\mathrm{Fro}^\mathbb{Z},\mathrm{Spec}^\mathrm{CS}\widetilde{\nabla}_{*,\Gamma,A}/\mathrm{Fro}^\mathbb{Z},\mathrm{Spec}^\mathrm{CS}\widetilde{\Phi}_{*,\Gamma,A}/\mathrm{Fro}^\mathbb{Z},\mathrm{Spec}^\mathrm{CS}\widetilde{\Delta}^+_{*,\Gamma,A}/\mathrm{Fro}^\mathbb{Z},\\
&\mathrm{Spec}^\mathrm{CS}\widetilde{\nabla}^+_{*,\Gamma,A}/\mathrm{Fro}^\mathbb{Z}, \mathrm{Spec}^\mathrm{CS}\widetilde{\Delta}^\dagger_{*,\Gamma,A}/\mathrm{Fro}^\mathbb{Z},\mathrm{Spec}^\mathrm{CS}\widetilde{\nabla}^\dagger_{*,\Gamma,A}/\mathrm{Fro}^\mathbb{Z},	\\
\end{align}
\begin{align}
&\mathrm{Spec}^\mathrm{CS}\breve{\Delta}_{*,\Gamma,A}/\mathrm{Fro}^\mathbb{Z},\breve{\nabla}_{*,\Gamma,A}/\mathrm{Fro}^\mathbb{Z},\mathrm{Spec}^\mathrm{CS}\breve{\Phi}_{*,\Gamma,A}/\mathrm{Fro}^\mathbb{Z},\mathrm{Spec}^\mathrm{CS}\breve{\Delta}^+_{*,\Gamma,A}/\mathrm{Fro}^\mathbb{Z},\\
&\mathrm{Spec}^\mathrm{CS}\breve{\nabla}^+_{*,\Gamma,A}/\mathrm{Fro}^\mathbb{Z}, \mathrm{Spec}^\mathrm{CS}\breve{\Delta}^\dagger_{*,\Gamma,A}/\mathrm{Fro}^\mathbb{Z},\mathrm{Spec}^\mathrm{CS}\breve{\nabla}^\dagger_{*,\Gamma,A}/\mathrm{Fro}^\mathbb{Z},	\\
\end{align}
\begin{align}
&\mathrm{Spec}^\mathrm{CS}{\Delta}_{*,\Gamma,A}/\mathrm{Fro}^\mathbb{Z},\mathrm{Spec}^\mathrm{CS}{\nabla}_{*,\Gamma,A}/\mathrm{Fro}^\mathbb{Z},\mathrm{Spec}^\mathrm{CS}{\Phi}_{*,\Gamma,A}/\mathrm{Fro}^\mathbb{Z},\mathrm{Spec}^\mathrm{CS}{\Delta}^+_{*,\Gamma,A}/\mathrm{Fro}^\mathbb{Z},\\
&\mathrm{Spec}^\mathrm{CS}{\nabla}^+_{*,\Gamma,A}/\mathrm{Fro}^\mathbb{Z}, \mathrm{Spec}^\mathrm{CS}{\Delta}^\dagger_{*,\Gamma,A}/\mathrm{Fro}^\mathbb{Z},\mathrm{Spec}^\mathrm{CS}{\nabla}^\dagger_{*,\Gamma,A}/\mathrm{Fro}^\mathbb{Z}.	
\end{align}
Here for those space with notations related to the radius and the corresponding interval we consider the total unions $\bigcap_r,\bigcup_I$ in order to achieve the whole spaces to achieve the analogues of the corresponding FF curves from \cite{10KL1}, \cite{10KL2}, \cite{10FF} for
\[
\xymatrix@R+0pc@C+0pc{
\underset{r}{\mathrm{homotopylimit}}~\mathrm{Spec}^\mathrm{CS}\widetilde{\Phi}^r_{*,\Gamma,A},\underset{I}{\mathrm{homotopycolimit}}~\mathrm{Spec}^\mathrm{CS}\widetilde{\Phi}^I_{*,\Gamma,A},	\\
}
\]
\[
\xymatrix@R+0pc@C+0pc{
\underset{r}{\mathrm{homotopylimit}}~\mathrm{Spec}^\mathrm{CS}\breve{\Phi}^r_{*,\Gamma,A},\underset{I}{\mathrm{homotopycolimit}}~\mathrm{Spec}^\mathrm{CS}\breve{\Phi}^I_{*,\Gamma,A},	\\
}
\]
\[
\xymatrix@R+0pc@C+0pc{
\underset{r}{\mathrm{homotopylimit}}~\mathrm{Spec}^\mathrm{CS}{\Phi}^r_{*,\Gamma,A},\underset{I}{\mathrm{homotopycolimit}}~\mathrm{Spec}^\mathrm{CS}{\Phi}^I_{*,\Gamma,A}.	
}
\]
\[ 
\xymatrix@R+0pc@C+0pc{
\underset{r}{\mathrm{homotopylimit}}~\mathrm{Spec}^\mathrm{CS}\widetilde{\Phi}^r_{*,\Gamma,A}/\mathrm{Fro}^\mathbb{Z},\underset{I}{\mathrm{homotopycolimit}}~\mathrm{Spec}^\mathrm{CS}\widetilde{\Phi}^I_{*,\Gamma,A}/\mathrm{Fro}^\mathbb{Z},	\\
}
\]
\[ 
\xymatrix@R+0pc@C+0pc{
\underset{r}{\mathrm{homotopylimit}}~\mathrm{Spec}^\mathrm{CS}\breve{\Phi}^r_{*,\Gamma,A}/\mathrm{Fro}^\mathbb{Z},\underset{I}{\mathrm{homotopycolimit}}~\breve{\Phi}^I_{*,\Gamma,A}/\mathrm{Fro}^\mathbb{Z},	\\
}
\]
\[ 
\xymatrix@R+0pc@C+0pc{
\underset{r}{\mathrm{homotopylimit}}~\mathrm{Spec}^\mathrm{CS}{\Phi}^r_{*,\Gamma,A}/\mathrm{Fro}^\mathbb{Z},\underset{I}{\mathrm{homotopycolimit}}~\mathrm{Spec}^\mathrm{CS}{\Phi}^I_{*,\Gamma,A}/\mathrm{Fro}^\mathbb{Z}.	
}
\]

\end{definition}

\

\begin{definition}
We then consider the corresponding quasipresheaves of the corresponding ind-Banach or monomorphic ind-Banach modules from \cite{10BBK}, \cite{10KKM}:
\begin{align}
\mathrm{Quasicoherentpresheaves,IndBanach}_{*}	
\end{align}
where $*$ is one of the following spaces:
\begin{align}
&\mathrm{Spec}^\mathrm{BK}\widetilde{\Phi}_{*,\Gamma,A}/\mathrm{Fro}^\mathbb{Z},	\\
\end{align}
\begin{align}
&\mathrm{Spec}^\mathrm{BK}\breve{\Phi}_{*,\Gamma,A}/\mathrm{Fro}^\mathbb{Z},	\\
\end{align}
\begin{align}
&\mathrm{Spec}^\mathrm{BK}{\Phi}_{*,\Gamma,A}/\mathrm{Fro}^\mathbb{Z}.	
\end{align}
Here for those space without notation related to the radius and the corresponding interval we consider the total unions $\bigcap_r,\bigcup_I$ in order to achieve the whole spaces to achieve the analogues of the corresponding FF curves from \cite{10KL1}, \cite{10KL2}, \cite{10FF} for
\[
\xymatrix@R+0pc@C+0pc{
\underset{r}{\mathrm{homotopylimit}}~\mathrm{Spec}^\mathrm{BK}\widetilde{\Phi}^r_{*,\Gamma,A},\underset{I}{\mathrm{homotopycolimit}}~\mathrm{Spec}^\mathrm{BK}\widetilde{\Phi}^I_{*,\Gamma,A},	\\
}
\]
\[
\xymatrix@R+0pc@C+0pc{
\underset{r}{\mathrm{homotopylimit}}~\mathrm{Spec}^\mathrm{BK}\breve{\Phi}^r_{*,\Gamma,A},\underset{I}{\mathrm{homotopycolimit}}~\mathrm{Spec}^\mathrm{BK}\breve{\Phi}^I_{*,\Gamma,A},	\\
}
\]
\[
\xymatrix@R+0pc@C+0pc{
\underset{r}{\mathrm{homotopylimit}}~\mathrm{Spec}^\mathrm{BK}{\Phi}^r_{*,\Gamma,A},\underset{I}{\mathrm{homotopycolimit}}~\mathrm{Spec}^\mathrm{BK}{\Phi}^I_{*,\Gamma,A}.	
}
\]
\[  
\xymatrix@R+0pc@C+0pc{
\underset{r}{\mathrm{homotopylimit}}~\mathrm{Spec}^\mathrm{BK}\widetilde{\Phi}^r_{*,\Gamma,A}/\mathrm{Fro}^\mathbb{Z},\underset{I}{\mathrm{homotopycolimit}}~\mathrm{Spec}^\mathrm{BK}\widetilde{\Phi}^I_{*,\Gamma,A}/\mathrm{Fro}^\mathbb{Z},	\\
}
\]
\[ 
\xymatrix@R+0pc@C+0pc{
\underset{r}{\mathrm{homotopylimit}}~\mathrm{Spec}^\mathrm{BK}\breve{\Phi}^r_{*,\Gamma,A}/\mathrm{Fro}^\mathbb{Z},\underset{I}{\mathrm{homotopycolimit}}~\mathrm{Spec}^\mathrm{BK}\breve{\Phi}^I_{*,\Gamma,A}/\mathrm{Fro}^\mathbb{Z},	\\
}
\]
\[ 
\xymatrix@R+0pc@C+0pc{
\underset{r}{\mathrm{homotopylimit}}~\mathrm{Spec}^\mathrm{BK}{\Phi}^r_{*,\Gamma,A}/\mathrm{Fro}^\mathbb{Z},\underset{I}{\mathrm{homotopycolimit}}~\mathrm{Spec}^\mathrm{BK}{\Phi}^I_{*,\Gamma,A}/\mathrm{Fro}^\mathbb{Z}.	
}
\]

\end{definition}

\begin{definition}
We then consider the corresponding quasisheaves of the corresponding condensed solid topological modules from \cite{10CS2}:
\begin{align}
\mathrm{Quasicoherentsheaves, Condensed}_{*}	
\end{align}
where $*$ is one of the following spaces:
\begin{align}
&\mathrm{Spec}^\mathrm{CS}\widetilde{\Delta}_{*,\Gamma,A}/\mathrm{Fro}^\mathbb{Z},\mathrm{Spec}^\mathrm{CS}\widetilde{\nabla}_{*,\Gamma,A}/\mathrm{Fro}^\mathbb{Z},\mathrm{Spec}^\mathrm{CS}\widetilde{\Phi}_{*,\Gamma,A}/\mathrm{Fro}^\mathbb{Z},\mathrm{Spec}^\mathrm{CS}\widetilde{\Delta}^+_{*,\Gamma,A}/\mathrm{Fro}^\mathbb{Z},\\
&\mathrm{Spec}^\mathrm{CS}\widetilde{\nabla}^+_{*,\Gamma,A}/\mathrm{Fro}^\mathbb{Z},\mathrm{Spec}^\mathrm{CS}\widetilde{\Delta}^\dagger_{*,\Gamma,A}/\mathrm{Fro}^\mathbb{Z},\mathrm{Spec}^\mathrm{CS}\widetilde{\nabla}^\dagger_{*,\Gamma,A}/\mathrm{Fro}^\mathbb{Z},	\\
\end{align}
\begin{align}
&\mathrm{Spec}^\mathrm{CS}\breve{\Delta}_{*,\Gamma,A}/\mathrm{Fro}^\mathbb{Z},\breve{\nabla}_{*,\Gamma,A}/\mathrm{Fro}^\mathbb{Z},\mathrm{Spec}^\mathrm{CS}\breve{\Phi}_{*,\Gamma,A}/\mathrm{Fro}^\mathbb{Z},\mathrm{Spec}^\mathrm{CS}\breve{\Delta}^+_{*,\Gamma,A}/\mathrm{Fro}^\mathbb{Z},\\
&\mathrm{Spec}^\mathrm{CS}\breve{\nabla}^+_{*,\Gamma,A}/\mathrm{Fro}^\mathbb{Z},\mathrm{Spec}^\mathrm{CS}\breve{\Delta}^\dagger_{*,\Gamma,A}/\mathrm{Fro}^\mathbb{Z},\mathrm{Spec}^\mathrm{CS}\breve{\nabla}^\dagger_{*,\Gamma,A}/\mathrm{Fro}^\mathbb{Z},	\\
\end{align}
\begin{align}
&\mathrm{Spec}^\mathrm{CS}{\Delta}_{*,\Gamma,A}/\mathrm{Fro}^\mathbb{Z},\mathrm{Spec}^\mathrm{CS}{\nabla}_{*,\Gamma,A}/\mathrm{Fro}^\mathbb{Z},\mathrm{Spec}^\mathrm{CS}{\Phi}_{*,\Gamma,A}/\mathrm{Fro}^\mathbb{Z},\mathrm{Spec}^\mathrm{CS}{\Delta}^+_{*,\Gamma,A}/\mathrm{Fro}^\mathbb{Z},\\
&\mathrm{Spec}^\mathrm{CS}{\nabla}^+_{*,\Gamma,A}/\mathrm{Fro}^\mathbb{Z}, \mathrm{Spec}^\mathrm{CS}{\Delta}^\dagger_{*,\Gamma,A}/\mathrm{Fro}^\mathbb{Z},\mathrm{Spec}^\mathrm{CS}{\nabla}^\dagger_{*,\Gamma,A}/\mathrm{Fro}^\mathbb{Z}.	
\end{align}
Here for those space with notations related to the radius and the corresponding interval we consider the total unions $\bigcap_r,\bigcup_I$ in order to achieve the whole spaces to achieve the analogues of the corresponding FF curves from \cite{10KL1}, \cite{10KL2}, \cite{10FF} for
\[
\xymatrix@R+0pc@C+0pc{
\underset{r}{\mathrm{homotopylimit}}~\mathrm{Spec}^\mathrm{CS}\widetilde{\Phi}^r_{*,\Gamma,A},\underset{I}{\mathrm{homotopycolimit}}~\mathrm{Spec}^\mathrm{CS}\widetilde{\Phi}^I_{*,\Gamma,A},	\\
}
\]
\[
\xymatrix@R+0pc@C+0pc{
\underset{r}{\mathrm{homotopylimit}}~\mathrm{Spec}^\mathrm{CS}\breve{\Phi}^r_{*,\Gamma,A},\underset{I}{\mathrm{homotopycolimit}}~\mathrm{Spec}^\mathrm{CS}\breve{\Phi}^I_{*,\Gamma,A},	\\
}
\]
\[
\xymatrix@R+0pc@C+0pc{
\underset{r}{\mathrm{homotopylimit}}~\mathrm{Spec}^\mathrm{CS}{\Phi}^r_{*,\Gamma,A},\underset{I}{\mathrm{homotopycolimit}}~\mathrm{Spec}^\mathrm{CS}{\Phi}^I_{*,\Gamma,A}.	
}
\]
\[ 
\xymatrix@R+0pc@C+0pc{
\underset{r}{\mathrm{homotopylimit}}~\mathrm{Spec}^\mathrm{CS}\widetilde{\Phi}^r_{*,\Gamma,A}/\mathrm{Fro}^\mathbb{Z},\underset{I}{\mathrm{homotopycolimit}}~\mathrm{Spec}^\mathrm{CS}\widetilde{\Phi}^I_{*,\Gamma,A}/\mathrm{Fro}^\mathbb{Z},	\\
}
\]
\[ 
\xymatrix@R+0pc@C+0pc{
\underset{r}{\mathrm{homotopylimit}}~\mathrm{Spec}^\mathrm{CS}\breve{\Phi}^r_{*,\Gamma,A}/\mathrm{Fro}^\mathbb{Z},\underset{I}{\mathrm{homotopycolimit}}~\breve{\Phi}^I_{*,\Gamma,A}/\mathrm{Fro}^\mathbb{Z},	\\
}
\]
\[ 
\xymatrix@R+0pc@C+0pc{
\underset{r}{\mathrm{homotopylimit}}~\mathrm{Spec}^\mathrm{CS}{\Phi}^r_{*,\Gamma,A}/\mathrm{Fro}^\mathbb{Z},\underset{I}{\mathrm{homotopycolimit}}~\mathrm{Spec}^\mathrm{CS}{\Phi}^I_{*,\Gamma,A}/\mathrm{Fro}^\mathbb{Z}.	
}
\]

\end{definition}

\

\begin{proposition}
There is a well-defined functor from the $\infty$-category 
\begin{align}
\mathrm{Quasicoherentpresheaves,Condensed}_{*}	
\end{align}
where $*$ is one of the following spaces:
\begin{align}
&\mathrm{Spec}^\mathrm{CS}\widetilde{\Phi}_{*,\Gamma,A}/\mathrm{Fro}^\mathbb{Z},	\\
\end{align}
\begin{align}
&\mathrm{Spec}^\mathrm{CS}\breve{\Phi}_{*,\Gamma,A}/\mathrm{Fro}^\mathbb{Z},	\\
\end{align}
\begin{align}
&\mathrm{Spec}^\mathrm{CS}{\Phi}_{*,\Gamma,A}/\mathrm{Fro}^\mathbb{Z},	
\end{align}
to the $\infty$-category of $\mathrm{Fro}$-equivariant quasicoherent presheaves over similar spaces above correspondingly without the $\mathrm{Fro}$-quotients, and to the $\infty$-category of $\mathrm{Fro}$-equivariant quasicoherent modules over global sections of the structure $\infty$-sheaves of the similar spaces above correspondingly without the $\mathrm{Fro}$-quotients. Here for those space without notation related to the radius and the corresponding interval we consider the total unions $\bigcap_r,\bigcup_I$ in order to achieve the whole spaces to achieve the analogues of the corresponding FF curves from \cite{10KL1}, \cite{10KL2}, \cite{10FF} for
\[
\xymatrix@R+0pc@C+0pc{
\underset{r}{\mathrm{homotopylimit}}~\mathrm{Spec}^\mathrm{CS}\widetilde{\Phi}^r_{*,\Gamma,A},\underset{I}{\mathrm{homotopycolimit}}~\mathrm{Spec}^\mathrm{CS}\widetilde{\Phi}^I_{*,\Gamma,A},	\\
}
\]
\[
\xymatrix@R+0pc@C+0pc{
\underset{r}{\mathrm{homotopylimit}}~\mathrm{Spec}^\mathrm{CS}\breve{\Phi}^r_{*,\Gamma,A},\underset{I}{\mathrm{homotopycolimit}}~\mathrm{Spec}^\mathrm{CS}\breve{\Phi}^I_{*,\Gamma,A},	\\
}
\]
\[
\xymatrix@R+0pc@C+0pc{
\underset{r}{\mathrm{homotopylimit}}~\mathrm{Spec}^\mathrm{CS}{\Phi}^r_{*,\Gamma,A},\underset{I}{\mathrm{homotopycolimit}}~\mathrm{Spec}^\mathrm{CS}{\Phi}^I_{*,\Gamma,A}.	
}
\]
\[ 
\xymatrix@R+0pc@C+0pc{
\underset{r}{\mathrm{homotopylimit}}~\mathrm{Spec}^\mathrm{CS}\widetilde{\Phi}^r_{*,\Gamma,A}/\mathrm{Fro}^\mathbb{Z},\underset{I}{\mathrm{homotopycolimit}}~\mathrm{Spec}^\mathrm{CS}\widetilde{\Phi}^I_{*,\Gamma,A}/\mathrm{Fro}^\mathbb{Z},	\\
}
\]
\[ 
\xymatrix@R+0pc@C+0pc{
\underset{r}{\mathrm{homotopylimit}}~\mathrm{Spec}^\mathrm{CS}\breve{\Phi}^r_{*,\Gamma,A}/\mathrm{Fro}^\mathbb{Z},\underset{I}{\mathrm{homotopycolimit}}~\breve{\Phi}^I_{*,\Gamma,A}/\mathrm{Fro}^\mathbb{Z},	\\
}
\]
\[ 
\xymatrix@R+0pc@C+0pc{
\underset{r}{\mathrm{homotopylimit}}~\mathrm{Spec}^\mathrm{CS}{\Phi}^r_{*,\Gamma,A}/\mathrm{Fro}^\mathbb{Z},\underset{I}{\mathrm{homotopycolimit}}~\mathrm{Spec}^\mathrm{CS}{\Phi}^I_{*,\Gamma,A}/\mathrm{Fro}^\mathbb{Z}.	
}
\]	
In this situation we will have the target category being family parametrized by $r$ or $I$ in compatible glueing sense as in \cite[Definition 5.4.10]{10KL2}. In this situation for modules parametrized by the intervals we have the equivalence of $\infty$-categories by using \cite[Proposition 13.8]{10CS2}. Here the corresponding quasicoherent Frobenius modules are defined to be the corresponding homotopy colimits and limits of Frobenius modules:
\begin{align}
\underset{r}{\mathrm{homotopycolimit}}~M_r,\\
\underset{I}{\mathrm{homotopylimit}}~M_I,	
\end{align}
where each $M_r$ is a Frobenius-equivariant module over the period ring with respect to some radius $r$ while each $M_I$ is a Frobenius-equivariant module over the period ring with respect to some interval $I$.\\
\end{proposition}

\begin{proposition}
Similar proposition holds for 
\begin{align}
\mathrm{Quasicoherentsheaves,IndBanach}_{*}.	
\end{align}	
\end{proposition}

\

\begin{definition}
We then consider the corresponding quasipresheaves of perfect complexes the corresponding ind-Banach or monomorphic ind-Banach modules from \cite{10BBK}, \cite{10KKM}:
\begin{align}
\mathrm{Quasicoherentpresheaves,Perfectcomplex,IndBanach}_{*}	
\end{align}
where $*$ is one of the following spaces:
\begin{align}
&\mathrm{Spec}^\mathrm{BK}\widetilde{\Phi}_{*,\Gamma,A}/\mathrm{Fro}^\mathbb{Z},	\\
\end{align}
\begin{align}
&\mathrm{Spec}^\mathrm{BK}\breve{\Phi}_{*,\Gamma,A}/\mathrm{Fro}^\mathbb{Z},	\\
\end{align}
\begin{align}
&\mathrm{Spec}^\mathrm{BK}{\Phi}_{*,\Gamma,A}/\mathrm{Fro}^\mathbb{Z}.	
\end{align}
Here for those space without notation related to the radius and the corresponding interval we consider the total unions $\bigcap_r,\bigcup_I$ in order to achieve the whole spaces to achieve the analogues of the corresponding FF curves from \cite{10KL1}, \cite{10KL2}, \cite{10FF} for
\[
\xymatrix@R+0pc@C+0pc{
\underset{r}{\mathrm{homotopylimit}}~\mathrm{Spec}^\mathrm{BK}\widetilde{\Phi}^r_{*,\Gamma,A},\underset{I}{\mathrm{homotopycolimit}}~\mathrm{Spec}^\mathrm{BK}\widetilde{\Phi}^I_{*,\Gamma,A},	\\
}
\]
\[
\xymatrix@R+0pc@C+0pc{
\underset{r}{\mathrm{homotopylimit}}~\mathrm{Spec}^\mathrm{BK}\breve{\Phi}^r_{*,\Gamma,A},\underset{I}{\mathrm{homotopycolimit}}~\mathrm{Spec}^\mathrm{BK}\breve{\Phi}^I_{*,\Gamma,A},	\\
}
\]
\[
\xymatrix@R+0pc@C+0pc{
\underset{r}{\mathrm{homotopylimit}}~\mathrm{Spec}^\mathrm{BK}{\Phi}^r_{*,\Gamma,A},\underset{I}{\mathrm{homotopycolimit}}~\mathrm{Spec}^\mathrm{BK}{\Phi}^I_{*,\Gamma,A}.	
}
\]
\[  
\xymatrix@R+0pc@C+0pc{
\underset{r}{\mathrm{homotopylimit}}~\mathrm{Spec}^\mathrm{BK}\widetilde{\Phi}^r_{*,\Gamma,A}/\mathrm{Fro}^\mathbb{Z},\underset{I}{\mathrm{homotopycolimit}}~\mathrm{Spec}^\mathrm{BK}\widetilde{\Phi}^I_{*,\Gamma,A}/\mathrm{Fro}^\mathbb{Z},	\\
}
\]
\[ 
\xymatrix@R+0pc@C+0pc{
\underset{r}{\mathrm{homotopylimit}}~\mathrm{Spec}^\mathrm{BK}\breve{\Phi}^r_{*,\Gamma,A}/\mathrm{Fro}^\mathbb{Z},\underset{I}{\mathrm{homotopycolimit}}~\mathrm{Spec}^\mathrm{BK}\breve{\Phi}^I_{*,\Gamma,A}/\mathrm{Fro}^\mathbb{Z},	\\
}
\]
\[ 
\xymatrix@R+0pc@C+0pc{
\underset{r}{\mathrm{homotopylimit}}~\mathrm{Spec}^\mathrm{BK}{\Phi}^r_{*,\Gamma,A}/\mathrm{Fro}^\mathbb{Z},\underset{I}{\mathrm{homotopycolimit}}~\mathrm{Spec}^\mathrm{BK}{\Phi}^I_{*,\Gamma,A}/\mathrm{Fro}^\mathbb{Z}.	
}
\]

\end{definition}

\begin{definition}
We then consider the corresponding quasisheaves of perfect complexes of the corresponding condensed solid topological modules from \cite{10CS2}:
\begin{align}
\mathrm{Quasicoherentsheaves, Perfectcomplex, Condensed}_{*}	
\end{align}
where $*$ is one of the following spaces:
\begin{align}
&\mathrm{Spec}^\mathrm{CS}\widetilde{\Delta}_{*,\Gamma,A}/\mathrm{Fro}^\mathbb{Z},\mathrm{Spec}^\mathrm{CS}\widetilde{\nabla}_{*,\Gamma,A}/\mathrm{Fro}^\mathbb{Z},\mathrm{Spec}^\mathrm{CS}\widetilde{\Phi}_{*,\Gamma,A}/\mathrm{Fro}^\mathbb{Z},\mathrm{Spec}^\mathrm{CS}\widetilde{\Delta}^+_{*,\Gamma,A}/\mathrm{Fro}^\mathbb{Z},\\
&\mathrm{Spec}^\mathrm{CS}\widetilde{\nabla}^+_{*,\Gamma,A}/\mathrm{Fro}^\mathbb{Z},\mathrm{Spec}^\mathrm{CS}\widetilde{\Delta}^\dagger_{*,\Gamma,A}/\mathrm{Fro}^\mathbb{Z},\mathrm{Spec}^\mathrm{CS}\widetilde{\nabla}^\dagger_{*,\Gamma,A}/\mathrm{Fro}^\mathbb{Z},	\\
\end{align}
\begin{align}
&\mathrm{Spec}^\mathrm{CS}\breve{\Delta}_{*,\Gamma,A}/\mathrm{Fro}^\mathbb{Z},\breve{\nabla}_{*,\Gamma,A}/\mathrm{Fro}^\mathbb{Z},\mathrm{Spec}^\mathrm{CS}\breve{\Phi}_{*,\Gamma,A}/\mathrm{Fro}^\mathbb{Z},\mathrm{Spec}^\mathrm{CS}\breve{\Delta}^+_{*,\Gamma,A}/\mathrm{Fro}^\mathbb{Z},\\
&\mathrm{Spec}^\mathrm{CS}\breve{\nabla}^+_{*,\Gamma,A}/\mathrm{Fro}^\mathbb{Z},\mathrm{Spec}^\mathrm{CS}\breve{\Delta}^\dagger_{*,\Gamma,A}/\mathrm{Fro}^\mathbb{Z},\mathrm{Spec}^\mathrm{CS}\breve{\nabla}^\dagger_{*,\Gamma,A}/\mathrm{Fro}^\mathbb{Z},	\\
\end{align}
\begin{align}
&\mathrm{Spec}^\mathrm{CS}{\Delta}_{*,\Gamma,A}/\mathrm{Fro}^\mathbb{Z},\mathrm{Spec}^\mathrm{CS}{\nabla}_{*,\Gamma,A}/\mathrm{Fro}^\mathbb{Z},\mathrm{Spec}^\mathrm{CS}{\Phi}_{*,\Gamma,A}/\mathrm{Fro}^\mathbb{Z},\mathrm{Spec}^\mathrm{CS}{\Delta}^+_{*,\Gamma,A}/\mathrm{Fro}^\mathbb{Z},\\
&\mathrm{Spec}^\mathrm{CS}{\nabla}^+_{*,\Gamma,A}/\mathrm{Fro}^\mathbb{Z}, \mathrm{Spec}^\mathrm{CS}{\Delta}^\dagger_{*,\Gamma,A}/\mathrm{Fro}^\mathbb{Z},\mathrm{Spec}^\mathrm{CS}{\nabla}^\dagger_{*,\Gamma,A}/\mathrm{Fro}^\mathbb{Z}.	
\end{align}
Here for those space with notations related to the radius and the corresponding interval we consider the total unions $\bigcap_r,\bigcup_I$ in order to achieve the whole spaces to achieve the analogues of the corresponding FF curves from \cite{10KL1}, \cite{10KL2}, \cite{10FF} for
\[
\xymatrix@R+0pc@C+0pc{
\underset{r}{\mathrm{homotopylimit}}~\mathrm{Spec}^\mathrm{CS}\widetilde{\Phi}^r_{*,\Gamma,A},\underset{I}{\mathrm{homotopycolimit}}~\mathrm{Spec}^\mathrm{CS}\widetilde{\Phi}^I_{*,\Gamma,A},	\\
}
\]
\[
\xymatrix@R+0pc@C+0pc{
\underset{r}{\mathrm{homotopylimit}}~\mathrm{Spec}^\mathrm{CS}\breve{\Phi}^r_{*,\Gamma,A},\underset{I}{\mathrm{homotopycolimit}}~\mathrm{Spec}^\mathrm{CS}\breve{\Phi}^I_{*,\Gamma,A},	\\
}
\]
\[
\xymatrix@R+0pc@C+0pc{
\underset{r}{\mathrm{homotopylimit}}~\mathrm{Spec}^\mathrm{CS}{\Phi}^r_{*,\Gamma,A},\underset{I}{\mathrm{homotopycolimit}}~\mathrm{Spec}^\mathrm{CS}{\Phi}^I_{*,\Gamma,A}.	
}
\]
\[ 
\xymatrix@R+0pc@C+0pc{
\underset{r}{\mathrm{homotopylimit}}~\mathrm{Spec}^\mathrm{CS}\widetilde{\Phi}^r_{*,\Gamma,A}/\mathrm{Fro}^\mathbb{Z},\underset{I}{\mathrm{homotopycolimit}}~\mathrm{Spec}^\mathrm{CS}\widetilde{\Phi}^I_{*,\Gamma,A}/\mathrm{Fro}^\mathbb{Z},	\\
}
\]
\[ 
\xymatrix@R+0pc@C+0pc{
\underset{r}{\mathrm{homotopylimit}}~\mathrm{Spec}^\mathrm{CS}\breve{\Phi}^r_{*,\Gamma,A}/\mathrm{Fro}^\mathbb{Z},\underset{I}{\mathrm{homotopycolimit}}~\breve{\Phi}^I_{*,\Gamma,A}/\mathrm{Fro}^\mathbb{Z},	\\
}
\]
\[ 
\xymatrix@R+0pc@C+0pc{
\underset{r}{\mathrm{homotopylimit}}~\mathrm{Spec}^\mathrm{CS}{\Phi}^r_{*,\Gamma,A}/\mathrm{Fro}^\mathbb{Z},\underset{I}{\mathrm{homotopycolimit}}~\mathrm{Spec}^\mathrm{CS}{\Phi}^I_{*,\Gamma,A}/\mathrm{Fro}^\mathbb{Z}.	
}
\]

\end{definition}

\begin{proposition}
There is a well-defined functor from the $\infty$-category 
\begin{align}
\mathrm{Quasicoherentpresheaves,Perfectcomplex,Condensed}_{*}	
\end{align}
where $*$ is one of the following spaces:
\begin{align}
&\mathrm{Spec}^\mathrm{CS}\widetilde{\Phi}_{*,\Gamma,A}/\mathrm{Fro}^\mathbb{Z},	\\
\end{align}
\begin{align}
&\mathrm{Spec}^\mathrm{CS}\breve{\Phi}_{*,\Gamma,A}/\mathrm{Fro}^\mathbb{Z},	\\
\end{align}
\begin{align}
&\mathrm{Spec}^\mathrm{CS}{\Phi}_{*,\Gamma,A}/\mathrm{Fro}^\mathbb{Z},	
\end{align}
to the $\infty$-category of $\mathrm{Fro}$-equivariant quasicoherent presheaves over similar spaces above correspondingly without the $\mathrm{Fro}$-quotients, and to the $\infty$-category of $\mathrm{Fro}$-equivariant quasicoherent modules over global sections of the structure $\infty$-sheaves of the similar spaces above correspondingly without the $\mathrm{Fro}$-quotients. Here for those space without notation related to the radius and the corresponding interval we consider the total unions $\bigcap_r,\bigcup_I$ in order to achieve the whole spaces to achieve the analogues of the corresponding FF curves from \cite{10KL1}, \cite{10KL2}, \cite{10FF} for
\[
\xymatrix@R+0pc@C+0pc{
\underset{r}{\mathrm{homotopylimit}}~\mathrm{Spec}^\mathrm{CS}\widetilde{\Phi}^r_{*,\Gamma,A},\underset{I}{\mathrm{homotopycolimit}}~\mathrm{Spec}^\mathrm{CS}\widetilde{\Phi}^I_{*,\Gamma,A},	\\
}
\]
\[
\xymatrix@R+0pc@C+0pc{
\underset{r}{\mathrm{homotopylimit}}~\mathrm{Spec}^\mathrm{CS}\breve{\Phi}^r_{*,\Gamma,A},\underset{I}{\mathrm{homotopycolimit}}~\mathrm{Spec}^\mathrm{CS}\breve{\Phi}^I_{*,\Gamma,A},	\\
}
\]
\[
\xymatrix@R+0pc@C+0pc{
\underset{r}{\mathrm{homotopylimit}}~\mathrm{Spec}^\mathrm{CS}{\Phi}^r_{*,\Gamma,A},\underset{I}{\mathrm{homotopycolimit}}~\mathrm{Spec}^\mathrm{CS}{\Phi}^I_{*,\Gamma,A}.	
}
\]
\[ 
\xymatrix@R+0pc@C+0pc{
\underset{r}{\mathrm{homotopylimit}}~\mathrm{Spec}^\mathrm{CS}\widetilde{\Phi}^r_{*,\Gamma,A}/\mathrm{Fro}^\mathbb{Z},\underset{I}{\mathrm{homotopycolimit}}~\mathrm{Spec}^\mathrm{CS}\widetilde{\Phi}^I_{*,\Gamma,A}/\mathrm{Fro}^\mathbb{Z},	\\
}
\]
\[ 
\xymatrix@R+0pc@C+0pc{
\underset{r}{\mathrm{homotopylimit}}~\mathrm{Spec}^\mathrm{CS}\breve{\Phi}^r_{*,\Gamma,A}/\mathrm{Fro}^\mathbb{Z},\underset{I}{\mathrm{homotopycolimit}}~\breve{\Phi}^I_{*,\Gamma,A}/\mathrm{Fro}^\mathbb{Z},	\\
}
\]
\[ 
\xymatrix@R+0pc@C+0pc{
\underset{r}{\mathrm{homotopylimit}}~\mathrm{Spec}^\mathrm{CS}{\Phi}^r_{*,\Gamma,A}/\mathrm{Fro}^\mathbb{Z},\underset{I}{\mathrm{homotopycolimit}}~\mathrm{Spec}^\mathrm{CS}{\Phi}^I_{*,\Gamma,A}/\mathrm{Fro}^\mathbb{Z}.	
}
\]	
In this situation we will have the target category being family parametrized by $r$ or $I$ in compatible glueing sense as in \cite[Definition 5.4.10]{10KL2}. In this situation for modules parametrized by the intervals we have the equivalence of $\infty$-categories by using \cite[Proposition 12.18]{10CS2}. Here the corresponding quasicoherent Frobenius modules are defined to be the corresponding homotopy colimits and limits of Frobenius modules:
\begin{align}
\underset{r}{\mathrm{homotopycolimit}}~M_r,\\
\underset{I}{\mathrm{homotopylimit}}~M_I,	
\end{align}
where each $M_r$ is a Frobenius-equivariant module over the period ring with respect to some radius $r$ while each $M_I$ is a Frobenius-equivariant module over the period ring with respect to some interval $I$.\\
\end{proposition}

\begin{proposition}
Similar proposition holds for 
\begin{align}
\mathrm{Quasicoherentsheaves,Perfectcomplex,IndBanach}_{*}.	
\end{align}	
\end{proposition}

\newpage
\subsection{Frobenius Quasicoherent Prestacks II: Deformation in Banach Rings}

\begin{definition}
We now consider the pro-\'etale site of $\mathrm{Spa}\mathbb{Q}_p\left<X_1^{\pm 1},...,X_k^{\pm 1}\right>$, denote that by $*$. To be more accurate we replace one component for $\Gamma$ with the pro-\'etale site of $\mathrm{Spa}\mathbb{Q}_p\left<X_1^{\pm 1},...,X_k^{\pm 1}\right>$. And we treat then all the functor to be prestacks for this site. Then from \cite{10KL1} and \cite[Definition 5.2.1]{10KL2} we have the following class of Kedlaya-Liu rings (with the following replacement: $\Delta$ stands for $A$, $\nabla$ stands for $B$, while $\Phi$ stands for $C$) by taking product in the sense of self $\Gamma$-th power:

\[
\xymatrix@R+0pc@C+0pc{
\widetilde{\Delta}_{*,\Gamma},\widetilde{\nabla}_{*,\Gamma},\widetilde{\Phi}_{*,\Gamma},\widetilde{\Delta}^+_{*,\Gamma},\widetilde{\nabla}^+_{*,\Gamma},\widetilde{\Delta}^\dagger_{*,\Gamma},\widetilde{\nabla}^\dagger_{*,\Gamma},\widetilde{\Phi}^r_{*,\Gamma},\widetilde{\Phi}^I_{*,\Gamma}, 
}
\]

\[
\xymatrix@R+0pc@C+0pc{
\breve{\Delta}_{*,\Gamma},\breve{\nabla}_{*,\Gamma},\breve{\Phi}_{*,\Gamma},\breve{\Delta}^+_{*,\Gamma},\breve{\nabla}^+_{*,\Gamma},\breve{\Delta}^\dagger_{*,\Gamma},\breve{\nabla}^\dagger_{*,\Gamma},\breve{\Phi}^r_{*,\Gamma},\breve{\Phi}^I_{*,\Gamma},	
}
\]

\[
\xymatrix@R+0pc@C+0pc{
{\Delta}_{*,\Gamma},{\nabla}_{*,\Gamma},{\Phi}_{*,\Gamma},{\Delta}^+_{*,\Gamma},{\nabla}^+_{*,\Gamma},{\Delta}^\dagger_{*,\Gamma},{\nabla}^\dagger_{*,\Gamma},{\Phi}^r_{*,\Gamma},{\Phi}^I_{*,\Gamma}.	
}
\]
We now consider the following rings with $-$ being any deforming Banach ring over $\mathbb{Q}_p$. Taking the product we have:
\[
\xymatrix@R+0pc@C+0pc{
\widetilde{\Phi}_{*,\Gamma,-},\widetilde{\Phi}^r_{*,\Gamma,-},\widetilde{\Phi}^I_{*,\Gamma,-},	
}
\]
\[
\xymatrix@R+0pc@C+0pc{
\breve{\Phi}_{*,\Gamma,-},\breve{\Phi}^r_{*,\Gamma,-},\breve{\Phi}^I_{*,\Gamma,-},	
}
\]
\[
\xymatrix@R+0pc@C+0pc{
{\Phi}_{*,\Gamma,-},{\Phi}^r_{*,\Gamma,-},{\Phi}^I_{*,\Gamma,-}.	
}
\]
They carry multi Frobenius action $\varphi_\Gamma$ and multi $\mathrm{Lie}_\Gamma:=\mathbb{Z}_p^{\times\Gamma}$ action. In our current situation after \cite{10CKZ} and \cite{10PZ} we consider the following $(\infty,1)$-categories of $(\infty,1)$-modules.\\
\end{definition}

\begin{definition}
First we consider the Bambozzi-Kremnizer spectrum $\mathrm{Spec}^\mathrm{BK}(*)$ attached to any of those in the above from \cite{10BK} by taking derived rational localization:
\begin{align}
&\mathrm{Spec}^\mathrm{BK}\widetilde{\Phi}_{*,\Gamma,-},\mathrm{Spec}^\mathrm{BK}\widetilde{\Phi}^r_{*,\Gamma,-},\mathrm{Spec}^\mathrm{BK}\widetilde{\Phi}^I_{*,\Gamma,-},	
\end{align}
\begin{align}
&\mathrm{Spec}^\mathrm{BK}\breve{\Phi}_{*,\Gamma,-},\mathrm{Spec}^\mathrm{BK}\breve{\Phi}^r_{*,\Gamma,-},\mathrm{Spec}^\mathrm{BK}\breve{\Phi}^I_{*,\Gamma,-},	
\end{align}
\begin{align}
&\mathrm{Spec}^\mathrm{BK}{\Phi}_{*,\Gamma,-},
\mathrm{Spec}^\mathrm{BK}{\Phi}^r_{*,\Gamma,-},\mathrm{Spec}^\mathrm{BK}{\Phi}^I_{*,\Gamma,-}.	
\end{align}

Then we take the corresponding quotients by using the corresponding Frobenius operators:
\begin{align}
&\mathrm{Spec}^\mathrm{BK}\widetilde{\Phi}_{*,\Gamma,-}/\mathrm{Fro}^\mathbb{Z},	\\
\end{align}
\begin{align}
&\mathrm{Spec}^\mathrm{BK}\breve{\Phi}_{*,\Gamma,-}/\mathrm{Fro}^\mathbb{Z},	\\
\end{align}
\begin{align}
&\mathrm{Spec}^\mathrm{BK}{\Phi}_{*,\Gamma,-}/\mathrm{Fro}^\mathbb{Z}.	
\end{align}
Here for those space without notation related to the radius and the corresponding interval we consider the total unions $\bigcap_r,\bigcup_I$ in order to achieve the whole spaces to achieve the analogues of the corresponding FF curves from \cite{10KL1}, \cite{10KL2}, \cite{10FF} for
\[
\xymatrix@R+0pc@C+0pc{
\underset{r}{\mathrm{homotopylimit}}~\mathrm{Spec}^\mathrm{BK}\widetilde{\Phi}^r_{*,\Gamma,-},\underset{I}{\mathrm{homotopycolimit}}~\mathrm{Spec}^\mathrm{BK}\widetilde{\Phi}^I_{*,\Gamma,-},	\\
}
\]
\[
\xymatrix@R+0pc@C+0pc{
\underset{r}{\mathrm{homotopylimit}}~\mathrm{Spec}^\mathrm{BK}\breve{\Phi}^r_{*,\Gamma,-},\underset{I}{\mathrm{homotopycolimit}}~\mathrm{Spec}^\mathrm{BK}\breve{\Phi}^I_{*,\Gamma,-},	\\
}
\]
\[
\xymatrix@R+0pc@C+0pc{
\underset{r}{\mathrm{homotopylimit}}~\mathrm{Spec}^\mathrm{BK}{\Phi}^r_{*,\Gamma,-},\underset{I}{\mathrm{homotopycolimit}}~\mathrm{Spec}^\mathrm{BK}{\Phi}^I_{*,\Gamma,-}.	
}
\]
\[  
\xymatrix@R+0pc@C+0pc{
\underset{r}{\mathrm{homotopylimit}}~\mathrm{Spec}^\mathrm{BK}\widetilde{\Phi}^r_{*,\Gamma,-}/\mathrm{Fro}^\mathbb{Z},\underset{I}{\mathrm{homotopycolimit}}~\mathrm{Spec}^\mathrm{BK}\widetilde{\Phi}^I_{*,\Gamma,-}/\mathrm{Fro}^\mathbb{Z},	\\
}
\]
\[ 
\xymatrix@R+0pc@C+0pc{
\underset{r}{\mathrm{homotopylimit}}~\mathrm{Spec}^\mathrm{BK}\breve{\Phi}^r_{*,\Gamma,-}/\mathrm{Fro}^\mathbb{Z},\underset{I}{\mathrm{homotopycolimit}}~\mathrm{Spec}^\mathrm{BK}\breve{\Phi}^I_{*,\Gamma,-}/\mathrm{Fro}^\mathbb{Z},	\\
}
\]
\[ 
\xymatrix@R+0pc@C+0pc{
\underset{r}{\mathrm{homotopylimit}}~\mathrm{Spec}^\mathrm{BK}{\Phi}^r_{*,\Gamma,-}/\mathrm{Fro}^\mathbb{Z},\underset{I}{\mathrm{homotopycolimit}}~\mathrm{Spec}^\mathrm{BK}{\Phi}^I_{*,\Gamma,-}/\mathrm{Fro}^\mathbb{Z}.	
}
\]

\end{definition}

\indent Meanwhile we have the corresponding Clausen-Scholze analytic stacks from \cite{10CS2}, therefore applying their construction we have:

\begin{definition}
Here we define the following products by using the solidified tensor product from \cite{10CS1} and \cite{10CS2}. Namely $A$ will still as above as a Banach ring over $\mathbb{Q}_p$. Then we take solidified tensor product $\overset{\blacksquare}{\otimes}$ of any of the following
\[
\xymatrix@R+0pc@C+0pc{
\widetilde{\Delta}_{*,\Gamma},\widetilde{\nabla}_{*,\Gamma},\widetilde{\Phi}_{*,\Gamma},\widetilde{\Delta}^+_{*,\Gamma},\widetilde{\nabla}^+_{*,\Gamma},\widetilde{\Delta}^\dagger_{*,\Gamma},\widetilde{\nabla}^\dagger_{*,\Gamma},\widetilde{\Phi}^r_{*,\Gamma},\widetilde{\Phi}^I_{*,\Gamma}, 
}
\]

\[
\xymatrix@R+0pc@C+0pc{
\breve{\Delta}_{*,\Gamma},\breve{\nabla}_{*,\Gamma},\breve{\Phi}_{*,\Gamma},\breve{\Delta}^+_{*,\Gamma},\breve{\nabla}^+_{*,\Gamma},\breve{\Delta}^\dagger_{*,\Gamma},\breve{\nabla}^\dagger_{*,\Gamma},\breve{\Phi}^r_{*,\Gamma},\breve{\Phi}^I_{*,\Gamma},	
}
\]

\[
\xymatrix@R+0pc@C+0pc{
{\Delta}_{*,\Gamma},{\nabla}_{*,\Gamma},{\Phi}_{*,\Gamma},{\Delta}^+_{*,\Gamma},{\nabla}^+_{*,\Gamma},{\Delta}^\dagger_{*,\Gamma},{\nabla}^\dagger_{*,\Gamma},{\Phi}^r_{*,\Gamma},{\Phi}^I_{*,\Gamma},	
}
\]  	
with $A$. Then we have the notations:
\[
\xymatrix@R+0pc@C+0pc{
\widetilde{\Delta}_{*,\Gamma,-},\widetilde{\nabla}_{*,\Gamma,-},\widetilde{\Phi}_{*,\Gamma,-},\widetilde{\Delta}^+_{*,\Gamma,-},\widetilde{\nabla}^+_{*,\Gamma,-},\widetilde{\Delta}^\dagger_{*,\Gamma,-},\widetilde{\nabla}^\dagger_{*,\Gamma,-},\widetilde{\Phi}^r_{*,\Gamma,-},\widetilde{\Phi}^I_{*,\Gamma,-}, 
}
\]

\[
\xymatrix@R+0pc@C+0pc{
\breve{\Delta}_{*,\Gamma,-},\breve{\nabla}_{*,\Gamma,-},\breve{\Phi}_{*,\Gamma,-},\breve{\Delta}^+_{*,\Gamma,-},\breve{\nabla}^+_{*,\Gamma,-},\breve{\Delta}^\dagger_{*,\Gamma,-},\breve{\nabla}^\dagger_{*,\Gamma,-},\breve{\Phi}^r_{*,\Gamma,-},\breve{\Phi}^I_{*,\Gamma,-},	
}
\]

\[
\xymatrix@R+0pc@C+0pc{
{\Delta}_{*,\Gamma,-},{\nabla}_{*,\Gamma,-},{\Phi}_{*,\Gamma,-},{\Delta}^+_{*,\Gamma,-},{\nabla}^+_{*,\Gamma,-},{\Delta}^\dagger_{*,\Gamma,-},{\nabla}^\dagger_{*,\Gamma,-},{\Phi}^r_{*,\Gamma,-},{\Phi}^I_{*,\Gamma,-}.	
}
\]
\end{definition}

\begin{definition}
First we consider the Clausen-Scholze spectrum $\mathrm{Spec}^\mathrm{CS}(*)$ attached to any of those in the above from \cite{10CS2} by taking derived rational localization:
\begin{align}
\mathrm{Spec}^\mathrm{CS}\widetilde{\Delta}_{*,\Gamma,-},\mathrm{Spec}^\mathrm{CS}\widetilde{\nabla}_{*,\Gamma,-},\mathrm{Spec}^\mathrm{CS}\widetilde{\Phi}_{*,\Gamma,-},\mathrm{Spec}^\mathrm{CS}\widetilde{\Delta}^+_{*,\Gamma,-},\mathrm{Spec}^\mathrm{CS}\widetilde{\nabla}^+_{*,\Gamma,-},\\
\mathrm{Spec}^\mathrm{CS}\widetilde{\Delta}^\dagger_{*,\Gamma,-},\mathrm{Spec}^\mathrm{CS}\widetilde{\nabla}^\dagger_{*,\Gamma,-},\mathrm{Spec}^\mathrm{CS}\widetilde{\Phi}^r_{*,\Gamma,-},\mathrm{Spec}^\mathrm{CS}\widetilde{\Phi}^I_{*,\Gamma,-},	\\
\end{align}
\begin{align}
\mathrm{Spec}^\mathrm{CS}\breve{\Delta}_{*,\Gamma,-},\breve{\nabla}_{*,\Gamma,-},\mathrm{Spec}^\mathrm{CS}\breve{\Phi}_{*,\Gamma,-},\mathrm{Spec}^\mathrm{CS}\breve{\Delta}^+_{*,\Gamma,-},\mathrm{Spec}^\mathrm{CS}\breve{\nabla}^+_{*,\Gamma,-},\\
\mathrm{Spec}^\mathrm{CS}\breve{\Delta}^\dagger_{*,\Gamma,-},\mathrm{Spec}^\mathrm{CS}\breve{\nabla}^\dagger_{*,\Gamma,-},\mathrm{Spec}^\mathrm{CS}\breve{\Phi}^r_{*,\Gamma,-},\breve{\Phi}^I_{*,\Gamma,-},	\\
\end{align}
\begin{align}
\mathrm{Spec}^\mathrm{CS}{\Delta}_{*,\Gamma,-},\mathrm{Spec}^\mathrm{CS}{\nabla}_{*,\Gamma,-},\mathrm{Spec}^\mathrm{CS}{\Phi}_{*,\Gamma,-},\mathrm{Spec}^\mathrm{CS}{\Delta}^+_{*,\Gamma,-},\mathrm{Spec}^\mathrm{CS}{\nabla}^+_{*,\Gamma,-},\\
\mathrm{Spec}^\mathrm{CS}{\Delta}^\dagger_{*,\Gamma,-},\mathrm{Spec}^\mathrm{CS}{\nabla}^\dagger_{*,\Gamma,-},\mathrm{Spec}^\mathrm{CS}{\Phi}^r_{*,\Gamma,-},\mathrm{Spec}^\mathrm{CS}{\Phi}^I_{*,\Gamma,-}.	
\end{align}

Then we take the corresponding quotients by using the corresponding Frobenius operators:
\begin{align}
&\mathrm{Spec}^\mathrm{CS}\widetilde{\Delta}_{*,\Gamma,-}/\mathrm{Fro}^\mathbb{Z},\mathrm{Spec}^\mathrm{CS}\widetilde{\nabla}_{*,\Gamma,-}/\mathrm{Fro}^\mathbb{Z},\mathrm{Spec}^\mathrm{CS}\widetilde{\Phi}_{*,\Gamma,-}/\mathrm{Fro}^\mathbb{Z},\mathrm{Spec}^\mathrm{CS}\widetilde{\Delta}^+_{*,\Gamma,-}/\mathrm{Fro}^\mathbb{Z},\\
&\mathrm{Spec}^\mathrm{CS}\widetilde{\nabla}^+_{*,\Gamma,-}/\mathrm{Fro}^\mathbb{Z}, \mathrm{Spec}^\mathrm{CS}\widetilde{\Delta}^\dagger_{*,\Gamma,-}/\mathrm{Fro}^\mathbb{Z},\mathrm{Spec}^\mathrm{CS}\widetilde{\nabla}^\dagger_{*,\Gamma,-}/\mathrm{Fro}^\mathbb{Z},	\\
\end{align}
\begin{align}
&\mathrm{Spec}^\mathrm{CS}\breve{\Delta}_{*,\Gamma,-}/\mathrm{Fro}^\mathbb{Z},\breve{\nabla}_{*,\Gamma,-}/\mathrm{Fro}^\mathbb{Z},\mathrm{Spec}^\mathrm{CS}\breve{\Phi}_{*,\Gamma,-}/\mathrm{Fro}^\mathbb{Z},\mathrm{Spec}^\mathrm{CS}\breve{\Delta}^+_{*,\Gamma,-}/\mathrm{Fro}^\mathbb{Z},\\
&\mathrm{Spec}^\mathrm{CS}\breve{\nabla}^+_{*,\Gamma,-}/\mathrm{Fro}^\mathbb{Z}, \mathrm{Spec}^\mathrm{CS}\breve{\Delta}^\dagger_{*,\Gamma,-}/\mathrm{Fro}^\mathbb{Z},\mathrm{Spec}^\mathrm{CS}\breve{\nabla}^\dagger_{*,\Gamma,-}/\mathrm{Fro}^\mathbb{Z},	\\
\end{align}
\begin{align}
&\mathrm{Spec}^\mathrm{CS}{\Delta}_{*,\Gamma,-}/\mathrm{Fro}^\mathbb{Z},\mathrm{Spec}^\mathrm{CS}{\nabla}_{*,\Gamma,-}/\mathrm{Fro}^\mathbb{Z},\mathrm{Spec}^\mathrm{CS}{\Phi}_{*,\Gamma,-}/\mathrm{Fro}^\mathbb{Z},\mathrm{Spec}^\mathrm{CS}{\Delta}^+_{*,\Gamma,-}/\mathrm{Fro}^\mathbb{Z},\\
&\mathrm{Spec}^\mathrm{CS}{\nabla}^+_{*,\Gamma,-}/\mathrm{Fro}^\mathbb{Z}, \mathrm{Spec}^\mathrm{CS}{\Delta}^\dagger_{*,\Gamma,-}/\mathrm{Fro}^\mathbb{Z},\mathrm{Spec}^\mathrm{CS}{\nabla}^\dagger_{*,\Gamma,-}/\mathrm{Fro}^\mathbb{Z}.	
\end{align}
Here for those space with notations related to the radius and the corresponding interval we consider the total unions $\bigcap_r,\bigcup_I$ in order to achieve the whole spaces to achieve the analogues of the corresponding FF curves from \cite{10KL1}, \cite{10KL2}, \cite{10FF} for
\[
\xymatrix@R+0pc@C+0pc{
\underset{r}{\mathrm{homotopylimit}}~\mathrm{Spec}^\mathrm{CS}\widetilde{\Phi}^r_{*,\Gamma,-},\underset{I}{\mathrm{homotopycolimit}}~\mathrm{Spec}^\mathrm{CS}\widetilde{\Phi}^I_{*,\Gamma,-},	\\
}
\]
\[
\xymatrix@R+0pc@C+0pc{
\underset{r}{\mathrm{homotopylimit}}~\mathrm{Spec}^\mathrm{CS}\breve{\Phi}^r_{*,\Gamma,-},\underset{I}{\mathrm{homotopycolimit}}~\mathrm{Spec}^\mathrm{CS}\breve{\Phi}^I_{*,\Gamma,-},	\\
}
\]
\[
\xymatrix@R+0pc@C+0pc{
\underset{r}{\mathrm{homotopylimit}}~\mathrm{Spec}^\mathrm{CS}{\Phi}^r_{*,\Gamma,-},\underset{I}{\mathrm{homotopycolimit}}~\mathrm{Spec}^\mathrm{CS}{\Phi}^I_{*,\Gamma,-}.	
}
\]
\[ 
\xymatrix@R+0pc@C+0pc{
\underset{r}{\mathrm{homotopylimit}}~\mathrm{Spec}^\mathrm{CS}\widetilde{\Phi}^r_{*,\Gamma,-}/\mathrm{Fro}^\mathbb{Z},\underset{I}{\mathrm{homotopycolimit}}~\mathrm{Spec}^\mathrm{CS}\widetilde{\Phi}^I_{*,\Gamma,-}/\mathrm{Fro}^\mathbb{Z},	\\
}
\]
\[ 
\xymatrix@R+0pc@C+0pc{
\underset{r}{\mathrm{homotopylimit}}~\mathrm{Spec}^\mathrm{CS}\breve{\Phi}^r_{*,\Gamma,-}/\mathrm{Fro}^\mathbb{Z},\underset{I}{\mathrm{homotopycolimit}}~\breve{\Phi}^I_{*,\Gamma,-}/\mathrm{Fro}^\mathbb{Z},	\\
}
\]
\[ 
\xymatrix@R+0pc@C+0pc{
\underset{r}{\mathrm{homotopylimit}}~\mathrm{Spec}^\mathrm{CS}{\Phi}^r_{*,\Gamma,-}/\mathrm{Fro}^\mathbb{Z},\underset{I}{\mathrm{homotopycolimit}}~\mathrm{Spec}^\mathrm{CS}{\Phi}^I_{*,\Gamma,-}/\mathrm{Fro}^\mathbb{Z}.	
}
\]

\end{definition}

\

\begin{definition}
We then consider the corresponding quasipresheaves of the corresponding ind-Banach or monomorphic ind-Banach modules from \cite{10BBK}, \cite{10KKM}:
\begin{align}
\mathrm{Quasicoherentpresheaves,IndBanach}_{*}	
\end{align}
where $*$ is one of the following spaces:
\begin{align}
&\mathrm{Spec}^\mathrm{BK}\widetilde{\Phi}_{*,\Gamma,-}/\mathrm{Fro}^\mathbb{Z},	\\
\end{align}
\begin{align}
&\mathrm{Spec}^\mathrm{BK}\breve{\Phi}_{*,\Gamma,-}/\mathrm{Fro}^\mathbb{Z},	\\
\end{align}
\begin{align}
&\mathrm{Spec}^\mathrm{BK}{\Phi}_{*,\Gamma,-}/\mathrm{Fro}^\mathbb{Z}.	
\end{align}
Here for those space without notation related to the radius and the corresponding interval we consider the total unions $\bigcap_r,\bigcup_I$ in order to achieve the whole spaces to achieve the analogues of the corresponding FF curves from \cite{10KL1}, \cite{10KL2}, \cite{10FF} for
\[
\xymatrix@R+0pc@C+0pc{
\underset{r}{\mathrm{homotopylimit}}~\mathrm{Spec}^\mathrm{BK}\widetilde{\Phi}^r_{*,\Gamma,-},\underset{I}{\mathrm{homotopycolimit}}~\mathrm{Spec}^\mathrm{BK}\widetilde{\Phi}^I_{*,\Gamma,-},	\\
}
\]
\[
\xymatrix@R+0pc@C+0pc{
\underset{r}{\mathrm{homotopylimit}}~\mathrm{Spec}^\mathrm{BK}\breve{\Phi}^r_{*,\Gamma,-},\underset{I}{\mathrm{homotopycolimit}}~\mathrm{Spec}^\mathrm{BK}\breve{\Phi}^I_{*,\Gamma,-},	\\
}
\]
\[
\xymatrix@R+0pc@C+0pc{
\underset{r}{\mathrm{homotopylimit}}~\mathrm{Spec}^\mathrm{BK}{\Phi}^r_{*,\Gamma,-},\underset{I}{\mathrm{homotopycolimit}}~\mathrm{Spec}^\mathrm{BK}{\Phi}^I_{*,\Gamma,-}.	
}
\]
\[  
\xymatrix@R+0pc@C+0pc{
\underset{r}{\mathrm{homotopylimit}}~\mathrm{Spec}^\mathrm{BK}\widetilde{\Phi}^r_{*,\Gamma,-}/\mathrm{Fro}^\mathbb{Z},\underset{I}{\mathrm{homotopycolimit}}~\mathrm{Spec}^\mathrm{BK}\widetilde{\Phi}^I_{*,\Gamma,-}/\mathrm{Fro}^\mathbb{Z},	\\
}
\]
\[ 
\xymatrix@R+0pc@C+0pc{
\underset{r}{\mathrm{homotopylimit}}~\mathrm{Spec}^\mathrm{BK}\breve{\Phi}^r_{*,\Gamma,-}/\mathrm{Fro}^\mathbb{Z},\underset{I}{\mathrm{homotopycolimit}}~\mathrm{Spec}^\mathrm{BK}\breve{\Phi}^I_{*,\Gamma,-}/\mathrm{Fro}^\mathbb{Z},	\\
}
\]
\[ 
\xymatrix@R+0pc@C+0pc{
\underset{r}{\mathrm{homotopylimit}}~\mathrm{Spec}^\mathrm{BK}{\Phi}^r_{*,\Gamma,-}/\mathrm{Fro}^\mathbb{Z},\underset{I}{\mathrm{homotopycolimit}}~\mathrm{Spec}^\mathrm{BK}{\Phi}^I_{*,\Gamma,-}/\mathrm{Fro}^\mathbb{Z}.	
}
\]

\end{definition}

\begin{definition}
We then consider the corresponding quasisheaves of the corresponding condensed solid topological modules from \cite{10CS2}:
\begin{align}
\mathrm{Quasicoherentsheaves, Condensed}_{*}	
\end{align}
where $*$ is one of the following spaces:
\begin{align}
&\mathrm{Spec}^\mathrm{CS}\widetilde{\Delta}_{*,\Gamma,-}/\mathrm{Fro}^\mathbb{Z},\mathrm{Spec}^\mathrm{CS}\widetilde{\nabla}_{*,\Gamma,-}/\mathrm{Fro}^\mathbb{Z},\mathrm{Spec}^\mathrm{CS}\widetilde{\Phi}_{*,\Gamma,-}/\mathrm{Fro}^\mathbb{Z},\mathrm{Spec}^\mathrm{CS}\widetilde{\Delta}^+_{*,\Gamma,-}/\mathrm{Fro}^\mathbb{Z},\\
&\mathrm{Spec}^\mathrm{CS}\widetilde{\nabla}^+_{*,\Gamma,-}/\mathrm{Fro}^\mathbb{Z},\mathrm{Spec}^\mathrm{CS}\widetilde{\Delta}^\dagger_{*,\Gamma,-}/\mathrm{Fro}^\mathbb{Z},\mathrm{Spec}^\mathrm{CS}\widetilde{\nabla}^\dagger_{*,\Gamma,-}/\mathrm{Fro}^\mathbb{Z},	\\
\end{align}
\begin{align}
&\mathrm{Spec}^\mathrm{CS}\breve{\Delta}_{*,\Gamma,-}/\mathrm{Fro}^\mathbb{Z},\breve{\nabla}_{*,\Gamma,-}/\mathrm{Fro}^\mathbb{Z},\mathrm{Spec}^\mathrm{CS}\breve{\Phi}_{*,\Gamma,-}/\mathrm{Fro}^\mathbb{Z},\mathrm{Spec}^\mathrm{CS}\breve{\Delta}^+_{*,\Gamma,-}/\mathrm{Fro}^\mathbb{Z},\\
&\mathrm{Spec}^\mathrm{CS}\breve{\nabla}^+_{*,\Gamma,-}/\mathrm{Fro}^\mathbb{Z},\mathrm{Spec}^\mathrm{CS}\breve{\Delta}^\dagger_{*,\Gamma,-}/\mathrm{Fro}^\mathbb{Z},\mathrm{Spec}^\mathrm{CS}\breve{\nabla}^\dagger_{*,\Gamma,-}/\mathrm{Fro}^\mathbb{Z},	\\
\end{align}
\begin{align}
&\mathrm{Spec}^\mathrm{CS}{\Delta}_{*,\Gamma,-}/\mathrm{Fro}^\mathbb{Z},\mathrm{Spec}^\mathrm{CS}{\nabla}_{*,\Gamma,-}/\mathrm{Fro}^\mathbb{Z},\mathrm{Spec}^\mathrm{CS}{\Phi}_{*,\Gamma,-}/\mathrm{Fro}^\mathbb{Z},\mathrm{Spec}^\mathrm{CS}{\Delta}^+_{*,\Gamma,-}/\mathrm{Fro}^\mathbb{Z},\\
&\mathrm{Spec}^\mathrm{CS}{\nabla}^+_{*,\Gamma,-}/\mathrm{Fro}^\mathbb{Z}, \mathrm{Spec}^\mathrm{CS}{\Delta}^\dagger_{*,\Gamma,-}/\mathrm{Fro}^\mathbb{Z},\mathrm{Spec}^\mathrm{CS}{\nabla}^\dagger_{*,\Gamma,-}/\mathrm{Fro}^\mathbb{Z}.	
\end{align}
Here for those space with notations related to the radius and the corresponding interval we consider the total unions $\bigcap_r,\bigcup_I$ in order to achieve the whole spaces to achieve the analogues of the corresponding FF curves from \cite{10KL1}, \cite{10KL2}, \cite{10FF} for
\[
\xymatrix@R+0pc@C+0pc{
\underset{r}{\mathrm{homotopylimit}}~\mathrm{Spec}^\mathrm{CS}\widetilde{\Phi}^r_{*,\Gamma,-},\underset{I}{\mathrm{homotopycolimit}}~\mathrm{Spec}^\mathrm{CS}\widetilde{\Phi}^I_{*,\Gamma,-},	\\
}
\]
\[
\xymatrix@R+0pc@C+0pc{
\underset{r}{\mathrm{homotopylimit}}~\mathrm{Spec}^\mathrm{CS}\breve{\Phi}^r_{*,\Gamma,-},\underset{I}{\mathrm{homotopycolimit}}~\mathrm{Spec}^\mathrm{CS}\breve{\Phi}^I_{*,\Gamma,-},	\\
}
\]
\[
\xymatrix@R+0pc@C+0pc{
\underset{r}{\mathrm{homotopylimit}}~\mathrm{Spec}^\mathrm{CS}{\Phi}^r_{*,\Gamma,-},\underset{I}{\mathrm{homotopycolimit}}~\mathrm{Spec}^\mathrm{CS}{\Phi}^I_{*,\Gamma,-}.	
}
\]
\[ 
\xymatrix@R+0pc@C+0pc{
\underset{r}{\mathrm{homotopylimit}}~\mathrm{Spec}^\mathrm{CS}\widetilde{\Phi}^r_{*,\Gamma,-}/\mathrm{Fro}^\mathbb{Z},\underset{I}{\mathrm{homotopycolimit}}~\mathrm{Spec}^\mathrm{CS}\widetilde{\Phi}^I_{*,\Gamma,-}/\mathrm{Fro}^\mathbb{Z},	\\
}
\]
\[ 
\xymatrix@R+0pc@C+0pc{
\underset{r}{\mathrm{homotopylimit}}~\mathrm{Spec}^\mathrm{CS}\breve{\Phi}^r_{*,\Gamma,-}/\mathrm{Fro}^\mathbb{Z},\underset{I}{\mathrm{homotopycolimit}}~\breve{\Phi}^I_{*,\Gamma,-}/\mathrm{Fro}^\mathbb{Z},	\\
}
\]
\[ 
\xymatrix@R+0pc@C+0pc{
\underset{r}{\mathrm{homotopylimit}}~\mathrm{Spec}^\mathrm{CS}{\Phi}^r_{*,\Gamma,-}/\mathrm{Fro}^\mathbb{Z},\underset{I}{\mathrm{homotopycolimit}}~\mathrm{Spec}^\mathrm{CS}{\Phi}^I_{*,\Gamma,-}/\mathrm{Fro}^\mathbb{Z}.	
}
\]

\end{definition}

\

\begin{proposition}
There is a well-defined functor from the $\infty$-category 
\begin{align}
\mathrm{Quasicoherentpresheaves,Condensed}_{*}	
\end{align}
where $*$ is one of the following spaces:
\begin{align}
&\mathrm{Spec}^\mathrm{CS}\widetilde{\Phi}_{*,\Gamma,-}/\mathrm{Fro}^\mathbb{Z},	\\
\end{align}
\begin{align}
&\mathrm{Spec}^\mathrm{CS}\breve{\Phi}_{*,\Gamma,-}/\mathrm{Fro}^\mathbb{Z},	\\
\end{align}
\begin{align}
&\mathrm{Spec}^\mathrm{CS}{\Phi}_{*,\Gamma,-}/\mathrm{Fro}^\mathbb{Z},	
\end{align}
to the $\infty$-category of $\mathrm{Fro}$-equivariant quasicoherent presheaves over similar spaces above correspondingly without the $\mathrm{Fro}$-quotients, and to the $\infty$-category of $\mathrm{Fro}$-equivariant quasicoherent modules over global sections of the structure $\infty$-sheaves of the similar spaces above correspondingly without the $\mathrm{Fro}$-quotients. Here for those space without notation related to the radius and the corresponding interval we consider the total unions $\bigcap_r,\bigcup_I$ in order to achieve the whole spaces to achieve the analogues of the corresponding FF curves from \cite{10KL1}, \cite{10KL2}, \cite{10FF} for
\[
\xymatrix@R+0pc@C+0pc{
\underset{r}{\mathrm{homotopylimit}}~\mathrm{Spec}^\mathrm{CS}\widetilde{\Phi}^r_{*,\Gamma,-},\underset{I}{\mathrm{homotopycolimit}}~\mathrm{Spec}^\mathrm{CS}\widetilde{\Phi}^I_{*,\Gamma,-},	\\
}
\]
\[
\xymatrix@R+0pc@C+0pc{
\underset{r}{\mathrm{homotopylimit}}~\mathrm{Spec}^\mathrm{CS}\breve{\Phi}^r_{*,\Gamma,-},\underset{I}{\mathrm{homotopycolimit}}~\mathrm{Spec}^\mathrm{CS}\breve{\Phi}^I_{*,\Gamma,-},	\\
}
\]
\[
\xymatrix@R+0pc@C+0pc{
\underset{r}{\mathrm{homotopylimit}}~\mathrm{Spec}^\mathrm{CS}{\Phi}^r_{*,\Gamma,-},\underset{I}{\mathrm{homotopycolimit}}~\mathrm{Spec}^\mathrm{CS}{\Phi}^I_{*,\Gamma,-}.	
}
\]
\[ 
\xymatrix@R+0pc@C+0pc{
\underset{r}{\mathrm{homotopylimit}}~\mathrm{Spec}^\mathrm{CS}\widetilde{\Phi}^r_{*,\Gamma,-}/\mathrm{Fro}^\mathbb{Z},\underset{I}{\mathrm{homotopycolimit}}~\mathrm{Spec}^\mathrm{CS}\widetilde{\Phi}^I_{*,\Gamma,-}/\mathrm{Fro}^\mathbb{Z},	\\
}
\]
\[ 
\xymatrix@R+0pc@C+0pc{
\underset{r}{\mathrm{homotopylimit}}~\mathrm{Spec}^\mathrm{CS}\breve{\Phi}^r_{*,\Gamma,-}/\mathrm{Fro}^\mathbb{Z},\underset{I}{\mathrm{homotopycolimit}}~\breve{\Phi}^I_{*,\Gamma,-}/\mathrm{Fro}^\mathbb{Z},	\\
}
\]
\[ 
\xymatrix@R+0pc@C+0pc{
\underset{r}{\mathrm{homotopylimit}}~\mathrm{Spec}^\mathrm{CS}{\Phi}^r_{*,\Gamma,-}/\mathrm{Fro}^\mathbb{Z},\underset{I}{\mathrm{homotopycolimit}}~\mathrm{Spec}^\mathrm{CS}{\Phi}^I_{*,\Gamma,-}/\mathrm{Fro}^\mathbb{Z}.	
}
\]	
In this situation we will have the target category being family parametrized by $r$ or $I$ in compatible glueing sense as in \cite[Definition 5.4.10]{10KL2}. In this situation for modules parametrized by the intervals we have the equivalence of $\infty$-categories by using \cite[Proposition 13.8]{10CS2}. Here the corresponding quasicoherent Frobenius modules are defined to be the corresponding homotopy colimits and limits of Frobenius modules:
\begin{align}
\underset{r}{\mathrm{homotopycolimit}}~M_r,\\
\underset{I}{\mathrm{homotopylimit}}~M_I,	
\end{align}
where each $M_r$ is a Frobenius-equivariant module over the period ring with respect to some radius $r$ while each $M_I$ is a Frobenius-equivariant module over the period ring with respect to some interval $I$.\\
\end{proposition}

\begin{proposition}
Similar proposition holds for 
\begin{align}
\mathrm{Quasicoherentsheaves,IndBanach}_{*}.	
\end{align}	
\end{proposition}

\

\begin{definition}
We then consider the corresponding quasipresheaves of perfect complexes the corresponding ind-Banach or monomorphic ind-Banach modules from \cite{10BBK}, \cite{10KKM}:
\begin{align}
\mathrm{Quasicoherentpresheaves,Perfectcomplex,IndBanach}_{*}	
\end{align}
where $*$ is one of the following spaces:
\begin{align}
&\mathrm{Spec}^\mathrm{BK}\widetilde{\Phi}_{*,\Gamma,-}/\mathrm{Fro}^\mathbb{Z},	\\
\end{align}
\begin{align}
&\mathrm{Spec}^\mathrm{BK}\breve{\Phi}_{*,\Gamma,-}/\mathrm{Fro}^\mathbb{Z},	\\
\end{align}
\begin{align}
&\mathrm{Spec}^\mathrm{BK}{\Phi}_{*,\Gamma,-}/\mathrm{Fro}^\mathbb{Z}.	
\end{align}
Here for those space without notation related to the radius and the corresponding interval we consider the total unions $\bigcap_r,\bigcup_I$ in order to achieve the whole spaces to achieve the analogues of the corresponding FF curves from \cite{10KL1}, \cite{10KL2}, \cite{10FF} for
\[
\xymatrix@R+0pc@C+0pc{
\underset{r}{\mathrm{homotopylimit}}~\mathrm{Spec}^\mathrm{BK}\widetilde{\Phi}^r_{*,\Gamma,-},\underset{I}{\mathrm{homotopycolimit}}~\mathrm{Spec}^\mathrm{BK}\widetilde{\Phi}^I_{*,\Gamma,-},	\\
}
\]
\[
\xymatrix@R+0pc@C+0pc{
\underset{r}{\mathrm{homotopylimit}}~\mathrm{Spec}^\mathrm{BK}\breve{\Phi}^r_{*,\Gamma,-},\underset{I}{\mathrm{homotopycolimit}}~\mathrm{Spec}^\mathrm{BK}\breve{\Phi}^I_{*,\Gamma,-},	\\
}
\]
\[
\xymatrix@R+0pc@C+0pc{
\underset{r}{\mathrm{homotopylimit}}~\mathrm{Spec}^\mathrm{BK}{\Phi}^r_{*,\Gamma,-},\underset{I}{\mathrm{homotopycolimit}}~\mathrm{Spec}^\mathrm{BK}{\Phi}^I_{*,\Gamma,-}.	
}
\]
\[  
\xymatrix@R+0pc@C+0pc{
\underset{r}{\mathrm{homotopylimit}}~\mathrm{Spec}^\mathrm{BK}\widetilde{\Phi}^r_{*,\Gamma,-}/\mathrm{Fro}^\mathbb{Z},\underset{I}{\mathrm{homotopycolimit}}~\mathrm{Spec}^\mathrm{BK}\widetilde{\Phi}^I_{*,\Gamma,-}/\mathrm{Fro}^\mathbb{Z},	\\
}
\]
\[ 
\xymatrix@R+0pc@C+0pc{
\underset{r}{\mathrm{homotopylimit}}~\mathrm{Spec}^\mathrm{BK}\breve{\Phi}^r_{*,\Gamma,-}/\mathrm{Fro}^\mathbb{Z},\underset{I}{\mathrm{homotopycolimit}}~\mathrm{Spec}^\mathrm{BK}\breve{\Phi}^I_{*,\Gamma,-}/\mathrm{Fro}^\mathbb{Z},	\\
}
\]
\[ 
\xymatrix@R+0pc@C+0pc{
\underset{r}{\mathrm{homotopylimit}}~\mathrm{Spec}^\mathrm{BK}{\Phi}^r_{*,\Gamma,-}/\mathrm{Fro}^\mathbb{Z},\underset{I}{\mathrm{homotopycolimit}}~\mathrm{Spec}^\mathrm{BK}{\Phi}^I_{*,\Gamma,-}/\mathrm{Fro}^\mathbb{Z}.	
}
\]

\end{definition}

\begin{definition}
We then consider the corresponding quasisheaves of perfect complexes of the corresponding condensed solid topological modules from \cite{10CS2}:
\begin{align}
\mathrm{Quasicoherentsheaves, Perfectcomplex, Condensed}_{*}	
\end{align}
where $*$ is one of the following spaces:
\begin{align}
&\mathrm{Spec}^\mathrm{CS}\widetilde{\Delta}_{*,\Gamma,-}/\mathrm{Fro}^\mathbb{Z},\mathrm{Spec}^\mathrm{CS}\widetilde{\nabla}_{*,\Gamma,-}/\mathrm{Fro}^\mathbb{Z},\mathrm{Spec}^\mathrm{CS}\widetilde{\Phi}_{*,\Gamma,-}/\mathrm{Fro}^\mathbb{Z},\mathrm{Spec}^\mathrm{CS}\widetilde{\Delta}^+_{*,\Gamma,-}/\mathrm{Fro}^\mathbb{Z},\\
&\mathrm{Spec}^\mathrm{CS}\widetilde{\nabla}^+_{*,\Gamma,-}/\mathrm{Fro}^\mathbb{Z},\mathrm{Spec}^\mathrm{CS}\widetilde{\Delta}^\dagger_{*,\Gamma,-}/\mathrm{Fro}^\mathbb{Z},\mathrm{Spec}^\mathrm{CS}\widetilde{\nabla}^\dagger_{*,\Gamma,-}/\mathrm{Fro}^\mathbb{Z},	\\
\end{align}
\begin{align}
&\mathrm{Spec}^\mathrm{CS}\breve{\Delta}_{*,\Gamma,-}/\mathrm{Fro}^\mathbb{Z},\breve{\nabla}_{*,\Gamma,-}/\mathrm{Fro}^\mathbb{Z},\mathrm{Spec}^\mathrm{CS}\breve{\Phi}_{*,\Gamma,-}/\mathrm{Fro}^\mathbb{Z},\mathrm{Spec}^\mathrm{CS}\breve{\Delta}^+_{*,\Gamma,-}/\mathrm{Fro}^\mathbb{Z},\\
&\mathrm{Spec}^\mathrm{CS}\breve{\nabla}^+_{*,\Gamma,-}/\mathrm{Fro}^\mathbb{Z},\mathrm{Spec}^\mathrm{CS}\breve{\Delta}^\dagger_{*,\Gamma,-}/\mathrm{Fro}^\mathbb{Z},\mathrm{Spec}^\mathrm{CS}\breve{\nabla}^\dagger_{*,\Gamma,-}/\mathrm{Fro}^\mathbb{Z},	\\
\end{align}
\begin{align}
&\mathrm{Spec}^\mathrm{CS}{\Delta}_{*,\Gamma,-}/\mathrm{Fro}^\mathbb{Z},\mathrm{Spec}^\mathrm{CS}{\nabla}_{*,\Gamma,-}/\mathrm{Fro}^\mathbb{Z},\mathrm{Spec}^\mathrm{CS}{\Phi}_{*,\Gamma,-}/\mathrm{Fro}^\mathbb{Z},\mathrm{Spec}^\mathrm{CS}{\Delta}^+_{*,\Gamma,-}/\mathrm{Fro}^\mathbb{Z},\\
&\mathrm{Spec}^\mathrm{CS}{\nabla}^+_{*,\Gamma,-}/\mathrm{Fro}^\mathbb{Z}, \mathrm{Spec}^\mathrm{CS}{\Delta}^\dagger_{*,\Gamma,-}/\mathrm{Fro}^\mathbb{Z},\mathrm{Spec}^\mathrm{CS}{\nabla}^\dagger_{*,\Gamma,-}/\mathrm{Fro}^\mathbb{Z}.	
\end{align}
Here for those space with notations related to the radius and the corresponding interval we consider the total unions $\bigcap_r,\bigcup_I$ in order to achieve the whole spaces to achieve the analogues of the corresponding FF curves from \cite{10KL1}, \cite{10KL2}, \cite{10FF} for
\[
\xymatrix@R+0pc@C+0pc{
\underset{r}{\mathrm{homotopylimit}}~\mathrm{Spec}^\mathrm{CS}\widetilde{\Phi}^r_{*,\Gamma,-},\underset{I}{\mathrm{homotopycolimit}}~\mathrm{Spec}^\mathrm{CS}\widetilde{\Phi}^I_{*,\Gamma,-},	\\
}
\]
\[
\xymatrix@R+0pc@C+0pc{
\underset{r}{\mathrm{homotopylimit}}~\mathrm{Spec}^\mathrm{CS}\breve{\Phi}^r_{*,\Gamma,-},\underset{I}{\mathrm{homotopycolimit}}~\mathrm{Spec}^\mathrm{CS}\breve{\Phi}^I_{*,\Gamma,-},	\\
}
\]
\[
\xymatrix@R+0pc@C+0pc{
\underset{r}{\mathrm{homotopylimit}}~\mathrm{Spec}^\mathrm{CS}{\Phi}^r_{*,\Gamma,-},\underset{I}{\mathrm{homotopycolimit}}~\mathrm{Spec}^\mathrm{CS}{\Phi}^I_{*,\Gamma,-}.	
}
\]
\[ 
\xymatrix@R+0pc@C+0pc{
\underset{r}{\mathrm{homotopylimit}}~\mathrm{Spec}^\mathrm{CS}\widetilde{\Phi}^r_{*,\Gamma,-}/\mathrm{Fro}^\mathbb{Z},\underset{I}{\mathrm{homotopycolimit}}~\mathrm{Spec}^\mathrm{CS}\widetilde{\Phi}^I_{*,\Gamma,-}/\mathrm{Fro}^\mathbb{Z},	\\
}
\]
\[ 
\xymatrix@R+0pc@C+0pc{
\underset{r}{\mathrm{homotopylimit}}~\mathrm{Spec}^\mathrm{CS}\breve{\Phi}^r_{*,\Gamma,-}/\mathrm{Fro}^\mathbb{Z},\underset{I}{\mathrm{homotopycolimit}}~\breve{\Phi}^I_{*,\Gamma,-}/\mathrm{Fro}^\mathbb{Z},	\\
}
\]
\[ 
\xymatrix@R+0pc@C+0pc{
\underset{r}{\mathrm{homotopylimit}}~\mathrm{Spec}^\mathrm{CS}{\Phi}^r_{*,\Gamma,-}/\mathrm{Fro}^\mathbb{Z},\underset{I}{\mathrm{homotopycolimit}}~\mathrm{Spec}^\mathrm{CS}{\Phi}^I_{*,\Gamma,-}/\mathrm{Fro}^\mathbb{Z}.	
}
\]

\end{definition}

\begin{proposition}
There is a well-defined functor from the $\infty$-category 
\begin{align}
\mathrm{Quasicoherentpresheaves,Perfectcomplex,Condensed}_{*}	
\end{align}
where $*$ is one of the following spaces:
\begin{align}
&\mathrm{Spec}^\mathrm{CS}\widetilde{\Phi}_{*,\Gamma,-}/\mathrm{Fro}^\mathbb{Z},	\\
\end{align}
\begin{align}
&\mathrm{Spec}^\mathrm{CS}\breve{\Phi}_{*,\Gamma,-}/\mathrm{Fro}^\mathbb{Z},	\\
\end{align}
\begin{align}
&\mathrm{Spec}^\mathrm{CS}{\Phi}_{*,\Gamma,-}/\mathrm{Fro}^\mathbb{Z},	
\end{align}
to the $\infty$-category of $\mathrm{Fro}$-equivariant quasicoherent presheaves over similar spaces above correspondingly without the $\mathrm{Fro}$-quotients, and to the $\infty$-category of $\mathrm{Fro}$-equivariant quasicoherent modules over global sections of the structure $\infty$-sheaves of the similar spaces above correspondingly without the $\mathrm{Fro}$-quotients. Here for those space without notation related to the radius and the corresponding interval we consider the total unions $\bigcap_r,\bigcup_I$ in order to achieve the whole spaces to achieve the analogues of the corresponding FF curves from \cite{10KL1}, \cite{10KL2}, \cite{10FF} for
\[
\xymatrix@R+0pc@C+0pc{
\underset{r}{\mathrm{homotopylimit}}~\mathrm{Spec}^\mathrm{CS}\widetilde{\Phi}^r_{*,\Gamma,-},\underset{I}{\mathrm{homotopycolimit}}~\mathrm{Spec}^\mathrm{CS}\widetilde{\Phi}^I_{*,\Gamma,-},	\\
}
\]
\[
\xymatrix@R+0pc@C+0pc{
\underset{r}{\mathrm{homotopylimit}}~\mathrm{Spec}^\mathrm{CS}\breve{\Phi}^r_{*,\Gamma,-},\underset{I}{\mathrm{homotopycolimit}}~\mathrm{Spec}^\mathrm{CS}\breve{\Phi}^I_{*,\Gamma,-},	\\
}
\]
\[
\xymatrix@R+0pc@C+0pc{
\underset{r}{\mathrm{homotopylimit}}~\mathrm{Spec}^\mathrm{CS}{\Phi}^r_{*,\Gamma,-},\underset{I}{\mathrm{homotopycolimit}}~\mathrm{Spec}^\mathrm{CS}{\Phi}^I_{*,\Gamma,-}.	
}
\]
\[ 
\xymatrix@R+0pc@C+0pc{
\underset{r}{\mathrm{homotopylimit}}~\mathrm{Spec}^\mathrm{CS}\widetilde{\Phi}^r_{*,\Gamma,-}/\mathrm{Fro}^\mathbb{Z},\underset{I}{\mathrm{homotopycolimit}}~\mathrm{Spec}^\mathrm{CS}\widetilde{\Phi}^I_{*,\Gamma,-}/\mathrm{Fro}^\mathbb{Z},	\\
}
\]
\[ 
\xymatrix@R+0pc@C+0pc{
\underset{r}{\mathrm{homotopylimit}}~\mathrm{Spec}^\mathrm{CS}\breve{\Phi}^r_{*,\Gamma,-}/\mathrm{Fro}^\mathbb{Z},\underset{I}{\mathrm{homotopycolimit}}~\breve{\Phi}^I_{*,\Gamma,-}/\mathrm{Fro}^\mathbb{Z},	\\
}
\]
\[ 
\xymatrix@R+0pc@C+0pc{
\underset{r}{\mathrm{homotopylimit}}~\mathrm{Spec}^\mathrm{CS}{\Phi}^r_{*,\Gamma,-}/\mathrm{Fro}^\mathbb{Z},\underset{I}{\mathrm{homotopycolimit}}~\mathrm{Spec}^\mathrm{CS}{\Phi}^I_{*,\Gamma,-}/\mathrm{Fro}^\mathbb{Z}.	
}
\]	
In this situation we will have the target category being family parametrized by $r$ or $I$ in compatible glueing sense as in \cite[Definition 5.4.10]{10KL2}. In this situation for modules parametrized by the intervals we have the equivalence of $\infty$-categories by using \cite[Proposition 12.18]{10CS2}. Here the corresponding quasicoherent Frobenius modules are defined to be the corresponding homotopy colimits and limits of Frobenius modules:
\begin{align}
\underset{r}{\mathrm{homotopycolimit}}~M_r,\\
\underset{I}{\mathrm{homotopylimit}}~M_I,	
\end{align}
where each $M_r$ is a Frobenius-equivariant module over the period ring with respect to some radius $r$ while each $M_I$ is a Frobenius-equivariant module over the period ring with respect to some interval $I$.\\
\end{proposition}

\begin{proposition}
Similar proposition holds for 
\begin{align}
\mathrm{Quasicoherentsheaves,Perfectcomplex,IndBanach}_{*}.	
\end{align}	
\end{proposition}

\newpage
\subsection{Frobenius Quasicoherent Prestacks III: Deformation in $(\infty,1)$-Ind-Banach Rings}

\begin{definition}
We now consider the pro-\'etale site of $\mathrm{Spa}\mathbb{Q}_p\left<X_1^{\pm 1},...,X_k^{\pm 1}\right>$, denote that by $*$. To be more accurate we replace one component for $\Gamma$ with the pro-\'etale site of $\mathrm{Spa}\mathbb{Q}_p\left<X_1^{\pm 1},...,X_k^{\pm 1}\right>$. And we treat then all the functor to be prestacks for this site. Then from \cite{10KL1} and \cite[Definition 5.2.1]{10KL2} we have the following class of Kedlaya-Liu rings (with the following replacement: $\Delta$ stands for $A$, $\nabla$ stands for $B$, while $\Phi$ stands for $C$) by taking product in the sense of self $\Gamma$-th power:

\[
\xymatrix@R+0pc@C+0pc{
\widetilde{\Delta}_{*,\Gamma},\widetilde{\nabla}_{*,\Gamma},\widetilde{\Phi}_{*,\Gamma},\widetilde{\Delta}^+_{*,\Gamma},\widetilde{\nabla}^+_{*,\Gamma},\widetilde{\Delta}^\dagger_{*,\Gamma},\widetilde{\nabla}^\dagger_{*,\Gamma},\widetilde{\Phi}^r_{*,\Gamma},\widetilde{\Phi}^I_{*,\Gamma}, 
}
\]

\[
\xymatrix@R+0pc@C+0pc{
\breve{\Delta}_{*,\Gamma},\breve{\nabla}_{*,\Gamma},\breve{\Phi}_{*,\Gamma},\breve{\Delta}^+_{*,\Gamma},\breve{\nabla}^+_{*,\Gamma},\breve{\Delta}^\dagger_{*,\Gamma},\breve{\nabla}^\dagger_{*,\Gamma},\breve{\Phi}^r_{*,\Gamma},\breve{\Phi}^I_{*,\Gamma},	
}
\]

\[
\xymatrix@R+0pc@C+0pc{
{\Delta}_{*,\Gamma},{\nabla}_{*,\Gamma},{\Phi}_{*,\Gamma},{\Delta}^+_{*,\Gamma},{\nabla}^+_{*,\Gamma},{\Delta}^\dagger_{*,\Gamma},{\nabla}^\dagger_{*,\Gamma},{\Phi}^r_{*,\Gamma},{\Phi}^I_{*,\Gamma}.	
}
\]
We now consider the following rings with $\square$ being a homotopy colimit
\begin{align}
 \underset{i}{\mathrm{homotopycolimit}}\square_i
 \end{align}
 of $\mathbb{Q}_p\left<Y_1,...,Y_i\right>,i=1,2,...$ in $\infty$-categories of simplicial ind-Banach rings in \cite{10BBBK}
 \begin{align}
  \mathrm{SimplicialInd-BanachRings}_{\mathbb{Q}_p}
\end{align}  
or animated analytic condensed commutative algebras in \cite{10CS2} 
\begin{align}   
\mathrm{SimplicialAnalyticCondensed}_{\mathbb{Q}_p}.
\end{align}   
Taking the product we have:
\[
\xymatrix@R+0pc@C+0pc{
\widetilde{\Phi}_{*,\Gamma,\square},\widetilde{\Phi}^r_{*,\Gamma,\square},\widetilde{\Phi}^I_{*,\Gamma,\square},	
}
\]
\[
\xymatrix@R+0pc@C+0pc{
\breve{\Phi}_{*,\Gamma,\square},\breve{\Phi}^r_{*,\Gamma,\square},\breve{\Phi}^I_{*,\Gamma,\square},	
}
\]
\[
\xymatrix@R+0pc@C+0pc{
{\Phi}_{*,\Gamma,\square},{\Phi}^r_{*,\Gamma,\square},{\Phi}^I_{*,\Gamma,\square}.	
}
\]
They carry multi Frobenius action $\varphi_\Gamma$ and multi $\mathrm{Lie}_\Gamma:=\mathbb{Z}_p^{\times\Gamma}$ action. In our current situation after \cite{10CKZ} and \cite{10PZ} we consider the following $(\infty,1)$-categories of $(\infty,1)$-modules.\\
\end{definition}

\begin{definition}
First we consider the Bambozzi-Kremnizer spectrum $\mathrm{Spec}^\mathrm{BK}(*)$ attached to any of those in the above from \cite{10BK} by taking derived rational localization:
\begin{align}
&\mathrm{Spec}^\mathrm{BK}\widetilde{\Phi}_{*,\Gamma,\square},\mathrm{Spec}^\mathrm{BK}\widetilde{\Phi}^r_{*,\Gamma,\square},\mathrm{Spec}^\mathrm{BK}\widetilde{\Phi}^I_{*,\Gamma,\square},	
\end{align}
\begin{align}
&\mathrm{Spec}^\mathrm{BK}\breve{\Phi}_{*,\Gamma,\square},\mathrm{Spec}^\mathrm{BK}\breve{\Phi}^r_{*,\Gamma,\square},\mathrm{Spec}^\mathrm{BK}\breve{\Phi}^I_{*,\Gamma,\square},	
\end{align}
\begin{align}
&\mathrm{Spec}^\mathrm{BK}{\Phi}_{*,\Gamma,\square},
\mathrm{Spec}^\mathrm{BK}{\Phi}^r_{*,\Gamma,\square},\mathrm{Spec}^\mathrm{BK}{\Phi}^I_{*,\Gamma,\square}.	
\end{align}

Then we take the corresponding quotients by using the corresponding Frobenius operators:
\begin{align}
&\mathrm{Spec}^\mathrm{BK}\widetilde{\Phi}_{*,\Gamma,\square}/\mathrm{Fro}^\mathbb{Z},	\\
\end{align}
\begin{align}
&\mathrm{Spec}^\mathrm{BK}\breve{\Phi}_{*,\Gamma,\square}/\mathrm{Fro}^\mathbb{Z},	\\
\end{align}
\begin{align}
&\mathrm{Spec}^\mathrm{BK}{\Phi}_{*,\Gamma,\square}/\mathrm{Fro}^\mathbb{Z}.	
\end{align}
Here for those space without notation related to the radius and the corresponding interval we consider the total unions $\bigcap_r,\bigcup_I$ in order to achieve the whole spaces to achieve the analogues of the corresponding FF curves from \cite{10KL1}, \cite{10KL2}, \cite{10FF} for
\[
\xymatrix@R+0pc@C+0pc{
\underset{r}{\mathrm{homotopylimit}}~\mathrm{Spec}^\mathrm{BK}\widetilde{\Phi}^r_{*,\Gamma,\square},\underset{I}{\mathrm{homotopycolimit}}~\mathrm{Spec}^\mathrm{BK}\widetilde{\Phi}^I_{*,\Gamma,\square},	\\
}
\]
\[
\xymatrix@R+0pc@C+0pc{
\underset{r}{\mathrm{homotopylimit}}~\mathrm{Spec}^\mathrm{BK}\breve{\Phi}^r_{*,\Gamma,\square},\underset{I}{\mathrm{homotopycolimit}}~\mathrm{Spec}^\mathrm{BK}\breve{\Phi}^I_{*,\Gamma,\square},	\\
}
\]
\[
\xymatrix@R+0pc@C+0pc{
\underset{r}{\mathrm{homotopylimit}}~\mathrm{Spec}^\mathrm{BK}{\Phi}^r_{*,\Gamma,\square},\underset{I}{\mathrm{homotopycolimit}}~\mathrm{Spec}^\mathrm{BK}{\Phi}^I_{*,\Gamma,\square}.	
}
\]
\[  
\xymatrix@R+0pc@C+0pc{
\underset{r}{\mathrm{homotopylimit}}~\mathrm{Spec}^\mathrm{BK}\widetilde{\Phi}^r_{*,\Gamma,\square}/\mathrm{Fro}^\mathbb{Z},\underset{I}{\mathrm{homotopycolimit}}~\mathrm{Spec}^\mathrm{BK}\widetilde{\Phi}^I_{*,\Gamma,\square}/\mathrm{Fro}^\mathbb{Z},	\\
}
\]
\[ 
\xymatrix@R+0pc@C+0pc{
\underset{r}{\mathrm{homotopylimit}}~\mathrm{Spec}^\mathrm{BK}\breve{\Phi}^r_{*,\Gamma,\square}/\mathrm{Fro}^\mathbb{Z},\underset{I}{\mathrm{homotopycolimit}}~\mathrm{Spec}^\mathrm{BK}\breve{\Phi}^I_{*,\Gamma,\square}/\mathrm{Fro}^\mathbb{Z},	\\
}
\]
\[ 
\xymatrix@R+0pc@C+0pc{
\underset{r}{\mathrm{homotopylimit}}~\mathrm{Spec}^\mathrm{BK}{\Phi}^r_{*,\Gamma,\square}/\mathrm{Fro}^\mathbb{Z},\underset{I}{\mathrm{homotopycolimit}}~\mathrm{Spec}^\mathrm{BK}{\Phi}^I_{*,\Gamma,\square}/\mathrm{Fro}^\mathbb{Z}.	
}
\]

\end{definition}

\indent Meanwhile we have the corresponding Clausen-Scholze analytic stacks from \cite{10CS2}, therefore applying their construction we have:

\begin{definition}
Here we define the following products by using the solidified tensor product from \cite{10CS1} and \cite{10CS2}. Namely $A$ will still as above as a Banach ring over $\mathbb{Q}_p$. Then we take solidified tensor product $\overset{\blacksquare}{\otimes}$ of any of the following
\[
\xymatrix@R+0pc@C+0pc{
\widetilde{\Delta}_{*,\Gamma},\widetilde{\nabla}_{*,\Gamma},\widetilde{\Phi}_{*,\Gamma},\widetilde{\Delta}^+_{*,\Gamma},\widetilde{\nabla}^+_{*,\Gamma},\widetilde{\Delta}^\dagger_{*,\Gamma},\widetilde{\nabla}^\dagger_{*,\Gamma},\widetilde{\Phi}^r_{*,\Gamma},\widetilde{\Phi}^I_{*,\Gamma}, 
}
\]

\[
\xymatrix@R+0pc@C+0pc{
\breve{\Delta}_{*,\Gamma},\breve{\nabla}_{*,\Gamma},\breve{\Phi}_{*,\Gamma},\breve{\Delta}^+_{*,\Gamma},\breve{\nabla}^+_{*,\Gamma},\breve{\Delta}^\dagger_{*,\Gamma},\breve{\nabla}^\dagger_{*,\Gamma},\breve{\Phi}^r_{*,\Gamma},\breve{\Phi}^I_{*,\Gamma},	
}
\]

\[
\xymatrix@R+0pc@C+0pc{
{\Delta}_{*,\Gamma},{\nabla}_{*,\Gamma},{\Phi}_{*,\Gamma},{\Delta}^+_{*,\Gamma},{\nabla}^+_{*,\Gamma},{\Delta}^\dagger_{*,\Gamma},{\nabla}^\dagger_{*,\Gamma},{\Phi}^r_{*,\Gamma},{\Phi}^I_{*,\Gamma},	
}
\]  	
with $A$. Then we have the notations:
\[
\xymatrix@R+0pc@C+0pc{
\widetilde{\Delta}_{*,\Gamma,\square},\widetilde{\nabla}_{*,\Gamma,\square},\widetilde{\Phi}_{*,\Gamma,\square},\widetilde{\Delta}^+_{*,\Gamma,\square},\widetilde{\nabla}^+_{*,\Gamma,\square},\widetilde{\Delta}^\dagger_{*,\Gamma,\square},\widetilde{\nabla}^\dagger_{*,\Gamma,\square},\widetilde{\Phi}^r_{*,\Gamma,\square},\widetilde{\Phi}^I_{*,\Gamma,\square}, 
}
\]

\[
\xymatrix@R+0pc@C+0pc{
\breve{\Delta}_{*,\Gamma,\square},\breve{\nabla}_{*,\Gamma,\square},\breve{\Phi}_{*,\Gamma,\square},\breve{\Delta}^+_{*,\Gamma,\square},\breve{\nabla}^+_{*,\Gamma,\square},\breve{\Delta}^\dagger_{*,\Gamma,\square},\breve{\nabla}^\dagger_{*,\Gamma,\square},\breve{\Phi}^r_{*,\Gamma,\square},\breve{\Phi}^I_{*,\Gamma,\square},	
}
\]

\[
\xymatrix@R+0pc@C+0pc{
{\Delta}_{*,\Gamma,\square},{\nabla}_{*,\Gamma,\square},{\Phi}_{*,\Gamma,\square},{\Delta}^+_{*,\Gamma,\square},{\nabla}^+_{*,\Gamma,\square},{\Delta}^\dagger_{*,\Gamma,\square},{\nabla}^\dagger_{*,\Gamma,\square},{\Phi}^r_{*,\Gamma,\square},{\Phi}^I_{*,\Gamma,\square}.	
}
\]
\end{definition}

\begin{definition}
First we consider the Clausen-Scholze spectrum $\mathrm{Spec}^\mathrm{CS}(*)$ attached to any of those in the above from \cite{10CS2} by taking derived rational localization:
\begin{align}
\mathrm{Spec}^\mathrm{CS}\widetilde{\Delta}_{*,\Gamma,\square},\mathrm{Spec}^\mathrm{CS}\widetilde{\nabla}_{*,\Gamma,\square},\mathrm{Spec}^\mathrm{CS}\widetilde{\Phi}_{*,\Gamma,\square},\mathrm{Spec}^\mathrm{CS}\widetilde{\Delta}^+_{*,\Gamma,\square},\mathrm{Spec}^\mathrm{CS}\widetilde{\nabla}^+_{*,\Gamma,\square},\\
\mathrm{Spec}^\mathrm{CS}\widetilde{\Delta}^\dagger_{*,\Gamma,\square},\mathrm{Spec}^\mathrm{CS}\widetilde{\nabla}^\dagger_{*,\Gamma,\square},\mathrm{Spec}^\mathrm{CS}\widetilde{\Phi}^r_{*,\Gamma,\square},\mathrm{Spec}^\mathrm{CS}\widetilde{\Phi}^I_{*,\Gamma,\square},	\\
\end{align}
\begin{align}
\mathrm{Spec}^\mathrm{CS}\breve{\Delta}_{*,\Gamma,\square},\breve{\nabla}_{*,\Gamma,\square},\mathrm{Spec}^\mathrm{CS}\breve{\Phi}_{*,\Gamma,\square},\mathrm{Spec}^\mathrm{CS}\breve{\Delta}^+_{*,\Gamma,\square},\mathrm{Spec}^\mathrm{CS}\breve{\nabla}^+_{*,\Gamma,\square},\\
\mathrm{Spec}^\mathrm{CS}\breve{\Delta}^\dagger_{*,\Gamma,\square},\mathrm{Spec}^\mathrm{CS}\breve{\nabla}^\dagger_{*,\Gamma,\square},\mathrm{Spec}^\mathrm{CS}\breve{\Phi}^r_{*,\Gamma,\square},\breve{\Phi}^I_{*,\Gamma,\square},	\\
\end{align}
\begin{align}
\mathrm{Spec}^\mathrm{CS}{\Delta}_{*,\Gamma,\square},\mathrm{Spec}^\mathrm{CS}{\nabla}_{*,\Gamma,\square},\mathrm{Spec}^\mathrm{CS}{\Phi}_{*,\Gamma,\square},\mathrm{Spec}^\mathrm{CS}{\Delta}^+_{*,\Gamma,\square},\mathrm{Spec}^\mathrm{CS}{\nabla}^+_{*,\Gamma,\square},\\
\mathrm{Spec}^\mathrm{CS}{\Delta}^\dagger_{*,\Gamma,\square},\mathrm{Spec}^\mathrm{CS}{\nabla}^\dagger_{*,\Gamma,\square},\mathrm{Spec}^\mathrm{CS}{\Phi}^r_{*,\Gamma,\square},\mathrm{Spec}^\mathrm{CS}{\Phi}^I_{*,\Gamma,\square}.	
\end{align}

Then we take the corresponding quotients by using the corresponding Frobenius operators:
\begin{align}
&\mathrm{Spec}^\mathrm{CS}\widetilde{\Delta}_{*,\Gamma,\square}/\mathrm{Fro}^\mathbb{Z},\mathrm{Spec}^\mathrm{CS}\widetilde{\nabla}_{*,\Gamma,\square}/\mathrm{Fro}^\mathbb{Z},\mathrm{Spec}^\mathrm{CS}\widetilde{\Phi}_{*,\Gamma,\square}/\mathrm{Fro}^\mathbb{Z},\mathrm{Spec}^\mathrm{CS}\widetilde{\Delta}^+_{*,\Gamma,\square}/\mathrm{Fro}^\mathbb{Z},\\
&\mathrm{Spec}^\mathrm{CS}\widetilde{\nabla}^+_{*,\Gamma,\square}/\mathrm{Fro}^\mathbb{Z}, \mathrm{Spec}^\mathrm{CS}\widetilde{\Delta}^\dagger_{*,\Gamma,\square}/\mathrm{Fro}^\mathbb{Z},\mathrm{Spec}^\mathrm{CS}\widetilde{\nabla}^\dagger_{*,\Gamma,\square}/\mathrm{Fro}^\mathbb{Z},	\\
\end{align}
\begin{align}
&\mathrm{Spec}^\mathrm{CS}\breve{\Delta}_{*,\Gamma,\square}/\mathrm{Fro}^\mathbb{Z},\breve{\nabla}_{*,\Gamma,\square}/\mathrm{Fro}^\mathbb{Z},\mathrm{Spec}^\mathrm{CS}\breve{\Phi}_{*,\Gamma,\square}/\mathrm{Fro}^\mathbb{Z},\mathrm{Spec}^\mathrm{CS}\breve{\Delta}^+_{*,\Gamma,\square}/\mathrm{Fro}^\mathbb{Z},\\
&\mathrm{Spec}^\mathrm{CS}\breve{\nabla}^+_{*,\Gamma,\square}/\mathrm{Fro}^\mathbb{Z}, \mathrm{Spec}^\mathrm{CS}\breve{\Delta}^\dagger_{*,\Gamma,\square}/\mathrm{Fro}^\mathbb{Z},\mathrm{Spec}^\mathrm{CS}\breve{\nabla}^\dagger_{*,\Gamma,\square}/\mathrm{Fro}^\mathbb{Z},	\\
\end{align}
\begin{align}
&\mathrm{Spec}^\mathrm{CS}{\Delta}_{*,\Gamma,\square}/\mathrm{Fro}^\mathbb{Z},\mathrm{Spec}^\mathrm{CS}{\nabla}_{*,\Gamma,\square}/\mathrm{Fro}^\mathbb{Z},\mathrm{Spec}^\mathrm{CS}{\Phi}_{*,\Gamma,\square}/\mathrm{Fro}^\mathbb{Z},\mathrm{Spec}^\mathrm{CS}{\Delta}^+_{*,\Gamma,\square}/\mathrm{Fro}^\mathbb{Z},\\
&\mathrm{Spec}^\mathrm{CS}{\nabla}^+_{*,\Gamma,\square}/\mathrm{Fro}^\mathbb{Z}, \mathrm{Spec}^\mathrm{CS}{\Delta}^\dagger_{*,\Gamma,\square}/\mathrm{Fro}^\mathbb{Z},\mathrm{Spec}^\mathrm{CS}{\nabla}^\dagger_{*,\Gamma,\square}/\mathrm{Fro}^\mathbb{Z}.	
\end{align}
Here for those space with notations related to the radius and the corresponding interval we consider the total unions $\bigcap_r,\bigcup_I$ in order to achieve the whole spaces to achieve the analogues of the corresponding FF curves from \cite{10KL1}, \cite{10KL2}, \cite{10FF} for
\[
\xymatrix@R+0pc@C+0pc{
\underset{r}{\mathrm{homotopylimit}}~\mathrm{Spec}^\mathrm{CS}\widetilde{\Phi}^r_{*,\Gamma,\square},\underset{I}{\mathrm{homotopycolimit}}~\mathrm{Spec}^\mathrm{CS}\widetilde{\Phi}^I_{*,\Gamma,\square},	\\
}
\]
\[
\xymatrix@R+0pc@C+0pc{
\underset{r}{\mathrm{homotopylimit}}~\mathrm{Spec}^\mathrm{CS}\breve{\Phi}^r_{*,\Gamma,\square},\underset{I}{\mathrm{homotopycolimit}}~\mathrm{Spec}^\mathrm{CS}\breve{\Phi}^I_{*,\Gamma,\square},	\\
}
\]
\[
\xymatrix@R+0pc@C+0pc{
\underset{r}{\mathrm{homotopylimit}}~\mathrm{Spec}^\mathrm{CS}{\Phi}^r_{*,\Gamma,\square},\underset{I}{\mathrm{homotopycolimit}}~\mathrm{Spec}^\mathrm{CS}{\Phi}^I_{*,\Gamma,\square}.	
}
\]
\[ 
\xymatrix@R+0pc@C+0pc{
\underset{r}{\mathrm{homotopylimit}}~\mathrm{Spec}^\mathrm{CS}\widetilde{\Phi}^r_{*,\Gamma,\square}/\mathrm{Fro}^\mathbb{Z},\underset{I}{\mathrm{homotopycolimit}}~\mathrm{Spec}^\mathrm{CS}\widetilde{\Phi}^I_{*,\Gamma,\square}/\mathrm{Fro}^\mathbb{Z},	\\
}
\]
\[ 
\xymatrix@R+0pc@C+0pc{
\underset{r}{\mathrm{homotopylimit}}~\mathrm{Spec}^\mathrm{CS}\breve{\Phi}^r_{*,\Gamma,\square}/\mathrm{Fro}^\mathbb{Z},\underset{I}{\mathrm{homotopycolimit}}~\breve{\Phi}^I_{*,\Gamma,\square}/\mathrm{Fro}^\mathbb{Z},	\\
}
\]
\[ 
\xymatrix@R+0pc@C+0pc{
\underset{r}{\mathrm{homotopylimit}}~\mathrm{Spec}^\mathrm{CS}{\Phi}^r_{*,\Gamma,\square}/\mathrm{Fro}^\mathbb{Z},\underset{I}{\mathrm{homotopycolimit}}~\mathrm{Spec}^\mathrm{CS}{\Phi}^I_{*,\Gamma,\square}/\mathrm{Fro}^\mathbb{Z}.	
}
\]

\end{definition}

\

\begin{definition}
We then consider the corresponding quasipresheaves of the corresponding ind-Banach or monomorphic ind-Banach modules from \cite{10BBK}, \cite{10KKM}:
\begin{align}
\mathrm{Quasicoherentpresheaves,IndBanach}_{*}	
\end{align}
where $*$ is one of the following spaces:
\begin{align}
&\mathrm{Spec}^\mathrm{BK}\widetilde{\Phi}_{*,\Gamma,\square}/\mathrm{Fro}^\mathbb{Z},	\\
\end{align}
\begin{align}
&\mathrm{Spec}^\mathrm{BK}\breve{\Phi}_{*,\Gamma,\square}/\mathrm{Fro}^\mathbb{Z},	\\
\end{align}
\begin{align}
&\mathrm{Spec}^\mathrm{BK}{\Phi}_{*,\Gamma,\square}/\mathrm{Fro}^\mathbb{Z}.	
\end{align}
Here for those space without notation related to the radius and the corresponding interval we consider the total unions $\bigcap_r,\bigcup_I$ in order to achieve the whole spaces to achieve the analogues of the corresponding FF curves from \cite{10KL1}, \cite{10KL2}, \cite{10FF} for
\[
\xymatrix@R+0pc@C+0pc{
\underset{r}{\mathrm{homotopylimit}}~\mathrm{Spec}^\mathrm{BK}\widetilde{\Phi}^r_{*,\Gamma,\square},\underset{I}{\mathrm{homotopycolimit}}~\mathrm{Spec}^\mathrm{BK}\widetilde{\Phi}^I_{*,\Gamma,\square},	\\
}
\]
\[
\xymatrix@R+0pc@C+0pc{
\underset{r}{\mathrm{homotopylimit}}~\mathrm{Spec}^\mathrm{BK}\breve{\Phi}^r_{*,\Gamma,\square},\underset{I}{\mathrm{homotopycolimit}}~\mathrm{Spec}^\mathrm{BK}\breve{\Phi}^I_{*,\Gamma,\square},	\\
}
\]
\[
\xymatrix@R+0pc@C+0pc{
\underset{r}{\mathrm{homotopylimit}}~\mathrm{Spec}^\mathrm{BK}{\Phi}^r_{*,\Gamma,\square},\underset{I}{\mathrm{homotopycolimit}}~\mathrm{Spec}^\mathrm{BK}{\Phi}^I_{*,\Gamma,\square}.	
}
\]
\[  
\xymatrix@R+0pc@C+0pc{
\underset{r}{\mathrm{homotopylimit}}~\mathrm{Spec}^\mathrm{BK}\widetilde{\Phi}^r_{*,\Gamma,\square}/\mathrm{Fro}^\mathbb{Z},\underset{I}{\mathrm{homotopycolimit}}~\mathrm{Spec}^\mathrm{BK}\widetilde{\Phi}^I_{*,\Gamma,\square}/\mathrm{Fro}^\mathbb{Z},	\\
}
\]
\[ 
\xymatrix@R+0pc@C+0pc{
\underset{r}{\mathrm{homotopylimit}}~\mathrm{Spec}^\mathrm{BK}\breve{\Phi}^r_{*,\Gamma,\square}/\mathrm{Fro}^\mathbb{Z},\underset{I}{\mathrm{homotopycolimit}}~\mathrm{Spec}^\mathrm{BK}\breve{\Phi}^I_{*,\Gamma,\square}/\mathrm{Fro}^\mathbb{Z},	\\
}
\]
\[ 
\xymatrix@R+0pc@C+0pc{
\underset{r}{\mathrm{homotopylimit}}~\mathrm{Spec}^\mathrm{BK}{\Phi}^r_{*,\Gamma,\square}/\mathrm{Fro}^\mathbb{Z},\underset{I}{\mathrm{homotopycolimit}}~\mathrm{Spec}^\mathrm{BK}{\Phi}^I_{*,\Gamma,\square}/\mathrm{Fro}^\mathbb{Z}.	
}
\]

\end{definition}

\begin{definition}
We then consider the corresponding quasisheaves of the corresponding condensed solid topological modules from \cite{10CS2}:
\begin{align}
\mathrm{Quasicoherentsheaves, Condensed}_{*}	
\end{align}
where $*$ is one of the following spaces:
\begin{align}
&\mathrm{Spec}^\mathrm{CS}\widetilde{\Delta}_{*,\Gamma,\square}/\mathrm{Fro}^\mathbb{Z},\mathrm{Spec}^\mathrm{CS}\widetilde{\nabla}_{*,\Gamma,\square}/\mathrm{Fro}^\mathbb{Z},\mathrm{Spec}^\mathrm{CS}\widetilde{\Phi}_{*,\Gamma,\square}/\mathrm{Fro}^\mathbb{Z},\mathrm{Spec}^\mathrm{CS}\widetilde{\Delta}^+_{*,\Gamma,\square}/\mathrm{Fro}^\mathbb{Z},\\
&\mathrm{Spec}^\mathrm{CS}\widetilde{\nabla}^+_{*,\Gamma,\square}/\mathrm{Fro}^\mathbb{Z},\mathrm{Spec}^\mathrm{CS}\widetilde{\Delta}^\dagger_{*,\Gamma,\square}/\mathrm{Fro}^\mathbb{Z},\mathrm{Spec}^\mathrm{CS}\widetilde{\nabla}^\dagger_{*,\Gamma,\square}/\mathrm{Fro}^\mathbb{Z},	\\
\end{align}
\begin{align}
&\mathrm{Spec}^\mathrm{CS}\breve{\Delta}_{*,\Gamma,\square}/\mathrm{Fro}^\mathbb{Z},\breve{\nabla}_{*,\Gamma,\square}/\mathrm{Fro}^\mathbb{Z},\mathrm{Spec}^\mathrm{CS}\breve{\Phi}_{*,\Gamma,\square}/\mathrm{Fro}^\mathbb{Z},\mathrm{Spec}^\mathrm{CS}\breve{\Delta}^+_{*,\Gamma,\square}/\mathrm{Fro}^\mathbb{Z},\\
&\mathrm{Spec}^\mathrm{CS}\breve{\nabla}^+_{*,\Gamma,\square}/\mathrm{Fro}^\mathbb{Z},\mathrm{Spec}^\mathrm{CS}\breve{\Delta}^\dagger_{*,\Gamma,\square}/\mathrm{Fro}^\mathbb{Z},\mathrm{Spec}^\mathrm{CS}\breve{\nabla}^\dagger_{*,\Gamma,\square}/\mathrm{Fro}^\mathbb{Z},	\\
\end{align}
\begin{align}
&\mathrm{Spec}^\mathrm{CS}{\Delta}_{*,\Gamma,\square}/\mathrm{Fro}^\mathbb{Z},\mathrm{Spec}^\mathrm{CS}{\nabla}_{*,\Gamma,\square}/\mathrm{Fro}^\mathbb{Z},\mathrm{Spec}^\mathrm{CS}{\Phi}_{*,\Gamma,\square}/\mathrm{Fro}^\mathbb{Z},\mathrm{Spec}^\mathrm{CS}{\Delta}^+_{*,\Gamma,\square}/\mathrm{Fro}^\mathbb{Z},\\
&\mathrm{Spec}^\mathrm{CS}{\nabla}^+_{*,\Gamma,\square}/\mathrm{Fro}^\mathbb{Z}, \mathrm{Spec}^\mathrm{CS}{\Delta}^\dagger_{*,\Gamma,\square}/\mathrm{Fro}^\mathbb{Z},\mathrm{Spec}^\mathrm{CS}{\nabla}^\dagger_{*,\Gamma,\square}/\mathrm{Fro}^\mathbb{Z}.	
\end{align}
Here for those space with notations related to the radius and the corresponding interval we consider the total unions $\bigcap_r,\bigcup_I$ in order to achieve the whole spaces to achieve the analogues of the corresponding FF curves from \cite{10KL1}, \cite{10KL2}, \cite{10FF} for
\[
\xymatrix@R+0pc@C+0pc{
\underset{r}{\mathrm{homotopylimit}}~\mathrm{Spec}^\mathrm{CS}\widetilde{\Phi}^r_{*,\Gamma,\square},\underset{I}{\mathrm{homotopycolimit}}~\mathrm{Spec}^\mathrm{CS}\widetilde{\Phi}^I_{*,\Gamma,\square},	\\
}
\]
\[
\xymatrix@R+0pc@C+0pc{
\underset{r}{\mathrm{homotopylimit}}~\mathrm{Spec}^\mathrm{CS}\breve{\Phi}^r_{*,\Gamma,\square},\underset{I}{\mathrm{homotopycolimit}}~\mathrm{Spec}^\mathrm{CS}\breve{\Phi}^I_{*,\Gamma,\square},	\\
}
\]
\[
\xymatrix@R+0pc@C+0pc{
\underset{r}{\mathrm{homotopylimit}}~\mathrm{Spec}^\mathrm{CS}{\Phi}^r_{*,\Gamma,\square},\underset{I}{\mathrm{homotopycolimit}}~\mathrm{Spec}^\mathrm{CS}{\Phi}^I_{*,\Gamma,\square}.	
}
\]
\[ 
\xymatrix@R+0pc@C+0pc{
\underset{r}{\mathrm{homotopylimit}}~\mathrm{Spec}^\mathrm{CS}\widetilde{\Phi}^r_{*,\Gamma,\square}/\mathrm{Fro}^\mathbb{Z},\underset{I}{\mathrm{homotopycolimit}}~\mathrm{Spec}^\mathrm{CS}\widetilde{\Phi}^I_{*,\Gamma,\square}/\mathrm{Fro}^\mathbb{Z},	\\
}
\]
\[ 
\xymatrix@R+0pc@C+0pc{
\underset{r}{\mathrm{homotopylimit}}~\mathrm{Spec}^\mathrm{CS}\breve{\Phi}^r_{*,\Gamma,\square}/\mathrm{Fro}^\mathbb{Z},\underset{I}{\mathrm{homotopycolimit}}~\breve{\Phi}^I_{*,\Gamma,\square}/\mathrm{Fro}^\mathbb{Z},	\\
}
\]
\[ 
\xymatrix@R+0pc@C+0pc{
\underset{r}{\mathrm{homotopylimit}}~\mathrm{Spec}^\mathrm{CS}{\Phi}^r_{*,\Gamma,\square}/\mathrm{Fro}^\mathbb{Z},\underset{I}{\mathrm{homotopycolimit}}~\mathrm{Spec}^\mathrm{CS}{\Phi}^I_{*,\Gamma,\square}/\mathrm{Fro}^\mathbb{Z}.	
}
\]

\end{definition}

\

\begin{proposition}
There is a well-defined functor from the $\infty$-category 
\begin{align}
\mathrm{Quasicoherentpresheaves,Condensed}_{*}	
\end{align}
where $*$ is one of the following spaces:
\begin{align}
&\mathrm{Spec}^\mathrm{CS}\widetilde{\Phi}_{*,\Gamma,\square}/\mathrm{Fro}^\mathbb{Z},	\\
\end{align}
\begin{align}
&\mathrm{Spec}^\mathrm{CS}\breve{\Phi}_{*,\Gamma,\square}/\mathrm{Fro}^\mathbb{Z},	\\
\end{align}
\begin{align}
&\mathrm{Spec}^\mathrm{CS}{\Phi}_{*,\Gamma,\square}/\mathrm{Fro}^\mathbb{Z},	
\end{align}
to the $\infty$-category of $\mathrm{Fro}$-equivariant quasicoherent presheaves over similar spaces above correspondingly without the $\mathrm{Fro}$-quotients, and to the $\infty$-category of $\mathrm{Fro}$-equivariant quasicoherent modules over global sections of the structure $\infty$-sheaves of the similar spaces above correspondingly without the $\mathrm{Fro}$-quotients. Here for those space without notation related to the radius and the corresponding interval we consider the total unions $\bigcap_r,\bigcup_I$ in order to achieve the whole spaces to achieve the analogues of the corresponding FF curves from \cite{10KL1}, \cite{10KL2}, \cite{10FF} for
\[
\xymatrix@R+0pc@C+0pc{
\underset{r}{\mathrm{homotopylimit}}~\mathrm{Spec}^\mathrm{CS}\widetilde{\Phi}^r_{*,\Gamma,\square},\underset{I}{\mathrm{homotopycolimit}}~\mathrm{Spec}^\mathrm{CS}\widetilde{\Phi}^I_{*,\Gamma,\square},	\\
}
\]
\[
\xymatrix@R+0pc@C+0pc{
\underset{r}{\mathrm{homotopylimit}}~\mathrm{Spec}^\mathrm{CS}\breve{\Phi}^r_{*,\Gamma,\square},\underset{I}{\mathrm{homotopycolimit}}~\mathrm{Spec}^\mathrm{CS}\breve{\Phi}^I_{*,\Gamma,\square},	\\
}
\]
\[
\xymatrix@R+0pc@C+0pc{
\underset{r}{\mathrm{homotopylimit}}~\mathrm{Spec}^\mathrm{CS}{\Phi}^r_{*,\Gamma,\square},\underset{I}{\mathrm{homotopycolimit}}~\mathrm{Spec}^\mathrm{CS}{\Phi}^I_{*,\Gamma,\square}.	
}
\]
\[ 
\xymatrix@R+0pc@C+0pc{
\underset{r}{\mathrm{homotopylimit}}~\mathrm{Spec}^\mathrm{CS}\widetilde{\Phi}^r_{*,\Gamma,\square}/\mathrm{Fro}^\mathbb{Z},\underset{I}{\mathrm{homotopycolimit}}~\mathrm{Spec}^\mathrm{CS}\widetilde{\Phi}^I_{*,\Gamma,\square}/\mathrm{Fro}^\mathbb{Z},	\\
}
\]
\[ 
\xymatrix@R+0pc@C+0pc{
\underset{r}{\mathrm{homotopylimit}}~\mathrm{Spec}^\mathrm{CS}\breve{\Phi}^r_{*,\Gamma,\square}/\mathrm{Fro}^\mathbb{Z},\underset{I}{\mathrm{homotopycolimit}}~\breve{\Phi}^I_{*,\Gamma,\square}/\mathrm{Fro}^\mathbb{Z},	\\
}
\]
\[ 
\xymatrix@R+0pc@C+0pc{
\underset{r}{\mathrm{homotopylimit}}~\mathrm{Spec}^\mathrm{CS}{\Phi}^r_{*,\Gamma,\square}/\mathrm{Fro}^\mathbb{Z},\underset{I}{\mathrm{homotopycolimit}}~\mathrm{Spec}^\mathrm{CS}{\Phi}^I_{*,\Gamma,\square}/\mathrm{Fro}^\mathbb{Z}.	
}
\]	
In this situation we will have the target category being family parametrized by $r$ or $I$ in compatible glueing sense as in \cite[Definition 5.4.10]{10KL2}. In this situation for modules parametrized by the intervals we have the equivalence of $\infty$-categories by using \cite[Proposition 13.8]{10CS2}. Here the corresponding quasicoherent Frobenius modules are defined to be the corresponding homotopy colimits and limits of Frobenius modules:
\begin{align}
\underset{r}{\mathrm{homotopycolimit}}~M_r,\\
\underset{I}{\mathrm{homotopylimit}}~M_I,	
\end{align}
where each $M_r$ is a Frobenius-equivariant module over the period ring with respect to some radius $r$ while each $M_I$ is a Frobenius-equivariant module over the period ring with respect to some interval $I$.\\
\end{proposition}

\begin{proposition}
Similar proposition holds for 
\begin{align}
\mathrm{Quasicoherentsheaves,IndBanach}_{*}.	
\end{align}	
\end{proposition}

\

\begin{definition}
We then consider the corresponding quasipresheaves of perfect complexes the corresponding ind-Banach or monomorphic ind-Banach modules from \cite{10BBK}, \cite{10KKM}:
\begin{align}
\mathrm{Quasicoherentpresheaves,Perfectcomplex,IndBanach}_{*}	
\end{align}
where $*$ is one of the following spaces:
\begin{align}
&\mathrm{Spec}^\mathrm{BK}\widetilde{\Phi}_{*,\Gamma,\square}/\mathrm{Fro}^\mathbb{Z},	\\
\end{align}
\begin{align}
&\mathrm{Spec}^\mathrm{BK}\breve{\Phi}_{*,\Gamma,\square}/\mathrm{Fro}^\mathbb{Z},	\\
\end{align}
\begin{align}
&\mathrm{Spec}^\mathrm{BK}{\Phi}_{*,\Gamma,\square}/\mathrm{Fro}^\mathbb{Z}.	
\end{align}
Here for those space without notation related to the radius and the corresponding interval we consider the total unions $\bigcap_r,\bigcup_I$ in order to achieve the whole spaces to achieve the analogues of the corresponding FF curves from \cite{10KL1}, \cite{10KL2}, \cite{10FF} for
\[
\xymatrix@R+0pc@C+0pc{
\underset{r}{\mathrm{homotopylimit}}~\mathrm{Spec}^\mathrm{BK}\widetilde{\Phi}^r_{*,\Gamma,\square},\underset{I}{\mathrm{homotopycolimit}}~\mathrm{Spec}^\mathrm{BK}\widetilde{\Phi}^I_{*,\Gamma,\square},	\\
}
\]
\[
\xymatrix@R+0pc@C+0pc{
\underset{r}{\mathrm{homotopylimit}}~\mathrm{Spec}^\mathrm{BK}\breve{\Phi}^r_{*,\Gamma,\square},\underset{I}{\mathrm{homotopycolimit}}~\mathrm{Spec}^\mathrm{BK}\breve{\Phi}^I_{*,\Gamma,\square},	\\
}
\]
\[
\xymatrix@R+0pc@C+0pc{
\underset{r}{\mathrm{homotopylimit}}~\mathrm{Spec}^\mathrm{BK}{\Phi}^r_{*,\Gamma,\square},\underset{I}{\mathrm{homotopycolimit}}~\mathrm{Spec}^\mathrm{BK}{\Phi}^I_{*,\Gamma,\square}.	
}
\]
\[  
\xymatrix@R+0pc@C+0pc{
\underset{r}{\mathrm{homotopylimit}}~\mathrm{Spec}^\mathrm{BK}\widetilde{\Phi}^r_{*,\Gamma,\square}/\mathrm{Fro}^\mathbb{Z},\underset{I}{\mathrm{homotopycolimit}}~\mathrm{Spec}^\mathrm{BK}\widetilde{\Phi}^I_{*,\Gamma,\square}/\mathrm{Fro}^\mathbb{Z},	\\
}
\]
\[ 
\xymatrix@R+0pc@C+0pc{
\underset{r}{\mathrm{homotopylimit}}~\mathrm{Spec}^\mathrm{BK}\breve{\Phi}^r_{*,\Gamma,\square}/\mathrm{Fro}^\mathbb{Z},\underset{I}{\mathrm{homotopycolimit}}~\mathrm{Spec}^\mathrm{BK}\breve{\Phi}^I_{*,\Gamma,\square}/\mathrm{Fro}^\mathbb{Z},	\\
}
\]
\[ 
\xymatrix@R+0pc@C+0pc{
\underset{r}{\mathrm{homotopylimit}}~\mathrm{Spec}^\mathrm{BK}{\Phi}^r_{*,\Gamma,\square}/\mathrm{Fro}^\mathbb{Z},\underset{I}{\mathrm{homotopycolimit}}~\mathrm{Spec}^\mathrm{BK}{\Phi}^I_{*,\Gamma,\square}/\mathrm{Fro}^\mathbb{Z}.	
}
\]

\end{definition}

\begin{definition}
We then consider the corresponding quasisheaves of perfect complexes of the corresponding condensed solid topological modules from \cite{10CS2}:
\begin{align}
\mathrm{Quasicoherentsheaves, Perfectcomplex, Condensed}_{*}	
\end{align}
where $*$ is one of the following spaces:
\begin{align}
&\mathrm{Spec}^\mathrm{CS}\widetilde{\Delta}_{*,\Gamma,\square}/\mathrm{Fro}^\mathbb{Z},\mathrm{Spec}^\mathrm{CS}\widetilde{\nabla}_{*,\Gamma,\square}/\mathrm{Fro}^\mathbb{Z},\mathrm{Spec}^\mathrm{CS}\widetilde{\Phi}_{*,\Gamma,\square}/\mathrm{Fro}^\mathbb{Z},\mathrm{Spec}^\mathrm{CS}\widetilde{\Delta}^+_{*,\Gamma,\square}/\mathrm{Fro}^\mathbb{Z},\\
&\mathrm{Spec}^\mathrm{CS}\widetilde{\nabla}^+_{*,\Gamma,\square}/\mathrm{Fro}^\mathbb{Z},\mathrm{Spec}^\mathrm{CS}\widetilde{\Delta}^\dagger_{*,\Gamma,\square}/\mathrm{Fro}^\mathbb{Z},\mathrm{Spec}^\mathrm{CS}\widetilde{\nabla}^\dagger_{*,\Gamma,\square}/\mathrm{Fro}^\mathbb{Z},	\\
\end{align}
\begin{align}
&\mathrm{Spec}^\mathrm{CS}\breve{\Delta}_{*,\Gamma,\square}/\mathrm{Fro}^\mathbb{Z},\breve{\nabla}_{*,\Gamma,\square}/\mathrm{Fro}^\mathbb{Z},\mathrm{Spec}^\mathrm{CS}\breve{\Phi}_{*,\Gamma,\square}/\mathrm{Fro}^\mathbb{Z},\mathrm{Spec}^\mathrm{CS}\breve{\Delta}^+_{*,\Gamma,\square}/\mathrm{Fro}^\mathbb{Z},\\
&\mathrm{Spec}^\mathrm{CS}\breve{\nabla}^+_{*,\Gamma,\square}/\mathrm{Fro}^\mathbb{Z},\mathrm{Spec}^\mathrm{CS}\breve{\Delta}^\dagger_{*,\Gamma,\square}/\mathrm{Fro}^\mathbb{Z},\mathrm{Spec}^\mathrm{CS}\breve{\nabla}^\dagger_{*,\Gamma,\square}/\mathrm{Fro}^\mathbb{Z},	\\
\end{align}
\begin{align}
&\mathrm{Spec}^\mathrm{CS}{\Delta}_{*,\Gamma,\square}/\mathrm{Fro}^\mathbb{Z},\mathrm{Spec}^\mathrm{CS}{\nabla}_{*,\Gamma,\square}/\mathrm{Fro}^\mathbb{Z},\mathrm{Spec}^\mathrm{CS}{\Phi}_{*,\Gamma,\square}/\mathrm{Fro}^\mathbb{Z},\mathrm{Spec}^\mathrm{CS}{\Delta}^+_{*,\Gamma,\square}/\mathrm{Fro}^\mathbb{Z},\\
&\mathrm{Spec}^\mathrm{CS}{\nabla}^+_{*,\Gamma,\square}/\mathrm{Fro}^\mathbb{Z}, \mathrm{Spec}^\mathrm{CS}{\Delta}^\dagger_{*,\Gamma,\square}/\mathrm{Fro}^\mathbb{Z},\mathrm{Spec}^\mathrm{CS}{\nabla}^\dagger_{*,\Gamma,\square}/\mathrm{Fro}^\mathbb{Z}.	
\end{align}
Here for those space with notations related to the radius and the corresponding interval we consider the total unions $\bigcap_r,\bigcup_I$ in order to achieve the whole spaces to achieve the analogues of the corresponding FF curves from \cite{10KL1}, \cite{10KL2}, \cite{10FF} for
\[
\xymatrix@R+0pc@C+0pc{
\underset{r}{\mathrm{homotopylimit}}~\mathrm{Spec}^\mathrm{CS}\widetilde{\Phi}^r_{*,\Gamma,\square},\underset{I}{\mathrm{homotopycolimit}}~\mathrm{Spec}^\mathrm{CS}\widetilde{\Phi}^I_{*,\Gamma,\square},	\\
}
\]
\[
\xymatrix@R+0pc@C+0pc{
\underset{r}{\mathrm{homotopylimit}}~\mathrm{Spec}^\mathrm{CS}\breve{\Phi}^r_{*,\Gamma,\square},\underset{I}{\mathrm{homotopycolimit}}~\mathrm{Spec}^\mathrm{CS}\breve{\Phi}^I_{*,\Gamma,\square},	\\
}
\]
\[
\xymatrix@R+0pc@C+0pc{
\underset{r}{\mathrm{homotopylimit}}~\mathrm{Spec}^\mathrm{CS}{\Phi}^r_{*,\Gamma,\square},\underset{I}{\mathrm{homotopycolimit}}~\mathrm{Spec}^\mathrm{CS}{\Phi}^I_{*,\Gamma,\square}.	
}
\]
\[ 
\xymatrix@R+0pc@C+0pc{
\underset{r}{\mathrm{homotopylimit}}~\mathrm{Spec}^\mathrm{CS}\widetilde{\Phi}^r_{*,\Gamma,\square}/\mathrm{Fro}^\mathbb{Z},\underset{I}{\mathrm{homotopycolimit}}~\mathrm{Spec}^\mathrm{CS}\widetilde{\Phi}^I_{*,\Gamma,\square}/\mathrm{Fro}^\mathbb{Z},	\\
}
\]
\[ 
\xymatrix@R+0pc@C+0pc{
\underset{r}{\mathrm{homotopylimit}}~\mathrm{Spec}^\mathrm{CS}\breve{\Phi}^r_{*,\Gamma,\square}/\mathrm{Fro}^\mathbb{Z},\underset{I}{\mathrm{homotopycolimit}}~\breve{\Phi}^I_{*,\Gamma,\square}/\mathrm{Fro}^\mathbb{Z},	\\
}
\]
\[ 
\xymatrix@R+0pc@C+0pc{
\underset{r}{\mathrm{homotopylimit}}~\mathrm{Spec}^\mathrm{CS}{\Phi}^r_{*,\Gamma,\square}/\mathrm{Fro}^\mathbb{Z},\underset{I}{\mathrm{homotopycolimit}}~\mathrm{Spec}^\mathrm{CS}{\Phi}^I_{*,\Gamma,\square}/\mathrm{Fro}^\mathbb{Z}.	
}
\]

\end{definition}

\begin{proposition}
There is a well-defined functor from the $\infty$-category 
\begin{align}
\mathrm{Quasicoherentpresheaves,Perfectcomplex,Condensed}_{*}	
\end{align}
where $*$ is one of the following spaces:
\begin{align}
&\mathrm{Spec}^\mathrm{CS}\widetilde{\Phi}_{*,\Gamma,\square}/\mathrm{Fro}^\mathbb{Z},	\\
\end{align}
\begin{align}
&\mathrm{Spec}^\mathrm{CS}\breve{\Phi}_{*,\Gamma,\square}/\mathrm{Fro}^\mathbb{Z},	\\
\end{align}
\begin{align}
&\mathrm{Spec}^\mathrm{CS}{\Phi}_{*,\Gamma,\square}/\mathrm{Fro}^\mathbb{Z},	
\end{align}
to the $\infty$-category of $\mathrm{Fro}$-equivariant quasicoherent presheaves over similar spaces above correspondingly without the $\mathrm{Fro}$-quotients, and to the $\infty$-category of $\mathrm{Fro}$-equivariant quasicoherent modules over global sections of the structure $\infty$-sheaves of the similar spaces above correspondingly without the $\mathrm{Fro}$-quotients. Here for those space without notation related to the radius and the corresponding interval we consider the total unions $\bigcap_r,\bigcup_I$ in order to achieve the whole spaces to achieve the analogues of the corresponding FF curves from \cite{10KL1}, \cite{10KL2}, \cite{10FF} for
\[
\xymatrix@R+0pc@C+0pc{
\underset{r}{\mathrm{homotopylimit}}~\mathrm{Spec}^\mathrm{CS}\widetilde{\Phi}^r_{*,\Gamma,\square},\underset{I}{\mathrm{homotopycolimit}}~\mathrm{Spec}^\mathrm{CS}\widetilde{\Phi}^I_{*,\Gamma,\square},	\\
}
\]
\[
\xymatrix@R+0pc@C+0pc{
\underset{r}{\mathrm{homotopylimit}}~\mathrm{Spec}^\mathrm{CS}\breve{\Phi}^r_{*,\Gamma,\square},\underset{I}{\mathrm{homotopycolimit}}~\mathrm{Spec}^\mathrm{CS}\breve{\Phi}^I_{*,\Gamma,\square},	\\
}
\]
\[
\xymatrix@R+0pc@C+0pc{
\underset{r}{\mathrm{homotopylimit}}~\mathrm{Spec}^\mathrm{CS}{\Phi}^r_{*,\Gamma,\square},\underset{I}{\mathrm{homotopycolimit}}~\mathrm{Spec}^\mathrm{CS}{\Phi}^I_{*,\Gamma,\square}.	
}
\]
\[ 
\xymatrix@R+0pc@C+0pc{
\underset{r}{\mathrm{homotopylimit}}~\mathrm{Spec}^\mathrm{CS}\widetilde{\Phi}^r_{*,\Gamma,\square}/\mathrm{Fro}^\mathbb{Z},\underset{I}{\mathrm{homotopycolimit}}~\mathrm{Spec}^\mathrm{CS}\widetilde{\Phi}^I_{*,\Gamma,\square}/\mathrm{Fro}^\mathbb{Z},	\\
}
\]
\[ 
\xymatrix@R+0pc@C+0pc{
\underset{r}{\mathrm{homotopylimit}}~\mathrm{Spec}^\mathrm{CS}\breve{\Phi}^r_{*,\Gamma,\square}/\mathrm{Fro}^\mathbb{Z},\underset{I}{\mathrm{homotopycolimit}}~\breve{\Phi}^I_{*,\Gamma,\square}/\mathrm{Fro}^\mathbb{Z},	\\
}
\]
\[ 
\xymatrix@R+0pc@C+0pc{
\underset{r}{\mathrm{homotopylimit}}~\mathrm{Spec}^\mathrm{CS}{\Phi}^r_{*,\Gamma,\square}/\mathrm{Fro}^\mathbb{Z},\underset{I}{\mathrm{homotopycolimit}}~\mathrm{Spec}^\mathrm{CS}{\Phi}^I_{*,\Gamma,\square}/\mathrm{Fro}^\mathbb{Z}.	
}
\]	
In this situation we will have the target category being family parametrized by $r$ or $I$ in compatible glueing sense as in \cite[Definition 5.4.10]{10KL2}. In this situation for modules parametrized by the intervals we have the equivalence of $\infty$-categories by using \cite[Proposition 12.18]{10CS2}. Here the corresponding quasicoherent Frobenius modules are defined to be the corresponding homotopy colimits and limits of Frobenius modules:
\begin{align}
\underset{r}{\mathrm{homotopycolimit}}~M_r,\\
\underset{I}{\mathrm{homotopylimit}}~M_I,	
\end{align}
where each $M_r$ is a Frobenius-equivariant module over the period ring with respect to some radius $r$ while each $M_I$ is a Frobenius-equivariant module over the period ring with respect to some interval $I$.\\
\end{proposition}

\begin{proposition}
Similar proposition holds for 
\begin{align}
\mathrm{Quasicoherentsheaves,Perfectcomplex,IndBanach}_{*}.	
\end{align}	
\end{proposition}

\section{Univariate Hodge Iwasawa Prestacks}

This chapter follows closely \cite{10T1}, \cite{10T2}, \cite{10T3}, \cite{10KPX}, \cite{10KP}, \cite{10KL1}, \cite{10KL2}, \cite{10BK}, \cite{10BBBK}, \cite{10BBM}, \cite{10KKM}, \cite{10CS1}, \cite{10CS2}, \cite{10LBV}.

\subsection{Frobenius Quasicoherent Prestacks I}

\begin{definition}
We now consider the pro-\'etale site of $\mathrm{Spa}\mathbb{Q}_p\left<X_1^{\pm 1},...,X_k^{\pm 1}\right>$, denote that by $*$. To be more accurate we replace one component for $\Gamma$ with the pro-\'etale site of $\mathrm{Spa}\mathbb{Q}_p\left<X_1^{\pm 1},...,X_k^{\pm 1}\right>$. And we treat then all the functor to be prestacks for this site. Then from \cite{10KL1} and \cite[Definition 5.2.1]{10KL2} we have the following class of Kedlaya-Liu rings (with the following replacement: $\Delta$ stands for $A$, $\nabla$ stands for $B$, while $\Phi$ stands for $C$) by taking product in the sense of self $\Gamma$-th power\footnote{Here $|\Gamma|=1$.}:

\[
\xymatrix@R+0pc@C+0pc{
\widetilde{\Delta}_{*},\widetilde{\nabla}_{*},\widetilde{\Phi}_{*},\widetilde{\Delta}^+_{*},\widetilde{\nabla}^+_{*},\widetilde{\Delta}^\dagger_{*},\widetilde{\nabla}^\dagger_{*},\widetilde{\Phi}^r_{*},\widetilde{\Phi}^I_{*}, 
}
\]

\[
\xymatrix@R+0pc@C+0pc{
\breve{\Delta}_{*},\breve{\nabla}_{*},\breve{\Phi}_{*},\breve{\Delta}^+_{*},\breve{\nabla}^+_{*},\breve{\Delta}^\dagger_{*},\breve{\nabla}^\dagger_{*},\breve{\Phi}^r_{*},\breve{\Phi}^I_{*},	
}
\]

\[
\xymatrix@R+0pc@C+0pc{
{\Delta}_{*},{\nabla}_{*},{\Phi}_{*},{\Delta}^+_{*},{\nabla}^+_{*},{\Delta}^\dagger_{*},{\nabla}^\dagger_{*},{\Phi}^r_{*},{\Phi}^I_{*}.	
}
\]
We now consider the following rings with $A$ being a Banach ring over $\mathbb{Q}_p$. Taking the product we have:
\[
\xymatrix@R+0pc@C+0pc{
\widetilde{\Phi}_{*,A},\widetilde{\Phi}^r_{*,A},\widetilde{\Phi}^I_{*,A},	
}
\]
\[
\xymatrix@R+0pc@C+0pc{
\breve{\Phi}_{*,A},\breve{\Phi}^r_{*,A},\breve{\Phi}^I_{*,A},	
}
\]
\[
\xymatrix@R+0pc@C+0pc{
{\Phi}_{*,A},{\Phi}^r_{*,A},{\Phi}^I_{*,A}.	
}
\]
They carry multi Frobenius action $\varphi_\Gamma$ and multi $\mathrm{Lie}_\Gamma:=\mathbb{Z}_p^{\times\Gamma}$ action. In our current situation after \cite{10CKZ} and \cite{10PZ} we consider the following $(\infty,1)$-categories of $(\infty,1)$-modules.\\
\end{definition}

\begin{definition}
First we consider the Bambozzi-Kremnizer spectrum $\mathrm{Spec}^\mathrm{BK}(*)$ attached to any of those in the above from \cite{10BK} by taking derived rational localization:
\begin{align}
&\mathrm{Spec}^\mathrm{BK}\widetilde{\Phi}_{*,A},\mathrm{Spec}^\mathrm{BK}\widetilde{\Phi}^r_{*,A},\mathrm{Spec}^\mathrm{BK}\widetilde{\Phi}^I_{*,A},	
\end{align}
\begin{align}
&\mathrm{Spec}^\mathrm{BK}\breve{\Phi}_{*,A},\mathrm{Spec}^\mathrm{BK}\breve{\Phi}^r_{*,A},\mathrm{Spec}^\mathrm{BK}\breve{\Phi}^I_{*,A},	
\end{align}
\begin{align}
&\mathrm{Spec}^\mathrm{BK}{\Phi}_{*,A},
\mathrm{Spec}^\mathrm{BK}{\Phi}^r_{*,A},\mathrm{Spec}^\mathrm{BK}{\Phi}^I_{*,A}.	
\end{align}

Then we take the corresponding quotients by using the corresponding Frobenius operators:
\begin{align}
&\mathrm{Spec}^\mathrm{BK}\widetilde{\Phi}_{*,A}/\mathrm{Fro}^\mathbb{Z},	\\
\end{align}
\begin{align}
&\mathrm{Spec}^\mathrm{BK}\breve{\Phi}_{*,A}/\mathrm{Fro}^\mathbb{Z},	\\
\end{align}
\begin{align}
&\mathrm{Spec}^\mathrm{BK}{\Phi}_{*,A}/\mathrm{Fro}^\mathbb{Z}.	
\end{align}
Here for those space without notation related to the radius and the corresponding interval we consider the total unions $\bigcap_r,\bigcup_I$ in order to achieve the whole spaces to achieve the analogues of the corresponding FF curves from \cite{10KL1}, \cite{10KL2}, \cite{10FF} for
\[
\xymatrix@R+0pc@C+0pc{
\underset{r}{\mathrm{homotopylimit}}~\mathrm{Spec}^\mathrm{BK}\widetilde{\Phi}^r_{*,A},\underset{I}{\mathrm{homotopycolimit}}~\mathrm{Spec}^\mathrm{BK}\widetilde{\Phi}^I_{*,A},	\\
}
\]
\[
\xymatrix@R+0pc@C+0pc{
\underset{r}{\mathrm{homotopylimit}}~\mathrm{Spec}^\mathrm{BK}\breve{\Phi}^r_{*,A},\underset{I}{\mathrm{homotopycolimit}}~\mathrm{Spec}^\mathrm{BK}\breve{\Phi}^I_{*,A},	\\
}
\]
\[
\xymatrix@R+0pc@C+0pc{
\underset{r}{\mathrm{homotopylimit}}~\mathrm{Spec}^\mathrm{BK}{\Phi}^r_{*,A},\underset{I}{\mathrm{homotopycolimit}}~\mathrm{Spec}^\mathrm{BK}{\Phi}^I_{*,A}.	
}
\]
\[  
\xymatrix@R+0pc@C+0pc{
\underset{r}{\mathrm{homotopylimit}}~\mathrm{Spec}^\mathrm{BK}\widetilde{\Phi}^r_{*,A}/\mathrm{Fro}^\mathbb{Z},\underset{I}{\mathrm{homotopycolimit}}~\mathrm{Spec}^\mathrm{BK}\widetilde{\Phi}^I_{*,A}/\mathrm{Fro}^\mathbb{Z},	\\
}
\]
\[ 
\xymatrix@R+0pc@C+0pc{
\underset{r}{\mathrm{homotopylimit}}~\mathrm{Spec}^\mathrm{BK}\breve{\Phi}^r_{*,A}/\mathrm{Fro}^\mathbb{Z},\underset{I}{\mathrm{homotopycolimit}}~\mathrm{Spec}^\mathrm{BK}\breve{\Phi}^I_{*,A}/\mathrm{Fro}^\mathbb{Z},	\\
}
\]
\[ 
\xymatrix@R+0pc@C+0pc{
\underset{r}{\mathrm{homotopylimit}}~\mathrm{Spec}^\mathrm{BK}{\Phi}^r_{*,A}/\mathrm{Fro}^\mathbb{Z},\underset{I}{\mathrm{homotopycolimit}}~\mathrm{Spec}^\mathrm{BK}{\Phi}^I_{*,A}/\mathrm{Fro}^\mathbb{Z}.	
}
\]

\end{definition}

\indent Meanwhile we have the corresponding Clausen-Scholze analytic stacks from \cite{10CS2}, therefore applying their construction we have:

\begin{definition}
Here we define the following products by using the solidified tensor product from \cite{10CS1} and \cite{10CS2}. Namely $A$ will still as above as a Banach ring over $\mathbb{Q}_p$. Then we take solidified tensor product $\overset{\blacksquare}{\otimes}$ of any of the following
\[
\xymatrix@R+0pc@C+0pc{
\widetilde{\Delta}_{*},\widetilde{\nabla}_{*},\widetilde{\Phi}_{*},\widetilde{\Delta}^+_{*},\widetilde{\nabla}^+_{*},\widetilde{\Delta}^\dagger_{*},\widetilde{\nabla}^\dagger_{*},\widetilde{\Phi}^r_{*},\widetilde{\Phi}^I_{*}, 
}
\]

\[
\xymatrix@R+0pc@C+0pc{
\breve{\Delta}_{*},\breve{\nabla}_{*},\breve{\Phi}_{*},\breve{\Delta}^+_{*},\breve{\nabla}^+_{*},\breve{\Delta}^\dagger_{*},\breve{\nabla}^\dagger_{*},\breve{\Phi}^r_{*},\breve{\Phi}^I_{*},	
}
\]

\[
\xymatrix@R+0pc@C+0pc{
{\Delta}_{*},{\nabla}_{*},{\Phi}_{*},{\Delta}^+_{*},{\nabla}^+_{*},{\Delta}^\dagger_{*},{\nabla}^\dagger_{*},{\Phi}^r_{*},{\Phi}^I_{*},	
}
\]  	
with $A$. Then we have the notations:
\[
\xymatrix@R+0pc@C+0pc{
\widetilde{\Delta}_{*,A},\widetilde{\nabla}_{*,A},\widetilde{\Phi}_{*,A},\widetilde{\Delta}^+_{*,A},\widetilde{\nabla}^+_{*,A},\widetilde{\Delta}^\dagger_{*,A},\widetilde{\nabla}^\dagger_{*,A},\widetilde{\Phi}^r_{*,A},\widetilde{\Phi}^I_{*,A}, 
}
\]

\[
\xymatrix@R+0pc@C+0pc{
\breve{\Delta}_{*,A},\breve{\nabla}_{*,A},\breve{\Phi}_{*,A},\breve{\Delta}^+_{*,A},\breve{\nabla}^+_{*,A},\breve{\Delta}^\dagger_{*,A},\breve{\nabla}^\dagger_{*,A},\breve{\Phi}^r_{*,A},\breve{\Phi}^I_{*,A},	
}
\]

\[
\xymatrix@R+0pc@C+0pc{
{\Delta}_{*,A},{\nabla}_{*,A},{\Phi}_{*,A},{\Delta}^+_{*,A},{\nabla}^+_{*,A},{\Delta}^\dagger_{*,A},{\nabla}^\dagger_{*,A},{\Phi}^r_{*,A},{\Phi}^I_{*,A}.	
}
\]
\end{definition}

\begin{definition}
First we consider the Clausen-Scholze spectrum $\mathrm{Spec}^\mathrm{CS}(*)$ attached to any of those in the above from \cite{10CS2} by taking derived rational localization:
\begin{align}
\mathrm{Spec}^\mathrm{CS}\widetilde{\Delta}_{*,A},\mathrm{Spec}^\mathrm{CS}\widetilde{\nabla}_{*,A},\mathrm{Spec}^\mathrm{CS}\widetilde{\Phi}_{*,A},\mathrm{Spec}^\mathrm{CS}\widetilde{\Delta}^+_{*,A},\mathrm{Spec}^\mathrm{CS}\widetilde{\nabla}^+_{*,A},\\
\mathrm{Spec}^\mathrm{CS}\widetilde{\Delta}^\dagger_{*,A},\mathrm{Spec}^\mathrm{CS}\widetilde{\nabla}^\dagger_{*,A},\mathrm{Spec}^\mathrm{CS}\widetilde{\Phi}^r_{*,A},\mathrm{Spec}^\mathrm{CS}\widetilde{\Phi}^I_{*,A},	\\
\end{align}
\begin{align}
\mathrm{Spec}^\mathrm{CS}\breve{\Delta}_{*,A},\breve{\nabla}_{*,A},\mathrm{Spec}^\mathrm{CS}\breve{\Phi}_{*,A},\mathrm{Spec}^\mathrm{CS}\breve{\Delta}^+_{*,A},\mathrm{Spec}^\mathrm{CS}\breve{\nabla}^+_{*,A},\\
\mathrm{Spec}^\mathrm{CS}\breve{\Delta}^\dagger_{*,A},\mathrm{Spec}^\mathrm{CS}\breve{\nabla}^\dagger_{*,A},\mathrm{Spec}^\mathrm{CS}\breve{\Phi}^r_{*,A},\breve{\Phi}^I_{*,A},	\\
\end{align}
\begin{align}
\mathrm{Spec}^\mathrm{CS}{\Delta}_{*,A},\mathrm{Spec}^\mathrm{CS}{\nabla}_{*,A},\mathrm{Spec}^\mathrm{CS}{\Phi}_{*,A},\mathrm{Spec}^\mathrm{CS}{\Delta}^+_{*,A},\mathrm{Spec}^\mathrm{CS}{\nabla}^+_{*,A},\\
\mathrm{Spec}^\mathrm{CS}{\Delta}^\dagger_{*,A},\mathrm{Spec}^\mathrm{CS}{\nabla}^\dagger_{*,A},\mathrm{Spec}^\mathrm{CS}{\Phi}^r_{*,A},\mathrm{Spec}^\mathrm{CS}{\Phi}^I_{*,A}.	
\end{align}

Then we take the corresponding quotients by using the corresponding Frobenius operators:
\begin{align}
&\mathrm{Spec}^\mathrm{CS}\widetilde{\Delta}_{*,A}/\mathrm{Fro}^\mathbb{Z},\mathrm{Spec}^\mathrm{CS}\widetilde{\nabla}_{*,A}/\mathrm{Fro}^\mathbb{Z},\mathrm{Spec}^\mathrm{CS}\widetilde{\Phi}_{*,A}/\mathrm{Fro}^\mathbb{Z},\mathrm{Spec}^\mathrm{CS}\widetilde{\Delta}^+_{*,A}/\mathrm{Fro}^\mathbb{Z},\\
&\mathrm{Spec}^\mathrm{CS}\widetilde{\nabla}^+_{*,A}/\mathrm{Fro}^\mathbb{Z}, \mathrm{Spec}^\mathrm{CS}\widetilde{\Delta}^\dagger_{*,A}/\mathrm{Fro}^\mathbb{Z},\mathrm{Spec}^\mathrm{CS}\widetilde{\nabla}^\dagger_{*,A}/\mathrm{Fro}^\mathbb{Z},	\\
\end{align}
\begin{align}
&\mathrm{Spec}^\mathrm{CS}\breve{\Delta}_{*,A}/\mathrm{Fro}^\mathbb{Z},\breve{\nabla}_{*,A}/\mathrm{Fro}^\mathbb{Z},\mathrm{Spec}^\mathrm{CS}\breve{\Phi}_{*,A}/\mathrm{Fro}^\mathbb{Z},\mathrm{Spec}^\mathrm{CS}\breve{\Delta}^+_{*,A}/\mathrm{Fro}^\mathbb{Z},\\
&\mathrm{Spec}^\mathrm{CS}\breve{\nabla}^+_{*,A}/\mathrm{Fro}^\mathbb{Z}, \mathrm{Spec}^\mathrm{CS}\breve{\Delta}^\dagger_{*,A}/\mathrm{Fro}^\mathbb{Z},\mathrm{Spec}^\mathrm{CS}\breve{\nabla}^\dagger_{*,A}/\mathrm{Fro}^\mathbb{Z},	\\
\end{align}
\begin{align}
&\mathrm{Spec}^\mathrm{CS}{\Delta}_{*,A}/\mathrm{Fro}^\mathbb{Z},\mathrm{Spec}^\mathrm{CS}{\nabla}_{*,A}/\mathrm{Fro}^\mathbb{Z},\mathrm{Spec}^\mathrm{CS}{\Phi}_{*,A}/\mathrm{Fro}^\mathbb{Z},\mathrm{Spec}^\mathrm{CS}{\Delta}^+_{*,A}/\mathrm{Fro}^\mathbb{Z},\\
&\mathrm{Spec}^\mathrm{CS}{\nabla}^+_{*,A}/\mathrm{Fro}^\mathbb{Z}, \mathrm{Spec}^\mathrm{CS}{\Delta}^\dagger_{*,A}/\mathrm{Fro}^\mathbb{Z},\mathrm{Spec}^\mathrm{CS}{\nabla}^\dagger_{*,A}/\mathrm{Fro}^\mathbb{Z}.	
\end{align}
Here for those space with notations related to the radius and the corresponding interval we consider the total unions $\bigcap_r,\bigcup_I$ in order to achieve the whole spaces to achieve the analogues of the corresponding FF curves from \cite{10KL1}, \cite{10KL2}, \cite{10FF} for
\[
\xymatrix@R+0pc@C+0pc{
\underset{r}{\mathrm{homotopylimit}}~\mathrm{Spec}^\mathrm{CS}\widetilde{\Phi}^r_{*,A},\underset{I}{\mathrm{homotopycolimit}}~\mathrm{Spec}^\mathrm{CS}\widetilde{\Phi}^I_{*,A},	\\
}
\]
\[
\xymatrix@R+0pc@C+0pc{
\underset{r}{\mathrm{homotopylimit}}~\mathrm{Spec}^\mathrm{CS}\breve{\Phi}^r_{*,A},\underset{I}{\mathrm{homotopycolimit}}~\mathrm{Spec}^\mathrm{CS}\breve{\Phi}^I_{*,A},	\\
}
\]
\[
\xymatrix@R+0pc@C+0pc{
\underset{r}{\mathrm{homotopylimit}}~\mathrm{Spec}^\mathrm{CS}{\Phi}^r_{*,A},\underset{I}{\mathrm{homotopycolimit}}~\mathrm{Spec}^\mathrm{CS}{\Phi}^I_{*,A}.	
}
\]
\[ 
\xymatrix@R+0pc@C+0pc{
\underset{r}{\mathrm{homotopylimit}}~\mathrm{Spec}^\mathrm{CS}\widetilde{\Phi}^r_{*,A}/\mathrm{Fro}^\mathbb{Z},\underset{I}{\mathrm{homotopycolimit}}~\mathrm{Spec}^\mathrm{CS}\widetilde{\Phi}^I_{*,A}/\mathrm{Fro}^\mathbb{Z},	\\
}
\]
\[ 
\xymatrix@R+0pc@C+0pc{
\underset{r}{\mathrm{homotopylimit}}~\mathrm{Spec}^\mathrm{CS}\breve{\Phi}^r_{*,A}/\mathrm{Fro}^\mathbb{Z},\underset{I}{\mathrm{homotopycolimit}}~\breve{\Phi}^I_{*,A}/\mathrm{Fro}^\mathbb{Z},	\\
}
\]
\[ 
\xymatrix@R+0pc@C+0pc{
\underset{r}{\mathrm{homotopylimit}}~\mathrm{Spec}^\mathrm{CS}{\Phi}^r_{*,A}/\mathrm{Fro}^\mathbb{Z},\underset{I}{\mathrm{homotopycolimit}}~\mathrm{Spec}^\mathrm{CS}{\Phi}^I_{*,A}/\mathrm{Fro}^\mathbb{Z}.	
}
\]

\end{definition}

\

\begin{definition}
We then consider the corresponding quasipresheaves of the corresponding ind-Banach or monomorphic ind-Banach modules from \cite{10BBK}, \cite{10KKM}:
\begin{align}
\mathrm{Quasicoherentpresheaves,IndBanach}_{*}	
\end{align}
where $*$ is one of the following spaces:
\begin{align}
&\mathrm{Spec}^\mathrm{BK}\widetilde{\Phi}_{*,A}/\mathrm{Fro}^\mathbb{Z},	\\
\end{align}
\begin{align}
&\mathrm{Spec}^\mathrm{BK}\breve{\Phi}_{*,A}/\mathrm{Fro}^\mathbb{Z},	\\
\end{align}
\begin{align}
&\mathrm{Spec}^\mathrm{BK}{\Phi}_{*,A}/\mathrm{Fro}^\mathbb{Z}.	
\end{align}
Here for those space without notation related to the radius and the corresponding interval we consider the total unions $\bigcap_r,\bigcup_I$ in order to achieve the whole spaces to achieve the analogues of the corresponding FF curves from \cite{10KL1}, \cite{10KL2}, \cite{10FF} for
\[
\xymatrix@R+0pc@C+0pc{
\underset{r}{\mathrm{homotopylimit}}~\mathrm{Spec}^\mathrm{BK}\widetilde{\Phi}^r_{*,A},\underset{I}{\mathrm{homotopycolimit}}~\mathrm{Spec}^\mathrm{BK}\widetilde{\Phi}^I_{*,A},	\\
}
\]
\[
\xymatrix@R+0pc@C+0pc{
\underset{r}{\mathrm{homotopylimit}}~\mathrm{Spec}^\mathrm{BK}\breve{\Phi}^r_{*,A},\underset{I}{\mathrm{homotopycolimit}}~\mathrm{Spec}^\mathrm{BK}\breve{\Phi}^I_{*,A},	\\
}
\]
\[
\xymatrix@R+0pc@C+0pc{
\underset{r}{\mathrm{homotopylimit}}~\mathrm{Spec}^\mathrm{BK}{\Phi}^r_{*,A},\underset{I}{\mathrm{homotopycolimit}}~\mathrm{Spec}^\mathrm{BK}{\Phi}^I_{*,A}.	
}
\]
\[  
\xymatrix@R+0pc@C+0pc{
\underset{r}{\mathrm{homotopylimit}}~\mathrm{Spec}^\mathrm{BK}\widetilde{\Phi}^r_{*,A}/\mathrm{Fro}^\mathbb{Z},\underset{I}{\mathrm{homotopycolimit}}~\mathrm{Spec}^\mathrm{BK}\widetilde{\Phi}^I_{*,A}/\mathrm{Fro}^\mathbb{Z},	\\
}
\]
\[ 
\xymatrix@R+0pc@C+0pc{
\underset{r}{\mathrm{homotopylimit}}~\mathrm{Spec}^\mathrm{BK}\breve{\Phi}^r_{*,A}/\mathrm{Fro}^\mathbb{Z},\underset{I}{\mathrm{homotopycolimit}}~\mathrm{Spec}^\mathrm{BK}\breve{\Phi}^I_{*,A}/\mathrm{Fro}^\mathbb{Z},	\\
}
\]
\[ 
\xymatrix@R+0pc@C+0pc{
\underset{r}{\mathrm{homotopylimit}}~\mathrm{Spec}^\mathrm{BK}{\Phi}^r_{*,A}/\mathrm{Fro}^\mathbb{Z},\underset{I}{\mathrm{homotopycolimit}}~\mathrm{Spec}^\mathrm{BK}{\Phi}^I_{*,A}/\mathrm{Fro}^\mathbb{Z}.	
}
\]

\end{definition}

\begin{definition}
We then consider the corresponding quasisheaves of the corresponding condensed solid topological modules from \cite{10CS2}:
\begin{align}
\mathrm{Quasicoherentsheaves, Condensed}_{*}	
\end{align}
where $*$ is one of the following spaces:
\begin{align}
&\mathrm{Spec}^\mathrm{CS}\widetilde{\Delta}_{*,A}/\mathrm{Fro}^\mathbb{Z},\mathrm{Spec}^\mathrm{CS}\widetilde{\nabla}_{*,A}/\mathrm{Fro}^\mathbb{Z},\mathrm{Spec}^\mathrm{CS}\widetilde{\Phi}_{*,A}/\mathrm{Fro}^\mathbb{Z},\mathrm{Spec}^\mathrm{CS}\widetilde{\Delta}^+_{*,A}/\mathrm{Fro}^\mathbb{Z},\\
&\mathrm{Spec}^\mathrm{CS}\widetilde{\nabla}^+_{*,A}/\mathrm{Fro}^\mathbb{Z},\mathrm{Spec}^\mathrm{CS}\widetilde{\Delta}^\dagger_{*,A}/\mathrm{Fro}^\mathbb{Z},\mathrm{Spec}^\mathrm{CS}\widetilde{\nabla}^\dagger_{*,A}/\mathrm{Fro}^\mathbb{Z},	\\
\end{align}
\begin{align}
&\mathrm{Spec}^\mathrm{CS}\breve{\Delta}_{*,A}/\mathrm{Fro}^\mathbb{Z},\breve{\nabla}_{*,A}/\mathrm{Fro}^\mathbb{Z},\mathrm{Spec}^\mathrm{CS}\breve{\Phi}_{*,A}/\mathrm{Fro}^\mathbb{Z},\mathrm{Spec}^\mathrm{CS}\breve{\Delta}^+_{*,A}/\mathrm{Fro}^\mathbb{Z},\\
&\mathrm{Spec}^\mathrm{CS}\breve{\nabla}^+_{*,A}/\mathrm{Fro}^\mathbb{Z},\mathrm{Spec}^\mathrm{CS}\breve{\Delta}^\dagger_{*,A}/\mathrm{Fro}^\mathbb{Z},\mathrm{Spec}^\mathrm{CS}\breve{\nabla}^\dagger_{*,A}/\mathrm{Fro}^\mathbb{Z},	\\
\end{align}
\begin{align}
&\mathrm{Spec}^\mathrm{CS}{\Delta}_{*,A}/\mathrm{Fro}^\mathbb{Z},\mathrm{Spec}^\mathrm{CS}{\nabla}_{*,A}/\mathrm{Fro}^\mathbb{Z},\mathrm{Spec}^\mathrm{CS}{\Phi}_{*,A}/\mathrm{Fro}^\mathbb{Z},\mathrm{Spec}^\mathrm{CS}{\Delta}^+_{*,A}/\mathrm{Fro}^\mathbb{Z},\\
&\mathrm{Spec}^\mathrm{CS}{\nabla}^+_{*,A}/\mathrm{Fro}^\mathbb{Z}, \mathrm{Spec}^\mathrm{CS}{\Delta}^\dagger_{*,A}/\mathrm{Fro}^\mathbb{Z},\mathrm{Spec}^\mathrm{CS}{\nabla}^\dagger_{*,A}/\mathrm{Fro}^\mathbb{Z}.	
\end{align}
Here for those space with notations related to the radius and the corresponding interval we consider the total unions $\bigcap_r,\bigcup_I$ in order to achieve the whole spaces to achieve the analogues of the corresponding FF curves from \cite{10KL1}, \cite{10KL2}, \cite{10FF} for
\[
\xymatrix@R+0pc@C+0pc{
\underset{r}{\mathrm{homotopylimit}}~\mathrm{Spec}^\mathrm{CS}\widetilde{\Phi}^r_{*,A},\underset{I}{\mathrm{homotopycolimit}}~\mathrm{Spec}^\mathrm{CS}\widetilde{\Phi}^I_{*,A},	\\
}
\]
\[
\xymatrix@R+0pc@C+0pc{
\underset{r}{\mathrm{homotopylimit}}~\mathrm{Spec}^\mathrm{CS}\breve{\Phi}^r_{*,A},\underset{I}{\mathrm{homotopycolimit}}~\mathrm{Spec}^\mathrm{CS}\breve{\Phi}^I_{*,A},	\\
}
\]
\[
\xymatrix@R+0pc@C+0pc{
\underset{r}{\mathrm{homotopylimit}}~\mathrm{Spec}^\mathrm{CS}{\Phi}^r_{*,A},\underset{I}{\mathrm{homotopycolimit}}~\mathrm{Spec}^\mathrm{CS}{\Phi}^I_{*,A}.	
}
\]
\[ 
\xymatrix@R+0pc@C+0pc{
\underset{r}{\mathrm{homotopylimit}}~\mathrm{Spec}^\mathrm{CS}\widetilde{\Phi}^r_{*,A}/\mathrm{Fro}^\mathbb{Z},\underset{I}{\mathrm{homotopycolimit}}~\mathrm{Spec}^\mathrm{CS}\widetilde{\Phi}^I_{*,A}/\mathrm{Fro}^\mathbb{Z},	\\
}
\]
\[ 
\xymatrix@R+0pc@C+0pc{
\underset{r}{\mathrm{homotopylimit}}~\mathrm{Spec}^\mathrm{CS}\breve{\Phi}^r_{*,A}/\mathrm{Fro}^\mathbb{Z},\underset{I}{\mathrm{homotopycolimit}}~\breve{\Phi}^I_{*,A}/\mathrm{Fro}^\mathbb{Z},	\\
}
\]
\[ 
\xymatrix@R+0pc@C+0pc{
\underset{r}{\mathrm{homotopylimit}}~\mathrm{Spec}^\mathrm{CS}{\Phi}^r_{*,A}/\mathrm{Fro}^\mathbb{Z},\underset{I}{\mathrm{homotopycolimit}}~\mathrm{Spec}^\mathrm{CS}{\Phi}^I_{*,A}/\mathrm{Fro}^\mathbb{Z}.	
}
\]

\end{definition}

\

\begin{proposition}
There is a well-defined functor from the $\infty$-category 
\begin{align}
\mathrm{Quasicoherentpresheaves,Condensed}_{*}	
\end{align}
where $*$ is one of the following spaces:
\begin{align}
&\mathrm{Spec}^\mathrm{CS}\widetilde{\Phi}_{*,A}/\mathrm{Fro}^\mathbb{Z},	\\
\end{align}
\begin{align}
&\mathrm{Spec}^\mathrm{CS}\breve{\Phi}_{*,A}/\mathrm{Fro}^\mathbb{Z},	\\
\end{align}
\begin{align}
&\mathrm{Spec}^\mathrm{CS}{\Phi}_{*,A}/\mathrm{Fro}^\mathbb{Z},	
\end{align}
to the $\infty$-category of $\mathrm{Fro}$-equivariant quasicoherent presheaves over similar spaces above correspondingly without the $\mathrm{Fro}$-quotients, and to the $\infty$-category of $\mathrm{Fro}$-equivariant quasicoherent modules over global sections of the structure $\infty$-sheaves of the similar spaces above correspondingly without the $\mathrm{Fro}$-quotients. Here for those space without notation related to the radius and the corresponding interval we consider the total unions $\bigcap_r,\bigcup_I$ in order to achieve the whole spaces to achieve the analogues of the corresponding FF curves from \cite{10KL1}, \cite{10KL2}, \cite{10FF} for
\[
\xymatrix@R+0pc@C+0pc{
\underset{r}{\mathrm{homotopylimit}}~\mathrm{Spec}^\mathrm{CS}\widetilde{\Phi}^r_{*,A},\underset{I}{\mathrm{homotopycolimit}}~\mathrm{Spec}^\mathrm{CS}\widetilde{\Phi}^I_{*,A},	\\
}
\]
\[
\xymatrix@R+0pc@C+0pc{
\underset{r}{\mathrm{homotopylimit}}~\mathrm{Spec}^\mathrm{CS}\breve{\Phi}^r_{*,A},\underset{I}{\mathrm{homotopycolimit}}~\mathrm{Spec}^\mathrm{CS}\breve{\Phi}^I_{*,A},	\\
}
\]
\[
\xymatrix@R+0pc@C+0pc{
\underset{r}{\mathrm{homotopylimit}}~\mathrm{Spec}^\mathrm{CS}{\Phi}^r_{*,A},\underset{I}{\mathrm{homotopycolimit}}~\mathrm{Spec}^\mathrm{CS}{\Phi}^I_{*,A}.	
}
\]
\[ 
\xymatrix@R+0pc@C+0pc{
\underset{r}{\mathrm{homotopylimit}}~\mathrm{Spec}^\mathrm{CS}\widetilde{\Phi}^r_{*,A}/\mathrm{Fro}^\mathbb{Z},\underset{I}{\mathrm{homotopycolimit}}~\mathrm{Spec}^\mathrm{CS}\widetilde{\Phi}^I_{*,A}/\mathrm{Fro}^\mathbb{Z},	\\
}
\]
\[ 
\xymatrix@R+0pc@C+0pc{
\underset{r}{\mathrm{homotopylimit}}~\mathrm{Spec}^\mathrm{CS}\breve{\Phi}^r_{*,A}/\mathrm{Fro}^\mathbb{Z},\underset{I}{\mathrm{homotopycolimit}}~\breve{\Phi}^I_{*,A}/\mathrm{Fro}^\mathbb{Z},	\\
}
\]
\[ 
\xymatrix@R+0pc@C+0pc{
\underset{r}{\mathrm{homotopylimit}}~\mathrm{Spec}^\mathrm{CS}{\Phi}^r_{*,A}/\mathrm{Fro}^\mathbb{Z},\underset{I}{\mathrm{homotopycolimit}}~\mathrm{Spec}^\mathrm{CS}{\Phi}^I_{*,A}/\mathrm{Fro}^\mathbb{Z}.	
}
\]	
In this situation we will have the target category being family parametrized by $r$ or $I$ in compatible glueing sense as in \cite[Definition 5.4.10]{10KL2}. In this situation for modules parametrized by the intervals we have the equivalence of $\infty$-categories by using \cite[Proposition 13.8]{10CS2}. Here the corresponding quasicoherent Frobenius modules are defined to be the corresponding homotopy colimits and limits of Frobenius modules:
\begin{align}
\underset{r}{\mathrm{homotopycolimit}}~M_r,\\
\underset{I}{\mathrm{homotopylimit}}~M_I,	
\end{align}
where each $M_r$ is a Frobenius-equivariant module over the period ring with respect to some radius $r$ while each $M_I$ is a Frobenius-equivariant module over the period ring with respect to some interval $I$.\\
\end{proposition}

\begin{proposition}
Similar proposition holds for 
\begin{align}
\mathrm{Quasicoherentsheaves,IndBanach}_{*}.	
\end{align}	
\end{proposition}

\

\begin{definition}
We then consider the corresponding quasipresheaves of perfect complexes the corresponding ind-Banach or monomorphic ind-Banach modules from \cite{10BBK}, \cite{10KKM}:
\begin{align}
\mathrm{Quasicoherentpresheaves,Perfectcomplex,IndBanach}_{*}	
\end{align}
where $*$ is one of the following spaces:
\begin{align}
&\mathrm{Spec}^\mathrm{BK}\widetilde{\Phi}_{*,A}/\mathrm{Fro}^\mathbb{Z},	\\
\end{align}
\begin{align}
&\mathrm{Spec}^\mathrm{BK}\breve{\Phi}_{*,A}/\mathrm{Fro}^\mathbb{Z},	\\
\end{align}
\begin{align}
&\mathrm{Spec}^\mathrm{BK}{\Phi}_{*,A}/\mathrm{Fro}^\mathbb{Z}.	
\end{align}
Here for those space without notation related to the radius and the corresponding interval we consider the total unions $\bigcap_r,\bigcup_I$ in order to achieve the whole spaces to achieve the analogues of the corresponding FF curves from \cite{10KL1}, \cite{10KL2}, \cite{10FF} for
\[
\xymatrix@R+0pc@C+0pc{
\underset{r}{\mathrm{homotopylimit}}~\mathrm{Spec}^\mathrm{BK}\widetilde{\Phi}^r_{*,A},\underset{I}{\mathrm{homotopycolimit}}~\mathrm{Spec}^\mathrm{BK}\widetilde{\Phi}^I_{*,A},	\\
}
\]
\[
\xymatrix@R+0pc@C+0pc{
\underset{r}{\mathrm{homotopylimit}}~\mathrm{Spec}^\mathrm{BK}\breve{\Phi}^r_{*,A},\underset{I}{\mathrm{homotopycolimit}}~\mathrm{Spec}^\mathrm{BK}\breve{\Phi}^I_{*,A},	\\
}
\]
\[
\xymatrix@R+0pc@C+0pc{
\underset{r}{\mathrm{homotopylimit}}~\mathrm{Spec}^\mathrm{BK}{\Phi}^r_{*,A},\underset{I}{\mathrm{homotopycolimit}}~\mathrm{Spec}^\mathrm{BK}{\Phi}^I_{*,A}.	
}
\]
\[  
\xymatrix@R+0pc@C+0pc{
\underset{r}{\mathrm{homotopylimit}}~\mathrm{Spec}^\mathrm{BK}\widetilde{\Phi}^r_{*,A}/\mathrm{Fro}^\mathbb{Z},\underset{I}{\mathrm{homotopycolimit}}~\mathrm{Spec}^\mathrm{BK}\widetilde{\Phi}^I_{*,A}/\mathrm{Fro}^\mathbb{Z},	\\
}
\]
\[ 
\xymatrix@R+0pc@C+0pc{
\underset{r}{\mathrm{homotopylimit}}~\mathrm{Spec}^\mathrm{BK}\breve{\Phi}^r_{*,A}/\mathrm{Fro}^\mathbb{Z},\underset{I}{\mathrm{homotopycolimit}}~\mathrm{Spec}^\mathrm{BK}\breve{\Phi}^I_{*,A}/\mathrm{Fro}^\mathbb{Z},	\\
}
\]
\[ 
\xymatrix@R+0pc@C+0pc{
\underset{r}{\mathrm{homotopylimit}}~\mathrm{Spec}^\mathrm{BK}{\Phi}^r_{*,A}/\mathrm{Fro}^\mathbb{Z},\underset{I}{\mathrm{homotopycolimit}}~\mathrm{Spec}^\mathrm{BK}{\Phi}^I_{*,A}/\mathrm{Fro}^\mathbb{Z}.	
}
\]

\end{definition}

\begin{definition}
We then consider the corresponding quasisheaves of perfect complexes of the corresponding condensed solid topological modules from \cite{10CS2}:
\begin{align}
\mathrm{Quasicoherentsheaves, Perfectcomplex, Condensed}_{*}	
\end{align}
where $*$ is one of the following spaces:
\begin{align}
&\mathrm{Spec}^\mathrm{CS}\widetilde{\Delta}_{*,A}/\mathrm{Fro}^\mathbb{Z},\mathrm{Spec}^\mathrm{CS}\widetilde{\nabla}_{*,A}/\mathrm{Fro}^\mathbb{Z},\mathrm{Spec}^\mathrm{CS}\widetilde{\Phi}_{*,A}/\mathrm{Fro}^\mathbb{Z},\mathrm{Spec}^\mathrm{CS}\widetilde{\Delta}^+_{*,A}/\mathrm{Fro}^\mathbb{Z},\\
&\mathrm{Spec}^\mathrm{CS}\widetilde{\nabla}^+_{*,A}/\mathrm{Fro}^\mathbb{Z},\mathrm{Spec}^\mathrm{CS}\widetilde{\Delta}^\dagger_{*,A}/\mathrm{Fro}^\mathbb{Z},\mathrm{Spec}^\mathrm{CS}\widetilde{\nabla}^\dagger_{*,A}/\mathrm{Fro}^\mathbb{Z},	\\
\end{align}
\begin{align}
&\mathrm{Spec}^\mathrm{CS}\breve{\Delta}_{*,A}/\mathrm{Fro}^\mathbb{Z},\breve{\nabla}_{*,A}/\mathrm{Fro}^\mathbb{Z},\mathrm{Spec}^\mathrm{CS}\breve{\Phi}_{*,A}/\mathrm{Fro}^\mathbb{Z},\mathrm{Spec}^\mathrm{CS}\breve{\Delta}^+_{*,A}/\mathrm{Fro}^\mathbb{Z},\\
&\mathrm{Spec}^\mathrm{CS}\breve{\nabla}^+_{*,A}/\mathrm{Fro}^\mathbb{Z},\mathrm{Spec}^\mathrm{CS}\breve{\Delta}^\dagger_{*,A}/\mathrm{Fro}^\mathbb{Z},\mathrm{Spec}^\mathrm{CS}\breve{\nabla}^\dagger_{*,A}/\mathrm{Fro}^\mathbb{Z},	\\
\end{align}
\begin{align}
&\mathrm{Spec}^\mathrm{CS}{\Delta}_{*,A}/\mathrm{Fro}^\mathbb{Z},\mathrm{Spec}^\mathrm{CS}{\nabla}_{*,A}/\mathrm{Fro}^\mathbb{Z},\mathrm{Spec}^\mathrm{CS}{\Phi}_{*,A}/\mathrm{Fro}^\mathbb{Z},\mathrm{Spec}^\mathrm{CS}{\Delta}^+_{*,A}/\mathrm{Fro}^\mathbb{Z},\\
&\mathrm{Spec}^\mathrm{CS}{\nabla}^+_{*,A}/\mathrm{Fro}^\mathbb{Z}, \mathrm{Spec}^\mathrm{CS}{\Delta}^\dagger_{*,A}/\mathrm{Fro}^\mathbb{Z},\mathrm{Spec}^\mathrm{CS}{\nabla}^\dagger_{*,A}/\mathrm{Fro}^\mathbb{Z}.	
\end{align}
Here for those space with notations related to the radius and the corresponding interval we consider the total unions $\bigcap_r,\bigcup_I$ in order to achieve the whole spaces to achieve the analogues of the corresponding FF curves from \cite{10KL1}, \cite{10KL2}, \cite{10FF} for
\[
\xymatrix@R+0pc@C+0pc{
\underset{r}{\mathrm{homotopylimit}}~\mathrm{Spec}^\mathrm{CS}\widetilde{\Phi}^r_{*,A},\underset{I}{\mathrm{homotopycolimit}}~\mathrm{Spec}^\mathrm{CS}\widetilde{\Phi}^I_{*,A},	\\
}
\]
\[
\xymatrix@R+0pc@C+0pc{
\underset{r}{\mathrm{homotopylimit}}~\mathrm{Spec}^\mathrm{CS}\breve{\Phi}^r_{*,A},\underset{I}{\mathrm{homotopycolimit}}~\mathrm{Spec}^\mathrm{CS}\breve{\Phi}^I_{*,A},	\\
}
\]
\[
\xymatrix@R+0pc@C+0pc{
\underset{r}{\mathrm{homotopylimit}}~\mathrm{Spec}^\mathrm{CS}{\Phi}^r_{*,A},\underset{I}{\mathrm{homotopycolimit}}~\mathrm{Spec}^\mathrm{CS}{\Phi}^I_{*,A}.	
}
\]
\[ 
\xymatrix@R+0pc@C+0pc{
\underset{r}{\mathrm{homotopylimit}}~\mathrm{Spec}^\mathrm{CS}\widetilde{\Phi}^r_{*,A}/\mathrm{Fro}^\mathbb{Z},\underset{I}{\mathrm{homotopycolimit}}~\mathrm{Spec}^\mathrm{CS}\widetilde{\Phi}^I_{*,A}/\mathrm{Fro}^\mathbb{Z},	\\
}
\]
\[ 
\xymatrix@R+0pc@C+0pc{
\underset{r}{\mathrm{homotopylimit}}~\mathrm{Spec}^\mathrm{CS}\breve{\Phi}^r_{*,A}/\mathrm{Fro}^\mathbb{Z},\underset{I}{\mathrm{homotopycolimit}}~\breve{\Phi}^I_{*,A}/\mathrm{Fro}^\mathbb{Z},	\\
}
\]
\[ 
\xymatrix@R+0pc@C+0pc{
\underset{r}{\mathrm{homotopylimit}}~\mathrm{Spec}^\mathrm{CS}{\Phi}^r_{*,A}/\mathrm{Fro}^\mathbb{Z},\underset{I}{\mathrm{homotopycolimit}}~\mathrm{Spec}^\mathrm{CS}{\Phi}^I_{*,A}/\mathrm{Fro}^\mathbb{Z}.	
}
\]

\end{definition}

\begin{proposition}
There is a well-defined functor from the $\infty$-category 
\begin{align}
\mathrm{Quasicoherentpresheaves,Perfectcomplex,Condensed}_{*}	
\end{align}
where $*$ is one of the following spaces:
\begin{align}
&\mathrm{Spec}^\mathrm{CS}\widetilde{\Phi}_{*,A}/\mathrm{Fro}^\mathbb{Z},	\\
\end{align}
\begin{align}
&\mathrm{Spec}^\mathrm{CS}\breve{\Phi}_{*,A}/\mathrm{Fro}^\mathbb{Z},	\\
\end{align}
\begin{align}
&\mathrm{Spec}^\mathrm{CS}{\Phi}_{*,A}/\mathrm{Fro}^\mathbb{Z},	
\end{align}
to the $\infty$-category of $\mathrm{Fro}$-equivariant quasicoherent presheaves over similar spaces above correspondingly without the $\mathrm{Fro}$-quotients, and to the $\infty$-category of $\mathrm{Fro}$-equivariant quasicoherent modules over global sections of the structure $\infty$-sheaves of the similar spaces above correspondingly without the $\mathrm{Fro}$-quotients. Here for those space without notation related to the radius and the corresponding interval we consider the total unions $\bigcap_r,\bigcup_I$ in order to achieve the whole spaces to achieve the analogues of the corresponding FF curves from \cite{10KL1}, \cite{10KL2}, \cite{10FF} for
\[
\xymatrix@R+0pc@C+0pc{
\underset{r}{\mathrm{homotopylimit}}~\mathrm{Spec}^\mathrm{CS}\widetilde{\Phi}^r_{*,A},\underset{I}{\mathrm{homotopycolimit}}~\mathrm{Spec}^\mathrm{CS}\widetilde{\Phi}^I_{*,A},	\\
}
\]
\[
\xymatrix@R+0pc@C+0pc{
\underset{r}{\mathrm{homotopylimit}}~\mathrm{Spec}^\mathrm{CS}\breve{\Phi}^r_{*,A},\underset{I}{\mathrm{homotopycolimit}}~\mathrm{Spec}^\mathrm{CS}\breve{\Phi}^I_{*,A},	\\
}
\]
\[
\xymatrix@R+0pc@C+0pc{
\underset{r}{\mathrm{homotopylimit}}~\mathrm{Spec}^\mathrm{CS}{\Phi}^r_{*,A},\underset{I}{\mathrm{homotopycolimit}}~\mathrm{Spec}^\mathrm{CS}{\Phi}^I_{*,A}.	
}
\]
\[ 
\xymatrix@R+0pc@C+0pc{
\underset{r}{\mathrm{homotopylimit}}~\mathrm{Spec}^\mathrm{CS}\widetilde{\Phi}^r_{*,A}/\mathrm{Fro}^\mathbb{Z},\underset{I}{\mathrm{homotopycolimit}}~\mathrm{Spec}^\mathrm{CS}\widetilde{\Phi}^I_{*,A}/\mathrm{Fro}^\mathbb{Z},	\\
}
\]
\[ 
\xymatrix@R+0pc@C+0pc{
\underset{r}{\mathrm{homotopylimit}}~\mathrm{Spec}^\mathrm{CS}\breve{\Phi}^r_{*,A}/\mathrm{Fro}^\mathbb{Z},\underset{I}{\mathrm{homotopycolimit}}~\breve{\Phi}^I_{*,A}/\mathrm{Fro}^\mathbb{Z},	\\
}
\]
\[ 
\xymatrix@R+0pc@C+0pc{
\underset{r}{\mathrm{homotopylimit}}~\mathrm{Spec}^\mathrm{CS}{\Phi}^r_{*,A}/\mathrm{Fro}^\mathbb{Z},\underset{I}{\mathrm{homotopycolimit}}~\mathrm{Spec}^\mathrm{CS}{\Phi}^I_{*,A}/\mathrm{Fro}^\mathbb{Z}.	
}
\]	
In this situation we will have the target category being family parametrized by $r$ or $I$ in compatible glueing sense as in \cite[Definition 5.4.10]{10KL2}. In this situation for modules parametrized by the intervals we have the equivalence of $\infty$-categories by using \cite[Proposition 12.18]{10CS2}. Here the corresponding quasicoherent Frobenius modules are defined to be the corresponding homotopy colimits and limits of Frobenius modules:
\begin{align}
\underset{r}{\mathrm{homotopycolimit}}~M_r,\\
\underset{I}{\mathrm{homotopylimit}}~M_I,	
\end{align}
where each $M_r$ is a Frobenius-equivariant module over the period ring with respect to some radius $r$ while each $M_I$ is a Frobenius-equivariant module over the period ring with respect to some interval $I$.\\
\end{proposition}

\begin{proposition}
Similar proposition holds for 
\begin{align}
\mathrm{Quasicoherentsheaves,Perfectcomplex,IndBanach}_{*}.	
\end{align}	
\end{proposition}

\newpage
\subsection{Frobenius Quasicoherent Prestacks II: Deformation in Banach Rings}

\begin{definition}
We now consider the pro-\'etale site of $\mathrm{Spa}\mathbb{Q}_p\left<X_1^{\pm 1},...,X_k^{\pm 1}\right>$, denote that by $*$. To be more accurate we replace one component for $\Gamma$ with the pro-\'etale site of $\mathrm{Spa}\mathbb{Q}_p\left<X_1^{\pm 1},...,X_k^{\pm 1}\right>$. And we treat then all the functor to be prestacks for this site. Then from \cite{10KL1} and \cite[Definition 5.2.1]{10KL2} we have the following class of Kedlaya-Liu rings (with the following replacement: $\Delta$ stands for $A$, $\nabla$ stands for $B$, while $\Phi$ stands for $C$) by taking product in the sense of self $\Gamma$-th power\footnote{Here $|\Gamma|=1$.}:

\[
\xymatrix@R+0pc@C+0pc{
\widetilde{\Delta}_{*},\widetilde{\nabla}_{*},\widetilde{\Phi}_{*},\widetilde{\Delta}^+_{*},\widetilde{\nabla}^+_{*},\widetilde{\Delta}^\dagger_{*},\widetilde{\nabla}^\dagger_{*},\widetilde{\Phi}^r_{*},\widetilde{\Phi}^I_{*}, 
}
\]

\[
\xymatrix@R+0pc@C+0pc{
\breve{\Delta}_{*},\breve{\nabla}_{*},\breve{\Phi}_{*},\breve{\Delta}^+_{*},\breve{\nabla}^+_{*},\breve{\Delta}^\dagger_{*},\breve{\nabla}^\dagger_{*},\breve{\Phi}^r_{*},\breve{\Phi}^I_{*},	
}
\]

\[
\xymatrix@R+0pc@C+0pc{
{\Delta}_{*},{\nabla}_{*},{\Phi}_{*},{\Delta}^+_{*},{\nabla}^+_{*},{\Delta}^\dagger_{*},{\nabla}^\dagger_{*},{\Phi}^r_{*},{\Phi}^I_{*}.	
}
\]
We now consider the following rings with $-$ being any deforming Banach ring over $\mathbb{Q}_p$. Taking the product we have:
\[
\xymatrix@R+0pc@C+0pc{
\widetilde{\Phi}_{*,-},\widetilde{\Phi}^r_{*,-},\widetilde{\Phi}^I_{*,-},	
}
\]
\[
\xymatrix@R+0pc@C+0pc{
\breve{\Phi}_{*,-},\breve{\Phi}^r_{*,-},\breve{\Phi}^I_{*,-},	
}
\]
\[
\xymatrix@R+0pc@C+0pc{
{\Phi}_{*,-},{\Phi}^r_{*,-},{\Phi}^I_{*,-}.	
}
\]
They carry multi Frobenius action $\varphi_\Gamma$ and multi $\mathrm{Lie}_\Gamma:=\mathbb{Z}_p^{\times\Gamma}$ action. In our current situation after \cite{10CKZ} and \cite{10PZ} we consider the following $(\infty,1)$-categories of $(\infty,1)$-modules.\\
\end{definition}

\begin{definition}
First we consider the Bambozzi-Kremnizer spectrum $\mathrm{Spec}^\mathrm{BK}(*)$ attached to any of those in the above from \cite{10BK} by taking derived rational localization:
\begin{align}
&\mathrm{Spec}^\mathrm{BK}\widetilde{\Phi}_{*,-},\mathrm{Spec}^\mathrm{BK}\widetilde{\Phi}^r_{*,-},\mathrm{Spec}^\mathrm{BK}\widetilde{\Phi}^I_{*,-},	
\end{align}
\begin{align}
&\mathrm{Spec}^\mathrm{BK}\breve{\Phi}_{*,-},\mathrm{Spec}^\mathrm{BK}\breve{\Phi}^r_{*,-},\mathrm{Spec}^\mathrm{BK}\breve{\Phi}^I_{*,-},	
\end{align}
\begin{align}
&\mathrm{Spec}^\mathrm{BK}{\Phi}_{*,-},
\mathrm{Spec}^\mathrm{BK}{\Phi}^r_{*,-},\mathrm{Spec}^\mathrm{BK}{\Phi}^I_{*,-}.	
\end{align}

Then we take the corresponding quotients by using the corresponding Frobenius operators:
\begin{align}
&\mathrm{Spec}^\mathrm{BK}\widetilde{\Phi}_{*,-}/\mathrm{Fro}^\mathbb{Z},	\\
\end{align}
\begin{align}
&\mathrm{Spec}^\mathrm{BK}\breve{\Phi}_{*,-}/\mathrm{Fro}^\mathbb{Z},	\\
\end{align}
\begin{align}
&\mathrm{Spec}^\mathrm{BK}{\Phi}_{*,-}/\mathrm{Fro}^\mathbb{Z}.	
\end{align}
Here for those space without notation related to the radius and the corresponding interval we consider the total unions $\bigcap_r,\bigcup_I$ in order to achieve the whole spaces to achieve the analogues of the corresponding FF curves from \cite{10KL1}, \cite{10KL2}, \cite{10FF} for
\[
\xymatrix@R+0pc@C+0pc{
\underset{r}{\mathrm{homotopylimit}}~\mathrm{Spec}^\mathrm{BK}\widetilde{\Phi}^r_{*,-},\underset{I}{\mathrm{homotopycolimit}}~\mathrm{Spec}^\mathrm{BK}\widetilde{\Phi}^I_{*,-},	\\
}
\]
\[
\xymatrix@R+0pc@C+0pc{
\underset{r}{\mathrm{homotopylimit}}~\mathrm{Spec}^\mathrm{BK}\breve{\Phi}^r_{*,-},\underset{I}{\mathrm{homotopycolimit}}~\mathrm{Spec}^\mathrm{BK}\breve{\Phi}^I_{*,-},	\\
}
\]
\[
\xymatrix@R+0pc@C+0pc{
\underset{r}{\mathrm{homotopylimit}}~\mathrm{Spec}^\mathrm{BK}{\Phi}^r_{*,-},\underset{I}{\mathrm{homotopycolimit}}~\mathrm{Spec}^\mathrm{BK}{\Phi}^I_{*,-}.	
}
\]
\[  
\xymatrix@R+0pc@C+0pc{
\underset{r}{\mathrm{homotopylimit}}~\mathrm{Spec}^\mathrm{BK}\widetilde{\Phi}^r_{*,-}/\mathrm{Fro}^\mathbb{Z},\underset{I}{\mathrm{homotopycolimit}}~\mathrm{Spec}^\mathrm{BK}\widetilde{\Phi}^I_{*,-}/\mathrm{Fro}^\mathbb{Z},	\\
}
\]
\[ 
\xymatrix@R+0pc@C+0pc{
\underset{r}{\mathrm{homotopylimit}}~\mathrm{Spec}^\mathrm{BK}\breve{\Phi}^r_{*,-}/\mathrm{Fro}^\mathbb{Z},\underset{I}{\mathrm{homotopycolimit}}~\mathrm{Spec}^\mathrm{BK}\breve{\Phi}^I_{*,-}/\mathrm{Fro}^\mathbb{Z},	\\
}
\]
\[ 
\xymatrix@R+0pc@C+0pc{
\underset{r}{\mathrm{homotopylimit}}~\mathrm{Spec}^\mathrm{BK}{\Phi}^r_{*,-}/\mathrm{Fro}^\mathbb{Z},\underset{I}{\mathrm{homotopycolimit}}~\mathrm{Spec}^\mathrm{BK}{\Phi}^I_{*,-}/\mathrm{Fro}^\mathbb{Z}.	
}
\]

\end{definition}

\indent Meanwhile we have the corresponding Clausen-Scholze analytic stacks from \cite{10CS2}, therefore applying their construction we have:

\begin{definition}
Here we define the following products by using the solidified tensor product from \cite{10CS1} and \cite{10CS2}. Namely $A$ will still as above as a Banach ring over $\mathbb{Q}_p$. Then we take solidified tensor product $\overset{\blacksquare}{\otimes}$ of any of the following
\[
\xymatrix@R+0pc@C+0pc{
\widetilde{\Delta}_{*},\widetilde{\nabla}_{*},\widetilde{\Phi}_{*},\widetilde{\Delta}^+_{*},\widetilde{\nabla}^+_{*},\widetilde{\Delta}^\dagger_{*},\widetilde{\nabla}^\dagger_{*},\widetilde{\Phi}^r_{*},\widetilde{\Phi}^I_{*}, 
}
\]

\[
\xymatrix@R+0pc@C+0pc{
\breve{\Delta}_{*},\breve{\nabla}_{*},\breve{\Phi}_{*},\breve{\Delta}^+_{*},\breve{\nabla}^+_{*},\breve{\Delta}^\dagger_{*},\breve{\nabla}^\dagger_{*},\breve{\Phi}^r_{*},\breve{\Phi}^I_{*},	
}
\]

\[
\xymatrix@R+0pc@C+0pc{
{\Delta}_{*},{\nabla}_{*},{\Phi}_{*},{\Delta}^+_{*},{\nabla}^+_{*},{\Delta}^\dagger_{*},{\nabla}^\dagger_{*},{\Phi}^r_{*},{\Phi}^I_{*},	
}
\]  	
with $A$. Then we have the notations:
\[
\xymatrix@R+0pc@C+0pc{
\widetilde{\Delta}_{*,-},\widetilde{\nabla}_{*,-},\widetilde{\Phi}_{*,-},\widetilde{\Delta}^+_{*,-},\widetilde{\nabla}^+_{*,-},\widetilde{\Delta}^\dagger_{*,-},\widetilde{\nabla}^\dagger_{*,-},\widetilde{\Phi}^r_{*,-},\widetilde{\Phi}^I_{*,-}, 
}
\]

\[
\xymatrix@R+0pc@C+0pc{
\breve{\Delta}_{*,-},\breve{\nabla}_{*,-},\breve{\Phi}_{*,-},\breve{\Delta}^+_{*,-},\breve{\nabla}^+_{*,-},\breve{\Delta}^\dagger_{*,-},\breve{\nabla}^\dagger_{*,-},\breve{\Phi}^r_{*,-},\breve{\Phi}^I_{*,-},	
}
\]

\[
\xymatrix@R+0pc@C+0pc{
{\Delta}_{*,-},{\nabla}_{*,-},{\Phi}_{*,-},{\Delta}^+_{*,-},{\nabla}^+_{*,-},{\Delta}^\dagger_{*,-},{\nabla}^\dagger_{*,-},{\Phi}^r_{*,-},{\Phi}^I_{*,-}.	
}
\]
\end{definition}

\begin{definition}
First we consider the Clausen-Scholze spectrum $\mathrm{Spec}^\mathrm{CS}(*)$ attached to any of those in the above from \cite{10CS2} by taking derived rational localization:
\begin{align}
\mathrm{Spec}^\mathrm{CS}\widetilde{\Delta}_{*,-},\mathrm{Spec}^\mathrm{CS}\widetilde{\nabla}_{*,-},\mathrm{Spec}^\mathrm{CS}\widetilde{\Phi}_{*,-},\mathrm{Spec}^\mathrm{CS}\widetilde{\Delta}^+_{*,-},\mathrm{Spec}^\mathrm{CS}\widetilde{\nabla}^+_{*,-},\\
\mathrm{Spec}^\mathrm{CS}\widetilde{\Delta}^\dagger_{*,-},\mathrm{Spec}^\mathrm{CS}\widetilde{\nabla}^\dagger_{*,-},\mathrm{Spec}^\mathrm{CS}\widetilde{\Phi}^r_{*,-},\mathrm{Spec}^\mathrm{CS}\widetilde{\Phi}^I_{*,-},	\\
\end{align}
\begin{align}
\mathrm{Spec}^\mathrm{CS}\breve{\Delta}_{*,-},\breve{\nabla}_{*,-},\mathrm{Spec}^\mathrm{CS}\breve{\Phi}_{*,-},\mathrm{Spec}^\mathrm{CS}\breve{\Delta}^+_{*,-},\mathrm{Spec}^\mathrm{CS}\breve{\nabla}^+_{*,-},\\
\mathrm{Spec}^\mathrm{CS}\breve{\Delta}^\dagger_{*,-},\mathrm{Spec}^\mathrm{CS}\breve{\nabla}^\dagger_{*,-},\mathrm{Spec}^\mathrm{CS}\breve{\Phi}^r_{*,-},\breve{\Phi}^I_{*,-},	\\
\end{align}
\begin{align}
\mathrm{Spec}^\mathrm{CS}{\Delta}_{*,-},\mathrm{Spec}^\mathrm{CS}{\nabla}_{*,-},\mathrm{Spec}^\mathrm{CS}{\Phi}_{*,-},\mathrm{Spec}^\mathrm{CS}{\Delta}^+_{*,-},\mathrm{Spec}^\mathrm{CS}{\nabla}^+_{*,-},\\
\mathrm{Spec}^\mathrm{CS}{\Delta}^\dagger_{*,-},\mathrm{Spec}^\mathrm{CS}{\nabla}^\dagger_{*,-},\mathrm{Spec}^\mathrm{CS}{\Phi}^r_{*,-},\mathrm{Spec}^\mathrm{CS}{\Phi}^I_{*,-}.	
\end{align}

Then we take the corresponding quotients by using the corresponding Frobenius operators:
\begin{align}
&\mathrm{Spec}^\mathrm{CS}\widetilde{\Delta}_{*,-}/\mathrm{Fro}^\mathbb{Z},\mathrm{Spec}^\mathrm{CS}\widetilde{\nabla}_{*,-}/\mathrm{Fro}^\mathbb{Z},\mathrm{Spec}^\mathrm{CS}\widetilde{\Phi}_{*,-}/\mathrm{Fro}^\mathbb{Z},\mathrm{Spec}^\mathrm{CS}\widetilde{\Delta}^+_{*,-}/\mathrm{Fro}^\mathbb{Z},\\
&\mathrm{Spec}^\mathrm{CS}\widetilde{\nabla}^+_{*,-}/\mathrm{Fro}^\mathbb{Z}, \mathrm{Spec}^\mathrm{CS}\widetilde{\Delta}^\dagger_{*,-}/\mathrm{Fro}^\mathbb{Z},\mathrm{Spec}^\mathrm{CS}\widetilde{\nabla}^\dagger_{*,-}/\mathrm{Fro}^\mathbb{Z},	\\
\end{align}
\begin{align}
&\mathrm{Spec}^\mathrm{CS}\breve{\Delta}_{*,-}/\mathrm{Fro}^\mathbb{Z},\breve{\nabla}_{*,-}/\mathrm{Fro}^\mathbb{Z},\mathrm{Spec}^\mathrm{CS}\breve{\Phi}_{*,-}/\mathrm{Fro}^\mathbb{Z},\mathrm{Spec}^\mathrm{CS}\breve{\Delta}^+_{*,-}/\mathrm{Fro}^\mathbb{Z},\\
&\mathrm{Spec}^\mathrm{CS}\breve{\nabla}^+_{*,-}/\mathrm{Fro}^\mathbb{Z}, \mathrm{Spec}^\mathrm{CS}\breve{\Delta}^\dagger_{*,-}/\mathrm{Fro}^\mathbb{Z},\mathrm{Spec}^\mathrm{CS}\breve{\nabla}^\dagger_{*,-}/\mathrm{Fro}^\mathbb{Z},	\\
\end{align}
\begin{align}
&\mathrm{Spec}^\mathrm{CS}{\Delta}_{*,-}/\mathrm{Fro}^\mathbb{Z},\mathrm{Spec}^\mathrm{CS}{\nabla}_{*,-}/\mathrm{Fro}^\mathbb{Z},\mathrm{Spec}^\mathrm{CS}{\Phi}_{*,-}/\mathrm{Fro}^\mathbb{Z},\mathrm{Spec}^\mathrm{CS}{\Delta}^+_{*,-}/\mathrm{Fro}^\mathbb{Z},\\
&\mathrm{Spec}^\mathrm{CS}{\nabla}^+_{*,-}/\mathrm{Fro}^\mathbb{Z}, \mathrm{Spec}^\mathrm{CS}{\Delta}^\dagger_{*,-}/\mathrm{Fro}^\mathbb{Z},\mathrm{Spec}^\mathrm{CS}{\nabla}^\dagger_{*,-}/\mathrm{Fro}^\mathbb{Z}.	
\end{align}
Here for those space with notations related to the radius and the corresponding interval we consider the total unions $\bigcap_r,\bigcup_I$ in order to achieve the whole spaces to achieve the analogues of the corresponding FF curves from \cite{10KL1}, \cite{10KL2}, \cite{10FF} for
\[
\xymatrix@R+0pc@C+0pc{
\underset{r}{\mathrm{homotopylimit}}~\mathrm{Spec}^\mathrm{CS}\widetilde{\Phi}^r_{*,-},\underset{I}{\mathrm{homotopycolimit}}~\mathrm{Spec}^\mathrm{CS}\widetilde{\Phi}^I_{*,-},	\\
}
\]
\[
\xymatrix@R+0pc@C+0pc{
\underset{r}{\mathrm{homotopylimit}}~\mathrm{Spec}^\mathrm{CS}\breve{\Phi}^r_{*,-},\underset{I}{\mathrm{homotopycolimit}}~\mathrm{Spec}^\mathrm{CS}\breve{\Phi}^I_{*,-},	\\
}
\]
\[
\xymatrix@R+0pc@C+0pc{
\underset{r}{\mathrm{homotopylimit}}~\mathrm{Spec}^\mathrm{CS}{\Phi}^r_{*,-},\underset{I}{\mathrm{homotopycolimit}}~\mathrm{Spec}^\mathrm{CS}{\Phi}^I_{*,-}.	
}
\]
\[ 
\xymatrix@R+0pc@C+0pc{
\underset{r}{\mathrm{homotopylimit}}~\mathrm{Spec}^\mathrm{CS}\widetilde{\Phi}^r_{*,-}/\mathrm{Fro}^\mathbb{Z},\underset{I}{\mathrm{homotopycolimit}}~\mathrm{Spec}^\mathrm{CS}\widetilde{\Phi}^I_{*,-}/\mathrm{Fro}^\mathbb{Z},	\\
}
\]
\[ 
\xymatrix@R+0pc@C+0pc{
\underset{r}{\mathrm{homotopylimit}}~\mathrm{Spec}^\mathrm{CS}\breve{\Phi}^r_{*,-}/\mathrm{Fro}^\mathbb{Z},\underset{I}{\mathrm{homotopycolimit}}~\breve{\Phi}^I_{*,-}/\mathrm{Fro}^\mathbb{Z},	\\
}
\]
\[ 
\xymatrix@R+0pc@C+0pc{
\underset{r}{\mathrm{homotopylimit}}~\mathrm{Spec}^\mathrm{CS}{\Phi}^r_{*,-}/\mathrm{Fro}^\mathbb{Z},\underset{I}{\mathrm{homotopycolimit}}~\mathrm{Spec}^\mathrm{CS}{\Phi}^I_{*,-}/\mathrm{Fro}^\mathbb{Z}.	
}
\]

\end{definition}

\

\begin{definition}
We then consider the corresponding quasipresheaves of the corresponding ind-Banach or monomorphic ind-Banach modules from \cite{10BBK}, \cite{10KKM}:
\begin{align}
\mathrm{Quasicoherentpresheaves,IndBanach}_{*}	
\end{align}
where $*$ is one of the following spaces:
\begin{align}
&\mathrm{Spec}^\mathrm{BK}\widetilde{\Phi}_{*,-}/\mathrm{Fro}^\mathbb{Z},	\\
\end{align}
\begin{align}
&\mathrm{Spec}^\mathrm{BK}\breve{\Phi}_{*,-}/\mathrm{Fro}^\mathbb{Z},	\\
\end{align}
\begin{align}
&\mathrm{Spec}^\mathrm{BK}{\Phi}_{*,-}/\mathrm{Fro}^\mathbb{Z}.	
\end{align}
Here for those space without notation related to the radius and the corresponding interval we consider the total unions $\bigcap_r,\bigcup_I$ in order to achieve the whole spaces to achieve the analogues of the corresponding FF curves from \cite{10KL1}, \cite{10KL2}, \cite{10FF} for
\[
\xymatrix@R+0pc@C+0pc{
\underset{r}{\mathrm{homotopylimit}}~\mathrm{Spec}^\mathrm{BK}\widetilde{\Phi}^r_{*,-},\underset{I}{\mathrm{homotopycolimit}}~\mathrm{Spec}^\mathrm{BK}\widetilde{\Phi}^I_{*,-},	\\
}
\]
\[
\xymatrix@R+0pc@C+0pc{
\underset{r}{\mathrm{homotopylimit}}~\mathrm{Spec}^\mathrm{BK}\breve{\Phi}^r_{*,-},\underset{I}{\mathrm{homotopycolimit}}~\mathrm{Spec}^\mathrm{BK}\breve{\Phi}^I_{*,-},	\\
}
\]
\[
\xymatrix@R+0pc@C+0pc{
\underset{r}{\mathrm{homotopylimit}}~\mathrm{Spec}^\mathrm{BK}{\Phi}^r_{*,-},\underset{I}{\mathrm{homotopycolimit}}~\mathrm{Spec}^\mathrm{BK}{\Phi}^I_{*,-}.	
}
\]
\[  
\xymatrix@R+0pc@C+0pc{
\underset{r}{\mathrm{homotopylimit}}~\mathrm{Spec}^\mathrm{BK}\widetilde{\Phi}^r_{*,-}/\mathrm{Fro}^\mathbb{Z},\underset{I}{\mathrm{homotopycolimit}}~\mathrm{Spec}^\mathrm{BK}\widetilde{\Phi}^I_{*,-}/\mathrm{Fro}^\mathbb{Z},	\\
}
\]
\[ 
\xymatrix@R+0pc@C+0pc{
\underset{r}{\mathrm{homotopylimit}}~\mathrm{Spec}^\mathrm{BK}\breve{\Phi}^r_{*,-}/\mathrm{Fro}^\mathbb{Z},\underset{I}{\mathrm{homotopycolimit}}~\mathrm{Spec}^\mathrm{BK}\breve{\Phi}^I_{*,-}/\mathrm{Fro}^\mathbb{Z},	\\
}
\]
\[ 
\xymatrix@R+0pc@C+0pc{
\underset{r}{\mathrm{homotopylimit}}~\mathrm{Spec}^\mathrm{BK}{\Phi}^r_{*,-}/\mathrm{Fro}^\mathbb{Z},\underset{I}{\mathrm{homotopycolimit}}~\mathrm{Spec}^\mathrm{BK}{\Phi}^I_{*,-}/\mathrm{Fro}^\mathbb{Z}.	
}
\]

\end{definition}

\begin{definition}
We then consider the corresponding quasisheaves of the corresponding condensed solid topological modules from \cite{10CS2}:
\begin{align}
\mathrm{Quasicoherentsheaves, Condensed}_{*}	
\end{align}
where $*$ is one of the following spaces:
\begin{align}
&\mathrm{Spec}^\mathrm{CS}\widetilde{\Delta}_{*,-}/\mathrm{Fro}^\mathbb{Z},\mathrm{Spec}^\mathrm{CS}\widetilde{\nabla}_{*,-}/\mathrm{Fro}^\mathbb{Z},\mathrm{Spec}^\mathrm{CS}\widetilde{\Phi}_{*,-}/\mathrm{Fro}^\mathbb{Z},\mathrm{Spec}^\mathrm{CS}\widetilde{\Delta}^+_{*,-}/\mathrm{Fro}^\mathbb{Z},\\
&\mathrm{Spec}^\mathrm{CS}\widetilde{\nabla}^+_{*,-}/\mathrm{Fro}^\mathbb{Z},\mathrm{Spec}^\mathrm{CS}\widetilde{\Delta}^\dagger_{*,-}/\mathrm{Fro}^\mathbb{Z},\mathrm{Spec}^\mathrm{CS}\widetilde{\nabla}^\dagger_{*,-}/\mathrm{Fro}^\mathbb{Z},	\\
\end{align}
\begin{align}
&\mathrm{Spec}^\mathrm{CS}\breve{\Delta}_{*,-}/\mathrm{Fro}^\mathbb{Z},\breve{\nabla}_{*,-}/\mathrm{Fro}^\mathbb{Z},\mathrm{Spec}^\mathrm{CS}\breve{\Phi}_{*,-}/\mathrm{Fro}^\mathbb{Z},\mathrm{Spec}^\mathrm{CS}\breve{\Delta}^+_{*,-}/\mathrm{Fro}^\mathbb{Z},\\
&\mathrm{Spec}^\mathrm{CS}\breve{\nabla}^+_{*,-}/\mathrm{Fro}^\mathbb{Z},\mathrm{Spec}^\mathrm{CS}\breve{\Delta}^\dagger_{*,-}/\mathrm{Fro}^\mathbb{Z},\mathrm{Spec}^\mathrm{CS}\breve{\nabla}^\dagger_{*,-}/\mathrm{Fro}^\mathbb{Z},	\\
\end{align}
\begin{align}
&\mathrm{Spec}^\mathrm{CS}{\Delta}_{*,-}/\mathrm{Fro}^\mathbb{Z},\mathrm{Spec}^\mathrm{CS}{\nabla}_{*,-}/\mathrm{Fro}^\mathbb{Z},\mathrm{Spec}^\mathrm{CS}{\Phi}_{*,-}/\mathrm{Fro}^\mathbb{Z},\mathrm{Spec}^\mathrm{CS}{\Delta}^+_{*,-}/\mathrm{Fro}^\mathbb{Z},\\
&\mathrm{Spec}^\mathrm{CS}{\nabla}^+_{*,-}/\mathrm{Fro}^\mathbb{Z}, \mathrm{Spec}^\mathrm{CS}{\Delta}^\dagger_{*,-}/\mathrm{Fro}^\mathbb{Z},\mathrm{Spec}^\mathrm{CS}{\nabla}^\dagger_{*,-}/\mathrm{Fro}^\mathbb{Z}.	
\end{align}
Here for those space with notations related to the radius and the corresponding interval we consider the total unions $\bigcap_r,\bigcup_I$ in order to achieve the whole spaces to achieve the analogues of the corresponding FF curves from \cite{10KL1}, \cite{10KL2}, \cite{10FF} for
\[
\xymatrix@R+0pc@C+0pc{
\underset{r}{\mathrm{homotopylimit}}~\mathrm{Spec}^\mathrm{CS}\widetilde{\Phi}^r_{*,-},\underset{I}{\mathrm{homotopycolimit}}~\mathrm{Spec}^\mathrm{CS}\widetilde{\Phi}^I_{*,-},	\\
}
\]
\[
\xymatrix@R+0pc@C+0pc{
\underset{r}{\mathrm{homotopylimit}}~\mathrm{Spec}^\mathrm{CS}\breve{\Phi}^r_{*,-},\underset{I}{\mathrm{homotopycolimit}}~\mathrm{Spec}^\mathrm{CS}\breve{\Phi}^I_{*,-},	\\
}
\]
\[
\xymatrix@R+0pc@C+0pc{
\underset{r}{\mathrm{homotopylimit}}~\mathrm{Spec}^\mathrm{CS}{\Phi}^r_{*,-},\underset{I}{\mathrm{homotopycolimit}}~\mathrm{Spec}^\mathrm{CS}{\Phi}^I_{*,-}.	
}
\]
\[ 
\xymatrix@R+0pc@C+0pc{
\underset{r}{\mathrm{homotopylimit}}~\mathrm{Spec}^\mathrm{CS}\widetilde{\Phi}^r_{*,-}/\mathrm{Fro}^\mathbb{Z},\underset{I}{\mathrm{homotopycolimit}}~\mathrm{Spec}^\mathrm{CS}\widetilde{\Phi}^I_{*,-}/\mathrm{Fro}^\mathbb{Z},	\\
}
\]
\[ 
\xymatrix@R+0pc@C+0pc{
\underset{r}{\mathrm{homotopylimit}}~\mathrm{Spec}^\mathrm{CS}\breve{\Phi}^r_{*,-}/\mathrm{Fro}^\mathbb{Z},\underset{I}{\mathrm{homotopycolimit}}~\breve{\Phi}^I_{*,-}/\mathrm{Fro}^\mathbb{Z},	\\
}
\]
\[ 
\xymatrix@R+0pc@C+0pc{
\underset{r}{\mathrm{homotopylimit}}~\mathrm{Spec}^\mathrm{CS}{\Phi}^r_{*,-}/\mathrm{Fro}^\mathbb{Z},\underset{I}{\mathrm{homotopycolimit}}~\mathrm{Spec}^\mathrm{CS}{\Phi}^I_{*,-}/\mathrm{Fro}^\mathbb{Z}.	
}
\]

\end{definition}

\

\begin{proposition}
There is a well-defined functor from the $\infty$-category 
\begin{align}
\mathrm{Quasicoherentpresheaves,Condensed}_{*}	
\end{align}
where $*$ is one of the following spaces:
\begin{align}
&\mathrm{Spec}^\mathrm{CS}\widetilde{\Phi}_{*,-}/\mathrm{Fro}^\mathbb{Z},	\\
\end{align}
\begin{align}
&\mathrm{Spec}^\mathrm{CS}\breve{\Phi}_{*,-}/\mathrm{Fro}^\mathbb{Z},	\\
\end{align}
\begin{align}
&\mathrm{Spec}^\mathrm{CS}{\Phi}_{*,-}/\mathrm{Fro}^\mathbb{Z},	
\end{align}
to the $\infty$-category of $\mathrm{Fro}$-equivariant quasicoherent presheaves over similar spaces above correspondingly without the $\mathrm{Fro}$-quotients, and to the $\infty$-category of $\mathrm{Fro}$-equivariant quasicoherent modules over global sections of the structure $\infty$-sheaves of the similar spaces above correspondingly without the $\mathrm{Fro}$-quotients. Here for those space without notation related to the radius and the corresponding interval we consider the total unions $\bigcap_r,\bigcup_I$ in order to achieve the whole spaces to achieve the analogues of the corresponding FF curves from \cite{10KL1}, \cite{10KL2}, \cite{10FF} for
\[
\xymatrix@R+0pc@C+0pc{
\underset{r}{\mathrm{homotopylimit}}~\mathrm{Spec}^\mathrm{CS}\widetilde{\Phi}^r_{*,-},\underset{I}{\mathrm{homotopycolimit}}~\mathrm{Spec}^\mathrm{CS}\widetilde{\Phi}^I_{*,-},	\\
}
\]
\[
\xymatrix@R+0pc@C+0pc{
\underset{r}{\mathrm{homotopylimit}}~\mathrm{Spec}^\mathrm{CS}\breve{\Phi}^r_{*,-},\underset{I}{\mathrm{homotopycolimit}}~\mathrm{Spec}^\mathrm{CS}\breve{\Phi}^I_{*,-},	\\
}
\]
\[
\xymatrix@R+0pc@C+0pc{
\underset{r}{\mathrm{homotopylimit}}~\mathrm{Spec}^\mathrm{CS}{\Phi}^r_{*,-},\underset{I}{\mathrm{homotopycolimit}}~\mathrm{Spec}^\mathrm{CS}{\Phi}^I_{*,-}.	
}
\]
\[ 
\xymatrix@R+0pc@C+0pc{
\underset{r}{\mathrm{homotopylimit}}~\mathrm{Spec}^\mathrm{CS}\widetilde{\Phi}^r_{*,-}/\mathrm{Fro}^\mathbb{Z},\underset{I}{\mathrm{homotopycolimit}}~\mathrm{Spec}^\mathrm{CS}\widetilde{\Phi}^I_{*,-}/\mathrm{Fro}^\mathbb{Z},	\\
}
\]
\[ 
\xymatrix@R+0pc@C+0pc{
\underset{r}{\mathrm{homotopylimit}}~\mathrm{Spec}^\mathrm{CS}\breve{\Phi}^r_{*,-}/\mathrm{Fro}^\mathbb{Z},\underset{I}{\mathrm{homotopycolimit}}~\breve{\Phi}^I_{*,-}/\mathrm{Fro}^\mathbb{Z},	\\
}
\]
\[ 
\xymatrix@R+0pc@C+0pc{
\underset{r}{\mathrm{homotopylimit}}~\mathrm{Spec}^\mathrm{CS}{\Phi}^r_{*,-}/\mathrm{Fro}^\mathbb{Z},\underset{I}{\mathrm{homotopycolimit}}~\mathrm{Spec}^\mathrm{CS}{\Phi}^I_{*,-}/\mathrm{Fro}^\mathbb{Z}.	
}
\]	
In this situation we will have the target category being family parametrized by $r$ or $I$ in compatible glueing sense as in \cite[Definition 5.4.10]{10KL2}. In this situation for modules parametrized by the intervals we have the equivalence of $\infty$-categories by using \cite[Proposition 13.8]{10CS2}. Here the corresponding quasicoherent Frobenius modules are defined to be the corresponding homotopy colimits and limits of Frobenius modules:
\begin{align}
\underset{r}{\mathrm{homotopycolimit}}~M_r,\\
\underset{I}{\mathrm{homotopylimit}}~M_I,	
\end{align}
where each $M_r$ is a Frobenius-equivariant module over the period ring with respect to some radius $r$ while each $M_I$ is a Frobenius-equivariant module over the period ring with respect to some interval $I$.\\
\end{proposition}

\begin{proposition}
Similar proposition holds for 
\begin{align}
\mathrm{Quasicoherentsheaves,IndBanach}_{*}.	
\end{align}	
\end{proposition}

\

\begin{definition}
We then consider the corresponding quasipresheaves of perfect complexes the corresponding ind-Banach or monomorphic ind-Banach modules from \cite{10BBK}, \cite{10KKM}:
\begin{align}
\mathrm{Quasicoherentpresheaves,Perfectcomplex,IndBanach}_{*}	
\end{align}
where $*$ is one of the following spaces:
\begin{align}
&\mathrm{Spec}^\mathrm{BK}\widetilde{\Phi}_{*,-}/\mathrm{Fro}^\mathbb{Z},	\\
\end{align}
\begin{align}
&\mathrm{Spec}^\mathrm{BK}\breve{\Phi}_{*,-}/\mathrm{Fro}^\mathbb{Z},	\\
\end{align}
\begin{align}
&\mathrm{Spec}^\mathrm{BK}{\Phi}_{*,-}/\mathrm{Fro}^\mathbb{Z}.	
\end{align}
Here for those space without notation related to the radius and the corresponding interval we consider the total unions $\bigcap_r,\bigcup_I$ in order to achieve the whole spaces to achieve the analogues of the corresponding FF curves from \cite{10KL1}, \cite{10KL2}, \cite{10FF} for
\[
\xymatrix@R+0pc@C+0pc{
\underset{r}{\mathrm{homotopylimit}}~\mathrm{Spec}^\mathrm{BK}\widetilde{\Phi}^r_{*,-},\underset{I}{\mathrm{homotopycolimit}}~\mathrm{Spec}^\mathrm{BK}\widetilde{\Phi}^I_{*,-},	\\
}
\]
\[
\xymatrix@R+0pc@C+0pc{
\underset{r}{\mathrm{homotopylimit}}~\mathrm{Spec}^\mathrm{BK}\breve{\Phi}^r_{*,-},\underset{I}{\mathrm{homotopycolimit}}~\mathrm{Spec}^\mathrm{BK}\breve{\Phi}^I_{*,-},	\\
}
\]
\[
\xymatrix@R+0pc@C+0pc{
\underset{r}{\mathrm{homotopylimit}}~\mathrm{Spec}^\mathrm{BK}{\Phi}^r_{*,-},\underset{I}{\mathrm{homotopycolimit}}~\mathrm{Spec}^\mathrm{BK}{\Phi}^I_{*,-}.	
}
\]
\[  
\xymatrix@R+0pc@C+0pc{
\underset{r}{\mathrm{homotopylimit}}~\mathrm{Spec}^\mathrm{BK}\widetilde{\Phi}^r_{*,-}/\mathrm{Fro}^\mathbb{Z},\underset{I}{\mathrm{homotopycolimit}}~\mathrm{Spec}^\mathrm{BK}\widetilde{\Phi}^I_{*,-}/\mathrm{Fro}^\mathbb{Z},	\\
}
\]
\[ 
\xymatrix@R+0pc@C+0pc{
\underset{r}{\mathrm{homotopylimit}}~\mathrm{Spec}^\mathrm{BK}\breve{\Phi}^r_{*,-}/\mathrm{Fro}^\mathbb{Z},\underset{I}{\mathrm{homotopycolimit}}~\mathrm{Spec}^\mathrm{BK}\breve{\Phi}^I_{*,-}/\mathrm{Fro}^\mathbb{Z},	\\
}
\]
\[ 
\xymatrix@R+0pc@C+0pc{
\underset{r}{\mathrm{homotopylimit}}~\mathrm{Spec}^\mathrm{BK}{\Phi}^r_{*,-}/\mathrm{Fro}^\mathbb{Z},\underset{I}{\mathrm{homotopycolimit}}~\mathrm{Spec}^\mathrm{BK}{\Phi}^I_{*,-}/\mathrm{Fro}^\mathbb{Z}.	
}
\]

\end{definition}

\begin{definition}
We then consider the corresponding quasisheaves of perfect complexes of the corresponding condensed solid topological modules from \cite{10CS2}:
\begin{align}
\mathrm{Quasicoherentsheaves, Perfectcomplex, Condensed}_{*}	
\end{align}
where $*$ is one of the following spaces:
\begin{align}
&\mathrm{Spec}^\mathrm{CS}\widetilde{\Delta}_{*,-}/\mathrm{Fro}^\mathbb{Z},\mathrm{Spec}^\mathrm{CS}\widetilde{\nabla}_{*,-}/\mathrm{Fro}^\mathbb{Z},\mathrm{Spec}^\mathrm{CS}\widetilde{\Phi}_{*,-}/\mathrm{Fro}^\mathbb{Z},\mathrm{Spec}^\mathrm{CS}\widetilde{\Delta}^+_{*,-}/\mathrm{Fro}^\mathbb{Z},\\
&\mathrm{Spec}^\mathrm{CS}\widetilde{\nabla}^+_{*,-}/\mathrm{Fro}^\mathbb{Z},\mathrm{Spec}^\mathrm{CS}\widetilde{\Delta}^\dagger_{*,-}/\mathrm{Fro}^\mathbb{Z},\mathrm{Spec}^\mathrm{CS}\widetilde{\nabla}^\dagger_{*,-}/\mathrm{Fro}^\mathbb{Z},	\\
\end{align}
\begin{align}
&\mathrm{Spec}^\mathrm{CS}\breve{\Delta}_{*,-}/\mathrm{Fro}^\mathbb{Z},\breve{\nabla}_{*,-}/\mathrm{Fro}^\mathbb{Z},\mathrm{Spec}^\mathrm{CS}\breve{\Phi}_{*,-}/\mathrm{Fro}^\mathbb{Z},\mathrm{Spec}^\mathrm{CS}\breve{\Delta}^+_{*,-}/\mathrm{Fro}^\mathbb{Z},\\
&\mathrm{Spec}^\mathrm{CS}\breve{\nabla}^+_{*,-}/\mathrm{Fro}^\mathbb{Z},\mathrm{Spec}^\mathrm{CS}\breve{\Delta}^\dagger_{*,-}/\mathrm{Fro}^\mathbb{Z},\mathrm{Spec}^\mathrm{CS}\breve{\nabla}^\dagger_{*,-}/\mathrm{Fro}^\mathbb{Z},	\\
\end{align}
\begin{align}
&\mathrm{Spec}^\mathrm{CS}{\Delta}_{*,-}/\mathrm{Fro}^\mathbb{Z},\mathrm{Spec}^\mathrm{CS}{\nabla}_{*,-}/\mathrm{Fro}^\mathbb{Z},\mathrm{Spec}^\mathrm{CS}{\Phi}_{*,-}/\mathrm{Fro}^\mathbb{Z},\mathrm{Spec}^\mathrm{CS}{\Delta}^+_{*,-}/\mathrm{Fro}^\mathbb{Z},\\
&\mathrm{Spec}^\mathrm{CS}{\nabla}^+_{*,-}/\mathrm{Fro}^\mathbb{Z}, \mathrm{Spec}^\mathrm{CS}{\Delta}^\dagger_{*,-}/\mathrm{Fro}^\mathbb{Z},\mathrm{Spec}^\mathrm{CS}{\nabla}^\dagger_{*,-}/\mathrm{Fro}^\mathbb{Z}.	
\end{align}
Here for those space with notations related to the radius and the corresponding interval we consider the total unions $\bigcap_r,\bigcup_I$ in order to achieve the whole spaces to achieve the analogues of the corresponding FF curves from \cite{10KL1}, \cite{10KL2}, \cite{10FF} for
\[
\xymatrix@R+0pc@C+0pc{
\underset{r}{\mathrm{homotopylimit}}~\mathrm{Spec}^\mathrm{CS}\widetilde{\Phi}^r_{*,-},\underset{I}{\mathrm{homotopycolimit}}~\mathrm{Spec}^\mathrm{CS}\widetilde{\Phi}^I_{*,-},	\\
}
\]
\[
\xymatrix@R+0pc@C+0pc{
\underset{r}{\mathrm{homotopylimit}}~\mathrm{Spec}^\mathrm{CS}\breve{\Phi}^r_{*,-},\underset{I}{\mathrm{homotopycolimit}}~\mathrm{Spec}^\mathrm{CS}\breve{\Phi}^I_{*,-},	\\
}
\]
\[
\xymatrix@R+0pc@C+0pc{
\underset{r}{\mathrm{homotopylimit}}~\mathrm{Spec}^\mathrm{CS}{\Phi}^r_{*,-},\underset{I}{\mathrm{homotopycolimit}}~\mathrm{Spec}^\mathrm{CS}{\Phi}^I_{*,-}.	
}
\]
\[ 
\xymatrix@R+0pc@C+0pc{
\underset{r}{\mathrm{homotopylimit}}~\mathrm{Spec}^\mathrm{CS}\widetilde{\Phi}^r_{*,-}/\mathrm{Fro}^\mathbb{Z},\underset{I}{\mathrm{homotopycolimit}}~\mathrm{Spec}^\mathrm{CS}\widetilde{\Phi}^I_{*,-}/\mathrm{Fro}^\mathbb{Z},	\\
}
\]
\[ 
\xymatrix@R+0pc@C+0pc{
\underset{r}{\mathrm{homotopylimit}}~\mathrm{Spec}^\mathrm{CS}\breve{\Phi}^r_{*,-}/\mathrm{Fro}^\mathbb{Z},\underset{I}{\mathrm{homotopycolimit}}~\breve{\Phi}^I_{*,-}/\mathrm{Fro}^\mathbb{Z},	\\
}
\]
\[ 
\xymatrix@R+0pc@C+0pc{
\underset{r}{\mathrm{homotopylimit}}~\mathrm{Spec}^\mathrm{CS}{\Phi}^r_{*,-}/\mathrm{Fro}^\mathbb{Z},\underset{I}{\mathrm{homotopycolimit}}~\mathrm{Spec}^\mathrm{CS}{\Phi}^I_{*,-}/\mathrm{Fro}^\mathbb{Z}.	
}
\]

\end{definition}

\begin{proposition}
There is a well-defined functor from the $\infty$-category 
\begin{align}
\mathrm{Quasicoherentpresheaves,Perfectcomplex,Condensed}_{*}	
\end{align}
where $*$ is one of the following spaces:
\begin{align}
&\mathrm{Spec}^\mathrm{CS}\widetilde{\Phi}_{*,-}/\mathrm{Fro}^\mathbb{Z},	\\
\end{align}
\begin{align}
&\mathrm{Spec}^\mathrm{CS}\breve{\Phi}_{*,-}/\mathrm{Fro}^\mathbb{Z},	\\
\end{align}
\begin{align}
&\mathrm{Spec}^\mathrm{CS}{\Phi}_{*,-}/\mathrm{Fro}^\mathbb{Z},	
\end{align}
to the $\infty$-category of $\mathrm{Fro}$-equivariant quasicoherent presheaves over similar spaces above correspondingly without the $\mathrm{Fro}$-quotients, and to the $\infty$-category of $\mathrm{Fro}$-equivariant quasicoherent modules over global sections of the structure $\infty$-sheaves of the similar spaces above correspondingly without the $\mathrm{Fro}$-quotients. Here for those space without notation related to the radius and the corresponding interval we consider the total unions $\bigcap_r,\bigcup_I$ in order to achieve the whole spaces to achieve the analogues of the corresponding FF curves from \cite{10KL1}, \cite{10KL2}, \cite{10FF} for
\[
\xymatrix@R+0pc@C+0pc{
\underset{r}{\mathrm{homotopylimit}}~\mathrm{Spec}^\mathrm{CS}\widetilde{\Phi}^r_{*,-},\underset{I}{\mathrm{homotopycolimit}}~\mathrm{Spec}^\mathrm{CS}\widetilde{\Phi}^I_{*,-},	\\
}
\]
\[
\xymatrix@R+0pc@C+0pc{
\underset{r}{\mathrm{homotopylimit}}~\mathrm{Spec}^\mathrm{CS}\breve{\Phi}^r_{*,-},\underset{I}{\mathrm{homotopycolimit}}~\mathrm{Spec}^\mathrm{CS}\breve{\Phi}^I_{*,-},	\\
}
\]
\[
\xymatrix@R+0pc@C+0pc{
\underset{r}{\mathrm{homotopylimit}}~\mathrm{Spec}^\mathrm{CS}{\Phi}^r_{*,-},\underset{I}{\mathrm{homotopycolimit}}~\mathrm{Spec}^\mathrm{CS}{\Phi}^I_{*,-}.	
}
\]
\[ 
\xymatrix@R+0pc@C+0pc{
\underset{r}{\mathrm{homotopylimit}}~\mathrm{Spec}^\mathrm{CS}\widetilde{\Phi}^r_{*,-}/\mathrm{Fro}^\mathbb{Z},\underset{I}{\mathrm{homotopycolimit}}~\mathrm{Spec}^\mathrm{CS}\widetilde{\Phi}^I_{*,-}/\mathrm{Fro}^\mathbb{Z},	\\
}
\]
\[ 
\xymatrix@R+0pc@C+0pc{
\underset{r}{\mathrm{homotopylimit}}~\mathrm{Spec}^\mathrm{CS}\breve{\Phi}^r_{*,-}/\mathrm{Fro}^\mathbb{Z},\underset{I}{\mathrm{homotopycolimit}}~\breve{\Phi}^I_{*,-}/\mathrm{Fro}^\mathbb{Z},	\\
}
\]
\[ 
\xymatrix@R+0pc@C+0pc{
\underset{r}{\mathrm{homotopylimit}}~\mathrm{Spec}^\mathrm{CS}{\Phi}^r_{*,-}/\mathrm{Fro}^\mathbb{Z},\underset{I}{\mathrm{homotopycolimit}}~\mathrm{Spec}^\mathrm{CS}{\Phi}^I_{*,-}/\mathrm{Fro}^\mathbb{Z}.	
}
\]	
In this situation we will have the target category being family parametrized by $r$ or $I$ in compatible glueing sense as in \cite[Definition 5.4.10]{10KL2}. In this situation for modules parametrized by the intervals we have the equivalence of $\infty$-categories by using \cite[Proposition 12.18]{10CS2}. Here the corresponding quasicoherent Frobenius modules are defined to be the corresponding homotopy colimits and limits of Frobenius modules:
\begin{align}
\underset{r}{\mathrm{homotopycolimit}}~M_r,\\
\underset{I}{\mathrm{homotopylimit}}~M_I,	
\end{align}
where each $M_r$ is a Frobenius-equivariant module over the period ring with respect to some radius $r$ while each $M_I$ is a Frobenius-equivariant module over the period ring with respect to some interval $I$.\\
\end{proposition}

\begin{proposition}
Similar proposition holds for 
\begin{align}
\mathrm{Quasicoherentsheaves,Perfectcomplex,IndBanach}_{*}.	
\end{align}	
\end{proposition}

\chapter{Deformation}

\section{Multivariate Hodge Iwasawa Modules by Deformation}

This chapter follows closely \cite{10T1}, \cite{10T2}, \cite{10T3}, \cite{10T4}, \cite{10T5}, \cite{10T6}, \cite{10KPX}, \cite{10KP}, \cite{10KL1}, \cite{10KL2}, \cite{10BK}, \cite{10BBBK}, \cite{10BBM}, \cite{10KKM}, \cite{10CS1}, \cite{10CS2}, \cite{10CKZ}, \cite{10PZ}, \cite{10BCM}, \cite{10LBV}.

\begin{remark}
In the following chapters, we remind the readers of the fact that the notations for the deformation $A,-,\square$ in our following discussion will mean different thing, the deformation with respect to $A,-,\square$ will happen along the structure sheaves $\mathcal{O}$ directly. The $\infty$-descent results in \cite{10BK}, \cite{10BBBK}, \cite{10BBM}, \cite{10KKM}, \cite{10CS1}, \cite{10CS2} will guarantee that the deformed sheaves are still quasicoherent sheaves over $\mathcal{O}$.
\end{remark}

\subsection{Frobenius Quasicoherent Modules I}

\begin{definition}
First we consider the Bambozzi-Kremnizer spectrum $\underset{\mathrm{Spec}}{\mathcal{O}}^\mathrm{BK}(*)$ attached to any of those in the above from \cite{10BK} by taking derived rational localization:
\begin{align}
&\underset{\mathrm{Spec}}{\mathcal{O}}^\mathrm{BK}\widetilde{\Phi}_{\psi,\Gamma,A},\underset{\mathrm{Spec}}{\mathcal{O}}^\mathrm{BK}\widetilde{\Phi}^r_{\psi,\Gamma,A},\underset{\mathrm{Spec}}{\mathcal{O}}^\mathrm{BK}\widetilde{\Phi}^I_{\psi,\Gamma,A},	
\end{align}
\begin{align}
&\underset{\mathrm{Spec}}{\mathcal{O}}^\mathrm{BK}\breve{\Phi}_{\psi,\Gamma,A},\underset{\mathrm{Spec}}{\mathcal{O}}^\mathrm{BK}\breve{\Phi}^r_{\psi,\Gamma,A},\underset{\mathrm{Spec}}{\mathcal{O}}^\mathrm{BK}\breve{\Phi}^I_{\psi,\Gamma,A},	
\end{align}
\begin{align}
&\underset{\mathrm{Spec}}{\mathcal{O}}^\mathrm{BK}{\Phi}_{\psi,\Gamma,A},
\underset{\mathrm{Spec}}{\mathcal{O}}^\mathrm{BK}{\Phi}^r_{\psi,\Gamma,A},\underset{\mathrm{Spec}}{\mathcal{O}}^\mathrm{BK}{\Phi}^I_{\psi,\Gamma,A}.	
\end{align}

Then we take the corresponding quotients by using the corresponding Frobenius operators:
\begin{align}
&\underset{\mathrm{Spec}}{\mathcal{O}}^\mathrm{BK}\widetilde{\Phi}_{\psi,\Gamma,A}/\mathrm{Fro}^\mathbb{Z},	\\
\end{align}
\begin{align}
&\underset{\mathrm{Spec}}{\mathcal{O}}^\mathrm{BK}\breve{\Phi}_{\psi,\Gamma,A}/\mathrm{Fro}^\mathbb{Z},	\\
\end{align}
\begin{align}
&\underset{\mathrm{Spec}}{\mathcal{O}}^\mathrm{BK}{\Phi}_{\psi,\Gamma,A}/\mathrm{Fro}^\mathbb{Z}.	
\end{align}
Here for those space without notation related to the radius and the corresponding interval we consider the total unions $\bigcap_r,\bigcup_I$ in order to achieve the whole spaces to achieve the analogues of the corresponding FF curves from \cite{10KL1}, \cite{10KL2}, \cite{10FF} for
\[
\xymatrix@R+0pc@C+0pc{
\underset{r}{\mathrm{homotopycolimit}}~\underset{\mathrm{Spec}}{\mathcal{O}}^\mathrm{BK}\widetilde{\Phi}^r_{\psi,\Gamma,A},\underset{I}{\mathrm{homotopylimit}}~\underset{\mathrm{Spec}}{\mathcal{O}}^\mathrm{BK}\widetilde{\Phi}^I_{\psi,\Gamma,A},	\\
}
\]
\[
\xymatrix@R+0pc@C+0pc{
\underset{r}{\mathrm{homotopycolimit}}~\underset{\mathrm{Spec}}{\mathcal{O}}^\mathrm{BK}\breve{\Phi}^r_{\psi,\Gamma,A},\underset{I}{\mathrm{homotopylimit}}~\underset{\mathrm{Spec}}{\mathcal{O}}^\mathrm{BK}\breve{\Phi}^I_{\psi,\Gamma,A},	\\
}
\]
\[
\xymatrix@R+0pc@C+0pc{
\underset{r}{\mathrm{homotopycolimit}}~\underset{\mathrm{Spec}}{\mathcal{O}}^\mathrm{BK}{\Phi}^r_{\psi,\Gamma,A},\underset{I}{\mathrm{homotopylimit}}~\underset{\mathrm{Spec}}{\mathcal{O}}^\mathrm{BK}{\Phi}^I_{\psi,\Gamma,A}.	
}
\]
\[  
\xymatrix@R+0pc@C+0pc{
\underset{r}{\mathrm{homotopycolimit}}~\underset{\mathrm{Spec}}{\mathcal{O}}^\mathrm{BK}\widetilde{\Phi}^r_{\psi,\Gamma,A}/\mathrm{Fro}^\mathbb{Z},\underset{I}{\mathrm{homotopylimit}}~\underset{\mathrm{Spec}}{\mathcal{O}}^\mathrm{BK}\widetilde{\Phi}^I_{\psi,\Gamma,A}/\mathrm{Fro}^\mathbb{Z},	\\
}
\]
\[ 
\xymatrix@R+0pc@C+0pc{
\underset{r}{\mathrm{homotopycolimit}}~\underset{\mathrm{Spec}}{\mathcal{O}}^\mathrm{BK}\breve{\Phi}^r_{\psi,\Gamma,A}/\mathrm{Fro}^\mathbb{Z},\underset{I}{\mathrm{homotopylimit}}~\underset{\mathrm{Spec}}{\mathcal{O}}^\mathrm{BK}\breve{\Phi}^I_{\psi,\Gamma,A}/\mathrm{Fro}^\mathbb{Z},	\\
}
\]
\[ 
\xymatrix@R+0pc@C+0pc{
\underset{r}{\mathrm{homotopycolimit}}~\underset{\mathrm{Spec}}{\mathcal{O}}^\mathrm{BK}{\Phi}^r_{\psi,\Gamma,A}/\mathrm{Fro}^\mathbb{Z},\underset{I}{\mathrm{homotopylimit}}~\underset{\mathrm{Spec}}{\mathcal{O}}^\mathrm{BK}{\Phi}^I_{\psi,\Gamma,A}/\mathrm{Fro}^\mathbb{Z}.	
}
\]

\end{definition}

\indent Meanwhile we have the corresponding Clausen-Scholze analytic stacks from \cite{10CS2}, therefore applying their construction we have:

\begin{definition}
Here we define the following products by using the solidified tensor product from \cite{10CS1} and \cite{10CS2}. Namely $A$ will still as above as a Banach ring over $\mathbb{Q}_p$. Then we take solidified tensor product $\overset{\blacksquare}{\otimes}$ of any of the following
\[
\xymatrix@R+0pc@C+0pc{
\widetilde{\Delta}_{\psi,\Gamma},\widetilde{\nabla}_{\psi,\Gamma},\widetilde{\Phi}_{\psi,\Gamma},\widetilde{\Delta}^+_{\psi,\Gamma},\widetilde{\nabla}^+_{\psi,\Gamma},\widetilde{\Delta}^\dagger_{\psi,\Gamma},\widetilde{\nabla}^\dagger_{\psi,\Gamma},\widetilde{\Phi}^r_{\psi,\Gamma},\widetilde{\Phi}^I_{\psi,\Gamma}, 
}
\]

\[
\xymatrix@R+0pc@C+0pc{
\breve{\Delta}_{\psi,\Gamma},\breve{\nabla}_{\psi,\Gamma},\breve{\Phi}_{\psi,\Gamma},\breve{\Delta}^+_{\psi,\Gamma},\breve{\nabla}^+_{\psi,\Gamma},\breve{\Delta}^\dagger_{\psi,\Gamma},\breve{\nabla}^\dagger_{\psi,\Gamma},\breve{\Phi}^r_{\psi,\Gamma},\breve{\Phi}^I_{\psi,\Gamma},	
}
\]

\[
\xymatrix@R+0pc@C+0pc{
{\Delta}_{\psi,\Gamma},{\nabla}_{\psi,\Gamma},{\Phi}_{\psi,\Gamma},{\Delta}^+_{\psi,\Gamma},{\nabla}^+_{\psi,\Gamma},{\Delta}^\dagger_{\psi,\Gamma},{\nabla}^\dagger_{\psi,\Gamma},{\Phi}^r_{\psi,\Gamma},{\Phi}^I_{\psi,\Gamma},	
}
\]  	
with $A$. Then we have the notations:
\[
\xymatrix@R+0pc@C+0pc{
\widetilde{\Delta}_{\psi,\Gamma,A},\widetilde{\nabla}_{\psi,\Gamma,A},\widetilde{\Phi}_{\psi,\Gamma,A},\widetilde{\Delta}^+_{\psi,\Gamma,A},\widetilde{\nabla}^+_{\psi,\Gamma,A},\widetilde{\Delta}^\dagger_{\psi,\Gamma,A},\widetilde{\nabla}^\dagger_{\psi,\Gamma,A},\widetilde{\Phi}^r_{\psi,\Gamma,A},\widetilde{\Phi}^I_{\psi,\Gamma,A}, 
}
\]

\[
\xymatrix@R+0pc@C+0pc{
\breve{\Delta}_{\psi,\Gamma,A},\breve{\nabla}_{\psi,\Gamma,A},\breve{\Phi}_{\psi,\Gamma,A},\breve{\Delta}^+_{\psi,\Gamma,A},\breve{\nabla}^+_{\psi,\Gamma,A},\breve{\Delta}^\dagger_{\psi,\Gamma,A},\breve{\nabla}^\dagger_{\psi,\Gamma,A},\breve{\Phi}^r_{\psi,\Gamma,A},\breve{\Phi}^I_{\psi,\Gamma,A},	
}
\]

\[
\xymatrix@R+0pc@C+0pc{
{\Delta}_{\psi,\Gamma,A},{\nabla}_{\psi,\Gamma,A},{\Phi}_{\psi,\Gamma,A},{\Delta}^+_{\psi,\Gamma,A},{\nabla}^+_{\psi,\Gamma,A},{\Delta}^\dagger_{\psi,\Gamma,A},{\nabla}^\dagger_{\psi,\Gamma,A},{\Phi}^r_{\psi,\Gamma,A},{\Phi}^I_{\psi,\Gamma,A}.	
}
\]
\end{definition}

\begin{definition}
First we consider the Clausen-Scholze spectrum $\underset{\mathrm{Spec}}{\mathcal{O}}^\mathrm{CS}(*)$ attached to any of those in the above from \cite{10CS2} by taking derived rational localization:
\begin{align}
\underset{\mathrm{Spec}}{\mathcal{O}}^\mathrm{CS}\widetilde{\Delta}_{\psi,\Gamma,A},\underset{\mathrm{Spec}}{\mathcal{O}}^\mathrm{CS}\widetilde{\nabla}_{\psi,\Gamma,A},\underset{\mathrm{Spec}}{\mathcal{O}}^\mathrm{CS}\widetilde{\Phi}_{\psi,\Gamma,A},\underset{\mathrm{Spec}}{\mathcal{O}}^\mathrm{CS}\widetilde{\Delta}^+_{\psi,\Gamma,A},\underset{\mathrm{Spec}}{\mathcal{O}}^\mathrm{CS}\widetilde{\nabla}^+_{\psi,\Gamma,A},\\
\underset{\mathrm{Spec}}{\mathcal{O}}^\mathrm{CS}\widetilde{\Delta}^\dagger_{\psi,\Gamma,A},\underset{\mathrm{Spec}}{\mathcal{O}}^\mathrm{CS}\widetilde{\nabla}^\dagger_{\psi,\Gamma,A},\underset{\mathrm{Spec}}{\mathcal{O}}^\mathrm{CS}\widetilde{\Phi}^r_{\psi,\Gamma,A},\underset{\mathrm{Spec}}{\mathcal{O}}^\mathrm{CS}\widetilde{\Phi}^I_{\psi,\Gamma,A},	\\
\end{align}
\begin{align}
\underset{\mathrm{Spec}}{\mathcal{O}}^\mathrm{CS}\breve{\Delta}_{\psi,\Gamma,A},\breve{\nabla}_{\psi,\Gamma,A},\underset{\mathrm{Spec}}{\mathcal{O}}^\mathrm{CS}\breve{\Phi}_{\psi,\Gamma,A},\underset{\mathrm{Spec}}{\mathcal{O}}^\mathrm{CS}\breve{\Delta}^+_{\psi,\Gamma,A},\underset{\mathrm{Spec}}{\mathcal{O}}^\mathrm{CS}\breve{\nabla}^+_{\psi,\Gamma,A},\\
\underset{\mathrm{Spec}}{\mathcal{O}}^\mathrm{CS}\breve{\Delta}^\dagger_{\psi,\Gamma,A},\underset{\mathrm{Spec}}{\mathcal{O}}^\mathrm{CS}\breve{\nabla}^\dagger_{\psi,\Gamma,A},\underset{\mathrm{Spec}}{\mathcal{O}}^\mathrm{CS}\breve{\Phi}^r_{\psi,\Gamma,A},\breve{\Phi}^I_{\psi,\Gamma,A},	\\
\end{align}
\begin{align}
\underset{\mathrm{Spec}}{\mathcal{O}}^\mathrm{CS}{\Delta}_{\psi,\Gamma,A},\underset{\mathrm{Spec}}{\mathcal{O}}^\mathrm{CS}{\nabla}_{\psi,\Gamma,A},\underset{\mathrm{Spec}}{\mathcal{O}}^\mathrm{CS}{\Phi}_{\psi,\Gamma,A},\underset{\mathrm{Spec}}{\mathcal{O}}^\mathrm{CS}{\Delta}^+_{\psi,\Gamma,A},\underset{\mathrm{Spec}}{\mathcal{O}}^\mathrm{CS}{\nabla}^+_{\psi,\Gamma,A},\\
\underset{\mathrm{Spec}}{\mathcal{O}}^\mathrm{CS}{\Delta}^\dagger_{\psi,\Gamma,A},\underset{\mathrm{Spec}}{\mathcal{O}}^\mathrm{CS}{\nabla}^\dagger_{\psi,\Gamma,A},\underset{\mathrm{Spec}}{\mathcal{O}}^\mathrm{CS}{\Phi}^r_{\psi,\Gamma,A},\underset{\mathrm{Spec}}{\mathcal{O}}^\mathrm{CS}{\Phi}^I_{\psi,\Gamma,A}.	
\end{align}

Then we take the corresponding quotients by using the corresponding Frobenius operators:
\begin{align}
&\underset{\mathrm{Spec}}{\mathcal{O}}^\mathrm{CS}\widetilde{\Delta}_{\psi,\Gamma,A}/\mathrm{Fro}^\mathbb{Z},\underset{\mathrm{Spec}}{\mathcal{O}}^\mathrm{CS}\widetilde{\nabla}_{\psi,\Gamma,A}/\mathrm{Fro}^\mathbb{Z},\underset{\mathrm{Spec}}{\mathcal{O}}^\mathrm{CS}\widetilde{\Phi}_{\psi,\Gamma,A}/\mathrm{Fro}^\mathbb{Z},\underset{\mathrm{Spec}}{\mathcal{O}}^\mathrm{CS}\widetilde{\Delta}^+_{\psi,\Gamma,A}/\mathrm{Fro}^\mathbb{Z},\\
&\underset{\mathrm{Spec}}{\mathcal{O}}^\mathrm{CS}\widetilde{\nabla}^+_{\psi,\Gamma,A}/\mathrm{Fro}^\mathbb{Z}, \underset{\mathrm{Spec}}{\mathcal{O}}^\mathrm{CS}\widetilde{\Delta}^\dagger_{\psi,\Gamma,A}/\mathrm{Fro}^\mathbb{Z},\underset{\mathrm{Spec}}{\mathcal{O}}^\mathrm{CS}\widetilde{\nabla}^\dagger_{\psi,\Gamma,A}/\mathrm{Fro}^\mathbb{Z},	\\
\end{align}
\begin{align}
&\underset{\mathrm{Spec}}{\mathcal{O}}^\mathrm{CS}\breve{\Delta}_{\psi,\Gamma,A}/\mathrm{Fro}^\mathbb{Z},\breve{\nabla}_{\psi,\Gamma,A}/\mathrm{Fro}^\mathbb{Z},\underset{\mathrm{Spec}}{\mathcal{O}}^\mathrm{CS}\breve{\Phi}_{\psi,\Gamma,A}/\mathrm{Fro}^\mathbb{Z},\underset{\mathrm{Spec}}{\mathcal{O}}^\mathrm{CS}\breve{\Delta}^+_{\psi,\Gamma,A}/\mathrm{Fro}^\mathbb{Z},\\
&\underset{\mathrm{Spec}}{\mathcal{O}}^\mathrm{CS}\breve{\nabla}^+_{\psi,\Gamma,A}/\mathrm{Fro}^\mathbb{Z}, \underset{\mathrm{Spec}}{\mathcal{O}}^\mathrm{CS}\breve{\Delta}^\dagger_{\psi,\Gamma,A}/\mathrm{Fro}^\mathbb{Z},\underset{\mathrm{Spec}}{\mathcal{O}}^\mathrm{CS}\breve{\nabla}^\dagger_{\psi,\Gamma,A}/\mathrm{Fro}^\mathbb{Z},	\\
\end{align}
\begin{align}
&\underset{\mathrm{Spec}}{\mathcal{O}}^\mathrm{CS}{\Delta}_{\psi,\Gamma,A}/\mathrm{Fro}^\mathbb{Z},\underset{\mathrm{Spec}}{\mathcal{O}}^\mathrm{CS}{\nabla}_{\psi,\Gamma,A}/\mathrm{Fro}^\mathbb{Z},\underset{\mathrm{Spec}}{\mathcal{O}}^\mathrm{CS}{\Phi}_{\psi,\Gamma,A}/\mathrm{Fro}^\mathbb{Z},\underset{\mathrm{Spec}}{\mathcal{O}}^\mathrm{CS}{\Delta}^+_{\psi,\Gamma,A}/\mathrm{Fro}^\mathbb{Z},\\
&\underset{\mathrm{Spec}}{\mathcal{O}}^\mathrm{CS}{\nabla}^+_{\psi,\Gamma,A}/\mathrm{Fro}^\mathbb{Z}, \underset{\mathrm{Spec}}{\mathcal{O}}^\mathrm{CS}{\Delta}^\dagger_{\psi,\Gamma,A}/\mathrm{Fro}^\mathbb{Z},\underset{\mathrm{Spec}}{\mathcal{O}}^\mathrm{CS}{\nabla}^\dagger_{\psi,\Gamma,A}/\mathrm{Fro}^\mathbb{Z}.	
\end{align}
Here for those space with notations related to the radius and the corresponding interval we consider the total unions $\bigcap_r,\bigcup_I$ in order to achieve the whole spaces to achieve the analogues of the corresponding FF curves from \cite{10KL1}, \cite{10KL2}, \cite{10FF} for
\[
\xymatrix@R+0pc@C+0pc{
\underset{r}{\mathrm{homotopycolimit}}~\underset{\mathrm{Spec}}{\mathcal{O}}^\mathrm{CS}\widetilde{\Phi}^r_{\psi,\Gamma,A},\underset{I}{\mathrm{homotopylimit}}~\underset{\mathrm{Spec}}{\mathcal{O}}^\mathrm{CS}\widetilde{\Phi}^I_{\psi,\Gamma,A},	\\
}
\]
\[
\xymatrix@R+0pc@C+0pc{
\underset{r}{\mathrm{homotopycolimit}}~\underset{\mathrm{Spec}}{\mathcal{O}}^\mathrm{CS}\breve{\Phi}^r_{\psi,\Gamma,A},\underset{I}{\mathrm{homotopylimit}}~\underset{\mathrm{Spec}}{\mathcal{O}}^\mathrm{CS}\breve{\Phi}^I_{\psi,\Gamma,A},	\\
}
\]
\[
\xymatrix@R+0pc@C+0pc{
\underset{r}{\mathrm{homotopycolimit}}~\underset{\mathrm{Spec}}{\mathcal{O}}^\mathrm{CS}{\Phi}^r_{\psi,\Gamma,A},\underset{I}{\mathrm{homotopylimit}}~\underset{\mathrm{Spec}}{\mathcal{O}}^\mathrm{CS}{\Phi}^I_{\psi,\Gamma,A}.	
}
\]
\[ 
\xymatrix@R+0pc@C+0pc{
\underset{r}{\mathrm{homotopycolimit}}~\underset{\mathrm{Spec}}{\mathcal{O}}^\mathrm{CS}\widetilde{\Phi}^r_{\psi,\Gamma,A}/\mathrm{Fro}^\mathbb{Z},\underset{I}{\mathrm{homotopylimit}}~\underset{\mathrm{Spec}}{\mathcal{O}}^\mathrm{CS}\widetilde{\Phi}^I_{\psi,\Gamma,A}/\mathrm{Fro}^\mathbb{Z},	\\
}
\]
\[ 
\xymatrix@R+0pc@C+0pc{
\underset{r}{\mathrm{homotopycolimit}}~\underset{\mathrm{Spec}}{\mathcal{O}}^\mathrm{CS}\breve{\Phi}^r_{\psi,\Gamma,A}/\mathrm{Fro}^\mathbb{Z},\underset{I}{\mathrm{homotopylimit}}~\breve{\Phi}^I_{\psi,\Gamma,A}/\mathrm{Fro}^\mathbb{Z},	\\
}
\]
\[ 
\xymatrix@R+0pc@C+0pc{
\underset{r}{\mathrm{homotopycolimit}}~\underset{\mathrm{Spec}}{\mathcal{O}}^\mathrm{CS}{\Phi}^r_{\psi,\Gamma,A}/\mathrm{Fro}^\mathbb{Z},\underset{I}{\mathrm{homotopylimit}}~\underset{\mathrm{Spec}}{\mathcal{O}}^\mathrm{CS}{\Phi}^I_{\psi,\Gamma,A}/\mathrm{Fro}^\mathbb{Z}.	
}
\]

\end{definition}

\

\begin{definition}
We then consider the corresponding quasipresheaves of the corresponding ind-Banach or monomorphic ind-Banach modules from \cite{10BBK}, \cite{10KKM}:
\begin{align}
\mathrm{Quasicoherentpresheaves,IndBanach}_{*}	
\end{align}
where $*$ is one of the following spaces:
\begin{align}
&\underset{\mathrm{Spec}}{\mathcal{O}}^\mathrm{BK}\widetilde{\Phi}_{\psi,\Gamma,A}/\mathrm{Fro}^\mathbb{Z},	\\
\end{align}
\begin{align}
&\underset{\mathrm{Spec}}{\mathcal{O}}^\mathrm{BK}\breve{\Phi}_{\psi,\Gamma,A}/\mathrm{Fro}^\mathbb{Z},	\\
\end{align}
\begin{align}
&\underset{\mathrm{Spec}}{\mathcal{O}}^\mathrm{BK}{\Phi}_{\psi,\Gamma,A}/\mathrm{Fro}^\mathbb{Z}.	
\end{align}
Here for those space without notation related to the radius and the corresponding interval we consider the total unions $\bigcap_r,\bigcup_I$ in order to achieve the whole spaces to achieve the analogues of the corresponding FF curves from \cite{10KL1}, \cite{10KL2}, \cite{10FF} for
\[
\xymatrix@R+0pc@C+0pc{
\underset{r}{\mathrm{homotopycolimit}}~\underset{\mathrm{Spec}}{\mathcal{O}}^\mathrm{BK}\widetilde{\Phi}^r_{\psi,\Gamma,A},\underset{I}{\mathrm{homotopylimit}}~\underset{\mathrm{Spec}}{\mathcal{O}}^\mathrm{BK}\widetilde{\Phi}^I_{\psi,\Gamma,A},	\\
}
\]
\[
\xymatrix@R+0pc@C+0pc{
\underset{r}{\mathrm{homotopycolimit}}~\underset{\mathrm{Spec}}{\mathcal{O}}^\mathrm{BK}\breve{\Phi}^r_{\psi,\Gamma,A},\underset{I}{\mathrm{homotopylimit}}~\underset{\mathrm{Spec}}{\mathcal{O}}^\mathrm{BK}\breve{\Phi}^I_{\psi,\Gamma,A},	\\
}
\]
\[
\xymatrix@R+0pc@C+0pc{
\underset{r}{\mathrm{homotopycolimit}}~\underset{\mathrm{Spec}}{\mathcal{O}}^\mathrm{BK}{\Phi}^r_{\psi,\Gamma,A},\underset{I}{\mathrm{homotopylimit}}~\underset{\mathrm{Spec}}{\mathcal{O}}^\mathrm{BK}{\Phi}^I_{\psi,\Gamma,A}.	
}
\]
\[  
\xymatrix@R+0pc@C+0pc{
\underset{r}{\mathrm{homotopycolimit}}~\underset{\mathrm{Spec}}{\mathcal{O}}^\mathrm{BK}\widetilde{\Phi}^r_{\psi,\Gamma,A}/\mathrm{Fro}^\mathbb{Z},\underset{I}{\mathrm{homotopylimit}}~\underset{\mathrm{Spec}}{\mathcal{O}}^\mathrm{BK}\widetilde{\Phi}^I_{\psi,\Gamma,A}/\mathrm{Fro}^\mathbb{Z},	\\
}
\]
\[ 
\xymatrix@R+0pc@C+0pc{
\underset{r}{\mathrm{homotopycolimit}}~\underset{\mathrm{Spec}}{\mathcal{O}}^\mathrm{BK}\breve{\Phi}^r_{\psi,\Gamma,A}/\mathrm{Fro}^\mathbb{Z},\underset{I}{\mathrm{homotopylimit}}~\underset{\mathrm{Spec}}{\mathcal{O}}^\mathrm{BK}\breve{\Phi}^I_{\psi,\Gamma,A}/\mathrm{Fro}^\mathbb{Z},	\\
}
\]
\[ 
\xymatrix@R+0pc@C+0pc{
\underset{r}{\mathrm{homotopycolimit}}~\underset{\mathrm{Spec}}{\mathcal{O}}^\mathrm{BK}{\Phi}^r_{\psi,\Gamma,A}/\mathrm{Fro}^\mathbb{Z},\underset{I}{\mathrm{homotopylimit}}~\underset{\mathrm{Spec}}{\mathcal{O}}^\mathrm{BK}{\Phi}^I_{\psi,\Gamma,A}/\mathrm{Fro}^\mathbb{Z}.	
}
\]

\end{definition}

\begin{definition}
We then consider the corresponding quasisheaves of the corresponding condensed solid topological modules from \cite{10CS2}:
\begin{align}
\mathrm{Quasicoherentsheaves, Condensed}_{*}	
\end{align}
where $*$ is one of the following spaces:
\begin{align}
&\underset{\mathrm{Spec}}{\mathcal{O}}^\mathrm{CS}\widetilde{\Delta}_{\psi,\Gamma,A}/\mathrm{Fro}^\mathbb{Z},\underset{\mathrm{Spec}}{\mathcal{O}}^\mathrm{CS}\widetilde{\nabla}_{\psi,\Gamma,A}/\mathrm{Fro}^\mathbb{Z},\underset{\mathrm{Spec}}{\mathcal{O}}^\mathrm{CS}\widetilde{\Phi}_{\psi,\Gamma,A}/\mathrm{Fro}^\mathbb{Z},\underset{\mathrm{Spec}}{\mathcal{O}}^\mathrm{CS}\widetilde{\Delta}^+_{\psi,\Gamma,A}/\mathrm{Fro}^\mathbb{Z},\\
&\underset{\mathrm{Spec}}{\mathcal{O}}^\mathrm{CS}\widetilde{\nabla}^+_{\psi,\Gamma,A}/\mathrm{Fro}^\mathbb{Z},\underset{\mathrm{Spec}}{\mathcal{O}}^\mathrm{CS}\widetilde{\Delta}^\dagger_{\psi,\Gamma,A}/\mathrm{Fro}^\mathbb{Z},\underset{\mathrm{Spec}}{\mathcal{O}}^\mathrm{CS}\widetilde{\nabla}^\dagger_{\psi,\Gamma,A}/\mathrm{Fro}^\mathbb{Z},	\\
\end{align}
\begin{align}
&\underset{\mathrm{Spec}}{\mathcal{O}}^\mathrm{CS}\breve{\Delta}_{\psi,\Gamma,A}/\mathrm{Fro}^\mathbb{Z},\breve{\nabla}_{\psi,\Gamma,A}/\mathrm{Fro}^\mathbb{Z},\underset{\mathrm{Spec}}{\mathcal{O}}^\mathrm{CS}\breve{\Phi}_{\psi,\Gamma,A}/\mathrm{Fro}^\mathbb{Z},\underset{\mathrm{Spec}}{\mathcal{O}}^\mathrm{CS}\breve{\Delta}^+_{\psi,\Gamma,A}/\mathrm{Fro}^\mathbb{Z},\\
&\underset{\mathrm{Spec}}{\mathcal{O}}^\mathrm{CS}\breve{\nabla}^+_{\psi,\Gamma,A}/\mathrm{Fro}^\mathbb{Z},\underset{\mathrm{Spec}}{\mathcal{O}}^\mathrm{CS}\breve{\Delta}^\dagger_{\psi,\Gamma,A}/\mathrm{Fro}^\mathbb{Z},\underset{\mathrm{Spec}}{\mathcal{O}}^\mathrm{CS}\breve{\nabla}^\dagger_{\psi,\Gamma,A}/\mathrm{Fro}^\mathbb{Z},	\\
\end{align}
\begin{align}
&\underset{\mathrm{Spec}}{\mathcal{O}}^\mathrm{CS}{\Delta}_{\psi,\Gamma,A}/\mathrm{Fro}^\mathbb{Z},\underset{\mathrm{Spec}}{\mathcal{O}}^\mathrm{CS}{\nabla}_{\psi,\Gamma,A}/\mathrm{Fro}^\mathbb{Z},\underset{\mathrm{Spec}}{\mathcal{O}}^\mathrm{CS}{\Phi}_{\psi,\Gamma,A}/\mathrm{Fro}^\mathbb{Z},\underset{\mathrm{Spec}}{\mathcal{O}}^\mathrm{CS}{\Delta}^+_{\psi,\Gamma,A}/\mathrm{Fro}^\mathbb{Z},\\
&\underset{\mathrm{Spec}}{\mathcal{O}}^\mathrm{CS}{\nabla}^+_{\psi,\Gamma,A}/\mathrm{Fro}^\mathbb{Z}, \underset{\mathrm{Spec}}{\mathcal{O}}^\mathrm{CS}{\Delta}^\dagger_{\psi,\Gamma,A}/\mathrm{Fro}^\mathbb{Z},\underset{\mathrm{Spec}}{\mathcal{O}}^\mathrm{CS}{\nabla}^\dagger_{\psi,\Gamma,A}/\mathrm{Fro}^\mathbb{Z}.	
\end{align}
Here for those space with notations related to the radius and the corresponding interval we consider the total unions $\bigcap_r,\bigcup_I$ in order to achieve the whole spaces to achieve the analogues of the corresponding FF curves from \cite{10KL1}, \cite{10KL2}, \cite{10FF} for
\[
\xymatrix@R+0pc@C+0pc{
\underset{r}{\mathrm{homotopycolimit}}~\underset{\mathrm{Spec}}{\mathcal{O}}^\mathrm{CS}\widetilde{\Phi}^r_{\psi,\Gamma,A},\underset{I}{\mathrm{homotopylimit}}~\underset{\mathrm{Spec}}{\mathcal{O}}^\mathrm{CS}\widetilde{\Phi}^I_{\psi,\Gamma,A},	\\
}
\]
\[
\xymatrix@R+0pc@C+0pc{
\underset{r}{\mathrm{homotopycolimit}}~\underset{\mathrm{Spec}}{\mathcal{O}}^\mathrm{CS}\breve{\Phi}^r_{\psi,\Gamma,A},\underset{I}{\mathrm{homotopylimit}}~\underset{\mathrm{Spec}}{\mathcal{O}}^\mathrm{CS}\breve{\Phi}^I_{\psi,\Gamma,A},	\\
}
\]
\[
\xymatrix@R+0pc@C+0pc{
\underset{r}{\mathrm{homotopycolimit}}~\underset{\mathrm{Spec}}{\mathcal{O}}^\mathrm{CS}{\Phi}^r_{\psi,\Gamma,A},\underset{I}{\mathrm{homotopylimit}}~\underset{\mathrm{Spec}}{\mathcal{O}}^\mathrm{CS}{\Phi}^I_{\psi,\Gamma,A}.	
}
\]
\[ 
\xymatrix@R+0pc@C+0pc{
\underset{r}{\mathrm{homotopycolimit}}~\underset{\mathrm{Spec}}{\mathcal{O}}^\mathrm{CS}\widetilde{\Phi}^r_{\psi,\Gamma,A}/\mathrm{Fro}^\mathbb{Z},\underset{I}{\mathrm{homotopylimit}}~\underset{\mathrm{Spec}}{\mathcal{O}}^\mathrm{CS}\widetilde{\Phi}^I_{\psi,\Gamma,A}/\mathrm{Fro}^\mathbb{Z},	\\
}
\]
\[ 
\xymatrix@R+0pc@C+0pc{
\underset{r}{\mathrm{homotopycolimit}}~\underset{\mathrm{Spec}}{\mathcal{O}}^\mathrm{CS}\breve{\Phi}^r_{\psi,\Gamma,A}/\mathrm{Fro}^\mathbb{Z},\underset{I}{\mathrm{homotopylimit}}~\breve{\Phi}^I_{\psi,\Gamma,A}/\mathrm{Fro}^\mathbb{Z},	\\
}
\]
\[ 
\xymatrix@R+0pc@C+0pc{
\underset{r}{\mathrm{homotopycolimit}}~\underset{\mathrm{Spec}}{\mathcal{O}}^\mathrm{CS}{\Phi}^r_{\psi,\Gamma,A}/\mathrm{Fro}^\mathbb{Z},\underset{I}{\mathrm{homotopylimit}}~\underset{\mathrm{Spec}}{\mathcal{O}}^\mathrm{CS}{\Phi}^I_{\psi,\Gamma,A}/\mathrm{Fro}^\mathbb{Z}.	
}
\]

\end{definition}

\

\begin{proposition}
There is a well-defined functor from the $\infty$-category 
\begin{align}
\mathrm{Quasicoherentpresheaves,Condensed}_{*}	
\end{align}
where $*$ is one of the following spaces:
\begin{align}
&\underset{\mathrm{Spec}}{\mathcal{O}}^\mathrm{CS}\widetilde{\Phi}_{\psi,\Gamma,A}/\mathrm{Fro}^\mathbb{Z},	\\
\end{align}
\begin{align}
&\underset{\mathrm{Spec}}{\mathcal{O}}^\mathrm{CS}\breve{\Phi}_{\psi,\Gamma,A}/\mathrm{Fro}^\mathbb{Z},	\\
\end{align}
\begin{align}
&\underset{\mathrm{Spec}}{\mathcal{O}}^\mathrm{CS}{\Phi}_{\psi,\Gamma,A}/\mathrm{Fro}^\mathbb{Z},	
\end{align}
to the $\infty$-category of $\mathrm{Fro}$-equivariant quasicoherent presheaves over similar spaces above correspondingly without the $\mathrm{Fro}$-quotients, and to the $\infty$-category of $\mathrm{Fro}$-equivariant quasicoherent modules over global sections of the structure $\infty$-sheaves of the similar spaces above correspondingly without the $\mathrm{Fro}$-quotients. Here for those space without notation related to the radius and the corresponding interval we consider the total unions $\bigcap_r,\bigcup_I$ in order to achieve the whole spaces to achieve the analogues of the corresponding FF curves from \cite{10KL1}, \cite{10KL2}, \cite{10FF} for
\[
\xymatrix@R+0pc@C+0pc{
\underset{r}{\mathrm{homotopycolimit}}~\underset{\mathrm{Spec}}{\mathcal{O}}^\mathrm{CS}\widetilde{\Phi}^r_{\psi,\Gamma,A},\underset{I}{\mathrm{homotopylimit}}~\underset{\mathrm{Spec}}{\mathcal{O}}^\mathrm{CS}\widetilde{\Phi}^I_{\psi,\Gamma,A},	\\
}
\]
\[
\xymatrix@R+0pc@C+0pc{
\underset{r}{\mathrm{homotopycolimit}}~\underset{\mathrm{Spec}}{\mathcal{O}}^\mathrm{CS}\breve{\Phi}^r_{\psi,\Gamma,A},\underset{I}{\mathrm{homotopylimit}}~\underset{\mathrm{Spec}}{\mathcal{O}}^\mathrm{CS}\breve{\Phi}^I_{\psi,\Gamma,A},	\\
}
\]
\[
\xymatrix@R+0pc@C+0pc{
\underset{r}{\mathrm{homotopycolimit}}~\underset{\mathrm{Spec}}{\mathcal{O}}^\mathrm{CS}{\Phi}^r_{\psi,\Gamma,A},\underset{I}{\mathrm{homotopylimit}}~\underset{\mathrm{Spec}}{\mathcal{O}}^\mathrm{CS}{\Phi}^I_{\psi,\Gamma,A}.	
}
\]
\[ 
\xymatrix@R+0pc@C+0pc{
\underset{r}{\mathrm{homotopycolimit}}~\underset{\mathrm{Spec}}{\mathcal{O}}^\mathrm{CS}\widetilde{\Phi}^r_{\psi,\Gamma,A}/\mathrm{Fro}^\mathbb{Z},\underset{I}{\mathrm{homotopylimit}}~\underset{\mathrm{Spec}}{\mathcal{O}}^\mathrm{CS}\widetilde{\Phi}^I_{\psi,\Gamma,A}/\mathrm{Fro}^\mathbb{Z},	\\
}
\]
\[ 
\xymatrix@R+0pc@C+0pc{
\underset{r}{\mathrm{homotopycolimit}}~\underset{\mathrm{Spec}}{\mathcal{O}}^\mathrm{CS}\breve{\Phi}^r_{\psi,\Gamma,A}/\mathrm{Fro}^\mathbb{Z},\underset{I}{\mathrm{homotopylimit}}~\breve{\Phi}^I_{\psi,\Gamma,A}/\mathrm{Fro}^\mathbb{Z},	\\
}
\]
\[ 
\xymatrix@R+0pc@C+0pc{
\underset{r}{\mathrm{homotopycolimit}}~\underset{\mathrm{Spec}}{\mathcal{O}}^\mathrm{CS}{\Phi}^r_{\psi,\Gamma,A}/\mathrm{Fro}^\mathbb{Z},\underset{I}{\mathrm{homotopylimit}}~\underset{\mathrm{Spec}}{\mathcal{O}}^\mathrm{CS}{\Phi}^I_{\psi,\Gamma,A}/\mathrm{Fro}^\mathbb{Z}.	
}
\]	
In this situation we will have the target category being family parametrized by $r$ or $I$ in compatible glueing sense as in \cite[Definition 5.4.10]{10KL2}. In this situation for modules parametrized by the intervals we have the equivalence of $\infty$-categories by using \cite[Proposition 13.8]{10CS2}. Here the corresponding quasicoherent Frobenius modules are defined to be the corresponding homotopy colimits and limits of Frobenius modules:
\begin{align}
\underset{r}{\mathrm{homotopycolimit}}~M_r,\\
\underset{I}{\mathrm{homotopylimit}}~M_I,	
\end{align}
where each $M_r$ is a Frobenius-equivariant module over the period ring with respect to some radius $r$ while each $M_I$ is a Frobenius-equivariant module over the period ring with respect to some interval $I$.\\
\end{proposition}

\begin{proposition}
Similar proposition holds for 
\begin{align}
\mathrm{Quasicoherentsheaves,IndBanach}_{*}.	
\end{align}	
\end{proposition}

\

\begin{definition}
We then consider the corresponding quasipresheaves of perfect complexes the corresponding ind-Banach or monomorphic ind-Banach modules from \cite{10BBK}, \cite{10KKM}:
\begin{align}
\mathrm{Quasicoherentpresheaves,Perfectcomplex,IndBanach}_{*}	
\end{align}
where $*$ is one of the following spaces:
\begin{align}
&\underset{\mathrm{Spec}}{\mathcal{O}}^\mathrm{BK}\widetilde{\Phi}_{\psi,\Gamma,A}/\mathrm{Fro}^\mathbb{Z},	\\
\end{align}
\begin{align}
&\underset{\mathrm{Spec}}{\mathcal{O}}^\mathrm{BK}\breve{\Phi}_{\psi,\Gamma,A}/\mathrm{Fro}^\mathbb{Z},	\\
\end{align}
\begin{align}
&\underset{\mathrm{Spec}}{\mathcal{O}}^\mathrm{BK}{\Phi}_{\psi,\Gamma,A}/\mathrm{Fro}^\mathbb{Z}.	
\end{align}
Here for those space without notation related to the radius and the corresponding interval we consider the total unions $\bigcap_r,\bigcup_I$ in order to achieve the whole spaces to achieve the analogues of the corresponding FF curves from \cite{10KL1}, \cite{10KL2}, \cite{10FF} for
\[
\xymatrix@R+0pc@C+0pc{
\underset{r}{\mathrm{homotopycolimit}}~\underset{\mathrm{Spec}}{\mathcal{O}}^\mathrm{BK}\widetilde{\Phi}^r_{\psi,\Gamma,A},\underset{I}{\mathrm{homotopylimit}}~\underset{\mathrm{Spec}}{\mathcal{O}}^\mathrm{BK}\widetilde{\Phi}^I_{\psi,\Gamma,A},	\\
}
\]
\[
\xymatrix@R+0pc@C+0pc{
\underset{r}{\mathrm{homotopycolimit}}~\underset{\mathrm{Spec}}{\mathcal{O}}^\mathrm{BK}\breve{\Phi}^r_{\psi,\Gamma,A},\underset{I}{\mathrm{homotopylimit}}~\underset{\mathrm{Spec}}{\mathcal{O}}^\mathrm{BK}\breve{\Phi}^I_{\psi,\Gamma,A},	\\
}
\]
\[
\xymatrix@R+0pc@C+0pc{
\underset{r}{\mathrm{homotopycolimit}}~\underset{\mathrm{Spec}}{\mathcal{O}}^\mathrm{BK}{\Phi}^r_{\psi,\Gamma,A},\underset{I}{\mathrm{homotopylimit}}~\underset{\mathrm{Spec}}{\mathcal{O}}^\mathrm{BK}{\Phi}^I_{\psi,\Gamma,A}.	
}
\]
\[  
\xymatrix@R+0pc@C+0pc{
\underset{r}{\mathrm{homotopycolimit}}~\underset{\mathrm{Spec}}{\mathcal{O}}^\mathrm{BK}\widetilde{\Phi}^r_{\psi,\Gamma,A}/\mathrm{Fro}^\mathbb{Z},\underset{I}{\mathrm{homotopylimit}}~\underset{\mathrm{Spec}}{\mathcal{O}}^\mathrm{BK}\widetilde{\Phi}^I_{\psi,\Gamma,A}/\mathrm{Fro}^\mathbb{Z},	\\
}
\]
\[ 
\xymatrix@R+0pc@C+0pc{
\underset{r}{\mathrm{homotopycolimit}}~\underset{\mathrm{Spec}}{\mathcal{O}}^\mathrm{BK}\breve{\Phi}^r_{\psi,\Gamma,A}/\mathrm{Fro}^\mathbb{Z},\underset{I}{\mathrm{homotopylimit}}~\underset{\mathrm{Spec}}{\mathcal{O}}^\mathrm{BK}\breve{\Phi}^I_{\psi,\Gamma,A}/\mathrm{Fro}^\mathbb{Z},	\\
}
\]
\[ 
\xymatrix@R+0pc@C+0pc{
\underset{r}{\mathrm{homotopycolimit}}~\underset{\mathrm{Spec}}{\mathcal{O}}^\mathrm{BK}{\Phi}^r_{\psi,\Gamma,A}/\mathrm{Fro}^\mathbb{Z},\underset{I}{\mathrm{homotopylimit}}~\underset{\mathrm{Spec}}{\mathcal{O}}^\mathrm{BK}{\Phi}^I_{\psi,\Gamma,A}/\mathrm{Fro}^\mathbb{Z}.	
}
\]

\end{definition}

\begin{definition}
We then consider the corresponding quasisheaves of perfect complexes of the corresponding condensed solid topological modules from \cite{10CS2}:
\begin{align}
\mathrm{Quasicoherentsheaves, Perfectcomplex, Condensed}_{*}	
\end{align}
where $*$ is one of the following spaces:
\begin{align}
&\underset{\mathrm{Spec}}{\mathcal{O}}^\mathrm{CS}\widetilde{\Delta}_{\psi,\Gamma,A}/\mathrm{Fro}^\mathbb{Z},\underset{\mathrm{Spec}}{\mathcal{O}}^\mathrm{CS}\widetilde{\nabla}_{\psi,\Gamma,A}/\mathrm{Fro}^\mathbb{Z},\underset{\mathrm{Spec}}{\mathcal{O}}^\mathrm{CS}\widetilde{\Phi}_{\psi,\Gamma,A}/\mathrm{Fro}^\mathbb{Z},\underset{\mathrm{Spec}}{\mathcal{O}}^\mathrm{CS}\widetilde{\Delta}^+_{\psi,\Gamma,A}/\mathrm{Fro}^\mathbb{Z},\\
&\underset{\mathrm{Spec}}{\mathcal{O}}^\mathrm{CS}\widetilde{\nabla}^+_{\psi,\Gamma,A}/\mathrm{Fro}^\mathbb{Z},\underset{\mathrm{Spec}}{\mathcal{O}}^\mathrm{CS}\widetilde{\Delta}^\dagger_{\psi,\Gamma,A}/\mathrm{Fro}^\mathbb{Z},\underset{\mathrm{Spec}}{\mathcal{O}}^\mathrm{CS}\widetilde{\nabla}^\dagger_{\psi,\Gamma,A}/\mathrm{Fro}^\mathbb{Z},	\\
\end{align}
\begin{align}
&\underset{\mathrm{Spec}}{\mathcal{O}}^\mathrm{CS}\breve{\Delta}_{\psi,\Gamma,A}/\mathrm{Fro}^\mathbb{Z},\breve{\nabla}_{\psi,\Gamma,A}/\mathrm{Fro}^\mathbb{Z},\underset{\mathrm{Spec}}{\mathcal{O}}^\mathrm{CS}\breve{\Phi}_{\psi,\Gamma,A}/\mathrm{Fro}^\mathbb{Z},\underset{\mathrm{Spec}}{\mathcal{O}}^\mathrm{CS}\breve{\Delta}^+_{\psi,\Gamma,A}/\mathrm{Fro}^\mathbb{Z},\\
&\underset{\mathrm{Spec}}{\mathcal{O}}^\mathrm{CS}\breve{\nabla}^+_{\psi,\Gamma,A}/\mathrm{Fro}^\mathbb{Z},\underset{\mathrm{Spec}}{\mathcal{O}}^\mathrm{CS}\breve{\Delta}^\dagger_{\psi,\Gamma,A}/\mathrm{Fro}^\mathbb{Z},\underset{\mathrm{Spec}}{\mathcal{O}}^\mathrm{CS}\breve{\nabla}^\dagger_{\psi,\Gamma,A}/\mathrm{Fro}^\mathbb{Z},	\\
\end{align}
\begin{align}
&\underset{\mathrm{Spec}}{\mathcal{O}}^\mathrm{CS}{\Delta}_{\psi,\Gamma,A}/\mathrm{Fro}^\mathbb{Z},\underset{\mathrm{Spec}}{\mathcal{O}}^\mathrm{CS}{\nabla}_{\psi,\Gamma,A}/\mathrm{Fro}^\mathbb{Z},\underset{\mathrm{Spec}}{\mathcal{O}}^\mathrm{CS}{\Phi}_{\psi,\Gamma,A}/\mathrm{Fro}^\mathbb{Z},\underset{\mathrm{Spec}}{\mathcal{O}}^\mathrm{CS}{\Delta}^+_{\psi,\Gamma,A}/\mathrm{Fro}^\mathbb{Z},\\
&\underset{\mathrm{Spec}}{\mathcal{O}}^\mathrm{CS}{\nabla}^+_{\psi,\Gamma,A}/\mathrm{Fro}^\mathbb{Z}, \underset{\mathrm{Spec}}{\mathcal{O}}^\mathrm{CS}{\Delta}^\dagger_{\psi,\Gamma,A}/\mathrm{Fro}^\mathbb{Z},\underset{\mathrm{Spec}}{\mathcal{O}}^\mathrm{CS}{\nabla}^\dagger_{\psi,\Gamma,A}/\mathrm{Fro}^\mathbb{Z}.	
\end{align}
Here for those space with notations related to the radius and the corresponding interval we consider the total unions $\bigcap_r,\bigcup_I$ in order to achieve the whole spaces to achieve the analogues of the corresponding FF curves from \cite{10KL1}, \cite{10KL2}, \cite{10FF} for
\[
\xymatrix@R+0pc@C+0pc{
\underset{r}{\mathrm{homotopycolimit}}~\underset{\mathrm{Spec}}{\mathcal{O}}^\mathrm{CS}\widetilde{\Phi}^r_{\psi,\Gamma,A},\underset{I}{\mathrm{homotopylimit}}~\underset{\mathrm{Spec}}{\mathcal{O}}^\mathrm{CS}\widetilde{\Phi}^I_{\psi,\Gamma,A},	\\
}
\]
\[
\xymatrix@R+0pc@C+0pc{
\underset{r}{\mathrm{homotopycolimit}}~\underset{\mathrm{Spec}}{\mathcal{O}}^\mathrm{CS}\breve{\Phi}^r_{\psi,\Gamma,A},\underset{I}{\mathrm{homotopylimit}}~\underset{\mathrm{Spec}}{\mathcal{O}}^\mathrm{CS}\breve{\Phi}^I_{\psi,\Gamma,A},	\\
}
\]
\[
\xymatrix@R+0pc@C+0pc{
\underset{r}{\mathrm{homotopycolimit}}~\underset{\mathrm{Spec}}{\mathcal{O}}^\mathrm{CS}{\Phi}^r_{\psi,\Gamma,A},\underset{I}{\mathrm{homotopylimit}}~\underset{\mathrm{Spec}}{\mathcal{O}}^\mathrm{CS}{\Phi}^I_{\psi,\Gamma,A}.	
}
\]
\[ 
\xymatrix@R+0pc@C+0pc{
\underset{r}{\mathrm{homotopycolimit}}~\underset{\mathrm{Spec}}{\mathcal{O}}^\mathrm{CS}\widetilde{\Phi}^r_{\psi,\Gamma,A}/\mathrm{Fro}^\mathbb{Z},\underset{I}{\mathrm{homotopylimit}}~\underset{\mathrm{Spec}}{\mathcal{O}}^\mathrm{CS}\widetilde{\Phi}^I_{\psi,\Gamma,A}/\mathrm{Fro}^\mathbb{Z},	\\
}
\]
\[ 
\xymatrix@R+0pc@C+0pc{
\underset{r}{\mathrm{homotopycolimit}}~\underset{\mathrm{Spec}}{\mathcal{O}}^\mathrm{CS}\breve{\Phi}^r_{\psi,\Gamma,A}/\mathrm{Fro}^\mathbb{Z},\underset{I}{\mathrm{homotopylimit}}~\breve{\Phi}^I_{\psi,\Gamma,A}/\mathrm{Fro}^\mathbb{Z},	\\
}
\]
\[ 
\xymatrix@R+0pc@C+0pc{
\underset{r}{\mathrm{homotopycolimit}}~\underset{\mathrm{Spec}}{\mathcal{O}}^\mathrm{CS}{\Phi}^r_{\psi,\Gamma,A}/\mathrm{Fro}^\mathbb{Z},\underset{I}{\mathrm{homotopylimit}}~\underset{\mathrm{Spec}}{\mathcal{O}}^\mathrm{CS}{\Phi}^I_{\psi,\Gamma,A}/\mathrm{Fro}^\mathbb{Z}.	
}
\]

\end{definition}

\begin{proposition}
There is a well-defined functor from the $\infty$-category 
\begin{align}
\mathrm{Quasicoherentpresheaves,Perfectcomplex,Condensed}_{*}	
\end{align}
where $*$ is one of the following spaces:
\begin{align}
&\underset{\mathrm{Spec}}{\mathcal{O}}^\mathrm{CS}\widetilde{\Phi}_{\psi,\Gamma,A}/\mathrm{Fro}^\mathbb{Z},	\\
\end{align}
\begin{align}
&\underset{\mathrm{Spec}}{\mathcal{O}}^\mathrm{CS}\breve{\Phi}_{\psi,\Gamma,A}/\mathrm{Fro}^\mathbb{Z},	\\
\end{align}
\begin{align}
&\underset{\mathrm{Spec}}{\mathcal{O}}^\mathrm{CS}{\Phi}_{\psi,\Gamma,A}/\mathrm{Fro}^\mathbb{Z},	
\end{align}
to the $\infty$-category of $\mathrm{Fro}$-equivariant quasicoherent presheaves over similar spaces above correspondingly without the $\mathrm{Fro}$-quotients, and to the $\infty$-category of $\mathrm{Fro}$-equivariant quasicoherent modules over global sections of the structure $\infty$-sheaves of the similar spaces above correspondingly without the $\mathrm{Fro}$-quotients. Here for those space without notation related to the radius and the corresponding interval we consider the total unions $\bigcap_r,\bigcup_I$ in order to achieve the whole spaces to achieve the analogues of the corresponding FF curves from \cite{10KL1}, \cite{10KL2}, \cite{10FF} for
\[
\xymatrix@R+0pc@C+0pc{
\underset{r}{\mathrm{homotopycolimit}}~\underset{\mathrm{Spec}}{\mathcal{O}}^\mathrm{CS}\widetilde{\Phi}^r_{\psi,\Gamma,A},\underset{I}{\mathrm{homotopylimit}}~\underset{\mathrm{Spec}}{\mathcal{O}}^\mathrm{CS}\widetilde{\Phi}^I_{\psi,\Gamma,A},	\\
}
\]
\[
\xymatrix@R+0pc@C+0pc{
\underset{r}{\mathrm{homotopycolimit}}~\underset{\mathrm{Spec}}{\mathcal{O}}^\mathrm{CS}\breve{\Phi}^r_{\psi,\Gamma,A},\underset{I}{\mathrm{homotopylimit}}~\underset{\mathrm{Spec}}{\mathcal{O}}^\mathrm{CS}\breve{\Phi}^I_{\psi,\Gamma,A},	\\
}
\]
\[
\xymatrix@R+0pc@C+0pc{
\underset{r}{\mathrm{homotopycolimit}}~\underset{\mathrm{Spec}}{\mathcal{O}}^\mathrm{CS}{\Phi}^r_{\psi,\Gamma,A},\underset{I}{\mathrm{homotopylimit}}~\underset{\mathrm{Spec}}{\mathcal{O}}^\mathrm{CS}{\Phi}^I_{\psi,\Gamma,A}.	
}
\]
\[ 
\xymatrix@R+0pc@C+0pc{
\underset{r}{\mathrm{homotopycolimit}}~\underset{\mathrm{Spec}}{\mathcal{O}}^\mathrm{CS}\widetilde{\Phi}^r_{\psi,\Gamma,A}/\mathrm{Fro}^\mathbb{Z},\underset{I}{\mathrm{homotopylimit}}~\underset{\mathrm{Spec}}{\mathcal{O}}^\mathrm{CS}\widetilde{\Phi}^I_{\psi,\Gamma,A}/\mathrm{Fro}^\mathbb{Z},	\\
}
\]
\[ 
\xymatrix@R+0pc@C+0pc{
\underset{r}{\mathrm{homotopycolimit}}~\underset{\mathrm{Spec}}{\mathcal{O}}^\mathrm{CS}\breve{\Phi}^r_{\psi,\Gamma,A}/\mathrm{Fro}^\mathbb{Z},\underset{I}{\mathrm{homotopylimit}}~\breve{\Phi}^I_{\psi,\Gamma,A}/\mathrm{Fro}^\mathbb{Z},	\\
}
\]
\[ 
\xymatrix@R+0pc@C+0pc{
\underset{r}{\mathrm{homotopycolimit}}~\underset{\mathrm{Spec}}{\mathcal{O}}^\mathrm{CS}{\Phi}^r_{\psi,\Gamma,A}/\mathrm{Fro}^\mathbb{Z},\underset{I}{\mathrm{homotopylimit}}~\underset{\mathrm{Spec}}{\mathcal{O}}^\mathrm{CS}{\Phi}^I_{\psi,\Gamma,A}/\mathrm{Fro}^\mathbb{Z}.	
}
\]	
In this situation we will have the target category being family parametrized by $r$ or $I$ in compatible glueing sense as in \cite[Definition 5.4.10]{10KL2}. In this situation for modules parametrized by the intervals we have the equivalence of $\infty$-categories by using \cite[Proposition 12.18]{10CS2}. Here the corresponding quasicoherent Frobenius modules are defined to be the corresponding homotopy colimits and limits of Frobenius modules:
\begin{align}
\underset{r}{\mathrm{homotopycolimit}}~M_r,\\
\underset{I}{\mathrm{homotopylimit}}~M_I,	
\end{align}
where each $M_r$ is a Frobenius-equivariant module over the period ring with respect to some radius $r$ while each $M_I$ is a Frobenius-equivariant module over the period ring with respect to some interval $I$.\\
\end{proposition}

\begin{proposition}
Similar proposition holds for 
\begin{align}
\mathrm{Quasicoherentsheaves,Perfectcomplex,IndBanach}_{*}.	
\end{align}	
\end{proposition}

\newpage
\subsection{Frobenius Quasicoherent Modules II: Deformation in Banach Rings}

\begin{definition}
Let $\psi$ be a toric tower over $\mathbb{Q}_p$ as in \cite[Chapter 7]{10KL2} with base $\mathbb{Q}_p\left<X_1^{\pm 1},...,X_k^{\pm 1}\right>$. Then from \cite{10KL1} and \cite[Definition 5.2.1]{10KL2} we have the following class of Kedlaya-Liu rings (with the following replacement: $\Delta$ stands for $A$, $\nabla$ stands for $B$, while $\Phi$ stands for $C$) by taking product in the sense of self $\Gamma$-th power:

\[
\xymatrix@R+0pc@C+0pc{
\widetilde{\Delta}_{\psi,\Gamma},\widetilde{\nabla}_{\psi,\Gamma},\widetilde{\Phi}_{\psi,\Gamma},\widetilde{\Delta}^+_{\psi,\Gamma},\widetilde{\nabla}^+_{\psi,\Gamma},\widetilde{\Delta}^\dagger_{\psi,\Gamma},\widetilde{\nabla}^\dagger_{\psi,\Gamma},\widetilde{\Phi}^r_{\psi,\Gamma},\widetilde{\Phi}^I_{\psi,\Gamma}, 
}
\]

\[
\xymatrix@R+0pc@C+0pc{
\breve{\Delta}_{\psi,\Gamma},\breve{\nabla}_{\psi,\Gamma},\breve{\Phi}_{\psi,\Gamma},\breve{\Delta}^+_{\psi,\Gamma},\breve{\nabla}^+_{\psi,\Gamma},\breve{\Delta}^\dagger_{\psi,\Gamma},\breve{\nabla}^\dagger_{\psi,\Gamma},\breve{\Phi}^r_{\psi,\Gamma},\breve{\Phi}^I_{\psi,\Gamma},	
}
\]

\[
\xymatrix@R+0pc@C+0pc{
{\Delta}_{\psi,\Gamma},{\nabla}_{\psi,\Gamma},{\Phi}_{\psi,\Gamma},{\Delta}^+_{\psi,\Gamma},{\nabla}^+_{\psi,\Gamma},{\Delta}^\dagger_{\psi,\Gamma},{\nabla}^\dagger_{\psi,\Gamma},{\Phi}^r_{\psi,\Gamma},{\Phi}^I_{\psi,\Gamma}.	
}
\]
We now consider the following rings with $-$ being any deforming Banach ring over $\mathbb{Q}_p$. Taking the product we have:
\[
\xymatrix@R+0pc@C+0pc{
\widetilde{\Phi}_{\psi,\Gamma,-},\widetilde{\Phi}^r_{\psi,\Gamma,-},\widetilde{\Phi}^I_{\psi,\Gamma,-},	
}
\]
\[
\xymatrix@R+0pc@C+0pc{
\breve{\Phi}_{\psi,\Gamma,-},\breve{\Phi}^r_{\psi,\Gamma,-},\breve{\Phi}^I_{\psi,\Gamma,-},	
}
\]
\[
\xymatrix@R+0pc@C+0pc{
{\Phi}_{\psi,\Gamma,-},{\Phi}^r_{\psi,\Gamma,-},{\Phi}^I_{\psi,\Gamma,-}.	
}
\]
They carry multi Frobenius action $\varphi_\Gamma$ and multi $\mathrm{Lie}_\Gamma:=\mathbb{Z}_p^{\times\Gamma}$ action. In our current situation after \cite{10CKZ} and \cite{10PZ} we consider the following $(\infty,1)$-categories of $(\infty,1)$-modules.\\
\end{definition}

\begin{definition}
First we consider the Bambozzi-Kremnizer spectrum $\underset{\mathrm{Spec}}{\mathcal{O}}^\mathrm{BK}(*)$ attached to any of those in the above from \cite{10BK} by taking derived rational localization:
\begin{align}
&\underset{\mathrm{Spec}}{\mathcal{O}}^\mathrm{BK}\widetilde{\Phi}_{\psi,\Gamma,-},\underset{\mathrm{Spec}}{\mathcal{O}}^\mathrm{BK}\widetilde{\Phi}^r_{\psi,\Gamma,-},\underset{\mathrm{Spec}}{\mathcal{O}}^\mathrm{BK}\widetilde{\Phi}^I_{\psi,\Gamma,-},	
\end{align}
\begin{align}
&\underset{\mathrm{Spec}}{\mathcal{O}}^\mathrm{BK}\breve{\Phi}_{\psi,\Gamma,-},\underset{\mathrm{Spec}}{\mathcal{O}}^\mathrm{BK}\breve{\Phi}^r_{\psi,\Gamma,-},\underset{\mathrm{Spec}}{\mathcal{O}}^\mathrm{BK}\breve{\Phi}^I_{\psi,\Gamma,-},	
\end{align}
\begin{align}
&\underset{\mathrm{Spec}}{\mathcal{O}}^\mathrm{BK}{\Phi}_{\psi,\Gamma,-},
\underset{\mathrm{Spec}}{\mathcal{O}}^\mathrm{BK}{\Phi}^r_{\psi,\Gamma,-},\underset{\mathrm{Spec}}{\mathcal{O}}^\mathrm{BK}{\Phi}^I_{\psi,\Gamma,-}.	
\end{align}

Then we take the corresponding quotients by using the corresponding Frobenius operators:
\begin{align}
&\underset{\mathrm{Spec}}{\mathcal{O}}^\mathrm{BK}\widetilde{\Phi}_{\psi,\Gamma,-}/\mathrm{Fro}^\mathbb{Z},	\\
\end{align}
\begin{align}
&\underset{\mathrm{Spec}}{\mathcal{O}}^\mathrm{BK}\breve{\Phi}_{\psi,\Gamma,-}/\mathrm{Fro}^\mathbb{Z},	\\
\end{align}
\begin{align}
&\underset{\mathrm{Spec}}{\mathcal{O}}^\mathrm{BK}{\Phi}_{\psi,\Gamma,-}/\mathrm{Fro}^\mathbb{Z}.	
\end{align}
Here for those space without notation related to the radius and the corresponding interval we consider the total unions $\bigcap_r,\bigcup_I$ in order to achieve the whole spaces to achieve the analogues of the corresponding FF curves from \cite{10KL1}, \cite{10KL2}, \cite{10FF} for
\[
\xymatrix@R+0pc@C+0pc{
\underset{r}{\mathrm{homotopycolimit}}~\underset{\mathrm{Spec}}{\mathcal{O}}^\mathrm{BK}\widetilde{\Phi}^r_{\psi,\Gamma,-},\underset{I}{\mathrm{homotopylimit}}~\underset{\mathrm{Spec}}{\mathcal{O}}^\mathrm{BK}\widetilde{\Phi}^I_{\psi,\Gamma,-},	\\
}
\]
\[
\xymatrix@R+0pc@C+0pc{
\underset{r}{\mathrm{homotopycolimit}}~\underset{\mathrm{Spec}}{\mathcal{O}}^\mathrm{BK}\breve{\Phi}^r_{\psi,\Gamma,-},\underset{I}{\mathrm{homotopylimit}}~\underset{\mathrm{Spec}}{\mathcal{O}}^\mathrm{BK}\breve{\Phi}^I_{\psi,\Gamma,-},	\\
}
\]
\[
\xymatrix@R+0pc@C+0pc{
\underset{r}{\mathrm{homotopycolimit}}~\underset{\mathrm{Spec}}{\mathcal{O}}^\mathrm{BK}{\Phi}^r_{\psi,\Gamma,-},\underset{I}{\mathrm{homotopylimit}}~\underset{\mathrm{Spec}}{\mathcal{O}}^\mathrm{BK}{\Phi}^I_{\psi,\Gamma,-}.	
}
\]
\[  
\xymatrix@R+0pc@C+0pc{
\underset{r}{\mathrm{homotopycolimit}}~\underset{\mathrm{Spec}}{\mathcal{O}}^\mathrm{BK}\widetilde{\Phi}^r_{\psi,\Gamma,-}/\mathrm{Fro}^\mathbb{Z},\underset{I}{\mathrm{homotopylimit}}~\underset{\mathrm{Spec}}{\mathcal{O}}^\mathrm{BK}\widetilde{\Phi}^I_{\psi,\Gamma,-}/\mathrm{Fro}^\mathbb{Z},	\\
}
\]
\[ 
\xymatrix@R+0pc@C+0pc{
\underset{r}{\mathrm{homotopycolimit}}~\underset{\mathrm{Spec}}{\mathcal{O}}^\mathrm{BK}\breve{\Phi}^r_{\psi,\Gamma,-}/\mathrm{Fro}^\mathbb{Z},\underset{I}{\mathrm{homotopylimit}}~\underset{\mathrm{Spec}}{\mathcal{O}}^\mathrm{BK}\breve{\Phi}^I_{\psi,\Gamma,-}/\mathrm{Fro}^\mathbb{Z},	\\
}
\]
\[ 
\xymatrix@R+0pc@C+0pc{
\underset{r}{\mathrm{homotopycolimit}}~\underset{\mathrm{Spec}}{\mathcal{O}}^\mathrm{BK}{\Phi}^r_{\psi,\Gamma,-}/\mathrm{Fro}^\mathbb{Z},\underset{I}{\mathrm{homotopylimit}}~\underset{\mathrm{Spec}}{\mathcal{O}}^\mathrm{BK}{\Phi}^I_{\psi,\Gamma,-}/\mathrm{Fro}^\mathbb{Z}.	
}
\]

\end{definition}

\indent Meanwhile we have the corresponding Clausen-Scholze analytic stacks from \cite{10CS2}, therefore applying their construction we have:

\begin{definition}
Here we define the following products by using the solidified tensor product from \cite{10CS1} and \cite{10CS2}. Namely $A$ will still as above as a Banach ring over $\mathbb{Q}_p$. Then we take solidified tensor product $\overset{\blacksquare}{\otimes}$ of any of the following
\[
\xymatrix@R+0pc@C+0pc{
\widetilde{\Delta}_{\psi,\Gamma},\widetilde{\nabla}_{\psi,\Gamma},\widetilde{\Phi}_{\psi,\Gamma},\widetilde{\Delta}^+_{\psi,\Gamma},\widetilde{\nabla}^+_{\psi,\Gamma},\widetilde{\Delta}^\dagger_{\psi,\Gamma},\widetilde{\nabla}^\dagger_{\psi,\Gamma},\widetilde{\Phi}^r_{\psi,\Gamma},\widetilde{\Phi}^I_{\psi,\Gamma}, 
}
\]

\[
\xymatrix@R+0pc@C+0pc{
\breve{\Delta}_{\psi,\Gamma},\breve{\nabla}_{\psi,\Gamma},\breve{\Phi}_{\psi,\Gamma},\breve{\Delta}^+_{\psi,\Gamma},\breve{\nabla}^+_{\psi,\Gamma},\breve{\Delta}^\dagger_{\psi,\Gamma},\breve{\nabla}^\dagger_{\psi,\Gamma},\breve{\Phi}^r_{\psi,\Gamma},\breve{\Phi}^I_{\psi,\Gamma},	
}
\]

\[
\xymatrix@R+0pc@C+0pc{
{\Delta}_{\psi,\Gamma},{\nabla}_{\psi,\Gamma},{\Phi}_{\psi,\Gamma},{\Delta}^+_{\psi,\Gamma},{\nabla}^+_{\psi,\Gamma},{\Delta}^\dagger_{\psi,\Gamma},{\nabla}^\dagger_{\psi,\Gamma},{\Phi}^r_{\psi,\Gamma},{\Phi}^I_{\psi,\Gamma},	
}
\]  	
with $A$. Then we have the notations:
\[
\xymatrix@R+0pc@C+0pc{
\widetilde{\Delta}_{\psi,\Gamma,-},\widetilde{\nabla}_{\psi,\Gamma,-},\widetilde{\Phi}_{\psi,\Gamma,-},\widetilde{\Delta}^+_{\psi,\Gamma,-},\widetilde{\nabla}^+_{\psi,\Gamma,-},\widetilde{\Delta}^\dagger_{\psi,\Gamma,-},\widetilde{\nabla}^\dagger_{\psi,\Gamma,-},\widetilde{\Phi}^r_{\psi,\Gamma,-},\widetilde{\Phi}^I_{\psi,\Gamma,-}, 
}
\]

\[
\xymatrix@R+0pc@C+0pc{
\breve{\Delta}_{\psi,\Gamma,-},\breve{\nabla}_{\psi,\Gamma,-},\breve{\Phi}_{\psi,\Gamma,-},\breve{\Delta}^+_{\psi,\Gamma,-},\breve{\nabla}^+_{\psi,\Gamma,-},\breve{\Delta}^\dagger_{\psi,\Gamma,-},\breve{\nabla}^\dagger_{\psi,\Gamma,-},\breve{\Phi}^r_{\psi,\Gamma,-},\breve{\Phi}^I_{\psi,\Gamma,-},	
}
\]

\[
\xymatrix@R+0pc@C+0pc{
{\Delta}_{\psi,\Gamma,-},{\nabla}_{\psi,\Gamma,-},{\Phi}_{\psi,\Gamma,-},{\Delta}^+_{\psi,\Gamma,-},{\nabla}^+_{\psi,\Gamma,-},{\Delta}^\dagger_{\psi,\Gamma,-},{\nabla}^\dagger_{\psi,\Gamma,-},{\Phi}^r_{\psi,\Gamma,-},{\Phi}^I_{\psi,\Gamma,-}.	
}
\]
\end{definition}

\begin{definition}
First we consider the Clausen-Scholze spectrum $\underset{\mathrm{Spec}}{\mathcal{O}}^\mathrm{CS}(*)$ attached to any of those in the above from \cite{10CS2} by taking derived rational localization:
\begin{align}
\underset{\mathrm{Spec}}{\mathcal{O}}^\mathrm{CS}\widetilde{\Delta}_{\psi,\Gamma,-},\underset{\mathrm{Spec}}{\mathcal{O}}^\mathrm{CS}\widetilde{\nabla}_{\psi,\Gamma,-},\underset{\mathrm{Spec}}{\mathcal{O}}^\mathrm{CS}\widetilde{\Phi}_{\psi,\Gamma,-},\underset{\mathrm{Spec}}{\mathcal{O}}^\mathrm{CS}\widetilde{\Delta}^+_{\psi,\Gamma,-},\underset{\mathrm{Spec}}{\mathcal{O}}^\mathrm{CS}\widetilde{\nabla}^+_{\psi,\Gamma,-},\\
\underset{\mathrm{Spec}}{\mathcal{O}}^\mathrm{CS}\widetilde{\Delta}^\dagger_{\psi,\Gamma,-},\underset{\mathrm{Spec}}{\mathcal{O}}^\mathrm{CS}\widetilde{\nabla}^\dagger_{\psi,\Gamma,-},\underset{\mathrm{Spec}}{\mathcal{O}}^\mathrm{CS}\widetilde{\Phi}^r_{\psi,\Gamma,-},\underset{\mathrm{Spec}}{\mathcal{O}}^\mathrm{CS}\widetilde{\Phi}^I_{\psi,\Gamma,-},	\\
\end{align}
\begin{align}
\underset{\mathrm{Spec}}{\mathcal{O}}^\mathrm{CS}\breve{\Delta}_{\psi,\Gamma,-},\breve{\nabla}_{\psi,\Gamma,-},\underset{\mathrm{Spec}}{\mathcal{O}}^\mathrm{CS}\breve{\Phi}_{\psi,\Gamma,-},\underset{\mathrm{Spec}}{\mathcal{O}}^\mathrm{CS}\breve{\Delta}^+_{\psi,\Gamma,-},\underset{\mathrm{Spec}}{\mathcal{O}}^\mathrm{CS}\breve{\nabla}^+_{\psi,\Gamma,-},\\
\underset{\mathrm{Spec}}{\mathcal{O}}^\mathrm{CS}\breve{\Delta}^\dagger_{\psi,\Gamma,-},\underset{\mathrm{Spec}}{\mathcal{O}}^\mathrm{CS}\breve{\nabla}^\dagger_{\psi,\Gamma,-},\underset{\mathrm{Spec}}{\mathcal{O}}^\mathrm{CS}\breve{\Phi}^r_{\psi,\Gamma,-},\breve{\Phi}^I_{\psi,\Gamma,-},	\\
\end{align}
\begin{align}
\underset{\mathrm{Spec}}{\mathcal{O}}^\mathrm{CS}{\Delta}_{\psi,\Gamma,-},\underset{\mathrm{Spec}}{\mathcal{O}}^\mathrm{CS}{\nabla}_{\psi,\Gamma,-},\underset{\mathrm{Spec}}{\mathcal{O}}^\mathrm{CS}{\Phi}_{\psi,\Gamma,-},\underset{\mathrm{Spec}}{\mathcal{O}}^\mathrm{CS}{\Delta}^+_{\psi,\Gamma,-},\underset{\mathrm{Spec}}{\mathcal{O}}^\mathrm{CS}{\nabla}^+_{\psi,\Gamma,-},\\
\underset{\mathrm{Spec}}{\mathcal{O}}^\mathrm{CS}{\Delta}^\dagger_{\psi,\Gamma,-},\underset{\mathrm{Spec}}{\mathcal{O}}^\mathrm{CS}{\nabla}^\dagger_{\psi,\Gamma,-},\underset{\mathrm{Spec}}{\mathcal{O}}^\mathrm{CS}{\Phi}^r_{\psi,\Gamma,-},\underset{\mathrm{Spec}}{\mathcal{O}}^\mathrm{CS}{\Phi}^I_{\psi,\Gamma,-}.	
\end{align}

Then we take the corresponding quotients by using the corresponding Frobenius operators:
\begin{align}
&\underset{\mathrm{Spec}}{\mathcal{O}}^\mathrm{CS}\widetilde{\Delta}_{\psi,\Gamma,-}/\mathrm{Fro}^\mathbb{Z},\underset{\mathrm{Spec}}{\mathcal{O}}^\mathrm{CS}\widetilde{\nabla}_{\psi,\Gamma,-}/\mathrm{Fro}^\mathbb{Z},\underset{\mathrm{Spec}}{\mathcal{O}}^\mathrm{CS}\widetilde{\Phi}_{\psi,\Gamma,-}/\mathrm{Fro}^\mathbb{Z},\underset{\mathrm{Spec}}{\mathcal{O}}^\mathrm{CS}\widetilde{\Delta}^+_{\psi,\Gamma,-}/\mathrm{Fro}^\mathbb{Z},\\
&\underset{\mathrm{Spec}}{\mathcal{O}}^\mathrm{CS}\widetilde{\nabla}^+_{\psi,\Gamma,-}/\mathrm{Fro}^\mathbb{Z}, \underset{\mathrm{Spec}}{\mathcal{O}}^\mathrm{CS}\widetilde{\Delta}^\dagger_{\psi,\Gamma,-}/\mathrm{Fro}^\mathbb{Z},\underset{\mathrm{Spec}}{\mathcal{O}}^\mathrm{CS}\widetilde{\nabla}^\dagger_{\psi,\Gamma,-}/\mathrm{Fro}^\mathbb{Z},	\\
\end{align}
\begin{align}
&\underset{\mathrm{Spec}}{\mathcal{O}}^\mathrm{CS}\breve{\Delta}_{\psi,\Gamma,-}/\mathrm{Fro}^\mathbb{Z},\breve{\nabla}_{\psi,\Gamma,-}/\mathrm{Fro}^\mathbb{Z},\underset{\mathrm{Spec}}{\mathcal{O}}^\mathrm{CS}\breve{\Phi}_{\psi,\Gamma,-}/\mathrm{Fro}^\mathbb{Z},\underset{\mathrm{Spec}}{\mathcal{O}}^\mathrm{CS}\breve{\Delta}^+_{\psi,\Gamma,-}/\mathrm{Fro}^\mathbb{Z},\\
&\underset{\mathrm{Spec}}{\mathcal{O}}^\mathrm{CS}\breve{\nabla}^+_{\psi,\Gamma,-}/\mathrm{Fro}^\mathbb{Z}, \underset{\mathrm{Spec}}{\mathcal{O}}^\mathrm{CS}\breve{\Delta}^\dagger_{\psi,\Gamma,-}/\mathrm{Fro}^\mathbb{Z},\underset{\mathrm{Spec}}{\mathcal{O}}^\mathrm{CS}\breve{\nabla}^\dagger_{\psi,\Gamma,-}/\mathrm{Fro}^\mathbb{Z},	\\
\end{align}
\begin{align}
&\underset{\mathrm{Spec}}{\mathcal{O}}^\mathrm{CS}{\Delta}_{\psi,\Gamma,-}/\mathrm{Fro}^\mathbb{Z},\underset{\mathrm{Spec}}{\mathcal{O}}^\mathrm{CS}{\nabla}_{\psi,\Gamma,-}/\mathrm{Fro}^\mathbb{Z},\underset{\mathrm{Spec}}{\mathcal{O}}^\mathrm{CS}{\Phi}_{\psi,\Gamma,-}/\mathrm{Fro}^\mathbb{Z},\underset{\mathrm{Spec}}{\mathcal{O}}^\mathrm{CS}{\Delta}^+_{\psi,\Gamma,-}/\mathrm{Fro}^\mathbb{Z},\\
&\underset{\mathrm{Spec}}{\mathcal{O}}^\mathrm{CS}{\nabla}^+_{\psi,\Gamma,-}/\mathrm{Fro}^\mathbb{Z}, \underset{\mathrm{Spec}}{\mathcal{O}}^\mathrm{CS}{\Delta}^\dagger_{\psi,\Gamma,-}/\mathrm{Fro}^\mathbb{Z},\underset{\mathrm{Spec}}{\mathcal{O}}^\mathrm{CS}{\nabla}^\dagger_{\psi,\Gamma,-}/\mathrm{Fro}^\mathbb{Z}.	
\end{align}
Here for those space with notations related to the radius and the corresponding interval we consider the total unions $\bigcap_r,\bigcup_I$ in order to achieve the whole spaces to achieve the analogues of the corresponding FF curves from \cite{10KL1}, \cite{10KL2}, \cite{10FF} for
\[
\xymatrix@R+0pc@C+0pc{
\underset{r}{\mathrm{homotopycolimit}}~\underset{\mathrm{Spec}}{\mathcal{O}}^\mathrm{CS}\widetilde{\Phi}^r_{\psi,\Gamma,-},\underset{I}{\mathrm{homotopylimit}}~\underset{\mathrm{Spec}}{\mathcal{O}}^\mathrm{CS}\widetilde{\Phi}^I_{\psi,\Gamma,-},	\\
}
\]
\[
\xymatrix@R+0pc@C+0pc{
\underset{r}{\mathrm{homotopycolimit}}~\underset{\mathrm{Spec}}{\mathcal{O}}^\mathrm{CS}\breve{\Phi}^r_{\psi,\Gamma,-},\underset{I}{\mathrm{homotopylimit}}~\underset{\mathrm{Spec}}{\mathcal{O}}^\mathrm{CS}\breve{\Phi}^I_{\psi,\Gamma,-},	\\
}
\]
\[
\xymatrix@R+0pc@C+0pc{
\underset{r}{\mathrm{homotopycolimit}}~\underset{\mathrm{Spec}}{\mathcal{O}}^\mathrm{CS}{\Phi}^r_{\psi,\Gamma,-},\underset{I}{\mathrm{homotopylimit}}~\underset{\mathrm{Spec}}{\mathcal{O}}^\mathrm{CS}{\Phi}^I_{\psi,\Gamma,-}.	
}
\]
\[ 
\xymatrix@R+0pc@C+0pc{
\underset{r}{\mathrm{homotopycolimit}}~\underset{\mathrm{Spec}}{\mathcal{O}}^\mathrm{CS}\widetilde{\Phi}^r_{\psi,\Gamma,-}/\mathrm{Fro}^\mathbb{Z},\underset{I}{\mathrm{homotopylimit}}~\underset{\mathrm{Spec}}{\mathcal{O}}^\mathrm{CS}\widetilde{\Phi}^I_{\psi,\Gamma,-}/\mathrm{Fro}^\mathbb{Z},	\\
}
\]
\[ 
\xymatrix@R+0pc@C+0pc{
\underset{r}{\mathrm{homotopycolimit}}~\underset{\mathrm{Spec}}{\mathcal{O}}^\mathrm{CS}\breve{\Phi}^r_{\psi,\Gamma,-}/\mathrm{Fro}^\mathbb{Z},\underset{I}{\mathrm{homotopylimit}}~\breve{\Phi}^I_{\psi,\Gamma,-}/\mathrm{Fro}^\mathbb{Z},	\\
}
\]
\[ 
\xymatrix@R+0pc@C+0pc{
\underset{r}{\mathrm{homotopycolimit}}~\underset{\mathrm{Spec}}{\mathcal{O}}^\mathrm{CS}{\Phi}^r_{\psi,\Gamma,-}/\mathrm{Fro}^\mathbb{Z},\underset{I}{\mathrm{homotopylimit}}~\underset{\mathrm{Spec}}{\mathcal{O}}^\mathrm{CS}{\Phi}^I_{\psi,\Gamma,-}/\mathrm{Fro}^\mathbb{Z}.	
}
\]

\end{definition}

\

\begin{definition}
We then consider the corresponding quasipresheaves of the corresponding ind-Banach or monomorphic ind-Banach modules from \cite{10BBK}, \cite{10KKM}:
\begin{align}
\mathrm{Quasicoherentpresheaves,IndBanach}_{*}	
\end{align}
where $*$ is one of the following spaces:
\begin{align}
&\underset{\mathrm{Spec}}{\mathcal{O}}^\mathrm{BK}\widetilde{\Phi}_{\psi,\Gamma,-}/\mathrm{Fro}^\mathbb{Z},	\\
\end{align}
\begin{align}
&\underset{\mathrm{Spec}}{\mathcal{O}}^\mathrm{BK}\breve{\Phi}_{\psi,\Gamma,-}/\mathrm{Fro}^\mathbb{Z},	\\
\end{align}
\begin{align}
&\underset{\mathrm{Spec}}{\mathcal{O}}^\mathrm{BK}{\Phi}_{\psi,\Gamma,-}/\mathrm{Fro}^\mathbb{Z}.	
\end{align}
Here for those space without notation related to the radius and the corresponding interval we consider the total unions $\bigcap_r,\bigcup_I$ in order to achieve the whole spaces to achieve the analogues of the corresponding FF curves from \cite{10KL1}, \cite{10KL2}, \cite{10FF} for
\[
\xymatrix@R+0pc@C+0pc{
\underset{r}{\mathrm{homotopycolimit}}~\underset{\mathrm{Spec}}{\mathcal{O}}^\mathrm{BK}\widetilde{\Phi}^r_{\psi,\Gamma,-},\underset{I}{\mathrm{homotopylimit}}~\underset{\mathrm{Spec}}{\mathcal{O}}^\mathrm{BK}\widetilde{\Phi}^I_{\psi,\Gamma,-},	\\
}
\]
\[
\xymatrix@R+0pc@C+0pc{
\underset{r}{\mathrm{homotopycolimit}}~\underset{\mathrm{Spec}}{\mathcal{O}}^\mathrm{BK}\breve{\Phi}^r_{\psi,\Gamma,-},\underset{I}{\mathrm{homotopylimit}}~\underset{\mathrm{Spec}}{\mathcal{O}}^\mathrm{BK}\breve{\Phi}^I_{\psi,\Gamma,-},	\\
}
\]
\[
\xymatrix@R+0pc@C+0pc{
\underset{r}{\mathrm{homotopycolimit}}~\underset{\mathrm{Spec}}{\mathcal{O}}^\mathrm{BK}{\Phi}^r_{\psi,\Gamma,-},\underset{I}{\mathrm{homotopylimit}}~\underset{\mathrm{Spec}}{\mathcal{O}}^\mathrm{BK}{\Phi}^I_{\psi,\Gamma,-}.	
}
\]
\[  
\xymatrix@R+0pc@C+0pc{
\underset{r}{\mathrm{homotopycolimit}}~\underset{\mathrm{Spec}}{\mathcal{O}}^\mathrm{BK}\widetilde{\Phi}^r_{\psi,\Gamma,-}/\mathrm{Fro}^\mathbb{Z},\underset{I}{\mathrm{homotopylimit}}~\underset{\mathrm{Spec}}{\mathcal{O}}^\mathrm{BK}\widetilde{\Phi}^I_{\psi,\Gamma,-}/\mathrm{Fro}^\mathbb{Z},	\\
}
\]
\[ 
\xymatrix@R+0pc@C+0pc{
\underset{r}{\mathrm{homotopycolimit}}~\underset{\mathrm{Spec}}{\mathcal{O}}^\mathrm{BK}\breve{\Phi}^r_{\psi,\Gamma,-}/\mathrm{Fro}^\mathbb{Z},\underset{I}{\mathrm{homotopylimit}}~\underset{\mathrm{Spec}}{\mathcal{O}}^\mathrm{BK}\breve{\Phi}^I_{\psi,\Gamma,-}/\mathrm{Fro}^\mathbb{Z},	\\
}
\]
\[ 
\xymatrix@R+0pc@C+0pc{
\underset{r}{\mathrm{homotopycolimit}}~\underset{\mathrm{Spec}}{\mathcal{O}}^\mathrm{BK}{\Phi}^r_{\psi,\Gamma,-}/\mathrm{Fro}^\mathbb{Z},\underset{I}{\mathrm{homotopylimit}}~\underset{\mathrm{Spec}}{\mathcal{O}}^\mathrm{BK}{\Phi}^I_{\psi,\Gamma,-}/\mathrm{Fro}^\mathbb{Z}.	
}
\]

\end{definition}

\begin{definition}
We then consider the corresponding quasisheaves of the corresponding condensed solid topological modules from \cite{10CS2}:
\begin{align}
\mathrm{Quasicoherentsheaves, Condensed}_{*}	
\end{align}
where $*$ is one of the following spaces:
\begin{align}
&\underset{\mathrm{Spec}}{\mathcal{O}}^\mathrm{CS}\widetilde{\Delta}_{\psi,\Gamma,-}/\mathrm{Fro}^\mathbb{Z},\underset{\mathrm{Spec}}{\mathcal{O}}^\mathrm{CS}\widetilde{\nabla}_{\psi,\Gamma,-}/\mathrm{Fro}^\mathbb{Z},\underset{\mathrm{Spec}}{\mathcal{O}}^\mathrm{CS}\widetilde{\Phi}_{\psi,\Gamma,-}/\mathrm{Fro}^\mathbb{Z},\underset{\mathrm{Spec}}{\mathcal{O}}^\mathrm{CS}\widetilde{\Delta}^+_{\psi,\Gamma,-}/\mathrm{Fro}^\mathbb{Z},\\
&\underset{\mathrm{Spec}}{\mathcal{O}}^\mathrm{CS}\widetilde{\nabla}^+_{\psi,\Gamma,-}/\mathrm{Fro}^\mathbb{Z},\underset{\mathrm{Spec}}{\mathcal{O}}^\mathrm{CS}\widetilde{\Delta}^\dagger_{\psi,\Gamma,-}/\mathrm{Fro}^\mathbb{Z},\underset{\mathrm{Spec}}{\mathcal{O}}^\mathrm{CS}\widetilde{\nabla}^\dagger_{\psi,\Gamma,-}/\mathrm{Fro}^\mathbb{Z},	\\
\end{align}
\begin{align}
&\underset{\mathrm{Spec}}{\mathcal{O}}^\mathrm{CS}\breve{\Delta}_{\psi,\Gamma,-}/\mathrm{Fro}^\mathbb{Z},\breve{\nabla}_{\psi,\Gamma,-}/\mathrm{Fro}^\mathbb{Z},\underset{\mathrm{Spec}}{\mathcal{O}}^\mathrm{CS}\breve{\Phi}_{\psi,\Gamma,-}/\mathrm{Fro}^\mathbb{Z},\underset{\mathrm{Spec}}{\mathcal{O}}^\mathrm{CS}\breve{\Delta}^+_{\psi,\Gamma,-}/\mathrm{Fro}^\mathbb{Z},\\
&\underset{\mathrm{Spec}}{\mathcal{O}}^\mathrm{CS}\breve{\nabla}^+_{\psi,\Gamma,-}/\mathrm{Fro}^\mathbb{Z},\underset{\mathrm{Spec}}{\mathcal{O}}^\mathrm{CS}\breve{\Delta}^\dagger_{\psi,\Gamma,-}/\mathrm{Fro}^\mathbb{Z},\underset{\mathrm{Spec}}{\mathcal{O}}^\mathrm{CS}\breve{\nabla}^\dagger_{\psi,\Gamma,-}/\mathrm{Fro}^\mathbb{Z},	\\
\end{align}
\begin{align}
&\underset{\mathrm{Spec}}{\mathcal{O}}^\mathrm{CS}{\Delta}_{\psi,\Gamma,-}/\mathrm{Fro}^\mathbb{Z},\underset{\mathrm{Spec}}{\mathcal{O}}^\mathrm{CS}{\nabla}_{\psi,\Gamma,-}/\mathrm{Fro}^\mathbb{Z},\underset{\mathrm{Spec}}{\mathcal{O}}^\mathrm{CS}{\Phi}_{\psi,\Gamma,-}/\mathrm{Fro}^\mathbb{Z},\underset{\mathrm{Spec}}{\mathcal{O}}^\mathrm{CS}{\Delta}^+_{\psi,\Gamma,-}/\mathrm{Fro}^\mathbb{Z},\\
&\underset{\mathrm{Spec}}{\mathcal{O}}^\mathrm{CS}{\nabla}^+_{\psi,\Gamma,-}/\mathrm{Fro}^\mathbb{Z}, \underset{\mathrm{Spec}}{\mathcal{O}}^\mathrm{CS}{\Delta}^\dagger_{\psi,\Gamma,-}/\mathrm{Fro}^\mathbb{Z},\underset{\mathrm{Spec}}{\mathcal{O}}^\mathrm{CS}{\nabla}^\dagger_{\psi,\Gamma,-}/\mathrm{Fro}^\mathbb{Z}.	
\end{align}
Here for those space with notations related to the radius and the corresponding interval we consider the total unions $\bigcap_r,\bigcup_I$ in order to achieve the whole spaces to achieve the analogues of the corresponding FF curves from \cite{10KL1}, \cite{10KL2}, \cite{10FF} for
\[
\xymatrix@R+0pc@C+0pc{
\underset{r}{\mathrm{homotopycolimit}}~\underset{\mathrm{Spec}}{\mathcal{O}}^\mathrm{CS}\widetilde{\Phi}^r_{\psi,\Gamma,-},\underset{I}{\mathrm{homotopylimit}}~\underset{\mathrm{Spec}}{\mathcal{O}}^\mathrm{CS}\widetilde{\Phi}^I_{\psi,\Gamma,-},	\\
}
\]
\[
\xymatrix@R+0pc@C+0pc{
\underset{r}{\mathrm{homotopycolimit}}~\underset{\mathrm{Spec}}{\mathcal{O}}^\mathrm{CS}\breve{\Phi}^r_{\psi,\Gamma,-},\underset{I}{\mathrm{homotopylimit}}~\underset{\mathrm{Spec}}{\mathcal{O}}^\mathrm{CS}\breve{\Phi}^I_{\psi,\Gamma,-},	\\
}
\]
\[
\xymatrix@R+0pc@C+0pc{
\underset{r}{\mathrm{homotopycolimit}}~\underset{\mathrm{Spec}}{\mathcal{O}}^\mathrm{CS}{\Phi}^r_{\psi,\Gamma,-},\underset{I}{\mathrm{homotopylimit}}~\underset{\mathrm{Spec}}{\mathcal{O}}^\mathrm{CS}{\Phi}^I_{\psi,\Gamma,-}.	
}
\]
\[ 
\xymatrix@R+0pc@C+0pc{
\underset{r}{\mathrm{homotopycolimit}}~\underset{\mathrm{Spec}}{\mathcal{O}}^\mathrm{CS}\widetilde{\Phi}^r_{\psi,\Gamma,-}/\mathrm{Fro}^\mathbb{Z},\underset{I}{\mathrm{homotopylimit}}~\underset{\mathrm{Spec}}{\mathcal{O}}^\mathrm{CS}\widetilde{\Phi}^I_{\psi,\Gamma,-}/\mathrm{Fro}^\mathbb{Z},	\\
}
\]
\[ 
\xymatrix@R+0pc@C+0pc{
\underset{r}{\mathrm{homotopycolimit}}~\underset{\mathrm{Spec}}{\mathcal{O}}^\mathrm{CS}\breve{\Phi}^r_{\psi,\Gamma,-}/\mathrm{Fro}^\mathbb{Z},\underset{I}{\mathrm{homotopylimit}}~\breve{\Phi}^I_{\psi,\Gamma,-}/\mathrm{Fro}^\mathbb{Z},	\\
}
\]
\[ 
\xymatrix@R+0pc@C+0pc{
\underset{r}{\mathrm{homotopycolimit}}~\underset{\mathrm{Spec}}{\mathcal{O}}^\mathrm{CS}{\Phi}^r_{\psi,\Gamma,-}/\mathrm{Fro}^\mathbb{Z},\underset{I}{\mathrm{homotopylimit}}~\underset{\mathrm{Spec}}{\mathcal{O}}^\mathrm{CS}{\Phi}^I_{\psi,\Gamma,-}/\mathrm{Fro}^\mathbb{Z}.	
}
\]

\end{definition}

\

\begin{proposition}
There is a well-defined functor from the $\infty$-category 
\begin{align}
\mathrm{Quasicoherentpresheaves,Condensed}_{*}	
\end{align}
where $*$ is one of the following spaces:
\begin{align}
&\underset{\mathrm{Spec}}{\mathcal{O}}^\mathrm{CS}\widetilde{\Phi}_{\psi,\Gamma,-}/\mathrm{Fro}^\mathbb{Z},	\\
\end{align}
\begin{align}
&\underset{\mathrm{Spec}}{\mathcal{O}}^\mathrm{CS}\breve{\Phi}_{\psi,\Gamma,-}/\mathrm{Fro}^\mathbb{Z},	\\
\end{align}
\begin{align}
&\underset{\mathrm{Spec}}{\mathcal{O}}^\mathrm{CS}{\Phi}_{\psi,\Gamma,-}/\mathrm{Fro}^\mathbb{Z},	
\end{align}
to the $\infty$-category of $\mathrm{Fro}$-equivariant quasicoherent presheaves over similar spaces above correspondingly without the $\mathrm{Fro}$-quotients, and to the $\infty$-category of $\mathrm{Fro}$-equivariant quasicoherent modules over global sections of the structure $\infty$-sheaves of the similar spaces above correspondingly without the $\mathrm{Fro}$-quotients. Here for those space without notation related to the radius and the corresponding interval we consider the total unions $\bigcap_r,\bigcup_I$ in order to achieve the whole spaces to achieve the analogues of the corresponding FF curves from \cite{10KL1}, \cite{10KL2}, \cite{10FF} for
\[
\xymatrix@R+0pc@C+0pc{
\underset{r}{\mathrm{homotopycolimit}}~\underset{\mathrm{Spec}}{\mathcal{O}}^\mathrm{CS}\widetilde{\Phi}^r_{\psi,\Gamma,-},\underset{I}{\mathrm{homotopylimit}}~\underset{\mathrm{Spec}}{\mathcal{O}}^\mathrm{CS}\widetilde{\Phi}^I_{\psi,\Gamma,-},	\\
}
\]
\[
\xymatrix@R+0pc@C+0pc{
\underset{r}{\mathrm{homotopycolimit}}~\underset{\mathrm{Spec}}{\mathcal{O}}^\mathrm{CS}\breve{\Phi}^r_{\psi,\Gamma,-},\underset{I}{\mathrm{homotopylimit}}~\underset{\mathrm{Spec}}{\mathcal{O}}^\mathrm{CS}\breve{\Phi}^I_{\psi,\Gamma,-},	\\
}
\]
\[
\xymatrix@R+0pc@C+0pc{
\underset{r}{\mathrm{homotopycolimit}}~\underset{\mathrm{Spec}}{\mathcal{O}}^\mathrm{CS}{\Phi}^r_{\psi,\Gamma,-},\underset{I}{\mathrm{homotopylimit}}~\underset{\mathrm{Spec}}{\mathcal{O}}^\mathrm{CS}{\Phi}^I_{\psi,\Gamma,-}.	
}
\]
\[ 
\xymatrix@R+0pc@C+0pc{
\underset{r}{\mathrm{homotopycolimit}}~\underset{\mathrm{Spec}}{\mathcal{O}}^\mathrm{CS}\widetilde{\Phi}^r_{\psi,\Gamma,-}/\mathrm{Fro}^\mathbb{Z},\underset{I}{\mathrm{homotopylimit}}~\underset{\mathrm{Spec}}{\mathcal{O}}^\mathrm{CS}\widetilde{\Phi}^I_{\psi,\Gamma,-}/\mathrm{Fro}^\mathbb{Z},	\\
}
\]
\[ 
\xymatrix@R+0pc@C+0pc{
\underset{r}{\mathrm{homotopycolimit}}~\underset{\mathrm{Spec}}{\mathcal{O}}^\mathrm{CS}\breve{\Phi}^r_{\psi,\Gamma,-}/\mathrm{Fro}^\mathbb{Z},\underset{I}{\mathrm{homotopylimit}}~\breve{\Phi}^I_{\psi,\Gamma,-}/\mathrm{Fro}^\mathbb{Z},	\\
}
\]
\[ 
\xymatrix@R+0pc@C+0pc{
\underset{r}{\mathrm{homotopycolimit}}~\underset{\mathrm{Spec}}{\mathcal{O}}^\mathrm{CS}{\Phi}^r_{\psi,\Gamma,-}/\mathrm{Fro}^\mathbb{Z},\underset{I}{\mathrm{homotopylimit}}~\underset{\mathrm{Spec}}{\mathcal{O}}^\mathrm{CS}{\Phi}^I_{\psi,\Gamma,-}/\mathrm{Fro}^\mathbb{Z}.	
}
\]	
In this situation we will have the target category being family parametrized by $r$ or $I$ in compatible glueing sense as in \cite[Definition 5.4.10]{10KL2}. In this situation for modules parametrized by the intervals we have the equivalence of $\infty$-categories by using \cite[Proposition 13.8]{10CS2}. Here the corresponding quasicoherent Frobenius modules are defined to be the corresponding homotopy colimits and limits of Frobenius modules:
\begin{align}
\underset{r}{\mathrm{homotopycolimit}}~M_r,\\
\underset{I}{\mathrm{homotopylimit}}~M_I,	
\end{align}
where each $M_r$ is a Frobenius-equivariant module over the period ring with respect to some radius $r$ while each $M_I$ is a Frobenius-equivariant module over the period ring with respect to some interval $I$.\\
\end{proposition}

\begin{proposition}
Similar proposition holds for 
\begin{align}
\mathrm{Quasicoherentsheaves,IndBanach}_{*}.	
\end{align}	
\end{proposition}

\

\begin{definition}
We then consider the corresponding quasipresheaves of perfect complexes the corresponding ind-Banach or monomorphic ind-Banach modules from \cite{10BBK}, \cite{10KKM}:
\begin{align}
\mathrm{Quasicoherentpresheaves,Perfectcomplex,IndBanach}_{*}	
\end{align}
where $*$ is one of the following spaces:
\begin{align}
&\underset{\mathrm{Spec}}{\mathcal{O}}^\mathrm{BK}\widetilde{\Phi}_{\psi,\Gamma,-}/\mathrm{Fro}^\mathbb{Z},	\\
\end{align}
\begin{align}
&\underset{\mathrm{Spec}}{\mathcal{O}}^\mathrm{BK}\breve{\Phi}_{\psi,\Gamma,-}/\mathrm{Fro}^\mathbb{Z},	\\
\end{align}
\begin{align}
&\underset{\mathrm{Spec}}{\mathcal{O}}^\mathrm{BK}{\Phi}_{\psi,\Gamma,-}/\mathrm{Fro}^\mathbb{Z}.	
\end{align}
Here for those space without notation related to the radius and the corresponding interval we consider the total unions $\bigcap_r,\bigcup_I$ in order to achieve the whole spaces to achieve the analogues of the corresponding FF curves from \cite{10KL1}, \cite{10KL2}, \cite{10FF} for
\[
\xymatrix@R+0pc@C+0pc{
\underset{r}{\mathrm{homotopycolimit}}~\underset{\mathrm{Spec}}{\mathcal{O}}^\mathrm{BK}\widetilde{\Phi}^r_{\psi,\Gamma,-},\underset{I}{\mathrm{homotopylimit}}~\underset{\mathrm{Spec}}{\mathcal{O}}^\mathrm{BK}\widetilde{\Phi}^I_{\psi,\Gamma,-},	\\
}
\]
\[
\xymatrix@R+0pc@C+0pc{
\underset{r}{\mathrm{homotopycolimit}}~\underset{\mathrm{Spec}}{\mathcal{O}}^\mathrm{BK}\breve{\Phi}^r_{\psi,\Gamma,-},\underset{I}{\mathrm{homotopylimit}}~\underset{\mathrm{Spec}}{\mathcal{O}}^\mathrm{BK}\breve{\Phi}^I_{\psi,\Gamma,-},	\\
}
\]
\[
\xymatrix@R+0pc@C+0pc{
\underset{r}{\mathrm{homotopycolimit}}~\underset{\mathrm{Spec}}{\mathcal{O}}^\mathrm{BK}{\Phi}^r_{\psi,\Gamma,-},\underset{I}{\mathrm{homotopylimit}}~\underset{\mathrm{Spec}}{\mathcal{O}}^\mathrm{BK}{\Phi}^I_{\psi,\Gamma,-}.	
}
\]
\[  
\xymatrix@R+0pc@C+0pc{
\underset{r}{\mathrm{homotopycolimit}}~\underset{\mathrm{Spec}}{\mathcal{O}}^\mathrm{BK}\widetilde{\Phi}^r_{\psi,\Gamma,-}/\mathrm{Fro}^\mathbb{Z},\underset{I}{\mathrm{homotopylimit}}~\underset{\mathrm{Spec}}{\mathcal{O}}^\mathrm{BK}\widetilde{\Phi}^I_{\psi,\Gamma,-}/\mathrm{Fro}^\mathbb{Z},	\\
}
\]
\[ 
\xymatrix@R+0pc@C+0pc{
\underset{r}{\mathrm{homotopycolimit}}~\underset{\mathrm{Spec}}{\mathcal{O}}^\mathrm{BK}\breve{\Phi}^r_{\psi,\Gamma,-}/\mathrm{Fro}^\mathbb{Z},\underset{I}{\mathrm{homotopylimit}}~\underset{\mathrm{Spec}}{\mathcal{O}}^\mathrm{BK}\breve{\Phi}^I_{\psi,\Gamma,-}/\mathrm{Fro}^\mathbb{Z},	\\
}
\]
\[ 
\xymatrix@R+0pc@C+0pc{
\underset{r}{\mathrm{homotopycolimit}}~\underset{\mathrm{Spec}}{\mathcal{O}}^\mathrm{BK}{\Phi}^r_{\psi,\Gamma,-}/\mathrm{Fro}^\mathbb{Z},\underset{I}{\mathrm{homotopylimit}}~\underset{\mathrm{Spec}}{\mathcal{O}}^\mathrm{BK}{\Phi}^I_{\psi,\Gamma,-}/\mathrm{Fro}^\mathbb{Z}.	
}
\]

\end{definition}

\begin{definition}
We then consider the corresponding quasisheaves of perfect complexes of the corresponding condensed solid topological modules from \cite{10CS2}:
\begin{align}
\mathrm{Quasicoherentsheaves, Perfectcomplex, Condensed}_{*}	
\end{align}
where $*$ is one of the following spaces:
\begin{align}
&\underset{\mathrm{Spec}}{\mathcal{O}}^\mathrm{CS}\widetilde{\Delta}_{\psi,\Gamma,-}/\mathrm{Fro}^\mathbb{Z},\underset{\mathrm{Spec}}{\mathcal{O}}^\mathrm{CS}\widetilde{\nabla}_{\psi,\Gamma,-}/\mathrm{Fro}^\mathbb{Z},\underset{\mathrm{Spec}}{\mathcal{O}}^\mathrm{CS}\widetilde{\Phi}_{\psi,\Gamma,-}/\mathrm{Fro}^\mathbb{Z},\underset{\mathrm{Spec}}{\mathcal{O}}^\mathrm{CS}\widetilde{\Delta}^+_{\psi,\Gamma,-}/\mathrm{Fro}^\mathbb{Z},\\
&\underset{\mathrm{Spec}}{\mathcal{O}}^\mathrm{CS}\widetilde{\nabla}^+_{\psi,\Gamma,-}/\mathrm{Fro}^\mathbb{Z},\underset{\mathrm{Spec}}{\mathcal{O}}^\mathrm{CS}\widetilde{\Delta}^\dagger_{\psi,\Gamma,-}/\mathrm{Fro}^\mathbb{Z},\underset{\mathrm{Spec}}{\mathcal{O}}^\mathrm{CS}\widetilde{\nabla}^\dagger_{\psi,\Gamma,-}/\mathrm{Fro}^\mathbb{Z},	\\
\end{align}
\begin{align}
&\underset{\mathrm{Spec}}{\mathcal{O}}^\mathrm{CS}\breve{\Delta}_{\psi,\Gamma,-}/\mathrm{Fro}^\mathbb{Z},\breve{\nabla}_{\psi,\Gamma,-}/\mathrm{Fro}^\mathbb{Z},\underset{\mathrm{Spec}}{\mathcal{O}}^\mathrm{CS}\breve{\Phi}_{\psi,\Gamma,-}/\mathrm{Fro}^\mathbb{Z},\underset{\mathrm{Spec}}{\mathcal{O}}^\mathrm{CS}\breve{\Delta}^+_{\psi,\Gamma,-}/\mathrm{Fro}^\mathbb{Z},\\
&\underset{\mathrm{Spec}}{\mathcal{O}}^\mathrm{CS}\breve{\nabla}^+_{\psi,\Gamma,-}/\mathrm{Fro}^\mathbb{Z},\underset{\mathrm{Spec}}{\mathcal{O}}^\mathrm{CS}\breve{\Delta}^\dagger_{\psi,\Gamma,-}/\mathrm{Fro}^\mathbb{Z},\underset{\mathrm{Spec}}{\mathcal{O}}^\mathrm{CS}\breve{\nabla}^\dagger_{\psi,\Gamma,-}/\mathrm{Fro}^\mathbb{Z},	\\
\end{align}
\begin{align}
&\underset{\mathrm{Spec}}{\mathcal{O}}^\mathrm{CS}{\Delta}_{\psi,\Gamma,-}/\mathrm{Fro}^\mathbb{Z},\underset{\mathrm{Spec}}{\mathcal{O}}^\mathrm{CS}{\nabla}_{\psi,\Gamma,-}/\mathrm{Fro}^\mathbb{Z},\underset{\mathrm{Spec}}{\mathcal{O}}^\mathrm{CS}{\Phi}_{\psi,\Gamma,-}/\mathrm{Fro}^\mathbb{Z},\underset{\mathrm{Spec}}{\mathcal{O}}^\mathrm{CS}{\Delta}^+_{\psi,\Gamma,-}/\mathrm{Fro}^\mathbb{Z},\\
&\underset{\mathrm{Spec}}{\mathcal{O}}^\mathrm{CS}{\nabla}^+_{\psi,\Gamma,-}/\mathrm{Fro}^\mathbb{Z}, \underset{\mathrm{Spec}}{\mathcal{O}}^\mathrm{CS}{\Delta}^\dagger_{\psi,\Gamma,-}/\mathrm{Fro}^\mathbb{Z},\underset{\mathrm{Spec}}{\mathcal{O}}^\mathrm{CS}{\nabla}^\dagger_{\psi,\Gamma,-}/\mathrm{Fro}^\mathbb{Z}.	
\end{align}
Here for those space with notations related to the radius and the corresponding interval we consider the total unions $\bigcap_r,\bigcup_I$ in order to achieve the whole spaces to achieve the analogues of the corresponding FF curves from \cite{10KL1}, \cite{10KL2}, \cite{10FF} for
\[
\xymatrix@R+0pc@C+0pc{
\underset{r}{\mathrm{homotopycolimit}}~\underset{\mathrm{Spec}}{\mathcal{O}}^\mathrm{CS}\widetilde{\Phi}^r_{\psi,\Gamma,-},\underset{I}{\mathrm{homotopylimit}}~\underset{\mathrm{Spec}}{\mathcal{O}}^\mathrm{CS}\widetilde{\Phi}^I_{\psi,\Gamma,-},	\\
}
\]
\[
\xymatrix@R+0pc@C+0pc{
\underset{r}{\mathrm{homotopycolimit}}~\underset{\mathrm{Spec}}{\mathcal{O}}^\mathrm{CS}\breve{\Phi}^r_{\psi,\Gamma,-},\underset{I}{\mathrm{homotopylimit}}~\underset{\mathrm{Spec}}{\mathcal{O}}^\mathrm{CS}\breve{\Phi}^I_{\psi,\Gamma,-},	\\
}
\]
\[
\xymatrix@R+0pc@C+0pc{
\underset{r}{\mathrm{homotopycolimit}}~\underset{\mathrm{Spec}}{\mathcal{O}}^\mathrm{CS}{\Phi}^r_{\psi,\Gamma,-},\underset{I}{\mathrm{homotopylimit}}~\underset{\mathrm{Spec}}{\mathcal{O}}^\mathrm{CS}{\Phi}^I_{\psi,\Gamma,-}.	
}
\]
\[ 
\xymatrix@R+0pc@C+0pc{
\underset{r}{\mathrm{homotopycolimit}}~\underset{\mathrm{Spec}}{\mathcal{O}}^\mathrm{CS}\widetilde{\Phi}^r_{\psi,\Gamma,-}/\mathrm{Fro}^\mathbb{Z},\underset{I}{\mathrm{homotopylimit}}~\underset{\mathrm{Spec}}{\mathcal{O}}^\mathrm{CS}\widetilde{\Phi}^I_{\psi,\Gamma,-}/\mathrm{Fro}^\mathbb{Z},	\\
}
\]
\[ 
\xymatrix@R+0pc@C+0pc{
\underset{r}{\mathrm{homotopycolimit}}~\underset{\mathrm{Spec}}{\mathcal{O}}^\mathrm{CS}\breve{\Phi}^r_{\psi,\Gamma,-}/\mathrm{Fro}^\mathbb{Z},\underset{I}{\mathrm{homotopylimit}}~\breve{\Phi}^I_{\psi,\Gamma,-}/\mathrm{Fro}^\mathbb{Z},	\\
}
\]
\[ 
\xymatrix@R+0pc@C+0pc{
\underset{r}{\mathrm{homotopycolimit}}~\underset{\mathrm{Spec}}{\mathcal{O}}^\mathrm{CS}{\Phi}^r_{\psi,\Gamma,-}/\mathrm{Fro}^\mathbb{Z},\underset{I}{\mathrm{homotopylimit}}~\underset{\mathrm{Spec}}{\mathcal{O}}^\mathrm{CS}{\Phi}^I_{\psi,\Gamma,-}/\mathrm{Fro}^\mathbb{Z}.	
}
\]

\end{definition}

\begin{proposition}
There is a well-defined functor from the $\infty$-category 
\begin{align}
\mathrm{Quasicoherentpresheaves,Perfectcomplex,Condensed}_{*}	
\end{align}
where $*$ is one of the following spaces:
\begin{align}
&\underset{\mathrm{Spec}}{\mathcal{O}}^\mathrm{CS}\widetilde{\Phi}_{\psi,\Gamma,-}/\mathrm{Fro}^\mathbb{Z},	\\
\end{align}
\begin{align}
&\underset{\mathrm{Spec}}{\mathcal{O}}^\mathrm{CS}\breve{\Phi}_{\psi,\Gamma,-}/\mathrm{Fro}^\mathbb{Z},	\\
\end{align}
\begin{align}
&\underset{\mathrm{Spec}}{\mathcal{O}}^\mathrm{CS}{\Phi}_{\psi,\Gamma,-}/\mathrm{Fro}^\mathbb{Z},	
\end{align}
to the $\infty$-category of $\mathrm{Fro}$-equivariant quasicoherent presheaves over similar spaces above correspondingly without the $\mathrm{Fro}$-quotients, and to the $\infty$-category of $\mathrm{Fro}$-equivariant quasicoherent modules over global sections of the structure $\infty$-sheaves of the similar spaces above correspondingly without the $\mathrm{Fro}$-quotients. Here for those space without notation related to the radius and the corresponding interval we consider the total unions $\bigcap_r,\bigcup_I$ in order to achieve the whole spaces to achieve the analogues of the corresponding FF curves from \cite{10KL1}, \cite{10KL2}, \cite{10FF} for
\[
\xymatrix@R+0pc@C+0pc{
\underset{r}{\mathrm{homotopycolimit}}~\underset{\mathrm{Spec}}{\mathcal{O}}^\mathrm{CS}\widetilde{\Phi}^r_{\psi,\Gamma,-},\underset{I}{\mathrm{homotopylimit}}~\underset{\mathrm{Spec}}{\mathcal{O}}^\mathrm{CS}\widetilde{\Phi}^I_{\psi,\Gamma,-},	\\
}
\]
\[
\xymatrix@R+0pc@C+0pc{
\underset{r}{\mathrm{homotopycolimit}}~\underset{\mathrm{Spec}}{\mathcal{O}}^\mathrm{CS}\breve{\Phi}^r_{\psi,\Gamma,-},\underset{I}{\mathrm{homotopylimit}}~\underset{\mathrm{Spec}}{\mathcal{O}}^\mathrm{CS}\breve{\Phi}^I_{\psi,\Gamma,-},	\\
}
\]
\[
\xymatrix@R+0pc@C+0pc{
\underset{r}{\mathrm{homotopycolimit}}~\underset{\mathrm{Spec}}{\mathcal{O}}^\mathrm{CS}{\Phi}^r_{\psi,\Gamma,-},\underset{I}{\mathrm{homotopylimit}}~\underset{\mathrm{Spec}}{\mathcal{O}}^\mathrm{CS}{\Phi}^I_{\psi,\Gamma,-}.	
}
\]
\[ 
\xymatrix@R+0pc@C+0pc{
\underset{r}{\mathrm{homotopycolimit}}~\underset{\mathrm{Spec}}{\mathcal{O}}^\mathrm{CS}\widetilde{\Phi}^r_{\psi,\Gamma,-}/\mathrm{Fro}^\mathbb{Z},\underset{I}{\mathrm{homotopylimit}}~\underset{\mathrm{Spec}}{\mathcal{O}}^\mathrm{CS}\widetilde{\Phi}^I_{\psi,\Gamma,-}/\mathrm{Fro}^\mathbb{Z},	\\
}
\]
\[ 
\xymatrix@R+0pc@C+0pc{
\underset{r}{\mathrm{homotopycolimit}}~\underset{\mathrm{Spec}}{\mathcal{O}}^\mathrm{CS}\breve{\Phi}^r_{\psi,\Gamma,-}/\mathrm{Fro}^\mathbb{Z},\underset{I}{\mathrm{homotopylimit}}~\breve{\Phi}^I_{\psi,\Gamma,-}/\mathrm{Fro}^\mathbb{Z},	\\
}
\]
\[ 
\xymatrix@R+0pc@C+0pc{
\underset{r}{\mathrm{homotopycolimit}}~\underset{\mathrm{Spec}}{\mathcal{O}}^\mathrm{CS}{\Phi}^r_{\psi,\Gamma,-}/\mathrm{Fro}^\mathbb{Z},\underset{I}{\mathrm{homotopylimit}}~\underset{\mathrm{Spec}}{\mathcal{O}}^\mathrm{CS}{\Phi}^I_{\psi,\Gamma,-}/\mathrm{Fro}^\mathbb{Z}.	
}
\]	
In this situation we will have the target category being family parametrized by $r$ or $I$ in compatible glueing sense as in \cite[Definition 5.4.10]{10KL2}. In this situation for modules parametrized by the intervals we have the equivalence of $\infty$-categories by using \cite[Proposition 12.18]{10CS2}. Here the corresponding quasicoherent Frobenius modules are defined to be the corresponding homotopy colimits and limits of Frobenius modules:
\begin{align}
\underset{r}{\mathrm{homotopycolimit}}~M_r,\\
\underset{I}{\mathrm{homotopylimit}}~M_I,	
\end{align}
where each $M_r$ is a Frobenius-equivariant module over the period ring with respect to some radius $r$ while each $M_I$ is a Frobenius-equivariant module over the period ring with respect to some interval $I$.\\
\end{proposition}

\begin{proposition}
Similar proposition holds for 
\begin{align}
\mathrm{Quasicoherentsheaves,Perfectcomplex,IndBanach}_{*}.	
\end{align}	
\end{proposition}

\newpage
\subsection{Frobenius Quasicoherent Modules III: Deformation in $(\infty,1)$-Ind-Banach Rings}

\begin{definition}
Let $\psi$ be a toric tower over $\mathbb{Q}_p$ as in \cite[Chapter 7]{10KL2} with base $\mathbb{Q}_p\left<X_1^{\pm 1},...,X_k^{\pm 1}\right>$. Then from \cite{10KL1} and \cite[Definition 5.2.1]{10KL2} we have the following class of Kedlaya-Liu rings (with the following replacement: $\Delta$ stands for $A$, $\nabla$ stands for $B$, while $\Phi$ stands for $C$) by taking product in the sense of self $\Gamma$-th power:

\[
\xymatrix@R+0pc@C+0pc{
\widetilde{\Delta}_{\psi,\Gamma},\widetilde{\nabla}_{\psi,\Gamma},\widetilde{\Phi}_{\psi,\Gamma},\widetilde{\Delta}^+_{\psi,\Gamma},\widetilde{\nabla}^+_{\psi,\Gamma},\widetilde{\Delta}^\dagger_{\psi,\Gamma},\widetilde{\nabla}^\dagger_{\psi,\Gamma},\widetilde{\Phi}^r_{\psi,\Gamma},\widetilde{\Phi}^I_{\psi,\Gamma}, 
}
\]

\[
\xymatrix@R+0pc@C+0pc{
\breve{\Delta}_{\psi,\Gamma},\breve{\nabla}_{\psi,\Gamma},\breve{\Phi}_{\psi,\Gamma},\breve{\Delta}^+_{\psi,\Gamma},\breve{\nabla}^+_{\psi,\Gamma},\breve{\Delta}^\dagger_{\psi,\Gamma},\breve{\nabla}^\dagger_{\psi,\Gamma},\breve{\Phi}^r_{\psi,\Gamma},\breve{\Phi}^I_{\psi,\Gamma},	
}
\]

\[
\xymatrix@R+0pc@C+0pc{
{\Delta}_{\psi,\Gamma},{\nabla}_{\psi,\Gamma},{\Phi}_{\psi,\Gamma},{\Delta}^+_{\psi,\Gamma},{\nabla}^+_{\psi,\Gamma},{\Delta}^\dagger_{\psi,\Gamma},{\nabla}^\dagger_{\psi,\Gamma},{\Phi}^r_{\psi,\Gamma},{\Phi}^I_{\psi,\Gamma}.	
}
\]
We now consider the following rings with $\square$ being a homotopy colimit
\begin{align}
 \underset{I}{\mathrm{homotopylimit}}\square_i
 \end{align}
 of $\mathbb{Q}_p\left<Y_1,...,Y_i\right>,i=1,2,...$ in $\infty$-categories of simplicial ind-Banach rings in \cite{10BBBK}
 \begin{align}
  \mathrm{SimplicialInd-BanachRings}_{\mathbb{Q}_p}
\end{align}  
or animated analytic condensed commutative algebras in \cite{10CS2} 
\begin{align}   
\mathrm{SimplicialAnalyticCondensed}_{\mathbb{Q}_p}.
\end{align}   
Taking the product we have:
\[
\xymatrix@R+0pc@C+0pc{
\widetilde{\Phi}_{\psi,\Gamma,\square},\widetilde{\Phi}^r_{\psi,\Gamma,\square},\widetilde{\Phi}^I_{\psi,\Gamma,\square},	
}
\]
\[
\xymatrix@R+0pc@C+0pc{
\breve{\Phi}_{\psi,\Gamma,\square},\breve{\Phi}^r_{\psi,\Gamma,\square},\breve{\Phi}^I_{\psi,\Gamma,\square},	
}
\]
\[
\xymatrix@R+0pc@C+0pc{
{\Phi}_{\psi,\Gamma,\square},{\Phi}^r_{\psi,\Gamma,\square},{\Phi}^I_{\psi,\Gamma,\square}.	
}
\]
They carry multi Frobenius action $\varphi_\Gamma$ and multi $\mathrm{Lie}_\Gamma:=\mathbb{Z}_p^{\times\Gamma}$ action. In our current situation after \cite{10CKZ} and \cite{10PZ} we consider the following $(\infty,1)$-categories of $(\infty,1)$-modules.\\
\end{definition}

\begin{definition}
First we consider the Bambozzi-Kremnizer spectrum $\underset{\mathrm{Spec}}{\mathcal{O}}^\mathrm{BK}(*)$ attached to any of those in the above from \cite{10BK} by taking derived rational localization:
\begin{align}
&\underset{\mathrm{Spec}}{\mathcal{O}}^\mathrm{BK}\widetilde{\Phi}_{\psi,\Gamma,\square},\underset{\mathrm{Spec}}{\mathcal{O}}^\mathrm{BK}\widetilde{\Phi}^r_{\psi,\Gamma,\square},\underset{\mathrm{Spec}}{\mathcal{O}}^\mathrm{BK}\widetilde{\Phi}^I_{\psi,\Gamma,\square},	
\end{align}
\begin{align}
&\underset{\mathrm{Spec}}{\mathcal{O}}^\mathrm{BK}\breve{\Phi}_{\psi,\Gamma,\square},\underset{\mathrm{Spec}}{\mathcal{O}}^\mathrm{BK}\breve{\Phi}^r_{\psi,\Gamma,\square},\underset{\mathrm{Spec}}{\mathcal{O}}^\mathrm{BK}\breve{\Phi}^I_{\psi,\Gamma,\square},	
\end{align}
\begin{align}
&\underset{\mathrm{Spec}}{\mathcal{O}}^\mathrm{BK}{\Phi}_{\psi,\Gamma,\square},
\underset{\mathrm{Spec}}{\mathcal{O}}^\mathrm{BK}{\Phi}^r_{\psi,\Gamma,\square},\underset{\mathrm{Spec}}{\mathcal{O}}^\mathrm{BK}{\Phi}^I_{\psi,\Gamma,\square}.	
\end{align}

Then we take the corresponding quotients by using the corresponding Frobenius operators:
\begin{align}
&\underset{\mathrm{Spec}}{\mathcal{O}}^\mathrm{BK}\widetilde{\Phi}_{\psi,\Gamma,\square}/\mathrm{Fro}^\mathbb{Z},	\\
\end{align}
\begin{align}
&\underset{\mathrm{Spec}}{\mathcal{O}}^\mathrm{BK}\breve{\Phi}_{\psi,\Gamma,\square}/\mathrm{Fro}^\mathbb{Z},	\\
\end{align}
\begin{align}
&\underset{\mathrm{Spec}}{\mathcal{O}}^\mathrm{BK}{\Phi}_{\psi,\Gamma,\square}/\mathrm{Fro}^\mathbb{Z}.	
\end{align}
Here for those space without notation related to the radius and the corresponding interval we consider the total unions $\bigcap_r,\bigcup_I$ in order to achieve the whole spaces to achieve the analogues of the corresponding FF curves from \cite{10KL1}, \cite{10KL2}, \cite{10FF} for
\[
\xymatrix@R+0pc@C+0pc{
\underset{r}{\mathrm{homotopycolimit}}~\underset{\mathrm{Spec}}{\mathcal{O}}^\mathrm{BK}\widetilde{\Phi}^r_{\psi,\Gamma,\square},\underset{I}{\mathrm{homotopylimit}}~\underset{\mathrm{Spec}}{\mathcal{O}}^\mathrm{BK}\widetilde{\Phi}^I_{\psi,\Gamma,\square},	\\
}
\]
\[
\xymatrix@R+0pc@C+0pc{
\underset{r}{\mathrm{homotopycolimit}}~\underset{\mathrm{Spec}}{\mathcal{O}}^\mathrm{BK}\breve{\Phi}^r_{\psi,\Gamma,\square},\underset{I}{\mathrm{homotopylimit}}~\underset{\mathrm{Spec}}{\mathcal{O}}^\mathrm{BK}\breve{\Phi}^I_{\psi,\Gamma,\square},	\\
}
\]
\[
\xymatrix@R+0pc@C+0pc{
\underset{r}{\mathrm{homotopycolimit}}~\underset{\mathrm{Spec}}{\mathcal{O}}^\mathrm{BK}{\Phi}^r_{\psi,\Gamma,\square},\underset{I}{\mathrm{homotopylimit}}~\underset{\mathrm{Spec}}{\mathcal{O}}^\mathrm{BK}{\Phi}^I_{\psi,\Gamma,\square}.	
}
\]
\[  
\xymatrix@R+0pc@C+0pc{
\underset{r}{\mathrm{homotopycolimit}}~\underset{\mathrm{Spec}}{\mathcal{O}}^\mathrm{BK}\widetilde{\Phi}^r_{\psi,\Gamma,\square}/\mathrm{Fro}^\mathbb{Z},\underset{I}{\mathrm{homotopylimit}}~\underset{\mathrm{Spec}}{\mathcal{O}}^\mathrm{BK}\widetilde{\Phi}^I_{\psi,\Gamma,\square}/\mathrm{Fro}^\mathbb{Z},	\\
}
\]
\[ 
\xymatrix@R+0pc@C+0pc{
\underset{r}{\mathrm{homotopycolimit}}~\underset{\mathrm{Spec}}{\mathcal{O}}^\mathrm{BK}\breve{\Phi}^r_{\psi,\Gamma,\square}/\mathrm{Fro}^\mathbb{Z},\underset{I}{\mathrm{homotopylimit}}~\underset{\mathrm{Spec}}{\mathcal{O}}^\mathrm{BK}\breve{\Phi}^I_{\psi,\Gamma,\square}/\mathrm{Fro}^\mathbb{Z},	\\
}
\]
\[ 
\xymatrix@R+0pc@C+0pc{
\underset{r}{\mathrm{homotopycolimit}}~\underset{\mathrm{Spec}}{\mathcal{O}}^\mathrm{BK}{\Phi}^r_{\psi,\Gamma,\square}/\mathrm{Fro}^\mathbb{Z},\underset{I}{\mathrm{homotopylimit}}~\underset{\mathrm{Spec}}{\mathcal{O}}^\mathrm{BK}{\Phi}^I_{\psi,\Gamma,\square}/\mathrm{Fro}^\mathbb{Z}.	
}
\]

\end{definition}

\indent Meanwhile we have the corresponding Clausen-Scholze analytic stacks from \cite{10CS2}, therefore applying their construction we have:

\begin{definition}
Here we define the following products by using the solidified tensor product from \cite{10CS1} and \cite{10CS2}. Namely $A$ will still as above as a Banach ring over $\mathbb{Q}_p$. Then we take solidified tensor product $\overset{\blacksquare}{\otimes}$ of any of the following
\[
\xymatrix@R+0pc@C+0pc{
\widetilde{\Delta}_{\psi,\Gamma},\widetilde{\nabla}_{\psi,\Gamma},\widetilde{\Phi}_{\psi,\Gamma},\widetilde{\Delta}^+_{\psi,\Gamma},\widetilde{\nabla}^+_{\psi,\Gamma},\widetilde{\Delta}^\dagger_{\psi,\Gamma},\widetilde{\nabla}^\dagger_{\psi,\Gamma},\widetilde{\Phi}^r_{\psi,\Gamma},\widetilde{\Phi}^I_{\psi,\Gamma}, 
}
\]

\[
\xymatrix@R+0pc@C+0pc{
\breve{\Delta}_{\psi,\Gamma},\breve{\nabla}_{\psi,\Gamma},\breve{\Phi}_{\psi,\Gamma},\breve{\Delta}^+_{\psi,\Gamma},\breve{\nabla}^+_{\psi,\Gamma},\breve{\Delta}^\dagger_{\psi,\Gamma},\breve{\nabla}^\dagger_{\psi,\Gamma},\breve{\Phi}^r_{\psi,\Gamma},\breve{\Phi}^I_{\psi,\Gamma},	
}
\]

\[
\xymatrix@R+0pc@C+0pc{
{\Delta}_{\psi,\Gamma},{\nabla}_{\psi,\Gamma},{\Phi}_{\psi,\Gamma},{\Delta}^+_{\psi,\Gamma},{\nabla}^+_{\psi,\Gamma},{\Delta}^\dagger_{\psi,\Gamma},{\nabla}^\dagger_{\psi,\Gamma},{\Phi}^r_{\psi,\Gamma},{\Phi}^I_{\psi,\Gamma},	
}
\]  	
with $A$. Then we have the notations:
\[
\xymatrix@R+0pc@C+0pc{
\widetilde{\Delta}_{\psi,\Gamma,\square},\widetilde{\nabla}_{\psi,\Gamma,\square},\widetilde{\Phi}_{\psi,\Gamma,\square},\widetilde{\Delta}^+_{\psi,\Gamma,\square},\widetilde{\nabla}^+_{\psi,\Gamma,\square},\widetilde{\Delta}^\dagger_{\psi,\Gamma,\square},\widetilde{\nabla}^\dagger_{\psi,\Gamma,\square},\widetilde{\Phi}^r_{\psi,\Gamma,\square},\widetilde{\Phi}^I_{\psi,\Gamma,\square}, 
}
\]

\[
\xymatrix@R+0pc@C+0pc{
\breve{\Delta}_{\psi,\Gamma,\square},\breve{\nabla}_{\psi,\Gamma,\square},\breve{\Phi}_{\psi,\Gamma,\square},\breve{\Delta}^+_{\psi,\Gamma,\square},\breve{\nabla}^+_{\psi,\Gamma,\square},\breve{\Delta}^\dagger_{\psi,\Gamma,\square},\breve{\nabla}^\dagger_{\psi,\Gamma,\square},\breve{\Phi}^r_{\psi,\Gamma,\square},\breve{\Phi}^I_{\psi,\Gamma,\square},	
}
\]

\[
\xymatrix@R+0pc@C+0pc{
{\Delta}_{\psi,\Gamma,\square},{\nabla}_{\psi,\Gamma,\square},{\Phi}_{\psi,\Gamma,\square},{\Delta}^+_{\psi,\Gamma,\square},{\nabla}^+_{\psi,\Gamma,\square},{\Delta}^\dagger_{\psi,\Gamma,\square},{\nabla}^\dagger_{\psi,\Gamma,\square},{\Phi}^r_{\psi,\Gamma,\square},{\Phi}^I_{\psi,\Gamma,\square}.	
}
\]
\end{definition}

\begin{definition}
First we consider the Clausen-Scholze spectrum $\underset{\mathrm{Spec}}{\mathcal{O}}^\mathrm{CS}(*)$ attached to any of those in the above from \cite{10CS2} by taking derived rational localization:
\begin{align}
\underset{\mathrm{Spec}}{\mathcal{O}}^\mathrm{CS}\widetilde{\Delta}_{\psi,\Gamma,\square},\underset{\mathrm{Spec}}{\mathcal{O}}^\mathrm{CS}\widetilde{\nabla}_{\psi,\Gamma,\square},\underset{\mathrm{Spec}}{\mathcal{O}}^\mathrm{CS}\widetilde{\Phi}_{\psi,\Gamma,\square},\underset{\mathrm{Spec}}{\mathcal{O}}^\mathrm{CS}\widetilde{\Delta}^+_{\psi,\Gamma,\square},\underset{\mathrm{Spec}}{\mathcal{O}}^\mathrm{CS}\widetilde{\nabla}^+_{\psi,\Gamma,\square},\\
\underset{\mathrm{Spec}}{\mathcal{O}}^\mathrm{CS}\widetilde{\Delta}^\dagger_{\psi,\Gamma,\square},\underset{\mathrm{Spec}}{\mathcal{O}}^\mathrm{CS}\widetilde{\nabla}^\dagger_{\psi,\Gamma,\square},\underset{\mathrm{Spec}}{\mathcal{O}}^\mathrm{CS}\widetilde{\Phi}^r_{\psi,\Gamma,\square},\underset{\mathrm{Spec}}{\mathcal{O}}^\mathrm{CS}\widetilde{\Phi}^I_{\psi,\Gamma,\square},	\\
\end{align}
\begin{align}
\underset{\mathrm{Spec}}{\mathcal{O}}^\mathrm{CS}\breve{\Delta}_{\psi,\Gamma,\square},\breve{\nabla}_{\psi,\Gamma,\square},\underset{\mathrm{Spec}}{\mathcal{O}}^\mathrm{CS}\breve{\Phi}_{\psi,\Gamma,\square},\underset{\mathrm{Spec}}{\mathcal{O}}^\mathrm{CS}\breve{\Delta}^+_{\psi,\Gamma,\square},\underset{\mathrm{Spec}}{\mathcal{O}}^\mathrm{CS}\breve{\nabla}^+_{\psi,\Gamma,\square},\\
\underset{\mathrm{Spec}}{\mathcal{O}}^\mathrm{CS}\breve{\Delta}^\dagger_{\psi,\Gamma,\square},\underset{\mathrm{Spec}}{\mathcal{O}}^\mathrm{CS}\breve{\nabla}^\dagger_{\psi,\Gamma,\square},\underset{\mathrm{Spec}}{\mathcal{O}}^\mathrm{CS}\breve{\Phi}^r_{\psi,\Gamma,\square},\breve{\Phi}^I_{\psi,\Gamma,\square},	\\
\end{align}
\begin{align}
\underset{\mathrm{Spec}}{\mathcal{O}}^\mathrm{CS}{\Delta}_{\psi,\Gamma,\square},\underset{\mathrm{Spec}}{\mathcal{O}}^\mathrm{CS}{\nabla}_{\psi,\Gamma,\square},\underset{\mathrm{Spec}}{\mathcal{O}}^\mathrm{CS}{\Phi}_{\psi,\Gamma,\square},\underset{\mathrm{Spec}}{\mathcal{O}}^\mathrm{CS}{\Delta}^+_{\psi,\Gamma,\square},\underset{\mathrm{Spec}}{\mathcal{O}}^\mathrm{CS}{\nabla}^+_{\psi,\Gamma,\square},\\
\underset{\mathrm{Spec}}{\mathcal{O}}^\mathrm{CS}{\Delta}^\dagger_{\psi,\Gamma,\square},\underset{\mathrm{Spec}}{\mathcal{O}}^\mathrm{CS}{\nabla}^\dagger_{\psi,\Gamma,\square},\underset{\mathrm{Spec}}{\mathcal{O}}^\mathrm{CS}{\Phi}^r_{\psi,\Gamma,\square},\underset{\mathrm{Spec}}{\mathcal{O}}^\mathrm{CS}{\Phi}^I_{\psi,\Gamma,\square}.	
\end{align}

Then we take the corresponding quotients by using the corresponding Frobenius operators:
\begin{align}
&\underset{\mathrm{Spec}}{\mathcal{O}}^\mathrm{CS}\widetilde{\Delta}_{\psi,\Gamma,\square}/\mathrm{Fro}^\mathbb{Z},\underset{\mathrm{Spec}}{\mathcal{O}}^\mathrm{CS}\widetilde{\nabla}_{\psi,\Gamma,\square}/\mathrm{Fro}^\mathbb{Z},\underset{\mathrm{Spec}}{\mathcal{O}}^\mathrm{CS}\widetilde{\Phi}_{\psi,\Gamma,\square}/\mathrm{Fro}^\mathbb{Z},\underset{\mathrm{Spec}}{\mathcal{O}}^\mathrm{CS}\widetilde{\Delta}^+_{\psi,\Gamma,\square}/\mathrm{Fro}^\mathbb{Z},\\
&\underset{\mathrm{Spec}}{\mathcal{O}}^\mathrm{CS}\widetilde{\nabla}^+_{\psi,\Gamma,\square}/\mathrm{Fro}^\mathbb{Z}, \underset{\mathrm{Spec}}{\mathcal{O}}^\mathrm{CS}\widetilde{\Delta}^\dagger_{\psi,\Gamma,\square}/\mathrm{Fro}^\mathbb{Z},\underset{\mathrm{Spec}}{\mathcal{O}}^\mathrm{CS}\widetilde{\nabla}^\dagger_{\psi,\Gamma,\square}/\mathrm{Fro}^\mathbb{Z},	\\
\end{align}
\begin{align}
&\underset{\mathrm{Spec}}{\mathcal{O}}^\mathrm{CS}\breve{\Delta}_{\psi,\Gamma,\square}/\mathrm{Fro}^\mathbb{Z},\breve{\nabla}_{\psi,\Gamma,\square}/\mathrm{Fro}^\mathbb{Z},\underset{\mathrm{Spec}}{\mathcal{O}}^\mathrm{CS}\breve{\Phi}_{\psi,\Gamma,\square}/\mathrm{Fro}^\mathbb{Z},\underset{\mathrm{Spec}}{\mathcal{O}}^\mathrm{CS}\breve{\Delta}^+_{\psi,\Gamma,\square}/\mathrm{Fro}^\mathbb{Z},\\
&\underset{\mathrm{Spec}}{\mathcal{O}}^\mathrm{CS}\breve{\nabla}^+_{\psi,\Gamma,\square}/\mathrm{Fro}^\mathbb{Z}, \underset{\mathrm{Spec}}{\mathcal{O}}^\mathrm{CS}\breve{\Delta}^\dagger_{\psi,\Gamma,\square}/\mathrm{Fro}^\mathbb{Z},\underset{\mathrm{Spec}}{\mathcal{O}}^\mathrm{CS}\breve{\nabla}^\dagger_{\psi,\Gamma,\square}/\mathrm{Fro}^\mathbb{Z},	\\
\end{align}
\begin{align}
&\underset{\mathrm{Spec}}{\mathcal{O}}^\mathrm{CS}{\Delta}_{\psi,\Gamma,\square}/\mathrm{Fro}^\mathbb{Z},\underset{\mathrm{Spec}}{\mathcal{O}}^\mathrm{CS}{\nabla}_{\psi,\Gamma,\square}/\mathrm{Fro}^\mathbb{Z},\underset{\mathrm{Spec}}{\mathcal{O}}^\mathrm{CS}{\Phi}_{\psi,\Gamma,\square}/\mathrm{Fro}^\mathbb{Z},\underset{\mathrm{Spec}}{\mathcal{O}}^\mathrm{CS}{\Delta}^+_{\psi,\Gamma,\square}/\mathrm{Fro}^\mathbb{Z},\\
&\underset{\mathrm{Spec}}{\mathcal{O}}^\mathrm{CS}{\nabla}^+_{\psi,\Gamma,\square}/\mathrm{Fro}^\mathbb{Z}, \underset{\mathrm{Spec}}{\mathcal{O}}^\mathrm{CS}{\Delta}^\dagger_{\psi,\Gamma,\square}/\mathrm{Fro}^\mathbb{Z},\underset{\mathrm{Spec}}{\mathcal{O}}^\mathrm{CS}{\nabla}^\dagger_{\psi,\Gamma,\square}/\mathrm{Fro}^\mathbb{Z}.	
\end{align}
Here for those space with notations related to the radius and the corresponding interval we consider the total unions $\bigcap_r,\bigcup_I$ in order to achieve the whole spaces to achieve the analogues of the corresponding FF curves from \cite{10KL1}, \cite{10KL2}, \cite{10FF} for
\[
\xymatrix@R+0pc@C+0pc{
\underset{r}{\mathrm{homotopycolimit}}~\underset{\mathrm{Spec}}{\mathcal{O}}^\mathrm{CS}\widetilde{\Phi}^r_{\psi,\Gamma,\square},\underset{I}{\mathrm{homotopylimit}}~\underset{\mathrm{Spec}}{\mathcal{O}}^\mathrm{CS}\widetilde{\Phi}^I_{\psi,\Gamma,\square},	\\
}
\]
\[
\xymatrix@R+0pc@C+0pc{
\underset{r}{\mathrm{homotopycolimit}}~\underset{\mathrm{Spec}}{\mathcal{O}}^\mathrm{CS}\breve{\Phi}^r_{\psi,\Gamma,\square},\underset{I}{\mathrm{homotopylimit}}~\underset{\mathrm{Spec}}{\mathcal{O}}^\mathrm{CS}\breve{\Phi}^I_{\psi,\Gamma,\square},	\\
}
\]
\[
\xymatrix@R+0pc@C+0pc{
\underset{r}{\mathrm{homotopycolimit}}~\underset{\mathrm{Spec}}{\mathcal{O}}^\mathrm{CS}{\Phi}^r_{\psi,\Gamma,\square},\underset{I}{\mathrm{homotopylimit}}~\underset{\mathrm{Spec}}{\mathcal{O}}^\mathrm{CS}{\Phi}^I_{\psi,\Gamma,\square}.	
}
\]
\[ 
\xymatrix@R+0pc@C+0pc{
\underset{r}{\mathrm{homotopycolimit}}~\underset{\mathrm{Spec}}{\mathcal{O}}^\mathrm{CS}\widetilde{\Phi}^r_{\psi,\Gamma,\square}/\mathrm{Fro}^\mathbb{Z},\underset{I}{\mathrm{homotopylimit}}~\underset{\mathrm{Spec}}{\mathcal{O}}^\mathrm{CS}\widetilde{\Phi}^I_{\psi,\Gamma,\square}/\mathrm{Fro}^\mathbb{Z},	\\
}
\]
\[ 
\xymatrix@R+0pc@C+0pc{
\underset{r}{\mathrm{homotopycolimit}}~\underset{\mathrm{Spec}}{\mathcal{O}}^\mathrm{CS}\breve{\Phi}^r_{\psi,\Gamma,\square}/\mathrm{Fro}^\mathbb{Z},\underset{I}{\mathrm{homotopylimit}}~\breve{\Phi}^I_{\psi,\Gamma,\square}/\mathrm{Fro}^\mathbb{Z},	\\
}
\]
\[ 
\xymatrix@R+0pc@C+0pc{
\underset{r}{\mathrm{homotopycolimit}}~\underset{\mathrm{Spec}}{\mathcal{O}}^\mathrm{CS}{\Phi}^r_{\psi,\Gamma,\square}/\mathrm{Fro}^\mathbb{Z},\underset{I}{\mathrm{homotopylimit}}~\underset{\mathrm{Spec}}{\mathcal{O}}^\mathrm{CS}{\Phi}^I_{\psi,\Gamma,\square}/\mathrm{Fro}^\mathbb{Z}.	
}
\]

\end{definition}

\

\begin{definition}
We then consider the corresponding quasipresheaves of the corresponding ind-Banach or monomorphic ind-Banach modules from \cite{10BBK}, \cite{10KKM}\footnote{Here the categories are defined to be the corresponding homotopy colimits of the corresponding categories with respect to each $\square_i$.}:
\begin{align}
\mathrm{Quasicoherentpresheaves,IndBanach}_{*}	
\end{align}
where $*$ is one of the following spaces:
\begin{align}
&\underset{\mathrm{Spec}}{\mathcal{O}}^\mathrm{BK}\widetilde{\Phi}_{\psi,\Gamma,\square}/\mathrm{Fro}^\mathbb{Z},	\\
\end{align}
\begin{align}
&\underset{\mathrm{Spec}}{\mathcal{O}}^\mathrm{BK}\breve{\Phi}_{\psi,\Gamma,\square}/\mathrm{Fro}^\mathbb{Z},	\\
\end{align}
\begin{align}
&\underset{\mathrm{Spec}}{\mathcal{O}}^\mathrm{BK}{\Phi}_{\psi,\Gamma,\square}/\mathrm{Fro}^\mathbb{Z}.	
\end{align}
Here for those space without notation related to the radius and the corresponding interval we consider the total unions $\bigcap_r,\bigcup_I$ in order to achieve the whole spaces to achieve the analogues of the corresponding FF curves from \cite{10KL1}, \cite{10KL2}, \cite{10FF} for
\[
\xymatrix@R+0pc@C+0pc{
\underset{r}{\mathrm{homotopycolimit}}~\underset{\mathrm{Spec}}{\mathcal{O}}^\mathrm{BK}\widetilde{\Phi}^r_{\psi,\Gamma,\square},\underset{I}{\mathrm{homotopylimit}}~\underset{\mathrm{Spec}}{\mathcal{O}}^\mathrm{BK}\widetilde{\Phi}^I_{\psi,\Gamma,\square},	\\
}
\]
\[
\xymatrix@R+0pc@C+0pc{
\underset{r}{\mathrm{homotopycolimit}}~\underset{\mathrm{Spec}}{\mathcal{O}}^\mathrm{BK}\breve{\Phi}^r_{\psi,\Gamma,\square},\underset{I}{\mathrm{homotopylimit}}~\underset{\mathrm{Spec}}{\mathcal{O}}^\mathrm{BK}\breve{\Phi}^I_{\psi,\Gamma,\square},	\\
}
\]
\[
\xymatrix@R+0pc@C+0pc{
\underset{r}{\mathrm{homotopycolimit}}~\underset{\mathrm{Spec}}{\mathcal{O}}^\mathrm{BK}{\Phi}^r_{\psi,\Gamma,\square},\underset{I}{\mathrm{homotopylimit}}~\underset{\mathrm{Spec}}{\mathcal{O}}^\mathrm{BK}{\Phi}^I_{\psi,\Gamma,\square}.	
}
\]
\[  
\xymatrix@R+0pc@C+0pc{
\underset{r}{\mathrm{homotopycolimit}}~\underset{\mathrm{Spec}}{\mathcal{O}}^\mathrm{BK}\widetilde{\Phi}^r_{\psi,\Gamma,\square}/\mathrm{Fro}^\mathbb{Z},\underset{I}{\mathrm{homotopylimit}}~\underset{\mathrm{Spec}}{\mathcal{O}}^\mathrm{BK}\widetilde{\Phi}^I_{\psi,\Gamma,\square}/\mathrm{Fro}^\mathbb{Z},	\\
}
\]
\[ 
\xymatrix@R+0pc@C+0pc{
\underset{r}{\mathrm{homotopycolimit}}~\underset{\mathrm{Spec}}{\mathcal{O}}^\mathrm{BK}\breve{\Phi}^r_{\psi,\Gamma,\square}/\mathrm{Fro}^\mathbb{Z},\underset{I}{\mathrm{homotopylimit}}~\underset{\mathrm{Spec}}{\mathcal{O}}^\mathrm{BK}\breve{\Phi}^I_{\psi,\Gamma,\square}/\mathrm{Fro}^\mathbb{Z},	\\
}
\]
\[ 
\xymatrix@R+0pc@C+0pc{
\underset{r}{\mathrm{homotopycolimit}}~\underset{\mathrm{Spec}}{\mathcal{O}}^\mathrm{BK}{\Phi}^r_{\psi,\Gamma,\square}/\mathrm{Fro}^\mathbb{Z},\underset{I}{\mathrm{homotopylimit}}~\underset{\mathrm{Spec}}{\mathcal{O}}^\mathrm{BK}{\Phi}^I_{\psi,\Gamma,\square}/\mathrm{Fro}^\mathbb{Z}.	
}
\]

\end{definition}

\begin{definition}
We then consider the corresponding quasisheaves of the corresponding condensed solid topological modules from \cite{10CS2}:
\begin{align}
\mathrm{Quasicoherentsheaves, Condensed}_{*}	
\end{align}
where $*$ is one of the following spaces:
\begin{align}
&\underset{\mathrm{Spec}}{\mathcal{O}}^\mathrm{CS}\widetilde{\Delta}_{\psi,\Gamma,\square}/\mathrm{Fro}^\mathbb{Z},\underset{\mathrm{Spec}}{\mathcal{O}}^\mathrm{CS}\widetilde{\nabla}_{\psi,\Gamma,\square}/\mathrm{Fro}^\mathbb{Z},\underset{\mathrm{Spec}}{\mathcal{O}}^\mathrm{CS}\widetilde{\Phi}_{\psi,\Gamma,\square}/\mathrm{Fro}^\mathbb{Z},\underset{\mathrm{Spec}}{\mathcal{O}}^\mathrm{CS}\widetilde{\Delta}^+_{\psi,\Gamma,\square}/\mathrm{Fro}^\mathbb{Z},\\
&\underset{\mathrm{Spec}}{\mathcal{O}}^\mathrm{CS}\widetilde{\nabla}^+_{\psi,\Gamma,\square}/\mathrm{Fro}^\mathbb{Z},\underset{\mathrm{Spec}}{\mathcal{O}}^\mathrm{CS}\widetilde{\Delta}^\dagger_{\psi,\Gamma,\square}/\mathrm{Fro}^\mathbb{Z},\underset{\mathrm{Spec}}{\mathcal{O}}^\mathrm{CS}\widetilde{\nabla}^\dagger_{\psi,\Gamma,\square}/\mathrm{Fro}^\mathbb{Z},	\\
\end{align}
\begin{align}
&\underset{\mathrm{Spec}}{\mathcal{O}}^\mathrm{CS}\breve{\Delta}_{\psi,\Gamma,\square}/\mathrm{Fro}^\mathbb{Z},\breve{\nabla}_{\psi,\Gamma,\square}/\mathrm{Fro}^\mathbb{Z},\underset{\mathrm{Spec}}{\mathcal{O}}^\mathrm{CS}\breve{\Phi}_{\psi,\Gamma,\square}/\mathrm{Fro}^\mathbb{Z},\underset{\mathrm{Spec}}{\mathcal{O}}^\mathrm{CS}\breve{\Delta}^+_{\psi,\Gamma,\square}/\mathrm{Fro}^\mathbb{Z},\\
&\underset{\mathrm{Spec}}{\mathcal{O}}^\mathrm{CS}\breve{\nabla}^+_{\psi,\Gamma,\square}/\mathrm{Fro}^\mathbb{Z},\underset{\mathrm{Spec}}{\mathcal{O}}^\mathrm{CS}\breve{\Delta}^\dagger_{\psi,\Gamma,\square}/\mathrm{Fro}^\mathbb{Z},\underset{\mathrm{Spec}}{\mathcal{O}}^\mathrm{CS}\breve{\nabla}^\dagger_{\psi,\Gamma,\square}/\mathrm{Fro}^\mathbb{Z},	\\
\end{align}
\begin{align}
&\underset{\mathrm{Spec}}{\mathcal{O}}^\mathrm{CS}{\Delta}_{\psi,\Gamma,\square}/\mathrm{Fro}^\mathbb{Z},\underset{\mathrm{Spec}}{\mathcal{O}}^\mathrm{CS}{\nabla}_{\psi,\Gamma,\square}/\mathrm{Fro}^\mathbb{Z},\underset{\mathrm{Spec}}{\mathcal{O}}^\mathrm{CS}{\Phi}_{\psi,\Gamma,\square}/\mathrm{Fro}^\mathbb{Z},\underset{\mathrm{Spec}}{\mathcal{O}}^\mathrm{CS}{\Delta}^+_{\psi,\Gamma,\square}/\mathrm{Fro}^\mathbb{Z},\\
&\underset{\mathrm{Spec}}{\mathcal{O}}^\mathrm{CS}{\nabla}^+_{\psi,\Gamma,\square}/\mathrm{Fro}^\mathbb{Z}, \underset{\mathrm{Spec}}{\mathcal{O}}^\mathrm{CS}{\Delta}^\dagger_{\psi,\Gamma,\square}/\mathrm{Fro}^\mathbb{Z},\underset{\mathrm{Spec}}{\mathcal{O}}^\mathrm{CS}{\nabla}^\dagger_{\psi,\Gamma,\square}/\mathrm{Fro}^\mathbb{Z}.	
\end{align}
Here for those space with notations related to the radius and the corresponding interval we consider the total unions $\bigcap_r,\bigcup_I$ in order to achieve the whole spaces to achieve the analogues of the corresponding FF curves from \cite{10KL1}, \cite{10KL2}, \cite{10FF} for
\[
\xymatrix@R+0pc@C+0pc{
\underset{r}{\mathrm{homotopycolimit}}~\underset{\mathrm{Spec}}{\mathcal{O}}^\mathrm{CS}\widetilde{\Phi}^r_{\psi,\Gamma,\square},\underset{I}{\mathrm{homotopylimit}}~\underset{\mathrm{Spec}}{\mathcal{O}}^\mathrm{CS}\widetilde{\Phi}^I_{\psi,\Gamma,\square},	\\
}
\]
\[
\xymatrix@R+0pc@C+0pc{
\underset{r}{\mathrm{homotopycolimit}}~\underset{\mathrm{Spec}}{\mathcal{O}}^\mathrm{CS}\breve{\Phi}^r_{\psi,\Gamma,\square},\underset{I}{\mathrm{homotopylimit}}~\underset{\mathrm{Spec}}{\mathcal{O}}^\mathrm{CS}\breve{\Phi}^I_{\psi,\Gamma,\square},	\\
}
\]
\[
\xymatrix@R+0pc@C+0pc{
\underset{r}{\mathrm{homotopycolimit}}~\underset{\mathrm{Spec}}{\mathcal{O}}^\mathrm{CS}{\Phi}^r_{\psi,\Gamma,\square},\underset{I}{\mathrm{homotopylimit}}~\underset{\mathrm{Spec}}{\mathcal{O}}^\mathrm{CS}{\Phi}^I_{\psi,\Gamma,\square}.	
}
\]
\[ 
\xymatrix@R+0pc@C+0pc{
\underset{r}{\mathrm{homotopycolimit}}~\underset{\mathrm{Spec}}{\mathcal{O}}^\mathrm{CS}\widetilde{\Phi}^r_{\psi,\Gamma,\square}/\mathrm{Fro}^\mathbb{Z},\underset{I}{\mathrm{homotopylimit}}~\underset{\mathrm{Spec}}{\mathcal{O}}^\mathrm{CS}\widetilde{\Phi}^I_{\psi,\Gamma,\square}/\mathrm{Fro}^\mathbb{Z},	\\
}
\]
\[ 
\xymatrix@R+0pc@C+0pc{
\underset{r}{\mathrm{homotopycolimit}}~\underset{\mathrm{Spec}}{\mathcal{O}}^\mathrm{CS}\breve{\Phi}^r_{\psi,\Gamma,\square}/\mathrm{Fro}^\mathbb{Z},\underset{I}{\mathrm{homotopylimit}}~\breve{\Phi}^I_{\psi,\Gamma,\square}/\mathrm{Fro}^\mathbb{Z},	\\
}
\]
\[ 
\xymatrix@R+0pc@C+0pc{
\underset{r}{\mathrm{homotopycolimit}}~\underset{\mathrm{Spec}}{\mathcal{O}}^\mathrm{CS}{\Phi}^r_{\psi,\Gamma,\square}/\mathrm{Fro}^\mathbb{Z},\underset{I}{\mathrm{homotopylimit}}~\underset{\mathrm{Spec}}{\mathcal{O}}^\mathrm{CS}{\Phi}^I_{\psi,\Gamma,\square}/\mathrm{Fro}^\mathbb{Z}.	
}
\]

\end{definition}

\

\begin{proposition}
There is a well-defined functor from the $\infty$-category 
\begin{align}
\mathrm{Quasicoherentpresheaves,Condensed}_{*}	
\end{align}
where $*$ is one of the following spaces:
\begin{align}
&\underset{\mathrm{Spec}}{\mathcal{O}}^\mathrm{CS}\widetilde{\Phi}_{\psi,\Gamma,\square}/\mathrm{Fro}^\mathbb{Z},	\\
\end{align}
\begin{align}
&\underset{\mathrm{Spec}}{\mathcal{O}}^\mathrm{CS}\breve{\Phi}_{\psi,\Gamma,\square}/\mathrm{Fro}^\mathbb{Z},	\\
\end{align}
\begin{align}
&\underset{\mathrm{Spec}}{\mathcal{O}}^\mathrm{CS}{\Phi}_{\psi,\Gamma,\square}/\mathrm{Fro}^\mathbb{Z},	
\end{align}
to the $\infty$-category of $\mathrm{Fro}$-equivariant quasicoherent presheaves over similar spaces above correspondingly without the $\mathrm{Fro}$-quotients, and to the $\infty$-category of $\mathrm{Fro}$-equivariant quasicoherent modules over global sections of the structure $\infty$-sheaves of the similar spaces above correspondingly without the $\mathrm{Fro}$-quotients. Here for those space without notation related to the radius and the corresponding interval we consider the total unions $\bigcap_r,\bigcup_I$ in order to achieve the whole spaces to achieve the analogues of the corresponding FF curves from \cite{10KL1}, \cite{10KL2}, \cite{10FF} for
\[
\xymatrix@R+0pc@C+0pc{
\underset{r}{\mathrm{homotopycolimit}}~\underset{\mathrm{Spec}}{\mathcal{O}}^\mathrm{CS}\widetilde{\Phi}^r_{\psi,\Gamma,\square},\underset{I}{\mathrm{homotopylimit}}~\underset{\mathrm{Spec}}{\mathcal{O}}^\mathrm{CS}\widetilde{\Phi}^I_{\psi,\Gamma,\square},	\\
}
\]
\[
\xymatrix@R+0pc@C+0pc{
\underset{r}{\mathrm{homotopycolimit}}~\underset{\mathrm{Spec}}{\mathcal{O}}^\mathrm{CS}\breve{\Phi}^r_{\psi,\Gamma,\square},\underset{I}{\mathrm{homotopylimit}}~\underset{\mathrm{Spec}}{\mathcal{O}}^\mathrm{CS}\breve{\Phi}^I_{\psi,\Gamma,\square},	\\
}
\]
\[
\xymatrix@R+0pc@C+0pc{
\underset{r}{\mathrm{homotopycolimit}}~\underset{\mathrm{Spec}}{\mathcal{O}}^\mathrm{CS}{\Phi}^r_{\psi,\Gamma,\square},\underset{I}{\mathrm{homotopylimit}}~\underset{\mathrm{Spec}}{\mathcal{O}}^\mathrm{CS}{\Phi}^I_{\psi,\Gamma,\square}.	
}
\]
\[ 
\xymatrix@R+0pc@C+0pc{
\underset{r}{\mathrm{homotopycolimit}}~\underset{\mathrm{Spec}}{\mathcal{O}}^\mathrm{CS}\widetilde{\Phi}^r_{\psi,\Gamma,\square}/\mathrm{Fro}^\mathbb{Z},\underset{I}{\mathrm{homotopylimit}}~\underset{\mathrm{Spec}}{\mathcal{O}}^\mathrm{CS}\widetilde{\Phi}^I_{\psi,\Gamma,\square}/\mathrm{Fro}^\mathbb{Z},	\\
}
\]
\[ 
\xymatrix@R+0pc@C+0pc{
\underset{r}{\mathrm{homotopycolimit}}~\underset{\mathrm{Spec}}{\mathcal{O}}^\mathrm{CS}\breve{\Phi}^r_{\psi,\Gamma,\square}/\mathrm{Fro}^\mathbb{Z},\underset{I}{\mathrm{homotopylimit}}~\breve{\Phi}^I_{\psi,\Gamma,\square}/\mathrm{Fro}^\mathbb{Z},	\\
}
\]
\[ 
\xymatrix@R+0pc@C+0pc{
\underset{r}{\mathrm{homotopycolimit}}~\underset{\mathrm{Spec}}{\mathcal{O}}^\mathrm{CS}{\Phi}^r_{\psi,\Gamma,\square}/\mathrm{Fro}^\mathbb{Z},\underset{I}{\mathrm{homotopylimit}}~\underset{\mathrm{Spec}}{\mathcal{O}}^\mathrm{CS}{\Phi}^I_{\psi,\Gamma,\square}/\mathrm{Fro}^\mathbb{Z}.	
}
\]	
In this situation we will have the target category being family parametrized by $r$ or $I$ in compatible glueing sense as in \cite[Definition 5.4.10]{10KL2}. In this situation for modules parametrized by the intervals we have the equivalence of $\infty$-categories by using \cite[Proposition 13.8]{10CS2}. Here the corresponding quasicoherent Frobenius modules are defined to be the corresponding homotopy colimits and limits of Frobenius modules:
\begin{align}
\underset{r}{\mathrm{homotopycolimit}}~M_r,\\
\underset{I}{\mathrm{homotopylimit}}~M_I,	
\end{align}
where each $M_r$ is a Frobenius-equivariant module over the period ring with respect to some radius $r$ while each $M_I$ is a Frobenius-equivariant module over the period ring with respect to some interval $I$.\\
\end{proposition}

\begin{proposition}
Similar proposition holds for 
\begin{align}
\mathrm{Quasicoherentsheaves,IndBanach}_{*}.	
\end{align}	
\end{proposition}

\

\begin{definition}
We then consider the corresponding quasipresheaves of perfect complexes the corresponding ind-Banach or monomorphic ind-Banach modules from \cite{10BBK}, \cite{10KKM}:
\begin{align}
\mathrm{Quasicoherentpresheaves,Perfectcomplex,IndBanach}_{*}	
\end{align}
where $*$ is one of the following spaces:
\begin{align}
&\underset{\mathrm{Spec}}{\mathcal{O}}^\mathrm{BK}\widetilde{\Phi}_{\psi,\Gamma,\square}/\mathrm{Fro}^\mathbb{Z},	\\
\end{align}
\begin{align}
&\underset{\mathrm{Spec}}{\mathcal{O}}^\mathrm{BK}\breve{\Phi}_{\psi,\Gamma,\square}/\mathrm{Fro}^\mathbb{Z},	\\
\end{align}
\begin{align}
&\underset{\mathrm{Spec}}{\mathcal{O}}^\mathrm{BK}{\Phi}_{\psi,\Gamma,\square}/\mathrm{Fro}^\mathbb{Z}.	
\end{align}
Here for those space without notation related to the radius and the corresponding interval we consider the total unions $\bigcap_r,\bigcup_I$ in order to achieve the whole spaces to achieve the analogues of the corresponding FF curves from \cite{10KL1}, \cite{10KL2}, \cite{10FF} for
\[
\xymatrix@R+0pc@C+0pc{
\underset{r}{\mathrm{homotopycolimit}}~\underset{\mathrm{Spec}}{\mathcal{O}}^\mathrm{BK}\widetilde{\Phi}^r_{\psi,\Gamma,\square},\underset{I}{\mathrm{homotopylimit}}~\underset{\mathrm{Spec}}{\mathcal{O}}^\mathrm{BK}\widetilde{\Phi}^I_{\psi,\Gamma,\square},	\\
}
\]
\[
\xymatrix@R+0pc@C+0pc{
\underset{r}{\mathrm{homotopycolimit}}~\underset{\mathrm{Spec}}{\mathcal{O}}^\mathrm{BK}\breve{\Phi}^r_{\psi,\Gamma,\square},\underset{I}{\mathrm{homotopylimit}}~\underset{\mathrm{Spec}}{\mathcal{O}}^\mathrm{BK}\breve{\Phi}^I_{\psi,\Gamma,\square},	\\
}
\]
\[
\xymatrix@R+0pc@C+0pc{
\underset{r}{\mathrm{homotopycolimit}}~\underset{\mathrm{Spec}}{\mathcal{O}}^\mathrm{BK}{\Phi}^r_{\psi,\Gamma,\square},\underset{I}{\mathrm{homotopylimit}}~\underset{\mathrm{Spec}}{\mathcal{O}}^\mathrm{BK}{\Phi}^I_{\psi,\Gamma,\square}.	
}
\]
\[  
\xymatrix@R+0pc@C+0pc{
\underset{r}{\mathrm{homotopycolimit}}~\underset{\mathrm{Spec}}{\mathcal{O}}^\mathrm{BK}\widetilde{\Phi}^r_{\psi,\Gamma,\square}/\mathrm{Fro}^\mathbb{Z},\underset{I}{\mathrm{homotopylimit}}~\underset{\mathrm{Spec}}{\mathcal{O}}^\mathrm{BK}\widetilde{\Phi}^I_{\psi,\Gamma,\square}/\mathrm{Fro}^\mathbb{Z},	\\
}
\]
\[ 
\xymatrix@R+0pc@C+0pc{
\underset{r}{\mathrm{homotopycolimit}}~\underset{\mathrm{Spec}}{\mathcal{O}}^\mathrm{BK}\breve{\Phi}^r_{\psi,\Gamma,\square}/\mathrm{Fro}^\mathbb{Z},\underset{I}{\mathrm{homotopylimit}}~\underset{\mathrm{Spec}}{\mathcal{O}}^\mathrm{BK}\breve{\Phi}^I_{\psi,\Gamma,\square}/\mathrm{Fro}^\mathbb{Z},	\\
}
\]
\[ 
\xymatrix@R+0pc@C+0pc{
\underset{r}{\mathrm{homotopycolimit}}~\underset{\mathrm{Spec}}{\mathcal{O}}^\mathrm{BK}{\Phi}^r_{\psi,\Gamma,\square}/\mathrm{Fro}^\mathbb{Z},\underset{I}{\mathrm{homotopylimit}}~\underset{\mathrm{Spec}}{\mathcal{O}}^\mathrm{BK}{\Phi}^I_{\psi,\Gamma,\square}/\mathrm{Fro}^\mathbb{Z}.	
}
\]

\end{definition}

\begin{definition}
We then consider the corresponding quasisheaves of perfect complexes of the corresponding condensed solid topological modules from \cite{10CS2}:
\begin{align}
\mathrm{Quasicoherentsheaves, Perfectcomplex, Condensed}_{*}	
\end{align}
where $*$ is one of the following spaces:
\begin{align}
&\underset{\mathrm{Spec}}{\mathcal{O}}^\mathrm{CS}\widetilde{\Delta}_{\psi,\Gamma,\square}/\mathrm{Fro}^\mathbb{Z},\underset{\mathrm{Spec}}{\mathcal{O}}^\mathrm{CS}\widetilde{\nabla}_{\psi,\Gamma,\square}/\mathrm{Fro}^\mathbb{Z},\underset{\mathrm{Spec}}{\mathcal{O}}^\mathrm{CS}\widetilde{\Phi}_{\psi,\Gamma,\square}/\mathrm{Fro}^\mathbb{Z},\underset{\mathrm{Spec}}{\mathcal{O}}^\mathrm{CS}\widetilde{\Delta}^+_{\psi,\Gamma,\square}/\mathrm{Fro}^\mathbb{Z},\\
&\underset{\mathrm{Spec}}{\mathcal{O}}^\mathrm{CS}\widetilde{\nabla}^+_{\psi,\Gamma,\square}/\mathrm{Fro}^\mathbb{Z},\underset{\mathrm{Spec}}{\mathcal{O}}^\mathrm{CS}\widetilde{\Delta}^\dagger_{\psi,\Gamma,\square}/\mathrm{Fro}^\mathbb{Z},\underset{\mathrm{Spec}}{\mathcal{O}}^\mathrm{CS}\widetilde{\nabla}^\dagger_{\psi,\Gamma,\square}/\mathrm{Fro}^\mathbb{Z},	\\
\end{align}
\begin{align}
&\underset{\mathrm{Spec}}{\mathcal{O}}^\mathrm{CS}\breve{\Delta}_{\psi,\Gamma,\square}/\mathrm{Fro}^\mathbb{Z},\breve{\nabla}_{\psi,\Gamma,\square}/\mathrm{Fro}^\mathbb{Z},\underset{\mathrm{Spec}}{\mathcal{O}}^\mathrm{CS}\breve{\Phi}_{\psi,\Gamma,\square}/\mathrm{Fro}^\mathbb{Z},\underset{\mathrm{Spec}}{\mathcal{O}}^\mathrm{CS}\breve{\Delta}^+_{\psi,\Gamma,\square}/\mathrm{Fro}^\mathbb{Z},\\
&\underset{\mathrm{Spec}}{\mathcal{O}}^\mathrm{CS}\breve{\nabla}^+_{\psi,\Gamma,\square}/\mathrm{Fro}^\mathbb{Z},\underset{\mathrm{Spec}}{\mathcal{O}}^\mathrm{CS}\breve{\Delta}^\dagger_{\psi,\Gamma,\square}/\mathrm{Fro}^\mathbb{Z},\underset{\mathrm{Spec}}{\mathcal{O}}^\mathrm{CS}\breve{\nabla}^\dagger_{\psi,\Gamma,\square}/\mathrm{Fro}^\mathbb{Z},	\\
\end{align}
\begin{align}
&\underset{\mathrm{Spec}}{\mathcal{O}}^\mathrm{CS}{\Delta}_{\psi,\Gamma,\square}/\mathrm{Fro}^\mathbb{Z},\underset{\mathrm{Spec}}{\mathcal{O}}^\mathrm{CS}{\nabla}_{\psi,\Gamma,\square}/\mathrm{Fro}^\mathbb{Z},\underset{\mathrm{Spec}}{\mathcal{O}}^\mathrm{CS}{\Phi}_{\psi,\Gamma,\square}/\mathrm{Fro}^\mathbb{Z},\underset{\mathrm{Spec}}{\mathcal{O}}^\mathrm{CS}{\Delta}^+_{\psi,\Gamma,\square}/\mathrm{Fro}^\mathbb{Z},\\
&\underset{\mathrm{Spec}}{\mathcal{O}}^\mathrm{CS}{\nabla}^+_{\psi,\Gamma,\square}/\mathrm{Fro}^\mathbb{Z}, \underset{\mathrm{Spec}}{\mathcal{O}}^\mathrm{CS}{\Delta}^\dagger_{\psi,\Gamma,\square}/\mathrm{Fro}^\mathbb{Z},\underset{\mathrm{Spec}}{\mathcal{O}}^\mathrm{CS}{\nabla}^\dagger_{\psi,\Gamma,\square}/\mathrm{Fro}^\mathbb{Z}.	
\end{align}
Here for those space with notations related to the radius and the corresponding interval we consider the total unions $\bigcap_r,\bigcup_I$ in order to achieve the whole spaces to achieve the analogues of the corresponding FF curves from \cite{10KL1}, \cite{10KL2}, \cite{10FF} for
\[
\xymatrix@R+0pc@C+0pc{
\underset{r}{\mathrm{homotopycolimit}}~\underset{\mathrm{Spec}}{\mathcal{O}}^\mathrm{CS}\widetilde{\Phi}^r_{\psi,\Gamma,\square},\underset{I}{\mathrm{homotopylimit}}~\underset{\mathrm{Spec}}{\mathcal{O}}^\mathrm{CS}\widetilde{\Phi}^I_{\psi,\Gamma,\square},	\\
}
\]
\[
\xymatrix@R+0pc@C+0pc{
\underset{r}{\mathrm{homotopycolimit}}~\underset{\mathrm{Spec}}{\mathcal{O}}^\mathrm{CS}\breve{\Phi}^r_{\psi,\Gamma,\square},\underset{I}{\mathrm{homotopylimit}}~\underset{\mathrm{Spec}}{\mathcal{O}}^\mathrm{CS}\breve{\Phi}^I_{\psi,\Gamma,\square},	\\
}
\]
\[
\xymatrix@R+0pc@C+0pc{
\underset{r}{\mathrm{homotopycolimit}}~\underset{\mathrm{Spec}}{\mathcal{O}}^\mathrm{CS}{\Phi}^r_{\psi,\Gamma,\square},\underset{I}{\mathrm{homotopylimit}}~\underset{\mathrm{Spec}}{\mathcal{O}}^\mathrm{CS}{\Phi}^I_{\psi,\Gamma,\square}.	
}
\]
\[ 
\xymatrix@R+0pc@C+0pc{
\underset{r}{\mathrm{homotopycolimit}}~\underset{\mathrm{Spec}}{\mathcal{O}}^\mathrm{CS}\widetilde{\Phi}^r_{\psi,\Gamma,\square}/\mathrm{Fro}^\mathbb{Z},\underset{I}{\mathrm{homotopylimit}}~\underset{\mathrm{Spec}}{\mathcal{O}}^\mathrm{CS}\widetilde{\Phi}^I_{\psi,\Gamma,\square}/\mathrm{Fro}^\mathbb{Z},	\\
}
\]
\[ 
\xymatrix@R+0pc@C+0pc{
\underset{r}{\mathrm{homotopycolimit}}~\underset{\mathrm{Spec}}{\mathcal{O}}^\mathrm{CS}\breve{\Phi}^r_{\psi,\Gamma,\square}/\mathrm{Fro}^\mathbb{Z},\underset{I}{\mathrm{homotopylimit}}~\breve{\Phi}^I_{\psi,\Gamma,\square}/\mathrm{Fro}^\mathbb{Z},	\\
}
\]
\[ 
\xymatrix@R+0pc@C+0pc{
\underset{r}{\mathrm{homotopycolimit}}~\underset{\mathrm{Spec}}{\mathcal{O}}^\mathrm{CS}{\Phi}^r_{\psi,\Gamma,\square}/\mathrm{Fro}^\mathbb{Z},\underset{I}{\mathrm{homotopylimit}}~\underset{\mathrm{Spec}}{\mathcal{O}}^\mathrm{CS}{\Phi}^I_{\psi,\Gamma,\square}/\mathrm{Fro}^\mathbb{Z}.	
}
\]

\end{definition}

\begin{proposition}
There is a well-defined functor from the $\infty$-category 
\begin{align}
\mathrm{Quasicoherentpresheaves,Perfectcomplex,Condensed}_{*}	
\end{align}
where $*$ is one of the following spaces:
\begin{align}
&\underset{\mathrm{Spec}}{\mathcal{O}}^\mathrm{CS}\widetilde{\Phi}_{\psi,\Gamma,\square}/\mathrm{Fro}^\mathbb{Z},	\\
\end{align}
\begin{align}
&\underset{\mathrm{Spec}}{\mathcal{O}}^\mathrm{CS}\breve{\Phi}_{\psi,\Gamma,\square}/\mathrm{Fro}^\mathbb{Z},	\\
\end{align}
\begin{align}
&\underset{\mathrm{Spec}}{\mathcal{O}}^\mathrm{CS}{\Phi}_{\psi,\Gamma,\square}/\mathrm{Fro}^\mathbb{Z},	
\end{align}
to the $\infty$-category of $\mathrm{Fro}$-equivariant quasicoherent presheaves over similar spaces above correspondingly without the $\mathrm{Fro}$-quotients, and to the $\infty$-category of $\mathrm{Fro}$-equivariant quasicoherent modules over global sections of the structure $\infty$-sheaves of the similar spaces above correspondingly without the $\mathrm{Fro}$-quotients. Here for those space without notation related to the radius and the corresponding interval we consider the total unions $\bigcap_r,\bigcup_I$ in order to achieve the whole spaces to achieve the analogues of the corresponding FF curves from \cite{10KL1}, \cite{10KL2}, \cite{10FF} for
\[
\xymatrix@R+0pc@C+0pc{
\underset{r}{\mathrm{homotopycolimit}}~\underset{\mathrm{Spec}}{\mathcal{O}}^\mathrm{CS}\widetilde{\Phi}^r_{\psi,\Gamma,\square},\underset{I}{\mathrm{homotopylimit}}~\underset{\mathrm{Spec}}{\mathcal{O}}^\mathrm{CS}\widetilde{\Phi}^I_{\psi,\Gamma,\square},	\\
}
\]
\[
\xymatrix@R+0pc@C+0pc{
\underset{r}{\mathrm{homotopycolimit}}~\underset{\mathrm{Spec}}{\mathcal{O}}^\mathrm{CS}\breve{\Phi}^r_{\psi,\Gamma,\square},\underset{I}{\mathrm{homotopylimit}}~\underset{\mathrm{Spec}}{\mathcal{O}}^\mathrm{CS}\breve{\Phi}^I_{\psi,\Gamma,\square},	\\
}
\]
\[
\xymatrix@R+0pc@C+0pc{
\underset{r}{\mathrm{homotopycolimit}}~\underset{\mathrm{Spec}}{\mathcal{O}}^\mathrm{CS}{\Phi}^r_{\psi,\Gamma,\square},\underset{I}{\mathrm{homotopylimit}}~\underset{\mathrm{Spec}}{\mathcal{O}}^\mathrm{CS}{\Phi}^I_{\psi,\Gamma,\square}.	
}
\]
\[ 
\xymatrix@R+0pc@C+0pc{
\underset{r}{\mathrm{homotopycolimit}}~\underset{\mathrm{Spec}}{\mathcal{O}}^\mathrm{CS}\widetilde{\Phi}^r_{\psi,\Gamma,\square}/\mathrm{Fro}^\mathbb{Z},\underset{I}{\mathrm{homotopylimit}}~\underset{\mathrm{Spec}}{\mathcal{O}}^\mathrm{CS}\widetilde{\Phi}^I_{\psi,\Gamma,\square}/\mathrm{Fro}^\mathbb{Z},	\\
}
\]
\[ 
\xymatrix@R+0pc@C+0pc{
\underset{r}{\mathrm{homotopycolimit}}~\underset{\mathrm{Spec}}{\mathcal{O}}^\mathrm{CS}\breve{\Phi}^r_{\psi,\Gamma,\square}/\mathrm{Fro}^\mathbb{Z},\underset{I}{\mathrm{homotopylimit}}~\breve{\Phi}^I_{\psi,\Gamma,\square}/\mathrm{Fro}^\mathbb{Z},	\\
}
\]
\[ 
\xymatrix@R+0pc@C+0pc{
\underset{r}{\mathrm{homotopycolimit}}~\underset{\mathrm{Spec}}{\mathcal{O}}^\mathrm{CS}{\Phi}^r_{\psi,\Gamma,\square}/\mathrm{Fro}^\mathbb{Z},\underset{I}{\mathrm{homotopylimit}}~\underset{\mathrm{Spec}}{\mathcal{O}}^\mathrm{CS}{\Phi}^I_{\psi,\Gamma,\square}/\mathrm{Fro}^\mathbb{Z}.	
}
\]	
In this situation we will have the target category being family parametrized by $r$ or $I$ in compatible glueing sense as in \cite[Definition 5.4.10]{10KL2}. In this situation for modules parametrized by the intervals we have the equivalence of $\infty$-categories by using \cite[Proposition 12.18]{10CS2}. Here the corresponding quasicoherent Frobenius modules are defined to be the corresponding homotopy colimits and limits of Frobenius modules:
\begin{align}
\underset{r}{\mathrm{homotopycolimit}}~M_r,\\
\underset{I}{\mathrm{homotopylimit}}~M_I,	
\end{align}
where each $M_r$ is a Frobenius-equivariant module over the period ring with respect to some radius $r$ while each $M_I$ is a Frobenius-equivariant module over the period ring with respect to some interval $I$.\\
\end{proposition}

\begin{proposition}
Similar proposition holds for 
\begin{align}
\mathrm{Quasicoherentsheaves,Perfectcomplex,IndBanach}_{*}.	
\end{align}	
\end{proposition}

\section{Univariate Hodge Iwasawa Modules by Deformation}

This chapter follows closely \cite{10T1}, \cite{10T2}, \cite{10T3}, \cite{10T4}, \cite{10T5}, \cite{10T6}, \cite{10KPX}, \cite{10KP}, \cite{10KL1}, \cite{10KL2}, \cite{10BK}, \cite{10BBBK}, \cite{10BBM}, \cite{10KKM}, \cite{10CS1}, \cite{10CS2}, \cite{10LBV}.

\subsection{Frobenius Quasicoherent Modules I}

\begin{definition}
First we consider the Bambozzi-Kremnizer spectrum $\underset{\mathrm{Spec}}{\mathcal{O}}^\mathrm{BK}(*)$ attached to any of those in the above from \cite{10BK} by taking derived rational localization:
\begin{align}
&\underset{\mathrm{Spec}}{\mathcal{O}}^\mathrm{BK}\widetilde{\Phi}_{\psi,A},\underset{\mathrm{Spec}}{\mathcal{O}}^\mathrm{BK}\widetilde{\Phi}^r_{\psi,A},\underset{\mathrm{Spec}}{\mathcal{O}}^\mathrm{BK}\widetilde{\Phi}^I_{\psi,A},	
\end{align}
\begin{align}
&\underset{\mathrm{Spec}}{\mathcal{O}}^\mathrm{BK}\breve{\Phi}_{\psi,A},\underset{\mathrm{Spec}}{\mathcal{O}}^\mathrm{BK}\breve{\Phi}^r_{\psi,A},\underset{\mathrm{Spec}}{\mathcal{O}}^\mathrm{BK}\breve{\Phi}^I_{\psi,A},	
\end{align}
\begin{align}
&\underset{\mathrm{Spec}}{\mathcal{O}}^\mathrm{BK}{\Phi}_{\psi,A},
\underset{\mathrm{Spec}}{\mathcal{O}}^\mathrm{BK}{\Phi}^r_{\psi,A},\underset{\mathrm{Spec}}{\mathcal{O}}^\mathrm{BK}{\Phi}^I_{\psi,A}.	
\end{align}

Then we take the corresponding quotients by using the corresponding Frobenius operators:
\begin{align}
&\underset{\mathrm{Spec}}{\mathcal{O}}^\mathrm{BK}\widetilde{\Phi}_{\psi,A}/\mathrm{Fro}^\mathbb{Z},	\\
\end{align}
\begin{align}
&\underset{\mathrm{Spec}}{\mathcal{O}}^\mathrm{BK}\breve{\Phi}_{\psi,A}/\mathrm{Fro}^\mathbb{Z},	\\
\end{align}
\begin{align}
&\underset{\mathrm{Spec}}{\mathcal{O}}^\mathrm{BK}{\Phi}_{\psi,A}/\mathrm{Fro}^\mathbb{Z}.	
\end{align}
Here for those space without notation related to the radius and the corresponding interval we consider the total unions $\bigcap_r,\bigcup_I$ in order to achieve the whole spaces to achieve the analogues of the corresponding FF curves from \cite{10KL1}, \cite{10KL2}, \cite{10FF} for
\[
\xymatrix@R+0pc@C+0pc{
\underset{r}{\mathrm{homotopycolimit}}~\underset{\mathrm{Spec}}{\mathcal{O}}^\mathrm{BK}\widetilde{\Phi}^r_{\psi,A},\underset{I}{\mathrm{homotopylimit}}~\underset{\mathrm{Spec}}{\mathcal{O}}^\mathrm{BK}\widetilde{\Phi}^I_{\psi,A},	\\
}
\]
\[
\xymatrix@R+0pc@C+0pc{
\underset{r}{\mathrm{homotopycolimit}}~\underset{\mathrm{Spec}}{\mathcal{O}}^\mathrm{BK}\breve{\Phi}^r_{\psi,A},\underset{I}{\mathrm{homotopylimit}}~\underset{\mathrm{Spec}}{\mathcal{O}}^\mathrm{BK}\breve{\Phi}^I_{\psi,A},	\\
}
\]
\[
\xymatrix@R+0pc@C+0pc{
\underset{r}{\mathrm{homotopycolimit}}~\underset{\mathrm{Spec}}{\mathcal{O}}^\mathrm{BK}{\Phi}^r_{\psi,A},\underset{I}{\mathrm{homotopylimit}}~\underset{\mathrm{Spec}}{\mathcal{O}}^\mathrm{BK}{\Phi}^I_{\psi,A}.	
}
\]
\[  
\xymatrix@R+0pc@C+0pc{
\underset{r}{\mathrm{homotopycolimit}}~\underset{\mathrm{Spec}}{\mathcal{O}}^\mathrm{BK}\widetilde{\Phi}^r_{\psi,A}/\mathrm{Fro}^\mathbb{Z},\underset{I}{\mathrm{homotopylimit}}~\underset{\mathrm{Spec}}{\mathcal{O}}^\mathrm{BK}\widetilde{\Phi}^I_{\psi,A}/\mathrm{Fro}^\mathbb{Z},	\\
}
\]
\[ 
\xymatrix@R+0pc@C+0pc{
\underset{r}{\mathrm{homotopycolimit}}~\underset{\mathrm{Spec}}{\mathcal{O}}^\mathrm{BK}\breve{\Phi}^r_{\psi,A}/\mathrm{Fro}^\mathbb{Z},\underset{I}{\mathrm{homotopylimit}}~\underset{\mathrm{Spec}}{\mathcal{O}}^\mathrm{BK}\breve{\Phi}^I_{\psi,A}/\mathrm{Fro}^\mathbb{Z},	\\
}
\]
\[ 
\xymatrix@R+0pc@C+0pc{
\underset{r}{\mathrm{homotopycolimit}}~\underset{\mathrm{Spec}}{\mathcal{O}}^\mathrm{BK}{\Phi}^r_{\psi,A}/\mathrm{Fro}^\mathbb{Z},\underset{I}{\mathrm{homotopylimit}}~\underset{\mathrm{Spec}}{\mathcal{O}}^\mathrm{BK}{\Phi}^I_{\psi,A}/\mathrm{Fro}^\mathbb{Z}.	
}
\]

\end{definition}

\indent Meanwhile we have the corresponding Clausen-Scholze analytic stacks from \cite{10CS2}, therefore applying their construction we have:

\begin{definition}
Here we define the following products by using the solidified tensor product from \cite{10CS1} and \cite{10CS2}. Namely $A$ will still as above as a Banach ring over $\mathbb{Q}_p$. Then we take solidified tensor product $\overset{\blacksquare}{\otimes}$ of any of the following
\[
\xymatrix@R+0pc@C+0pc{
\widetilde{\Delta}_{\psi},\widetilde{\nabla}_{\psi},\widetilde{\Phi}_{\psi},\widetilde{\Delta}^+_{\psi},\widetilde{\nabla}^+_{\psi},\widetilde{\Delta}^\dagger_{\psi},\widetilde{\nabla}^\dagger_{\psi},\widetilde{\Phi}^r_{\psi},\widetilde{\Phi}^I_{\psi}, 
}
\]

\[
\xymatrix@R+0pc@C+0pc{
\breve{\Delta}_{\psi},\breve{\nabla}_{\psi},\breve{\Phi}_{\psi},\breve{\Delta}^+_{\psi},\breve{\nabla}^+_{\psi},\breve{\Delta}^\dagger_{\psi},\breve{\nabla}^\dagger_{\psi},\breve{\Phi}^r_{\psi},\breve{\Phi}^I_{\psi},	
}
\]

\[
\xymatrix@R+0pc@C+0pc{
{\Delta}_{\psi},{\nabla}_{\psi},{\Phi}_{\psi},{\Delta}^+_{\psi},{\nabla}^+_{\psi},{\Delta}^\dagger_{\psi},{\nabla}^\dagger_{\psi},{\Phi}^r_{\psi},{\Phi}^I_{\psi},	
}
\]  	
with $A$. Then we have the notations:
\[
\xymatrix@R+0pc@C+0pc{
\widetilde{\Delta}_{\psi,A},\widetilde{\nabla}_{\psi,A},\widetilde{\Phi}_{\psi,A},\widetilde{\Delta}^+_{\psi,A},\widetilde{\nabla}^+_{\psi,A},\widetilde{\Delta}^\dagger_{\psi,A},\widetilde{\nabla}^\dagger_{\psi,A},\widetilde{\Phi}^r_{\psi,A},\widetilde{\Phi}^I_{\psi,A}, 
}
\]

\[
\xymatrix@R+0pc@C+0pc{
\breve{\Delta}_{\psi,A},\breve{\nabla}_{\psi,A},\breve{\Phi}_{\psi,A},\breve{\Delta}^+_{\psi,A},\breve{\nabla}^+_{\psi,A},\breve{\Delta}^\dagger_{\psi,A},\breve{\nabla}^\dagger_{\psi,A},\breve{\Phi}^r_{\psi,A},\breve{\Phi}^I_{\psi,A},	
}
\]

\[
\xymatrix@R+0pc@C+0pc{
{\Delta}_{\psi,A},{\nabla}_{\psi,A},{\Phi}_{\psi,A},{\Delta}^+_{\psi,A},{\nabla}^+_{\psi,A},{\Delta}^\dagger_{\psi,A},{\nabla}^\dagger_{\psi,A},{\Phi}^r_{\psi,A},{\Phi}^I_{\psi,A}.	
}
\]
\end{definition}

\begin{definition}
First we consider the Clausen-Scholze spectrum $\underset{\mathrm{Spec}}{\mathcal{O}}^\mathrm{CS}(*)$ attached to any of those in the above from \cite{10CS2} by taking derived rational localization:
\begin{align}
\underset{\mathrm{Spec}}{\mathcal{O}}^\mathrm{CS}\widetilde{\Delta}_{\psi,A},\underset{\mathrm{Spec}}{\mathcal{O}}^\mathrm{CS}\widetilde{\nabla}_{\psi,A},\underset{\mathrm{Spec}}{\mathcal{O}}^\mathrm{CS}\widetilde{\Phi}_{\psi,A},\underset{\mathrm{Spec}}{\mathcal{O}}^\mathrm{CS}\widetilde{\Delta}^+_{\psi,A},\underset{\mathrm{Spec}}{\mathcal{O}}^\mathrm{CS}\widetilde{\nabla}^+_{\psi,A},\\
\underset{\mathrm{Spec}}{\mathcal{O}}^\mathrm{CS}\widetilde{\Delta}^\dagger_{\psi,A},\underset{\mathrm{Spec}}{\mathcal{O}}^\mathrm{CS}\widetilde{\nabla}^\dagger_{\psi,A},\underset{\mathrm{Spec}}{\mathcal{O}}^\mathrm{CS}\widetilde{\Phi}^r_{\psi,A},\underset{\mathrm{Spec}}{\mathcal{O}}^\mathrm{CS}\widetilde{\Phi}^I_{\psi,A},	\\
\end{align}
\begin{align}
\underset{\mathrm{Spec}}{\mathcal{O}}^\mathrm{CS}\breve{\Delta}_{\psi,A},\breve{\nabla}_{\psi,A},\underset{\mathrm{Spec}}{\mathcal{O}}^\mathrm{CS}\breve{\Phi}_{\psi,A},\underset{\mathrm{Spec}}{\mathcal{O}}^\mathrm{CS}\breve{\Delta}^+_{\psi,A},\underset{\mathrm{Spec}}{\mathcal{O}}^\mathrm{CS}\breve{\nabla}^+_{\psi,A},\\
\underset{\mathrm{Spec}}{\mathcal{O}}^\mathrm{CS}\breve{\Delta}^\dagger_{\psi,A},\underset{\mathrm{Spec}}{\mathcal{O}}^\mathrm{CS}\breve{\nabla}^\dagger_{\psi,A},\underset{\mathrm{Spec}}{\mathcal{O}}^\mathrm{CS}\breve{\Phi}^r_{\psi,A},\breve{\Phi}^I_{\psi,A},	\\
\end{align}
\begin{align}
\underset{\mathrm{Spec}}{\mathcal{O}}^\mathrm{CS}{\Delta}_{\psi,A},\underset{\mathrm{Spec}}{\mathcal{O}}^\mathrm{CS}{\nabla}_{\psi,A},\underset{\mathrm{Spec}}{\mathcal{O}}^\mathrm{CS}{\Phi}_{\psi,A},\underset{\mathrm{Spec}}{\mathcal{O}}^\mathrm{CS}{\Delta}^+_{\psi,A},\underset{\mathrm{Spec}}{\mathcal{O}}^\mathrm{CS}{\nabla}^+_{\psi,A},\\
\underset{\mathrm{Spec}}{\mathcal{O}}^\mathrm{CS}{\Delta}^\dagger_{\psi,A},\underset{\mathrm{Spec}}{\mathcal{O}}^\mathrm{CS}{\nabla}^\dagger_{\psi,A},\underset{\mathrm{Spec}}{\mathcal{O}}^\mathrm{CS}{\Phi}^r_{\psi,A},\underset{\mathrm{Spec}}{\mathcal{O}}^\mathrm{CS}{\Phi}^I_{\psi,A}.	
\end{align}

Then we take the corresponding quotients by using the corresponding Frobenius operators:
\begin{align}
&\underset{\mathrm{Spec}}{\mathcal{O}}^\mathrm{CS}\widetilde{\Delta}_{\psi,A}/\mathrm{Fro}^\mathbb{Z},\underset{\mathrm{Spec}}{\mathcal{O}}^\mathrm{CS}\widetilde{\nabla}_{\psi,A}/\mathrm{Fro}^\mathbb{Z},\underset{\mathrm{Spec}}{\mathcal{O}}^\mathrm{CS}\widetilde{\Phi}_{\psi,A}/\mathrm{Fro}^\mathbb{Z},\underset{\mathrm{Spec}}{\mathcal{O}}^\mathrm{CS}\widetilde{\Delta}^+_{\psi,A}/\mathrm{Fro}^\mathbb{Z},\\
&\underset{\mathrm{Spec}}{\mathcal{O}}^\mathrm{CS}\widetilde{\nabla}^+_{\psi,A}/\mathrm{Fro}^\mathbb{Z}, \underset{\mathrm{Spec}}{\mathcal{O}}^\mathrm{CS}\widetilde{\Delta}^\dagger_{\psi,A}/\mathrm{Fro}^\mathbb{Z},\underset{\mathrm{Spec}}{\mathcal{O}}^\mathrm{CS}\widetilde{\nabla}^\dagger_{\psi,A}/\mathrm{Fro}^\mathbb{Z},	\\
\end{align}
\begin{align}
&\underset{\mathrm{Spec}}{\mathcal{O}}^\mathrm{CS}\breve{\Delta}_{\psi,A}/\mathrm{Fro}^\mathbb{Z},\breve{\nabla}_{\psi,A}/\mathrm{Fro}^\mathbb{Z},\underset{\mathrm{Spec}}{\mathcal{O}}^\mathrm{CS}\breve{\Phi}_{\psi,A}/\mathrm{Fro}^\mathbb{Z},\underset{\mathrm{Spec}}{\mathcal{O}}^\mathrm{CS}\breve{\Delta}^+_{\psi,A}/\mathrm{Fro}^\mathbb{Z},\\
&\underset{\mathrm{Spec}}{\mathcal{O}}^\mathrm{CS}\breve{\nabla}^+_{\psi,A}/\mathrm{Fro}^\mathbb{Z}, \underset{\mathrm{Spec}}{\mathcal{O}}^\mathrm{CS}\breve{\Delta}^\dagger_{\psi,A}/\mathrm{Fro}^\mathbb{Z},\underset{\mathrm{Spec}}{\mathcal{O}}^\mathrm{CS}\breve{\nabla}^\dagger_{\psi,A}/\mathrm{Fro}^\mathbb{Z},	\\
\end{align}
\begin{align}
&\underset{\mathrm{Spec}}{\mathcal{O}}^\mathrm{CS}{\Delta}_{\psi,A}/\mathrm{Fro}^\mathbb{Z},\underset{\mathrm{Spec}}{\mathcal{O}}^\mathrm{CS}{\nabla}_{\psi,A}/\mathrm{Fro}^\mathbb{Z},\underset{\mathrm{Spec}}{\mathcal{O}}^\mathrm{CS}{\Phi}_{\psi,A}/\mathrm{Fro}^\mathbb{Z},\underset{\mathrm{Spec}}{\mathcal{O}}^\mathrm{CS}{\Delta}^+_{\psi,A}/\mathrm{Fro}^\mathbb{Z},\\
&\underset{\mathrm{Spec}}{\mathcal{O}}^\mathrm{CS}{\nabla}^+_{\psi,A}/\mathrm{Fro}^\mathbb{Z}, \underset{\mathrm{Spec}}{\mathcal{O}}^\mathrm{CS}{\Delta}^\dagger_{\psi,A}/\mathrm{Fro}^\mathbb{Z},\underset{\mathrm{Spec}}{\mathcal{O}}^\mathrm{CS}{\nabla}^\dagger_{\psi,A}/\mathrm{Fro}^\mathbb{Z}.	
\end{align}
Here for those space with notations related to the radius and the corresponding interval we consider the total unions $\bigcap_r,\bigcup_I$ in order to achieve the whole spaces to achieve the analogues of the corresponding FF curves from \cite{10KL1}, \cite{10KL2}, \cite{10FF} for
\[
\xymatrix@R+0pc@C+0pc{
\underset{r}{\mathrm{homotopycolimit}}~\underset{\mathrm{Spec}}{\mathcal{O}}^\mathrm{CS}\widetilde{\Phi}^r_{\psi,A},\underset{I}{\mathrm{homotopylimit}}~\underset{\mathrm{Spec}}{\mathcal{O}}^\mathrm{CS}\widetilde{\Phi}^I_{\psi,A},	\\
}
\]
\[
\xymatrix@R+0pc@C+0pc{
\underset{r}{\mathrm{homotopycolimit}}~\underset{\mathrm{Spec}}{\mathcal{O}}^\mathrm{CS}\breve{\Phi}^r_{\psi,A},\underset{I}{\mathrm{homotopylimit}}~\underset{\mathrm{Spec}}{\mathcal{O}}^\mathrm{CS}\breve{\Phi}^I_{\psi,A},	\\
}
\]
\[
\xymatrix@R+0pc@C+0pc{
\underset{r}{\mathrm{homotopycolimit}}~\underset{\mathrm{Spec}}{\mathcal{O}}^\mathrm{CS}{\Phi}^r_{\psi,A},\underset{I}{\mathrm{homotopylimit}}~\underset{\mathrm{Spec}}{\mathcal{O}}^\mathrm{CS}{\Phi}^I_{\psi,A}.	
}
\]
\[ 
\xymatrix@R+0pc@C+0pc{
\underset{r}{\mathrm{homotopycolimit}}~\underset{\mathrm{Spec}}{\mathcal{O}}^\mathrm{CS}\widetilde{\Phi}^r_{\psi,A}/\mathrm{Fro}^\mathbb{Z},\underset{I}{\mathrm{homotopylimit}}~\underset{\mathrm{Spec}}{\mathcal{O}}^\mathrm{CS}\widetilde{\Phi}^I_{\psi,A}/\mathrm{Fro}^\mathbb{Z},	\\
}
\]
\[ 
\xymatrix@R+0pc@C+0pc{
\underset{r}{\mathrm{homotopycolimit}}~\underset{\mathrm{Spec}}{\mathcal{O}}^\mathrm{CS}\breve{\Phi}^r_{\psi,A}/\mathrm{Fro}^\mathbb{Z},\underset{I}{\mathrm{homotopylimit}}~\breve{\Phi}^I_{\psi,A}/\mathrm{Fro}^\mathbb{Z},	\\
}
\]
\[ 
\xymatrix@R+0pc@C+0pc{
\underset{r}{\mathrm{homotopycolimit}}~\underset{\mathrm{Spec}}{\mathcal{O}}^\mathrm{CS}{\Phi}^r_{\psi,A}/\mathrm{Fro}^\mathbb{Z},\underset{I}{\mathrm{homotopylimit}}~\underset{\mathrm{Spec}}{\mathcal{O}}^\mathrm{CS}{\Phi}^I_{\psi,A}/\mathrm{Fro}^\mathbb{Z}.	
}
\]

\end{definition}

\

\begin{definition}
We then consider the corresponding quasipresheaves of the corresponding ind-Banach or monomorphic ind-Banach modules from \cite{10BBK}, \cite{10KKM}:
\begin{align}
\mathrm{Quasicoherentpresheaves,IndBanach}_{*}	
\end{align}
where $*$ is one of the following spaces:
\begin{align}
&\underset{\mathrm{Spec}}{\mathcal{O}}^\mathrm{BK}\widetilde{\Phi}_{\psi,A}/\mathrm{Fro}^\mathbb{Z},	\\
\end{align}
\begin{align}
&\underset{\mathrm{Spec}}{\mathcal{O}}^\mathrm{BK}\breve{\Phi}_{\psi,A}/\mathrm{Fro}^\mathbb{Z},	\\
\end{align}
\begin{align}
&\underset{\mathrm{Spec}}{\mathcal{O}}^\mathrm{BK}{\Phi}_{\psi,A}/\mathrm{Fro}^\mathbb{Z}.	
\end{align}
Here for those space without notation related to the radius and the corresponding interval we consider the total unions $\bigcap_r,\bigcup_I$ in order to achieve the whole spaces to achieve the analogues of the corresponding FF curves from \cite{10KL1}, \cite{10KL2}, \cite{10FF} for
\[
\xymatrix@R+0pc@C+0pc{
\underset{r}{\mathrm{homotopycolimit}}~\underset{\mathrm{Spec}}{\mathcal{O}}^\mathrm{BK}\widetilde{\Phi}^r_{\psi,A},\underset{I}{\mathrm{homotopylimit}}~\underset{\mathrm{Spec}}{\mathcal{O}}^\mathrm{BK}\widetilde{\Phi}^I_{\psi,A},	\\
}
\]
\[
\xymatrix@R+0pc@C+0pc{
\underset{r}{\mathrm{homotopycolimit}}~\underset{\mathrm{Spec}}{\mathcal{O}}^\mathrm{BK}\breve{\Phi}^r_{\psi,A},\underset{I}{\mathrm{homotopylimit}}~\underset{\mathrm{Spec}}{\mathcal{O}}^\mathrm{BK}\breve{\Phi}^I_{\psi,A},	\\
}
\]
\[
\xymatrix@R+0pc@C+0pc{
\underset{r}{\mathrm{homotopycolimit}}~\underset{\mathrm{Spec}}{\mathcal{O}}^\mathrm{BK}{\Phi}^r_{\psi,A},\underset{I}{\mathrm{homotopylimit}}~\underset{\mathrm{Spec}}{\mathcal{O}}^\mathrm{BK}{\Phi}^I_{\psi,A}.	
}
\]
\[  
\xymatrix@R+0pc@C+0pc{
\underset{r}{\mathrm{homotopycolimit}}~\underset{\mathrm{Spec}}{\mathcal{O}}^\mathrm{BK}\widetilde{\Phi}^r_{\psi,A}/\mathrm{Fro}^\mathbb{Z},\underset{I}{\mathrm{homotopylimit}}~\underset{\mathrm{Spec}}{\mathcal{O}}^\mathrm{BK}\widetilde{\Phi}^I_{\psi,A}/\mathrm{Fro}^\mathbb{Z},	\\
}
\]
\[ 
\xymatrix@R+0pc@C+0pc{
\underset{r}{\mathrm{homotopycolimit}}~\underset{\mathrm{Spec}}{\mathcal{O}}^\mathrm{BK}\breve{\Phi}^r_{\psi,A}/\mathrm{Fro}^\mathbb{Z},\underset{I}{\mathrm{homotopylimit}}~\underset{\mathrm{Spec}}{\mathcal{O}}^\mathrm{BK}\breve{\Phi}^I_{\psi,A}/\mathrm{Fro}^\mathbb{Z},	\\
}
\]
\[ 
\xymatrix@R+0pc@C+0pc{
\underset{r}{\mathrm{homotopycolimit}}~\underset{\mathrm{Spec}}{\mathcal{O}}^\mathrm{BK}{\Phi}^r_{\psi,A}/\mathrm{Fro}^\mathbb{Z},\underset{I}{\mathrm{homotopylimit}}~\underset{\mathrm{Spec}}{\mathcal{O}}^\mathrm{BK}{\Phi}^I_{\psi,A}/\mathrm{Fro}^\mathbb{Z}.	
}
\]

\end{definition}

\begin{definition}
We then consider the corresponding quasisheaves of the corresponding condensed solid topological modules from \cite{10CS2}:
\begin{align}
\mathrm{Quasicoherentsheaves, Condensed}_{*}	
\end{align}
where $*$ is one of the following spaces:
\begin{align}
&\underset{\mathrm{Spec}}{\mathcal{O}}^\mathrm{CS}\widetilde{\Delta}_{\psi,A}/\mathrm{Fro}^\mathbb{Z},\underset{\mathrm{Spec}}{\mathcal{O}}^\mathrm{CS}\widetilde{\nabla}_{\psi,A}/\mathrm{Fro}^\mathbb{Z},\underset{\mathrm{Spec}}{\mathcal{O}}^\mathrm{CS}\widetilde{\Phi}_{\psi,A}/\mathrm{Fro}^\mathbb{Z},\underset{\mathrm{Spec}}{\mathcal{O}}^\mathrm{CS}\widetilde{\Delta}^+_{\psi,A}/\mathrm{Fro}^\mathbb{Z},\\
&\underset{\mathrm{Spec}}{\mathcal{O}}^\mathrm{CS}\widetilde{\nabla}^+_{\psi,A}/\mathrm{Fro}^\mathbb{Z},\underset{\mathrm{Spec}}{\mathcal{O}}^\mathrm{CS}\widetilde{\Delta}^\dagger_{\psi,A}/\mathrm{Fro}^\mathbb{Z},\underset{\mathrm{Spec}}{\mathcal{O}}^\mathrm{CS}\widetilde{\nabla}^\dagger_{\psi,A}/\mathrm{Fro}^\mathbb{Z},	\\
\end{align}
\begin{align}
&\underset{\mathrm{Spec}}{\mathcal{O}}^\mathrm{CS}\breve{\Delta}_{\psi,A}/\mathrm{Fro}^\mathbb{Z},\breve{\nabla}_{\psi,A}/\mathrm{Fro}^\mathbb{Z},\underset{\mathrm{Spec}}{\mathcal{O}}^\mathrm{CS}\breve{\Phi}_{\psi,A}/\mathrm{Fro}^\mathbb{Z},\underset{\mathrm{Spec}}{\mathcal{O}}^\mathrm{CS}\breve{\Delta}^+_{\psi,A}/\mathrm{Fro}^\mathbb{Z},\\
&\underset{\mathrm{Spec}}{\mathcal{O}}^\mathrm{CS}\breve{\nabla}^+_{\psi,A}/\mathrm{Fro}^\mathbb{Z},\underset{\mathrm{Spec}}{\mathcal{O}}^\mathrm{CS}\breve{\Delta}^\dagger_{\psi,A}/\mathrm{Fro}^\mathbb{Z},\underset{\mathrm{Spec}}{\mathcal{O}}^\mathrm{CS}\breve{\nabla}^\dagger_{\psi,A}/\mathrm{Fro}^\mathbb{Z},	\\
\end{align}
\begin{align}
&\underset{\mathrm{Spec}}{\mathcal{O}}^\mathrm{CS}{\Delta}_{\psi,A}/\mathrm{Fro}^\mathbb{Z},\underset{\mathrm{Spec}}{\mathcal{O}}^\mathrm{CS}{\nabla}_{\psi,A}/\mathrm{Fro}^\mathbb{Z},\underset{\mathrm{Spec}}{\mathcal{O}}^\mathrm{CS}{\Phi}_{\psi,A}/\mathrm{Fro}^\mathbb{Z},\underset{\mathrm{Spec}}{\mathcal{O}}^\mathrm{CS}{\Delta}^+_{\psi,A}/\mathrm{Fro}^\mathbb{Z},\\
&\underset{\mathrm{Spec}}{\mathcal{O}}^\mathrm{CS}{\nabla}^+_{\psi,A}/\mathrm{Fro}^\mathbb{Z}, \underset{\mathrm{Spec}}{\mathcal{O}}^\mathrm{CS}{\Delta}^\dagger_{\psi,A}/\mathrm{Fro}^\mathbb{Z},\underset{\mathrm{Spec}}{\mathcal{O}}^\mathrm{CS}{\nabla}^\dagger_{\psi,A}/\mathrm{Fro}^\mathbb{Z}.	
\end{align}
Here for those space with notations related to the radius and the corresponding interval we consider the total unions $\bigcap_r,\bigcup_I$ in order to achieve the whole spaces to achieve the analogues of the corresponding FF curves from \cite{10KL1}, \cite{10KL2}, \cite{10FF} for
\[
\xymatrix@R+0pc@C+0pc{
\underset{r}{\mathrm{homotopycolimit}}~\underset{\mathrm{Spec}}{\mathcal{O}}^\mathrm{CS}\widetilde{\Phi}^r_{\psi,A},\underset{I}{\mathrm{homotopylimit}}~\underset{\mathrm{Spec}}{\mathcal{O}}^\mathrm{CS}\widetilde{\Phi}^I_{\psi,A},	\\
}
\]
\[
\xymatrix@R+0pc@C+0pc{
\underset{r}{\mathrm{homotopycolimit}}~\underset{\mathrm{Spec}}{\mathcal{O}}^\mathrm{CS}\breve{\Phi}^r_{\psi,A},\underset{I}{\mathrm{homotopylimit}}~\underset{\mathrm{Spec}}{\mathcal{O}}^\mathrm{CS}\breve{\Phi}^I_{\psi,A},	\\
}
\]
\[
\xymatrix@R+0pc@C+0pc{
\underset{r}{\mathrm{homotopycolimit}}~\underset{\mathrm{Spec}}{\mathcal{O}}^\mathrm{CS}{\Phi}^r_{\psi,A},\underset{I}{\mathrm{homotopylimit}}~\underset{\mathrm{Spec}}{\mathcal{O}}^\mathrm{CS}{\Phi}^I_{\psi,A}.	
}
\]
\[ 
\xymatrix@R+0pc@C+0pc{
\underset{r}{\mathrm{homotopycolimit}}~\underset{\mathrm{Spec}}{\mathcal{O}}^\mathrm{CS}\widetilde{\Phi}^r_{\psi,A}/\mathrm{Fro}^\mathbb{Z},\underset{I}{\mathrm{homotopylimit}}~\underset{\mathrm{Spec}}{\mathcal{O}}^\mathrm{CS}\widetilde{\Phi}^I_{\psi,A}/\mathrm{Fro}^\mathbb{Z},	\\
}
\]
\[ 
\xymatrix@R+0pc@C+0pc{
\underset{r}{\mathrm{homotopycolimit}}~\underset{\mathrm{Spec}}{\mathcal{O}}^\mathrm{CS}\breve{\Phi}^r_{\psi,A}/\mathrm{Fro}^\mathbb{Z},\underset{I}{\mathrm{homotopylimit}}~\breve{\Phi}^I_{\psi,A}/\mathrm{Fro}^\mathbb{Z},	\\
}
\]
\[ 
\xymatrix@R+0pc@C+0pc{
\underset{r}{\mathrm{homotopycolimit}}~\underset{\mathrm{Spec}}{\mathcal{O}}^\mathrm{CS}{\Phi}^r_{\psi,A}/\mathrm{Fro}^\mathbb{Z},\underset{I}{\mathrm{homotopylimit}}~\underset{\mathrm{Spec}}{\mathcal{O}}^\mathrm{CS}{\Phi}^I_{\psi,A}/\mathrm{Fro}^\mathbb{Z}.	
}
\]

\end{definition}

\

\begin{proposition}
There is a well-defined functor from the $\infty$-category 
\begin{align}
\mathrm{Quasicoherentpresheaves,Condensed}_{*}	
\end{align}
where $*$ is one of the following spaces:
\begin{align}
&\underset{\mathrm{Spec}}{\mathcal{O}}^\mathrm{CS}\widetilde{\Phi}_{\psi,A}/\mathrm{Fro}^\mathbb{Z},	\\
\end{align}
\begin{align}
&\underset{\mathrm{Spec}}{\mathcal{O}}^\mathrm{CS}\breve{\Phi}_{\psi,A}/\mathrm{Fro}^\mathbb{Z},	\\
\end{align}
\begin{align}
&\underset{\mathrm{Spec}}{\mathcal{O}}^\mathrm{CS}{\Phi}_{\psi,A}/\mathrm{Fro}^\mathbb{Z},	
\end{align}
to the $\infty$-category of $\mathrm{Fro}$-equivariant quasicoherent presheaves over similar spaces above correspondingly without the $\mathrm{Fro}$-quotients, and to the $\infty$-category of $\mathrm{Fro}$-equivariant quasicoherent modules over global sections of the structure $\infty$-sheaves of the similar spaces above correspondingly without the $\mathrm{Fro}$-quotients. Here for those space without notation related to the radius and the corresponding interval we consider the total unions $\bigcap_r,\bigcup_I$ in order to achieve the whole spaces to achieve the analogues of the corresponding FF curves from \cite{10KL1}, \cite{10KL2}, \cite{10FF} for
\[
\xymatrix@R+0pc@C+0pc{
\underset{r}{\mathrm{homotopycolimit}}~\underset{\mathrm{Spec}}{\mathcal{O}}^\mathrm{CS}\widetilde{\Phi}^r_{\psi,A},\underset{I}{\mathrm{homotopylimit}}~\underset{\mathrm{Spec}}{\mathcal{O}}^\mathrm{CS}\widetilde{\Phi}^I_{\psi,A},	\\
}
\]
\[
\xymatrix@R+0pc@C+0pc{
\underset{r}{\mathrm{homotopycolimit}}~\underset{\mathrm{Spec}}{\mathcal{O}}^\mathrm{CS}\breve{\Phi}^r_{\psi,A},\underset{I}{\mathrm{homotopylimit}}~\underset{\mathrm{Spec}}{\mathcal{O}}^\mathrm{CS}\breve{\Phi}^I_{\psi,A},	\\
}
\]
\[
\xymatrix@R+0pc@C+0pc{
\underset{r}{\mathrm{homotopycolimit}}~\underset{\mathrm{Spec}}{\mathcal{O}}^\mathrm{CS}{\Phi}^r_{\psi,A},\underset{I}{\mathrm{homotopylimit}}~\underset{\mathrm{Spec}}{\mathcal{O}}^\mathrm{CS}{\Phi}^I_{\psi,A}.	
}
\]
\[ 
\xymatrix@R+0pc@C+0pc{
\underset{r}{\mathrm{homotopycolimit}}~\underset{\mathrm{Spec}}{\mathcal{O}}^\mathrm{CS}\widetilde{\Phi}^r_{\psi,A}/\mathrm{Fro}^\mathbb{Z},\underset{I}{\mathrm{homotopylimit}}~\underset{\mathrm{Spec}}{\mathcal{O}}^\mathrm{CS}\widetilde{\Phi}^I_{\psi,A}/\mathrm{Fro}^\mathbb{Z},	\\
}
\]
\[ 
\xymatrix@R+0pc@C+0pc{
\underset{r}{\mathrm{homotopycolimit}}~\underset{\mathrm{Spec}}{\mathcal{O}}^\mathrm{CS}\breve{\Phi}^r_{\psi,A}/\mathrm{Fro}^\mathbb{Z},\underset{I}{\mathrm{homotopylimit}}~\breve{\Phi}^I_{\psi,A}/\mathrm{Fro}^\mathbb{Z},	\\
}
\]
\[ 
\xymatrix@R+0pc@C+0pc{
\underset{r}{\mathrm{homotopycolimit}}~\underset{\mathrm{Spec}}{\mathcal{O}}^\mathrm{CS}{\Phi}^r_{\psi,A}/\mathrm{Fro}^\mathbb{Z},\underset{I}{\mathrm{homotopylimit}}~\underset{\mathrm{Spec}}{\mathcal{O}}^\mathrm{CS}{\Phi}^I_{\psi,A}/\mathrm{Fro}^\mathbb{Z}.	
}
\]	
In this situation we will have the target category being family parametrized by $r$ or $I$ in compatible glueing sense as in \cite[Definition 5.4.10]{10KL2}. In this situation for modules parametrized by the intervals we have the equivalence of $\infty$-categories by using \cite[Proposition 13.8]{10CS2}. Here the corresponding quasicoherent Frobenius modules are defined to be the corresponding homotopy colimits and limits of Frobenius modules:
\begin{align}
\underset{r}{\mathrm{homotopycolimit}}~M_r,\\
\underset{I}{\mathrm{homotopylimit}}~M_I,	
\end{align}
where each $M_r$ is a Frobenius-equivariant module over the period ring with respect to some radius $r$ while each $M_I$ is a Frobenius-equivariant module over the period ring with respect to some interval $I$.\\
\end{proposition}

\begin{proposition}
Similar proposition holds for 
\begin{align}
\mathrm{Quasicoherentsheaves,IndBanach}_{*}.	
\end{align}	
\end{proposition}

\

\begin{definition}
We then consider the corresponding quasipresheaves of perfect complexes the corresponding ind-Banach or monomorphic ind-Banach modules from \cite{10BBK}, \cite{10KKM}:
\begin{align}
\mathrm{Quasicoherentpresheaves,Perfectcomplex,IndBanach}_{*}	
\end{align}
where $*$ is one of the following spaces:
\begin{align}
&\underset{\mathrm{Spec}}{\mathcal{O}}^\mathrm{BK}\widetilde{\Phi}_{\psi,A}/\mathrm{Fro}^\mathbb{Z},	\\
\end{align}
\begin{align}
&\underset{\mathrm{Spec}}{\mathcal{O}}^\mathrm{BK}\breve{\Phi}_{\psi,A}/\mathrm{Fro}^\mathbb{Z},	\\
\end{align}
\begin{align}
&\underset{\mathrm{Spec}}{\mathcal{O}}^\mathrm{BK}{\Phi}_{\psi,A}/\mathrm{Fro}^\mathbb{Z}.	
\end{align}
Here for those space without notation related to the radius and the corresponding interval we consider the total unions $\bigcap_r,\bigcup_I$ in order to achieve the whole spaces to achieve the analogues of the corresponding FF curves from \cite{10KL1}, \cite{10KL2}, \cite{10FF} for
\[
\xymatrix@R+0pc@C+0pc{
\underset{r}{\mathrm{homotopycolimit}}~\underset{\mathrm{Spec}}{\mathcal{O}}^\mathrm{BK}\widetilde{\Phi}^r_{\psi,A},\underset{I}{\mathrm{homotopylimit}}~\underset{\mathrm{Spec}}{\mathcal{O}}^\mathrm{BK}\widetilde{\Phi}^I_{\psi,A},	\\
}
\]
\[
\xymatrix@R+0pc@C+0pc{
\underset{r}{\mathrm{homotopycolimit}}~\underset{\mathrm{Spec}}{\mathcal{O}}^\mathrm{BK}\breve{\Phi}^r_{\psi,A},\underset{I}{\mathrm{homotopylimit}}~\underset{\mathrm{Spec}}{\mathcal{O}}^\mathrm{BK}\breve{\Phi}^I_{\psi,A},	\\
}
\]
\[
\xymatrix@R+0pc@C+0pc{
\underset{r}{\mathrm{homotopycolimit}}~\underset{\mathrm{Spec}}{\mathcal{O}}^\mathrm{BK}{\Phi}^r_{\psi,A},\underset{I}{\mathrm{homotopylimit}}~\underset{\mathrm{Spec}}{\mathcal{O}}^\mathrm{BK}{\Phi}^I_{\psi,A}.	
}
\]
\[  
\xymatrix@R+0pc@C+0pc{
\underset{r}{\mathrm{homotopycolimit}}~\underset{\mathrm{Spec}}{\mathcal{O}}^\mathrm{BK}\widetilde{\Phi}^r_{\psi,A}/\mathrm{Fro}^\mathbb{Z},\underset{I}{\mathrm{homotopylimit}}~\underset{\mathrm{Spec}}{\mathcal{O}}^\mathrm{BK}\widetilde{\Phi}^I_{\psi,A}/\mathrm{Fro}^\mathbb{Z},	\\
}
\]
\[ 
\xymatrix@R+0pc@C+0pc{
\underset{r}{\mathrm{homotopycolimit}}~\underset{\mathrm{Spec}}{\mathcal{O}}^\mathrm{BK}\breve{\Phi}^r_{\psi,A}/\mathrm{Fro}^\mathbb{Z},\underset{I}{\mathrm{homotopylimit}}~\underset{\mathrm{Spec}}{\mathcal{O}}^\mathrm{BK}\breve{\Phi}^I_{\psi,A}/\mathrm{Fro}^\mathbb{Z},	\\
}
\]
\[ 
\xymatrix@R+0pc@C+0pc{
\underset{r}{\mathrm{homotopycolimit}}~\underset{\mathrm{Spec}}{\mathcal{O}}^\mathrm{BK}{\Phi}^r_{\psi,A}/\mathrm{Fro}^\mathbb{Z},\underset{I}{\mathrm{homotopylimit}}~\underset{\mathrm{Spec}}{\mathcal{O}}^\mathrm{BK}{\Phi}^I_{\psi,A}/\mathrm{Fro}^\mathbb{Z}.	
}
\]

\end{definition}

\begin{definition}
We then consider the corresponding quasisheaves of perfect complexes of the corresponding condensed solid topological modules from \cite{10CS2}:
\begin{align}
\mathrm{Quasicoherentsheaves, Perfectcomplex, Condensed}_{*}	
\end{align}
where $*$ is one of the following spaces:
\begin{align}
&\underset{\mathrm{Spec}}{\mathcal{O}}^\mathrm{CS}\widetilde{\Delta}_{\psi,A}/\mathrm{Fro}^\mathbb{Z},\underset{\mathrm{Spec}}{\mathcal{O}}^\mathrm{CS}\widetilde{\nabla}_{\psi,A}/\mathrm{Fro}^\mathbb{Z},\underset{\mathrm{Spec}}{\mathcal{O}}^\mathrm{CS}\widetilde{\Phi}_{\psi,A}/\mathrm{Fro}^\mathbb{Z},\underset{\mathrm{Spec}}{\mathcal{O}}^\mathrm{CS}\widetilde{\Delta}^+_{\psi,A}/\mathrm{Fro}^\mathbb{Z},\\
&\underset{\mathrm{Spec}}{\mathcal{O}}^\mathrm{CS}\widetilde{\nabla}^+_{\psi,A}/\mathrm{Fro}^\mathbb{Z},\underset{\mathrm{Spec}}{\mathcal{O}}^\mathrm{CS}\widetilde{\Delta}^\dagger_{\psi,A}/\mathrm{Fro}^\mathbb{Z},\underset{\mathrm{Spec}}{\mathcal{O}}^\mathrm{CS}\widetilde{\nabla}^\dagger_{\psi,A}/\mathrm{Fro}^\mathbb{Z},	\\
\end{align}
\begin{align}
&\underset{\mathrm{Spec}}{\mathcal{O}}^\mathrm{CS}\breve{\Delta}_{\psi,A}/\mathrm{Fro}^\mathbb{Z},\breve{\nabla}_{\psi,A}/\mathrm{Fro}^\mathbb{Z},\underset{\mathrm{Spec}}{\mathcal{O}}^\mathrm{CS}\breve{\Phi}_{\psi,A}/\mathrm{Fro}^\mathbb{Z},\underset{\mathrm{Spec}}{\mathcal{O}}^\mathrm{CS}\breve{\Delta}^+_{\psi,A}/\mathrm{Fro}^\mathbb{Z},\\
&\underset{\mathrm{Spec}}{\mathcal{O}}^\mathrm{CS}\breve{\nabla}^+_{\psi,A}/\mathrm{Fro}^\mathbb{Z},\underset{\mathrm{Spec}}{\mathcal{O}}^\mathrm{CS}\breve{\Delta}^\dagger_{\psi,A}/\mathrm{Fro}^\mathbb{Z},\underset{\mathrm{Spec}}{\mathcal{O}}^\mathrm{CS}\breve{\nabla}^\dagger_{\psi,A}/\mathrm{Fro}^\mathbb{Z},	\\
\end{align}
\begin{align}
&\underset{\mathrm{Spec}}{\mathcal{O}}^\mathrm{CS}{\Delta}_{\psi,A}/\mathrm{Fro}^\mathbb{Z},\underset{\mathrm{Spec}}{\mathcal{O}}^\mathrm{CS}{\nabla}_{\psi,A}/\mathrm{Fro}^\mathbb{Z},\underset{\mathrm{Spec}}{\mathcal{O}}^\mathrm{CS}{\Phi}_{\psi,A}/\mathrm{Fro}^\mathbb{Z},\underset{\mathrm{Spec}}{\mathcal{O}}^\mathrm{CS}{\Delta}^+_{\psi,A}/\mathrm{Fro}^\mathbb{Z},\\
&\underset{\mathrm{Spec}}{\mathcal{O}}^\mathrm{CS}{\nabla}^+_{\psi,A}/\mathrm{Fro}^\mathbb{Z}, \underset{\mathrm{Spec}}{\mathcal{O}}^\mathrm{CS}{\Delta}^\dagger_{\psi,A}/\mathrm{Fro}^\mathbb{Z},\underset{\mathrm{Spec}}{\mathcal{O}}^\mathrm{CS}{\nabla}^\dagger_{\psi,A}/\mathrm{Fro}^\mathbb{Z}.	
\end{align}
Here for those space with notations related to the radius and the corresponding interval we consider the total unions $\bigcap_r,\bigcup_I$ in order to achieve the whole spaces to achieve the analogues of the corresponding FF curves from \cite{10KL1}, \cite{10KL2}, \cite{10FF} for
\[
\xymatrix@R+0pc@C+0pc{
\underset{r}{\mathrm{homotopycolimit}}~\underset{\mathrm{Spec}}{\mathcal{O}}^\mathrm{CS}\widetilde{\Phi}^r_{\psi,A},\underset{I}{\mathrm{homotopylimit}}~\underset{\mathrm{Spec}}{\mathcal{O}}^\mathrm{CS}\widetilde{\Phi}^I_{\psi,A},	\\
}
\]
\[
\xymatrix@R+0pc@C+0pc{
\underset{r}{\mathrm{homotopycolimit}}~\underset{\mathrm{Spec}}{\mathcal{O}}^\mathrm{CS}\breve{\Phi}^r_{\psi,A},\underset{I}{\mathrm{homotopylimit}}~\underset{\mathrm{Spec}}{\mathcal{O}}^\mathrm{CS}\breve{\Phi}^I_{\psi,A},	\\
}
\]
\[
\xymatrix@R+0pc@C+0pc{
\underset{r}{\mathrm{homotopycolimit}}~\underset{\mathrm{Spec}}{\mathcal{O}}^\mathrm{CS}{\Phi}^r_{\psi,A},\underset{I}{\mathrm{homotopylimit}}~\underset{\mathrm{Spec}}{\mathcal{O}}^\mathrm{CS}{\Phi}^I_{\psi,A}.	
}
\]
\[ 
\xymatrix@R+0pc@C+0pc{
\underset{r}{\mathrm{homotopycolimit}}~\underset{\mathrm{Spec}}{\mathcal{O}}^\mathrm{CS}\widetilde{\Phi}^r_{\psi,A}/\mathrm{Fro}^\mathbb{Z},\underset{I}{\mathrm{homotopylimit}}~\underset{\mathrm{Spec}}{\mathcal{O}}^\mathrm{CS}\widetilde{\Phi}^I_{\psi,A}/\mathrm{Fro}^\mathbb{Z},	\\
}
\]
\[ 
\xymatrix@R+0pc@C+0pc{
\underset{r}{\mathrm{homotopycolimit}}~\underset{\mathrm{Spec}}{\mathcal{O}}^\mathrm{CS}\breve{\Phi}^r_{\psi,A}/\mathrm{Fro}^\mathbb{Z},\underset{I}{\mathrm{homotopylimit}}~\breve{\Phi}^I_{\psi,A}/\mathrm{Fro}^\mathbb{Z},	\\
}
\]
\[ 
\xymatrix@R+0pc@C+0pc{
\underset{r}{\mathrm{homotopycolimit}}~\underset{\mathrm{Spec}}{\mathcal{O}}^\mathrm{CS}{\Phi}^r_{\psi,A}/\mathrm{Fro}^\mathbb{Z},\underset{I}{\mathrm{homotopylimit}}~\underset{\mathrm{Spec}}{\mathcal{O}}^\mathrm{CS}{\Phi}^I_{\psi,A}/\mathrm{Fro}^\mathbb{Z}.	
}
\]

\end{definition}

\begin{proposition}
There is a well-defined functor from the $\infty$-category 
\begin{align}
\mathrm{Quasicoherentpresheaves,Perfectcomplex,Condensed}_{*}	
\end{align}
where $*$ is one of the following spaces:
\begin{align}
&\underset{\mathrm{Spec}}{\mathcal{O}}^\mathrm{CS}\widetilde{\Phi}_{\psi,A}/\mathrm{Fro}^\mathbb{Z},	\\
\end{align}
\begin{align}
&\underset{\mathrm{Spec}}{\mathcal{O}}^\mathrm{CS}\breve{\Phi}_{\psi,A}/\mathrm{Fro}^\mathbb{Z},	\\
\end{align}
\begin{align}
&\underset{\mathrm{Spec}}{\mathcal{O}}^\mathrm{CS}{\Phi}_{\psi,A}/\mathrm{Fro}^\mathbb{Z},	
\end{align}
to the $\infty$-category of $\mathrm{Fro}$-equivariant quasicoherent presheaves over similar spaces above correspondingly without the $\mathrm{Fro}$-quotients, and to the $\infty$-category of $\mathrm{Fro}$-equivariant quasicoherent modules over global sections of the structure $\infty$-sheaves of the similar spaces above correspondingly without the $\mathrm{Fro}$-quotients. Here for those space without notation related to the radius and the corresponding interval we consider the total unions $\bigcap_r,\bigcup_I$ in order to achieve the whole spaces to achieve the analogues of the corresponding FF curves from \cite{10KL1}, \cite{10KL2}, \cite{10FF} for
\[
\xymatrix@R+0pc@C+0pc{
\underset{r}{\mathrm{homotopycolimit}}~\underset{\mathrm{Spec}}{\mathcal{O}}^\mathrm{CS}\widetilde{\Phi}^r_{\psi,A},\underset{I}{\mathrm{homotopylimit}}~\underset{\mathrm{Spec}}{\mathcal{O}}^\mathrm{CS}\widetilde{\Phi}^I_{\psi,A},	\\
}
\]
\[
\xymatrix@R+0pc@C+0pc{
\underset{r}{\mathrm{homotopycolimit}}~\underset{\mathrm{Spec}}{\mathcal{O}}^\mathrm{CS}\breve{\Phi}^r_{\psi,A},\underset{I}{\mathrm{homotopylimit}}~\underset{\mathrm{Spec}}{\mathcal{O}}^\mathrm{CS}\breve{\Phi}^I_{\psi,A},	\\
}
\]
\[
\xymatrix@R+0pc@C+0pc{
\underset{r}{\mathrm{homotopycolimit}}~\underset{\mathrm{Spec}}{\mathcal{O}}^\mathrm{CS}{\Phi}^r_{\psi,A},\underset{I}{\mathrm{homotopylimit}}~\underset{\mathrm{Spec}}{\mathcal{O}}^\mathrm{CS}{\Phi}^I_{\psi,A}.	
}
\]
\[ 
\xymatrix@R+0pc@C+0pc{
\underset{r}{\mathrm{homotopycolimit}}~\underset{\mathrm{Spec}}{\mathcal{O}}^\mathrm{CS}\widetilde{\Phi}^r_{\psi,A}/\mathrm{Fro}^\mathbb{Z},\underset{I}{\mathrm{homotopylimit}}~\underset{\mathrm{Spec}}{\mathcal{O}}^\mathrm{CS}\widetilde{\Phi}^I_{\psi,A}/\mathrm{Fro}^\mathbb{Z},	\\
}
\]
\[ 
\xymatrix@R+0pc@C+0pc{
\underset{r}{\mathrm{homotopycolimit}}~\underset{\mathrm{Spec}}{\mathcal{O}}^\mathrm{CS}\breve{\Phi}^r_{\psi,A}/\mathrm{Fro}^\mathbb{Z},\underset{I}{\mathrm{homotopylimit}}~\breve{\Phi}^I_{\psi,A}/\mathrm{Fro}^\mathbb{Z},	\\
}
\]
\[ 
\xymatrix@R+0pc@C+0pc{
\underset{r}{\mathrm{homotopycolimit}}~\underset{\mathrm{Spec}}{\mathcal{O}}^\mathrm{CS}{\Phi}^r_{\psi,A}/\mathrm{Fro}^\mathbb{Z},\underset{I}{\mathrm{homotopylimit}}~\underset{\mathrm{Spec}}{\mathcal{O}}^\mathrm{CS}{\Phi}^I_{\psi,A}/\mathrm{Fro}^\mathbb{Z}.	
}
\]	
In this situation we will have the target category being family parametrized by $r$ or $I$ in compatible glueing sense as in \cite[Definition 5.4.10]{10KL2}. In this situation for modules parametrized by the intervals we have the equivalence of $\infty$-categories by using \cite[Proposition 12.18]{10CS2}. Here the corresponding quasicoherent Frobenius modules are defined to be the corresponding homotopy colimits and limits of Frobenius modules:
\begin{align}
\underset{r}{\mathrm{homotopycolimit}}~M_r,\\
\underset{I}{\mathrm{homotopylimit}}~M_I,	
\end{align}
where each $M_r$ is a Frobenius-equivariant module over the period ring with respect to some radius $r$ while each $M_I$ is a Frobenius-equivariant module over the period ring with respect to some interval $I$.\\
\end{proposition}

\begin{proposition}
Similar proposition holds for 
\begin{align}
\mathrm{Quasicoherentsheaves,Perfectcomplex,IndBanach}_{*}.	
\end{align}	
\end{proposition}

\newpage
\subsection{Frobenius Quasicoherent Modules II: Deformation in Banach Rings}

\begin{definition}
Let $\psi$ be a toric tower over $\mathbb{Q}_p$ as in \cite[Chapter 7]{10KL2} with base $\mathbb{Q}_p\left<X_1^{\pm 1},...,X_k^{\pm 1}\right>$. Then from \cite{10KL1} and \cite[Definition 5.2.1]{10KL2} we have the following class of Kedlaya-Liu rings (with the following replacement: $\Delta$ stands for $A$, $\nabla$ stands for $B$, while $\Phi$ stands for $C$) by taking product in the sense of self $\Gamma$-th power\footnote{Here $|\Gamma|=1$.}:

\[
\xymatrix@R+0pc@C+0pc{
\widetilde{\Delta}_{\psi},\widetilde{\nabla}_{\psi},\widetilde{\Phi}_{\psi},\widetilde{\Delta}^+_{\psi},\widetilde{\nabla}^+_{\psi},\widetilde{\Delta}^\dagger_{\psi},\widetilde{\nabla}^\dagger_{\psi},\widetilde{\Phi}^r_{\psi},\widetilde{\Phi}^I_{\psi}, 
}
\]

\[
\xymatrix@R+0pc@C+0pc{
\breve{\Delta}_{\psi},\breve{\nabla}_{\psi},\breve{\Phi}_{\psi},\breve{\Delta}^+_{\psi},\breve{\nabla}^+_{\psi},\breve{\Delta}^\dagger_{\psi},\breve{\nabla}^\dagger_{\psi},\breve{\Phi}^r_{\psi},\breve{\Phi}^I_{\psi},	
}
\]

\[
\xymatrix@R+0pc@C+0pc{
{\Delta}_{\psi},{\nabla}_{\psi},{\Phi}_{\psi},{\Delta}^+_{\psi},{\nabla}^+_{\psi},{\Delta}^\dagger_{\psi},{\nabla}^\dagger_{\psi},{\Phi}^r_{\psi},{\Phi}^I_{\psi}.	
}
\]
We now consider the following rings with $-$ being any deforming Banach ring over $\mathbb{Q}_p$. Taking the product we have:
\[
\xymatrix@R+0pc@C+0pc{
\widetilde{\Phi}_{\psi,-},\widetilde{\Phi}^r_{\psi,-},\widetilde{\Phi}^I_{\psi,-},	
}
\]
\[
\xymatrix@R+0pc@C+0pc{
\breve{\Phi}_{\psi,-},\breve{\Phi}^r_{\psi,-},\breve{\Phi}^I_{\psi,-},	
}
\]
\[
\xymatrix@R+0pc@C+0pc{
{\Phi}_{\psi,-},{\Phi}^r_{\psi,-},{\Phi}^I_{\psi,-}.	
}
\]
They carry multi Frobenius action $\varphi_\Gamma$ and multi $\mathrm{Lie}_\Gamma:=\mathbb{Z}_p^{\times\Gamma}$ action. In our current situation after \cite{10CKZ} and \cite{10PZ} we consider the following $(\infty,1)$-categories of $(\infty,1)$-modules.\\
\end{definition}

\begin{definition}
First we consider the Bambozzi-Kremnizer spectrum $\underset{\mathrm{Spec}}{\mathcal{O}}^\mathrm{BK}(*)$ attached to any of those in the above from \cite{10BK} by taking derived rational localization:
\begin{align}
&\underset{\mathrm{Spec}}{\mathcal{O}}^\mathrm{BK}\widetilde{\Phi}_{\psi,-},\underset{\mathrm{Spec}}{\mathcal{O}}^\mathrm{BK}\widetilde{\Phi}^r_{\psi,-},\underset{\mathrm{Spec}}{\mathcal{O}}^\mathrm{BK}\widetilde{\Phi}^I_{\psi,-},	
\end{align}
\begin{align}
&\underset{\mathrm{Spec}}{\mathcal{O}}^\mathrm{BK}\breve{\Phi}_{\psi,-},\underset{\mathrm{Spec}}{\mathcal{O}}^\mathrm{BK}\breve{\Phi}^r_{\psi,-},\underset{\mathrm{Spec}}{\mathcal{O}}^\mathrm{BK}\breve{\Phi}^I_{\psi,-},	
\end{align}
\begin{align}
&\underset{\mathrm{Spec}}{\mathcal{O}}^\mathrm{BK}{\Phi}_{\psi,-},
\underset{\mathrm{Spec}}{\mathcal{O}}^\mathrm{BK}{\Phi}^r_{\psi,-},\underset{\mathrm{Spec}}{\mathcal{O}}^\mathrm{BK}{\Phi}^I_{\psi,-}.	
\end{align}

Then we take the corresponding quotients by using the corresponding Frobenius operators:
\begin{align}
&\underset{\mathrm{Spec}}{\mathcal{O}}^\mathrm{BK}\widetilde{\Phi}_{\psi,-}/\mathrm{Fro}^\mathbb{Z},	\\
\end{align}
\begin{align}
&\underset{\mathrm{Spec}}{\mathcal{O}}^\mathrm{BK}\breve{\Phi}_{\psi,-}/\mathrm{Fro}^\mathbb{Z},	\\
\end{align}
\begin{align}
&\underset{\mathrm{Spec}}{\mathcal{O}}^\mathrm{BK}{\Phi}_{\psi,-}/\mathrm{Fro}^\mathbb{Z}.	
\end{align}
Here for those space without notation related to the radius and the corresponding interval we consider the total unions $\bigcap_r,\bigcup_I$ in order to achieve the whole spaces to achieve the analogues of the corresponding FF curves from \cite{10KL1}, \cite{10KL2}, \cite{10FF} for
\[
\xymatrix@R+0pc@C+0pc{
\underset{r}{\mathrm{homotopycolimit}}~\underset{\mathrm{Spec}}{\mathcal{O}}^\mathrm{BK}\widetilde{\Phi}^r_{\psi,-},\underset{I}{\mathrm{homotopylimit}}~\underset{\mathrm{Spec}}{\mathcal{O}}^\mathrm{BK}\widetilde{\Phi}^I_{\psi,-},	\\
}
\]
\[
\xymatrix@R+0pc@C+0pc{
\underset{r}{\mathrm{homotopycolimit}}~\underset{\mathrm{Spec}}{\mathcal{O}}^\mathrm{BK}\breve{\Phi}^r_{\psi,-},\underset{I}{\mathrm{homotopylimit}}~\underset{\mathrm{Spec}}{\mathcal{O}}^\mathrm{BK}\breve{\Phi}^I_{\psi,-},	\\
}
\]
\[
\xymatrix@R+0pc@C+0pc{
\underset{r}{\mathrm{homotopycolimit}}~\underset{\mathrm{Spec}}{\mathcal{O}}^\mathrm{BK}{\Phi}^r_{\psi,-},\underset{I}{\mathrm{homotopylimit}}~\underset{\mathrm{Spec}}{\mathcal{O}}^\mathrm{BK}{\Phi}^I_{\psi,-}.	
}
\]
\[  
\xymatrix@R+0pc@C+0pc{
\underset{r}{\mathrm{homotopycolimit}}~\underset{\mathrm{Spec}}{\mathcal{O}}^\mathrm{BK}\widetilde{\Phi}^r_{\psi,-}/\mathrm{Fro}^\mathbb{Z},\underset{I}{\mathrm{homotopylimit}}~\underset{\mathrm{Spec}}{\mathcal{O}}^\mathrm{BK}\widetilde{\Phi}^I_{\psi,-}/\mathrm{Fro}^\mathbb{Z},	\\
}
\]
\[ 
\xymatrix@R+0pc@C+0pc{
\underset{r}{\mathrm{homotopycolimit}}~\underset{\mathrm{Spec}}{\mathcal{O}}^\mathrm{BK}\breve{\Phi}^r_{\psi,-}/\mathrm{Fro}^\mathbb{Z},\underset{I}{\mathrm{homotopylimit}}~\underset{\mathrm{Spec}}{\mathcal{O}}^\mathrm{BK}\breve{\Phi}^I_{\psi,-}/\mathrm{Fro}^\mathbb{Z},	\\
}
\]
\[ 
\xymatrix@R+0pc@C+0pc{
\underset{r}{\mathrm{homotopycolimit}}~\underset{\mathrm{Spec}}{\mathcal{O}}^\mathrm{BK}{\Phi}^r_{\psi,-}/\mathrm{Fro}^\mathbb{Z},\underset{I}{\mathrm{homotopylimit}}~\underset{\mathrm{Spec}}{\mathcal{O}}^\mathrm{BK}{\Phi}^I_{\psi,-}/\mathrm{Fro}^\mathbb{Z}.	
}
\]

\end{definition}

\indent Meanwhile we have the corresponding Clausen-Scholze analytic stacks from \cite{10CS2}, therefore applying their construction we have:

\begin{definition}
Here we define the following products by using the solidified tensor product from \cite{10CS1} and \cite{10CS2}. Namely $A$ will still as above as a Banach ring over $\mathbb{Q}_p$. Then we take solidified tensor product $\overset{\blacksquare}{\otimes}$ of any of the following
\[
\xymatrix@R+0pc@C+0pc{
\widetilde{\Delta}_{\psi},\widetilde{\nabla}_{\psi},\widetilde{\Phi}_{\psi},\widetilde{\Delta}^+_{\psi},\widetilde{\nabla}^+_{\psi},\widetilde{\Delta}^\dagger_{\psi},\widetilde{\nabla}^\dagger_{\psi},\widetilde{\Phi}^r_{\psi},\widetilde{\Phi}^I_{\psi}, 
}
\]

\[
\xymatrix@R+0pc@C+0pc{
\breve{\Delta}_{\psi},\breve{\nabla}_{\psi},\breve{\Phi}_{\psi},\breve{\Delta}^+_{\psi},\breve{\nabla}^+_{\psi},\breve{\Delta}^\dagger_{\psi},\breve{\nabla}^\dagger_{\psi},\breve{\Phi}^r_{\psi},\breve{\Phi}^I_{\psi},	
}
\]

\[
\xymatrix@R+0pc@C+0pc{
{\Delta}_{\psi},{\nabla}_{\psi},{\Phi}_{\psi},{\Delta}^+_{\psi},{\nabla}^+_{\psi},{\Delta}^\dagger_{\psi},{\nabla}^\dagger_{\psi},{\Phi}^r_{\psi},{\Phi}^I_{\psi},	
}
\]  	
with $A$. Then we have the notations:
\[
\xymatrix@R+0pc@C+0pc{
\widetilde{\Delta}_{\psi,-},\widetilde{\nabla}_{\psi,-},\widetilde{\Phi}_{\psi,-},\widetilde{\Delta}^+_{\psi,-},\widetilde{\nabla}^+_{\psi,-},\widetilde{\Delta}^\dagger_{\psi,-},\widetilde{\nabla}^\dagger_{\psi,-},\widetilde{\Phi}^r_{\psi,-},\widetilde{\Phi}^I_{\psi,-}, 
}
\]

\[
\xymatrix@R+0pc@C+0pc{
\breve{\Delta}_{\psi,-},\breve{\nabla}_{\psi,-},\breve{\Phi}_{\psi,-},\breve{\Delta}^+_{\psi,-},\breve{\nabla}^+_{\psi,-},\breve{\Delta}^\dagger_{\psi,-},\breve{\nabla}^\dagger_{\psi,-},\breve{\Phi}^r_{\psi,-},\breve{\Phi}^I_{\psi,-},	
}
\]

\[
\xymatrix@R+0pc@C+0pc{
{\Delta}_{\psi,-},{\nabla}_{\psi,-},{\Phi}_{\psi,-},{\Delta}^+_{\psi,-},{\nabla}^+_{\psi,-},{\Delta}^\dagger_{\psi,-},{\nabla}^\dagger_{\psi,-},{\Phi}^r_{\psi,-},{\Phi}^I_{\psi,-}.	
}
\]
\end{definition}

\begin{definition}
First we consider the Clausen-Scholze spectrum $\underset{\mathrm{Spec}}{\mathcal{O}}^\mathrm{CS}(*)$ attached to any of those in the above from \cite{10CS2} by taking derived rational localization:
\begin{align}
\underset{\mathrm{Spec}}{\mathcal{O}}^\mathrm{CS}\widetilde{\Delta}_{\psi,-},\underset{\mathrm{Spec}}{\mathcal{O}}^\mathrm{CS}\widetilde{\nabla}_{\psi,-},\underset{\mathrm{Spec}}{\mathcal{O}}^\mathrm{CS}\widetilde{\Phi}_{\psi,-},\underset{\mathrm{Spec}}{\mathcal{O}}^\mathrm{CS}\widetilde{\Delta}^+_{\psi,-},\underset{\mathrm{Spec}}{\mathcal{O}}^\mathrm{CS}\widetilde{\nabla}^+_{\psi,-},\\
\underset{\mathrm{Spec}}{\mathcal{O}}^\mathrm{CS}\widetilde{\Delta}^\dagger_{\psi,-},\underset{\mathrm{Spec}}{\mathcal{O}}^\mathrm{CS}\widetilde{\nabla}^\dagger_{\psi,-},\underset{\mathrm{Spec}}{\mathcal{O}}^\mathrm{CS}\widetilde{\Phi}^r_{\psi,-},\underset{\mathrm{Spec}}{\mathcal{O}}^\mathrm{CS}\widetilde{\Phi}^I_{\psi,-},	\\
\end{align}
\begin{align}
\underset{\mathrm{Spec}}{\mathcal{O}}^\mathrm{CS}\breve{\Delta}_{\psi,-},\breve{\nabla}_{\psi,-},\underset{\mathrm{Spec}}{\mathcal{O}}^\mathrm{CS}\breve{\Phi}_{\psi,-},\underset{\mathrm{Spec}}{\mathcal{O}}^\mathrm{CS}\breve{\Delta}^+_{\psi,-},\underset{\mathrm{Spec}}{\mathcal{O}}^\mathrm{CS}\breve{\nabla}^+_{\psi,-},\\
\underset{\mathrm{Spec}}{\mathcal{O}}^\mathrm{CS}\breve{\Delta}^\dagger_{\psi,-},\underset{\mathrm{Spec}}{\mathcal{O}}^\mathrm{CS}\breve{\nabla}^\dagger_{\psi,-},\underset{\mathrm{Spec}}{\mathcal{O}}^\mathrm{CS}\breve{\Phi}^r_{\psi,-},\breve{\Phi}^I_{\psi,-},	\\
\end{align}
\begin{align}
\underset{\mathrm{Spec}}{\mathcal{O}}^\mathrm{CS}{\Delta}_{\psi,-},\underset{\mathrm{Spec}}{\mathcal{O}}^\mathrm{CS}{\nabla}_{\psi,-},\underset{\mathrm{Spec}}{\mathcal{O}}^\mathrm{CS}{\Phi}_{\psi,-},\underset{\mathrm{Spec}}{\mathcal{O}}^\mathrm{CS}{\Delta}^+_{\psi,-},\underset{\mathrm{Spec}}{\mathcal{O}}^\mathrm{CS}{\nabla}^+_{\psi,-},\\
\underset{\mathrm{Spec}}{\mathcal{O}}^\mathrm{CS}{\Delta}^\dagger_{\psi,-},\underset{\mathrm{Spec}}{\mathcal{O}}^\mathrm{CS}{\nabla}^\dagger_{\psi,-},\underset{\mathrm{Spec}}{\mathcal{O}}^\mathrm{CS}{\Phi}^r_{\psi,-},\underset{\mathrm{Spec}}{\mathcal{O}}^\mathrm{CS}{\Phi}^I_{\psi,-}.	
\end{align}

Then we take the corresponding quotients by using the corresponding Frobenius operators:
\begin{align}
&\underset{\mathrm{Spec}}{\mathcal{O}}^\mathrm{CS}\widetilde{\Delta}_{\psi,-}/\mathrm{Fro}^\mathbb{Z},\underset{\mathrm{Spec}}{\mathcal{O}}^\mathrm{CS}\widetilde{\nabla}_{\psi,-}/\mathrm{Fro}^\mathbb{Z},\underset{\mathrm{Spec}}{\mathcal{O}}^\mathrm{CS}\widetilde{\Phi}_{\psi,-}/\mathrm{Fro}^\mathbb{Z},\underset{\mathrm{Spec}}{\mathcal{O}}^\mathrm{CS}\widetilde{\Delta}^+_{\psi,-}/\mathrm{Fro}^\mathbb{Z},\\
&\underset{\mathrm{Spec}}{\mathcal{O}}^\mathrm{CS}\widetilde{\nabla}^+_{\psi,-}/\mathrm{Fro}^\mathbb{Z}, \underset{\mathrm{Spec}}{\mathcal{O}}^\mathrm{CS}\widetilde{\Delta}^\dagger_{\psi,-}/\mathrm{Fro}^\mathbb{Z},\underset{\mathrm{Spec}}{\mathcal{O}}^\mathrm{CS}\widetilde{\nabla}^\dagger_{\psi,-}/\mathrm{Fro}^\mathbb{Z},	\\
\end{align}
\begin{align}
&\underset{\mathrm{Spec}}{\mathcal{O}}^\mathrm{CS}\breve{\Delta}_{\psi,-}/\mathrm{Fro}^\mathbb{Z},\breve{\nabla}_{\psi,-}/\mathrm{Fro}^\mathbb{Z},\underset{\mathrm{Spec}}{\mathcal{O}}^\mathrm{CS}\breve{\Phi}_{\psi,-}/\mathrm{Fro}^\mathbb{Z},\underset{\mathrm{Spec}}{\mathcal{O}}^\mathrm{CS}\breve{\Delta}^+_{\psi,-}/\mathrm{Fro}^\mathbb{Z},\\
&\underset{\mathrm{Spec}}{\mathcal{O}}^\mathrm{CS}\breve{\nabla}^+_{\psi,-}/\mathrm{Fro}^\mathbb{Z}, \underset{\mathrm{Spec}}{\mathcal{O}}^\mathrm{CS}\breve{\Delta}^\dagger_{\psi,-}/\mathrm{Fro}^\mathbb{Z},\underset{\mathrm{Spec}}{\mathcal{O}}^\mathrm{CS}\breve{\nabla}^\dagger_{\psi,-}/\mathrm{Fro}^\mathbb{Z},	\\
\end{align}
\begin{align}
&\underset{\mathrm{Spec}}{\mathcal{O}}^\mathrm{CS}{\Delta}_{\psi,-}/\mathrm{Fro}^\mathbb{Z},\underset{\mathrm{Spec}}{\mathcal{O}}^\mathrm{CS}{\nabla}_{\psi,-}/\mathrm{Fro}^\mathbb{Z},\underset{\mathrm{Spec}}{\mathcal{O}}^\mathrm{CS}{\Phi}_{\psi,-}/\mathrm{Fro}^\mathbb{Z},\underset{\mathrm{Spec}}{\mathcal{O}}^\mathrm{CS}{\Delta}^+_{\psi,-}/\mathrm{Fro}^\mathbb{Z},\\
&\underset{\mathrm{Spec}}{\mathcal{O}}^\mathrm{CS}{\nabla}^+_{\psi,-}/\mathrm{Fro}^\mathbb{Z}, \underset{\mathrm{Spec}}{\mathcal{O}}^\mathrm{CS}{\Delta}^\dagger_{\psi,-}/\mathrm{Fro}^\mathbb{Z},\underset{\mathrm{Spec}}{\mathcal{O}}^\mathrm{CS}{\nabla}^\dagger_{\psi,-}/\mathrm{Fro}^\mathbb{Z}.	
\end{align}
Here for those space with notations related to the radius and the corresponding interval we consider the total unions $\bigcap_r,\bigcup_I$ in order to achieve the whole spaces to achieve the analogues of the corresponding FF curves from \cite{10KL1}, \cite{10KL2}, \cite{10FF} for
\[
\xymatrix@R+0pc@C+0pc{
\underset{r}{\mathrm{homotopycolimit}}~\underset{\mathrm{Spec}}{\mathcal{O}}^\mathrm{CS}\widetilde{\Phi}^r_{\psi,-},\underset{I}{\mathrm{homotopylimit}}~\underset{\mathrm{Spec}}{\mathcal{O}}^\mathrm{CS}\widetilde{\Phi}^I_{\psi,-},	\\
}
\]
\[
\xymatrix@R+0pc@C+0pc{
\underset{r}{\mathrm{homotopycolimit}}~\underset{\mathrm{Spec}}{\mathcal{O}}^\mathrm{CS}\breve{\Phi}^r_{\psi,-},\underset{I}{\mathrm{homotopylimit}}~\underset{\mathrm{Spec}}{\mathcal{O}}^\mathrm{CS}\breve{\Phi}^I_{\psi,-},	\\
}
\]
\[
\xymatrix@R+0pc@C+0pc{
\underset{r}{\mathrm{homotopycolimit}}~\underset{\mathrm{Spec}}{\mathcal{O}}^\mathrm{CS}{\Phi}^r_{\psi,-},\underset{I}{\mathrm{homotopylimit}}~\underset{\mathrm{Spec}}{\mathcal{O}}^\mathrm{CS}{\Phi}^I_{\psi,-}.	
}
\]
\[ 
\xymatrix@R+0pc@C+0pc{
\underset{r}{\mathrm{homotopycolimit}}~\underset{\mathrm{Spec}}{\mathcal{O}}^\mathrm{CS}\widetilde{\Phi}^r_{\psi,-}/\mathrm{Fro}^\mathbb{Z},\underset{I}{\mathrm{homotopylimit}}~\underset{\mathrm{Spec}}{\mathcal{O}}^\mathrm{CS}\widetilde{\Phi}^I_{\psi,-}/\mathrm{Fro}^\mathbb{Z},	\\
}
\]
\[ 
\xymatrix@R+0pc@C+0pc{
\underset{r}{\mathrm{homotopycolimit}}~\underset{\mathrm{Spec}}{\mathcal{O}}^\mathrm{CS}\breve{\Phi}^r_{\psi,-}/\mathrm{Fro}^\mathbb{Z},\underset{I}{\mathrm{homotopylimit}}~\breve{\Phi}^I_{\psi,-}/\mathrm{Fro}^\mathbb{Z},	\\
}
\]
\[ 
\xymatrix@R+0pc@C+0pc{
\underset{r}{\mathrm{homotopycolimit}}~\underset{\mathrm{Spec}}{\mathcal{O}}^\mathrm{CS}{\Phi}^r_{\psi,-}/\mathrm{Fro}^\mathbb{Z},\underset{I}{\mathrm{homotopylimit}}~\underset{\mathrm{Spec}}{\mathcal{O}}^\mathrm{CS}{\Phi}^I_{\psi,-}/\mathrm{Fro}^\mathbb{Z}.	
}
\]

\end{definition}

\

\begin{definition}
We then consider the corresponding quasipresheaves of the corresponding ind-Banach or monomorphic ind-Banach modules from \cite{10BBK}, \cite{10KKM}:
\begin{align}
\mathrm{Quasicoherentpresheaves,IndBanach}_{*}	
\end{align}
where $*$ is one of the following spaces:
\begin{align}
&\underset{\mathrm{Spec}}{\mathcal{O}}^\mathrm{BK}\widetilde{\Phi}_{\psi,-}/\mathrm{Fro}^\mathbb{Z},	\\
\end{align}
\begin{align}
&\underset{\mathrm{Spec}}{\mathcal{O}}^\mathrm{BK}\breve{\Phi}_{\psi,-}/\mathrm{Fro}^\mathbb{Z},	\\
\end{align}
\begin{align}
&\underset{\mathrm{Spec}}{\mathcal{O}}^\mathrm{BK}{\Phi}_{\psi,-}/\mathrm{Fro}^\mathbb{Z}.	
\end{align}
Here for those space without notation related to the radius and the corresponding interval we consider the total unions $\bigcap_r,\bigcup_I$ in order to achieve the whole spaces to achieve the analogues of the corresponding FF curves from \cite{10KL1}, \cite{10KL2}, \cite{10FF} for
\[
\xymatrix@R+0pc@C+0pc{
\underset{r}{\mathrm{homotopycolimit}}~\underset{\mathrm{Spec}}{\mathcal{O}}^\mathrm{BK}\widetilde{\Phi}^r_{\psi,-},\underset{I}{\mathrm{homotopylimit}}~\underset{\mathrm{Spec}}{\mathcal{O}}^\mathrm{BK}\widetilde{\Phi}^I_{\psi,-},	\\
}
\]
\[
\xymatrix@R+0pc@C+0pc{
\underset{r}{\mathrm{homotopycolimit}}~\underset{\mathrm{Spec}}{\mathcal{O}}^\mathrm{BK}\breve{\Phi}^r_{\psi,-},\underset{I}{\mathrm{homotopylimit}}~\underset{\mathrm{Spec}}{\mathcal{O}}^\mathrm{BK}\breve{\Phi}^I_{\psi,-},	\\
}
\]
\[
\xymatrix@R+0pc@C+0pc{
\underset{r}{\mathrm{homotopycolimit}}~\underset{\mathrm{Spec}}{\mathcal{O}}^\mathrm{BK}{\Phi}^r_{\psi,-},\underset{I}{\mathrm{homotopylimit}}~\underset{\mathrm{Spec}}{\mathcal{O}}^\mathrm{BK}{\Phi}^I_{\psi,-}.	
}
\]
\[  
\xymatrix@R+0pc@C+0pc{
\underset{r}{\mathrm{homotopycolimit}}~\underset{\mathrm{Spec}}{\mathcal{O}}^\mathrm{BK}\widetilde{\Phi}^r_{\psi,-}/\mathrm{Fro}^\mathbb{Z},\underset{I}{\mathrm{homotopylimit}}~\underset{\mathrm{Spec}}{\mathcal{O}}^\mathrm{BK}\widetilde{\Phi}^I_{\psi,-}/\mathrm{Fro}^\mathbb{Z},	\\
}
\]
\[ 
\xymatrix@R+0pc@C+0pc{
\underset{r}{\mathrm{homotopycolimit}}~\underset{\mathrm{Spec}}{\mathcal{O}}^\mathrm{BK}\breve{\Phi}^r_{\psi,-}/\mathrm{Fro}^\mathbb{Z},\underset{I}{\mathrm{homotopylimit}}~\underset{\mathrm{Spec}}{\mathcal{O}}^\mathrm{BK}\breve{\Phi}^I_{\psi,-}/\mathrm{Fro}^\mathbb{Z},	\\
}
\]
\[ 
\xymatrix@R+0pc@C+0pc{
\underset{r}{\mathrm{homotopycolimit}}~\underset{\mathrm{Spec}}{\mathcal{O}}^\mathrm{BK}{\Phi}^r_{\psi,-}/\mathrm{Fro}^\mathbb{Z},\underset{I}{\mathrm{homotopylimit}}~\underset{\mathrm{Spec}}{\mathcal{O}}^\mathrm{BK}{\Phi}^I_{\psi,-}/\mathrm{Fro}^\mathbb{Z}.	
}
\]

\end{definition}

\begin{definition}
We then consider the corresponding quasisheaves of the corresponding condensed solid topological modules from \cite{10CS2}:
\begin{align}
\mathrm{Quasicoherentsheaves, Condensed}_{*}	
\end{align}
where $*$ is one of the following spaces:
\begin{align}
&\underset{\mathrm{Spec}}{\mathcal{O}}^\mathrm{CS}\widetilde{\Delta}_{\psi,-}/\mathrm{Fro}^\mathbb{Z},\underset{\mathrm{Spec}}{\mathcal{O}}^\mathrm{CS}\widetilde{\nabla}_{\psi,-}/\mathrm{Fro}^\mathbb{Z},\underset{\mathrm{Spec}}{\mathcal{O}}^\mathrm{CS}\widetilde{\Phi}_{\psi,-}/\mathrm{Fro}^\mathbb{Z},\underset{\mathrm{Spec}}{\mathcal{O}}^\mathrm{CS}\widetilde{\Delta}^+_{\psi,-}/\mathrm{Fro}^\mathbb{Z},\\
&\underset{\mathrm{Spec}}{\mathcal{O}}^\mathrm{CS}\widetilde{\nabla}^+_{\psi,-}/\mathrm{Fro}^\mathbb{Z},\underset{\mathrm{Spec}}{\mathcal{O}}^\mathrm{CS}\widetilde{\Delta}^\dagger_{\psi,-}/\mathrm{Fro}^\mathbb{Z},\underset{\mathrm{Spec}}{\mathcal{O}}^\mathrm{CS}\widetilde{\nabla}^\dagger_{\psi,-}/\mathrm{Fro}^\mathbb{Z},	\\
\end{align}
\begin{align}
&\underset{\mathrm{Spec}}{\mathcal{O}}^\mathrm{CS}\breve{\Delta}_{\psi,-}/\mathrm{Fro}^\mathbb{Z},\breve{\nabla}_{\psi,-}/\mathrm{Fro}^\mathbb{Z},\underset{\mathrm{Spec}}{\mathcal{O}}^\mathrm{CS}\breve{\Phi}_{\psi,-}/\mathrm{Fro}^\mathbb{Z},\underset{\mathrm{Spec}}{\mathcal{O}}^\mathrm{CS}\breve{\Delta}^+_{\psi,-}/\mathrm{Fro}^\mathbb{Z},\\
&\underset{\mathrm{Spec}}{\mathcal{O}}^\mathrm{CS}\breve{\nabla}^+_{\psi,-}/\mathrm{Fro}^\mathbb{Z},\underset{\mathrm{Spec}}{\mathcal{O}}^\mathrm{CS}\breve{\Delta}^\dagger_{\psi,-}/\mathrm{Fro}^\mathbb{Z},\underset{\mathrm{Spec}}{\mathcal{O}}^\mathrm{CS}\breve{\nabla}^\dagger_{\psi,-}/\mathrm{Fro}^\mathbb{Z},	\\
\end{align}
\begin{align}
&\underset{\mathrm{Spec}}{\mathcal{O}}^\mathrm{CS}{\Delta}_{\psi,-}/\mathrm{Fro}^\mathbb{Z},\underset{\mathrm{Spec}}{\mathcal{O}}^\mathrm{CS}{\nabla}_{\psi,-}/\mathrm{Fro}^\mathbb{Z},\underset{\mathrm{Spec}}{\mathcal{O}}^\mathrm{CS}{\Phi}_{\psi,-}/\mathrm{Fro}^\mathbb{Z},\underset{\mathrm{Spec}}{\mathcal{O}}^\mathrm{CS}{\Delta}^+_{\psi,-}/\mathrm{Fro}^\mathbb{Z},\\
&\underset{\mathrm{Spec}}{\mathcal{O}}^\mathrm{CS}{\nabla}^+_{\psi,-}/\mathrm{Fro}^\mathbb{Z}, \underset{\mathrm{Spec}}{\mathcal{O}}^\mathrm{CS}{\Delta}^\dagger_{\psi,-}/\mathrm{Fro}^\mathbb{Z},\underset{\mathrm{Spec}}{\mathcal{O}}^\mathrm{CS}{\nabla}^\dagger_{\psi,-}/\mathrm{Fro}^\mathbb{Z}.	
\end{align}
Here for those space with notations related to the radius and the corresponding interval we consider the total unions $\bigcap_r,\bigcup_I$ in order to achieve the whole spaces to achieve the analogues of the corresponding FF curves from \cite{10KL1}, \cite{10KL2}, \cite{10FF} for
\[
\xymatrix@R+0pc@C+0pc{
\underset{r}{\mathrm{homotopycolimit}}~\underset{\mathrm{Spec}}{\mathcal{O}}^\mathrm{CS}\widetilde{\Phi}^r_{\psi,-},\underset{I}{\mathrm{homotopylimit}}~\underset{\mathrm{Spec}}{\mathcal{O}}^\mathrm{CS}\widetilde{\Phi}^I_{\psi,-},	\\
}
\]
\[
\xymatrix@R+0pc@C+0pc{
\underset{r}{\mathrm{homotopycolimit}}~\underset{\mathrm{Spec}}{\mathcal{O}}^\mathrm{CS}\breve{\Phi}^r_{\psi,-},\underset{I}{\mathrm{homotopylimit}}~\underset{\mathrm{Spec}}{\mathcal{O}}^\mathrm{CS}\breve{\Phi}^I_{\psi,-},	\\
}
\]
\[
\xymatrix@R+0pc@C+0pc{
\underset{r}{\mathrm{homotopycolimit}}~\underset{\mathrm{Spec}}{\mathcal{O}}^\mathrm{CS}{\Phi}^r_{\psi,-},\underset{I}{\mathrm{homotopylimit}}~\underset{\mathrm{Spec}}{\mathcal{O}}^\mathrm{CS}{\Phi}^I_{\psi,-}.	
}
\]
\[ 
\xymatrix@R+0pc@C+0pc{
\underset{r}{\mathrm{homotopycolimit}}~\underset{\mathrm{Spec}}{\mathcal{O}}^\mathrm{CS}\widetilde{\Phi}^r_{\psi,-}/\mathrm{Fro}^\mathbb{Z},\underset{I}{\mathrm{homotopylimit}}~\underset{\mathrm{Spec}}{\mathcal{O}}^\mathrm{CS}\widetilde{\Phi}^I_{\psi,-}/\mathrm{Fro}^\mathbb{Z},	\\
}
\]
\[ 
\xymatrix@R+0pc@C+0pc{
\underset{r}{\mathrm{homotopycolimit}}~\underset{\mathrm{Spec}}{\mathcal{O}}^\mathrm{CS}\breve{\Phi}^r_{\psi,-}/\mathrm{Fro}^\mathbb{Z},\underset{I}{\mathrm{homotopylimit}}~\breve{\Phi}^I_{\psi,-}/\mathrm{Fro}^\mathbb{Z},	\\
}
\]
\[ 
\xymatrix@R+0pc@C+0pc{
\underset{r}{\mathrm{homotopycolimit}}~\underset{\mathrm{Spec}}{\mathcal{O}}^\mathrm{CS}{\Phi}^r_{\psi,-}/\mathrm{Fro}^\mathbb{Z},\underset{I}{\mathrm{homotopylimit}}~\underset{\mathrm{Spec}}{\mathcal{O}}^\mathrm{CS}{\Phi}^I_{\psi,-}/\mathrm{Fro}^\mathbb{Z}.	
}
\]

\end{definition}

\

\begin{proposition}
There is a well-defined functor from the $\infty$-category 
\begin{align}
\mathrm{Quasicoherentpresheaves,Condensed}_{*}	
\end{align}
where $*$ is one of the following spaces:
\begin{align}
&\underset{\mathrm{Spec}}{\mathcal{O}}^\mathrm{CS}\widetilde{\Phi}_{\psi,-}/\mathrm{Fro}^\mathbb{Z},	\\
\end{align}
\begin{align}
&\underset{\mathrm{Spec}}{\mathcal{O}}^\mathrm{CS}\breve{\Phi}_{\psi,-}/\mathrm{Fro}^\mathbb{Z},	\\
\end{align}
\begin{align}
&\underset{\mathrm{Spec}}{\mathcal{O}}^\mathrm{CS}{\Phi}_{\psi,-}/\mathrm{Fro}^\mathbb{Z},	
\end{align}
to the $\infty$-category of $\mathrm{Fro}$-equivariant quasicoherent presheaves over similar spaces above correspondingly without the $\mathrm{Fro}$-quotients, and to the $\infty$-category of $\mathrm{Fro}$-equivariant quasicoherent modules over global sections of the structure $\infty$-sheaves of the similar spaces above correspondingly without the $\mathrm{Fro}$-quotients. Here for those space without notation related to the radius and the corresponding interval we consider the total unions $\bigcap_r,\bigcup_I$ in order to achieve the whole spaces to achieve the analogues of the corresponding FF curves from \cite{10KL1}, \cite{10KL2}, \cite{10FF} for
\[
\xymatrix@R+0pc@C+0pc{
\underset{r}{\mathrm{homotopycolimit}}~\underset{\mathrm{Spec}}{\mathcal{O}}^\mathrm{CS}\widetilde{\Phi}^r_{\psi,-},\underset{I}{\mathrm{homotopylimit}}~\underset{\mathrm{Spec}}{\mathcal{O}}^\mathrm{CS}\widetilde{\Phi}^I_{\psi,-},	\\
}
\]
\[
\xymatrix@R+0pc@C+0pc{
\underset{r}{\mathrm{homotopycolimit}}~\underset{\mathrm{Spec}}{\mathcal{O}}^\mathrm{CS}\breve{\Phi}^r_{\psi,-},\underset{I}{\mathrm{homotopylimit}}~\underset{\mathrm{Spec}}{\mathcal{O}}^\mathrm{CS}\breve{\Phi}^I_{\psi,-},	\\
}
\]
\[
\xymatrix@R+0pc@C+0pc{
\underset{r}{\mathrm{homotopycolimit}}~\underset{\mathrm{Spec}}{\mathcal{O}}^\mathrm{CS}{\Phi}^r_{\psi,-},\underset{I}{\mathrm{homotopylimit}}~\underset{\mathrm{Spec}}{\mathcal{O}}^\mathrm{CS}{\Phi}^I_{\psi,-}.	
}
\]
\[ 
\xymatrix@R+0pc@C+0pc{
\underset{r}{\mathrm{homotopycolimit}}~\underset{\mathrm{Spec}}{\mathcal{O}}^\mathrm{CS}\widetilde{\Phi}^r_{\psi,-}/\mathrm{Fro}^\mathbb{Z},\underset{I}{\mathrm{homotopylimit}}~\underset{\mathrm{Spec}}{\mathcal{O}}^\mathrm{CS}\widetilde{\Phi}^I_{\psi,-}/\mathrm{Fro}^\mathbb{Z},	\\
}
\]
\[ 
\xymatrix@R+0pc@C+0pc{
\underset{r}{\mathrm{homotopycolimit}}~\underset{\mathrm{Spec}}{\mathcal{O}}^\mathrm{CS}\breve{\Phi}^r_{\psi,-}/\mathrm{Fro}^\mathbb{Z},\underset{I}{\mathrm{homotopylimit}}~\breve{\Phi}^I_{\psi,-}/\mathrm{Fro}^\mathbb{Z},	\\
}
\]
\[ 
\xymatrix@R+0pc@C+0pc{
\underset{r}{\mathrm{homotopycolimit}}~\underset{\mathrm{Spec}}{\mathcal{O}}^\mathrm{CS}{\Phi}^r_{\psi,-}/\mathrm{Fro}^\mathbb{Z},\underset{I}{\mathrm{homotopylimit}}~\underset{\mathrm{Spec}}{\mathcal{O}}^\mathrm{CS}{\Phi}^I_{\psi,-}/\mathrm{Fro}^\mathbb{Z}.	
}
\]	
In this situation we will have the target category being family parametrized by $r$ or $I$ in compatible glueing sense as in \cite[Definition 5.4.10]{10KL2}. In this situation for modules parametrized by the intervals we have the equivalence of $\infty$-categories by using \cite[Proposition 13.8]{10CS2}. Here the corresponding quasicoherent Frobenius modules are defined to be the corresponding homotopy colimits and limits of Frobenius modules:
\begin{align}
\underset{r}{\mathrm{homotopycolimit}}~M_r,\\
\underset{I}{\mathrm{homotopylimit}}~M_I,	
\end{align}
where each $M_r$ is a Frobenius-equivariant module over the period ring with respect to some radius $r$ while each $M_I$ is a Frobenius-equivariant module over the period ring with respect to some interval $I$.\\
\end{proposition}

\begin{proposition}
Similar proposition holds for 
\begin{align}
\mathrm{Quasicoherentsheaves,IndBanach}_{*}.	
\end{align}	
\end{proposition}

\

\begin{definition}
We then consider the corresponding quasipresheaves of perfect complexes the corresponding ind-Banach or monomorphic ind-Banach modules from \cite{10BBK}, \cite{10KKM}:
\begin{align}
\mathrm{Quasicoherentpresheaves,Perfectcomplex,IndBanach}_{*}	
\end{align}
where $*$ is one of the following spaces:
\begin{align}
&\underset{\mathrm{Spec}}{\mathcal{O}}^\mathrm{BK}\widetilde{\Phi}_{\psi,-}/\mathrm{Fro}^\mathbb{Z},	\\
\end{align}
\begin{align}
&\underset{\mathrm{Spec}}{\mathcal{O}}^\mathrm{BK}\breve{\Phi}_{\psi,-}/\mathrm{Fro}^\mathbb{Z},	\\
\end{align}
\begin{align}
&\underset{\mathrm{Spec}}{\mathcal{O}}^\mathrm{BK}{\Phi}_{\psi,-}/\mathrm{Fro}^\mathbb{Z}.	
\end{align}
Here for those space without notation related to the radius and the corresponding interval we consider the total unions $\bigcap_r,\bigcup_I$ in order to achieve the whole spaces to achieve the analogues of the corresponding FF curves from \cite{10KL1}, \cite{10KL2}, \cite{10FF} for
\[
\xymatrix@R+0pc@C+0pc{
\underset{r}{\mathrm{homotopycolimit}}~\underset{\mathrm{Spec}}{\mathcal{O}}^\mathrm{BK}\widetilde{\Phi}^r_{\psi,-},\underset{I}{\mathrm{homotopylimit}}~\underset{\mathrm{Spec}}{\mathcal{O}}^\mathrm{BK}\widetilde{\Phi}^I_{\psi,-},	\\
}
\]
\[
\xymatrix@R+0pc@C+0pc{
\underset{r}{\mathrm{homotopycolimit}}~\underset{\mathrm{Spec}}{\mathcal{O}}^\mathrm{BK}\breve{\Phi}^r_{\psi,-},\underset{I}{\mathrm{homotopylimit}}~\underset{\mathrm{Spec}}{\mathcal{O}}^\mathrm{BK}\breve{\Phi}^I_{\psi,-},	\\
}
\]
\[
\xymatrix@R+0pc@C+0pc{
\underset{r}{\mathrm{homotopycolimit}}~\underset{\mathrm{Spec}}{\mathcal{O}}^\mathrm{BK}{\Phi}^r_{\psi,-},\underset{I}{\mathrm{homotopylimit}}~\underset{\mathrm{Spec}}{\mathcal{O}}^\mathrm{BK}{\Phi}^I_{\psi,-}.	
}
\]
\[  
\xymatrix@R+0pc@C+0pc{
\underset{r}{\mathrm{homotopycolimit}}~\underset{\mathrm{Spec}}{\mathcal{O}}^\mathrm{BK}\widetilde{\Phi}^r_{\psi,-}/\mathrm{Fro}^\mathbb{Z},\underset{I}{\mathrm{homotopylimit}}~\underset{\mathrm{Spec}}{\mathcal{O}}^\mathrm{BK}\widetilde{\Phi}^I_{\psi,-}/\mathrm{Fro}^\mathbb{Z},	\\
}
\]
\[ 
\xymatrix@R+0pc@C+0pc{
\underset{r}{\mathrm{homotopycolimit}}~\underset{\mathrm{Spec}}{\mathcal{O}}^\mathrm{BK}\breve{\Phi}^r_{\psi,-}/\mathrm{Fro}^\mathbb{Z},\underset{I}{\mathrm{homotopylimit}}~\underset{\mathrm{Spec}}{\mathcal{O}}^\mathrm{BK}\breve{\Phi}^I_{\psi,-}/\mathrm{Fro}^\mathbb{Z},	\\
}
\]
\[ 
\xymatrix@R+0pc@C+0pc{
\underset{r}{\mathrm{homotopycolimit}}~\underset{\mathrm{Spec}}{\mathcal{O}}^\mathrm{BK}{\Phi}^r_{\psi,-}/\mathrm{Fro}^\mathbb{Z},\underset{I}{\mathrm{homotopylimit}}~\underset{\mathrm{Spec}}{\mathcal{O}}^\mathrm{BK}{\Phi}^I_{\psi,-}/\mathrm{Fro}^\mathbb{Z}.	
}
\]

\end{definition}

\begin{definition}
We then consider the corresponding quasisheaves of perfect complexes of the corresponding condensed solid topological modules from \cite{10CS2}:
\begin{align}
\mathrm{Quasicoherentsheaves, Perfectcomplex, Condensed}_{*}	
\end{align}
where $*$ is one of the following spaces:
\begin{align}
&\underset{\mathrm{Spec}}{\mathcal{O}}^\mathrm{CS}\widetilde{\Delta}_{\psi,-}/\mathrm{Fro}^\mathbb{Z},\underset{\mathrm{Spec}}{\mathcal{O}}^\mathrm{CS}\widetilde{\nabla}_{\psi,-}/\mathrm{Fro}^\mathbb{Z},\underset{\mathrm{Spec}}{\mathcal{O}}^\mathrm{CS}\widetilde{\Phi}_{\psi,-}/\mathrm{Fro}^\mathbb{Z},\underset{\mathrm{Spec}}{\mathcal{O}}^\mathrm{CS}\widetilde{\Delta}^+_{\psi,-}/\mathrm{Fro}^\mathbb{Z},\\
&\underset{\mathrm{Spec}}{\mathcal{O}}^\mathrm{CS}\widetilde{\nabla}^+_{\psi,-}/\mathrm{Fro}^\mathbb{Z},\underset{\mathrm{Spec}}{\mathcal{O}}^\mathrm{CS}\widetilde{\Delta}^\dagger_{\psi,-}/\mathrm{Fro}^\mathbb{Z},\underset{\mathrm{Spec}}{\mathcal{O}}^\mathrm{CS}\widetilde{\nabla}^\dagger_{\psi,-}/\mathrm{Fro}^\mathbb{Z},	\\
\end{align}
\begin{align}
&\underset{\mathrm{Spec}}{\mathcal{O}}^\mathrm{CS}\breve{\Delta}_{\psi,-}/\mathrm{Fro}^\mathbb{Z},\breve{\nabla}_{\psi,-}/\mathrm{Fro}^\mathbb{Z},\underset{\mathrm{Spec}}{\mathcal{O}}^\mathrm{CS}\breve{\Phi}_{\psi,-}/\mathrm{Fro}^\mathbb{Z},\underset{\mathrm{Spec}}{\mathcal{O}}^\mathrm{CS}\breve{\Delta}^+_{\psi,-}/\mathrm{Fro}^\mathbb{Z},\\
&\underset{\mathrm{Spec}}{\mathcal{O}}^\mathrm{CS}\breve{\nabla}^+_{\psi,-}/\mathrm{Fro}^\mathbb{Z},\underset{\mathrm{Spec}}{\mathcal{O}}^\mathrm{CS}\breve{\Delta}^\dagger_{\psi,-}/\mathrm{Fro}^\mathbb{Z},\underset{\mathrm{Spec}}{\mathcal{O}}^\mathrm{CS}\breve{\nabla}^\dagger_{\psi,-}/\mathrm{Fro}^\mathbb{Z},	\\
\end{align}
\begin{align}
&\underset{\mathrm{Spec}}{\mathcal{O}}^\mathrm{CS}{\Delta}_{\psi,-}/\mathrm{Fro}^\mathbb{Z},\underset{\mathrm{Spec}}{\mathcal{O}}^\mathrm{CS}{\nabla}_{\psi,-}/\mathrm{Fro}^\mathbb{Z},\underset{\mathrm{Spec}}{\mathcal{O}}^\mathrm{CS}{\Phi}_{\psi,-}/\mathrm{Fro}^\mathbb{Z},\underset{\mathrm{Spec}}{\mathcal{O}}^\mathrm{CS}{\Delta}^+_{\psi,-}/\mathrm{Fro}^\mathbb{Z},\\
&\underset{\mathrm{Spec}}{\mathcal{O}}^\mathrm{CS}{\nabla}^+_{\psi,-}/\mathrm{Fro}^\mathbb{Z}, \underset{\mathrm{Spec}}{\mathcal{O}}^\mathrm{CS}{\Delta}^\dagger_{\psi,-}/\mathrm{Fro}^\mathbb{Z},\underset{\mathrm{Spec}}{\mathcal{O}}^\mathrm{CS}{\nabla}^\dagger_{\psi,-}/\mathrm{Fro}^\mathbb{Z}.	
\end{align}
Here for those space with notations related to the radius and the corresponding interval we consider the total unions $\bigcap_r,\bigcup_I$ in order to achieve the whole spaces to achieve the analogues of the corresponding FF curves from \cite{10KL1}, \cite{10KL2}, \cite{10FF} for
\[
\xymatrix@R+0pc@C+0pc{
\underset{r}{\mathrm{homotopycolimit}}~\underset{\mathrm{Spec}}{\mathcal{O}}^\mathrm{CS}\widetilde{\Phi}^r_{\psi,-},\underset{I}{\mathrm{homotopylimit}}~\underset{\mathrm{Spec}}{\mathcal{O}}^\mathrm{CS}\widetilde{\Phi}^I_{\psi,-},	\\
}
\]
\[
\xymatrix@R+0pc@C+0pc{
\underset{r}{\mathrm{homotopycolimit}}~\underset{\mathrm{Spec}}{\mathcal{O}}^\mathrm{CS}\breve{\Phi}^r_{\psi,-},\underset{I}{\mathrm{homotopylimit}}~\underset{\mathrm{Spec}}{\mathcal{O}}^\mathrm{CS}\breve{\Phi}^I_{\psi,-},	\\
}
\]
\[
\xymatrix@R+0pc@C+0pc{
\underset{r}{\mathrm{homotopycolimit}}~\underset{\mathrm{Spec}}{\mathcal{O}}^\mathrm{CS}{\Phi}^r_{\psi,-},\underset{I}{\mathrm{homotopylimit}}~\underset{\mathrm{Spec}}{\mathcal{O}}^\mathrm{CS}{\Phi}^I_{\psi,-}.	
}
\]
\[ 
\xymatrix@R+0pc@C+0pc{
\underset{r}{\mathrm{homotopycolimit}}~\underset{\mathrm{Spec}}{\mathcal{O}}^\mathrm{CS}\widetilde{\Phi}^r_{\psi,-}/\mathrm{Fro}^\mathbb{Z},\underset{I}{\mathrm{homotopylimit}}~\underset{\mathrm{Spec}}{\mathcal{O}}^\mathrm{CS}\widetilde{\Phi}^I_{\psi,-}/\mathrm{Fro}^\mathbb{Z},	\\
}
\]
\[ 
\xymatrix@R+0pc@C+0pc{
\underset{r}{\mathrm{homotopycolimit}}~\underset{\mathrm{Spec}}{\mathcal{O}}^\mathrm{CS}\breve{\Phi}^r_{\psi,-}/\mathrm{Fro}^\mathbb{Z},\underset{I}{\mathrm{homotopylimit}}~\breve{\Phi}^I_{\psi,-}/\mathrm{Fro}^\mathbb{Z},	\\
}
\]
\[ 
\xymatrix@R+0pc@C+0pc{
\underset{r}{\mathrm{homotopycolimit}}~\underset{\mathrm{Spec}}{\mathcal{O}}^\mathrm{CS}{\Phi}^r_{\psi,-}/\mathrm{Fro}^\mathbb{Z},\underset{I}{\mathrm{homotopylimit}}~\underset{\mathrm{Spec}}{\mathcal{O}}^\mathrm{CS}{\Phi}^I_{\psi,-}/\mathrm{Fro}^\mathbb{Z}.	
}
\]

\end{definition}

\begin{proposition}
There is a well-defined functor from the $\infty$-category 
\begin{align}
\mathrm{Quasicoherentpresheaves,Perfectcomplex,Condensed}_{*}	
\end{align}
where $*$ is one of the following spaces:
\begin{align}
&\underset{\mathrm{Spec}}{\mathcal{O}}^\mathrm{CS}\widetilde{\Phi}_{\psi,-}/\mathrm{Fro}^\mathbb{Z},	\\
\end{align}
\begin{align}
&\underset{\mathrm{Spec}}{\mathcal{O}}^\mathrm{CS}\breve{\Phi}_{\psi,-}/\mathrm{Fro}^\mathbb{Z},	\\
\end{align}
\begin{align}
&\underset{\mathrm{Spec}}{\mathcal{O}}^\mathrm{CS}{\Phi}_{\psi,-}/\mathrm{Fro}^\mathbb{Z},	
\end{align}
to the $\infty$-category of $\mathrm{Fro}$-equivariant quasicoherent presheaves over similar spaces above correspondingly without the $\mathrm{Fro}$-quotients, and to the $\infty$-category of $\mathrm{Fro}$-equivariant quasicoherent modules over global sections of the structure $\infty$-sheaves of the similar spaces above correspondingly without the $\mathrm{Fro}$-quotients. Here for those space without notation related to the radius and the corresponding interval we consider the total unions $\bigcap_r,\bigcup_I$ in order to achieve the whole spaces to achieve the analogues of the corresponding FF curves from \cite{10KL1}, \cite{10KL2}, \cite{10FF} for
\[
\xymatrix@R+0pc@C+0pc{
\underset{r}{\mathrm{homotopycolimit}}~\underset{\mathrm{Spec}}{\mathcal{O}}^\mathrm{CS}\widetilde{\Phi}^r_{\psi,-},\underset{I}{\mathrm{homotopylimit}}~\underset{\mathrm{Spec}}{\mathcal{O}}^\mathrm{CS}\widetilde{\Phi}^I_{\psi,-},	\\
}
\]
\[
\xymatrix@R+0pc@C+0pc{
\underset{r}{\mathrm{homotopycolimit}}~\underset{\mathrm{Spec}}{\mathcal{O}}^\mathrm{CS}\breve{\Phi}^r_{\psi,-},\underset{I}{\mathrm{homotopylimit}}~\underset{\mathrm{Spec}}{\mathcal{O}}^\mathrm{CS}\breve{\Phi}^I_{\psi,-},	\\
}
\]
\[
\xymatrix@R+0pc@C+0pc{
\underset{r}{\mathrm{homotopycolimit}}~\underset{\mathrm{Spec}}{\mathcal{O}}^\mathrm{CS}{\Phi}^r_{\psi,-},\underset{I}{\mathrm{homotopylimit}}~\underset{\mathrm{Spec}}{\mathcal{O}}^\mathrm{CS}{\Phi}^I_{\psi,-}.	
}
\]
\[ 
\xymatrix@R+0pc@C+0pc{
\underset{r}{\mathrm{homotopycolimit}}~\underset{\mathrm{Spec}}{\mathcal{O}}^\mathrm{CS}\widetilde{\Phi}^r_{\psi,-}/\mathrm{Fro}^\mathbb{Z},\underset{I}{\mathrm{homotopylimit}}~\underset{\mathrm{Spec}}{\mathcal{O}}^\mathrm{CS}\widetilde{\Phi}^I_{\psi,-}/\mathrm{Fro}^\mathbb{Z},	\\
}
\]
\[ 
\xymatrix@R+0pc@C+0pc{
\underset{r}{\mathrm{homotopycolimit}}~\underset{\mathrm{Spec}}{\mathcal{O}}^\mathrm{CS}\breve{\Phi}^r_{\psi,-}/\mathrm{Fro}^\mathbb{Z},\underset{I}{\mathrm{homotopylimit}}~\breve{\Phi}^I_{\psi,-}/\mathrm{Fro}^\mathbb{Z},	\\
}
\]
\[ 
\xymatrix@R+0pc@C+0pc{
\underset{r}{\mathrm{homotopycolimit}}~\underset{\mathrm{Spec}}{\mathcal{O}}^\mathrm{CS}{\Phi}^r_{\psi,-}/\mathrm{Fro}^\mathbb{Z},\underset{I}{\mathrm{homotopylimit}}~\underset{\mathrm{Spec}}{\mathcal{O}}^\mathrm{CS}{\Phi}^I_{\psi,-}/\mathrm{Fro}^\mathbb{Z}.	
}
\]	
In this situation we will have the target category being family parametrized by $r$ or $I$ in compatible glueing sense as in \cite[Definition 5.4.10]{10KL2}. In this situation for modules parametrized by the intervals we have the equivalence of $\infty$-categories by using \cite[Proposition 12.18]{10CS2}. Here the corresponding quasicoherent Frobenius modules are defined to be the corresponding homotopy colimits and limits of Frobenius modules:
\begin{align}
\underset{r}{\mathrm{homotopycolimit}}~M_r,\\
\underset{I}{\mathrm{homotopylimit}}~M_I,	
\end{align}
where each $M_r$ is a Frobenius-equivariant module over the period ring with respect to some radius $r$ while each $M_I$ is a Frobenius-equivariant module over the period ring with respect to some interval $I$.\\
\end{proposition}

\begin{proposition}
Similar proposition holds for 
\begin{align}
\mathrm{Quasicoherentsheaves,Perfectcomplex,IndBanach}_{*}.	
\end{align}	
\end{proposition}

\section{Multivariate Hodge Iwasawa Prestacks by Deformation}

This chapter follows closely \cite{10T1}, \cite{10T2}, \cite{10T3}, \cite{10T4}, \cite{10T5}, \cite{10T6}, \cite{10KPX}, \cite{10KP}, \cite{10KL1}, \cite{10KL2}, \cite{10BK}, \cite{10BBBK}, \cite{10BBM}, \cite{10KKM}, \cite{10CS1}, \cite{10CS2}, \cite{10CKZ}, \cite{10PZ}, \cite{10BCM}, \cite{10LBV}.

\subsection{Frobenius Quasicoherent Prestacks I}

\begin{definition}
First we consider the Bambozzi-Kremnizer spectrum $\underset{\mathrm{Spec}}{\mathcal{O}}^\mathrm{BK}(*)$ attached to any of those in the above from \cite{10BK} by taking derived rational localization:
\begin{align}
&\underset{\mathrm{Spec}}{\mathcal{O}}^\mathrm{BK}\widetilde{\Phi}_{*,\Gamma,A},\underset{\mathrm{Spec}}{\mathcal{O}}^\mathrm{BK}\widetilde{\Phi}^r_{*,\Gamma,A},\underset{\mathrm{Spec}}{\mathcal{O}}^\mathrm{BK}\widetilde{\Phi}^I_{*,\Gamma,A},	
\end{align}
\begin{align}
&\underset{\mathrm{Spec}}{\mathcal{O}}^\mathrm{BK}\breve{\Phi}_{*,\Gamma,A},\underset{\mathrm{Spec}}{\mathcal{O}}^\mathrm{BK}\breve{\Phi}^r_{*,\Gamma,A},\underset{\mathrm{Spec}}{\mathcal{O}}^\mathrm{BK}\breve{\Phi}^I_{*,\Gamma,A},	
\end{align}
\begin{align}
&\underset{\mathrm{Spec}}{\mathcal{O}}^\mathrm{BK}{\Phi}_{*,\Gamma,A},
\underset{\mathrm{Spec}}{\mathcal{O}}^\mathrm{BK}{\Phi}^r_{*,\Gamma,A},\underset{\mathrm{Spec}}{\mathcal{O}}^\mathrm{BK}{\Phi}^I_{*,\Gamma,A}.	
\end{align}

Then we take the corresponding quotients by using the corresponding Frobenius operators:
\begin{align}
&\underset{\mathrm{Spec}}{\mathcal{O}}^\mathrm{BK}\widetilde{\Phi}_{*,\Gamma,A}/\mathrm{Fro}^\mathbb{Z},	\\
\end{align}
\begin{align}
&\underset{\mathrm{Spec}}{\mathcal{O}}^\mathrm{BK}\breve{\Phi}_{*,\Gamma,A}/\mathrm{Fro}^\mathbb{Z},	\\
\end{align}
\begin{align}
&\underset{\mathrm{Spec}}{\mathcal{O}}^\mathrm{BK}{\Phi}_{*,\Gamma,A}/\mathrm{Fro}^\mathbb{Z}.	
\end{align}
Here for those space without notation related to the radius and the corresponding interval we consider the total unions $\bigcap_r,\bigcup_I$ in order to achieve the whole spaces to achieve the analogues of the corresponding FF curves from \cite{10KL1}, \cite{10KL2}, \cite{10FF} for
\[
\xymatrix@R+0pc@C+0pc{
\underset{r}{\mathrm{homotopycolimit}}~\underset{\mathrm{Spec}}{\mathcal{O}}^\mathrm{BK}\widetilde{\Phi}^r_{*,\Gamma,A},\underset{I}{\mathrm{homotopylimit}}~\underset{\mathrm{Spec}}{\mathcal{O}}^\mathrm{BK}\widetilde{\Phi}^I_{*,\Gamma,A},	\\
}
\]
\[
\xymatrix@R+0pc@C+0pc{
\underset{r}{\mathrm{homotopycolimit}}~\underset{\mathrm{Spec}}{\mathcal{O}}^\mathrm{BK}\breve{\Phi}^r_{*,\Gamma,A},\underset{I}{\mathrm{homotopylimit}}~\underset{\mathrm{Spec}}{\mathcal{O}}^\mathrm{BK}\breve{\Phi}^I_{*,\Gamma,A},	\\
}
\]
\[
\xymatrix@R+0pc@C+0pc{
\underset{r}{\mathrm{homotopycolimit}}~\underset{\mathrm{Spec}}{\mathcal{O}}^\mathrm{BK}{\Phi}^r_{*,\Gamma,A},\underset{I}{\mathrm{homotopylimit}}~\underset{\mathrm{Spec}}{\mathcal{O}}^\mathrm{BK}{\Phi}^I_{*,\Gamma,A}.	
}
\]
\[  
\xymatrix@R+0pc@C+0pc{
\underset{r}{\mathrm{homotopycolimit}}~\underset{\mathrm{Spec}}{\mathcal{O}}^\mathrm{BK}\widetilde{\Phi}^r_{*,\Gamma,A}/\mathrm{Fro}^\mathbb{Z},\underset{I}{\mathrm{homotopylimit}}~\underset{\mathrm{Spec}}{\mathcal{O}}^\mathrm{BK}\widetilde{\Phi}^I_{*,\Gamma,A}/\mathrm{Fro}^\mathbb{Z},	\\
}
\]
\[ 
\xymatrix@R+0pc@C+0pc{
\underset{r}{\mathrm{homotopycolimit}}~\underset{\mathrm{Spec}}{\mathcal{O}}^\mathrm{BK}\breve{\Phi}^r_{*,\Gamma,A}/\mathrm{Fro}^\mathbb{Z},\underset{I}{\mathrm{homotopylimit}}~\underset{\mathrm{Spec}}{\mathcal{O}}^\mathrm{BK}\breve{\Phi}^I_{*,\Gamma,A}/\mathrm{Fro}^\mathbb{Z},	\\
}
\]
\[ 
\xymatrix@R+0pc@C+0pc{
\underset{r}{\mathrm{homotopycolimit}}~\underset{\mathrm{Spec}}{\mathcal{O}}^\mathrm{BK}{\Phi}^r_{*,\Gamma,A}/\mathrm{Fro}^\mathbb{Z},\underset{I}{\mathrm{homotopylimit}}~\underset{\mathrm{Spec}}{\mathcal{O}}^\mathrm{BK}{\Phi}^I_{*,\Gamma,A}/\mathrm{Fro}^\mathbb{Z}.	
}
\]

\end{definition}

\indent Meanwhile we have the corresponding Clausen-Scholze analytic stacks from \cite{10CS2}, therefore applying their construction we have:

\begin{definition}
Here we define the following products by using the solidified tensor product from \cite{10CS1} and \cite{10CS2}. Namely $A$ will still as above as a Banach ring over $\mathbb{Q}_p$. Then we take solidified tensor product $\overset{\blacksquare}{\otimes}$ of any of the following
\[
\xymatrix@R+0pc@C+0pc{
\widetilde{\Delta}_{*,\Gamma},\widetilde{\nabla}_{*,\Gamma},\widetilde{\Phi}_{*,\Gamma},\widetilde{\Delta}^+_{*,\Gamma},\widetilde{\nabla}^+_{*,\Gamma},\widetilde{\Delta}^\dagger_{*,\Gamma},\widetilde{\nabla}^\dagger_{*,\Gamma},\widetilde{\Phi}^r_{*,\Gamma},\widetilde{\Phi}^I_{*,\Gamma}, 
}
\]

\[
\xymatrix@R+0pc@C+0pc{
\breve{\Delta}_{*,\Gamma},\breve{\nabla}_{*,\Gamma},\breve{\Phi}_{*,\Gamma},\breve{\Delta}^+_{*,\Gamma},\breve{\nabla}^+_{*,\Gamma},\breve{\Delta}^\dagger_{*,\Gamma},\breve{\nabla}^\dagger_{*,\Gamma},\breve{\Phi}^r_{*,\Gamma},\breve{\Phi}^I_{*,\Gamma},	
}
\]

\[
\xymatrix@R+0pc@C+0pc{
{\Delta}_{*,\Gamma},{\nabla}_{*,\Gamma},{\Phi}_{*,\Gamma},{\Delta}^+_{*,\Gamma},{\nabla}^+_{*,\Gamma},{\Delta}^\dagger_{*,\Gamma},{\nabla}^\dagger_{*,\Gamma},{\Phi}^r_{*,\Gamma},{\Phi}^I_{*,\Gamma},	
}
\]  	
with $A$. Then we have the notations:
\[
\xymatrix@R+0pc@C+0pc{
\widetilde{\Delta}_{*,\Gamma,A},\widetilde{\nabla}_{*,\Gamma,A},\widetilde{\Phi}_{*,\Gamma,A},\widetilde{\Delta}^+_{*,\Gamma,A},\widetilde{\nabla}^+_{*,\Gamma,A},\widetilde{\Delta}^\dagger_{*,\Gamma,A},\widetilde{\nabla}^\dagger_{*,\Gamma,A},\widetilde{\Phi}^r_{*,\Gamma,A},\widetilde{\Phi}^I_{*,\Gamma,A}, 
}
\]

\[
\xymatrix@R+0pc@C+0pc{
\breve{\Delta}_{*,\Gamma,A},\breve{\nabla}_{*,\Gamma,A},\breve{\Phi}_{*,\Gamma,A},\breve{\Delta}^+_{*,\Gamma,A},\breve{\nabla}^+_{*,\Gamma,A},\breve{\Delta}^\dagger_{*,\Gamma,A},\breve{\nabla}^\dagger_{*,\Gamma,A},\breve{\Phi}^r_{*,\Gamma,A},\breve{\Phi}^I_{*,\Gamma,A},	
}
\]

\[
\xymatrix@R+0pc@C+0pc{
{\Delta}_{*,\Gamma,A},{\nabla}_{*,\Gamma,A},{\Phi}_{*,\Gamma,A},{\Delta}^+_{*,\Gamma,A},{\nabla}^+_{*,\Gamma,A},{\Delta}^\dagger_{*,\Gamma,A},{\nabla}^\dagger_{*,\Gamma,A},{\Phi}^r_{*,\Gamma,A},{\Phi}^I_{*,\Gamma,A}.	
}
\]
\end{definition}

\begin{definition}
First we consider the Clausen-Scholze spectrum $\underset{\mathrm{Spec}}{\mathcal{O}}^\mathrm{CS}(*)$ attached to any of those in the above from \cite{10CS2} by taking derived rational localization:
\begin{align}
\underset{\mathrm{Spec}}{\mathcal{O}}^\mathrm{CS}\widetilde{\Delta}_{*,\Gamma,A},\underset{\mathrm{Spec}}{\mathcal{O}}^\mathrm{CS}\widetilde{\nabla}_{*,\Gamma,A},\underset{\mathrm{Spec}}{\mathcal{O}}^\mathrm{CS}\widetilde{\Phi}_{*,\Gamma,A},\underset{\mathrm{Spec}}{\mathcal{O}}^\mathrm{CS}\widetilde{\Delta}^+_{*,\Gamma,A},\underset{\mathrm{Spec}}{\mathcal{O}}^\mathrm{CS}\widetilde{\nabla}^+_{*,\Gamma,A},\\
\underset{\mathrm{Spec}}{\mathcal{O}}^\mathrm{CS}\widetilde{\Delta}^\dagger_{*,\Gamma,A},\underset{\mathrm{Spec}}{\mathcal{O}}^\mathrm{CS}\widetilde{\nabla}^\dagger_{*,\Gamma,A},\underset{\mathrm{Spec}}{\mathcal{O}}^\mathrm{CS}\widetilde{\Phi}^r_{*,\Gamma,A},\underset{\mathrm{Spec}}{\mathcal{O}}^\mathrm{CS}\widetilde{\Phi}^I_{*,\Gamma,A},	\\
\end{align}
\begin{align}
\underset{\mathrm{Spec}}{\mathcal{O}}^\mathrm{CS}\breve{\Delta}_{*,\Gamma,A},\breve{\nabla}_{*,\Gamma,A},\underset{\mathrm{Spec}}{\mathcal{O}}^\mathrm{CS}\breve{\Phi}_{*,\Gamma,A},\underset{\mathrm{Spec}}{\mathcal{O}}^\mathrm{CS}\breve{\Delta}^+_{*,\Gamma,A},\underset{\mathrm{Spec}}{\mathcal{O}}^\mathrm{CS}\breve{\nabla}^+_{*,\Gamma,A},\\
\underset{\mathrm{Spec}}{\mathcal{O}}^\mathrm{CS}\breve{\Delta}^\dagger_{*,\Gamma,A},\underset{\mathrm{Spec}}{\mathcal{O}}^\mathrm{CS}\breve{\nabla}^\dagger_{*,\Gamma,A},\underset{\mathrm{Spec}}{\mathcal{O}}^\mathrm{CS}\breve{\Phi}^r_{*,\Gamma,A},\breve{\Phi}^I_{*,\Gamma,A},	\\
\end{align}
\begin{align}
\underset{\mathrm{Spec}}{\mathcal{O}}^\mathrm{CS}{\Delta}_{*,\Gamma,A},\underset{\mathrm{Spec}}{\mathcal{O}}^\mathrm{CS}{\nabla}_{*,\Gamma,A},\underset{\mathrm{Spec}}{\mathcal{O}}^\mathrm{CS}{\Phi}_{*,\Gamma,A},\underset{\mathrm{Spec}}{\mathcal{O}}^\mathrm{CS}{\Delta}^+_{*,\Gamma,A},\underset{\mathrm{Spec}}{\mathcal{O}}^\mathrm{CS}{\nabla}^+_{*,\Gamma,A},\\
\underset{\mathrm{Spec}}{\mathcal{O}}^\mathrm{CS}{\Delta}^\dagger_{*,\Gamma,A},\underset{\mathrm{Spec}}{\mathcal{O}}^\mathrm{CS}{\nabla}^\dagger_{*,\Gamma,A},\underset{\mathrm{Spec}}{\mathcal{O}}^\mathrm{CS}{\Phi}^r_{*,\Gamma,A},\underset{\mathrm{Spec}}{\mathcal{O}}^\mathrm{CS}{\Phi}^I_{*,\Gamma,A}.	
\end{align}

Then we take the corresponding quotients by using the corresponding Frobenius operators:
\begin{align}
&\underset{\mathrm{Spec}}{\mathcal{O}}^\mathrm{CS}\widetilde{\Delta}_{*,\Gamma,A}/\mathrm{Fro}^\mathbb{Z},\underset{\mathrm{Spec}}{\mathcal{O}}^\mathrm{CS}\widetilde{\nabla}_{*,\Gamma,A}/\mathrm{Fro}^\mathbb{Z},\underset{\mathrm{Spec}}{\mathcal{O}}^\mathrm{CS}\widetilde{\Phi}_{*,\Gamma,A}/\mathrm{Fro}^\mathbb{Z},\underset{\mathrm{Spec}}{\mathcal{O}}^\mathrm{CS}\widetilde{\Delta}^+_{*,\Gamma,A}/\mathrm{Fro}^\mathbb{Z},\\
&\underset{\mathrm{Spec}}{\mathcal{O}}^\mathrm{CS}\widetilde{\nabla}^+_{*,\Gamma,A}/\mathrm{Fro}^\mathbb{Z}, \underset{\mathrm{Spec}}{\mathcal{O}}^\mathrm{CS}\widetilde{\Delta}^\dagger_{*,\Gamma,A}/\mathrm{Fro}^\mathbb{Z},\underset{\mathrm{Spec}}{\mathcal{O}}^\mathrm{CS}\widetilde{\nabla}^\dagger_{*,\Gamma,A}/\mathrm{Fro}^\mathbb{Z},	\\
\end{align}
\begin{align}
&\underset{\mathrm{Spec}}{\mathcal{O}}^\mathrm{CS}\breve{\Delta}_{*,\Gamma,A}/\mathrm{Fro}^\mathbb{Z},\breve{\nabla}_{*,\Gamma,A}/\mathrm{Fro}^\mathbb{Z},\underset{\mathrm{Spec}}{\mathcal{O}}^\mathrm{CS}\breve{\Phi}_{*,\Gamma,A}/\mathrm{Fro}^\mathbb{Z},\underset{\mathrm{Spec}}{\mathcal{O}}^\mathrm{CS}\breve{\Delta}^+_{*,\Gamma,A}/\mathrm{Fro}^\mathbb{Z},\\
&\underset{\mathrm{Spec}}{\mathcal{O}}^\mathrm{CS}\breve{\nabla}^+_{*,\Gamma,A}/\mathrm{Fro}^\mathbb{Z}, \underset{\mathrm{Spec}}{\mathcal{O}}^\mathrm{CS}\breve{\Delta}^\dagger_{*,\Gamma,A}/\mathrm{Fro}^\mathbb{Z},\underset{\mathrm{Spec}}{\mathcal{O}}^\mathrm{CS}\breve{\nabla}^\dagger_{*,\Gamma,A}/\mathrm{Fro}^\mathbb{Z},	\\
\end{align}
\begin{align}
&\underset{\mathrm{Spec}}{\mathcal{O}}^\mathrm{CS}{\Delta}_{*,\Gamma,A}/\mathrm{Fro}^\mathbb{Z},\underset{\mathrm{Spec}}{\mathcal{O}}^\mathrm{CS}{\nabla}_{*,\Gamma,A}/\mathrm{Fro}^\mathbb{Z},\underset{\mathrm{Spec}}{\mathcal{O}}^\mathrm{CS}{\Phi}_{*,\Gamma,A}/\mathrm{Fro}^\mathbb{Z},\underset{\mathrm{Spec}}{\mathcal{O}}^\mathrm{CS}{\Delta}^+_{*,\Gamma,A}/\mathrm{Fro}^\mathbb{Z},\\
&\underset{\mathrm{Spec}}{\mathcal{O}}^\mathrm{CS}{\nabla}^+_{*,\Gamma,A}/\mathrm{Fro}^\mathbb{Z}, \underset{\mathrm{Spec}}{\mathcal{O}}^\mathrm{CS}{\Delta}^\dagger_{*,\Gamma,A}/\mathrm{Fro}^\mathbb{Z},\underset{\mathrm{Spec}}{\mathcal{O}}^\mathrm{CS}{\nabla}^\dagger_{*,\Gamma,A}/\mathrm{Fro}^\mathbb{Z}.	
\end{align}
Here for those space with notations related to the radius and the corresponding interval we consider the total unions $\bigcap_r,\bigcup_I$ in order to achieve the whole spaces to achieve the analogues of the corresponding FF curves from \cite{10KL1}, \cite{10KL2}, \cite{10FF} for
\[
\xymatrix@R+0pc@C+0pc{
\underset{r}{\mathrm{homotopycolimit}}~\underset{\mathrm{Spec}}{\mathcal{O}}^\mathrm{CS}\widetilde{\Phi}^r_{*,\Gamma,A},\underset{I}{\mathrm{homotopylimit}}~\underset{\mathrm{Spec}}{\mathcal{O}}^\mathrm{CS}\widetilde{\Phi}^I_{*,\Gamma,A},	\\
}
\]
\[
\xymatrix@R+0pc@C+0pc{
\underset{r}{\mathrm{homotopycolimit}}~\underset{\mathrm{Spec}}{\mathcal{O}}^\mathrm{CS}\breve{\Phi}^r_{*,\Gamma,A},\underset{I}{\mathrm{homotopylimit}}~\underset{\mathrm{Spec}}{\mathcal{O}}^\mathrm{CS}\breve{\Phi}^I_{*,\Gamma,A},	\\
}
\]
\[
\xymatrix@R+0pc@C+0pc{
\underset{r}{\mathrm{homotopycolimit}}~\underset{\mathrm{Spec}}{\mathcal{O}}^\mathrm{CS}{\Phi}^r_{*,\Gamma,A},\underset{I}{\mathrm{homotopylimit}}~\underset{\mathrm{Spec}}{\mathcal{O}}^\mathrm{CS}{\Phi}^I_{*,\Gamma,A}.	
}
\]
\[ 
\xymatrix@R+0pc@C+0pc{
\underset{r}{\mathrm{homotopycolimit}}~\underset{\mathrm{Spec}}{\mathcal{O}}^\mathrm{CS}\widetilde{\Phi}^r_{*,\Gamma,A}/\mathrm{Fro}^\mathbb{Z},\underset{I}{\mathrm{homotopylimit}}~\underset{\mathrm{Spec}}{\mathcal{O}}^\mathrm{CS}\widetilde{\Phi}^I_{*,\Gamma,A}/\mathrm{Fro}^\mathbb{Z},	\\
}
\]
\[ 
\xymatrix@R+0pc@C+0pc{
\underset{r}{\mathrm{homotopycolimit}}~\underset{\mathrm{Spec}}{\mathcal{O}}^\mathrm{CS}\breve{\Phi}^r_{*,\Gamma,A}/\mathrm{Fro}^\mathbb{Z},\underset{I}{\mathrm{homotopylimit}}~\breve{\Phi}^I_{*,\Gamma,A}/\mathrm{Fro}^\mathbb{Z},	\\
}
\]
\[ 
\xymatrix@R+0pc@C+0pc{
\underset{r}{\mathrm{homotopycolimit}}~\underset{\mathrm{Spec}}{\mathcal{O}}^\mathrm{CS}{\Phi}^r_{*,\Gamma,A}/\mathrm{Fro}^\mathbb{Z},\underset{I}{\mathrm{homotopylimit}}~\underset{\mathrm{Spec}}{\mathcal{O}}^\mathrm{CS}{\Phi}^I_{*,\Gamma,A}/\mathrm{Fro}^\mathbb{Z}.	
}
\]

\end{definition}

\

\begin{definition}
We then consider the corresponding quasipresheaves of the corresponding ind-Banach or monomorphic ind-Banach modules from \cite{10BBK}, \cite{10KKM}:
\begin{align}
\mathrm{Quasicoherentpresheaves,IndBanach}_{*}	
\end{align}
where $*$ is one of the following spaces:
\begin{align}
&\underset{\mathrm{Spec}}{\mathcal{O}}^\mathrm{BK}\widetilde{\Phi}_{*,\Gamma,A}/\mathrm{Fro}^\mathbb{Z},	\\
\end{align}
\begin{align}
&\underset{\mathrm{Spec}}{\mathcal{O}}^\mathrm{BK}\breve{\Phi}_{*,\Gamma,A}/\mathrm{Fro}^\mathbb{Z},	\\
\end{align}
\begin{align}
&\underset{\mathrm{Spec}}{\mathcal{O}}^\mathrm{BK}{\Phi}_{*,\Gamma,A}/\mathrm{Fro}^\mathbb{Z}.	
\end{align}
Here for those space without notation related to the radius and the corresponding interval we consider the total unions $\bigcap_r,\bigcup_I$ in order to achieve the whole spaces to achieve the analogues of the corresponding FF curves from \cite{10KL1}, \cite{10KL2}, \cite{10FF} for
\[
\xymatrix@R+0pc@C+0pc{
\underset{r}{\mathrm{homotopycolimit}}~\underset{\mathrm{Spec}}{\mathcal{O}}^\mathrm{BK}\widetilde{\Phi}^r_{*,\Gamma,A},\underset{I}{\mathrm{homotopylimit}}~\underset{\mathrm{Spec}}{\mathcal{O}}^\mathrm{BK}\widetilde{\Phi}^I_{*,\Gamma,A},	\\
}
\]
\[
\xymatrix@R+0pc@C+0pc{
\underset{r}{\mathrm{homotopycolimit}}~\underset{\mathrm{Spec}}{\mathcal{O}}^\mathrm{BK}\breve{\Phi}^r_{*,\Gamma,A},\underset{I}{\mathrm{homotopylimit}}~\underset{\mathrm{Spec}}{\mathcal{O}}^\mathrm{BK}\breve{\Phi}^I_{*,\Gamma,A},	\\
}
\]
\[
\xymatrix@R+0pc@C+0pc{
\underset{r}{\mathrm{homotopycolimit}}~\underset{\mathrm{Spec}}{\mathcal{O}}^\mathrm{BK}{\Phi}^r_{*,\Gamma,A},\underset{I}{\mathrm{homotopylimit}}~\underset{\mathrm{Spec}}{\mathcal{O}}^\mathrm{BK}{\Phi}^I_{*,\Gamma,A}.	
}
\]
\[  
\xymatrix@R+0pc@C+0pc{
\underset{r}{\mathrm{homotopycolimit}}~\underset{\mathrm{Spec}}{\mathcal{O}}^\mathrm{BK}\widetilde{\Phi}^r_{*,\Gamma,A}/\mathrm{Fro}^\mathbb{Z},\underset{I}{\mathrm{homotopylimit}}~\underset{\mathrm{Spec}}{\mathcal{O}}^\mathrm{BK}\widetilde{\Phi}^I_{*,\Gamma,A}/\mathrm{Fro}^\mathbb{Z},	\\
}
\]
\[ 
\xymatrix@R+0pc@C+0pc{
\underset{r}{\mathrm{homotopycolimit}}~\underset{\mathrm{Spec}}{\mathcal{O}}^\mathrm{BK}\breve{\Phi}^r_{*,\Gamma,A}/\mathrm{Fro}^\mathbb{Z},\underset{I}{\mathrm{homotopylimit}}~\underset{\mathrm{Spec}}{\mathcal{O}}^\mathrm{BK}\breve{\Phi}^I_{*,\Gamma,A}/\mathrm{Fro}^\mathbb{Z},	\\
}
\]
\[ 
\xymatrix@R+0pc@C+0pc{
\underset{r}{\mathrm{homotopycolimit}}~\underset{\mathrm{Spec}}{\mathcal{O}}^\mathrm{BK}{\Phi}^r_{*,\Gamma,A}/\mathrm{Fro}^\mathbb{Z},\underset{I}{\mathrm{homotopylimit}}~\underset{\mathrm{Spec}}{\mathcal{O}}^\mathrm{BK}{\Phi}^I_{*,\Gamma,A}/\mathrm{Fro}^\mathbb{Z}.	
}
\]

\end{definition}

\begin{definition}
We then consider the corresponding quasisheaves of the corresponding condensed solid topological modules from \cite{10CS2}:
\begin{align}
\mathrm{Quasicoherentsheaves, Condensed}_{*}	
\end{align}
where $*$ is one of the following spaces:
\begin{align}
&\underset{\mathrm{Spec}}{\mathcal{O}}^\mathrm{CS}\widetilde{\Delta}_{*,\Gamma,A}/\mathrm{Fro}^\mathbb{Z},\underset{\mathrm{Spec}}{\mathcal{O}}^\mathrm{CS}\widetilde{\nabla}_{*,\Gamma,A}/\mathrm{Fro}^\mathbb{Z},\underset{\mathrm{Spec}}{\mathcal{O}}^\mathrm{CS}\widetilde{\Phi}_{*,\Gamma,A}/\mathrm{Fro}^\mathbb{Z},\underset{\mathrm{Spec}}{\mathcal{O}}^\mathrm{CS}\widetilde{\Delta}^+_{*,\Gamma,A}/\mathrm{Fro}^\mathbb{Z},\\
&\underset{\mathrm{Spec}}{\mathcal{O}}^\mathrm{CS}\widetilde{\nabla}^+_{*,\Gamma,A}/\mathrm{Fro}^\mathbb{Z},\underset{\mathrm{Spec}}{\mathcal{O}}^\mathrm{CS}\widetilde{\Delta}^\dagger_{*,\Gamma,A}/\mathrm{Fro}^\mathbb{Z},\underset{\mathrm{Spec}}{\mathcal{O}}^\mathrm{CS}\widetilde{\nabla}^\dagger_{*,\Gamma,A}/\mathrm{Fro}^\mathbb{Z},	\\
\end{align}
\begin{align}
&\underset{\mathrm{Spec}}{\mathcal{O}}^\mathrm{CS}\breve{\Delta}_{*,\Gamma,A}/\mathrm{Fro}^\mathbb{Z},\breve{\nabla}_{*,\Gamma,A}/\mathrm{Fro}^\mathbb{Z},\underset{\mathrm{Spec}}{\mathcal{O}}^\mathrm{CS}\breve{\Phi}_{*,\Gamma,A}/\mathrm{Fro}^\mathbb{Z},\underset{\mathrm{Spec}}{\mathcal{O}}^\mathrm{CS}\breve{\Delta}^+_{*,\Gamma,A}/\mathrm{Fro}^\mathbb{Z},\\
&\underset{\mathrm{Spec}}{\mathcal{O}}^\mathrm{CS}\breve{\nabla}^+_{*,\Gamma,A}/\mathrm{Fro}^\mathbb{Z},\underset{\mathrm{Spec}}{\mathcal{O}}^\mathrm{CS}\breve{\Delta}^\dagger_{*,\Gamma,A}/\mathrm{Fro}^\mathbb{Z},\underset{\mathrm{Spec}}{\mathcal{O}}^\mathrm{CS}\breve{\nabla}^\dagger_{*,\Gamma,A}/\mathrm{Fro}^\mathbb{Z},	\\
\end{align}
\begin{align}
&\underset{\mathrm{Spec}}{\mathcal{O}}^\mathrm{CS}{\Delta}_{*,\Gamma,A}/\mathrm{Fro}^\mathbb{Z},\underset{\mathrm{Spec}}{\mathcal{O}}^\mathrm{CS}{\nabla}_{*,\Gamma,A}/\mathrm{Fro}^\mathbb{Z},\underset{\mathrm{Spec}}{\mathcal{O}}^\mathrm{CS}{\Phi}_{*,\Gamma,A}/\mathrm{Fro}^\mathbb{Z},\underset{\mathrm{Spec}}{\mathcal{O}}^\mathrm{CS}{\Delta}^+_{*,\Gamma,A}/\mathrm{Fro}^\mathbb{Z},\\
&\underset{\mathrm{Spec}}{\mathcal{O}}^\mathrm{CS}{\nabla}^+_{*,\Gamma,A}/\mathrm{Fro}^\mathbb{Z}, \underset{\mathrm{Spec}}{\mathcal{O}}^\mathrm{CS}{\Delta}^\dagger_{*,\Gamma,A}/\mathrm{Fro}^\mathbb{Z},\underset{\mathrm{Spec}}{\mathcal{O}}^\mathrm{CS}{\nabla}^\dagger_{*,\Gamma,A}/\mathrm{Fro}^\mathbb{Z}.	
\end{align}
Here for those space with notations related to the radius and the corresponding interval we consider the total unions $\bigcap_r,\bigcup_I$ in order to achieve the whole spaces to achieve the analogues of the corresponding FF curves from \cite{10KL1}, \cite{10KL2}, \cite{10FF} for
\[
\xymatrix@R+0pc@C+0pc{
\underset{r}{\mathrm{homotopycolimit}}~\underset{\mathrm{Spec}}{\mathcal{O}}^\mathrm{CS}\widetilde{\Phi}^r_{*,\Gamma,A},\underset{I}{\mathrm{homotopylimit}}~\underset{\mathrm{Spec}}{\mathcal{O}}^\mathrm{CS}\widetilde{\Phi}^I_{*,\Gamma,A},	\\
}
\]
\[
\xymatrix@R+0pc@C+0pc{
\underset{r}{\mathrm{homotopycolimit}}~\underset{\mathrm{Spec}}{\mathcal{O}}^\mathrm{CS}\breve{\Phi}^r_{*,\Gamma,A},\underset{I}{\mathrm{homotopylimit}}~\underset{\mathrm{Spec}}{\mathcal{O}}^\mathrm{CS}\breve{\Phi}^I_{*,\Gamma,A},	\\
}
\]
\[
\xymatrix@R+0pc@C+0pc{
\underset{r}{\mathrm{homotopycolimit}}~\underset{\mathrm{Spec}}{\mathcal{O}}^\mathrm{CS}{\Phi}^r_{*,\Gamma,A},\underset{I}{\mathrm{homotopylimit}}~\underset{\mathrm{Spec}}{\mathcal{O}}^\mathrm{CS}{\Phi}^I_{*,\Gamma,A}.	
}
\]
\[ 
\xymatrix@R+0pc@C+0pc{
\underset{r}{\mathrm{homotopycolimit}}~\underset{\mathrm{Spec}}{\mathcal{O}}^\mathrm{CS}\widetilde{\Phi}^r_{*,\Gamma,A}/\mathrm{Fro}^\mathbb{Z},\underset{I}{\mathrm{homotopylimit}}~\underset{\mathrm{Spec}}{\mathcal{O}}^\mathrm{CS}\widetilde{\Phi}^I_{*,\Gamma,A}/\mathrm{Fro}^\mathbb{Z},	\\
}
\]
\[ 
\xymatrix@R+0pc@C+0pc{
\underset{r}{\mathrm{homotopycolimit}}~\underset{\mathrm{Spec}}{\mathcal{O}}^\mathrm{CS}\breve{\Phi}^r_{*,\Gamma,A}/\mathrm{Fro}^\mathbb{Z},\underset{I}{\mathrm{homotopylimit}}~\breve{\Phi}^I_{*,\Gamma,A}/\mathrm{Fro}^\mathbb{Z},	\\
}
\]
\[ 
\xymatrix@R+0pc@C+0pc{
\underset{r}{\mathrm{homotopycolimit}}~\underset{\mathrm{Spec}}{\mathcal{O}}^\mathrm{CS}{\Phi}^r_{*,\Gamma,A}/\mathrm{Fro}^\mathbb{Z},\underset{I}{\mathrm{homotopylimit}}~\underset{\mathrm{Spec}}{\mathcal{O}}^\mathrm{CS}{\Phi}^I_{*,\Gamma,A}/\mathrm{Fro}^\mathbb{Z}.	
}
\]

\end{definition}

\

\begin{proposition}
There is a well-defined functor from the $\infty$-category 
\begin{align}
\mathrm{Quasicoherentpresheaves,Condensed}_{*}	
\end{align}
where $*$ is one of the following spaces:
\begin{align}
&\underset{\mathrm{Spec}}{\mathcal{O}}^\mathrm{CS}\widetilde{\Phi}_{*,\Gamma,A}/\mathrm{Fro}^\mathbb{Z},	\\
\end{align}
\begin{align}
&\underset{\mathrm{Spec}}{\mathcal{O}}^\mathrm{CS}\breve{\Phi}_{*,\Gamma,A}/\mathrm{Fro}^\mathbb{Z},	\\
\end{align}
\begin{align}
&\underset{\mathrm{Spec}}{\mathcal{O}}^\mathrm{CS}{\Phi}_{*,\Gamma,A}/\mathrm{Fro}^\mathbb{Z},	
\end{align}
to the $\infty$-category of $\mathrm{Fro}$-equivariant quasicoherent presheaves over similar spaces above correspondingly without the $\mathrm{Fro}$-quotients, and to the $\infty$-category of $\mathrm{Fro}$-equivariant quasicoherent modules over global sections of the structure $\infty$-sheaves of the similar spaces above correspondingly without the $\mathrm{Fro}$-quotients. Here for those space without notation related to the radius and the corresponding interval we consider the total unions $\bigcap_r,\bigcup_I$ in order to achieve the whole spaces to achieve the analogues of the corresponding FF curves from \cite{10KL1}, \cite{10KL2}, \cite{10FF} for
\[
\xymatrix@R+0pc@C+0pc{
\underset{r}{\mathrm{homotopycolimit}}~\underset{\mathrm{Spec}}{\mathcal{O}}^\mathrm{CS}\widetilde{\Phi}^r_{*,\Gamma,A},\underset{I}{\mathrm{homotopylimit}}~\underset{\mathrm{Spec}}{\mathcal{O}}^\mathrm{CS}\widetilde{\Phi}^I_{*,\Gamma,A},	\\
}
\]
\[
\xymatrix@R+0pc@C+0pc{
\underset{r}{\mathrm{homotopycolimit}}~\underset{\mathrm{Spec}}{\mathcal{O}}^\mathrm{CS}\breve{\Phi}^r_{*,\Gamma,A},\underset{I}{\mathrm{homotopylimit}}~\underset{\mathrm{Spec}}{\mathcal{O}}^\mathrm{CS}\breve{\Phi}^I_{*,\Gamma,A},	\\
}
\]
\[
\xymatrix@R+0pc@C+0pc{
\underset{r}{\mathrm{homotopycolimit}}~\underset{\mathrm{Spec}}{\mathcal{O}}^\mathrm{CS}{\Phi}^r_{*,\Gamma,A},\underset{I}{\mathrm{homotopylimit}}~\underset{\mathrm{Spec}}{\mathcal{O}}^\mathrm{CS}{\Phi}^I_{*,\Gamma,A}.	
}
\]
\[ 
\xymatrix@R+0pc@C+0pc{
\underset{r}{\mathrm{homotopycolimit}}~\underset{\mathrm{Spec}}{\mathcal{O}}^\mathrm{CS}\widetilde{\Phi}^r_{*,\Gamma,A}/\mathrm{Fro}^\mathbb{Z},\underset{I}{\mathrm{homotopylimit}}~\underset{\mathrm{Spec}}{\mathcal{O}}^\mathrm{CS}\widetilde{\Phi}^I_{*,\Gamma,A}/\mathrm{Fro}^\mathbb{Z},	\\
}
\]
\[ 
\xymatrix@R+0pc@C+0pc{
\underset{r}{\mathrm{homotopycolimit}}~\underset{\mathrm{Spec}}{\mathcal{O}}^\mathrm{CS}\breve{\Phi}^r_{*,\Gamma,A}/\mathrm{Fro}^\mathbb{Z},\underset{I}{\mathrm{homotopylimit}}~\breve{\Phi}^I_{*,\Gamma,A}/\mathrm{Fro}^\mathbb{Z},	\\
}
\]
\[ 
\xymatrix@R+0pc@C+0pc{
\underset{r}{\mathrm{homotopycolimit}}~\underset{\mathrm{Spec}}{\mathcal{O}}^\mathrm{CS}{\Phi}^r_{*,\Gamma,A}/\mathrm{Fro}^\mathbb{Z},\underset{I}{\mathrm{homotopylimit}}~\underset{\mathrm{Spec}}{\mathcal{O}}^\mathrm{CS}{\Phi}^I_{*,\Gamma,A}/\mathrm{Fro}^\mathbb{Z}.	
}
\]	
In this situation we will have the target category being family parametrized by $r$ or $I$ in compatible glueing sense as in \cite[Definition 5.4.10]{10KL2}. In this situation for modules parametrized by the intervals we have the equivalence of $\infty$-categories by using \cite[Proposition 13.8]{10CS2}. Here the corresponding quasicoherent Frobenius modules are defined to be the corresponding homotopy colimits and limits of Frobenius modules:
\begin{align}
\underset{r}{\mathrm{homotopycolimit}}~M_r,\\
\underset{I}{\mathrm{homotopylimit}}~M_I,	
\end{align}
where each $M_r$ is a Frobenius-equivariant module over the period ring with respect to some radius $r$ while each $M_I$ is a Frobenius-equivariant module over the period ring with respect to some interval $I$.\\
\end{proposition}

\begin{proposition}
Similar proposition holds for 
\begin{align}
\mathrm{Quasicoherentsheaves,IndBanach}_{*}.	
\end{align}	
\end{proposition}

\

\begin{definition}
We then consider the corresponding quasipresheaves of perfect complexes the corresponding ind-Banach or monomorphic ind-Banach modules from \cite{10BBK}, \cite{10KKM}:
\begin{align}
\mathrm{Quasicoherentpresheaves,Perfectcomplex,IndBanach}_{*}	
\end{align}
where $*$ is one of the following spaces:
\begin{align}
&\underset{\mathrm{Spec}}{\mathcal{O}}^\mathrm{BK}\widetilde{\Phi}_{*,\Gamma,A}/\mathrm{Fro}^\mathbb{Z},	\\
\end{align}
\begin{align}
&\underset{\mathrm{Spec}}{\mathcal{O}}^\mathrm{BK}\breve{\Phi}_{*,\Gamma,A}/\mathrm{Fro}^\mathbb{Z},	\\
\end{align}
\begin{align}
&\underset{\mathrm{Spec}}{\mathcal{O}}^\mathrm{BK}{\Phi}_{*,\Gamma,A}/\mathrm{Fro}^\mathbb{Z}.	
\end{align}
Here for those space without notation related to the radius and the corresponding interval we consider the total unions $\bigcap_r,\bigcup_I$ in order to achieve the whole spaces to achieve the analogues of the corresponding FF curves from \cite{10KL1}, \cite{10KL2}, \cite{10FF} for
\[
\xymatrix@R+0pc@C+0pc{
\underset{r}{\mathrm{homotopycolimit}}~\underset{\mathrm{Spec}}{\mathcal{O}}^\mathrm{BK}\widetilde{\Phi}^r_{*,\Gamma,A},\underset{I}{\mathrm{homotopylimit}}~\underset{\mathrm{Spec}}{\mathcal{O}}^\mathrm{BK}\widetilde{\Phi}^I_{*,\Gamma,A},	\\
}
\]
\[
\xymatrix@R+0pc@C+0pc{
\underset{r}{\mathrm{homotopycolimit}}~\underset{\mathrm{Spec}}{\mathcal{O}}^\mathrm{BK}\breve{\Phi}^r_{*,\Gamma,A},\underset{I}{\mathrm{homotopylimit}}~\underset{\mathrm{Spec}}{\mathcal{O}}^\mathrm{BK}\breve{\Phi}^I_{*,\Gamma,A},	\\
}
\]
\[
\xymatrix@R+0pc@C+0pc{
\underset{r}{\mathrm{homotopycolimit}}~\underset{\mathrm{Spec}}{\mathcal{O}}^\mathrm{BK}{\Phi}^r_{*,\Gamma,A},\underset{I}{\mathrm{homotopylimit}}~\underset{\mathrm{Spec}}{\mathcal{O}}^\mathrm{BK}{\Phi}^I_{*,\Gamma,A}.	
}
\]
\[  
\xymatrix@R+0pc@C+0pc{
\underset{r}{\mathrm{homotopycolimit}}~\underset{\mathrm{Spec}}{\mathcal{O}}^\mathrm{BK}\widetilde{\Phi}^r_{*,\Gamma,A}/\mathrm{Fro}^\mathbb{Z},\underset{I}{\mathrm{homotopylimit}}~\underset{\mathrm{Spec}}{\mathcal{O}}^\mathrm{BK}\widetilde{\Phi}^I_{*,\Gamma,A}/\mathrm{Fro}^\mathbb{Z},	\\
}
\]
\[ 
\xymatrix@R+0pc@C+0pc{
\underset{r}{\mathrm{homotopycolimit}}~\underset{\mathrm{Spec}}{\mathcal{O}}^\mathrm{BK}\breve{\Phi}^r_{*,\Gamma,A}/\mathrm{Fro}^\mathbb{Z},\underset{I}{\mathrm{homotopylimit}}~\underset{\mathrm{Spec}}{\mathcal{O}}^\mathrm{BK}\breve{\Phi}^I_{*,\Gamma,A}/\mathrm{Fro}^\mathbb{Z},	\\
}
\]
\[ 
\xymatrix@R+0pc@C+0pc{
\underset{r}{\mathrm{homotopycolimit}}~\underset{\mathrm{Spec}}{\mathcal{O}}^\mathrm{BK}{\Phi}^r_{*,\Gamma,A}/\mathrm{Fro}^\mathbb{Z},\underset{I}{\mathrm{homotopylimit}}~\underset{\mathrm{Spec}}{\mathcal{O}}^\mathrm{BK}{\Phi}^I_{*,\Gamma,A}/\mathrm{Fro}^\mathbb{Z}.	
}
\]

\end{definition}

\begin{definition}
We then consider the corresponding quasisheaves of perfect complexes of the corresponding condensed solid topological modules from \cite{10CS2}:
\begin{align}
\mathrm{Quasicoherentsheaves, Perfectcomplex, Condensed}_{*}	
\end{align}
where $*$ is one of the following spaces:
\begin{align}
&\underset{\mathrm{Spec}}{\mathcal{O}}^\mathrm{CS}\widetilde{\Delta}_{*,\Gamma,A}/\mathrm{Fro}^\mathbb{Z},\underset{\mathrm{Spec}}{\mathcal{O}}^\mathrm{CS}\widetilde{\nabla}_{*,\Gamma,A}/\mathrm{Fro}^\mathbb{Z},\underset{\mathrm{Spec}}{\mathcal{O}}^\mathrm{CS}\widetilde{\Phi}_{*,\Gamma,A}/\mathrm{Fro}^\mathbb{Z},\underset{\mathrm{Spec}}{\mathcal{O}}^\mathrm{CS}\widetilde{\Delta}^+_{*,\Gamma,A}/\mathrm{Fro}^\mathbb{Z},\\
&\underset{\mathrm{Spec}}{\mathcal{O}}^\mathrm{CS}\widetilde{\nabla}^+_{*,\Gamma,A}/\mathrm{Fro}^\mathbb{Z},\underset{\mathrm{Spec}}{\mathcal{O}}^\mathrm{CS}\widetilde{\Delta}^\dagger_{*,\Gamma,A}/\mathrm{Fro}^\mathbb{Z},\underset{\mathrm{Spec}}{\mathcal{O}}^\mathrm{CS}\widetilde{\nabla}^\dagger_{*,\Gamma,A}/\mathrm{Fro}^\mathbb{Z},	\\
\end{align}
\begin{align}
&\underset{\mathrm{Spec}}{\mathcal{O}}^\mathrm{CS}\breve{\Delta}_{*,\Gamma,A}/\mathrm{Fro}^\mathbb{Z},\breve{\nabla}_{*,\Gamma,A}/\mathrm{Fro}^\mathbb{Z},\underset{\mathrm{Spec}}{\mathcal{O}}^\mathrm{CS}\breve{\Phi}_{*,\Gamma,A}/\mathrm{Fro}^\mathbb{Z},\underset{\mathrm{Spec}}{\mathcal{O}}^\mathrm{CS}\breve{\Delta}^+_{*,\Gamma,A}/\mathrm{Fro}^\mathbb{Z},\\
&\underset{\mathrm{Spec}}{\mathcal{O}}^\mathrm{CS}\breve{\nabla}^+_{*,\Gamma,A}/\mathrm{Fro}^\mathbb{Z},\underset{\mathrm{Spec}}{\mathcal{O}}^\mathrm{CS}\breve{\Delta}^\dagger_{*,\Gamma,A}/\mathrm{Fro}^\mathbb{Z},\underset{\mathrm{Spec}}{\mathcal{O}}^\mathrm{CS}\breve{\nabla}^\dagger_{*,\Gamma,A}/\mathrm{Fro}^\mathbb{Z},	\\
\end{align}
\begin{align}
&\underset{\mathrm{Spec}}{\mathcal{O}}^\mathrm{CS}{\Delta}_{*,\Gamma,A}/\mathrm{Fro}^\mathbb{Z},\underset{\mathrm{Spec}}{\mathcal{O}}^\mathrm{CS}{\nabla}_{*,\Gamma,A}/\mathrm{Fro}^\mathbb{Z},\underset{\mathrm{Spec}}{\mathcal{O}}^\mathrm{CS}{\Phi}_{*,\Gamma,A}/\mathrm{Fro}^\mathbb{Z},\underset{\mathrm{Spec}}{\mathcal{O}}^\mathrm{CS}{\Delta}^+_{*,\Gamma,A}/\mathrm{Fro}^\mathbb{Z},\\
&\underset{\mathrm{Spec}}{\mathcal{O}}^\mathrm{CS}{\nabla}^+_{*,\Gamma,A}/\mathrm{Fro}^\mathbb{Z}, \underset{\mathrm{Spec}}{\mathcal{O}}^\mathrm{CS}{\Delta}^\dagger_{*,\Gamma,A}/\mathrm{Fro}^\mathbb{Z},\underset{\mathrm{Spec}}{\mathcal{O}}^\mathrm{CS}{\nabla}^\dagger_{*,\Gamma,A}/\mathrm{Fro}^\mathbb{Z}.	
\end{align}
Here for those space with notations related to the radius and the corresponding interval we consider the total unions $\bigcap_r,\bigcup_I$ in order to achieve the whole spaces to achieve the analogues of the corresponding FF curves from \cite{10KL1}, \cite{10KL2}, \cite{10FF} for
\[
\xymatrix@R+0pc@C+0pc{
\underset{r}{\mathrm{homotopycolimit}}~\underset{\mathrm{Spec}}{\mathcal{O}}^\mathrm{CS}\widetilde{\Phi}^r_{*,\Gamma,A},\underset{I}{\mathrm{homotopylimit}}~\underset{\mathrm{Spec}}{\mathcal{O}}^\mathrm{CS}\widetilde{\Phi}^I_{*,\Gamma,A},	\\
}
\]
\[
\xymatrix@R+0pc@C+0pc{
\underset{r}{\mathrm{homotopycolimit}}~\underset{\mathrm{Spec}}{\mathcal{O}}^\mathrm{CS}\breve{\Phi}^r_{*,\Gamma,A},\underset{I}{\mathrm{homotopylimit}}~\underset{\mathrm{Spec}}{\mathcal{O}}^\mathrm{CS}\breve{\Phi}^I_{*,\Gamma,A},	\\
}
\]
\[
\xymatrix@R+0pc@C+0pc{
\underset{r}{\mathrm{homotopycolimit}}~\underset{\mathrm{Spec}}{\mathcal{O}}^\mathrm{CS}{\Phi}^r_{*,\Gamma,A},\underset{I}{\mathrm{homotopylimit}}~\underset{\mathrm{Spec}}{\mathcal{O}}^\mathrm{CS}{\Phi}^I_{*,\Gamma,A}.	
}
\]
\[ 
\xymatrix@R+0pc@C+0pc{
\underset{r}{\mathrm{homotopycolimit}}~\underset{\mathrm{Spec}}{\mathcal{O}}^\mathrm{CS}\widetilde{\Phi}^r_{*,\Gamma,A}/\mathrm{Fro}^\mathbb{Z},\underset{I}{\mathrm{homotopylimit}}~\underset{\mathrm{Spec}}{\mathcal{O}}^\mathrm{CS}\widetilde{\Phi}^I_{*,\Gamma,A}/\mathrm{Fro}^\mathbb{Z},	\\
}
\]
\[ 
\xymatrix@R+0pc@C+0pc{
\underset{r}{\mathrm{homotopycolimit}}~\underset{\mathrm{Spec}}{\mathcal{O}}^\mathrm{CS}\breve{\Phi}^r_{*,\Gamma,A}/\mathrm{Fro}^\mathbb{Z},\underset{I}{\mathrm{homotopylimit}}~\breve{\Phi}^I_{*,\Gamma,A}/\mathrm{Fro}^\mathbb{Z},	\\
}
\]
\[ 
\xymatrix@R+0pc@C+0pc{
\underset{r}{\mathrm{homotopycolimit}}~\underset{\mathrm{Spec}}{\mathcal{O}}^\mathrm{CS}{\Phi}^r_{*,\Gamma,A}/\mathrm{Fro}^\mathbb{Z},\underset{I}{\mathrm{homotopylimit}}~\underset{\mathrm{Spec}}{\mathcal{O}}^\mathrm{CS}{\Phi}^I_{*,\Gamma,A}/\mathrm{Fro}^\mathbb{Z}.	
}
\]

\end{definition}

\begin{proposition}
There is a well-defined functor from the $\infty$-category 
\begin{align}
\mathrm{Quasicoherentpresheaves,Perfectcomplex,Condensed}_{*}	
\end{align}
where $*$ is one of the following spaces:
\begin{align}
&\underset{\mathrm{Spec}}{\mathcal{O}}^\mathrm{CS}\widetilde{\Phi}_{*,\Gamma,A}/\mathrm{Fro}^\mathbb{Z},	\\
\end{align}
\begin{align}
&\underset{\mathrm{Spec}}{\mathcal{O}}^\mathrm{CS}\breve{\Phi}_{*,\Gamma,A}/\mathrm{Fro}^\mathbb{Z},	\\
\end{align}
\begin{align}
&\underset{\mathrm{Spec}}{\mathcal{O}}^\mathrm{CS}{\Phi}_{*,\Gamma,A}/\mathrm{Fro}^\mathbb{Z},	
\end{align}
to the $\infty$-category of $\mathrm{Fro}$-equivariant quasicoherent presheaves over similar spaces above correspondingly without the $\mathrm{Fro}$-quotients, and to the $\infty$-category of $\mathrm{Fro}$-equivariant quasicoherent modules over global sections of the structure $\infty$-sheaves of the similar spaces above correspondingly without the $\mathrm{Fro}$-quotients. Here for those space without notation related to the radius and the corresponding interval we consider the total unions $\bigcap_r,\bigcup_I$ in order to achieve the whole spaces to achieve the analogues of the corresponding FF curves from \cite{10KL1}, \cite{10KL2}, \cite{10FF} for
\[
\xymatrix@R+0pc@C+0pc{
\underset{r}{\mathrm{homotopycolimit}}~\underset{\mathrm{Spec}}{\mathcal{O}}^\mathrm{CS}\widetilde{\Phi}^r_{*,\Gamma,A},\underset{I}{\mathrm{homotopylimit}}~\underset{\mathrm{Spec}}{\mathcal{O}}^\mathrm{CS}\widetilde{\Phi}^I_{*,\Gamma,A},	\\
}
\]
\[
\xymatrix@R+0pc@C+0pc{
\underset{r}{\mathrm{homotopycolimit}}~\underset{\mathrm{Spec}}{\mathcal{O}}^\mathrm{CS}\breve{\Phi}^r_{*,\Gamma,A},\underset{I}{\mathrm{homotopylimit}}~\underset{\mathrm{Spec}}{\mathcal{O}}^\mathrm{CS}\breve{\Phi}^I_{*,\Gamma,A},	\\
}
\]
\[
\xymatrix@R+0pc@C+0pc{
\underset{r}{\mathrm{homotopycolimit}}~\underset{\mathrm{Spec}}{\mathcal{O}}^\mathrm{CS}{\Phi}^r_{*,\Gamma,A},\underset{I}{\mathrm{homotopylimit}}~\underset{\mathrm{Spec}}{\mathcal{O}}^\mathrm{CS}{\Phi}^I_{*,\Gamma,A}.	
}
\]
\[ 
\xymatrix@R+0pc@C+0pc{
\underset{r}{\mathrm{homotopycolimit}}~\underset{\mathrm{Spec}}{\mathcal{O}}^\mathrm{CS}\widetilde{\Phi}^r_{*,\Gamma,A}/\mathrm{Fro}^\mathbb{Z},\underset{I}{\mathrm{homotopylimit}}~\underset{\mathrm{Spec}}{\mathcal{O}}^\mathrm{CS}\widetilde{\Phi}^I_{*,\Gamma,A}/\mathrm{Fro}^\mathbb{Z},	\\
}
\]
\[ 
\xymatrix@R+0pc@C+0pc{
\underset{r}{\mathrm{homotopycolimit}}~\underset{\mathrm{Spec}}{\mathcal{O}}^\mathrm{CS}\breve{\Phi}^r_{*,\Gamma,A}/\mathrm{Fro}^\mathbb{Z},\underset{I}{\mathrm{homotopylimit}}~\breve{\Phi}^I_{*,\Gamma,A}/\mathrm{Fro}^\mathbb{Z},	\\
}
\]
\[ 
\xymatrix@R+0pc@C+0pc{
\underset{r}{\mathrm{homotopycolimit}}~\underset{\mathrm{Spec}}{\mathcal{O}}^\mathrm{CS}{\Phi}^r_{*,\Gamma,A}/\mathrm{Fro}^\mathbb{Z},\underset{I}{\mathrm{homotopylimit}}~\underset{\mathrm{Spec}}{\mathcal{O}}^\mathrm{CS}{\Phi}^I_{*,\Gamma,A}/\mathrm{Fro}^\mathbb{Z}.	
}
\]	
In this situation we will have the target category being family parametrized by $r$ or $I$ in compatible glueing sense as in \cite[Definition 5.4.10]{10KL2}. In this situation for modules parametrized by the intervals we have the equivalence of $\infty$-categories by using \cite[Proposition 12.18]{10CS2}. Here the corresponding quasicoherent Frobenius modules are defined to be the corresponding homotopy colimits and limits of Frobenius modules:
\begin{align}
\underset{r}{\mathrm{homotopycolimit}}~M_r,\\
\underset{I}{\mathrm{homotopylimit}}~M_I,	
\end{align}
where each $M_r$ is a Frobenius-equivariant module over the period ring with respect to some radius $r$ while each $M_I$ is a Frobenius-equivariant module over the period ring with respect to some interval $I$.\\
\end{proposition}

\begin{proposition}
Similar proposition holds for 
\begin{align}
\mathrm{Quasicoherentsheaves,Perfectcomplex,IndBanach}_{*}.	
\end{align}	
\end{proposition}

\newpage
\subsection{Frobenius Quasicoherent Prestacks II: Deformation in Banach Rings}

\begin{definition}
We now consider the pro-\'etale site of $\mathrm{Spa}\mathbb{Q}_p\left<X_1^{\pm 1},...,X_k^{\pm 1}\right>$, denote that by $*$. To be more accurate we replace one component for $\Gamma$ with the pro-\'etale site of $\mathrm{Spa}\mathbb{Q}_p\left<X_1^{\pm 1},...,X_k^{\pm 1}\right>$. And we treat then all the functor to be prestacks for this site. Then from \cite{10KL1} and \cite[Definition 5.2.1]{10KL2} we have the following class of Kedlaya-Liu rings (with the following replacement: $\Delta$ stands for $A$, $\nabla$ stands for $B$, while $\Phi$ stands for $C$) by taking product in the sense of self $\Gamma$-th power:

\[
\xymatrix@R+0pc@C+0pc{
\widetilde{\Delta}_{*,\Gamma},\widetilde{\nabla}_{*,\Gamma},\widetilde{\Phi}_{*,\Gamma},\widetilde{\Delta}^+_{*,\Gamma},\widetilde{\nabla}^+_{*,\Gamma},\widetilde{\Delta}^\dagger_{*,\Gamma},\widetilde{\nabla}^\dagger_{*,\Gamma},\widetilde{\Phi}^r_{*,\Gamma},\widetilde{\Phi}^I_{*,\Gamma}, 
}
\]

\[
\xymatrix@R+0pc@C+0pc{
\breve{\Delta}_{*,\Gamma},\breve{\nabla}_{*,\Gamma},\breve{\Phi}_{*,\Gamma},\breve{\Delta}^+_{*,\Gamma},\breve{\nabla}^+_{*,\Gamma},\breve{\Delta}^\dagger_{*,\Gamma},\breve{\nabla}^\dagger_{*,\Gamma},\breve{\Phi}^r_{*,\Gamma},\breve{\Phi}^I_{*,\Gamma},	
}
\]

\[
\xymatrix@R+0pc@C+0pc{
{\Delta}_{*,\Gamma},{\nabla}_{*,\Gamma},{\Phi}_{*,\Gamma},{\Delta}^+_{*,\Gamma},{\nabla}^+_{*,\Gamma},{\Delta}^\dagger_{*,\Gamma},{\nabla}^\dagger_{*,\Gamma},{\Phi}^r_{*,\Gamma},{\Phi}^I_{*,\Gamma}.	
}
\]
We now consider the following rings with $-$ being any deforming Banach ring over $\mathbb{Q}_p$. Taking the product we have:
\[
\xymatrix@R+0pc@C+0pc{
\widetilde{\Phi}_{*,\Gamma,-},\widetilde{\Phi}^r_{*,\Gamma,-},\widetilde{\Phi}^I_{*,\Gamma,-},	
}
\]
\[
\xymatrix@R+0pc@C+0pc{
\breve{\Phi}_{*,\Gamma,-},\breve{\Phi}^r_{*,\Gamma,-},\breve{\Phi}^I_{*,\Gamma,-},	
}
\]
\[
\xymatrix@R+0pc@C+0pc{
{\Phi}_{*,\Gamma,-},{\Phi}^r_{*,\Gamma,-},{\Phi}^I_{*,\Gamma,-}.	
}
\]
They carry multi Frobenius action $\varphi_\Gamma$ and multi $\mathrm{Lie}_\Gamma:=\mathbb{Z}_p^{\times\Gamma}$ action. In our current situation after \cite{10CKZ} and \cite{10PZ} we consider the following $(\infty,1)$-categories of $(\infty,1)$-modules.\\
\end{definition}

\begin{definition}
First we consider the Bambozzi-Kremnizer spectrum $\underset{\mathrm{Spec}}{\mathcal{O}}^\mathrm{BK}(*)$ attached to any of those in the above from \cite{10BK} by taking derived rational localization:
\begin{align}
&\underset{\mathrm{Spec}}{\mathcal{O}}^\mathrm{BK}\widetilde{\Phi}_{*,\Gamma,-},\underset{\mathrm{Spec}}{\mathcal{O}}^\mathrm{BK}\widetilde{\Phi}^r_{*,\Gamma,-},\underset{\mathrm{Spec}}{\mathcal{O}}^\mathrm{BK}\widetilde{\Phi}^I_{*,\Gamma,-},	
\end{align}
\begin{align}
&\underset{\mathrm{Spec}}{\mathcal{O}}^\mathrm{BK}\breve{\Phi}_{*,\Gamma,-},\underset{\mathrm{Spec}}{\mathcal{O}}^\mathrm{BK}\breve{\Phi}^r_{*,\Gamma,-},\underset{\mathrm{Spec}}{\mathcal{O}}^\mathrm{BK}\breve{\Phi}^I_{*,\Gamma,-},	
\end{align}
\begin{align}
&\underset{\mathrm{Spec}}{\mathcal{O}}^\mathrm{BK}{\Phi}_{*,\Gamma,-},
\underset{\mathrm{Spec}}{\mathcal{O}}^\mathrm{BK}{\Phi}^r_{*,\Gamma,-},\underset{\mathrm{Spec}}{\mathcal{O}}^\mathrm{BK}{\Phi}^I_{*,\Gamma,-}.	
\end{align}

Then we take the corresponding quotients by using the corresponding Frobenius operators:
\begin{align}
&\underset{\mathrm{Spec}}{\mathcal{O}}^\mathrm{BK}\widetilde{\Phi}_{*,\Gamma,-}/\mathrm{Fro}^\mathbb{Z},	\\
\end{align}
\begin{align}
&\underset{\mathrm{Spec}}{\mathcal{O}}^\mathrm{BK}\breve{\Phi}_{*,\Gamma,-}/\mathrm{Fro}^\mathbb{Z},	\\
\end{align}
\begin{align}
&\underset{\mathrm{Spec}}{\mathcal{O}}^\mathrm{BK}{\Phi}_{*,\Gamma,-}/\mathrm{Fro}^\mathbb{Z}.	
\end{align}
Here for those space without notation related to the radius and the corresponding interval we consider the total unions $\bigcap_r,\bigcup_I$ in order to achieve the whole spaces to achieve the analogues of the corresponding FF curves from \cite{10KL1}, \cite{10KL2}, \cite{10FF} for
\[
\xymatrix@R+0pc@C+0pc{
\underset{r}{\mathrm{homotopycolimit}}~\underset{\mathrm{Spec}}{\mathcal{O}}^\mathrm{BK}\widetilde{\Phi}^r_{*,\Gamma,-},\underset{I}{\mathrm{homotopylimit}}~\underset{\mathrm{Spec}}{\mathcal{O}}^\mathrm{BK}\widetilde{\Phi}^I_{*,\Gamma,-},	\\
}
\]
\[
\xymatrix@R+0pc@C+0pc{
\underset{r}{\mathrm{homotopycolimit}}~\underset{\mathrm{Spec}}{\mathcal{O}}^\mathrm{BK}\breve{\Phi}^r_{*,\Gamma,-},\underset{I}{\mathrm{homotopylimit}}~\underset{\mathrm{Spec}}{\mathcal{O}}^\mathrm{BK}\breve{\Phi}^I_{*,\Gamma,-},	\\
}
\]
\[
\xymatrix@R+0pc@C+0pc{
\underset{r}{\mathrm{homotopycolimit}}~\underset{\mathrm{Spec}}{\mathcal{O}}^\mathrm{BK}{\Phi}^r_{*,\Gamma,-},\underset{I}{\mathrm{homotopylimit}}~\underset{\mathrm{Spec}}{\mathcal{O}}^\mathrm{BK}{\Phi}^I_{*,\Gamma,-}.	
}
\]
\[  
\xymatrix@R+0pc@C+0pc{
\underset{r}{\mathrm{homotopycolimit}}~\underset{\mathrm{Spec}}{\mathcal{O}}^\mathrm{BK}\widetilde{\Phi}^r_{*,\Gamma,-}/\mathrm{Fro}^\mathbb{Z},\underset{I}{\mathrm{homotopylimit}}~\underset{\mathrm{Spec}}{\mathcal{O}}^\mathrm{BK}\widetilde{\Phi}^I_{*,\Gamma,-}/\mathrm{Fro}^\mathbb{Z},	\\
}
\]
\[ 
\xymatrix@R+0pc@C+0pc{
\underset{r}{\mathrm{homotopycolimit}}~\underset{\mathrm{Spec}}{\mathcal{O}}^\mathrm{BK}\breve{\Phi}^r_{*,\Gamma,-}/\mathrm{Fro}^\mathbb{Z},\underset{I}{\mathrm{homotopylimit}}~\underset{\mathrm{Spec}}{\mathcal{O}}^\mathrm{BK}\breve{\Phi}^I_{*,\Gamma,-}/\mathrm{Fro}^\mathbb{Z},	\\
}
\]
\[ 
\xymatrix@R+0pc@C+0pc{
\underset{r}{\mathrm{homotopycolimit}}~\underset{\mathrm{Spec}}{\mathcal{O}}^\mathrm{BK}{\Phi}^r_{*,\Gamma,-}/\mathrm{Fro}^\mathbb{Z},\underset{I}{\mathrm{homotopylimit}}~\underset{\mathrm{Spec}}{\mathcal{O}}^\mathrm{BK}{\Phi}^I_{*,\Gamma,-}/\mathrm{Fro}^\mathbb{Z}.	
}
\]

\end{definition}

\indent Meanwhile we have the corresponding Clausen-Scholze analytic stacks from \cite{10CS2}, therefore applying their construction we have:

\begin{definition}
Here we define the following products by using the solidified tensor product from \cite{10CS1} and \cite{10CS2}. Namely $A$ will still as above as a Banach ring over $\mathbb{Q}_p$. Then we take solidified tensor product $\overset{\blacksquare}{\otimes}$ of any of the following
\[
\xymatrix@R+0pc@C+0pc{
\widetilde{\Delta}_{*,\Gamma},\widetilde{\nabla}_{*,\Gamma},\widetilde{\Phi}_{*,\Gamma},\widetilde{\Delta}^+_{*,\Gamma},\widetilde{\nabla}^+_{*,\Gamma},\widetilde{\Delta}^\dagger_{*,\Gamma},\widetilde{\nabla}^\dagger_{*,\Gamma},\widetilde{\Phi}^r_{*,\Gamma},\widetilde{\Phi}^I_{*,\Gamma}, 
}
\]

\[
\xymatrix@R+0pc@C+0pc{
\breve{\Delta}_{*,\Gamma},\breve{\nabla}_{*,\Gamma},\breve{\Phi}_{*,\Gamma},\breve{\Delta}^+_{*,\Gamma},\breve{\nabla}^+_{*,\Gamma},\breve{\Delta}^\dagger_{*,\Gamma},\breve{\nabla}^\dagger_{*,\Gamma},\breve{\Phi}^r_{*,\Gamma},\breve{\Phi}^I_{*,\Gamma},	
}
\]

\[
\xymatrix@R+0pc@C+0pc{
{\Delta}_{*,\Gamma},{\nabla}_{*,\Gamma},{\Phi}_{*,\Gamma},{\Delta}^+_{*,\Gamma},{\nabla}^+_{*,\Gamma},{\Delta}^\dagger_{*,\Gamma},{\nabla}^\dagger_{*,\Gamma},{\Phi}^r_{*,\Gamma},{\Phi}^I_{*,\Gamma},	
}
\]  	
with $A$. Then we have the notations:
\[
\xymatrix@R+0pc@C+0pc{
\widetilde{\Delta}_{*,\Gamma,-},\widetilde{\nabla}_{*,\Gamma,-},\widetilde{\Phi}_{*,\Gamma,-},\widetilde{\Delta}^+_{*,\Gamma,-},\widetilde{\nabla}^+_{*,\Gamma,-},\widetilde{\Delta}^\dagger_{*,\Gamma,-},\widetilde{\nabla}^\dagger_{*,\Gamma,-},\widetilde{\Phi}^r_{*,\Gamma,-},\widetilde{\Phi}^I_{*,\Gamma,-}, 
}
\]

\[
\xymatrix@R+0pc@C+0pc{
\breve{\Delta}_{*,\Gamma,-},\breve{\nabla}_{*,\Gamma,-},\breve{\Phi}_{*,\Gamma,-},\breve{\Delta}^+_{*,\Gamma,-},\breve{\nabla}^+_{*,\Gamma,-},\breve{\Delta}^\dagger_{*,\Gamma,-},\breve{\nabla}^\dagger_{*,\Gamma,-},\breve{\Phi}^r_{*,\Gamma,-},\breve{\Phi}^I_{*,\Gamma,-},	
}
\]

\[
\xymatrix@R+0pc@C+0pc{
{\Delta}_{*,\Gamma,-},{\nabla}_{*,\Gamma,-},{\Phi}_{*,\Gamma,-},{\Delta}^+_{*,\Gamma,-},{\nabla}^+_{*,\Gamma,-},{\Delta}^\dagger_{*,\Gamma,-},{\nabla}^\dagger_{*,\Gamma,-},{\Phi}^r_{*,\Gamma,-},{\Phi}^I_{*,\Gamma,-}.	
}
\]
\end{definition}

\begin{definition}
First we consider the Clausen-Scholze spectrum $\underset{\mathrm{Spec}}{\mathcal{O}}^\mathrm{CS}(*)$ attached to any of those in the above from \cite{10CS2} by taking derived rational localization:
\begin{align}
\underset{\mathrm{Spec}}{\mathcal{O}}^\mathrm{CS}\widetilde{\Delta}_{*,\Gamma,-},\underset{\mathrm{Spec}}{\mathcal{O}}^\mathrm{CS}\widetilde{\nabla}_{*,\Gamma,-},\underset{\mathrm{Spec}}{\mathcal{O}}^\mathrm{CS}\widetilde{\Phi}_{*,\Gamma,-},\underset{\mathrm{Spec}}{\mathcal{O}}^\mathrm{CS}\widetilde{\Delta}^+_{*,\Gamma,-},\underset{\mathrm{Spec}}{\mathcal{O}}^\mathrm{CS}\widetilde{\nabla}^+_{*,\Gamma,-},\\
\underset{\mathrm{Spec}}{\mathcal{O}}^\mathrm{CS}\widetilde{\Delta}^\dagger_{*,\Gamma,-},\underset{\mathrm{Spec}}{\mathcal{O}}^\mathrm{CS}\widetilde{\nabla}^\dagger_{*,\Gamma,-},\underset{\mathrm{Spec}}{\mathcal{O}}^\mathrm{CS}\widetilde{\Phi}^r_{*,\Gamma,-},\underset{\mathrm{Spec}}{\mathcal{O}}^\mathrm{CS}\widetilde{\Phi}^I_{*,\Gamma,-},	\\
\end{align}
\begin{align}
\underset{\mathrm{Spec}}{\mathcal{O}}^\mathrm{CS}\breve{\Delta}_{*,\Gamma,-},\breve{\nabla}_{*,\Gamma,-},\underset{\mathrm{Spec}}{\mathcal{O}}^\mathrm{CS}\breve{\Phi}_{*,\Gamma,-},\underset{\mathrm{Spec}}{\mathcal{O}}^\mathrm{CS}\breve{\Delta}^+_{*,\Gamma,-},\underset{\mathrm{Spec}}{\mathcal{O}}^\mathrm{CS}\breve{\nabla}^+_{*,\Gamma,-},\\
\underset{\mathrm{Spec}}{\mathcal{O}}^\mathrm{CS}\breve{\Delta}^\dagger_{*,\Gamma,-},\underset{\mathrm{Spec}}{\mathcal{O}}^\mathrm{CS}\breve{\nabla}^\dagger_{*,\Gamma,-},\underset{\mathrm{Spec}}{\mathcal{O}}^\mathrm{CS}\breve{\Phi}^r_{*,\Gamma,-},\breve{\Phi}^I_{*,\Gamma,-},	\\
\end{align}
\begin{align}
\underset{\mathrm{Spec}}{\mathcal{O}}^\mathrm{CS}{\Delta}_{*,\Gamma,-},\underset{\mathrm{Spec}}{\mathcal{O}}^\mathrm{CS}{\nabla}_{*,\Gamma,-},\underset{\mathrm{Spec}}{\mathcal{O}}^\mathrm{CS}{\Phi}_{*,\Gamma,-},\underset{\mathrm{Spec}}{\mathcal{O}}^\mathrm{CS}{\Delta}^+_{*,\Gamma,-},\underset{\mathrm{Spec}}{\mathcal{O}}^\mathrm{CS}{\nabla}^+_{*,\Gamma,-},\\
\underset{\mathrm{Spec}}{\mathcal{O}}^\mathrm{CS}{\Delta}^\dagger_{*,\Gamma,-},\underset{\mathrm{Spec}}{\mathcal{O}}^\mathrm{CS}{\nabla}^\dagger_{*,\Gamma,-},\underset{\mathrm{Spec}}{\mathcal{O}}^\mathrm{CS}{\Phi}^r_{*,\Gamma,-},\underset{\mathrm{Spec}}{\mathcal{O}}^\mathrm{CS}{\Phi}^I_{*,\Gamma,-}.	
\end{align}

Then we take the corresponding quotients by using the corresponding Frobenius operators:
\begin{align}
&\underset{\mathrm{Spec}}{\mathcal{O}}^\mathrm{CS}\widetilde{\Delta}_{*,\Gamma,-}/\mathrm{Fro}^\mathbb{Z},\underset{\mathrm{Spec}}{\mathcal{O}}^\mathrm{CS}\widetilde{\nabla}_{*,\Gamma,-}/\mathrm{Fro}^\mathbb{Z},\underset{\mathrm{Spec}}{\mathcal{O}}^\mathrm{CS}\widetilde{\Phi}_{*,\Gamma,-}/\mathrm{Fro}^\mathbb{Z},\underset{\mathrm{Spec}}{\mathcal{O}}^\mathrm{CS}\widetilde{\Delta}^+_{*,\Gamma,-}/\mathrm{Fro}^\mathbb{Z},\\
&\underset{\mathrm{Spec}}{\mathcal{O}}^\mathrm{CS}\widetilde{\nabla}^+_{*,\Gamma,-}/\mathrm{Fro}^\mathbb{Z}, \underset{\mathrm{Spec}}{\mathcal{O}}^\mathrm{CS}\widetilde{\Delta}^\dagger_{*,\Gamma,-}/\mathrm{Fro}^\mathbb{Z},\underset{\mathrm{Spec}}{\mathcal{O}}^\mathrm{CS}\widetilde{\nabla}^\dagger_{*,\Gamma,-}/\mathrm{Fro}^\mathbb{Z},	\\
\end{align}
\begin{align}
&\underset{\mathrm{Spec}}{\mathcal{O}}^\mathrm{CS}\breve{\Delta}_{*,\Gamma,-}/\mathrm{Fro}^\mathbb{Z},\breve{\nabla}_{*,\Gamma,-}/\mathrm{Fro}^\mathbb{Z},\underset{\mathrm{Spec}}{\mathcal{O}}^\mathrm{CS}\breve{\Phi}_{*,\Gamma,-}/\mathrm{Fro}^\mathbb{Z},\underset{\mathrm{Spec}}{\mathcal{O}}^\mathrm{CS}\breve{\Delta}^+_{*,\Gamma,-}/\mathrm{Fro}^\mathbb{Z},\\
&\underset{\mathrm{Spec}}{\mathcal{O}}^\mathrm{CS}\breve{\nabla}^+_{*,\Gamma,-}/\mathrm{Fro}^\mathbb{Z}, \underset{\mathrm{Spec}}{\mathcal{O}}^\mathrm{CS}\breve{\Delta}^\dagger_{*,\Gamma,-}/\mathrm{Fro}^\mathbb{Z},\underset{\mathrm{Spec}}{\mathcal{O}}^\mathrm{CS}\breve{\nabla}^\dagger_{*,\Gamma,-}/\mathrm{Fro}^\mathbb{Z},	\\
\end{align}
\begin{align}
&\underset{\mathrm{Spec}}{\mathcal{O}}^\mathrm{CS}{\Delta}_{*,\Gamma,-}/\mathrm{Fro}^\mathbb{Z},\underset{\mathrm{Spec}}{\mathcal{O}}^\mathrm{CS}{\nabla}_{*,\Gamma,-}/\mathrm{Fro}^\mathbb{Z},\underset{\mathrm{Spec}}{\mathcal{O}}^\mathrm{CS}{\Phi}_{*,\Gamma,-}/\mathrm{Fro}^\mathbb{Z},\underset{\mathrm{Spec}}{\mathcal{O}}^\mathrm{CS}{\Delta}^+_{*,\Gamma,-}/\mathrm{Fro}^\mathbb{Z},\\
&\underset{\mathrm{Spec}}{\mathcal{O}}^\mathrm{CS}{\nabla}^+_{*,\Gamma,-}/\mathrm{Fro}^\mathbb{Z}, \underset{\mathrm{Spec}}{\mathcal{O}}^\mathrm{CS}{\Delta}^\dagger_{*,\Gamma,-}/\mathrm{Fro}^\mathbb{Z},\underset{\mathrm{Spec}}{\mathcal{O}}^\mathrm{CS}{\nabla}^\dagger_{*,\Gamma,-}/\mathrm{Fro}^\mathbb{Z}.	
\end{align}
Here for those space with notations related to the radius and the corresponding interval we consider the total unions $\bigcap_r,\bigcup_I$ in order to achieve the whole spaces to achieve the analogues of the corresponding FF curves from \cite{10KL1}, \cite{10KL2}, \cite{10FF} for
\[
\xymatrix@R+0pc@C+0pc{
\underset{r}{\mathrm{homotopycolimit}}~\underset{\mathrm{Spec}}{\mathcal{O}}^\mathrm{CS}\widetilde{\Phi}^r_{*,\Gamma,-},\underset{I}{\mathrm{homotopylimit}}~\underset{\mathrm{Spec}}{\mathcal{O}}^\mathrm{CS}\widetilde{\Phi}^I_{*,\Gamma,-},	\\
}
\]
\[
\xymatrix@R+0pc@C+0pc{
\underset{r}{\mathrm{homotopycolimit}}~\underset{\mathrm{Spec}}{\mathcal{O}}^\mathrm{CS}\breve{\Phi}^r_{*,\Gamma,-},\underset{I}{\mathrm{homotopylimit}}~\underset{\mathrm{Spec}}{\mathcal{O}}^\mathrm{CS}\breve{\Phi}^I_{*,\Gamma,-},	\\
}
\]
\[
\xymatrix@R+0pc@C+0pc{
\underset{r}{\mathrm{homotopycolimit}}~\underset{\mathrm{Spec}}{\mathcal{O}}^\mathrm{CS}{\Phi}^r_{*,\Gamma,-},\underset{I}{\mathrm{homotopylimit}}~\underset{\mathrm{Spec}}{\mathcal{O}}^\mathrm{CS}{\Phi}^I_{*,\Gamma,-}.	
}
\]
\[ 
\xymatrix@R+0pc@C+0pc{
\underset{r}{\mathrm{homotopycolimit}}~\underset{\mathrm{Spec}}{\mathcal{O}}^\mathrm{CS}\widetilde{\Phi}^r_{*,\Gamma,-}/\mathrm{Fro}^\mathbb{Z},\underset{I}{\mathrm{homotopylimit}}~\underset{\mathrm{Spec}}{\mathcal{O}}^\mathrm{CS}\widetilde{\Phi}^I_{*,\Gamma,-}/\mathrm{Fro}^\mathbb{Z},	\\
}
\]
\[ 
\xymatrix@R+0pc@C+0pc{
\underset{r}{\mathrm{homotopycolimit}}~\underset{\mathrm{Spec}}{\mathcal{O}}^\mathrm{CS}\breve{\Phi}^r_{*,\Gamma,-}/\mathrm{Fro}^\mathbb{Z},\underset{I}{\mathrm{homotopylimit}}~\breve{\Phi}^I_{*,\Gamma,-}/\mathrm{Fro}^\mathbb{Z},	\\
}
\]
\[ 
\xymatrix@R+0pc@C+0pc{
\underset{r}{\mathrm{homotopycolimit}}~\underset{\mathrm{Spec}}{\mathcal{O}}^\mathrm{CS}{\Phi}^r_{*,\Gamma,-}/\mathrm{Fro}^\mathbb{Z},\underset{I}{\mathrm{homotopylimit}}~\underset{\mathrm{Spec}}{\mathcal{O}}^\mathrm{CS}{\Phi}^I_{*,\Gamma,-}/\mathrm{Fro}^\mathbb{Z}.	
}
\]

\end{definition}

\

\begin{definition}
We then consider the corresponding quasipresheaves of the corresponding ind-Banach or monomorphic ind-Banach modules from \cite{10BBK}, \cite{10KKM}:
\begin{align}
\mathrm{Quasicoherentpresheaves,IndBanach}_{*}	
\end{align}
where $*$ is one of the following spaces:
\begin{align}
&\underset{\mathrm{Spec}}{\mathcal{O}}^\mathrm{BK}\widetilde{\Phi}_{*,\Gamma,-}/\mathrm{Fro}^\mathbb{Z},	\\
\end{align}
\begin{align}
&\underset{\mathrm{Spec}}{\mathcal{O}}^\mathrm{BK}\breve{\Phi}_{*,\Gamma,-}/\mathrm{Fro}^\mathbb{Z},	\\
\end{align}
\begin{align}
&\underset{\mathrm{Spec}}{\mathcal{O}}^\mathrm{BK}{\Phi}_{*,\Gamma,-}/\mathrm{Fro}^\mathbb{Z}.	
\end{align}
Here for those space without notation related to the radius and the corresponding interval we consider the total unions $\bigcap_r,\bigcup_I$ in order to achieve the whole spaces to achieve the analogues of the corresponding FF curves from \cite{10KL1}, \cite{10KL2}, \cite{10FF} for
\[
\xymatrix@R+0pc@C+0pc{
\underset{r}{\mathrm{homotopycolimit}}~\underset{\mathrm{Spec}}{\mathcal{O}}^\mathrm{BK}\widetilde{\Phi}^r_{*,\Gamma,-},\underset{I}{\mathrm{homotopylimit}}~\underset{\mathrm{Spec}}{\mathcal{O}}^\mathrm{BK}\widetilde{\Phi}^I_{*,\Gamma,-},	\\
}
\]
\[
\xymatrix@R+0pc@C+0pc{
\underset{r}{\mathrm{homotopycolimit}}~\underset{\mathrm{Spec}}{\mathcal{O}}^\mathrm{BK}\breve{\Phi}^r_{*,\Gamma,-},\underset{I}{\mathrm{homotopylimit}}~\underset{\mathrm{Spec}}{\mathcal{O}}^\mathrm{BK}\breve{\Phi}^I_{*,\Gamma,-},	\\
}
\]
\[
\xymatrix@R+0pc@C+0pc{
\underset{r}{\mathrm{homotopycolimit}}~\underset{\mathrm{Spec}}{\mathcal{O}}^\mathrm{BK}{\Phi}^r_{*,\Gamma,-},\underset{I}{\mathrm{homotopylimit}}~\underset{\mathrm{Spec}}{\mathcal{O}}^\mathrm{BK}{\Phi}^I_{*,\Gamma,-}.	
}
\]
\[  
\xymatrix@R+0pc@C+0pc{
\underset{r}{\mathrm{homotopycolimit}}~\underset{\mathrm{Spec}}{\mathcal{O}}^\mathrm{BK}\widetilde{\Phi}^r_{*,\Gamma,-}/\mathrm{Fro}^\mathbb{Z},\underset{I}{\mathrm{homotopylimit}}~\underset{\mathrm{Spec}}{\mathcal{O}}^\mathrm{BK}\widetilde{\Phi}^I_{*,\Gamma,-}/\mathrm{Fro}^\mathbb{Z},	\\
}
\]
\[ 
\xymatrix@R+0pc@C+0pc{
\underset{r}{\mathrm{homotopycolimit}}~\underset{\mathrm{Spec}}{\mathcal{O}}^\mathrm{BK}\breve{\Phi}^r_{*,\Gamma,-}/\mathrm{Fro}^\mathbb{Z},\underset{I}{\mathrm{homotopylimit}}~\underset{\mathrm{Spec}}{\mathcal{O}}^\mathrm{BK}\breve{\Phi}^I_{*,\Gamma,-}/\mathrm{Fro}^\mathbb{Z},	\\
}
\]
\[ 
\xymatrix@R+0pc@C+0pc{
\underset{r}{\mathrm{homotopycolimit}}~\underset{\mathrm{Spec}}{\mathcal{O}}^\mathrm{BK}{\Phi}^r_{*,\Gamma,-}/\mathrm{Fro}^\mathbb{Z},\underset{I}{\mathrm{homotopylimit}}~\underset{\mathrm{Spec}}{\mathcal{O}}^\mathrm{BK}{\Phi}^I_{*,\Gamma,-}/\mathrm{Fro}^\mathbb{Z}.	
}
\]

\end{definition}

\begin{definition}
We then consider the corresponding quasisheaves of the corresponding condensed solid topological modules from \cite{10CS2}:
\begin{align}
\mathrm{Quasicoherentsheaves, Condensed}_{*}	
\end{align}
where $*$ is one of the following spaces:
\begin{align}
&\underset{\mathrm{Spec}}{\mathcal{O}}^\mathrm{CS}\widetilde{\Delta}_{*,\Gamma,-}/\mathrm{Fro}^\mathbb{Z},\underset{\mathrm{Spec}}{\mathcal{O}}^\mathrm{CS}\widetilde{\nabla}_{*,\Gamma,-}/\mathrm{Fro}^\mathbb{Z},\underset{\mathrm{Spec}}{\mathcal{O}}^\mathrm{CS}\widetilde{\Phi}_{*,\Gamma,-}/\mathrm{Fro}^\mathbb{Z},\underset{\mathrm{Spec}}{\mathcal{O}}^\mathrm{CS}\widetilde{\Delta}^+_{*,\Gamma,-}/\mathrm{Fro}^\mathbb{Z},\\
&\underset{\mathrm{Spec}}{\mathcal{O}}^\mathrm{CS}\widetilde{\nabla}^+_{*,\Gamma,-}/\mathrm{Fro}^\mathbb{Z},\underset{\mathrm{Spec}}{\mathcal{O}}^\mathrm{CS}\widetilde{\Delta}^\dagger_{*,\Gamma,-}/\mathrm{Fro}^\mathbb{Z},\underset{\mathrm{Spec}}{\mathcal{O}}^\mathrm{CS}\widetilde{\nabla}^\dagger_{*,\Gamma,-}/\mathrm{Fro}^\mathbb{Z},	\\
\end{align}
\begin{align}
&\underset{\mathrm{Spec}}{\mathcal{O}}^\mathrm{CS}\breve{\Delta}_{*,\Gamma,-}/\mathrm{Fro}^\mathbb{Z},\breve{\nabla}_{*,\Gamma,-}/\mathrm{Fro}^\mathbb{Z},\underset{\mathrm{Spec}}{\mathcal{O}}^\mathrm{CS}\breve{\Phi}_{*,\Gamma,-}/\mathrm{Fro}^\mathbb{Z},\underset{\mathrm{Spec}}{\mathcal{O}}^\mathrm{CS}\breve{\Delta}^+_{*,\Gamma,-}/\mathrm{Fro}^\mathbb{Z},\\
&\underset{\mathrm{Spec}}{\mathcal{O}}^\mathrm{CS}\breve{\nabla}^+_{*,\Gamma,-}/\mathrm{Fro}^\mathbb{Z},\underset{\mathrm{Spec}}{\mathcal{O}}^\mathrm{CS}\breve{\Delta}^\dagger_{*,\Gamma,-}/\mathrm{Fro}^\mathbb{Z},\underset{\mathrm{Spec}}{\mathcal{O}}^\mathrm{CS}\breve{\nabla}^\dagger_{*,\Gamma,-}/\mathrm{Fro}^\mathbb{Z},	\\
\end{align}
\begin{align}
&\underset{\mathrm{Spec}}{\mathcal{O}}^\mathrm{CS}{\Delta}_{*,\Gamma,-}/\mathrm{Fro}^\mathbb{Z},\underset{\mathrm{Spec}}{\mathcal{O}}^\mathrm{CS}{\nabla}_{*,\Gamma,-}/\mathrm{Fro}^\mathbb{Z},\underset{\mathrm{Spec}}{\mathcal{O}}^\mathrm{CS}{\Phi}_{*,\Gamma,-}/\mathrm{Fro}^\mathbb{Z},\underset{\mathrm{Spec}}{\mathcal{O}}^\mathrm{CS}{\Delta}^+_{*,\Gamma,-}/\mathrm{Fro}^\mathbb{Z},\\
&\underset{\mathrm{Spec}}{\mathcal{O}}^\mathrm{CS}{\nabla}^+_{*,\Gamma,-}/\mathrm{Fro}^\mathbb{Z}, \underset{\mathrm{Spec}}{\mathcal{O}}^\mathrm{CS}{\Delta}^\dagger_{*,\Gamma,-}/\mathrm{Fro}^\mathbb{Z},\underset{\mathrm{Spec}}{\mathcal{O}}^\mathrm{CS}{\nabla}^\dagger_{*,\Gamma,-}/\mathrm{Fro}^\mathbb{Z}.	
\end{align}
Here for those space with notations related to the radius and the corresponding interval we consider the total unions $\bigcap_r,\bigcup_I$ in order to achieve the whole spaces to achieve the analogues of the corresponding FF curves from \cite{10KL1}, \cite{10KL2}, \cite{10FF} for
\[
\xymatrix@R+0pc@C+0pc{
\underset{r}{\mathrm{homotopycolimit}}~\underset{\mathrm{Spec}}{\mathcal{O}}^\mathrm{CS}\widetilde{\Phi}^r_{*,\Gamma,-},\underset{I}{\mathrm{homotopylimit}}~\underset{\mathrm{Spec}}{\mathcal{O}}^\mathrm{CS}\widetilde{\Phi}^I_{*,\Gamma,-},	\\
}
\]
\[
\xymatrix@R+0pc@C+0pc{
\underset{r}{\mathrm{homotopycolimit}}~\underset{\mathrm{Spec}}{\mathcal{O}}^\mathrm{CS}\breve{\Phi}^r_{*,\Gamma,-},\underset{I}{\mathrm{homotopylimit}}~\underset{\mathrm{Spec}}{\mathcal{O}}^\mathrm{CS}\breve{\Phi}^I_{*,\Gamma,-},	\\
}
\]
\[
\xymatrix@R+0pc@C+0pc{
\underset{r}{\mathrm{homotopycolimit}}~\underset{\mathrm{Spec}}{\mathcal{O}}^\mathrm{CS}{\Phi}^r_{*,\Gamma,-},\underset{I}{\mathrm{homotopylimit}}~\underset{\mathrm{Spec}}{\mathcal{O}}^\mathrm{CS}{\Phi}^I_{*,\Gamma,-}.	
}
\]
\[ 
\xymatrix@R+0pc@C+0pc{
\underset{r}{\mathrm{homotopycolimit}}~\underset{\mathrm{Spec}}{\mathcal{O}}^\mathrm{CS}\widetilde{\Phi}^r_{*,\Gamma,-}/\mathrm{Fro}^\mathbb{Z},\underset{I}{\mathrm{homotopylimit}}~\underset{\mathrm{Spec}}{\mathcal{O}}^\mathrm{CS}\widetilde{\Phi}^I_{*,\Gamma,-}/\mathrm{Fro}^\mathbb{Z},	\\
}
\]
\[ 
\xymatrix@R+0pc@C+0pc{
\underset{r}{\mathrm{homotopycolimit}}~\underset{\mathrm{Spec}}{\mathcal{O}}^\mathrm{CS}\breve{\Phi}^r_{*,\Gamma,-}/\mathrm{Fro}^\mathbb{Z},\underset{I}{\mathrm{homotopylimit}}~\breve{\Phi}^I_{*,\Gamma,-}/\mathrm{Fro}^\mathbb{Z},	\\
}
\]
\[ 
\xymatrix@R+0pc@C+0pc{
\underset{r}{\mathrm{homotopycolimit}}~\underset{\mathrm{Spec}}{\mathcal{O}}^\mathrm{CS}{\Phi}^r_{*,\Gamma,-}/\mathrm{Fro}^\mathbb{Z},\underset{I}{\mathrm{homotopylimit}}~\underset{\mathrm{Spec}}{\mathcal{O}}^\mathrm{CS}{\Phi}^I_{*,\Gamma,-}/\mathrm{Fro}^\mathbb{Z}.	
}
\]

\end{definition}

\

\begin{proposition}
There is a well-defined functor from the $\infty$-category 
\begin{align}
\mathrm{Quasicoherentpresheaves,Condensed}_{*}	
\end{align}
where $*$ is one of the following spaces:
\begin{align}
&\underset{\mathrm{Spec}}{\mathcal{O}}^\mathrm{CS}\widetilde{\Phi}_{*,\Gamma,-}/\mathrm{Fro}^\mathbb{Z},	\\
\end{align}
\begin{align}
&\underset{\mathrm{Spec}}{\mathcal{O}}^\mathrm{CS}\breve{\Phi}_{*,\Gamma,-}/\mathrm{Fro}^\mathbb{Z},	\\
\end{align}
\begin{align}
&\underset{\mathrm{Spec}}{\mathcal{O}}^\mathrm{CS}{\Phi}_{*,\Gamma,-}/\mathrm{Fro}^\mathbb{Z},	
\end{align}
to the $\infty$-category of $\mathrm{Fro}$-equivariant quasicoherent presheaves over similar spaces above correspondingly without the $\mathrm{Fro}$-quotients, and to the $\infty$-category of $\mathrm{Fro}$-equivariant quasicoherent modules over global sections of the structure $\infty$-sheaves of the similar spaces above correspondingly without the $\mathrm{Fro}$-quotients. Here for those space without notation related to the radius and the corresponding interval we consider the total unions $\bigcap_r,\bigcup_I$ in order to achieve the whole spaces to achieve the analogues of the corresponding FF curves from \cite{10KL1}, \cite{10KL2}, \cite{10FF} for
\[
\xymatrix@R+0pc@C+0pc{
\underset{r}{\mathrm{homotopycolimit}}~\underset{\mathrm{Spec}}{\mathcal{O}}^\mathrm{CS}\widetilde{\Phi}^r_{*,\Gamma,-},\underset{I}{\mathrm{homotopylimit}}~\underset{\mathrm{Spec}}{\mathcal{O}}^\mathrm{CS}\widetilde{\Phi}^I_{*,\Gamma,-},	\\
}
\]
\[
\xymatrix@R+0pc@C+0pc{
\underset{r}{\mathrm{homotopycolimit}}~\underset{\mathrm{Spec}}{\mathcal{O}}^\mathrm{CS}\breve{\Phi}^r_{*,\Gamma,-},\underset{I}{\mathrm{homotopylimit}}~\underset{\mathrm{Spec}}{\mathcal{O}}^\mathrm{CS}\breve{\Phi}^I_{*,\Gamma,-},	\\
}
\]
\[
\xymatrix@R+0pc@C+0pc{
\underset{r}{\mathrm{homotopycolimit}}~\underset{\mathrm{Spec}}{\mathcal{O}}^\mathrm{CS}{\Phi}^r_{*,\Gamma,-},\underset{I}{\mathrm{homotopylimit}}~\underset{\mathrm{Spec}}{\mathcal{O}}^\mathrm{CS}{\Phi}^I_{*,\Gamma,-}.	
}
\]
\[ 
\xymatrix@R+0pc@C+0pc{
\underset{r}{\mathrm{homotopycolimit}}~\underset{\mathrm{Spec}}{\mathcal{O}}^\mathrm{CS}\widetilde{\Phi}^r_{*,\Gamma,-}/\mathrm{Fro}^\mathbb{Z},\underset{I}{\mathrm{homotopylimit}}~\underset{\mathrm{Spec}}{\mathcal{O}}^\mathrm{CS}\widetilde{\Phi}^I_{*,\Gamma,-}/\mathrm{Fro}^\mathbb{Z},	\\
}
\]
\[ 
\xymatrix@R+0pc@C+0pc{
\underset{r}{\mathrm{homotopycolimit}}~\underset{\mathrm{Spec}}{\mathcal{O}}^\mathrm{CS}\breve{\Phi}^r_{*,\Gamma,-}/\mathrm{Fro}^\mathbb{Z},\underset{I}{\mathrm{homotopylimit}}~\breve{\Phi}^I_{*,\Gamma,-}/\mathrm{Fro}^\mathbb{Z},	\\
}
\]
\[ 
\xymatrix@R+0pc@C+0pc{
\underset{r}{\mathrm{homotopycolimit}}~\underset{\mathrm{Spec}}{\mathcal{O}}^\mathrm{CS}{\Phi}^r_{*,\Gamma,-}/\mathrm{Fro}^\mathbb{Z},\underset{I}{\mathrm{homotopylimit}}~\underset{\mathrm{Spec}}{\mathcal{O}}^\mathrm{CS}{\Phi}^I_{*,\Gamma,-}/\mathrm{Fro}^\mathbb{Z}.	
}
\]	
In this situation we will have the target category being family parametrized by $r$ or $I$ in compatible glueing sense as in \cite[Definition 5.4.10]{10KL2}. In this situation for modules parametrized by the intervals we have the equivalence of $\infty$-categories by using \cite[Proposition 13.8]{10CS2}. Here the corresponding quasicoherent Frobenius modules are defined to be the corresponding homotopy colimits and limits of Frobenius modules:
\begin{align}
\underset{r}{\mathrm{homotopycolimit}}~M_r,\\
\underset{I}{\mathrm{homotopylimit}}~M_I,	
\end{align}
where each $M_r$ is a Frobenius-equivariant module over the period ring with respect to some radius $r$ while each $M_I$ is a Frobenius-equivariant module over the period ring with respect to some interval $I$.\\
\end{proposition}

\begin{proposition}
Similar proposition holds for 
\begin{align}
\mathrm{Quasicoherentsheaves,IndBanach}_{*}.	
\end{align}	
\end{proposition}

\

\begin{definition}
We then consider the corresponding quasipresheaves of perfect complexes the corresponding ind-Banach or monomorphic ind-Banach modules from \cite{10BBK}, \cite{10KKM}:
\begin{align}
\mathrm{Quasicoherentpresheaves,Perfectcomplex,IndBanach}_{*}	
\end{align}
where $*$ is one of the following spaces:
\begin{align}
&\underset{\mathrm{Spec}}{\mathcal{O}}^\mathrm{BK}\widetilde{\Phi}_{*,\Gamma,-}/\mathrm{Fro}^\mathbb{Z},	\\
\end{align}
\begin{align}
&\underset{\mathrm{Spec}}{\mathcal{O}}^\mathrm{BK}\breve{\Phi}_{*,\Gamma,-}/\mathrm{Fro}^\mathbb{Z},	\\
\end{align}
\begin{align}
&\underset{\mathrm{Spec}}{\mathcal{O}}^\mathrm{BK}{\Phi}_{*,\Gamma,-}/\mathrm{Fro}^\mathbb{Z}.	
\end{align}
Here for those space without notation related to the radius and the corresponding interval we consider the total unions $\bigcap_r,\bigcup_I$ in order to achieve the whole spaces to achieve the analogues of the corresponding FF curves from \cite{10KL1}, \cite{10KL2}, \cite{10FF} for
\[
\xymatrix@R+0pc@C+0pc{
\underset{r}{\mathrm{homotopycolimit}}~\underset{\mathrm{Spec}}{\mathcal{O}}^\mathrm{BK}\widetilde{\Phi}^r_{*,\Gamma,-},\underset{I}{\mathrm{homotopylimit}}~\underset{\mathrm{Spec}}{\mathcal{O}}^\mathrm{BK}\widetilde{\Phi}^I_{*,\Gamma,-},	\\
}
\]
\[
\xymatrix@R+0pc@C+0pc{
\underset{r}{\mathrm{homotopycolimit}}~\underset{\mathrm{Spec}}{\mathcal{O}}^\mathrm{BK}\breve{\Phi}^r_{*,\Gamma,-},\underset{I}{\mathrm{homotopylimit}}~\underset{\mathrm{Spec}}{\mathcal{O}}^\mathrm{BK}\breve{\Phi}^I_{*,\Gamma,-},	\\
}
\]
\[
\xymatrix@R+0pc@C+0pc{
\underset{r}{\mathrm{homotopycolimit}}~\underset{\mathrm{Spec}}{\mathcal{O}}^\mathrm{BK}{\Phi}^r_{*,\Gamma,-},\underset{I}{\mathrm{homotopylimit}}~\underset{\mathrm{Spec}}{\mathcal{O}}^\mathrm{BK}{\Phi}^I_{*,\Gamma,-}.	
}
\]
\[  
\xymatrix@R+0pc@C+0pc{
\underset{r}{\mathrm{homotopycolimit}}~\underset{\mathrm{Spec}}{\mathcal{O}}^\mathrm{BK}\widetilde{\Phi}^r_{*,\Gamma,-}/\mathrm{Fro}^\mathbb{Z},\underset{I}{\mathrm{homotopylimit}}~\underset{\mathrm{Spec}}{\mathcal{O}}^\mathrm{BK}\widetilde{\Phi}^I_{*,\Gamma,-}/\mathrm{Fro}^\mathbb{Z},	\\
}
\]
\[ 
\xymatrix@R+0pc@C+0pc{
\underset{r}{\mathrm{homotopycolimit}}~\underset{\mathrm{Spec}}{\mathcal{O}}^\mathrm{BK}\breve{\Phi}^r_{*,\Gamma,-}/\mathrm{Fro}^\mathbb{Z},\underset{I}{\mathrm{homotopylimit}}~\underset{\mathrm{Spec}}{\mathcal{O}}^\mathrm{BK}\breve{\Phi}^I_{*,\Gamma,-}/\mathrm{Fro}^\mathbb{Z},	\\
}
\]
\[ 
\xymatrix@R+0pc@C+0pc{
\underset{r}{\mathrm{homotopycolimit}}~\underset{\mathrm{Spec}}{\mathcal{O}}^\mathrm{BK}{\Phi}^r_{*,\Gamma,-}/\mathrm{Fro}^\mathbb{Z},\underset{I}{\mathrm{homotopylimit}}~\underset{\mathrm{Spec}}{\mathcal{O}}^\mathrm{BK}{\Phi}^I_{*,\Gamma,-}/\mathrm{Fro}^\mathbb{Z}.	
}
\]

\end{definition}

\begin{definition}
We then consider the corresponding quasisheaves of perfect complexes of the corresponding condensed solid topological modules from \cite{10CS2}:
\begin{align}
\mathrm{Quasicoherentsheaves, Perfectcomplex, Condensed}_{*}	
\end{align}
where $*$ is one of the following spaces:
\begin{align}
&\underset{\mathrm{Spec}}{\mathcal{O}}^\mathrm{CS}\widetilde{\Delta}_{*,\Gamma,-}/\mathrm{Fro}^\mathbb{Z},\underset{\mathrm{Spec}}{\mathcal{O}}^\mathrm{CS}\widetilde{\nabla}_{*,\Gamma,-}/\mathrm{Fro}^\mathbb{Z},\underset{\mathrm{Spec}}{\mathcal{O}}^\mathrm{CS}\widetilde{\Phi}_{*,\Gamma,-}/\mathrm{Fro}^\mathbb{Z},\underset{\mathrm{Spec}}{\mathcal{O}}^\mathrm{CS}\widetilde{\Delta}^+_{*,\Gamma,-}/\mathrm{Fro}^\mathbb{Z},\\
&\underset{\mathrm{Spec}}{\mathcal{O}}^\mathrm{CS}\widetilde{\nabla}^+_{*,\Gamma,-}/\mathrm{Fro}^\mathbb{Z},\underset{\mathrm{Spec}}{\mathcal{O}}^\mathrm{CS}\widetilde{\Delta}^\dagger_{*,\Gamma,-}/\mathrm{Fro}^\mathbb{Z},\underset{\mathrm{Spec}}{\mathcal{O}}^\mathrm{CS}\widetilde{\nabla}^\dagger_{*,\Gamma,-}/\mathrm{Fro}^\mathbb{Z},	\\
\end{align}
\begin{align}
&\underset{\mathrm{Spec}}{\mathcal{O}}^\mathrm{CS}\breve{\Delta}_{*,\Gamma,-}/\mathrm{Fro}^\mathbb{Z},\breve{\nabla}_{*,\Gamma,-}/\mathrm{Fro}^\mathbb{Z},\underset{\mathrm{Spec}}{\mathcal{O}}^\mathrm{CS}\breve{\Phi}_{*,\Gamma,-}/\mathrm{Fro}^\mathbb{Z},\underset{\mathrm{Spec}}{\mathcal{O}}^\mathrm{CS}\breve{\Delta}^+_{*,\Gamma,-}/\mathrm{Fro}^\mathbb{Z},\\
&\underset{\mathrm{Spec}}{\mathcal{O}}^\mathrm{CS}\breve{\nabla}^+_{*,\Gamma,-}/\mathrm{Fro}^\mathbb{Z},\underset{\mathrm{Spec}}{\mathcal{O}}^\mathrm{CS}\breve{\Delta}^\dagger_{*,\Gamma,-}/\mathrm{Fro}^\mathbb{Z},\underset{\mathrm{Spec}}{\mathcal{O}}^\mathrm{CS}\breve{\nabla}^\dagger_{*,\Gamma,-}/\mathrm{Fro}^\mathbb{Z},	\\
\end{align}
\begin{align}
&\underset{\mathrm{Spec}}{\mathcal{O}}^\mathrm{CS}{\Delta}_{*,\Gamma,-}/\mathrm{Fro}^\mathbb{Z},\underset{\mathrm{Spec}}{\mathcal{O}}^\mathrm{CS}{\nabla}_{*,\Gamma,-}/\mathrm{Fro}^\mathbb{Z},\underset{\mathrm{Spec}}{\mathcal{O}}^\mathrm{CS}{\Phi}_{*,\Gamma,-}/\mathrm{Fro}^\mathbb{Z},\underset{\mathrm{Spec}}{\mathcal{O}}^\mathrm{CS}{\Delta}^+_{*,\Gamma,-}/\mathrm{Fro}^\mathbb{Z},\\
&\underset{\mathrm{Spec}}{\mathcal{O}}^\mathrm{CS}{\nabla}^+_{*,\Gamma,-}/\mathrm{Fro}^\mathbb{Z}, \underset{\mathrm{Spec}}{\mathcal{O}}^\mathrm{CS}{\Delta}^\dagger_{*,\Gamma,-}/\mathrm{Fro}^\mathbb{Z},\underset{\mathrm{Spec}}{\mathcal{O}}^\mathrm{CS}{\nabla}^\dagger_{*,\Gamma,-}/\mathrm{Fro}^\mathbb{Z}.	
\end{align}
Here for those space with notations related to the radius and the corresponding interval we consider the total unions $\bigcap_r,\bigcup_I$ in order to achieve the whole spaces to achieve the analogues of the corresponding FF curves from \cite{10KL1}, \cite{10KL2}, \cite{10FF} for
\[
\xymatrix@R+0pc@C+0pc{
\underset{r}{\mathrm{homotopycolimit}}~\underset{\mathrm{Spec}}{\mathcal{O}}^\mathrm{CS}\widetilde{\Phi}^r_{*,\Gamma,-},\underset{I}{\mathrm{homotopylimit}}~\underset{\mathrm{Spec}}{\mathcal{O}}^\mathrm{CS}\widetilde{\Phi}^I_{*,\Gamma,-},	\\
}
\]
\[
\xymatrix@R+0pc@C+0pc{
\underset{r}{\mathrm{homotopycolimit}}~\underset{\mathrm{Spec}}{\mathcal{O}}^\mathrm{CS}\breve{\Phi}^r_{*,\Gamma,-},\underset{I}{\mathrm{homotopylimit}}~\underset{\mathrm{Spec}}{\mathcal{O}}^\mathrm{CS}\breve{\Phi}^I_{*,\Gamma,-},	\\
}
\]
\[
\xymatrix@R+0pc@C+0pc{
\underset{r}{\mathrm{homotopycolimit}}~\underset{\mathrm{Spec}}{\mathcal{O}}^\mathrm{CS}{\Phi}^r_{*,\Gamma,-},\underset{I}{\mathrm{homotopylimit}}~\underset{\mathrm{Spec}}{\mathcal{O}}^\mathrm{CS}{\Phi}^I_{*,\Gamma,-}.	
}
\]
\[ 
\xymatrix@R+0pc@C+0pc{
\underset{r}{\mathrm{homotopycolimit}}~\underset{\mathrm{Spec}}{\mathcal{O}}^\mathrm{CS}\widetilde{\Phi}^r_{*,\Gamma,-}/\mathrm{Fro}^\mathbb{Z},\underset{I}{\mathrm{homotopylimit}}~\underset{\mathrm{Spec}}{\mathcal{O}}^\mathrm{CS}\widetilde{\Phi}^I_{*,\Gamma,-}/\mathrm{Fro}^\mathbb{Z},	\\
}
\]
\[ 
\xymatrix@R+0pc@C+0pc{
\underset{r}{\mathrm{homotopycolimit}}~\underset{\mathrm{Spec}}{\mathcal{O}}^\mathrm{CS}\breve{\Phi}^r_{*,\Gamma,-}/\mathrm{Fro}^\mathbb{Z},\underset{I}{\mathrm{homotopylimit}}~\breve{\Phi}^I_{*,\Gamma,-}/\mathrm{Fro}^\mathbb{Z},	\\
}
\]
\[ 
\xymatrix@R+0pc@C+0pc{
\underset{r}{\mathrm{homotopycolimit}}~\underset{\mathrm{Spec}}{\mathcal{O}}^\mathrm{CS}{\Phi}^r_{*,\Gamma,-}/\mathrm{Fro}^\mathbb{Z},\underset{I}{\mathrm{homotopylimit}}~\underset{\mathrm{Spec}}{\mathcal{O}}^\mathrm{CS}{\Phi}^I_{*,\Gamma,-}/\mathrm{Fro}^\mathbb{Z}.	
}
\]

\end{definition}

\begin{proposition}
There is a well-defined functor from the $\infty$-category 
\begin{align}
\mathrm{Quasicoherentpresheaves,Perfectcomplex,Condensed}_{*}	
\end{align}
where $*$ is one of the following spaces:
\begin{align}
&\underset{\mathrm{Spec}}{\mathcal{O}}^\mathrm{CS}\widetilde{\Phi}_{*,\Gamma,-}/\mathrm{Fro}^\mathbb{Z},	\\
\end{align}
\begin{align}
&\underset{\mathrm{Spec}}{\mathcal{O}}^\mathrm{CS}\breve{\Phi}_{*,\Gamma,-}/\mathrm{Fro}^\mathbb{Z},	\\
\end{align}
\begin{align}
&\underset{\mathrm{Spec}}{\mathcal{O}}^\mathrm{CS}{\Phi}_{*,\Gamma,-}/\mathrm{Fro}^\mathbb{Z},	
\end{align}
to the $\infty$-category of $\mathrm{Fro}$-equivariant quasicoherent presheaves over similar spaces above correspondingly without the $\mathrm{Fro}$-quotients, and to the $\infty$-category of $\mathrm{Fro}$-equivariant quasicoherent modules over global sections of the structure $\infty$-sheaves of the similar spaces above correspondingly without the $\mathrm{Fro}$-quotients. Here for those space without notation related to the radius and the corresponding interval we consider the total unions $\bigcap_r,\bigcup_I$ in order to achieve the whole spaces to achieve the analogues of the corresponding FF curves from \cite{10KL1}, \cite{10KL2}, \cite{10FF} for
\[
\xymatrix@R+0pc@C+0pc{
\underset{r}{\mathrm{homotopycolimit}}~\underset{\mathrm{Spec}}{\mathcal{O}}^\mathrm{CS}\widetilde{\Phi}^r_{*,\Gamma,-},\underset{I}{\mathrm{homotopylimit}}~\underset{\mathrm{Spec}}{\mathcal{O}}^\mathrm{CS}\widetilde{\Phi}^I_{*,\Gamma,-},	\\
}
\]
\[
\xymatrix@R+0pc@C+0pc{
\underset{r}{\mathrm{homotopycolimit}}~\underset{\mathrm{Spec}}{\mathcal{O}}^\mathrm{CS}\breve{\Phi}^r_{*,\Gamma,-},\underset{I}{\mathrm{homotopylimit}}~\underset{\mathrm{Spec}}{\mathcal{O}}^\mathrm{CS}\breve{\Phi}^I_{*,\Gamma,-},	\\
}
\]
\[
\xymatrix@R+0pc@C+0pc{
\underset{r}{\mathrm{homotopycolimit}}~\underset{\mathrm{Spec}}{\mathcal{O}}^\mathrm{CS}{\Phi}^r_{*,\Gamma,-},\underset{I}{\mathrm{homotopylimit}}~\underset{\mathrm{Spec}}{\mathcal{O}}^\mathrm{CS}{\Phi}^I_{*,\Gamma,-}.	
}
\]
\[ 
\xymatrix@R+0pc@C+0pc{
\underset{r}{\mathrm{homotopycolimit}}~\underset{\mathrm{Spec}}{\mathcal{O}}^\mathrm{CS}\widetilde{\Phi}^r_{*,\Gamma,-}/\mathrm{Fro}^\mathbb{Z},\underset{I}{\mathrm{homotopylimit}}~\underset{\mathrm{Spec}}{\mathcal{O}}^\mathrm{CS}\widetilde{\Phi}^I_{*,\Gamma,-}/\mathrm{Fro}^\mathbb{Z},	\\
}
\]
\[ 
\xymatrix@R+0pc@C+0pc{
\underset{r}{\mathrm{homotopycolimit}}~\underset{\mathrm{Spec}}{\mathcal{O}}^\mathrm{CS}\breve{\Phi}^r_{*,\Gamma,-}/\mathrm{Fro}^\mathbb{Z},\underset{I}{\mathrm{homotopylimit}}~\breve{\Phi}^I_{*,\Gamma,-}/\mathrm{Fro}^\mathbb{Z},	\\
}
\]
\[ 
\xymatrix@R+0pc@C+0pc{
\underset{r}{\mathrm{homotopycolimit}}~\underset{\mathrm{Spec}}{\mathcal{O}}^\mathrm{CS}{\Phi}^r_{*,\Gamma,-}/\mathrm{Fro}^\mathbb{Z},\underset{I}{\mathrm{homotopylimit}}~\underset{\mathrm{Spec}}{\mathcal{O}}^\mathrm{CS}{\Phi}^I_{*,\Gamma,-}/\mathrm{Fro}^\mathbb{Z}.	
}
\]	
In this situation we will have the target category being family parametrized by $r$ or $I$ in compatible glueing sense as in \cite[Definition 5.4.10]{10KL2}. In this situation for modules parametrized by the intervals we have the equivalence of $\infty$-categories by using \cite[Proposition 12.18]{10CS2}. Here the corresponding quasicoherent Frobenius modules are defined to be the corresponding homotopy colimits and limits of Frobenius modules:
\begin{align}
\underset{r}{\mathrm{homotopycolimit}}~M_r,\\
\underset{I}{\mathrm{homotopylimit}}~M_I,	
\end{align}
where each $M_r$ is a Frobenius-equivariant module over the period ring with respect to some radius $r$ while each $M_I$ is a Frobenius-equivariant module over the period ring with respect to some interval $I$.\\
\end{proposition}

\begin{proposition}
Similar proposition holds for 
\begin{align}
\mathrm{Quasicoherentsheaves,Perfectcomplex,IndBanach}_{*}.	
\end{align}	
\end{proposition}

\newpage
\subsection{Frobenius Quasicoherent Prestacks III: Deformation in $(\infty,1)$-Ind-Banach Rings}

\begin{definition}
We now consider the pro-\'etale site of $\mathrm{Spa}\mathbb{Q}_p\left<X_1^{\pm 1},...,X_k^{\pm 1}\right>$, denote that by $*$. To be more accurate we replace one component for $\Gamma$ with the pro-\'etale site of $\mathrm{Spa}\mathbb{Q}_p\left<X_1^{\pm 1},...,X_k^{\pm 1}\right>$. And we treat then all the functor to be prestacks for this site. Then from \cite{10KL1} and \cite[Definition 5.2.1]{10KL2} we have the following class of Kedlaya-Liu rings (with the following replacement: $\Delta$ stands for $A$, $\nabla$ stands for $B$, while $\Phi$ stands for $C$) by taking product in the sense of self $\Gamma$-th power:

\[
\xymatrix@R+0pc@C+0pc{
\widetilde{\Delta}_{*,\Gamma},\widetilde{\nabla}_{*,\Gamma},\widetilde{\Phi}_{*,\Gamma},\widetilde{\Delta}^+_{*,\Gamma},\widetilde{\nabla}^+_{*,\Gamma},\widetilde{\Delta}^\dagger_{*,\Gamma},\widetilde{\nabla}^\dagger_{*,\Gamma},\widetilde{\Phi}^r_{*,\Gamma},\widetilde{\Phi}^I_{*,\Gamma}, 
}
\]

\[
\xymatrix@R+0pc@C+0pc{
\breve{\Delta}_{*,\Gamma},\breve{\nabla}_{*,\Gamma},\breve{\Phi}_{*,\Gamma},\breve{\Delta}^+_{*,\Gamma},\breve{\nabla}^+_{*,\Gamma},\breve{\Delta}^\dagger_{*,\Gamma},\breve{\nabla}^\dagger_{*,\Gamma},\breve{\Phi}^r_{*,\Gamma},\breve{\Phi}^I_{*,\Gamma},	
}
\]

\[
\xymatrix@R+0pc@C+0pc{
{\Delta}_{*,\Gamma},{\nabla}_{*,\Gamma},{\Phi}_{*,\Gamma},{\Delta}^+_{*,\Gamma},{\nabla}^+_{*,\Gamma},{\Delta}^\dagger_{*,\Gamma},{\nabla}^\dagger_{*,\Gamma},{\Phi}^r_{*,\Gamma},{\Phi}^I_{*,\Gamma}.	
}
\]
We now consider the following rings with $\square$ being a homotopy colimit
\begin{align}
 \underset{I}{\mathrm{homotopylimit}}\square_i
 \end{align}
 of $\mathbb{Q}_p\left<Y_1,...,Y_i\right>,i=1,2,...$ in $\infty$-categories of simplicial ind-Banach rings in \cite{10BBBK}
 \begin{align}
  \mathrm{SimplicialInd-BanachRings}_{\mathbb{Q}_p}
\end{align}  
or animated analytic condensed commutative algebras in \cite{10CS2} 
\begin{align}   
\mathrm{SimplicialAnalyticCondensed}_{\mathbb{Q}_p}.
\end{align}   
Taking the product we have:
\[
\xymatrix@R+0pc@C+0pc{
\widetilde{\Phi}_{*,\Gamma,\square},\widetilde{\Phi}^r_{*,\Gamma,\square},\widetilde{\Phi}^I_{*,\Gamma,\square},	
}
\]
\[
\xymatrix@R+0pc@C+0pc{
\breve{\Phi}_{*,\Gamma,\square},\breve{\Phi}^r_{*,\Gamma,\square},\breve{\Phi}^I_{*,\Gamma,\square},	
}
\]
\[
\xymatrix@R+0pc@C+0pc{
{\Phi}_{*,\Gamma,\square},{\Phi}^r_{*,\Gamma,\square},{\Phi}^I_{*,\Gamma,\square}.	
}
\]
They carry multi Frobenius action $\varphi_\Gamma$ and multi $\mathrm{Lie}_\Gamma:=\mathbb{Z}_p^{\times\Gamma}$ action. In our current situation after \cite{10CKZ} and \cite{10PZ} we consider the following $(\infty,1)$-categories of $(\infty,1)$-modules.\\
\end{definition}

\begin{definition}
First we consider the Bambozzi-Kremnizer spectrum $\underset{\mathrm{Spec}}{\mathcal{O}}^\mathrm{BK}(*)$ attached to any of those in the above from \cite{10BK} by taking derived rational localization:
\begin{align}
&\underset{\mathrm{Spec}}{\mathcal{O}}^\mathrm{BK}\widetilde{\Phi}_{*,\Gamma,\square},\underset{\mathrm{Spec}}{\mathcal{O}}^\mathrm{BK}\widetilde{\Phi}^r_{*,\Gamma,\square},\underset{\mathrm{Spec}}{\mathcal{O}}^\mathrm{BK}\widetilde{\Phi}^I_{*,\Gamma,\square},	
\end{align}
\begin{align}
&\underset{\mathrm{Spec}}{\mathcal{O}}^\mathrm{BK}\breve{\Phi}_{*,\Gamma,\square},\underset{\mathrm{Spec}}{\mathcal{O}}^\mathrm{BK}\breve{\Phi}^r_{*,\Gamma,\square},\underset{\mathrm{Spec}}{\mathcal{O}}^\mathrm{BK}\breve{\Phi}^I_{*,\Gamma,\square},	
\end{align}
\begin{align}
&\underset{\mathrm{Spec}}{\mathcal{O}}^\mathrm{BK}{\Phi}_{*,\Gamma,\square},
\underset{\mathrm{Spec}}{\mathcal{O}}^\mathrm{BK}{\Phi}^r_{*,\Gamma,\square},\underset{\mathrm{Spec}}{\mathcal{O}}^\mathrm{BK}{\Phi}^I_{*,\Gamma,\square}.	
\end{align}

Then we take the corresponding quotients by using the corresponding Frobenius operators:
\begin{align}
&\underset{\mathrm{Spec}}{\mathcal{O}}^\mathrm{BK}\widetilde{\Phi}_{*,\Gamma,\square}/\mathrm{Fro}^\mathbb{Z},	\\
\end{align}
\begin{align}
&\underset{\mathrm{Spec}}{\mathcal{O}}^\mathrm{BK}\breve{\Phi}_{*,\Gamma,\square}/\mathrm{Fro}^\mathbb{Z},	\\
\end{align}
\begin{align}
&\underset{\mathrm{Spec}}{\mathcal{O}}^\mathrm{BK}{\Phi}_{*,\Gamma,\square}/\mathrm{Fro}^\mathbb{Z}.	
\end{align}
Here for those space without notation related to the radius and the corresponding interval we consider the total unions $\bigcap_r,\bigcup_I$ in order to achieve the whole spaces to achieve the analogues of the corresponding FF curves from \cite{10KL1}, \cite{10KL2}, \cite{10FF} for
\[
\xymatrix@R+0pc@C+0pc{
\underset{r}{\mathrm{homotopycolimit}}~\underset{\mathrm{Spec}}{\mathcal{O}}^\mathrm{BK}\widetilde{\Phi}^r_{*,\Gamma,\square},\underset{I}{\mathrm{homotopylimit}}~\underset{\mathrm{Spec}}{\mathcal{O}}^\mathrm{BK}\widetilde{\Phi}^I_{*,\Gamma,\square},	\\
}
\]
\[
\xymatrix@R+0pc@C+0pc{
\underset{r}{\mathrm{homotopycolimit}}~\underset{\mathrm{Spec}}{\mathcal{O}}^\mathrm{BK}\breve{\Phi}^r_{*,\Gamma,\square},\underset{I}{\mathrm{homotopylimit}}~\underset{\mathrm{Spec}}{\mathcal{O}}^\mathrm{BK}\breve{\Phi}^I_{*,\Gamma,\square},	\\
}
\]
\[
\xymatrix@R+0pc@C+0pc{
\underset{r}{\mathrm{homotopycolimit}}~\underset{\mathrm{Spec}}{\mathcal{O}}^\mathrm{BK}{\Phi}^r_{*,\Gamma,\square},\underset{I}{\mathrm{homotopylimit}}~\underset{\mathrm{Spec}}{\mathcal{O}}^\mathrm{BK}{\Phi}^I_{*,\Gamma,\square}.	
}
\]
\[  
\xymatrix@R+0pc@C+0pc{
\underset{r}{\mathrm{homotopycolimit}}~\underset{\mathrm{Spec}}{\mathcal{O}}^\mathrm{BK}\widetilde{\Phi}^r_{*,\Gamma,\square}/\mathrm{Fro}^\mathbb{Z},\underset{I}{\mathrm{homotopylimit}}~\underset{\mathrm{Spec}}{\mathcal{O}}^\mathrm{BK}\widetilde{\Phi}^I_{*,\Gamma,\square}/\mathrm{Fro}^\mathbb{Z},	\\
}
\]
\[ 
\xymatrix@R+0pc@C+0pc{
\underset{r}{\mathrm{homotopycolimit}}~\underset{\mathrm{Spec}}{\mathcal{O}}^\mathrm{BK}\breve{\Phi}^r_{*,\Gamma,\square}/\mathrm{Fro}^\mathbb{Z},\underset{I}{\mathrm{homotopylimit}}~\underset{\mathrm{Spec}}{\mathcal{O}}^\mathrm{BK}\breve{\Phi}^I_{*,\Gamma,\square}/\mathrm{Fro}^\mathbb{Z},	\\
}
\]
\[ 
\xymatrix@R+0pc@C+0pc{
\underset{r}{\mathrm{homotopycolimit}}~\underset{\mathrm{Spec}}{\mathcal{O}}^\mathrm{BK}{\Phi}^r_{*,\Gamma,\square}/\mathrm{Fro}^\mathbb{Z},\underset{I}{\mathrm{homotopylimit}}~\underset{\mathrm{Spec}}{\mathcal{O}}^\mathrm{BK}{\Phi}^I_{*,\Gamma,\square}/\mathrm{Fro}^\mathbb{Z}.	
}
\]

\end{definition}

\indent Meanwhile we have the corresponding Clausen-Scholze analytic stacks from \cite{10CS2}, therefore applying their construction we have:

\begin{definition}
Here we define the following products by using the solidified tensor product from \cite{10CS1} and \cite{10CS2}. Namely $A$ will still as above as a Banach ring over $\mathbb{Q}_p$. Then we take solidified tensor product $\overset{\blacksquare}{\otimes}$ of any of the following
\[
\xymatrix@R+0pc@C+0pc{
\widetilde{\Delta}_{*,\Gamma},\widetilde{\nabla}_{*,\Gamma},\widetilde{\Phi}_{*,\Gamma},\widetilde{\Delta}^+_{*,\Gamma},\widetilde{\nabla}^+_{*,\Gamma},\widetilde{\Delta}^\dagger_{*,\Gamma},\widetilde{\nabla}^\dagger_{*,\Gamma},\widetilde{\Phi}^r_{*,\Gamma},\widetilde{\Phi}^I_{*,\Gamma}, 
}
\]

\[
\xymatrix@R+0pc@C+0pc{
\breve{\Delta}_{*,\Gamma},\breve{\nabla}_{*,\Gamma},\breve{\Phi}_{*,\Gamma},\breve{\Delta}^+_{*,\Gamma},\breve{\nabla}^+_{*,\Gamma},\breve{\Delta}^\dagger_{*,\Gamma},\breve{\nabla}^\dagger_{*,\Gamma},\breve{\Phi}^r_{*,\Gamma},\breve{\Phi}^I_{*,\Gamma},	
}
\]

\[
\xymatrix@R+0pc@C+0pc{
{\Delta}_{*,\Gamma},{\nabla}_{*,\Gamma},{\Phi}_{*,\Gamma},{\Delta}^+_{*,\Gamma},{\nabla}^+_{*,\Gamma},{\Delta}^\dagger_{*,\Gamma},{\nabla}^\dagger_{*,\Gamma},{\Phi}^r_{*,\Gamma},{\Phi}^I_{*,\Gamma},	
}
\]  	
with $A$. Then we have the notations:
\[
\xymatrix@R+0pc@C+0pc{
\widetilde{\Delta}_{*,\Gamma,\square},\widetilde{\nabla}_{*,\Gamma,\square},\widetilde{\Phi}_{*,\Gamma,\square},\widetilde{\Delta}^+_{*,\Gamma,\square},\widetilde{\nabla}^+_{*,\Gamma,\square},\widetilde{\Delta}^\dagger_{*,\Gamma,\square},\widetilde{\nabla}^\dagger_{*,\Gamma,\square},\widetilde{\Phi}^r_{*,\Gamma,\square},\widetilde{\Phi}^I_{*,\Gamma,\square}, 
}
\]

\[
\xymatrix@R+0pc@C+0pc{
\breve{\Delta}_{*,\Gamma,\square},\breve{\nabla}_{*,\Gamma,\square},\breve{\Phi}_{*,\Gamma,\square},\breve{\Delta}^+_{*,\Gamma,\square},\breve{\nabla}^+_{*,\Gamma,\square},\breve{\Delta}^\dagger_{*,\Gamma,\square},\breve{\nabla}^\dagger_{*,\Gamma,\square},\breve{\Phi}^r_{*,\Gamma,\square},\breve{\Phi}^I_{*,\Gamma,\square},	
}
\]

\[
\xymatrix@R+0pc@C+0pc{
{\Delta}_{*,\Gamma,\square},{\nabla}_{*,\Gamma,\square},{\Phi}_{*,\Gamma,\square},{\Delta}^+_{*,\Gamma,\square},{\nabla}^+_{*,\Gamma,\square},{\Delta}^\dagger_{*,\Gamma,\square},{\nabla}^\dagger_{*,\Gamma,\square},{\Phi}^r_{*,\Gamma,\square},{\Phi}^I_{*,\Gamma,\square}.	
}
\]
\end{definition}

\begin{definition}
First we consider the Clausen-Scholze spectrum $\underset{\mathrm{Spec}}{\mathcal{O}}^\mathrm{CS}(*)$ attached to any of those in the above from \cite{10CS2} by taking derived rational localization:
\begin{align}
\underset{\mathrm{Spec}}{\mathcal{O}}^\mathrm{CS}\widetilde{\Delta}_{*,\Gamma,\square},\underset{\mathrm{Spec}}{\mathcal{O}}^\mathrm{CS}\widetilde{\nabla}_{*,\Gamma,\square},\underset{\mathrm{Spec}}{\mathcal{O}}^\mathrm{CS}\widetilde{\Phi}_{*,\Gamma,\square},\underset{\mathrm{Spec}}{\mathcal{O}}^\mathrm{CS}\widetilde{\Delta}^+_{*,\Gamma,\square},\underset{\mathrm{Spec}}{\mathcal{O}}^\mathrm{CS}\widetilde{\nabla}^+_{*,\Gamma,\square},\\
\underset{\mathrm{Spec}}{\mathcal{O}}^\mathrm{CS}\widetilde{\Delta}^\dagger_{*,\Gamma,\square},\underset{\mathrm{Spec}}{\mathcal{O}}^\mathrm{CS}\widetilde{\nabla}^\dagger_{*,\Gamma,\square},\underset{\mathrm{Spec}}{\mathcal{O}}^\mathrm{CS}\widetilde{\Phi}^r_{*,\Gamma,\square},\underset{\mathrm{Spec}}{\mathcal{O}}^\mathrm{CS}\widetilde{\Phi}^I_{*,\Gamma,\square},	\\
\end{align}
\begin{align}
\underset{\mathrm{Spec}}{\mathcal{O}}^\mathrm{CS}\breve{\Delta}_{*,\Gamma,\square},\breve{\nabla}_{*,\Gamma,\square},\underset{\mathrm{Spec}}{\mathcal{O}}^\mathrm{CS}\breve{\Phi}_{*,\Gamma,\square},\underset{\mathrm{Spec}}{\mathcal{O}}^\mathrm{CS}\breve{\Delta}^+_{*,\Gamma,\square},\underset{\mathrm{Spec}}{\mathcal{O}}^\mathrm{CS}\breve{\nabla}^+_{*,\Gamma,\square},\\
\underset{\mathrm{Spec}}{\mathcal{O}}^\mathrm{CS}\breve{\Delta}^\dagger_{*,\Gamma,\square},\underset{\mathrm{Spec}}{\mathcal{O}}^\mathrm{CS}\breve{\nabla}^\dagger_{*,\Gamma,\square},\underset{\mathrm{Spec}}{\mathcal{O}}^\mathrm{CS}\breve{\Phi}^r_{*,\Gamma,\square},\breve{\Phi}^I_{*,\Gamma,\square},	\\
\end{align}
\begin{align}
\underset{\mathrm{Spec}}{\mathcal{O}}^\mathrm{CS}{\Delta}_{*,\Gamma,\square},\underset{\mathrm{Spec}}{\mathcal{O}}^\mathrm{CS}{\nabla}_{*,\Gamma,\square},\underset{\mathrm{Spec}}{\mathcal{O}}^\mathrm{CS}{\Phi}_{*,\Gamma,\square},\underset{\mathrm{Spec}}{\mathcal{O}}^\mathrm{CS}{\Delta}^+_{*,\Gamma,\square},\underset{\mathrm{Spec}}{\mathcal{O}}^\mathrm{CS}{\nabla}^+_{*,\Gamma,\square},\\
\underset{\mathrm{Spec}}{\mathcal{O}}^\mathrm{CS}{\Delta}^\dagger_{*,\Gamma,\square},\underset{\mathrm{Spec}}{\mathcal{O}}^\mathrm{CS}{\nabla}^\dagger_{*,\Gamma,\square},\underset{\mathrm{Spec}}{\mathcal{O}}^\mathrm{CS}{\Phi}^r_{*,\Gamma,\square},\underset{\mathrm{Spec}}{\mathcal{O}}^\mathrm{CS}{\Phi}^I_{*,\Gamma,\square}.	
\end{align}

Then we take the corresponding quotients by using the corresponding Frobenius operators:
\begin{align}
&\underset{\mathrm{Spec}}{\mathcal{O}}^\mathrm{CS}\widetilde{\Delta}_{*,\Gamma,\square}/\mathrm{Fro}^\mathbb{Z},\underset{\mathrm{Spec}}{\mathcal{O}}^\mathrm{CS}\widetilde{\nabla}_{*,\Gamma,\square}/\mathrm{Fro}^\mathbb{Z},\underset{\mathrm{Spec}}{\mathcal{O}}^\mathrm{CS}\widetilde{\Phi}_{*,\Gamma,\square}/\mathrm{Fro}^\mathbb{Z},\underset{\mathrm{Spec}}{\mathcal{O}}^\mathrm{CS}\widetilde{\Delta}^+_{*,\Gamma,\square}/\mathrm{Fro}^\mathbb{Z},\\
&\underset{\mathrm{Spec}}{\mathcal{O}}^\mathrm{CS}\widetilde{\nabla}^+_{*,\Gamma,\square}/\mathrm{Fro}^\mathbb{Z}, \underset{\mathrm{Spec}}{\mathcal{O}}^\mathrm{CS}\widetilde{\Delta}^\dagger_{*,\Gamma,\square}/\mathrm{Fro}^\mathbb{Z},\underset{\mathrm{Spec}}{\mathcal{O}}^\mathrm{CS}\widetilde{\nabla}^\dagger_{*,\Gamma,\square}/\mathrm{Fro}^\mathbb{Z},	\\
\end{align}
\begin{align}
&\underset{\mathrm{Spec}}{\mathcal{O}}^\mathrm{CS}\breve{\Delta}_{*,\Gamma,\square}/\mathrm{Fro}^\mathbb{Z},\breve{\nabla}_{*,\Gamma,\square}/\mathrm{Fro}^\mathbb{Z},\underset{\mathrm{Spec}}{\mathcal{O}}^\mathrm{CS}\breve{\Phi}_{*,\Gamma,\square}/\mathrm{Fro}^\mathbb{Z},\underset{\mathrm{Spec}}{\mathcal{O}}^\mathrm{CS}\breve{\Delta}^+_{*,\Gamma,\square}/\mathrm{Fro}^\mathbb{Z},\\
&\underset{\mathrm{Spec}}{\mathcal{O}}^\mathrm{CS}\breve{\nabla}^+_{*,\Gamma,\square}/\mathrm{Fro}^\mathbb{Z}, \underset{\mathrm{Spec}}{\mathcal{O}}^\mathrm{CS}\breve{\Delta}^\dagger_{*,\Gamma,\square}/\mathrm{Fro}^\mathbb{Z},\underset{\mathrm{Spec}}{\mathcal{O}}^\mathrm{CS}\breve{\nabla}^\dagger_{*,\Gamma,\square}/\mathrm{Fro}^\mathbb{Z},	\\
\end{align}
\begin{align}
&\underset{\mathrm{Spec}}{\mathcal{O}}^\mathrm{CS}{\Delta}_{*,\Gamma,\square}/\mathrm{Fro}^\mathbb{Z},\underset{\mathrm{Spec}}{\mathcal{O}}^\mathrm{CS}{\nabla}_{*,\Gamma,\square}/\mathrm{Fro}^\mathbb{Z},\underset{\mathrm{Spec}}{\mathcal{O}}^\mathrm{CS}{\Phi}_{*,\Gamma,\square}/\mathrm{Fro}^\mathbb{Z},\underset{\mathrm{Spec}}{\mathcal{O}}^\mathrm{CS}{\Delta}^+_{*,\Gamma,\square}/\mathrm{Fro}^\mathbb{Z},\\
&\underset{\mathrm{Spec}}{\mathcal{O}}^\mathrm{CS}{\nabla}^+_{*,\Gamma,\square}/\mathrm{Fro}^\mathbb{Z}, \underset{\mathrm{Spec}}{\mathcal{O}}^\mathrm{CS}{\Delta}^\dagger_{*,\Gamma,\square}/\mathrm{Fro}^\mathbb{Z},\underset{\mathrm{Spec}}{\mathcal{O}}^\mathrm{CS}{\nabla}^\dagger_{*,\Gamma,\square}/\mathrm{Fro}^\mathbb{Z}.	
\end{align}
Here for those space with notations related to the radius and the corresponding interval we consider the total unions $\bigcap_r,\bigcup_I$ in order to achieve the whole spaces to achieve the analogues of the corresponding FF curves from \cite{10KL1}, \cite{10KL2}, \cite{10FF} for
\[
\xymatrix@R+0pc@C+0pc{
\underset{r}{\mathrm{homotopycolimit}}~\underset{\mathrm{Spec}}{\mathcal{O}}^\mathrm{CS}\widetilde{\Phi}^r_{*,\Gamma,\square},\underset{I}{\mathrm{homotopylimit}}~\underset{\mathrm{Spec}}{\mathcal{O}}^\mathrm{CS}\widetilde{\Phi}^I_{*,\Gamma,\square},	\\
}
\]
\[
\xymatrix@R+0pc@C+0pc{
\underset{r}{\mathrm{homotopycolimit}}~\underset{\mathrm{Spec}}{\mathcal{O}}^\mathrm{CS}\breve{\Phi}^r_{*,\Gamma,\square},\underset{I}{\mathrm{homotopylimit}}~\underset{\mathrm{Spec}}{\mathcal{O}}^\mathrm{CS}\breve{\Phi}^I_{*,\Gamma,\square},	\\
}
\]
\[
\xymatrix@R+0pc@C+0pc{
\underset{r}{\mathrm{homotopycolimit}}~\underset{\mathrm{Spec}}{\mathcal{O}}^\mathrm{CS}{\Phi}^r_{*,\Gamma,\square},\underset{I}{\mathrm{homotopylimit}}~\underset{\mathrm{Spec}}{\mathcal{O}}^\mathrm{CS}{\Phi}^I_{*,\Gamma,\square}.	
}
\]
\[ 
\xymatrix@R+0pc@C+0pc{
\underset{r}{\mathrm{homotopycolimit}}~\underset{\mathrm{Spec}}{\mathcal{O}}^\mathrm{CS}\widetilde{\Phi}^r_{*,\Gamma,\square}/\mathrm{Fro}^\mathbb{Z},\underset{I}{\mathrm{homotopylimit}}~\underset{\mathrm{Spec}}{\mathcal{O}}^\mathrm{CS}\widetilde{\Phi}^I_{*,\Gamma,\square}/\mathrm{Fro}^\mathbb{Z},	\\
}
\]
\[ 
\xymatrix@R+0pc@C+0pc{
\underset{r}{\mathrm{homotopycolimit}}~\underset{\mathrm{Spec}}{\mathcal{O}}^\mathrm{CS}\breve{\Phi}^r_{*,\Gamma,\square}/\mathrm{Fro}^\mathbb{Z},\underset{I}{\mathrm{homotopylimit}}~\breve{\Phi}^I_{*,\Gamma,\square}/\mathrm{Fro}^\mathbb{Z},	\\
}
\]
\[ 
\xymatrix@R+0pc@C+0pc{
\underset{r}{\mathrm{homotopycolimit}}~\underset{\mathrm{Spec}}{\mathcal{O}}^\mathrm{CS}{\Phi}^r_{*,\Gamma,\square}/\mathrm{Fro}^\mathbb{Z},\underset{I}{\mathrm{homotopylimit}}~\underset{\mathrm{Spec}}{\mathcal{O}}^\mathrm{CS}{\Phi}^I_{*,\Gamma,\square}/\mathrm{Fro}^\mathbb{Z}.	
}
\]

\end{definition}

\

\begin{definition}
We then consider the corresponding quasipresheaves of the corresponding ind-Banach or monomorphic ind-Banach modules from \cite{10BBK}, \cite{10KKM}:
\begin{align}
\mathrm{Quasicoherentpresheaves,IndBanach}_{*}	
\end{align}
where $*$ is one of the following spaces:
\begin{align}
&\underset{\mathrm{Spec}}{\mathcal{O}}^\mathrm{BK}\widetilde{\Phi}_{*,\Gamma,\square}/\mathrm{Fro}^\mathbb{Z},	\\
\end{align}
\begin{align}
&\underset{\mathrm{Spec}}{\mathcal{O}}^\mathrm{BK}\breve{\Phi}_{*,\Gamma,\square}/\mathrm{Fro}^\mathbb{Z},	\\
\end{align}
\begin{align}
&\underset{\mathrm{Spec}}{\mathcal{O}}^\mathrm{BK}{\Phi}_{*,\Gamma,\square}/\mathrm{Fro}^\mathbb{Z}.	
\end{align}
Here for those space without notation related to the radius and the corresponding interval we consider the total unions $\bigcap_r,\bigcup_I$ in order to achieve the whole spaces to achieve the analogues of the corresponding FF curves from \cite{10KL1}, \cite{10KL2}, \cite{10FF} for
\[
\xymatrix@R+0pc@C+0pc{
\underset{r}{\mathrm{homotopycolimit}}~\underset{\mathrm{Spec}}{\mathcal{O}}^\mathrm{BK}\widetilde{\Phi}^r_{*,\Gamma,\square},\underset{I}{\mathrm{homotopylimit}}~\underset{\mathrm{Spec}}{\mathcal{O}}^\mathrm{BK}\widetilde{\Phi}^I_{*,\Gamma,\square},	\\
}
\]
\[
\xymatrix@R+0pc@C+0pc{
\underset{r}{\mathrm{homotopycolimit}}~\underset{\mathrm{Spec}}{\mathcal{O}}^\mathrm{BK}\breve{\Phi}^r_{*,\Gamma,\square},\underset{I}{\mathrm{homotopylimit}}~\underset{\mathrm{Spec}}{\mathcal{O}}^\mathrm{BK}\breve{\Phi}^I_{*,\Gamma,\square},	\\
}
\]
\[
\xymatrix@R+0pc@C+0pc{
\underset{r}{\mathrm{homotopycolimit}}~\underset{\mathrm{Spec}}{\mathcal{O}}^\mathrm{BK}{\Phi}^r_{*,\Gamma,\square},\underset{I}{\mathrm{homotopylimit}}~\underset{\mathrm{Spec}}{\mathcal{O}}^\mathrm{BK}{\Phi}^I_{*,\Gamma,\square}.	
}
\]
\[  
\xymatrix@R+0pc@C+0pc{
\underset{r}{\mathrm{homotopycolimit}}~\underset{\mathrm{Spec}}{\mathcal{O}}^\mathrm{BK}\widetilde{\Phi}^r_{*,\Gamma,\square}/\mathrm{Fro}^\mathbb{Z},\underset{I}{\mathrm{homotopylimit}}~\underset{\mathrm{Spec}}{\mathcal{O}}^\mathrm{BK}\widetilde{\Phi}^I_{*,\Gamma,\square}/\mathrm{Fro}^\mathbb{Z},	\\
}
\]
\[ 
\xymatrix@R+0pc@C+0pc{
\underset{r}{\mathrm{homotopycolimit}}~\underset{\mathrm{Spec}}{\mathcal{O}}^\mathrm{BK}\breve{\Phi}^r_{*,\Gamma,\square}/\mathrm{Fro}^\mathbb{Z},\underset{I}{\mathrm{homotopylimit}}~\underset{\mathrm{Spec}}{\mathcal{O}}^\mathrm{BK}\breve{\Phi}^I_{*,\Gamma,\square}/\mathrm{Fro}^\mathbb{Z},	\\
}
\]
\[ 
\xymatrix@R+0pc@C+0pc{
\underset{r}{\mathrm{homotopycolimit}}~\underset{\mathrm{Spec}}{\mathcal{O}}^\mathrm{BK}{\Phi}^r_{*,\Gamma,\square}/\mathrm{Fro}^\mathbb{Z},\underset{I}{\mathrm{homotopylimit}}~\underset{\mathrm{Spec}}{\mathcal{O}}^\mathrm{BK}{\Phi}^I_{*,\Gamma,\square}/\mathrm{Fro}^\mathbb{Z}.	
}
\]

\end{definition}

\begin{definition}
We then consider the corresponding quasisheaves of the corresponding condensed solid topological modules from \cite{10CS2}:
\begin{align}
\mathrm{Quasicoherentsheaves, Condensed}_{*}	
\end{align}
where $*$ is one of the following spaces:
\begin{align}
&\underset{\mathrm{Spec}}{\mathcal{O}}^\mathrm{CS}\widetilde{\Delta}_{*,\Gamma,\square}/\mathrm{Fro}^\mathbb{Z},\underset{\mathrm{Spec}}{\mathcal{O}}^\mathrm{CS}\widetilde{\nabla}_{*,\Gamma,\square}/\mathrm{Fro}^\mathbb{Z},\underset{\mathrm{Spec}}{\mathcal{O}}^\mathrm{CS}\widetilde{\Phi}_{*,\Gamma,\square}/\mathrm{Fro}^\mathbb{Z},\underset{\mathrm{Spec}}{\mathcal{O}}^\mathrm{CS}\widetilde{\Delta}^+_{*,\Gamma,\square}/\mathrm{Fro}^\mathbb{Z},\\
&\underset{\mathrm{Spec}}{\mathcal{O}}^\mathrm{CS}\widetilde{\nabla}^+_{*,\Gamma,\square}/\mathrm{Fro}^\mathbb{Z},\underset{\mathrm{Spec}}{\mathcal{O}}^\mathrm{CS}\widetilde{\Delta}^\dagger_{*,\Gamma,\square}/\mathrm{Fro}^\mathbb{Z},\underset{\mathrm{Spec}}{\mathcal{O}}^\mathrm{CS}\widetilde{\nabla}^\dagger_{*,\Gamma,\square}/\mathrm{Fro}^\mathbb{Z},	\\
\end{align}
\begin{align}
&\underset{\mathrm{Spec}}{\mathcal{O}}^\mathrm{CS}\breve{\Delta}_{*,\Gamma,\square}/\mathrm{Fro}^\mathbb{Z},\breve{\nabla}_{*,\Gamma,\square}/\mathrm{Fro}^\mathbb{Z},\underset{\mathrm{Spec}}{\mathcal{O}}^\mathrm{CS}\breve{\Phi}_{*,\Gamma,\square}/\mathrm{Fro}^\mathbb{Z},\underset{\mathrm{Spec}}{\mathcal{O}}^\mathrm{CS}\breve{\Delta}^+_{*,\Gamma,\square}/\mathrm{Fro}^\mathbb{Z},\\
&\underset{\mathrm{Spec}}{\mathcal{O}}^\mathrm{CS}\breve{\nabla}^+_{*,\Gamma,\square}/\mathrm{Fro}^\mathbb{Z},\underset{\mathrm{Spec}}{\mathcal{O}}^\mathrm{CS}\breve{\Delta}^\dagger_{*,\Gamma,\square}/\mathrm{Fro}^\mathbb{Z},\underset{\mathrm{Spec}}{\mathcal{O}}^\mathrm{CS}\breve{\nabla}^\dagger_{*,\Gamma,\square}/\mathrm{Fro}^\mathbb{Z},	\\
\end{align}
\begin{align}
&\underset{\mathrm{Spec}}{\mathcal{O}}^\mathrm{CS}{\Delta}_{*,\Gamma,\square}/\mathrm{Fro}^\mathbb{Z},\underset{\mathrm{Spec}}{\mathcal{O}}^\mathrm{CS}{\nabla}_{*,\Gamma,\square}/\mathrm{Fro}^\mathbb{Z},\underset{\mathrm{Spec}}{\mathcal{O}}^\mathrm{CS}{\Phi}_{*,\Gamma,\square}/\mathrm{Fro}^\mathbb{Z},\underset{\mathrm{Spec}}{\mathcal{O}}^\mathrm{CS}{\Delta}^+_{*,\Gamma,\square}/\mathrm{Fro}^\mathbb{Z},\\
&\underset{\mathrm{Spec}}{\mathcal{O}}^\mathrm{CS}{\nabla}^+_{*,\Gamma,\square}/\mathrm{Fro}^\mathbb{Z}, \underset{\mathrm{Spec}}{\mathcal{O}}^\mathrm{CS}{\Delta}^\dagger_{*,\Gamma,\square}/\mathrm{Fro}^\mathbb{Z},\underset{\mathrm{Spec}}{\mathcal{O}}^\mathrm{CS}{\nabla}^\dagger_{*,\Gamma,\square}/\mathrm{Fro}^\mathbb{Z}.	
\end{align}
Here for those space with notations related to the radius and the corresponding interval we consider the total unions $\bigcap_r,\bigcup_I$ in order to achieve the whole spaces to achieve the analogues of the corresponding FF curves from \cite{10KL1}, \cite{10KL2}, \cite{10FF} for
\[
\xymatrix@R+0pc@C+0pc{
\underset{r}{\mathrm{homotopycolimit}}~\underset{\mathrm{Spec}}{\mathcal{O}}^\mathrm{CS}\widetilde{\Phi}^r_{*,\Gamma,\square},\underset{I}{\mathrm{homotopylimit}}~\underset{\mathrm{Spec}}{\mathcal{O}}^\mathrm{CS}\widetilde{\Phi}^I_{*,\Gamma,\square},	\\
}
\]
\[
\xymatrix@R+0pc@C+0pc{
\underset{r}{\mathrm{homotopycolimit}}~\underset{\mathrm{Spec}}{\mathcal{O}}^\mathrm{CS}\breve{\Phi}^r_{*,\Gamma,\square},\underset{I}{\mathrm{homotopylimit}}~\underset{\mathrm{Spec}}{\mathcal{O}}^\mathrm{CS}\breve{\Phi}^I_{*,\Gamma,\square},	\\
}
\]
\[
\xymatrix@R+0pc@C+0pc{
\underset{r}{\mathrm{homotopycolimit}}~\underset{\mathrm{Spec}}{\mathcal{O}}^\mathrm{CS}{\Phi}^r_{*,\Gamma,\square},\underset{I}{\mathrm{homotopylimit}}~\underset{\mathrm{Spec}}{\mathcal{O}}^\mathrm{CS}{\Phi}^I_{*,\Gamma,\square}.	
}
\]
\[ 
\xymatrix@R+0pc@C+0pc{
\underset{r}{\mathrm{homotopycolimit}}~\underset{\mathrm{Spec}}{\mathcal{O}}^\mathrm{CS}\widetilde{\Phi}^r_{*,\Gamma,\square}/\mathrm{Fro}^\mathbb{Z},\underset{I}{\mathrm{homotopylimit}}~\underset{\mathrm{Spec}}{\mathcal{O}}^\mathrm{CS}\widetilde{\Phi}^I_{*,\Gamma,\square}/\mathrm{Fro}^\mathbb{Z},	\\
}
\]
\[ 
\xymatrix@R+0pc@C+0pc{
\underset{r}{\mathrm{homotopycolimit}}~\underset{\mathrm{Spec}}{\mathcal{O}}^\mathrm{CS}\breve{\Phi}^r_{*,\Gamma,\square}/\mathrm{Fro}^\mathbb{Z},\underset{I}{\mathrm{homotopylimit}}~\breve{\Phi}^I_{*,\Gamma,\square}/\mathrm{Fro}^\mathbb{Z},	\\
}
\]
\[ 
\xymatrix@R+0pc@C+0pc{
\underset{r}{\mathrm{homotopycolimit}}~\underset{\mathrm{Spec}}{\mathcal{O}}^\mathrm{CS}{\Phi}^r_{*,\Gamma,\square}/\mathrm{Fro}^\mathbb{Z},\underset{I}{\mathrm{homotopylimit}}~\underset{\mathrm{Spec}}{\mathcal{O}}^\mathrm{CS}{\Phi}^I_{*,\Gamma,\square}/\mathrm{Fro}^\mathbb{Z}.	
}
\]

\end{definition}

\

\begin{proposition}
There is a well-defined functor from the $\infty$-category 
\begin{align}
\mathrm{Quasicoherentpresheaves,Condensed}_{*}	
\end{align}
where $*$ is one of the following spaces:
\begin{align}
&\underset{\mathrm{Spec}}{\mathcal{O}}^\mathrm{CS}\widetilde{\Phi}_{*,\Gamma,\square}/\mathrm{Fro}^\mathbb{Z},	\\
\end{align}
\begin{align}
&\underset{\mathrm{Spec}}{\mathcal{O}}^\mathrm{CS}\breve{\Phi}_{*,\Gamma,\square}/\mathrm{Fro}^\mathbb{Z},	\\
\end{align}
\begin{align}
&\underset{\mathrm{Spec}}{\mathcal{O}}^\mathrm{CS}{\Phi}_{*,\Gamma,\square}/\mathrm{Fro}^\mathbb{Z},	
\end{align}
to the $\infty$-category of $\mathrm{Fro}$-equivariant quasicoherent presheaves over similar spaces above correspondingly without the $\mathrm{Fro}$-quotients, and to the $\infty$-category of $\mathrm{Fro}$-equivariant quasicoherent modules over global sections of the structure $\infty$-sheaves of the similar spaces above correspondingly without the $\mathrm{Fro}$-quotients. Here for those space without notation related to the radius and the corresponding interval we consider the total unions $\bigcap_r,\bigcup_I$ in order to achieve the whole spaces to achieve the analogues of the corresponding FF curves from \cite{10KL1}, \cite{10KL2}, \cite{10FF} for
\[
\xymatrix@R+0pc@C+0pc{
\underset{r}{\mathrm{homotopycolimit}}~\underset{\mathrm{Spec}}{\mathcal{O}}^\mathrm{CS}\widetilde{\Phi}^r_{*,\Gamma,\square},\underset{I}{\mathrm{homotopylimit}}~\underset{\mathrm{Spec}}{\mathcal{O}}^\mathrm{CS}\widetilde{\Phi}^I_{*,\Gamma,\square},	\\
}
\]
\[
\xymatrix@R+0pc@C+0pc{
\underset{r}{\mathrm{homotopycolimit}}~\underset{\mathrm{Spec}}{\mathcal{O}}^\mathrm{CS}\breve{\Phi}^r_{*,\Gamma,\square},\underset{I}{\mathrm{homotopylimit}}~\underset{\mathrm{Spec}}{\mathcal{O}}^\mathrm{CS}\breve{\Phi}^I_{*,\Gamma,\square},	\\
}
\]
\[
\xymatrix@R+0pc@C+0pc{
\underset{r}{\mathrm{homotopycolimit}}~\underset{\mathrm{Spec}}{\mathcal{O}}^\mathrm{CS}{\Phi}^r_{*,\Gamma,\square},\underset{I}{\mathrm{homotopylimit}}~\underset{\mathrm{Spec}}{\mathcal{O}}^\mathrm{CS}{\Phi}^I_{*,\Gamma,\square}.	
}
\]
\[ 
\xymatrix@R+0pc@C+0pc{
\underset{r}{\mathrm{homotopycolimit}}~\underset{\mathrm{Spec}}{\mathcal{O}}^\mathrm{CS}\widetilde{\Phi}^r_{*,\Gamma,\square}/\mathrm{Fro}^\mathbb{Z},\underset{I}{\mathrm{homotopylimit}}~\underset{\mathrm{Spec}}{\mathcal{O}}^\mathrm{CS}\widetilde{\Phi}^I_{*,\Gamma,\square}/\mathrm{Fro}^\mathbb{Z},	\\
}
\]
\[ 
\xymatrix@R+0pc@C+0pc{
\underset{r}{\mathrm{homotopycolimit}}~\underset{\mathrm{Spec}}{\mathcal{O}}^\mathrm{CS}\breve{\Phi}^r_{*,\Gamma,\square}/\mathrm{Fro}^\mathbb{Z},\underset{I}{\mathrm{homotopylimit}}~\breve{\Phi}^I_{*,\Gamma,\square}/\mathrm{Fro}^\mathbb{Z},	\\
}
\]
\[ 
\xymatrix@R+0pc@C+0pc{
\underset{r}{\mathrm{homotopycolimit}}~\underset{\mathrm{Spec}}{\mathcal{O}}^\mathrm{CS}{\Phi}^r_{*,\Gamma,\square}/\mathrm{Fro}^\mathbb{Z},\underset{I}{\mathrm{homotopylimit}}~\underset{\mathrm{Spec}}{\mathcal{O}}^\mathrm{CS}{\Phi}^I_{*,\Gamma,\square}/\mathrm{Fro}^\mathbb{Z}.	
}
\]	
In this situation we will have the target category being family parametrized by $r$ or $I$ in compatible glueing sense as in \cite[Definition 5.4.10]{10KL2}. In this situation for modules parametrized by the intervals we have the equivalence of $\infty$-categories by using \cite[Proposition 13.8]{10CS2}. Here the corresponding quasicoherent Frobenius modules are defined to be the corresponding homotopy colimits and limits of Frobenius modules:
\begin{align}
\underset{r}{\mathrm{homotopycolimit}}~M_r,\\
\underset{I}{\mathrm{homotopylimit}}~M_I,	
\end{align}
where each $M_r$ is a Frobenius-equivariant module over the period ring with respect to some radius $r$ while each $M_I$ is a Frobenius-equivariant module over the period ring with respect to some interval $I$.\\
\end{proposition}

\begin{proposition}
Similar proposition holds for 
\begin{align}
\mathrm{Quasicoherentsheaves,IndBanach}_{*}.	
\end{align}	
\end{proposition}

\

\begin{definition}
We then consider the corresponding quasipresheaves of perfect complexes the corresponding ind-Banach or monomorphic ind-Banach modules from \cite{10BBK}, \cite{10KKM}:
\begin{align}
\mathrm{Quasicoherentpresheaves,Perfectcomplex,IndBanach}_{*}	
\end{align}
where $*$ is one of the following spaces:
\begin{align}
&\underset{\mathrm{Spec}}{\mathcal{O}}^\mathrm{BK}\widetilde{\Phi}_{*,\Gamma,\square}/\mathrm{Fro}^\mathbb{Z},	\\
\end{align}
\begin{align}
&\underset{\mathrm{Spec}}{\mathcal{O}}^\mathrm{BK}\breve{\Phi}_{*,\Gamma,\square}/\mathrm{Fro}^\mathbb{Z},	\\
\end{align}
\begin{align}
&\underset{\mathrm{Spec}}{\mathcal{O}}^\mathrm{BK}{\Phi}_{*,\Gamma,\square}/\mathrm{Fro}^\mathbb{Z}.	
\end{align}
Here for those space without notation related to the radius and the corresponding interval we consider the total unions $\bigcap_r,\bigcup_I$ in order to achieve the whole spaces to achieve the analogues of the corresponding FF curves from \cite{10KL1}, \cite{10KL2}, \cite{10FF} for
\[
\xymatrix@R+0pc@C+0pc{
\underset{r}{\mathrm{homotopycolimit}}~\underset{\mathrm{Spec}}{\mathcal{O}}^\mathrm{BK}\widetilde{\Phi}^r_{*,\Gamma,\square},\underset{I}{\mathrm{homotopylimit}}~\underset{\mathrm{Spec}}{\mathcal{O}}^\mathrm{BK}\widetilde{\Phi}^I_{*,\Gamma,\square},	\\
}
\]
\[
\xymatrix@R+0pc@C+0pc{
\underset{r}{\mathrm{homotopycolimit}}~\underset{\mathrm{Spec}}{\mathcal{O}}^\mathrm{BK}\breve{\Phi}^r_{*,\Gamma,\square},\underset{I}{\mathrm{homotopylimit}}~\underset{\mathrm{Spec}}{\mathcal{O}}^\mathrm{BK}\breve{\Phi}^I_{*,\Gamma,\square},	\\
}
\]
\[
\xymatrix@R+0pc@C+0pc{
\underset{r}{\mathrm{homotopycolimit}}~\underset{\mathrm{Spec}}{\mathcal{O}}^\mathrm{BK}{\Phi}^r_{*,\Gamma,\square},\underset{I}{\mathrm{homotopylimit}}~\underset{\mathrm{Spec}}{\mathcal{O}}^\mathrm{BK}{\Phi}^I_{*,\Gamma,\square}.	
}
\]
\[  
\xymatrix@R+0pc@C+0pc{
\underset{r}{\mathrm{homotopycolimit}}~\underset{\mathrm{Spec}}{\mathcal{O}}^\mathrm{BK}\widetilde{\Phi}^r_{*,\Gamma,\square}/\mathrm{Fro}^\mathbb{Z},\underset{I}{\mathrm{homotopylimit}}~\underset{\mathrm{Spec}}{\mathcal{O}}^\mathrm{BK}\widetilde{\Phi}^I_{*,\Gamma,\square}/\mathrm{Fro}^\mathbb{Z},	\\
}
\]
\[ 
\xymatrix@R+0pc@C+0pc{
\underset{r}{\mathrm{homotopycolimit}}~\underset{\mathrm{Spec}}{\mathcal{O}}^\mathrm{BK}\breve{\Phi}^r_{*,\Gamma,\square}/\mathrm{Fro}^\mathbb{Z},\underset{I}{\mathrm{homotopylimit}}~\underset{\mathrm{Spec}}{\mathcal{O}}^\mathrm{BK}\breve{\Phi}^I_{*,\Gamma,\square}/\mathrm{Fro}^\mathbb{Z},	\\
}
\]
\[ 
\xymatrix@R+0pc@C+0pc{
\underset{r}{\mathrm{homotopycolimit}}~\underset{\mathrm{Spec}}{\mathcal{O}}^\mathrm{BK}{\Phi}^r_{*,\Gamma,\square}/\mathrm{Fro}^\mathbb{Z},\underset{I}{\mathrm{homotopylimit}}~\underset{\mathrm{Spec}}{\mathcal{O}}^\mathrm{BK}{\Phi}^I_{*,\Gamma,\square}/\mathrm{Fro}^\mathbb{Z}.	
}
\]

\end{definition}

\begin{definition}
We then consider the corresponding quasisheaves of perfect complexes of the corresponding condensed solid topological modules from \cite{10CS2}:
\begin{align}
\mathrm{Quasicoherentsheaves, Perfectcomplex, Condensed}_{*}	
\end{align}
where $*$ is one of the following spaces:
\begin{align}
&\underset{\mathrm{Spec}}{\mathcal{O}}^\mathrm{CS}\widetilde{\Delta}_{*,\Gamma,\square}/\mathrm{Fro}^\mathbb{Z},\underset{\mathrm{Spec}}{\mathcal{O}}^\mathrm{CS}\widetilde{\nabla}_{*,\Gamma,\square}/\mathrm{Fro}^\mathbb{Z},\underset{\mathrm{Spec}}{\mathcal{O}}^\mathrm{CS}\widetilde{\Phi}_{*,\Gamma,\square}/\mathrm{Fro}^\mathbb{Z},\underset{\mathrm{Spec}}{\mathcal{O}}^\mathrm{CS}\widetilde{\Delta}^+_{*,\Gamma,\square}/\mathrm{Fro}^\mathbb{Z},\\
&\underset{\mathrm{Spec}}{\mathcal{O}}^\mathrm{CS}\widetilde{\nabla}^+_{*,\Gamma,\square}/\mathrm{Fro}^\mathbb{Z},\underset{\mathrm{Spec}}{\mathcal{O}}^\mathrm{CS}\widetilde{\Delta}^\dagger_{*,\Gamma,\square}/\mathrm{Fro}^\mathbb{Z},\underset{\mathrm{Spec}}{\mathcal{O}}^\mathrm{CS}\widetilde{\nabla}^\dagger_{*,\Gamma,\square}/\mathrm{Fro}^\mathbb{Z},	\\
\end{align}
\begin{align}
&\underset{\mathrm{Spec}}{\mathcal{O}}^\mathrm{CS}\breve{\Delta}_{*,\Gamma,\square}/\mathrm{Fro}^\mathbb{Z},\breve{\nabla}_{*,\Gamma,\square}/\mathrm{Fro}^\mathbb{Z},\underset{\mathrm{Spec}}{\mathcal{O}}^\mathrm{CS}\breve{\Phi}_{*,\Gamma,\square}/\mathrm{Fro}^\mathbb{Z},\underset{\mathrm{Spec}}{\mathcal{O}}^\mathrm{CS}\breve{\Delta}^+_{*,\Gamma,\square}/\mathrm{Fro}^\mathbb{Z},\\
&\underset{\mathrm{Spec}}{\mathcal{O}}^\mathrm{CS}\breve{\nabla}^+_{*,\Gamma,\square}/\mathrm{Fro}^\mathbb{Z},\underset{\mathrm{Spec}}{\mathcal{O}}^\mathrm{CS}\breve{\Delta}^\dagger_{*,\Gamma,\square}/\mathrm{Fro}^\mathbb{Z},\underset{\mathrm{Spec}}{\mathcal{O}}^\mathrm{CS}\breve{\nabla}^\dagger_{*,\Gamma,\square}/\mathrm{Fro}^\mathbb{Z},	\\
\end{align}
\begin{align}
&\underset{\mathrm{Spec}}{\mathcal{O}}^\mathrm{CS}{\Delta}_{*,\Gamma,\square}/\mathrm{Fro}^\mathbb{Z},\underset{\mathrm{Spec}}{\mathcal{O}}^\mathrm{CS}{\nabla}_{*,\Gamma,\square}/\mathrm{Fro}^\mathbb{Z},\underset{\mathrm{Spec}}{\mathcal{O}}^\mathrm{CS}{\Phi}_{*,\Gamma,\square}/\mathrm{Fro}^\mathbb{Z},\underset{\mathrm{Spec}}{\mathcal{O}}^\mathrm{CS}{\Delta}^+_{*,\Gamma,\square}/\mathrm{Fro}^\mathbb{Z},\\
&\underset{\mathrm{Spec}}{\mathcal{O}}^\mathrm{CS}{\nabla}^+_{*,\Gamma,\square}/\mathrm{Fro}^\mathbb{Z}, \underset{\mathrm{Spec}}{\mathcal{O}}^\mathrm{CS}{\Delta}^\dagger_{*,\Gamma,\square}/\mathrm{Fro}^\mathbb{Z},\underset{\mathrm{Spec}}{\mathcal{O}}^\mathrm{CS}{\nabla}^\dagger_{*,\Gamma,\square}/\mathrm{Fro}^\mathbb{Z}.	
\end{align}
Here for those space with notations related to the radius and the corresponding interval we consider the total unions $\bigcap_r,\bigcup_I$ in order to achieve the whole spaces to achieve the analogues of the corresponding FF curves from \cite{10KL1}, \cite{10KL2}, \cite{10FF} for
\[
\xymatrix@R+0pc@C+0pc{
\underset{r}{\mathrm{homotopycolimit}}~\underset{\mathrm{Spec}}{\mathcal{O}}^\mathrm{CS}\widetilde{\Phi}^r_{*,\Gamma,\square},\underset{I}{\mathrm{homotopylimit}}~\underset{\mathrm{Spec}}{\mathcal{O}}^\mathrm{CS}\widetilde{\Phi}^I_{*,\Gamma,\square},	\\
}
\]
\[
\xymatrix@R+0pc@C+0pc{
\underset{r}{\mathrm{homotopycolimit}}~\underset{\mathrm{Spec}}{\mathcal{O}}^\mathrm{CS}\breve{\Phi}^r_{*,\Gamma,\square},\underset{I}{\mathrm{homotopylimit}}~\underset{\mathrm{Spec}}{\mathcal{O}}^\mathrm{CS}\breve{\Phi}^I_{*,\Gamma,\square},	\\
}
\]
\[
\xymatrix@R+0pc@C+0pc{
\underset{r}{\mathrm{homotopycolimit}}~\underset{\mathrm{Spec}}{\mathcal{O}}^\mathrm{CS}{\Phi}^r_{*,\Gamma,\square},\underset{I}{\mathrm{homotopylimit}}~\underset{\mathrm{Spec}}{\mathcal{O}}^\mathrm{CS}{\Phi}^I_{*,\Gamma,\square}.	
}
\]
\[ 
\xymatrix@R+0pc@C+0pc{
\underset{r}{\mathrm{homotopycolimit}}~\underset{\mathrm{Spec}}{\mathcal{O}}^\mathrm{CS}\widetilde{\Phi}^r_{*,\Gamma,\square}/\mathrm{Fro}^\mathbb{Z},\underset{I}{\mathrm{homotopylimit}}~\underset{\mathrm{Spec}}{\mathcal{O}}^\mathrm{CS}\widetilde{\Phi}^I_{*,\Gamma,\square}/\mathrm{Fro}^\mathbb{Z},	\\
}
\]
\[ 
\xymatrix@R+0pc@C+0pc{
\underset{r}{\mathrm{homotopycolimit}}~\underset{\mathrm{Spec}}{\mathcal{O}}^\mathrm{CS}\breve{\Phi}^r_{*,\Gamma,\square}/\mathrm{Fro}^\mathbb{Z},\underset{I}{\mathrm{homotopylimit}}~\breve{\Phi}^I_{*,\Gamma,\square}/\mathrm{Fro}^\mathbb{Z},	\\
}
\]
\[ 
\xymatrix@R+0pc@C+0pc{
\underset{r}{\mathrm{homotopycolimit}}~\underset{\mathrm{Spec}}{\mathcal{O}}^\mathrm{CS}{\Phi}^r_{*,\Gamma,\square}/\mathrm{Fro}^\mathbb{Z},\underset{I}{\mathrm{homotopylimit}}~\underset{\mathrm{Spec}}{\mathcal{O}}^\mathrm{CS}{\Phi}^I_{*,\Gamma,\square}/\mathrm{Fro}^\mathbb{Z}.	
}
\]

\end{definition}

\begin{proposition}
There is a well-defined functor from the $\infty$-category 
\begin{align}
\mathrm{Quasicoherentpresheaves,Perfectcomplex,Condensed}_{*}	
\end{align}
where $*$ is one of the following spaces:
\begin{align}
&\underset{\mathrm{Spec}}{\mathcal{O}}^\mathrm{CS}\widetilde{\Phi}_{*,\Gamma,\square}/\mathrm{Fro}^\mathbb{Z},	\\
\end{align}
\begin{align}
&\underset{\mathrm{Spec}}{\mathcal{O}}^\mathrm{CS}\breve{\Phi}_{*,\Gamma,\square}/\mathrm{Fro}^\mathbb{Z},	\\
\end{align}
\begin{align}
&\underset{\mathrm{Spec}}{\mathcal{O}}^\mathrm{CS}{\Phi}_{*,\Gamma,\square}/\mathrm{Fro}^\mathbb{Z},	
\end{align}
to the $\infty$-category of $\mathrm{Fro}$-equivariant quasicoherent presheaves over similar spaces above correspondingly without the $\mathrm{Fro}$-quotients, and to the $\infty$-category of $\mathrm{Fro}$-equivariant quasicoherent modules over global sections of the structure $\infty$-sheaves of the similar spaces above correspondingly without the $\mathrm{Fro}$-quotients. Here for those space without notation related to the radius and the corresponding interval we consider the total unions $\bigcap_r,\bigcup_I$ in order to achieve the whole spaces to achieve the analogues of the corresponding FF curves from \cite{10KL1}, \cite{10KL2}, \cite{10FF} for
\[
\xymatrix@R+0pc@C+0pc{
\underset{r}{\mathrm{homotopycolimit}}~\underset{\mathrm{Spec}}{\mathcal{O}}^\mathrm{CS}\widetilde{\Phi}^r_{*,\Gamma,\square},\underset{I}{\mathrm{homotopylimit}}~\underset{\mathrm{Spec}}{\mathcal{O}}^\mathrm{CS}\widetilde{\Phi}^I_{*,\Gamma,\square},	\\
}
\]
\[
\xymatrix@R+0pc@C+0pc{
\underset{r}{\mathrm{homotopycolimit}}~\underset{\mathrm{Spec}}{\mathcal{O}}^\mathrm{CS}\breve{\Phi}^r_{*,\Gamma,\square},\underset{I}{\mathrm{homotopylimit}}~\underset{\mathrm{Spec}}{\mathcal{O}}^\mathrm{CS}\breve{\Phi}^I_{*,\Gamma,\square},	\\
}
\]
\[
\xymatrix@R+0pc@C+0pc{
\underset{r}{\mathrm{homotopycolimit}}~\underset{\mathrm{Spec}}{\mathcal{O}}^\mathrm{CS}{\Phi}^r_{*,\Gamma,\square},\underset{I}{\mathrm{homotopylimit}}~\underset{\mathrm{Spec}}{\mathcal{O}}^\mathrm{CS}{\Phi}^I_{*,\Gamma,\square}.	
}
\]
\[ 
\xymatrix@R+0pc@C+0pc{
\underset{r}{\mathrm{homotopycolimit}}~\underset{\mathrm{Spec}}{\mathcal{O}}^\mathrm{CS}\widetilde{\Phi}^r_{*,\Gamma,\square}/\mathrm{Fro}^\mathbb{Z},\underset{I}{\mathrm{homotopylimit}}~\underset{\mathrm{Spec}}{\mathcal{O}}^\mathrm{CS}\widetilde{\Phi}^I_{*,\Gamma,\square}/\mathrm{Fro}^\mathbb{Z},	\\
}
\]
\[ 
\xymatrix@R+0pc@C+0pc{
\underset{r}{\mathrm{homotopycolimit}}~\underset{\mathrm{Spec}}{\mathcal{O}}^\mathrm{CS}\breve{\Phi}^r_{*,\Gamma,\square}/\mathrm{Fro}^\mathbb{Z},\underset{I}{\mathrm{homotopylimit}}~\breve{\Phi}^I_{*,\Gamma,\square}/\mathrm{Fro}^\mathbb{Z},	\\
}
\]
\[ 
\xymatrix@R+0pc@C+0pc{
\underset{r}{\mathrm{homotopycolimit}}~\underset{\mathrm{Spec}}{\mathcal{O}}^\mathrm{CS}{\Phi}^r_{*,\Gamma,\square}/\mathrm{Fro}^\mathbb{Z},\underset{I}{\mathrm{homotopylimit}}~\underset{\mathrm{Spec}}{\mathcal{O}}^\mathrm{CS}{\Phi}^I_{*,\Gamma,\square}/\mathrm{Fro}^\mathbb{Z}.	
}
\]	
In this situation we will have the target category being family parametrized by $r$ or $I$ in compatible glueing sense as in \cite[Definition 5.4.10]{10KL2}. In this situation for modules parametrized by the intervals we have the equivalence of $\infty$-categories by using \cite[Proposition 12.18]{10CS2}. Here the corresponding quasicoherent Frobenius modules are defined to be the corresponding homotopy colimits and limits of Frobenius modules:
\begin{align}
\underset{r}{\mathrm{homotopycolimit}}~M_r,\\
\underset{I}{\mathrm{homotopylimit}}~M_I,	
\end{align}
where each $M_r$ is a Frobenius-equivariant module over the period ring with respect to some radius $r$ while each $M_I$ is a Frobenius-equivariant module over the period ring with respect to some interval $I$.\\
\end{proposition}

\begin{proposition}
Similar proposition holds for 
\begin{align}
\mathrm{Quasicoherentsheaves,Perfectcomplex,IndBanach}_{*}.	
\end{align}	
\end{proposition}

\section{Univariate Hodge Iwasawa Prestacks by Deformation}

This chapter follows closely \cite{10T1}, \cite{10T2}, \cite{10T3}, \cite{10T4}, \cite{10T5}, \cite{10T6}, \cite{10KPX}, \cite{10KP}, \cite{10KL1}, \cite{10KL2}, \cite{10BK}, \cite{10BBBK}, \cite{10BBM}, \cite{10KKM}, \cite{10CS1}, \cite{10CS2}, \cite{10LBV}.

\subsection{Frobenius Quasicoherent Prestacks I}

\begin{definition}
First we consider the Bambozzi-Kremnizer spectrum $\underset{\mathrm{Spec}}{\mathcal{O}}^\mathrm{BK}(*)$ attached to any of those in the above from \cite{10BK} by taking derived rational localization:
\begin{align}
&\underset{\mathrm{Spec}}{\mathcal{O}}^\mathrm{BK}\widetilde{\Phi}_{*,A},\underset{\mathrm{Spec}}{\mathcal{O}}^\mathrm{BK}\widetilde{\Phi}^r_{*,A},\underset{\mathrm{Spec}}{\mathcal{O}}^\mathrm{BK}\widetilde{\Phi}^I_{*,A},	
\end{align}
\begin{align}
&\underset{\mathrm{Spec}}{\mathcal{O}}^\mathrm{BK}\breve{\Phi}_{*,A},\underset{\mathrm{Spec}}{\mathcal{O}}^\mathrm{BK}\breve{\Phi}^r_{*,A},\underset{\mathrm{Spec}}{\mathcal{O}}^\mathrm{BK}\breve{\Phi}^I_{*,A},	
\end{align}
\begin{align}
&\underset{\mathrm{Spec}}{\mathcal{O}}^\mathrm{BK}{\Phi}_{*,A},
\underset{\mathrm{Spec}}{\mathcal{O}}^\mathrm{BK}{\Phi}^r_{*,A},\underset{\mathrm{Spec}}{\mathcal{O}}^\mathrm{BK}{\Phi}^I_{*,A}.	
\end{align}

Then we take the corresponding quotients by using the corresponding Frobenius operators:
\begin{align}
&\underset{\mathrm{Spec}}{\mathcal{O}}^\mathrm{BK}\widetilde{\Phi}_{*,A}/\mathrm{Fro}^\mathbb{Z},	\\
\end{align}
\begin{align}
&\underset{\mathrm{Spec}}{\mathcal{O}}^\mathrm{BK}\breve{\Phi}_{*,A}/\mathrm{Fro}^\mathbb{Z},	\\
\end{align}
\begin{align}
&\underset{\mathrm{Spec}}{\mathcal{O}}^\mathrm{BK}{\Phi}_{*,A}/\mathrm{Fro}^\mathbb{Z}.	
\end{align}
Here for those space without notation related to the radius and the corresponding interval we consider the total unions $\bigcap_r,\bigcup_I$ in order to achieve the whole spaces to achieve the analogues of the corresponding FF curves from \cite{10KL1}, \cite{10KL2}, \cite{10FF} for
\[
\xymatrix@R+0pc@C+0pc{
\underset{r}{\mathrm{homotopycolimit}}~\underset{\mathrm{Spec}}{\mathcal{O}}^\mathrm{BK}\widetilde{\Phi}^r_{*,A},\underset{I}{\mathrm{homotopylimit}}~\underset{\mathrm{Spec}}{\mathcal{O}}^\mathrm{BK}\widetilde{\Phi}^I_{*,A},	\\
}
\]
\[
\xymatrix@R+0pc@C+0pc{
\underset{r}{\mathrm{homotopycolimit}}~\underset{\mathrm{Spec}}{\mathcal{O}}^\mathrm{BK}\breve{\Phi}^r_{*,A},\underset{I}{\mathrm{homotopylimit}}~\underset{\mathrm{Spec}}{\mathcal{O}}^\mathrm{BK}\breve{\Phi}^I_{*,A},	\\
}
\]
\[
\xymatrix@R+0pc@C+0pc{
\underset{r}{\mathrm{homotopycolimit}}~\underset{\mathrm{Spec}}{\mathcal{O}}^\mathrm{BK}{\Phi}^r_{*,A},\underset{I}{\mathrm{homotopylimit}}~\underset{\mathrm{Spec}}{\mathcal{O}}^\mathrm{BK}{\Phi}^I_{*,A}.	
}
\]
\[  
\xymatrix@R+0pc@C+0pc{
\underset{r}{\mathrm{homotopycolimit}}~\underset{\mathrm{Spec}}{\mathcal{O}}^\mathrm{BK}\widetilde{\Phi}^r_{*,A}/\mathrm{Fro}^\mathbb{Z},\underset{I}{\mathrm{homotopylimit}}~\underset{\mathrm{Spec}}{\mathcal{O}}^\mathrm{BK}\widetilde{\Phi}^I_{*,A}/\mathrm{Fro}^\mathbb{Z},	\\
}
\]
\[ 
\xymatrix@R+0pc@C+0pc{
\underset{r}{\mathrm{homotopycolimit}}~\underset{\mathrm{Spec}}{\mathcal{O}}^\mathrm{BK}\breve{\Phi}^r_{*,A}/\mathrm{Fro}^\mathbb{Z},\underset{I}{\mathrm{homotopylimit}}~\underset{\mathrm{Spec}}{\mathcal{O}}^\mathrm{BK}\breve{\Phi}^I_{*,A}/\mathrm{Fro}^\mathbb{Z},	\\
}
\]
\[ 
\xymatrix@R+0pc@C+0pc{
\underset{r}{\mathrm{homotopycolimit}}~\underset{\mathrm{Spec}}{\mathcal{O}}^\mathrm{BK}{\Phi}^r_{*,A}/\mathrm{Fro}^\mathbb{Z},\underset{I}{\mathrm{homotopylimit}}~\underset{\mathrm{Spec}}{\mathcal{O}}^\mathrm{BK}{\Phi}^I_{*,A}/\mathrm{Fro}^\mathbb{Z}.	
}
\]

\end{definition}

\indent Meanwhile we have the corresponding Clausen-Scholze analytic stacks from \cite{10CS2}, therefore applying their construction we have:

\begin{definition}
Here we define the following products by using the solidified tensor product from \cite{10CS1} and \cite{10CS2}. Namely $A$ will still as above as a Banach ring over $\mathbb{Q}_p$. Then we take solidified tensor product $\overset{\blacksquare}{\otimes}$ of any of the following
\[
\xymatrix@R+0pc@C+0pc{
\widetilde{\Delta}_{*},\widetilde{\nabla}_{*},\widetilde{\Phi}_{*},\widetilde{\Delta}^+_{*},\widetilde{\nabla}^+_{*},\widetilde{\Delta}^\dagger_{*},\widetilde{\nabla}^\dagger_{*},\widetilde{\Phi}^r_{*},\widetilde{\Phi}^I_{*}, 
}
\]

\[
\xymatrix@R+0pc@C+0pc{
\breve{\Delta}_{*},\breve{\nabla}_{*},\breve{\Phi}_{*},\breve{\Delta}^+_{*},\breve{\nabla}^+_{*},\breve{\Delta}^\dagger_{*},\breve{\nabla}^\dagger_{*},\breve{\Phi}^r_{*},\breve{\Phi}^I_{*},	
}
\]

\[
\xymatrix@R+0pc@C+0pc{
{\Delta}_{*},{\nabla}_{*},{\Phi}_{*},{\Delta}^+_{*},{\nabla}^+_{*},{\Delta}^\dagger_{*},{\nabla}^\dagger_{*},{\Phi}^r_{*},{\Phi}^I_{*},	
}
\]  	
with $A$. Then we have the notations:
\[
\xymatrix@R+0pc@C+0pc{
\widetilde{\Delta}_{*,A},\widetilde{\nabla}_{*,A},\widetilde{\Phi}_{*,A},\widetilde{\Delta}^+_{*,A},\widetilde{\nabla}^+_{*,A},\widetilde{\Delta}^\dagger_{*,A},\widetilde{\nabla}^\dagger_{*,A},\widetilde{\Phi}^r_{*,A},\widetilde{\Phi}^I_{*,A}, 
}
\]

\[
\xymatrix@R+0pc@C+0pc{
\breve{\Delta}_{*,A},\breve{\nabla}_{*,A},\breve{\Phi}_{*,A},\breve{\Delta}^+_{*,A},\breve{\nabla}^+_{*,A},\breve{\Delta}^\dagger_{*,A},\breve{\nabla}^\dagger_{*,A},\breve{\Phi}^r_{*,A},\breve{\Phi}^I_{*,A},	
}
\]

\[
\xymatrix@R+0pc@C+0pc{
{\Delta}_{*,A},{\nabla}_{*,A},{\Phi}_{*,A},{\Delta}^+_{*,A},{\nabla}^+_{*,A},{\Delta}^\dagger_{*,A},{\nabla}^\dagger_{*,A},{\Phi}^r_{*,A},{\Phi}^I_{*,A}.	
}
\]
\end{definition}

\begin{definition}
First we consider the Clausen-Scholze spectrum $\underset{\mathrm{Spec}}{\mathcal{O}}^\mathrm{CS}(*)$ attached to any of those in the above from \cite{10CS2} by taking derived rational localization:
\begin{align}
\underset{\mathrm{Spec}}{\mathcal{O}}^\mathrm{CS}\widetilde{\Delta}_{*,A},\underset{\mathrm{Spec}}{\mathcal{O}}^\mathrm{CS}\widetilde{\nabla}_{*,A},\underset{\mathrm{Spec}}{\mathcal{O}}^\mathrm{CS}\widetilde{\Phi}_{*,A},\underset{\mathrm{Spec}}{\mathcal{O}}^\mathrm{CS}\widetilde{\Delta}^+_{*,A},\underset{\mathrm{Spec}}{\mathcal{O}}^\mathrm{CS}\widetilde{\nabla}^+_{*,A},\\
\underset{\mathrm{Spec}}{\mathcal{O}}^\mathrm{CS}\widetilde{\Delta}^\dagger_{*,A},\underset{\mathrm{Spec}}{\mathcal{O}}^\mathrm{CS}\widetilde{\nabla}^\dagger_{*,A},\underset{\mathrm{Spec}}{\mathcal{O}}^\mathrm{CS}\widetilde{\Phi}^r_{*,A},\underset{\mathrm{Spec}}{\mathcal{O}}^\mathrm{CS}\widetilde{\Phi}^I_{*,A},	\\
\end{align}
\begin{align}
\underset{\mathrm{Spec}}{\mathcal{O}}^\mathrm{CS}\breve{\Delta}_{*,A},\breve{\nabla}_{*,A},\underset{\mathrm{Spec}}{\mathcal{O}}^\mathrm{CS}\breve{\Phi}_{*,A},\underset{\mathrm{Spec}}{\mathcal{O}}^\mathrm{CS}\breve{\Delta}^+_{*,A},\underset{\mathrm{Spec}}{\mathcal{O}}^\mathrm{CS}\breve{\nabla}^+_{*,A},\\
\underset{\mathrm{Spec}}{\mathcal{O}}^\mathrm{CS}\breve{\Delta}^\dagger_{*,A},\underset{\mathrm{Spec}}{\mathcal{O}}^\mathrm{CS}\breve{\nabla}^\dagger_{*,A},\underset{\mathrm{Spec}}{\mathcal{O}}^\mathrm{CS}\breve{\Phi}^r_{*,A},\breve{\Phi}^I_{*,A},	\\
\end{align}
\begin{align}
\underset{\mathrm{Spec}}{\mathcal{O}}^\mathrm{CS}{\Delta}_{*,A},\underset{\mathrm{Spec}}{\mathcal{O}}^\mathrm{CS}{\nabla}_{*,A},\underset{\mathrm{Spec}}{\mathcal{O}}^\mathrm{CS}{\Phi}_{*,A},\underset{\mathrm{Spec}}{\mathcal{O}}^\mathrm{CS}{\Delta}^+_{*,A},\underset{\mathrm{Spec}}{\mathcal{O}}^\mathrm{CS}{\nabla}^+_{*,A},\\
\underset{\mathrm{Spec}}{\mathcal{O}}^\mathrm{CS}{\Delta}^\dagger_{*,A},\underset{\mathrm{Spec}}{\mathcal{O}}^\mathrm{CS}{\nabla}^\dagger_{*,A},\underset{\mathrm{Spec}}{\mathcal{O}}^\mathrm{CS}{\Phi}^r_{*,A},\underset{\mathrm{Spec}}{\mathcal{O}}^\mathrm{CS}{\Phi}^I_{*,A}.	
\end{align}

Then we take the corresponding quotients by using the corresponding Frobenius operators:
\begin{align}
&\underset{\mathrm{Spec}}{\mathcal{O}}^\mathrm{CS}\widetilde{\Delta}_{*,A}/\mathrm{Fro}^\mathbb{Z},\underset{\mathrm{Spec}}{\mathcal{O}}^\mathrm{CS}\widetilde{\nabla}_{*,A}/\mathrm{Fro}^\mathbb{Z},\underset{\mathrm{Spec}}{\mathcal{O}}^\mathrm{CS}\widetilde{\Phi}_{*,A}/\mathrm{Fro}^\mathbb{Z},\underset{\mathrm{Spec}}{\mathcal{O}}^\mathrm{CS}\widetilde{\Delta}^+_{*,A}/\mathrm{Fro}^\mathbb{Z},\\
&\underset{\mathrm{Spec}}{\mathcal{O}}^\mathrm{CS}\widetilde{\nabla}^+_{*,A}/\mathrm{Fro}^\mathbb{Z}, \underset{\mathrm{Spec}}{\mathcal{O}}^\mathrm{CS}\widetilde{\Delta}^\dagger_{*,A}/\mathrm{Fro}^\mathbb{Z},\underset{\mathrm{Spec}}{\mathcal{O}}^\mathrm{CS}\widetilde{\nabla}^\dagger_{*,A}/\mathrm{Fro}^\mathbb{Z},	\\
\end{align}
\begin{align}
&\underset{\mathrm{Spec}}{\mathcal{O}}^\mathrm{CS}\breve{\Delta}_{*,A}/\mathrm{Fro}^\mathbb{Z},\breve{\nabla}_{*,A}/\mathrm{Fro}^\mathbb{Z},\underset{\mathrm{Spec}}{\mathcal{O}}^\mathrm{CS}\breve{\Phi}_{*,A}/\mathrm{Fro}^\mathbb{Z},\underset{\mathrm{Spec}}{\mathcal{O}}^\mathrm{CS}\breve{\Delta}^+_{*,A}/\mathrm{Fro}^\mathbb{Z},\\
&\underset{\mathrm{Spec}}{\mathcal{O}}^\mathrm{CS}\breve{\nabla}^+_{*,A}/\mathrm{Fro}^\mathbb{Z}, \underset{\mathrm{Spec}}{\mathcal{O}}^\mathrm{CS}\breve{\Delta}^\dagger_{*,A}/\mathrm{Fro}^\mathbb{Z},\underset{\mathrm{Spec}}{\mathcal{O}}^\mathrm{CS}\breve{\nabla}^\dagger_{*,A}/\mathrm{Fro}^\mathbb{Z},	\\
\end{align}
\begin{align}
&\underset{\mathrm{Spec}}{\mathcal{O}}^\mathrm{CS}{\Delta}_{*,A}/\mathrm{Fro}^\mathbb{Z},\underset{\mathrm{Spec}}{\mathcal{O}}^\mathrm{CS}{\nabla}_{*,A}/\mathrm{Fro}^\mathbb{Z},\underset{\mathrm{Spec}}{\mathcal{O}}^\mathrm{CS}{\Phi}_{*,A}/\mathrm{Fro}^\mathbb{Z},\underset{\mathrm{Spec}}{\mathcal{O}}^\mathrm{CS}{\Delta}^+_{*,A}/\mathrm{Fro}^\mathbb{Z},\\
&\underset{\mathrm{Spec}}{\mathcal{O}}^\mathrm{CS}{\nabla}^+_{*,A}/\mathrm{Fro}^\mathbb{Z}, \underset{\mathrm{Spec}}{\mathcal{O}}^\mathrm{CS}{\Delta}^\dagger_{*,A}/\mathrm{Fro}^\mathbb{Z},\underset{\mathrm{Spec}}{\mathcal{O}}^\mathrm{CS}{\nabla}^\dagger_{*,A}/\mathrm{Fro}^\mathbb{Z}.	
\end{align}
Here for those space with notations related to the radius and the corresponding interval we consider the total unions $\bigcap_r,\bigcup_I$ in order to achieve the whole spaces to achieve the analogues of the corresponding FF curves from \cite{10KL1}, \cite{10KL2}, \cite{10FF} for
\[
\xymatrix@R+0pc@C+0pc{
\underset{r}{\mathrm{homotopycolimit}}~\underset{\mathrm{Spec}}{\mathcal{O}}^\mathrm{CS}\widetilde{\Phi}^r_{*,A},\underset{I}{\mathrm{homotopylimit}}~\underset{\mathrm{Spec}}{\mathcal{O}}^\mathrm{CS}\widetilde{\Phi}^I_{*,A},	\\
}
\]
\[
\xymatrix@R+0pc@C+0pc{
\underset{r}{\mathrm{homotopycolimit}}~\underset{\mathrm{Spec}}{\mathcal{O}}^\mathrm{CS}\breve{\Phi}^r_{*,A},\underset{I}{\mathrm{homotopylimit}}~\underset{\mathrm{Spec}}{\mathcal{O}}^\mathrm{CS}\breve{\Phi}^I_{*,A},	\\
}
\]
\[
\xymatrix@R+0pc@C+0pc{
\underset{r}{\mathrm{homotopycolimit}}~\underset{\mathrm{Spec}}{\mathcal{O}}^\mathrm{CS}{\Phi}^r_{*,A},\underset{I}{\mathrm{homotopylimit}}~\underset{\mathrm{Spec}}{\mathcal{O}}^\mathrm{CS}{\Phi}^I_{*,A}.	
}
\]
\[ 
\xymatrix@R+0pc@C+0pc{
\underset{r}{\mathrm{homotopycolimit}}~\underset{\mathrm{Spec}}{\mathcal{O}}^\mathrm{CS}\widetilde{\Phi}^r_{*,A}/\mathrm{Fro}^\mathbb{Z},\underset{I}{\mathrm{homotopylimit}}~\underset{\mathrm{Spec}}{\mathcal{O}}^\mathrm{CS}\widetilde{\Phi}^I_{*,A}/\mathrm{Fro}^\mathbb{Z},	\\
}
\]
\[ 
\xymatrix@R+0pc@C+0pc{
\underset{r}{\mathrm{homotopycolimit}}~\underset{\mathrm{Spec}}{\mathcal{O}}^\mathrm{CS}\breve{\Phi}^r_{*,A}/\mathrm{Fro}^\mathbb{Z},\underset{I}{\mathrm{homotopylimit}}~\breve{\Phi}^I_{*,A}/\mathrm{Fro}^\mathbb{Z},	\\
}
\]
\[ 
\xymatrix@R+0pc@C+0pc{
\underset{r}{\mathrm{homotopycolimit}}~\underset{\mathrm{Spec}}{\mathcal{O}}^\mathrm{CS}{\Phi}^r_{*,A}/\mathrm{Fro}^\mathbb{Z},\underset{I}{\mathrm{homotopylimit}}~\underset{\mathrm{Spec}}{\mathcal{O}}^\mathrm{CS}{\Phi}^I_{*,A}/\mathrm{Fro}^\mathbb{Z}.	
}
\]

\end{definition}

\

\begin{definition}
We then consider the corresponding quasipresheaves of the corresponding ind-Banach or monomorphic ind-Banach modules from \cite{10BBK}, \cite{10KKM}:
\begin{align}
\mathrm{Quasicoherentpresheaves,IndBanach}_{*}	
\end{align}
where $*$ is one of the following spaces:
\begin{align}
&\underset{\mathrm{Spec}}{\mathcal{O}}^\mathrm{BK}\widetilde{\Phi}_{*,A}/\mathrm{Fro}^\mathbb{Z},	\\
\end{align}
\begin{align}
&\underset{\mathrm{Spec}}{\mathcal{O}}^\mathrm{BK}\breve{\Phi}_{*,A}/\mathrm{Fro}^\mathbb{Z},	\\
\end{align}
\begin{align}
&\underset{\mathrm{Spec}}{\mathcal{O}}^\mathrm{BK}{\Phi}_{*,A}/\mathrm{Fro}^\mathbb{Z}.	
\end{align}
Here for those space without notation related to the radius and the corresponding interval we consider the total unions $\bigcap_r,\bigcup_I$ in order to achieve the whole spaces to achieve the analogues of the corresponding FF curves from \cite{10KL1}, \cite{10KL2}, \cite{10FF} for
\[
\xymatrix@R+0pc@C+0pc{
\underset{r}{\mathrm{homotopycolimit}}~\underset{\mathrm{Spec}}{\mathcal{O}}^\mathrm{BK}\widetilde{\Phi}^r_{*,A},\underset{I}{\mathrm{homotopylimit}}~\underset{\mathrm{Spec}}{\mathcal{O}}^\mathrm{BK}\widetilde{\Phi}^I_{*,A},	\\
}
\]
\[
\xymatrix@R+0pc@C+0pc{
\underset{r}{\mathrm{homotopycolimit}}~\underset{\mathrm{Spec}}{\mathcal{O}}^\mathrm{BK}\breve{\Phi}^r_{*,A},\underset{I}{\mathrm{homotopylimit}}~\underset{\mathrm{Spec}}{\mathcal{O}}^\mathrm{BK}\breve{\Phi}^I_{*,A},	\\
}
\]
\[
\xymatrix@R+0pc@C+0pc{
\underset{r}{\mathrm{homotopycolimit}}~\underset{\mathrm{Spec}}{\mathcal{O}}^\mathrm{BK}{\Phi}^r_{*,A},\underset{I}{\mathrm{homotopylimit}}~\underset{\mathrm{Spec}}{\mathcal{O}}^\mathrm{BK}{\Phi}^I_{*,A}.	
}
\]
\[  
\xymatrix@R+0pc@C+0pc{
\underset{r}{\mathrm{homotopycolimit}}~\underset{\mathrm{Spec}}{\mathcal{O}}^\mathrm{BK}\widetilde{\Phi}^r_{*,A}/\mathrm{Fro}^\mathbb{Z},\underset{I}{\mathrm{homotopylimit}}~\underset{\mathrm{Spec}}{\mathcal{O}}^\mathrm{BK}\widetilde{\Phi}^I_{*,A}/\mathrm{Fro}^\mathbb{Z},	\\
}
\]
\[ 
\xymatrix@R+0pc@C+0pc{
\underset{r}{\mathrm{homotopycolimit}}~\underset{\mathrm{Spec}}{\mathcal{O}}^\mathrm{BK}\breve{\Phi}^r_{*,A}/\mathrm{Fro}^\mathbb{Z},\underset{I}{\mathrm{homotopylimit}}~\underset{\mathrm{Spec}}{\mathcal{O}}^\mathrm{BK}\breve{\Phi}^I_{*,A}/\mathrm{Fro}^\mathbb{Z},	\\
}
\]
\[ 
\xymatrix@R+0pc@C+0pc{
\underset{r}{\mathrm{homotopycolimit}}~\underset{\mathrm{Spec}}{\mathcal{O}}^\mathrm{BK}{\Phi}^r_{*,A}/\mathrm{Fro}^\mathbb{Z},\underset{I}{\mathrm{homotopylimit}}~\underset{\mathrm{Spec}}{\mathcal{O}}^\mathrm{BK}{\Phi}^I_{*,A}/\mathrm{Fro}^\mathbb{Z}.	
}
\]

\end{definition}

\begin{definition}
We then consider the corresponding quasisheaves of the corresponding condensed solid topological modules from \cite{10CS2}:
\begin{align}
\mathrm{Quasicoherentsheaves, Condensed}_{*}	
\end{align}
where $*$ is one of the following spaces:
\begin{align}
&\underset{\mathrm{Spec}}{\mathcal{O}}^\mathrm{CS}\widetilde{\Delta}_{*,A}/\mathrm{Fro}^\mathbb{Z},\underset{\mathrm{Spec}}{\mathcal{O}}^\mathrm{CS}\widetilde{\nabla}_{*,A}/\mathrm{Fro}^\mathbb{Z},\underset{\mathrm{Spec}}{\mathcal{O}}^\mathrm{CS}\widetilde{\Phi}_{*,A}/\mathrm{Fro}^\mathbb{Z},\underset{\mathrm{Spec}}{\mathcal{O}}^\mathrm{CS}\widetilde{\Delta}^+_{*,A}/\mathrm{Fro}^\mathbb{Z},\\
&\underset{\mathrm{Spec}}{\mathcal{O}}^\mathrm{CS}\widetilde{\nabla}^+_{*,A}/\mathrm{Fro}^\mathbb{Z},\underset{\mathrm{Spec}}{\mathcal{O}}^\mathrm{CS}\widetilde{\Delta}^\dagger_{*,A}/\mathrm{Fro}^\mathbb{Z},\underset{\mathrm{Spec}}{\mathcal{O}}^\mathrm{CS}\widetilde{\nabla}^\dagger_{*,A}/\mathrm{Fro}^\mathbb{Z},	\\
\end{align}
\begin{align}
&\underset{\mathrm{Spec}}{\mathcal{O}}^\mathrm{CS}\breve{\Delta}_{*,A}/\mathrm{Fro}^\mathbb{Z},\breve{\nabla}_{*,A}/\mathrm{Fro}^\mathbb{Z},\underset{\mathrm{Spec}}{\mathcal{O}}^\mathrm{CS}\breve{\Phi}_{*,A}/\mathrm{Fro}^\mathbb{Z},\underset{\mathrm{Spec}}{\mathcal{O}}^\mathrm{CS}\breve{\Delta}^+_{*,A}/\mathrm{Fro}^\mathbb{Z},\\
&\underset{\mathrm{Spec}}{\mathcal{O}}^\mathrm{CS}\breve{\nabla}^+_{*,A}/\mathrm{Fro}^\mathbb{Z},\underset{\mathrm{Spec}}{\mathcal{O}}^\mathrm{CS}\breve{\Delta}^\dagger_{*,A}/\mathrm{Fro}^\mathbb{Z},\underset{\mathrm{Spec}}{\mathcal{O}}^\mathrm{CS}\breve{\nabla}^\dagger_{*,A}/\mathrm{Fro}^\mathbb{Z},	\\
\end{align}
\begin{align}
&\underset{\mathrm{Spec}}{\mathcal{O}}^\mathrm{CS}{\Delta}_{*,A}/\mathrm{Fro}^\mathbb{Z},\underset{\mathrm{Spec}}{\mathcal{O}}^\mathrm{CS}{\nabla}_{*,A}/\mathrm{Fro}^\mathbb{Z},\underset{\mathrm{Spec}}{\mathcal{O}}^\mathrm{CS}{\Phi}_{*,A}/\mathrm{Fro}^\mathbb{Z},\underset{\mathrm{Spec}}{\mathcal{O}}^\mathrm{CS}{\Delta}^+_{*,A}/\mathrm{Fro}^\mathbb{Z},\\
&\underset{\mathrm{Spec}}{\mathcal{O}}^\mathrm{CS}{\nabla}^+_{*,A}/\mathrm{Fro}^\mathbb{Z}, \underset{\mathrm{Spec}}{\mathcal{O}}^\mathrm{CS}{\Delta}^\dagger_{*,A}/\mathrm{Fro}^\mathbb{Z},\underset{\mathrm{Spec}}{\mathcal{O}}^\mathrm{CS}{\nabla}^\dagger_{*,A}/\mathrm{Fro}^\mathbb{Z}.	
\end{align}
Here for those space with notations related to the radius and the corresponding interval we consider the total unions $\bigcap_r,\bigcup_I$ in order to achieve the whole spaces to achieve the analogues of the corresponding FF curves from \cite{10KL1}, \cite{10KL2}, \cite{10FF} for
\[
\xymatrix@R+0pc@C+0pc{
\underset{r}{\mathrm{homotopycolimit}}~\underset{\mathrm{Spec}}{\mathcal{O}}^\mathrm{CS}\widetilde{\Phi}^r_{*,A},\underset{I}{\mathrm{homotopylimit}}~\underset{\mathrm{Spec}}{\mathcal{O}}^\mathrm{CS}\widetilde{\Phi}^I_{*,A},	\\
}
\]
\[
\xymatrix@R+0pc@C+0pc{
\underset{r}{\mathrm{homotopycolimit}}~\underset{\mathrm{Spec}}{\mathcal{O}}^\mathrm{CS}\breve{\Phi}^r_{*,A},\underset{I}{\mathrm{homotopylimit}}~\underset{\mathrm{Spec}}{\mathcal{O}}^\mathrm{CS}\breve{\Phi}^I_{*,A},	\\
}
\]
\[
\xymatrix@R+0pc@C+0pc{
\underset{r}{\mathrm{homotopycolimit}}~\underset{\mathrm{Spec}}{\mathcal{O}}^\mathrm{CS}{\Phi}^r_{*,A},\underset{I}{\mathrm{homotopylimit}}~\underset{\mathrm{Spec}}{\mathcal{O}}^\mathrm{CS}{\Phi}^I_{*,A}.	
}
\]
\[ 
\xymatrix@R+0pc@C+0pc{
\underset{r}{\mathrm{homotopycolimit}}~\underset{\mathrm{Spec}}{\mathcal{O}}^\mathrm{CS}\widetilde{\Phi}^r_{*,A}/\mathrm{Fro}^\mathbb{Z},\underset{I}{\mathrm{homotopylimit}}~\underset{\mathrm{Spec}}{\mathcal{O}}^\mathrm{CS}\widetilde{\Phi}^I_{*,A}/\mathrm{Fro}^\mathbb{Z},	\\
}
\]
\[ 
\xymatrix@R+0pc@C+0pc{
\underset{r}{\mathrm{homotopycolimit}}~\underset{\mathrm{Spec}}{\mathcal{O}}^\mathrm{CS}\breve{\Phi}^r_{*,A}/\mathrm{Fro}^\mathbb{Z},\underset{I}{\mathrm{homotopylimit}}~\breve{\Phi}^I_{*,A}/\mathrm{Fro}^\mathbb{Z},	\\
}
\]
\[ 
\xymatrix@R+0pc@C+0pc{
\underset{r}{\mathrm{homotopycolimit}}~\underset{\mathrm{Spec}}{\mathcal{O}}^\mathrm{CS}{\Phi}^r_{*,A}/\mathrm{Fro}^\mathbb{Z},\underset{I}{\mathrm{homotopylimit}}~\underset{\mathrm{Spec}}{\mathcal{O}}^\mathrm{CS}{\Phi}^I_{*,A}/\mathrm{Fro}^\mathbb{Z}.	
}
\]

\end{definition}

\

\begin{proposition}
There is a well-defined functor from the $\infty$-category 
\begin{align}
\mathrm{Quasicoherentpresheaves,Condensed}_{*}	
\end{align}
where $*$ is one of the following spaces:
\begin{align}
&\underset{\mathrm{Spec}}{\mathcal{O}}^\mathrm{CS}\widetilde{\Phi}_{*,A}/\mathrm{Fro}^\mathbb{Z},	\\
\end{align}
\begin{align}
&\underset{\mathrm{Spec}}{\mathcal{O}}^\mathrm{CS}\breve{\Phi}_{*,A}/\mathrm{Fro}^\mathbb{Z},	\\
\end{align}
\begin{align}
&\underset{\mathrm{Spec}}{\mathcal{O}}^\mathrm{CS}{\Phi}_{*,A}/\mathrm{Fro}^\mathbb{Z},	
\end{align}
to the $\infty$-category of $\mathrm{Fro}$-equivariant quasicoherent presheaves over similar spaces above correspondingly without the $\mathrm{Fro}$-quotients, and to the $\infty$-category of $\mathrm{Fro}$-equivariant quasicoherent modules over global sections of the structure $\infty$-sheaves of the similar spaces above correspondingly without the $\mathrm{Fro}$-quotients. Here for those space without notation related to the radius and the corresponding interval we consider the total unions $\bigcap_r,\bigcup_I$ in order to achieve the whole spaces to achieve the analogues of the corresponding FF curves from \cite{10KL1}, \cite{10KL2}, \cite{10FF} for
\[
\xymatrix@R+0pc@C+0pc{
\underset{r}{\mathrm{homotopycolimit}}~\underset{\mathrm{Spec}}{\mathcal{O}}^\mathrm{CS}\widetilde{\Phi}^r_{*,A},\underset{I}{\mathrm{homotopylimit}}~\underset{\mathrm{Spec}}{\mathcal{O}}^\mathrm{CS}\widetilde{\Phi}^I_{*,A},	\\
}
\]
\[
\xymatrix@R+0pc@C+0pc{
\underset{r}{\mathrm{homotopycolimit}}~\underset{\mathrm{Spec}}{\mathcal{O}}^\mathrm{CS}\breve{\Phi}^r_{*,A},\underset{I}{\mathrm{homotopylimit}}~\underset{\mathrm{Spec}}{\mathcal{O}}^\mathrm{CS}\breve{\Phi}^I_{*,A},	\\
}
\]
\[
\xymatrix@R+0pc@C+0pc{
\underset{r}{\mathrm{homotopycolimit}}~\underset{\mathrm{Spec}}{\mathcal{O}}^\mathrm{CS}{\Phi}^r_{*,A},\underset{I}{\mathrm{homotopylimit}}~\underset{\mathrm{Spec}}{\mathcal{O}}^\mathrm{CS}{\Phi}^I_{*,A}.	
}
\]
\[ 
\xymatrix@R+0pc@C+0pc{
\underset{r}{\mathrm{homotopycolimit}}~\underset{\mathrm{Spec}}{\mathcal{O}}^\mathrm{CS}\widetilde{\Phi}^r_{*,A}/\mathrm{Fro}^\mathbb{Z},\underset{I}{\mathrm{homotopylimit}}~\underset{\mathrm{Spec}}{\mathcal{O}}^\mathrm{CS}\widetilde{\Phi}^I_{*,A}/\mathrm{Fro}^\mathbb{Z},	\\
}
\]
\[ 
\xymatrix@R+0pc@C+0pc{
\underset{r}{\mathrm{homotopycolimit}}~\underset{\mathrm{Spec}}{\mathcal{O}}^\mathrm{CS}\breve{\Phi}^r_{*,A}/\mathrm{Fro}^\mathbb{Z},\underset{I}{\mathrm{homotopylimit}}~\breve{\Phi}^I_{*,A}/\mathrm{Fro}^\mathbb{Z},	\\
}
\]
\[ 
\xymatrix@R+0pc@C+0pc{
\underset{r}{\mathrm{homotopycolimit}}~\underset{\mathrm{Spec}}{\mathcal{O}}^\mathrm{CS}{\Phi}^r_{*,A}/\mathrm{Fro}^\mathbb{Z},\underset{I}{\mathrm{homotopylimit}}~\underset{\mathrm{Spec}}{\mathcal{O}}^\mathrm{CS}{\Phi}^I_{*,A}/\mathrm{Fro}^\mathbb{Z}.	
}
\]	
In this situation we will have the target category being family parametrized by $r$ or $I$ in compatible glueing sense as in \cite[Definition 5.4.10]{10KL2}. In this situation for modules parametrized by the intervals we have the equivalence of $\infty$-categories by using \cite[Proposition 13.8]{10CS2}. Here the corresponding quasicoherent Frobenius modules are defined to be the corresponding homotopy colimits and limits of Frobenius modules:
\begin{align}
\underset{r}{\mathrm{homotopycolimit}}~M_r,\\
\underset{I}{\mathrm{homotopylimit}}~M_I,	
\end{align}
where each $M_r$ is a Frobenius-equivariant module over the period ring with respect to some radius $r$ while each $M_I$ is a Frobenius-equivariant module over the period ring with respect to some interval $I$.\\
\end{proposition}

\begin{proposition}
Similar proposition holds for 
\begin{align}
\mathrm{Quasicoherentsheaves,IndBanach}_{*}.	
\end{align}	
\end{proposition}

\

\begin{definition}
We then consider the corresponding quasipresheaves of perfect complexes the corresponding ind-Banach or monomorphic ind-Banach modules from \cite{10BBK}, \cite{10KKM}:
\begin{align}
\mathrm{Quasicoherentpresheaves,Perfectcomplex,IndBanach}_{*}	
\end{align}
where $*$ is one of the following spaces:
\begin{align}
&\underset{\mathrm{Spec}}{\mathcal{O}}^\mathrm{BK}\widetilde{\Phi}_{*,A}/\mathrm{Fro}^\mathbb{Z},	\\
\end{align}
\begin{align}
&\underset{\mathrm{Spec}}{\mathcal{O}}^\mathrm{BK}\breve{\Phi}_{*,A}/\mathrm{Fro}^\mathbb{Z},	\\
\end{align}
\begin{align}
&\underset{\mathrm{Spec}}{\mathcal{O}}^\mathrm{BK}{\Phi}_{*,A}/\mathrm{Fro}^\mathbb{Z}.	
\end{align}
Here for those space without notation related to the radius and the corresponding interval we consider the total unions $\bigcap_r,\bigcup_I$ in order to achieve the whole spaces to achieve the analogues of the corresponding FF curves from \cite{10KL1}, \cite{10KL2}, \cite{10FF} for
\[
\xymatrix@R+0pc@C+0pc{
\underset{r}{\mathrm{homotopycolimit}}~\underset{\mathrm{Spec}}{\mathcal{O}}^\mathrm{BK}\widetilde{\Phi}^r_{*,A},\underset{I}{\mathrm{homotopylimit}}~\underset{\mathrm{Spec}}{\mathcal{O}}^\mathrm{BK}\widetilde{\Phi}^I_{*,A},	\\
}
\]
\[
\xymatrix@R+0pc@C+0pc{
\underset{r}{\mathrm{homotopycolimit}}~\underset{\mathrm{Spec}}{\mathcal{O}}^\mathrm{BK}\breve{\Phi}^r_{*,A},\underset{I}{\mathrm{homotopylimit}}~\underset{\mathrm{Spec}}{\mathcal{O}}^\mathrm{BK}\breve{\Phi}^I_{*,A},	\\
}
\]
\[
\xymatrix@R+0pc@C+0pc{
\underset{r}{\mathrm{homotopycolimit}}~\underset{\mathrm{Spec}}{\mathcal{O}}^\mathrm{BK}{\Phi}^r_{*,A},\underset{I}{\mathrm{homotopylimit}}~\underset{\mathrm{Spec}}{\mathcal{O}}^\mathrm{BK}{\Phi}^I_{*,A}.	
}
\]
\[  
\xymatrix@R+0pc@C+0pc{
\underset{r}{\mathrm{homotopycolimit}}~\underset{\mathrm{Spec}}{\mathcal{O}}^\mathrm{BK}\widetilde{\Phi}^r_{*,A}/\mathrm{Fro}^\mathbb{Z},\underset{I}{\mathrm{homotopylimit}}~\underset{\mathrm{Spec}}{\mathcal{O}}^\mathrm{BK}\widetilde{\Phi}^I_{*,A}/\mathrm{Fro}^\mathbb{Z},	\\
}
\]
\[ 
\xymatrix@R+0pc@C+0pc{
\underset{r}{\mathrm{homotopycolimit}}~\underset{\mathrm{Spec}}{\mathcal{O}}^\mathrm{BK}\breve{\Phi}^r_{*,A}/\mathrm{Fro}^\mathbb{Z},\underset{I}{\mathrm{homotopylimit}}~\underset{\mathrm{Spec}}{\mathcal{O}}^\mathrm{BK}\breve{\Phi}^I_{*,A}/\mathrm{Fro}^\mathbb{Z},	\\
}
\]
\[ 
\xymatrix@R+0pc@C+0pc{
\underset{r}{\mathrm{homotopycolimit}}~\underset{\mathrm{Spec}}{\mathcal{O}}^\mathrm{BK}{\Phi}^r_{*,A}/\mathrm{Fro}^\mathbb{Z},\underset{I}{\mathrm{homotopylimit}}~\underset{\mathrm{Spec}}{\mathcal{O}}^\mathrm{BK}{\Phi}^I_{*,A}/\mathrm{Fro}^\mathbb{Z}.	
}
\]

\end{definition}

\begin{definition}
We then consider the corresponding quasisheaves of perfect complexes of the corresponding condensed solid topological modules from \cite{10CS2}:
\begin{align}
\mathrm{Quasicoherentsheaves, Perfectcomplex, Condensed}_{*}	
\end{align}
where $*$ is one of the following spaces:
\begin{align}
&\underset{\mathrm{Spec}}{\mathcal{O}}^\mathrm{CS}\widetilde{\Delta}_{*,A}/\mathrm{Fro}^\mathbb{Z},\underset{\mathrm{Spec}}{\mathcal{O}}^\mathrm{CS}\widetilde{\nabla}_{*,A}/\mathrm{Fro}^\mathbb{Z},\underset{\mathrm{Spec}}{\mathcal{O}}^\mathrm{CS}\widetilde{\Phi}_{*,A}/\mathrm{Fro}^\mathbb{Z},\underset{\mathrm{Spec}}{\mathcal{O}}^\mathrm{CS}\widetilde{\Delta}^+_{*,A}/\mathrm{Fro}^\mathbb{Z},\\
&\underset{\mathrm{Spec}}{\mathcal{O}}^\mathrm{CS}\widetilde{\nabla}^+_{*,A}/\mathrm{Fro}^\mathbb{Z},\underset{\mathrm{Spec}}{\mathcal{O}}^\mathrm{CS}\widetilde{\Delta}^\dagger_{*,A}/\mathrm{Fro}^\mathbb{Z},\underset{\mathrm{Spec}}{\mathcal{O}}^\mathrm{CS}\widetilde{\nabla}^\dagger_{*,A}/\mathrm{Fro}^\mathbb{Z},	\\
\end{align}
\begin{align}
&\underset{\mathrm{Spec}}{\mathcal{O}}^\mathrm{CS}\breve{\Delta}_{*,A}/\mathrm{Fro}^\mathbb{Z},\breve{\nabla}_{*,A}/\mathrm{Fro}^\mathbb{Z},\underset{\mathrm{Spec}}{\mathcal{O}}^\mathrm{CS}\breve{\Phi}_{*,A}/\mathrm{Fro}^\mathbb{Z},\underset{\mathrm{Spec}}{\mathcal{O}}^\mathrm{CS}\breve{\Delta}^+_{*,A}/\mathrm{Fro}^\mathbb{Z},\\
&\underset{\mathrm{Spec}}{\mathcal{O}}^\mathrm{CS}\breve{\nabla}^+_{*,A}/\mathrm{Fro}^\mathbb{Z},\underset{\mathrm{Spec}}{\mathcal{O}}^\mathrm{CS}\breve{\Delta}^\dagger_{*,A}/\mathrm{Fro}^\mathbb{Z},\underset{\mathrm{Spec}}{\mathcal{O}}^\mathrm{CS}\breve{\nabla}^\dagger_{*,A}/\mathrm{Fro}^\mathbb{Z},	\\
\end{align}
\begin{align}
&\underset{\mathrm{Spec}}{\mathcal{O}}^\mathrm{CS}{\Delta}_{*,A}/\mathrm{Fro}^\mathbb{Z},\underset{\mathrm{Spec}}{\mathcal{O}}^\mathrm{CS}{\nabla}_{*,A}/\mathrm{Fro}^\mathbb{Z},\underset{\mathrm{Spec}}{\mathcal{O}}^\mathrm{CS}{\Phi}_{*,A}/\mathrm{Fro}^\mathbb{Z},\underset{\mathrm{Spec}}{\mathcal{O}}^\mathrm{CS}{\Delta}^+_{*,A}/\mathrm{Fro}^\mathbb{Z},\\
&\underset{\mathrm{Spec}}{\mathcal{O}}^\mathrm{CS}{\nabla}^+_{*,A}/\mathrm{Fro}^\mathbb{Z}, \underset{\mathrm{Spec}}{\mathcal{O}}^\mathrm{CS}{\Delta}^\dagger_{*,A}/\mathrm{Fro}^\mathbb{Z},\underset{\mathrm{Spec}}{\mathcal{O}}^\mathrm{CS}{\nabla}^\dagger_{*,A}/\mathrm{Fro}^\mathbb{Z}.	
\end{align}
Here for those space with notations related to the radius and the corresponding interval we consider the total unions $\bigcap_r,\bigcup_I$ in order to achieve the whole spaces to achieve the analogues of the corresponding FF curves from \cite{10KL1}, \cite{10KL2}, \cite{10FF} for
\[
\xymatrix@R+0pc@C+0pc{
\underset{r}{\mathrm{homotopycolimit}}~\underset{\mathrm{Spec}}{\mathcal{O}}^\mathrm{CS}\widetilde{\Phi}^r_{*,A},\underset{I}{\mathrm{homotopylimit}}~\underset{\mathrm{Spec}}{\mathcal{O}}^\mathrm{CS}\widetilde{\Phi}^I_{*,A},	\\
}
\]
\[
\xymatrix@R+0pc@C+0pc{
\underset{r}{\mathrm{homotopycolimit}}~\underset{\mathrm{Spec}}{\mathcal{O}}^\mathrm{CS}\breve{\Phi}^r_{*,A},\underset{I}{\mathrm{homotopylimit}}~\underset{\mathrm{Spec}}{\mathcal{O}}^\mathrm{CS}\breve{\Phi}^I_{*,A},	\\
}
\]
\[
\xymatrix@R+0pc@C+0pc{
\underset{r}{\mathrm{homotopycolimit}}~\underset{\mathrm{Spec}}{\mathcal{O}}^\mathrm{CS}{\Phi}^r_{*,A},\underset{I}{\mathrm{homotopylimit}}~\underset{\mathrm{Spec}}{\mathcal{O}}^\mathrm{CS}{\Phi}^I_{*,A}.	
}
\]
\[ 
\xymatrix@R+0pc@C+0pc{
\underset{r}{\mathrm{homotopycolimit}}~\underset{\mathrm{Spec}}{\mathcal{O}}^\mathrm{CS}\widetilde{\Phi}^r_{*,A}/\mathrm{Fro}^\mathbb{Z},\underset{I}{\mathrm{homotopylimit}}~\underset{\mathrm{Spec}}{\mathcal{O}}^\mathrm{CS}\widetilde{\Phi}^I_{*,A}/\mathrm{Fro}^\mathbb{Z},	\\
}
\]
\[ 
\xymatrix@R+0pc@C+0pc{
\underset{r}{\mathrm{homotopycolimit}}~\underset{\mathrm{Spec}}{\mathcal{O}}^\mathrm{CS}\breve{\Phi}^r_{*,A}/\mathrm{Fro}^\mathbb{Z},\underset{I}{\mathrm{homotopylimit}}~\breve{\Phi}^I_{*,A}/\mathrm{Fro}^\mathbb{Z},	\\
}
\]
\[ 
\xymatrix@R+0pc@C+0pc{
\underset{r}{\mathrm{homotopycolimit}}~\underset{\mathrm{Spec}}{\mathcal{O}}^\mathrm{CS}{\Phi}^r_{*,A}/\mathrm{Fro}^\mathbb{Z},\underset{I}{\mathrm{homotopylimit}}~\underset{\mathrm{Spec}}{\mathcal{O}}^\mathrm{CS}{\Phi}^I_{*,A}/\mathrm{Fro}^\mathbb{Z}.	
}
\]

\end{definition}

\begin{proposition}
There is a well-defined functor from the $\infty$-category 
\begin{align}
\mathrm{Quasicoherentpresheaves,Perfectcomplex,Condensed}_{*}	
\end{align}
where $*$ is one of the following spaces:
\begin{align}
&\underset{\mathrm{Spec}}{\mathcal{O}}^\mathrm{CS}\widetilde{\Phi}_{*,A}/\mathrm{Fro}^\mathbb{Z},	\\
\end{align}
\begin{align}
&\underset{\mathrm{Spec}}{\mathcal{O}}^\mathrm{CS}\breve{\Phi}_{*,A}/\mathrm{Fro}^\mathbb{Z},	\\
\end{align}
\begin{align}
&\underset{\mathrm{Spec}}{\mathcal{O}}^\mathrm{CS}{\Phi}_{*,A}/\mathrm{Fro}^\mathbb{Z},	
\end{align}
to the $\infty$-category of $\mathrm{Fro}$-equivariant quasicoherent presheaves over similar spaces above correspondingly without the $\mathrm{Fro}$-quotients, and to the $\infty$-category of $\mathrm{Fro}$-equivariant quasicoherent modules over global sections of the structure $\infty$-sheaves of the similar spaces above correspondingly without the $\mathrm{Fro}$-quotients. Here for those space without notation related to the radius and the corresponding interval we consider the total unions $\bigcap_r,\bigcup_I$ in order to achieve the whole spaces to achieve the analogues of the corresponding FF curves from \cite{10KL1}, \cite{10KL2}, \cite{10FF} for
\[
\xymatrix@R+0pc@C+0pc{
\underset{r}{\mathrm{homotopycolimit}}~\underset{\mathrm{Spec}}{\mathcal{O}}^\mathrm{CS}\widetilde{\Phi}^r_{*,A},\underset{I}{\mathrm{homotopylimit}}~\underset{\mathrm{Spec}}{\mathcal{O}}^\mathrm{CS}\widetilde{\Phi}^I_{*,A},	\\
}
\]
\[
\xymatrix@R+0pc@C+0pc{
\underset{r}{\mathrm{homotopycolimit}}~\underset{\mathrm{Spec}}{\mathcal{O}}^\mathrm{CS}\breve{\Phi}^r_{*,A},\underset{I}{\mathrm{homotopylimit}}~\underset{\mathrm{Spec}}{\mathcal{O}}^\mathrm{CS}\breve{\Phi}^I_{*,A},	\\
}
\]
\[
\xymatrix@R+0pc@C+0pc{
\underset{r}{\mathrm{homotopycolimit}}~\underset{\mathrm{Spec}}{\mathcal{O}}^\mathrm{CS}{\Phi}^r_{*,A},\underset{I}{\mathrm{homotopylimit}}~\underset{\mathrm{Spec}}{\mathcal{O}}^\mathrm{CS}{\Phi}^I_{*,A}.	
}
\]
\[ 
\xymatrix@R+0pc@C+0pc{
\underset{r}{\mathrm{homotopycolimit}}~\underset{\mathrm{Spec}}{\mathcal{O}}^\mathrm{CS}\widetilde{\Phi}^r_{*,A}/\mathrm{Fro}^\mathbb{Z},\underset{I}{\mathrm{homotopylimit}}~\underset{\mathrm{Spec}}{\mathcal{O}}^\mathrm{CS}\widetilde{\Phi}^I_{*,A}/\mathrm{Fro}^\mathbb{Z},	\\
}
\]
\[ 
\xymatrix@R+0pc@C+0pc{
\underset{r}{\mathrm{homotopycolimit}}~\underset{\mathrm{Spec}}{\mathcal{O}}^\mathrm{CS}\breve{\Phi}^r_{*,A}/\mathrm{Fro}^\mathbb{Z},\underset{I}{\mathrm{homotopylimit}}~\breve{\Phi}^I_{*,A}/\mathrm{Fro}^\mathbb{Z},	\\
}
\]
\[ 
\xymatrix@R+0pc@C+0pc{
\underset{r}{\mathrm{homotopycolimit}}~\underset{\mathrm{Spec}}{\mathcal{O}}^\mathrm{CS}{\Phi}^r_{*,A}/\mathrm{Fro}^\mathbb{Z},\underset{I}{\mathrm{homotopylimit}}~\underset{\mathrm{Spec}}{\mathcal{O}}^\mathrm{CS}{\Phi}^I_{*,A}/\mathrm{Fro}^\mathbb{Z}.	
}
\]	
In this situation we will have the target category being family parametrized by $r$ or $I$ in compatible glueing sense as in \cite[Definition 5.4.10]{10KL2}. In this situation for modules parametrized by the intervals we have the equivalence of $\infty$-categories by using \cite[Proposition 12.18]{10CS2}. Here the corresponding quasicoherent Frobenius modules are defined to be the corresponding homotopy colimits and limits of Frobenius modules:
\begin{align}
\underset{r}{\mathrm{homotopycolimit}}~M_r,\\
\underset{I}{\mathrm{homotopylimit}}~M_I,	
\end{align}
where each $M_r$ is a Frobenius-equivariant module over the period ring with respect to some radius $r$ while each $M_I$ is a Frobenius-equivariant module over the period ring with respect to some interval $I$.\\
\end{proposition}

\begin{proposition}
Similar proposition holds for 
\begin{align}
\mathrm{Quasicoherentsheaves,Perfectcomplex,IndBanach}_{*}.	
\end{align}	
\end{proposition}

\newpage
\subsection{Frobenius Quasicoherent Prestacks II: Deformation in Banach Rings}

\begin{definition}
We now consider the pro-\'etale site of $\mathrm{Spa}\mathbb{Q}_p\left<X_1^{\pm 1},...,X_k^{\pm 1}\right>$, denote that by $*$. To be more accurate we replace one component for $\Gamma$ with the pro-\'etale site of $\mathrm{Spa}\mathbb{Q}_p\left<X_1^{\pm 1},...,X_k^{\pm 1}\right>$. And we treat then all the functor to be prestacks for this site. Then from \cite{10KL1} and \cite[Definition 5.2.1]{10KL2} we have the following class of Kedlaya-Liu rings (with the following replacement: $\Delta$ stands for $A$, $\nabla$ stands for $B$, while $\Phi$ stands for $C$) by taking product in the sense of self $\Gamma$-th power\footnote{Here $|\Gamma|=1$.}:

\[
\xymatrix@R+0pc@C+0pc{
\widetilde{\Delta}_{*},\widetilde{\nabla}_{*},\widetilde{\Phi}_{*},\widetilde{\Delta}^+_{*},\widetilde{\nabla}^+_{*},\widetilde{\Delta}^\dagger_{*},\widetilde{\nabla}^\dagger_{*},\widetilde{\Phi}^r_{*},\widetilde{\Phi}^I_{*}, 
}
\]

\[
\xymatrix@R+0pc@C+0pc{
\breve{\Delta}_{*},\breve{\nabla}_{*},\breve{\Phi}_{*},\breve{\Delta}^+_{*},\breve{\nabla}^+_{*},\breve{\Delta}^\dagger_{*},\breve{\nabla}^\dagger_{*},\breve{\Phi}^r_{*},\breve{\Phi}^I_{*},	
}
\]

\[
\xymatrix@R+0pc@C+0pc{
{\Delta}_{*},{\nabla}_{*},{\Phi}_{*},{\Delta}^+_{*},{\nabla}^+_{*},{\Delta}^\dagger_{*},{\nabla}^\dagger_{*},{\Phi}^r_{*},{\Phi}^I_{*}.	
}
\]
We now consider the following rings with $-$ being any deforming Banach ring over $\mathbb{Q}_p$. Taking the product we have:
\[
\xymatrix@R+0pc@C+0pc{
\widetilde{\Phi}_{*,-},\widetilde{\Phi}^r_{*,-},\widetilde{\Phi}^I_{*,-},	
}
\]
\[
\xymatrix@R+0pc@C+0pc{
\breve{\Phi}_{*,-},\breve{\Phi}^r_{*,-},\breve{\Phi}^I_{*,-},	
}
\]
\[
\xymatrix@R+0pc@C+0pc{
{\Phi}_{*,-},{\Phi}^r_{*,-},{\Phi}^I_{*,-}.	
}
\]
They carry multi Frobenius action $\varphi_\Gamma$ and multi $\mathrm{Lie}_\Gamma:=\mathbb{Z}_p^{\times\Gamma}$ action. In our current situation after \cite{10CKZ} and \cite{10PZ} we consider the following $(\infty,1)$-categories of $(\infty,1)$-modules.\\
\end{definition}

\begin{definition}
First we consider the Bambozzi-Kremnizer spectrum $\underset{\mathrm{Spec}}{\mathcal{O}}^\mathrm{BK}(*)$ attached to any of those in the above from \cite{10BK} by taking derived rational localization:
\begin{align}
&\underset{\mathrm{Spec}}{\mathcal{O}}^\mathrm{BK}\widetilde{\Phi}_{*,-},\underset{\mathrm{Spec}}{\mathcal{O}}^\mathrm{BK}\widetilde{\Phi}^r_{*,-},\underset{\mathrm{Spec}}{\mathcal{O}}^\mathrm{BK}\widetilde{\Phi}^I_{*,-},	
\end{align}
\begin{align}
&\underset{\mathrm{Spec}}{\mathcal{O}}^\mathrm{BK}\breve{\Phi}_{*,-},\underset{\mathrm{Spec}}{\mathcal{O}}^\mathrm{BK}\breve{\Phi}^r_{*,-},\underset{\mathrm{Spec}}{\mathcal{O}}^\mathrm{BK}\breve{\Phi}^I_{*,-},	
\end{align}
\begin{align}
&\underset{\mathrm{Spec}}{\mathcal{O}}^\mathrm{BK}{\Phi}_{*,-},
\underset{\mathrm{Spec}}{\mathcal{O}}^\mathrm{BK}{\Phi}^r_{*,-},\underset{\mathrm{Spec}}{\mathcal{O}}^\mathrm{BK}{\Phi}^I_{*,-}.	
\end{align}

Then we take the corresponding quotients by using the corresponding Frobenius operators:
\begin{align}
&\underset{\mathrm{Spec}}{\mathcal{O}}^\mathrm{BK}\widetilde{\Phi}_{*,-}/\mathrm{Fro}^\mathbb{Z},	\\
\end{align}
\begin{align}
&\underset{\mathrm{Spec}}{\mathcal{O}}^\mathrm{BK}\breve{\Phi}_{*,-}/\mathrm{Fro}^\mathbb{Z},	\\
\end{align}
\begin{align}
&\underset{\mathrm{Spec}}{\mathcal{O}}^\mathrm{BK}{\Phi}_{*,-}/\mathrm{Fro}^\mathbb{Z}.	
\end{align}
Here for those space without notation related to the radius and the corresponding interval we consider the total unions $\bigcap_r,\bigcup_I$ in order to achieve the whole spaces to achieve the analogues of the corresponding FF curves from \cite{10KL1}, \cite{10KL2}, \cite{10FF} for
\[
\xymatrix@R+0pc@C+0pc{
\underset{r}{\mathrm{homotopycolimit}}~\underset{\mathrm{Spec}}{\mathcal{O}}^\mathrm{BK}\widetilde{\Phi}^r_{*,-},\underset{I}{\mathrm{homotopylimit}}~\underset{\mathrm{Spec}}{\mathcal{O}}^\mathrm{BK}\widetilde{\Phi}^I_{*,-},	\\
}
\]
\[
\xymatrix@R+0pc@C+0pc{
\underset{r}{\mathrm{homotopycolimit}}~\underset{\mathrm{Spec}}{\mathcal{O}}^\mathrm{BK}\breve{\Phi}^r_{*,-},\underset{I}{\mathrm{homotopylimit}}~\underset{\mathrm{Spec}}{\mathcal{O}}^\mathrm{BK}\breve{\Phi}^I_{*,-},	\\
}
\]
\[
\xymatrix@R+0pc@C+0pc{
\underset{r}{\mathrm{homotopycolimit}}~\underset{\mathrm{Spec}}{\mathcal{O}}^\mathrm{BK}{\Phi}^r_{*,-},\underset{I}{\mathrm{homotopylimit}}~\underset{\mathrm{Spec}}{\mathcal{O}}^\mathrm{BK}{\Phi}^I_{*,-}.	
}
\]
\[  
\xymatrix@R+0pc@C+0pc{
\underset{r}{\mathrm{homotopycolimit}}~\underset{\mathrm{Spec}}{\mathcal{O}}^\mathrm{BK}\widetilde{\Phi}^r_{*,-}/\mathrm{Fro}^\mathbb{Z},\underset{I}{\mathrm{homotopylimit}}~\underset{\mathrm{Spec}}{\mathcal{O}}^\mathrm{BK}\widetilde{\Phi}^I_{*,-}/\mathrm{Fro}^\mathbb{Z},	\\
}
\]
\[ 
\xymatrix@R+0pc@C+0pc{
\underset{r}{\mathrm{homotopycolimit}}~\underset{\mathrm{Spec}}{\mathcal{O}}^\mathrm{BK}\breve{\Phi}^r_{*,-}/\mathrm{Fro}^\mathbb{Z},\underset{I}{\mathrm{homotopylimit}}~\underset{\mathrm{Spec}}{\mathcal{O}}^\mathrm{BK}\breve{\Phi}^I_{*,-}/\mathrm{Fro}^\mathbb{Z},	\\
}
\]
\[ 
\xymatrix@R+0pc@C+0pc{
\underset{r}{\mathrm{homotopycolimit}}~\underset{\mathrm{Spec}}{\mathcal{O}}^\mathrm{BK}{\Phi}^r_{*,-}/\mathrm{Fro}^\mathbb{Z},\underset{I}{\mathrm{homotopylimit}}~\underset{\mathrm{Spec}}{\mathcal{O}}^\mathrm{BK}{\Phi}^I_{*,-}/\mathrm{Fro}^\mathbb{Z}.	
}
\]

\end{definition}

\indent Meanwhile we have the corresponding Clausen-Scholze analytic stacks from \cite{10CS2}, therefore applying their construction we have:

\begin{definition}
Here we define the following products by using the solidified tensor product from \cite{10CS1} and \cite{10CS2}. Namely $A$ will still as above as a Banach ring over $\mathbb{Q}_p$. Then we take solidified tensor product $\overset{\blacksquare}{\otimes}$ of any of the following
\[
\xymatrix@R+0pc@C+0pc{
\widetilde{\Delta}_{*},\widetilde{\nabla}_{*},\widetilde{\Phi}_{*},\widetilde{\Delta}^+_{*},\widetilde{\nabla}^+_{*},\widetilde{\Delta}^\dagger_{*},\widetilde{\nabla}^\dagger_{*},\widetilde{\Phi}^r_{*},\widetilde{\Phi}^I_{*}, 
}
\]

\[
\xymatrix@R+0pc@C+0pc{
\breve{\Delta}_{*},\breve{\nabla}_{*},\breve{\Phi}_{*},\breve{\Delta}^+_{*},\breve{\nabla}^+_{*},\breve{\Delta}^\dagger_{*},\breve{\nabla}^\dagger_{*},\breve{\Phi}^r_{*},\breve{\Phi}^I_{*},	
}
\]

\[
\xymatrix@R+0pc@C+0pc{
{\Delta}_{*},{\nabla}_{*},{\Phi}_{*},{\Delta}^+_{*},{\nabla}^+_{*},{\Delta}^\dagger_{*},{\nabla}^\dagger_{*},{\Phi}^r_{*},{\Phi}^I_{*},	
}
\]  	
with $A$. Then we have the notations:
\[
\xymatrix@R+0pc@C+0pc{
\widetilde{\Delta}_{*,-},\widetilde{\nabla}_{*,-},\widetilde{\Phi}_{*,-},\widetilde{\Delta}^+_{*,-},\widetilde{\nabla}^+_{*,-},\widetilde{\Delta}^\dagger_{*,-},\widetilde{\nabla}^\dagger_{*,-},\widetilde{\Phi}^r_{*,-},\widetilde{\Phi}^I_{*,-}, 
}
\]

\[
\xymatrix@R+0pc@C+0pc{
\breve{\Delta}_{*,-},\breve{\nabla}_{*,-},\breve{\Phi}_{*,-},\breve{\Delta}^+_{*,-},\breve{\nabla}^+_{*,-},\breve{\Delta}^\dagger_{*,-},\breve{\nabla}^\dagger_{*,-},\breve{\Phi}^r_{*,-},\breve{\Phi}^I_{*,-},	
}
\]

\[
\xymatrix@R+0pc@C+0pc{
{\Delta}_{*,-},{\nabla}_{*,-},{\Phi}_{*,-},{\Delta}^+_{*,-},{\nabla}^+_{*,-},{\Delta}^\dagger_{*,-},{\nabla}^\dagger_{*,-},{\Phi}^r_{*,-},{\Phi}^I_{*,-}.	
}
\]
\end{definition}

\begin{definition}
First we consider the Clausen-Scholze spectrum $\underset{\mathrm{Spec}}{\mathcal{O}}^\mathrm{CS}(*)$ attached to any of those in the above from \cite{10CS2} by taking derived rational localization:
\begin{align}
\underset{\mathrm{Spec}}{\mathcal{O}}^\mathrm{CS}\widetilde{\Delta}_{*,-},\underset{\mathrm{Spec}}{\mathcal{O}}^\mathrm{CS}\widetilde{\nabla}_{*,-},\underset{\mathrm{Spec}}{\mathcal{O}}^\mathrm{CS}\widetilde{\Phi}_{*,-},\underset{\mathrm{Spec}}{\mathcal{O}}^\mathrm{CS}\widetilde{\Delta}^+_{*,-},\underset{\mathrm{Spec}}{\mathcal{O}}^\mathrm{CS}\widetilde{\nabla}^+_{*,-},\\
\underset{\mathrm{Spec}}{\mathcal{O}}^\mathrm{CS}\widetilde{\Delta}^\dagger_{*,-},\underset{\mathrm{Spec}}{\mathcal{O}}^\mathrm{CS}\widetilde{\nabla}^\dagger_{*,-},\underset{\mathrm{Spec}}{\mathcal{O}}^\mathrm{CS}\widetilde{\Phi}^r_{*,-},\underset{\mathrm{Spec}}{\mathcal{O}}^\mathrm{CS}\widetilde{\Phi}^I_{*,-},	\\
\end{align}
\begin{align}
\underset{\mathrm{Spec}}{\mathcal{O}}^\mathrm{CS}\breve{\Delta}_{*,-},\breve{\nabla}_{*,-},\underset{\mathrm{Spec}}{\mathcal{O}}^\mathrm{CS}\breve{\Phi}_{*,-},\underset{\mathrm{Spec}}{\mathcal{O}}^\mathrm{CS}\breve{\Delta}^+_{*,-},\underset{\mathrm{Spec}}{\mathcal{O}}^\mathrm{CS}\breve{\nabla}^+_{*,-},\\
\underset{\mathrm{Spec}}{\mathcal{O}}^\mathrm{CS}\breve{\Delta}^\dagger_{*,-},\underset{\mathrm{Spec}}{\mathcal{O}}^\mathrm{CS}\breve{\nabla}^\dagger_{*,-},\underset{\mathrm{Spec}}{\mathcal{O}}^\mathrm{CS}\breve{\Phi}^r_{*,-},\breve{\Phi}^I_{*,-},	\\
\end{align}
\begin{align}
\underset{\mathrm{Spec}}{\mathcal{O}}^\mathrm{CS}{\Delta}_{*,-},\underset{\mathrm{Spec}}{\mathcal{O}}^\mathrm{CS}{\nabla}_{*,-},\underset{\mathrm{Spec}}{\mathcal{O}}^\mathrm{CS}{\Phi}_{*,-},\underset{\mathrm{Spec}}{\mathcal{O}}^\mathrm{CS}{\Delta}^+_{*,-},\underset{\mathrm{Spec}}{\mathcal{O}}^\mathrm{CS}{\nabla}^+_{*,-},\\
\underset{\mathrm{Spec}}{\mathcal{O}}^\mathrm{CS}{\Delta}^\dagger_{*,-},\underset{\mathrm{Spec}}{\mathcal{O}}^\mathrm{CS}{\nabla}^\dagger_{*,-},\underset{\mathrm{Spec}}{\mathcal{O}}^\mathrm{CS}{\Phi}^r_{*,-},\underset{\mathrm{Spec}}{\mathcal{O}}^\mathrm{CS}{\Phi}^I_{*,-}.	
\end{align}

Then we take the corresponding quotients by using the corresponding Frobenius operators:
\begin{align}
&\underset{\mathrm{Spec}}{\mathcal{O}}^\mathrm{CS}\widetilde{\Delta}_{*,-}/\mathrm{Fro}^\mathbb{Z},\underset{\mathrm{Spec}}{\mathcal{O}}^\mathrm{CS}\widetilde{\nabla}_{*,-}/\mathrm{Fro}^\mathbb{Z},\underset{\mathrm{Spec}}{\mathcal{O}}^\mathrm{CS}\widetilde{\Phi}_{*,-}/\mathrm{Fro}^\mathbb{Z},\underset{\mathrm{Spec}}{\mathcal{O}}^\mathrm{CS}\widetilde{\Delta}^+_{*,-}/\mathrm{Fro}^\mathbb{Z},\\
&\underset{\mathrm{Spec}}{\mathcal{O}}^\mathrm{CS}\widetilde{\nabla}^+_{*,-}/\mathrm{Fro}^\mathbb{Z}, \underset{\mathrm{Spec}}{\mathcal{O}}^\mathrm{CS}\widetilde{\Delta}^\dagger_{*,-}/\mathrm{Fro}^\mathbb{Z},\underset{\mathrm{Spec}}{\mathcal{O}}^\mathrm{CS}\widetilde{\nabla}^\dagger_{*,-}/\mathrm{Fro}^\mathbb{Z},	\\
\end{align}
\begin{align}
&\underset{\mathrm{Spec}}{\mathcal{O}}^\mathrm{CS}\breve{\Delta}_{*,-}/\mathrm{Fro}^\mathbb{Z},\breve{\nabla}_{*,-}/\mathrm{Fro}^\mathbb{Z},\underset{\mathrm{Spec}}{\mathcal{O}}^\mathrm{CS}\breve{\Phi}_{*,-}/\mathrm{Fro}^\mathbb{Z},\underset{\mathrm{Spec}}{\mathcal{O}}^\mathrm{CS}\breve{\Delta}^+_{*,-}/\mathrm{Fro}^\mathbb{Z},\\
&\underset{\mathrm{Spec}}{\mathcal{O}}^\mathrm{CS}\breve{\nabla}^+_{*,-}/\mathrm{Fro}^\mathbb{Z}, \underset{\mathrm{Spec}}{\mathcal{O}}^\mathrm{CS}\breve{\Delta}^\dagger_{*,-}/\mathrm{Fro}^\mathbb{Z},\underset{\mathrm{Spec}}{\mathcal{O}}^\mathrm{CS}\breve{\nabla}^\dagger_{*,-}/\mathrm{Fro}^\mathbb{Z},	\\
\end{align}
\begin{align}
&\underset{\mathrm{Spec}}{\mathcal{O}}^\mathrm{CS}{\Delta}_{*,-}/\mathrm{Fro}^\mathbb{Z},\underset{\mathrm{Spec}}{\mathcal{O}}^\mathrm{CS}{\nabla}_{*,-}/\mathrm{Fro}^\mathbb{Z},\underset{\mathrm{Spec}}{\mathcal{O}}^\mathrm{CS}{\Phi}_{*,-}/\mathrm{Fro}^\mathbb{Z},\underset{\mathrm{Spec}}{\mathcal{O}}^\mathrm{CS}{\Delta}^+_{*,-}/\mathrm{Fro}^\mathbb{Z},\\
&\underset{\mathrm{Spec}}{\mathcal{O}}^\mathrm{CS}{\nabla}^+_{*,-}/\mathrm{Fro}^\mathbb{Z}, \underset{\mathrm{Spec}}{\mathcal{O}}^\mathrm{CS}{\Delta}^\dagger_{*,-}/\mathrm{Fro}^\mathbb{Z},\underset{\mathrm{Spec}}{\mathcal{O}}^\mathrm{CS}{\nabla}^\dagger_{*,-}/\mathrm{Fro}^\mathbb{Z}.	
\end{align}
Here for those space with notations related to the radius and the corresponding interval we consider the total unions $\bigcap_r,\bigcup_I$ in order to achieve the whole spaces to achieve the analogues of the corresponding FF curves from \cite{10KL1}, \cite{10KL2}, \cite{10FF} for
\[
\xymatrix@R+0pc@C+0pc{
\underset{r}{\mathrm{homotopycolimit}}~\underset{\mathrm{Spec}}{\mathcal{O}}^\mathrm{CS}\widetilde{\Phi}^r_{*,-},\underset{I}{\mathrm{homotopylimit}}~\underset{\mathrm{Spec}}{\mathcal{O}}^\mathrm{CS}\widetilde{\Phi}^I_{*,-},	\\
}
\]
\[
\xymatrix@R+0pc@C+0pc{
\underset{r}{\mathrm{homotopycolimit}}~\underset{\mathrm{Spec}}{\mathcal{O}}^\mathrm{CS}\breve{\Phi}^r_{*,-},\underset{I}{\mathrm{homotopylimit}}~\underset{\mathrm{Spec}}{\mathcal{O}}^\mathrm{CS}\breve{\Phi}^I_{*,-},	\\
}
\]
\[
\xymatrix@R+0pc@C+0pc{
\underset{r}{\mathrm{homotopycolimit}}~\underset{\mathrm{Spec}}{\mathcal{O}}^\mathrm{CS}{\Phi}^r_{*,-},\underset{I}{\mathrm{homotopylimit}}~\underset{\mathrm{Spec}}{\mathcal{O}}^\mathrm{CS}{\Phi}^I_{*,-}.	
}
\]
\[ 
\xymatrix@R+0pc@C+0pc{
\underset{r}{\mathrm{homotopycolimit}}~\underset{\mathrm{Spec}}{\mathcal{O}}^\mathrm{CS}\widetilde{\Phi}^r_{*,-}/\mathrm{Fro}^\mathbb{Z},\underset{I}{\mathrm{homotopylimit}}~\underset{\mathrm{Spec}}{\mathcal{O}}^\mathrm{CS}\widetilde{\Phi}^I_{*,-}/\mathrm{Fro}^\mathbb{Z},	\\
}
\]
\[ 
\xymatrix@R+0pc@C+0pc{
\underset{r}{\mathrm{homotopycolimit}}~\underset{\mathrm{Spec}}{\mathcal{O}}^\mathrm{CS}\breve{\Phi}^r_{*,-}/\mathrm{Fro}^\mathbb{Z},\underset{I}{\mathrm{homotopylimit}}~\breve{\Phi}^I_{*,-}/\mathrm{Fro}^\mathbb{Z},	\\
}
\]
\[ 
\xymatrix@R+0pc@C+0pc{
\underset{r}{\mathrm{homotopycolimit}}~\underset{\mathrm{Spec}}{\mathcal{O}}^\mathrm{CS}{\Phi}^r_{*,-}/\mathrm{Fro}^\mathbb{Z},\underset{I}{\mathrm{homotopylimit}}~\underset{\mathrm{Spec}}{\mathcal{O}}^\mathrm{CS}{\Phi}^I_{*,-}/\mathrm{Fro}^\mathbb{Z}.	
}
\]

\end{definition}

\

\begin{definition}
We then consider the corresponding quasipresheaves of the corresponding ind-Banach or monomorphic ind-Banach modules from \cite{10BBK}, \cite{10KKM}:
\begin{align}
\mathrm{Quasicoherentpresheaves,IndBanach}_{*}	
\end{align}
where $*$ is one of the following spaces:
\begin{align}
&\underset{\mathrm{Spec}}{\mathcal{O}}^\mathrm{BK}\widetilde{\Phi}_{*,-}/\mathrm{Fro}^\mathbb{Z},	\\
\end{align}
\begin{align}
&\underset{\mathrm{Spec}}{\mathcal{O}}^\mathrm{BK}\breve{\Phi}_{*,-}/\mathrm{Fro}^\mathbb{Z},	\\
\end{align}
\begin{align}
&\underset{\mathrm{Spec}}{\mathcal{O}}^\mathrm{BK}{\Phi}_{*,-}/\mathrm{Fro}^\mathbb{Z}.	
\end{align}
Here for those space without notation related to the radius and the corresponding interval we consider the total unions $\bigcap_r,\bigcup_I$ in order to achieve the whole spaces to achieve the analogues of the corresponding FF curves from \cite{10KL1}, \cite{10KL2}, \cite{10FF} for
\[
\xymatrix@R+0pc@C+0pc{
\underset{r}{\mathrm{homotopycolimit}}~\underset{\mathrm{Spec}}{\mathcal{O}}^\mathrm{BK}\widetilde{\Phi}^r_{*,-},\underset{I}{\mathrm{homotopylimit}}~\underset{\mathrm{Spec}}{\mathcal{O}}^\mathrm{BK}\widetilde{\Phi}^I_{*,-},	\\
}
\]
\[
\xymatrix@R+0pc@C+0pc{
\underset{r}{\mathrm{homotopycolimit}}~\underset{\mathrm{Spec}}{\mathcal{O}}^\mathrm{BK}\breve{\Phi}^r_{*,-},\underset{I}{\mathrm{homotopylimit}}~\underset{\mathrm{Spec}}{\mathcal{O}}^\mathrm{BK}\breve{\Phi}^I_{*,-},	\\
}
\]
\[
\xymatrix@R+0pc@C+0pc{
\underset{r}{\mathrm{homotopycolimit}}~\underset{\mathrm{Spec}}{\mathcal{O}}^\mathrm{BK}{\Phi}^r_{*,-},\underset{I}{\mathrm{homotopylimit}}~\underset{\mathrm{Spec}}{\mathcal{O}}^\mathrm{BK}{\Phi}^I_{*,-}.	
}
\]
\[  
\xymatrix@R+0pc@C+0pc{
\underset{r}{\mathrm{homotopycolimit}}~\underset{\mathrm{Spec}}{\mathcal{O}}^\mathrm{BK}\widetilde{\Phi}^r_{*,-}/\mathrm{Fro}^\mathbb{Z},\underset{I}{\mathrm{homotopylimit}}~\underset{\mathrm{Spec}}{\mathcal{O}}^\mathrm{BK}\widetilde{\Phi}^I_{*,-}/\mathrm{Fro}^\mathbb{Z},	\\
}
\]
\[ 
\xymatrix@R+0pc@C+0pc{
\underset{r}{\mathrm{homotopycolimit}}~\underset{\mathrm{Spec}}{\mathcal{O}}^\mathrm{BK}\breve{\Phi}^r_{*,-}/\mathrm{Fro}^\mathbb{Z},\underset{I}{\mathrm{homotopylimit}}~\underset{\mathrm{Spec}}{\mathcal{O}}^\mathrm{BK}\breve{\Phi}^I_{*,-}/\mathrm{Fro}^\mathbb{Z},	\\
}
\]
\[ 
\xymatrix@R+0pc@C+0pc{
\underset{r}{\mathrm{homotopycolimit}}~\underset{\mathrm{Spec}}{\mathcal{O}}^\mathrm{BK}{\Phi}^r_{*,-}/\mathrm{Fro}^\mathbb{Z},\underset{I}{\mathrm{homotopylimit}}~\underset{\mathrm{Spec}}{\mathcal{O}}^\mathrm{BK}{\Phi}^I_{*,-}/\mathrm{Fro}^\mathbb{Z}.	
}
\]

\end{definition}

\begin{definition}
We then consider the corresponding quasisheaves of the corresponding condensed solid topological modules from \cite{10CS2}:
\begin{align}
\mathrm{Quasicoherentsheaves, Condensed}_{*}	
\end{align}
where $*$ is one of the following spaces:
\begin{align}
&\underset{\mathrm{Spec}}{\mathcal{O}}^\mathrm{CS}\widetilde{\Delta}_{*,-}/\mathrm{Fro}^\mathbb{Z},\underset{\mathrm{Spec}}{\mathcal{O}}^\mathrm{CS}\widetilde{\nabla}_{*,-}/\mathrm{Fro}^\mathbb{Z},\underset{\mathrm{Spec}}{\mathcal{O}}^\mathrm{CS}\widetilde{\Phi}_{*,-}/\mathrm{Fro}^\mathbb{Z},\underset{\mathrm{Spec}}{\mathcal{O}}^\mathrm{CS}\widetilde{\Delta}^+_{*,-}/\mathrm{Fro}^\mathbb{Z},\\
&\underset{\mathrm{Spec}}{\mathcal{O}}^\mathrm{CS}\widetilde{\nabla}^+_{*,-}/\mathrm{Fro}^\mathbb{Z},\underset{\mathrm{Spec}}{\mathcal{O}}^\mathrm{CS}\widetilde{\Delta}^\dagger_{*,-}/\mathrm{Fro}^\mathbb{Z},\underset{\mathrm{Spec}}{\mathcal{O}}^\mathrm{CS}\widetilde{\nabla}^\dagger_{*,-}/\mathrm{Fro}^\mathbb{Z},	\\
\end{align}
\begin{align}
&\underset{\mathrm{Spec}}{\mathcal{O}}^\mathrm{CS}\breve{\Delta}_{*,-}/\mathrm{Fro}^\mathbb{Z},\breve{\nabla}_{*,-}/\mathrm{Fro}^\mathbb{Z},\underset{\mathrm{Spec}}{\mathcal{O}}^\mathrm{CS}\breve{\Phi}_{*,-}/\mathrm{Fro}^\mathbb{Z},\underset{\mathrm{Spec}}{\mathcal{O}}^\mathrm{CS}\breve{\Delta}^+_{*,-}/\mathrm{Fro}^\mathbb{Z},\\
&\underset{\mathrm{Spec}}{\mathcal{O}}^\mathrm{CS}\breve{\nabla}^+_{*,-}/\mathrm{Fro}^\mathbb{Z},\underset{\mathrm{Spec}}{\mathcal{O}}^\mathrm{CS}\breve{\Delta}^\dagger_{*,-}/\mathrm{Fro}^\mathbb{Z},\underset{\mathrm{Spec}}{\mathcal{O}}^\mathrm{CS}\breve{\nabla}^\dagger_{*,-}/\mathrm{Fro}^\mathbb{Z},	\\
\end{align}
\begin{align}
&\underset{\mathrm{Spec}}{\mathcal{O}}^\mathrm{CS}{\Delta}_{*,-}/\mathrm{Fro}^\mathbb{Z},\underset{\mathrm{Spec}}{\mathcal{O}}^\mathrm{CS}{\nabla}_{*,-}/\mathrm{Fro}^\mathbb{Z},\underset{\mathrm{Spec}}{\mathcal{O}}^\mathrm{CS}{\Phi}_{*,-}/\mathrm{Fro}^\mathbb{Z},\underset{\mathrm{Spec}}{\mathcal{O}}^\mathrm{CS}{\Delta}^+_{*,-}/\mathrm{Fro}^\mathbb{Z},\\
&\underset{\mathrm{Spec}}{\mathcal{O}}^\mathrm{CS}{\nabla}^+_{*,-}/\mathrm{Fro}^\mathbb{Z}, \underset{\mathrm{Spec}}{\mathcal{O}}^\mathrm{CS}{\Delta}^\dagger_{*,-}/\mathrm{Fro}^\mathbb{Z},\underset{\mathrm{Spec}}{\mathcal{O}}^\mathrm{CS}{\nabla}^\dagger_{*,-}/\mathrm{Fro}^\mathbb{Z}.	
\end{align}
Here for those space with notations related to the radius and the corresponding interval we consider the total unions $\bigcap_r,\bigcup_I$ in order to achieve the whole spaces to achieve the analogues of the corresponding FF curves from \cite{10KL1}, \cite{10KL2}, \cite{10FF} for
\[
\xymatrix@R+0pc@C+0pc{
\underset{r}{\mathrm{homotopycolimit}}~\underset{\mathrm{Spec}}{\mathcal{O}}^\mathrm{CS}\widetilde{\Phi}^r_{*,-},\underset{I}{\mathrm{homotopylimit}}~\underset{\mathrm{Spec}}{\mathcal{O}}^\mathrm{CS}\widetilde{\Phi}^I_{*,-},	\\
}
\]
\[
\xymatrix@R+0pc@C+0pc{
\underset{r}{\mathrm{homotopycolimit}}~\underset{\mathrm{Spec}}{\mathcal{O}}^\mathrm{CS}\breve{\Phi}^r_{*,-},\underset{I}{\mathrm{homotopylimit}}~\underset{\mathrm{Spec}}{\mathcal{O}}^\mathrm{CS}\breve{\Phi}^I_{*,-},	\\
}
\]
\[
\xymatrix@R+0pc@C+0pc{
\underset{r}{\mathrm{homotopycolimit}}~\underset{\mathrm{Spec}}{\mathcal{O}}^\mathrm{CS}{\Phi}^r_{*,-},\underset{I}{\mathrm{homotopylimit}}~\underset{\mathrm{Spec}}{\mathcal{O}}^\mathrm{CS}{\Phi}^I_{*,-}.	
}
\]
\[ 
\xymatrix@R+0pc@C+0pc{
\underset{r}{\mathrm{homotopycolimit}}~\underset{\mathrm{Spec}}{\mathcal{O}}^\mathrm{CS}\widetilde{\Phi}^r_{*,-}/\mathrm{Fro}^\mathbb{Z},\underset{I}{\mathrm{homotopylimit}}~\underset{\mathrm{Spec}}{\mathcal{O}}^\mathrm{CS}\widetilde{\Phi}^I_{*,-}/\mathrm{Fro}^\mathbb{Z},	\\
}
\]
\[ 
\xymatrix@R+0pc@C+0pc{
\underset{r}{\mathrm{homotopycolimit}}~\underset{\mathrm{Spec}}{\mathcal{O}}^\mathrm{CS}\breve{\Phi}^r_{*,-}/\mathrm{Fro}^\mathbb{Z},\underset{I}{\mathrm{homotopylimit}}~\breve{\Phi}^I_{*,-}/\mathrm{Fro}^\mathbb{Z},	\\
}
\]
\[ 
\xymatrix@R+0pc@C+0pc{
\underset{r}{\mathrm{homotopycolimit}}~\underset{\mathrm{Spec}}{\mathcal{O}}^\mathrm{CS}{\Phi}^r_{*,-}/\mathrm{Fro}^\mathbb{Z},\underset{I}{\mathrm{homotopylimit}}~\underset{\mathrm{Spec}}{\mathcal{O}}^\mathrm{CS}{\Phi}^I_{*,-}/\mathrm{Fro}^\mathbb{Z}.	
}
\]

\end{definition}

\

\begin{proposition}
There is a well-defined functor from the $\infty$-category 
\begin{align}
\mathrm{Quasicoherentpresheaves,Condensed}_{*}	
\end{align}
where $*$ is one of the following spaces:
\begin{align}
&\underset{\mathrm{Spec}}{\mathcal{O}}^\mathrm{CS}\widetilde{\Phi}_{*,-}/\mathrm{Fro}^\mathbb{Z},	\\
\end{align}
\begin{align}
&\underset{\mathrm{Spec}}{\mathcal{O}}^\mathrm{CS}\breve{\Phi}_{*,-}/\mathrm{Fro}^\mathbb{Z},	\\
\end{align}
\begin{align}
&\underset{\mathrm{Spec}}{\mathcal{O}}^\mathrm{CS}{\Phi}_{*,-}/\mathrm{Fro}^\mathbb{Z},	
\end{align}
to the $\infty$-category of $\mathrm{Fro}$-equivariant quasicoherent presheaves over similar spaces above correspondingly without the $\mathrm{Fro}$-quotients, and to the $\infty$-category of $\mathrm{Fro}$-equivariant quasicoherent modules over global sections of the structure $\infty$-sheaves of the similar spaces above correspondingly without the $\mathrm{Fro}$-quotients. Here for those space without notation related to the radius and the corresponding interval we consider the total unions $\bigcap_r,\bigcup_I$ in order to achieve the whole spaces to achieve the analogues of the corresponding FF curves from \cite{10KL1}, \cite{10KL2}, \cite{10FF} for
\[
\xymatrix@R+0pc@C+0pc{
\underset{r}{\mathrm{homotopycolimit}}~\underset{\mathrm{Spec}}{\mathcal{O}}^\mathrm{CS}\widetilde{\Phi}^r_{*,-},\underset{I}{\mathrm{homotopylimit}}~\underset{\mathrm{Spec}}{\mathcal{O}}^\mathrm{CS}\widetilde{\Phi}^I_{*,-},	\\
}
\]
\[
\xymatrix@R+0pc@C+0pc{
\underset{r}{\mathrm{homotopycolimit}}~\underset{\mathrm{Spec}}{\mathcal{O}}^\mathrm{CS}\breve{\Phi}^r_{*,-},\underset{I}{\mathrm{homotopylimit}}~\underset{\mathrm{Spec}}{\mathcal{O}}^\mathrm{CS}\breve{\Phi}^I_{*,-},	\\
}
\]
\[
\xymatrix@R+0pc@C+0pc{
\underset{r}{\mathrm{homotopycolimit}}~\underset{\mathrm{Spec}}{\mathcal{O}}^\mathrm{CS}{\Phi}^r_{*,-},\underset{I}{\mathrm{homotopylimit}}~\underset{\mathrm{Spec}}{\mathcal{O}}^\mathrm{CS}{\Phi}^I_{*,-}.	
}
\]
\[ 
\xymatrix@R+0pc@C+0pc{
\underset{r}{\mathrm{homotopycolimit}}~\underset{\mathrm{Spec}}{\mathcal{O}}^\mathrm{CS}\widetilde{\Phi}^r_{*,-}/\mathrm{Fro}^\mathbb{Z},\underset{I}{\mathrm{homotopylimit}}~\underset{\mathrm{Spec}}{\mathcal{O}}^\mathrm{CS}\widetilde{\Phi}^I_{*,-}/\mathrm{Fro}^\mathbb{Z},	\\
}
\]
\[ 
\xymatrix@R+0pc@C+0pc{
\underset{r}{\mathrm{homotopycolimit}}~\underset{\mathrm{Spec}}{\mathcal{O}}^\mathrm{CS}\breve{\Phi}^r_{*,-}/\mathrm{Fro}^\mathbb{Z},\underset{I}{\mathrm{homotopylimit}}~\breve{\Phi}^I_{*,-}/\mathrm{Fro}^\mathbb{Z},	\\
}
\]
\[ 
\xymatrix@R+0pc@C+0pc{
\underset{r}{\mathrm{homotopycolimit}}~\underset{\mathrm{Spec}}{\mathcal{O}}^\mathrm{CS}{\Phi}^r_{*,-}/\mathrm{Fro}^\mathbb{Z},\underset{I}{\mathrm{homotopylimit}}~\underset{\mathrm{Spec}}{\mathcal{O}}^\mathrm{CS}{\Phi}^I_{*,-}/\mathrm{Fro}^\mathbb{Z}.	
}
\]	
In this situation we will have the target category being family parametrized by $r$ or $I$ in compatible glueing sense as in \cite[Definition 5.4.10]{10KL2}. In this situation for modules parametrized by the intervals we have the equivalence of $\infty$-categories by using \cite[Proposition 13.8]{10CS2}. Here the corresponding quasicoherent Frobenius modules are defined to be the corresponding homotopy colimits and limits of Frobenius modules:
\begin{align}
\underset{r}{\mathrm{homotopycolimit}}~M_r,\\
\underset{I}{\mathrm{homotopylimit}}~M_I,	
\end{align}
where each $M_r$ is a Frobenius-equivariant module over the period ring with respect to some radius $r$ while each $M_I$ is a Frobenius-equivariant module over the period ring with respect to some interval $I$.\\
\end{proposition}

\begin{proposition}
Similar proposition holds for 
\begin{align}
\mathrm{Quasicoherentsheaves,IndBanach}_{*}.	
\end{align}	
\end{proposition}

\

\begin{definition}
We then consider the corresponding quasipresheaves of perfect complexes the corresponding ind-Banach or monomorphic ind-Banach modules from \cite{10BBK}, \cite{10KKM}:
\begin{align}
\mathrm{Quasicoherentpresheaves,Perfectcomplex,IndBanach}_{*}	
\end{align}
where $*$ is one of the following spaces:
\begin{align}
&\underset{\mathrm{Spec}}{\mathcal{O}}^\mathrm{BK}\widetilde{\Phi}_{*,-}/\mathrm{Fro}^\mathbb{Z},	\\
\end{align}
\begin{align}
&\underset{\mathrm{Spec}}{\mathcal{O}}^\mathrm{BK}\breve{\Phi}_{*,-}/\mathrm{Fro}^\mathbb{Z},	\\
\end{align}
\begin{align}
&\underset{\mathrm{Spec}}{\mathcal{O}}^\mathrm{BK}{\Phi}_{*,-}/\mathrm{Fro}^\mathbb{Z}.	
\end{align}
Here for those space without notation related to the radius and the corresponding interval we consider the total unions $\bigcap_r,\bigcup_I$ in order to achieve the whole spaces to achieve the analogues of the corresponding FF curves from \cite{10KL1}, \cite{10KL2}, \cite{10FF} for
\[
\xymatrix@R+0pc@C+0pc{
\underset{r}{\mathrm{homotopycolimit}}~\underset{\mathrm{Spec}}{\mathcal{O}}^\mathrm{BK}\widetilde{\Phi}^r_{*,-},\underset{I}{\mathrm{homotopylimit}}~\underset{\mathrm{Spec}}{\mathcal{O}}^\mathrm{BK}\widetilde{\Phi}^I_{*,-},	\\
}
\]
\[
\xymatrix@R+0pc@C+0pc{
\underset{r}{\mathrm{homotopycolimit}}~\underset{\mathrm{Spec}}{\mathcal{O}}^\mathrm{BK}\breve{\Phi}^r_{*,-},\underset{I}{\mathrm{homotopylimit}}~\underset{\mathrm{Spec}}{\mathcal{O}}^\mathrm{BK}\breve{\Phi}^I_{*,-},	\\
}
\]
\[
\xymatrix@R+0pc@C+0pc{
\underset{r}{\mathrm{homotopycolimit}}~\underset{\mathrm{Spec}}{\mathcal{O}}^\mathrm{BK}{\Phi}^r_{*,-},\underset{I}{\mathrm{homotopylimit}}~\underset{\mathrm{Spec}}{\mathcal{O}}^\mathrm{BK}{\Phi}^I_{*,-}.	
}
\]
\[  
\xymatrix@R+0pc@C+0pc{
\underset{r}{\mathrm{homotopycolimit}}~\underset{\mathrm{Spec}}{\mathcal{O}}^\mathrm{BK}\widetilde{\Phi}^r_{*,-}/\mathrm{Fro}^\mathbb{Z},\underset{I}{\mathrm{homotopylimit}}~\underset{\mathrm{Spec}}{\mathcal{O}}^\mathrm{BK}\widetilde{\Phi}^I_{*,-}/\mathrm{Fro}^\mathbb{Z},	\\
}
\]
\[ 
\xymatrix@R+0pc@C+0pc{
\underset{r}{\mathrm{homotopycolimit}}~\underset{\mathrm{Spec}}{\mathcal{O}}^\mathrm{BK}\breve{\Phi}^r_{*,-}/\mathrm{Fro}^\mathbb{Z},\underset{I}{\mathrm{homotopylimit}}~\underset{\mathrm{Spec}}{\mathcal{O}}^\mathrm{BK}\breve{\Phi}^I_{*,-}/\mathrm{Fro}^\mathbb{Z},	\\
}
\]
\[ 
\xymatrix@R+0pc@C+0pc{
\underset{r}{\mathrm{homotopycolimit}}~\underset{\mathrm{Spec}}{\mathcal{O}}^\mathrm{BK}{\Phi}^r_{*,-}/\mathrm{Fro}^\mathbb{Z},\underset{I}{\mathrm{homotopylimit}}~\underset{\mathrm{Spec}}{\mathcal{O}}^\mathrm{BK}{\Phi}^I_{*,-}/\mathrm{Fro}^\mathbb{Z}.	
}
\]

\end{definition}

\begin{definition}
We then consider the corresponding quasisheaves of perfect complexes of the corresponding condensed solid topological modules from \cite{10CS2}:
\begin{align}
\mathrm{Quasicoherentsheaves, Perfectcomplex, Condensed}_{*}	
\end{align}
where $*$ is one of the following spaces:
\begin{align}
&\underset{\mathrm{Spec}}{\mathcal{O}}^\mathrm{CS}\widetilde{\Delta}_{*,-}/\mathrm{Fro}^\mathbb{Z},\underset{\mathrm{Spec}}{\mathcal{O}}^\mathrm{CS}\widetilde{\nabla}_{*,-}/\mathrm{Fro}^\mathbb{Z},\underset{\mathrm{Spec}}{\mathcal{O}}^\mathrm{CS}\widetilde{\Phi}_{*,-}/\mathrm{Fro}^\mathbb{Z},\underset{\mathrm{Spec}}{\mathcal{O}}^\mathrm{CS}\widetilde{\Delta}^+_{*,-}/\mathrm{Fro}^\mathbb{Z},\\
&\underset{\mathrm{Spec}}{\mathcal{O}}^\mathrm{CS}\widetilde{\nabla}^+_{*,-}/\mathrm{Fro}^\mathbb{Z},\underset{\mathrm{Spec}}{\mathcal{O}}^\mathrm{CS}\widetilde{\Delta}^\dagger_{*,-}/\mathrm{Fro}^\mathbb{Z},\underset{\mathrm{Spec}}{\mathcal{O}}^\mathrm{CS}\widetilde{\nabla}^\dagger_{*,-}/\mathrm{Fro}^\mathbb{Z},	\\
\end{align}
\begin{align}
&\underset{\mathrm{Spec}}{\mathcal{O}}^\mathrm{CS}\breve{\Delta}_{*,-}/\mathrm{Fro}^\mathbb{Z},\breve{\nabla}_{*,-}/\mathrm{Fro}^\mathbb{Z},\underset{\mathrm{Spec}}{\mathcal{O}}^\mathrm{CS}\breve{\Phi}_{*,-}/\mathrm{Fro}^\mathbb{Z},\underset{\mathrm{Spec}}{\mathcal{O}}^\mathrm{CS}\breve{\Delta}^+_{*,-}/\mathrm{Fro}^\mathbb{Z},\\
&\underset{\mathrm{Spec}}{\mathcal{O}}^\mathrm{CS}\breve{\nabla}^+_{*,-}/\mathrm{Fro}^\mathbb{Z},\underset{\mathrm{Spec}}{\mathcal{O}}^\mathrm{CS}\breve{\Delta}^\dagger_{*,-}/\mathrm{Fro}^\mathbb{Z},\underset{\mathrm{Spec}}{\mathcal{O}}^\mathrm{CS}\breve{\nabla}^\dagger_{*,-}/\mathrm{Fro}^\mathbb{Z},	\\
\end{align}
\begin{align}
&\underset{\mathrm{Spec}}{\mathcal{O}}^\mathrm{CS}{\Delta}_{*,-}/\mathrm{Fro}^\mathbb{Z},\underset{\mathrm{Spec}}{\mathcal{O}}^\mathrm{CS}{\nabla}_{*,-}/\mathrm{Fro}^\mathbb{Z},\underset{\mathrm{Spec}}{\mathcal{O}}^\mathrm{CS}{\Phi}_{*,-}/\mathrm{Fro}^\mathbb{Z},\underset{\mathrm{Spec}}{\mathcal{O}}^\mathrm{CS}{\Delta}^+_{*,-}/\mathrm{Fro}^\mathbb{Z},\\
&\underset{\mathrm{Spec}}{\mathcal{O}}^\mathrm{CS}{\nabla}^+_{*,-}/\mathrm{Fro}^\mathbb{Z}, \underset{\mathrm{Spec}}{\mathcal{O}}^\mathrm{CS}{\Delta}^\dagger_{*,-}/\mathrm{Fro}^\mathbb{Z},\underset{\mathrm{Spec}}{\mathcal{O}}^\mathrm{CS}{\nabla}^\dagger_{*,-}/\mathrm{Fro}^\mathbb{Z}.	
\end{align}
Here for those space with notations related to the radius and the corresponding interval we consider the total unions $\bigcap_r,\bigcup_I$ in order to achieve the whole spaces to achieve the analogues of the corresponding FF curves from \cite{10KL1}, \cite{10KL2}, \cite{10FF} for
\[
\xymatrix@R+0pc@C+0pc{
\underset{r}{\mathrm{homotopycolimit}}~\underset{\mathrm{Spec}}{\mathcal{O}}^\mathrm{CS}\widetilde{\Phi}^r_{*,-},\underset{I}{\mathrm{homotopylimit}}~\underset{\mathrm{Spec}}{\mathcal{O}}^\mathrm{CS}\widetilde{\Phi}^I_{*,-},	\\
}
\]
\[
\xymatrix@R+0pc@C+0pc{
\underset{r}{\mathrm{homotopycolimit}}~\underset{\mathrm{Spec}}{\mathcal{O}}^\mathrm{CS}\breve{\Phi}^r_{*,-},\underset{I}{\mathrm{homotopylimit}}~\underset{\mathrm{Spec}}{\mathcal{O}}^\mathrm{CS}\breve{\Phi}^I_{*,-},	\\
}
\]
\[
\xymatrix@R+0pc@C+0pc{
\underset{r}{\mathrm{homotopycolimit}}~\underset{\mathrm{Spec}}{\mathcal{O}}^\mathrm{CS}{\Phi}^r_{*,-},\underset{I}{\mathrm{homotopylimit}}~\underset{\mathrm{Spec}}{\mathcal{O}}^\mathrm{CS}{\Phi}^I_{*,-}.	
}
\]
\[ 
\xymatrix@R+0pc@C+0pc{
\underset{r}{\mathrm{homotopycolimit}}~\underset{\mathrm{Spec}}{\mathcal{O}}^\mathrm{CS}\widetilde{\Phi}^r_{*,-}/\mathrm{Fro}^\mathbb{Z},\underset{I}{\mathrm{homotopylimit}}~\underset{\mathrm{Spec}}{\mathcal{O}}^\mathrm{CS}\widetilde{\Phi}^I_{*,-}/\mathrm{Fro}^\mathbb{Z},	\\
}
\]
\[ 
\xymatrix@R+0pc@C+0pc{
\underset{r}{\mathrm{homotopycolimit}}~\underset{\mathrm{Spec}}{\mathcal{O}}^\mathrm{CS}\breve{\Phi}^r_{*,-}/\mathrm{Fro}^\mathbb{Z},\underset{I}{\mathrm{homotopylimit}}~\breve{\Phi}^I_{*,-}/\mathrm{Fro}^\mathbb{Z},	\\
}
\]
\[ 
\xymatrix@R+0pc@C+0pc{
\underset{r}{\mathrm{homotopycolimit}}~\underset{\mathrm{Spec}}{\mathcal{O}}^\mathrm{CS}{\Phi}^r_{*,-}/\mathrm{Fro}^\mathbb{Z},\underset{I}{\mathrm{homotopylimit}}~\underset{\mathrm{Spec}}{\mathcal{O}}^\mathrm{CS}{\Phi}^I_{*,-}/\mathrm{Fro}^\mathbb{Z}.	
}
\]

\end{definition}

\begin{proposition}
There is a well-defined functor from the $\infty$-category 
\begin{align}
\mathrm{Quasicoherentpresheaves,Perfectcomplex,Condensed}_{*}	
\end{align}
where $*$ is one of the following spaces:
\begin{align}
&\underset{\mathrm{Spec}}{\mathcal{O}}^\mathrm{CS}\widetilde{\Phi}_{*,-}/\mathrm{Fro}^\mathbb{Z},	\\
\end{align}
\begin{align}
&\underset{\mathrm{Spec}}{\mathcal{O}}^\mathrm{CS}\breve{\Phi}_{*,-}/\mathrm{Fro}^\mathbb{Z},	\\
\end{align}
\begin{align}
&\underset{\mathrm{Spec}}{\mathcal{O}}^\mathrm{CS}{\Phi}_{*,-}/\mathrm{Fro}^\mathbb{Z},	
\end{align}
to the $\infty$-category of $\mathrm{Fro}$-equivariant quasicoherent presheaves over similar spaces above correspondingly without the $\mathrm{Fro}$-quotients, and to the $\infty$-category of $\mathrm{Fro}$-equivariant quasicoherent modules over global sections of the structure $\infty$-sheaves of the similar spaces above correspondingly without the $\mathrm{Fro}$-quotients. Here for those space without notation related to the radius and the corresponding interval we consider the total unions $\bigcap_r,\bigcup_I$ in order to achieve the whole spaces to achieve the analogues of the corresponding FF curves from \cite{10KL1}, \cite{10KL2}, \cite{10FF} for
\[
\xymatrix@R+0pc@C+0pc{
\underset{r}{\mathrm{homotopycolimit}}~\underset{\mathrm{Spec}}{\mathcal{O}}^\mathrm{CS}\widetilde{\Phi}^r_{*,-},\underset{I}{\mathrm{homotopylimit}}~\underset{\mathrm{Spec}}{\mathcal{O}}^\mathrm{CS}\widetilde{\Phi}^I_{*,-},	\\
}
\]
\[
\xymatrix@R+0pc@C+0pc{
\underset{r}{\mathrm{homotopycolimit}}~\underset{\mathrm{Spec}}{\mathcal{O}}^\mathrm{CS}\breve{\Phi}^r_{*,-},\underset{I}{\mathrm{homotopylimit}}~\underset{\mathrm{Spec}}{\mathcal{O}}^\mathrm{CS}\breve{\Phi}^I_{*,-},	\\
}
\]
\[
\xymatrix@R+0pc@C+0pc{
\underset{r}{\mathrm{homotopycolimit}}~\underset{\mathrm{Spec}}{\mathcal{O}}^\mathrm{CS}{\Phi}^r_{*,-},\underset{I}{\mathrm{homotopylimit}}~\underset{\mathrm{Spec}}{\mathcal{O}}^\mathrm{CS}{\Phi}^I_{*,-}.	
}
\]
\[ 
\xymatrix@R+0pc@C+0pc{
\underset{r}{\mathrm{homotopycolimit}}~\underset{\mathrm{Spec}}{\mathcal{O}}^\mathrm{CS}\widetilde{\Phi}^r_{*,-}/\mathrm{Fro}^\mathbb{Z},\underset{I}{\mathrm{homotopylimit}}~\underset{\mathrm{Spec}}{\mathcal{O}}^\mathrm{CS}\widetilde{\Phi}^I_{*,-}/\mathrm{Fro}^\mathbb{Z},	\\
}
\]
\[ 
\xymatrix@R+0pc@C+0pc{
\underset{r}{\mathrm{homotopycolimit}}~\underset{\mathrm{Spec}}{\mathcal{O}}^\mathrm{CS}\breve{\Phi}^r_{*,-}/\mathrm{Fro}^\mathbb{Z},\underset{I}{\mathrm{homotopylimit}}~\breve{\Phi}^I_{*,-}/\mathrm{Fro}^\mathbb{Z},	\\
}
\]
\[ 
\xymatrix@R+0pc@C+0pc{
\underset{r}{\mathrm{homotopycolimit}}~\underset{\mathrm{Spec}}{\mathcal{O}}^\mathrm{CS}{\Phi}^r_{*,-}/\mathrm{Fro}^\mathbb{Z},\underset{I}{\mathrm{homotopylimit}}~\underset{\mathrm{Spec}}{\mathcal{O}}^\mathrm{CS}{\Phi}^I_{*,-}/\mathrm{Fro}^\mathbb{Z}.	
}
\]	
In this situation we will have the target category being family parametrized by $r$ or $I$ in compatible glueing sense as in \cite[Definition 5.4.10]{10KL2}. In this situation for modules parametrized by the intervals we have the equivalence of $\infty$-categories by using \cite[Proposition 12.18]{10CS2}. Here the corresponding quasicoherent Frobenius modules are defined to be the corresponding homotopy colimits and limits of Frobenius modules:
\begin{align}
\underset{r}{\mathrm{homotopycolimit}}~M_r,\\
\underset{I}{\mathrm{homotopylimit}}~M_I,	
\end{align}
where each $M_r$ is a Frobenius-equivariant module over the period ring with respect to some radius $r$ while each $M_I$ is a Frobenius-equivariant module over the period ring with respect to some interval $I$.\\
\end{proposition}

\begin{proposition}
Similar proposition holds for 
\begin{align}
\mathrm{Quasicoherentsheaves,Perfectcomplex,IndBanach}_{*}.	
\end{align}	
\end{proposition}

\chapter{Over General Stacks}

\section{Over Preadic Spaces}

This chapter follows closely \cite{10T1}, \cite{10T2}, \cite{10T3}, \cite{10KPX}, \cite{10KP}, \cite{10KL1}, \cite{10KL2}, \cite{10BK}, \cite{10BBBK}, \cite{10BBM}, \cite{10KKM}, \cite{10CS1}, \cite{10CS2}, \cite{10CKZ}, \cite{10PZ}, \cite{10BCM}, \cite{10LBV}. All the preadic spaces will be those defined as in \cite{10Hu1}, \cite{10Hu2}, \cite{10KL1}, \cite{10KL2}, \cite{10SW}, while we regard them as the corresponding $\infty$-ringed $\infty$-toposes in the sense of Bambozzi-Kremnizer and Clausen-Scholze as in \cite{10BK}, \cite{10BBBK}, \cite{10BBM}, \cite{10KKM}, \cite{10CS1}, \cite{10CS2} by directly animating the corresponding presheaves in the previous categories to reach the enhancement.

\subsection{Multivariate Hodge Iwasawa Modules}

\subsubsection{Frobenius Quasicoherent Modules I}

\begin{definition}
Let $\psi$ be a toric tower over $\mathbb{Q}_p$ as in \cite[Chapter 7]{10KL2} with base $\mathbb{Q}_p\left<X_1^{\pm 1},...,X_k^{\pm 1}\right>$. Then from \cite{10KL1} and \cite[Definition 5.2.1]{10KL2} we have the following class of Kedlaya-Liu rings (with the following replacement: $\Delta$ stands for $A$, $\nabla$ stands for $B$, while $\Phi$ stands for $C$) by taking product in the sense of self $\Gamma$-th power:

\[
\xymatrix@R+0pc@C+0pc{
\widetilde{\Delta}_{\psi,\Gamma},\widetilde{\nabla}_{\psi,\Gamma},\widetilde{\Phi}_{\psi,\Gamma},\widetilde{\Delta}^+_{\psi,\Gamma},\widetilde{\nabla}^+_{\psi,\Gamma},\widetilde{\Delta}^\dagger_{\psi,\Gamma},\widetilde{\nabla}^\dagger_{\psi,\Gamma},\widetilde{\Phi}^r_{\psi,\Gamma},\widetilde{\Phi}^I_{\psi,\Gamma}, 
}
\]

\[
\xymatrix@R+0pc@C+0pc{
\breve{\Delta}_{\psi,\Gamma},\breve{\nabla}_{\psi,\Gamma},\breve{\Phi}_{\psi,\Gamma},\breve{\Delta}^+_{\psi,\Gamma},\breve{\nabla}^+_{\psi,\Gamma},\breve{\Delta}^\dagger_{\psi,\Gamma},\breve{\nabla}^\dagger_{\psi,\Gamma},\breve{\Phi}^r_{\psi,\Gamma},\breve{\Phi}^I_{\psi,\Gamma},	
}
\]

\[
\xymatrix@R+0pc@C+0pc{
{\Delta}_{\psi,\Gamma},{\nabla}_{\psi,\Gamma},{\Phi}_{\psi,\Gamma},{\Delta}^+_{\psi,\Gamma},{\nabla}^+_{\psi,\Gamma},{\Delta}^\dagger_{\psi,\Gamma},{\nabla}^\dagger_{\psi,\Gamma},{\Phi}^r_{\psi,\Gamma},{\Phi}^I_{\psi,\Gamma}.	
}
\]
Here in the following we have $X$ a preadic space over $\mathbb{Q}_p$. Taking the product we have:
\[
\xymatrix@R+0pc@C+0pc{
\widetilde{\Phi}_{\psi,\Gamma,X},\widetilde{\Phi}^r_{\psi,\Gamma,X},\widetilde{\Phi}^I_{\psi,\Gamma,X},	
}
\]
\[
\xymatrix@R+0pc@C+0pc{
\breve{\Phi}_{\psi,\Gamma,X},\breve{\Phi}^r_{\psi,\Gamma,X},\breve{\Phi}^I_{\psi,\Gamma,X},	
}
\]
\[
\xymatrix@R+0pc@C+0pc{
{\Phi}_{\psi,\Gamma,X},{\Phi}^r_{\psi,\Gamma,X},{\Phi}^I_{\psi,\Gamma,X}.	
}
\]
They carry multi Frobenius action $\varphi_\Gamma$ and multi $\mathrm{Lie}_\Gamma:=\mathbb{Z}_p^{\times\Gamma}$ action. In our current situation after \cite{10CKZ} and \cite{10PZ} we consider the following $(\infty,1)$-categories of $(\infty,1)$-modules.\\
\end{definition}

\begin{definition}
First we consider the Bambozzi-Kremnizer spectrum $\mathrm{Spec}^\mathrm{BK}(*)$ attached to any of those in the above from \cite{10BK} by taking derived rational localization:
\begin{align}
&\mathrm{Spec}^\mathrm{BK}\widetilde{\Phi}_{\psi,\Gamma,X},\mathrm{Spec}^\mathrm{BK}\widetilde{\Phi}^r_{\psi,\Gamma,X},\mathrm{Spec}^\mathrm{BK}\widetilde{\Phi}^I_{\psi,\Gamma,X},	
\end{align}
\begin{align}
&\mathrm{Spec}^\mathrm{BK}\breve{\Phi}_{\psi,\Gamma,X},\mathrm{Spec}^\mathrm{BK}\breve{\Phi}^r_{\psi,\Gamma,X},\mathrm{Spec}^\mathrm{BK}\breve{\Phi}^I_{\psi,\Gamma,X},	
\end{align}
\begin{align}
&\mathrm{Spec}^\mathrm{BK}{\Phi}_{\psi,\Gamma,X},
\mathrm{Spec}^\mathrm{BK}{\Phi}^r_{\psi,\Gamma,X},\mathrm{Spec}^\mathrm{BK}{\Phi}^I_{\psi,\Gamma,X}.	
\end{align}

Then we take the corresponding quotients by using the corresponding Frobenius operators:
\begin{align}
&\mathrm{Spec}^\mathrm{BK}\widetilde{\Phi}_{\psi,\Gamma,X}/\mathrm{Fro}^\mathbb{Z},	\\
\end{align}
\begin{align}
&\mathrm{Spec}^\mathrm{BK}\breve{\Phi}_{\psi,\Gamma,X}/\mathrm{Fro}^\mathbb{Z},	\\
\end{align}
\begin{align}
&\mathrm{Spec}^\mathrm{BK}{\Phi}_{\psi,\Gamma,X}/\mathrm{Fro}^\mathbb{Z}.	
\end{align}
Here for those space without notation related to the radius and the corresponding interval we consider the total unions $\bigcap_r,\bigcup_I$ in order to achieve the whole spaces to achieve the analogues of the corresponding FF curves from \cite{10KL1}, \cite{10KL2}, \cite{10FF} for
\[
\xymatrix@R+0pc@C+0pc{
\underset{r}{\mathrm{homotopylimit}}~\mathrm{Spec}^\mathrm{BK}\widetilde{\Phi}^r_{\psi,\Gamma,X},\underset{I}{\mathrm{homotopycolimit}}~\mathrm{Spec}^\mathrm{BK}\widetilde{\Phi}^I_{\psi,\Gamma,X},	\\
}
\]
\[
\xymatrix@R+0pc@C+0pc{
\underset{r}{\mathrm{homotopylimit}}~\mathrm{Spec}^\mathrm{BK}\breve{\Phi}^r_{\psi,\Gamma,X},\underset{I}{\mathrm{homotopycolimit}}~\mathrm{Spec}^\mathrm{BK}\breve{\Phi}^I_{\psi,\Gamma,X},	\\
}
\]
\[
\xymatrix@R+0pc@C+0pc{
\underset{r}{\mathrm{homotopylimit}}~\mathrm{Spec}^\mathrm{BK}{\Phi}^r_{\psi,\Gamma,X},\underset{I}{\mathrm{homotopycolimit}}~\mathrm{Spec}^\mathrm{BK}{\Phi}^I_{\psi,\Gamma,X}.	
}
\]
\[  
\xymatrix@R+0pc@C+0pc{
\underset{r}{\mathrm{homotopylimit}}~\mathrm{Spec}^\mathrm{BK}\widetilde{\Phi}^r_{\psi,\Gamma,X}/\mathrm{Fro}^\mathbb{Z},\underset{I}{\mathrm{homotopycolimit}}~\mathrm{Spec}^\mathrm{BK}\widetilde{\Phi}^I_{\psi,\Gamma,X}/\mathrm{Fro}^\mathbb{Z},	\\
}
\]
\[ 
\xymatrix@R+0pc@C+0pc{
\underset{r}{\mathrm{homotopylimit}}~\mathrm{Spec}^\mathrm{BK}\breve{\Phi}^r_{\psi,\Gamma,X}/\mathrm{Fro}^\mathbb{Z},\underset{I}{\mathrm{homotopycolimit}}~\mathrm{Spec}^\mathrm{BK}\breve{\Phi}^I_{\psi,\Gamma,X}/\mathrm{Fro}^\mathbb{Z},	\\
}
\]
\[ 
\xymatrix@R+0pc@C+0pc{
\underset{r}{\mathrm{homotopylimit}}~\mathrm{Spec}^\mathrm{BK}{\Phi}^r_{\psi,\Gamma,X}/\mathrm{Fro}^\mathbb{Z},\underset{I}{\mathrm{homotopycolimit}}~\mathrm{Spec}^\mathrm{BK}{\Phi}^I_{\psi,\Gamma,X}/\mathrm{Fro}^\mathbb{Z}.	
}
\]

\end{definition}

\indent Meanwhile we have the corresponding Clausen-Scholze analytic stacks from \cite{10CS2}, therefore applying their construction we have:

\begin{definition}
Here we define the following products by using the solidified tensor product from \cite{10CS1} and \cite{10CS2}. Then we take solidified tensor product $\overset{\blacksquare}{\otimes}$ of any of the following
\[
\xymatrix@R+0pc@C+0pc{
\widetilde{\Delta}_{\psi,\Gamma},\widetilde{\nabla}_{\psi,\Gamma},\widetilde{\Phi}_{\psi,\Gamma},\widetilde{\Delta}^+_{\psi,\Gamma},\widetilde{\nabla}^+_{\psi,\Gamma},\widetilde{\Delta}^\dagger_{\psi,\Gamma},\widetilde{\nabla}^\dagger_{\psi,\Gamma},\widetilde{\Phi}^r_{\psi,\Gamma},\widetilde{\Phi}^I_{\psi,\Gamma}, 
}
\]

\[
\xymatrix@R+0pc@C+0pc{
\breve{\Delta}_{\psi,\Gamma},\breve{\nabla}_{\psi,\Gamma},\breve{\Phi}_{\psi,\Gamma},\breve{\Delta}^+_{\psi,\Gamma},\breve{\nabla}^+_{\psi,\Gamma},\breve{\Delta}^\dagger_{\psi,\Gamma},\breve{\nabla}^\dagger_{\psi,\Gamma},\breve{\Phi}^r_{\psi,\Gamma},\breve{\Phi}^I_{\psi,\Gamma},	
}
\]

\[
\xymatrix@R+0pc@C+0pc{
{\Delta}_{\psi,\Gamma},{\nabla}_{\psi,\Gamma},{\Phi}_{\psi,\Gamma},{\Delta}^+_{\psi,\Gamma},{\nabla}^+_{\psi,\Gamma},{\Delta}^\dagger_{\psi,\Gamma},{\nabla}^\dagger_{\psi,\Gamma},{\Phi}^r_{\psi,\Gamma},{\Phi}^I_{\psi,\Gamma},	
}
\]  	
with $X$. Then we have the notations:
\[
\xymatrix@R+0pc@C+0pc{
\widetilde{\Delta}_{\psi,\Gamma,X},\widetilde{\nabla}_{\psi,\Gamma,X},\widetilde{\Phi}_{\psi,\Gamma,X},\widetilde{\Delta}^+_{\psi,\Gamma,X},\widetilde{\nabla}^+_{\psi,\Gamma,X},\widetilde{\Delta}^\dagger_{\psi,\Gamma,X},\widetilde{\nabla}^\dagger_{\psi,\Gamma,X},\widetilde{\Phi}^r_{\psi,\Gamma,X},\widetilde{\Phi}^I_{\psi,\Gamma,X}, 
}
\]

\[
\xymatrix@R+0pc@C+0pc{
\breve{\Delta}_{\psi,\Gamma,X},\breve{\nabla}_{\psi,\Gamma,X},\breve{\Phi}_{\psi,\Gamma,X},\breve{\Delta}^+_{\psi,\Gamma,X},\breve{\nabla}^+_{\psi,\Gamma,X},\breve{\Delta}^\dagger_{\psi,\Gamma,X},\breve{\nabla}^\dagger_{\psi,\Gamma,X},\breve{\Phi}^r_{\psi,\Gamma,X},\breve{\Phi}^I_{\psi,\Gamma,X},	
}
\]

\[
\xymatrix@R+0pc@C+0pc{
{\Delta}_{\psi,\Gamma,X},{\nabla}_{\psi,\Gamma,X},{\Phi}_{\psi,\Gamma,X},{\Delta}^+_{\psi,\Gamma,X},{\nabla}^+_{\psi,\Gamma,X},{\Delta}^\dagger_{\psi,\Gamma,X},{\nabla}^\dagger_{\psi,\Gamma,X},{\Phi}^r_{\psi,\Gamma,X},{\Phi}^I_{\psi,\Gamma,X}.	
}
\]
\end{definition}

\begin{definition}
First we consider the Clausen-Scholze spectrum $\mathrm{Spec}^\mathrm{CS}(*)$ attached to any of those in the above from \cite{10CS2} by taking derived rational localization:
\begin{align}
\mathrm{Spec}^\mathrm{CS}\widetilde{\Delta}_{\psi,\Gamma,X},\mathrm{Spec}^\mathrm{CS}\widetilde{\nabla}_{\psi,\Gamma,X},\mathrm{Spec}^\mathrm{CS}\widetilde{\Phi}_{\psi,\Gamma,X},\mathrm{Spec}^\mathrm{CS}\widetilde{\Delta}^+_{\psi,\Gamma,X},\mathrm{Spec}^\mathrm{CS}\widetilde{\nabla}^+_{\psi,\Gamma,X},\\
\mathrm{Spec}^\mathrm{CS}\widetilde{\Delta}^\dagger_{\psi,\Gamma,X},\mathrm{Spec}^\mathrm{CS}\widetilde{\nabla}^\dagger_{\psi,\Gamma,X},\mathrm{Spec}^\mathrm{CS}\widetilde{\Phi}^r_{\psi,\Gamma,X},\mathrm{Spec}^\mathrm{CS}\widetilde{\Phi}^I_{\psi,\Gamma,X},	\\
\end{align}
\begin{align}
\mathrm{Spec}^\mathrm{CS}\breve{\Delta}_{\psi,\Gamma,X},\breve{\nabla}_{\psi,\Gamma,X},\mathrm{Spec}^\mathrm{CS}\breve{\Phi}_{\psi,\Gamma,X},\mathrm{Spec}^\mathrm{CS}\breve{\Delta}^+_{\psi,\Gamma,X},\mathrm{Spec}^\mathrm{CS}\breve{\nabla}^+_{\psi,\Gamma,X},\\
\mathrm{Spec}^\mathrm{CS}\breve{\Delta}^\dagger_{\psi,\Gamma,X},\mathrm{Spec}^\mathrm{CS}\breve{\nabla}^\dagger_{\psi,\Gamma,X},\mathrm{Spec}^\mathrm{CS}\breve{\Phi}^r_{\psi,\Gamma,X},\breve{\Phi}^I_{\psi,\Gamma,X},	\\
\end{align}
\begin{align}
\mathrm{Spec}^\mathrm{CS}{\Delta}_{\psi,\Gamma,X},\mathrm{Spec}^\mathrm{CS}{\nabla}_{\psi,\Gamma,X},\mathrm{Spec}^\mathrm{CS}{\Phi}_{\psi,\Gamma,X},\mathrm{Spec}^\mathrm{CS}{\Delta}^+_{\psi,\Gamma,X},\mathrm{Spec}^\mathrm{CS}{\nabla}^+_{\psi,\Gamma,X},\\
\mathrm{Spec}^\mathrm{CS}{\Delta}^\dagger_{\psi,\Gamma,X},\mathrm{Spec}^\mathrm{CS}{\nabla}^\dagger_{\psi,\Gamma,X},\mathrm{Spec}^\mathrm{CS}{\Phi}^r_{\psi,\Gamma,X},\mathrm{Spec}^\mathrm{CS}{\Phi}^I_{\psi,\Gamma,X}.	
\end{align}

Then we take the corresponding quotients by using the corresponding Frobenius operators:
\begin{align}
&\mathrm{Spec}^\mathrm{CS}\widetilde{\Delta}_{\psi,\Gamma,X}/\mathrm{Fro}^\mathbb{Z},\mathrm{Spec}^\mathrm{CS}\widetilde{\nabla}_{\psi,\Gamma,X}/\mathrm{Fro}^\mathbb{Z},\mathrm{Spec}^\mathrm{CS}\widetilde{\Phi}_{\psi,\Gamma,X}/\mathrm{Fro}^\mathbb{Z},\mathrm{Spec}^\mathrm{CS}\widetilde{\Delta}^+_{\psi,\Gamma,X}/\mathrm{Fro}^\mathbb{Z},\\
&\mathrm{Spec}^\mathrm{CS}\widetilde{\nabla}^+_{\psi,\Gamma,X}/\mathrm{Fro}^\mathbb{Z}, \mathrm{Spec}^\mathrm{CS}\widetilde{\Delta}^\dagger_{\psi,\Gamma,X}/\mathrm{Fro}^\mathbb{Z},\mathrm{Spec}^\mathrm{CS}\widetilde{\nabla}^\dagger_{\psi,\Gamma,X}/\mathrm{Fro}^\mathbb{Z},	\\
\end{align}
\begin{align}
&\mathrm{Spec}^\mathrm{CS}\breve{\Delta}_{\psi,\Gamma,X}/\mathrm{Fro}^\mathbb{Z},\breve{\nabla}_{\psi,\Gamma,X}/\mathrm{Fro}^\mathbb{Z},\mathrm{Spec}^\mathrm{CS}\breve{\Phi}_{\psi,\Gamma,X}/\mathrm{Fro}^\mathbb{Z},\mathrm{Spec}^\mathrm{CS}\breve{\Delta}^+_{\psi,\Gamma,X}/\mathrm{Fro}^\mathbb{Z},\\
&\mathrm{Spec}^\mathrm{CS}\breve{\nabla}^+_{\psi,\Gamma,X}/\mathrm{Fro}^\mathbb{Z}, \mathrm{Spec}^\mathrm{CS}\breve{\Delta}^\dagger_{\psi,\Gamma,X}/\mathrm{Fro}^\mathbb{Z},\mathrm{Spec}^\mathrm{CS}\breve{\nabla}^\dagger_{\psi,\Gamma,X}/\mathrm{Fro}^\mathbb{Z},	\\
\end{align}
\begin{align}
&\mathrm{Spec}^\mathrm{CS}{\Delta}_{\psi,\Gamma,X}/\mathrm{Fro}^\mathbb{Z},\mathrm{Spec}^\mathrm{CS}{\nabla}_{\psi,\Gamma,X}/\mathrm{Fro}^\mathbb{Z},\mathrm{Spec}^\mathrm{CS}{\Phi}_{\psi,\Gamma,X}/\mathrm{Fro}^\mathbb{Z},\mathrm{Spec}^\mathrm{CS}{\Delta}^+_{\psi,\Gamma,X}/\mathrm{Fro}^\mathbb{Z},\\
&\mathrm{Spec}^\mathrm{CS}{\nabla}^+_{\psi,\Gamma,X}/\mathrm{Fro}^\mathbb{Z}, \mathrm{Spec}^\mathrm{CS}{\Delta}^\dagger_{\psi,\Gamma,X}/\mathrm{Fro}^\mathbb{Z},\mathrm{Spec}^\mathrm{CS}{\nabla}^\dagger_{\psi,\Gamma,X}/\mathrm{Fro}^\mathbb{Z}.	
\end{align}
Here for those space with notations related to the radius and the corresponding interval we consider the total unions $\bigcap_r,\bigcup_I$ in order to achieve the whole spaces to achieve the analogues of the corresponding FF curves from \cite{10KL1}, \cite{10KL2}, \cite{10FF} for
\[
\xymatrix@R+0pc@C+0pc{
\underset{r}{\mathrm{homotopylimit}}~\mathrm{Spec}^\mathrm{CS}\widetilde{\Phi}^r_{\psi,\Gamma,X},\underset{I}{\mathrm{homotopycolimit}}~\mathrm{Spec}^\mathrm{CS}\widetilde{\Phi}^I_{\psi,\Gamma,X},	\\
}
\]
\[
\xymatrix@R+0pc@C+0pc{
\underset{r}{\mathrm{homotopylimit}}~\mathrm{Spec}^\mathrm{CS}\breve{\Phi}^r_{\psi,\Gamma,X},\underset{I}{\mathrm{homotopycolimit}}~\mathrm{Spec}^\mathrm{CS}\breve{\Phi}^I_{\psi,\Gamma,X},	\\
}
\]
\[
\xymatrix@R+0pc@C+0pc{
\underset{r}{\mathrm{homotopylimit}}~\mathrm{Spec}^\mathrm{CS}{\Phi}^r_{\psi,\Gamma,X},\underset{I}{\mathrm{homotopycolimit}}~\mathrm{Spec}^\mathrm{CS}{\Phi}^I_{\psi,\Gamma,X}.	
}
\]
\[ 
\xymatrix@R+0pc@C+0pc{
\underset{r}{\mathrm{homotopylimit}}~\mathrm{Spec}^\mathrm{CS}\widetilde{\Phi}^r_{\psi,\Gamma,X}/\mathrm{Fro}^\mathbb{Z},\underset{I}{\mathrm{homotopycolimit}}~\mathrm{Spec}^\mathrm{CS}\widetilde{\Phi}^I_{\psi,\Gamma,X}/\mathrm{Fro}^\mathbb{Z},	\\
}
\]
\[ 
\xymatrix@R+0pc@C+0pc{
\underset{r}{\mathrm{homotopylimit}}~\mathrm{Spec}^\mathrm{CS}\breve{\Phi}^r_{\psi,\Gamma,X}/\mathrm{Fro}^\mathbb{Z},\underset{I}{\mathrm{homotopycolimit}}~\breve{\Phi}^I_{\psi,\Gamma,X}/\mathrm{Fro}^\mathbb{Z},	\\
}
\]
\[ 
\xymatrix@R+0pc@C+0pc{
\underset{r}{\mathrm{homotopylimit}}~\mathrm{Spec}^\mathrm{CS}{\Phi}^r_{\psi,\Gamma,X}/\mathrm{Fro}^\mathbb{Z},\underset{I}{\mathrm{homotopycolimit}}~\mathrm{Spec}^\mathrm{CS}{\Phi}^I_{\psi,\Gamma,X}/\mathrm{Fro}^\mathbb{Z}.	
}
\]

\end{definition}

\

\begin{definition}
We then consider the corresponding quasipresheaves of the corresponding ind-Banach or monomorphic ind-Banach modules from \cite{10BBK}, \cite{10KKM}:
\begin{align}
\mathrm{Quasicoherentpresheaves,IndBanach}_{*}	
\end{align}
where $*$ is one of the following spaces:
\begin{align}
&\mathrm{Spec}^\mathrm{BK}\widetilde{\Phi}_{\psi,\Gamma,X}/\mathrm{Fro}^\mathbb{Z},	\\
\end{align}
\begin{align}
&\mathrm{Spec}^\mathrm{BK}\breve{\Phi}_{\psi,\Gamma,X}/\mathrm{Fro}^\mathbb{Z},	\\
\end{align}
\begin{align}
&\mathrm{Spec}^\mathrm{BK}{\Phi}_{\psi,\Gamma,X}/\mathrm{Fro}^\mathbb{Z}.	
\end{align}
Here for those space without notation related to the radius and the corresponding interval we consider the total unions $\bigcap_r,\bigcup_I$ in order to achieve the whole spaces to achieve the analogues of the corresponding FF curves from \cite{10KL1}, \cite{10KL2}, \cite{10FF} for
\[
\xymatrix@R+0pc@C+0pc{
\underset{r}{\mathrm{homotopylimit}}~\mathrm{Spec}^\mathrm{BK}\widetilde{\Phi}^r_{\psi,\Gamma,X},\underset{I}{\mathrm{homotopycolimit}}~\mathrm{Spec}^\mathrm{BK}\widetilde{\Phi}^I_{\psi,\Gamma,X},	\\
}
\]
\[
\xymatrix@R+0pc@C+0pc{
\underset{r}{\mathrm{homotopylimit}}~\mathrm{Spec}^\mathrm{BK}\breve{\Phi}^r_{\psi,\Gamma,X},\underset{I}{\mathrm{homotopycolimit}}~\mathrm{Spec}^\mathrm{BK}\breve{\Phi}^I_{\psi,\Gamma,X},	\\
}
\]
\[
\xymatrix@R+0pc@C+0pc{
\underset{r}{\mathrm{homotopylimit}}~\mathrm{Spec}^\mathrm{BK}{\Phi}^r_{\psi,\Gamma,X},\underset{I}{\mathrm{homotopycolimit}}~\mathrm{Spec}^\mathrm{BK}{\Phi}^I_{\psi,\Gamma,X}.	
}
\]
\[  
\xymatrix@R+0pc@C+0pc{
\underset{r}{\mathrm{homotopylimit}}~\mathrm{Spec}^\mathrm{BK}\widetilde{\Phi}^r_{\psi,\Gamma,X}/\mathrm{Fro}^\mathbb{Z},\underset{I}{\mathrm{homotopycolimit}}~\mathrm{Spec}^\mathrm{BK}\widetilde{\Phi}^I_{\psi,\Gamma,X}/\mathrm{Fro}^\mathbb{Z},	\\
}
\]
\[ 
\xymatrix@R+0pc@C+0pc{
\underset{r}{\mathrm{homotopylimit}}~\mathrm{Spec}^\mathrm{BK}\breve{\Phi}^r_{\psi,\Gamma,X}/\mathrm{Fro}^\mathbb{Z},\underset{I}{\mathrm{homotopycolimit}}~\mathrm{Spec}^\mathrm{BK}\breve{\Phi}^I_{\psi,\Gamma,X}/\mathrm{Fro}^\mathbb{Z},	\\
}
\]
\[ 
\xymatrix@R+0pc@C+0pc{
\underset{r}{\mathrm{homotopylimit}}~\mathrm{Spec}^\mathrm{BK}{\Phi}^r_{\psi,\Gamma,X}/\mathrm{Fro}^\mathbb{Z},\underset{I}{\mathrm{homotopycolimit}}~\mathrm{Spec}^\mathrm{BK}{\Phi}^I_{\psi,\Gamma,X}/\mathrm{Fro}^\mathbb{Z}.	
}
\]

\end{definition}

\begin{definition}
We then consider the corresponding quasisheaves of the corresponding condensed solid topological modules from \cite{10CS2}:
\begin{align}
\mathrm{Quasicoherentsheaves, Condensed}_{*}	
\end{align}
where $*$ is one of the following spaces:
\begin{align}
&\mathrm{Spec}^\mathrm{CS}\widetilde{\Delta}_{\psi,\Gamma,X}/\mathrm{Fro}^\mathbb{Z},\mathrm{Spec}^\mathrm{CS}\widetilde{\nabla}_{\psi,\Gamma,X}/\mathrm{Fro}^\mathbb{Z},\mathrm{Spec}^\mathrm{CS}\widetilde{\Phi}_{\psi,\Gamma,X}/\mathrm{Fro}^\mathbb{Z},\mathrm{Spec}^\mathrm{CS}\widetilde{\Delta}^+_{\psi,\Gamma,X}/\mathrm{Fro}^\mathbb{Z},\\
&\mathrm{Spec}^\mathrm{CS}\widetilde{\nabla}^+_{\psi,\Gamma,X}/\mathrm{Fro}^\mathbb{Z},\mathrm{Spec}^\mathrm{CS}\widetilde{\Delta}^\dagger_{\psi,\Gamma,X}/\mathrm{Fro}^\mathbb{Z},\mathrm{Spec}^\mathrm{CS}\widetilde{\nabla}^\dagger_{\psi,\Gamma,X}/\mathrm{Fro}^\mathbb{Z},	\\
\end{align}
\begin{align}
&\mathrm{Spec}^\mathrm{CS}\breve{\Delta}_{\psi,\Gamma,X}/\mathrm{Fro}^\mathbb{Z},\breve{\nabla}_{\psi,\Gamma,X}/\mathrm{Fro}^\mathbb{Z},\mathrm{Spec}^\mathrm{CS}\breve{\Phi}_{\psi,\Gamma,X}/\mathrm{Fro}^\mathbb{Z},\mathrm{Spec}^\mathrm{CS}\breve{\Delta}^+_{\psi,\Gamma,X}/\mathrm{Fro}^\mathbb{Z},\\
&\mathrm{Spec}^\mathrm{CS}\breve{\nabla}^+_{\psi,\Gamma,X}/\mathrm{Fro}^\mathbb{Z},\mathrm{Spec}^\mathrm{CS}\breve{\Delta}^\dagger_{\psi,\Gamma,X}/\mathrm{Fro}^\mathbb{Z},\mathrm{Spec}^\mathrm{CS}\breve{\nabla}^\dagger_{\psi,\Gamma,X}/\mathrm{Fro}^\mathbb{Z},	\\
\end{align}
\begin{align}
&\mathrm{Spec}^\mathrm{CS}{\Delta}_{\psi,\Gamma,X}/\mathrm{Fro}^\mathbb{Z},\mathrm{Spec}^\mathrm{CS}{\nabla}_{\psi,\Gamma,X}/\mathrm{Fro}^\mathbb{Z},\mathrm{Spec}^\mathrm{CS}{\Phi}_{\psi,\Gamma,X}/\mathrm{Fro}^\mathbb{Z},\mathrm{Spec}^\mathrm{CS}{\Delta}^+_{\psi,\Gamma,X}/\mathrm{Fro}^\mathbb{Z},\\
&\mathrm{Spec}^\mathrm{CS}{\nabla}^+_{\psi,\Gamma,X}/\mathrm{Fro}^\mathbb{Z}, \mathrm{Spec}^\mathrm{CS}{\Delta}^\dagger_{\psi,\Gamma,X}/\mathrm{Fro}^\mathbb{Z},\mathrm{Spec}^\mathrm{CS}{\nabla}^\dagger_{\psi,\Gamma,X}/\mathrm{Fro}^\mathbb{Z}.	
\end{align}
Here for those space with notations related to the radius and the corresponding interval we consider the total unions $\bigcap_r,\bigcup_I$ in order to achieve the whole spaces to achieve the analogues of the corresponding FF curves from \cite{10KL1}, \cite{10KL2}, \cite{10FF} for
\[
\xymatrix@R+0pc@C+0pc{
\underset{r}{\mathrm{homotopylimit}}~\mathrm{Spec}^\mathrm{CS}\widetilde{\Phi}^r_{\psi,\Gamma,X},\underset{I}{\mathrm{homotopycolimit}}~\mathrm{Spec}^\mathrm{CS}\widetilde{\Phi}^I_{\psi,\Gamma,X},	\\
}
\]
\[
\xymatrix@R+0pc@C+0pc{
\underset{r}{\mathrm{homotopylimit}}~\mathrm{Spec}^\mathrm{CS}\breve{\Phi}^r_{\psi,\Gamma,X},\underset{I}{\mathrm{homotopycolimit}}~\mathrm{Spec}^\mathrm{CS}\breve{\Phi}^I_{\psi,\Gamma,X},	\\
}
\]
\[
\xymatrix@R+0pc@C+0pc{
\underset{r}{\mathrm{homotopylimit}}~\mathrm{Spec}^\mathrm{CS}{\Phi}^r_{\psi,\Gamma,X},\underset{I}{\mathrm{homotopycolimit}}~\mathrm{Spec}^\mathrm{CS}{\Phi}^I_{\psi,\Gamma,X}.	
}
\]
\[ 
\xymatrix@R+0pc@C+0pc{
\underset{r}{\mathrm{homotopylimit}}~\mathrm{Spec}^\mathrm{CS}\widetilde{\Phi}^r_{\psi,\Gamma,X}/\mathrm{Fro}^\mathbb{Z},\underset{I}{\mathrm{homotopycolimit}}~\mathrm{Spec}^\mathrm{CS}\widetilde{\Phi}^I_{\psi,\Gamma,X}/\mathrm{Fro}^\mathbb{Z},	\\
}
\]
\[ 
\xymatrix@R+0pc@C+0pc{
\underset{r}{\mathrm{homotopylimit}}~\mathrm{Spec}^\mathrm{CS}\breve{\Phi}^r_{\psi,\Gamma,X}/\mathrm{Fro}^\mathbb{Z},\underset{I}{\mathrm{homotopycolimit}}~\breve{\Phi}^I_{\psi,\Gamma,X}/\mathrm{Fro}^\mathbb{Z},	\\
}
\]
\[ 
\xymatrix@R+0pc@C+0pc{
\underset{r}{\mathrm{homotopylimit}}~\mathrm{Spec}^\mathrm{CS}{\Phi}^r_{\psi,\Gamma,X}/\mathrm{Fro}^\mathbb{Z},\underset{I}{\mathrm{homotopycolimit}}~\mathrm{Spec}^\mathrm{CS}{\Phi}^I_{\psi,\Gamma,X}/\mathrm{Fro}^\mathbb{Z}.	
}
\]

\end{definition}

\

\begin{proposition}
There is a well-defined functor from the $\infty$-category 
\begin{align}
\mathrm{Quasicoherentpresheaves,Condensed}_{*}	
\end{align}
where $*$ is one of the following spaces:
\begin{align}
&\mathrm{Spec}^\mathrm{CS}\widetilde{\Phi}_{\psi,\Gamma,X}/\mathrm{Fro}^\mathbb{Z},	\\
\end{align}
\begin{align}
&\mathrm{Spec}^\mathrm{CS}\breve{\Phi}_{\psi,\Gamma,X}/\mathrm{Fro}^\mathbb{Z},	\\
\end{align}
\begin{align}
&\mathrm{Spec}^\mathrm{CS}{\Phi}_{\psi,\Gamma,X}/\mathrm{Fro}^\mathbb{Z},	
\end{align}
to the $\infty$-category of $\mathrm{Fro}$-equivariant quasicoherent presheaves over similar spaces above correspondingly without the $\mathrm{Fro}$-quotients, and to the $\infty$-category of $\mathrm{Fro}$-equivariant quasicoherent modules over global sections of the structure $\infty$-sheaves of the similar spaces above correspondingly without the $\mathrm{Fro}$-quotients\footnote{When we are considering at current generality the corresponding deformation over general preadic stacks, we will further consider the corresponding Frobenius modules over period sheaves (not rings though we talk about rings if it is not too confusing) carrying coefficients in preadic stacks.}. Here for those space without notation related to the radius and the corresponding interval we consider the total unions $\bigcap_r,\bigcup_I$ in order to achieve the whole spaces to achieve the analogues of the corresponding FF curves from \cite{10KL1}, \cite{10KL2}, \cite{10FF} for
\[
\xymatrix@R+0pc@C+0pc{
\underset{r}{\mathrm{homotopylimit}}~\mathrm{Spec}^\mathrm{CS}\widetilde{\Phi}^r_{\psi,\Gamma,X},\underset{I}{\mathrm{homotopycolimit}}~\mathrm{Spec}^\mathrm{CS}\widetilde{\Phi}^I_{\psi,\Gamma,X},	\\
}
\]
\[
\xymatrix@R+0pc@C+0pc{
\underset{r}{\mathrm{homotopylimit}}~\mathrm{Spec}^\mathrm{CS}\breve{\Phi}^r_{\psi,\Gamma,X},\underset{I}{\mathrm{homotopycolimit}}~\mathrm{Spec}^\mathrm{CS}\breve{\Phi}^I_{\psi,\Gamma,X},	\\
}
\]
\[
\xymatrix@R+0pc@C+0pc{
\underset{r}{\mathrm{homotopylimit}}~\mathrm{Spec}^\mathrm{CS}{\Phi}^r_{\psi,\Gamma,X},\underset{I}{\mathrm{homotopycolimit}}~\mathrm{Spec}^\mathrm{CS}{\Phi}^I_{\psi,\Gamma,X}.	
}
\]
\[ 
\xymatrix@R+0pc@C+0pc{
\underset{r}{\mathrm{homotopylimit}}~\mathrm{Spec}^\mathrm{CS}\widetilde{\Phi}^r_{\psi,\Gamma,X}/\mathrm{Fro}^\mathbb{Z},\underset{I}{\mathrm{homotopycolimit}}~\mathrm{Spec}^\mathrm{CS}\widetilde{\Phi}^I_{\psi,\Gamma,X}/\mathrm{Fro}^\mathbb{Z},	\\
}
\]
\[ 
\xymatrix@R+0pc@C+0pc{
\underset{r}{\mathrm{homotopylimit}}~\mathrm{Spec}^\mathrm{CS}\breve{\Phi}^r_{\psi,\Gamma,X}/\mathrm{Fro}^\mathbb{Z},\underset{I}{\mathrm{homotopycolimit}}~\breve{\Phi}^I_{\psi,\Gamma,X}/\mathrm{Fro}^\mathbb{Z},	\\
}
\]
\[ 
\xymatrix@R+0pc@C+0pc{
\underset{r}{\mathrm{homotopylimit}}~\mathrm{Spec}^\mathrm{CS}{\Phi}^r_{\psi,\Gamma,X}/\mathrm{Fro}^\mathbb{Z},\underset{I}{\mathrm{homotopycolimit}}~\mathrm{Spec}^\mathrm{CS}{\Phi}^I_{\psi,\Gamma,X}/\mathrm{Fro}^\mathbb{Z}.	
}
\]	
In this situation we will have the target category being family parametrized by $r$ or $I$ in compatible glueing sense as in \cite[Definition 5.4.10]{10KL2}. In this situation for modules parametrized by the intervals we have the equivalence of $\infty$-categories by using \cite[Proposition 13.8]{10CS2}. Here the corresponding quasicoherent Frobenius modules are defined to be the corresponding homotopy colimits and limits of Frobenius modules:
\begin{align}
\underset{r}{\mathrm{homotopycolimit}}~M_r,\\
\underset{I}{\mathrm{homotopylimit}}~M_I,	
\end{align}
where each $M_r$ is a Frobenius-equivariant module over the period ring with respect to some radius $r$ while each $M_I$ is a Frobenius-equivariant module over the period ring with respect to some interval $I$.\\
\end{proposition}

\begin{proposition}
Similar proposition holds for 
\begin{align}
\mathrm{Quasicoherentsheaves,IndBanach}_{*}.	
\end{align}	
\end{proposition}

\

\begin{definition}
We then consider the corresponding quasipresheaves of perfect complexes the corresponding ind-Banach or monomorphic ind-Banach modules from \cite{10BBK}, \cite{10KKM}:
\begin{align}
\mathrm{Quasicoherentpresheaves,Perfectcomplex,IndBanach}_{*}	
\end{align}
where $*$ is one of the following spaces:
\begin{align}
&\mathrm{Spec}^\mathrm{BK}\widetilde{\Phi}_{\psi,\Gamma,X}/\mathrm{Fro}^\mathbb{Z},	\\
\end{align}
\begin{align}
&\mathrm{Spec}^\mathrm{BK}\breve{\Phi}_{\psi,\Gamma,X}/\mathrm{Fro}^\mathbb{Z},	\\
\end{align}
\begin{align}
&\mathrm{Spec}^\mathrm{BK}{\Phi}_{\psi,\Gamma,X}/\mathrm{Fro}^\mathbb{Z}.	
\end{align}
Here for those space without notation related to the radius and the corresponding interval we consider the total unions $\bigcap_r,\bigcup_I$ in order to achieve the whole spaces to achieve the analogues of the corresponding FF curves from \cite{10KL1}, \cite{10KL2}, \cite{10FF} for
\[
\xymatrix@R+0pc@C+0pc{
\underset{r}{\mathrm{homotopylimit}}~\mathrm{Spec}^\mathrm{BK}\widetilde{\Phi}^r_{\psi,\Gamma,X},\underset{I}{\mathrm{homotopycolimit}}~\mathrm{Spec}^\mathrm{BK}\widetilde{\Phi}^I_{\psi,\Gamma,X},	\\
}
\]
\[
\xymatrix@R+0pc@C+0pc{
\underset{r}{\mathrm{homotopylimit}}~\mathrm{Spec}^\mathrm{BK}\breve{\Phi}^r_{\psi,\Gamma,X},\underset{I}{\mathrm{homotopycolimit}}~\mathrm{Spec}^\mathrm{BK}\breve{\Phi}^I_{\psi,\Gamma,X},	\\
}
\]
\[
\xymatrix@R+0pc@C+0pc{
\underset{r}{\mathrm{homotopylimit}}~\mathrm{Spec}^\mathrm{BK}{\Phi}^r_{\psi,\Gamma,X},\underset{I}{\mathrm{homotopycolimit}}~\mathrm{Spec}^\mathrm{BK}{\Phi}^I_{\psi,\Gamma,X}.	
}
\]
\[  
\xymatrix@R+0pc@C+0pc{
\underset{r}{\mathrm{homotopylimit}}~\mathrm{Spec}^\mathrm{BK}\widetilde{\Phi}^r_{\psi,\Gamma,X}/\mathrm{Fro}^\mathbb{Z},\underset{I}{\mathrm{homotopycolimit}}~\mathrm{Spec}^\mathrm{BK}\widetilde{\Phi}^I_{\psi,\Gamma,X}/\mathrm{Fro}^\mathbb{Z},	\\
}
\]
\[ 
\xymatrix@R+0pc@C+0pc{
\underset{r}{\mathrm{homotopylimit}}~\mathrm{Spec}^\mathrm{BK}\breve{\Phi}^r_{\psi,\Gamma,X}/\mathrm{Fro}^\mathbb{Z},\underset{I}{\mathrm{homotopycolimit}}~\mathrm{Spec}^\mathrm{BK}\breve{\Phi}^I_{\psi,\Gamma,X}/\mathrm{Fro}^\mathbb{Z},	\\
}
\]
\[ 
\xymatrix@R+0pc@C+0pc{
\underset{r}{\mathrm{homotopylimit}}~\mathrm{Spec}^\mathrm{BK}{\Phi}^r_{\psi,\Gamma,X}/\mathrm{Fro}^\mathbb{Z},\underset{I}{\mathrm{homotopycolimit}}~\mathrm{Spec}^\mathrm{BK}{\Phi}^I_{\psi,\Gamma,X}/\mathrm{Fro}^\mathbb{Z}.	
}
\]

\end{definition}

\begin{definition}
We then consider the corresponding quasisheaves of perfect complexes of the corresponding condensed solid topological modules from \cite{10CS2}:
\begin{align}
\mathrm{Quasicoherentsheaves, Perfectcomplex, Condensed}_{*}	
\end{align}
where $*$ is one of the following spaces:
\begin{align}
&\mathrm{Spec}^\mathrm{CS}\widetilde{\Delta}_{\psi,\Gamma,X}/\mathrm{Fro}^\mathbb{Z},\mathrm{Spec}^\mathrm{CS}\widetilde{\nabla}_{\psi,\Gamma,X}/\mathrm{Fro}^\mathbb{Z},\mathrm{Spec}^\mathrm{CS}\widetilde{\Phi}_{\psi,\Gamma,X}/\mathrm{Fro}^\mathbb{Z},\mathrm{Spec}^\mathrm{CS}\widetilde{\Delta}^+_{\psi,\Gamma,X}/\mathrm{Fro}^\mathbb{Z},\\
&\mathrm{Spec}^\mathrm{CS}\widetilde{\nabla}^+_{\psi,\Gamma,X}/\mathrm{Fro}^\mathbb{Z},\mathrm{Spec}^\mathrm{CS}\widetilde{\Delta}^\dagger_{\psi,\Gamma,X}/\mathrm{Fro}^\mathbb{Z},\mathrm{Spec}^\mathrm{CS}\widetilde{\nabla}^\dagger_{\psi,\Gamma,X}/\mathrm{Fro}^\mathbb{Z},	\\
\end{align}
\begin{align}
&\mathrm{Spec}^\mathrm{CS}\breve{\Delta}_{\psi,\Gamma,X}/\mathrm{Fro}^\mathbb{Z},\breve{\nabla}_{\psi,\Gamma,X}/\mathrm{Fro}^\mathbb{Z},\mathrm{Spec}^\mathrm{CS}\breve{\Phi}_{\psi,\Gamma,X}/\mathrm{Fro}^\mathbb{Z},\mathrm{Spec}^\mathrm{CS}\breve{\Delta}^+_{\psi,\Gamma,X}/\mathrm{Fro}^\mathbb{Z},\\
&\mathrm{Spec}^\mathrm{CS}\breve{\nabla}^+_{\psi,\Gamma,X}/\mathrm{Fro}^\mathbb{Z},\mathrm{Spec}^\mathrm{CS}\breve{\Delta}^\dagger_{\psi,\Gamma,X}/\mathrm{Fro}^\mathbb{Z},\mathrm{Spec}^\mathrm{CS}\breve{\nabla}^\dagger_{\psi,\Gamma,X}/\mathrm{Fro}^\mathbb{Z},	\\
\end{align}
\begin{align}
&\mathrm{Spec}^\mathrm{CS}{\Delta}_{\psi,\Gamma,X}/\mathrm{Fro}^\mathbb{Z},\mathrm{Spec}^\mathrm{CS}{\nabla}_{\psi,\Gamma,X}/\mathrm{Fro}^\mathbb{Z},\mathrm{Spec}^\mathrm{CS}{\Phi}_{\psi,\Gamma,X}/\mathrm{Fro}^\mathbb{Z},\mathrm{Spec}^\mathrm{CS}{\Delta}^+_{\psi,\Gamma,X}/\mathrm{Fro}^\mathbb{Z},\\
&\mathrm{Spec}^\mathrm{CS}{\nabla}^+_{\psi,\Gamma,X}/\mathrm{Fro}^\mathbb{Z}, \mathrm{Spec}^\mathrm{CS}{\Delta}^\dagger_{\psi,\Gamma,X}/\mathrm{Fro}^\mathbb{Z},\mathrm{Spec}^\mathrm{CS}{\nabla}^\dagger_{\psi,\Gamma,X}/\mathrm{Fro}^\mathbb{Z}.	
\end{align}
Here for those space with notations related to the radius and the corresponding interval we consider the total unions $\bigcap_r,\bigcup_I$ in order to achieve the whole spaces to achieve the analogues of the corresponding FF curves from \cite{10KL1}, \cite{10KL2}, \cite{10FF} for
\[
\xymatrix@R+0pc@C+0pc{
\underset{r}{\mathrm{homotopylimit}}~\mathrm{Spec}^\mathrm{CS}\widetilde{\Phi}^r_{\psi,\Gamma,X},\underset{I}{\mathrm{homotopycolimit}}~\mathrm{Spec}^\mathrm{CS}\widetilde{\Phi}^I_{\psi,\Gamma,X},	\\
}
\]
\[
\xymatrix@R+0pc@C+0pc{
\underset{r}{\mathrm{homotopylimit}}~\mathrm{Spec}^\mathrm{CS}\breve{\Phi}^r_{\psi,\Gamma,X},\underset{I}{\mathrm{homotopycolimit}}~\mathrm{Spec}^\mathrm{CS}\breve{\Phi}^I_{\psi,\Gamma,X},	\\
}
\]
\[
\xymatrix@R+0pc@C+0pc{
\underset{r}{\mathrm{homotopylimit}}~\mathrm{Spec}^\mathrm{CS}{\Phi}^r_{\psi,\Gamma,X},\underset{I}{\mathrm{homotopycolimit}}~\mathrm{Spec}^\mathrm{CS}{\Phi}^I_{\psi,\Gamma,X}.	
}
\]
\[ 
\xymatrix@R+0pc@C+0pc{
\underset{r}{\mathrm{homotopylimit}}~\mathrm{Spec}^\mathrm{CS}\widetilde{\Phi}^r_{\psi,\Gamma,X}/\mathrm{Fro}^\mathbb{Z},\underset{I}{\mathrm{homotopycolimit}}~\mathrm{Spec}^\mathrm{CS}\widetilde{\Phi}^I_{\psi,\Gamma,X}/\mathrm{Fro}^\mathbb{Z},	\\
}
\]
\[ 
\xymatrix@R+0pc@C+0pc{
\underset{r}{\mathrm{homotopylimit}}~\mathrm{Spec}^\mathrm{CS}\breve{\Phi}^r_{\psi,\Gamma,X}/\mathrm{Fro}^\mathbb{Z},\underset{I}{\mathrm{homotopycolimit}}~\breve{\Phi}^I_{\psi,\Gamma,X}/\mathrm{Fro}^\mathbb{Z},	\\
}
\]
\[ 
\xymatrix@R+0pc@C+0pc{
\underset{r}{\mathrm{homotopylimit}}~\mathrm{Spec}^\mathrm{CS}{\Phi}^r_{\psi,\Gamma,X}/\mathrm{Fro}^\mathbb{Z},\underset{I}{\mathrm{homotopycolimit}}~\mathrm{Spec}^\mathrm{CS}{\Phi}^I_{\psi,\Gamma,X}/\mathrm{Fro}^\mathbb{Z}.	
}
\]

\end{definition}

\begin{proposition}
There is a well-defined functor from the $\infty$-category 
\begin{align}
\mathrm{Quasicoherentpresheaves,Perfectcomplex,Condensed}_{*}	
\end{align}
where $*$ is one of the following spaces:
\begin{align}
&\mathrm{Spec}^\mathrm{CS}\widetilde{\Phi}_{\psi,\Gamma,X}/\mathrm{Fro}^\mathbb{Z},	\\
\end{align}
\begin{align}
&\mathrm{Spec}^\mathrm{CS}\breve{\Phi}_{\psi,\Gamma,X}/\mathrm{Fro}^\mathbb{Z},	\\
\end{align}
\begin{align}
&\mathrm{Spec}^\mathrm{CS}{\Phi}_{\psi,\Gamma,X}/\mathrm{Fro}^\mathbb{Z},	
\end{align}
to the $\infty$-category of $\mathrm{Fro}$-equivariant quasicoherent presheaves over similar spaces above correspondingly without the $\mathrm{Fro}$-quotients, and to the $\infty$-category of $\mathrm{Fro}$-equivariant quasicoherent modules over global sections of the structure $\infty$-sheaves of the similar spaces above correspondingly without the $\mathrm{Fro}$-quotients. Here for those space without notation related to the radius and the corresponding interval we consider the total unions $\bigcap_r,\bigcup_I$ in order to achieve the whole spaces to achieve the analogues of the corresponding FF curves from \cite{10KL1}, \cite{10KL2}, \cite{10FF} for
\[
\xymatrix@R+0pc@C+0pc{
\underset{r}{\mathrm{homotopylimit}}~\mathrm{Spec}^\mathrm{CS}\widetilde{\Phi}^r_{\psi,\Gamma,X},\underset{I}{\mathrm{homotopycolimit}}~\mathrm{Spec}^\mathrm{CS}\widetilde{\Phi}^I_{\psi,\Gamma,X},	\\
}
\]
\[
\xymatrix@R+0pc@C+0pc{
\underset{r}{\mathrm{homotopylimit}}~\mathrm{Spec}^\mathrm{CS}\breve{\Phi}^r_{\psi,\Gamma,X},\underset{I}{\mathrm{homotopycolimit}}~\mathrm{Spec}^\mathrm{CS}\breve{\Phi}^I_{\psi,\Gamma,X},	\\
}
\]
\[
\xymatrix@R+0pc@C+0pc{
\underset{r}{\mathrm{homotopylimit}}~\mathrm{Spec}^\mathrm{CS}{\Phi}^r_{\psi,\Gamma,X},\underset{I}{\mathrm{homotopycolimit}}~\mathrm{Spec}^\mathrm{CS}{\Phi}^I_{\psi,\Gamma,X}.	
}
\]
\[ 
\xymatrix@R+0pc@C+0pc{
\underset{r}{\mathrm{homotopylimit}}~\mathrm{Spec}^\mathrm{CS}\widetilde{\Phi}^r_{\psi,\Gamma,X}/\mathrm{Fro}^\mathbb{Z},\underset{I}{\mathrm{homotopycolimit}}~\mathrm{Spec}^\mathrm{CS}\widetilde{\Phi}^I_{\psi,\Gamma,X}/\mathrm{Fro}^\mathbb{Z},	\\
}
\]
\[ 
\xymatrix@R+0pc@C+0pc{
\underset{r}{\mathrm{homotopylimit}}~\mathrm{Spec}^\mathrm{CS}\breve{\Phi}^r_{\psi,\Gamma,X}/\mathrm{Fro}^\mathbb{Z},\underset{I}{\mathrm{homotopycolimit}}~\breve{\Phi}^I_{\psi,\Gamma,X}/\mathrm{Fro}^\mathbb{Z},	\\
}
\]
\[ 
\xymatrix@R+0pc@C+0pc{
\underset{r}{\mathrm{homotopylimit}}~\mathrm{Spec}^\mathrm{CS}{\Phi}^r_{\psi,\Gamma,X}/\mathrm{Fro}^\mathbb{Z},\underset{I}{\mathrm{homotopycolimit}}~\mathrm{Spec}^\mathrm{CS}{\Phi}^I_{\psi,\Gamma,X}/\mathrm{Fro}^\mathbb{Z}.	
}
\]	
In this situation we will have the target category being family parametrized by $r$ or $I$ in compatible glueing sense as in \cite[Definition 5.4.10]{10KL2}. In this situation for modules parametrized by the intervals we have the equivalence of $\infty$-categories by using \cite[Proposition 12.18]{10CS2}. Here the corresponding quasicoherent Frobenius modules are defined to be the corresponding homotopy colimits and limits of Frobenius modules:
\begin{align}
\underset{r}{\mathrm{homotopycolimit}}~M_r,\\
\underset{I}{\mathrm{homotopylimit}}~M_I,	
\end{align}
where each $M_r$ is a Frobenius-equivariant module over the period ring with respect to some radius $r$ while each $M_I$ is a Frobenius-equivariant module over the period ring with respect to some interval $I$.\\
\end{proposition}

\begin{proposition}
Similar proposition holds for 
\begin{align}
\mathrm{Quasicoherentsheaves,Perfectcomplex,IndBanach}_{*}.	
\end{align}	
\end{proposition}

\subsubsection{Frobenius Quasicoherent Modules II: Deformation in Preadic Spaces}

\begin{definition}
Let $\psi$ be a toric tower over $\mathbb{Q}_p$ as in \cite[Chapter 7]{10KL2} with base $\mathbb{Q}_p\left<X_1^{\pm 1},...,X_k^{\pm 1}\right>$. Then from \cite{10KL1} and \cite[Definition 5.2.1]{10KL2} we have the following class of Kedlaya-Liu rings (with the following replacement: $\Delta$ stands for $A$, $\nabla$ stands for $B$, while $\Phi$ stands for $C$) by taking product in the sense of self $\Gamma$-th power:

\[
\xymatrix@R+0pc@C+0pc{
\widetilde{\Delta}_{\psi,\Gamma},\widetilde{\nabla}_{\psi,\Gamma},\widetilde{\Phi}_{\psi,\Gamma},\widetilde{\Delta}^+_{\psi,\Gamma},\widetilde{\nabla}^+_{\psi,\Gamma},\widetilde{\Delta}^\dagger_{\psi,\Gamma},\widetilde{\nabla}^\dagger_{\psi,\Gamma},\widetilde{\Phi}^r_{\psi,\Gamma},\widetilde{\Phi}^I_{\psi,\Gamma}, 
}
\]

\[
\xymatrix@R+0pc@C+0pc{
\breve{\Delta}_{\psi,\Gamma},\breve{\nabla}_{\psi,\Gamma},\breve{\Phi}_{\psi,\Gamma},\breve{\Delta}^+_{\psi,\Gamma},\breve{\nabla}^+_{\psi,\Gamma},\breve{\Delta}^\dagger_{\psi,\Gamma},\breve{\nabla}^\dagger_{\psi,\Gamma},\breve{\Phi}^r_{\psi,\Gamma},\breve{\Phi}^I_{\psi,\Gamma},	
}
\]

\[
\xymatrix@R+0pc@C+0pc{
{\Delta}_{\psi,\Gamma},{\nabla}_{\psi,\Gamma},{\Phi}_{\psi,\Gamma},{\Delta}^+_{\psi,\Gamma},{\nabla}^+_{\psi,\Gamma},{\Delta}^\dagger_{\psi,\Gamma},{\nabla}^\dagger_{\psi,\Gamma},{\Phi}^r_{\psi,\Gamma},{\Phi}^I_{\psi,\Gamma}.	
}
\]
We now consider $\circ$ being deforming preadic space over $\mathbb{Q}_p$. Taking the product we have:
\[
\xymatrix@R+0pc@C+0pc{
\widetilde{\Phi}_{\psi,\Gamma,\circ},\widetilde{\Phi}^r_{\psi,\Gamma,\circ},\widetilde{\Phi}^I_{\psi,\Gamma,\circ},	
}
\]
\[
\xymatrix@R+0pc@C+0pc{
\breve{\Phi}_{\psi,\Gamma,\circ},\breve{\Phi}^r_{\psi,\Gamma,\circ},\breve{\Phi}^I_{\psi,\Gamma,\circ},	
}
\]
\[
\xymatrix@R+0pc@C+0pc{
{\Phi}_{\psi,\Gamma,\circ},{\Phi}^r_{\psi,\Gamma,\circ},{\Phi}^I_{\psi,\Gamma,\circ}.	
}
\]
They carry multi Frobenius action $\varphi_\Gamma$ and multi $\mathrm{Lie}_\Gamma:=\mathbb{Z}_p^{\times\Gamma}$ action. In our current situation after \cite{10CKZ} and \cite{10PZ} we consider the following $(\infty,1)$-categories of $(\infty,1)$-modules.\\
\end{definition}

\begin{definition}
First we consider the Bambozzi-Kremnizer spectrum $\mathrm{Spec}^\mathrm{BK}(*)$ attached to any of those in the above from \cite{10BK} by taking derived rational localization:
\begin{align}
&\mathrm{Spec}^\mathrm{BK}\widetilde{\Phi}_{\psi,\Gamma,\circ},\mathrm{Spec}^\mathrm{BK}\widetilde{\Phi}^r_{\psi,\Gamma,\circ},\mathrm{Spec}^\mathrm{BK}\widetilde{\Phi}^I_{\psi,\Gamma,\circ},	
\end{align}
\begin{align}
&\mathrm{Spec}^\mathrm{BK}\breve{\Phi}_{\psi,\Gamma,\circ},\mathrm{Spec}^\mathrm{BK}\breve{\Phi}^r_{\psi,\Gamma,\circ},\mathrm{Spec}^\mathrm{BK}\breve{\Phi}^I_{\psi,\Gamma,\circ},	
\end{align}
\begin{align}
&\mathrm{Spec}^\mathrm{BK}{\Phi}_{\psi,\Gamma,\circ},
\mathrm{Spec}^\mathrm{BK}{\Phi}^r_{\psi,\Gamma,\circ},\mathrm{Spec}^\mathrm{BK}{\Phi}^I_{\psi,\Gamma,\circ}.	
\end{align}

Then we take the corresponding quotients by using the corresponding Frobenius operators:
\begin{align}
&\mathrm{Spec}^\mathrm{BK}\widetilde{\Phi}_{\psi,\Gamma,\circ}/\mathrm{Fro}^\mathbb{Z},	\\
\end{align}
\begin{align}
&\mathrm{Spec}^\mathrm{BK}\breve{\Phi}_{\psi,\Gamma,\circ}/\mathrm{Fro}^\mathbb{Z},	\\
\end{align}
\begin{align}
&\mathrm{Spec}^\mathrm{BK}{\Phi}_{\psi,\Gamma,\circ}/\mathrm{Fro}^\mathbb{Z}.	
\end{align}
Here for those space without notation related to the radius and the corresponding interval we consider the total unions $\bigcap_r,\bigcup_I$ in order to achieve the whole spaces to achieve the analogues of the corresponding FF curves from \cite{10KL1}, \cite{10KL2}, \cite{10FF} for
\[
\xymatrix@R+0pc@C+0pc{
\underset{r}{\mathrm{homotopylimit}}~\mathrm{Spec}^\mathrm{BK}\widetilde{\Phi}^r_{\psi,\Gamma,\circ},\underset{I}{\mathrm{homotopycolimit}}~\mathrm{Spec}^\mathrm{BK}\widetilde{\Phi}^I_{\psi,\Gamma,\circ},	\\
}
\]
\[
\xymatrix@R+0pc@C+0pc{
\underset{r}{\mathrm{homotopylimit}}~\mathrm{Spec}^\mathrm{BK}\breve{\Phi}^r_{\psi,\Gamma,\circ},\underset{I}{\mathrm{homotopycolimit}}~\mathrm{Spec}^\mathrm{BK}\breve{\Phi}^I_{\psi,\Gamma,\circ},	\\
}
\]
\[
\xymatrix@R+0pc@C+0pc{
\underset{r}{\mathrm{homotopylimit}}~\mathrm{Spec}^\mathrm{BK}{\Phi}^r_{\psi,\Gamma,\circ},\underset{I}{\mathrm{homotopycolimit}}~\mathrm{Spec}^\mathrm{BK}{\Phi}^I_{\psi,\Gamma,\circ}.	
}
\]
\[  
\xymatrix@R+0pc@C+0pc{
\underset{r}{\mathrm{homotopylimit}}~\mathrm{Spec}^\mathrm{BK}\widetilde{\Phi}^r_{\psi,\Gamma,\circ}/\mathrm{Fro}^\mathbb{Z},\underset{I}{\mathrm{homotopycolimit}}~\mathrm{Spec}^\mathrm{BK}\widetilde{\Phi}^I_{\psi,\Gamma,\circ}/\mathrm{Fro}^\mathbb{Z},	\\
}
\]
\[ 
\xymatrix@R+0pc@C+0pc{
\underset{r}{\mathrm{homotopylimit}}~\mathrm{Spec}^\mathrm{BK}\breve{\Phi}^r_{\psi,\Gamma,\circ}/\mathrm{Fro}^\mathbb{Z},\underset{I}{\mathrm{homotopycolimit}}~\mathrm{Spec}^\mathrm{BK}\breve{\Phi}^I_{\psi,\Gamma,\circ}/\mathrm{Fro}^\mathbb{Z},	\\
}
\]
\[ 
\xymatrix@R+0pc@C+0pc{
\underset{r}{\mathrm{homotopylimit}}~\mathrm{Spec}^\mathrm{BK}{\Phi}^r_{\psi,\Gamma,\circ}/\mathrm{Fro}^\mathbb{Z},\underset{I}{\mathrm{homotopycolimit}}~\mathrm{Spec}^\mathrm{BK}{\Phi}^I_{\psi,\Gamma,\circ}/\mathrm{Fro}^\mathbb{Z}.	
}
\]

\end{definition}

\indent Meanwhile we have the corresponding Clausen-Scholze analytic stacks from \cite{10CS2}, therefore applying their construction we have:

\begin{definition}
Here we define the following products by using the solidified tensor product from \cite{10CS1} and \cite{10CS2}. Then we take solidified tensor product $\overset{\blacksquare}{\otimes}$ of any of the following
\[
\xymatrix@R+0pc@C+0pc{
\widetilde{\Delta}_{\psi,\Gamma},\widetilde{\nabla}_{\psi,\Gamma},\widetilde{\Phi}_{\psi,\Gamma},\widetilde{\Delta}^+_{\psi,\Gamma},\widetilde{\nabla}^+_{\psi,\Gamma},\widetilde{\Delta}^\dagger_{\psi,\Gamma},\widetilde{\nabla}^\dagger_{\psi,\Gamma},\widetilde{\Phi}^r_{\psi,\Gamma},\widetilde{\Phi}^I_{\psi,\Gamma}, 
}
\]

\[
\xymatrix@R+0pc@C+0pc{
\breve{\Delta}_{\psi,\Gamma},\breve{\nabla}_{\psi,\Gamma},\breve{\Phi}_{\psi,\Gamma},\breve{\Delta}^+_{\psi,\Gamma},\breve{\nabla}^+_{\psi,\Gamma},\breve{\Delta}^\dagger_{\psi,\Gamma},\breve{\nabla}^\dagger_{\psi,\Gamma},\breve{\Phi}^r_{\psi,\Gamma},\breve{\Phi}^I_{\psi,\Gamma},	
}
\]

\[
\xymatrix@R+0pc@C+0pc{
{\Delta}_{\psi,\Gamma},{\nabla}_{\psi,\Gamma},{\Phi}_{\psi,\Gamma},{\Delta}^+_{\psi,\Gamma},{\nabla}^+_{\psi,\Gamma},{\Delta}^\dagger_{\psi,\Gamma},{\nabla}^\dagger_{\psi,\Gamma},{\Phi}^r_{\psi,\Gamma},{\Phi}^I_{\psi,\Gamma},	
}
\]  	
with $\circ$. Then we have the notations:
\[
\xymatrix@R+0pc@C+0pc{
\widetilde{\Delta}_{\psi,\Gamma,\circ},\widetilde{\nabla}_{\psi,\Gamma,\circ},\widetilde{\Phi}_{\psi,\Gamma,\circ},\widetilde{\Delta}^+_{\psi,\Gamma,\circ},\widetilde{\nabla}^+_{\psi,\Gamma,\circ},\widetilde{\Delta}^\dagger_{\psi,\Gamma,\circ},\widetilde{\nabla}^\dagger_{\psi,\Gamma,\circ},\widetilde{\Phi}^r_{\psi,\Gamma,\circ},\widetilde{\Phi}^I_{\psi,\Gamma,\circ}, 
}
\]

\[
\xymatrix@R+0pc@C+0pc{
\breve{\Delta}_{\psi,\Gamma,\circ},\breve{\nabla}_{\psi,\Gamma,\circ},\breve{\Phi}_{\psi,\Gamma,\circ},\breve{\Delta}^+_{\psi,\Gamma,\circ},\breve{\nabla}^+_{\psi,\Gamma,\circ},\breve{\Delta}^\dagger_{\psi,\Gamma,\circ},\breve{\nabla}^\dagger_{\psi,\Gamma,\circ},\breve{\Phi}^r_{\psi,\Gamma,\circ},\breve{\Phi}^I_{\psi,\Gamma,\circ},	
}
\]

\[
\xymatrix@R+0pc@C+0pc{
{\Delta}_{\psi,\Gamma,\circ},{\nabla}_{\psi,\Gamma,\circ},{\Phi}_{\psi,\Gamma,\circ},{\Delta}^+_{\psi,\Gamma,\circ},{\nabla}^+_{\psi,\Gamma,\circ},{\Delta}^\dagger_{\psi,\Gamma,\circ},{\nabla}^\dagger_{\psi,\Gamma,\circ},{\Phi}^r_{\psi,\Gamma,\circ},{\Phi}^I_{\psi,\Gamma,\circ}.	
}
\]
\end{definition}

\begin{definition}
First we consider the Clausen-Scholze spectrum $\mathrm{Spec}^\mathrm{CS}(*)$ attached to any of those in the above from \cite{10CS2} by taking derived rational localization:
\begin{align}
\mathrm{Spec}^\mathrm{CS}\widetilde{\Delta}_{\psi,\Gamma,\circ},\mathrm{Spec}^\mathrm{CS}\widetilde{\nabla}_{\psi,\Gamma,\circ},\mathrm{Spec}^\mathrm{CS}\widetilde{\Phi}_{\psi,\Gamma,\circ},\mathrm{Spec}^\mathrm{CS}\widetilde{\Delta}^+_{\psi,\Gamma,\circ},\mathrm{Spec}^\mathrm{CS}\widetilde{\nabla}^+_{\psi,\Gamma,\circ},\\
\mathrm{Spec}^\mathrm{CS}\widetilde{\Delta}^\dagger_{\psi,\Gamma,\circ},\mathrm{Spec}^\mathrm{CS}\widetilde{\nabla}^\dagger_{\psi,\Gamma,\circ},\mathrm{Spec}^\mathrm{CS}\widetilde{\Phi}^r_{\psi,\Gamma,\circ},\mathrm{Spec}^\mathrm{CS}\widetilde{\Phi}^I_{\psi,\Gamma,\circ},	\\
\end{align}
\begin{align}
\mathrm{Spec}^\mathrm{CS}\breve{\Delta}_{\psi,\Gamma,\circ},\breve{\nabla}_{\psi,\Gamma,\circ},\mathrm{Spec}^\mathrm{CS}\breve{\Phi}_{\psi,\Gamma,\circ},\mathrm{Spec}^\mathrm{CS}\breve{\Delta}^+_{\psi,\Gamma,\circ},\mathrm{Spec}^\mathrm{CS}\breve{\nabla}^+_{\psi,\Gamma,\circ},\\
\mathrm{Spec}^\mathrm{CS}\breve{\Delta}^\dagger_{\psi,\Gamma,\circ},\mathrm{Spec}^\mathrm{CS}\breve{\nabla}^\dagger_{\psi,\Gamma,\circ},\mathrm{Spec}^\mathrm{CS}\breve{\Phi}^r_{\psi,\Gamma,\circ},\breve{\Phi}^I_{\psi,\Gamma,\circ},	\\
\end{align}
\begin{align}
\mathrm{Spec}^\mathrm{CS}{\Delta}_{\psi,\Gamma,\circ},\mathrm{Spec}^\mathrm{CS}{\nabla}_{\psi,\Gamma,\circ},\mathrm{Spec}^\mathrm{CS}{\Phi}_{\psi,\Gamma,\circ},\mathrm{Spec}^\mathrm{CS}{\Delta}^+_{\psi,\Gamma,\circ},\mathrm{Spec}^\mathrm{CS}{\nabla}^+_{\psi,\Gamma,\circ},\\
\mathrm{Spec}^\mathrm{CS}{\Delta}^\dagger_{\psi,\Gamma,\circ},\mathrm{Spec}^\mathrm{CS}{\nabla}^\dagger_{\psi,\Gamma,\circ},\mathrm{Spec}^\mathrm{CS}{\Phi}^r_{\psi,\Gamma,\circ},\mathrm{Spec}^\mathrm{CS}{\Phi}^I_{\psi,\Gamma,\circ}.	
\end{align}

Then we take the corresponding quotients by using the corresponding Frobenius operators:
\begin{align}
&\mathrm{Spec}^\mathrm{CS}\widetilde{\Delta}_{\psi,\Gamma,\circ}/\mathrm{Fro}^\mathbb{Z},\mathrm{Spec}^\mathrm{CS}\widetilde{\nabla}_{\psi,\Gamma,\circ}/\mathrm{Fro}^\mathbb{Z},\mathrm{Spec}^\mathrm{CS}\widetilde{\Phi}_{\psi,\Gamma,\circ}/\mathrm{Fro}^\mathbb{Z},\mathrm{Spec}^\mathrm{CS}\widetilde{\Delta}^+_{\psi,\Gamma,\circ}/\mathrm{Fro}^\mathbb{Z},\\
&\mathrm{Spec}^\mathrm{CS}\widetilde{\nabla}^+_{\psi,\Gamma,\circ}/\mathrm{Fro}^\mathbb{Z}, \mathrm{Spec}^\mathrm{CS}\widetilde{\Delta}^\dagger_{\psi,\Gamma,\circ}/\mathrm{Fro}^\mathbb{Z},\mathrm{Spec}^\mathrm{CS}\widetilde{\nabla}^\dagger_{\psi,\Gamma,\circ}/\mathrm{Fro}^\mathbb{Z},	\\
\end{align}
\begin{align}
&\mathrm{Spec}^\mathrm{CS}\breve{\Delta}_{\psi,\Gamma,\circ}/\mathrm{Fro}^\mathbb{Z},\breve{\nabla}_{\psi,\Gamma,\circ}/\mathrm{Fro}^\mathbb{Z},\mathrm{Spec}^\mathrm{CS}\breve{\Phi}_{\psi,\Gamma,\circ}/\mathrm{Fro}^\mathbb{Z},\mathrm{Spec}^\mathrm{CS}\breve{\Delta}^+_{\psi,\Gamma,\circ}/\mathrm{Fro}^\mathbb{Z},\\
&\mathrm{Spec}^\mathrm{CS}\breve{\nabla}^+_{\psi,\Gamma,\circ}/\mathrm{Fro}^\mathbb{Z}, \mathrm{Spec}^\mathrm{CS}\breve{\Delta}^\dagger_{\psi,\Gamma,\circ}/\mathrm{Fro}^\mathbb{Z},\mathrm{Spec}^\mathrm{CS}\breve{\nabla}^\dagger_{\psi,\Gamma,\circ}/\mathrm{Fro}^\mathbb{Z},	\\
\end{align}
\begin{align}
&\mathrm{Spec}^\mathrm{CS}{\Delta}_{\psi,\Gamma,\circ}/\mathrm{Fro}^\mathbb{Z},\mathrm{Spec}^\mathrm{CS}{\nabla}_{\psi,\Gamma,\circ}/\mathrm{Fro}^\mathbb{Z},\mathrm{Spec}^\mathrm{CS}{\Phi}_{\psi,\Gamma,\circ}/\mathrm{Fro}^\mathbb{Z},\mathrm{Spec}^\mathrm{CS}{\Delta}^+_{\psi,\Gamma,\circ}/\mathrm{Fro}^\mathbb{Z},\\
&\mathrm{Spec}^\mathrm{CS}{\nabla}^+_{\psi,\Gamma,\circ}/\mathrm{Fro}^\mathbb{Z}, \mathrm{Spec}^\mathrm{CS}{\Delta}^\dagger_{\psi,\Gamma,\circ}/\mathrm{Fro}^\mathbb{Z},\mathrm{Spec}^\mathrm{CS}{\nabla}^\dagger_{\psi,\Gamma,\circ}/\mathrm{Fro}^\mathbb{Z}.	
\end{align}
Here for those space with notations related to the radius and the corresponding interval we consider the total unions $\bigcap_r,\bigcup_I$ in order to achieve the whole spaces to achieve the analogues of the corresponding FF curves from \cite{10KL1}, \cite{10KL2}, \cite{10FF} for
\[
\xymatrix@R+0pc@C+0pc{
\underset{r}{\mathrm{homotopylimit}}~\mathrm{Spec}^\mathrm{CS}\widetilde{\Phi}^r_{\psi,\Gamma,\circ},\underset{I}{\mathrm{homotopycolimit}}~\mathrm{Spec}^\mathrm{CS}\widetilde{\Phi}^I_{\psi,\Gamma,\circ},	\\
}
\]
\[
\xymatrix@R+0pc@C+0pc{
\underset{r}{\mathrm{homotopylimit}}~\mathrm{Spec}^\mathrm{CS}\breve{\Phi}^r_{\psi,\Gamma,\circ},\underset{I}{\mathrm{homotopycolimit}}~\mathrm{Spec}^\mathrm{CS}\breve{\Phi}^I_{\psi,\Gamma,\circ},	\\
}
\]
\[
\xymatrix@R+0pc@C+0pc{
\underset{r}{\mathrm{homotopylimit}}~\mathrm{Spec}^\mathrm{CS}{\Phi}^r_{\psi,\Gamma,\circ},\underset{I}{\mathrm{homotopycolimit}}~\mathrm{Spec}^\mathrm{CS}{\Phi}^I_{\psi,\Gamma,\circ}.	
}
\]
\[ 
\xymatrix@R+0pc@C+0pc{
\underset{r}{\mathrm{homotopylimit}}~\mathrm{Spec}^\mathrm{CS}\widetilde{\Phi}^r_{\psi,\Gamma,\circ}/\mathrm{Fro}^\mathbb{Z},\underset{I}{\mathrm{homotopycolimit}}~\mathrm{Spec}^\mathrm{CS}\widetilde{\Phi}^I_{\psi,\Gamma,\circ}/\mathrm{Fro}^\mathbb{Z},	\\
}
\]
\[ 
\xymatrix@R+0pc@C+0pc{
\underset{r}{\mathrm{homotopylimit}}~\mathrm{Spec}^\mathrm{CS}\breve{\Phi}^r_{\psi,\Gamma,\circ}/\mathrm{Fro}^\mathbb{Z},\underset{I}{\mathrm{homotopycolimit}}~\breve{\Phi}^I_{\psi,\Gamma,\circ}/\mathrm{Fro}^\mathbb{Z},	\\
}
\]
\[ 
\xymatrix@R+0pc@C+0pc{
\underset{r}{\mathrm{homotopylimit}}~\mathrm{Spec}^\mathrm{CS}{\Phi}^r_{\psi,\Gamma,\circ}/\mathrm{Fro}^\mathbb{Z},\underset{I}{\mathrm{homotopycolimit}}~\mathrm{Spec}^\mathrm{CS}{\Phi}^I_{\psi,\Gamma,\circ}/\mathrm{Fro}^\mathbb{Z}.	
}
\]

\end{definition}

\

\begin{definition}
We then consider the corresponding quasipresheaves of the corresponding ind-Banach or monomorphic ind-Banach modules from \cite{10BBK}, \cite{10KKM}:
\begin{align}
\mathrm{Quasicoherentpresheaves,IndBanach}_{*}	
\end{align}
where $*$ is one of the following spaces:
\begin{align}
&\mathrm{Spec}^\mathrm{BK}\widetilde{\Phi}_{\psi,\Gamma,\circ}/\mathrm{Fro}^\mathbb{Z},	\\
\end{align}
\begin{align}
&\mathrm{Spec}^\mathrm{BK}\breve{\Phi}_{\psi,\Gamma,\circ}/\mathrm{Fro}^\mathbb{Z},	\\
\end{align}
\begin{align}
&\mathrm{Spec}^\mathrm{BK}{\Phi}_{\psi,\Gamma,\circ}/\mathrm{Fro}^\mathbb{Z}.	
\end{align}
Here for those space without notation related to the radius and the corresponding interval we consider the total unions $\bigcap_r,\bigcup_I$ in order to achieve the whole spaces to achieve the analogues of the corresponding FF curves from \cite{10KL1}, \cite{10KL2}, \cite{10FF} for
\[
\xymatrix@R+0pc@C+0pc{
\underset{r}{\mathrm{homotopylimit}}~\mathrm{Spec}^\mathrm{BK}\widetilde{\Phi}^r_{\psi,\Gamma,\circ},\underset{I}{\mathrm{homotopycolimit}}~\mathrm{Spec}^\mathrm{BK}\widetilde{\Phi}^I_{\psi,\Gamma,\circ},	\\
}
\]
\[
\xymatrix@R+0pc@C+0pc{
\underset{r}{\mathrm{homotopylimit}}~\mathrm{Spec}^\mathrm{BK}\breve{\Phi}^r_{\psi,\Gamma,\circ},\underset{I}{\mathrm{homotopycolimit}}~\mathrm{Spec}^\mathrm{BK}\breve{\Phi}^I_{\psi,\Gamma,\circ},	\\
}
\]
\[
\xymatrix@R+0pc@C+0pc{
\underset{r}{\mathrm{homotopylimit}}~\mathrm{Spec}^\mathrm{BK}{\Phi}^r_{\psi,\Gamma,\circ},\underset{I}{\mathrm{homotopycolimit}}~\mathrm{Spec}^\mathrm{BK}{\Phi}^I_{\psi,\Gamma,\circ}.	
}
\]
\[  
\xymatrix@R+0pc@C+0pc{
\underset{r}{\mathrm{homotopylimit}}~\mathrm{Spec}^\mathrm{BK}\widetilde{\Phi}^r_{\psi,\Gamma,\circ}/\mathrm{Fro}^\mathbb{Z},\underset{I}{\mathrm{homotopycolimit}}~\mathrm{Spec}^\mathrm{BK}\widetilde{\Phi}^I_{\psi,\Gamma,\circ}/\mathrm{Fro}^\mathbb{Z},	\\
}
\]
\[ 
\xymatrix@R+0pc@C+0pc{
\underset{r}{\mathrm{homotopylimit}}~\mathrm{Spec}^\mathrm{BK}\breve{\Phi}^r_{\psi,\Gamma,\circ}/\mathrm{Fro}^\mathbb{Z},\underset{I}{\mathrm{homotopycolimit}}~\mathrm{Spec}^\mathrm{BK}\breve{\Phi}^I_{\psi,\Gamma,\circ}/\mathrm{Fro}^\mathbb{Z},	\\
}
\]
\[ 
\xymatrix@R+0pc@C+0pc{
\underset{r}{\mathrm{homotopylimit}}~\mathrm{Spec}^\mathrm{BK}{\Phi}^r_{\psi,\Gamma,\circ}/\mathrm{Fro}^\mathbb{Z},\underset{I}{\mathrm{homotopycolimit}}~\mathrm{Spec}^\mathrm{BK}{\Phi}^I_{\psi,\Gamma,\circ}/\mathrm{Fro}^\mathbb{Z}.	
}
\]

\end{definition}

\begin{definition}
We then consider the corresponding quasisheaves of the corresponding condensed solid topological modules from \cite{10CS2}:
\begin{align}
\mathrm{Quasicoherentsheaves, Condensed}_{*}	
\end{align}
where $*$ is one of the following spaces:
\begin{align}
&\mathrm{Spec}^\mathrm{CS}\widetilde{\Delta}_{\psi,\Gamma,\circ}/\mathrm{Fro}^\mathbb{Z},\mathrm{Spec}^\mathrm{CS}\widetilde{\nabla}_{\psi,\Gamma,\circ}/\mathrm{Fro}^\mathbb{Z},\mathrm{Spec}^\mathrm{CS}\widetilde{\Phi}_{\psi,\Gamma,\circ}/\mathrm{Fro}^\mathbb{Z},\mathrm{Spec}^\mathrm{CS}\widetilde{\Delta}^+_{\psi,\Gamma,\circ}/\mathrm{Fro}^\mathbb{Z},\\
&\mathrm{Spec}^\mathrm{CS}\widetilde{\nabla}^+_{\psi,\Gamma,\circ}/\mathrm{Fro}^\mathbb{Z},\mathrm{Spec}^\mathrm{CS}\widetilde{\Delta}^\dagger_{\psi,\Gamma,\circ}/\mathrm{Fro}^\mathbb{Z},\mathrm{Spec}^\mathrm{CS}\widetilde{\nabla}^\dagger_{\psi,\Gamma,\circ}/\mathrm{Fro}^\mathbb{Z},	\\
\end{align}
\begin{align}
&\mathrm{Spec}^\mathrm{CS}\breve{\Delta}_{\psi,\Gamma,\circ}/\mathrm{Fro}^\mathbb{Z},\breve{\nabla}_{\psi,\Gamma,\circ}/\mathrm{Fro}^\mathbb{Z},\mathrm{Spec}^\mathrm{CS}\breve{\Phi}_{\psi,\Gamma,\circ}/\mathrm{Fro}^\mathbb{Z},\mathrm{Spec}^\mathrm{CS}\breve{\Delta}^+_{\psi,\Gamma,\circ}/\mathrm{Fro}^\mathbb{Z},\\
&\mathrm{Spec}^\mathrm{CS}\breve{\nabla}^+_{\psi,\Gamma,\circ}/\mathrm{Fro}^\mathbb{Z},\mathrm{Spec}^\mathrm{CS}\breve{\Delta}^\dagger_{\psi,\Gamma,\circ}/\mathrm{Fro}^\mathbb{Z},\mathrm{Spec}^\mathrm{CS}\breve{\nabla}^\dagger_{\psi,\Gamma,\circ}/\mathrm{Fro}^\mathbb{Z},	\\
\end{align}
\begin{align}
&\mathrm{Spec}^\mathrm{CS}{\Delta}_{\psi,\Gamma,\circ}/\mathrm{Fro}^\mathbb{Z},\mathrm{Spec}^\mathrm{CS}{\nabla}_{\psi,\Gamma,\circ}/\mathrm{Fro}^\mathbb{Z},\mathrm{Spec}^\mathrm{CS}{\Phi}_{\psi,\Gamma,\circ}/\mathrm{Fro}^\mathbb{Z},\mathrm{Spec}^\mathrm{CS}{\Delta}^+_{\psi,\Gamma,\circ}/\mathrm{Fro}^\mathbb{Z},\\
&\mathrm{Spec}^\mathrm{CS}{\nabla}^+_{\psi,\Gamma,\circ}/\mathrm{Fro}^\mathbb{Z}, \mathrm{Spec}^\mathrm{CS}{\Delta}^\dagger_{\psi,\Gamma,\circ}/\mathrm{Fro}^\mathbb{Z},\mathrm{Spec}^\mathrm{CS}{\nabla}^\dagger_{\psi,\Gamma,\circ}/\mathrm{Fro}^\mathbb{Z}.	
\end{align}
Here for those space with notations related to the radius and the corresponding interval we consider the total unions $\bigcap_r,\bigcup_I$ in order to achieve the whole spaces to achieve the analogues of the corresponding FF curves from \cite{10KL1}, \cite{10KL2}, \cite{10FF} for
\[
\xymatrix@R+0pc@C+0pc{
\underset{r}{\mathrm{homotopylimit}}~\mathrm{Spec}^\mathrm{CS}\widetilde{\Phi}^r_{\psi,\Gamma,\circ},\underset{I}{\mathrm{homotopycolimit}}~\mathrm{Spec}^\mathrm{CS}\widetilde{\Phi}^I_{\psi,\Gamma,\circ},	\\
}
\]
\[
\xymatrix@R+0pc@C+0pc{
\underset{r}{\mathrm{homotopylimit}}~\mathrm{Spec}^\mathrm{CS}\breve{\Phi}^r_{\psi,\Gamma,\circ},\underset{I}{\mathrm{homotopycolimit}}~\mathrm{Spec}^\mathrm{CS}\breve{\Phi}^I_{\psi,\Gamma,\circ},	\\
}
\]
\[
\xymatrix@R+0pc@C+0pc{
\underset{r}{\mathrm{homotopylimit}}~\mathrm{Spec}^\mathrm{CS}{\Phi}^r_{\psi,\Gamma,\circ},\underset{I}{\mathrm{homotopycolimit}}~\mathrm{Spec}^\mathrm{CS}{\Phi}^I_{\psi,\Gamma,\circ}.	
}
\]
\[ 
\xymatrix@R+0pc@C+0pc{
\underset{r}{\mathrm{homotopylimit}}~\mathrm{Spec}^\mathrm{CS}\widetilde{\Phi}^r_{\psi,\Gamma,\circ}/\mathrm{Fro}^\mathbb{Z},\underset{I}{\mathrm{homotopycolimit}}~\mathrm{Spec}^\mathrm{CS}\widetilde{\Phi}^I_{\psi,\Gamma,\circ}/\mathrm{Fro}^\mathbb{Z},	\\
}
\]
\[ 
\xymatrix@R+0pc@C+0pc{
\underset{r}{\mathrm{homotopylimit}}~\mathrm{Spec}^\mathrm{CS}\breve{\Phi}^r_{\psi,\Gamma,\circ}/\mathrm{Fro}^\mathbb{Z},\underset{I}{\mathrm{homotopycolimit}}~\breve{\Phi}^I_{\psi,\Gamma,\circ}/\mathrm{Fro}^\mathbb{Z},	\\
}
\]
\[ 
\xymatrix@R+0pc@C+0pc{
\underset{r}{\mathrm{homotopylimit}}~\mathrm{Spec}^\mathrm{CS}{\Phi}^r_{\psi,\Gamma,\circ}/\mathrm{Fro}^\mathbb{Z},\underset{I}{\mathrm{homotopycolimit}}~\mathrm{Spec}^\mathrm{CS}{\Phi}^I_{\psi,\Gamma,\circ}/\mathrm{Fro}^\mathbb{Z}.	
}
\]

\end{definition}

\

\begin{proposition}
There is a well-defined functor from the $\infty$-category 
\begin{align}
\mathrm{Quasicoherentpresheaves,Condensed}_{*}	
\end{align}
where $*$ is one of the following spaces:
\begin{align}
&\mathrm{Spec}^\mathrm{CS}\widetilde{\Phi}_{\psi,\Gamma,\circ}/\mathrm{Fro}^\mathbb{Z},	\\
\end{align}
\begin{align}
&\mathrm{Spec}^\mathrm{CS}\breve{\Phi}_{\psi,\Gamma,\circ}/\mathrm{Fro}^\mathbb{Z},	\\
\end{align}
\begin{align}
&\mathrm{Spec}^\mathrm{CS}{\Phi}_{\psi,\Gamma,\circ}/\mathrm{Fro}^\mathbb{Z},	
\end{align}
to the $\infty$-category of $\mathrm{Fro}$-equivariant quasicoherent presheaves over similar spaces above correspondingly without the $\mathrm{Fro}$-quotients, and to the $\infty$-category of $\mathrm{Fro}$-equivariant quasicoherent modules over global sections of the structure $\infty$-sheaves of the similar spaces above correspondingly without the $\mathrm{Fro}$-quotients. Here for those space without notation related to the radius and the corresponding interval we consider the total unions $\bigcap_r,\bigcup_I$ in order to achieve the whole spaces to achieve the analogues of the corresponding FF curves from \cite{10KL1}, \cite{10KL2}, \cite{10FF} for
\[
\xymatrix@R+0pc@C+0pc{
\underset{r}{\mathrm{homotopylimit}}~\mathrm{Spec}^\mathrm{CS}\widetilde{\Phi}^r_{\psi,\Gamma,\circ},\underset{I}{\mathrm{homotopycolimit}}~\mathrm{Spec}^\mathrm{CS}\widetilde{\Phi}^I_{\psi,\Gamma,\circ},	\\
}
\]
\[
\xymatrix@R+0pc@C+0pc{
\underset{r}{\mathrm{homotopylimit}}~\mathrm{Spec}^\mathrm{CS}\breve{\Phi}^r_{\psi,\Gamma,\circ},\underset{I}{\mathrm{homotopycolimit}}~\mathrm{Spec}^\mathrm{CS}\breve{\Phi}^I_{\psi,\Gamma,\circ},	\\
}
\]
\[
\xymatrix@R+0pc@C+0pc{
\underset{r}{\mathrm{homotopylimit}}~\mathrm{Spec}^\mathrm{CS}{\Phi}^r_{\psi,\Gamma,\circ},\underset{I}{\mathrm{homotopycolimit}}~\mathrm{Spec}^\mathrm{CS}{\Phi}^I_{\psi,\Gamma,\circ}.	
}
\]
\[ 
\xymatrix@R+0pc@C+0pc{
\underset{r}{\mathrm{homotopylimit}}~\mathrm{Spec}^\mathrm{CS}\widetilde{\Phi}^r_{\psi,\Gamma,\circ}/\mathrm{Fro}^\mathbb{Z},\underset{I}{\mathrm{homotopycolimit}}~\mathrm{Spec}^\mathrm{CS}\widetilde{\Phi}^I_{\psi,\Gamma,\circ}/\mathrm{Fro}^\mathbb{Z},	\\
}
\]
\[ 
\xymatrix@R+0pc@C+0pc{
\underset{r}{\mathrm{homotopylimit}}~\mathrm{Spec}^\mathrm{CS}\breve{\Phi}^r_{\psi,\Gamma,\circ}/\mathrm{Fro}^\mathbb{Z},\underset{I}{\mathrm{homotopycolimit}}~\breve{\Phi}^I_{\psi,\Gamma,\circ}/\mathrm{Fro}^\mathbb{Z},	\\
}
\]
\[ 
\xymatrix@R+0pc@C+0pc{
\underset{r}{\mathrm{homotopylimit}}~\mathrm{Spec}^\mathrm{CS}{\Phi}^r_{\psi,\Gamma,\circ}/\mathrm{Fro}^\mathbb{Z},\underset{I}{\mathrm{homotopycolimit}}~\mathrm{Spec}^\mathrm{CS}{\Phi}^I_{\psi,\Gamma,\circ}/\mathrm{Fro}^\mathbb{Z}.	
}
\]	
In this situation we will have the target category being family parametrized by $r$ or $I$ in compatible glueing sense as in \cite[Definition 5.4.10]{10KL2}. In this situation for modules parametrized by the intervals we have the equivalence of $\infty$-categories by using \cite[Proposition 13.8]{10CS2}. Here the corresponding quasicoherent Frobenius modules are defined to be the corresponding homotopy colimits and limits of Frobenius modules:
\begin{align}
\underset{r}{\mathrm{homotopycolimit}}~M_r,\\
\underset{I}{\mathrm{homotopylimit}}~M_I,	
\end{align}
where each $M_r$ is a Frobenius-equivariant module over the period ring with respect to some radius $r$ while each $M_I$ is a Frobenius-equivariant module over the period ring with respect to some interval $I$.\\
\end{proposition}

\begin{proposition}
Similar proposition holds for 
\begin{align}
\mathrm{Quasicoherentsheaves,IndBanach}_{*}.	
\end{align}	
\end{proposition}

\

\begin{definition}
We then consider the corresponding quasipresheaves of perfect complexes the corresponding ind-Banach or monomorphic ind-Banach modules from \cite{10BBK}, \cite{10KKM}:
\begin{align}
\mathrm{Quasicoherentpresheaves,Perfectcomplex,IndBanach}_{*}	
\end{align}
where $*$ is one of the following spaces:
\begin{align}
&\mathrm{Spec}^\mathrm{BK}\widetilde{\Phi}_{\psi,\Gamma,\circ}/\mathrm{Fro}^\mathbb{Z},	\\
\end{align}
\begin{align}
&\mathrm{Spec}^\mathrm{BK}\breve{\Phi}_{\psi,\Gamma,\circ}/\mathrm{Fro}^\mathbb{Z},	\\
\end{align}
\begin{align}
&\mathrm{Spec}^\mathrm{BK}{\Phi}_{\psi,\Gamma,\circ}/\mathrm{Fro}^\mathbb{Z}.	
\end{align}
Here for those space without notation related to the radius and the corresponding interval we consider the total unions $\bigcap_r,\bigcup_I$ in order to achieve the whole spaces to achieve the analogues of the corresponding FF curves from \cite{10KL1}, \cite{10KL2}, \cite{10FF} for
\[
\xymatrix@R+0pc@C+0pc{
\underset{r}{\mathrm{homotopylimit}}~\mathrm{Spec}^\mathrm{BK}\widetilde{\Phi}^r_{\psi,\Gamma,\circ},\underset{I}{\mathrm{homotopycolimit}}~\mathrm{Spec}^\mathrm{BK}\widetilde{\Phi}^I_{\psi,\Gamma,\circ},	\\
}
\]
\[
\xymatrix@R+0pc@C+0pc{
\underset{r}{\mathrm{homotopylimit}}~\mathrm{Spec}^\mathrm{BK}\breve{\Phi}^r_{\psi,\Gamma,\circ},\underset{I}{\mathrm{homotopycolimit}}~\mathrm{Spec}^\mathrm{BK}\breve{\Phi}^I_{\psi,\Gamma,\circ},	\\
}
\]
\[
\xymatrix@R+0pc@C+0pc{
\underset{r}{\mathrm{homotopylimit}}~\mathrm{Spec}^\mathrm{BK}{\Phi}^r_{\psi,\Gamma,\circ},\underset{I}{\mathrm{homotopycolimit}}~\mathrm{Spec}^\mathrm{BK}{\Phi}^I_{\psi,\Gamma,\circ}.	
}
\]
\[  
\xymatrix@R+0pc@C+0pc{
\underset{r}{\mathrm{homotopylimit}}~\mathrm{Spec}^\mathrm{BK}\widetilde{\Phi}^r_{\psi,\Gamma,\circ}/\mathrm{Fro}^\mathbb{Z},\underset{I}{\mathrm{homotopycolimit}}~\mathrm{Spec}^\mathrm{BK}\widetilde{\Phi}^I_{\psi,\Gamma,\circ}/\mathrm{Fro}^\mathbb{Z},	\\
}
\]
\[ 
\xymatrix@R+0pc@C+0pc{
\underset{r}{\mathrm{homotopylimit}}~\mathrm{Spec}^\mathrm{BK}\breve{\Phi}^r_{\psi,\Gamma,\circ}/\mathrm{Fro}^\mathbb{Z},\underset{I}{\mathrm{homotopycolimit}}~\mathrm{Spec}^\mathrm{BK}\breve{\Phi}^I_{\psi,\Gamma,\circ}/\mathrm{Fro}^\mathbb{Z},	\\
}
\]
\[ 
\xymatrix@R+0pc@C+0pc{
\underset{r}{\mathrm{homotopylimit}}~\mathrm{Spec}^\mathrm{BK}{\Phi}^r_{\psi,\Gamma,\circ}/\mathrm{Fro}^\mathbb{Z},\underset{I}{\mathrm{homotopycolimit}}~\mathrm{Spec}^\mathrm{BK}{\Phi}^I_{\psi,\Gamma,\circ}/\mathrm{Fro}^\mathbb{Z}.	
}
\]

\end{definition}

\begin{definition}
We then consider the corresponding quasisheaves of perfect complexes of the corresponding condensed solid topological modules from \cite{10CS2}:
\begin{align}
\mathrm{Quasicoherentsheaves, Perfectcomplex, Condensed}_{*}	
\end{align}
where $*$ is one of the following spaces:
\begin{align}
&\mathrm{Spec}^\mathrm{CS}\widetilde{\Delta}_{\psi,\Gamma,\circ}/\mathrm{Fro}^\mathbb{Z},\mathrm{Spec}^\mathrm{CS}\widetilde{\nabla}_{\psi,\Gamma,\circ}/\mathrm{Fro}^\mathbb{Z},\mathrm{Spec}^\mathrm{CS}\widetilde{\Phi}_{\psi,\Gamma,\circ}/\mathrm{Fro}^\mathbb{Z},\mathrm{Spec}^\mathrm{CS}\widetilde{\Delta}^+_{\psi,\Gamma,\circ}/\mathrm{Fro}^\mathbb{Z},\\
&\mathrm{Spec}^\mathrm{CS}\widetilde{\nabla}^+_{\psi,\Gamma,\circ}/\mathrm{Fro}^\mathbb{Z},\mathrm{Spec}^\mathrm{CS}\widetilde{\Delta}^\dagger_{\psi,\Gamma,\circ}/\mathrm{Fro}^\mathbb{Z},\mathrm{Spec}^\mathrm{CS}\widetilde{\nabla}^\dagger_{\psi,\Gamma,\circ}/\mathrm{Fro}^\mathbb{Z},	\\
\end{align}
\begin{align}
&\mathrm{Spec}^\mathrm{CS}\breve{\Delta}_{\psi,\Gamma,\circ}/\mathrm{Fro}^\mathbb{Z},\breve{\nabla}_{\psi,\Gamma,\circ}/\mathrm{Fro}^\mathbb{Z},\mathrm{Spec}^\mathrm{CS}\breve{\Phi}_{\psi,\Gamma,\circ}/\mathrm{Fro}^\mathbb{Z},\mathrm{Spec}^\mathrm{CS}\breve{\Delta}^+_{\psi,\Gamma,\circ}/\mathrm{Fro}^\mathbb{Z},\\
&\mathrm{Spec}^\mathrm{CS}\breve{\nabla}^+_{\psi,\Gamma,\circ}/\mathrm{Fro}^\mathbb{Z},\mathrm{Spec}^\mathrm{CS}\breve{\Delta}^\dagger_{\psi,\Gamma,\circ}/\mathrm{Fro}^\mathbb{Z},\mathrm{Spec}^\mathrm{CS}\breve{\nabla}^\dagger_{\psi,\Gamma,\circ}/\mathrm{Fro}^\mathbb{Z},	\\
\end{align}
\begin{align}
&\mathrm{Spec}^\mathrm{CS}{\Delta}_{\psi,\Gamma,\circ}/\mathrm{Fro}^\mathbb{Z},\mathrm{Spec}^\mathrm{CS}{\nabla}_{\psi,\Gamma,\circ}/\mathrm{Fro}^\mathbb{Z},\mathrm{Spec}^\mathrm{CS}{\Phi}_{\psi,\Gamma,\circ}/\mathrm{Fro}^\mathbb{Z},\mathrm{Spec}^\mathrm{CS}{\Delta}^+_{\psi,\Gamma,\circ}/\mathrm{Fro}^\mathbb{Z},\\
&\mathrm{Spec}^\mathrm{CS}{\nabla}^+_{\psi,\Gamma,\circ}/\mathrm{Fro}^\mathbb{Z}, \mathrm{Spec}^\mathrm{CS}{\Delta}^\dagger_{\psi,\Gamma,\circ}/\mathrm{Fro}^\mathbb{Z},\mathrm{Spec}^\mathrm{CS}{\nabla}^\dagger_{\psi,\Gamma,\circ}/\mathrm{Fro}^\mathbb{Z}.	
\end{align}
Here for those space with notations related to the radius and the corresponding interval we consider the total unions $\bigcap_r,\bigcup_I$ in order to achieve the whole spaces to achieve the analogues of the corresponding FF curves from \cite{10KL1}, \cite{10KL2}, \cite{10FF} for
\[
\xymatrix@R+0pc@C+0pc{
\underset{r}{\mathrm{homotopylimit}}~\mathrm{Spec}^\mathrm{CS}\widetilde{\Phi}^r_{\psi,\Gamma,\circ},\underset{I}{\mathrm{homotopycolimit}}~\mathrm{Spec}^\mathrm{CS}\widetilde{\Phi}^I_{\psi,\Gamma,\circ},	\\
}
\]
\[
\xymatrix@R+0pc@C+0pc{
\underset{r}{\mathrm{homotopylimit}}~\mathrm{Spec}^\mathrm{CS}\breve{\Phi}^r_{\psi,\Gamma,\circ},\underset{I}{\mathrm{homotopycolimit}}~\mathrm{Spec}^\mathrm{CS}\breve{\Phi}^I_{\psi,\Gamma,\circ},	\\
}
\]
\[
\xymatrix@R+0pc@C+0pc{
\underset{r}{\mathrm{homotopylimit}}~\mathrm{Spec}^\mathrm{CS}{\Phi}^r_{\psi,\Gamma,\circ},\underset{I}{\mathrm{homotopycolimit}}~\mathrm{Spec}^\mathrm{CS}{\Phi}^I_{\psi,\Gamma,\circ}.	
}
\]
\[ 
\xymatrix@R+0pc@C+0pc{
\underset{r}{\mathrm{homotopylimit}}~\mathrm{Spec}^\mathrm{CS}\widetilde{\Phi}^r_{\psi,\Gamma,\circ}/\mathrm{Fro}^\mathbb{Z},\underset{I}{\mathrm{homotopycolimit}}~\mathrm{Spec}^\mathrm{CS}\widetilde{\Phi}^I_{\psi,\Gamma,\circ}/\mathrm{Fro}^\mathbb{Z},	\\
}
\]
\[ 
\xymatrix@R+0pc@C+0pc{
\underset{r}{\mathrm{homotopylimit}}~\mathrm{Spec}^\mathrm{CS}\breve{\Phi}^r_{\psi,\Gamma,\circ}/\mathrm{Fro}^\mathbb{Z},\underset{I}{\mathrm{homotopycolimit}}~\breve{\Phi}^I_{\psi,\Gamma,\circ}/\mathrm{Fro}^\mathbb{Z},	\\
}
\]
\[ 
\xymatrix@R+0pc@C+0pc{
\underset{r}{\mathrm{homotopylimit}}~\mathrm{Spec}^\mathrm{CS}{\Phi}^r_{\psi,\Gamma,\circ}/\mathrm{Fro}^\mathbb{Z},\underset{I}{\mathrm{homotopycolimit}}~\mathrm{Spec}^\mathrm{CS}{\Phi}^I_{\psi,\Gamma,\circ}/\mathrm{Fro}^\mathbb{Z}.	
}
\]

\end{definition}

\begin{proposition}
There is a well-defined functor from the $\infty$-category 
\begin{align}
\mathrm{Quasicoherentpresheaves,Perfectcomplex,Condensed}_{*}	
\end{align}
where $*$ is one of the following spaces:
\begin{align}
&\mathrm{Spec}^\mathrm{CS}\widetilde{\Phi}_{\psi,\Gamma,\circ}/\mathrm{Fro}^\mathbb{Z},	\\
\end{align}
\begin{align}
&\mathrm{Spec}^\mathrm{CS}\breve{\Phi}_{\psi,\Gamma,\circ}/\mathrm{Fro}^\mathbb{Z},	\\
\end{align}
\begin{align}
&\mathrm{Spec}^\mathrm{CS}{\Phi}_{\psi,\Gamma,\circ}/\mathrm{Fro}^\mathbb{Z},	
\end{align}
to the $\infty$-category of $\mathrm{Fro}$-equivariant quasicoherent presheaves over similar spaces above correspondingly without the $\mathrm{Fro}$-quotients, and to the $\infty$-category of $\mathrm{Fro}$-equivariant quasicoherent modules over global sections of the structure $\infty$-sheaves of the similar spaces above correspondingly without the $\mathrm{Fro}$-quotients. Here for those space without notation related to the radius and the corresponding interval we consider the total unions $\bigcap_r,\bigcup_I$ in order to achieve the whole spaces to achieve the analogues of the corresponding FF curves from \cite{10KL1}, \cite{10KL2}, \cite{10FF} for
\[
\xymatrix@R+0pc@C+0pc{
\underset{r}{\mathrm{homotopylimit}}~\mathrm{Spec}^\mathrm{CS}\widetilde{\Phi}^r_{\psi,\Gamma,\circ},\underset{I}{\mathrm{homotopycolimit}}~\mathrm{Spec}^\mathrm{CS}\widetilde{\Phi}^I_{\psi,\Gamma,\circ},	\\
}
\]
\[
\xymatrix@R+0pc@C+0pc{
\underset{r}{\mathrm{homotopylimit}}~\mathrm{Spec}^\mathrm{CS}\breve{\Phi}^r_{\psi,\Gamma,\circ},\underset{I}{\mathrm{homotopycolimit}}~\mathrm{Spec}^\mathrm{CS}\breve{\Phi}^I_{\psi,\Gamma,\circ},	\\
}
\]
\[
\xymatrix@R+0pc@C+0pc{
\underset{r}{\mathrm{homotopylimit}}~\mathrm{Spec}^\mathrm{CS}{\Phi}^r_{\psi,\Gamma,\circ},\underset{I}{\mathrm{homotopycolimit}}~\mathrm{Spec}^\mathrm{CS}{\Phi}^I_{\psi,\Gamma,\circ}.	
}
\]
\[ 
\xymatrix@R+0pc@C+0pc{
\underset{r}{\mathrm{homotopylimit}}~\mathrm{Spec}^\mathrm{CS}\widetilde{\Phi}^r_{\psi,\Gamma,\circ}/\mathrm{Fro}^\mathbb{Z},\underset{I}{\mathrm{homotopycolimit}}~\mathrm{Spec}^\mathrm{CS}\widetilde{\Phi}^I_{\psi,\Gamma,\circ}/\mathrm{Fro}^\mathbb{Z},	\\
}
\]
\[ 
\xymatrix@R+0pc@C+0pc{
\underset{r}{\mathrm{homotopylimit}}~\mathrm{Spec}^\mathrm{CS}\breve{\Phi}^r_{\psi,\Gamma,\circ}/\mathrm{Fro}^\mathbb{Z},\underset{I}{\mathrm{homotopycolimit}}~\breve{\Phi}^I_{\psi,\Gamma,\circ}/\mathrm{Fro}^\mathbb{Z},	\\
}
\]
\[ 
\xymatrix@R+0pc@C+0pc{
\underset{r}{\mathrm{homotopylimit}}~\mathrm{Spec}^\mathrm{CS}{\Phi}^r_{\psi,\Gamma,\circ}/\mathrm{Fro}^\mathbb{Z},\underset{I}{\mathrm{homotopycolimit}}~\mathrm{Spec}^\mathrm{CS}{\Phi}^I_{\psi,\Gamma,\circ}/\mathrm{Fro}^\mathbb{Z}.	
}
\]	
In this situation we will have the target category being family parametrized by $r$ or $I$ in compatible glueing sense as in \cite[Definition 5.4.10]{10KL2}. In this situation for modules parametrized by the intervals we have the equivalence of $\infty$-categories by using \cite[Proposition 12.18]{10CS2}. Here the corresponding quasicoherent Frobenius modules are defined to be the corresponding homotopy colimits and limits of Frobenius modules:
\begin{align}
\underset{r}{\mathrm{homotopycolimit}}~M_r,\\
\underset{I}{\mathrm{homotopylimit}}~M_I,	
\end{align}
where each $M_r$ is a Frobenius-equivariant module over the period ring with respect to some radius $r$ while each $M_I$ is a Frobenius-equivariant module over the period ring with respect to some interval $I$.\\
\end{proposition}

\begin{proposition}
Similar proposition holds for 
\begin{align}
\mathrm{Quasicoherentsheaves,Perfectcomplex,IndBanach}_{*}.	
\end{align}	
\end{proposition}

\subsubsection{Frobenius Quasicoherent Modules III: Deformation in $(\infty,1)$-Ind-Preadic Spaces}

\begin{definition}
Let $\psi$ be a toric tower over $\mathbb{Q}_p$ as in \cite[Chapter 7]{10KL2} with base $\mathbb{Q}_p\left<X_1^{\pm 1},...,X_k^{\pm 1}\right>$. Then from \cite{10KL1} and \cite[Definition 5.2.1]{10KL2} we have the following class of Kedlaya-Liu rings (with the following replacement: $\Delta$ stands for $A$, $\nabla$ stands for $B$, while $\Phi$ stands for $C$) by taking product in the sense of self $\Gamma$-th power:

\[
\xymatrix@R+0pc@C+0pc{
\widetilde{\Delta}_{\psi,\Gamma},\widetilde{\nabla}_{\psi,\Gamma},\widetilde{\Phi}_{\psi,\Gamma},\widetilde{\Delta}^+_{\psi,\Gamma},\widetilde{\nabla}^+_{\psi,\Gamma},\widetilde{\Delta}^\dagger_{\psi,\Gamma},\widetilde{\nabla}^\dagger_{\psi,\Gamma},\widetilde{\Phi}^r_{\psi,\Gamma},\widetilde{\Phi}^I_{\psi,\Gamma}, 
}
\]

\[
\xymatrix@R+0pc@C+0pc{
\breve{\Delta}_{\psi,\Gamma},\breve{\nabla}_{\psi,\Gamma},\breve{\Phi}_{\psi,\Gamma},\breve{\Delta}^+_{\psi,\Gamma},\breve{\nabla}^+_{\psi,\Gamma},\breve{\Delta}^\dagger_{\psi,\Gamma},\breve{\nabla}^\dagger_{\psi,\Gamma},\breve{\Phi}^r_{\psi,\Gamma},\breve{\Phi}^I_{\psi,\Gamma},	
}
\]

\[
\xymatrix@R+0pc@C+0pc{
{\Delta}_{\psi,\Gamma},{\nabla}_{\psi,\Gamma},{\Phi}_{\psi,\Gamma},{\Delta}^+_{\psi,\Gamma},{\nabla}^+_{\psi,\Gamma},{\Delta}^\dagger_{\psi,\Gamma},{\nabla}^\dagger_{\psi,\Gamma},{\Phi}^r_{\psi,\Gamma},{\Phi}^I_{\psi,\Gamma}.	
}
\]
We now consider the following rings with $X_\square$ being a homotopy colimit
\begin{align}
 \underset{i}{\mathrm{homotopycolimit}}X_{\square_i}
 \end{align}
in the $\infty$-categories of analytic stacks from \cite{10BBBK} and \cite{10CS2}.   
Taking the product we have:
\[
\xymatrix@R+0pc@C+0pc{
\widetilde{\Phi}_{\psi,\Gamma,X_\square},\widetilde{\Phi}^r_{\psi,\Gamma,X_\square},\widetilde{\Phi}^I_{\psi,\Gamma,X_\square},	
}
\]
\[
\xymatrix@R+0pc@C+0pc{
\breve{\Phi}_{\psi,\Gamma,X_\square},\breve{\Phi}^r_{\psi,\Gamma,X_\square},\breve{\Phi}^I_{\psi,\Gamma,X_\square},	
}
\]
\[
\xymatrix@R+0pc@C+0pc{
{\Phi}_{\psi,\Gamma,X_\square},{\Phi}^r_{\psi,\Gamma,X_\square},{\Phi}^I_{\psi,\Gamma,X_\square}.	
}
\]
They carry multi Frobenius action $\varphi_\Gamma$ and multi $\mathrm{Lie}_\Gamma:=\mathbb{Z}_p^{\times\Gamma}$ action. In our current situation after \cite{10CKZ} and \cite{10PZ} we consider the following $(\infty,1)$-categories of $(\infty,1)$-modules.\\
\end{definition}

\begin{definition}
First we consider the Bambozzi-Kremnizer spectrum $\mathrm{Spec}^\mathrm{BK}(*)$ attached to any of those in the above from \cite{10BK} by taking derived rational localization:
\begin{align}
&\mathrm{Spec}^\mathrm{BK}\widetilde{\Phi}_{\psi,\Gamma,X_\square},\mathrm{Spec}^\mathrm{BK}\widetilde{\Phi}^r_{\psi,\Gamma,X_\square},\mathrm{Spec}^\mathrm{BK}\widetilde{\Phi}^I_{\psi,\Gamma,X_\square},	
\end{align}
\begin{align}
&\mathrm{Spec}^\mathrm{BK}\breve{\Phi}_{\psi,\Gamma,X_\square},\mathrm{Spec}^\mathrm{BK}\breve{\Phi}^r_{\psi,\Gamma,X_\square},\mathrm{Spec}^\mathrm{BK}\breve{\Phi}^I_{\psi,\Gamma,X_\square},	
\end{align}
\begin{align}
&\mathrm{Spec}^\mathrm{BK}{\Phi}_{\psi,\Gamma,X_\square},
\mathrm{Spec}^\mathrm{BK}{\Phi}^r_{\psi,\Gamma,X_\square},\mathrm{Spec}^\mathrm{BK}{\Phi}^I_{\psi,\Gamma,X_\square}.	
\end{align}

Then we take the corresponding quotients by using the corresponding Frobenius operators:
\begin{align}
&\mathrm{Spec}^\mathrm{BK}\widetilde{\Phi}_{\psi,\Gamma,X_\square}/\mathrm{Fro}^\mathbb{Z},	\\
\end{align}
\begin{align}
&\mathrm{Spec}^\mathrm{BK}\breve{\Phi}_{\psi,\Gamma,X_\square}/\mathrm{Fro}^\mathbb{Z},	\\
\end{align}
\begin{align}
&\mathrm{Spec}^\mathrm{BK}{\Phi}_{\psi,\Gamma,X_\square}/\mathrm{Fro}^\mathbb{Z}.	
\end{align}
Here for those space without notation related to the radius and the corresponding interval we consider the total unions $\bigcap_r,\bigcup_I$ in order to achieve the whole spaces to achieve the analogues of the corresponding FF curves from \cite{10KL1}, \cite{10KL2}, \cite{10FF} for
\[
\xymatrix@R+0pc@C+0pc{
\underset{r}{\mathrm{homotopylimit}}~\mathrm{Spec}^\mathrm{BK}\widetilde{\Phi}^r_{\psi,\Gamma,X_\square},\underset{I}{\mathrm{homotopycolimit}}~\mathrm{Spec}^\mathrm{BK}\widetilde{\Phi}^I_{\psi,\Gamma,X_\square},	\\
}
\]
\[
\xymatrix@R+0pc@C+0pc{
\underset{r}{\mathrm{homotopylimit}}~\mathrm{Spec}^\mathrm{BK}\breve{\Phi}^r_{\psi,\Gamma,X_\square},\underset{I}{\mathrm{homotopycolimit}}~\mathrm{Spec}^\mathrm{BK}\breve{\Phi}^I_{\psi,\Gamma,X_\square},	\\
}
\]
\[
\xymatrix@R+0pc@C+0pc{
\underset{r}{\mathrm{homotopylimit}}~\mathrm{Spec}^\mathrm{BK}{\Phi}^r_{\psi,\Gamma,X_\square},\underset{I}{\mathrm{homotopycolimit}}~\mathrm{Spec}^\mathrm{BK}{\Phi}^I_{\psi,\Gamma,X_\square}.	
}
\]
\[  
\xymatrix@R+0pc@C+0pc{
\underset{r}{\mathrm{homotopylimit}}~\mathrm{Spec}^\mathrm{BK}\widetilde{\Phi}^r_{\psi,\Gamma,X_\square}/\mathrm{Fro}^\mathbb{Z},\underset{I}{\mathrm{homotopycolimit}}~\mathrm{Spec}^\mathrm{BK}\widetilde{\Phi}^I_{\psi,\Gamma,X_\square}/\mathrm{Fro}^\mathbb{Z},	\\
}
\]
\[ 
\xymatrix@R+0pc@C+0pc{
\underset{r}{\mathrm{homotopylimit}}~\mathrm{Spec}^\mathrm{BK}\breve{\Phi}^r_{\psi,\Gamma,X_\square}/\mathrm{Fro}^\mathbb{Z},\underset{I}{\mathrm{homotopycolimit}}~\mathrm{Spec}^\mathrm{BK}\breve{\Phi}^I_{\psi,\Gamma,X_\square}/\mathrm{Fro}^\mathbb{Z},	\\
}
\]
\[ 
\xymatrix@R+0pc@C+0pc{
\underset{r}{\mathrm{homotopylimit}}~\mathrm{Spec}^\mathrm{BK}{\Phi}^r_{\psi,\Gamma,X_\square}/\mathrm{Fro}^\mathbb{Z},\underset{I}{\mathrm{homotopycolimit}}~\mathrm{Spec}^\mathrm{BK}{\Phi}^I_{\psi,\Gamma,X_\square}/\mathrm{Fro}^\mathbb{Z}.	
}
\]

\end{definition}

\indent Meanwhile we have the corresponding Clausen-Scholze analytic stacks from \cite{10CS2}, therefore applying their construction we have:

\begin{definition}
Here we define the following products by using the solidified tensor product from \cite{10CS1} and \cite{10CS2}. Then we take solidified tensor product $\overset{\blacksquare}{\otimes}$ of any of the following
\[
\xymatrix@R+0pc@C+0pc{
\widetilde{\Delta}_{\psi,\Gamma},\widetilde{\nabla}_{\psi,\Gamma},\widetilde{\Phi}_{\psi,\Gamma},\widetilde{\Delta}^+_{\psi,\Gamma},\widetilde{\nabla}^+_{\psi,\Gamma},\widetilde{\Delta}^\dagger_{\psi,\Gamma},\widetilde{\nabla}^\dagger_{\psi,\Gamma},\widetilde{\Phi}^r_{\psi,\Gamma},\widetilde{\Phi}^I_{\psi,\Gamma}, 
}
\]

\[
\xymatrix@R+0pc@C+0pc{
\breve{\Delta}_{\psi,\Gamma},\breve{\nabla}_{\psi,\Gamma},\breve{\Phi}_{\psi,\Gamma},\breve{\Delta}^+_{\psi,\Gamma},\breve{\nabla}^+_{\psi,\Gamma},\breve{\Delta}^\dagger_{\psi,\Gamma},\breve{\nabla}^\dagger_{\psi,\Gamma},\breve{\Phi}^r_{\psi,\Gamma},\breve{\Phi}^I_{\psi,\Gamma},	
}
\]

\[
\xymatrix@R+0pc@C+0pc{
{\Delta}_{\psi,\Gamma},{\nabla}_{\psi,\Gamma},{\Phi}_{\psi,\Gamma},{\Delta}^+_{\psi,\Gamma},{\nabla}^+_{\psi,\Gamma},{\Delta}^\dagger_{\psi,\Gamma},{\nabla}^\dagger_{\psi,\Gamma},{\Phi}^r_{\psi,\Gamma},{\Phi}^I_{\psi,\Gamma},	
}
\]  	
with $X_\square$. Then we have the notations:
\[
\xymatrix@R+0pc@C+0pc{
\widetilde{\Delta}_{\psi,\Gamma,X_\square},\widetilde{\nabla}_{\psi,\Gamma,X_\square},\widetilde{\Phi}_{\psi,\Gamma,X_\square},\widetilde{\Delta}^+_{\psi,\Gamma,X_\square},\widetilde{\nabla}^+_{\psi,\Gamma,X_\square},\widetilde{\Delta}^\dagger_{\psi,\Gamma,X_\square},\widetilde{\nabla}^\dagger_{\psi,\Gamma,X_\square},\widetilde{\Phi}^r_{\psi,\Gamma,X_\square},\widetilde{\Phi}^I_{\psi,\Gamma,X_\square}, 
}
\]

\[
\xymatrix@R+0pc@C+0pc{
\breve{\Delta}_{\psi,\Gamma,X_\square},\breve{\nabla}_{\psi,\Gamma,X_\square},\breve{\Phi}_{\psi,\Gamma,X_\square},\breve{\Delta}^+_{\psi,\Gamma,X_\square},\breve{\nabla}^+_{\psi,\Gamma,X_\square},\breve{\Delta}^\dagger_{\psi,\Gamma,X_\square},\breve{\nabla}^\dagger_{\psi,\Gamma,X_\square},\breve{\Phi}^r_{\psi,\Gamma,X_\square},\breve{\Phi}^I_{\psi,\Gamma,X_\square},	
}
\]

\[
\xymatrix@R+0pc@C+0pc{
{\Delta}_{\psi,\Gamma,X_\square},{\nabla}_{\psi,\Gamma,X_\square},{\Phi}_{\psi,\Gamma,X_\square},{\Delta}^+_{\psi,\Gamma,X_\square},{\nabla}^+_{\psi,\Gamma,X_\square},{\Delta}^\dagger_{\psi,\Gamma,X_\square},{\nabla}^\dagger_{\psi,\Gamma,X_\square},{\Phi}^r_{\psi,\Gamma,X_\square},{\Phi}^I_{\psi,\Gamma,X_\square}.	
}
\]
\end{definition}

\begin{definition}
First we consider the Clausen-Scholze spectrum $\mathrm{Spec}^\mathrm{CS}(*)$ attached to any of those in the above from \cite{10CS2} by taking derived rational localization:
\begin{align}
\mathrm{Spec}^\mathrm{CS}\widetilde{\Delta}_{\psi,\Gamma,X_\square},\mathrm{Spec}^\mathrm{CS}\widetilde{\nabla}_{\psi,\Gamma,X_\square},\mathrm{Spec}^\mathrm{CS}\widetilde{\Phi}_{\psi,\Gamma,X_\square},\mathrm{Spec}^\mathrm{CS}\widetilde{\Delta}^+_{\psi,\Gamma,X_\square},\mathrm{Spec}^\mathrm{CS}\widetilde{\nabla}^+_{\psi,\Gamma,X_\square},\\
\mathrm{Spec}^\mathrm{CS}\widetilde{\Delta}^\dagger_{\psi,\Gamma,X_\square},\mathrm{Spec}^\mathrm{CS}\widetilde{\nabla}^\dagger_{\psi,\Gamma,X_\square},\mathrm{Spec}^\mathrm{CS}\widetilde{\Phi}^r_{\psi,\Gamma,X_\square},\mathrm{Spec}^\mathrm{CS}\widetilde{\Phi}^I_{\psi,\Gamma,X_\square},	\\
\end{align}
\begin{align}
\mathrm{Spec}^\mathrm{CS}\breve{\Delta}_{\psi,\Gamma,X_\square},\breve{\nabla}_{\psi,\Gamma,X_\square},\mathrm{Spec}^\mathrm{CS}\breve{\Phi}_{\psi,\Gamma,X_\square},\mathrm{Spec}^\mathrm{CS}\breve{\Delta}^+_{\psi,\Gamma,X_\square},\mathrm{Spec}^\mathrm{CS}\breve{\nabla}^+_{\psi,\Gamma,X_\square},\\
\mathrm{Spec}^\mathrm{CS}\breve{\Delta}^\dagger_{\psi,\Gamma,X_\square},\mathrm{Spec}^\mathrm{CS}\breve{\nabla}^\dagger_{\psi,\Gamma,X_\square},\mathrm{Spec}^\mathrm{CS}\breve{\Phi}^r_{\psi,\Gamma,X_\square},\breve{\Phi}^I_{\psi,\Gamma,X_\square},	\\
\end{align}
\begin{align}
\mathrm{Spec}^\mathrm{CS}{\Delta}_{\psi,\Gamma,X_\square},\mathrm{Spec}^\mathrm{CS}{\nabla}_{\psi,\Gamma,X_\square},\mathrm{Spec}^\mathrm{CS}{\Phi}_{\psi,\Gamma,X_\square},\mathrm{Spec}^\mathrm{CS}{\Delta}^+_{\psi,\Gamma,X_\square},\mathrm{Spec}^\mathrm{CS}{\nabla}^+_{\psi,\Gamma,X_\square},\\
\mathrm{Spec}^\mathrm{CS}{\Delta}^\dagger_{\psi,\Gamma,X_\square},\mathrm{Spec}^\mathrm{CS}{\nabla}^\dagger_{\psi,\Gamma,X_\square},\mathrm{Spec}^\mathrm{CS}{\Phi}^r_{\psi,\Gamma,X_\square},\mathrm{Spec}^\mathrm{CS}{\Phi}^I_{\psi,\Gamma,X_\square}.	
\end{align}

Then we take the corresponding quotients by using the corresponding Frobenius operators:
\begin{align}
&\mathrm{Spec}^\mathrm{CS}\widetilde{\Delta}_{\psi,\Gamma,X_\square}/\mathrm{Fro}^\mathbb{Z},\mathrm{Spec}^\mathrm{CS}\widetilde{\nabla}_{\psi,\Gamma,X_\square}/\mathrm{Fro}^\mathbb{Z},\mathrm{Spec}^\mathrm{CS}\widetilde{\Phi}_{\psi,\Gamma,X_\square}/\mathrm{Fro}^\mathbb{Z},\mathrm{Spec}^\mathrm{CS}\widetilde{\Delta}^+_{\psi,\Gamma,X_\square}/\mathrm{Fro}^\mathbb{Z},\\
&\mathrm{Spec}^\mathrm{CS}\widetilde{\nabla}^+_{\psi,\Gamma,X_\square}/\mathrm{Fro}^\mathbb{Z}, \mathrm{Spec}^\mathrm{CS}\widetilde{\Delta}^\dagger_{\psi,\Gamma,X_\square}/\mathrm{Fro}^\mathbb{Z},\mathrm{Spec}^\mathrm{CS}\widetilde{\nabla}^\dagger_{\psi,\Gamma,X_\square}/\mathrm{Fro}^\mathbb{Z},	\\
\end{align}
\begin{align}
&\mathrm{Spec}^\mathrm{CS}\breve{\Delta}_{\psi,\Gamma,X_\square}/\mathrm{Fro}^\mathbb{Z},\breve{\nabla}_{\psi,\Gamma,X_\square}/\mathrm{Fro}^\mathbb{Z},\mathrm{Spec}^\mathrm{CS}\breve{\Phi}_{\psi,\Gamma,X_\square}/\mathrm{Fro}^\mathbb{Z},\mathrm{Spec}^\mathrm{CS}\breve{\Delta}^+_{\psi,\Gamma,X_\square}/\mathrm{Fro}^\mathbb{Z},\\
&\mathrm{Spec}^\mathrm{CS}\breve{\nabla}^+_{\psi,\Gamma,X_\square}/\mathrm{Fro}^\mathbb{Z}, \mathrm{Spec}^\mathrm{CS}\breve{\Delta}^\dagger_{\psi,\Gamma,X_\square}/\mathrm{Fro}^\mathbb{Z},\mathrm{Spec}^\mathrm{CS}\breve{\nabla}^\dagger_{\psi,\Gamma,X_\square}/\mathrm{Fro}^\mathbb{Z},	\\
\end{align}
\begin{align}
&\mathrm{Spec}^\mathrm{CS}{\Delta}_{\psi,\Gamma,X_\square}/\mathrm{Fro}^\mathbb{Z},\mathrm{Spec}^\mathrm{CS}{\nabla}_{\psi,\Gamma,X_\square}/\mathrm{Fro}^\mathbb{Z},\mathrm{Spec}^\mathrm{CS}{\Phi}_{\psi,\Gamma,X_\square}/\mathrm{Fro}^\mathbb{Z},\mathrm{Spec}^\mathrm{CS}{\Delta}^+_{\psi,\Gamma,X_\square}/\mathrm{Fro}^\mathbb{Z},\\
&\mathrm{Spec}^\mathrm{CS}{\nabla}^+_{\psi,\Gamma,X_\square}/\mathrm{Fro}^\mathbb{Z}, \mathrm{Spec}^\mathrm{CS}{\Delta}^\dagger_{\psi,\Gamma,X_\square}/\mathrm{Fro}^\mathbb{Z},\mathrm{Spec}^\mathrm{CS}{\nabla}^\dagger_{\psi,\Gamma,X_\square}/\mathrm{Fro}^\mathbb{Z}.	
\end{align}
Here for those space with notations related to the radius and the corresponding interval we consider the total unions $\bigcap_r,\bigcup_I$ in order to achieve the whole spaces to achieve the analogues of the corresponding FF curves from \cite{10KL1}, \cite{10KL2}, \cite{10FF} for
\[
\xymatrix@R+0pc@C+0pc{
\underset{r}{\mathrm{homotopylimit}}~\mathrm{Spec}^\mathrm{CS}\widetilde{\Phi}^r_{\psi,\Gamma,X_\square},\underset{I}{\mathrm{homotopycolimit}}~\mathrm{Spec}^\mathrm{CS}\widetilde{\Phi}^I_{\psi,\Gamma,X_\square},	\\
}
\]
\[
\xymatrix@R+0pc@C+0pc{
\underset{r}{\mathrm{homotopylimit}}~\mathrm{Spec}^\mathrm{CS}\breve{\Phi}^r_{\psi,\Gamma,X_\square},\underset{I}{\mathrm{homotopycolimit}}~\mathrm{Spec}^\mathrm{CS}\breve{\Phi}^I_{\psi,\Gamma,X_\square},	\\
}
\]
\[
\xymatrix@R+0pc@C+0pc{
\underset{r}{\mathrm{homotopylimit}}~\mathrm{Spec}^\mathrm{CS}{\Phi}^r_{\psi,\Gamma,X_\square},\underset{I}{\mathrm{homotopycolimit}}~\mathrm{Spec}^\mathrm{CS}{\Phi}^I_{\psi,\Gamma,X_\square}.	
}
\]
\[ 
\xymatrix@R+0pc@C+0pc{
\underset{r}{\mathrm{homotopylimit}}~\mathrm{Spec}^\mathrm{CS}\widetilde{\Phi}^r_{\psi,\Gamma,X_\square}/\mathrm{Fro}^\mathbb{Z},\underset{I}{\mathrm{homotopycolimit}}~\mathrm{Spec}^\mathrm{CS}\widetilde{\Phi}^I_{\psi,\Gamma,X_\square}/\mathrm{Fro}^\mathbb{Z},	\\
}
\]
\[ 
\xymatrix@R+0pc@C+0pc{
\underset{r}{\mathrm{homotopylimit}}~\mathrm{Spec}^\mathrm{CS}\breve{\Phi}^r_{\psi,\Gamma,X_\square}/\mathrm{Fro}^\mathbb{Z},\underset{I}{\mathrm{homotopycolimit}}~\breve{\Phi}^I_{\psi,\Gamma,X_\square}/\mathrm{Fro}^\mathbb{Z},	\\
}
\]
\[ 
\xymatrix@R+0pc@C+0pc{
\underset{r}{\mathrm{homotopylimit}}~\mathrm{Spec}^\mathrm{CS}{\Phi}^r_{\psi,\Gamma,X_\square}/\mathrm{Fro}^\mathbb{Z},\underset{I}{\mathrm{homotopycolimit}}~\mathrm{Spec}^\mathrm{CS}{\Phi}^I_{\psi,\Gamma,X_\square}/\mathrm{Fro}^\mathbb{Z}.	
}
\]

\end{definition}

\

\begin{definition}
We then consider the corresponding quasipresheaves of the corresponding ind-Banach or monomorphic ind-Banach modules from \cite{10BBK}, \cite{10KKM}\footnote{Here the categories are defined to be the corresponding homotopy colimits of the corresponding categories with respect to each $\square_i$.}:
\begin{align}
\mathrm{Quasicoherentpresheaves,IndBanach}_{*}	
\end{align}
where $*$ is one of the following spaces:
\begin{align}
&\mathrm{Spec}^\mathrm{BK}\widetilde{\Phi}_{\psi,\Gamma,X_\square}/\mathrm{Fro}^\mathbb{Z},	\\
\end{align}
\begin{align}
&\mathrm{Spec}^\mathrm{BK}\breve{\Phi}_{\psi,\Gamma,X_\square}/\mathrm{Fro}^\mathbb{Z},	\\
\end{align}
\begin{align}
&\mathrm{Spec}^\mathrm{BK}{\Phi}_{\psi,\Gamma,X_\square}/\mathrm{Fro}^\mathbb{Z}.	
\end{align}
Here for those space without notation related to the radius and the corresponding interval we consider the total unions $\bigcap_r,\bigcup_I$ in order to achieve the whole spaces to achieve the analogues of the corresponding FF curves from \cite{10KL1}, \cite{10KL2}, \cite{10FF} for
\[
\xymatrix@R+0pc@C+0pc{
\underset{r}{\mathrm{homotopylimit}}~\mathrm{Spec}^\mathrm{BK}\widetilde{\Phi}^r_{\psi,\Gamma,X_\square},\underset{I}{\mathrm{homotopycolimit}}~\mathrm{Spec}^\mathrm{BK}\widetilde{\Phi}^I_{\psi,\Gamma,X_\square},	\\
}
\]
\[
\xymatrix@R+0pc@C+0pc{
\underset{r}{\mathrm{homotopylimit}}~\mathrm{Spec}^\mathrm{BK}\breve{\Phi}^r_{\psi,\Gamma,X_\square},\underset{I}{\mathrm{homotopycolimit}}~\mathrm{Spec}^\mathrm{BK}\breve{\Phi}^I_{\psi,\Gamma,X_\square},	\\
}
\]
\[
\xymatrix@R+0pc@C+0pc{
\underset{r}{\mathrm{homotopylimit}}~\mathrm{Spec}^\mathrm{BK}{\Phi}^r_{\psi,\Gamma,X_\square},\underset{I}{\mathrm{homotopycolimit}}~\mathrm{Spec}^\mathrm{BK}{\Phi}^I_{\psi,\Gamma,X_\square}.	
}
\]
\[  
\xymatrix@R+0pc@C+0pc{
\underset{r}{\mathrm{homotopylimit}}~\mathrm{Spec}^\mathrm{BK}\widetilde{\Phi}^r_{\psi,\Gamma,X_\square}/\mathrm{Fro}^\mathbb{Z},\underset{I}{\mathrm{homotopycolimit}}~\mathrm{Spec}^\mathrm{BK}\widetilde{\Phi}^I_{\psi,\Gamma,X_\square}/\mathrm{Fro}^\mathbb{Z},	\\
}
\]
\[ 
\xymatrix@R+0pc@C+0pc{
\underset{r}{\mathrm{homotopylimit}}~\mathrm{Spec}^\mathrm{BK}\breve{\Phi}^r_{\psi,\Gamma,X_\square}/\mathrm{Fro}^\mathbb{Z},\underset{I}{\mathrm{homotopycolimit}}~\mathrm{Spec}^\mathrm{BK}\breve{\Phi}^I_{\psi,\Gamma,X_\square}/\mathrm{Fro}^\mathbb{Z},	\\
}
\]
\[ 
\xymatrix@R+0pc@C+0pc{
\underset{r}{\mathrm{homotopylimit}}~\mathrm{Spec}^\mathrm{BK}{\Phi}^r_{\psi,\Gamma,X_\square}/\mathrm{Fro}^\mathbb{Z},\underset{I}{\mathrm{homotopycolimit}}~\mathrm{Spec}^\mathrm{BK}{\Phi}^I_{\psi,\Gamma,X_\square}/\mathrm{Fro}^\mathbb{Z}.	
}
\]

\end{definition}

\begin{definition}
We then consider the corresponding quasisheaves of the corresponding condensed solid topological modules from \cite{10CS2}:
\begin{align}
\mathrm{Quasicoherentsheaves, Condensed}_{*}	
\end{align}
where $*$ is one of the following spaces:
\begin{align}
&\mathrm{Spec}^\mathrm{CS}\widetilde{\Delta}_{\psi,\Gamma,X_\square}/\mathrm{Fro}^\mathbb{Z},\mathrm{Spec}^\mathrm{CS}\widetilde{\nabla}_{\psi,\Gamma,X_\square}/\mathrm{Fro}^\mathbb{Z},\mathrm{Spec}^\mathrm{CS}\widetilde{\Phi}_{\psi,\Gamma,X_\square}/\mathrm{Fro}^\mathbb{Z},\mathrm{Spec}^\mathrm{CS}\widetilde{\Delta}^+_{\psi,\Gamma,X_\square}/\mathrm{Fro}^\mathbb{Z},\\
&\mathrm{Spec}^\mathrm{CS}\widetilde{\nabla}^+_{\psi,\Gamma,X_\square}/\mathrm{Fro}^\mathbb{Z},\mathrm{Spec}^\mathrm{CS}\widetilde{\Delta}^\dagger_{\psi,\Gamma,X_\square}/\mathrm{Fro}^\mathbb{Z},\mathrm{Spec}^\mathrm{CS}\widetilde{\nabla}^\dagger_{\psi,\Gamma,X_\square}/\mathrm{Fro}^\mathbb{Z},	\\
\end{align}
\begin{align}
&\mathrm{Spec}^\mathrm{CS}\breve{\Delta}_{\psi,\Gamma,X_\square}/\mathrm{Fro}^\mathbb{Z},\breve{\nabla}_{\psi,\Gamma,X_\square}/\mathrm{Fro}^\mathbb{Z},\mathrm{Spec}^\mathrm{CS}\breve{\Phi}_{\psi,\Gamma,X_\square}/\mathrm{Fro}^\mathbb{Z},\mathrm{Spec}^\mathrm{CS}\breve{\Delta}^+_{\psi,\Gamma,X_\square}/\mathrm{Fro}^\mathbb{Z},\\
&\mathrm{Spec}^\mathrm{CS}\breve{\nabla}^+_{\psi,\Gamma,X_\square}/\mathrm{Fro}^\mathbb{Z},\mathrm{Spec}^\mathrm{CS}\breve{\Delta}^\dagger_{\psi,\Gamma,X_\square}/\mathrm{Fro}^\mathbb{Z},\mathrm{Spec}^\mathrm{CS}\breve{\nabla}^\dagger_{\psi,\Gamma,X_\square}/\mathrm{Fro}^\mathbb{Z},	\\
\end{align}
\begin{align}
&\mathrm{Spec}^\mathrm{CS}{\Delta}_{\psi,\Gamma,X_\square}/\mathrm{Fro}^\mathbb{Z},\mathrm{Spec}^\mathrm{CS}{\nabla}_{\psi,\Gamma,X_\square}/\mathrm{Fro}^\mathbb{Z},\mathrm{Spec}^\mathrm{CS}{\Phi}_{\psi,\Gamma,X_\square}/\mathrm{Fro}^\mathbb{Z},\mathrm{Spec}^\mathrm{CS}{\Delta}^+_{\psi,\Gamma,X_\square}/\mathrm{Fro}^\mathbb{Z},\\
&\mathrm{Spec}^\mathrm{CS}{\nabla}^+_{\psi,\Gamma,X_\square}/\mathrm{Fro}^\mathbb{Z}, \mathrm{Spec}^\mathrm{CS}{\Delta}^\dagger_{\psi,\Gamma,X_\square}/\mathrm{Fro}^\mathbb{Z},\mathrm{Spec}^\mathrm{CS}{\nabla}^\dagger_{\psi,\Gamma,X_\square}/\mathrm{Fro}^\mathbb{Z}.	
\end{align}
Here for those space with notations related to the radius and the corresponding interval we consider the total unions $\bigcap_r,\bigcup_I$ in order to achieve the whole spaces to achieve the analogues of the corresponding FF curves from \cite{10KL1}, \cite{10KL2}, \cite{10FF} for
\[
\xymatrix@R+0pc@C+0pc{
\underset{r}{\mathrm{homotopylimit}}~\mathrm{Spec}^\mathrm{CS}\widetilde{\Phi}^r_{\psi,\Gamma,X_\square},\underset{I}{\mathrm{homotopycolimit}}~\mathrm{Spec}^\mathrm{CS}\widetilde{\Phi}^I_{\psi,\Gamma,X_\square},	\\
}
\]
\[
\xymatrix@R+0pc@C+0pc{
\underset{r}{\mathrm{homotopylimit}}~\mathrm{Spec}^\mathrm{CS}\breve{\Phi}^r_{\psi,\Gamma,X_\square},\underset{I}{\mathrm{homotopycolimit}}~\mathrm{Spec}^\mathrm{CS}\breve{\Phi}^I_{\psi,\Gamma,X_\square},	\\
}
\]
\[
\xymatrix@R+0pc@C+0pc{
\underset{r}{\mathrm{homotopylimit}}~\mathrm{Spec}^\mathrm{CS}{\Phi}^r_{\psi,\Gamma,X_\square},\underset{I}{\mathrm{homotopycolimit}}~\mathrm{Spec}^\mathrm{CS}{\Phi}^I_{\psi,\Gamma,X_\square}.	
}
\]
\[ 
\xymatrix@R+0pc@C+0pc{
\underset{r}{\mathrm{homotopylimit}}~\mathrm{Spec}^\mathrm{CS}\widetilde{\Phi}^r_{\psi,\Gamma,X_\square}/\mathrm{Fro}^\mathbb{Z},\underset{I}{\mathrm{homotopycolimit}}~\mathrm{Spec}^\mathrm{CS}\widetilde{\Phi}^I_{\psi,\Gamma,X_\square}/\mathrm{Fro}^\mathbb{Z},	\\
}
\]
\[ 
\xymatrix@R+0pc@C+0pc{
\underset{r}{\mathrm{homotopylimit}}~\mathrm{Spec}^\mathrm{CS}\breve{\Phi}^r_{\psi,\Gamma,X_\square}/\mathrm{Fro}^\mathbb{Z},\underset{I}{\mathrm{homotopycolimit}}~\breve{\Phi}^I_{\psi,\Gamma,X_\square}/\mathrm{Fro}^\mathbb{Z},	\\
}
\]
\[ 
\xymatrix@R+0pc@C+0pc{
\underset{r}{\mathrm{homotopylimit}}~\mathrm{Spec}^\mathrm{CS}{\Phi}^r_{\psi,\Gamma,X_\square}/\mathrm{Fro}^\mathbb{Z},\underset{I}{\mathrm{homotopycolimit}}~\mathrm{Spec}^\mathrm{CS}{\Phi}^I_{\psi,\Gamma,X_\square}/\mathrm{Fro}^\mathbb{Z}.	
}
\]

\end{definition}

\

\begin{proposition}
There is a well-defined functor from the $\infty$-category 
\begin{align}
\mathrm{Quasicoherentpresheaves,Condensed}_{*}	
\end{align}
where $*$ is one of the following spaces:
\begin{align}
&\mathrm{Spec}^\mathrm{CS}\widetilde{\Phi}_{\psi,\Gamma,X_\square}/\mathrm{Fro}^\mathbb{Z},	\\
\end{align}
\begin{align}
&\mathrm{Spec}^\mathrm{CS}\breve{\Phi}_{\psi,\Gamma,X_\square}/\mathrm{Fro}^\mathbb{Z},	\\
\end{align}
\begin{align}
&\mathrm{Spec}^\mathrm{CS}{\Phi}_{\psi,\Gamma,X_\square}/\mathrm{Fro}^\mathbb{Z},	
\end{align}
to the $\infty$-category of $\mathrm{Fro}$-equivariant quasicoherent presheaves over similar spaces above correspondingly without the $\mathrm{Fro}$-quotients, and to the $\infty$-category of $\mathrm{Fro}$-equivariant quasicoherent modules over global sections of the structure $\infty$-sheaves of the similar spaces above correspondingly without the $\mathrm{Fro}$-quotients. Here for those space without notation related to the radius and the corresponding interval we consider the total unions $\bigcap_r,\bigcup_I$ in order to achieve the whole spaces to achieve the analogues of the corresponding FF curves from \cite{10KL1}, \cite{10KL2}, \cite{10FF} for
\[
\xymatrix@R+0pc@C+0pc{
\underset{r}{\mathrm{homotopylimit}}~\mathrm{Spec}^\mathrm{CS}\widetilde{\Phi}^r_{\psi,\Gamma,X_\square},\underset{I}{\mathrm{homotopycolimit}}~\mathrm{Spec}^\mathrm{CS}\widetilde{\Phi}^I_{\psi,\Gamma,X_\square},	\\
}
\]
\[
\xymatrix@R+0pc@C+0pc{
\underset{r}{\mathrm{homotopylimit}}~\mathrm{Spec}^\mathrm{CS}\breve{\Phi}^r_{\psi,\Gamma,X_\square},\underset{I}{\mathrm{homotopycolimit}}~\mathrm{Spec}^\mathrm{CS}\breve{\Phi}^I_{\psi,\Gamma,X_\square},	\\
}
\]
\[
\xymatrix@R+0pc@C+0pc{
\underset{r}{\mathrm{homotopylimit}}~\mathrm{Spec}^\mathrm{CS}{\Phi}^r_{\psi,\Gamma,X_\square},\underset{I}{\mathrm{homotopycolimit}}~\mathrm{Spec}^\mathrm{CS}{\Phi}^I_{\psi,\Gamma,X_\square}.	
}
\]
\[ 
\xymatrix@R+0pc@C+0pc{
\underset{r}{\mathrm{homotopylimit}}~\mathrm{Spec}^\mathrm{CS}\widetilde{\Phi}^r_{\psi,\Gamma,X_\square}/\mathrm{Fro}^\mathbb{Z},\underset{I}{\mathrm{homotopycolimit}}~\mathrm{Spec}^\mathrm{CS}\widetilde{\Phi}^I_{\psi,\Gamma,X_\square}/\mathrm{Fro}^\mathbb{Z},	\\
}
\]
\[ 
\xymatrix@R+0pc@C+0pc{
\underset{r}{\mathrm{homotopylimit}}~\mathrm{Spec}^\mathrm{CS}\breve{\Phi}^r_{\psi,\Gamma,X_\square}/\mathrm{Fro}^\mathbb{Z},\underset{I}{\mathrm{homotopycolimit}}~\breve{\Phi}^I_{\psi,\Gamma,X_\square}/\mathrm{Fro}^\mathbb{Z},	\\
}
\]
\[ 
\xymatrix@R+0pc@C+0pc{
\underset{r}{\mathrm{homotopylimit}}~\mathrm{Spec}^\mathrm{CS}{\Phi}^r_{\psi,\Gamma,X_\square}/\mathrm{Fro}^\mathbb{Z},\underset{I}{\mathrm{homotopycolimit}}~\mathrm{Spec}^\mathrm{CS}{\Phi}^I_{\psi,\Gamma,X_\square}/\mathrm{Fro}^\mathbb{Z}.	
}
\]	
In this situation we will have the target category being family parametrized by $r$ or $I$ in compatible glueing sense as in \cite[Definition 5.4.10]{10KL2}. In this situation for modules parametrized by the intervals we have the equivalence of $\infty$-categories by using \cite[Proposition 13.8]{10CS2}. Here the corresponding quasicoherent Frobenius modules are defined to be the corresponding homotopy colimits and limits of Frobenius modules:
\begin{align}
\underset{r}{\mathrm{homotopycolimit}}~M_r,\\
\underset{I}{\mathrm{homotopylimit}}~M_I,	
\end{align}
where each $M_r$ is a Frobenius-equivariant module over the period ring with respect to some radius $r$ while each $M_I$ is a Frobenius-equivariant module over the period ring with respect to some interval $I$.\\
\end{proposition}

\begin{proposition}
Similar proposition holds for 
\begin{align}
\mathrm{Quasicoherentsheaves,IndBanach}_{*}.	
\end{align}	
\end{proposition}

\

\begin{definition}
We then consider the corresponding quasipresheaves of perfect complexes the corresponding ind-Banach or monomorphic ind-Banach modules from \cite{10BBK}, \cite{10KKM}:
\begin{align}
\mathrm{Quasicoherentpresheaves,Perfectcomplex,IndBanach}_{*}	
\end{align}
where $*$ is one of the following spaces:
\begin{align}
&\mathrm{Spec}^\mathrm{BK}\widetilde{\Phi}_{\psi,\Gamma,X_\square}/\mathrm{Fro}^\mathbb{Z},	\\
\end{align}
\begin{align}
&\mathrm{Spec}^\mathrm{BK}\breve{\Phi}_{\psi,\Gamma,X_\square}/\mathrm{Fro}^\mathbb{Z},	\\
\end{align}
\begin{align}
&\mathrm{Spec}^\mathrm{BK}{\Phi}_{\psi,\Gamma,X_\square}/\mathrm{Fro}^\mathbb{Z}.	
\end{align}
Here for those space without notation related to the radius and the corresponding interval we consider the total unions $\bigcap_r,\bigcup_I$ in order to achieve the whole spaces to achieve the analogues of the corresponding FF curves from \cite{10KL1}, \cite{10KL2}, \cite{10FF} for
\[
\xymatrix@R+0pc@C+0pc{
\underset{r}{\mathrm{homotopylimit}}~\mathrm{Spec}^\mathrm{BK}\widetilde{\Phi}^r_{\psi,\Gamma,X_\square},\underset{I}{\mathrm{homotopycolimit}}~\mathrm{Spec}^\mathrm{BK}\widetilde{\Phi}^I_{\psi,\Gamma,X_\square},	\\
}
\]
\[
\xymatrix@R+0pc@C+0pc{
\underset{r}{\mathrm{homotopylimit}}~\mathrm{Spec}^\mathrm{BK}\breve{\Phi}^r_{\psi,\Gamma,X_\square},\underset{I}{\mathrm{homotopycolimit}}~\mathrm{Spec}^\mathrm{BK}\breve{\Phi}^I_{\psi,\Gamma,X_\square},	\\
}
\]
\[
\xymatrix@R+0pc@C+0pc{
\underset{r}{\mathrm{homotopylimit}}~\mathrm{Spec}^\mathrm{BK}{\Phi}^r_{\psi,\Gamma,X_\square},\underset{I}{\mathrm{homotopycolimit}}~\mathrm{Spec}^\mathrm{BK}{\Phi}^I_{\psi,\Gamma,X_\square}.	
}
\]
\[  
\xymatrix@R+0pc@C+0pc{
\underset{r}{\mathrm{homotopylimit}}~\mathrm{Spec}^\mathrm{BK}\widetilde{\Phi}^r_{\psi,\Gamma,X_\square}/\mathrm{Fro}^\mathbb{Z},\underset{I}{\mathrm{homotopycolimit}}~\mathrm{Spec}^\mathrm{BK}\widetilde{\Phi}^I_{\psi,\Gamma,X_\square}/\mathrm{Fro}^\mathbb{Z},	\\
}
\]
\[ 
\xymatrix@R+0pc@C+0pc{
\underset{r}{\mathrm{homotopylimit}}~\mathrm{Spec}^\mathrm{BK}\breve{\Phi}^r_{\psi,\Gamma,X_\square}/\mathrm{Fro}^\mathbb{Z},\underset{I}{\mathrm{homotopycolimit}}~\mathrm{Spec}^\mathrm{BK}\breve{\Phi}^I_{\psi,\Gamma,X_\square}/\mathrm{Fro}^\mathbb{Z},	\\
}
\]
\[ 
\xymatrix@R+0pc@C+0pc{
\underset{r}{\mathrm{homotopylimit}}~\mathrm{Spec}^\mathrm{BK}{\Phi}^r_{\psi,\Gamma,X_\square}/\mathrm{Fro}^\mathbb{Z},\underset{I}{\mathrm{homotopycolimit}}~\mathrm{Spec}^\mathrm{BK}{\Phi}^I_{\psi,\Gamma,X_\square}/\mathrm{Fro}^\mathbb{Z}.	
}
\]

\end{definition}

\begin{definition}
We then consider the corresponding quasisheaves of perfect complexes of the corresponding condensed solid topological modules from \cite{10CS2}:
\begin{align}
\mathrm{Quasicoherentsheaves, Perfectcomplex, Condensed}_{*}	
\end{align}
where $*$ is one of the following spaces:
\begin{align}
&\mathrm{Spec}^\mathrm{CS}\widetilde{\Delta}_{\psi,\Gamma,X_\square}/\mathrm{Fro}^\mathbb{Z},\mathrm{Spec}^\mathrm{CS}\widetilde{\nabla}_{\psi,\Gamma,X_\square}/\mathrm{Fro}^\mathbb{Z},\mathrm{Spec}^\mathrm{CS}\widetilde{\Phi}_{\psi,\Gamma,X_\square}/\mathrm{Fro}^\mathbb{Z},\mathrm{Spec}^\mathrm{CS}\widetilde{\Delta}^+_{\psi,\Gamma,X_\square}/\mathrm{Fro}^\mathbb{Z},\\
&\mathrm{Spec}^\mathrm{CS}\widetilde{\nabla}^+_{\psi,\Gamma,X_\square}/\mathrm{Fro}^\mathbb{Z},\mathrm{Spec}^\mathrm{CS}\widetilde{\Delta}^\dagger_{\psi,\Gamma,X_\square}/\mathrm{Fro}^\mathbb{Z},\mathrm{Spec}^\mathrm{CS}\widetilde{\nabla}^\dagger_{\psi,\Gamma,X_\square}/\mathrm{Fro}^\mathbb{Z},	\\
\end{align}
\begin{align}
&\mathrm{Spec}^\mathrm{CS}\breve{\Delta}_{\psi,\Gamma,X_\square}/\mathrm{Fro}^\mathbb{Z},\breve{\nabla}_{\psi,\Gamma,X_\square}/\mathrm{Fro}^\mathbb{Z},\mathrm{Spec}^\mathrm{CS}\breve{\Phi}_{\psi,\Gamma,X_\square}/\mathrm{Fro}^\mathbb{Z},\mathrm{Spec}^\mathrm{CS}\breve{\Delta}^+_{\psi,\Gamma,X_\square}/\mathrm{Fro}^\mathbb{Z},\\
&\mathrm{Spec}^\mathrm{CS}\breve{\nabla}^+_{\psi,\Gamma,X_\square}/\mathrm{Fro}^\mathbb{Z},\mathrm{Spec}^\mathrm{CS}\breve{\Delta}^\dagger_{\psi,\Gamma,X_\square}/\mathrm{Fro}^\mathbb{Z},\mathrm{Spec}^\mathrm{CS}\breve{\nabla}^\dagger_{\psi,\Gamma,X_\square}/\mathrm{Fro}^\mathbb{Z},	\\
\end{align}
\begin{align}
&\mathrm{Spec}^\mathrm{CS}{\Delta}_{\psi,\Gamma,X_\square}/\mathrm{Fro}^\mathbb{Z},\mathrm{Spec}^\mathrm{CS}{\nabla}_{\psi,\Gamma,X_\square}/\mathrm{Fro}^\mathbb{Z},\mathrm{Spec}^\mathrm{CS}{\Phi}_{\psi,\Gamma,X_\square}/\mathrm{Fro}^\mathbb{Z},\mathrm{Spec}^\mathrm{CS}{\Delta}^+_{\psi,\Gamma,X_\square}/\mathrm{Fro}^\mathbb{Z},\\
&\mathrm{Spec}^\mathrm{CS}{\nabla}^+_{\psi,\Gamma,X_\square}/\mathrm{Fro}^\mathbb{Z}, \mathrm{Spec}^\mathrm{CS}{\Delta}^\dagger_{\psi,\Gamma,X_\square}/\mathrm{Fro}^\mathbb{Z},\mathrm{Spec}^\mathrm{CS}{\nabla}^\dagger_{\psi,\Gamma,X_\square}/\mathrm{Fro}^\mathbb{Z}.	
\end{align}
Here for those space with notations related to the radius and the corresponding interval we consider the total unions $\bigcap_r,\bigcup_I$ in order to achieve the whole spaces to achieve the analogues of the corresponding FF curves from \cite{10KL1}, \cite{10KL2}, \cite{10FF} for
\[
\xymatrix@R+0pc@C+0pc{
\underset{r}{\mathrm{homotopylimit}}~\mathrm{Spec}^\mathrm{CS}\widetilde{\Phi}^r_{\psi,\Gamma,X_\square},\underset{I}{\mathrm{homotopycolimit}}~\mathrm{Spec}^\mathrm{CS}\widetilde{\Phi}^I_{\psi,\Gamma,X_\square},	\\
}
\]
\[
\xymatrix@R+0pc@C+0pc{
\underset{r}{\mathrm{homotopylimit}}~\mathrm{Spec}^\mathrm{CS}\breve{\Phi}^r_{\psi,\Gamma,X_\square},\underset{I}{\mathrm{homotopycolimit}}~\mathrm{Spec}^\mathrm{CS}\breve{\Phi}^I_{\psi,\Gamma,X_\square},	\\
}
\]
\[
\xymatrix@R+0pc@C+0pc{
\underset{r}{\mathrm{homotopylimit}}~\mathrm{Spec}^\mathrm{CS}{\Phi}^r_{\psi,\Gamma,X_\square},\underset{I}{\mathrm{homotopycolimit}}~\mathrm{Spec}^\mathrm{CS}{\Phi}^I_{\psi,\Gamma,X_\square}.	
}
\]
\[ 
\xymatrix@R+0pc@C+0pc{
\underset{r}{\mathrm{homotopylimit}}~\mathrm{Spec}^\mathrm{CS}\widetilde{\Phi}^r_{\psi,\Gamma,X_\square}/\mathrm{Fro}^\mathbb{Z},\underset{I}{\mathrm{homotopycolimit}}~\mathrm{Spec}^\mathrm{CS}\widetilde{\Phi}^I_{\psi,\Gamma,X_\square}/\mathrm{Fro}^\mathbb{Z},	\\
}
\]
\[ 
\xymatrix@R+0pc@C+0pc{
\underset{r}{\mathrm{homotopylimit}}~\mathrm{Spec}^\mathrm{CS}\breve{\Phi}^r_{\psi,\Gamma,X_\square}/\mathrm{Fro}^\mathbb{Z},\underset{I}{\mathrm{homotopycolimit}}~\breve{\Phi}^I_{\psi,\Gamma,X_\square}/\mathrm{Fro}^\mathbb{Z},	\\
}
\]
\[ 
\xymatrix@R+0pc@C+0pc{
\underset{r}{\mathrm{homotopylimit}}~\mathrm{Spec}^\mathrm{CS}{\Phi}^r_{\psi,\Gamma,X_\square}/\mathrm{Fro}^\mathbb{Z},\underset{I}{\mathrm{homotopycolimit}}~\mathrm{Spec}^\mathrm{CS}{\Phi}^I_{\psi,\Gamma,X_\square}/\mathrm{Fro}^\mathbb{Z}.	
}
\]

\end{definition}

\begin{proposition}
There is a well-defined functor from the $\infty$-category 
\begin{align}
\mathrm{Quasicoherentpresheaves,Perfectcomplex,Condensed}_{*}	
\end{align}
where $*$ is one of the following spaces:
\begin{align}
&\mathrm{Spec}^\mathrm{CS}\widetilde{\Phi}_{\psi,\Gamma,X_\square}/\mathrm{Fro}^\mathbb{Z},	\\
\end{align}
\begin{align}
&\mathrm{Spec}^\mathrm{CS}\breve{\Phi}_{\psi,\Gamma,X_\square}/\mathrm{Fro}^\mathbb{Z},	\\
\end{align}
\begin{align}
&\mathrm{Spec}^\mathrm{CS}{\Phi}_{\psi,\Gamma,X_\square}/\mathrm{Fro}^\mathbb{Z},	
\end{align}
to the $\infty$-category of $\mathrm{Fro}$-equivariant quasicoherent presheaves over similar spaces above correspondingly without the $\mathrm{Fro}$-quotients, and to the $\infty$-category of $\mathrm{Fro}$-equivariant quasicoherent modules over global sections of the structure $\infty$-sheaves of the similar spaces above correspondingly without the $\mathrm{Fro}$-quotients. Here for those space without notation related to the radius and the corresponding interval we consider the total unions $\bigcap_r,\bigcup_I$ in order to achieve the whole spaces to achieve the analogues of the corresponding FF curves from \cite{10KL1}, \cite{10KL2}, \cite{10FF} for
\[
\xymatrix@R+0pc@C+0pc{
\underset{r}{\mathrm{homotopylimit}}~\mathrm{Spec}^\mathrm{CS}\widetilde{\Phi}^r_{\psi,\Gamma,X_\square},\underset{I}{\mathrm{homotopycolimit}}~\mathrm{Spec}^\mathrm{CS}\widetilde{\Phi}^I_{\psi,\Gamma,X_\square},	\\
}
\]
\[
\xymatrix@R+0pc@C+0pc{
\underset{r}{\mathrm{homotopylimit}}~\mathrm{Spec}^\mathrm{CS}\breve{\Phi}^r_{\psi,\Gamma,X_\square},\underset{I}{\mathrm{homotopycolimit}}~\mathrm{Spec}^\mathrm{CS}\breve{\Phi}^I_{\psi,\Gamma,X_\square},	\\
}
\]
\[
\xymatrix@R+0pc@C+0pc{
\underset{r}{\mathrm{homotopylimit}}~\mathrm{Spec}^\mathrm{CS}{\Phi}^r_{\psi,\Gamma,X_\square},\underset{I}{\mathrm{homotopycolimit}}~\mathrm{Spec}^\mathrm{CS}{\Phi}^I_{\psi,\Gamma,X_\square}.	
}
\]
\[ 
\xymatrix@R+0pc@C+0pc{
\underset{r}{\mathrm{homotopylimit}}~\mathrm{Spec}^\mathrm{CS}\widetilde{\Phi}^r_{\psi,\Gamma,X_\square}/\mathrm{Fro}^\mathbb{Z},\underset{I}{\mathrm{homotopycolimit}}~\mathrm{Spec}^\mathrm{CS}\widetilde{\Phi}^I_{\psi,\Gamma,X_\square}/\mathrm{Fro}^\mathbb{Z},	\\
}
\]
\[ 
\xymatrix@R+0pc@C+0pc{
\underset{r}{\mathrm{homotopylimit}}~\mathrm{Spec}^\mathrm{CS}\breve{\Phi}^r_{\psi,\Gamma,X_\square}/\mathrm{Fro}^\mathbb{Z},\underset{I}{\mathrm{homotopycolimit}}~\breve{\Phi}^I_{\psi,\Gamma,X_\square}/\mathrm{Fro}^\mathbb{Z},	\\
}
\]
\[ 
\xymatrix@R+0pc@C+0pc{
\underset{r}{\mathrm{homotopylimit}}~\mathrm{Spec}^\mathrm{CS}{\Phi}^r_{\psi,\Gamma,X_\square}/\mathrm{Fro}^\mathbb{Z},\underset{I}{\mathrm{homotopycolimit}}~\mathrm{Spec}^\mathrm{CS}{\Phi}^I_{\psi,\Gamma,X_\square}/\mathrm{Fro}^\mathbb{Z}.	
}
\]	
In this situation we will have the target category being family parametrized by $r$ or $I$ in compatible glueing sense as in \cite[Definition 5.4.10]{10KL2}. In this situation for modules parametrized by the intervals we have the equivalence of $\infty$-categories by using \cite[Proposition 12.18]{10CS2}. Here the corresponding quasicoherent Frobenius modules are defined to be the corresponding homotopy colimits and limits of Frobenius modules:
\begin{align}
\underset{r}{\mathrm{homotopycolimit}}~M_r,\\
\underset{I}{\mathrm{homotopylimit}}~M_I,	
\end{align}
where each $M_r$ is a Frobenius-equivariant module over the period ring with respect to some radius $r$ while each $M_I$ is a Frobenius-equivariant module over the period ring with respect to some interval $I$.\\
\end{proposition}

\begin{proposition}
Similar proposition holds for 
\begin{align}
\mathrm{Quasicoherentsheaves,Perfectcomplex,IndBanach}_{*}.	
\end{align}	
\end{proposition}

\subsection{Univariate Hodge Iwasawa Modules}

\subsubsection{Frobenius Quasicoherent Modules I}

\begin{definition}
Let $\psi$ be a toric tower over $\mathbb{Q}_p$ as in \cite[Chapter 7]{10KL2} with base $\mathbb{Q}_p\left<X_1^{\pm 1},...,X_k^{\pm 1}\right>$. Then from \cite{10KL1} and \cite[Definition 5.2.1]{10KL2} we have the following class of Kedlaya-Liu rings (with the following replacement: $\Delta$ stands for $A$, $\nabla$ stands for $B$, while $\Phi$ stands for $C$) by taking product in the sense of self $\Gamma$-th power\footnote{Here $|\Gamma|=1$.}:

\[
\xymatrix@R+0pc@C+0pc{
\widetilde{\Delta}_{\psi},\widetilde{\nabla}_{\psi},\widetilde{\Phi}_{\psi},\widetilde{\Delta}^+_{\psi},\widetilde{\nabla}^+_{\psi},\widetilde{\Delta}^\dagger_{\psi},\widetilde{\nabla}^\dagger_{\psi},\widetilde{\Phi}^r_{\psi},\widetilde{\Phi}^I_{\psi}, 
}
\]

\[
\xymatrix@R+0pc@C+0pc{
\breve{\Delta}_{\psi},\breve{\nabla}_{\psi},\breve{\Phi}_{\psi},\breve{\Delta}^+_{\psi},\breve{\nabla}^+_{\psi},\breve{\Delta}^\dagger_{\psi},\breve{\nabla}^\dagger_{\psi},\breve{\Phi}^r_{\psi},\breve{\Phi}^I_{\psi},	
}
\]

\[
\xymatrix@R+0pc@C+0pc{
{\Delta}_{\psi},{\nabla}_{\psi},{\Phi}_{\psi},{\Delta}^+_{\psi},{\nabla}^+_{\psi},{\Delta}^\dagger_{\psi},{\nabla}^\dagger_{\psi},{\Phi}^r_{\psi},{\Phi}^I_{\psi}.	
}
\]
Now consider $X$ being preadic space over $\mathbb{Q}_p$. Taking the product we have:
\[
\xymatrix@R+0pc@C+0pc{
\widetilde{\Phi}_{\psi,X},\widetilde{\Phi}^r_{\psi,X},\widetilde{\Phi}^I_{\psi,X},	
}
\]
\[
\xymatrix@R+0pc@C+0pc{
\breve{\Phi}_{\psi,X},\breve{\Phi}^r_{\psi,X},\breve{\Phi}^I_{\psi,X},	
}
\]
\[
\xymatrix@R+0pc@C+0pc{
{\Phi}_{\psi,X},{\Phi}^r_{\psi,X},{\Phi}^I_{\psi,X}.	
}
\]
They carry multi Frobenius action $\varphi_\Gamma$ and multi $\mathrm{Lie}_\Gamma:=\mathbb{Z}_p^{\times\Gamma}$ action. In our current situation after \cite{10CKZ} and \cite{10PZ} we consider the following $(\infty,1)$-categories of $(\infty,1)$-modules.\\
\end{definition}

\begin{definition}
First we consider the Bambozzi-Kremnizer spectrum $\mathrm{Spec}^\mathrm{BK}(*)$ attached to any of those in the above from \cite{10BK} by taking derived rational localization:
\begin{align}
&\mathrm{Spec}^\mathrm{BK}\widetilde{\Phi}_{\psi,X},\mathrm{Spec}^\mathrm{BK}\widetilde{\Phi}^r_{\psi,X},\mathrm{Spec}^\mathrm{BK}\widetilde{\Phi}^I_{\psi,X},	
\end{align}
\begin{align}
&\mathrm{Spec}^\mathrm{BK}\breve{\Phi}_{\psi,X},\mathrm{Spec}^\mathrm{BK}\breve{\Phi}^r_{\psi,X},\mathrm{Spec}^\mathrm{BK}\breve{\Phi}^I_{\psi,X},	
\end{align}
\begin{align}
&\mathrm{Spec}^\mathrm{BK}{\Phi}_{\psi,X},
\mathrm{Spec}^\mathrm{BK}{\Phi}^r_{\psi,X},\mathrm{Spec}^\mathrm{BK}{\Phi}^I_{\psi,X}.	
\end{align}

Then we take the corresponding quotients by using the corresponding Frobenius operators:
\begin{align}
&\mathrm{Spec}^\mathrm{BK}\widetilde{\Phi}_{\psi,X}/\mathrm{Fro}^\mathbb{Z},	\\
\end{align}
\begin{align}
&\mathrm{Spec}^\mathrm{BK}\breve{\Phi}_{\psi,X}/\mathrm{Fro}^\mathbb{Z},	\\
\end{align}
\begin{align}
&\mathrm{Spec}^\mathrm{BK}{\Phi}_{\psi,X}/\mathrm{Fro}^\mathbb{Z}.	
\end{align}
Here for those space without notation related to the radius and the corresponding interval we consider the total unions $\bigcap_r,\bigcup_I$ in order to achieve the whole spaces to achieve the analogues of the corresponding FF curves from \cite{10KL1}, \cite{10KL2}, \cite{10FF} for
\[
\xymatrix@R+0pc@C+0pc{
\underset{r}{\mathrm{homotopylimit}}~\mathrm{Spec}^\mathrm{BK}\widetilde{\Phi}^r_{\psi,X},\underset{I}{\mathrm{homotopycolimit}}~\mathrm{Spec}^\mathrm{BK}\widetilde{\Phi}^I_{\psi,X},	\\
}
\]
\[
\xymatrix@R+0pc@C+0pc{
\underset{r}{\mathrm{homotopylimit}}~\mathrm{Spec}^\mathrm{BK}\breve{\Phi}^r_{\psi,X},\underset{I}{\mathrm{homotopycolimit}}~\mathrm{Spec}^\mathrm{BK}\breve{\Phi}^I_{\psi,X},	\\
}
\]
\[
\xymatrix@R+0pc@C+0pc{
\underset{r}{\mathrm{homotopylimit}}~\mathrm{Spec}^\mathrm{BK}{\Phi}^r_{\psi,X},\underset{I}{\mathrm{homotopycolimit}}~\mathrm{Spec}^\mathrm{BK}{\Phi}^I_{\psi,X}.	
}
\]
\[  
\xymatrix@R+0pc@C+0pc{
\underset{r}{\mathrm{homotopylimit}}~\mathrm{Spec}^\mathrm{BK}\widetilde{\Phi}^r_{\psi,X}/\mathrm{Fro}^\mathbb{Z},\underset{I}{\mathrm{homotopycolimit}}~\mathrm{Spec}^\mathrm{BK}\widetilde{\Phi}^I_{\psi,X}/\mathrm{Fro}^\mathbb{Z},	\\
}
\]
\[ 
\xymatrix@R+0pc@C+0pc{
\underset{r}{\mathrm{homotopylimit}}~\mathrm{Spec}^\mathrm{BK}\breve{\Phi}^r_{\psi,X}/\mathrm{Fro}^\mathbb{Z},\underset{I}{\mathrm{homotopycolimit}}~\mathrm{Spec}^\mathrm{BK}\breve{\Phi}^I_{\psi,X}/\mathrm{Fro}^\mathbb{Z},	\\
}
\]
\[ 
\xymatrix@R+0pc@C+0pc{
\underset{r}{\mathrm{homotopylimit}}~\mathrm{Spec}^\mathrm{BK}{\Phi}^r_{\psi,X}/\mathrm{Fro}^\mathbb{Z},\underset{I}{\mathrm{homotopycolimit}}~\mathrm{Spec}^\mathrm{BK}{\Phi}^I_{\psi,X}/\mathrm{Fro}^\mathbb{Z}.	
}
\]

\end{definition}

\indent Meanwhile we have the corresponding Clausen-Scholze analytic stacks from \cite{10CS2}, therefore applying their construction we have:

\begin{definition}
Here we define the following products by using the solidified tensor product from \cite{10CS1} and \cite{10CS2}. Then we take solidified tensor product $\overset{\blacksquare}{\otimes}$ of any of the following
\[
\xymatrix@R+0pc@C+0pc{
\widetilde{\Delta}_{\psi},\widetilde{\nabla}_{\psi},\widetilde{\Phi}_{\psi},\widetilde{\Delta}^+_{\psi},\widetilde{\nabla}^+_{\psi},\widetilde{\Delta}^\dagger_{\psi},\widetilde{\nabla}^\dagger_{\psi},\widetilde{\Phi}^r_{\psi},\widetilde{\Phi}^I_{\psi}, 
}
\]

\[
\xymatrix@R+0pc@C+0pc{
\breve{\Delta}_{\psi},\breve{\nabla}_{\psi},\breve{\Phi}_{\psi},\breve{\Delta}^+_{\psi},\breve{\nabla}^+_{\psi},\breve{\Delta}^\dagger_{\psi},\breve{\nabla}^\dagger_{\psi},\breve{\Phi}^r_{\psi},\breve{\Phi}^I_{\psi},	
}
\]

\[
\xymatrix@R+0pc@C+0pc{
{\Delta}_{\psi},{\nabla}_{\psi},{\Phi}_{\psi},{\Delta}^+_{\psi},{\nabla}^+_{\psi},{\Delta}^\dagger_{\psi},{\nabla}^\dagger_{\psi},{\Phi}^r_{\psi},{\Phi}^I_{\psi},	
}
\]  	
with $X$. Then we have the notations:
\[
\xymatrix@R+0pc@C+0pc{
\widetilde{\Delta}_{\psi,X},\widetilde{\nabla}_{\psi,X},\widetilde{\Phi}_{\psi,X},\widetilde{\Delta}^+_{\psi,X},\widetilde{\nabla}^+_{\psi,X},\widetilde{\Delta}^\dagger_{\psi,X},\widetilde{\nabla}^\dagger_{\psi,X},\widetilde{\Phi}^r_{\psi,X},\widetilde{\Phi}^I_{\psi,X}, 
}
\]

\[
\xymatrix@R+0pc@C+0pc{
\breve{\Delta}_{\psi,X},\breve{\nabla}_{\psi,X},\breve{\Phi}_{\psi,X},\breve{\Delta}^+_{\psi,X},\breve{\nabla}^+_{\psi,X},\breve{\Delta}^\dagger_{\psi,X},\breve{\nabla}^\dagger_{\psi,X},\breve{\Phi}^r_{\psi,X},\breve{\Phi}^I_{\psi,X},	
}
\]

\[
\xymatrix@R+0pc@C+0pc{
{\Delta}_{\psi,X},{\nabla}_{\psi,X},{\Phi}_{\psi,X},{\Delta}^+_{\psi,X},{\nabla}^+_{\psi,X},{\Delta}^\dagger_{\psi,X},{\nabla}^\dagger_{\psi,X},{\Phi}^r_{\psi,X},{\Phi}^I_{\psi,X}.	
}
\]
\end{definition}

\begin{definition}
First we consider the Clausen-Scholze spectrum $\mathrm{Spec}^\mathrm{CS}(*)$ attached to any of those in the above from \cite{10CS2} by taking derived rational localization:
\begin{align}
\mathrm{Spec}^\mathrm{CS}\widetilde{\Delta}_{\psi,X},\mathrm{Spec}^\mathrm{CS}\widetilde{\nabla}_{\psi,X},\mathrm{Spec}^\mathrm{CS}\widetilde{\Phi}_{\psi,X},\mathrm{Spec}^\mathrm{CS}\widetilde{\Delta}^+_{\psi,X},\mathrm{Spec}^\mathrm{CS}\widetilde{\nabla}^+_{\psi,X},\\
\mathrm{Spec}^\mathrm{CS}\widetilde{\Delta}^\dagger_{\psi,X},\mathrm{Spec}^\mathrm{CS}\widetilde{\nabla}^\dagger_{\psi,X},\mathrm{Spec}^\mathrm{CS}\widetilde{\Phi}^r_{\psi,X},\mathrm{Spec}^\mathrm{CS}\widetilde{\Phi}^I_{\psi,X},	\\
\end{align}
\begin{align}
\mathrm{Spec}^\mathrm{CS}\breve{\Delta}_{\psi,X},\breve{\nabla}_{\psi,X},\mathrm{Spec}^\mathrm{CS}\breve{\Phi}_{\psi,X},\mathrm{Spec}^\mathrm{CS}\breve{\Delta}^+_{\psi,X},\mathrm{Spec}^\mathrm{CS}\breve{\nabla}^+_{\psi,X},\\
\mathrm{Spec}^\mathrm{CS}\breve{\Delta}^\dagger_{\psi,X},\mathrm{Spec}^\mathrm{CS}\breve{\nabla}^\dagger_{\psi,X},\mathrm{Spec}^\mathrm{CS}\breve{\Phi}^r_{\psi,X},\breve{\Phi}^I_{\psi,X},	\\
\end{align}
\begin{align}
\mathrm{Spec}^\mathrm{CS}{\Delta}_{\psi,X},\mathrm{Spec}^\mathrm{CS}{\nabla}_{\psi,X},\mathrm{Spec}^\mathrm{CS}{\Phi}_{\psi,X},\mathrm{Spec}^\mathrm{CS}{\Delta}^+_{\psi,X},\mathrm{Spec}^\mathrm{CS}{\nabla}^+_{\psi,X},\\
\mathrm{Spec}^\mathrm{CS}{\Delta}^\dagger_{\psi,X},\mathrm{Spec}^\mathrm{CS}{\nabla}^\dagger_{\psi,X},\mathrm{Spec}^\mathrm{CS}{\Phi}^r_{\psi,X},\mathrm{Spec}^\mathrm{CS}{\Phi}^I_{\psi,X}.	
\end{align}

Then we take the corresponding quotients by using the corresponding Frobenius operators:
\begin{align}
&\mathrm{Spec}^\mathrm{CS}\widetilde{\Delta}_{\psi,X}/\mathrm{Fro}^\mathbb{Z},\mathrm{Spec}^\mathrm{CS}\widetilde{\nabla}_{\psi,X}/\mathrm{Fro}^\mathbb{Z},\mathrm{Spec}^\mathrm{CS}\widetilde{\Phi}_{\psi,X}/\mathrm{Fro}^\mathbb{Z},\mathrm{Spec}^\mathrm{CS}\widetilde{\Delta}^+_{\psi,X}/\mathrm{Fro}^\mathbb{Z},\\
&\mathrm{Spec}^\mathrm{CS}\widetilde{\nabla}^+_{\psi,X}/\mathrm{Fro}^\mathbb{Z}, \mathrm{Spec}^\mathrm{CS}\widetilde{\Delta}^\dagger_{\psi,X}/\mathrm{Fro}^\mathbb{Z},\mathrm{Spec}^\mathrm{CS}\widetilde{\nabla}^\dagger_{\psi,X}/\mathrm{Fro}^\mathbb{Z},	\\
\end{align}
\begin{align}
&\mathrm{Spec}^\mathrm{CS}\breve{\Delta}_{\psi,X}/\mathrm{Fro}^\mathbb{Z},\breve{\nabla}_{\psi,X}/\mathrm{Fro}^\mathbb{Z},\mathrm{Spec}^\mathrm{CS}\breve{\Phi}_{\psi,X}/\mathrm{Fro}^\mathbb{Z},\mathrm{Spec}^\mathrm{CS}\breve{\Delta}^+_{\psi,X}/\mathrm{Fro}^\mathbb{Z},\\
&\mathrm{Spec}^\mathrm{CS}\breve{\nabla}^+_{\psi,X}/\mathrm{Fro}^\mathbb{Z}, \mathrm{Spec}^\mathrm{CS}\breve{\Delta}^\dagger_{\psi,X}/\mathrm{Fro}^\mathbb{Z},\mathrm{Spec}^\mathrm{CS}\breve{\nabla}^\dagger_{\psi,X}/\mathrm{Fro}^\mathbb{Z},	\\
\end{align}
\begin{align}
&\mathrm{Spec}^\mathrm{CS}{\Delta}_{\psi,X}/\mathrm{Fro}^\mathbb{Z},\mathrm{Spec}^\mathrm{CS}{\nabla}_{\psi,X}/\mathrm{Fro}^\mathbb{Z},\mathrm{Spec}^\mathrm{CS}{\Phi}_{\psi,X}/\mathrm{Fro}^\mathbb{Z},\mathrm{Spec}^\mathrm{CS}{\Delta}^+_{\psi,X}/\mathrm{Fro}^\mathbb{Z},\\
&\mathrm{Spec}^\mathrm{CS}{\nabla}^+_{\psi,X}/\mathrm{Fro}^\mathbb{Z}, \mathrm{Spec}^\mathrm{CS}{\Delta}^\dagger_{\psi,X}/\mathrm{Fro}^\mathbb{Z},\mathrm{Spec}^\mathrm{CS}{\nabla}^\dagger_{\psi,X}/\mathrm{Fro}^\mathbb{Z}.	
\end{align}
Here for those space with notations related to the radius and the corresponding interval we consider the total unions $\bigcap_r,\bigcup_I$ in order to achieve the whole spaces to achieve the analogues of the corresponding FF curves from \cite{10KL1}, \cite{10KL2}, \cite{10FF} for
\[
\xymatrix@R+0pc@C+0pc{
\underset{r}{\mathrm{homotopylimit}}~\mathrm{Spec}^\mathrm{CS}\widetilde{\Phi}^r_{\psi,X},\underset{I}{\mathrm{homotopycolimit}}~\mathrm{Spec}^\mathrm{CS}\widetilde{\Phi}^I_{\psi,X},	\\
}
\]
\[
\xymatrix@R+0pc@C+0pc{
\underset{r}{\mathrm{homotopylimit}}~\mathrm{Spec}^\mathrm{CS}\breve{\Phi}^r_{\psi,X},\underset{I}{\mathrm{homotopycolimit}}~\mathrm{Spec}^\mathrm{CS}\breve{\Phi}^I_{\psi,X},	\\
}
\]
\[
\xymatrix@R+0pc@C+0pc{
\underset{r}{\mathrm{homotopylimit}}~\mathrm{Spec}^\mathrm{CS}{\Phi}^r_{\psi,X},\underset{I}{\mathrm{homotopycolimit}}~\mathrm{Spec}^\mathrm{CS}{\Phi}^I_{\psi,X}.	
}
\]
\[ 
\xymatrix@R+0pc@C+0pc{
\underset{r}{\mathrm{homotopylimit}}~\mathrm{Spec}^\mathrm{CS}\widetilde{\Phi}^r_{\psi,X}/\mathrm{Fro}^\mathbb{Z},\underset{I}{\mathrm{homotopycolimit}}~\mathrm{Spec}^\mathrm{CS}\widetilde{\Phi}^I_{\psi,X}/\mathrm{Fro}^\mathbb{Z},	\\
}
\]
\[ 
\xymatrix@R+0pc@C+0pc{
\underset{r}{\mathrm{homotopylimit}}~\mathrm{Spec}^\mathrm{CS}\breve{\Phi}^r_{\psi,X}/\mathrm{Fro}^\mathbb{Z},\underset{I}{\mathrm{homotopycolimit}}~\breve{\Phi}^I_{\psi,X}/\mathrm{Fro}^\mathbb{Z},	\\
}
\]
\[ 
\xymatrix@R+0pc@C+0pc{
\underset{r}{\mathrm{homotopylimit}}~\mathrm{Spec}^\mathrm{CS}{\Phi}^r_{\psi,X}/\mathrm{Fro}^\mathbb{Z},\underset{I}{\mathrm{homotopycolimit}}~\mathrm{Spec}^\mathrm{CS}{\Phi}^I_{\psi,X}/\mathrm{Fro}^\mathbb{Z}.	
}
\]

\end{definition}

\

\begin{definition}
We then consider the corresponding quasipresheaves of the corresponding ind-Banach or monomorphic ind-Banach modules from \cite{10BBK}, \cite{10KKM}:
\begin{align}
\mathrm{Quasicoherentpresheaves,IndBanach}_{*}	
\end{align}
where $*$ is one of the following spaces:
\begin{align}
&\mathrm{Spec}^\mathrm{BK}\widetilde{\Phi}_{\psi,X}/\mathrm{Fro}^\mathbb{Z},	\\
\end{align}
\begin{align}
&\mathrm{Spec}^\mathrm{BK}\breve{\Phi}_{\psi,X}/\mathrm{Fro}^\mathbb{Z},	\\
\end{align}
\begin{align}
&\mathrm{Spec}^\mathrm{BK}{\Phi}_{\psi,X}/\mathrm{Fro}^\mathbb{Z}.	
\end{align}
Here for those space without notation related to the radius and the corresponding interval we consider the total unions $\bigcap_r,\bigcup_I$ in order to achieve the whole spaces to achieve the analogues of the corresponding FF curves from \cite{10KL1}, \cite{10KL2}, \cite{10FF} for
\[
\xymatrix@R+0pc@C+0pc{
\underset{r}{\mathrm{homotopylimit}}~\mathrm{Spec}^\mathrm{BK}\widetilde{\Phi}^r_{\psi,X},\underset{I}{\mathrm{homotopycolimit}}~\mathrm{Spec}^\mathrm{BK}\widetilde{\Phi}^I_{\psi,X},	\\
}
\]
\[
\xymatrix@R+0pc@C+0pc{
\underset{r}{\mathrm{homotopylimit}}~\mathrm{Spec}^\mathrm{BK}\breve{\Phi}^r_{\psi,X},\underset{I}{\mathrm{homotopycolimit}}~\mathrm{Spec}^\mathrm{BK}\breve{\Phi}^I_{\psi,X},	\\
}
\]
\[
\xymatrix@R+0pc@C+0pc{
\underset{r}{\mathrm{homotopylimit}}~\mathrm{Spec}^\mathrm{BK}{\Phi}^r_{\psi,X},\underset{I}{\mathrm{homotopycolimit}}~\mathrm{Spec}^\mathrm{BK}{\Phi}^I_{\psi,X}.	
}
\]
\[  
\xymatrix@R+0pc@C+0pc{
\underset{r}{\mathrm{homotopylimit}}~\mathrm{Spec}^\mathrm{BK}\widetilde{\Phi}^r_{\psi,X}/\mathrm{Fro}^\mathbb{Z},\underset{I}{\mathrm{homotopycolimit}}~\mathrm{Spec}^\mathrm{BK}\widetilde{\Phi}^I_{\psi,X}/\mathrm{Fro}^\mathbb{Z},	\\
}
\]
\[ 
\xymatrix@R+0pc@C+0pc{
\underset{r}{\mathrm{homotopylimit}}~\mathrm{Spec}^\mathrm{BK}\breve{\Phi}^r_{\psi,X}/\mathrm{Fro}^\mathbb{Z},\underset{I}{\mathrm{homotopycolimit}}~\mathrm{Spec}^\mathrm{BK}\breve{\Phi}^I_{\psi,X}/\mathrm{Fro}^\mathbb{Z},	\\
}
\]
\[ 
\xymatrix@R+0pc@C+0pc{
\underset{r}{\mathrm{homotopylimit}}~\mathrm{Spec}^\mathrm{BK}{\Phi}^r_{\psi,X}/\mathrm{Fro}^\mathbb{Z},\underset{I}{\mathrm{homotopycolimit}}~\mathrm{Spec}^\mathrm{BK}{\Phi}^I_{\psi,X}/\mathrm{Fro}^\mathbb{Z}.	
}
\]

\end{definition}

\begin{definition}
We then consider the corresponding quasisheaves of the corresponding condensed solid topological modules from \cite{10CS2}:
\begin{align}
\mathrm{Quasicoherentsheaves, Condensed}_{*}	
\end{align}
where $*$ is one of the following spaces:
\begin{align}
&\mathrm{Spec}^\mathrm{CS}\widetilde{\Delta}_{\psi,X}/\mathrm{Fro}^\mathbb{Z},\mathrm{Spec}^\mathrm{CS}\widetilde{\nabla}_{\psi,X}/\mathrm{Fro}^\mathbb{Z},\mathrm{Spec}^\mathrm{CS}\widetilde{\Phi}_{\psi,X}/\mathrm{Fro}^\mathbb{Z},\mathrm{Spec}^\mathrm{CS}\widetilde{\Delta}^+_{\psi,X}/\mathrm{Fro}^\mathbb{Z},\\
&\mathrm{Spec}^\mathrm{CS}\widetilde{\nabla}^+_{\psi,X}/\mathrm{Fro}^\mathbb{Z},\mathrm{Spec}^\mathrm{CS}\widetilde{\Delta}^\dagger_{\psi,X}/\mathrm{Fro}^\mathbb{Z},\mathrm{Spec}^\mathrm{CS}\widetilde{\nabla}^\dagger_{\psi,X}/\mathrm{Fro}^\mathbb{Z},	\\
\end{align}
\begin{align}
&\mathrm{Spec}^\mathrm{CS}\breve{\Delta}_{\psi,X}/\mathrm{Fro}^\mathbb{Z},\breve{\nabla}_{\psi,X}/\mathrm{Fro}^\mathbb{Z},\mathrm{Spec}^\mathrm{CS}\breve{\Phi}_{\psi,X}/\mathrm{Fro}^\mathbb{Z},\mathrm{Spec}^\mathrm{CS}\breve{\Delta}^+_{\psi,X}/\mathrm{Fro}^\mathbb{Z},\\
&\mathrm{Spec}^\mathrm{CS}\breve{\nabla}^+_{\psi,X}/\mathrm{Fro}^\mathbb{Z},\mathrm{Spec}^\mathrm{CS}\breve{\Delta}^\dagger_{\psi,X}/\mathrm{Fro}^\mathbb{Z},\mathrm{Spec}^\mathrm{CS}\breve{\nabla}^\dagger_{\psi,X}/\mathrm{Fro}^\mathbb{Z},	\\
\end{align}
\begin{align}
&\mathrm{Spec}^\mathrm{CS}{\Delta}_{\psi,X}/\mathrm{Fro}^\mathbb{Z},\mathrm{Spec}^\mathrm{CS}{\nabla}_{\psi,X}/\mathrm{Fro}^\mathbb{Z},\mathrm{Spec}^\mathrm{CS}{\Phi}_{\psi,X}/\mathrm{Fro}^\mathbb{Z},\mathrm{Spec}^\mathrm{CS}{\Delta}^+_{\psi,X}/\mathrm{Fro}^\mathbb{Z},\\
&\mathrm{Spec}^\mathrm{CS}{\nabla}^+_{\psi,X}/\mathrm{Fro}^\mathbb{Z}, \mathrm{Spec}^\mathrm{CS}{\Delta}^\dagger_{\psi,X}/\mathrm{Fro}^\mathbb{Z},\mathrm{Spec}^\mathrm{CS}{\nabla}^\dagger_{\psi,X}/\mathrm{Fro}^\mathbb{Z}.	
\end{align}
Here for those space with notations related to the radius and the corresponding interval we consider the total unions $\bigcap_r,\bigcup_I$ in order to achieve the whole spaces to achieve the analogues of the corresponding FF curves from \cite{10KL1}, \cite{10KL2}, \cite{10FF} for
\[
\xymatrix@R+0pc@C+0pc{
\underset{r}{\mathrm{homotopylimit}}~\mathrm{Spec}^\mathrm{CS}\widetilde{\Phi}^r_{\psi,X},\underset{I}{\mathrm{homotopycolimit}}~\mathrm{Spec}^\mathrm{CS}\widetilde{\Phi}^I_{\psi,X},	\\
}
\]
\[
\xymatrix@R+0pc@C+0pc{
\underset{r}{\mathrm{homotopylimit}}~\mathrm{Spec}^\mathrm{CS}\breve{\Phi}^r_{\psi,X},\underset{I}{\mathrm{homotopycolimit}}~\mathrm{Spec}^\mathrm{CS}\breve{\Phi}^I_{\psi,X},	\\
}
\]
\[
\xymatrix@R+0pc@C+0pc{
\underset{r}{\mathrm{homotopylimit}}~\mathrm{Spec}^\mathrm{CS}{\Phi}^r_{\psi,X},\underset{I}{\mathrm{homotopycolimit}}~\mathrm{Spec}^\mathrm{CS}{\Phi}^I_{\psi,X}.	
}
\]
\[ 
\xymatrix@R+0pc@C+0pc{
\underset{r}{\mathrm{homotopylimit}}~\mathrm{Spec}^\mathrm{CS}\widetilde{\Phi}^r_{\psi,X}/\mathrm{Fro}^\mathbb{Z},\underset{I}{\mathrm{homotopycolimit}}~\mathrm{Spec}^\mathrm{CS}\widetilde{\Phi}^I_{\psi,X}/\mathrm{Fro}^\mathbb{Z},	\\
}
\]
\[ 
\xymatrix@R+0pc@C+0pc{
\underset{r}{\mathrm{homotopylimit}}~\mathrm{Spec}^\mathrm{CS}\breve{\Phi}^r_{\psi,X}/\mathrm{Fro}^\mathbb{Z},\underset{I}{\mathrm{homotopycolimit}}~\breve{\Phi}^I_{\psi,X}/\mathrm{Fro}^\mathbb{Z},	\\
}
\]
\[ 
\xymatrix@R+0pc@C+0pc{
\underset{r}{\mathrm{homotopylimit}}~\mathrm{Spec}^\mathrm{CS}{\Phi}^r_{\psi,X}/\mathrm{Fro}^\mathbb{Z},\underset{I}{\mathrm{homotopycolimit}}~\mathrm{Spec}^\mathrm{CS}{\Phi}^I_{\psi,X}/\mathrm{Fro}^\mathbb{Z}.	
}
\]

\end{definition}

\

\begin{proposition}
There is a well-defined functor from the $\infty$-category 
\begin{align}
\mathrm{Quasicoherentpresheaves,Condensed}_{*}	
\end{align}
where $*$ is one of the following spaces:
\begin{align}
&\mathrm{Spec}^\mathrm{CS}\widetilde{\Phi}_{\psi,X}/\mathrm{Fro}^\mathbb{Z},	\\
\end{align}
\begin{align}
&\mathrm{Spec}^\mathrm{CS}\breve{\Phi}_{\psi,X}/\mathrm{Fro}^\mathbb{Z},	\\
\end{align}
\begin{align}
&\mathrm{Spec}^\mathrm{CS}{\Phi}_{\psi,X}/\mathrm{Fro}^\mathbb{Z},	
\end{align}
to the $\infty$-category of $\mathrm{Fro}$-equivariant quasicoherent presheaves over similar spaces above correspondingly without the $\mathrm{Fro}$-quotients, and to the $\infty$-category of $\mathrm{Fro}$-equivariant quasicoherent modules over global sections of the structure $\infty$-sheaves of the similar spaces above correspondingly without the $\mathrm{Fro}$-quotients. Here for those space without notation related to the radius and the corresponding interval we consider the total unions $\bigcap_r,\bigcup_I$ in order to achieve the whole spaces to achieve the analogues of the corresponding FF curves from \cite{10KL1}, \cite{10KL2}, \cite{10FF} for
\[
\xymatrix@R+0pc@C+0pc{
\underset{r}{\mathrm{homotopylimit}}~\mathrm{Spec}^\mathrm{CS}\widetilde{\Phi}^r_{\psi,X},\underset{I}{\mathrm{homotopycolimit}}~\mathrm{Spec}^\mathrm{CS}\widetilde{\Phi}^I_{\psi,X},	\\
}
\]
\[
\xymatrix@R+0pc@C+0pc{
\underset{r}{\mathrm{homotopylimit}}~\mathrm{Spec}^\mathrm{CS}\breve{\Phi}^r_{\psi,X},\underset{I}{\mathrm{homotopycolimit}}~\mathrm{Spec}^\mathrm{CS}\breve{\Phi}^I_{\psi,X},	\\
}
\]
\[
\xymatrix@R+0pc@C+0pc{
\underset{r}{\mathrm{homotopylimit}}~\mathrm{Spec}^\mathrm{CS}{\Phi}^r_{\psi,X},\underset{I}{\mathrm{homotopycolimit}}~\mathrm{Spec}^\mathrm{CS}{\Phi}^I_{\psi,X}.	
}
\]
\[ 
\xymatrix@R+0pc@C+0pc{
\underset{r}{\mathrm{homotopylimit}}~\mathrm{Spec}^\mathrm{CS}\widetilde{\Phi}^r_{\psi,X}/\mathrm{Fro}^\mathbb{Z},\underset{I}{\mathrm{homotopycolimit}}~\mathrm{Spec}^\mathrm{CS}\widetilde{\Phi}^I_{\psi,X}/\mathrm{Fro}^\mathbb{Z},	\\
}
\]
\[ 
\xymatrix@R+0pc@C+0pc{
\underset{r}{\mathrm{homotopylimit}}~\mathrm{Spec}^\mathrm{CS}\breve{\Phi}^r_{\psi,X}/\mathrm{Fro}^\mathbb{Z},\underset{I}{\mathrm{homotopycolimit}}~\breve{\Phi}^I_{\psi,X}/\mathrm{Fro}^\mathbb{Z},	\\
}
\]
\[ 
\xymatrix@R+0pc@C+0pc{
\underset{r}{\mathrm{homotopylimit}}~\mathrm{Spec}^\mathrm{CS}{\Phi}^r_{\psi,X}/\mathrm{Fro}^\mathbb{Z},\underset{I}{\mathrm{homotopycolimit}}~\mathrm{Spec}^\mathrm{CS}{\Phi}^I_{\psi,X}/\mathrm{Fro}^\mathbb{Z}.	
}
\]	
In this situation we will have the target category being family parametrized by $r$ or $I$ in compatible glueing sense as in \cite[Definition 5.4.10]{10KL2}. In this situation for modules parametrized by the intervals we have the equivalence of $\infty$-categories by using \cite[Proposition 13.8]{10CS2}. Here the corresponding quasicoherent Frobenius modules are defined to be the corresponding homotopy colimits and limits of Frobenius modules:
\begin{align}
\underset{r}{\mathrm{homotopycolimit}}~M_r,\\
\underset{I}{\mathrm{homotopylimit}}~M_I,	
\end{align}
where each $M_r$ is a Frobenius-equivariant module over the period ring with respect to some radius $r$ while each $M_I$ is a Frobenius-equivariant module over the period ring with respect to some interval $I$.\\
\end{proposition}

\begin{proposition}
Similar proposition holds for 
\begin{align}
\mathrm{Quasicoherentsheaves,IndBanach}_{*}.	
\end{align}	
\end{proposition}

\

\begin{definition}
We then consider the corresponding quasipresheaves of perfect complexes the corresponding ind-Banach or monomorphic ind-Banach modules from \cite{10BBK}, \cite{10KKM}:
\begin{align}
\mathrm{Quasicoherentpresheaves,Perfectcomplex,IndBanach}_{*}	
\end{align}
where $*$ is one of the following spaces:
\begin{align}
&\mathrm{Spec}^\mathrm{BK}\widetilde{\Phi}_{\psi,X}/\mathrm{Fro}^\mathbb{Z},	\\
\end{align}
\begin{align}
&\mathrm{Spec}^\mathrm{BK}\breve{\Phi}_{\psi,X}/\mathrm{Fro}^\mathbb{Z},	\\
\end{align}
\begin{align}
&\mathrm{Spec}^\mathrm{BK}{\Phi}_{\psi,X}/\mathrm{Fro}^\mathbb{Z}.	
\end{align}
Here for those space without notation related to the radius and the corresponding interval we consider the total unions $\bigcap_r,\bigcup_I$ in order to achieve the whole spaces to achieve the analogues of the corresponding FF curves from \cite{10KL1}, \cite{10KL2}, \cite{10FF} for
\[
\xymatrix@R+0pc@C+0pc{
\underset{r}{\mathrm{homotopylimit}}~\mathrm{Spec}^\mathrm{BK}\widetilde{\Phi}^r_{\psi,X},\underset{I}{\mathrm{homotopycolimit}}~\mathrm{Spec}^\mathrm{BK}\widetilde{\Phi}^I_{\psi,X},	\\
}
\]
\[
\xymatrix@R+0pc@C+0pc{
\underset{r}{\mathrm{homotopylimit}}~\mathrm{Spec}^\mathrm{BK}\breve{\Phi}^r_{\psi,X},\underset{I}{\mathrm{homotopycolimit}}~\mathrm{Spec}^\mathrm{BK}\breve{\Phi}^I_{\psi,X},	\\
}
\]
\[
\xymatrix@R+0pc@C+0pc{
\underset{r}{\mathrm{homotopylimit}}~\mathrm{Spec}^\mathrm{BK}{\Phi}^r_{\psi,X},\underset{I}{\mathrm{homotopycolimit}}~\mathrm{Spec}^\mathrm{BK}{\Phi}^I_{\psi,X}.	
}
\]
\[  
\xymatrix@R+0pc@C+0pc{
\underset{r}{\mathrm{homotopylimit}}~\mathrm{Spec}^\mathrm{BK}\widetilde{\Phi}^r_{\psi,X}/\mathrm{Fro}^\mathbb{Z},\underset{I}{\mathrm{homotopycolimit}}~\mathrm{Spec}^\mathrm{BK}\widetilde{\Phi}^I_{\psi,X}/\mathrm{Fro}^\mathbb{Z},	\\
}
\]
\[ 
\xymatrix@R+0pc@C+0pc{
\underset{r}{\mathrm{homotopylimit}}~\mathrm{Spec}^\mathrm{BK}\breve{\Phi}^r_{\psi,X}/\mathrm{Fro}^\mathbb{Z},\underset{I}{\mathrm{homotopycolimit}}~\mathrm{Spec}^\mathrm{BK}\breve{\Phi}^I_{\psi,X}/\mathrm{Fro}^\mathbb{Z},	\\
}
\]
\[ 
\xymatrix@R+0pc@C+0pc{
\underset{r}{\mathrm{homotopylimit}}~\mathrm{Spec}^\mathrm{BK}{\Phi}^r_{\psi,X}/\mathrm{Fro}^\mathbb{Z},\underset{I}{\mathrm{homotopycolimit}}~\mathrm{Spec}^\mathrm{BK}{\Phi}^I_{\psi,X}/\mathrm{Fro}^\mathbb{Z}.	
}
\]

\end{definition}

\begin{definition}
We then consider the corresponding quasisheaves of perfect complexes of the corresponding condensed solid topological modules from \cite{10CS2}:
\begin{align}
\mathrm{Quasicoherentsheaves, Perfectcomplex, Condensed}_{*}	
\end{align}
where $*$ is one of the following spaces:
\begin{align}
&\mathrm{Spec}^\mathrm{CS}\widetilde{\Delta}_{\psi,X}/\mathrm{Fro}^\mathbb{Z},\mathrm{Spec}^\mathrm{CS}\widetilde{\nabla}_{\psi,X}/\mathrm{Fro}^\mathbb{Z},\mathrm{Spec}^\mathrm{CS}\widetilde{\Phi}_{\psi,X}/\mathrm{Fro}^\mathbb{Z},\mathrm{Spec}^\mathrm{CS}\widetilde{\Delta}^+_{\psi,X}/\mathrm{Fro}^\mathbb{Z},\\
&\mathrm{Spec}^\mathrm{CS}\widetilde{\nabla}^+_{\psi,X}/\mathrm{Fro}^\mathbb{Z},\mathrm{Spec}^\mathrm{CS}\widetilde{\Delta}^\dagger_{\psi,X}/\mathrm{Fro}^\mathbb{Z},\mathrm{Spec}^\mathrm{CS}\widetilde{\nabla}^\dagger_{\psi,X}/\mathrm{Fro}^\mathbb{Z},	\\
\end{align}
\begin{align}
&\mathrm{Spec}^\mathrm{CS}\breve{\Delta}_{\psi,X}/\mathrm{Fro}^\mathbb{Z},\breve{\nabla}_{\psi,X}/\mathrm{Fro}^\mathbb{Z},\mathrm{Spec}^\mathrm{CS}\breve{\Phi}_{\psi,X}/\mathrm{Fro}^\mathbb{Z},\mathrm{Spec}^\mathrm{CS}\breve{\Delta}^+_{\psi,X}/\mathrm{Fro}^\mathbb{Z},\\
&\mathrm{Spec}^\mathrm{CS}\breve{\nabla}^+_{\psi,X}/\mathrm{Fro}^\mathbb{Z},\mathrm{Spec}^\mathrm{CS}\breve{\Delta}^\dagger_{\psi,X}/\mathrm{Fro}^\mathbb{Z},\mathrm{Spec}^\mathrm{CS}\breve{\nabla}^\dagger_{\psi,X}/\mathrm{Fro}^\mathbb{Z},	\\
\end{align}
\begin{align}
&\mathrm{Spec}^\mathrm{CS}{\Delta}_{\psi,X}/\mathrm{Fro}^\mathbb{Z},\mathrm{Spec}^\mathrm{CS}{\nabla}_{\psi,X}/\mathrm{Fro}^\mathbb{Z},\mathrm{Spec}^\mathrm{CS}{\Phi}_{\psi,X}/\mathrm{Fro}^\mathbb{Z},\mathrm{Spec}^\mathrm{CS}{\Delta}^+_{\psi,X}/\mathrm{Fro}^\mathbb{Z},\\
&\mathrm{Spec}^\mathrm{CS}{\nabla}^+_{\psi,X}/\mathrm{Fro}^\mathbb{Z}, \mathrm{Spec}^\mathrm{CS}{\Delta}^\dagger_{\psi,X}/\mathrm{Fro}^\mathbb{Z},\mathrm{Spec}^\mathrm{CS}{\nabla}^\dagger_{\psi,X}/\mathrm{Fro}^\mathbb{Z}.	
\end{align}
Here for those space with notations related to the radius and the corresponding interval we consider the total unions $\bigcap_r,\bigcup_I$ in order to achieve the whole spaces to achieve the analogues of the corresponding FF curves from \cite{10KL1}, \cite{10KL2}, \cite{10FF} for
\[
\xymatrix@R+0pc@C+0pc{
\underset{r}{\mathrm{homotopylimit}}~\mathrm{Spec}^\mathrm{CS}\widetilde{\Phi}^r_{\psi,X},\underset{I}{\mathrm{homotopycolimit}}~\mathrm{Spec}^\mathrm{CS}\widetilde{\Phi}^I_{\psi,X},	\\
}
\]
\[
\xymatrix@R+0pc@C+0pc{
\underset{r}{\mathrm{homotopylimit}}~\mathrm{Spec}^\mathrm{CS}\breve{\Phi}^r_{\psi,X},\underset{I}{\mathrm{homotopycolimit}}~\mathrm{Spec}^\mathrm{CS}\breve{\Phi}^I_{\psi,X},	\\
}
\]
\[
\xymatrix@R+0pc@C+0pc{
\underset{r}{\mathrm{homotopylimit}}~\mathrm{Spec}^\mathrm{CS}{\Phi}^r_{\psi,X},\underset{I}{\mathrm{homotopycolimit}}~\mathrm{Spec}^\mathrm{CS}{\Phi}^I_{\psi,X}.	
}
\]
\[ 
\xymatrix@R+0pc@C+0pc{
\underset{r}{\mathrm{homotopylimit}}~\mathrm{Spec}^\mathrm{CS}\widetilde{\Phi}^r_{\psi,X}/\mathrm{Fro}^\mathbb{Z},\underset{I}{\mathrm{homotopycolimit}}~\mathrm{Spec}^\mathrm{CS}\widetilde{\Phi}^I_{\psi,X}/\mathrm{Fro}^\mathbb{Z},	\\
}
\]
\[ 
\xymatrix@R+0pc@C+0pc{
\underset{r}{\mathrm{homotopylimit}}~\mathrm{Spec}^\mathrm{CS}\breve{\Phi}^r_{\psi,X}/\mathrm{Fro}^\mathbb{Z},\underset{I}{\mathrm{homotopycolimit}}~\breve{\Phi}^I_{\psi,X}/\mathrm{Fro}^\mathbb{Z},	\\
}
\]
\[ 
\xymatrix@R+0pc@C+0pc{
\underset{r}{\mathrm{homotopylimit}}~\mathrm{Spec}^\mathrm{CS}{\Phi}^r_{\psi,X}/\mathrm{Fro}^\mathbb{Z},\underset{I}{\mathrm{homotopycolimit}}~\mathrm{Spec}^\mathrm{CS}{\Phi}^I_{\psi,X}/\mathrm{Fro}^\mathbb{Z}.	
}
\]

\end{definition}

\begin{proposition}
There is a well-defined functor from the $\infty$-category 
\begin{align}
\mathrm{Quasicoherentpresheaves,Perfectcomplex,Condensed}_{*}	
\end{align}
where $*$ is one of the following spaces:
\begin{align}
&\mathrm{Spec}^\mathrm{CS}\widetilde{\Phi}_{\psi,X}/\mathrm{Fro}^\mathbb{Z},	\\
\end{align}
\begin{align}
&\mathrm{Spec}^\mathrm{CS}\breve{\Phi}_{\psi,X}/\mathrm{Fro}^\mathbb{Z},	\\
\end{align}
\begin{align}
&\mathrm{Spec}^\mathrm{CS}{\Phi}_{\psi,X}/\mathrm{Fro}^\mathbb{Z},	
\end{align}
to the $\infty$-category of $\mathrm{Fro}$-equivariant quasicoherent presheaves over similar spaces above correspondingly without the $\mathrm{Fro}$-quotients, and to the $\infty$-category of $\mathrm{Fro}$-equivariant quasicoherent modules over global sections of the structure $\infty$-sheaves of the similar spaces above correspondingly without the $\mathrm{Fro}$-quotients. Here for those space without notation related to the radius and the corresponding interval we consider the total unions $\bigcap_r,\bigcup_I$ in order to achieve the whole spaces to achieve the analogues of the corresponding FF curves from \cite{10KL1}, \cite{10KL2}, \cite{10FF} for
\[
\xymatrix@R+0pc@C+0pc{
\underset{r}{\mathrm{homotopylimit}}~\mathrm{Spec}^\mathrm{CS}\widetilde{\Phi}^r_{\psi,X},\underset{I}{\mathrm{homotopycolimit}}~\mathrm{Spec}^\mathrm{CS}\widetilde{\Phi}^I_{\psi,X},	\\
}
\]
\[
\xymatrix@R+0pc@C+0pc{
\underset{r}{\mathrm{homotopylimit}}~\mathrm{Spec}^\mathrm{CS}\breve{\Phi}^r_{\psi,X},\underset{I}{\mathrm{homotopycolimit}}~\mathrm{Spec}^\mathrm{CS}\breve{\Phi}^I_{\psi,X},	\\
}
\]
\[
\xymatrix@R+0pc@C+0pc{
\underset{r}{\mathrm{homotopylimit}}~\mathrm{Spec}^\mathrm{CS}{\Phi}^r_{\psi,X},\underset{I}{\mathrm{homotopycolimit}}~\mathrm{Spec}^\mathrm{CS}{\Phi}^I_{\psi,X}.	
}
\]
\[ 
\xymatrix@R+0pc@C+0pc{
\underset{r}{\mathrm{homotopylimit}}~\mathrm{Spec}^\mathrm{CS}\widetilde{\Phi}^r_{\psi,X}/\mathrm{Fro}^\mathbb{Z},\underset{I}{\mathrm{homotopycolimit}}~\mathrm{Spec}^\mathrm{CS}\widetilde{\Phi}^I_{\psi,X}/\mathrm{Fro}^\mathbb{Z},	\\
}
\]
\[ 
\xymatrix@R+0pc@C+0pc{
\underset{r}{\mathrm{homotopylimit}}~\mathrm{Spec}^\mathrm{CS}\breve{\Phi}^r_{\psi,X}/\mathrm{Fro}^\mathbb{Z},\underset{I}{\mathrm{homotopycolimit}}~\breve{\Phi}^I_{\psi,X}/\mathrm{Fro}^\mathbb{Z},	\\
}
\]
\[ 
\xymatrix@R+0pc@C+0pc{
\underset{r}{\mathrm{homotopylimit}}~\mathrm{Spec}^\mathrm{CS}{\Phi}^r_{\psi,X}/\mathrm{Fro}^\mathbb{Z},\underset{I}{\mathrm{homotopycolimit}}~\mathrm{Spec}^\mathrm{CS}{\Phi}^I_{\psi,X}/\mathrm{Fro}^\mathbb{Z}.	
}
\]	
In this situation we will have the target category being family parametrized by $r$ or $I$ in compatible glueing sense as in \cite[Definition 5.4.10]{10KL2}. In this situation for modules parametrized by the intervals we have the equivalence of $\infty$-categories by using \cite[Proposition 12.18]{10CS2}. Here the corresponding quasicoherent Frobenius modules are defined to be the corresponding homotopy colimits and limits of Frobenius modules:
\begin{align}
\underset{r}{\mathrm{homotopycolimit}}~M_r,\\
\underset{I}{\mathrm{homotopylimit}}~M_I,	
\end{align}
where each $M_r$ is a Frobenius-equivariant module over the period ring with respect to some radius $r$ while each $M_I$ is a Frobenius-equivariant module over the period ring with respect to some interval $I$.\\
\end{proposition}

\begin{proposition}
Similar proposition holds for 
\begin{align}
\mathrm{Quasicoherentsheaves,Perfectcomplex,IndBanach}_{*}.	
\end{align}	
\end{proposition}

\subsubsection{Frobenius Quasicoherent Modules II: Deformation in Preadic Spaces}

\begin{definition}
Let $\psi$ be a toric tower over $\mathbb{Q}_p$ as in \cite[Chapter 7]{10KL2} with base $\mathbb{Q}_p\left<X_1^{\pm 1},...,X_k^{\pm 1}\right>$. Then from \cite{10KL1} and \cite[Definition 5.2.1]{10KL2} we have the following class of Kedlaya-Liu rings (with the following replacement: $\Delta$ stands for $A$, $\nabla$ stands for $B$, while $\Phi$ stands for $C$) by taking product in the sense of self $\Gamma$-th power\footnote{Here $|\Gamma|=1$.}:

\[
\xymatrix@R+0pc@C+0pc{
\widetilde{\Delta}_{\psi},\widetilde{\nabla}_{\psi},\widetilde{\Phi}_{\psi},\widetilde{\Delta}^+_{\psi},\widetilde{\nabla}^+_{\psi},\widetilde{\Delta}^\dagger_{\psi},\widetilde{\nabla}^\dagger_{\psi},\widetilde{\Phi}^r_{\psi},\widetilde{\Phi}^I_{\psi}, 
}
\]

\[
\xymatrix@R+0pc@C+0pc{
\breve{\Delta}_{\psi},\breve{\nabla}_{\psi},\breve{\Phi}_{\psi},\breve{\Delta}^+_{\psi},\breve{\nabla}^+_{\psi},\breve{\Delta}^\dagger_{\psi},\breve{\nabla}^\dagger_{\psi},\breve{\Phi}^r_{\psi},\breve{\Phi}^I_{\psi},	
}
\]

\[
\xymatrix@R+0pc@C+0pc{
{\Delta}_{\psi},{\nabla}_{\psi},{\Phi}_{\psi},{\Delta}^+_{\psi},{\nabla}^+_{\psi},{\Delta}^\dagger_{\psi},{\nabla}^\dagger_{\psi},{\Phi}^r_{\psi},{\Phi}^I_{\psi}.	
}
\]
We now consider $\circ$ being deforming preadic space over $\mathbb{Q}_p$. Taking the product we have:
\[
\xymatrix@R+0pc@C+0pc{
\widetilde{\Phi}_{\psi,\circ},\widetilde{\Phi}^r_{\psi,\circ},\widetilde{\Phi}^I_{\psi,\circ},	
}
\]
\[
\xymatrix@R+0pc@C+0pc{
\breve{\Phi}_{\psi,\circ},\breve{\Phi}^r_{\psi,\circ},\breve{\Phi}^I_{\psi,\circ},	
}
\]
\[
\xymatrix@R+0pc@C+0pc{
{\Phi}_{\psi,\circ},{\Phi}^r_{\psi,\circ},{\Phi}^I_{\psi,\circ}.	
}
\]
They carry multi Frobenius action $\varphi_\Gamma$ and multi $\mathrm{Lie}_\Gamma:=\mathbb{Z}_p^{\times\Gamma}$ action. In our current situation after \cite{10CKZ} and \cite{10PZ} we consider the following $(\infty,1)$-categories of $(\infty,1)$-modules.\\
\end{definition}

\begin{definition}
First we consider the Bambozzi-Kremnizer spectrum $\mathrm{Spec}^\mathrm{BK}(*)$ attached to any of those in the above from \cite{10BK} by taking derived rational localization:
\begin{align}
&\mathrm{Spec}^\mathrm{BK}\widetilde{\Phi}_{\psi,\circ},\mathrm{Spec}^\mathrm{BK}\widetilde{\Phi}^r_{\psi,\circ},\mathrm{Spec}^\mathrm{BK}\widetilde{\Phi}^I_{\psi,\circ},	
\end{align}
\begin{align}
&\mathrm{Spec}^\mathrm{BK}\breve{\Phi}_{\psi,\circ},\mathrm{Spec}^\mathrm{BK}\breve{\Phi}^r_{\psi,\circ},\mathrm{Spec}^\mathrm{BK}\breve{\Phi}^I_{\psi,\circ},	
\end{align}
\begin{align}
&\mathrm{Spec}^\mathrm{BK}{\Phi}_{\psi,\circ},
\mathrm{Spec}^\mathrm{BK}{\Phi}^r_{\psi,\circ},\mathrm{Spec}^\mathrm{BK}{\Phi}^I_{\psi,\circ}.	
\end{align}

Then we take the corresponding quotients by using the corresponding Frobenius operators:
\begin{align}
&\mathrm{Spec}^\mathrm{BK}\widetilde{\Phi}_{\psi,\circ}/\mathrm{Fro}^\mathbb{Z},	\\
\end{align}
\begin{align}
&\mathrm{Spec}^\mathrm{BK}\breve{\Phi}_{\psi,\circ}/\mathrm{Fro}^\mathbb{Z},	\\
\end{align}
\begin{align}
&\mathrm{Spec}^\mathrm{BK}{\Phi}_{\psi,\circ}/\mathrm{Fro}^\mathbb{Z}.	
\end{align}
Here for those space without notation related to the radius and the corresponding interval we consider the total unions $\bigcap_r,\bigcup_I$ in order to achieve the whole spaces to achieve the analogues of the corresponding FF curves from \cite{10KL1}, \cite{10KL2}, \cite{10FF} for
\[
\xymatrix@R+0pc@C+0pc{
\underset{r}{\mathrm{homotopylimit}}~\mathrm{Spec}^\mathrm{BK}\widetilde{\Phi}^r_{\psi,\circ},\underset{I}{\mathrm{homotopycolimit}}~\mathrm{Spec}^\mathrm{BK}\widetilde{\Phi}^I_{\psi,\circ},	\\
}
\]
\[
\xymatrix@R+0pc@C+0pc{
\underset{r}{\mathrm{homotopylimit}}~\mathrm{Spec}^\mathrm{BK}\breve{\Phi}^r_{\psi,\circ},\underset{I}{\mathrm{homotopycolimit}}~\mathrm{Spec}^\mathrm{BK}\breve{\Phi}^I_{\psi,\circ},	\\
}
\]
\[
\xymatrix@R+0pc@C+0pc{
\underset{r}{\mathrm{homotopylimit}}~\mathrm{Spec}^\mathrm{BK}{\Phi}^r_{\psi,\circ},\underset{I}{\mathrm{homotopycolimit}}~\mathrm{Spec}^\mathrm{BK}{\Phi}^I_{\psi,\circ}.	
}
\]
\[  
\xymatrix@R+0pc@C+0pc{
\underset{r}{\mathrm{homotopylimit}}~\mathrm{Spec}^\mathrm{BK}\widetilde{\Phi}^r_{\psi,\circ}/\mathrm{Fro}^\mathbb{Z},\underset{I}{\mathrm{homotopycolimit}}~\mathrm{Spec}^\mathrm{BK}\widetilde{\Phi}^I_{\psi,\circ}/\mathrm{Fro}^\mathbb{Z},	\\
}
\]
\[ 
\xymatrix@R+0pc@C+0pc{
\underset{r}{\mathrm{homotopylimit}}~\mathrm{Spec}^\mathrm{BK}\breve{\Phi}^r_{\psi,\circ}/\mathrm{Fro}^\mathbb{Z},\underset{I}{\mathrm{homotopycolimit}}~\mathrm{Spec}^\mathrm{BK}\breve{\Phi}^I_{\psi,\circ}/\mathrm{Fro}^\mathbb{Z},	\\
}
\]
\[ 
\xymatrix@R+0pc@C+0pc{
\underset{r}{\mathrm{homotopylimit}}~\mathrm{Spec}^\mathrm{BK}{\Phi}^r_{\psi,\circ}/\mathrm{Fro}^\mathbb{Z},\underset{I}{\mathrm{homotopycolimit}}~\mathrm{Spec}^\mathrm{BK}{\Phi}^I_{\psi,\circ}/\mathrm{Fro}^\mathbb{Z}.	
}
\]

\end{definition}

\indent Meanwhile we have the corresponding Clausen-Scholze analytic stacks from \cite{10CS2}, therefore applying their construction we have:

\begin{definition}
Here we define the following products by using the solidified tensor product from \cite{10CS1} and \cite{10CS2}. Then we take solidified tensor product $\overset{\blacksquare}{\otimes}$ of any of the following
\[
\xymatrix@R+0pc@C+0pc{
\widetilde{\Delta}_{\psi},\widetilde{\nabla}_{\psi},\widetilde{\Phi}_{\psi},\widetilde{\Delta}^+_{\psi},\widetilde{\nabla}^+_{\psi},\widetilde{\Delta}^\dagger_{\psi},\widetilde{\nabla}^\dagger_{\psi},\widetilde{\Phi}^r_{\psi},\widetilde{\Phi}^I_{\psi}, 
}
\]

\[
\xymatrix@R+0pc@C+0pc{
\breve{\Delta}_{\psi},\breve{\nabla}_{\psi},\breve{\Phi}_{\psi},\breve{\Delta}^+_{\psi},\breve{\nabla}^+_{\psi},\breve{\Delta}^\dagger_{\psi},\breve{\nabla}^\dagger_{\psi},\breve{\Phi}^r_{\psi},\breve{\Phi}^I_{\psi},	
}
\]

\[
\xymatrix@R+0pc@C+0pc{
{\Delta}_{\psi},{\nabla}_{\psi},{\Phi}_{\psi},{\Delta}^+_{\psi},{\nabla}^+_{\psi},{\Delta}^\dagger_{\psi},{\nabla}^\dagger_{\psi},{\Phi}^r_{\psi},{\Phi}^I_{\psi},	
}
\]  	
with $\circ$. Then we have the notations:
\[
\xymatrix@R+0pc@C+0pc{
\widetilde{\Delta}_{\psi,\circ},\widetilde{\nabla}_{\psi,\circ},\widetilde{\Phi}_{\psi,\circ},\widetilde{\Delta}^+_{\psi,\circ},\widetilde{\nabla}^+_{\psi,\circ},\widetilde{\Delta}^\dagger_{\psi,\circ},\widetilde{\nabla}^\dagger_{\psi,\circ},\widetilde{\Phi}^r_{\psi,\circ},\widetilde{\Phi}^I_{\psi,\circ}, 
}
\]

\[
\xymatrix@R+0pc@C+0pc{
\breve{\Delta}_{\psi,\circ},\breve{\nabla}_{\psi,\circ},\breve{\Phi}_{\psi,\circ},\breve{\Delta}^+_{\psi,\circ},\breve{\nabla}^+_{\psi,\circ},\breve{\Delta}^\dagger_{\psi,\circ},\breve{\nabla}^\dagger_{\psi,\circ},\breve{\Phi}^r_{\psi,\circ},\breve{\Phi}^I_{\psi,\circ},	
}
\]

\[
\xymatrix@R+0pc@C+0pc{
{\Delta}_{\psi,\circ},{\nabla}_{\psi,\circ},{\Phi}_{\psi,\circ},{\Delta}^+_{\psi,\circ},{\nabla}^+_{\psi,\circ},{\Delta}^\dagger_{\psi,\circ},{\nabla}^\dagger_{\psi,\circ},{\Phi}^r_{\psi,\circ},{\Phi}^I_{\psi,\circ}.	
}
\]
\end{definition}

\begin{definition}
First we consider the Clausen-Scholze spectrum $\mathrm{Spec}^\mathrm{CS}(*)$ attached to any of those in the above from \cite{10CS2} by taking derived rational localization:
\begin{align}
\mathrm{Spec}^\mathrm{CS}\widetilde{\Delta}_{\psi,\circ},\mathrm{Spec}^\mathrm{CS}\widetilde{\nabla}_{\psi,\circ},\mathrm{Spec}^\mathrm{CS}\widetilde{\Phi}_{\psi,\circ},\mathrm{Spec}^\mathrm{CS}\widetilde{\Delta}^+_{\psi,\circ},\mathrm{Spec}^\mathrm{CS}\widetilde{\nabla}^+_{\psi,\circ},\\
\mathrm{Spec}^\mathrm{CS}\widetilde{\Delta}^\dagger_{\psi,\circ},\mathrm{Spec}^\mathrm{CS}\widetilde{\nabla}^\dagger_{\psi,\circ},\mathrm{Spec}^\mathrm{CS}\widetilde{\Phi}^r_{\psi,\circ},\mathrm{Spec}^\mathrm{CS}\widetilde{\Phi}^I_{\psi,\circ},	\\
\end{align}
\begin{align}
\mathrm{Spec}^\mathrm{CS}\breve{\Delta}_{\psi,\circ},\breve{\nabla}_{\psi,\circ},\mathrm{Spec}^\mathrm{CS}\breve{\Phi}_{\psi,\circ},\mathrm{Spec}^\mathrm{CS}\breve{\Delta}^+_{\psi,\circ},\mathrm{Spec}^\mathrm{CS}\breve{\nabla}^+_{\psi,\circ},\\
\mathrm{Spec}^\mathrm{CS}\breve{\Delta}^\dagger_{\psi,\circ},\mathrm{Spec}^\mathrm{CS}\breve{\nabla}^\dagger_{\psi,\circ},\mathrm{Spec}^\mathrm{CS}\breve{\Phi}^r_{\psi,\circ},\breve{\Phi}^I_{\psi,\circ},	\\
\end{align}
\begin{align}
\mathrm{Spec}^\mathrm{CS}{\Delta}_{\psi,\circ},\mathrm{Spec}^\mathrm{CS}{\nabla}_{\psi,\circ},\mathrm{Spec}^\mathrm{CS}{\Phi}_{\psi,\circ},\mathrm{Spec}^\mathrm{CS}{\Delta}^+_{\psi,\circ},\mathrm{Spec}^\mathrm{CS}{\nabla}^+_{\psi,\circ},\\
\mathrm{Spec}^\mathrm{CS}{\Delta}^\dagger_{\psi,\circ},\mathrm{Spec}^\mathrm{CS}{\nabla}^\dagger_{\psi,\circ},\mathrm{Spec}^\mathrm{CS}{\Phi}^r_{\psi,\circ},\mathrm{Spec}^\mathrm{CS}{\Phi}^I_{\psi,\circ}.	
\end{align}

Then we take the corresponding quotients by using the corresponding Frobenius operators:
\begin{align}
&\mathrm{Spec}^\mathrm{CS}\widetilde{\Delta}_{\psi,\circ}/\mathrm{Fro}^\mathbb{Z},\mathrm{Spec}^\mathrm{CS}\widetilde{\nabla}_{\psi,\circ}/\mathrm{Fro}^\mathbb{Z},\mathrm{Spec}^\mathrm{CS}\widetilde{\Phi}_{\psi,\circ}/\mathrm{Fro}^\mathbb{Z},\mathrm{Spec}^\mathrm{CS}\widetilde{\Delta}^+_{\psi,\circ}/\mathrm{Fro}^\mathbb{Z},\\
&\mathrm{Spec}^\mathrm{CS}\widetilde{\nabla}^+_{\psi,\circ}/\mathrm{Fro}^\mathbb{Z}, \mathrm{Spec}^\mathrm{CS}\widetilde{\Delta}^\dagger_{\psi,\circ}/\mathrm{Fro}^\mathbb{Z},\mathrm{Spec}^\mathrm{CS}\widetilde{\nabla}^\dagger_{\psi,\circ}/\mathrm{Fro}^\mathbb{Z},	\\
\end{align}
\begin{align}
&\mathrm{Spec}^\mathrm{CS}\breve{\Delta}_{\psi,\circ}/\mathrm{Fro}^\mathbb{Z},\breve{\nabla}_{\psi,\circ}/\mathrm{Fro}^\mathbb{Z},\mathrm{Spec}^\mathrm{CS}\breve{\Phi}_{\psi,\circ}/\mathrm{Fro}^\mathbb{Z},\mathrm{Spec}^\mathrm{CS}\breve{\Delta}^+_{\psi,\circ}/\mathrm{Fro}^\mathbb{Z},\\
&\mathrm{Spec}^\mathrm{CS}\breve{\nabla}^+_{\psi,\circ}/\mathrm{Fro}^\mathbb{Z}, \mathrm{Spec}^\mathrm{CS}\breve{\Delta}^\dagger_{\psi,\circ}/\mathrm{Fro}^\mathbb{Z},\mathrm{Spec}^\mathrm{CS}\breve{\nabla}^\dagger_{\psi,\circ}/\mathrm{Fro}^\mathbb{Z},	\\
\end{align}
\begin{align}
&\mathrm{Spec}^\mathrm{CS}{\Delta}_{\psi,\circ}/\mathrm{Fro}^\mathbb{Z},\mathrm{Spec}^\mathrm{CS}{\nabla}_{\psi,\circ}/\mathrm{Fro}^\mathbb{Z},\mathrm{Spec}^\mathrm{CS}{\Phi}_{\psi,\circ}/\mathrm{Fro}^\mathbb{Z},\mathrm{Spec}^\mathrm{CS}{\Delta}^+_{\psi,\circ}/\mathrm{Fro}^\mathbb{Z},\\
&\mathrm{Spec}^\mathrm{CS}{\nabla}^+_{\psi,\circ}/\mathrm{Fro}^\mathbb{Z}, \mathrm{Spec}^\mathrm{CS}{\Delta}^\dagger_{\psi,\circ}/\mathrm{Fro}^\mathbb{Z},\mathrm{Spec}^\mathrm{CS}{\nabla}^\dagger_{\psi,\circ}/\mathrm{Fro}^\mathbb{Z}.	
\end{align}
Here for those space with notations related to the radius and the corresponding interval we consider the total unions $\bigcap_r,\bigcup_I$ in order to achieve the whole spaces to achieve the analogues of the corresponding FF curves from \cite{10KL1}, \cite{10KL2}, \cite{10FF} for
\[
\xymatrix@R+0pc@C+0pc{
\underset{r}{\mathrm{homotopylimit}}~\mathrm{Spec}^\mathrm{CS}\widetilde{\Phi}^r_{\psi,\circ},\underset{I}{\mathrm{homotopycolimit}}~\mathrm{Spec}^\mathrm{CS}\widetilde{\Phi}^I_{\psi,\circ},	\\
}
\]
\[
\xymatrix@R+0pc@C+0pc{
\underset{r}{\mathrm{homotopylimit}}~\mathrm{Spec}^\mathrm{CS}\breve{\Phi}^r_{\psi,\circ},\underset{I}{\mathrm{homotopycolimit}}~\mathrm{Spec}^\mathrm{CS}\breve{\Phi}^I_{\psi,\circ},	\\
}
\]
\[
\xymatrix@R+0pc@C+0pc{
\underset{r}{\mathrm{homotopylimit}}~\mathrm{Spec}^\mathrm{CS}{\Phi}^r_{\psi,\circ},\underset{I}{\mathrm{homotopycolimit}}~\mathrm{Spec}^\mathrm{CS}{\Phi}^I_{\psi,\circ}.	
}
\]
\[ 
\xymatrix@R+0pc@C+0pc{
\underset{r}{\mathrm{homotopylimit}}~\mathrm{Spec}^\mathrm{CS}\widetilde{\Phi}^r_{\psi,\circ}/\mathrm{Fro}^\mathbb{Z},\underset{I}{\mathrm{homotopycolimit}}~\mathrm{Spec}^\mathrm{CS}\widetilde{\Phi}^I_{\psi,\circ}/\mathrm{Fro}^\mathbb{Z},	\\
}
\]
\[ 
\xymatrix@R+0pc@C+0pc{
\underset{r}{\mathrm{homotopylimit}}~\mathrm{Spec}^\mathrm{CS}\breve{\Phi}^r_{\psi,\circ}/\mathrm{Fro}^\mathbb{Z},\underset{I}{\mathrm{homotopycolimit}}~\breve{\Phi}^I_{\psi,\circ}/\mathrm{Fro}^\mathbb{Z},	\\
}
\]
\[ 
\xymatrix@R+0pc@C+0pc{
\underset{r}{\mathrm{homotopylimit}}~\mathrm{Spec}^\mathrm{CS}{\Phi}^r_{\psi,\circ}/\mathrm{Fro}^\mathbb{Z},\underset{I}{\mathrm{homotopycolimit}}~\mathrm{Spec}^\mathrm{CS}{\Phi}^I_{\psi,\circ}/\mathrm{Fro}^\mathbb{Z}.	
}
\]

\end{definition}

\

\begin{definition}
We then consider the corresponding quasipresheaves of the corresponding ind-Banach or monomorphic ind-Banach modules from \cite{10BBK}, \cite{10KKM}:
\begin{align}
\mathrm{Quasicoherentpresheaves,IndBanach}_{*}	
\end{align}
where $*$ is one of the following spaces:
\begin{align}
&\mathrm{Spec}^\mathrm{BK}\widetilde{\Phi}_{\psi,\circ}/\mathrm{Fro}^\mathbb{Z},	\\
\end{align}
\begin{align}
&\mathrm{Spec}^\mathrm{BK}\breve{\Phi}_{\psi,\circ}/\mathrm{Fro}^\mathbb{Z},	\\
\end{align}
\begin{align}
&\mathrm{Spec}^\mathrm{BK}{\Phi}_{\psi,\circ}/\mathrm{Fro}^\mathbb{Z}.	
\end{align}
Here for those space without notation related to the radius and the corresponding interval we consider the total unions $\bigcap_r,\bigcup_I$ in order to achieve the whole spaces to achieve the analogues of the corresponding FF curves from \cite{10KL1}, \cite{10KL2}, \cite{10FF} for
\[
\xymatrix@R+0pc@C+0pc{
\underset{r}{\mathrm{homotopylimit}}~\mathrm{Spec}^\mathrm{BK}\widetilde{\Phi}^r_{\psi,\circ},\underset{I}{\mathrm{homotopycolimit}}~\mathrm{Spec}^\mathrm{BK}\widetilde{\Phi}^I_{\psi,\circ},	\\
}
\]
\[
\xymatrix@R+0pc@C+0pc{
\underset{r}{\mathrm{homotopylimit}}~\mathrm{Spec}^\mathrm{BK}\breve{\Phi}^r_{\psi,\circ},\underset{I}{\mathrm{homotopycolimit}}~\mathrm{Spec}^\mathrm{BK}\breve{\Phi}^I_{\psi,\circ},	\\
}
\]
\[
\xymatrix@R+0pc@C+0pc{
\underset{r}{\mathrm{homotopylimit}}~\mathrm{Spec}^\mathrm{BK}{\Phi}^r_{\psi,\circ},\underset{I}{\mathrm{homotopycolimit}}~\mathrm{Spec}^\mathrm{BK}{\Phi}^I_{\psi,\circ}.	
}
\]
\[  
\xymatrix@R+0pc@C+0pc{
\underset{r}{\mathrm{homotopylimit}}~\mathrm{Spec}^\mathrm{BK}\widetilde{\Phi}^r_{\psi,\circ}/\mathrm{Fro}^\mathbb{Z},\underset{I}{\mathrm{homotopycolimit}}~\mathrm{Spec}^\mathrm{BK}\widetilde{\Phi}^I_{\psi,\circ}/\mathrm{Fro}^\mathbb{Z},	\\
}
\]
\[ 
\xymatrix@R+0pc@C+0pc{
\underset{r}{\mathrm{homotopylimit}}~\mathrm{Spec}^\mathrm{BK}\breve{\Phi}^r_{\psi,\circ}/\mathrm{Fro}^\mathbb{Z},\underset{I}{\mathrm{homotopycolimit}}~\mathrm{Spec}^\mathrm{BK}\breve{\Phi}^I_{\psi,\circ}/\mathrm{Fro}^\mathbb{Z},	\\
}
\]
\[ 
\xymatrix@R+0pc@C+0pc{
\underset{r}{\mathrm{homotopylimit}}~\mathrm{Spec}^\mathrm{BK}{\Phi}^r_{\psi,\circ}/\mathrm{Fro}^\mathbb{Z},\underset{I}{\mathrm{homotopycolimit}}~\mathrm{Spec}^\mathrm{BK}{\Phi}^I_{\psi,\circ}/\mathrm{Fro}^\mathbb{Z}.	
}
\]

\end{definition}

\begin{definition}
We then consider the corresponding quasisheaves of the corresponding condensed solid topological modules from \cite{10CS2}:
\begin{align}
\mathrm{Quasicoherentsheaves, Condensed}_{*}	
\end{align}
where $*$ is one of the following spaces:
\begin{align}
&\mathrm{Spec}^\mathrm{CS}\widetilde{\Delta}_{\psi,\circ}/\mathrm{Fro}^\mathbb{Z},\mathrm{Spec}^\mathrm{CS}\widetilde{\nabla}_{\psi,\circ}/\mathrm{Fro}^\mathbb{Z},\mathrm{Spec}^\mathrm{CS}\widetilde{\Phi}_{\psi,\circ}/\mathrm{Fro}^\mathbb{Z},\mathrm{Spec}^\mathrm{CS}\widetilde{\Delta}^+_{\psi,\circ}/\mathrm{Fro}^\mathbb{Z},\\
&\mathrm{Spec}^\mathrm{CS}\widetilde{\nabla}^+_{\psi,\circ}/\mathrm{Fro}^\mathbb{Z},\mathrm{Spec}^\mathrm{CS}\widetilde{\Delta}^\dagger_{\psi,\circ}/\mathrm{Fro}^\mathbb{Z},\mathrm{Spec}^\mathrm{CS}\widetilde{\nabla}^\dagger_{\psi,\circ}/\mathrm{Fro}^\mathbb{Z},	\\
\end{align}
\begin{align}
&\mathrm{Spec}^\mathrm{CS}\breve{\Delta}_{\psi,\circ}/\mathrm{Fro}^\mathbb{Z},\breve{\nabla}_{\psi,\circ}/\mathrm{Fro}^\mathbb{Z},\mathrm{Spec}^\mathrm{CS}\breve{\Phi}_{\psi,\circ}/\mathrm{Fro}^\mathbb{Z},\mathrm{Spec}^\mathrm{CS}\breve{\Delta}^+_{\psi,\circ}/\mathrm{Fro}^\mathbb{Z},\\
&\mathrm{Spec}^\mathrm{CS}\breve{\nabla}^+_{\psi,\circ}/\mathrm{Fro}^\mathbb{Z},\mathrm{Spec}^\mathrm{CS}\breve{\Delta}^\dagger_{\psi,\circ}/\mathrm{Fro}^\mathbb{Z},\mathrm{Spec}^\mathrm{CS}\breve{\nabla}^\dagger_{\psi,\circ}/\mathrm{Fro}^\mathbb{Z},	\\
\end{align}
\begin{align}
&\mathrm{Spec}^\mathrm{CS}{\Delta}_{\psi,\circ}/\mathrm{Fro}^\mathbb{Z},\mathrm{Spec}^\mathrm{CS}{\nabla}_{\psi,\circ}/\mathrm{Fro}^\mathbb{Z},\mathrm{Spec}^\mathrm{CS}{\Phi}_{\psi,\circ}/\mathrm{Fro}^\mathbb{Z},\mathrm{Spec}^\mathrm{CS}{\Delta}^+_{\psi,\circ}/\mathrm{Fro}^\mathbb{Z},\\
&\mathrm{Spec}^\mathrm{CS}{\nabla}^+_{\psi,\circ}/\mathrm{Fro}^\mathbb{Z}, \mathrm{Spec}^\mathrm{CS}{\Delta}^\dagger_{\psi,\circ}/\mathrm{Fro}^\mathbb{Z},\mathrm{Spec}^\mathrm{CS}{\nabla}^\dagger_{\psi,\circ}/\mathrm{Fro}^\mathbb{Z}.	
\end{align}
Here for those space with notations related to the radius and the corresponding interval we consider the total unions $\bigcap_r,\bigcup_I$ in order to achieve the whole spaces to achieve the analogues of the corresponding FF curves from \cite{10KL1}, \cite{10KL2}, \cite{10FF} for
\[
\xymatrix@R+0pc@C+0pc{
\underset{r}{\mathrm{homotopylimit}}~\mathrm{Spec}^\mathrm{CS}\widetilde{\Phi}^r_{\psi,\circ},\underset{I}{\mathrm{homotopycolimit}}~\mathrm{Spec}^\mathrm{CS}\widetilde{\Phi}^I_{\psi,\circ},	\\
}
\]
\[
\xymatrix@R+0pc@C+0pc{
\underset{r}{\mathrm{homotopylimit}}~\mathrm{Spec}^\mathrm{CS}\breve{\Phi}^r_{\psi,\circ},\underset{I}{\mathrm{homotopycolimit}}~\mathrm{Spec}^\mathrm{CS}\breve{\Phi}^I_{\psi,\circ},	\\
}
\]
\[
\xymatrix@R+0pc@C+0pc{
\underset{r}{\mathrm{homotopylimit}}~\mathrm{Spec}^\mathrm{CS}{\Phi}^r_{\psi,\circ},\underset{I}{\mathrm{homotopycolimit}}~\mathrm{Spec}^\mathrm{CS}{\Phi}^I_{\psi,\circ}.	
}
\]
\[ 
\xymatrix@R+0pc@C+0pc{
\underset{r}{\mathrm{homotopylimit}}~\mathrm{Spec}^\mathrm{CS}\widetilde{\Phi}^r_{\psi,\circ}/\mathrm{Fro}^\mathbb{Z},\underset{I}{\mathrm{homotopycolimit}}~\mathrm{Spec}^\mathrm{CS}\widetilde{\Phi}^I_{\psi,\circ}/\mathrm{Fro}^\mathbb{Z},	\\
}
\]
\[ 
\xymatrix@R+0pc@C+0pc{
\underset{r}{\mathrm{homotopylimit}}~\mathrm{Spec}^\mathrm{CS}\breve{\Phi}^r_{\psi,\circ}/\mathrm{Fro}^\mathbb{Z},\underset{I}{\mathrm{homotopycolimit}}~\breve{\Phi}^I_{\psi,\circ}/\mathrm{Fro}^\mathbb{Z},	\\
}
\]
\[ 
\xymatrix@R+0pc@C+0pc{
\underset{r}{\mathrm{homotopylimit}}~\mathrm{Spec}^\mathrm{CS}{\Phi}^r_{\psi,\circ}/\mathrm{Fro}^\mathbb{Z},\underset{I}{\mathrm{homotopycolimit}}~\mathrm{Spec}^\mathrm{CS}{\Phi}^I_{\psi,\circ}/\mathrm{Fro}^\mathbb{Z}.	
}
\]

\end{definition}

\

\begin{proposition}
There is a well-defined functor from the $\infty$-category 
\begin{align}
\mathrm{Quasicoherentpresheaves,Condensed}_{*}	
\end{align}
where $*$ is one of the following spaces:
\begin{align}
&\mathrm{Spec}^\mathrm{CS}\widetilde{\Phi}_{\psi,\circ}/\mathrm{Fro}^\mathbb{Z},	\\
\end{align}
\begin{align}
&\mathrm{Spec}^\mathrm{CS}\breve{\Phi}_{\psi,\circ}/\mathrm{Fro}^\mathbb{Z},	\\
\end{align}
\begin{align}
&\mathrm{Spec}^\mathrm{CS}{\Phi}_{\psi,\circ}/\mathrm{Fro}^\mathbb{Z},	
\end{align}
to the $\infty$-category of $\mathrm{Fro}$-equivariant quasicoherent presheaves over similar spaces above correspondingly without the $\mathrm{Fro}$-quotients, and to the $\infty$-category of $\mathrm{Fro}$-equivariant quasicoherent modules over global sections of the structure $\infty$-sheaves of the similar spaces above correspondingly without the $\mathrm{Fro}$-quotients. Here for those space without notation related to the radius and the corresponding interval we consider the total unions $\bigcap_r,\bigcup_I$ in order to achieve the whole spaces to achieve the analogues of the corresponding FF curves from \cite{10KL1}, \cite{10KL2}, \cite{10FF} for
\[
\xymatrix@R+0pc@C+0pc{
\underset{r}{\mathrm{homotopylimit}}~\mathrm{Spec}^\mathrm{CS}\widetilde{\Phi}^r_{\psi,\circ},\underset{I}{\mathrm{homotopycolimit}}~\mathrm{Spec}^\mathrm{CS}\widetilde{\Phi}^I_{\psi,\circ},	\\
}
\]
\[
\xymatrix@R+0pc@C+0pc{
\underset{r}{\mathrm{homotopylimit}}~\mathrm{Spec}^\mathrm{CS}\breve{\Phi}^r_{\psi,\circ},\underset{I}{\mathrm{homotopycolimit}}~\mathrm{Spec}^\mathrm{CS}\breve{\Phi}^I_{\psi,\circ},	\\
}
\]
\[
\xymatrix@R+0pc@C+0pc{
\underset{r}{\mathrm{homotopylimit}}~\mathrm{Spec}^\mathrm{CS}{\Phi}^r_{\psi,\circ},\underset{I}{\mathrm{homotopycolimit}}~\mathrm{Spec}^\mathrm{CS}{\Phi}^I_{\psi,\circ}.	
}
\]
\[ 
\xymatrix@R+0pc@C+0pc{
\underset{r}{\mathrm{homotopylimit}}~\mathrm{Spec}^\mathrm{CS}\widetilde{\Phi}^r_{\psi,\circ}/\mathrm{Fro}^\mathbb{Z},\underset{I}{\mathrm{homotopycolimit}}~\mathrm{Spec}^\mathrm{CS}\widetilde{\Phi}^I_{\psi,\circ}/\mathrm{Fro}^\mathbb{Z},	\\
}
\]
\[ 
\xymatrix@R+0pc@C+0pc{
\underset{r}{\mathrm{homotopylimit}}~\mathrm{Spec}^\mathrm{CS}\breve{\Phi}^r_{\psi,\circ}/\mathrm{Fro}^\mathbb{Z},\underset{I}{\mathrm{homotopycolimit}}~\breve{\Phi}^I_{\psi,\circ}/\mathrm{Fro}^\mathbb{Z},	\\
}
\]
\[ 
\xymatrix@R+0pc@C+0pc{
\underset{r}{\mathrm{homotopylimit}}~\mathrm{Spec}^\mathrm{CS}{\Phi}^r_{\psi,\circ}/\mathrm{Fro}^\mathbb{Z},\underset{I}{\mathrm{homotopycolimit}}~\mathrm{Spec}^\mathrm{CS}{\Phi}^I_{\psi,\circ}/\mathrm{Fro}^\mathbb{Z}.	
}
\]	
In this situation we will have the target category being family parametrized by $r$ or $I$ in compatible glueing sense as in \cite[Definition 5.4.10]{10KL2}. In this situation for modules parametrized by the intervals we have the equivalence of $\infty$-categories by using \cite[Proposition 13.8]{10CS2}. Here the corresponding quasicoherent Frobenius modules are defined to be the corresponding homotopy colimits and limits of Frobenius modules:
\begin{align}
\underset{r}{\mathrm{homotopycolimit}}~M_r,\\
\underset{I}{\mathrm{homotopylimit}}~M_I,	
\end{align}
where each $M_r$ is a Frobenius-equivariant module over the period ring with respect to some radius $r$ while each $M_I$ is a Frobenius-equivariant module over the period ring with respect to some interval $I$.\\
\end{proposition}

\begin{proposition}
Similar proposition holds for 
\begin{align}
\mathrm{Quasicoherentsheaves,IndBanach}_{*}.	
\end{align}	
\end{proposition}

\

\begin{definition}
We then consider the corresponding quasipresheaves of perfect complexes the corresponding ind-Banach or monomorphic ind-Banach modules from \cite{10BBK}, \cite{10KKM}:
\begin{align}
\mathrm{Quasicoherentpresheaves,Perfectcomplex,IndBanach}_{*}	
\end{align}
where $*$ is one of the following spaces:
\begin{align}
&\mathrm{Spec}^\mathrm{BK}\widetilde{\Phi}_{\psi,\circ}/\mathrm{Fro}^\mathbb{Z},	\\
\end{align}
\begin{align}
&\mathrm{Spec}^\mathrm{BK}\breve{\Phi}_{\psi,\circ}/\mathrm{Fro}^\mathbb{Z},	\\
\end{align}
\begin{align}
&\mathrm{Spec}^\mathrm{BK}{\Phi}_{\psi,\circ}/\mathrm{Fro}^\mathbb{Z}.	
\end{align}
Here for those space without notation related to the radius and the corresponding interval we consider the total unions $\bigcap_r,\bigcup_I$ in order to achieve the whole spaces to achieve the analogues of the corresponding FF curves from \cite{10KL1}, \cite{10KL2}, \cite{10FF} for
\[
\xymatrix@R+0pc@C+0pc{
\underset{r}{\mathrm{homotopylimit}}~\mathrm{Spec}^\mathrm{BK}\widetilde{\Phi}^r_{\psi,\circ},\underset{I}{\mathrm{homotopycolimit}}~\mathrm{Spec}^\mathrm{BK}\widetilde{\Phi}^I_{\psi,\circ},	\\
}
\]
\[
\xymatrix@R+0pc@C+0pc{
\underset{r}{\mathrm{homotopylimit}}~\mathrm{Spec}^\mathrm{BK}\breve{\Phi}^r_{\psi,\circ},\underset{I}{\mathrm{homotopycolimit}}~\mathrm{Spec}^\mathrm{BK}\breve{\Phi}^I_{\psi,\circ},	\\
}
\]
\[
\xymatrix@R+0pc@C+0pc{
\underset{r}{\mathrm{homotopylimit}}~\mathrm{Spec}^\mathrm{BK}{\Phi}^r_{\psi,\circ},\underset{I}{\mathrm{homotopycolimit}}~\mathrm{Spec}^\mathrm{BK}{\Phi}^I_{\psi,\circ}.	
}
\]
\[  
\xymatrix@R+0pc@C+0pc{
\underset{r}{\mathrm{homotopylimit}}~\mathrm{Spec}^\mathrm{BK}\widetilde{\Phi}^r_{\psi,\circ}/\mathrm{Fro}^\mathbb{Z},\underset{I}{\mathrm{homotopycolimit}}~\mathrm{Spec}^\mathrm{BK}\widetilde{\Phi}^I_{\psi,\circ}/\mathrm{Fro}^\mathbb{Z},	\\
}
\]
\[ 
\xymatrix@R+0pc@C+0pc{
\underset{r}{\mathrm{homotopylimit}}~\mathrm{Spec}^\mathrm{BK}\breve{\Phi}^r_{\psi,\circ}/\mathrm{Fro}^\mathbb{Z},\underset{I}{\mathrm{homotopycolimit}}~\mathrm{Spec}^\mathrm{BK}\breve{\Phi}^I_{\psi,\circ}/\mathrm{Fro}^\mathbb{Z},	\\
}
\]
\[ 
\xymatrix@R+0pc@C+0pc{
\underset{r}{\mathrm{homotopylimit}}~\mathrm{Spec}^\mathrm{BK}{\Phi}^r_{\psi,\circ}/\mathrm{Fro}^\mathbb{Z},\underset{I}{\mathrm{homotopycolimit}}~\mathrm{Spec}^\mathrm{BK}{\Phi}^I_{\psi,\circ}/\mathrm{Fro}^\mathbb{Z}.	
}
\]

\end{definition}

\begin{definition}
We then consider the corresponding quasisheaves of perfect complexes of the corresponding condensed solid topological modules from \cite{10CS2}:
\begin{align}
\mathrm{Quasicoherentsheaves, Perfectcomplex, Condensed}_{*}	
\end{align}
where $*$ is one of the following spaces:
\begin{align}
&\mathrm{Spec}^\mathrm{CS}\widetilde{\Delta}_{\psi,\circ}/\mathrm{Fro}^\mathbb{Z},\mathrm{Spec}^\mathrm{CS}\widetilde{\nabla}_{\psi,\circ}/\mathrm{Fro}^\mathbb{Z},\mathrm{Spec}^\mathrm{CS}\widetilde{\Phi}_{\psi,\circ}/\mathrm{Fro}^\mathbb{Z},\mathrm{Spec}^\mathrm{CS}\widetilde{\Delta}^+_{\psi,\circ}/\mathrm{Fro}^\mathbb{Z},\\
&\mathrm{Spec}^\mathrm{CS}\widetilde{\nabla}^+_{\psi,\circ}/\mathrm{Fro}^\mathbb{Z},\mathrm{Spec}^\mathrm{CS}\widetilde{\Delta}^\dagger_{\psi,\circ}/\mathrm{Fro}^\mathbb{Z},\mathrm{Spec}^\mathrm{CS}\widetilde{\nabla}^\dagger_{\psi,\circ}/\mathrm{Fro}^\mathbb{Z},	\\
\end{align}
\begin{align}
&\mathrm{Spec}^\mathrm{CS}\breve{\Delta}_{\psi,\circ}/\mathrm{Fro}^\mathbb{Z},\breve{\nabla}_{\psi,\circ}/\mathrm{Fro}^\mathbb{Z},\mathrm{Spec}^\mathrm{CS}\breve{\Phi}_{\psi,\circ}/\mathrm{Fro}^\mathbb{Z},\mathrm{Spec}^\mathrm{CS}\breve{\Delta}^+_{\psi,\circ}/\mathrm{Fro}^\mathbb{Z},\\
&\mathrm{Spec}^\mathrm{CS}\breve{\nabla}^+_{\psi,\circ}/\mathrm{Fro}^\mathbb{Z},\mathrm{Spec}^\mathrm{CS}\breve{\Delta}^\dagger_{\psi,\circ}/\mathrm{Fro}^\mathbb{Z},\mathrm{Spec}^\mathrm{CS}\breve{\nabla}^\dagger_{\psi,\circ}/\mathrm{Fro}^\mathbb{Z},	\\
\end{align}
\begin{align}
&\mathrm{Spec}^\mathrm{CS}{\Delta}_{\psi,\circ}/\mathrm{Fro}^\mathbb{Z},\mathrm{Spec}^\mathrm{CS}{\nabla}_{\psi,\circ}/\mathrm{Fro}^\mathbb{Z},\mathrm{Spec}^\mathrm{CS}{\Phi}_{\psi,\circ}/\mathrm{Fro}^\mathbb{Z},\mathrm{Spec}^\mathrm{CS}{\Delta}^+_{\psi,\circ}/\mathrm{Fro}^\mathbb{Z},\\
&\mathrm{Spec}^\mathrm{CS}{\nabla}^+_{\psi,\circ}/\mathrm{Fro}^\mathbb{Z}, \mathrm{Spec}^\mathrm{CS}{\Delta}^\dagger_{\psi,\circ}/\mathrm{Fro}^\mathbb{Z},\mathrm{Spec}^\mathrm{CS}{\nabla}^\dagger_{\psi,\circ}/\mathrm{Fro}^\mathbb{Z}.	
\end{align}
Here for those space with notations related to the radius and the corresponding interval we consider the total unions $\bigcap_r,\bigcup_I$ in order to achieve the whole spaces to achieve the analogues of the corresponding FF curves from \cite{10KL1}, \cite{10KL2}, \cite{10FF} for
\[
\xymatrix@R+0pc@C+0pc{
\underset{r}{\mathrm{homotopylimit}}~\mathrm{Spec}^\mathrm{CS}\widetilde{\Phi}^r_{\psi,\circ},\underset{I}{\mathrm{homotopycolimit}}~\mathrm{Spec}^\mathrm{CS}\widetilde{\Phi}^I_{\psi,\circ},	\\
}
\]
\[
\xymatrix@R+0pc@C+0pc{
\underset{r}{\mathrm{homotopylimit}}~\mathrm{Spec}^\mathrm{CS}\breve{\Phi}^r_{\psi,\circ},\underset{I}{\mathrm{homotopycolimit}}~\mathrm{Spec}^\mathrm{CS}\breve{\Phi}^I_{\psi,\circ},	\\
}
\]
\[
\xymatrix@R+0pc@C+0pc{
\underset{r}{\mathrm{homotopylimit}}~\mathrm{Spec}^\mathrm{CS}{\Phi}^r_{\psi,\circ},\underset{I}{\mathrm{homotopycolimit}}~\mathrm{Spec}^\mathrm{CS}{\Phi}^I_{\psi,\circ}.	
}
\]
\[ 
\xymatrix@R+0pc@C+0pc{
\underset{r}{\mathrm{homotopylimit}}~\mathrm{Spec}^\mathrm{CS}\widetilde{\Phi}^r_{\psi,\circ}/\mathrm{Fro}^\mathbb{Z},\underset{I}{\mathrm{homotopycolimit}}~\mathrm{Spec}^\mathrm{CS}\widetilde{\Phi}^I_{\psi,\circ}/\mathrm{Fro}^\mathbb{Z},	\\
}
\]
\[ 
\xymatrix@R+0pc@C+0pc{
\underset{r}{\mathrm{homotopylimit}}~\mathrm{Spec}^\mathrm{CS}\breve{\Phi}^r_{\psi,\circ}/\mathrm{Fro}^\mathbb{Z},\underset{I}{\mathrm{homotopycolimit}}~\breve{\Phi}^I_{\psi,\circ}/\mathrm{Fro}^\mathbb{Z},	\\
}
\]
\[ 
\xymatrix@R+0pc@C+0pc{
\underset{r}{\mathrm{homotopylimit}}~\mathrm{Spec}^\mathrm{CS}{\Phi}^r_{\psi,\circ}/\mathrm{Fro}^\mathbb{Z},\underset{I}{\mathrm{homotopycolimit}}~\mathrm{Spec}^\mathrm{CS}{\Phi}^I_{\psi,\circ}/\mathrm{Fro}^\mathbb{Z}.	
}
\]

\end{definition}

\begin{proposition}
There is a well-defined functor from the $\infty$-category 
\begin{align}
\mathrm{Quasicoherentpresheaves,Perfectcomplex,Condensed}_{*}	
\end{align}
where $*$ is one of the following spaces:
\begin{align}
&\mathrm{Spec}^\mathrm{CS}\widetilde{\Phi}_{\psi,\circ}/\mathrm{Fro}^\mathbb{Z},	\\
\end{align}
\begin{align}
&\mathrm{Spec}^\mathrm{CS}\breve{\Phi}_{\psi,\circ}/\mathrm{Fro}^\mathbb{Z},	\\
\end{align}
\begin{align}
&\mathrm{Spec}^\mathrm{CS}{\Phi}_{\psi,\circ}/\mathrm{Fro}^\mathbb{Z},	
\end{align}
to the $\infty$-category of $\mathrm{Fro}$-equivariant quasicoherent presheaves over similar spaces above correspondingly without the $\mathrm{Fro}$-quotients, and to the $\infty$-category of $\mathrm{Fro}$-equivariant quasicoherent modules over global sections of the structure $\infty$-sheaves of the similar spaces above correspondingly without the $\mathrm{Fro}$-quotients. Here for those space without notation related to the radius and the corresponding interval we consider the total unions $\bigcap_r,\bigcup_I$ in order to achieve the whole spaces to achieve the analogues of the corresponding FF curves from \cite{10KL1}, \cite{10KL2}, \cite{10FF} for
\[
\xymatrix@R+0pc@C+0pc{
\underset{r}{\mathrm{homotopylimit}}~\mathrm{Spec}^\mathrm{CS}\widetilde{\Phi}^r_{\psi,\circ},\underset{I}{\mathrm{homotopycolimit}}~\mathrm{Spec}^\mathrm{CS}\widetilde{\Phi}^I_{\psi,\circ},	\\
}
\]
\[
\xymatrix@R+0pc@C+0pc{
\underset{r}{\mathrm{homotopylimit}}~\mathrm{Spec}^\mathrm{CS}\breve{\Phi}^r_{\psi,\circ},\underset{I}{\mathrm{homotopycolimit}}~\mathrm{Spec}^\mathrm{CS}\breve{\Phi}^I_{\psi,\circ},	\\
}
\]
\[
\xymatrix@R+0pc@C+0pc{
\underset{r}{\mathrm{homotopylimit}}~\mathrm{Spec}^\mathrm{CS}{\Phi}^r_{\psi,\circ},\underset{I}{\mathrm{homotopycolimit}}~\mathrm{Spec}^\mathrm{CS}{\Phi}^I_{\psi,\circ}.	
}
\]
\[ 
\xymatrix@R+0pc@C+0pc{
\underset{r}{\mathrm{homotopylimit}}~\mathrm{Spec}^\mathrm{CS}\widetilde{\Phi}^r_{\psi,\circ}/\mathrm{Fro}^\mathbb{Z},\underset{I}{\mathrm{homotopycolimit}}~\mathrm{Spec}^\mathrm{CS}\widetilde{\Phi}^I_{\psi,\circ}/\mathrm{Fro}^\mathbb{Z},	\\
}
\]
\[ 
\xymatrix@R+0pc@C+0pc{
\underset{r}{\mathrm{homotopylimit}}~\mathrm{Spec}^\mathrm{CS}\breve{\Phi}^r_{\psi,\circ}/\mathrm{Fro}^\mathbb{Z},\underset{I}{\mathrm{homotopycolimit}}~\breve{\Phi}^I_{\psi,\circ}/\mathrm{Fro}^\mathbb{Z},	\\
}
\]
\[ 
\xymatrix@R+0pc@C+0pc{
\underset{r}{\mathrm{homotopylimit}}~\mathrm{Spec}^\mathrm{CS}{\Phi}^r_{\psi,\circ}/\mathrm{Fro}^\mathbb{Z},\underset{I}{\mathrm{homotopycolimit}}~\mathrm{Spec}^\mathrm{CS}{\Phi}^I_{\psi,\circ}/\mathrm{Fro}^\mathbb{Z}.	
}
\]	
In this situation we will have the target category being family parametrized by $r$ or $I$ in compatible glueing sense as in \cite[Definition 5.4.10]{10KL2}. In this situation for modules parametrized by the intervals we have the equivalence of $\infty$-categories by using \cite[Proposition 12.18]{10CS2}. Here the corresponding quasicoherent Frobenius modules are defined to be the corresponding homotopy colimits and limits of Frobenius modules:
\begin{align}
\underset{r}{\mathrm{homotopycolimit}}~M_r,\\
\underset{I}{\mathrm{homotopylimit}}~M_I,	
\end{align}
where each $M_r$ is a Frobenius-equivariant module over the period ring with respect to some radius $r$ while each $M_I$ is a Frobenius-equivariant module over the period ring with respect to some interval $I$.\\
\end{proposition}

\begin{proposition}
Similar proposition holds for 
\begin{align}
\mathrm{Quasicoherentsheaves,Perfectcomplex,IndBanach}_{*}.	
\end{align}	
\end{proposition}

\subsection{Multivariate Hodge Iwasawa Prestacks}

\subsubsection{Frobenius Quasicoherent Prestacks I}

\begin{definition}
We now consider the pro-\'etale site of $\mathrm{Spa}\mathbb{Q}_p\left<X_1^{\pm 1},...,X_k^{\pm 1}\right>$ from \cite{10Sch}, denote that by $*$. To be more accurate we replace one component for $\Gamma$ with the pro-\'etale site of $\mathrm{Spa}\mathbb{Q}_p\left<X_1^{\pm 1},...,X_k^{\pm 1}\right>$. And we treat then all the functor to be prestacks for this site\footnote{Here for those imperfect rings, the notation will mean that the specific component forming the pro-\'etale site will be the perfect version of the corresponding ring. Certainly if we have $|\Gamma|=1$ then we have that all the rings are perfect in \cite{10KL1} and \cite{10KL2}.}. Then from \cite{10KL1} and \cite[Definition 5.2.1]{10KL2} we have the following class of Kedlaya-Liu rings (with the following replacement: $\Delta$ stands for $A$, $\nabla$ stands for $B$, while $\Phi$ stands for $C$) by taking product in the sense of self $\Gamma$-th power:

\[
\xymatrix@R+0pc@C+0pc{
\widetilde{\Delta}_{*,\Gamma},\widetilde{\nabla}_{*,\Gamma},\widetilde{\Phi}_{*,\Gamma},\widetilde{\Delta}^+_{*,\Gamma},\widetilde{\nabla}^+_{*,\Gamma},\widetilde{\Delta}^\dagger_{*,\Gamma},\widetilde{\nabla}^\dagger_{*,\Gamma},\widetilde{\Phi}^r_{*,\Gamma},\widetilde{\Phi}^I_{*,\Gamma}, 
}
\]

\[
\xymatrix@R+0pc@C+0pc{
\breve{\Delta}_{*,\Gamma},\breve{\nabla}_{*,\Gamma},\breve{\Phi}_{*,\Gamma},\breve{\Delta}^+_{*,\Gamma},\breve{\nabla}^+_{*,\Gamma},\breve{\Delta}^\dagger_{*,\Gamma},\breve{\nabla}^\dagger_{*,\Gamma},\breve{\Phi}^r_{*,\Gamma},\breve{\Phi}^I_{*,\Gamma},	
}
\]

\[
\xymatrix@R+0pc@C+0pc{
{\Delta}_{*,\Gamma},{\nabla}_{*,\Gamma},{\Phi}_{*,\Gamma},{\Delta}^+_{*,\Gamma},{\nabla}^+_{*,\Gamma},{\Delta}^\dagger_{*,\Gamma},{\nabla}^\dagger_{*,\Gamma},{\Phi}^r_{*,\Gamma},{\Phi}^I_{*,\Gamma}.	
}
\]
Taking the product we have:
\[
\xymatrix@R+0pc@C+0pc{
\widetilde{\Phi}_{*,\Gamma,X},\widetilde{\Phi}^r_{*,\Gamma,X},\widetilde{\Phi}^I_{*,\Gamma,X},	
}
\]
\[
\xymatrix@R+0pc@C+0pc{
\breve{\Phi}_{*,\Gamma,X},\breve{\Phi}^r_{*,\Gamma,X},\breve{\Phi}^I_{*,\Gamma,X},	
}
\]
\[
\xymatrix@R+0pc@C+0pc{
{\Phi}_{*,\Gamma,X},{\Phi}^r_{*,\Gamma,X},{\Phi}^I_{*,\Gamma,X}.	
}
\]
They carry multi Frobenius action $\varphi_\Gamma$ and multi $\mathrm{Lie}_\Gamma:=\mathbb{Z}_p^{\times\Gamma}$ action. In our current situation after \cite{10CKZ} and \cite{10PZ} we consider the following $(\infty,1)$-categories of $(\infty,1)$-modules.\\
\end{definition}

\begin{definition}
First we consider the Bambozzi-Kremnizer spectrum $\mathrm{Spec}^\mathrm{BK}(*)$ attached to any of those in the above from \cite{10BK} by taking derived rational localization:
\begin{align}
&\mathrm{Spec}^\mathrm{BK}\widetilde{\Phi}_{*,\Gamma,X},\mathrm{Spec}^\mathrm{BK}\widetilde{\Phi}^r_{*,\Gamma,X},\mathrm{Spec}^\mathrm{BK}\widetilde{\Phi}^I_{*,\Gamma,X},	
\end{align}
\begin{align}
&\mathrm{Spec}^\mathrm{BK}\breve{\Phi}_{*,\Gamma,X},\mathrm{Spec}^\mathrm{BK}\breve{\Phi}^r_{*,\Gamma,X},\mathrm{Spec}^\mathrm{BK}\breve{\Phi}^I_{*,\Gamma,X},	
\end{align}
\begin{align}
&\mathrm{Spec}^\mathrm{BK}{\Phi}_{*,\Gamma,X},
\mathrm{Spec}^\mathrm{BK}{\Phi}^r_{*,\Gamma,X},\mathrm{Spec}^\mathrm{BK}{\Phi}^I_{*,\Gamma,X}.	
\end{align}

Then we take the corresponding quotients by using the corresponding Frobenius operators:
\begin{align}
&\mathrm{Spec}^\mathrm{BK}\widetilde{\Phi}_{*,\Gamma,X}/\mathrm{Fro}^\mathbb{Z},	\\
\end{align}
\begin{align}
&\mathrm{Spec}^\mathrm{BK}\breve{\Phi}_{*,\Gamma,X}/\mathrm{Fro}^\mathbb{Z},	\\
\end{align}
\begin{align}
&\mathrm{Spec}^\mathrm{BK}{\Phi}_{*,\Gamma,X}/\mathrm{Fro}^\mathbb{Z}.	
\end{align}
Here for those space without notation related to the radius and the corresponding interval we consider the total unions $\bigcap_r,\bigcup_I$ in order to achieve the whole spaces to achieve the analogues of the corresponding FF curves from \cite{10KL1}, \cite{10KL2}, \cite{10FF} for
\[
\xymatrix@R+0pc@C+0pc{
\underset{r}{\mathrm{homotopylimit}}~\mathrm{Spec}^\mathrm{BK}\widetilde{\Phi}^r_{*,\Gamma,X},\underset{I}{\mathrm{homotopycolimit}}~\mathrm{Spec}^\mathrm{BK}\widetilde{\Phi}^I_{*,\Gamma,X},	\\
}
\]
\[
\xymatrix@R+0pc@C+0pc{
\underset{r}{\mathrm{homotopylimit}}~\mathrm{Spec}^\mathrm{BK}\breve{\Phi}^r_{*,\Gamma,X},\underset{I}{\mathrm{homotopycolimit}}~\mathrm{Spec}^\mathrm{BK}\breve{\Phi}^I_{*,\Gamma,X},	\\
}
\]
\[
\xymatrix@R+0pc@C+0pc{
\underset{r}{\mathrm{homotopylimit}}~\mathrm{Spec}^\mathrm{BK}{\Phi}^r_{*,\Gamma,X},\underset{I}{\mathrm{homotopycolimit}}~\mathrm{Spec}^\mathrm{BK}{\Phi}^I_{*,\Gamma,X}.	
}
\]
\[  
\xymatrix@R+0pc@C+0pc{
\underset{r}{\mathrm{homotopylimit}}~\mathrm{Spec}^\mathrm{BK}\widetilde{\Phi}^r_{*,\Gamma,X}/\mathrm{Fro}^\mathbb{Z},\underset{I}{\mathrm{homotopycolimit}}~\mathrm{Spec}^\mathrm{BK}\widetilde{\Phi}^I_{*,\Gamma,X}/\mathrm{Fro}^\mathbb{Z},	\\
}
\]
\[ 
\xymatrix@R+0pc@C+0pc{
\underset{r}{\mathrm{homotopylimit}}~\mathrm{Spec}^\mathrm{BK}\breve{\Phi}^r_{*,\Gamma,X}/\mathrm{Fro}^\mathbb{Z},\underset{I}{\mathrm{homotopycolimit}}~\mathrm{Spec}^\mathrm{BK}\breve{\Phi}^I_{*,\Gamma,X}/\mathrm{Fro}^\mathbb{Z},	\\
}
\]
\[ 
\xymatrix@R+0pc@C+0pc{
\underset{r}{\mathrm{homotopylimit}}~\mathrm{Spec}^\mathrm{BK}{\Phi}^r_{*,\Gamma,X}/\mathrm{Fro}^\mathbb{Z},\underset{I}{\mathrm{homotopycolimit}}~\mathrm{Spec}^\mathrm{BK}{\Phi}^I_{*,\Gamma,X}/\mathrm{Fro}^\mathbb{Z}.	
}
\]

\end{definition}

\indent Meanwhile we have the corresponding Clausen-Scholze analytic stacks from \cite{10CS2}, therefore applying their construction we have:

\begin{definition}
Here we define the following products by using the solidified tensor product from \cite{10CS1} and \cite{10CS2}. Then we take solidified tensor product $\overset{\blacksquare}{\otimes}$ of any of the following
\[
\xymatrix@R+0pc@C+0pc{
\widetilde{\Delta}_{*,\Gamma},\widetilde{\nabla}_{*,\Gamma},\widetilde{\Phi}_{*,\Gamma},\widetilde{\Delta}^+_{*,\Gamma},\widetilde{\nabla}^+_{*,\Gamma},\widetilde{\Delta}^\dagger_{*,\Gamma},\widetilde{\nabla}^\dagger_{*,\Gamma},\widetilde{\Phi}^r_{*,\Gamma},\widetilde{\Phi}^I_{*,\Gamma}, 
}
\]

\[
\xymatrix@R+0pc@C+0pc{
\breve{\Delta}_{*,\Gamma},\breve{\nabla}_{*,\Gamma},\breve{\Phi}_{*,\Gamma},\breve{\Delta}^+_{*,\Gamma},\breve{\nabla}^+_{*,\Gamma},\breve{\Delta}^\dagger_{*,\Gamma},\breve{\nabla}^\dagger_{*,\Gamma},\breve{\Phi}^r_{*,\Gamma},\breve{\Phi}^I_{*,\Gamma},	
}
\]

\[
\xymatrix@R+0pc@C+0pc{
{\Delta}_{*,\Gamma},{\nabla}_{*,\Gamma},{\Phi}_{*,\Gamma},{\Delta}^+_{*,\Gamma},{\nabla}^+_{*,\Gamma},{\Delta}^\dagger_{*,\Gamma},{\nabla}^\dagger_{*,\Gamma},{\Phi}^r_{*,\Gamma},{\Phi}^I_{*,\Gamma},	
}
\]  	
with $X$. Then we have the notations:
\[
\xymatrix@R+0pc@C+0pc{
\widetilde{\Delta}_{*,\Gamma,X},\widetilde{\nabla}_{*,\Gamma,X},\widetilde{\Phi}_{*,\Gamma,X},\widetilde{\Delta}^+_{*,\Gamma,X},\widetilde{\nabla}^+_{*,\Gamma,X},\widetilde{\Delta}^\dagger_{*,\Gamma,X},\widetilde{\nabla}^\dagger_{*,\Gamma,X},\widetilde{\Phi}^r_{*,\Gamma,X},\widetilde{\Phi}^I_{*,\Gamma,X}, 
}
\]

\[
\xymatrix@R+0pc@C+0pc{
\breve{\Delta}_{*,\Gamma,X},\breve{\nabla}_{*,\Gamma,X},\breve{\Phi}_{*,\Gamma,X},\breve{\Delta}^+_{*,\Gamma,X},\breve{\nabla}^+_{*,\Gamma,X},\breve{\Delta}^\dagger_{*,\Gamma,X},\breve{\nabla}^\dagger_{*,\Gamma,X},\breve{\Phi}^r_{*,\Gamma,X},\breve{\Phi}^I_{*,\Gamma,X},	
}
\]

\[
\xymatrix@R+0pc@C+0pc{
{\Delta}_{*,\Gamma,X},{\nabla}_{*,\Gamma,X},{\Phi}_{*,\Gamma,X},{\Delta}^+_{*,\Gamma,X},{\nabla}^+_{*,\Gamma,X},{\Delta}^\dagger_{*,\Gamma,X},{\nabla}^\dagger_{*,\Gamma,X},{\Phi}^r_{*,\Gamma,X},{\Phi}^I_{*,\Gamma,X}.	
}
\]
\end{definition}

\begin{definition}
First we consider the Clausen-Scholze spectrum $\mathrm{Spec}^\mathrm{CS}(*)$ attached to any of those in the above from \cite{10CS2} by taking derived rational localization:
\begin{align}
\mathrm{Spec}^\mathrm{CS}\widetilde{\Delta}_{*,\Gamma,X},\mathrm{Spec}^\mathrm{CS}\widetilde{\nabla}_{*,\Gamma,X},\mathrm{Spec}^\mathrm{CS}\widetilde{\Phi}_{*,\Gamma,X},\mathrm{Spec}^\mathrm{CS}\widetilde{\Delta}^+_{*,\Gamma,X},\mathrm{Spec}^\mathrm{CS}\widetilde{\nabla}^+_{*,\Gamma,X},\\
\mathrm{Spec}^\mathrm{CS}\widetilde{\Delta}^\dagger_{*,\Gamma,X},\mathrm{Spec}^\mathrm{CS}\widetilde{\nabla}^\dagger_{*,\Gamma,X},\mathrm{Spec}^\mathrm{CS}\widetilde{\Phi}^r_{*,\Gamma,X},\mathrm{Spec}^\mathrm{CS}\widetilde{\Phi}^I_{*,\Gamma,X},	\\
\end{align}
\begin{align}
\mathrm{Spec}^\mathrm{CS}\breve{\Delta}_{*,\Gamma,X},\breve{\nabla}_{*,\Gamma,X},\mathrm{Spec}^\mathrm{CS}\breve{\Phi}_{*,\Gamma,X},\mathrm{Spec}^\mathrm{CS}\breve{\Delta}^+_{*,\Gamma,X},\mathrm{Spec}^\mathrm{CS}\breve{\nabla}^+_{*,\Gamma,X},\\
\mathrm{Spec}^\mathrm{CS}\breve{\Delta}^\dagger_{*,\Gamma,X},\mathrm{Spec}^\mathrm{CS}\breve{\nabla}^\dagger_{*,\Gamma,X},\mathrm{Spec}^\mathrm{CS}\breve{\Phi}^r_{*,\Gamma,X},\breve{\Phi}^I_{*,\Gamma,X},	\\
\end{align}
\begin{align}
\mathrm{Spec}^\mathrm{CS}{\Delta}_{*,\Gamma,X},\mathrm{Spec}^\mathrm{CS}{\nabla}_{*,\Gamma,X},\mathrm{Spec}^\mathrm{CS}{\Phi}_{*,\Gamma,X},\mathrm{Spec}^\mathrm{CS}{\Delta}^+_{*,\Gamma,X},\mathrm{Spec}^\mathrm{CS}{\nabla}^+_{*,\Gamma,X},\\
\mathrm{Spec}^\mathrm{CS}{\Delta}^\dagger_{*,\Gamma,X},\mathrm{Spec}^\mathrm{CS}{\nabla}^\dagger_{*,\Gamma,X},\mathrm{Spec}^\mathrm{CS}{\Phi}^r_{*,\Gamma,X},\mathrm{Spec}^\mathrm{CS}{\Phi}^I_{*,\Gamma,X}.	
\end{align}

Then we take the corresponding quotients by using the corresponding Frobenius operators:
\begin{align}
&\mathrm{Spec}^\mathrm{CS}\widetilde{\Delta}_{*,\Gamma,X}/\mathrm{Fro}^\mathbb{Z},\mathrm{Spec}^\mathrm{CS}\widetilde{\nabla}_{*,\Gamma,X}/\mathrm{Fro}^\mathbb{Z},\mathrm{Spec}^\mathrm{CS}\widetilde{\Phi}_{*,\Gamma,X}/\mathrm{Fro}^\mathbb{Z},\mathrm{Spec}^\mathrm{CS}\widetilde{\Delta}^+_{*,\Gamma,X}/\mathrm{Fro}^\mathbb{Z},\\
&\mathrm{Spec}^\mathrm{CS}\widetilde{\nabla}^+_{*,\Gamma,X}/\mathrm{Fro}^\mathbb{Z}, \mathrm{Spec}^\mathrm{CS}\widetilde{\Delta}^\dagger_{*,\Gamma,X}/\mathrm{Fro}^\mathbb{Z},\mathrm{Spec}^\mathrm{CS}\widetilde{\nabla}^\dagger_{*,\Gamma,X}/\mathrm{Fro}^\mathbb{Z},	\\
\end{align}
\begin{align}
&\mathrm{Spec}^\mathrm{CS}\breve{\Delta}_{*,\Gamma,X}/\mathrm{Fro}^\mathbb{Z},\breve{\nabla}_{*,\Gamma,X}/\mathrm{Fro}^\mathbb{Z},\mathrm{Spec}^\mathrm{CS}\breve{\Phi}_{*,\Gamma,X}/\mathrm{Fro}^\mathbb{Z},\mathrm{Spec}^\mathrm{CS}\breve{\Delta}^+_{*,\Gamma,X}/\mathrm{Fro}^\mathbb{Z},\\
&\mathrm{Spec}^\mathrm{CS}\breve{\nabla}^+_{*,\Gamma,X}/\mathrm{Fro}^\mathbb{Z}, \mathrm{Spec}^\mathrm{CS}\breve{\Delta}^\dagger_{*,\Gamma,X}/\mathrm{Fro}^\mathbb{Z},\mathrm{Spec}^\mathrm{CS}\breve{\nabla}^\dagger_{*,\Gamma,X}/\mathrm{Fro}^\mathbb{Z},	\\
\end{align}
\begin{align}
&\mathrm{Spec}^\mathrm{CS}{\Delta}_{*,\Gamma,X}/\mathrm{Fro}^\mathbb{Z},\mathrm{Spec}^\mathrm{CS}{\nabla}_{*,\Gamma,X}/\mathrm{Fro}^\mathbb{Z},\mathrm{Spec}^\mathrm{CS}{\Phi}_{*,\Gamma,X}/\mathrm{Fro}^\mathbb{Z},\mathrm{Spec}^\mathrm{CS}{\Delta}^+_{*,\Gamma,X}/\mathrm{Fro}^\mathbb{Z},\\
&\mathrm{Spec}^\mathrm{CS}{\nabla}^+_{*,\Gamma,X}/\mathrm{Fro}^\mathbb{Z}, \mathrm{Spec}^\mathrm{CS}{\Delta}^\dagger_{*,\Gamma,X}/\mathrm{Fro}^\mathbb{Z},\mathrm{Spec}^\mathrm{CS}{\nabla}^\dagger_{*,\Gamma,X}/\mathrm{Fro}^\mathbb{Z}.	
\end{align}
Here for those space with notations related to the radius and the corresponding interval we consider the total unions $\bigcap_r,\bigcup_I$ in order to achieve the whole spaces to achieve the analogues of the corresponding FF curves from \cite{10KL1}, \cite{10KL2}, \cite{10FF} for
\[
\xymatrix@R+0pc@C+0pc{
\underset{r}{\mathrm{homotopylimit}}~\mathrm{Spec}^\mathrm{CS}\widetilde{\Phi}^r_{*,\Gamma,X},\underset{I}{\mathrm{homotopycolimit}}~\mathrm{Spec}^\mathrm{CS}\widetilde{\Phi}^I_{*,\Gamma,X},	\\
}
\]
\[
\xymatrix@R+0pc@C+0pc{
\underset{r}{\mathrm{homotopylimit}}~\mathrm{Spec}^\mathrm{CS}\breve{\Phi}^r_{*,\Gamma,X},\underset{I}{\mathrm{homotopycolimit}}~\mathrm{Spec}^\mathrm{CS}\breve{\Phi}^I_{*,\Gamma,X},	\\
}
\]
\[
\xymatrix@R+0pc@C+0pc{
\underset{r}{\mathrm{homotopylimit}}~\mathrm{Spec}^\mathrm{CS}{\Phi}^r_{*,\Gamma,X},\underset{I}{\mathrm{homotopycolimit}}~\mathrm{Spec}^\mathrm{CS}{\Phi}^I_{*,\Gamma,X}.	
}
\]
\[ 
\xymatrix@R+0pc@C+0pc{
\underset{r}{\mathrm{homotopylimit}}~\mathrm{Spec}^\mathrm{CS}\widetilde{\Phi}^r_{*,\Gamma,X}/\mathrm{Fro}^\mathbb{Z},\underset{I}{\mathrm{homotopycolimit}}~\mathrm{Spec}^\mathrm{CS}\widetilde{\Phi}^I_{*,\Gamma,X}/\mathrm{Fro}^\mathbb{Z},	\\
}
\]
\[ 
\xymatrix@R+0pc@C+0pc{
\underset{r}{\mathrm{homotopylimit}}~\mathrm{Spec}^\mathrm{CS}\breve{\Phi}^r_{*,\Gamma,X}/\mathrm{Fro}^\mathbb{Z},\underset{I}{\mathrm{homotopycolimit}}~\breve{\Phi}^I_{*,\Gamma,X}/\mathrm{Fro}^\mathbb{Z},	\\
}
\]
\[ 
\xymatrix@R+0pc@C+0pc{
\underset{r}{\mathrm{homotopylimit}}~\mathrm{Spec}^\mathrm{CS}{\Phi}^r_{*,\Gamma,X}/\mathrm{Fro}^\mathbb{Z},\underset{I}{\mathrm{homotopycolimit}}~\mathrm{Spec}^\mathrm{CS}{\Phi}^I_{*,\Gamma,X}/\mathrm{Fro}^\mathbb{Z}.	
}
\]

\end{definition}

\

\begin{definition}
We then consider the corresponding quasipresheaves of the corresponding ind-Banach or monomorphic ind-Banach modules from \cite{10BBK}, \cite{10KKM}:
\begin{align}
\mathrm{Quasicoherentpresheaves,IndBanach}_{*}	
\end{align}
where $*$ is one of the following spaces:
\begin{align}
&\mathrm{Spec}^\mathrm{BK}\widetilde{\Phi}_{*,\Gamma,X}/\mathrm{Fro}^\mathbb{Z},	\\
\end{align}
\begin{align}
&\mathrm{Spec}^\mathrm{BK}\breve{\Phi}_{*,\Gamma,X}/\mathrm{Fro}^\mathbb{Z},	\\
\end{align}
\begin{align}
&\mathrm{Spec}^\mathrm{BK}{\Phi}_{*,\Gamma,X}/\mathrm{Fro}^\mathbb{Z}.	
\end{align}
Here for those space without notation related to the radius and the corresponding interval we consider the total unions $\bigcap_r,\bigcup_I$ in order to achieve the whole spaces to achieve the analogues of the corresponding FF curves from \cite{10KL1}, \cite{10KL2}, \cite{10FF} for
\[
\xymatrix@R+0pc@C+0pc{
\underset{r}{\mathrm{homotopylimit}}~\mathrm{Spec}^\mathrm{BK}\widetilde{\Phi}^r_{*,\Gamma,X},\underset{I}{\mathrm{homotopycolimit}}~\mathrm{Spec}^\mathrm{BK}\widetilde{\Phi}^I_{*,\Gamma,X},	\\
}
\]
\[
\xymatrix@R+0pc@C+0pc{
\underset{r}{\mathrm{homotopylimit}}~\mathrm{Spec}^\mathrm{BK}\breve{\Phi}^r_{*,\Gamma,X},\underset{I}{\mathrm{homotopycolimit}}~\mathrm{Spec}^\mathrm{BK}\breve{\Phi}^I_{*,\Gamma,X},	\\
}
\]
\[
\xymatrix@R+0pc@C+0pc{
\underset{r}{\mathrm{homotopylimit}}~\mathrm{Spec}^\mathrm{BK}{\Phi}^r_{*,\Gamma,X},\underset{I}{\mathrm{homotopycolimit}}~\mathrm{Spec}^\mathrm{BK}{\Phi}^I_{*,\Gamma,X}.	
}
\]
\[  
\xymatrix@R+0pc@C+0pc{
\underset{r}{\mathrm{homotopylimit}}~\mathrm{Spec}^\mathrm{BK}\widetilde{\Phi}^r_{*,\Gamma,X}/\mathrm{Fro}^\mathbb{Z},\underset{I}{\mathrm{homotopycolimit}}~\mathrm{Spec}^\mathrm{BK}\widetilde{\Phi}^I_{*,\Gamma,X}/\mathrm{Fro}^\mathbb{Z},	\\
}
\]
\[ 
\xymatrix@R+0pc@C+0pc{
\underset{r}{\mathrm{homotopylimit}}~\mathrm{Spec}^\mathrm{BK}\breve{\Phi}^r_{*,\Gamma,X}/\mathrm{Fro}^\mathbb{Z},\underset{I}{\mathrm{homotopycolimit}}~\mathrm{Spec}^\mathrm{BK}\breve{\Phi}^I_{*,\Gamma,X}/\mathrm{Fro}^\mathbb{Z},	\\
}
\]
\[ 
\xymatrix@R+0pc@C+0pc{
\underset{r}{\mathrm{homotopylimit}}~\mathrm{Spec}^\mathrm{BK}{\Phi}^r_{*,\Gamma,X}/\mathrm{Fro}^\mathbb{Z},\underset{I}{\mathrm{homotopycolimit}}~\mathrm{Spec}^\mathrm{BK}{\Phi}^I_{*,\Gamma,X}/\mathrm{Fro}^\mathbb{Z}.	
}
\]

\end{definition}

\begin{definition}
We then consider the corresponding quasisheaves of the corresponding condensed solid topological modules from \cite{10CS2}:
\begin{align}
\mathrm{Quasicoherentsheaves, Condensed}_{*}	
\end{align}
where $*$ is one of the following spaces:
\begin{align}
&\mathrm{Spec}^\mathrm{CS}\widetilde{\Delta}_{*,\Gamma,X}/\mathrm{Fro}^\mathbb{Z},\mathrm{Spec}^\mathrm{CS}\widetilde{\nabla}_{*,\Gamma,X}/\mathrm{Fro}^\mathbb{Z},\mathrm{Spec}^\mathrm{CS}\widetilde{\Phi}_{*,\Gamma,X}/\mathrm{Fro}^\mathbb{Z},\mathrm{Spec}^\mathrm{CS}\widetilde{\Delta}^+_{*,\Gamma,X}/\mathrm{Fro}^\mathbb{Z},\\
&\mathrm{Spec}^\mathrm{CS}\widetilde{\nabla}^+_{*,\Gamma,X}/\mathrm{Fro}^\mathbb{Z},\mathrm{Spec}^\mathrm{CS}\widetilde{\Delta}^\dagger_{*,\Gamma,X}/\mathrm{Fro}^\mathbb{Z},\mathrm{Spec}^\mathrm{CS}\widetilde{\nabla}^\dagger_{*,\Gamma,X}/\mathrm{Fro}^\mathbb{Z},	\\
\end{align}
\begin{align}
&\mathrm{Spec}^\mathrm{CS}\breve{\Delta}_{*,\Gamma,X}/\mathrm{Fro}^\mathbb{Z},\breve{\nabla}_{*,\Gamma,X}/\mathrm{Fro}^\mathbb{Z},\mathrm{Spec}^\mathrm{CS}\breve{\Phi}_{*,\Gamma,X}/\mathrm{Fro}^\mathbb{Z},\mathrm{Spec}^\mathrm{CS}\breve{\Delta}^+_{*,\Gamma,X}/\mathrm{Fro}^\mathbb{Z},\\
&\mathrm{Spec}^\mathrm{CS}\breve{\nabla}^+_{*,\Gamma,X}/\mathrm{Fro}^\mathbb{Z},\mathrm{Spec}^\mathrm{CS}\breve{\Delta}^\dagger_{*,\Gamma,X}/\mathrm{Fro}^\mathbb{Z},\mathrm{Spec}^\mathrm{CS}\breve{\nabla}^\dagger_{*,\Gamma,X}/\mathrm{Fro}^\mathbb{Z},	\\
\end{align}
\begin{align}
&\mathrm{Spec}^\mathrm{CS}{\Delta}_{*,\Gamma,X}/\mathrm{Fro}^\mathbb{Z},\mathrm{Spec}^\mathrm{CS}{\nabla}_{*,\Gamma,X}/\mathrm{Fro}^\mathbb{Z},\mathrm{Spec}^\mathrm{CS}{\Phi}_{*,\Gamma,X}/\mathrm{Fro}^\mathbb{Z},\mathrm{Spec}^\mathrm{CS}{\Delta}^+_{*,\Gamma,X}/\mathrm{Fro}^\mathbb{Z},\\
&\mathrm{Spec}^\mathrm{CS}{\nabla}^+_{*,\Gamma,X}/\mathrm{Fro}^\mathbb{Z}, \mathrm{Spec}^\mathrm{CS}{\Delta}^\dagger_{*,\Gamma,X}/\mathrm{Fro}^\mathbb{Z},\mathrm{Spec}^\mathrm{CS}{\nabla}^\dagger_{*,\Gamma,X}/\mathrm{Fro}^\mathbb{Z}.	
\end{align}
Here for those space with notations related to the radius and the corresponding interval we consider the total unions $\bigcap_r,\bigcup_I$ in order to achieve the whole spaces to achieve the analogues of the corresponding FF curves from \cite{10KL1}, \cite{10KL2}, \cite{10FF} for
\[
\xymatrix@R+0pc@C+0pc{
\underset{r}{\mathrm{homotopylimit}}~\mathrm{Spec}^\mathrm{CS}\widetilde{\Phi}^r_{*,\Gamma,X},\underset{I}{\mathrm{homotopycolimit}}~\mathrm{Spec}^\mathrm{CS}\widetilde{\Phi}^I_{*,\Gamma,X},	\\
}
\]
\[
\xymatrix@R+0pc@C+0pc{
\underset{r}{\mathrm{homotopylimit}}~\mathrm{Spec}^\mathrm{CS}\breve{\Phi}^r_{*,\Gamma,X},\underset{I}{\mathrm{homotopycolimit}}~\mathrm{Spec}^\mathrm{CS}\breve{\Phi}^I_{*,\Gamma,X},	\\
}
\]
\[
\xymatrix@R+0pc@C+0pc{
\underset{r}{\mathrm{homotopylimit}}~\mathrm{Spec}^\mathrm{CS}{\Phi}^r_{*,\Gamma,X},\underset{I}{\mathrm{homotopycolimit}}~\mathrm{Spec}^\mathrm{CS}{\Phi}^I_{*,\Gamma,X}.	
}
\]
\[ 
\xymatrix@R+0pc@C+0pc{
\underset{r}{\mathrm{homotopylimit}}~\mathrm{Spec}^\mathrm{CS}\widetilde{\Phi}^r_{*,\Gamma,X}/\mathrm{Fro}^\mathbb{Z},\underset{I}{\mathrm{homotopycolimit}}~\mathrm{Spec}^\mathrm{CS}\widetilde{\Phi}^I_{*,\Gamma,X}/\mathrm{Fro}^\mathbb{Z},	\\
}
\]
\[ 
\xymatrix@R+0pc@C+0pc{
\underset{r}{\mathrm{homotopylimit}}~\mathrm{Spec}^\mathrm{CS}\breve{\Phi}^r_{*,\Gamma,X}/\mathrm{Fro}^\mathbb{Z},\underset{I}{\mathrm{homotopycolimit}}~\breve{\Phi}^I_{*,\Gamma,X}/\mathrm{Fro}^\mathbb{Z},	\\
}
\]
\[ 
\xymatrix@R+0pc@C+0pc{
\underset{r}{\mathrm{homotopylimit}}~\mathrm{Spec}^\mathrm{CS}{\Phi}^r_{*,\Gamma,X}/\mathrm{Fro}^\mathbb{Z},\underset{I}{\mathrm{homotopycolimit}}~\mathrm{Spec}^\mathrm{CS}{\Phi}^I_{*,\Gamma,X}/\mathrm{Fro}^\mathbb{Z}.	
}
\]

\end{definition}

\

\begin{proposition}
There is a well-defined functor from the $\infty$-category 
\begin{align}
\mathrm{Quasicoherentpresheaves,Condensed}_{*}	
\end{align}
where $*$ is one of the following spaces:
\begin{align}
&\mathrm{Spec}^\mathrm{CS}\widetilde{\Phi}_{*,\Gamma,X}/\mathrm{Fro}^\mathbb{Z},	\\
\end{align}
\begin{align}
&\mathrm{Spec}^\mathrm{CS}\breve{\Phi}_{*,\Gamma,X}/\mathrm{Fro}^\mathbb{Z},	\\
\end{align}
\begin{align}
&\mathrm{Spec}^\mathrm{CS}{\Phi}_{*,\Gamma,X}/\mathrm{Fro}^\mathbb{Z},	
\end{align}
to the $\infty$-category of $\mathrm{Fro}$-equivariant quasicoherent presheaves over similar spaces above correspondingly without the $\mathrm{Fro}$-quotients, and to the $\infty$-category of $\mathrm{Fro}$-equivariant quasicoherent modules over global sections of the structure $\infty$-sheaves of the similar spaces above correspondingly without the $\mathrm{Fro}$-quotients. Here for those space without notation related to the radius and the corresponding interval we consider the total unions $\bigcap_r,\bigcup_I$ in order to achieve the whole spaces to achieve the analogues of the corresponding FF curves from \cite{10KL1}, \cite{10KL2}, \cite{10FF} for
\[
\xymatrix@R+0pc@C+0pc{
\underset{r}{\mathrm{homotopylimit}}~\mathrm{Spec}^\mathrm{CS}\widetilde{\Phi}^r_{*,\Gamma,X},\underset{I}{\mathrm{homotopycolimit}}~\mathrm{Spec}^\mathrm{CS}\widetilde{\Phi}^I_{*,\Gamma,X},	\\
}
\]
\[
\xymatrix@R+0pc@C+0pc{
\underset{r}{\mathrm{homotopylimit}}~\mathrm{Spec}^\mathrm{CS}\breve{\Phi}^r_{*,\Gamma,X},\underset{I}{\mathrm{homotopycolimit}}~\mathrm{Spec}^\mathrm{CS}\breve{\Phi}^I_{*,\Gamma,X},	\\
}
\]
\[
\xymatrix@R+0pc@C+0pc{
\underset{r}{\mathrm{homotopylimit}}~\mathrm{Spec}^\mathrm{CS}{\Phi}^r_{*,\Gamma,X},\underset{I}{\mathrm{homotopycolimit}}~\mathrm{Spec}^\mathrm{CS}{\Phi}^I_{*,\Gamma,X}.	
}
\]
\[ 
\xymatrix@R+0pc@C+0pc{
\underset{r}{\mathrm{homotopylimit}}~\mathrm{Spec}^\mathrm{CS}\widetilde{\Phi}^r_{*,\Gamma,X}/\mathrm{Fro}^\mathbb{Z},\underset{I}{\mathrm{homotopycolimit}}~\mathrm{Spec}^\mathrm{CS}\widetilde{\Phi}^I_{*,\Gamma,X}/\mathrm{Fro}^\mathbb{Z},	\\
}
\]
\[ 
\xymatrix@R+0pc@C+0pc{
\underset{r}{\mathrm{homotopylimit}}~\mathrm{Spec}^\mathrm{CS}\breve{\Phi}^r_{*,\Gamma,X}/\mathrm{Fro}^\mathbb{Z},\underset{I}{\mathrm{homotopycolimit}}~\breve{\Phi}^I_{*,\Gamma,X}/\mathrm{Fro}^\mathbb{Z},	\\
}
\]
\[ 
\xymatrix@R+0pc@C+0pc{
\underset{r}{\mathrm{homotopylimit}}~\mathrm{Spec}^\mathrm{CS}{\Phi}^r_{*,\Gamma,X}/\mathrm{Fro}^\mathbb{Z},\underset{I}{\mathrm{homotopycolimit}}~\mathrm{Spec}^\mathrm{CS}{\Phi}^I_{*,\Gamma,X}/\mathrm{Fro}^\mathbb{Z}.	
}
\]	
In this situation we will have the target category being family parametrized by $r$ or $I$ in compatible glueing sense as in \cite[Definition 5.4.10]{10KL2}. In this situation for modules parametrized by the intervals we have the equivalence of $\infty$-categories by using \cite[Proposition 13.8]{10CS2}. Here the corresponding quasicoherent Frobenius modules are defined to be the corresponding homotopy colimits and limits of Frobenius modules:
\begin{align}
\underset{r}{\mathrm{homotopycolimit}}~M_r,\\
\underset{I}{\mathrm{homotopylimit}}~M_I,	
\end{align}
where each $M_r$ is a Frobenius-equivariant module over the period ring with respect to some radius $r$ while each $M_I$ is a Frobenius-equivariant module over the period ring with respect to some interval $I$.\\
\end{proposition}

\begin{proposition}
Similar proposition holds for 
\begin{align}
\mathrm{Quasicoherentsheaves,IndBanach}_{*}.	
\end{align}	
\end{proposition}

\

\begin{definition}
We then consider the corresponding quasipresheaves of perfect complexes the corresponding ind-Banach or monomorphic ind-Banach modules from \cite{10BBK}, \cite{10KKM}:
\begin{align}
\mathrm{Quasicoherentpresheaves,Perfectcomplex,IndBanach}_{*}	
\end{align}
where $*$ is one of the following spaces:
\begin{align}
&\mathrm{Spec}^\mathrm{BK}\widetilde{\Phi}_{*,\Gamma,X}/\mathrm{Fro}^\mathbb{Z},	\\
\end{align}
\begin{align}
&\mathrm{Spec}^\mathrm{BK}\breve{\Phi}_{*,\Gamma,X}/\mathrm{Fro}^\mathbb{Z},	\\
\end{align}
\begin{align}
&\mathrm{Spec}^\mathrm{BK}{\Phi}_{*,\Gamma,X}/\mathrm{Fro}^\mathbb{Z}.	
\end{align}
Here for those space without notation related to the radius and the corresponding interval we consider the total unions $\bigcap_r,\bigcup_I$ in order to achieve the whole spaces to achieve the analogues of the corresponding FF curves from \cite{10KL1}, \cite{10KL2}, \cite{10FF} for
\[
\xymatrix@R+0pc@C+0pc{
\underset{r}{\mathrm{homotopylimit}}~\mathrm{Spec}^\mathrm{BK}\widetilde{\Phi}^r_{*,\Gamma,X},\underset{I}{\mathrm{homotopycolimit}}~\mathrm{Spec}^\mathrm{BK}\widetilde{\Phi}^I_{*,\Gamma,X},	\\
}
\]
\[
\xymatrix@R+0pc@C+0pc{
\underset{r}{\mathrm{homotopylimit}}~\mathrm{Spec}^\mathrm{BK}\breve{\Phi}^r_{*,\Gamma,X},\underset{I}{\mathrm{homotopycolimit}}~\mathrm{Spec}^\mathrm{BK}\breve{\Phi}^I_{*,\Gamma,X},	\\
}
\]
\[
\xymatrix@R+0pc@C+0pc{
\underset{r}{\mathrm{homotopylimit}}~\mathrm{Spec}^\mathrm{BK}{\Phi}^r_{*,\Gamma,X},\underset{I}{\mathrm{homotopycolimit}}~\mathrm{Spec}^\mathrm{BK}{\Phi}^I_{*,\Gamma,X}.	
}
\]
\[  
\xymatrix@R+0pc@C+0pc{
\underset{r}{\mathrm{homotopylimit}}~\mathrm{Spec}^\mathrm{BK}\widetilde{\Phi}^r_{*,\Gamma,X}/\mathrm{Fro}^\mathbb{Z},\underset{I}{\mathrm{homotopycolimit}}~\mathrm{Spec}^\mathrm{BK}\widetilde{\Phi}^I_{*,\Gamma,X}/\mathrm{Fro}^\mathbb{Z},	\\
}
\]
\[ 
\xymatrix@R+0pc@C+0pc{
\underset{r}{\mathrm{homotopylimit}}~\mathrm{Spec}^\mathrm{BK}\breve{\Phi}^r_{*,\Gamma,X}/\mathrm{Fro}^\mathbb{Z},\underset{I}{\mathrm{homotopycolimit}}~\mathrm{Spec}^\mathrm{BK}\breve{\Phi}^I_{*,\Gamma,X}/\mathrm{Fro}^\mathbb{Z},	\\
}
\]
\[ 
\xymatrix@R+0pc@C+0pc{
\underset{r}{\mathrm{homotopylimit}}~\mathrm{Spec}^\mathrm{BK}{\Phi}^r_{*,\Gamma,X}/\mathrm{Fro}^\mathbb{Z},\underset{I}{\mathrm{homotopycolimit}}~\mathrm{Spec}^\mathrm{BK}{\Phi}^I_{*,\Gamma,X}/\mathrm{Fro}^\mathbb{Z}.	
}
\]

\end{definition}

\begin{definition}
We then consider the corresponding quasisheaves of perfect complexes of the corresponding condensed solid topological modules from \cite{10CS2}:
\begin{align}
\mathrm{Quasicoherentsheaves, Perfectcomplex, Condensed}_{*}	
\end{align}
where $*$ is one of the following spaces:
\begin{align}
&\mathrm{Spec}^\mathrm{CS}\widetilde{\Delta}_{*,\Gamma,X}/\mathrm{Fro}^\mathbb{Z},\mathrm{Spec}^\mathrm{CS}\widetilde{\nabla}_{*,\Gamma,X}/\mathrm{Fro}^\mathbb{Z},\mathrm{Spec}^\mathrm{CS}\widetilde{\Phi}_{*,\Gamma,X}/\mathrm{Fro}^\mathbb{Z},\mathrm{Spec}^\mathrm{CS}\widetilde{\Delta}^+_{*,\Gamma,X}/\mathrm{Fro}^\mathbb{Z},\\
&\mathrm{Spec}^\mathrm{CS}\widetilde{\nabla}^+_{*,\Gamma,X}/\mathrm{Fro}^\mathbb{Z},\mathrm{Spec}^\mathrm{CS}\widetilde{\Delta}^\dagger_{*,\Gamma,X}/\mathrm{Fro}^\mathbb{Z},\mathrm{Spec}^\mathrm{CS}\widetilde{\nabla}^\dagger_{*,\Gamma,X}/\mathrm{Fro}^\mathbb{Z},	\\
\end{align}
\begin{align}
&\mathrm{Spec}^\mathrm{CS}\breve{\Delta}_{*,\Gamma,X}/\mathrm{Fro}^\mathbb{Z},\breve{\nabla}_{*,\Gamma,X}/\mathrm{Fro}^\mathbb{Z},\mathrm{Spec}^\mathrm{CS}\breve{\Phi}_{*,\Gamma,X}/\mathrm{Fro}^\mathbb{Z},\mathrm{Spec}^\mathrm{CS}\breve{\Delta}^+_{*,\Gamma,X}/\mathrm{Fro}^\mathbb{Z},\\
&\mathrm{Spec}^\mathrm{CS}\breve{\nabla}^+_{*,\Gamma,X}/\mathrm{Fro}^\mathbb{Z},\mathrm{Spec}^\mathrm{CS}\breve{\Delta}^\dagger_{*,\Gamma,X}/\mathrm{Fro}^\mathbb{Z},\mathrm{Spec}^\mathrm{CS}\breve{\nabla}^\dagger_{*,\Gamma,X}/\mathrm{Fro}^\mathbb{Z},	\\
\end{align}
\begin{align}
&\mathrm{Spec}^\mathrm{CS}{\Delta}_{*,\Gamma,X}/\mathrm{Fro}^\mathbb{Z},\mathrm{Spec}^\mathrm{CS}{\nabla}_{*,\Gamma,X}/\mathrm{Fro}^\mathbb{Z},\mathrm{Spec}^\mathrm{CS}{\Phi}_{*,\Gamma,X}/\mathrm{Fro}^\mathbb{Z},\mathrm{Spec}^\mathrm{CS}{\Delta}^+_{*,\Gamma,X}/\mathrm{Fro}^\mathbb{Z},\\
&\mathrm{Spec}^\mathrm{CS}{\nabla}^+_{*,\Gamma,X}/\mathrm{Fro}^\mathbb{Z}, \mathrm{Spec}^\mathrm{CS}{\Delta}^\dagger_{*,\Gamma,X}/\mathrm{Fro}^\mathbb{Z},\mathrm{Spec}^\mathrm{CS}{\nabla}^\dagger_{*,\Gamma,X}/\mathrm{Fro}^\mathbb{Z}.	
\end{align}
Here for those space with notations related to the radius and the corresponding interval we consider the total unions $\bigcap_r,\bigcup_I$ in order to achieve the whole spaces to achieve the analogues of the corresponding FF curves from \cite{10KL1}, \cite{10KL2}, \cite{10FF} for
\[
\xymatrix@R+0pc@C+0pc{
\underset{r}{\mathrm{homotopylimit}}~\mathrm{Spec}^\mathrm{CS}\widetilde{\Phi}^r_{*,\Gamma,X},\underset{I}{\mathrm{homotopycolimit}}~\mathrm{Spec}^\mathrm{CS}\widetilde{\Phi}^I_{*,\Gamma,X},	\\
}
\]
\[
\xymatrix@R+0pc@C+0pc{
\underset{r}{\mathrm{homotopylimit}}~\mathrm{Spec}^\mathrm{CS}\breve{\Phi}^r_{*,\Gamma,X},\underset{I}{\mathrm{homotopycolimit}}~\mathrm{Spec}^\mathrm{CS}\breve{\Phi}^I_{*,\Gamma,X},	\\
}
\]
\[
\xymatrix@R+0pc@C+0pc{
\underset{r}{\mathrm{homotopylimit}}~\mathrm{Spec}^\mathrm{CS}{\Phi}^r_{*,\Gamma,X},\underset{I}{\mathrm{homotopycolimit}}~\mathrm{Spec}^\mathrm{CS}{\Phi}^I_{*,\Gamma,X}.	
}
\]
\[ 
\xymatrix@R+0pc@C+0pc{
\underset{r}{\mathrm{homotopylimit}}~\mathrm{Spec}^\mathrm{CS}\widetilde{\Phi}^r_{*,\Gamma,X}/\mathrm{Fro}^\mathbb{Z},\underset{I}{\mathrm{homotopycolimit}}~\mathrm{Spec}^\mathrm{CS}\widetilde{\Phi}^I_{*,\Gamma,X}/\mathrm{Fro}^\mathbb{Z},	\\
}
\]
\[ 
\xymatrix@R+0pc@C+0pc{
\underset{r}{\mathrm{homotopylimit}}~\mathrm{Spec}^\mathrm{CS}\breve{\Phi}^r_{*,\Gamma,X}/\mathrm{Fro}^\mathbb{Z},\underset{I}{\mathrm{homotopycolimit}}~\breve{\Phi}^I_{*,\Gamma,X}/\mathrm{Fro}^\mathbb{Z},	\\
}
\]
\[ 
\xymatrix@R+0pc@C+0pc{
\underset{r}{\mathrm{homotopylimit}}~\mathrm{Spec}^\mathrm{CS}{\Phi}^r_{*,\Gamma,X}/\mathrm{Fro}^\mathbb{Z},\underset{I}{\mathrm{homotopycolimit}}~\mathrm{Spec}^\mathrm{CS}{\Phi}^I_{*,\Gamma,X}/\mathrm{Fro}^\mathbb{Z}.	
}
\]

\end{definition}

\begin{proposition}
There is a well-defined functor from the $\infty$-category 
\begin{align}
\mathrm{Quasicoherentpresheaves,Perfectcomplex,Condensed}_{*}	
\end{align}
where $*$ is one of the following spaces:
\begin{align}
&\mathrm{Spec}^\mathrm{CS}\widetilde{\Phi}_{*,\Gamma,X}/\mathrm{Fro}^\mathbb{Z},	\\
\end{align}
\begin{align}
&\mathrm{Spec}^\mathrm{CS}\breve{\Phi}_{*,\Gamma,X}/\mathrm{Fro}^\mathbb{Z},	\\
\end{align}
\begin{align}
&\mathrm{Spec}^\mathrm{CS}{\Phi}_{*,\Gamma,X}/\mathrm{Fro}^\mathbb{Z},	
\end{align}
to the $\infty$-category of $\mathrm{Fro}$-equivariant quasicoherent presheaves over similar spaces above correspondingly without the $\mathrm{Fro}$-quotients, and to the $\infty$-category of $\mathrm{Fro}$-equivariant quasicoherent modules over global sections of the structure $\infty$-sheaves of the similar spaces above correspondingly without the $\mathrm{Fro}$-quotients. Here for those space without notation related to the radius and the corresponding interval we consider the total unions $\bigcap_r,\bigcup_I$ in order to achieve the whole spaces to achieve the analogues of the corresponding FF curves from \cite{10KL1}, \cite{10KL2}, \cite{10FF} for
\[
\xymatrix@R+0pc@C+0pc{
\underset{r}{\mathrm{homotopylimit}}~\mathrm{Spec}^\mathrm{CS}\widetilde{\Phi}^r_{*,\Gamma,X},\underset{I}{\mathrm{homotopycolimit}}~\mathrm{Spec}^\mathrm{CS}\widetilde{\Phi}^I_{*,\Gamma,X},	\\
}
\]
\[
\xymatrix@R+0pc@C+0pc{
\underset{r}{\mathrm{homotopylimit}}~\mathrm{Spec}^\mathrm{CS}\breve{\Phi}^r_{*,\Gamma,X},\underset{I}{\mathrm{homotopycolimit}}~\mathrm{Spec}^\mathrm{CS}\breve{\Phi}^I_{*,\Gamma,X},	\\
}
\]
\[
\xymatrix@R+0pc@C+0pc{
\underset{r}{\mathrm{homotopylimit}}~\mathrm{Spec}^\mathrm{CS}{\Phi}^r_{*,\Gamma,X},\underset{I}{\mathrm{homotopycolimit}}~\mathrm{Spec}^\mathrm{CS}{\Phi}^I_{*,\Gamma,X}.	
}
\]
\[ 
\xymatrix@R+0pc@C+0pc{
\underset{r}{\mathrm{homotopylimit}}~\mathrm{Spec}^\mathrm{CS}\widetilde{\Phi}^r_{*,\Gamma,X}/\mathrm{Fro}^\mathbb{Z},\underset{I}{\mathrm{homotopycolimit}}~\mathrm{Spec}^\mathrm{CS}\widetilde{\Phi}^I_{*,\Gamma,X}/\mathrm{Fro}^\mathbb{Z},	\\
}
\]
\[ 
\xymatrix@R+0pc@C+0pc{
\underset{r}{\mathrm{homotopylimit}}~\mathrm{Spec}^\mathrm{CS}\breve{\Phi}^r_{*,\Gamma,X}/\mathrm{Fro}^\mathbb{Z},\underset{I}{\mathrm{homotopycolimit}}~\breve{\Phi}^I_{*,\Gamma,X}/\mathrm{Fro}^\mathbb{Z},	\\
}
\]
\[ 
\xymatrix@R+0pc@C+0pc{
\underset{r}{\mathrm{homotopylimit}}~\mathrm{Spec}^\mathrm{CS}{\Phi}^r_{*,\Gamma,X}/\mathrm{Fro}^\mathbb{Z},\underset{I}{\mathrm{homotopycolimit}}~\mathrm{Spec}^\mathrm{CS}{\Phi}^I_{*,\Gamma,X}/\mathrm{Fro}^\mathbb{Z}.	
}
\]	
In this situation we will have the target category being family parametrized by $r$ or $I$ in compatible glueing sense as in \cite[Definition 5.4.10]{10KL2}. In this situation for modules parametrized by the intervals we have the equivalence of $\infty$-categories by using \cite[Proposition 12.18]{10CS2}. Here the corresponding quasicoherent Frobenius modules are defined to be the corresponding homotopy colimits and limits of Frobenius modules:
\begin{align}
\underset{r}{\mathrm{homotopycolimit}}~M_r,\\
\underset{I}{\mathrm{homotopylimit}}~M_I,	
\end{align}
where each $M_r$ is a Frobenius-equivariant module over the period ring with respect to some radius $r$ while each $M_I$ is a Frobenius-equivariant module over the period ring with respect to some interval $I$.\\
\end{proposition}

\begin{proposition}
Similar proposition holds for 
\begin{align}
\mathrm{Quasicoherentsheaves,Perfectcomplex,IndBanach}_{*}.	
\end{align}	
\end{proposition}

\subsubsection{Frobenius Quasicoherent Prestacks II: Deformation in Preadic Spaces}

\begin{definition}
We now consider the pro-\'etale site of $\mathrm{Spa}\mathbb{Q}_p\left<X_1^{\pm 1},...,X_k^{\pm 1}\right>$, denote that by $*$. To be more accurate we replace one component for $\Gamma$ with the pro-\'etale site of $\mathrm{Spa}\mathbb{Q}_p\left<X_1^{\pm 1},...,X_k^{\pm 1}\right>$. And we treat then all the functor to be prestacks for this site. Then from \cite{10KL1} and \cite[Definition 5.2.1]{10KL2} we have the following class of Kedlaya-Liu rings (with the following replacement: $\Delta$ stands for $A$, $\nabla$ stands for $B$, while $\Phi$ stands for $C$) by taking product in the sense of self $\Gamma$-th power:

\[
\xymatrix@R+0pc@C+0pc{
\widetilde{\Delta}_{*,\Gamma},\widetilde{\nabla}_{*,\Gamma},\widetilde{\Phi}_{*,\Gamma},\widetilde{\Delta}^+_{*,\Gamma},\widetilde{\nabla}^+_{*,\Gamma},\widetilde{\Delta}^\dagger_{*,\Gamma},\widetilde{\nabla}^\dagger_{*,\Gamma},\widetilde{\Phi}^r_{*,\Gamma},\widetilde{\Phi}^I_{*,\Gamma}, 
}
\]

\[
\xymatrix@R+0pc@C+0pc{
\breve{\Delta}_{*,\Gamma},\breve{\nabla}_{*,\Gamma},\breve{\Phi}_{*,\Gamma},\breve{\Delta}^+_{*,\Gamma},\breve{\nabla}^+_{*,\Gamma},\breve{\Delta}^\dagger_{*,\Gamma},\breve{\nabla}^\dagger_{*,\Gamma},\breve{\Phi}^r_{*,\Gamma},\breve{\Phi}^I_{*,\Gamma},	
}
\]

\[
\xymatrix@R+0pc@C+0pc{
{\Delta}_{*,\Gamma},{\nabla}_{*,\Gamma},{\Phi}_{*,\Gamma},{\Delta}^+_{*,\Gamma},{\nabla}^+_{*,\Gamma},{\Delta}^\dagger_{*,\Gamma},{\nabla}^\dagger_{*,\Gamma},{\Phi}^r_{*,\Gamma},{\Phi}^I_{*,\Gamma}.	
}
\]
Taking the product we have:
\[
\xymatrix@R+0pc@C+0pc{
\widetilde{\Phi}_{*,\Gamma,\circ},\widetilde{\Phi}^r_{*,\Gamma,\circ},\widetilde{\Phi}^I_{*,\Gamma,\circ},	
}
\]
\[
\xymatrix@R+0pc@C+0pc{
\breve{\Phi}_{*,\Gamma,\circ},\breve{\Phi}^r_{*,\Gamma,\circ},\breve{\Phi}^I_{*,\Gamma,\circ},	
}
\]
\[
\xymatrix@R+0pc@C+0pc{
{\Phi}_{*,\Gamma,\circ},{\Phi}^r_{*,\Gamma,\circ},{\Phi}^I_{*,\Gamma,\circ}.	
}
\]
They carry multi Frobenius action $\varphi_\Gamma$ and multi $\mathrm{Lie}_\Gamma:=\mathbb{Z}_p^{\times\Gamma}$ action. In our current situation after \cite{10CKZ} and \cite{10PZ} we consider the following $(\infty,1)$-categories of $(\infty,1)$-modules.\\
\end{definition}

\begin{definition}
First we consider the Bambozzi-Kremnizer spectrum $\mathrm{Spec}^\mathrm{BK}(*)$ attached to any of those in the above from \cite{10BK} by taking derived rational localization:
\begin{align}
&\mathrm{Spec}^\mathrm{BK}\widetilde{\Phi}_{*,\Gamma,\circ},\mathrm{Spec}^\mathrm{BK}\widetilde{\Phi}^r_{*,\Gamma,\circ},\mathrm{Spec}^\mathrm{BK}\widetilde{\Phi}^I_{*,\Gamma,\circ},	
\end{align}
\begin{align}
&\mathrm{Spec}^\mathrm{BK}\breve{\Phi}_{*,\Gamma,\circ},\mathrm{Spec}^\mathrm{BK}\breve{\Phi}^r_{*,\Gamma,\circ},\mathrm{Spec}^\mathrm{BK}\breve{\Phi}^I_{*,\Gamma,\circ},	
\end{align}
\begin{align}
&\mathrm{Spec}^\mathrm{BK}{\Phi}_{*,\Gamma,\circ},
\mathrm{Spec}^\mathrm{BK}{\Phi}^r_{*,\Gamma,\circ},\mathrm{Spec}^\mathrm{BK}{\Phi}^I_{*,\Gamma,\circ}.	
\end{align}

Then we take the corresponding quotients by using the corresponding Frobenius operators:
\begin{align}
&\mathrm{Spec}^\mathrm{BK}\widetilde{\Phi}_{*,\Gamma,\circ}/\mathrm{Fro}^\mathbb{Z},	\\
\end{align}
\begin{align}
&\mathrm{Spec}^\mathrm{BK}\breve{\Phi}_{*,\Gamma,\circ}/\mathrm{Fro}^\mathbb{Z},	\\
\end{align}
\begin{align}
&\mathrm{Spec}^\mathrm{BK}{\Phi}_{*,\Gamma,\circ}/\mathrm{Fro}^\mathbb{Z}.	
\end{align}
Here for those space without notation related to the radius and the corresponding interval we consider the total unions $\bigcap_r,\bigcup_I$ in order to achieve the whole spaces to achieve the analogues of the corresponding FF curves from \cite{10KL1}, \cite{10KL2}, \cite{10FF} for
\[
\xymatrix@R+0pc@C+0pc{
\underset{r}{\mathrm{homotopylimit}}~\mathrm{Spec}^\mathrm{BK}\widetilde{\Phi}^r_{*,\Gamma,\circ},\underset{I}{\mathrm{homotopycolimit}}~\mathrm{Spec}^\mathrm{BK}\widetilde{\Phi}^I_{*,\Gamma,\circ},	\\
}
\]
\[
\xymatrix@R+0pc@C+0pc{
\underset{r}{\mathrm{homotopylimit}}~\mathrm{Spec}^\mathrm{BK}\breve{\Phi}^r_{*,\Gamma,\circ},\underset{I}{\mathrm{homotopycolimit}}~\mathrm{Spec}^\mathrm{BK}\breve{\Phi}^I_{*,\Gamma,\circ},	\\
}
\]
\[
\xymatrix@R+0pc@C+0pc{
\underset{r}{\mathrm{homotopylimit}}~\mathrm{Spec}^\mathrm{BK}{\Phi}^r_{*,\Gamma,\circ},\underset{I}{\mathrm{homotopycolimit}}~\mathrm{Spec}^\mathrm{BK}{\Phi}^I_{*,\Gamma,\circ}.	
}
\]
\[  
\xymatrix@R+0pc@C+0pc{
\underset{r}{\mathrm{homotopylimit}}~\mathrm{Spec}^\mathrm{BK}\widetilde{\Phi}^r_{*,\Gamma,\circ}/\mathrm{Fro}^\mathbb{Z},\underset{I}{\mathrm{homotopycolimit}}~\mathrm{Spec}^\mathrm{BK}\widetilde{\Phi}^I_{*,\Gamma,\circ}/\mathrm{Fro}^\mathbb{Z},	\\
}
\]
\[ 
\xymatrix@R+0pc@C+0pc{
\underset{r}{\mathrm{homotopylimit}}~\mathrm{Spec}^\mathrm{BK}\breve{\Phi}^r_{*,\Gamma,\circ}/\mathrm{Fro}^\mathbb{Z},\underset{I}{\mathrm{homotopycolimit}}~\mathrm{Spec}^\mathrm{BK}\breve{\Phi}^I_{*,\Gamma,\circ}/\mathrm{Fro}^\mathbb{Z},	\\
}
\]
\[ 
\xymatrix@R+0pc@C+0pc{
\underset{r}{\mathrm{homotopylimit}}~\mathrm{Spec}^\mathrm{BK}{\Phi}^r_{*,\Gamma,\circ}/\mathrm{Fro}^\mathbb{Z},\underset{I}{\mathrm{homotopycolimit}}~\mathrm{Spec}^\mathrm{BK}{\Phi}^I_{*,\Gamma,\circ}/\mathrm{Fro}^\mathbb{Z}.	
}
\]

\end{definition}

\indent Meanwhile we have the corresponding Clausen-Scholze analytic stacks from \cite{10CS2}, therefore applying their construction we have:

\begin{definition}
Here we define the following products by using the solidified tensor product from \cite{10CS1} and \cite{10CS2}. Then we take solidified tensor product $\overset{\blacksquare}{\otimes}$ of any of the following
\[
\xymatrix@R+0pc@C+0pc{
\widetilde{\Delta}_{*,\Gamma},\widetilde{\nabla}_{*,\Gamma},\widetilde{\Phi}_{*,\Gamma},\widetilde{\Delta}^+_{*,\Gamma},\widetilde{\nabla}^+_{*,\Gamma},\widetilde{\Delta}^\dagger_{*,\Gamma},\widetilde{\nabla}^\dagger_{*,\Gamma},\widetilde{\Phi}^r_{*,\Gamma},\widetilde{\Phi}^I_{*,\Gamma}, 
}
\]

\[
\xymatrix@R+0pc@C+0pc{
\breve{\Delta}_{*,\Gamma},\breve{\nabla}_{*,\Gamma},\breve{\Phi}_{*,\Gamma},\breve{\Delta}^+_{*,\Gamma},\breve{\nabla}^+_{*,\Gamma},\breve{\Delta}^\dagger_{*,\Gamma},\breve{\nabla}^\dagger_{*,\Gamma},\breve{\Phi}^r_{*,\Gamma},\breve{\Phi}^I_{*,\Gamma},	
}
\]

\[
\xymatrix@R+0pc@C+0pc{
{\Delta}_{*,\Gamma},{\nabla}_{*,\Gamma},{\Phi}_{*,\Gamma},{\Delta}^+_{*,\Gamma},{\nabla}^+_{*,\Gamma},{\Delta}^\dagger_{*,\Gamma},{\nabla}^\dagger_{*,\Gamma},{\Phi}^r_{*,\Gamma},{\Phi}^I_{*,\Gamma},	
}
\]  	
with $\circ$. Then we have the notations:
\[
\xymatrix@R+0pc@C+0pc{
\widetilde{\Delta}_{*,\Gamma,\circ},\widetilde{\nabla}_{*,\Gamma,\circ},\widetilde{\Phi}_{*,\Gamma,\circ},\widetilde{\Delta}^+_{*,\Gamma,\circ},\widetilde{\nabla}^+_{*,\Gamma,\circ},\widetilde{\Delta}^\dagger_{*,\Gamma,\circ},\widetilde{\nabla}^\dagger_{*,\Gamma,\circ},\widetilde{\Phi}^r_{*,\Gamma,\circ},\widetilde{\Phi}^I_{*,\Gamma,\circ}, 
}
\]

\[
\xymatrix@R+0pc@C+0pc{
\breve{\Delta}_{*,\Gamma,\circ},\breve{\nabla}_{*,\Gamma,\circ},\breve{\Phi}_{*,\Gamma,\circ},\breve{\Delta}^+_{*,\Gamma,\circ},\breve{\nabla}^+_{*,\Gamma,\circ},\breve{\Delta}^\dagger_{*,\Gamma,\circ},\breve{\nabla}^\dagger_{*,\Gamma,\circ},\breve{\Phi}^r_{*,\Gamma,\circ},\breve{\Phi}^I_{*,\Gamma,\circ},	
}
\]

\[
\xymatrix@R+0pc@C+0pc{
{\Delta}_{*,\Gamma,\circ},{\nabla}_{*,\Gamma,\circ},{\Phi}_{*,\Gamma,\circ},{\Delta}^+_{*,\Gamma,\circ},{\nabla}^+_{*,\Gamma,\circ},{\Delta}^\dagger_{*,\Gamma,\circ},{\nabla}^\dagger_{*,\Gamma,\circ},{\Phi}^r_{*,\Gamma,\circ},{\Phi}^I_{*,\Gamma,\circ}.	
}
\]
\end{definition}

\begin{definition}
First we consider the Clausen-Scholze spectrum $\mathrm{Spec}^\mathrm{CS}(*)$ attached to any of those in the above from \cite{10CS2} by taking derived rational localization:
\begin{align}
\mathrm{Spec}^\mathrm{CS}\widetilde{\Delta}_{*,\Gamma,\circ},\mathrm{Spec}^\mathrm{CS}\widetilde{\nabla}_{*,\Gamma,\circ},\mathrm{Spec}^\mathrm{CS}\widetilde{\Phi}_{*,\Gamma,\circ},\mathrm{Spec}^\mathrm{CS}\widetilde{\Delta}^+_{*,\Gamma,\circ},\mathrm{Spec}^\mathrm{CS}\widetilde{\nabla}^+_{*,\Gamma,\circ},\\
\mathrm{Spec}^\mathrm{CS}\widetilde{\Delta}^\dagger_{*,\Gamma,\circ},\mathrm{Spec}^\mathrm{CS}\widetilde{\nabla}^\dagger_{*,\Gamma,\circ},\mathrm{Spec}^\mathrm{CS}\widetilde{\Phi}^r_{*,\Gamma,\circ},\mathrm{Spec}^\mathrm{CS}\widetilde{\Phi}^I_{*,\Gamma,\circ},	\\
\end{align}
\begin{align}
\mathrm{Spec}^\mathrm{CS}\breve{\Delta}_{*,\Gamma,\circ},\breve{\nabla}_{*,\Gamma,\circ},\mathrm{Spec}^\mathrm{CS}\breve{\Phi}_{*,\Gamma,\circ},\mathrm{Spec}^\mathrm{CS}\breve{\Delta}^+_{*,\Gamma,\circ},\mathrm{Spec}^\mathrm{CS}\breve{\nabla}^+_{*,\Gamma,\circ},\\
\mathrm{Spec}^\mathrm{CS}\breve{\Delta}^\dagger_{*,\Gamma,\circ},\mathrm{Spec}^\mathrm{CS}\breve{\nabla}^\dagger_{*,\Gamma,\circ},\mathrm{Spec}^\mathrm{CS}\breve{\Phi}^r_{*,\Gamma,\circ},\breve{\Phi}^I_{*,\Gamma,\circ},	\\
\end{align}
\begin{align}
\mathrm{Spec}^\mathrm{CS}{\Delta}_{*,\Gamma,\circ},\mathrm{Spec}^\mathrm{CS}{\nabla}_{*,\Gamma,\circ},\mathrm{Spec}^\mathrm{CS}{\Phi}_{*,\Gamma,\circ},\mathrm{Spec}^\mathrm{CS}{\Delta}^+_{*,\Gamma,\circ},\mathrm{Spec}^\mathrm{CS}{\nabla}^+_{*,\Gamma,\circ},\\
\mathrm{Spec}^\mathrm{CS}{\Delta}^\dagger_{*,\Gamma,\circ},\mathrm{Spec}^\mathrm{CS}{\nabla}^\dagger_{*,\Gamma,\circ},\mathrm{Spec}^\mathrm{CS}{\Phi}^r_{*,\Gamma,\circ},\mathrm{Spec}^\mathrm{CS}{\Phi}^I_{*,\Gamma,\circ}.	
\end{align}

Then we take the corresponding quotients by using the corresponding Frobenius operators:
\begin{align}
&\mathrm{Spec}^\mathrm{CS}\widetilde{\Delta}_{*,\Gamma,\circ}/\mathrm{Fro}^\mathbb{Z},\mathrm{Spec}^\mathrm{CS}\widetilde{\nabla}_{*,\Gamma,\circ}/\mathrm{Fro}^\mathbb{Z},\mathrm{Spec}^\mathrm{CS}\widetilde{\Phi}_{*,\Gamma,\circ}/\mathrm{Fro}^\mathbb{Z},\mathrm{Spec}^\mathrm{CS}\widetilde{\Delta}^+_{*,\Gamma,\circ}/\mathrm{Fro}^\mathbb{Z},\\
&\mathrm{Spec}^\mathrm{CS}\widetilde{\nabla}^+_{*,\Gamma,\circ}/\mathrm{Fro}^\mathbb{Z}, \mathrm{Spec}^\mathrm{CS}\widetilde{\Delta}^\dagger_{*,\Gamma,\circ}/\mathrm{Fro}^\mathbb{Z},\mathrm{Spec}^\mathrm{CS}\widetilde{\nabla}^\dagger_{*,\Gamma,\circ}/\mathrm{Fro}^\mathbb{Z},	\\
\end{align}
\begin{align}
&\mathrm{Spec}^\mathrm{CS}\breve{\Delta}_{*,\Gamma,\circ}/\mathrm{Fro}^\mathbb{Z},\breve{\nabla}_{*,\Gamma,\circ}/\mathrm{Fro}^\mathbb{Z},\mathrm{Spec}^\mathrm{CS}\breve{\Phi}_{*,\Gamma,\circ}/\mathrm{Fro}^\mathbb{Z},\mathrm{Spec}^\mathrm{CS}\breve{\Delta}^+_{*,\Gamma,\circ}/\mathrm{Fro}^\mathbb{Z},\\
&\mathrm{Spec}^\mathrm{CS}\breve{\nabla}^+_{*,\Gamma,\circ}/\mathrm{Fro}^\mathbb{Z}, \mathrm{Spec}^\mathrm{CS}\breve{\Delta}^\dagger_{*,\Gamma,\circ}/\mathrm{Fro}^\mathbb{Z},\mathrm{Spec}^\mathrm{CS}\breve{\nabla}^\dagger_{*,\Gamma,\circ}/\mathrm{Fro}^\mathbb{Z},	\\
\end{align}
\begin{align}
&\mathrm{Spec}^\mathrm{CS}{\Delta}_{*,\Gamma,\circ}/\mathrm{Fro}^\mathbb{Z},\mathrm{Spec}^\mathrm{CS}{\nabla}_{*,\Gamma,\circ}/\mathrm{Fro}^\mathbb{Z},\mathrm{Spec}^\mathrm{CS}{\Phi}_{*,\Gamma,\circ}/\mathrm{Fro}^\mathbb{Z},\mathrm{Spec}^\mathrm{CS}{\Delta}^+_{*,\Gamma,\circ}/\mathrm{Fro}^\mathbb{Z},\\
&\mathrm{Spec}^\mathrm{CS}{\nabla}^+_{*,\Gamma,\circ}/\mathrm{Fro}^\mathbb{Z}, \mathrm{Spec}^\mathrm{CS}{\Delta}^\dagger_{*,\Gamma,\circ}/\mathrm{Fro}^\mathbb{Z},\mathrm{Spec}^\mathrm{CS}{\nabla}^\dagger_{*,\Gamma,\circ}/\mathrm{Fro}^\mathbb{Z}.	
\end{align}
Here for those space with notations related to the radius and the corresponding interval we consider the total unions $\bigcap_r,\bigcup_I$ in order to achieve the whole spaces to achieve the analogues of the corresponding FF curves from \cite{10KL1}, \cite{10KL2}, \cite{10FF} for
\[
\xymatrix@R+0pc@C+0pc{
\underset{r}{\mathrm{homotopylimit}}~\mathrm{Spec}^\mathrm{CS}\widetilde{\Phi}^r_{*,\Gamma,\circ},\underset{I}{\mathrm{homotopycolimit}}~\mathrm{Spec}^\mathrm{CS}\widetilde{\Phi}^I_{*,\Gamma,\circ},	\\
}
\]
\[
\xymatrix@R+0pc@C+0pc{
\underset{r}{\mathrm{homotopylimit}}~\mathrm{Spec}^\mathrm{CS}\breve{\Phi}^r_{*,\Gamma,\circ},\underset{I}{\mathrm{homotopycolimit}}~\mathrm{Spec}^\mathrm{CS}\breve{\Phi}^I_{*,\Gamma,\circ},	\\
}
\]
\[
\xymatrix@R+0pc@C+0pc{
\underset{r}{\mathrm{homotopylimit}}~\mathrm{Spec}^\mathrm{CS}{\Phi}^r_{*,\Gamma,\circ},\underset{I}{\mathrm{homotopycolimit}}~\mathrm{Spec}^\mathrm{CS}{\Phi}^I_{*,\Gamma,\circ}.	
}
\]
\[ 
\xymatrix@R+0pc@C+0pc{
\underset{r}{\mathrm{homotopylimit}}~\mathrm{Spec}^\mathrm{CS}\widetilde{\Phi}^r_{*,\Gamma,\circ}/\mathrm{Fro}^\mathbb{Z},\underset{I}{\mathrm{homotopycolimit}}~\mathrm{Spec}^\mathrm{CS}\widetilde{\Phi}^I_{*,\Gamma,\circ}/\mathrm{Fro}^\mathbb{Z},	\\
}
\]
\[ 
\xymatrix@R+0pc@C+0pc{
\underset{r}{\mathrm{homotopylimit}}~\mathrm{Spec}^\mathrm{CS}\breve{\Phi}^r_{*,\Gamma,\circ}/\mathrm{Fro}^\mathbb{Z},\underset{I}{\mathrm{homotopycolimit}}~\breve{\Phi}^I_{*,\Gamma,\circ}/\mathrm{Fro}^\mathbb{Z},	\\
}
\]
\[ 
\xymatrix@R+0pc@C+0pc{
\underset{r}{\mathrm{homotopylimit}}~\mathrm{Spec}^\mathrm{CS}{\Phi}^r_{*,\Gamma,\circ}/\mathrm{Fro}^\mathbb{Z},\underset{I}{\mathrm{homotopycolimit}}~\mathrm{Spec}^\mathrm{CS}{\Phi}^I_{*,\Gamma,\circ}/\mathrm{Fro}^\mathbb{Z}.	
}
\]

\end{definition}

\

\begin{definition}
We then consider the corresponding quasipresheaves of the corresponding ind-Banach or monomorphic ind-Banach modules from \cite{10BBK}, \cite{10KKM}:
\begin{align}
\mathrm{Quasicoherentpresheaves,IndBanach}_{*}	
\end{align}
where $*$ is one of the following spaces:
\begin{align}
&\mathrm{Spec}^\mathrm{BK}\widetilde{\Phi}_{*,\Gamma,\circ}/\mathrm{Fro}^\mathbb{Z},	\\
\end{align}
\begin{align}
&\mathrm{Spec}^\mathrm{BK}\breve{\Phi}_{*,\Gamma,\circ}/\mathrm{Fro}^\mathbb{Z},	\\
\end{align}
\begin{align}
&\mathrm{Spec}^\mathrm{BK}{\Phi}_{*,\Gamma,\circ}/\mathrm{Fro}^\mathbb{Z}.	
\end{align}
Here for those space without notation related to the radius and the corresponding interval we consider the total unions $\bigcap_r,\bigcup_I$ in order to achieve the whole spaces to achieve the analogues of the corresponding FF curves from \cite{10KL1}, \cite{10KL2}, \cite{10FF} for
\[
\xymatrix@R+0pc@C+0pc{
\underset{r}{\mathrm{homotopylimit}}~\mathrm{Spec}^\mathrm{BK}\widetilde{\Phi}^r_{*,\Gamma,\circ},\underset{I}{\mathrm{homotopycolimit}}~\mathrm{Spec}^\mathrm{BK}\widetilde{\Phi}^I_{*,\Gamma,\circ},	\\
}
\]
\[
\xymatrix@R+0pc@C+0pc{
\underset{r}{\mathrm{homotopylimit}}~\mathrm{Spec}^\mathrm{BK}\breve{\Phi}^r_{*,\Gamma,\circ},\underset{I}{\mathrm{homotopycolimit}}~\mathrm{Spec}^\mathrm{BK}\breve{\Phi}^I_{*,\Gamma,\circ},	\\
}
\]
\[
\xymatrix@R+0pc@C+0pc{
\underset{r}{\mathrm{homotopylimit}}~\mathrm{Spec}^\mathrm{BK}{\Phi}^r_{*,\Gamma,\circ},\underset{I}{\mathrm{homotopycolimit}}~\mathrm{Spec}^\mathrm{BK}{\Phi}^I_{*,\Gamma,\circ}.	
}
\]
\[  
\xymatrix@R+0pc@C+0pc{
\underset{r}{\mathrm{homotopylimit}}~\mathrm{Spec}^\mathrm{BK}\widetilde{\Phi}^r_{*,\Gamma,\circ}/\mathrm{Fro}^\mathbb{Z},\underset{I}{\mathrm{homotopycolimit}}~\mathrm{Spec}^\mathrm{BK}\widetilde{\Phi}^I_{*,\Gamma,\circ}/\mathrm{Fro}^\mathbb{Z},	\\
}
\]
\[ 
\xymatrix@R+0pc@C+0pc{
\underset{r}{\mathrm{homotopylimit}}~\mathrm{Spec}^\mathrm{BK}\breve{\Phi}^r_{*,\Gamma,\circ}/\mathrm{Fro}^\mathbb{Z},\underset{I}{\mathrm{homotopycolimit}}~\mathrm{Spec}^\mathrm{BK}\breve{\Phi}^I_{*,\Gamma,\circ}/\mathrm{Fro}^\mathbb{Z},	\\
}
\]
\[ 
\xymatrix@R+0pc@C+0pc{
\underset{r}{\mathrm{homotopylimit}}~\mathrm{Spec}^\mathrm{BK}{\Phi}^r_{*,\Gamma,\circ}/\mathrm{Fro}^\mathbb{Z},\underset{I}{\mathrm{homotopycolimit}}~\mathrm{Spec}^\mathrm{BK}{\Phi}^I_{*,\Gamma,\circ}/\mathrm{Fro}^\mathbb{Z}.	
}
\]

\end{definition}

\begin{definition}
We then consider the corresponding quasisheaves of the corresponding condensed solid topological modules from \cite{10CS2}:
\begin{align}
\mathrm{Quasicoherentsheaves, Condensed}_{*}	
\end{align}
where $*$ is one of the following spaces:
\begin{align}
&\mathrm{Spec}^\mathrm{CS}\widetilde{\Delta}_{*,\Gamma,\circ}/\mathrm{Fro}^\mathbb{Z},\mathrm{Spec}^\mathrm{CS}\widetilde{\nabla}_{*,\Gamma,\circ}/\mathrm{Fro}^\mathbb{Z},\mathrm{Spec}^\mathrm{CS}\widetilde{\Phi}_{*,\Gamma,\circ}/\mathrm{Fro}^\mathbb{Z},\mathrm{Spec}^\mathrm{CS}\widetilde{\Delta}^+_{*,\Gamma,\circ}/\mathrm{Fro}^\mathbb{Z},\\
&\mathrm{Spec}^\mathrm{CS}\widetilde{\nabla}^+_{*,\Gamma,\circ}/\mathrm{Fro}^\mathbb{Z},\mathrm{Spec}^\mathrm{CS}\widetilde{\Delta}^\dagger_{*,\Gamma,\circ}/\mathrm{Fro}^\mathbb{Z},\mathrm{Spec}^\mathrm{CS}\widetilde{\nabla}^\dagger_{*,\Gamma,\circ}/\mathrm{Fro}^\mathbb{Z},	\\
\end{align}
\begin{align}
&\mathrm{Spec}^\mathrm{CS}\breve{\Delta}_{*,\Gamma,\circ}/\mathrm{Fro}^\mathbb{Z},\breve{\nabla}_{*,\Gamma,\circ}/\mathrm{Fro}^\mathbb{Z},\mathrm{Spec}^\mathrm{CS}\breve{\Phi}_{*,\Gamma,\circ}/\mathrm{Fro}^\mathbb{Z},\mathrm{Spec}^\mathrm{CS}\breve{\Delta}^+_{*,\Gamma,\circ}/\mathrm{Fro}^\mathbb{Z},\\
&\mathrm{Spec}^\mathrm{CS}\breve{\nabla}^+_{*,\Gamma,\circ}/\mathrm{Fro}^\mathbb{Z},\mathrm{Spec}^\mathrm{CS}\breve{\Delta}^\dagger_{*,\Gamma,\circ}/\mathrm{Fro}^\mathbb{Z},\mathrm{Spec}^\mathrm{CS}\breve{\nabla}^\dagger_{*,\Gamma,\circ}/\mathrm{Fro}^\mathbb{Z},	\\
\end{align}
\begin{align}
&\mathrm{Spec}^\mathrm{CS}{\Delta}_{*,\Gamma,\circ}/\mathrm{Fro}^\mathbb{Z},\mathrm{Spec}^\mathrm{CS}{\nabla}_{*,\Gamma,\circ}/\mathrm{Fro}^\mathbb{Z},\mathrm{Spec}^\mathrm{CS}{\Phi}_{*,\Gamma,\circ}/\mathrm{Fro}^\mathbb{Z},\mathrm{Spec}^\mathrm{CS}{\Delta}^+_{*,\Gamma,\circ}/\mathrm{Fro}^\mathbb{Z},\\
&\mathrm{Spec}^\mathrm{CS}{\nabla}^+_{*,\Gamma,\circ}/\mathrm{Fro}^\mathbb{Z}, \mathrm{Spec}^\mathrm{CS}{\Delta}^\dagger_{*,\Gamma,\circ}/\mathrm{Fro}^\mathbb{Z},\mathrm{Spec}^\mathrm{CS}{\nabla}^\dagger_{*,\Gamma,\circ}/\mathrm{Fro}^\mathbb{Z}.	
\end{align}
Here for those space with notations related to the radius and the corresponding interval we consider the total unions $\bigcap_r,\bigcup_I$ in order to achieve the whole spaces to achieve the analogues of the corresponding FF curves from \cite{10KL1}, \cite{10KL2}, \cite{10FF} for
\[
\xymatrix@R+0pc@C+0pc{
\underset{r}{\mathrm{homotopylimit}}~\mathrm{Spec}^\mathrm{CS}\widetilde{\Phi}^r_{*,\Gamma,\circ},\underset{I}{\mathrm{homotopycolimit}}~\mathrm{Spec}^\mathrm{CS}\widetilde{\Phi}^I_{*,\Gamma,\circ},	\\
}
\]
\[
\xymatrix@R+0pc@C+0pc{
\underset{r}{\mathrm{homotopylimit}}~\mathrm{Spec}^\mathrm{CS}\breve{\Phi}^r_{*,\Gamma,\circ},\underset{I}{\mathrm{homotopycolimit}}~\mathrm{Spec}^\mathrm{CS}\breve{\Phi}^I_{*,\Gamma,\circ},	\\
}
\]
\[
\xymatrix@R+0pc@C+0pc{
\underset{r}{\mathrm{homotopylimit}}~\mathrm{Spec}^\mathrm{CS}{\Phi}^r_{*,\Gamma,\circ},\underset{I}{\mathrm{homotopycolimit}}~\mathrm{Spec}^\mathrm{CS}{\Phi}^I_{*,\Gamma,\circ}.	
}
\]
\[ 
\xymatrix@R+0pc@C+0pc{
\underset{r}{\mathrm{homotopylimit}}~\mathrm{Spec}^\mathrm{CS}\widetilde{\Phi}^r_{*,\Gamma,\circ}/\mathrm{Fro}^\mathbb{Z},\underset{I}{\mathrm{homotopycolimit}}~\mathrm{Spec}^\mathrm{CS}\widetilde{\Phi}^I_{*,\Gamma,\circ}/\mathrm{Fro}^\mathbb{Z},	\\
}
\]
\[ 
\xymatrix@R+0pc@C+0pc{
\underset{r}{\mathrm{homotopylimit}}~\mathrm{Spec}^\mathrm{CS}\breve{\Phi}^r_{*,\Gamma,\circ}/\mathrm{Fro}^\mathbb{Z},\underset{I}{\mathrm{homotopycolimit}}~\breve{\Phi}^I_{*,\Gamma,\circ}/\mathrm{Fro}^\mathbb{Z},	\\
}
\]
\[ 
\xymatrix@R+0pc@C+0pc{
\underset{r}{\mathrm{homotopylimit}}~\mathrm{Spec}^\mathrm{CS}{\Phi}^r_{*,\Gamma,\circ}/\mathrm{Fro}^\mathbb{Z},\underset{I}{\mathrm{homotopycolimit}}~\mathrm{Spec}^\mathrm{CS}{\Phi}^I_{*,\Gamma,\circ}/\mathrm{Fro}^\mathbb{Z}.	
}
\]

\end{definition}

\

\begin{proposition}
There is a well-defined functor from the $\infty$-category 
\begin{align}
\mathrm{Quasicoherentpresheaves,Condensed}_{*}	
\end{align}
where $*$ is one of the following spaces:
\begin{align}
&\mathrm{Spec}^\mathrm{CS}\widetilde{\Phi}_{*,\Gamma,\circ}/\mathrm{Fro}^\mathbb{Z},	\\
\end{align}
\begin{align}
&\mathrm{Spec}^\mathrm{CS}\breve{\Phi}_{*,\Gamma,\circ}/\mathrm{Fro}^\mathbb{Z},	\\
\end{align}
\begin{align}
&\mathrm{Spec}^\mathrm{CS}{\Phi}_{*,\Gamma,\circ}/\mathrm{Fro}^\mathbb{Z},	
\end{align}
to the $\infty$-category of $\mathrm{Fro}$-equivariant quasicoherent presheaves over similar spaces above correspondingly without the $\mathrm{Fro}$-quotients, and to the $\infty$-category of $\mathrm{Fro}$-equivariant quasicoherent modules over global sections of the structure $\infty$-sheaves of the similar spaces above correspondingly without the $\mathrm{Fro}$-quotients. Here for those space without notation related to the radius and the corresponding interval we consider the total unions $\bigcap_r,\bigcup_I$ in order to achieve the whole spaces to achieve the analogues of the corresponding FF curves from \cite{10KL1}, \cite{10KL2}, \cite{10FF} for
\[
\xymatrix@R+0pc@C+0pc{
\underset{r}{\mathrm{homotopylimit}}~\mathrm{Spec}^\mathrm{CS}\widetilde{\Phi}^r_{*,\Gamma,\circ},\underset{I}{\mathrm{homotopycolimit}}~\mathrm{Spec}^\mathrm{CS}\widetilde{\Phi}^I_{*,\Gamma,\circ},	\\
}
\]
\[
\xymatrix@R+0pc@C+0pc{
\underset{r}{\mathrm{homotopylimit}}~\mathrm{Spec}^\mathrm{CS}\breve{\Phi}^r_{*,\Gamma,\circ},\underset{I}{\mathrm{homotopycolimit}}~\mathrm{Spec}^\mathrm{CS}\breve{\Phi}^I_{*,\Gamma,\circ},	\\
}
\]
\[
\xymatrix@R+0pc@C+0pc{
\underset{r}{\mathrm{homotopylimit}}~\mathrm{Spec}^\mathrm{CS}{\Phi}^r_{*,\Gamma,\circ},\underset{I}{\mathrm{homotopycolimit}}~\mathrm{Spec}^\mathrm{CS}{\Phi}^I_{*,\Gamma,\circ}.	
}
\]
\[ 
\xymatrix@R+0pc@C+0pc{
\underset{r}{\mathrm{homotopylimit}}~\mathrm{Spec}^\mathrm{CS}\widetilde{\Phi}^r_{*,\Gamma,\circ}/\mathrm{Fro}^\mathbb{Z},\underset{I}{\mathrm{homotopycolimit}}~\mathrm{Spec}^\mathrm{CS}\widetilde{\Phi}^I_{*,\Gamma,\circ}/\mathrm{Fro}^\mathbb{Z},	\\
}
\]
\[ 
\xymatrix@R+0pc@C+0pc{
\underset{r}{\mathrm{homotopylimit}}~\mathrm{Spec}^\mathrm{CS}\breve{\Phi}^r_{*,\Gamma,\circ}/\mathrm{Fro}^\mathbb{Z},\underset{I}{\mathrm{homotopycolimit}}~\breve{\Phi}^I_{*,\Gamma,\circ}/\mathrm{Fro}^\mathbb{Z},	\\
}
\]
\[ 
\xymatrix@R+0pc@C+0pc{
\underset{r}{\mathrm{homotopylimit}}~\mathrm{Spec}^\mathrm{CS}{\Phi}^r_{*,\Gamma,\circ}/\mathrm{Fro}^\mathbb{Z},\underset{I}{\mathrm{homotopycolimit}}~\mathrm{Spec}^\mathrm{CS}{\Phi}^I_{*,\Gamma,\circ}/\mathrm{Fro}^\mathbb{Z}.	
}
\]	
In this situation we will have the target category being family parametrized by $r$ or $I$ in compatible glueing sense as in \cite[Definition 5.4.10]{10KL2}. In this situation for modules parametrized by the intervals we have the equivalence of $\infty$-categories by using \cite[Proposition 13.8]{10CS2}. Here the corresponding quasicoherent Frobenius modules are defined to be the corresponding homotopy colimits and limits of Frobenius modules:
\begin{align}
\underset{r}{\mathrm{homotopycolimit}}~M_r,\\
\underset{I}{\mathrm{homotopylimit}}~M_I,	
\end{align}
where each $M_r$ is a Frobenius-equivariant module over the period ring with respect to some radius $r$ while each $M_I$ is a Frobenius-equivariant module over the period ring with respect to some interval $I$.\\
\end{proposition}

\begin{proposition}
Similar proposition holds for 
\begin{align}
\mathrm{Quasicoherentsheaves,IndBanach}_{*}.	
\end{align}	
\end{proposition}

\

\begin{definition}
We then consider the corresponding quasipresheaves of perfect complexes the corresponding ind-Banach or monomorphic ind-Banach modules from \cite{10BBK}, \cite{10KKM}:
\begin{align}
\mathrm{Quasicoherentpresheaves,Perfectcomplex,IndBanach}_{*}	
\end{align}
where $*$ is one of the following spaces:
\begin{align}
&\mathrm{Spec}^\mathrm{BK}\widetilde{\Phi}_{*,\Gamma,\circ}/\mathrm{Fro}^\mathbb{Z},	\\
\end{align}
\begin{align}
&\mathrm{Spec}^\mathrm{BK}\breve{\Phi}_{*,\Gamma,\circ}/\mathrm{Fro}^\mathbb{Z},	\\
\end{align}
\begin{align}
&\mathrm{Spec}^\mathrm{BK}{\Phi}_{*,\Gamma,\circ}/\mathrm{Fro}^\mathbb{Z}.	
\end{align}
Here for those space without notation related to the radius and the corresponding interval we consider the total unions $\bigcap_r,\bigcup_I$ in order to achieve the whole spaces to achieve the analogues of the corresponding FF curves from \cite{10KL1}, \cite{10KL2}, \cite{10FF} for
\[
\xymatrix@R+0pc@C+0pc{
\underset{r}{\mathrm{homotopylimit}}~\mathrm{Spec}^\mathrm{BK}\widetilde{\Phi}^r_{*,\Gamma,\circ},\underset{I}{\mathrm{homotopycolimit}}~\mathrm{Spec}^\mathrm{BK}\widetilde{\Phi}^I_{*,\Gamma,\circ},	\\
}
\]
\[
\xymatrix@R+0pc@C+0pc{
\underset{r}{\mathrm{homotopylimit}}~\mathrm{Spec}^\mathrm{BK}\breve{\Phi}^r_{*,\Gamma,\circ},\underset{I}{\mathrm{homotopycolimit}}~\mathrm{Spec}^\mathrm{BK}\breve{\Phi}^I_{*,\Gamma,\circ},	\\
}
\]
\[
\xymatrix@R+0pc@C+0pc{
\underset{r}{\mathrm{homotopylimit}}~\mathrm{Spec}^\mathrm{BK}{\Phi}^r_{*,\Gamma,\circ},\underset{I}{\mathrm{homotopycolimit}}~\mathrm{Spec}^\mathrm{BK}{\Phi}^I_{*,\Gamma,\circ}.	
}
\]
\[  
\xymatrix@R+0pc@C+0pc{
\underset{r}{\mathrm{homotopylimit}}~\mathrm{Spec}^\mathrm{BK}\widetilde{\Phi}^r_{*,\Gamma,\circ}/\mathrm{Fro}^\mathbb{Z},\underset{I}{\mathrm{homotopycolimit}}~\mathrm{Spec}^\mathrm{BK}\widetilde{\Phi}^I_{*,\Gamma,\circ}/\mathrm{Fro}^\mathbb{Z},	\\
}
\]
\[ 
\xymatrix@R+0pc@C+0pc{
\underset{r}{\mathrm{homotopylimit}}~\mathrm{Spec}^\mathrm{BK}\breve{\Phi}^r_{*,\Gamma,\circ}/\mathrm{Fro}^\mathbb{Z},\underset{I}{\mathrm{homotopycolimit}}~\mathrm{Spec}^\mathrm{BK}\breve{\Phi}^I_{*,\Gamma,\circ}/\mathrm{Fro}^\mathbb{Z},	\\
}
\]
\[ 
\xymatrix@R+0pc@C+0pc{
\underset{r}{\mathrm{homotopylimit}}~\mathrm{Spec}^\mathrm{BK}{\Phi}^r_{*,\Gamma,\circ}/\mathrm{Fro}^\mathbb{Z},\underset{I}{\mathrm{homotopycolimit}}~\mathrm{Spec}^\mathrm{BK}{\Phi}^I_{*,\Gamma,\circ}/\mathrm{Fro}^\mathbb{Z}.	
}
\]

\end{definition}

\begin{definition}
We then consider the corresponding quasisheaves of perfect complexes of the corresponding condensed solid topological modules from \cite{10CS2}:
\begin{align}
\mathrm{Quasicoherentsheaves, Perfectcomplex, Condensed}_{*}	
\end{align}
where $*$ is one of the following spaces:
\begin{align}
&\mathrm{Spec}^\mathrm{CS}\widetilde{\Delta}_{*,\Gamma,\circ}/\mathrm{Fro}^\mathbb{Z},\mathrm{Spec}^\mathrm{CS}\widetilde{\nabla}_{*,\Gamma,\circ}/\mathrm{Fro}^\mathbb{Z},\mathrm{Spec}^\mathrm{CS}\widetilde{\Phi}_{*,\Gamma,\circ}/\mathrm{Fro}^\mathbb{Z},\mathrm{Spec}^\mathrm{CS}\widetilde{\Delta}^+_{*,\Gamma,\circ}/\mathrm{Fro}^\mathbb{Z},\\
&\mathrm{Spec}^\mathrm{CS}\widetilde{\nabla}^+_{*,\Gamma,\circ}/\mathrm{Fro}^\mathbb{Z},\mathrm{Spec}^\mathrm{CS}\widetilde{\Delta}^\dagger_{*,\Gamma,\circ}/\mathrm{Fro}^\mathbb{Z},\mathrm{Spec}^\mathrm{CS}\widetilde{\nabla}^\dagger_{*,\Gamma,\circ}/\mathrm{Fro}^\mathbb{Z},	\\
\end{align}
\begin{align}
&\mathrm{Spec}^\mathrm{CS}\breve{\Delta}_{*,\Gamma,\circ}/\mathrm{Fro}^\mathbb{Z},\breve{\nabla}_{*,\Gamma,\circ}/\mathrm{Fro}^\mathbb{Z},\mathrm{Spec}^\mathrm{CS}\breve{\Phi}_{*,\Gamma,\circ}/\mathrm{Fro}^\mathbb{Z},\mathrm{Spec}^\mathrm{CS}\breve{\Delta}^+_{*,\Gamma,\circ}/\mathrm{Fro}^\mathbb{Z},\\
&\mathrm{Spec}^\mathrm{CS}\breve{\nabla}^+_{*,\Gamma,\circ}/\mathrm{Fro}^\mathbb{Z},\mathrm{Spec}^\mathrm{CS}\breve{\Delta}^\dagger_{*,\Gamma,\circ}/\mathrm{Fro}^\mathbb{Z},\mathrm{Spec}^\mathrm{CS}\breve{\nabla}^\dagger_{*,\Gamma,\circ}/\mathrm{Fro}^\mathbb{Z},	\\
\end{align}
\begin{align}
&\mathrm{Spec}^\mathrm{CS}{\Delta}_{*,\Gamma,\circ}/\mathrm{Fro}^\mathbb{Z},\mathrm{Spec}^\mathrm{CS}{\nabla}_{*,\Gamma,\circ}/\mathrm{Fro}^\mathbb{Z},\mathrm{Spec}^\mathrm{CS}{\Phi}_{*,\Gamma,\circ}/\mathrm{Fro}^\mathbb{Z},\mathrm{Spec}^\mathrm{CS}{\Delta}^+_{*,\Gamma,\circ}/\mathrm{Fro}^\mathbb{Z},\\
&\mathrm{Spec}^\mathrm{CS}{\nabla}^+_{*,\Gamma,\circ}/\mathrm{Fro}^\mathbb{Z}, \mathrm{Spec}^\mathrm{CS}{\Delta}^\dagger_{*,\Gamma,\circ}/\mathrm{Fro}^\mathbb{Z},\mathrm{Spec}^\mathrm{CS}{\nabla}^\dagger_{*,\Gamma,\circ}/\mathrm{Fro}^\mathbb{Z}.	
\end{align}
Here for those space with notations related to the radius and the corresponding interval we consider the total unions $\bigcap_r,\bigcup_I$ in order to achieve the whole spaces to achieve the analogues of the corresponding FF curves from \cite{10KL1}, \cite{10KL2}, \cite{10FF} for
\[
\xymatrix@R+0pc@C+0pc{
\underset{r}{\mathrm{homotopylimit}}~\mathrm{Spec}^\mathrm{CS}\widetilde{\Phi}^r_{*,\Gamma,\circ},\underset{I}{\mathrm{homotopycolimit}}~\mathrm{Spec}^\mathrm{CS}\widetilde{\Phi}^I_{*,\Gamma,\circ},	\\
}
\]
\[
\xymatrix@R+0pc@C+0pc{
\underset{r}{\mathrm{homotopylimit}}~\mathrm{Spec}^\mathrm{CS}\breve{\Phi}^r_{*,\Gamma,\circ},\underset{I}{\mathrm{homotopycolimit}}~\mathrm{Spec}^\mathrm{CS}\breve{\Phi}^I_{*,\Gamma,\circ},	\\
}
\]
\[
\xymatrix@R+0pc@C+0pc{
\underset{r}{\mathrm{homotopylimit}}~\mathrm{Spec}^\mathrm{CS}{\Phi}^r_{*,\Gamma,\circ},\underset{I}{\mathrm{homotopycolimit}}~\mathrm{Spec}^\mathrm{CS}{\Phi}^I_{*,\Gamma,\circ}.	
}
\]
\[ 
\xymatrix@R+0pc@C+0pc{
\underset{r}{\mathrm{homotopylimit}}~\mathrm{Spec}^\mathrm{CS}\widetilde{\Phi}^r_{*,\Gamma,\circ}/\mathrm{Fro}^\mathbb{Z},\underset{I}{\mathrm{homotopycolimit}}~\mathrm{Spec}^\mathrm{CS}\widetilde{\Phi}^I_{*,\Gamma,\circ}/\mathrm{Fro}^\mathbb{Z},	\\
}
\]
\[ 
\xymatrix@R+0pc@C+0pc{
\underset{r}{\mathrm{homotopylimit}}~\mathrm{Spec}^\mathrm{CS}\breve{\Phi}^r_{*,\Gamma,\circ}/\mathrm{Fro}^\mathbb{Z},\underset{I}{\mathrm{homotopycolimit}}~\breve{\Phi}^I_{*,\Gamma,\circ}/\mathrm{Fro}^\mathbb{Z},	\\
}
\]
\[ 
\xymatrix@R+0pc@C+0pc{
\underset{r}{\mathrm{homotopylimit}}~\mathrm{Spec}^\mathrm{CS}{\Phi}^r_{*,\Gamma,\circ}/\mathrm{Fro}^\mathbb{Z},\underset{I}{\mathrm{homotopycolimit}}~\mathrm{Spec}^\mathrm{CS}{\Phi}^I_{*,\Gamma,\circ}/\mathrm{Fro}^\mathbb{Z}.	
}
\]

\end{definition}

\begin{proposition}
There is a well-defined functor from the $\infty$-category 
\begin{align}
\mathrm{Quasicoherentpresheaves,Perfectcomplex,Condensed}_{*}	
\end{align}
where $*$ is one of the following spaces:
\begin{align}
&\mathrm{Spec}^\mathrm{CS}\widetilde{\Phi}_{*,\Gamma,\circ}/\mathrm{Fro}^\mathbb{Z},	\\
\end{align}
\begin{align}
&\mathrm{Spec}^\mathrm{CS}\breve{\Phi}_{*,\Gamma,\circ}/\mathrm{Fro}^\mathbb{Z},	\\
\end{align}
\begin{align}
&\mathrm{Spec}^\mathrm{CS}{\Phi}_{*,\Gamma,\circ}/\mathrm{Fro}^\mathbb{Z},	
\end{align}
to the $\infty$-category of $\mathrm{Fro}$-equivariant quasicoherent presheaves over similar spaces above correspondingly without the $\mathrm{Fro}$-quotients, and to the $\infty$-category of $\mathrm{Fro}$-equivariant quasicoherent modules over global sections of the structure $\infty$-sheaves of the similar spaces above correspondingly without the $\mathrm{Fro}$-quotients. Here for those space without notation related to the radius and the corresponding interval we consider the total unions $\bigcap_r,\bigcup_I$ in order to achieve the whole spaces to achieve the analogues of the corresponding FF curves from \cite{10KL1}, \cite{10KL2}, \cite{10FF} for
\[
\xymatrix@R+0pc@C+0pc{
\underset{r}{\mathrm{homotopylimit}}~\mathrm{Spec}^\mathrm{CS}\widetilde{\Phi}^r_{*,\Gamma,\circ},\underset{I}{\mathrm{homotopycolimit}}~\mathrm{Spec}^\mathrm{CS}\widetilde{\Phi}^I_{*,\Gamma,\circ},	\\
}
\]
\[
\xymatrix@R+0pc@C+0pc{
\underset{r}{\mathrm{homotopylimit}}~\mathrm{Spec}^\mathrm{CS}\breve{\Phi}^r_{*,\Gamma,\circ},\underset{I}{\mathrm{homotopycolimit}}~\mathrm{Spec}^\mathrm{CS}\breve{\Phi}^I_{*,\Gamma,\circ},	\\
}
\]
\[
\xymatrix@R+0pc@C+0pc{
\underset{r}{\mathrm{homotopylimit}}~\mathrm{Spec}^\mathrm{CS}{\Phi}^r_{*,\Gamma,\circ},\underset{I}{\mathrm{homotopycolimit}}~\mathrm{Spec}^\mathrm{CS}{\Phi}^I_{*,\Gamma,\circ}.	
}
\]
\[ 
\xymatrix@R+0pc@C+0pc{
\underset{r}{\mathrm{homotopylimit}}~\mathrm{Spec}^\mathrm{CS}\widetilde{\Phi}^r_{*,\Gamma,\circ}/\mathrm{Fro}^\mathbb{Z},\underset{I}{\mathrm{homotopycolimit}}~\mathrm{Spec}^\mathrm{CS}\widetilde{\Phi}^I_{*,\Gamma,\circ}/\mathrm{Fro}^\mathbb{Z},	\\
}
\]
\[ 
\xymatrix@R+0pc@C+0pc{
\underset{r}{\mathrm{homotopylimit}}~\mathrm{Spec}^\mathrm{CS}\breve{\Phi}^r_{*,\Gamma,\circ}/\mathrm{Fro}^\mathbb{Z},\underset{I}{\mathrm{homotopycolimit}}~\breve{\Phi}^I_{*,\Gamma,\circ}/\mathrm{Fro}^\mathbb{Z},	\\
}
\]
\[ 
\xymatrix@R+0pc@C+0pc{
\underset{r}{\mathrm{homotopylimit}}~\mathrm{Spec}^\mathrm{CS}{\Phi}^r_{*,\Gamma,\circ}/\mathrm{Fro}^\mathbb{Z},\underset{I}{\mathrm{homotopycolimit}}~\mathrm{Spec}^\mathrm{CS}{\Phi}^I_{*,\Gamma,\circ}/\mathrm{Fro}^\mathbb{Z}.	
}
\]	
In this situation we will have the target category being family parametrized by $r$ or $I$ in compatible glueing sense as in \cite[Definition 5.4.10]{10KL2}. In this situation for modules parametrized by the intervals we have the equivalence of $\infty$-categories by using \cite[Proposition 12.18]{10CS2}. Here the corresponding quasicoherent Frobenius modules are defined to be the corresponding homotopy colimits and limits of Frobenius modules:
\begin{align}
\underset{r}{\mathrm{homotopycolimit}}~M_r,\\
\underset{I}{\mathrm{homotopylimit}}~M_I,	
\end{align}
where each $M_r$ is a Frobenius-equivariant module over the period ring with respect to some radius $r$ while each $M_I$ is a Frobenius-equivariant module over the period ring with respect to some interval $I$.\\
\end{proposition}

\begin{proposition}
Similar proposition holds for 
\begin{align}
\mathrm{Quasicoherentsheaves,Perfectcomplex,IndBanach}_{*}.	
\end{align}	
\end{proposition}

\subsubsection{Frobenius Quasicoherent Prestacks III: Deformation in $(\infty,1)$-Ind-Preadic Spaces}

\begin{definition}
We now consider the pro-\'etale site of $\mathrm{Spa}\mathbb{Q}_p\left<X_1^{\pm 1},...,X_k^{\pm 1}\right>$, denote that by $*$. To be more accurate we replace one component for $\Gamma$ with the pro-\'etale site of $\mathrm{Spa}\mathbb{Q}_p\left<X_1^{\pm 1},...,X_k^{\pm 1}\right>$. And we treat then all the functor to be prestacks for this site. Then from \cite{10KL1} and \cite[Definition 5.2.1]{10KL2} we have the following class of Kedlaya-Liu rings (with the following replacement: $\Delta$ stands for $A$, $\nabla$ stands for $B$, while $\Phi$ stands for $C$) by taking product in the sense of self $\Gamma$-th power:

\[
\xymatrix@R+0pc@C+0pc{
\widetilde{\Delta}_{*,\Gamma},\widetilde{\nabla}_{*,\Gamma},\widetilde{\Phi}_{*,\Gamma},\widetilde{\Delta}^+_{*,\Gamma},\widetilde{\nabla}^+_{*,\Gamma},\widetilde{\Delta}^\dagger_{*,\Gamma},\widetilde{\nabla}^\dagger_{*,\Gamma},\widetilde{\Phi}^r_{*,\Gamma},\widetilde{\Phi}^I_{*,\Gamma}, 
}
\]

\[
\xymatrix@R+0pc@C+0pc{
\breve{\Delta}_{*,\Gamma},\breve{\nabla}_{*,\Gamma},\breve{\Phi}_{*,\Gamma},\breve{\Delta}^+_{*,\Gamma},\breve{\nabla}^+_{*,\Gamma},\breve{\Delta}^\dagger_{*,\Gamma},\breve{\nabla}^\dagger_{*,\Gamma},\breve{\Phi}^r_{*,\Gamma},\breve{\Phi}^I_{*,\Gamma},	
}
\]

\[
\xymatrix@R+0pc@C+0pc{
{\Delta}_{*,\Gamma},{\nabla}_{*,\Gamma},{\Phi}_{*,\Gamma},{\Delta}^+_{*,\Gamma},{\nabla}^+_{*,\Gamma},{\Delta}^\dagger_{*,\Gamma},{\nabla}^\dagger_{*,\Gamma},{\Phi}^r_{*,\Gamma},{\Phi}^I_{*,\Gamma}.	
}
\]
Taking the product we have:
\[
\xymatrix@R+0pc@C+0pc{
\widetilde{\Phi}_{*,\Gamma,X_\square},\widetilde{\Phi}^r_{*,\Gamma,X_\square},\widetilde{\Phi}^I_{*,\Gamma,X_\square},	
}
\]
\[
\xymatrix@R+0pc@C+0pc{
\breve{\Phi}_{*,\Gamma,X_\square},\breve{\Phi}^r_{*,\Gamma,X_\square},\breve{\Phi}^I_{*,\Gamma,X_\square},	
}
\]
\[
\xymatrix@R+0pc@C+0pc{
{\Phi}_{*,\Gamma,X_\square},{\Phi}^r_{*,\Gamma,X_\square},{\Phi}^I_{*,\Gamma,X_\square}.	
}
\]
They carry multi Frobenius action $\varphi_\Gamma$ and multi $\mathrm{Lie}_\Gamma:=\mathbb{Z}_p^{\times\Gamma}$ action. In our current situation after \cite{10CKZ} and \cite{10PZ} we consider the following $(\infty,1)$-categories of $(\infty,1)$-modules.\\
\end{definition}

\begin{definition}
First we consider the Bambozzi-Kremnizer spectrum $\mathrm{Spec}^\mathrm{BK}(*)$ attached to any of those in the above from \cite{10BK} by taking derived rational localization:
\begin{align}
&\mathrm{Spec}^\mathrm{BK}\widetilde{\Phi}_{*,\Gamma,X_\square},\mathrm{Spec}^\mathrm{BK}\widetilde{\Phi}^r_{*,\Gamma,X_\square},\mathrm{Spec}^\mathrm{BK}\widetilde{\Phi}^I_{*,\Gamma,X_\square},	
\end{align}
\begin{align}
&\mathrm{Spec}^\mathrm{BK}\breve{\Phi}_{*,\Gamma,X_\square},\mathrm{Spec}^\mathrm{BK}\breve{\Phi}^r_{*,\Gamma,X_\square},\mathrm{Spec}^\mathrm{BK}\breve{\Phi}^I_{*,\Gamma,X_\square},	
\end{align}
\begin{align}
&\mathrm{Spec}^\mathrm{BK}{\Phi}_{*,\Gamma,X_\square},
\mathrm{Spec}^\mathrm{BK}{\Phi}^r_{*,\Gamma,X_\square},\mathrm{Spec}^\mathrm{BK}{\Phi}^I_{*,\Gamma,X_\square}.	
\end{align}

Then we take the corresponding quotients by using the corresponding Frobenius operators:
\begin{align}
&\mathrm{Spec}^\mathrm{BK}\widetilde{\Phi}_{*,\Gamma,X_\square}/\mathrm{Fro}^\mathbb{Z},	\\
\end{align}
\begin{align}
&\mathrm{Spec}^\mathrm{BK}\breve{\Phi}_{*,\Gamma,X_\square}/\mathrm{Fro}^\mathbb{Z},	\\
\end{align}
\begin{align}
&\mathrm{Spec}^\mathrm{BK}{\Phi}_{*,\Gamma,X_\square}/\mathrm{Fro}^\mathbb{Z}.	
\end{align}
Here for those space without notation related to the radius and the corresponding interval we consider the total unions $\bigcap_r,\bigcup_I$ in order to achieve the whole spaces to achieve the analogues of the corresponding FF curves from \cite{10KL1}, \cite{10KL2}, \cite{10FF} for
\[
\xymatrix@R+0pc@C+0pc{
\underset{r}{\mathrm{homotopylimit}}~\mathrm{Spec}^\mathrm{BK}\widetilde{\Phi}^r_{*,\Gamma,X_\square},\underset{I}{\mathrm{homotopycolimit}}~\mathrm{Spec}^\mathrm{BK}\widetilde{\Phi}^I_{*,\Gamma,X_\square},	\\
}
\]
\[
\xymatrix@R+0pc@C+0pc{
\underset{r}{\mathrm{homotopylimit}}~\mathrm{Spec}^\mathrm{BK}\breve{\Phi}^r_{*,\Gamma,X_\square},\underset{I}{\mathrm{homotopycolimit}}~\mathrm{Spec}^\mathrm{BK}\breve{\Phi}^I_{*,\Gamma,X_\square},	\\
}
\]
\[
\xymatrix@R+0pc@C+0pc{
\underset{r}{\mathrm{homotopylimit}}~\mathrm{Spec}^\mathrm{BK}{\Phi}^r_{*,\Gamma,X_\square},\underset{I}{\mathrm{homotopycolimit}}~\mathrm{Spec}^\mathrm{BK}{\Phi}^I_{*,\Gamma,X_\square}.	
}
\]
\[  
\xymatrix@R+0pc@C+0pc{
\underset{r}{\mathrm{homotopylimit}}~\mathrm{Spec}^\mathrm{BK}\widetilde{\Phi}^r_{*,\Gamma,X_\square}/\mathrm{Fro}^\mathbb{Z},\underset{I}{\mathrm{homotopycolimit}}~\mathrm{Spec}^\mathrm{BK}\widetilde{\Phi}^I_{*,\Gamma,X_\square}/\mathrm{Fro}^\mathbb{Z},	\\
}
\]
\[ 
\xymatrix@R+0pc@C+0pc{
\underset{r}{\mathrm{homotopylimit}}~\mathrm{Spec}^\mathrm{BK}\breve{\Phi}^r_{*,\Gamma,X_\square}/\mathrm{Fro}^\mathbb{Z},\underset{I}{\mathrm{homotopycolimit}}~\mathrm{Spec}^\mathrm{BK}\breve{\Phi}^I_{*,\Gamma,X_\square}/\mathrm{Fro}^\mathbb{Z},	\\
}
\]
\[ 
\xymatrix@R+0pc@C+0pc{
\underset{r}{\mathrm{homotopylimit}}~\mathrm{Spec}^\mathrm{BK}{\Phi}^r_{*,\Gamma,X_\square}/\mathrm{Fro}^\mathbb{Z},\underset{I}{\mathrm{homotopycolimit}}~\mathrm{Spec}^\mathrm{BK}{\Phi}^I_{*,\Gamma,X_\square}/\mathrm{Fro}^\mathbb{Z}.	
}
\]

\end{definition}

\indent Meanwhile we have the corresponding Clausen-Scholze analytic stacks from \cite{10CS2}, therefore applying their construction we have:

\begin{definition}
Here we define the following products by using the solidified tensor product from \cite{10CS1} and \cite{10CS2}. Then we take solidified tensor product $\overset{\blacksquare}{\otimes}$ of any of the following
\[
\xymatrix@R+0pc@C+0pc{
\widetilde{\Delta}_{*,\Gamma},\widetilde{\nabla}_{*,\Gamma},\widetilde{\Phi}_{*,\Gamma},\widetilde{\Delta}^+_{*,\Gamma},\widetilde{\nabla}^+_{*,\Gamma},\widetilde{\Delta}^\dagger_{*,\Gamma},\widetilde{\nabla}^\dagger_{*,\Gamma},\widetilde{\Phi}^r_{*,\Gamma},\widetilde{\Phi}^I_{*,\Gamma}, 
}
\]

\[
\xymatrix@R+0pc@C+0pc{
\breve{\Delta}_{*,\Gamma},\breve{\nabla}_{*,\Gamma},\breve{\Phi}_{*,\Gamma},\breve{\Delta}^+_{*,\Gamma},\breve{\nabla}^+_{*,\Gamma},\breve{\Delta}^\dagger_{*,\Gamma},\breve{\nabla}^\dagger_{*,\Gamma},\breve{\Phi}^r_{*,\Gamma},\breve{\Phi}^I_{*,\Gamma},	
}
\]

\[
\xymatrix@R+0pc@C+0pc{
{\Delta}_{*,\Gamma},{\nabla}_{*,\Gamma},{\Phi}_{*,\Gamma},{\Delta}^+_{*,\Gamma},{\nabla}^+_{*,\Gamma},{\Delta}^\dagger_{*,\Gamma},{\nabla}^\dagger_{*,\Gamma},{\Phi}^r_{*,\Gamma},{\Phi}^I_{*,\Gamma},	
}
\]  	
with $X_\square$. Then we have the notations:
\[
\xymatrix@R+0pc@C+0pc{
\widetilde{\Delta}_{*,\Gamma,X_\square},\widetilde{\nabla}_{*,\Gamma,X_\square},\widetilde{\Phi}_{*,\Gamma,X_\square},\widetilde{\Delta}^+_{*,\Gamma,X_\square},\widetilde{\nabla}^+_{*,\Gamma,X_\square},\widetilde{\Delta}^\dagger_{*,\Gamma,X_\square},\widetilde{\nabla}^\dagger_{*,\Gamma,X_\square},\widetilde{\Phi}^r_{*,\Gamma,X_\square},\widetilde{\Phi}^I_{*,\Gamma,X_\square}, 
}
\]

\[
\xymatrix@R+0pc@C+0pc{
\breve{\Delta}_{*,\Gamma,X_\square},\breve{\nabla}_{*,\Gamma,X_\square},\breve{\Phi}_{*,\Gamma,X_\square},\breve{\Delta}^+_{*,\Gamma,X_\square},\breve{\nabla}^+_{*,\Gamma,X_\square},\breve{\Delta}^\dagger_{*,\Gamma,X_\square},\breve{\nabla}^\dagger_{*,\Gamma,X_\square},\breve{\Phi}^r_{*,\Gamma,X_\square},\breve{\Phi}^I_{*,\Gamma,X_\square},	
}
\]

\[
\xymatrix@R+0pc@C+0pc{
{\Delta}_{*,\Gamma,X_\square},{\nabla}_{*,\Gamma,X_\square},{\Phi}_{*,\Gamma,X_\square},{\Delta}^+_{*,\Gamma,X_\square},{\nabla}^+_{*,\Gamma,X_\square},{\Delta}^\dagger_{*,\Gamma,X_\square},{\nabla}^\dagger_{*,\Gamma,X_\square},{\Phi}^r_{*,\Gamma,X_\square},{\Phi}^I_{*,\Gamma,X_\square}.	
}
\]
\end{definition}

\begin{definition}
First we consider the Clausen-Scholze spectrum $\mathrm{Spec}^\mathrm{CS}(*)$ attached to any of those in the above from \cite{10CS2} by taking derived rational localization:
\begin{align}
\mathrm{Spec}^\mathrm{CS}\widetilde{\Delta}_{*,\Gamma,X_\square},\mathrm{Spec}^\mathrm{CS}\widetilde{\nabla}_{*,\Gamma,X_\square},\mathrm{Spec}^\mathrm{CS}\widetilde{\Phi}_{*,\Gamma,X_\square},\mathrm{Spec}^\mathrm{CS}\widetilde{\Delta}^+_{*,\Gamma,X_\square},\mathrm{Spec}^\mathrm{CS}\widetilde{\nabla}^+_{*,\Gamma,X_\square},\\
\mathrm{Spec}^\mathrm{CS}\widetilde{\Delta}^\dagger_{*,\Gamma,X_\square},\mathrm{Spec}^\mathrm{CS}\widetilde{\nabla}^\dagger_{*,\Gamma,X_\square},\mathrm{Spec}^\mathrm{CS}\widetilde{\Phi}^r_{*,\Gamma,X_\square},\mathrm{Spec}^\mathrm{CS}\widetilde{\Phi}^I_{*,\Gamma,X_\square},	\\
\end{align}
\begin{align}
\mathrm{Spec}^\mathrm{CS}\breve{\Delta}_{*,\Gamma,X_\square},\breve{\nabla}_{*,\Gamma,X_\square},\mathrm{Spec}^\mathrm{CS}\breve{\Phi}_{*,\Gamma,X_\square},\mathrm{Spec}^\mathrm{CS}\breve{\Delta}^+_{*,\Gamma,X_\square},\mathrm{Spec}^\mathrm{CS}\breve{\nabla}^+_{*,\Gamma,X_\square},\\
\mathrm{Spec}^\mathrm{CS}\breve{\Delta}^\dagger_{*,\Gamma,X_\square},\mathrm{Spec}^\mathrm{CS}\breve{\nabla}^\dagger_{*,\Gamma,X_\square},\mathrm{Spec}^\mathrm{CS}\breve{\Phi}^r_{*,\Gamma,X_\square},\breve{\Phi}^I_{*,\Gamma,X_\square},	\\
\end{align}
\begin{align}
\mathrm{Spec}^\mathrm{CS}{\Delta}_{*,\Gamma,X_\square},\mathrm{Spec}^\mathrm{CS}{\nabla}_{*,\Gamma,X_\square},\mathrm{Spec}^\mathrm{CS}{\Phi}_{*,\Gamma,X_\square},\mathrm{Spec}^\mathrm{CS}{\Delta}^+_{*,\Gamma,X_\square},\mathrm{Spec}^\mathrm{CS}{\nabla}^+_{*,\Gamma,X_\square},\\
\mathrm{Spec}^\mathrm{CS}{\Delta}^\dagger_{*,\Gamma,X_\square},\mathrm{Spec}^\mathrm{CS}{\nabla}^\dagger_{*,\Gamma,X_\square},\mathrm{Spec}^\mathrm{CS}{\Phi}^r_{*,\Gamma,X_\square},\mathrm{Spec}^\mathrm{CS}{\Phi}^I_{*,\Gamma,X_\square}.	
\end{align}

Then we take the corresponding quotients by using the corresponding Frobenius operators:
\begin{align}
&\mathrm{Spec}^\mathrm{CS}\widetilde{\Delta}_{*,\Gamma,X_\square}/\mathrm{Fro}^\mathbb{Z},\mathrm{Spec}^\mathrm{CS}\widetilde{\nabla}_{*,\Gamma,X_\square}/\mathrm{Fro}^\mathbb{Z},\mathrm{Spec}^\mathrm{CS}\widetilde{\Phi}_{*,\Gamma,X_\square}/\mathrm{Fro}^\mathbb{Z},\mathrm{Spec}^\mathrm{CS}\widetilde{\Delta}^+_{*,\Gamma,X_\square}/\mathrm{Fro}^\mathbb{Z},\\
&\mathrm{Spec}^\mathrm{CS}\widetilde{\nabla}^+_{*,\Gamma,X_\square}/\mathrm{Fro}^\mathbb{Z}, \mathrm{Spec}^\mathrm{CS}\widetilde{\Delta}^\dagger_{*,\Gamma,X_\square}/\mathrm{Fro}^\mathbb{Z},\mathrm{Spec}^\mathrm{CS}\widetilde{\nabla}^\dagger_{*,\Gamma,X_\square}/\mathrm{Fro}^\mathbb{Z},	\\
\end{align}
\begin{align}
&\mathrm{Spec}^\mathrm{CS}\breve{\Delta}_{*,\Gamma,X_\square}/\mathrm{Fro}^\mathbb{Z},\breve{\nabla}_{*,\Gamma,X_\square}/\mathrm{Fro}^\mathbb{Z},\mathrm{Spec}^\mathrm{CS}\breve{\Phi}_{*,\Gamma,X_\square}/\mathrm{Fro}^\mathbb{Z},\mathrm{Spec}^\mathrm{CS}\breve{\Delta}^+_{*,\Gamma,X_\square}/\mathrm{Fro}^\mathbb{Z},\\
&\mathrm{Spec}^\mathrm{CS}\breve{\nabla}^+_{*,\Gamma,X_\square}/\mathrm{Fro}^\mathbb{Z}, \mathrm{Spec}^\mathrm{CS}\breve{\Delta}^\dagger_{*,\Gamma,X_\square}/\mathrm{Fro}^\mathbb{Z},\mathrm{Spec}^\mathrm{CS}\breve{\nabla}^\dagger_{*,\Gamma,X_\square}/\mathrm{Fro}^\mathbb{Z},	\\
\end{align}
\begin{align}
&\mathrm{Spec}^\mathrm{CS}{\Delta}_{*,\Gamma,X_\square}/\mathrm{Fro}^\mathbb{Z},\mathrm{Spec}^\mathrm{CS}{\nabla}_{*,\Gamma,X_\square}/\mathrm{Fro}^\mathbb{Z},\mathrm{Spec}^\mathrm{CS}{\Phi}_{*,\Gamma,X_\square}/\mathrm{Fro}^\mathbb{Z},\mathrm{Spec}^\mathrm{CS}{\Delta}^+_{*,\Gamma,X_\square}/\mathrm{Fro}^\mathbb{Z},\\
&\mathrm{Spec}^\mathrm{CS}{\nabla}^+_{*,\Gamma,X_\square}/\mathrm{Fro}^\mathbb{Z}, \mathrm{Spec}^\mathrm{CS}{\Delta}^\dagger_{*,\Gamma,X_\square}/\mathrm{Fro}^\mathbb{Z},\mathrm{Spec}^\mathrm{CS}{\nabla}^\dagger_{*,\Gamma,X_\square}/\mathrm{Fro}^\mathbb{Z}.	
\end{align}
Here for those space with notations related to the radius and the corresponding interval we consider the total unions $\bigcap_r,\bigcup_I$ in order to achieve the whole spaces to achieve the analogues of the corresponding FF curves from \cite{10KL1}, \cite{10KL2}, \cite{10FF} for
\[
\xymatrix@R+0pc@C+0pc{
\underset{r}{\mathrm{homotopylimit}}~\mathrm{Spec}^\mathrm{CS}\widetilde{\Phi}^r_{*,\Gamma,X_\square},\underset{I}{\mathrm{homotopycolimit}}~\mathrm{Spec}^\mathrm{CS}\widetilde{\Phi}^I_{*,\Gamma,X_\square},	\\
}
\]
\[
\xymatrix@R+0pc@C+0pc{
\underset{r}{\mathrm{homotopylimit}}~\mathrm{Spec}^\mathrm{CS}\breve{\Phi}^r_{*,\Gamma,X_\square},\underset{I}{\mathrm{homotopycolimit}}~\mathrm{Spec}^\mathrm{CS}\breve{\Phi}^I_{*,\Gamma,X_\square},	\\
}
\]
\[
\xymatrix@R+0pc@C+0pc{
\underset{r}{\mathrm{homotopylimit}}~\mathrm{Spec}^\mathrm{CS}{\Phi}^r_{*,\Gamma,X_\square},\underset{I}{\mathrm{homotopycolimit}}~\mathrm{Spec}^\mathrm{CS}{\Phi}^I_{*,\Gamma,X_\square}.	
}
\]
\[ 
\xymatrix@R+0pc@C+0pc{
\underset{r}{\mathrm{homotopylimit}}~\mathrm{Spec}^\mathrm{CS}\widetilde{\Phi}^r_{*,\Gamma,X_\square}/\mathrm{Fro}^\mathbb{Z},\underset{I}{\mathrm{homotopycolimit}}~\mathrm{Spec}^\mathrm{CS}\widetilde{\Phi}^I_{*,\Gamma,X_\square}/\mathrm{Fro}^\mathbb{Z},	\\
}
\]
\[ 
\xymatrix@R+0pc@C+0pc{
\underset{r}{\mathrm{homotopylimit}}~\mathrm{Spec}^\mathrm{CS}\breve{\Phi}^r_{*,\Gamma,X_\square}/\mathrm{Fro}^\mathbb{Z},\underset{I}{\mathrm{homotopycolimit}}~\breve{\Phi}^I_{*,\Gamma,X_\square}/\mathrm{Fro}^\mathbb{Z},	\\
}
\]
\[ 
\xymatrix@R+0pc@C+0pc{
\underset{r}{\mathrm{homotopylimit}}~\mathrm{Spec}^\mathrm{CS}{\Phi}^r_{*,\Gamma,X_\square}/\mathrm{Fro}^\mathbb{Z},\underset{I}{\mathrm{homotopycolimit}}~\mathrm{Spec}^\mathrm{CS}{\Phi}^I_{*,\Gamma,X_\square}/\mathrm{Fro}^\mathbb{Z}.	
}
\]

\end{definition}

\

\begin{definition}
We then consider the corresponding quasipresheaves of the corresponding ind-Banach or monomorphic ind-Banach modules from \cite{10BBK}, \cite{10KKM}:
\begin{align}
\mathrm{Quasicoherentpresheaves,IndBanach}_{*}	
\end{align}
where $*$ is one of the following spaces:
\begin{align}
&\mathrm{Spec}^\mathrm{BK}\widetilde{\Phi}_{*,\Gamma,X_\square}/\mathrm{Fro}^\mathbb{Z},	\\
\end{align}
\begin{align}
&\mathrm{Spec}^\mathrm{BK}\breve{\Phi}_{*,\Gamma,X_\square}/\mathrm{Fro}^\mathbb{Z},	\\
\end{align}
\begin{align}
&\mathrm{Spec}^\mathrm{BK}{\Phi}_{*,\Gamma,X_\square}/\mathrm{Fro}^\mathbb{Z}.	
\end{align}
Here for those space without notation related to the radius and the corresponding interval we consider the total unions $\bigcap_r,\bigcup_I$ in order to achieve the whole spaces to achieve the analogues of the corresponding FF curves from \cite{10KL1}, \cite{10KL2}, \cite{10FF} for
\[
\xymatrix@R+0pc@C+0pc{
\underset{r}{\mathrm{homotopylimit}}~\mathrm{Spec}^\mathrm{BK}\widetilde{\Phi}^r_{*,\Gamma,X_\square},\underset{I}{\mathrm{homotopycolimit}}~\mathrm{Spec}^\mathrm{BK}\widetilde{\Phi}^I_{*,\Gamma,X_\square},	\\
}
\]
\[
\xymatrix@R+0pc@C+0pc{
\underset{r}{\mathrm{homotopylimit}}~\mathrm{Spec}^\mathrm{BK}\breve{\Phi}^r_{*,\Gamma,X_\square},\underset{I}{\mathrm{homotopycolimit}}~\mathrm{Spec}^\mathrm{BK}\breve{\Phi}^I_{*,\Gamma,X_\square},	\\
}
\]
\[
\xymatrix@R+0pc@C+0pc{
\underset{r}{\mathrm{homotopylimit}}~\mathrm{Spec}^\mathrm{BK}{\Phi}^r_{*,\Gamma,X_\square},\underset{I}{\mathrm{homotopycolimit}}~\mathrm{Spec}^\mathrm{BK}{\Phi}^I_{*,\Gamma,X_\square}.	
}
\]
\[  
\xymatrix@R+0pc@C+0pc{
\underset{r}{\mathrm{homotopylimit}}~\mathrm{Spec}^\mathrm{BK}\widetilde{\Phi}^r_{*,\Gamma,X_\square}/\mathrm{Fro}^\mathbb{Z},\underset{I}{\mathrm{homotopycolimit}}~\mathrm{Spec}^\mathrm{BK}\widetilde{\Phi}^I_{*,\Gamma,X_\square}/\mathrm{Fro}^\mathbb{Z},	\\
}
\]
\[ 
\xymatrix@R+0pc@C+0pc{
\underset{r}{\mathrm{homotopylimit}}~\mathrm{Spec}^\mathrm{BK}\breve{\Phi}^r_{*,\Gamma,X_\square}/\mathrm{Fro}^\mathbb{Z},\underset{I}{\mathrm{homotopycolimit}}~\mathrm{Spec}^\mathrm{BK}\breve{\Phi}^I_{*,\Gamma,X_\square}/\mathrm{Fro}^\mathbb{Z},	\\
}
\]
\[ 
\xymatrix@R+0pc@C+0pc{
\underset{r}{\mathrm{homotopylimit}}~\mathrm{Spec}^\mathrm{BK}{\Phi}^r_{*,\Gamma,X_\square}/\mathrm{Fro}^\mathbb{Z},\underset{I}{\mathrm{homotopycolimit}}~\mathrm{Spec}^\mathrm{BK}{\Phi}^I_{*,\Gamma,X_\square}/\mathrm{Fro}^\mathbb{Z}.	
}
\]

\end{definition}

\begin{definition}
We then consider the corresponding quasisheaves of the corresponding condensed solid topological modules from \cite{10CS2}:
\begin{align}
\mathrm{Quasicoherentsheaves, Condensed}_{*}	
\end{align}
where $*$ is one of the following spaces:
\begin{align}
&\mathrm{Spec}^\mathrm{CS}\widetilde{\Delta}_{*,\Gamma,X_\square}/\mathrm{Fro}^\mathbb{Z},\mathrm{Spec}^\mathrm{CS}\widetilde{\nabla}_{*,\Gamma,X_\square}/\mathrm{Fro}^\mathbb{Z},\mathrm{Spec}^\mathrm{CS}\widetilde{\Phi}_{*,\Gamma,X_\square}/\mathrm{Fro}^\mathbb{Z},\mathrm{Spec}^\mathrm{CS}\widetilde{\Delta}^+_{*,\Gamma,X_\square}/\mathrm{Fro}^\mathbb{Z},\\
&\mathrm{Spec}^\mathrm{CS}\widetilde{\nabla}^+_{*,\Gamma,X_\square}/\mathrm{Fro}^\mathbb{Z},\mathrm{Spec}^\mathrm{CS}\widetilde{\Delta}^\dagger_{*,\Gamma,X_\square}/\mathrm{Fro}^\mathbb{Z},\mathrm{Spec}^\mathrm{CS}\widetilde{\nabla}^\dagger_{*,\Gamma,X_\square}/\mathrm{Fro}^\mathbb{Z},	\\
\end{align}
\begin{align}
&\mathrm{Spec}^\mathrm{CS}\breve{\Delta}_{*,\Gamma,X_\square}/\mathrm{Fro}^\mathbb{Z},\breve{\nabla}_{*,\Gamma,X_\square}/\mathrm{Fro}^\mathbb{Z},\mathrm{Spec}^\mathrm{CS}\breve{\Phi}_{*,\Gamma,X_\square}/\mathrm{Fro}^\mathbb{Z},\mathrm{Spec}^\mathrm{CS}\breve{\Delta}^+_{*,\Gamma,X_\square}/\mathrm{Fro}^\mathbb{Z},\\
&\mathrm{Spec}^\mathrm{CS}\breve{\nabla}^+_{*,\Gamma,X_\square}/\mathrm{Fro}^\mathbb{Z},\mathrm{Spec}^\mathrm{CS}\breve{\Delta}^\dagger_{*,\Gamma,X_\square}/\mathrm{Fro}^\mathbb{Z},\mathrm{Spec}^\mathrm{CS}\breve{\nabla}^\dagger_{*,\Gamma,X_\square}/\mathrm{Fro}^\mathbb{Z},	\\
\end{align}
\begin{align}
&\mathrm{Spec}^\mathrm{CS}{\Delta}_{*,\Gamma,X_\square}/\mathrm{Fro}^\mathbb{Z},\mathrm{Spec}^\mathrm{CS}{\nabla}_{*,\Gamma,X_\square}/\mathrm{Fro}^\mathbb{Z},\mathrm{Spec}^\mathrm{CS}{\Phi}_{*,\Gamma,X_\square}/\mathrm{Fro}^\mathbb{Z},\mathrm{Spec}^\mathrm{CS}{\Delta}^+_{*,\Gamma,X_\square}/\mathrm{Fro}^\mathbb{Z},\\
&\mathrm{Spec}^\mathrm{CS}{\nabla}^+_{*,\Gamma,X_\square}/\mathrm{Fro}^\mathbb{Z}, \mathrm{Spec}^\mathrm{CS}{\Delta}^\dagger_{*,\Gamma,X_\square}/\mathrm{Fro}^\mathbb{Z},\mathrm{Spec}^\mathrm{CS}{\nabla}^\dagger_{*,\Gamma,X_\square}/\mathrm{Fro}^\mathbb{Z}.	
\end{align}
Here for those space with notations related to the radius and the corresponding interval we consider the total unions $\bigcap_r,\bigcup_I$ in order to achieve the whole spaces to achieve the analogues of the corresponding FF curves from \cite{10KL1}, \cite{10KL2}, \cite{10FF} for
\[
\xymatrix@R+0pc@C+0pc{
\underset{r}{\mathrm{homotopylimit}}~\mathrm{Spec}^\mathrm{CS}\widetilde{\Phi}^r_{*,\Gamma,X_\square},\underset{I}{\mathrm{homotopycolimit}}~\mathrm{Spec}^\mathrm{CS}\widetilde{\Phi}^I_{*,\Gamma,X_\square},	\\
}
\]
\[
\xymatrix@R+0pc@C+0pc{
\underset{r}{\mathrm{homotopylimit}}~\mathrm{Spec}^\mathrm{CS}\breve{\Phi}^r_{*,\Gamma,X_\square},\underset{I}{\mathrm{homotopycolimit}}~\mathrm{Spec}^\mathrm{CS}\breve{\Phi}^I_{*,\Gamma,X_\square},	\\
}
\]
\[
\xymatrix@R+0pc@C+0pc{
\underset{r}{\mathrm{homotopylimit}}~\mathrm{Spec}^\mathrm{CS}{\Phi}^r_{*,\Gamma,X_\square},\underset{I}{\mathrm{homotopycolimit}}~\mathrm{Spec}^\mathrm{CS}{\Phi}^I_{*,\Gamma,X_\square}.	
}
\]
\[ 
\xymatrix@R+0pc@C+0pc{
\underset{r}{\mathrm{homotopylimit}}~\mathrm{Spec}^\mathrm{CS}\widetilde{\Phi}^r_{*,\Gamma,X_\square}/\mathrm{Fro}^\mathbb{Z},\underset{I}{\mathrm{homotopycolimit}}~\mathrm{Spec}^\mathrm{CS}\widetilde{\Phi}^I_{*,\Gamma,X_\square}/\mathrm{Fro}^\mathbb{Z},	\\
}
\]
\[ 
\xymatrix@R+0pc@C+0pc{
\underset{r}{\mathrm{homotopylimit}}~\mathrm{Spec}^\mathrm{CS}\breve{\Phi}^r_{*,\Gamma,X_\square}/\mathrm{Fro}^\mathbb{Z},\underset{I}{\mathrm{homotopycolimit}}~\breve{\Phi}^I_{*,\Gamma,X_\square}/\mathrm{Fro}^\mathbb{Z},	\\
}
\]
\[ 
\xymatrix@R+0pc@C+0pc{
\underset{r}{\mathrm{homotopylimit}}~\mathrm{Spec}^\mathrm{CS}{\Phi}^r_{*,\Gamma,X_\square}/\mathrm{Fro}^\mathbb{Z},\underset{I}{\mathrm{homotopycolimit}}~\mathrm{Spec}^\mathrm{CS}{\Phi}^I_{*,\Gamma,X_\square}/\mathrm{Fro}^\mathbb{Z}.	
}
\]

\end{definition}

\

\begin{proposition}
There is a well-defined functor from the $\infty$-category 
\begin{align}
\mathrm{Quasicoherentpresheaves,Condensed}_{*}	
\end{align}
where $*$ is one of the following spaces:
\begin{align}
&\mathrm{Spec}^\mathrm{CS}\widetilde{\Phi}_{*,\Gamma,X_\square}/\mathrm{Fro}^\mathbb{Z},	\\
\end{align}
\begin{align}
&\mathrm{Spec}^\mathrm{CS}\breve{\Phi}_{*,\Gamma,X_\square}/\mathrm{Fro}^\mathbb{Z},	\\
\end{align}
\begin{align}
&\mathrm{Spec}^\mathrm{CS}{\Phi}_{*,\Gamma,X_\square}/\mathrm{Fro}^\mathbb{Z},	
\end{align}
to the $\infty$-category of $\mathrm{Fro}$-equivariant quasicoherent presheaves over similar spaces above correspondingly without the $\mathrm{Fro}$-quotients, and to the $\infty$-category of $\mathrm{Fro}$-equivariant quasicoherent modules over global sections of the structure $\infty$-sheaves of the similar spaces above correspondingly without the $\mathrm{Fro}$-quotients. Here for those space without notation related to the radius and the corresponding interval we consider the total unions $\bigcap_r,\bigcup_I$ in order to achieve the whole spaces to achieve the analogues of the corresponding FF curves from \cite{10KL1}, \cite{10KL2}, \cite{10FF} for
\[
\xymatrix@R+0pc@C+0pc{
\underset{r}{\mathrm{homotopylimit}}~\mathrm{Spec}^\mathrm{CS}\widetilde{\Phi}^r_{*,\Gamma,X_\square},\underset{I}{\mathrm{homotopycolimit}}~\mathrm{Spec}^\mathrm{CS}\widetilde{\Phi}^I_{*,\Gamma,X_\square},	\\
}
\]
\[
\xymatrix@R+0pc@C+0pc{
\underset{r}{\mathrm{homotopylimit}}~\mathrm{Spec}^\mathrm{CS}\breve{\Phi}^r_{*,\Gamma,X_\square},\underset{I}{\mathrm{homotopycolimit}}~\mathrm{Spec}^\mathrm{CS}\breve{\Phi}^I_{*,\Gamma,X_\square},	\\
}
\]
\[
\xymatrix@R+0pc@C+0pc{
\underset{r}{\mathrm{homotopylimit}}~\mathrm{Spec}^\mathrm{CS}{\Phi}^r_{*,\Gamma,X_\square},\underset{I}{\mathrm{homotopycolimit}}~\mathrm{Spec}^\mathrm{CS}{\Phi}^I_{*,\Gamma,X_\square}.	
}
\]
\[ 
\xymatrix@R+0pc@C+0pc{
\underset{r}{\mathrm{homotopylimit}}~\mathrm{Spec}^\mathrm{CS}\widetilde{\Phi}^r_{*,\Gamma,X_\square}/\mathrm{Fro}^\mathbb{Z},\underset{I}{\mathrm{homotopycolimit}}~\mathrm{Spec}^\mathrm{CS}\widetilde{\Phi}^I_{*,\Gamma,X_\square}/\mathrm{Fro}^\mathbb{Z},	\\
}
\]
\[ 
\xymatrix@R+0pc@C+0pc{
\underset{r}{\mathrm{homotopylimit}}~\mathrm{Spec}^\mathrm{CS}\breve{\Phi}^r_{*,\Gamma,X_\square}/\mathrm{Fro}^\mathbb{Z},\underset{I}{\mathrm{homotopycolimit}}~\breve{\Phi}^I_{*,\Gamma,X_\square}/\mathrm{Fro}^\mathbb{Z},	\\
}
\]
\[ 
\xymatrix@R+0pc@C+0pc{
\underset{r}{\mathrm{homotopylimit}}~\mathrm{Spec}^\mathrm{CS}{\Phi}^r_{*,\Gamma,X_\square}/\mathrm{Fro}^\mathbb{Z},\underset{I}{\mathrm{homotopycolimit}}~\mathrm{Spec}^\mathrm{CS}{\Phi}^I_{*,\Gamma,X_\square}/\mathrm{Fro}^\mathbb{Z}.	
}
\]	
In this situation we will have the target category being family parametrized by $r$ or $I$ in compatible glueing sense as in \cite[Definition 5.4.10]{10KL2}. In this situation for modules parametrized by the intervals we have the equivalence of $\infty$-categories by using \cite[Proposition 13.8]{10CS2}. Here the corresponding quasicoherent Frobenius modules are defined to be the corresponding homotopy colimits and limits of Frobenius modules:
\begin{align}
\underset{r}{\mathrm{homotopycolimit}}~M_r,\\
\underset{I}{\mathrm{homotopylimit}}~M_I,	
\end{align}
where each $M_r$ is a Frobenius-equivariant module over the period ring with respect to some radius $r$ while each $M_I$ is a Frobenius-equivariant module over the period ring with respect to some interval $I$.\\
\end{proposition}

\begin{proposition}
Similar proposition holds for 
\begin{align}
\mathrm{Quasicoherentsheaves,IndBanach}_{*}.	
\end{align}	
\end{proposition}

\

\begin{definition}
We then consider the corresponding quasipresheaves of perfect complexes the corresponding ind-Banach or monomorphic ind-Banach modules from \cite{10BBK}, \cite{10KKM}:
\begin{align}
\mathrm{Quasicoherentpresheaves,Perfectcomplex,IndBanach}_{*}	
\end{align}
where $*$ is one of the following spaces:
\begin{align}
&\mathrm{Spec}^\mathrm{BK}\widetilde{\Phi}_{*,\Gamma,X_\square}/\mathrm{Fro}^\mathbb{Z},	\\
\end{align}
\begin{align}
&\mathrm{Spec}^\mathrm{BK}\breve{\Phi}_{*,\Gamma,X_\square}/\mathrm{Fro}^\mathbb{Z},	\\
\end{align}
\begin{align}
&\mathrm{Spec}^\mathrm{BK}{\Phi}_{*,\Gamma,X_\square}/\mathrm{Fro}^\mathbb{Z}.	
\end{align}
Here for those space without notation related to the radius and the corresponding interval we consider the total unions $\bigcap_r,\bigcup_I$ in order to achieve the whole spaces to achieve the analogues of the corresponding FF curves from \cite{10KL1}, \cite{10KL2}, \cite{10FF} for
\[
\xymatrix@R+0pc@C+0pc{
\underset{r}{\mathrm{homotopylimit}}~\mathrm{Spec}^\mathrm{BK}\widetilde{\Phi}^r_{*,\Gamma,X_\square},\underset{I}{\mathrm{homotopycolimit}}~\mathrm{Spec}^\mathrm{BK}\widetilde{\Phi}^I_{*,\Gamma,X_\square},	\\
}
\]
\[
\xymatrix@R+0pc@C+0pc{
\underset{r}{\mathrm{homotopylimit}}~\mathrm{Spec}^\mathrm{BK}\breve{\Phi}^r_{*,\Gamma,X_\square},\underset{I}{\mathrm{homotopycolimit}}~\mathrm{Spec}^\mathrm{BK}\breve{\Phi}^I_{*,\Gamma,X_\square},	\\
}
\]
\[
\xymatrix@R+0pc@C+0pc{
\underset{r}{\mathrm{homotopylimit}}~\mathrm{Spec}^\mathrm{BK}{\Phi}^r_{*,\Gamma,X_\square},\underset{I}{\mathrm{homotopycolimit}}~\mathrm{Spec}^\mathrm{BK}{\Phi}^I_{*,\Gamma,X_\square}.	
}
\]
\[  
\xymatrix@R+0pc@C+0pc{
\underset{r}{\mathrm{homotopylimit}}~\mathrm{Spec}^\mathrm{BK}\widetilde{\Phi}^r_{*,\Gamma,X_\square}/\mathrm{Fro}^\mathbb{Z},\underset{I}{\mathrm{homotopycolimit}}~\mathrm{Spec}^\mathrm{BK}\widetilde{\Phi}^I_{*,\Gamma,X_\square}/\mathrm{Fro}^\mathbb{Z},	\\
}
\]
\[ 
\xymatrix@R+0pc@C+0pc{
\underset{r}{\mathrm{homotopylimit}}~\mathrm{Spec}^\mathrm{BK}\breve{\Phi}^r_{*,\Gamma,X_\square}/\mathrm{Fro}^\mathbb{Z},\underset{I}{\mathrm{homotopycolimit}}~\mathrm{Spec}^\mathrm{BK}\breve{\Phi}^I_{*,\Gamma,X_\square}/\mathrm{Fro}^\mathbb{Z},	\\
}
\]
\[ 
\xymatrix@R+0pc@C+0pc{
\underset{r}{\mathrm{homotopylimit}}~\mathrm{Spec}^\mathrm{BK}{\Phi}^r_{*,\Gamma,X_\square}/\mathrm{Fro}^\mathbb{Z},\underset{I}{\mathrm{homotopycolimit}}~\mathrm{Spec}^\mathrm{BK}{\Phi}^I_{*,\Gamma,X_\square}/\mathrm{Fro}^\mathbb{Z}.	
}
\]

\end{definition}

\begin{definition}
We then consider the corresponding quasisheaves of perfect complexes of the corresponding condensed solid topological modules from \cite{10CS2}:
\begin{align}
\mathrm{Quasicoherentsheaves, Perfectcomplex, Condensed}_{*}	
\end{align}
where $*$ is one of the following spaces:
\begin{align}
&\mathrm{Spec}^\mathrm{CS}\widetilde{\Delta}_{*,\Gamma,X_\square}/\mathrm{Fro}^\mathbb{Z},\mathrm{Spec}^\mathrm{CS}\widetilde{\nabla}_{*,\Gamma,X_\square}/\mathrm{Fro}^\mathbb{Z},\mathrm{Spec}^\mathrm{CS}\widetilde{\Phi}_{*,\Gamma,X_\square}/\mathrm{Fro}^\mathbb{Z},\mathrm{Spec}^\mathrm{CS}\widetilde{\Delta}^+_{*,\Gamma,X_\square}/\mathrm{Fro}^\mathbb{Z},\\
&\mathrm{Spec}^\mathrm{CS}\widetilde{\nabla}^+_{*,\Gamma,X_\square}/\mathrm{Fro}^\mathbb{Z},\mathrm{Spec}^\mathrm{CS}\widetilde{\Delta}^\dagger_{*,\Gamma,X_\square}/\mathrm{Fro}^\mathbb{Z},\mathrm{Spec}^\mathrm{CS}\widetilde{\nabla}^\dagger_{*,\Gamma,X_\square}/\mathrm{Fro}^\mathbb{Z},	\\
\end{align}
\begin{align}
&\mathrm{Spec}^\mathrm{CS}\breve{\Delta}_{*,\Gamma,X_\square}/\mathrm{Fro}^\mathbb{Z},\breve{\nabla}_{*,\Gamma,X_\square}/\mathrm{Fro}^\mathbb{Z},\mathrm{Spec}^\mathrm{CS}\breve{\Phi}_{*,\Gamma,X_\square}/\mathrm{Fro}^\mathbb{Z},\mathrm{Spec}^\mathrm{CS}\breve{\Delta}^+_{*,\Gamma,X_\square}/\mathrm{Fro}^\mathbb{Z},\\
&\mathrm{Spec}^\mathrm{CS}\breve{\nabla}^+_{*,\Gamma,X_\square}/\mathrm{Fro}^\mathbb{Z},\mathrm{Spec}^\mathrm{CS}\breve{\Delta}^\dagger_{*,\Gamma,X_\square}/\mathrm{Fro}^\mathbb{Z},\mathrm{Spec}^\mathrm{CS}\breve{\nabla}^\dagger_{*,\Gamma,X_\square}/\mathrm{Fro}^\mathbb{Z},	\\
\end{align}
\begin{align}
&\mathrm{Spec}^\mathrm{CS}{\Delta}_{*,\Gamma,X_\square}/\mathrm{Fro}^\mathbb{Z},\mathrm{Spec}^\mathrm{CS}{\nabla}_{*,\Gamma,X_\square}/\mathrm{Fro}^\mathbb{Z},\mathrm{Spec}^\mathrm{CS}{\Phi}_{*,\Gamma,X_\square}/\mathrm{Fro}^\mathbb{Z},\mathrm{Spec}^\mathrm{CS}{\Delta}^+_{*,\Gamma,X_\square}/\mathrm{Fro}^\mathbb{Z},\\
&\mathrm{Spec}^\mathrm{CS}{\nabla}^+_{*,\Gamma,X_\square}/\mathrm{Fro}^\mathbb{Z}, \mathrm{Spec}^\mathrm{CS}{\Delta}^\dagger_{*,\Gamma,X_\square}/\mathrm{Fro}^\mathbb{Z},\mathrm{Spec}^\mathrm{CS}{\nabla}^\dagger_{*,\Gamma,X_\square}/\mathrm{Fro}^\mathbb{Z}.	
\end{align}
Here for those space with notations related to the radius and the corresponding interval we consider the total unions $\bigcap_r,\bigcup_I$ in order to achieve the whole spaces to achieve the analogues of the corresponding FF curves from \cite{10KL1}, \cite{10KL2}, \cite{10FF} for
\[
\xymatrix@R+0pc@C+0pc{
\underset{r}{\mathrm{homotopylimit}}~\mathrm{Spec}^\mathrm{CS}\widetilde{\Phi}^r_{*,\Gamma,X_\square},\underset{I}{\mathrm{homotopycolimit}}~\mathrm{Spec}^\mathrm{CS}\widetilde{\Phi}^I_{*,\Gamma,X_\square},	\\
}
\]
\[
\xymatrix@R+0pc@C+0pc{
\underset{r}{\mathrm{homotopylimit}}~\mathrm{Spec}^\mathrm{CS}\breve{\Phi}^r_{*,\Gamma,X_\square},\underset{I}{\mathrm{homotopycolimit}}~\mathrm{Spec}^\mathrm{CS}\breve{\Phi}^I_{*,\Gamma,X_\square},	\\
}
\]
\[
\xymatrix@R+0pc@C+0pc{
\underset{r}{\mathrm{homotopylimit}}~\mathrm{Spec}^\mathrm{CS}{\Phi}^r_{*,\Gamma,X_\square},\underset{I}{\mathrm{homotopycolimit}}~\mathrm{Spec}^\mathrm{CS}{\Phi}^I_{*,\Gamma,X_\square}.	
}
\]
\[ 
\xymatrix@R+0pc@C+0pc{
\underset{r}{\mathrm{homotopylimit}}~\mathrm{Spec}^\mathrm{CS}\widetilde{\Phi}^r_{*,\Gamma,X_\square}/\mathrm{Fro}^\mathbb{Z},\underset{I}{\mathrm{homotopycolimit}}~\mathrm{Spec}^\mathrm{CS}\widetilde{\Phi}^I_{*,\Gamma,X_\square}/\mathrm{Fro}^\mathbb{Z},	\\
}
\]
\[ 
\xymatrix@R+0pc@C+0pc{
\underset{r}{\mathrm{homotopylimit}}~\mathrm{Spec}^\mathrm{CS}\breve{\Phi}^r_{*,\Gamma,X_\square}/\mathrm{Fro}^\mathbb{Z},\underset{I}{\mathrm{homotopycolimit}}~\breve{\Phi}^I_{*,\Gamma,X_\square}/\mathrm{Fro}^\mathbb{Z},	\\
}
\]
\[ 
\xymatrix@R+0pc@C+0pc{
\underset{r}{\mathrm{homotopylimit}}~\mathrm{Spec}^\mathrm{CS}{\Phi}^r_{*,\Gamma,X_\square}/\mathrm{Fro}^\mathbb{Z},\underset{I}{\mathrm{homotopycolimit}}~\mathrm{Spec}^\mathrm{CS}{\Phi}^I_{*,\Gamma,X_\square}/\mathrm{Fro}^\mathbb{Z}.	
}
\]

\end{definition}

\begin{proposition}
There is a well-defined functor from the $\infty$-category 
\begin{align}
\mathrm{Quasicoherentpresheaves,Perfectcomplex,Condensed}_{*}	
\end{align}
where $*$ is one of the following spaces:
\begin{align}
&\mathrm{Spec}^\mathrm{CS}\widetilde{\Phi}_{*,\Gamma,X_\square}/\mathrm{Fro}^\mathbb{Z},	\\
\end{align}
\begin{align}
&\mathrm{Spec}^\mathrm{CS}\breve{\Phi}_{*,\Gamma,X_\square}/\mathrm{Fro}^\mathbb{Z},	\\
\end{align}
\begin{align}
&\mathrm{Spec}^\mathrm{CS}{\Phi}_{*,\Gamma,X_\square}/\mathrm{Fro}^\mathbb{Z},	
\end{align}
to the $\infty$-category of $\mathrm{Fro}$-equivariant quasicoherent presheaves over similar spaces above correspondingly without the $\mathrm{Fro}$-quotients, and to the $\infty$-category of $\mathrm{Fro}$-equivariant quasicoherent modules over global sections of the structure $\infty$-sheaves of the similar spaces above correspondingly without the $\mathrm{Fro}$-quotients. Here for those space without notation related to the radius and the corresponding interval we consider the total unions $\bigcap_r,\bigcup_I$ in order to achieve the whole spaces to achieve the analogues of the corresponding FF curves from \cite{10KL1}, \cite{10KL2}, \cite{10FF} for
\[
\xymatrix@R+0pc@C+0pc{
\underset{r}{\mathrm{homotopylimit}}~\mathrm{Spec}^\mathrm{CS}\widetilde{\Phi}^r_{*,\Gamma,X_\square},\underset{I}{\mathrm{homotopycolimit}}~\mathrm{Spec}^\mathrm{CS}\widetilde{\Phi}^I_{*,\Gamma,X_\square},	\\
}
\]
\[
\xymatrix@R+0pc@C+0pc{
\underset{r}{\mathrm{homotopylimit}}~\mathrm{Spec}^\mathrm{CS}\breve{\Phi}^r_{*,\Gamma,X_\square},\underset{I}{\mathrm{homotopycolimit}}~\mathrm{Spec}^\mathrm{CS}\breve{\Phi}^I_{*,\Gamma,X_\square},	\\
}
\]
\[
\xymatrix@R+0pc@C+0pc{
\underset{r}{\mathrm{homotopylimit}}~\mathrm{Spec}^\mathrm{CS}{\Phi}^r_{*,\Gamma,X_\square},\underset{I}{\mathrm{homotopycolimit}}~\mathrm{Spec}^\mathrm{CS}{\Phi}^I_{*,\Gamma,X_\square}.	
}
\]
\[ 
\xymatrix@R+0pc@C+0pc{
\underset{r}{\mathrm{homotopylimit}}~\mathrm{Spec}^\mathrm{CS}\widetilde{\Phi}^r_{*,\Gamma,X_\square}/\mathrm{Fro}^\mathbb{Z},\underset{I}{\mathrm{homotopycolimit}}~\mathrm{Spec}^\mathrm{CS}\widetilde{\Phi}^I_{*,\Gamma,X_\square}/\mathrm{Fro}^\mathbb{Z},	\\
}
\]
\[ 
\xymatrix@R+0pc@C+0pc{
\underset{r}{\mathrm{homotopylimit}}~\mathrm{Spec}^\mathrm{CS}\breve{\Phi}^r_{*,\Gamma,X_\square}/\mathrm{Fro}^\mathbb{Z},\underset{I}{\mathrm{homotopycolimit}}~\breve{\Phi}^I_{*,\Gamma,X_\square}/\mathrm{Fro}^\mathbb{Z},	\\
}
\]
\[ 
\xymatrix@R+0pc@C+0pc{
\underset{r}{\mathrm{homotopylimit}}~\mathrm{Spec}^\mathrm{CS}{\Phi}^r_{*,\Gamma,X_\square}/\mathrm{Fro}^\mathbb{Z},\underset{I}{\mathrm{homotopycolimit}}~\mathrm{Spec}^\mathrm{CS}{\Phi}^I_{*,\Gamma,X_\square}/\mathrm{Fro}^\mathbb{Z}.	
}
\]	
In this situation we will have the target category being family parametrized by $r$ or $I$ in compatible glueing sense as in \cite[Definition 5.4.10]{10KL2}. In this situation for modules parametrized by the intervals we have the equivalence of $\infty$-categories by using \cite[Proposition 12.18]{10CS2}. Here the corresponding quasicoherent Frobenius modules are defined to be the corresponding homotopy colimits and limits of Frobenius modules:
\begin{align}
\underset{r}{\mathrm{homotopycolimit}}~M_r,\\
\underset{I}{\mathrm{homotopylimit}}~M_I,	
\end{align}
where each $M_r$ is a Frobenius-equivariant module over the period ring with respect to some radius $r$ while each $M_I$ is a Frobenius-equivariant module over the period ring with respect to some interval $I$.\\
\end{proposition}

\begin{proposition}
Similar proposition holds for 
\begin{align}
\mathrm{Quasicoherentsheaves,Perfectcomplex,IndBanach}_{*}.	
\end{align}	
\end{proposition}

\subsection{Univariate Hodge Iwasawa Prestacks}

\subsubsection{Frobenius Quasicoherent Prestacks I}

\begin{definition}
We now consider the pro-\'etale site of $\mathrm{Spa}\mathbb{Q}_p\left<X_1^{\pm 1},...,X_k^{\pm 1}\right>$, denote that by $*$. To be more accurate we replace one component for $\Gamma$ with the pro-\'etale site of $\mathrm{Spa}\mathbb{Q}_p\left<X_1^{\pm 1},...,X_k^{\pm 1}\right>$. And we treat then all the functor to be prestacks for this site. Then from \cite{10KL1} and \cite[Definition 5.2.1]{10KL2} we have the following class of Kedlaya-Liu rings (with the following replacement: $\Delta$ stands for $A$, $\nabla$ stands for $B$, while $\Phi$ stands for $C$) by taking product in the sense of self $\Gamma$-th power\footnote{Here $|\Gamma|=1$.}:

\[
\xymatrix@R+0pc@C+0pc{
\widetilde{\Delta}_{*},\widetilde{\nabla}_{*},\widetilde{\Phi}_{*},\widetilde{\Delta}^+_{*},\widetilde{\nabla}^+_{*},\widetilde{\Delta}^\dagger_{*},\widetilde{\nabla}^\dagger_{*},\widetilde{\Phi}^r_{*},\widetilde{\Phi}^I_{*}, 
}
\]

\[
\xymatrix@R+0pc@C+0pc{
\breve{\Delta}_{*},\breve{\nabla}_{*},\breve{\Phi}_{*},\breve{\Delta}^+_{*},\breve{\nabla}^+_{*},\breve{\Delta}^\dagger_{*},\breve{\nabla}^\dagger_{*},\breve{\Phi}^r_{*},\breve{\Phi}^I_{*},	
}
\]

\[
\xymatrix@R+0pc@C+0pc{
{\Delta}_{*},{\nabla}_{*},{\Phi}_{*},{\Delta}^+_{*},{\nabla}^+_{*},{\Delta}^\dagger_{*},{\nabla}^\dagger_{*},{\Phi}^r_{*},{\Phi}^I_{*}.	
}
\]
Taking the product we have:
\[
\xymatrix@R+0pc@C+0pc{
\widetilde{\Phi}_{*,X},\widetilde{\Phi}^r_{*,X},\widetilde{\Phi}^I_{*,X},	
}
\]
\[
\xymatrix@R+0pc@C+0pc{
\breve{\Phi}_{*,X},\breve{\Phi}^r_{*,X},\breve{\Phi}^I_{*,X},	
}
\]
\[
\xymatrix@R+0pc@C+0pc{
{\Phi}_{*,X},{\Phi}^r_{*,X},{\Phi}^I_{*,X}.	
}
\]
They carry multi Frobenius action $\varphi_\Gamma$ and multi $\mathrm{Lie}_\Gamma:=\mathbb{Z}_p^{\times\Gamma}$ action. In our current situation after \cite{10CKZ} and \cite{10PZ} we consider the following $(\infty,1)$-categories of $(\infty,1)$-modules.\\
\end{definition}

\begin{definition}
First we consider the Bambozzi-Kremnizer spectrum $\mathrm{Spec}^\mathrm{BK}(*)$ attached to any of those in the above from \cite{10BK} by taking derived rational localization:
\begin{align}
&\mathrm{Spec}^\mathrm{BK}\widetilde{\Phi}_{*,X},\mathrm{Spec}^\mathrm{BK}\widetilde{\Phi}^r_{*,X},\mathrm{Spec}^\mathrm{BK}\widetilde{\Phi}^I_{*,X},	
\end{align}
\begin{align}
&\mathrm{Spec}^\mathrm{BK}\breve{\Phi}_{*,X},\mathrm{Spec}^\mathrm{BK}\breve{\Phi}^r_{*,X},\mathrm{Spec}^\mathrm{BK}\breve{\Phi}^I_{*,X},	
\end{align}
\begin{align}
&\mathrm{Spec}^\mathrm{BK}{\Phi}_{*,X},
\mathrm{Spec}^\mathrm{BK}{\Phi}^r_{*,X},\mathrm{Spec}^\mathrm{BK}{\Phi}^I_{*,X}.	
\end{align}

Then we take the corresponding quotients by using the corresponding Frobenius operators:
\begin{align}
&\mathrm{Spec}^\mathrm{BK}\widetilde{\Phi}_{*,X}/\mathrm{Fro}^\mathbb{Z},	\\
\end{align}
\begin{align}
&\mathrm{Spec}^\mathrm{BK}\breve{\Phi}_{*,X}/\mathrm{Fro}^\mathbb{Z},	\\
\end{align}
\begin{align}
&\mathrm{Spec}^\mathrm{BK}{\Phi}_{*,X}/\mathrm{Fro}^\mathbb{Z}.	
\end{align}
Here for those space without notation related to the radius and the corresponding interval we consider the total unions $\bigcap_r,\bigcup_I$ in order to achieve the whole spaces to achieve the analogues of the corresponding FF curves from \cite{10KL1}, \cite{10KL2}, \cite{10FF} for
\[
\xymatrix@R+0pc@C+0pc{
\underset{r}{\mathrm{homotopylimit}}~\mathrm{Spec}^\mathrm{BK}\widetilde{\Phi}^r_{*,X},\underset{I}{\mathrm{homotopycolimit}}~\mathrm{Spec}^\mathrm{BK}\widetilde{\Phi}^I_{*,X},	\\
}
\]
\[
\xymatrix@R+0pc@C+0pc{
\underset{r}{\mathrm{homotopylimit}}~\mathrm{Spec}^\mathrm{BK}\breve{\Phi}^r_{*,X},\underset{I}{\mathrm{homotopycolimit}}~\mathrm{Spec}^\mathrm{BK}\breve{\Phi}^I_{*,X},	\\
}
\]
\[
\xymatrix@R+0pc@C+0pc{
\underset{r}{\mathrm{homotopylimit}}~\mathrm{Spec}^\mathrm{BK}{\Phi}^r_{*,X},\underset{I}{\mathrm{homotopycolimit}}~\mathrm{Spec}^\mathrm{BK}{\Phi}^I_{*,X}.	
}
\]
\[  
\xymatrix@R+0pc@C+0pc{
\underset{r}{\mathrm{homotopylimit}}~\mathrm{Spec}^\mathrm{BK}\widetilde{\Phi}^r_{*,X}/\mathrm{Fro}^\mathbb{Z},\underset{I}{\mathrm{homotopycolimit}}~\mathrm{Spec}^\mathrm{BK}\widetilde{\Phi}^I_{*,X}/\mathrm{Fro}^\mathbb{Z},	\\
}
\]
\[ 
\xymatrix@R+0pc@C+0pc{
\underset{r}{\mathrm{homotopylimit}}~\mathrm{Spec}^\mathrm{BK}\breve{\Phi}^r_{*,X}/\mathrm{Fro}^\mathbb{Z},\underset{I}{\mathrm{homotopycolimit}}~\mathrm{Spec}^\mathrm{BK}\breve{\Phi}^I_{*,X}/\mathrm{Fro}^\mathbb{Z},	\\
}
\]
\[ 
\xymatrix@R+0pc@C+0pc{
\underset{r}{\mathrm{homotopylimit}}~\mathrm{Spec}^\mathrm{BK}{\Phi}^r_{*,X}/\mathrm{Fro}^\mathbb{Z},\underset{I}{\mathrm{homotopycolimit}}~\mathrm{Spec}^\mathrm{BK}{\Phi}^I_{*,X}/\mathrm{Fro}^\mathbb{Z}.	
}
\]

\end{definition}

\indent Meanwhile we have the corresponding Clausen-Scholze analytic stacks from \cite{10CS2}, therefore applying their construction we have:

\begin{definition}
Here we define the following products by using the solidified tensor product from \cite{10CS1} and \cite{10CS2}. Then we take solidified tensor product $\overset{\blacksquare}{\otimes}$ of any of the following
\[
\xymatrix@R+0pc@C+0pc{
\widetilde{\Delta}_{*},\widetilde{\nabla}_{*},\widetilde{\Phi}_{*},\widetilde{\Delta}^+_{*},\widetilde{\nabla}^+_{*},\widetilde{\Delta}^\dagger_{*},\widetilde{\nabla}^\dagger_{*},\widetilde{\Phi}^r_{*},\widetilde{\Phi}^I_{*}, 
}
\]

\[
\xymatrix@R+0pc@C+0pc{
\breve{\Delta}_{*},\breve{\nabla}_{*},\breve{\Phi}_{*},\breve{\Delta}^+_{*},\breve{\nabla}^+_{*},\breve{\Delta}^\dagger_{*},\breve{\nabla}^\dagger_{*},\breve{\Phi}^r_{*},\breve{\Phi}^I_{*},	
}
\]

\[
\xymatrix@R+0pc@C+0pc{
{\Delta}_{*},{\nabla}_{*},{\Phi}_{*},{\Delta}^+_{*},{\nabla}^+_{*},{\Delta}^\dagger_{*},{\nabla}^\dagger_{*},{\Phi}^r_{*},{\Phi}^I_{*},	
}
\]  	
with $X$. Then we have the notations:
\[
\xymatrix@R+0pc@C+0pc{
\widetilde{\Delta}_{*,X},\widetilde{\nabla}_{*,X},\widetilde{\Phi}_{*,X},\widetilde{\Delta}^+_{*,X},\widetilde{\nabla}^+_{*,X},\widetilde{\Delta}^\dagger_{*,X},\widetilde{\nabla}^\dagger_{*,X},\widetilde{\Phi}^r_{*,X},\widetilde{\Phi}^I_{*,X}, 
}
\]

\[
\xymatrix@R+0pc@C+0pc{
\breve{\Delta}_{*,X},\breve{\nabla}_{*,X},\breve{\Phi}_{*,X},\breve{\Delta}^+_{*,X},\breve{\nabla}^+_{*,X},\breve{\Delta}^\dagger_{*,X},\breve{\nabla}^\dagger_{*,X},\breve{\Phi}^r_{*,X},\breve{\Phi}^I_{*,X},	
}
\]

\[
\xymatrix@R+0pc@C+0pc{
{\Delta}_{*,X},{\nabla}_{*,X},{\Phi}_{*,X},{\Delta}^+_{*,X},{\nabla}^+_{*,X},{\Delta}^\dagger_{*,X},{\nabla}^\dagger_{*,X},{\Phi}^r_{*,X},{\Phi}^I_{*,X}.	
}
\]
\end{definition}

\begin{definition}
First we consider the Clausen-Scholze spectrum $\mathrm{Spec}^\mathrm{CS}(*)$ attached to any of those in the above from \cite{10CS2} by taking derived rational localization:
\begin{align}
\mathrm{Spec}^\mathrm{CS}\widetilde{\Delta}_{*,X},\mathrm{Spec}^\mathrm{CS}\widetilde{\nabla}_{*,X},\mathrm{Spec}^\mathrm{CS}\widetilde{\Phi}_{*,X},\mathrm{Spec}^\mathrm{CS}\widetilde{\Delta}^+_{*,X},\mathrm{Spec}^\mathrm{CS}\widetilde{\nabla}^+_{*,X},\\
\mathrm{Spec}^\mathrm{CS}\widetilde{\Delta}^\dagger_{*,X},\mathrm{Spec}^\mathrm{CS}\widetilde{\nabla}^\dagger_{*,X},\mathrm{Spec}^\mathrm{CS}\widetilde{\Phi}^r_{*,X},\mathrm{Spec}^\mathrm{CS}\widetilde{\Phi}^I_{*,X},	\\
\end{align}
\begin{align}
\mathrm{Spec}^\mathrm{CS}\breve{\Delta}_{*,X},\breve{\nabla}_{*,X},\mathrm{Spec}^\mathrm{CS}\breve{\Phi}_{*,X},\mathrm{Spec}^\mathrm{CS}\breve{\Delta}^+_{*,X},\mathrm{Spec}^\mathrm{CS}\breve{\nabla}^+_{*,X},\\
\mathrm{Spec}^\mathrm{CS}\breve{\Delta}^\dagger_{*,X},\mathrm{Spec}^\mathrm{CS}\breve{\nabla}^\dagger_{*,X},\mathrm{Spec}^\mathrm{CS}\breve{\Phi}^r_{*,X},\breve{\Phi}^I_{*,X},	\\
\end{align}
\begin{align}
\mathrm{Spec}^\mathrm{CS}{\Delta}_{*,X},\mathrm{Spec}^\mathrm{CS}{\nabla}_{*,X},\mathrm{Spec}^\mathrm{CS}{\Phi}_{*,X},\mathrm{Spec}^\mathrm{CS}{\Delta}^+_{*,X},\mathrm{Spec}^\mathrm{CS}{\nabla}^+_{*,X},\\
\mathrm{Spec}^\mathrm{CS}{\Delta}^\dagger_{*,X},\mathrm{Spec}^\mathrm{CS}{\nabla}^\dagger_{*,X},\mathrm{Spec}^\mathrm{CS}{\Phi}^r_{*,X},\mathrm{Spec}^\mathrm{CS}{\Phi}^I_{*,X}.	
\end{align}

Then we take the corresponding quotients by using the corresponding Frobenius operators:
\begin{align}
&\mathrm{Spec}^\mathrm{CS}\widetilde{\Delta}_{*,X}/\mathrm{Fro}^\mathbb{Z},\mathrm{Spec}^\mathrm{CS}\widetilde{\nabla}_{*,X}/\mathrm{Fro}^\mathbb{Z},\mathrm{Spec}^\mathrm{CS}\widetilde{\Phi}_{*,X}/\mathrm{Fro}^\mathbb{Z},\mathrm{Spec}^\mathrm{CS}\widetilde{\Delta}^+_{*,X}/\mathrm{Fro}^\mathbb{Z},\\
&\mathrm{Spec}^\mathrm{CS}\widetilde{\nabla}^+_{*,X}/\mathrm{Fro}^\mathbb{Z}, \mathrm{Spec}^\mathrm{CS}\widetilde{\Delta}^\dagger_{*,X}/\mathrm{Fro}^\mathbb{Z},\mathrm{Spec}^\mathrm{CS}\widetilde{\nabla}^\dagger_{*,X}/\mathrm{Fro}^\mathbb{Z},	\\
\end{align}
\begin{align}
&\mathrm{Spec}^\mathrm{CS}\breve{\Delta}_{*,X}/\mathrm{Fro}^\mathbb{Z},\breve{\nabla}_{*,X}/\mathrm{Fro}^\mathbb{Z},\mathrm{Spec}^\mathrm{CS}\breve{\Phi}_{*,X}/\mathrm{Fro}^\mathbb{Z},\mathrm{Spec}^\mathrm{CS}\breve{\Delta}^+_{*,X}/\mathrm{Fro}^\mathbb{Z},\\
&\mathrm{Spec}^\mathrm{CS}\breve{\nabla}^+_{*,X}/\mathrm{Fro}^\mathbb{Z}, \mathrm{Spec}^\mathrm{CS}\breve{\Delta}^\dagger_{*,X}/\mathrm{Fro}^\mathbb{Z},\mathrm{Spec}^\mathrm{CS}\breve{\nabla}^\dagger_{*,X}/\mathrm{Fro}^\mathbb{Z},	\\
\end{align}
\begin{align}
&\mathrm{Spec}^\mathrm{CS}{\Delta}_{*,X}/\mathrm{Fro}^\mathbb{Z},\mathrm{Spec}^\mathrm{CS}{\nabla}_{*,X}/\mathrm{Fro}^\mathbb{Z},\mathrm{Spec}^\mathrm{CS}{\Phi}_{*,X}/\mathrm{Fro}^\mathbb{Z},\mathrm{Spec}^\mathrm{CS}{\Delta}^+_{*,X}/\mathrm{Fro}^\mathbb{Z},\\
&\mathrm{Spec}^\mathrm{CS}{\nabla}^+_{*,X}/\mathrm{Fro}^\mathbb{Z}, \mathrm{Spec}^\mathrm{CS}{\Delta}^\dagger_{*,X}/\mathrm{Fro}^\mathbb{Z},\mathrm{Spec}^\mathrm{CS}{\nabla}^\dagger_{*,X}/\mathrm{Fro}^\mathbb{Z}.	
\end{align}
Here for those space with notations related to the radius and the corresponding interval we consider the total unions $\bigcap_r,\bigcup_I$ in order to achieve the whole spaces to achieve the analogues of the corresponding FF curves from \cite{10KL1}, \cite{10KL2}, \cite{10FF} for
\[
\xymatrix@R+0pc@C+0pc{
\underset{r}{\mathrm{homotopylimit}}~\mathrm{Spec}^\mathrm{CS}\widetilde{\Phi}^r_{*,X},\underset{I}{\mathrm{homotopycolimit}}~\mathrm{Spec}^\mathrm{CS}\widetilde{\Phi}^I_{*,X},	\\
}
\]
\[
\xymatrix@R+0pc@C+0pc{
\underset{r}{\mathrm{homotopylimit}}~\mathrm{Spec}^\mathrm{CS}\breve{\Phi}^r_{*,X},\underset{I}{\mathrm{homotopycolimit}}~\mathrm{Spec}^\mathrm{CS}\breve{\Phi}^I_{*,X},	\\
}
\]
\[
\xymatrix@R+0pc@C+0pc{
\underset{r}{\mathrm{homotopylimit}}~\mathrm{Spec}^\mathrm{CS}{\Phi}^r_{*,X},\underset{I}{\mathrm{homotopycolimit}}~\mathrm{Spec}^\mathrm{CS}{\Phi}^I_{*,X}.	
}
\]
\[ 
\xymatrix@R+0pc@C+0pc{
\underset{r}{\mathrm{homotopylimit}}~\mathrm{Spec}^\mathrm{CS}\widetilde{\Phi}^r_{*,X}/\mathrm{Fro}^\mathbb{Z},\underset{I}{\mathrm{homotopycolimit}}~\mathrm{Spec}^\mathrm{CS}\widetilde{\Phi}^I_{*,X}/\mathrm{Fro}^\mathbb{Z},	\\
}
\]
\[ 
\xymatrix@R+0pc@C+0pc{
\underset{r}{\mathrm{homotopylimit}}~\mathrm{Spec}^\mathrm{CS}\breve{\Phi}^r_{*,X}/\mathrm{Fro}^\mathbb{Z},\underset{I}{\mathrm{homotopycolimit}}~\breve{\Phi}^I_{*,X}/\mathrm{Fro}^\mathbb{Z},	\\
}
\]
\[ 
\xymatrix@R+0pc@C+0pc{
\underset{r}{\mathrm{homotopylimit}}~\mathrm{Spec}^\mathrm{CS}{\Phi}^r_{*,X}/\mathrm{Fro}^\mathbb{Z},\underset{I}{\mathrm{homotopycolimit}}~\mathrm{Spec}^\mathrm{CS}{\Phi}^I_{*,X}/\mathrm{Fro}^\mathbb{Z}.	
}
\]

\end{definition}

\

\begin{definition}
We then consider the corresponding quasipresheaves of the corresponding ind-Banach or monomorphic ind-Banach modules from \cite{10BBK}, \cite{10KKM}:
\begin{align}
\mathrm{Quasicoherentpresheaves,IndBanach}_{*}	
\end{align}
where $*$ is one of the following spaces:
\begin{align}
&\mathrm{Spec}^\mathrm{BK}\widetilde{\Phi}_{*,X}/\mathrm{Fro}^\mathbb{Z},	\\
\end{align}
\begin{align}
&\mathrm{Spec}^\mathrm{BK}\breve{\Phi}_{*,X}/\mathrm{Fro}^\mathbb{Z},	\\
\end{align}
\begin{align}
&\mathrm{Spec}^\mathrm{BK}{\Phi}_{*,X}/\mathrm{Fro}^\mathbb{Z}.	
\end{align}
Here for those space without notation related to the radius and the corresponding interval we consider the total unions $\bigcap_r,\bigcup_I$ in order to achieve the whole spaces to achieve the analogues of the corresponding FF curves from \cite{10KL1}, \cite{10KL2}, \cite{10FF} for
\[
\xymatrix@R+0pc@C+0pc{
\underset{r}{\mathrm{homotopylimit}}~\mathrm{Spec}^\mathrm{BK}\widetilde{\Phi}^r_{*,X},\underset{I}{\mathrm{homotopycolimit}}~\mathrm{Spec}^\mathrm{BK}\widetilde{\Phi}^I_{*,X},	\\
}
\]
\[
\xymatrix@R+0pc@C+0pc{
\underset{r}{\mathrm{homotopylimit}}~\mathrm{Spec}^\mathrm{BK}\breve{\Phi}^r_{*,X},\underset{I}{\mathrm{homotopycolimit}}~\mathrm{Spec}^\mathrm{BK}\breve{\Phi}^I_{*,X},	\\
}
\]
\[
\xymatrix@R+0pc@C+0pc{
\underset{r}{\mathrm{homotopylimit}}~\mathrm{Spec}^\mathrm{BK}{\Phi}^r_{*,X},\underset{I}{\mathrm{homotopycolimit}}~\mathrm{Spec}^\mathrm{BK}{\Phi}^I_{*,X}.	
}
\]
\[  
\xymatrix@R+0pc@C+0pc{
\underset{r}{\mathrm{homotopylimit}}~\mathrm{Spec}^\mathrm{BK}\widetilde{\Phi}^r_{*,X}/\mathrm{Fro}^\mathbb{Z},\underset{I}{\mathrm{homotopycolimit}}~\mathrm{Spec}^\mathrm{BK}\widetilde{\Phi}^I_{*,X}/\mathrm{Fro}^\mathbb{Z},	\\
}
\]
\[ 
\xymatrix@R+0pc@C+0pc{
\underset{r}{\mathrm{homotopylimit}}~\mathrm{Spec}^\mathrm{BK}\breve{\Phi}^r_{*,X}/\mathrm{Fro}^\mathbb{Z},\underset{I}{\mathrm{homotopycolimit}}~\mathrm{Spec}^\mathrm{BK}\breve{\Phi}^I_{*,X}/\mathrm{Fro}^\mathbb{Z},	\\
}
\]
\[ 
\xymatrix@R+0pc@C+0pc{
\underset{r}{\mathrm{homotopylimit}}~\mathrm{Spec}^\mathrm{BK}{\Phi}^r_{*,X}/\mathrm{Fro}^\mathbb{Z},\underset{I}{\mathrm{homotopycolimit}}~\mathrm{Spec}^\mathrm{BK}{\Phi}^I_{*,X}/\mathrm{Fro}^\mathbb{Z}.	
}
\]

\end{definition}

\begin{definition}
We then consider the corresponding quasisheaves of the corresponding condensed solid topological modules from \cite{10CS2}:
\begin{align}
\mathrm{Quasicoherentsheaves, Condensed}_{*}	
\end{align}
where $*$ is one of the following spaces:
\begin{align}
&\mathrm{Spec}^\mathrm{CS}\widetilde{\Delta}_{*,X}/\mathrm{Fro}^\mathbb{Z},\mathrm{Spec}^\mathrm{CS}\widetilde{\nabla}_{*,X}/\mathrm{Fro}^\mathbb{Z},\mathrm{Spec}^\mathrm{CS}\widetilde{\Phi}_{*,X}/\mathrm{Fro}^\mathbb{Z},\mathrm{Spec}^\mathrm{CS}\widetilde{\Delta}^+_{*,X}/\mathrm{Fro}^\mathbb{Z},\\
&\mathrm{Spec}^\mathrm{CS}\widetilde{\nabla}^+_{*,X}/\mathrm{Fro}^\mathbb{Z},\mathrm{Spec}^\mathrm{CS}\widetilde{\Delta}^\dagger_{*,X}/\mathrm{Fro}^\mathbb{Z},\mathrm{Spec}^\mathrm{CS}\widetilde{\nabla}^\dagger_{*,X}/\mathrm{Fro}^\mathbb{Z},	\\
\end{align}
\begin{align}
&\mathrm{Spec}^\mathrm{CS}\breve{\Delta}_{*,X}/\mathrm{Fro}^\mathbb{Z},\breve{\nabla}_{*,X}/\mathrm{Fro}^\mathbb{Z},\mathrm{Spec}^\mathrm{CS}\breve{\Phi}_{*,X}/\mathrm{Fro}^\mathbb{Z},\mathrm{Spec}^\mathrm{CS}\breve{\Delta}^+_{*,X}/\mathrm{Fro}^\mathbb{Z},\\
&\mathrm{Spec}^\mathrm{CS}\breve{\nabla}^+_{*,X}/\mathrm{Fro}^\mathbb{Z},\mathrm{Spec}^\mathrm{CS}\breve{\Delta}^\dagger_{*,X}/\mathrm{Fro}^\mathbb{Z},\mathrm{Spec}^\mathrm{CS}\breve{\nabla}^\dagger_{*,X}/\mathrm{Fro}^\mathbb{Z},	\\
\end{align}
\begin{align}
&\mathrm{Spec}^\mathrm{CS}{\Delta}_{*,X}/\mathrm{Fro}^\mathbb{Z},\mathrm{Spec}^\mathrm{CS}{\nabla}_{*,X}/\mathrm{Fro}^\mathbb{Z},\mathrm{Spec}^\mathrm{CS}{\Phi}_{*,X}/\mathrm{Fro}^\mathbb{Z},\mathrm{Spec}^\mathrm{CS}{\Delta}^+_{*,X}/\mathrm{Fro}^\mathbb{Z},\\
&\mathrm{Spec}^\mathrm{CS}{\nabla}^+_{*,X}/\mathrm{Fro}^\mathbb{Z}, \mathrm{Spec}^\mathrm{CS}{\Delta}^\dagger_{*,X}/\mathrm{Fro}^\mathbb{Z},\mathrm{Spec}^\mathrm{CS}{\nabla}^\dagger_{*,X}/\mathrm{Fro}^\mathbb{Z}.	
\end{align}
Here for those space with notations related to the radius and the corresponding interval we consider the total unions $\bigcap_r,\bigcup_I$ in order to achieve the whole spaces to achieve the analogues of the corresponding FF curves from \cite{10KL1}, \cite{10KL2}, \cite{10FF} for
\[
\xymatrix@R+0pc@C+0pc{
\underset{r}{\mathrm{homotopylimit}}~\mathrm{Spec}^\mathrm{CS}\widetilde{\Phi}^r_{*,X},\underset{I}{\mathrm{homotopycolimit}}~\mathrm{Spec}^\mathrm{CS}\widetilde{\Phi}^I_{*,X},	\\
}
\]
\[
\xymatrix@R+0pc@C+0pc{
\underset{r}{\mathrm{homotopylimit}}~\mathrm{Spec}^\mathrm{CS}\breve{\Phi}^r_{*,X},\underset{I}{\mathrm{homotopycolimit}}~\mathrm{Spec}^\mathrm{CS}\breve{\Phi}^I_{*,X},	\\
}
\]
\[
\xymatrix@R+0pc@C+0pc{
\underset{r}{\mathrm{homotopylimit}}~\mathrm{Spec}^\mathrm{CS}{\Phi}^r_{*,X},\underset{I}{\mathrm{homotopycolimit}}~\mathrm{Spec}^\mathrm{CS}{\Phi}^I_{*,X}.	
}
\]
\[ 
\xymatrix@R+0pc@C+0pc{
\underset{r}{\mathrm{homotopylimit}}~\mathrm{Spec}^\mathrm{CS}\widetilde{\Phi}^r_{*,X}/\mathrm{Fro}^\mathbb{Z},\underset{I}{\mathrm{homotopycolimit}}~\mathrm{Spec}^\mathrm{CS}\widetilde{\Phi}^I_{*,X}/\mathrm{Fro}^\mathbb{Z},	\\
}
\]
\[ 
\xymatrix@R+0pc@C+0pc{
\underset{r}{\mathrm{homotopylimit}}~\mathrm{Spec}^\mathrm{CS}\breve{\Phi}^r_{*,X}/\mathrm{Fro}^\mathbb{Z},\underset{I}{\mathrm{homotopycolimit}}~\breve{\Phi}^I_{*,X}/\mathrm{Fro}^\mathbb{Z},	\\
}
\]
\[ 
\xymatrix@R+0pc@C+0pc{
\underset{r}{\mathrm{homotopylimit}}~\mathrm{Spec}^\mathrm{CS}{\Phi}^r_{*,X}/\mathrm{Fro}^\mathbb{Z},\underset{I}{\mathrm{homotopycolimit}}~\mathrm{Spec}^\mathrm{CS}{\Phi}^I_{*,X}/\mathrm{Fro}^\mathbb{Z}.	
}
\]

\end{definition}

\

\begin{proposition}
There is a well-defined functor from the $\infty$-category 
\begin{align}
\mathrm{Quasicoherentpresheaves,Condensed}_{*}	
\end{align}
where $*$ is one of the following spaces:
\begin{align}
&\mathrm{Spec}^\mathrm{CS}\widetilde{\Phi}_{*,X}/\mathrm{Fro}^\mathbb{Z},	\\
\end{align}
\begin{align}
&\mathrm{Spec}^\mathrm{CS}\breve{\Phi}_{*,X}/\mathrm{Fro}^\mathbb{Z},	\\
\end{align}
\begin{align}
&\mathrm{Spec}^\mathrm{CS}{\Phi}_{*,X}/\mathrm{Fro}^\mathbb{Z},	
\end{align}
to the $\infty$-category of $\mathrm{Fro}$-equivariant quasicoherent presheaves over similar spaces above correspondingly without the $\mathrm{Fro}$-quotients, and to the $\infty$-category of $\mathrm{Fro}$-equivariant quasicoherent modules over global sections of the structure $\infty$-sheaves of the similar spaces above correspondingly without the $\mathrm{Fro}$-quotients. Here for those space without notation related to the radius and the corresponding interval we consider the total unions $\bigcap_r,\bigcup_I$ in order to achieve the whole spaces to achieve the analogues of the corresponding FF curves from \cite{10KL1}, \cite{10KL2}, \cite{10FF} for
\[
\xymatrix@R+0pc@C+0pc{
\underset{r}{\mathrm{homotopylimit}}~\mathrm{Spec}^\mathrm{CS}\widetilde{\Phi}^r_{*,X},\underset{I}{\mathrm{homotopycolimit}}~\mathrm{Spec}^\mathrm{CS}\widetilde{\Phi}^I_{*,X},	\\
}
\]
\[
\xymatrix@R+0pc@C+0pc{
\underset{r}{\mathrm{homotopylimit}}~\mathrm{Spec}^\mathrm{CS}\breve{\Phi}^r_{*,X},\underset{I}{\mathrm{homotopycolimit}}~\mathrm{Spec}^\mathrm{CS}\breve{\Phi}^I_{*,X},	\\
}
\]
\[
\xymatrix@R+0pc@C+0pc{
\underset{r}{\mathrm{homotopylimit}}~\mathrm{Spec}^\mathrm{CS}{\Phi}^r_{*,X},\underset{I}{\mathrm{homotopycolimit}}~\mathrm{Spec}^\mathrm{CS}{\Phi}^I_{*,X}.	
}
\]
\[ 
\xymatrix@R+0pc@C+0pc{
\underset{r}{\mathrm{homotopylimit}}~\mathrm{Spec}^\mathrm{CS}\widetilde{\Phi}^r_{*,X}/\mathrm{Fro}^\mathbb{Z},\underset{I}{\mathrm{homotopycolimit}}~\mathrm{Spec}^\mathrm{CS}\widetilde{\Phi}^I_{*,X}/\mathrm{Fro}^\mathbb{Z},	\\
}
\]
\[ 
\xymatrix@R+0pc@C+0pc{
\underset{r}{\mathrm{homotopylimit}}~\mathrm{Spec}^\mathrm{CS}\breve{\Phi}^r_{*,X}/\mathrm{Fro}^\mathbb{Z},\underset{I}{\mathrm{homotopycolimit}}~\breve{\Phi}^I_{*,X}/\mathrm{Fro}^\mathbb{Z},	\\
}
\]
\[ 
\xymatrix@R+0pc@C+0pc{
\underset{r}{\mathrm{homotopylimit}}~\mathrm{Spec}^\mathrm{CS}{\Phi}^r_{*,X}/\mathrm{Fro}^\mathbb{Z},\underset{I}{\mathrm{homotopycolimit}}~\mathrm{Spec}^\mathrm{CS}{\Phi}^I_{*,X}/\mathrm{Fro}^\mathbb{Z}.	
}
\]	
In this situation we will have the target category being family parametrized by $r$ or $I$ in compatible glueing sense as in \cite[Definition 5.4.10]{10KL2}. In this situation for modules parametrized by the intervals we have the equivalence of $\infty$-categories by using \cite[Proposition 13.8]{10CS2}. Here the corresponding quasicoherent Frobenius modules are defined to be the corresponding homotopy colimits and limits of Frobenius modules:
\begin{align}
\underset{r}{\mathrm{homotopycolimit}}~M_r,\\
\underset{I}{\mathrm{homotopylimit}}~M_I,	
\end{align}
where each $M_r$ is a Frobenius-equivariant module over the period ring with respect to some radius $r$ while each $M_I$ is a Frobenius-equivariant module over the period ring with respect to some interval $I$.\\
\end{proposition}

\begin{proposition}
Similar proposition holds for 
\begin{align}
\mathrm{Quasicoherentsheaves,IndBanach}_{*}.	
\end{align}	
\end{proposition}

\

\begin{definition}
We then consider the corresponding quasipresheaves of perfect complexes the corresponding ind-Banach or monomorphic ind-Banach modules from \cite{10BBK}, \cite{10KKM}:
\begin{align}
\mathrm{Quasicoherentpresheaves,Perfectcomplex,IndBanach}_{*}	
\end{align}
where $*$ is one of the following spaces:
\begin{align}
&\mathrm{Spec}^\mathrm{BK}\widetilde{\Phi}_{*,X}/\mathrm{Fro}^\mathbb{Z},	\\
\end{align}
\begin{align}
&\mathrm{Spec}^\mathrm{BK}\breve{\Phi}_{*,X}/\mathrm{Fro}^\mathbb{Z},	\\
\end{align}
\begin{align}
&\mathrm{Spec}^\mathrm{BK}{\Phi}_{*,X}/\mathrm{Fro}^\mathbb{Z}.	
\end{align}
Here for those space without notation related to the radius and the corresponding interval we consider the total unions $\bigcap_r,\bigcup_I$ in order to achieve the whole spaces to achieve the analogues of the corresponding FF curves from \cite{10KL1}, \cite{10KL2}, \cite{10FF} for
\[
\xymatrix@R+0pc@C+0pc{
\underset{r}{\mathrm{homotopylimit}}~\mathrm{Spec}^\mathrm{BK}\widetilde{\Phi}^r_{*,X},\underset{I}{\mathrm{homotopycolimit}}~\mathrm{Spec}^\mathrm{BK}\widetilde{\Phi}^I_{*,X},	\\
}
\]
\[
\xymatrix@R+0pc@C+0pc{
\underset{r}{\mathrm{homotopylimit}}~\mathrm{Spec}^\mathrm{BK}\breve{\Phi}^r_{*,X},\underset{I}{\mathrm{homotopycolimit}}~\mathrm{Spec}^\mathrm{BK}\breve{\Phi}^I_{*,X},	\\
}
\]
\[
\xymatrix@R+0pc@C+0pc{
\underset{r}{\mathrm{homotopylimit}}~\mathrm{Spec}^\mathrm{BK}{\Phi}^r_{*,X},\underset{I}{\mathrm{homotopycolimit}}~\mathrm{Spec}^\mathrm{BK}{\Phi}^I_{*,X}.	
}
\]
\[  
\xymatrix@R+0pc@C+0pc{
\underset{r}{\mathrm{homotopylimit}}~\mathrm{Spec}^\mathrm{BK}\widetilde{\Phi}^r_{*,X}/\mathrm{Fro}^\mathbb{Z},\underset{I}{\mathrm{homotopycolimit}}~\mathrm{Spec}^\mathrm{BK}\widetilde{\Phi}^I_{*,X}/\mathrm{Fro}^\mathbb{Z},	\\
}
\]
\[ 
\xymatrix@R+0pc@C+0pc{
\underset{r}{\mathrm{homotopylimit}}~\mathrm{Spec}^\mathrm{BK}\breve{\Phi}^r_{*,X}/\mathrm{Fro}^\mathbb{Z},\underset{I}{\mathrm{homotopycolimit}}~\mathrm{Spec}^\mathrm{BK}\breve{\Phi}^I_{*,X}/\mathrm{Fro}^\mathbb{Z},	\\
}
\]
\[ 
\xymatrix@R+0pc@C+0pc{
\underset{r}{\mathrm{homotopylimit}}~\mathrm{Spec}^\mathrm{BK}{\Phi}^r_{*,X}/\mathrm{Fro}^\mathbb{Z},\underset{I}{\mathrm{homotopycolimit}}~\mathrm{Spec}^\mathrm{BK}{\Phi}^I_{*,X}/\mathrm{Fro}^\mathbb{Z}.	
}
\]

\end{definition}

\begin{definition}
We then consider the corresponding quasisheaves of perfect complexes of the corresponding condensed solid topological modules from \cite{10CS2}:
\begin{align}
\mathrm{Quasicoherentsheaves, Perfectcomplex, Condensed}_{*}	
\end{align}
where $*$ is one of the following spaces:
\begin{align}
&\mathrm{Spec}^\mathrm{CS}\widetilde{\Delta}_{*,X}/\mathrm{Fro}^\mathbb{Z},\mathrm{Spec}^\mathrm{CS}\widetilde{\nabla}_{*,X}/\mathrm{Fro}^\mathbb{Z},\mathrm{Spec}^\mathrm{CS}\widetilde{\Phi}_{*,X}/\mathrm{Fro}^\mathbb{Z},\mathrm{Spec}^\mathrm{CS}\widetilde{\Delta}^+_{*,X}/\mathrm{Fro}^\mathbb{Z},\\
&\mathrm{Spec}^\mathrm{CS}\widetilde{\nabla}^+_{*,X}/\mathrm{Fro}^\mathbb{Z},\mathrm{Spec}^\mathrm{CS}\widetilde{\Delta}^\dagger_{*,X}/\mathrm{Fro}^\mathbb{Z},\mathrm{Spec}^\mathrm{CS}\widetilde{\nabla}^\dagger_{*,X}/\mathrm{Fro}^\mathbb{Z},	\\
\end{align}
\begin{align}
&\mathrm{Spec}^\mathrm{CS}\breve{\Delta}_{*,X}/\mathrm{Fro}^\mathbb{Z},\breve{\nabla}_{*,X}/\mathrm{Fro}^\mathbb{Z},\mathrm{Spec}^\mathrm{CS}\breve{\Phi}_{*,X}/\mathrm{Fro}^\mathbb{Z},\mathrm{Spec}^\mathrm{CS}\breve{\Delta}^+_{*,X}/\mathrm{Fro}^\mathbb{Z},\\
&\mathrm{Spec}^\mathrm{CS}\breve{\nabla}^+_{*,X}/\mathrm{Fro}^\mathbb{Z},\mathrm{Spec}^\mathrm{CS}\breve{\Delta}^\dagger_{*,X}/\mathrm{Fro}^\mathbb{Z},\mathrm{Spec}^\mathrm{CS}\breve{\nabla}^\dagger_{*,X}/\mathrm{Fro}^\mathbb{Z},	\\
\end{align}
\begin{align}
&\mathrm{Spec}^\mathrm{CS}{\Delta}_{*,X}/\mathrm{Fro}^\mathbb{Z},\mathrm{Spec}^\mathrm{CS}{\nabla}_{*,X}/\mathrm{Fro}^\mathbb{Z},\mathrm{Spec}^\mathrm{CS}{\Phi}_{*,X}/\mathrm{Fro}^\mathbb{Z},\mathrm{Spec}^\mathrm{CS}{\Delta}^+_{*,X}/\mathrm{Fro}^\mathbb{Z},\\
&\mathrm{Spec}^\mathrm{CS}{\nabla}^+_{*,X}/\mathrm{Fro}^\mathbb{Z}, \mathrm{Spec}^\mathrm{CS}{\Delta}^\dagger_{*,X}/\mathrm{Fro}^\mathbb{Z},\mathrm{Spec}^\mathrm{CS}{\nabla}^\dagger_{*,X}/\mathrm{Fro}^\mathbb{Z}.	
\end{align}
Here for those space with notations related to the radius and the corresponding interval we consider the total unions $\bigcap_r,\bigcup_I$ in order to achieve the whole spaces to achieve the analogues of the corresponding FF curves from \cite{10KL1}, \cite{10KL2}, \cite{10FF} for
\[
\xymatrix@R+0pc@C+0pc{
\underset{r}{\mathrm{homotopylimit}}~\mathrm{Spec}^\mathrm{CS}\widetilde{\Phi}^r_{*,X},\underset{I}{\mathrm{homotopycolimit}}~\mathrm{Spec}^\mathrm{CS}\widetilde{\Phi}^I_{*,X},	\\
}
\]
\[
\xymatrix@R+0pc@C+0pc{
\underset{r}{\mathrm{homotopylimit}}~\mathrm{Spec}^\mathrm{CS}\breve{\Phi}^r_{*,X},\underset{I}{\mathrm{homotopycolimit}}~\mathrm{Spec}^\mathrm{CS}\breve{\Phi}^I_{*,X},	\\
}
\]
\[
\xymatrix@R+0pc@C+0pc{
\underset{r}{\mathrm{homotopylimit}}~\mathrm{Spec}^\mathrm{CS}{\Phi}^r_{*,X},\underset{I}{\mathrm{homotopycolimit}}~\mathrm{Spec}^\mathrm{CS}{\Phi}^I_{*,X}.	
}
\]
\[ 
\xymatrix@R+0pc@C+0pc{
\underset{r}{\mathrm{homotopylimit}}~\mathrm{Spec}^\mathrm{CS}\widetilde{\Phi}^r_{*,X}/\mathrm{Fro}^\mathbb{Z},\underset{I}{\mathrm{homotopycolimit}}~\mathrm{Spec}^\mathrm{CS}\widetilde{\Phi}^I_{*,X}/\mathrm{Fro}^\mathbb{Z},	\\
}
\]
\[ 
\xymatrix@R+0pc@C+0pc{
\underset{r}{\mathrm{homotopylimit}}~\mathrm{Spec}^\mathrm{CS}\breve{\Phi}^r_{*,X}/\mathrm{Fro}^\mathbb{Z},\underset{I}{\mathrm{homotopycolimit}}~\breve{\Phi}^I_{*,X}/\mathrm{Fro}^\mathbb{Z},	\\
}
\]
\[ 
\xymatrix@R+0pc@C+0pc{
\underset{r}{\mathrm{homotopylimit}}~\mathrm{Spec}^\mathrm{CS}{\Phi}^r_{*,X}/\mathrm{Fro}^\mathbb{Z},\underset{I}{\mathrm{homotopycolimit}}~\mathrm{Spec}^\mathrm{CS}{\Phi}^I_{*,X}/\mathrm{Fro}^\mathbb{Z}.	
}
\]

\end{definition}

\begin{proposition}
There is a well-defined functor from the $\infty$-category 
\begin{align}
\mathrm{Quasicoherentpresheaves,Perfectcomplex,Condensed}_{*}	
\end{align}
where $*$ is one of the following spaces:
\begin{align}
&\mathrm{Spec}^\mathrm{CS}\widetilde{\Phi}_{*,X}/\mathrm{Fro}^\mathbb{Z},	\\
\end{align}
\begin{align}
&\mathrm{Spec}^\mathrm{CS}\breve{\Phi}_{*,X}/\mathrm{Fro}^\mathbb{Z},	\\
\end{align}
\begin{align}
&\mathrm{Spec}^\mathrm{CS}{\Phi}_{*,X}/\mathrm{Fro}^\mathbb{Z},	
\end{align}
to the $\infty$-category of $\mathrm{Fro}$-equivariant quasicoherent presheaves over similar spaces above correspondingly without the $\mathrm{Fro}$-quotients, and to the $\infty$-category of $\mathrm{Fro}$-equivariant quasicoherent modules over global sections of the structure $\infty$-sheaves of the similar spaces above correspondingly without the $\mathrm{Fro}$-quotients. Here for those space without notation related to the radius and the corresponding interval we consider the total unions $\bigcap_r,\bigcup_I$ in order to achieve the whole spaces to achieve the analogues of the corresponding FF curves from \cite{10KL1}, \cite{10KL2}, \cite{10FF} for
\[
\xymatrix@R+0pc@C+0pc{
\underset{r}{\mathrm{homotopylimit}}~\mathrm{Spec}^\mathrm{CS}\widetilde{\Phi}^r_{*,X},\underset{I}{\mathrm{homotopycolimit}}~\mathrm{Spec}^\mathrm{CS}\widetilde{\Phi}^I_{*,X},	\\
}
\]
\[
\xymatrix@R+0pc@C+0pc{
\underset{r}{\mathrm{homotopylimit}}~\mathrm{Spec}^\mathrm{CS}\breve{\Phi}^r_{*,X},\underset{I}{\mathrm{homotopycolimit}}~\mathrm{Spec}^\mathrm{CS}\breve{\Phi}^I_{*,X},	\\
}
\]
\[
\xymatrix@R+0pc@C+0pc{
\underset{r}{\mathrm{homotopylimit}}~\mathrm{Spec}^\mathrm{CS}{\Phi}^r_{*,X},\underset{I}{\mathrm{homotopycolimit}}~\mathrm{Spec}^\mathrm{CS}{\Phi}^I_{*,X}.	
}
\]
\[ 
\xymatrix@R+0pc@C+0pc{
\underset{r}{\mathrm{homotopylimit}}~\mathrm{Spec}^\mathrm{CS}\widetilde{\Phi}^r_{*,X}/\mathrm{Fro}^\mathbb{Z},\underset{I}{\mathrm{homotopycolimit}}~\mathrm{Spec}^\mathrm{CS}\widetilde{\Phi}^I_{*,X}/\mathrm{Fro}^\mathbb{Z},	\\
}
\]
\[ 
\xymatrix@R+0pc@C+0pc{
\underset{r}{\mathrm{homotopylimit}}~\mathrm{Spec}^\mathrm{CS}\breve{\Phi}^r_{*,X}/\mathrm{Fro}^\mathbb{Z},\underset{I}{\mathrm{homotopycolimit}}~\breve{\Phi}^I_{*,X}/\mathrm{Fro}^\mathbb{Z},	\\
}
\]
\[ 
\xymatrix@R+0pc@C+0pc{
\underset{r}{\mathrm{homotopylimit}}~\mathrm{Spec}^\mathrm{CS}{\Phi}^r_{*,X}/\mathrm{Fro}^\mathbb{Z},\underset{I}{\mathrm{homotopycolimit}}~\mathrm{Spec}^\mathrm{CS}{\Phi}^I_{*,X}/\mathrm{Fro}^\mathbb{Z}.	
}
\]	
In this situation we will have the target category being family parametrized by $r$ or $I$ in compatible glueing sense as in \cite[Definition 5.4.10]{10KL2}. In this situation for modules parametrized by the intervals we have the equivalence of $\infty$-categories by using \cite[Proposition 12.18]{10CS2}. Here the corresponding quasicoherent Frobenius modules are defined to be the corresponding homotopy colimits and limits of Frobenius modules:
\begin{align}
\underset{r}{\mathrm{homotopycolimit}}~M_r,\\
\underset{I}{\mathrm{homotopylimit}}~M_I,	
\end{align}
where each $M_r$ is a Frobenius-equivariant module over the period ring with respect to some radius $r$ while each $M_I$ is a Frobenius-equivariant module over the period ring with respect to some interval $I$.\\
\end{proposition}

\begin{proposition}
Similar proposition holds for 
\begin{align}
\mathrm{Quasicoherentsheaves,Perfectcomplex,IndBanach}_{*}.	
\end{align}	
\end{proposition}

\subsubsection{Frobenius Quasicoherent Prestacks II: Deformation in Preadic Spaces}

\begin{definition}
We now consider the pro-\'etale site of $\mathrm{Spa}\mathbb{Q}_p\left<X_1^{\pm 1},...,X_k^{\pm 1}\right>$, denote that by $*$. To be more accurate we replace one component for $\Gamma$ with the pro-\'etale site of $\mathrm{Spa}\mathbb{Q}_p\left<X_1^{\pm 1},...,X_k^{\pm 1}\right>$. And we treat then all the functor to be prestacks for this site. Then from \cite{10KL1} and \cite[Definition 5.2.1]{10KL2} we have the following class of Kedlaya-Liu rings (with the following replacement: $\Delta$ stands for $A$, $\nabla$ stands for $B$, while $\Phi$ stands for $C$) by taking product in the sense of self $\Gamma$-th power\footnote{Here $|\Gamma|=1$.}:

\[
\xymatrix@R+0pc@C+0pc{
\widetilde{\Delta}_{*},\widetilde{\nabla}_{*},\widetilde{\Phi}_{*},\widetilde{\Delta}^+_{*},\widetilde{\nabla}^+_{*},\widetilde{\Delta}^\dagger_{*},\widetilde{\nabla}^\dagger_{*},\widetilde{\Phi}^r_{*},\widetilde{\Phi}^I_{*}, 
}
\]

\[
\xymatrix@R+0pc@C+0pc{
\breve{\Delta}_{*},\breve{\nabla}_{*},\breve{\Phi}_{*},\breve{\Delta}^+_{*},\breve{\nabla}^+_{*},\breve{\Delta}^\dagger_{*},\breve{\nabla}^\dagger_{*},\breve{\Phi}^r_{*},\breve{\Phi}^I_{*},	
}
\]

\[
\xymatrix@R+0pc@C+0pc{
{\Delta}_{*},{\nabla}_{*},{\Phi}_{*},{\Delta}^+_{*},{\nabla}^+_{*},{\Delta}^\dagger_{*},{\nabla}^\dagger_{*},{\Phi}^r_{*},{\Phi}^I_{*}.	
}
\]
Taking the product we have:
\[
\xymatrix@R+0pc@C+0pc{
\widetilde{\Phi}_{*,\circ},\widetilde{\Phi}^r_{*,\circ},\widetilde{\Phi}^I_{*,\circ},	
}
\]
\[
\xymatrix@R+0pc@C+0pc{
\breve{\Phi}_{*,\circ},\breve{\Phi}^r_{*,\circ},\breve{\Phi}^I_{*,\circ},	
}
\]
\[
\xymatrix@R+0pc@C+0pc{
{\Phi}_{*,\circ},{\Phi}^r_{*,\circ},{\Phi}^I_{*,\circ}.	
}
\]
They carry multi Frobenius action $\varphi_\Gamma$ and multi $\mathrm{Lie}_\Gamma:=\mathbb{Z}_p^{\times\Gamma}$ action. In our current situation after \cite{10CKZ} and \cite{10PZ} we consider the following $(\infty,1)$-categories of $(\infty,1)$-modules.\\
\end{definition}

\begin{definition}
First we consider the Bambozzi-Kremnizer spectrum $\mathrm{Spec}^\mathrm{BK}(*)$ attached to any of those in the above from \cite{10BK} by taking derived rational localization:
\begin{align}
&\mathrm{Spec}^\mathrm{BK}\widetilde{\Phi}_{*,\circ},\mathrm{Spec}^\mathrm{BK}\widetilde{\Phi}^r_{*,\circ},\mathrm{Spec}^\mathrm{BK}\widetilde{\Phi}^I_{*,\circ},	
\end{align}
\begin{align}
&\mathrm{Spec}^\mathrm{BK}\breve{\Phi}_{*,\circ},\mathrm{Spec}^\mathrm{BK}\breve{\Phi}^r_{*,\circ},\mathrm{Spec}^\mathrm{BK}\breve{\Phi}^I_{*,\circ},	
\end{align}
\begin{align}
&\mathrm{Spec}^\mathrm{BK}{\Phi}_{*,\circ},
\mathrm{Spec}^\mathrm{BK}{\Phi}^r_{*,\circ},\mathrm{Spec}^\mathrm{BK}{\Phi}^I_{*,\circ}.	
\end{align}

Then we take the corresponding quotients by using the corresponding Frobenius operators:
\begin{align}
&\mathrm{Spec}^\mathrm{BK}\widetilde{\Phi}_{*,\circ}/\mathrm{Fro}^\mathbb{Z},	\\
\end{align}
\begin{align}
&\mathrm{Spec}^\mathrm{BK}\breve{\Phi}_{*,\circ}/\mathrm{Fro}^\mathbb{Z},	\\
\end{align}
\begin{align}
&\mathrm{Spec}^\mathrm{BK}{\Phi}_{*,\circ}/\mathrm{Fro}^\mathbb{Z}.	
\end{align}
Here for those space without notation related to the radius and the corresponding interval we consider the total unions $\bigcap_r,\bigcup_I$ in order to achieve the whole spaces to achieve the analogues of the corresponding FF curves from \cite{10KL1}, \cite{10KL2}, \cite{10FF} for
\[
\xymatrix@R+0pc@C+0pc{
\underset{r}{\mathrm{homotopylimit}}~\mathrm{Spec}^\mathrm{BK}\widetilde{\Phi}^r_{*,\circ},\underset{I}{\mathrm{homotopycolimit}}~\mathrm{Spec}^\mathrm{BK}\widetilde{\Phi}^I_{*,\circ},	\\
}
\]
\[
\xymatrix@R+0pc@C+0pc{
\underset{r}{\mathrm{homotopylimit}}~\mathrm{Spec}^\mathrm{BK}\breve{\Phi}^r_{*,\circ},\underset{I}{\mathrm{homotopycolimit}}~\mathrm{Spec}^\mathrm{BK}\breve{\Phi}^I_{*,\circ},	\\
}
\]
\[
\xymatrix@R+0pc@C+0pc{
\underset{r}{\mathrm{homotopylimit}}~\mathrm{Spec}^\mathrm{BK}{\Phi}^r_{*,\circ},\underset{I}{\mathrm{homotopycolimit}}~\mathrm{Spec}^\mathrm{BK}{\Phi}^I_{*,\circ}.	
}
\]
\[  
\xymatrix@R+0pc@C+0pc{
\underset{r}{\mathrm{homotopylimit}}~\mathrm{Spec}^\mathrm{BK}\widetilde{\Phi}^r_{*,\circ}/\mathrm{Fro}^\mathbb{Z},\underset{I}{\mathrm{homotopycolimit}}~\mathrm{Spec}^\mathrm{BK}\widetilde{\Phi}^I_{*,\circ}/\mathrm{Fro}^\mathbb{Z},	\\
}
\]
\[ 
\xymatrix@R+0pc@C+0pc{
\underset{r}{\mathrm{homotopylimit}}~\mathrm{Spec}^\mathrm{BK}\breve{\Phi}^r_{*,\circ}/\mathrm{Fro}^\mathbb{Z},\underset{I}{\mathrm{homotopycolimit}}~\mathrm{Spec}^\mathrm{BK}\breve{\Phi}^I_{*,\circ}/\mathrm{Fro}^\mathbb{Z},	\\
}
\]
\[ 
\xymatrix@R+0pc@C+0pc{
\underset{r}{\mathrm{homotopylimit}}~\mathrm{Spec}^\mathrm{BK}{\Phi}^r_{*,\circ}/\mathrm{Fro}^\mathbb{Z},\underset{I}{\mathrm{homotopycolimit}}~\mathrm{Spec}^\mathrm{BK}{\Phi}^I_{*,\circ}/\mathrm{Fro}^\mathbb{Z}.	
}
\]

\end{definition}

\indent Meanwhile we have the corresponding Clausen-Scholze analytic stacks from \cite{10CS2}, therefore applying their construction we have:

\begin{definition}
Here we define the following products by using the solidified tensor product from \cite{10CS1} and \cite{10CS2}. Then we take solidified tensor product $\overset{\blacksquare}{\otimes}$ of any of the following
\[
\xymatrix@R+0pc@C+0pc{
\widetilde{\Delta}_{*},\widetilde{\nabla}_{*},\widetilde{\Phi}_{*},\widetilde{\Delta}^+_{*},\widetilde{\nabla}^+_{*},\widetilde{\Delta}^\dagger_{*},\widetilde{\nabla}^\dagger_{*},\widetilde{\Phi}^r_{*},\widetilde{\Phi}^I_{*}, 
}
\]

\[
\xymatrix@R+0pc@C+0pc{
\breve{\Delta}_{*},\breve{\nabla}_{*},\breve{\Phi}_{*},\breve{\Delta}^+_{*},\breve{\nabla}^+_{*},\breve{\Delta}^\dagger_{*},\breve{\nabla}^\dagger_{*},\breve{\Phi}^r_{*},\breve{\Phi}^I_{*},	
}
\]

\[
\xymatrix@R+0pc@C+0pc{
{\Delta}_{*},{\nabla}_{*},{\Phi}_{*},{\Delta}^+_{*},{\nabla}^+_{*},{\Delta}^\dagger_{*},{\nabla}^\dagger_{*},{\Phi}^r_{*},{\Phi}^I_{*},	
}
\]  	
with $\circ$. Then we have the notations:
\[
\xymatrix@R+0pc@C+0pc{
\widetilde{\Delta}_{*,\circ},\widetilde{\nabla}_{*,\circ},\widetilde{\Phi}_{*,\circ},\widetilde{\Delta}^+_{*,\circ},\widetilde{\nabla}^+_{*,\circ},\widetilde{\Delta}^\dagger_{*,\circ},\widetilde{\nabla}^\dagger_{*,\circ},\widetilde{\Phi}^r_{*,\circ},\widetilde{\Phi}^I_{*,\circ}, 
}
\]

\[
\xymatrix@R+0pc@C+0pc{
\breve{\Delta}_{*,\circ},\breve{\nabla}_{*,\circ},\breve{\Phi}_{*,\circ},\breve{\Delta}^+_{*,\circ},\breve{\nabla}^+_{*,\circ},\breve{\Delta}^\dagger_{*,\circ},\breve{\nabla}^\dagger_{*,\circ},\breve{\Phi}^r_{*,\circ},\breve{\Phi}^I_{*,\circ},	
}
\]

\[
\xymatrix@R+0pc@C+0pc{
{\Delta}_{*,\circ},{\nabla}_{*,\circ},{\Phi}_{*,\circ},{\Delta}^+_{*,\circ},{\nabla}^+_{*,\circ},{\Delta}^\dagger_{*,\circ},{\nabla}^\dagger_{*,\circ},{\Phi}^r_{*,\circ},{\Phi}^I_{*,\circ}.	
}
\]
\end{definition}

\begin{definition}
First we consider the Clausen-Scholze spectrum $\mathrm{Spec}^\mathrm{CS}(*)$ attached to any of those in the above from \cite{10CS2} by taking derived rational localization:
\begin{align}
\mathrm{Spec}^\mathrm{CS}\widetilde{\Delta}_{*,\circ},\mathrm{Spec}^\mathrm{CS}\widetilde{\nabla}_{*,\circ},\mathrm{Spec}^\mathrm{CS}\widetilde{\Phi}_{*,\circ},\mathrm{Spec}^\mathrm{CS}\widetilde{\Delta}^+_{*,\circ},\mathrm{Spec}^\mathrm{CS}\widetilde{\nabla}^+_{*,\circ},\\
\mathrm{Spec}^\mathrm{CS}\widetilde{\Delta}^\dagger_{*,\circ},\mathrm{Spec}^\mathrm{CS}\widetilde{\nabla}^\dagger_{*,\circ},\mathrm{Spec}^\mathrm{CS}\widetilde{\Phi}^r_{*,\circ},\mathrm{Spec}^\mathrm{CS}\widetilde{\Phi}^I_{*,\circ},	\\
\end{align}
\begin{align}
\mathrm{Spec}^\mathrm{CS}\breve{\Delta}_{*,\circ},\breve{\nabla}_{*,\circ},\mathrm{Spec}^\mathrm{CS}\breve{\Phi}_{*,\circ},\mathrm{Spec}^\mathrm{CS}\breve{\Delta}^+_{*,\circ},\mathrm{Spec}^\mathrm{CS}\breve{\nabla}^+_{*,\circ},\\
\mathrm{Spec}^\mathrm{CS}\breve{\Delta}^\dagger_{*,\circ},\mathrm{Spec}^\mathrm{CS}\breve{\nabla}^\dagger_{*,\circ},\mathrm{Spec}^\mathrm{CS}\breve{\Phi}^r_{*,\circ},\breve{\Phi}^I_{*,\circ},	\\
\end{align}
\begin{align}
\mathrm{Spec}^\mathrm{CS}{\Delta}_{*,\circ},\mathrm{Spec}^\mathrm{CS}{\nabla}_{*,\circ},\mathrm{Spec}^\mathrm{CS}{\Phi}_{*,\circ},\mathrm{Spec}^\mathrm{CS}{\Delta}^+_{*,\circ},\mathrm{Spec}^\mathrm{CS}{\nabla}^+_{*,\circ},\\
\mathrm{Spec}^\mathrm{CS}{\Delta}^\dagger_{*,\circ},\mathrm{Spec}^\mathrm{CS}{\nabla}^\dagger_{*,\circ},\mathrm{Spec}^\mathrm{CS}{\Phi}^r_{*,\circ},\mathrm{Spec}^\mathrm{CS}{\Phi}^I_{*,\circ}.	
\end{align}

Then we take the corresponding quotients by using the corresponding Frobenius operators:
\begin{align}
&\mathrm{Spec}^\mathrm{CS}\widetilde{\Delta}_{*,\circ}/\mathrm{Fro}^\mathbb{Z},\mathrm{Spec}^\mathrm{CS}\widetilde{\nabla}_{*,\circ}/\mathrm{Fro}^\mathbb{Z},\mathrm{Spec}^\mathrm{CS}\widetilde{\Phi}_{*,\circ}/\mathrm{Fro}^\mathbb{Z},\mathrm{Spec}^\mathrm{CS}\widetilde{\Delta}^+_{*,\circ}/\mathrm{Fro}^\mathbb{Z},\\
&\mathrm{Spec}^\mathrm{CS}\widetilde{\nabla}^+_{*,\circ}/\mathrm{Fro}^\mathbb{Z}, \mathrm{Spec}^\mathrm{CS}\widetilde{\Delta}^\dagger_{*,\circ}/\mathrm{Fro}^\mathbb{Z},\mathrm{Spec}^\mathrm{CS}\widetilde{\nabla}^\dagger_{*,\circ}/\mathrm{Fro}^\mathbb{Z},	\\
\end{align}
\begin{align}
&\mathrm{Spec}^\mathrm{CS}\breve{\Delta}_{*,\circ}/\mathrm{Fro}^\mathbb{Z},\breve{\nabla}_{*,\circ}/\mathrm{Fro}^\mathbb{Z},\mathrm{Spec}^\mathrm{CS}\breve{\Phi}_{*,\circ}/\mathrm{Fro}^\mathbb{Z},\mathrm{Spec}^\mathrm{CS}\breve{\Delta}^+_{*,\circ}/\mathrm{Fro}^\mathbb{Z},\\
&\mathrm{Spec}^\mathrm{CS}\breve{\nabla}^+_{*,\circ}/\mathrm{Fro}^\mathbb{Z}, \mathrm{Spec}^\mathrm{CS}\breve{\Delta}^\dagger_{*,\circ}/\mathrm{Fro}^\mathbb{Z},\mathrm{Spec}^\mathrm{CS}\breve{\nabla}^\dagger_{*,\circ}/\mathrm{Fro}^\mathbb{Z},	\\
\end{align}
\begin{align}
&\mathrm{Spec}^\mathrm{CS}{\Delta}_{*,\circ}/\mathrm{Fro}^\mathbb{Z},\mathrm{Spec}^\mathrm{CS}{\nabla}_{*,\circ}/\mathrm{Fro}^\mathbb{Z},\mathrm{Spec}^\mathrm{CS}{\Phi}_{*,\circ}/\mathrm{Fro}^\mathbb{Z},\mathrm{Spec}^\mathrm{CS}{\Delta}^+_{*,\circ}/\mathrm{Fro}^\mathbb{Z},\\
&\mathrm{Spec}^\mathrm{CS}{\nabla}^+_{*,\circ}/\mathrm{Fro}^\mathbb{Z}, \mathrm{Spec}^\mathrm{CS}{\Delta}^\dagger_{*,\circ}/\mathrm{Fro}^\mathbb{Z},\mathrm{Spec}^\mathrm{CS}{\nabla}^\dagger_{*,\circ}/\mathrm{Fro}^\mathbb{Z}.	
\end{align}
Here for those space with notations related to the radius and the corresponding interval we consider the total unions $\bigcap_r,\bigcup_I$ in order to achieve the whole spaces to achieve the analogues of the corresponding FF curves from \cite{10KL1}, \cite{10KL2}, \cite{10FF} for
\[
\xymatrix@R+0pc@C+0pc{
\underset{r}{\mathrm{homotopylimit}}~\mathrm{Spec}^\mathrm{CS}\widetilde{\Phi}^r_{*,\circ},\underset{I}{\mathrm{homotopycolimit}}~\mathrm{Spec}^\mathrm{CS}\widetilde{\Phi}^I_{*,\circ},	\\
}
\]
\[
\xymatrix@R+0pc@C+0pc{
\underset{r}{\mathrm{homotopylimit}}~\mathrm{Spec}^\mathrm{CS}\breve{\Phi}^r_{*,\circ},\underset{I}{\mathrm{homotopycolimit}}~\mathrm{Spec}^\mathrm{CS}\breve{\Phi}^I_{*,\circ},	\\
}
\]
\[
\xymatrix@R+0pc@C+0pc{
\underset{r}{\mathrm{homotopylimit}}~\mathrm{Spec}^\mathrm{CS}{\Phi}^r_{*,\circ},\underset{I}{\mathrm{homotopycolimit}}~\mathrm{Spec}^\mathrm{CS}{\Phi}^I_{*,\circ}.	
}
\]
\[ 
\xymatrix@R+0pc@C+0pc{
\underset{r}{\mathrm{homotopylimit}}~\mathrm{Spec}^\mathrm{CS}\widetilde{\Phi}^r_{*,\circ}/\mathrm{Fro}^\mathbb{Z},\underset{I}{\mathrm{homotopycolimit}}~\mathrm{Spec}^\mathrm{CS}\widetilde{\Phi}^I_{*,\circ}/\mathrm{Fro}^\mathbb{Z},	\\
}
\]
\[ 
\xymatrix@R+0pc@C+0pc{
\underset{r}{\mathrm{homotopylimit}}~\mathrm{Spec}^\mathrm{CS}\breve{\Phi}^r_{*,\circ}/\mathrm{Fro}^\mathbb{Z},\underset{I}{\mathrm{homotopycolimit}}~\breve{\Phi}^I_{*,\circ}/\mathrm{Fro}^\mathbb{Z},	\\
}
\]
\[ 
\xymatrix@R+0pc@C+0pc{
\underset{r}{\mathrm{homotopylimit}}~\mathrm{Spec}^\mathrm{CS}{\Phi}^r_{*,\circ}/\mathrm{Fro}^\mathbb{Z},\underset{I}{\mathrm{homotopycolimit}}~\mathrm{Spec}^\mathrm{CS}{\Phi}^I_{*,\circ}/\mathrm{Fro}^\mathbb{Z}.	
}
\]

\end{definition}

\

\begin{definition}
We then consider the corresponding quasipresheaves of the corresponding ind-Banach or monomorphic ind-Banach modules from \cite{10BBK}, \cite{10KKM}:
\begin{align}
\mathrm{Quasicoherentpresheaves,IndBanach}_{*}	
\end{align}
where $*$ is one of the following spaces:
\begin{align}
&\mathrm{Spec}^\mathrm{BK}\widetilde{\Phi}_{*,\circ}/\mathrm{Fro}^\mathbb{Z},	\\
\end{align}
\begin{align}
&\mathrm{Spec}^\mathrm{BK}\breve{\Phi}_{*,\circ}/\mathrm{Fro}^\mathbb{Z},	\\
\end{align}
\begin{align}
&\mathrm{Spec}^\mathrm{BK}{\Phi}_{*,\circ}/\mathrm{Fro}^\mathbb{Z}.	
\end{align}
Here for those space without notation related to the radius and the corresponding interval we consider the total unions $\bigcap_r,\bigcup_I$ in order to achieve the whole spaces to achieve the analogues of the corresponding FF curves from \cite{10KL1}, \cite{10KL2}, \cite{10FF} for
\[
\xymatrix@R+0pc@C+0pc{
\underset{r}{\mathrm{homotopylimit}}~\mathrm{Spec}^\mathrm{BK}\widetilde{\Phi}^r_{*,\circ},\underset{I}{\mathrm{homotopycolimit}}~\mathrm{Spec}^\mathrm{BK}\widetilde{\Phi}^I_{*,\circ},	\\
}
\]
\[
\xymatrix@R+0pc@C+0pc{
\underset{r}{\mathrm{homotopylimit}}~\mathrm{Spec}^\mathrm{BK}\breve{\Phi}^r_{*,\circ},\underset{I}{\mathrm{homotopycolimit}}~\mathrm{Spec}^\mathrm{BK}\breve{\Phi}^I_{*,\circ},	\\
}
\]
\[
\xymatrix@R+0pc@C+0pc{
\underset{r}{\mathrm{homotopylimit}}~\mathrm{Spec}^\mathrm{BK}{\Phi}^r_{*,\circ},\underset{I}{\mathrm{homotopycolimit}}~\mathrm{Spec}^\mathrm{BK}{\Phi}^I_{*,\circ}.	
}
\]
\[  
\xymatrix@R+0pc@C+0pc{
\underset{r}{\mathrm{homotopylimit}}~\mathrm{Spec}^\mathrm{BK}\widetilde{\Phi}^r_{*,\circ}/\mathrm{Fro}^\mathbb{Z},\underset{I}{\mathrm{homotopycolimit}}~\mathrm{Spec}^\mathrm{BK}\widetilde{\Phi}^I_{*,\circ}/\mathrm{Fro}^\mathbb{Z},	\\
}
\]
\[ 
\xymatrix@R+0pc@C+0pc{
\underset{r}{\mathrm{homotopylimit}}~\mathrm{Spec}^\mathrm{BK}\breve{\Phi}^r_{*,\circ}/\mathrm{Fro}^\mathbb{Z},\underset{I}{\mathrm{homotopycolimit}}~\mathrm{Spec}^\mathrm{BK}\breve{\Phi}^I_{*,\circ}/\mathrm{Fro}^\mathbb{Z},	\\
}
\]
\[ 
\xymatrix@R+0pc@C+0pc{
\underset{r}{\mathrm{homotopylimit}}~\mathrm{Spec}^\mathrm{BK}{\Phi}^r_{*,\circ}/\mathrm{Fro}^\mathbb{Z},\underset{I}{\mathrm{homotopycolimit}}~\mathrm{Spec}^\mathrm{BK}{\Phi}^I_{*,\circ}/\mathrm{Fro}^\mathbb{Z}.	
}
\]

\end{definition}

\begin{definition}
We then consider the corresponding quasisheaves of the corresponding condensed solid topological modules from \cite{10CS2}:
\begin{align}
\mathrm{Quasicoherentsheaves, Condensed}_{*}	
\end{align}
where $*$ is one of the following spaces:
\begin{align}
&\mathrm{Spec}^\mathrm{CS}\widetilde{\Delta}_{*,\circ}/\mathrm{Fro}^\mathbb{Z},\mathrm{Spec}^\mathrm{CS}\widetilde{\nabla}_{*,\circ}/\mathrm{Fro}^\mathbb{Z},\mathrm{Spec}^\mathrm{CS}\widetilde{\Phi}_{*,\circ}/\mathrm{Fro}^\mathbb{Z},\mathrm{Spec}^\mathrm{CS}\widetilde{\Delta}^+_{*,\circ}/\mathrm{Fro}^\mathbb{Z},\\
&\mathrm{Spec}^\mathrm{CS}\widetilde{\nabla}^+_{*,\circ}/\mathrm{Fro}^\mathbb{Z},\mathrm{Spec}^\mathrm{CS}\widetilde{\Delta}^\dagger_{*,\circ}/\mathrm{Fro}^\mathbb{Z},\mathrm{Spec}^\mathrm{CS}\widetilde{\nabla}^\dagger_{*,\circ}/\mathrm{Fro}^\mathbb{Z},	\\
\end{align}
\begin{align}
&\mathrm{Spec}^\mathrm{CS}\breve{\Delta}_{*,\circ}/\mathrm{Fro}^\mathbb{Z},\breve{\nabla}_{*,\circ}/\mathrm{Fro}^\mathbb{Z},\mathrm{Spec}^\mathrm{CS}\breve{\Phi}_{*,\circ}/\mathrm{Fro}^\mathbb{Z},\mathrm{Spec}^\mathrm{CS}\breve{\Delta}^+_{*,\circ}/\mathrm{Fro}^\mathbb{Z},\\
&\mathrm{Spec}^\mathrm{CS}\breve{\nabla}^+_{*,\circ}/\mathrm{Fro}^\mathbb{Z},\mathrm{Spec}^\mathrm{CS}\breve{\Delta}^\dagger_{*,\circ}/\mathrm{Fro}^\mathbb{Z},\mathrm{Spec}^\mathrm{CS}\breve{\nabla}^\dagger_{*,\circ}/\mathrm{Fro}^\mathbb{Z},	\\
\end{align}
\begin{align}
&\mathrm{Spec}^\mathrm{CS}{\Delta}_{*,\circ}/\mathrm{Fro}^\mathbb{Z},\mathrm{Spec}^\mathrm{CS}{\nabla}_{*,\circ}/\mathrm{Fro}^\mathbb{Z},\mathrm{Spec}^\mathrm{CS}{\Phi}_{*,\circ}/\mathrm{Fro}^\mathbb{Z},\mathrm{Spec}^\mathrm{CS}{\Delta}^+_{*,\circ}/\mathrm{Fro}^\mathbb{Z},\\
&\mathrm{Spec}^\mathrm{CS}{\nabla}^+_{*,\circ}/\mathrm{Fro}^\mathbb{Z}, \mathrm{Spec}^\mathrm{CS}{\Delta}^\dagger_{*,\circ}/\mathrm{Fro}^\mathbb{Z},\mathrm{Spec}^\mathrm{CS}{\nabla}^\dagger_{*,\circ}/\mathrm{Fro}^\mathbb{Z}.	
\end{align}
Here for those space with notations related to the radius and the corresponding interval we consider the total unions $\bigcap_r,\bigcup_I$ in order to achieve the whole spaces to achieve the analogues of the corresponding FF curves from \cite{10KL1}, \cite{10KL2}, \cite{10FF} for
\[
\xymatrix@R+0pc@C+0pc{
\underset{r}{\mathrm{homotopylimit}}~\mathrm{Spec}^\mathrm{CS}\widetilde{\Phi}^r_{*,\circ},\underset{I}{\mathrm{homotopycolimit}}~\mathrm{Spec}^\mathrm{CS}\widetilde{\Phi}^I_{*,\circ},	\\
}
\]
\[
\xymatrix@R+0pc@C+0pc{
\underset{r}{\mathrm{homotopylimit}}~\mathrm{Spec}^\mathrm{CS}\breve{\Phi}^r_{*,\circ},\underset{I}{\mathrm{homotopycolimit}}~\mathrm{Spec}^\mathrm{CS}\breve{\Phi}^I_{*,\circ},	\\
}
\]
\[
\xymatrix@R+0pc@C+0pc{
\underset{r}{\mathrm{homotopylimit}}~\mathrm{Spec}^\mathrm{CS}{\Phi}^r_{*,\circ},\underset{I}{\mathrm{homotopycolimit}}~\mathrm{Spec}^\mathrm{CS}{\Phi}^I_{*,\circ}.	
}
\]
\[ 
\xymatrix@R+0pc@C+0pc{
\underset{r}{\mathrm{homotopylimit}}~\mathrm{Spec}^\mathrm{CS}\widetilde{\Phi}^r_{*,\circ}/\mathrm{Fro}^\mathbb{Z},\underset{I}{\mathrm{homotopycolimit}}~\mathrm{Spec}^\mathrm{CS}\widetilde{\Phi}^I_{*,\circ}/\mathrm{Fro}^\mathbb{Z},	\\
}
\]
\[ 
\xymatrix@R+0pc@C+0pc{
\underset{r}{\mathrm{homotopylimit}}~\mathrm{Spec}^\mathrm{CS}\breve{\Phi}^r_{*,\circ}/\mathrm{Fro}^\mathbb{Z},\underset{I}{\mathrm{homotopycolimit}}~\breve{\Phi}^I_{*,\circ}/\mathrm{Fro}^\mathbb{Z},	\\
}
\]
\[ 
\xymatrix@R+0pc@C+0pc{
\underset{r}{\mathrm{homotopylimit}}~\mathrm{Spec}^\mathrm{CS}{\Phi}^r_{*,\circ}/\mathrm{Fro}^\mathbb{Z},\underset{I}{\mathrm{homotopycolimit}}~\mathrm{Spec}^\mathrm{CS}{\Phi}^I_{*,\circ}/\mathrm{Fro}^\mathbb{Z}.	
}
\]

\end{definition}

\

\begin{proposition}
There is a well-defined functor from the $\infty$-category 
\begin{align}
\mathrm{Quasicoherentpresheaves,Condensed}_{*}	
\end{align}
where $*$ is one of the following spaces:
\begin{align}
&\mathrm{Spec}^\mathrm{CS}\widetilde{\Phi}_{*,\circ}/\mathrm{Fro}^\mathbb{Z},	\\
\end{align}
\begin{align}
&\mathrm{Spec}^\mathrm{CS}\breve{\Phi}_{*,\circ}/\mathrm{Fro}^\mathbb{Z},	\\
\end{align}
\begin{align}
&\mathrm{Spec}^\mathrm{CS}{\Phi}_{*,\circ}/\mathrm{Fro}^\mathbb{Z},	
\end{align}
to the $\infty$-category of $\mathrm{Fro}$-equivariant quasicoherent presheaves over similar spaces above correspondingly without the $\mathrm{Fro}$-quotients, and to the $\infty$-category of $\mathrm{Fro}$-equivariant quasicoherent modules over global sections of the structure $\infty$-sheaves of the similar spaces above correspondingly without the $\mathrm{Fro}$-quotients. Here for those space without notation related to the radius and the corresponding interval we consider the total unions $\bigcap_r,\bigcup_I$ in order to achieve the whole spaces to achieve the analogues of the corresponding FF curves from \cite{10KL1}, \cite{10KL2}, \cite{10FF} for
\[
\xymatrix@R+0pc@C+0pc{
\underset{r}{\mathrm{homotopylimit}}~\mathrm{Spec}^\mathrm{CS}\widetilde{\Phi}^r_{*,\circ},\underset{I}{\mathrm{homotopycolimit}}~\mathrm{Spec}^\mathrm{CS}\widetilde{\Phi}^I_{*,\circ},	\\
}
\]
\[
\xymatrix@R+0pc@C+0pc{
\underset{r}{\mathrm{homotopylimit}}~\mathrm{Spec}^\mathrm{CS}\breve{\Phi}^r_{*,\circ},\underset{I}{\mathrm{homotopycolimit}}~\mathrm{Spec}^\mathrm{CS}\breve{\Phi}^I_{*,\circ},	\\
}
\]
\[
\xymatrix@R+0pc@C+0pc{
\underset{r}{\mathrm{homotopylimit}}~\mathrm{Spec}^\mathrm{CS}{\Phi}^r_{*,\circ},\underset{I}{\mathrm{homotopycolimit}}~\mathrm{Spec}^\mathrm{CS}{\Phi}^I_{*,\circ}.	
}
\]
\[ 
\xymatrix@R+0pc@C+0pc{
\underset{r}{\mathrm{homotopylimit}}~\mathrm{Spec}^\mathrm{CS}\widetilde{\Phi}^r_{*,\circ}/\mathrm{Fro}^\mathbb{Z},\underset{I}{\mathrm{homotopycolimit}}~\mathrm{Spec}^\mathrm{CS}\widetilde{\Phi}^I_{*,\circ}/\mathrm{Fro}^\mathbb{Z},	\\
}
\]
\[ 
\xymatrix@R+0pc@C+0pc{
\underset{r}{\mathrm{homotopylimit}}~\mathrm{Spec}^\mathrm{CS}\breve{\Phi}^r_{*,\circ}/\mathrm{Fro}^\mathbb{Z},\underset{I}{\mathrm{homotopycolimit}}~\breve{\Phi}^I_{*,\circ}/\mathrm{Fro}^\mathbb{Z},	\\
}
\]
\[ 
\xymatrix@R+0pc@C+0pc{
\underset{r}{\mathrm{homotopylimit}}~\mathrm{Spec}^\mathrm{CS}{\Phi}^r_{*,\circ}/\mathrm{Fro}^\mathbb{Z},\underset{I}{\mathrm{homotopycolimit}}~\mathrm{Spec}^\mathrm{CS}{\Phi}^I_{*,\circ}/\mathrm{Fro}^\mathbb{Z}.	
}
\]	
In this situation we will have the target category being family parametrized by $r$ or $I$ in compatible glueing sense as in \cite[Definition 5.4.10]{10KL2}. In this situation for modules parametrized by the intervals we have the equivalence of $\infty$-categories by using \cite[Proposition 13.8]{10CS2}. Here the corresponding quasicoherent Frobenius modules are defined to be the corresponding homotopy colimits and limits of Frobenius modules:
\begin{align}
\underset{r}{\mathrm{homotopycolimit}}~M_r,\\
\underset{I}{\mathrm{homotopylimit}}~M_I,	
\end{align}
where each $M_r$ is a Frobenius-equivariant module over the period ring with respect to some radius $r$ while each $M_I$ is a Frobenius-equivariant module over the period ring with respect to some interval $I$.\\
\end{proposition}

\begin{proposition}
Similar proposition holds for 
\begin{align}
\mathrm{Quasicoherentsheaves,IndBanach}_{*}.	
\end{align}	
\end{proposition}

\

\begin{definition}
We then consider the corresponding quasipresheaves of perfect complexes the corresponding ind-Banach or monomorphic ind-Banach modules from \cite{10BBK}, \cite{10KKM}:
\begin{align}
\mathrm{Quasicoherentpresheaves,Perfectcomplex,IndBanach}_{*}	
\end{align}
where $*$ is one of the following spaces:
\begin{align}
&\mathrm{Spec}^\mathrm{BK}\widetilde{\Phi}_{*,\circ}/\mathrm{Fro}^\mathbb{Z},	\\
\end{align}
\begin{align}
&\mathrm{Spec}^\mathrm{BK}\breve{\Phi}_{*,\circ}/\mathrm{Fro}^\mathbb{Z},	\\
\end{align}
\begin{align}
&\mathrm{Spec}^\mathrm{BK}{\Phi}_{*,\circ}/\mathrm{Fro}^\mathbb{Z}.	
\end{align}
Here for those space without notation related to the radius and the corresponding interval we consider the total unions $\bigcap_r,\bigcup_I$ in order to achieve the whole spaces to achieve the analogues of the corresponding FF curves from \cite{10KL1}, \cite{10KL2}, \cite{10FF} for
\[
\xymatrix@R+0pc@C+0pc{
\underset{r}{\mathrm{homotopylimit}}~\mathrm{Spec}^\mathrm{BK}\widetilde{\Phi}^r_{*,\circ},\underset{I}{\mathrm{homotopycolimit}}~\mathrm{Spec}^\mathrm{BK}\widetilde{\Phi}^I_{*,\circ},	\\
}
\]
\[
\xymatrix@R+0pc@C+0pc{
\underset{r}{\mathrm{homotopylimit}}~\mathrm{Spec}^\mathrm{BK}\breve{\Phi}^r_{*,\circ},\underset{I}{\mathrm{homotopycolimit}}~\mathrm{Spec}^\mathrm{BK}\breve{\Phi}^I_{*,\circ},	\\
}
\]
\[
\xymatrix@R+0pc@C+0pc{
\underset{r}{\mathrm{homotopylimit}}~\mathrm{Spec}^\mathrm{BK}{\Phi}^r_{*,\circ},\underset{I}{\mathrm{homotopycolimit}}~\mathrm{Spec}^\mathrm{BK}{\Phi}^I_{*,\circ}.	
}
\]
\[  
\xymatrix@R+0pc@C+0pc{
\underset{r}{\mathrm{homotopylimit}}~\mathrm{Spec}^\mathrm{BK}\widetilde{\Phi}^r_{*,\circ}/\mathrm{Fro}^\mathbb{Z},\underset{I}{\mathrm{homotopycolimit}}~\mathrm{Spec}^\mathrm{BK}\widetilde{\Phi}^I_{*,\circ}/\mathrm{Fro}^\mathbb{Z},	\\
}
\]
\[ 
\xymatrix@R+0pc@C+0pc{
\underset{r}{\mathrm{homotopylimit}}~\mathrm{Spec}^\mathrm{BK}\breve{\Phi}^r_{*,\circ}/\mathrm{Fro}^\mathbb{Z},\underset{I}{\mathrm{homotopycolimit}}~\mathrm{Spec}^\mathrm{BK}\breve{\Phi}^I_{*,\circ}/\mathrm{Fro}^\mathbb{Z},	\\
}
\]
\[ 
\xymatrix@R+0pc@C+0pc{
\underset{r}{\mathrm{homotopylimit}}~\mathrm{Spec}^\mathrm{BK}{\Phi}^r_{*,\circ}/\mathrm{Fro}^\mathbb{Z},\underset{I}{\mathrm{homotopycolimit}}~\mathrm{Spec}^\mathrm{BK}{\Phi}^I_{*,\circ}/\mathrm{Fro}^\mathbb{Z}.	
}
\]

\end{definition}

\begin{definition}
We then consider the corresponding quasisheaves of perfect complexes of the corresponding condensed solid topological modules from \cite{10CS2}:
\begin{align}
\mathrm{Quasicoherentsheaves, Perfectcomplex, Condensed}_{*}	
\end{align}
where $*$ is one of the following spaces:
\begin{align}
&\mathrm{Spec}^\mathrm{CS}\widetilde{\Delta}_{*,\circ}/\mathrm{Fro}^\mathbb{Z},\mathrm{Spec}^\mathrm{CS}\widetilde{\nabla}_{*,\circ}/\mathrm{Fro}^\mathbb{Z},\mathrm{Spec}^\mathrm{CS}\widetilde{\Phi}_{*,\circ}/\mathrm{Fro}^\mathbb{Z},\mathrm{Spec}^\mathrm{CS}\widetilde{\Delta}^+_{*,\circ}/\mathrm{Fro}^\mathbb{Z},\\
&\mathrm{Spec}^\mathrm{CS}\widetilde{\nabla}^+_{*,\circ}/\mathrm{Fro}^\mathbb{Z},\mathrm{Spec}^\mathrm{CS}\widetilde{\Delta}^\dagger_{*,\circ}/\mathrm{Fro}^\mathbb{Z},\mathrm{Spec}^\mathrm{CS}\widetilde{\nabla}^\dagger_{*,\circ}/\mathrm{Fro}^\mathbb{Z},	\\
\end{align}
\begin{align}
&\mathrm{Spec}^\mathrm{CS}\breve{\Delta}_{*,\circ}/\mathrm{Fro}^\mathbb{Z},\breve{\nabla}_{*,\circ}/\mathrm{Fro}^\mathbb{Z},\mathrm{Spec}^\mathrm{CS}\breve{\Phi}_{*,\circ}/\mathrm{Fro}^\mathbb{Z},\mathrm{Spec}^\mathrm{CS}\breve{\Delta}^+_{*,\circ}/\mathrm{Fro}^\mathbb{Z},\\
&\mathrm{Spec}^\mathrm{CS}\breve{\nabla}^+_{*,\circ}/\mathrm{Fro}^\mathbb{Z},\mathrm{Spec}^\mathrm{CS}\breve{\Delta}^\dagger_{*,\circ}/\mathrm{Fro}^\mathbb{Z},\mathrm{Spec}^\mathrm{CS}\breve{\nabla}^\dagger_{*,\circ}/\mathrm{Fro}^\mathbb{Z},	\\
\end{align}
\begin{align}
&\mathrm{Spec}^\mathrm{CS}{\Delta}_{*,\circ}/\mathrm{Fro}^\mathbb{Z},\mathrm{Spec}^\mathrm{CS}{\nabla}_{*,\circ}/\mathrm{Fro}^\mathbb{Z},\mathrm{Spec}^\mathrm{CS}{\Phi}_{*,\circ}/\mathrm{Fro}^\mathbb{Z},\mathrm{Spec}^\mathrm{CS}{\Delta}^+_{*,\circ}/\mathrm{Fro}^\mathbb{Z},\\
&\mathrm{Spec}^\mathrm{CS}{\nabla}^+_{*,\circ}/\mathrm{Fro}^\mathbb{Z}, \mathrm{Spec}^\mathrm{CS}{\Delta}^\dagger_{*,\circ}/\mathrm{Fro}^\mathbb{Z},\mathrm{Spec}^\mathrm{CS}{\nabla}^\dagger_{*,\circ}/\mathrm{Fro}^\mathbb{Z}.	
\end{align}
Here for those space with notations related to the radius and the corresponding interval we consider the total unions $\bigcap_r,\bigcup_I$ in order to achieve the whole spaces to achieve the analogues of the corresponding FF curves from \cite{10KL1}, \cite{10KL2}, \cite{10FF} for
\[
\xymatrix@R+0pc@C+0pc{
\underset{r}{\mathrm{homotopylimit}}~\mathrm{Spec}^\mathrm{CS}\widetilde{\Phi}^r_{*,\circ},\underset{I}{\mathrm{homotopycolimit}}~\mathrm{Spec}^\mathrm{CS}\widetilde{\Phi}^I_{*,\circ},	\\
}
\]
\[
\xymatrix@R+0pc@C+0pc{
\underset{r}{\mathrm{homotopylimit}}~\mathrm{Spec}^\mathrm{CS}\breve{\Phi}^r_{*,\circ},\underset{I}{\mathrm{homotopycolimit}}~\mathrm{Spec}^\mathrm{CS}\breve{\Phi}^I_{*,\circ},	\\
}
\]
\[
\xymatrix@R+0pc@C+0pc{
\underset{r}{\mathrm{homotopylimit}}~\mathrm{Spec}^\mathrm{CS}{\Phi}^r_{*,\circ},\underset{I}{\mathrm{homotopycolimit}}~\mathrm{Spec}^\mathrm{CS}{\Phi}^I_{*,\circ}.	
}
\]
\[ 
\xymatrix@R+0pc@C+0pc{
\underset{r}{\mathrm{homotopylimit}}~\mathrm{Spec}^\mathrm{CS}\widetilde{\Phi}^r_{*,\circ}/\mathrm{Fro}^\mathbb{Z},\underset{I}{\mathrm{homotopycolimit}}~\mathrm{Spec}^\mathrm{CS}\widetilde{\Phi}^I_{*,\circ}/\mathrm{Fro}^\mathbb{Z},	\\
}
\]
\[ 
\xymatrix@R+0pc@C+0pc{
\underset{r}{\mathrm{homotopylimit}}~\mathrm{Spec}^\mathrm{CS}\breve{\Phi}^r_{*,\circ}/\mathrm{Fro}^\mathbb{Z},\underset{I}{\mathrm{homotopycolimit}}~\breve{\Phi}^I_{*,\circ}/\mathrm{Fro}^\mathbb{Z},	\\
}
\]
\[ 
\xymatrix@R+0pc@C+0pc{
\underset{r}{\mathrm{homotopylimit}}~\mathrm{Spec}^\mathrm{CS}{\Phi}^r_{*,\circ}/\mathrm{Fro}^\mathbb{Z},\underset{I}{\mathrm{homotopycolimit}}~\mathrm{Spec}^\mathrm{CS}{\Phi}^I_{*,\circ}/\mathrm{Fro}^\mathbb{Z}.	
}
\]

\end{definition}

\begin{proposition}
There is a well-defined functor from the $\infty$-category 
\begin{align}
\mathrm{Quasicoherentpresheaves,Perfectcomplex,Condensed}_{*}	
\end{align}
where $*$ is one of the following spaces:
\begin{align}
&\mathrm{Spec}^\mathrm{CS}\widetilde{\Phi}_{*,\circ}/\mathrm{Fro}^\mathbb{Z},	\\
\end{align}
\begin{align}
&\mathrm{Spec}^\mathrm{CS}\breve{\Phi}_{*,\circ}/\mathrm{Fro}^\mathbb{Z},	\\
\end{align}
\begin{align}
&\mathrm{Spec}^\mathrm{CS}{\Phi}_{*,\circ}/\mathrm{Fro}^\mathbb{Z},	
\end{align}
to the $\infty$-category of $\mathrm{Fro}$-equivariant quasicoherent presheaves over similar spaces above correspondingly without the $\mathrm{Fro}$-quotients, and to the $\infty$-category of $\mathrm{Fro}$-equivariant quasicoherent modules over global sections of the structure $\infty$-sheaves of the similar spaces above correspondingly without the $\mathrm{Fro}$-quotients. Here for those space without notation related to the radius and the corresponding interval we consider the total unions $\bigcap_r,\bigcup_I$ in order to achieve the whole spaces to achieve the analogues of the corresponding FF curves from \cite{10KL1}, \cite{10KL2}, \cite{10FF} for
\[
\xymatrix@R+0pc@C+0pc{
\underset{r}{\mathrm{homotopylimit}}~\mathrm{Spec}^\mathrm{CS}\widetilde{\Phi}^r_{*,\circ},\underset{I}{\mathrm{homotopycolimit}}~\mathrm{Spec}^\mathrm{CS}\widetilde{\Phi}^I_{*,\circ},	\\
}
\]
\[
\xymatrix@R+0pc@C+0pc{
\underset{r}{\mathrm{homotopylimit}}~\mathrm{Spec}^\mathrm{CS}\breve{\Phi}^r_{*,\circ},\underset{I}{\mathrm{homotopycolimit}}~\mathrm{Spec}^\mathrm{CS}\breve{\Phi}^I_{*,\circ},	\\
}
\]
\[
\xymatrix@R+0pc@C+0pc{
\underset{r}{\mathrm{homotopylimit}}~\mathrm{Spec}^\mathrm{CS}{\Phi}^r_{*,\circ},\underset{I}{\mathrm{homotopycolimit}}~\mathrm{Spec}^\mathrm{CS}{\Phi}^I_{*,\circ}.	
}
\]
\[ 
\xymatrix@R+0pc@C+0pc{
\underset{r}{\mathrm{homotopylimit}}~\mathrm{Spec}^\mathrm{CS}\widetilde{\Phi}^r_{*,\circ}/\mathrm{Fro}^\mathbb{Z},\underset{I}{\mathrm{homotopycolimit}}~\mathrm{Spec}^\mathrm{CS}\widetilde{\Phi}^I_{*,\circ}/\mathrm{Fro}^\mathbb{Z},	\\
}
\]
\[ 
\xymatrix@R+0pc@C+0pc{
\underset{r}{\mathrm{homotopylimit}}~\mathrm{Spec}^\mathrm{CS}\breve{\Phi}^r_{*,\circ}/\mathrm{Fro}^\mathbb{Z},\underset{I}{\mathrm{homotopycolimit}}~\breve{\Phi}^I_{*,\circ}/\mathrm{Fro}^\mathbb{Z},	\\
}
\]
\[ 
\xymatrix@R+0pc@C+0pc{
\underset{r}{\mathrm{homotopylimit}}~\mathrm{Spec}^\mathrm{CS}{\Phi}^r_{*,\circ}/\mathrm{Fro}^\mathbb{Z},\underset{I}{\mathrm{homotopycolimit}}~\mathrm{Spec}^\mathrm{CS}{\Phi}^I_{*,\circ}/\mathrm{Fro}^\mathbb{Z}.	
}
\]	
In this situation we will have the target category being family parametrized by $r$ or $I$ in compatible glueing sense as in \cite[Definition 5.4.10]{10KL2}. In this situation for modules parametrized by the intervals we have the equivalence of $\infty$-categories by using \cite[Proposition 12.18]{10CS2}. Here the corresponding quasicoherent Frobenius modules are defined to be the corresponding homotopy colimits and limits of Frobenius modules:
\begin{align}
\underset{r}{\mathrm{homotopycolimit}}~M_r,\\
\underset{I}{\mathrm{homotopylimit}}~M_I,	
\end{align}
where each $M_r$ is a Frobenius-equivariant module over the period ring with respect to some radius $r$ while each $M_I$ is a Frobenius-equivariant module over the period ring with respect to some interval $I$.\\
\end{proposition}

\begin{proposition}
Similar proposition holds for 
\begin{align}
\mathrm{Quasicoherentsheaves,Perfectcomplex,IndBanach}_{*}.	
\end{align}	
\end{proposition}

\section{Over Affinoid Analytic Spaces}

This chapter follows closely \cite{10T1}, \cite{10T2}, \cite{10T3}, \cite{10KPX}, \cite{10KP}, \cite{10KL1}, \cite{10KL2}, \cite{10BK}, \cite{10BBBK}, \cite{10BBM}, \cite{10KKM}, \cite{10CS1}, \cite{10CS2}, \cite{10CKZ}, \cite{10PZ}, \cite{10BCM}, \cite{10LBV}, \cite{10T3}, \cite{10He}, \cite{10PR}, \cite{10SW}, \cite{10FS}, \cite{10RZ}, \cite{10Sch2}, where along one direction we will have the goal in mind to study the moduli stacks of Frobenius modules in some sense. \footnote{The consideration will be essentially after the work \cite{102EG}, the work \cite{10HHS} and the work \cite{10GEH}. See  \cite[Conjecture 5.1.18, Section 5.2, Theorem, 5.2.4]{10GEH} for the detail of Emerton-Gee-Hellmann conjecture on the moduli stack of $(\varphi,\Gamma)$-modules.} All the corresponding affinoid analytic spaces will be Clausen-Scholze spectra of analytic rings in \cite{10CS1} and \cite{10CS2}, in the notation of $X$, $\circ$, $X_\square$.

\subsubsection{Frobenius Quasicoherent Modules I}

\begin{definition}
Let $\psi$ be a toric tower over $\mathbb{Q}_p$ as in \cite[Chapter 7]{10KL2} with base $\mathbb{Q}_p\left<X_1^{\pm 1},...,X_k^{\pm 1}\right>$. Then from \cite{10KL1} and \cite[Definition 5.2.1]{10KL2} we have the following class of Kedlaya-Liu rings (with the following replacement: $\Delta$ stands for $A$, $\nabla$ stands for $B$, while $\Phi$ stands for $C$) by taking product in the sense of self $\Gamma$-th power:

\[
\xymatrix@R+0pc@C+0pc{
\widetilde{\Delta}_{\psi,\Gamma},\widetilde{\nabla}_{\psi,\Gamma},\widetilde{\Phi}_{\psi,\Gamma},\widetilde{\Delta}^+_{\psi,\Gamma},\widetilde{\nabla}^+_{\psi,\Gamma},\widetilde{\Delta}^\dagger_{\psi,\Gamma},\widetilde{\nabla}^\dagger_{\psi,\Gamma},\widetilde{\Phi}^r_{\psi,\Gamma},\widetilde{\Phi}^I_{\psi,\Gamma}, 
}
\]

\[
\xymatrix@R+0pc@C+0pc{
\breve{\Delta}_{\psi,\Gamma},\breve{\nabla}_{\psi,\Gamma},\breve{\Phi}_{\psi,\Gamma},\breve{\Delta}^+_{\psi,\Gamma},\breve{\nabla}^+_{\psi,\Gamma},\breve{\Delta}^\dagger_{\psi,\Gamma},\breve{\nabla}^\dagger_{\psi,\Gamma},\breve{\Phi}^r_{\psi,\Gamma},\breve{\Phi}^I_{\psi,\Gamma},	
}
\]

\[
\xymatrix@R+0pc@C+0pc{
{\Delta}_{\psi,\Gamma},{\nabla}_{\psi,\Gamma},{\Phi}_{\psi,\Gamma},{\Delta}^+_{\psi,\Gamma},{\nabla}^+_{\psi,\Gamma},{\Delta}^\dagger_{\psi,\Gamma},{\nabla}^\dagger_{\psi,\Gamma},{\Phi}^r_{\psi,\Gamma},{\Phi}^I_{\psi,\Gamma}.	
}
\]
 Taking the product we have:
\[
\xymatrix@R+0pc@C+0pc{
\widetilde{\Phi}_{\psi,\Gamma,X},\widetilde{\Phi}^r_{\psi,\Gamma,X},\widetilde{\Phi}^I_{\psi,\Gamma,X},	
}
\]
\[
\xymatrix@R+0pc@C+0pc{
\breve{\Phi}_{\psi,\Gamma,X},\breve{\Phi}^r_{\psi,\Gamma,X},\breve{\Phi}^I_{\psi,\Gamma,X},	
}
\]
\[
\xymatrix@R+0pc@C+0pc{
{\Phi}_{\psi,\Gamma,X},{\Phi}^r_{\psi,\Gamma,X},{\Phi}^I_{\psi,\Gamma,X}.	
}
\]
They carry multi Frobenius action $\varphi_\Gamma$ and multi $\mathrm{Lie}_\Gamma:=\mathbb{Z}_p^{\times\Gamma}$ action. In our current situation after \cite{10CKZ} and \cite{10PZ} we consider the following $(\infty,1)$-categories of $(\infty,1)$-modules.\\
\end{definition}

\indent Meanwhile we have the corresponding Clausen-Scholze analytic stacks from \cite{10CS2}, therefore applying their construction we have:

\begin{definition}
Here we define the following products by using the solidified tensor product from \cite{10CS1} and \cite{10CS2}. Then we take solidified tensor product $\overset{\blacksquare}{\otimes}$ of any of the following
\[
\xymatrix@R+0pc@C+0pc{
\widetilde{\Delta}_{\psi,\Gamma},\widetilde{\nabla}_{\psi,\Gamma},\widetilde{\Phi}_{\psi,\Gamma},\widetilde{\Delta}^+_{\psi,\Gamma},\widetilde{\nabla}^+_{\psi,\Gamma},\widetilde{\Delta}^\dagger_{\psi,\Gamma},\widetilde{\nabla}^\dagger_{\psi,\Gamma},\widetilde{\Phi}^r_{\psi,\Gamma},\widetilde{\Phi}^I_{\psi,\Gamma}, 
}
\]

\[
\xymatrix@R+0pc@C+0pc{
\breve{\Delta}_{\psi,\Gamma},\breve{\nabla}_{\psi,\Gamma},\breve{\Phi}_{\psi,\Gamma},\breve{\Delta}^+_{\psi,\Gamma},\breve{\nabla}^+_{\psi,\Gamma},\breve{\Delta}^\dagger_{\psi,\Gamma},\breve{\nabla}^\dagger_{\psi,\Gamma},\breve{\Phi}^r_{\psi,\Gamma},\breve{\Phi}^I_{\psi,\Gamma},	
}
\]

\[
\xymatrix@R+0pc@C+0pc{
{\Delta}_{\psi,\Gamma},{\nabla}_{\psi,\Gamma},{\Phi}_{\psi,\Gamma},{\Delta}^+_{\psi,\Gamma},{\nabla}^+_{\psi,\Gamma},{\Delta}^\dagger_{\psi,\Gamma},{\nabla}^\dagger_{\psi,\Gamma},{\Phi}^r_{\psi,\Gamma},{\Phi}^I_{\psi,\Gamma},	
}
\]  	
with $X$. Then we have the notations:
\[
\xymatrix@R+0pc@C+0pc{
\widetilde{\Delta}_{\psi,\Gamma,X},\widetilde{\nabla}_{\psi,\Gamma,X},\widetilde{\Phi}_{\psi,\Gamma,X},\widetilde{\Delta}^+_{\psi,\Gamma,X},\widetilde{\nabla}^+_{\psi,\Gamma,X},\widetilde{\Delta}^\dagger_{\psi,\Gamma,X},\widetilde{\nabla}^\dagger_{\psi,\Gamma,X},\widetilde{\Phi}^r_{\psi,\Gamma,X},\widetilde{\Phi}^I_{\psi,\Gamma,X}, 
}
\]

\[
\xymatrix@R+0pc@C+0pc{
\breve{\Delta}_{\psi,\Gamma,X},\breve{\nabla}_{\psi,\Gamma,X},\breve{\Phi}_{\psi,\Gamma,X},\breve{\Delta}^+_{\psi,\Gamma,X},\breve{\nabla}^+_{\psi,\Gamma,X},\breve{\Delta}^\dagger_{\psi,\Gamma,X},\breve{\nabla}^\dagger_{\psi,\Gamma,X},\breve{\Phi}^r_{\psi,\Gamma,X},\breve{\Phi}^I_{\psi,\Gamma,X},	
}
\]

\[
\xymatrix@R+0pc@C+0pc{
{\Delta}_{\psi,\Gamma,X},{\nabla}_{\psi,\Gamma,X},{\Phi}_{\psi,\Gamma,X},{\Delta}^+_{\psi,\Gamma,X},{\nabla}^+_{\psi,\Gamma,X},{\Delta}^\dagger_{\psi,\Gamma,X},{\nabla}^\dagger_{\psi,\Gamma,X},{\Phi}^r_{\psi,\Gamma,X},{\Phi}^I_{\psi,\Gamma,X}.	
}
\]
\end{definition}

\begin{definition}
First we consider the Clausen-Scholze spectrum $\mathrm{Spec}^\mathrm{CS}(*)$ attached to any of those in the above from \cite{10CS2} by taking derived rational localization:
\begin{align}
\mathrm{Spec}^\mathrm{CS}\widetilde{\Delta}_{\psi,\Gamma,X},\mathrm{Spec}^\mathrm{CS}\widetilde{\nabla}_{\psi,\Gamma,X},\mathrm{Spec}^\mathrm{CS}\widetilde{\Phi}_{\psi,\Gamma,X},\mathrm{Spec}^\mathrm{CS}\widetilde{\Delta}^+_{\psi,\Gamma,X},\mathrm{Spec}^\mathrm{CS}\widetilde{\nabla}^+_{\psi,\Gamma,X},\\
\mathrm{Spec}^\mathrm{CS}\widetilde{\Delta}^\dagger_{\psi,\Gamma,X},\mathrm{Spec}^\mathrm{CS}\widetilde{\nabla}^\dagger_{\psi,\Gamma,X},\mathrm{Spec}^\mathrm{CS}\widetilde{\Phi}^r_{\psi,\Gamma,X},\mathrm{Spec}^\mathrm{CS}\widetilde{\Phi}^I_{\psi,\Gamma,X},	\\
\end{align}
\begin{align}
\mathrm{Spec}^\mathrm{CS}\breve{\Delta}_{\psi,\Gamma,X},\breve{\nabla}_{\psi,\Gamma,X},\mathrm{Spec}^\mathrm{CS}\breve{\Phi}_{\psi,\Gamma,X},\mathrm{Spec}^\mathrm{CS}\breve{\Delta}^+_{\psi,\Gamma,X},\mathrm{Spec}^\mathrm{CS}\breve{\nabla}^+_{\psi,\Gamma,X},\\
\mathrm{Spec}^\mathrm{CS}\breve{\Delta}^\dagger_{\psi,\Gamma,X},\mathrm{Spec}^\mathrm{CS}\breve{\nabla}^\dagger_{\psi,\Gamma,X},\mathrm{Spec}^\mathrm{CS}\breve{\Phi}^r_{\psi,\Gamma,X},\breve{\Phi}^I_{\psi,\Gamma,X},	\\
\end{align}
\begin{align}
\mathrm{Spec}^\mathrm{CS}{\Delta}_{\psi,\Gamma,X},\mathrm{Spec}^\mathrm{CS}{\nabla}_{\psi,\Gamma,X},\mathrm{Spec}^\mathrm{CS}{\Phi}_{\psi,\Gamma,X},\mathrm{Spec}^\mathrm{CS}{\Delta}^+_{\psi,\Gamma,X},\mathrm{Spec}^\mathrm{CS}{\nabla}^+_{\psi,\Gamma,X},\\
\mathrm{Spec}^\mathrm{CS}{\Delta}^\dagger_{\psi,\Gamma,X},\mathrm{Spec}^\mathrm{CS}{\nabla}^\dagger_{\psi,\Gamma,X},\mathrm{Spec}^\mathrm{CS}{\Phi}^r_{\psi,\Gamma,X},\mathrm{Spec}^\mathrm{CS}{\Phi}^I_{\psi,\Gamma,X}.	
\end{align}

Then we take the corresponding quotients by using the corresponding Frobenius operators:
\begin{align}
&\mathrm{Spec}^\mathrm{CS}\widetilde{\Delta}_{\psi,\Gamma,X}/\mathrm{Fro}^\mathbb{Z},\mathrm{Spec}^\mathrm{CS}\widetilde{\nabla}_{\psi,\Gamma,X}/\mathrm{Fro}^\mathbb{Z},\mathrm{Spec}^\mathrm{CS}\widetilde{\Phi}_{\psi,\Gamma,X}/\mathrm{Fro}^\mathbb{Z},\mathrm{Spec}^\mathrm{CS}\widetilde{\Delta}^+_{\psi,\Gamma,X}/\mathrm{Fro}^\mathbb{Z},\\
&\mathrm{Spec}^\mathrm{CS}\widetilde{\nabla}^+_{\psi,\Gamma,X}/\mathrm{Fro}^\mathbb{Z}, \mathrm{Spec}^\mathrm{CS}\widetilde{\Delta}^\dagger_{\psi,\Gamma,X}/\mathrm{Fro}^\mathbb{Z},\mathrm{Spec}^\mathrm{CS}\widetilde{\nabla}^\dagger_{\psi,\Gamma,X}/\mathrm{Fro}^\mathbb{Z},	\\
\end{align}
\begin{align}
&\mathrm{Spec}^\mathrm{CS}\breve{\Delta}_{\psi,\Gamma,X}/\mathrm{Fro}^\mathbb{Z},\breve{\nabla}_{\psi,\Gamma,X}/\mathrm{Fro}^\mathbb{Z},\mathrm{Spec}^\mathrm{CS}\breve{\Phi}_{\psi,\Gamma,X}/\mathrm{Fro}^\mathbb{Z},\mathrm{Spec}^\mathrm{CS}\breve{\Delta}^+_{\psi,\Gamma,X}/\mathrm{Fro}^\mathbb{Z},\\
&\mathrm{Spec}^\mathrm{CS}\breve{\nabla}^+_{\psi,\Gamma,X}/\mathrm{Fro}^\mathbb{Z}, \mathrm{Spec}^\mathrm{CS}\breve{\Delta}^\dagger_{\psi,\Gamma,X}/\mathrm{Fro}^\mathbb{Z},\mathrm{Spec}^\mathrm{CS}\breve{\nabla}^\dagger_{\psi,\Gamma,X}/\mathrm{Fro}^\mathbb{Z},	\\
\end{align}
\begin{align}
&\mathrm{Spec}^\mathrm{CS}{\Delta}_{\psi,\Gamma,X}/\mathrm{Fro}^\mathbb{Z},\mathrm{Spec}^\mathrm{CS}{\nabla}_{\psi,\Gamma,X}/\mathrm{Fro}^\mathbb{Z},\mathrm{Spec}^\mathrm{CS}{\Phi}_{\psi,\Gamma,X}/\mathrm{Fro}^\mathbb{Z},\mathrm{Spec}^\mathrm{CS}{\Delta}^+_{\psi,\Gamma,X}/\mathrm{Fro}^\mathbb{Z},\\
&\mathrm{Spec}^\mathrm{CS}{\nabla}^+_{\psi,\Gamma,X}/\mathrm{Fro}^\mathbb{Z}, \mathrm{Spec}^\mathrm{CS}{\Delta}^\dagger_{\psi,\Gamma,X}/\mathrm{Fro}^\mathbb{Z},\mathrm{Spec}^\mathrm{CS}{\nabla}^\dagger_{\psi,\Gamma,X}/\mathrm{Fro}^\mathbb{Z}.	
\end{align}
Here for those space with notations related to the radius and the corresponding interval we consider the total unions $\bigcap_r,\bigcup_I$ in order to achieve the whole spaces to achieve the analogues of the corresponding FF curves from \cite{10KL1}, \cite{10KL2}, \cite{10FF} for
\[
\xymatrix@R+0pc@C+0pc{
\underset{r}{\mathrm{homotopylimit}}~\mathrm{Spec}^\mathrm{CS}\widetilde{\Phi}^r_{\psi,\Gamma,X},\underset{I}{\mathrm{homotopycolimit}}~\mathrm{Spec}^\mathrm{CS}\widetilde{\Phi}^I_{\psi,\Gamma,X},	\\
}
\]
\[
\xymatrix@R+0pc@C+0pc{
\underset{r}{\mathrm{homotopylimit}}~\mathrm{Spec}^\mathrm{CS}\breve{\Phi}^r_{\psi,\Gamma,X},\underset{I}{\mathrm{homotopycolimit}}~\mathrm{Spec}^\mathrm{CS}\breve{\Phi}^I_{\psi,\Gamma,X},	\\
}
\]
\[
\xymatrix@R+0pc@C+0pc{
\underset{r}{\mathrm{homotopylimit}}~\mathrm{Spec}^\mathrm{CS}{\Phi}^r_{\psi,\Gamma,X},\underset{I}{\mathrm{homotopycolimit}}~\mathrm{Spec}^\mathrm{CS}{\Phi}^I_{\psi,\Gamma,X}.	
}
\]
\[ 
\xymatrix@R+0pc@C+0pc{
\underset{r}{\mathrm{homotopylimit}}~\mathrm{Spec}^\mathrm{CS}\widetilde{\Phi}^r_{\psi,\Gamma,X}/\mathrm{Fro}^\mathbb{Z},\underset{I}{\mathrm{homotopycolimit}}~\mathrm{Spec}^\mathrm{CS}\widetilde{\Phi}^I_{\psi,\Gamma,X}/\mathrm{Fro}^\mathbb{Z},	\\
}
\]
\[ 
\xymatrix@R+0pc@C+0pc{
\underset{r}{\mathrm{homotopylimit}}~\mathrm{Spec}^\mathrm{CS}\breve{\Phi}^r_{\psi,\Gamma,X}/\mathrm{Fro}^\mathbb{Z},\underset{I}{\mathrm{homotopycolimit}}~\breve{\Phi}^I_{\psi,\Gamma,X}/\mathrm{Fro}^\mathbb{Z},	\\
}
\]
\[ 
\xymatrix@R+0pc@C+0pc{
\underset{r}{\mathrm{homotopylimit}}~\mathrm{Spec}^\mathrm{CS}{\Phi}^r_{\psi,\Gamma,X}/\mathrm{Fro}^\mathbb{Z},\underset{I}{\mathrm{homotopycolimit}}~\mathrm{Spec}^\mathrm{CS}{\Phi}^I_{\psi,\Gamma,X}/\mathrm{Fro}^\mathbb{Z}.	
}
\]

\end{definition}

\

\begin{definition}
We then consider the corresponding quasisheaves of the corresponding condensed solid topological modules from \cite{10CS2}:
\begin{align}
\mathrm{Quasicoherentsheaves, Condensed}_{*}	
\end{align}
where $*$ is one of the following spaces:
\begin{align}
&\mathrm{Spec}^\mathrm{CS}\widetilde{\Delta}_{\psi,\Gamma,X}/\mathrm{Fro}^\mathbb{Z},\mathrm{Spec}^\mathrm{CS}\widetilde{\nabla}_{\psi,\Gamma,X}/\mathrm{Fro}^\mathbb{Z},\mathrm{Spec}^\mathrm{CS}\widetilde{\Phi}_{\psi,\Gamma,X}/\mathrm{Fro}^\mathbb{Z},\mathrm{Spec}^\mathrm{CS}\widetilde{\Delta}^+_{\psi,\Gamma,X}/\mathrm{Fro}^\mathbb{Z},\\
&\mathrm{Spec}^\mathrm{CS}\widetilde{\nabla}^+_{\psi,\Gamma,X}/\mathrm{Fro}^\mathbb{Z},\mathrm{Spec}^\mathrm{CS}\widetilde{\Delta}^\dagger_{\psi,\Gamma,X}/\mathrm{Fro}^\mathbb{Z},\mathrm{Spec}^\mathrm{CS}\widetilde{\nabla}^\dagger_{\psi,\Gamma,X}/\mathrm{Fro}^\mathbb{Z},	\\
\end{align}
\begin{align}
&\mathrm{Spec}^\mathrm{CS}\breve{\Delta}_{\psi,\Gamma,X}/\mathrm{Fro}^\mathbb{Z},\breve{\nabla}_{\psi,\Gamma,X}/\mathrm{Fro}^\mathbb{Z},\mathrm{Spec}^\mathrm{CS}\breve{\Phi}_{\psi,\Gamma,X}/\mathrm{Fro}^\mathbb{Z},\mathrm{Spec}^\mathrm{CS}\breve{\Delta}^+_{\psi,\Gamma,X}/\mathrm{Fro}^\mathbb{Z},\\
&\mathrm{Spec}^\mathrm{CS}\breve{\nabla}^+_{\psi,\Gamma,X}/\mathrm{Fro}^\mathbb{Z},\mathrm{Spec}^\mathrm{CS}\breve{\Delta}^\dagger_{\psi,\Gamma,X}/\mathrm{Fro}^\mathbb{Z},\mathrm{Spec}^\mathrm{CS}\breve{\nabla}^\dagger_{\psi,\Gamma,X}/\mathrm{Fro}^\mathbb{Z},	\\
\end{align}
\begin{align}
&\mathrm{Spec}^\mathrm{CS}{\Delta}_{\psi,\Gamma,X}/\mathrm{Fro}^\mathbb{Z},\mathrm{Spec}^\mathrm{CS}{\nabla}_{\psi,\Gamma,X}/\mathrm{Fro}^\mathbb{Z},\mathrm{Spec}^\mathrm{CS}{\Phi}_{\psi,\Gamma,X}/\mathrm{Fro}^\mathbb{Z},\mathrm{Spec}^\mathrm{CS}{\Delta}^+_{\psi,\Gamma,X}/\mathrm{Fro}^\mathbb{Z},\\
&\mathrm{Spec}^\mathrm{CS}{\nabla}^+_{\psi,\Gamma,X}/\mathrm{Fro}^\mathbb{Z}, \mathrm{Spec}^\mathrm{CS}{\Delta}^\dagger_{\psi,\Gamma,X}/\mathrm{Fro}^\mathbb{Z},\mathrm{Spec}^\mathrm{CS}{\nabla}^\dagger_{\psi,\Gamma,X}/\mathrm{Fro}^\mathbb{Z}.	
\end{align}
Here for those space with notations related to the radius and the corresponding interval we consider the total unions $\bigcap_r,\bigcup_I$ in order to achieve the whole spaces to achieve the analogues of the corresponding FF curves from \cite{10KL1}, \cite{10KL2}, \cite{10FF} for
\[
\xymatrix@R+0pc@C+0pc{
\underset{r}{\mathrm{homotopylimit}}~\mathrm{Spec}^\mathrm{CS}\widetilde{\Phi}^r_{\psi,\Gamma,X},\underset{I}{\mathrm{homotopycolimit}}~\mathrm{Spec}^\mathrm{CS}\widetilde{\Phi}^I_{\psi,\Gamma,X},	\\
}
\]
\[
\xymatrix@R+0pc@C+0pc{
\underset{r}{\mathrm{homotopylimit}}~\mathrm{Spec}^\mathrm{CS}\breve{\Phi}^r_{\psi,\Gamma,X},\underset{I}{\mathrm{homotopycolimit}}~\mathrm{Spec}^\mathrm{CS}\breve{\Phi}^I_{\psi,\Gamma,X},	\\
}
\]
\[
\xymatrix@R+0pc@C+0pc{
\underset{r}{\mathrm{homotopylimit}}~\mathrm{Spec}^\mathrm{CS}{\Phi}^r_{\psi,\Gamma,X},\underset{I}{\mathrm{homotopycolimit}}~\mathrm{Spec}^\mathrm{CS}{\Phi}^I_{\psi,\Gamma,X}.	
}
\]
\[ 
\xymatrix@R+0pc@C+0pc{
\underset{r}{\mathrm{homotopylimit}}~\mathrm{Spec}^\mathrm{CS}\widetilde{\Phi}^r_{\psi,\Gamma,X}/\mathrm{Fro}^\mathbb{Z},\underset{I}{\mathrm{homotopycolimit}}~\mathrm{Spec}^\mathrm{CS}\widetilde{\Phi}^I_{\psi,\Gamma,X}/\mathrm{Fro}^\mathbb{Z},	\\
}
\]
\[ 
\xymatrix@R+0pc@C+0pc{
\underset{r}{\mathrm{homotopylimit}}~\mathrm{Spec}^\mathrm{CS}\breve{\Phi}^r_{\psi,\Gamma,X}/\mathrm{Fro}^\mathbb{Z},\underset{I}{\mathrm{homotopycolimit}}~\breve{\Phi}^I_{\psi,\Gamma,X}/\mathrm{Fro}^\mathbb{Z},	\\
}
\]
\[ 
\xymatrix@R+0pc@C+0pc{
\underset{r}{\mathrm{homotopylimit}}~\mathrm{Spec}^\mathrm{CS}{\Phi}^r_{\psi,\Gamma,X}/\mathrm{Fro}^\mathbb{Z},\underset{I}{\mathrm{homotopycolimit}}~\mathrm{Spec}^\mathrm{CS}{\Phi}^I_{\psi,\Gamma,X}/\mathrm{Fro}^\mathbb{Z}.	
}
\]

\end{definition}

\

\begin{proposition}
There is a well-defined functor from the $\infty$-category 
\begin{align}
\mathrm{Quasicoherentpresheaves,Condensed}_{*}	
\end{align}
where $*$ is one of the following spaces:
\begin{align}
&\mathrm{Spec}^\mathrm{CS}\widetilde{\Phi}_{\psi,\Gamma,X}/\mathrm{Fro}^\mathbb{Z},	\\
\end{align}
\begin{align}
&\mathrm{Spec}^\mathrm{CS}\breve{\Phi}_{\psi,\Gamma,X}/\mathrm{Fro}^\mathbb{Z},	\\
\end{align}
\begin{align}
&\mathrm{Spec}^\mathrm{CS}{\Phi}_{\psi,\Gamma,X}/\mathrm{Fro}^\mathbb{Z},	
\end{align}
to the $\infty$-category of $\mathrm{Fro}$-equivariant quasicoherent presheaves over similar spaces above correspondingly without the $\mathrm{Fro}$-quotients, and to the $\infty$-category of $\mathrm{Fro}$-equivariant quasicoherent modules over global sections of the structure $\infty$-sheaves of the similar spaces above correspondingly without the $\mathrm{Fro}$-quotients. Here for those space without notation related to the radius and the corresponding interval we consider the total unions $\bigcap_r,\bigcup_I$ in order to achieve the whole spaces to achieve the analogues of the corresponding FF curves from \cite{10KL1}, \cite{10KL2}, \cite{10FF} for
\[
\xymatrix@R+0pc@C+0pc{
\underset{r}{\mathrm{homotopylimit}}~\mathrm{Spec}^\mathrm{CS}\widetilde{\Phi}^r_{\psi,\Gamma,X},\underset{I}{\mathrm{homotopycolimit}}~\mathrm{Spec}^\mathrm{CS}\widetilde{\Phi}^I_{\psi,\Gamma,X},	\\
}
\]
\[
\xymatrix@R+0pc@C+0pc{
\underset{r}{\mathrm{homotopylimit}}~\mathrm{Spec}^\mathrm{CS}\breve{\Phi}^r_{\psi,\Gamma,X},\underset{I}{\mathrm{homotopycolimit}}~\mathrm{Spec}^\mathrm{CS}\breve{\Phi}^I_{\psi,\Gamma,X},	\\
}
\]
\[
\xymatrix@R+0pc@C+0pc{
\underset{r}{\mathrm{homotopylimit}}~\mathrm{Spec}^\mathrm{CS}{\Phi}^r_{\psi,\Gamma,X},\underset{I}{\mathrm{homotopycolimit}}~\mathrm{Spec}^\mathrm{CS}{\Phi}^I_{\psi,\Gamma,X}.	
}
\]
\[ 
\xymatrix@R+0pc@C+0pc{
\underset{r}{\mathrm{homotopylimit}}~\mathrm{Spec}^\mathrm{CS}\widetilde{\Phi}^r_{\psi,\Gamma,X}/\mathrm{Fro}^\mathbb{Z},\underset{I}{\mathrm{homotopycolimit}}~\mathrm{Spec}^\mathrm{CS}\widetilde{\Phi}^I_{\psi,\Gamma,X}/\mathrm{Fro}^\mathbb{Z},	\\
}
\]
\[ 
\xymatrix@R+0pc@C+0pc{
\underset{r}{\mathrm{homotopylimit}}~\mathrm{Spec}^\mathrm{CS}\breve{\Phi}^r_{\psi,\Gamma,X}/\mathrm{Fro}^\mathbb{Z},\underset{I}{\mathrm{homotopycolimit}}~\breve{\Phi}^I_{\psi,\Gamma,X}/\mathrm{Fro}^\mathbb{Z},	\\
}
\]
\[ 
\xymatrix@R+0pc@C+0pc{
\underset{r}{\mathrm{homotopylimit}}~\mathrm{Spec}^\mathrm{CS}{\Phi}^r_{\psi,\Gamma,X}/\mathrm{Fro}^\mathbb{Z},\underset{I}{\mathrm{homotopycolimit}}~\mathrm{Spec}^\mathrm{CS}{\Phi}^I_{\psi,\Gamma,X}/\mathrm{Fro}^\mathbb{Z}.	
}
\]	
In this situation we will have the target category being family parametrized by $r$ or $I$ in compatible glueing sense as in \cite[Definition 5.4.10]{10KL2}. In this situation for modules parametrized by the intervals we have the equivalence of $\infty$-categories by using \cite[Proposition 13.8]{10CS2}. Here the corresponding quasicoherent Frobenius modules are defined to be the corresponding homotopy colimits and limits of Frobenius modules:
\begin{align}
\underset{r}{\mathrm{homotopycolimit}}~M_r,\\
\underset{I}{\mathrm{homotopylimit}}~M_I,	
\end{align}
where each $M_r$ is a Frobenius-equivariant module over the period ring with respect to some radius $r$ while each $M_I$ is a Frobenius-equivariant module over the period ring with respect to some interval $I$.\\
\end{proposition}

\

\begin{definition}
We then consider the corresponding quasisheaves of perfect complexes of the corresponding condensed solid topological modules from \cite{10CS2}:
\begin{align}
\mathrm{Quasicoherentsheaves, Perfectcomplex, Condensed}_{*}	
\end{align}
where $*$ is one of the following spaces:
\begin{align}
&\mathrm{Spec}^\mathrm{CS}\widetilde{\Delta}_{\psi,\Gamma,X}/\mathrm{Fro}^\mathbb{Z},\mathrm{Spec}^\mathrm{CS}\widetilde{\nabla}_{\psi,\Gamma,X}/\mathrm{Fro}^\mathbb{Z},\mathrm{Spec}^\mathrm{CS}\widetilde{\Phi}_{\psi,\Gamma,X}/\mathrm{Fro}^\mathbb{Z},\mathrm{Spec}^\mathrm{CS}\widetilde{\Delta}^+_{\psi,\Gamma,X}/\mathrm{Fro}^\mathbb{Z},\\
&\mathrm{Spec}^\mathrm{CS}\widetilde{\nabla}^+_{\psi,\Gamma,X}/\mathrm{Fro}^\mathbb{Z},\mathrm{Spec}^\mathrm{CS}\widetilde{\Delta}^\dagger_{\psi,\Gamma,X}/\mathrm{Fro}^\mathbb{Z},\mathrm{Spec}^\mathrm{CS}\widetilde{\nabla}^\dagger_{\psi,\Gamma,X}/\mathrm{Fro}^\mathbb{Z},	\\
\end{align}
\begin{align}
&\mathrm{Spec}^\mathrm{CS}\breve{\Delta}_{\psi,\Gamma,X}/\mathrm{Fro}^\mathbb{Z},\breve{\nabla}_{\psi,\Gamma,X}/\mathrm{Fro}^\mathbb{Z},\mathrm{Spec}^\mathrm{CS}\breve{\Phi}_{\psi,\Gamma,X}/\mathrm{Fro}^\mathbb{Z},\mathrm{Spec}^\mathrm{CS}\breve{\Delta}^+_{\psi,\Gamma,X}/\mathrm{Fro}^\mathbb{Z},\\
&\mathrm{Spec}^\mathrm{CS}\breve{\nabla}^+_{\psi,\Gamma,X}/\mathrm{Fro}^\mathbb{Z},\mathrm{Spec}^\mathrm{CS}\breve{\Delta}^\dagger_{\psi,\Gamma,X}/\mathrm{Fro}^\mathbb{Z},\mathrm{Spec}^\mathrm{CS}\breve{\nabla}^\dagger_{\psi,\Gamma,X}/\mathrm{Fro}^\mathbb{Z},	\\
\end{align}
\begin{align}
&\mathrm{Spec}^\mathrm{CS}{\Delta}_{\psi,\Gamma,X}/\mathrm{Fro}^\mathbb{Z},\mathrm{Spec}^\mathrm{CS}{\nabla}_{\psi,\Gamma,X}/\mathrm{Fro}^\mathbb{Z},\mathrm{Spec}^\mathrm{CS}{\Phi}_{\psi,\Gamma,X}/\mathrm{Fro}^\mathbb{Z},\mathrm{Spec}^\mathrm{CS}{\Delta}^+_{\psi,\Gamma,X}/\mathrm{Fro}^\mathbb{Z},\\
&\mathrm{Spec}^\mathrm{CS}{\nabla}^+_{\psi,\Gamma,X}/\mathrm{Fro}^\mathbb{Z}, \mathrm{Spec}^\mathrm{CS}{\Delta}^\dagger_{\psi,\Gamma,X}/\mathrm{Fro}^\mathbb{Z},\mathrm{Spec}^\mathrm{CS}{\nabla}^\dagger_{\psi,\Gamma,X}/\mathrm{Fro}^\mathbb{Z}.	
\end{align}
Here for those space with notations related to the radius and the corresponding interval we consider the total unions $\bigcap_r,\bigcup_I$ in order to achieve the whole spaces to achieve the analogues of the corresponding FF curves from \cite{10KL1}, \cite{10KL2}, \cite{10FF} for
\[
\xymatrix@R+0pc@C+0pc{
\underset{r}{\mathrm{homotopylimit}}~\mathrm{Spec}^\mathrm{CS}\widetilde{\Phi}^r_{\psi,\Gamma,X},\underset{I}{\mathrm{homotopycolimit}}~\mathrm{Spec}^\mathrm{CS}\widetilde{\Phi}^I_{\psi,\Gamma,X},	\\
}
\]
\[
\xymatrix@R+0pc@C+0pc{
\underset{r}{\mathrm{homotopylimit}}~\mathrm{Spec}^\mathrm{CS}\breve{\Phi}^r_{\psi,\Gamma,X},\underset{I}{\mathrm{homotopycolimit}}~\mathrm{Spec}^\mathrm{CS}\breve{\Phi}^I_{\psi,\Gamma,X},	\\
}
\]
\[
\xymatrix@R+0pc@C+0pc{
\underset{r}{\mathrm{homotopylimit}}~\mathrm{Spec}^\mathrm{CS}{\Phi}^r_{\psi,\Gamma,X},\underset{I}{\mathrm{homotopycolimit}}~\mathrm{Spec}^\mathrm{CS}{\Phi}^I_{\psi,\Gamma,X}.	
}
\]
\[ 
\xymatrix@R+0pc@C+0pc{
\underset{r}{\mathrm{homotopylimit}}~\mathrm{Spec}^\mathrm{CS}\widetilde{\Phi}^r_{\psi,\Gamma,X}/\mathrm{Fro}^\mathbb{Z},\underset{I}{\mathrm{homotopycolimit}}~\mathrm{Spec}^\mathrm{CS}\widetilde{\Phi}^I_{\psi,\Gamma,X}/\mathrm{Fro}^\mathbb{Z},	\\
}
\]
\[ 
\xymatrix@R+0pc@C+0pc{
\underset{r}{\mathrm{homotopylimit}}~\mathrm{Spec}^\mathrm{CS}\breve{\Phi}^r_{\psi,\Gamma,X}/\mathrm{Fro}^\mathbb{Z},\underset{I}{\mathrm{homotopycolimit}}~\breve{\Phi}^I_{\psi,\Gamma,X}/\mathrm{Fro}^\mathbb{Z},	\\
}
\]
\[ 
\xymatrix@R+0pc@C+0pc{
\underset{r}{\mathrm{homotopylimit}}~\mathrm{Spec}^\mathrm{CS}{\Phi}^r_{\psi,\Gamma,X}/\mathrm{Fro}^\mathbb{Z},\underset{I}{\mathrm{homotopycolimit}}~\mathrm{Spec}^\mathrm{CS}{\Phi}^I_{\psi,\Gamma,X}/\mathrm{Fro}^\mathbb{Z}.	
}
\]

\end{definition}

\begin{proposition}
There is a well-defined functor from the $\infty$-category 
\begin{align}
\mathrm{Quasicoherentpresheaves,Perfectcomplex,Condensed}_{*}	
\end{align}
where $*$ is one of the following spaces:
\begin{align}
&\mathrm{Spec}^\mathrm{CS}\widetilde{\Phi}_{\psi,\Gamma,X}/\mathrm{Fro}^\mathbb{Z},	\\
\end{align}
\begin{align}
&\mathrm{Spec}^\mathrm{CS}\breve{\Phi}_{\psi,\Gamma,X}/\mathrm{Fro}^\mathbb{Z},	\\
\end{align}
\begin{align}
&\mathrm{Spec}^\mathrm{CS}{\Phi}_{\psi,\Gamma,X}/\mathrm{Fro}^\mathbb{Z},	
\end{align}
to the $\infty$-category of $\mathrm{Fro}$-equivariant quasicoherent presheaves over similar spaces above correspondingly without the $\mathrm{Fro}$-quotients, and to the $\infty$-category of $\mathrm{Fro}$-equivariant quasicoherent modules over global sections of the structure $\infty$-sheaves of the similar spaces above correspondingly without the $\mathrm{Fro}$-quotients. Here for those space without notation related to the radius and the corresponding interval we consider the total unions $\bigcap_r,\bigcup_I$ in order to achieve the whole spaces to achieve the analogues of the corresponding FF curves from \cite{10KL1}, \cite{10KL2}, \cite{10FF} for
\[
\xymatrix@R+0pc@C+0pc{
\underset{r}{\mathrm{homotopylimit}}~\mathrm{Spec}^\mathrm{CS}\widetilde{\Phi}^r_{\psi,\Gamma,X},\underset{I}{\mathrm{homotopycolimit}}~\mathrm{Spec}^\mathrm{CS}\widetilde{\Phi}^I_{\psi,\Gamma,X},	\\
}
\]
\[
\xymatrix@R+0pc@C+0pc{
\underset{r}{\mathrm{homotopylimit}}~\mathrm{Spec}^\mathrm{CS}\breve{\Phi}^r_{\psi,\Gamma,X},\underset{I}{\mathrm{homotopycolimit}}~\mathrm{Spec}^\mathrm{CS}\breve{\Phi}^I_{\psi,\Gamma,X},	\\
}
\]
\[
\xymatrix@R+0pc@C+0pc{
\underset{r}{\mathrm{homotopylimit}}~\mathrm{Spec}^\mathrm{CS}{\Phi}^r_{\psi,\Gamma,X},\underset{I}{\mathrm{homotopycolimit}}~\mathrm{Spec}^\mathrm{CS}{\Phi}^I_{\psi,\Gamma,X}.	
}
\]
\[ 
\xymatrix@R+0pc@C+0pc{
\underset{r}{\mathrm{homotopylimit}}~\mathrm{Spec}^\mathrm{CS}\widetilde{\Phi}^r_{\psi,\Gamma,X}/\mathrm{Fro}^\mathbb{Z},\underset{I}{\mathrm{homotopycolimit}}~\mathrm{Spec}^\mathrm{CS}\widetilde{\Phi}^I_{\psi,\Gamma,X}/\mathrm{Fro}^\mathbb{Z},	\\
}
\]
\[ 
\xymatrix@R+0pc@C+0pc{
\underset{r}{\mathrm{homotopylimit}}~\mathrm{Spec}^\mathrm{CS}\breve{\Phi}^r_{\psi,\Gamma,X}/\mathrm{Fro}^\mathbb{Z},\underset{I}{\mathrm{homotopycolimit}}~\breve{\Phi}^I_{\psi,\Gamma,X}/\mathrm{Fro}^\mathbb{Z},	\\
}
\]
\[ 
\xymatrix@R+0pc@C+0pc{
\underset{r}{\mathrm{homotopylimit}}~\mathrm{Spec}^\mathrm{CS}{\Phi}^r_{\psi,\Gamma,X}/\mathrm{Fro}^\mathbb{Z},\underset{I}{\mathrm{homotopycolimit}}~\mathrm{Spec}^\mathrm{CS}{\Phi}^I_{\psi,\Gamma,X}/\mathrm{Fro}^\mathbb{Z}.	
}
\]	
In this situation we will have the target category being family parametrized by $r$ or $I$ in compatible glueing sense as in \cite[Definition 5.4.10]{10KL2}. Here the corresponding quasicoherent Frobenius modules are defined to be the corresponding homotopy colimits and limits of Frobenius modules:
\begin{align}
\underset{r}{\mathrm{homotopycolimit}}~M_r,\\
\underset{I}{\mathrm{homotopylimit}}~M_I,	
\end{align}
where each $M_r$ is a Frobenius-equivariant module over the period ring with respect to some radius $r$ while each $M_I$ is a Frobenius-equivariant module over the period ring with respect to some interval $I$.\\
\end{proposition}

\subsubsection{Frobenius Quasicoherent Modules II: Moduli Stacks of Frobenius Modules}

\begin{definition}
Let $\psi$ be a toric tower over $\mathbb{Q}_p$ as in \cite[Chapter 7]{10KL2} with base $\mathbb{Q}_p\left<X_1^{\pm 1},...,X_k^{\pm 1}\right>$. Then from \cite{10KL1} and \cite[Definition 5.2.1]{10KL2} we have the following class of Kedlaya-Liu rings (with the following replacement: $\Delta$ stands for $A$, $\nabla$ stands for $B$, while $\Phi$ stands for $C$) by taking product in the sense of self $\Gamma$-th power:

\[
\xymatrix@R+0pc@C+0pc{
\widetilde{\Delta}_{\psi,\Gamma},\widetilde{\nabla}_{\psi,\Gamma},\widetilde{\Phi}_{\psi,\Gamma},\widetilde{\Delta}^+_{\psi,\Gamma},\widetilde{\nabla}^+_{\psi,\Gamma},\widetilde{\Delta}^\dagger_{\psi,\Gamma},\widetilde{\nabla}^\dagger_{\psi,\Gamma},\widetilde{\Phi}^r_{\psi,\Gamma},\widetilde{\Phi}^I_{\psi,\Gamma}, 
}
\]

\[
\xymatrix@R+0pc@C+0pc{
\breve{\Delta}_{\psi,\Gamma},\breve{\nabla}_{\psi,\Gamma},\breve{\Phi}_{\psi,\Gamma},\breve{\Delta}^+_{\psi,\Gamma},\breve{\nabla}^+_{\psi,\Gamma},\breve{\Delta}^\dagger_{\psi,\Gamma},\breve{\nabla}^\dagger_{\psi,\Gamma},\breve{\Phi}^r_{\psi,\Gamma},\breve{\Phi}^I_{\psi,\Gamma},	
}
\]

\[
\xymatrix@R+0pc@C+0pc{
{\Delta}_{\psi,\Gamma},{\nabla}_{\psi,\Gamma},{\Phi}_{\psi,\Gamma},{\Delta}^+_{\psi,\Gamma},{\nabla}^+_{\psi,\Gamma},{\Delta}^\dagger_{\psi,\Gamma},{\nabla}^\dagger_{\psi,\Gamma},{\Phi}^r_{\psi,\Gamma},{\Phi}^I_{\psi,\Gamma}.	
}
\]
Taking the product we have:
\[
\xymatrix@R+0pc@C+0pc{
\widetilde{\Phi}_{\psi,\Gamma,\circ},\widetilde{\Phi}^r_{\psi,\Gamma,\circ},\widetilde{\Phi}^I_{\psi,\Gamma,\circ},	
}
\]
\[
\xymatrix@R+0pc@C+0pc{
\breve{\Phi}_{\psi,\Gamma,\circ},\breve{\Phi}^r_{\psi,\Gamma,\circ},\breve{\Phi}^I_{\psi,\Gamma,\circ},	
}
\]
\[
\xymatrix@R+0pc@C+0pc{
{\Phi}_{\psi,\Gamma,\circ},{\Phi}^r_{\psi,\Gamma,\circ},{\Phi}^I_{\psi,\Gamma,\circ}.	
}
\]
They carry multi Frobenius action $\varphi_\Gamma$ and multi $\mathrm{Lie}_\Gamma:=\mathbb{Z}_p^{\times\Gamma}$ action. In our current situation after \cite{10CKZ} and \cite{10PZ} we consider the following $(\infty,1)$-categories of $(\infty,1)$-modules.\\
\end{definition}

\indent Meanwhile we have the corresponding Clausen-Scholze analytic stacks from \cite{10CS2}, therefore applying their construction we have:

\begin{definition}
Here we define the following products by using the solidified tensor product from \cite{10CS1} and \cite{10CS2}. Then we take solidified tensor product $\overset{\blacksquare}{\otimes}$ of any of the following
\[
\xymatrix@R+0pc@C+0pc{
\widetilde{\Delta}_{\psi,\Gamma},\widetilde{\nabla}_{\psi,\Gamma},\widetilde{\Phi}_{\psi,\Gamma},\widetilde{\Delta}^+_{\psi,\Gamma},\widetilde{\nabla}^+_{\psi,\Gamma},\widetilde{\Delta}^\dagger_{\psi,\Gamma},\widetilde{\nabla}^\dagger_{\psi,\Gamma},\widetilde{\Phi}^r_{\psi,\Gamma},\widetilde{\Phi}^I_{\psi,\Gamma}, 
}
\]

\[
\xymatrix@R+0pc@C+0pc{
\breve{\Delta}_{\psi,\Gamma},\breve{\nabla}_{\psi,\Gamma},\breve{\Phi}_{\psi,\Gamma},\breve{\Delta}^+_{\psi,\Gamma},\breve{\nabla}^+_{\psi,\Gamma},\breve{\Delta}^\dagger_{\psi,\Gamma},\breve{\nabla}^\dagger_{\psi,\Gamma},\breve{\Phi}^r_{\psi,\Gamma},\breve{\Phi}^I_{\psi,\Gamma},	
}
\]

\[
\xymatrix@R+0pc@C+0pc{
{\Delta}_{\psi,\Gamma},{\nabla}_{\psi,\Gamma},{\Phi}_{\psi,\Gamma},{\Delta}^+_{\psi,\Gamma},{\nabla}^+_{\psi,\Gamma},{\Delta}^\dagger_{\psi,\Gamma},{\nabla}^\dagger_{\psi,\Gamma},{\Phi}^r_{\psi,\Gamma},{\Phi}^I_{\psi,\Gamma},	
}
\]  	
with $\circ$. Then we have the notations:
\[
\xymatrix@R+0pc@C+0pc{
\widetilde{\Delta}_{\psi,\Gamma,\circ},\widetilde{\nabla}_{\psi,\Gamma,\circ},\widetilde{\Phi}_{\psi,\Gamma,\circ},\widetilde{\Delta}^+_{\psi,\Gamma,\circ},\widetilde{\nabla}^+_{\psi,\Gamma,\circ},\widetilde{\Delta}^\dagger_{\psi,\Gamma,\circ},\widetilde{\nabla}^\dagger_{\psi,\Gamma,\circ},\widetilde{\Phi}^r_{\psi,\Gamma,\circ},\widetilde{\Phi}^I_{\psi,\Gamma,\circ}, 
}
\]

\[
\xymatrix@R+0pc@C+0pc{
\breve{\Delta}_{\psi,\Gamma,\circ},\breve{\nabla}_{\psi,\Gamma,\circ},\breve{\Phi}_{\psi,\Gamma,\circ},\breve{\Delta}^+_{\psi,\Gamma,\circ},\breve{\nabla}^+_{\psi,\Gamma,\circ},\breve{\Delta}^\dagger_{\psi,\Gamma,\circ},\breve{\nabla}^\dagger_{\psi,\Gamma,\circ},\breve{\Phi}^r_{\psi,\Gamma,\circ},\breve{\Phi}^I_{\psi,\Gamma,\circ},	
}
\]

\[
\xymatrix@R+0pc@C+0pc{
{\Delta}_{\psi,\Gamma,\circ},{\nabla}_{\psi,\Gamma,\circ},{\Phi}_{\psi,\Gamma,\circ},{\Delta}^+_{\psi,\Gamma,\circ},{\nabla}^+_{\psi,\Gamma,\circ},{\Delta}^\dagger_{\psi,\Gamma,\circ},{\nabla}^\dagger_{\psi,\Gamma,\circ},{\Phi}^r_{\psi,\Gamma,\circ},{\Phi}^I_{\psi,\Gamma,\circ}.	
}
\]
\end{definition}

\begin{definition}
First we consider the Clausen-Scholze spectrum $\mathrm{Spec}^\mathrm{CS}(*)$ attached to any of those in the above from \cite{10CS2} by taking derived rational localization:
\begin{align}
\mathrm{Spec}^\mathrm{CS}\widetilde{\Delta}_{\psi,\Gamma,\circ},\mathrm{Spec}^\mathrm{CS}\widetilde{\nabla}_{\psi,\Gamma,\circ},\mathrm{Spec}^\mathrm{CS}\widetilde{\Phi}_{\psi,\Gamma,\circ},\mathrm{Spec}^\mathrm{CS}\widetilde{\Delta}^+_{\psi,\Gamma,\circ},\mathrm{Spec}^\mathrm{CS}\widetilde{\nabla}^+_{\psi,\Gamma,\circ},\\
\mathrm{Spec}^\mathrm{CS}\widetilde{\Delta}^\dagger_{\psi,\Gamma,\circ},\mathrm{Spec}^\mathrm{CS}\widetilde{\nabla}^\dagger_{\psi,\Gamma,\circ},\mathrm{Spec}^\mathrm{CS}\widetilde{\Phi}^r_{\psi,\Gamma,\circ},\mathrm{Spec}^\mathrm{CS}\widetilde{\Phi}^I_{\psi,\Gamma,\circ},	\\
\end{align}
\begin{align}
\mathrm{Spec}^\mathrm{CS}\breve{\Delta}_{\psi,\Gamma,\circ},\breve{\nabla}_{\psi,\Gamma,\circ},\mathrm{Spec}^\mathrm{CS}\breve{\Phi}_{\psi,\Gamma,\circ},\mathrm{Spec}^\mathrm{CS}\breve{\Delta}^+_{\psi,\Gamma,\circ},\mathrm{Spec}^\mathrm{CS}\breve{\nabla}^+_{\psi,\Gamma,\circ},\\
\mathrm{Spec}^\mathrm{CS}\breve{\Delta}^\dagger_{\psi,\Gamma,\circ},\mathrm{Spec}^\mathrm{CS}\breve{\nabla}^\dagger_{\psi,\Gamma,\circ},\mathrm{Spec}^\mathrm{CS}\breve{\Phi}^r_{\psi,\Gamma,\circ},\breve{\Phi}^I_{\psi,\Gamma,\circ},	\\
\end{align}
\begin{align}
\mathrm{Spec}^\mathrm{CS}{\Delta}_{\psi,\Gamma,\circ},\mathrm{Spec}^\mathrm{CS}{\nabla}_{\psi,\Gamma,\circ},\mathrm{Spec}^\mathrm{CS}{\Phi}_{\psi,\Gamma,\circ},\mathrm{Spec}^\mathrm{CS}{\Delta}^+_{\psi,\Gamma,\circ},\mathrm{Spec}^\mathrm{CS}{\nabla}^+_{\psi,\Gamma,\circ},\\
\mathrm{Spec}^\mathrm{CS}{\Delta}^\dagger_{\psi,\Gamma,\circ},\mathrm{Spec}^\mathrm{CS}{\nabla}^\dagger_{\psi,\Gamma,\circ},\mathrm{Spec}^\mathrm{CS}{\Phi}^r_{\psi,\Gamma,\circ},\mathrm{Spec}^\mathrm{CS}{\Phi}^I_{\psi,\Gamma,\circ}.	
\end{align}

Then we take the corresponding quotients by using the corresponding Frobenius operators:
\begin{align}
&\mathrm{Spec}^\mathrm{CS}\widetilde{\Delta}_{\psi,\Gamma,\circ}/\mathrm{Fro}^\mathbb{Z},\mathrm{Spec}^\mathrm{CS}\widetilde{\nabla}_{\psi,\Gamma,\circ}/\mathrm{Fro}^\mathbb{Z},\mathrm{Spec}^\mathrm{CS}\widetilde{\Phi}_{\psi,\Gamma,\circ}/\mathrm{Fro}^\mathbb{Z},\mathrm{Spec}^\mathrm{CS}\widetilde{\Delta}^+_{\psi,\Gamma,\circ}/\mathrm{Fro}^\mathbb{Z},\\
&\mathrm{Spec}^\mathrm{CS}\widetilde{\nabla}^+_{\psi,\Gamma,\circ}/\mathrm{Fro}^\mathbb{Z}, \mathrm{Spec}^\mathrm{CS}\widetilde{\Delta}^\dagger_{\psi,\Gamma,\circ}/\mathrm{Fro}^\mathbb{Z},\mathrm{Spec}^\mathrm{CS}\widetilde{\nabla}^\dagger_{\psi,\Gamma,\circ}/\mathrm{Fro}^\mathbb{Z},	\\
\end{align}
\begin{align}
&\mathrm{Spec}^\mathrm{CS}\breve{\Delta}_{\psi,\Gamma,\circ}/\mathrm{Fro}^\mathbb{Z},\breve{\nabla}_{\psi,\Gamma,\circ}/\mathrm{Fro}^\mathbb{Z},\mathrm{Spec}^\mathrm{CS}\breve{\Phi}_{\psi,\Gamma,\circ}/\mathrm{Fro}^\mathbb{Z},\mathrm{Spec}^\mathrm{CS}\breve{\Delta}^+_{\psi,\Gamma,\circ}/\mathrm{Fro}^\mathbb{Z},\\
&\mathrm{Spec}^\mathrm{CS}\breve{\nabla}^+_{\psi,\Gamma,\circ}/\mathrm{Fro}^\mathbb{Z}, \mathrm{Spec}^\mathrm{CS}\breve{\Delta}^\dagger_{\psi,\Gamma,\circ}/\mathrm{Fro}^\mathbb{Z},\mathrm{Spec}^\mathrm{CS}\breve{\nabla}^\dagger_{\psi,\Gamma,\circ}/\mathrm{Fro}^\mathbb{Z},	\\
\end{align}
\begin{align}
&\mathrm{Spec}^\mathrm{CS}{\Delta}_{\psi,\Gamma,\circ}/\mathrm{Fro}^\mathbb{Z},\mathrm{Spec}^\mathrm{CS}{\nabla}_{\psi,\Gamma,\circ}/\mathrm{Fro}^\mathbb{Z},\mathrm{Spec}^\mathrm{CS}{\Phi}_{\psi,\Gamma,\circ}/\mathrm{Fro}^\mathbb{Z},\mathrm{Spec}^\mathrm{CS}{\Delta}^+_{\psi,\Gamma,\circ}/\mathrm{Fro}^\mathbb{Z},\\
&\mathrm{Spec}^\mathrm{CS}{\nabla}^+_{\psi,\Gamma,\circ}/\mathrm{Fro}^\mathbb{Z}, \mathrm{Spec}^\mathrm{CS}{\Delta}^\dagger_{\psi,\Gamma,\circ}/\mathrm{Fro}^\mathbb{Z},\mathrm{Spec}^\mathrm{CS}{\nabla}^\dagger_{\psi,\Gamma,\circ}/\mathrm{Fro}^\mathbb{Z}.	
\end{align}
Here for those space with notations related to the radius and the corresponding interval we consider the total unions $\bigcap_r,\bigcup_I$ in order to achieve the whole spaces to achieve the analogues of the corresponding FF curves from \cite{10KL1}, \cite{10KL2}, \cite{10FF} for
\[
\xymatrix@R+0pc@C+0pc{
\underset{r}{\mathrm{homotopylimit}}~\mathrm{Spec}^\mathrm{CS}\widetilde{\Phi}^r_{\psi,\Gamma,\circ},\underset{I}{\mathrm{homotopycolimit}}~\mathrm{Spec}^\mathrm{CS}\widetilde{\Phi}^I_{\psi,\Gamma,\circ},	\\
}
\]
\[
\xymatrix@R+0pc@C+0pc{
\underset{r}{\mathrm{homotopylimit}}~\mathrm{Spec}^\mathrm{CS}\breve{\Phi}^r_{\psi,\Gamma,\circ},\underset{I}{\mathrm{homotopycolimit}}~\mathrm{Spec}^\mathrm{CS}\breve{\Phi}^I_{\psi,\Gamma,\circ},	\\
}
\]
\[
\xymatrix@R+0pc@C+0pc{
\underset{r}{\mathrm{homotopylimit}}~\mathrm{Spec}^\mathrm{CS}{\Phi}^r_{\psi,\Gamma,\circ},\underset{I}{\mathrm{homotopycolimit}}~\mathrm{Spec}^\mathrm{CS}{\Phi}^I_{\psi,\Gamma,\circ}.	
}
\]
\[ 
\xymatrix@R+0pc@C+0pc{
\underset{r}{\mathrm{homotopylimit}}~\mathrm{Spec}^\mathrm{CS}\widetilde{\Phi}^r_{\psi,\Gamma,\circ}/\mathrm{Fro}^\mathbb{Z},\underset{I}{\mathrm{homotopycolimit}}~\mathrm{Spec}^\mathrm{CS}\widetilde{\Phi}^I_{\psi,\Gamma,\circ}/\mathrm{Fro}^\mathbb{Z},	\\
}
\]
\[ 
\xymatrix@R+0pc@C+0pc{
\underset{r}{\mathrm{homotopylimit}}~\mathrm{Spec}^\mathrm{CS}\breve{\Phi}^r_{\psi,\Gamma,\circ}/\mathrm{Fro}^\mathbb{Z},\underset{I}{\mathrm{homotopycolimit}}~\breve{\Phi}^I_{\psi,\Gamma,\circ}/\mathrm{Fro}^\mathbb{Z},	\\
}
\]
\[ 
\xymatrix@R+0pc@C+0pc{
\underset{r}{\mathrm{homotopylimit}}~\mathrm{Spec}^\mathrm{CS}{\Phi}^r_{\psi,\Gamma,\circ}/\mathrm{Fro}^\mathbb{Z},\underset{I}{\mathrm{homotopycolimit}}~\mathrm{Spec}^\mathrm{CS}{\Phi}^I_{\psi,\Gamma,\circ}/\mathrm{Fro}^\mathbb{Z}.	
}
\]

\end{definition}

\

\begin{definition}
We then consider the corresponding quasisheaves of the corresponding condensed solid topological modules from \cite{10CS2}:
\begin{align}
\mathrm{Quasicoherentsheaves, Condensed}_{*}	
\end{align}
where $*$ is one of the following spaces:
\begin{align}
&\mathrm{Spec}^\mathrm{CS}\widetilde{\Delta}_{\psi,\Gamma,\circ}/\mathrm{Fro}^\mathbb{Z},\mathrm{Spec}^\mathrm{CS}\widetilde{\nabla}_{\psi,\Gamma,\circ}/\mathrm{Fro}^\mathbb{Z},\mathrm{Spec}^\mathrm{CS}\widetilde{\Phi}_{\psi,\Gamma,\circ}/\mathrm{Fro}^\mathbb{Z},\mathrm{Spec}^\mathrm{CS}\widetilde{\Delta}^+_{\psi,\Gamma,\circ}/\mathrm{Fro}^\mathbb{Z},\\
&\mathrm{Spec}^\mathrm{CS}\widetilde{\nabla}^+_{\psi,\Gamma,\circ}/\mathrm{Fro}^\mathbb{Z},\mathrm{Spec}^\mathrm{CS}\widetilde{\Delta}^\dagger_{\psi,\Gamma,\circ}/\mathrm{Fro}^\mathbb{Z},\mathrm{Spec}^\mathrm{CS}\widetilde{\nabla}^\dagger_{\psi,\Gamma,\circ}/\mathrm{Fro}^\mathbb{Z},	\\
\end{align}
\begin{align}
&\mathrm{Spec}^\mathrm{CS}\breve{\Delta}_{\psi,\Gamma,\circ}/\mathrm{Fro}^\mathbb{Z},\breve{\nabla}_{\psi,\Gamma,\circ}/\mathrm{Fro}^\mathbb{Z},\mathrm{Spec}^\mathrm{CS}\breve{\Phi}_{\psi,\Gamma,\circ}/\mathrm{Fro}^\mathbb{Z},\mathrm{Spec}^\mathrm{CS}\breve{\Delta}^+_{\psi,\Gamma,\circ}/\mathrm{Fro}^\mathbb{Z},\\
&\mathrm{Spec}^\mathrm{CS}\breve{\nabla}^+_{\psi,\Gamma,\circ}/\mathrm{Fro}^\mathbb{Z},\mathrm{Spec}^\mathrm{CS}\breve{\Delta}^\dagger_{\psi,\Gamma,\circ}/\mathrm{Fro}^\mathbb{Z},\mathrm{Spec}^\mathrm{CS}\breve{\nabla}^\dagger_{\psi,\Gamma,\circ}/\mathrm{Fro}^\mathbb{Z},	\\
\end{align}
\begin{align}
&\mathrm{Spec}^\mathrm{CS}{\Delta}_{\psi,\Gamma,\circ}/\mathrm{Fro}^\mathbb{Z},\mathrm{Spec}^\mathrm{CS}{\nabla}_{\psi,\Gamma,\circ}/\mathrm{Fro}^\mathbb{Z},\mathrm{Spec}^\mathrm{CS}{\Phi}_{\psi,\Gamma,\circ}/\mathrm{Fro}^\mathbb{Z},\mathrm{Spec}^\mathrm{CS}{\Delta}^+_{\psi,\Gamma,\circ}/\mathrm{Fro}^\mathbb{Z},\\
&\mathrm{Spec}^\mathrm{CS}{\nabla}^+_{\psi,\Gamma,\circ}/\mathrm{Fro}^\mathbb{Z}, \mathrm{Spec}^\mathrm{CS}{\Delta}^\dagger_{\psi,\Gamma,\circ}/\mathrm{Fro}^\mathbb{Z},\mathrm{Spec}^\mathrm{CS}{\nabla}^\dagger_{\psi,\Gamma,\circ}/\mathrm{Fro}^\mathbb{Z}.	
\end{align}
Here for those space with notations related to the radius and the corresponding interval we consider the total unions $\bigcap_r,\bigcup_I$ in order to achieve the whole spaces to achieve the analogues of the corresponding FF curves from \cite{10KL1}, \cite{10KL2}, \cite{10FF} for
\[
\xymatrix@R+0pc@C+0pc{
\underset{r}{\mathrm{homotopylimit}}~\mathrm{Spec}^\mathrm{CS}\widetilde{\Phi}^r_{\psi,\Gamma,\circ},\underset{I}{\mathrm{homotopycolimit}}~\mathrm{Spec}^\mathrm{CS}\widetilde{\Phi}^I_{\psi,\Gamma,\circ},	\\
}
\]
\[
\xymatrix@R+0pc@C+0pc{
\underset{r}{\mathrm{homotopylimit}}~\mathrm{Spec}^\mathrm{CS}\breve{\Phi}^r_{\psi,\Gamma,\circ},\underset{I}{\mathrm{homotopycolimit}}~\mathrm{Spec}^\mathrm{CS}\breve{\Phi}^I_{\psi,\Gamma,\circ},	\\
}
\]
\[
\xymatrix@R+0pc@C+0pc{
\underset{r}{\mathrm{homotopylimit}}~\mathrm{Spec}^\mathrm{CS}{\Phi}^r_{\psi,\Gamma,\circ},\underset{I}{\mathrm{homotopycolimit}}~\mathrm{Spec}^\mathrm{CS}{\Phi}^I_{\psi,\Gamma,\circ}.	
}
\]
\[ 
\xymatrix@R+0pc@C+0pc{
\underset{r}{\mathrm{homotopylimit}}~\mathrm{Spec}^\mathrm{CS}\widetilde{\Phi}^r_{\psi,\Gamma,\circ}/\mathrm{Fro}^\mathbb{Z},\underset{I}{\mathrm{homotopycolimit}}~\mathrm{Spec}^\mathrm{CS}\widetilde{\Phi}^I_{\psi,\Gamma,\circ}/\mathrm{Fro}^\mathbb{Z},	\\
}
\]
\[ 
\xymatrix@R+0pc@C+0pc{
\underset{r}{\mathrm{homotopylimit}}~\mathrm{Spec}^\mathrm{CS}\breve{\Phi}^r_{\psi,\Gamma,\circ}/\mathrm{Fro}^\mathbb{Z},\underset{I}{\mathrm{homotopycolimit}}~\breve{\Phi}^I_{\psi,\Gamma,\circ}/\mathrm{Fro}^\mathbb{Z},	\\
}
\]
\[ 
\xymatrix@R+0pc@C+0pc{
\underset{r}{\mathrm{homotopylimit}}~\mathrm{Spec}^\mathrm{CS}{\Phi}^r_{\psi,\Gamma,\circ}/\mathrm{Fro}^\mathbb{Z},\underset{I}{\mathrm{homotopycolimit}}~\mathrm{Spec}^\mathrm{CS}{\Phi}^I_{\psi,\Gamma,\circ}/\mathrm{Fro}^\mathbb{Z}.	
}
\]

\end{definition}

\

\begin{proposition}
There is a well-defined functor from the $\infty$-category 
\begin{align}
\mathrm{Quasicoherentpresheaves,Condensed}_{*}	
\end{align}
where $*$ is one of the following spaces:
\begin{align}
&\mathrm{Spec}^\mathrm{CS}\widetilde{\Phi}_{\psi,\Gamma,\circ}/\mathrm{Fro}^\mathbb{Z},	\\
\end{align}
\begin{align}
&\mathrm{Spec}^\mathrm{CS}\breve{\Phi}_{\psi,\Gamma,\circ}/\mathrm{Fro}^\mathbb{Z},	\\
\end{align}
\begin{align}
&\mathrm{Spec}^\mathrm{CS}{\Phi}_{\psi,\Gamma,\circ}/\mathrm{Fro}^\mathbb{Z},	
\end{align}
to the $\infty$-category of $\mathrm{Fro}$-equivariant quasicoherent presheaves over similar spaces above correspondingly without the $\mathrm{Fro}$-quotients, and to the $\infty$-category of $\mathrm{Fro}$-equivariant quasicoherent modules over global sections of the structure $\infty$-sheaves of the similar spaces above correspondingly without the $\mathrm{Fro}$-quotients. Here for those space without notation related to the radius and the corresponding interval we consider the total unions $\bigcap_r,\bigcup_I$ in order to achieve the whole spaces to achieve the analogues of the corresponding FF curves from \cite{10KL1}, \cite{10KL2}, \cite{10FF} for
\[
\xymatrix@R+0pc@C+0pc{
\underset{r}{\mathrm{homotopylimit}}~\mathrm{Spec}^\mathrm{CS}\widetilde{\Phi}^r_{\psi,\Gamma,\circ},\underset{I}{\mathrm{homotopycolimit}}~\mathrm{Spec}^\mathrm{CS}\widetilde{\Phi}^I_{\psi,\Gamma,\circ},	\\
}
\]
\[
\xymatrix@R+0pc@C+0pc{
\underset{r}{\mathrm{homotopylimit}}~\mathrm{Spec}^\mathrm{CS}\breve{\Phi}^r_{\psi,\Gamma,\circ},\underset{I}{\mathrm{homotopycolimit}}~\mathrm{Spec}^\mathrm{CS}\breve{\Phi}^I_{\psi,\Gamma,\circ},	\\
}
\]
\[
\xymatrix@R+0pc@C+0pc{
\underset{r}{\mathrm{homotopylimit}}~\mathrm{Spec}^\mathrm{CS}{\Phi}^r_{\psi,\Gamma,\circ},\underset{I}{\mathrm{homotopycolimit}}~\mathrm{Spec}^\mathrm{CS}{\Phi}^I_{\psi,\Gamma,\circ}.	
}
\]
\[ 
\xymatrix@R+0pc@C+0pc{
\underset{r}{\mathrm{homotopylimit}}~\mathrm{Spec}^\mathrm{CS}\widetilde{\Phi}^r_{\psi,\Gamma,\circ}/\mathrm{Fro}^\mathbb{Z},\underset{I}{\mathrm{homotopycolimit}}~\mathrm{Spec}^\mathrm{CS}\widetilde{\Phi}^I_{\psi,\Gamma,\circ}/\mathrm{Fro}^\mathbb{Z},	\\
}
\]
\[ 
\xymatrix@R+0pc@C+0pc{
\underset{r}{\mathrm{homotopylimit}}~\mathrm{Spec}^\mathrm{CS}\breve{\Phi}^r_{\psi,\Gamma,\circ}/\mathrm{Fro}^\mathbb{Z},\underset{I}{\mathrm{homotopycolimit}}~\breve{\Phi}^I_{\psi,\Gamma,\circ}/\mathrm{Fro}^\mathbb{Z},	\\
}
\]
\[ 
\xymatrix@R+0pc@C+0pc{
\underset{r}{\mathrm{homotopylimit}}~\mathrm{Spec}^\mathrm{CS}{\Phi}^r_{\psi,\Gamma,\circ}/\mathrm{Fro}^\mathbb{Z},\underset{I}{\mathrm{homotopycolimit}}~\mathrm{Spec}^\mathrm{CS}{\Phi}^I_{\psi,\Gamma,\circ}/\mathrm{Fro}^\mathbb{Z}.	
}
\]	
In this situation we will have the target category being family parametrized by $r$ or $I$ in compatible glueing sense as in \cite[Definition 5.4.10]{10KL2}. In this situation for modules parametrized by the intervals we have the equivalence of $\infty$-categories by using \cite[Proposition 13.8]{10CS2}. Here the corresponding quasicoherent Frobenius modules are defined to be the corresponding homotopy colimits and limits of Frobenius modules:
\begin{align}
\underset{r}{\mathrm{homotopycolimit}}~M_r,\\
\underset{I}{\mathrm{homotopylimit}}~M_I,	
\end{align}
where each $M_r$ is a Frobenius-equivariant module over the period ring with respect to some radius $r$ while each $M_I$ is a Frobenius-equivariant module over the period ring with respect to some interval $I$.\\
\end{proposition}

\

\begin{definition}
We then consider the corresponding quasisheaves of perfect complexes of the corresponding condensed solid topological modules from \cite{10CS2}:
\begin{align}
\mathrm{Quasicoherentsheaves, Perfectcomplex, Condensed}_{*}	
\end{align}
where $*$ is one of the following spaces:
\begin{align}
&\mathrm{Spec}^\mathrm{CS}\widetilde{\Delta}_{\psi,\Gamma,\circ}/\mathrm{Fro}^\mathbb{Z},\mathrm{Spec}^\mathrm{CS}\widetilde{\nabla}_{\psi,\Gamma,\circ}/\mathrm{Fro}^\mathbb{Z},\mathrm{Spec}^\mathrm{CS}\widetilde{\Phi}_{\psi,\Gamma,\circ}/\mathrm{Fro}^\mathbb{Z},\mathrm{Spec}^\mathrm{CS}\widetilde{\Delta}^+_{\psi,\Gamma,\circ}/\mathrm{Fro}^\mathbb{Z},\\
&\mathrm{Spec}^\mathrm{CS}\widetilde{\nabla}^+_{\psi,\Gamma,\circ}/\mathrm{Fro}^\mathbb{Z},\mathrm{Spec}^\mathrm{CS}\widetilde{\Delta}^\dagger_{\psi,\Gamma,\circ}/\mathrm{Fro}^\mathbb{Z},\mathrm{Spec}^\mathrm{CS}\widetilde{\nabla}^\dagger_{\psi,\Gamma,\circ}/\mathrm{Fro}^\mathbb{Z},	\\
\end{align}
\begin{align}
&\mathrm{Spec}^\mathrm{CS}\breve{\Delta}_{\psi,\Gamma,\circ}/\mathrm{Fro}^\mathbb{Z},\breve{\nabla}_{\psi,\Gamma,\circ}/\mathrm{Fro}^\mathbb{Z},\mathrm{Spec}^\mathrm{CS}\breve{\Phi}_{\psi,\Gamma,\circ}/\mathrm{Fro}^\mathbb{Z},\mathrm{Spec}^\mathrm{CS}\breve{\Delta}^+_{\psi,\Gamma,\circ}/\mathrm{Fro}^\mathbb{Z},\\
&\mathrm{Spec}^\mathrm{CS}\breve{\nabla}^+_{\psi,\Gamma,\circ}/\mathrm{Fro}^\mathbb{Z},\mathrm{Spec}^\mathrm{CS}\breve{\Delta}^\dagger_{\psi,\Gamma,\circ}/\mathrm{Fro}^\mathbb{Z},\mathrm{Spec}^\mathrm{CS}\breve{\nabla}^\dagger_{\psi,\Gamma,\circ}/\mathrm{Fro}^\mathbb{Z},	\\
\end{align}
\begin{align}
&\mathrm{Spec}^\mathrm{CS}{\Delta}_{\psi,\Gamma,\circ}/\mathrm{Fro}^\mathbb{Z},\mathrm{Spec}^\mathrm{CS}{\nabla}_{\psi,\Gamma,\circ}/\mathrm{Fro}^\mathbb{Z},\mathrm{Spec}^\mathrm{CS}{\Phi}_{\psi,\Gamma,\circ}/\mathrm{Fro}^\mathbb{Z},\mathrm{Spec}^\mathrm{CS}{\Delta}^+_{\psi,\Gamma,\circ}/\mathrm{Fro}^\mathbb{Z},\\
&\mathrm{Spec}^\mathrm{CS}{\nabla}^+_{\psi,\Gamma,\circ}/\mathrm{Fro}^\mathbb{Z}, \mathrm{Spec}^\mathrm{CS}{\Delta}^\dagger_{\psi,\Gamma,\circ}/\mathrm{Fro}^\mathbb{Z},\mathrm{Spec}^\mathrm{CS}{\nabla}^\dagger_{\psi,\Gamma,\circ}/\mathrm{Fro}^\mathbb{Z}.	
\end{align}
Here for those space with notations related to the radius and the corresponding interval we consider the total unions $\bigcap_r,\bigcup_I$ in order to achieve the whole spaces to achieve the analogues of the corresponding FF curves from \cite{10KL1}, \cite{10KL2}, \cite{10FF} for
\[
\xymatrix@R+0pc@C+0pc{
\underset{r}{\mathrm{homotopylimit}}~\mathrm{Spec}^\mathrm{CS}\widetilde{\Phi}^r_{\psi,\Gamma,\circ},\underset{I}{\mathrm{homotopycolimit}}~\mathrm{Spec}^\mathrm{CS}\widetilde{\Phi}^I_{\psi,\Gamma,\circ},	\\
}
\]
\[
\xymatrix@R+0pc@C+0pc{
\underset{r}{\mathrm{homotopylimit}}~\mathrm{Spec}^\mathrm{CS}\breve{\Phi}^r_{\psi,\Gamma,\circ},\underset{I}{\mathrm{homotopycolimit}}~\mathrm{Spec}^\mathrm{CS}\breve{\Phi}^I_{\psi,\Gamma,\circ},	\\
}
\]
\[
\xymatrix@R+0pc@C+0pc{
\underset{r}{\mathrm{homotopylimit}}~\mathrm{Spec}^\mathrm{CS}{\Phi}^r_{\psi,\Gamma,\circ},\underset{I}{\mathrm{homotopycolimit}}~\mathrm{Spec}^\mathrm{CS}{\Phi}^I_{\psi,\Gamma,\circ}.	
}
\]
\[ 
\xymatrix@R+0pc@C+0pc{
\underset{r}{\mathrm{homotopylimit}}~\mathrm{Spec}^\mathrm{CS}\widetilde{\Phi}^r_{\psi,\Gamma,\circ}/\mathrm{Fro}^\mathbb{Z},\underset{I}{\mathrm{homotopycolimit}}~\mathrm{Spec}^\mathrm{CS}\widetilde{\Phi}^I_{\psi,\Gamma,\circ}/\mathrm{Fro}^\mathbb{Z},	\\
}
\]
\[ 
\xymatrix@R+0pc@C+0pc{
\underset{r}{\mathrm{homotopylimit}}~\mathrm{Spec}^\mathrm{CS}\breve{\Phi}^r_{\psi,\Gamma,\circ}/\mathrm{Fro}^\mathbb{Z},\underset{I}{\mathrm{homotopycolimit}}~\breve{\Phi}^I_{\psi,\Gamma,\circ}/\mathrm{Fro}^\mathbb{Z},	\\
}
\]
\[ 
\xymatrix@R+0pc@C+0pc{
\underset{r}{\mathrm{homotopylimit}}~\mathrm{Spec}^\mathrm{CS}{\Phi}^r_{\psi,\Gamma,\circ}/\mathrm{Fro}^\mathbb{Z},\underset{I}{\mathrm{homotopycolimit}}~\mathrm{Spec}^\mathrm{CS}{\Phi}^I_{\psi,\Gamma,\circ}/\mathrm{Fro}^\mathbb{Z}.	
}
\]

\end{definition}

\begin{proposition}
There is a well-defined functor from the $\infty$-category 
\begin{align}
\mathrm{Quasicoherentpresheaves,Perfectcomplex,Condensed}_{*}	
\end{align}
where $*$ is one of the following spaces:
\begin{align}
&\mathrm{Spec}^\mathrm{CS}\widetilde{\Phi}_{\psi,\Gamma,\circ}/\mathrm{Fro}^\mathbb{Z},	\\
\end{align}
\begin{align}
&\mathrm{Spec}^\mathrm{CS}\breve{\Phi}_{\psi,\Gamma,\circ}/\mathrm{Fro}^\mathbb{Z},	\\
\end{align}
\begin{align}
&\mathrm{Spec}^\mathrm{CS}{\Phi}_{\psi,\Gamma,\circ}/\mathrm{Fro}^\mathbb{Z},	
\end{align}
to the $\infty$-category of $\mathrm{Fro}$-equivariant quasicoherent presheaves over similar spaces above correspondingly without the $\mathrm{Fro}$-quotients, and to the $\infty$-category of $\mathrm{Fro}$-equivariant quasicoherent modules over global sections of the structure $\infty$-sheaves of the similar spaces above correspondingly without the $\mathrm{Fro}$-quotients. Here for those space without notation related to the radius and the corresponding interval we consider the total unions $\bigcap_r,\bigcup_I$ in order to achieve the whole spaces to achieve the analogues of the corresponding FF curves from \cite{10KL1}, \cite{10KL2}, \cite{10FF} for
\[
\xymatrix@R+0pc@C+0pc{
\underset{r}{\mathrm{homotopylimit}}~\mathrm{Spec}^\mathrm{CS}\widetilde{\Phi}^r_{\psi,\Gamma,\circ},\underset{I}{\mathrm{homotopycolimit}}~\mathrm{Spec}^\mathrm{CS}\widetilde{\Phi}^I_{\psi,\Gamma,\circ},	\\
}
\]
\[
\xymatrix@R+0pc@C+0pc{
\underset{r}{\mathrm{homotopylimit}}~\mathrm{Spec}^\mathrm{CS}\breve{\Phi}^r_{\psi,\Gamma,\circ},\underset{I}{\mathrm{homotopycolimit}}~\mathrm{Spec}^\mathrm{CS}\breve{\Phi}^I_{\psi,\Gamma,\circ},	\\
}
\]
\[
\xymatrix@R+0pc@C+0pc{
\underset{r}{\mathrm{homotopylimit}}~\mathrm{Spec}^\mathrm{CS}{\Phi}^r_{\psi,\Gamma,\circ},\underset{I}{\mathrm{homotopycolimit}}~\mathrm{Spec}^\mathrm{CS}{\Phi}^I_{\psi,\Gamma,\circ}.	
}
\]
\[ 
\xymatrix@R+0pc@C+0pc{
\underset{r}{\mathrm{homotopylimit}}~\mathrm{Spec}^\mathrm{CS}\widetilde{\Phi}^r_{\psi,\Gamma,\circ}/\mathrm{Fro}^\mathbb{Z},\underset{I}{\mathrm{homotopycolimit}}~\mathrm{Spec}^\mathrm{CS}\widetilde{\Phi}^I_{\psi,\Gamma,\circ}/\mathrm{Fro}^\mathbb{Z},	\\
}
\]
\[ 
\xymatrix@R+0pc@C+0pc{
\underset{r}{\mathrm{homotopylimit}}~\mathrm{Spec}^\mathrm{CS}\breve{\Phi}^r_{\psi,\Gamma,\circ}/\mathrm{Fro}^\mathbb{Z},\underset{I}{\mathrm{homotopycolimit}}~\breve{\Phi}^I_{\psi,\Gamma,\circ}/\mathrm{Fro}^\mathbb{Z},	\\
}
\]
\[ 
\xymatrix@R+0pc@C+0pc{
\underset{r}{\mathrm{homotopylimit}}~\mathrm{Spec}^\mathrm{CS}{\Phi}^r_{\psi,\Gamma,\circ}/\mathrm{Fro}^\mathbb{Z},\underset{I}{\mathrm{homotopycolimit}}~\mathrm{Spec}^\mathrm{CS}{\Phi}^I_{\psi,\Gamma,\circ}/\mathrm{Fro}^\mathbb{Z}.	
}
\]	
In this situation we will have the target category being family parametrized by $r$ or $I$ in compatible glueing sense as in \cite[Definition 5.4.10]{10KL2}. Here the corresponding quasicoherent Frobenius modules are defined to be the corresponding homotopy colimits and limits of Frobenius modules:
\begin{align}
\underset{r}{\mathrm{homotopycolimit}}~M_r,\\
\underset{I}{\mathrm{homotopylimit}}~M_I,	
\end{align}
where each $M_r$ is a Frobenius-equivariant module over the period ring with respect to some radius $r$ while each $M_I$ is a Frobenius-equivariant module over the period ring with respect to some interval $I$.\\
\end{proposition}

\subsubsection{Frobenius Quasicoherent Modules III: Deformation in $(\infty,1)$-Analytic Spaces}

\begin{definition}
Let $\psi$ be a toric tower over $\mathbb{Q}_p$ as in \cite[Chapter 7]{10KL2} with base $\mathbb{Q}_p\left<X_1^{\pm 1},...,X_k^{\pm 1}\right>$. Then from \cite{10KL1} and \cite[Definition 5.2.1]{10KL2} we have the following class of Kedlaya-Liu rings (with the following replacement: $\Delta$ stands for $A$, $\nabla$ stands for $B$, while $\Phi$ stands for $C$) by taking product in the sense of self $\Gamma$-th power:

\[
\xymatrix@R+0pc@C+0pc{
\widetilde{\Delta}_{\psi,\Gamma},\widetilde{\nabla}_{\psi,\Gamma},\widetilde{\Phi}_{\psi,\Gamma},\widetilde{\Delta}^+_{\psi,\Gamma},\widetilde{\nabla}^+_{\psi,\Gamma},\widetilde{\Delta}^\dagger_{\psi,\Gamma},\widetilde{\nabla}^\dagger_{\psi,\Gamma},\widetilde{\Phi}^r_{\psi,\Gamma},\widetilde{\Phi}^I_{\psi,\Gamma}, 
}
\]

\[
\xymatrix@R+0pc@C+0pc{
\breve{\Delta}_{\psi,\Gamma},\breve{\nabla}_{\psi,\Gamma},\breve{\Phi}_{\psi,\Gamma},\breve{\Delta}^+_{\psi,\Gamma},\breve{\nabla}^+_{\psi,\Gamma},\breve{\Delta}^\dagger_{\psi,\Gamma},\breve{\nabla}^\dagger_{\psi,\Gamma},\breve{\Phi}^r_{\psi,\Gamma},\breve{\Phi}^I_{\psi,\Gamma},	
}
\]

\[
\xymatrix@R+0pc@C+0pc{
{\Delta}_{\psi,\Gamma},{\nabla}_{\psi,\Gamma},{\Phi}_{\psi,\Gamma},{\Delta}^+_{\psi,\Gamma},{\nabla}^+_{\psi,\Gamma},{\Delta}^\dagger_{\psi,\Gamma},{\nabla}^\dagger_{\psi,\Gamma},{\Phi}^r_{\psi,\Gamma},{\Phi}^I_{\psi,\Gamma}.	
}
\]
   
Taking the product we have:
\[
\xymatrix@R+0pc@C+0pc{
\widetilde{\Phi}_{\psi,\Gamma,X_\square},\widetilde{\Phi}^r_{\psi,\Gamma,X_\square},\widetilde{\Phi}^I_{\psi,\Gamma,X_\square},	
}
\]
\[
\xymatrix@R+0pc@C+0pc{
\breve{\Phi}_{\psi,\Gamma,X_\square},\breve{\Phi}^r_{\psi,\Gamma,X_\square},\breve{\Phi}^I_{\psi,\Gamma,X_\square},	
}
\]
\[
\xymatrix@R+0pc@C+0pc{
{\Phi}_{\psi,\Gamma,X_\square},{\Phi}^r_{\psi,\Gamma,X_\square},{\Phi}^I_{\psi,\Gamma,X_\square}.	
}
\]
They carry multi Frobenius action $\varphi_\Gamma$ and multi $\mathrm{Lie}_\Gamma:=\mathbb{Z}_p^{\times\Gamma}$ action. In our current situation after \cite{10CKZ} and \cite{10PZ} we consider the following $(\infty,1)$-categories of $(\infty,1)$-modules.\\
\end{definition}

\indent Meanwhile we have the corresponding Clausen-Scholze analytic stacks from \cite{10CS2}, therefore applying their construction we have:

\begin{definition}
Here we define the following products by using the solidified tensor product from \cite{10CS1} and \cite{10CS2}. Then we take solidified tensor product $\overset{\blacksquare}{\otimes}$ of any of the following
\[
\xymatrix@R+0pc@C+0pc{
\widetilde{\Delta}_{\psi,\Gamma},\widetilde{\nabla}_{\psi,\Gamma},\widetilde{\Phi}_{\psi,\Gamma},\widetilde{\Delta}^+_{\psi,\Gamma},\widetilde{\nabla}^+_{\psi,\Gamma},\widetilde{\Delta}^\dagger_{\psi,\Gamma},\widetilde{\nabla}^\dagger_{\psi,\Gamma},\widetilde{\Phi}^r_{\psi,\Gamma},\widetilde{\Phi}^I_{\psi,\Gamma}, 
}
\]

\[
\xymatrix@R+0pc@C+0pc{
\breve{\Delta}_{\psi,\Gamma},\breve{\nabla}_{\psi,\Gamma},\breve{\Phi}_{\psi,\Gamma},\breve{\Delta}^+_{\psi,\Gamma},\breve{\nabla}^+_{\psi,\Gamma},\breve{\Delta}^\dagger_{\psi,\Gamma},\breve{\nabla}^\dagger_{\psi,\Gamma},\breve{\Phi}^r_{\psi,\Gamma},\breve{\Phi}^I_{\psi,\Gamma},	
}
\]

\[
\xymatrix@R+0pc@C+0pc{
{\Delta}_{\psi,\Gamma},{\nabla}_{\psi,\Gamma},{\Phi}_{\psi,\Gamma},{\Delta}^+_{\psi,\Gamma},{\nabla}^+_{\psi,\Gamma},{\Delta}^\dagger_{\psi,\Gamma},{\nabla}^\dagger_{\psi,\Gamma},{\Phi}^r_{\psi,\Gamma},{\Phi}^I_{\psi,\Gamma},	
}
\]  	
with $X_\square$. Then we have the notations:
\[
\xymatrix@R+0pc@C+0pc{
\widetilde{\Delta}_{\psi,\Gamma,X_\square},\widetilde{\nabla}_{\psi,\Gamma,X_\square},\widetilde{\Phi}_{\psi,\Gamma,X_\square},\widetilde{\Delta}^+_{\psi,\Gamma,X_\square},\widetilde{\nabla}^+_{\psi,\Gamma,X_\square},\widetilde{\Delta}^\dagger_{\psi,\Gamma,X_\square},\widetilde{\nabla}^\dagger_{\psi,\Gamma,X_\square},\widetilde{\Phi}^r_{\psi,\Gamma,X_\square},\widetilde{\Phi}^I_{\psi,\Gamma,X_\square}, 
}
\]

\[
\xymatrix@R+0pc@C+0pc{
\breve{\Delta}_{\psi,\Gamma,X_\square},\breve{\nabla}_{\psi,\Gamma,X_\square},\breve{\Phi}_{\psi,\Gamma,X_\square},\breve{\Delta}^+_{\psi,\Gamma,X_\square},\breve{\nabla}^+_{\psi,\Gamma,X_\square},\breve{\Delta}^\dagger_{\psi,\Gamma,X_\square},\breve{\nabla}^\dagger_{\psi,\Gamma,X_\square},\breve{\Phi}^r_{\psi,\Gamma,X_\square},\breve{\Phi}^I_{\psi,\Gamma,X_\square},	
}
\]

\[
\xymatrix@R+0pc@C+0pc{
{\Delta}_{\psi,\Gamma,X_\square},{\nabla}_{\psi,\Gamma,X_\square},{\Phi}_{\psi,\Gamma,X_\square},{\Delta}^+_{\psi,\Gamma,X_\square},{\nabla}^+_{\psi,\Gamma,X_\square},{\Delta}^\dagger_{\psi,\Gamma,X_\square},{\nabla}^\dagger_{\psi,\Gamma,X_\square},{\Phi}^r_{\psi,\Gamma,X_\square},{\Phi}^I_{\psi,\Gamma,X_\square}.	
}
\]
\end{definition}

\begin{definition}
First we consider the Clausen-Scholze spectrum $\mathrm{Spec}^\mathrm{CS}(*)$ attached to any of those in the above from \cite{10CS2} by taking derived rational localization:
\begin{align}
\mathrm{Spec}^\mathrm{CS}\widetilde{\Delta}_{\psi,\Gamma,X_\square},\mathrm{Spec}^\mathrm{CS}\widetilde{\nabla}_{\psi,\Gamma,X_\square},\mathrm{Spec}^\mathrm{CS}\widetilde{\Phi}_{\psi,\Gamma,X_\square},\mathrm{Spec}^\mathrm{CS}\widetilde{\Delta}^+_{\psi,\Gamma,X_\square},\mathrm{Spec}^\mathrm{CS}\widetilde{\nabla}^+_{\psi,\Gamma,X_\square},\\
\mathrm{Spec}^\mathrm{CS}\widetilde{\Delta}^\dagger_{\psi,\Gamma,X_\square},\mathrm{Spec}^\mathrm{CS}\widetilde{\nabla}^\dagger_{\psi,\Gamma,X_\square},\mathrm{Spec}^\mathrm{CS}\widetilde{\Phi}^r_{\psi,\Gamma,X_\square},\mathrm{Spec}^\mathrm{CS}\widetilde{\Phi}^I_{\psi,\Gamma,X_\square},	\\
\end{align}
\begin{align}
\mathrm{Spec}^\mathrm{CS}\breve{\Delta}_{\psi,\Gamma,X_\square},\breve{\nabla}_{\psi,\Gamma,X_\square},\mathrm{Spec}^\mathrm{CS}\breve{\Phi}_{\psi,\Gamma,X_\square},\mathrm{Spec}^\mathrm{CS}\breve{\Delta}^+_{\psi,\Gamma,X_\square},\mathrm{Spec}^\mathrm{CS}\breve{\nabla}^+_{\psi,\Gamma,X_\square},\\
\mathrm{Spec}^\mathrm{CS}\breve{\Delta}^\dagger_{\psi,\Gamma,X_\square},\mathrm{Spec}^\mathrm{CS}\breve{\nabla}^\dagger_{\psi,\Gamma,X_\square},\mathrm{Spec}^\mathrm{CS}\breve{\Phi}^r_{\psi,\Gamma,X_\square},\breve{\Phi}^I_{\psi,\Gamma,X_\square},	\\
\end{align}
\begin{align}
\mathrm{Spec}^\mathrm{CS}{\Delta}_{\psi,\Gamma,X_\square},\mathrm{Spec}^\mathrm{CS}{\nabla}_{\psi,\Gamma,X_\square},\mathrm{Spec}^\mathrm{CS}{\Phi}_{\psi,\Gamma,X_\square},\mathrm{Spec}^\mathrm{CS}{\Delta}^+_{\psi,\Gamma,X_\square},\mathrm{Spec}^\mathrm{CS}{\nabla}^+_{\psi,\Gamma,X_\square},\\
\mathrm{Spec}^\mathrm{CS}{\Delta}^\dagger_{\psi,\Gamma,X_\square},\mathrm{Spec}^\mathrm{CS}{\nabla}^\dagger_{\psi,\Gamma,X_\square},\mathrm{Spec}^\mathrm{CS}{\Phi}^r_{\psi,\Gamma,X_\square},\mathrm{Spec}^\mathrm{CS}{\Phi}^I_{\psi,\Gamma,X_\square}.	
\end{align}

Then we take the corresponding quotients by using the corresponding Frobenius operators:
\begin{align}
&\mathrm{Spec}^\mathrm{CS}\widetilde{\Delta}_{\psi,\Gamma,X_\square}/\mathrm{Fro}^\mathbb{Z},\mathrm{Spec}^\mathrm{CS}\widetilde{\nabla}_{\psi,\Gamma,X_\square}/\mathrm{Fro}^\mathbb{Z},\mathrm{Spec}^\mathrm{CS}\widetilde{\Phi}_{\psi,\Gamma,X_\square}/\mathrm{Fro}^\mathbb{Z},\mathrm{Spec}^\mathrm{CS}\widetilde{\Delta}^+_{\psi,\Gamma,X_\square}/\mathrm{Fro}^\mathbb{Z},\\
&\mathrm{Spec}^\mathrm{CS}\widetilde{\nabla}^+_{\psi,\Gamma,X_\square}/\mathrm{Fro}^\mathbb{Z}, \mathrm{Spec}^\mathrm{CS}\widetilde{\Delta}^\dagger_{\psi,\Gamma,X_\square}/\mathrm{Fro}^\mathbb{Z},\mathrm{Spec}^\mathrm{CS}\widetilde{\nabla}^\dagger_{\psi,\Gamma,X_\square}/\mathrm{Fro}^\mathbb{Z},	\\
\end{align}
\begin{align}
&\mathrm{Spec}^\mathrm{CS}\breve{\Delta}_{\psi,\Gamma,X_\square}/\mathrm{Fro}^\mathbb{Z},\breve{\nabla}_{\psi,\Gamma,X_\square}/\mathrm{Fro}^\mathbb{Z},\mathrm{Spec}^\mathrm{CS}\breve{\Phi}_{\psi,\Gamma,X_\square}/\mathrm{Fro}^\mathbb{Z},\mathrm{Spec}^\mathrm{CS}\breve{\Delta}^+_{\psi,\Gamma,X_\square}/\mathrm{Fro}^\mathbb{Z},\\
&\mathrm{Spec}^\mathrm{CS}\breve{\nabla}^+_{\psi,\Gamma,X_\square}/\mathrm{Fro}^\mathbb{Z}, \mathrm{Spec}^\mathrm{CS}\breve{\Delta}^\dagger_{\psi,\Gamma,X_\square}/\mathrm{Fro}^\mathbb{Z},\mathrm{Spec}^\mathrm{CS}\breve{\nabla}^\dagger_{\psi,\Gamma,X_\square}/\mathrm{Fro}^\mathbb{Z},	\\
\end{align}
\begin{align}
&\mathrm{Spec}^\mathrm{CS}{\Delta}_{\psi,\Gamma,X_\square}/\mathrm{Fro}^\mathbb{Z},\mathrm{Spec}^\mathrm{CS}{\nabla}_{\psi,\Gamma,X_\square}/\mathrm{Fro}^\mathbb{Z},\mathrm{Spec}^\mathrm{CS}{\Phi}_{\psi,\Gamma,X_\square}/\mathrm{Fro}^\mathbb{Z},\mathrm{Spec}^\mathrm{CS}{\Delta}^+_{\psi,\Gamma,X_\square}/\mathrm{Fro}^\mathbb{Z},\\
&\mathrm{Spec}^\mathrm{CS}{\nabla}^+_{\psi,\Gamma,X_\square}/\mathrm{Fro}^\mathbb{Z}, \mathrm{Spec}^\mathrm{CS}{\Delta}^\dagger_{\psi,\Gamma,X_\square}/\mathrm{Fro}^\mathbb{Z},\mathrm{Spec}^\mathrm{CS}{\nabla}^\dagger_{\psi,\Gamma,X_\square}/\mathrm{Fro}^\mathbb{Z}.	
\end{align}
Here for those space with notations related to the radius and the corresponding interval we consider the total unions $\bigcap_r,\bigcup_I$ in order to achieve the whole spaces to achieve the analogues of the corresponding FF curves from \cite{10KL1}, \cite{10KL2}, \cite{10FF} for
\[
\xymatrix@R+0pc@C+0pc{
\underset{r}{\mathrm{homotopylimit}}~\mathrm{Spec}^\mathrm{CS}\widetilde{\Phi}^r_{\psi,\Gamma,X_\square},\underset{I}{\mathrm{homotopycolimit}}~\mathrm{Spec}^\mathrm{CS}\widetilde{\Phi}^I_{\psi,\Gamma,X_\square},	\\
}
\]
\[
\xymatrix@R+0pc@C+0pc{
\underset{r}{\mathrm{homotopylimit}}~\mathrm{Spec}^\mathrm{CS}\breve{\Phi}^r_{\psi,\Gamma,X_\square},\underset{I}{\mathrm{homotopycolimit}}~\mathrm{Spec}^\mathrm{CS}\breve{\Phi}^I_{\psi,\Gamma,X_\square},	\\
}
\]
\[
\xymatrix@R+0pc@C+0pc{
\underset{r}{\mathrm{homotopylimit}}~\mathrm{Spec}^\mathrm{CS}{\Phi}^r_{\psi,\Gamma,X_\square},\underset{I}{\mathrm{homotopycolimit}}~\mathrm{Spec}^\mathrm{CS}{\Phi}^I_{\psi,\Gamma,X_\square}.	
}
\]
\[ 
\xymatrix@R+0pc@C+0pc{
\underset{r}{\mathrm{homotopylimit}}~\mathrm{Spec}^\mathrm{CS}\widetilde{\Phi}^r_{\psi,\Gamma,X_\square}/\mathrm{Fro}^\mathbb{Z},\underset{I}{\mathrm{homotopycolimit}}~\mathrm{Spec}^\mathrm{CS}\widetilde{\Phi}^I_{\psi,\Gamma,X_\square}/\mathrm{Fro}^\mathbb{Z},	\\
}
\]
\[ 
\xymatrix@R+0pc@C+0pc{
\underset{r}{\mathrm{homotopylimit}}~\mathrm{Spec}^\mathrm{CS}\breve{\Phi}^r_{\psi,\Gamma,X_\square}/\mathrm{Fro}^\mathbb{Z},\underset{I}{\mathrm{homotopycolimit}}~\breve{\Phi}^I_{\psi,\Gamma,X_\square}/\mathrm{Fro}^\mathbb{Z},	\\
}
\]
\[ 
\xymatrix@R+0pc@C+0pc{
\underset{r}{\mathrm{homotopylimit}}~\mathrm{Spec}^\mathrm{CS}{\Phi}^r_{\psi,\Gamma,X_\square}/\mathrm{Fro}^\mathbb{Z},\underset{I}{\mathrm{homotopycolimit}}~\mathrm{Spec}^\mathrm{CS}{\Phi}^I_{\psi,\Gamma,X_\square}/\mathrm{Fro}^\mathbb{Z}.	
}
\]

\end{definition}

\

\begin{definition}
We then consider the corresponding quasisheaves of the corresponding condensed solid topological modules from \cite{10CS2}:
\begin{align}
\mathrm{Quasicoherentsheaves, Condensed}_{*}	
\end{align}
where $*$ is one of the following spaces:
\begin{align}
&\mathrm{Spec}^\mathrm{CS}\widetilde{\Delta}_{\psi,\Gamma,X_\square}/\mathrm{Fro}^\mathbb{Z},\mathrm{Spec}^\mathrm{CS}\widetilde{\nabla}_{\psi,\Gamma,X_\square}/\mathrm{Fro}^\mathbb{Z},\mathrm{Spec}^\mathrm{CS}\widetilde{\Phi}_{\psi,\Gamma,X_\square}/\mathrm{Fro}^\mathbb{Z},\mathrm{Spec}^\mathrm{CS}\widetilde{\Delta}^+_{\psi,\Gamma,X_\square}/\mathrm{Fro}^\mathbb{Z},\\
&\mathrm{Spec}^\mathrm{CS}\widetilde{\nabla}^+_{\psi,\Gamma,X_\square}/\mathrm{Fro}^\mathbb{Z},\mathrm{Spec}^\mathrm{CS}\widetilde{\Delta}^\dagger_{\psi,\Gamma,X_\square}/\mathrm{Fro}^\mathbb{Z},\mathrm{Spec}^\mathrm{CS}\widetilde{\nabla}^\dagger_{\psi,\Gamma,X_\square}/\mathrm{Fro}^\mathbb{Z},	\\
\end{align}
\begin{align}
&\mathrm{Spec}^\mathrm{CS}\breve{\Delta}_{\psi,\Gamma,X_\square}/\mathrm{Fro}^\mathbb{Z},\breve{\nabla}_{\psi,\Gamma,X_\square}/\mathrm{Fro}^\mathbb{Z},\mathrm{Spec}^\mathrm{CS}\breve{\Phi}_{\psi,\Gamma,X_\square}/\mathrm{Fro}^\mathbb{Z},\mathrm{Spec}^\mathrm{CS}\breve{\Delta}^+_{\psi,\Gamma,X_\square}/\mathrm{Fro}^\mathbb{Z},\\
&\mathrm{Spec}^\mathrm{CS}\breve{\nabla}^+_{\psi,\Gamma,X_\square}/\mathrm{Fro}^\mathbb{Z},\mathrm{Spec}^\mathrm{CS}\breve{\Delta}^\dagger_{\psi,\Gamma,X_\square}/\mathrm{Fro}^\mathbb{Z},\mathrm{Spec}^\mathrm{CS}\breve{\nabla}^\dagger_{\psi,\Gamma,X_\square}/\mathrm{Fro}^\mathbb{Z},	\\
\end{align}
\begin{align}
&\mathrm{Spec}^\mathrm{CS}{\Delta}_{\psi,\Gamma,X_\square}/\mathrm{Fro}^\mathbb{Z},\mathrm{Spec}^\mathrm{CS}{\nabla}_{\psi,\Gamma,X_\square}/\mathrm{Fro}^\mathbb{Z},\mathrm{Spec}^\mathrm{CS}{\Phi}_{\psi,\Gamma,X_\square}/\mathrm{Fro}^\mathbb{Z},\mathrm{Spec}^\mathrm{CS}{\Delta}^+_{\psi,\Gamma,X_\square}/\mathrm{Fro}^\mathbb{Z},\\
&\mathrm{Spec}^\mathrm{CS}{\nabla}^+_{\psi,\Gamma,X_\square}/\mathrm{Fro}^\mathbb{Z}, \mathrm{Spec}^\mathrm{CS}{\Delta}^\dagger_{\psi,\Gamma,X_\square}/\mathrm{Fro}^\mathbb{Z},\mathrm{Spec}^\mathrm{CS}{\nabla}^\dagger_{\psi,\Gamma,X_\square}/\mathrm{Fro}^\mathbb{Z}.	
\end{align}
Here for those space with notations related to the radius and the corresponding interval we consider the total unions $\bigcap_r,\bigcup_I$ in order to achieve the whole spaces to achieve the analogues of the corresponding FF curves from \cite{10KL1}, \cite{10KL2}, \cite{10FF} for
\[
\xymatrix@R+0pc@C+0pc{
\underset{r}{\mathrm{homotopylimit}}~\mathrm{Spec}^\mathrm{CS}\widetilde{\Phi}^r_{\psi,\Gamma,X_\square},\underset{I}{\mathrm{homotopycolimit}}~\mathrm{Spec}^\mathrm{CS}\widetilde{\Phi}^I_{\psi,\Gamma,X_\square},	\\
}
\]
\[
\xymatrix@R+0pc@C+0pc{
\underset{r}{\mathrm{homotopylimit}}~\mathrm{Spec}^\mathrm{CS}\breve{\Phi}^r_{\psi,\Gamma,X_\square},\underset{I}{\mathrm{homotopycolimit}}~\mathrm{Spec}^\mathrm{CS}\breve{\Phi}^I_{\psi,\Gamma,X_\square},	\\
}
\]
\[
\xymatrix@R+0pc@C+0pc{
\underset{r}{\mathrm{homotopylimit}}~\mathrm{Spec}^\mathrm{CS}{\Phi}^r_{\psi,\Gamma,X_\square},\underset{I}{\mathrm{homotopycolimit}}~\mathrm{Spec}^\mathrm{CS}{\Phi}^I_{\psi,\Gamma,X_\square}.	
}
\]
\[ 
\xymatrix@R+0pc@C+0pc{
\underset{r}{\mathrm{homotopylimit}}~\mathrm{Spec}^\mathrm{CS}\widetilde{\Phi}^r_{\psi,\Gamma,X_\square}/\mathrm{Fro}^\mathbb{Z},\underset{I}{\mathrm{homotopycolimit}}~\mathrm{Spec}^\mathrm{CS}\widetilde{\Phi}^I_{\psi,\Gamma,X_\square}/\mathrm{Fro}^\mathbb{Z},	\\
}
\]
\[ 
\xymatrix@R+0pc@C+0pc{
\underset{r}{\mathrm{homotopylimit}}~\mathrm{Spec}^\mathrm{CS}\breve{\Phi}^r_{\psi,\Gamma,X_\square}/\mathrm{Fro}^\mathbb{Z},\underset{I}{\mathrm{homotopycolimit}}~\breve{\Phi}^I_{\psi,\Gamma,X_\square}/\mathrm{Fro}^\mathbb{Z},	\\
}
\]
\[ 
\xymatrix@R+0pc@C+0pc{
\underset{r}{\mathrm{homotopylimit}}~\mathrm{Spec}^\mathrm{CS}{\Phi}^r_{\psi,\Gamma,X_\square}/\mathrm{Fro}^\mathbb{Z},\underset{I}{\mathrm{homotopycolimit}}~\mathrm{Spec}^\mathrm{CS}{\Phi}^I_{\psi,\Gamma,X_\square}/\mathrm{Fro}^\mathbb{Z}.	
}
\]

\end{definition}

\

\begin{proposition}
There is a well-defined functor from the $\infty$-category 
\begin{align}
\mathrm{Quasicoherentpresheaves,Condensed}_{*}	
\end{align}
where $*$ is one of the following spaces:
\begin{align}
&\mathrm{Spec}^\mathrm{CS}\widetilde{\Phi}_{\psi,\Gamma,X_\square}/\mathrm{Fro}^\mathbb{Z},	\\
\end{align}
\begin{align}
&\mathrm{Spec}^\mathrm{CS}\breve{\Phi}_{\psi,\Gamma,X_\square}/\mathrm{Fro}^\mathbb{Z},	\\
\end{align}
\begin{align}
&\mathrm{Spec}^\mathrm{CS}{\Phi}_{\psi,\Gamma,X_\square}/\mathrm{Fro}^\mathbb{Z},	
\end{align}
to the $\infty$-category of $\mathrm{Fro}$-equivariant quasicoherent presheaves over similar spaces above correspondingly without the $\mathrm{Fro}$-quotients, and to the $\infty$-category of $\mathrm{Fro}$-equivariant quasicoherent modules over global sections of the structure $\infty$-sheaves of the similar spaces above correspondingly without the $\mathrm{Fro}$-quotients. Here for those space without notation related to the radius and the corresponding interval we consider the total unions $\bigcap_r,\bigcup_I$ in order to achieve the whole spaces to achieve the analogues of the corresponding FF curves from \cite{10KL1}, \cite{10KL2}, \cite{10FF} for
\[
\xymatrix@R+0pc@C+0pc{
\underset{r}{\mathrm{homotopylimit}}~\mathrm{Spec}^\mathrm{CS}\widetilde{\Phi}^r_{\psi,\Gamma,X_\square},\underset{I}{\mathrm{homotopycolimit}}~\mathrm{Spec}^\mathrm{CS}\widetilde{\Phi}^I_{\psi,\Gamma,X_\square},	\\
}
\]
\[
\xymatrix@R+0pc@C+0pc{
\underset{r}{\mathrm{homotopylimit}}~\mathrm{Spec}^\mathrm{CS}\breve{\Phi}^r_{\psi,\Gamma,X_\square},\underset{I}{\mathrm{homotopycolimit}}~\mathrm{Spec}^\mathrm{CS}\breve{\Phi}^I_{\psi,\Gamma,X_\square},	\\
}
\]
\[
\xymatrix@R+0pc@C+0pc{
\underset{r}{\mathrm{homotopylimit}}~\mathrm{Spec}^\mathrm{CS}{\Phi}^r_{\psi,\Gamma,X_\square},\underset{I}{\mathrm{homotopycolimit}}~\mathrm{Spec}^\mathrm{CS}{\Phi}^I_{\psi,\Gamma,X_\square}.	
}
\]
\[ 
\xymatrix@R+0pc@C+0pc{
\underset{r}{\mathrm{homotopylimit}}~\mathrm{Spec}^\mathrm{CS}\widetilde{\Phi}^r_{\psi,\Gamma,X_\square}/\mathrm{Fro}^\mathbb{Z},\underset{I}{\mathrm{homotopycolimit}}~\mathrm{Spec}^\mathrm{CS}\widetilde{\Phi}^I_{\psi,\Gamma,X_\square}/\mathrm{Fro}^\mathbb{Z},	\\
}
\]
\[ 
\xymatrix@R+0pc@C+0pc{
\underset{r}{\mathrm{homotopylimit}}~\mathrm{Spec}^\mathrm{CS}\breve{\Phi}^r_{\psi,\Gamma,X_\square}/\mathrm{Fro}^\mathbb{Z},\underset{I}{\mathrm{homotopycolimit}}~\breve{\Phi}^I_{\psi,\Gamma,X_\square}/\mathrm{Fro}^\mathbb{Z},	\\
}
\]
\[ 
\xymatrix@R+0pc@C+0pc{
\underset{r}{\mathrm{homotopylimit}}~\mathrm{Spec}^\mathrm{CS}{\Phi}^r_{\psi,\Gamma,X_\square}/\mathrm{Fro}^\mathbb{Z},\underset{I}{\mathrm{homotopycolimit}}~\mathrm{Spec}^\mathrm{CS}{\Phi}^I_{\psi,\Gamma,X_\square}/\mathrm{Fro}^\mathbb{Z}.	
}
\]	
In this situation we will have the target category being family parametrized by $r$ or $I$ in compatible glueing sense as in \cite[Definition 5.4.10]{10KL2}. In this situation for modules parametrized by the intervals we have the equivalence of $\infty$-categories by using \cite[Proposition 13.8]{10CS2}. Here the corresponding quasicoherent Frobenius modules are defined to be the corresponding homotopy colimits and limits of Frobenius modules:
\begin{align}
\underset{r}{\mathrm{homotopycolimit}}~M_r,\\
\underset{I}{\mathrm{homotopylimit}}~M_I,	
\end{align}
where each $M_r$ is a Frobenius-equivariant module over the period ring with respect to some radius $r$ while each $M_I$ is a Frobenius-equivariant module over the period ring with respect to some interval $I$.\\
\end{proposition}

\

\begin{definition}
We then consider the corresponding quasisheaves of perfect complexes of the corresponding condensed solid topological modules from \cite{10CS2}:
\begin{align}
\mathrm{Quasicoherentsheaves, Perfectcomplex, Condensed}_{*}	
\end{align}
where $*$ is one of the following spaces:
\begin{align}
&\mathrm{Spec}^\mathrm{CS}\widetilde{\Delta}_{\psi,\Gamma,X_\square}/\mathrm{Fro}^\mathbb{Z},\mathrm{Spec}^\mathrm{CS}\widetilde{\nabla}_{\psi,\Gamma,X_\square}/\mathrm{Fro}^\mathbb{Z},\mathrm{Spec}^\mathrm{CS}\widetilde{\Phi}_{\psi,\Gamma,X_\square}/\mathrm{Fro}^\mathbb{Z},\mathrm{Spec}^\mathrm{CS}\widetilde{\Delta}^+_{\psi,\Gamma,X_\square}/\mathrm{Fro}^\mathbb{Z},\\
&\mathrm{Spec}^\mathrm{CS}\widetilde{\nabla}^+_{\psi,\Gamma,X_\square}/\mathrm{Fro}^\mathbb{Z},\mathrm{Spec}^\mathrm{CS}\widetilde{\Delta}^\dagger_{\psi,\Gamma,X_\square}/\mathrm{Fro}^\mathbb{Z},\mathrm{Spec}^\mathrm{CS}\widetilde{\nabla}^\dagger_{\psi,\Gamma,X_\square}/\mathrm{Fro}^\mathbb{Z},	\\
\end{align}
\begin{align}
&\mathrm{Spec}^\mathrm{CS}\breve{\Delta}_{\psi,\Gamma,X_\square}/\mathrm{Fro}^\mathbb{Z},\breve{\nabla}_{\psi,\Gamma,X_\square}/\mathrm{Fro}^\mathbb{Z},\mathrm{Spec}^\mathrm{CS}\breve{\Phi}_{\psi,\Gamma,X_\square}/\mathrm{Fro}^\mathbb{Z},\mathrm{Spec}^\mathrm{CS}\breve{\Delta}^+_{\psi,\Gamma,X_\square}/\mathrm{Fro}^\mathbb{Z},\\
&\mathrm{Spec}^\mathrm{CS}\breve{\nabla}^+_{\psi,\Gamma,X_\square}/\mathrm{Fro}^\mathbb{Z},\mathrm{Spec}^\mathrm{CS}\breve{\Delta}^\dagger_{\psi,\Gamma,X_\square}/\mathrm{Fro}^\mathbb{Z},\mathrm{Spec}^\mathrm{CS}\breve{\nabla}^\dagger_{\psi,\Gamma,X_\square}/\mathrm{Fro}^\mathbb{Z},	\\
\end{align}
\begin{align}
&\mathrm{Spec}^\mathrm{CS}{\Delta}_{\psi,\Gamma,X_\square}/\mathrm{Fro}^\mathbb{Z},\mathrm{Spec}^\mathrm{CS}{\nabla}_{\psi,\Gamma,X_\square}/\mathrm{Fro}^\mathbb{Z},\mathrm{Spec}^\mathrm{CS}{\Phi}_{\psi,\Gamma,X_\square}/\mathrm{Fro}^\mathbb{Z},\mathrm{Spec}^\mathrm{CS}{\Delta}^+_{\psi,\Gamma,X_\square}/\mathrm{Fro}^\mathbb{Z},\\
&\mathrm{Spec}^\mathrm{CS}{\nabla}^+_{\psi,\Gamma,X_\square}/\mathrm{Fro}^\mathbb{Z}, \mathrm{Spec}^\mathrm{CS}{\Delta}^\dagger_{\psi,\Gamma,X_\square}/\mathrm{Fro}^\mathbb{Z},\mathrm{Spec}^\mathrm{CS}{\nabla}^\dagger_{\psi,\Gamma,X_\square}/\mathrm{Fro}^\mathbb{Z}.	
\end{align}
Here for those space with notations related to the radius and the corresponding interval we consider the total unions $\bigcap_r,\bigcup_I$ in order to achieve the whole spaces to achieve the analogues of the corresponding FF curves from \cite{10KL1}, \cite{10KL2}, \cite{10FF} for
\[
\xymatrix@R+0pc@C+0pc{
\underset{r}{\mathrm{homotopylimit}}~\mathrm{Spec}^\mathrm{CS}\widetilde{\Phi}^r_{\psi,\Gamma,X_\square},\underset{I}{\mathrm{homotopycolimit}}~\mathrm{Spec}^\mathrm{CS}\widetilde{\Phi}^I_{\psi,\Gamma,X_\square},	\\
}
\]
\[
\xymatrix@R+0pc@C+0pc{
\underset{r}{\mathrm{homotopylimit}}~\mathrm{Spec}^\mathrm{CS}\breve{\Phi}^r_{\psi,\Gamma,X_\square},\underset{I}{\mathrm{homotopycolimit}}~\mathrm{Spec}^\mathrm{CS}\breve{\Phi}^I_{\psi,\Gamma,X_\square},	\\
}
\]
\[
\xymatrix@R+0pc@C+0pc{
\underset{r}{\mathrm{homotopylimit}}~\mathrm{Spec}^\mathrm{CS}{\Phi}^r_{\psi,\Gamma,X_\square},\underset{I}{\mathrm{homotopycolimit}}~\mathrm{Spec}^\mathrm{CS}{\Phi}^I_{\psi,\Gamma,X_\square}.	
}
\]
\[ 
\xymatrix@R+0pc@C+0pc{
\underset{r}{\mathrm{homotopylimit}}~\mathrm{Spec}^\mathrm{CS}\widetilde{\Phi}^r_{\psi,\Gamma,X_\square}/\mathrm{Fro}^\mathbb{Z},\underset{I}{\mathrm{homotopycolimit}}~\mathrm{Spec}^\mathrm{CS}\widetilde{\Phi}^I_{\psi,\Gamma,X_\square}/\mathrm{Fro}^\mathbb{Z},	\\
}
\]
\[ 
\xymatrix@R+0pc@C+0pc{
\underset{r}{\mathrm{homotopylimit}}~\mathrm{Spec}^\mathrm{CS}\breve{\Phi}^r_{\psi,\Gamma,X_\square}/\mathrm{Fro}^\mathbb{Z},\underset{I}{\mathrm{homotopycolimit}}~\breve{\Phi}^I_{\psi,\Gamma,X_\square}/\mathrm{Fro}^\mathbb{Z},	\\
}
\]
\[ 
\xymatrix@R+0pc@C+0pc{
\underset{r}{\mathrm{homotopylimit}}~\mathrm{Spec}^\mathrm{CS}{\Phi}^r_{\psi,\Gamma,X_\square}/\mathrm{Fro}^\mathbb{Z},\underset{I}{\mathrm{homotopycolimit}}~\mathrm{Spec}^\mathrm{CS}{\Phi}^I_{\psi,\Gamma,X_\square}/\mathrm{Fro}^\mathbb{Z}.	
}
\]

\end{definition}

\begin{proposition}
There is a well-defined functor from the $\infty$-category 
\begin{align}
\mathrm{Quasicoherentpresheaves,Perfectcomplex,Condensed}_{*}	
\end{align}
where $*$ is one of the following spaces:
\begin{align}
&\mathrm{Spec}^\mathrm{CS}\widetilde{\Phi}_{\psi,\Gamma,X_\square}/\mathrm{Fro}^\mathbb{Z},	\\
\end{align}
\begin{align}
&\mathrm{Spec}^\mathrm{CS}\breve{\Phi}_{\psi,\Gamma,X_\square}/\mathrm{Fro}^\mathbb{Z},	\\
\end{align}
\begin{align}
&\mathrm{Spec}^\mathrm{CS}{\Phi}_{\psi,\Gamma,X_\square}/\mathrm{Fro}^\mathbb{Z},	
\end{align}
to the $\infty$-category of $\mathrm{Fro}$-equivariant quasicoherent presheaves over similar spaces above correspondingly without the $\mathrm{Fro}$-quotients, and to the $\infty$-category of $\mathrm{Fro}$-equivariant quasicoherent modules over global sections of the structure $\infty$-sheaves of the similar spaces above correspondingly without the $\mathrm{Fro}$-quotients. Here for those space without notation related to the radius and the corresponding interval we consider the total unions $\bigcap_r,\bigcup_I$ in order to achieve the whole spaces to achieve the analogues of the corresponding FF curves from \cite{10KL1}, \cite{10KL2}, \cite{10FF} for
\[
\xymatrix@R+0pc@C+0pc{
\underset{r}{\mathrm{homotopylimit}}~\mathrm{Spec}^\mathrm{CS}\widetilde{\Phi}^r_{\psi,\Gamma,X_\square},\underset{I}{\mathrm{homotopycolimit}}~\mathrm{Spec}^\mathrm{CS}\widetilde{\Phi}^I_{\psi,\Gamma,X_\square},	\\
}
\]
\[
\xymatrix@R+0pc@C+0pc{
\underset{r}{\mathrm{homotopylimit}}~\mathrm{Spec}^\mathrm{CS}\breve{\Phi}^r_{\psi,\Gamma,X_\square},\underset{I}{\mathrm{homotopycolimit}}~\mathrm{Spec}^\mathrm{CS}\breve{\Phi}^I_{\psi,\Gamma,X_\square},	\\
}
\]
\[
\xymatrix@R+0pc@C+0pc{
\underset{r}{\mathrm{homotopylimit}}~\mathrm{Spec}^\mathrm{CS}{\Phi}^r_{\psi,\Gamma,X_\square},\underset{I}{\mathrm{homotopycolimit}}~\mathrm{Spec}^\mathrm{CS}{\Phi}^I_{\psi,\Gamma,X_\square}.	
}
\]
\[ 
\xymatrix@R+0pc@C+0pc{
\underset{r}{\mathrm{homotopylimit}}~\mathrm{Spec}^\mathrm{CS}\widetilde{\Phi}^r_{\psi,\Gamma,X_\square}/\mathrm{Fro}^\mathbb{Z},\underset{I}{\mathrm{homotopycolimit}}~\mathrm{Spec}^\mathrm{CS}\widetilde{\Phi}^I_{\psi,\Gamma,X_\square}/\mathrm{Fro}^\mathbb{Z},	\\
}
\]
\[ 
\xymatrix@R+0pc@C+0pc{
\underset{r}{\mathrm{homotopylimit}}~\mathrm{Spec}^\mathrm{CS}\breve{\Phi}^r_{\psi,\Gamma,X_\square}/\mathrm{Fro}^\mathbb{Z},\underset{I}{\mathrm{homotopycolimit}}~\breve{\Phi}^I_{\psi,\Gamma,X_\square}/\mathrm{Fro}^\mathbb{Z},	\\
}
\]
\[ 
\xymatrix@R+0pc@C+0pc{
\underset{r}{\mathrm{homotopylimit}}~\mathrm{Spec}^\mathrm{CS}{\Phi}^r_{\psi,\Gamma,X_\square}/\mathrm{Fro}^\mathbb{Z},\underset{I}{\mathrm{homotopycolimit}}~\mathrm{Spec}^\mathrm{CS}{\Phi}^I_{\psi,\Gamma,X_\square}/\mathrm{Fro}^\mathbb{Z}.	
}
\]	
In this situation we will have the target category being family parametrized by $r$ or $I$ in compatible glueing sense as in \cite[Definition 5.4.10]{10KL2}. Here the corresponding quasicoherent Frobenius modules are defined to be the corresponding homotopy colimits and limits of Frobenius modules:
\begin{align}
\underset{r}{\mathrm{homotopycolimit}}~M_r,\\
\underset{I}{\mathrm{homotopylimit}}~M_I,	
\end{align}
where each $M_r$ is a Frobenius-equivariant module over the period ring with respect to some radius $r$ while each $M_I$ is a Frobenius-equivariant module over the period ring with respect to some interval $I$.\\
\end{proposition}

\newpage

\subsection*{Acknowledgements}
From the discussion with Professor Kedlaya, we understand the importance of the application of \cite{10BK}, \cite{10BBBK}, \cite{10BBK}, \cite{10KKM}, \cite{10BBM}, \cite{10CS1}, \cite{10CS2} to many problems in our setting. We thank Professor Kedlaya for the related discussion on these topics. We also have the chance to learn much from the work \cite{10GEH} through the discussion with the authors.

\bibliographystyle{splncs}

\end{document}